\documentclass{amsbook}
\usepackage{amssymb,amscd}
\usepackage[cp1251]{inputenc}
\usepackage[T2A]{fontenc}
\usepackage[dvips]{graphicx}
\usepackage{pstricks,pst-grad,pst-node}%picture
\usepackage[all]{xy}

\input xypic

\linethickness{0.5pt}
\newtheorem{theorem}{Theorem}[section]
\newtheorem{proposition}[theorem]{Proposition}
\newtheorem{corollary}[theorem]{Corollary}
\newtheorem{lemma}[theorem]{Lemma}
\newtheorem{conjecture}[theorem]{Conjecture}
\theoremstyle{definition}
\newtheorem{definition}[theorem]{Definition}
\newtheorem{example}[theorem]{Example}
\newtheorem{construction}[theorem]{Construction}
\newtheorem{exercise}[theorem]{}
\newtheorem{problem}[theorem]{Problem}
\theoremstyle{remark}
\newtheorem*{claim}{Claim}
\newtheorem*{remark}{Remark}

\numberwithin{equation}{chapter}
\numberwithin{figure}{chapter}
\numberwithin{section}{chapter}

\def\C{\mathbb C}
\def\D{\mathbb D}
\def\I{\mathbb I}
\def\N{\mathbb N}
\def\Q{\mathbb Q}
\def\R{\mathbb R}
\def\T{\mathbb T}
\def\Z{\mathbb Z}

\def\sB{\mathcal B}
\def\sC{\mathcal C}
\def\sK{\mathcal K}
\def\sL{\mathcal L}
\def\sS{\mathcal S}

\def\k{\mathbf k}

\def\F{F}

\def\phi{\varphi}

\newcommand{\mb}[1]{{\textbf {\textit#1}}}
\newcommand{\bl}{\lambda\kern-0.53em\lambda}% bold lambda
\newcommand{\bmu}{\mu\kern-0.55em\mu}% bold mu
\newcommand{\bnu}{\nu\kern-0.51em\nu}% bold nu

\renewcommand{\u}[1]{ u_{{}_{\!#1}} }
\newcommand{\va}[1]{ v_{{}_{\!#1}} }
\newcommand{\ta}[1]{\tau_{{}_{\!#1}} }
\newcommand{\n}[1]{\nu_{{}_{\!#1}} }

\newcommand{\sbr}[2]{{\textstyle\genfrac{[}{]}{}{}{#1}{#2}}}

\newcommand{\bin}[2]{{\textstyle\binom{#1}{#2}}}

\renewcommand{\ge}{\geqslant}
\renewcommand{\le}{\leqslant}

\newcommand{\llra}{\relbar\joinrel\hspace{-1pt}\longrightarrow}
\newcommand{\lllra}{\relbar\joinrel\hspace{-1pt}\llra}

\def\dbs{/\!\!/}
\newcommand{\id}{\mathrm{id}}
\newcommand{\ind}{\mathop{\rm ind}\nolimits}
\def\scs{\raisebox{0.07\baselineskip}{$\scriptstyle\mathbin{\widetilde\#}$}}

\def\Ann{\mathop{\mathrm{Ann}}}
\newcommand{\bideg}{\mathop{\rm bideg}}
\newcommand{\codim}{\mathop{\rm codim}}
\newcommand{\depth}{\mathop{\mathrm{depth}}\nolimits}
\newcommand{\Ext}{\mathop{\rm Ext}\nolimits}
\newcommand{\Hom}{\mathop{\mathrm{Hom}}\nolimits}
\renewcommand{\Im}{\mathop{\mathrm{Im}}\nolimits}
\newcommand{\Ker}{\mathop{\rm Ker}}
\newcommand{\pdim}{\mathop{\mathrm{pdim}}\nolimits}
\newcommand{\rank}{\mathop{\mathrm{rank}}}
\renewcommand{\Re}{\mathop{\mathrm{Re}}}
\newcommand{\Tor}{\mathop{\rm Tor}\nolimits}

\newcommand{\cl}{\mathop\mathit{CL}\nolimits}
\newcommand{\fl}{\mathop\mathit{FL}\nolimits}

\def\Spec{\mathop{\mathrm{Spec}}}

\newcommand{\cat}[1]{\mbox{\sc #1}}
\newcommand{\fcat}[2]{\mbox{\raisebox{1pt}{\rm\scriptsize [}{\sc #1},\hspace{1pt}{\sc #2}\raisebox{1pt}{\rm\scriptsize ]}}\hspace{1pt}}
\newcommand{\scat}[1]{\mbox{\scriptsize{\sc #1}}}

\def\ca{{\text{\sc cat}}}
\def\cga{\text{\sc cga}}
\def\del{\mbox{\footnotesize{$\Delta$}}}
\def\top{\text{\sc top}}

\def\colim{\mathop\mathrm{colim}\nolimits}
\def\lim{\mathop\mathrm{lim}\nolimits}
\def\hocolim{\mathop\mathrm{hocolim}\nolimits}
\def\holim{\mathop\mathrm{holim}\nolimits}

\def\under{\!\downarrow\!}

\def\Ho{\mathop\mathit{Ho}}

\newcommand{\llongrightarrow}{\relbar\joinrel\longrightarrow}
\newcommand{\lllongrightarrow}{\relbar\joinrel\llongrightarrow}
\newcommand{\llllongrightarrow}{\relbar\joinrel\lllongrightarrow}
\newcommand{\lllllongrightarrow}{\relbar\joinrel\llllongrightarrow}

\newcommand{\bm}{\mathop{\mbox{bm}}\nolimits}%bistellar move
\newcommand{\cc}{\mathop{\rm cc}}%cubical complex
\def\cone{\mathop{\mathrm{cone}}}
\newcommand{\core}{\mathop{\rm core}}
\newcommand{\cub}{\mathop{\rm cub}}%cubical complex
\def\hatzero{\hat0}
\newcommand{\lk}{\mathop{\rm lk}\nolimits}%link
\newcommand{\ord}{\mathop{\rm ord}}%order complex
\def\ss{\mathop{\rm ss}\nolimits}%stellar subdivision
\newcommand{\st}{\mathop{\rm st}\nolimits}%star

\def\conv{\mathop{\mathrm{conv}}}
\def\vt{\mathop{\mathrm{vt}}}

\newcommand{\APL}{\mbox{$A_{P\!L}$}}
\newcommand{\cf}{\mbox{\it{c\!f\hspace{1pt}}}}
\def\ch{\mathop{\mathrm{ch}}\nolimits}
\newcommand{\cs}{\mathbin{\#}}
\def\pt{\mathit{pt}}
\newcommand{\sign}{\mathop{\rm sign}}
\def\sk{\mathop\mathrm{sk}\!{}}
\newcommand{\td}{\mathop{\rm td}}
\def\Th{\mathop{\mbox{\textit{Th}}}}

\def\MO{\mathop{\mbox{\it MO\/}}\nolimits}
\def\BO{\mathop{\mbox{\it BO\/}}}
\def\EO{\mathop{\mbox{\it EO\/}}}
\def\MU{\mathop{\mbox{\it MU\/}}\nolimits}
\def\BU{\mathop{\mbox{\it BU\/}}}
\def\EU{\mathop{\mbox{\it EU\/}}}

\def\BF{\mbox{\it BF\/}}
\newcommand{\djs}{\mathop{\mbox{\it DJ\/}}}

\newcommand{\rk}{\mathcal R_{\mathcal K}}
\newcommand{\zk}{\mathcal Z_{\mathcal K}}
\newcommand{\zp}{\mathcal Z_P}
\newcommand{\zs}{\mathcal Z_{\mathcal S}}
\def\ftr{\mathop{\mathrm{ftr}}} %free toral rank
\def\atr{\mathop{\mathrm{atr}}} %almost free toral rank

\newcommand{\<}{\langle}
\renewcommand{\>}{\rangle}

\def\Amap{A}

\begin{document}
\frontmatter
\title{Toric Topology}
\author{Victor M. Buchstaber}
\address{Steklov Mathematical Institute,
Russian Academy of Sciences, Gubkina Street 8, 119991 Moscow,
Russia}
\email{buchstab@mi.ras.ru}

\author{Taras E. Panov}
\address{Department of Mathematics and Mechanics, Moscow
State University, Leninskie Gory, 119991 Moscow, Russia}
\address{Institute for Information Transmission Problems, Russian Academy of Sciences, Moscow}
\address{Institute of Theoretical and Experimental Physics, Moscow}
\email{tpanov@mech.math.msu.su}

%\subjclass[2000]{52B70, 57Q15, 57R19, 14M25, 52B05, 13F55, 52C35}

\begin{abstract}
Toric topology emerged in the end of the 1990s on the borders of
equivariant topology, algebraic and symplectic geometry,
combinatorics and commutative algebra. It has quickly grown up
into a very active area with many interdisciplinary links and
applications, and continues to attract experts from different
fields.

The key players in toric topology are moment-angle manifolds, a
family of manifolds with torus actions defined in combinatorial
terms. Their construction links to combinatorial geometry and
algebraic geometry of toric varieties via the related notion of a
quasitoric manifold. Discovery of remarkable geometric structures
on moment-angle manifolds led to seminal connections with the
classical and modern areas of symplectic, Lagrangian and
non-K\"ahler complex geometry. A related categorical construction
of moment-angle complexes and their generalisations, polyhedral
products, provides a universal framework for many fundamental
constructions of homotopical topology. The study of polyhedral
products is now evolving into a separate area of homotopy theory,
with strong links to other areas of toric topology. A new
perspective on torus action has also contributed to the
development of classical areas of algebraic topology, such as
complex cobordism.

The book contains lots of open problems and is addressed to
experts interested in new ideas linking all the subjects involved,
as well as to graduate students and young researchers ready to
enter into a beautiful new area.
\end{abstract}

\maketitle
\tableofcontents

\mainmatter

%Сюжеты, не вошедшие в книгу:
%% tropical geometry
%% Riemannian metrics with special curvature properties

\chapter*{Introduction}
Traditionally, the study of torus actions on topological spaces
has been considered as a classical field of algebraic topology.
Specific problems connected with torus actions arise in different
areas of mathematics and mathematical physics, which results in
permanent interest in the theory, constant source of new
applications and penetration of new ideas into topology.

Since the 1970s, algebraic and symplectic viewpoints on torus
actions have enriched the subject with new combinatorial ideas and
methods, largely based on the convex-geometric concept of
polytopes.

The study of algebraic torus actions on algebraic varieties has
quickly developed into a much successful branch of algebraic
geometry, known as \emph{toric geometry}. It gives a bijection
between, on the one hand, \emph{toric varieties}, which are
complex algebraic varieties equipped with an action of an
algebraic torus with a dense orbit, and on the other hand,
\emph{fans}, which are combinatorial objects. The fan allows one
to completely translate various algebraic-geometric notions into
combinatorics. Projective toric varieties correspond to fans which
arise from convex polytopes. A valuable aspect of this theory is
that it provides many explicit examples of algebraic varieties,
leading to applications in deep subjects such as the singularity
theory and mirror symmetry.
%It also points the way toward
%the deep link between the study of these varieties and the
%associated combinatorics of their fans.

In symplectic geometry, since the early 1980s there has been much
activity in the field of Hamiltonian group actions on symplectic
manifolds. Such an action defines the \emph{moment map} from the
manifold to a Euclidean space (more precisely, the dual Lie
algebra of the torus) whose image is a convex polytope. If the
torus has half the dimension of the manifold the moment map image
determines the manifold up to equivariant symplectomorphism. The
class of polytopes which can arise as the images of moment maps
can be described explicitly, together with an effective procedure
of recovering a symplectic manifold from such a polytope. In
symplectic geometry, as in algebraic geometry, one translates
various geometric constructions into the language of convex
polytopes and combinatorics.

There is a tight relationship between the algebraic and the
symplectic pictures: a projective embedding of a toric manifold
determines a symplectic form and a moment map. The image of the
moment map is a convex polytope that is dual to the fan. In both
the smooth algebraic-geometric and the symplectic situations, the
compact torus action is locally isomorphic to the standard action
of $(S^1)^n$ on $\C^n$ by rotation of the coordinates. Thus the
quotient of the manifold by this action is naturally a manifold
with corners, stratified according to the dimension of the
stabilisers, and each stratum can be equipped with data that
encodes the isotropy torus action along that stratum. Not only
does this structure of the quotient provide a powerful means of
investigating the action, but some of its subtler combinatorial
properties may also be illuminated by a careful study of the
equivariant topology of the manifold. Thus, it should come as no
surprise that since the beginning of the 1990s, the ideas and
methodology of toric varieties and Hamiltonian torus actions have
started penetrating back into algebraic topology.

By 2000, several constructions of topological analogues of toric
varieties and symplectic toric manifolds have appeared in the
literature, together with different seemingly unrelated
realisations of what later has become known as the moment-angle
manifolds. We tried to systematise both known and emerging links
between torus actions and combinatorics in our 2000
paper~\cite{bu-pa00-2} in Russian Mathematical Surveys, where the
terms `moment-angle manifold' and `moment-angle complex' first
appeared. Two years later it grew up into a book `Torus Actions
and Their Applications in Topology and
Combinatorics'~\cite{bu-pa02} published by the AMS in 2002 (the
extended Russian edition~\cite{bu-pa04-2} appeared in 2004). The
title `Toric Topology' coined by our colleague Nigel Ray became
official after the 2006 Osaka conference under the same name. Its
proceedings volume~\cite{tori08} contained many important
contributions to the subject, as well as the introductory survey
`An invitation to toric topology: vertex four of a remarkable
tetrahedron' by Buchstaber and Ray. The vertices of the `toric
tetrahedron' are topology, combinatorics, algebraic and symplectic
geometry, and it symbolised much strengthened links between these
subjects.  With many young researchers entering the subject and
conferences held around the world every year, toric topology has
definitely grown up into a mature area. Its various aspects are
presented in this monograph, with an intention to consolidate the
foundations and stimulate further applications.

\section*{Chapter guide}
\begin{figure}[h]
\begin{center}
\begin{picture}(125,75)
\put(0,60){\line(1,0){20}}
\put(20,60){\line(0,1){15}}
\put(20,75){\line(-1,0){20}}
\put(0,75){\line(0,-1){15}}
\put(22,67.5){\vector(1,0){11}}
\put(9.5,70){1}
\put(4.5,65){\scriptsize{Polytopes}}
\put(35,60){\line(1,0){20}}
\put(55,60){\line(0,1){15}}
\put(55,75){\line(-1,0){20}}
\put(35,75){\line(0,-1){15}}
\put(57,67.5){\vector(1,0){11}}
\put(45,58){\vector(0,-1){11}}
\put(44.5,70){2}
\put(36.5,65){\scriptsize{Combinatorial}}
\put(39.5,62.5){\scriptsize{structures}}
\put(70,60){\line(1,0){20}}
\put(90,60){\line(0,1){15}}
\put(90,75){\line(-1,0){20}}
\put(70,75){\line(0,-1){15}}
\put(80,58){\vector(0,-1){11}}
\put(79.5,70){3}
\put(77.5,65){\scriptsize{Face}}
\put(77.3,62.5){\scriptsize{rings}}
%
%\put(0,30){\line(1,0){20}} \put(20,30){\line(0,1){15}}
%\put(20,45){\line(-1,0){20}} \put(0,45){\line(0,-1){15}}
%\put(33,37.5){\vector(-1,0){11}} \put(45,28){\vector(0,-1){11}}
%\put(9.5,40){8} \put(6.7,35){\scriptsize{Torus}}
%\put(4.5,32.5){\scriptsize{manifolds}}
%
\put(70,30){\line(1,0){20}}
\put(90,30){\line(0,1){15}}
\put(90,45){\line(-1,0){20}}
\put(70,45){\line(0,-1){15}}
\put(68,28){\vector(-4,-3){15}}
\put(92,37.5){\vector(1,0){11}}
\put(80,28){\vector(0,-1){11}}
\put(79.5,40){4}
\put(71.5,35){\scriptsize{Moment-angle}}
\put(74.2,32.5){\scriptsize{complexes}}
\put(35,30){\line(1,0){20}}
\put(55,30){\line(0,1){15}}
\put(55,45){\line(-1,0){20}}
\put(35,45){\line(0,-1){15}}
\put(57,28){\vector(4,-3){15}}
\put(45,28){\vector(0,-1){11}}
\put(44.5,40){5}
\put(42,35){\scriptsize{Toric}}
\put(40.2,32.5){\scriptsize{varieties}}
\put(35,0){\line(1,0){20}}
\put(55,0){\line(0,1){15}}
\put(55,15){\line(-1,0){20}}
\put(35,15){\line(0,-1){15}}
%\put(57,7.5){\vector(1,0){11}}
\put(44.5,10){6}
\put(36.5,5){\scriptsize{Moment-angle}}
\put(39.5,2.5){\scriptsize{manifolds}}
\put(70,00){\line(1,0){20}} \put(90,00){\line(0,1){15}}
\put(90,15){\line(-1,0){20}} \put(70,15){\line(0,-1){15}}
\put(92,7.5){\vector(1,0){11}} \put(79.5,10){7}
\put(75,5){\scriptsize{Half-dim}} \put(72.5,2.5){\scriptsize{torus
actions}}
\put(105,30){\line(1,0){20}}
\put(125,30){\line(0,1){15}}
\put(125,45){\line(-1,0){20}}
\put(105,45){\line(0,-1){15}}
\put(114.5,40){8}
\put(109,35){\scriptsize{Homotopy}}
\put(111.5,32.5){\scriptsize{theory}}
\put(105,0){\line(1,0){20}} \put(125,0){\line(0,1){15}}
\put(125,15){\line(-1,0){20}} \put(105,15){\line(0,-1){15}}
\put(113.5,10){9} \put(109,5){\scriptsize{Cobordism}}
\end{picture}
%\caption{Chapter dependence scheme.} \label{cont-scheme}
\end{center}
\end{figure}

At the heart of toric topology lies a class of torus actions whose
orbit spaces are highly structured in combinatorial terms, that
is, have lots of orbit types tied together in a nice combinatorial
way. We use the generic terms \emph{toric space}\label{toricspace}
and \emph{toric object} to refer to a topological space with a
nice torus action, or to a space produced from a torus action via
different standard topological or categorical constructions.
Examples of toric spaces include toric varieties, toric and
quasitoric manifolds and their generalisations, moment-angle
manifolds, moment-angle complexes and their Borel constructions,
polyhedral products, complements of coordinate subspace
arrangements, intersections of real or Hermitian quadrics, etc.

Each chapter and most sections have their own introductions with
more specific information about the contents.

In Chapter~1 we collect background material related to convex
polytopes, including basic convex-geometric constructions and the
combinatorial theory of face vectors. The famous $g$-theorem
describing integer sequences that can be the face vectors of
simple (or simplicial) polytopes was one of the most striking
applications of toric geometry to combinatorics. The concept of
Gale duality and Gale diagrams are important tools for the study
of moment-angle manifolds via intersections of quadrics. In the
additional sections we describe several combinatorial
constructions providing families of simple polytopes, including
nestohedra, graph associahedra, flagtopes and truncated cubes. The
classical series of permutahedra and associahedra (Stasheff
polytopes) are particular examples. The construction of nestohedra
takes its origin in singularity and representation theory. We
develop a differential algebraic formalism which links the
generating series of nestohedra to classical partial differential
equations. The potential of truncated cubes in toric topology is
yet to be fully exploited, as they provide an immense source of
explicitly constructed toric spaces.

In Chapter 2 we describe systematically combinatorial structures
that appear in the orbit spaces of toric objects. Besides convex
polytopes, these include fans, simplicial and cubical complexes,
simplicial posets. All these structures are objects of independent
interest for combinatorialists, and we emphasised the aspects of
their combinatorial theory most relevant to subsequent topological
applications.

The subject of Chapter~3 is the algebraic theory of face rings
(also known as Stanley--Reisner rings) of simplicial complexes,
and their generalisations to simplicial posets. With the
appearance of the face rings in the beginning of the 1970s in the
work of Reisner and Stanley many combinatorial problems were
translated into the language of commutative algebra, which paved
the way for their solution using the extensive machinery of
algebraic and homological methods. Algebraic tools used for
attacking combinatorial problems included regular sequences,
Cohen--Macaulay and Gorenstein rings, $\Tor$-algebras, local
cohomology, etc. A whole new thriving field appeared on the
borders of combinatorics and algebra, which has since become known
as \emph{combinatorial commutative algebra}.
%The basic reference here is Stanley's monograph~\cite{stan96}.

Chapter~4 is the first `toric' chapter of the book; it links the
combinatorial and algebraic constructions of the previous chapters
to the world of toric spaces. The concept of the moment-angle
complex $\zk$ is introduced as a functor from the category of
simplicial complexes $\sK$ to the category of topological spaces
with torus actions and equivariant maps. When $\sK$ is a
triangulated manifold, the moment-angle complex $\zk$ consists a
free orbit $\mathcal Z_{\varnothing}$ consisting of singular
points. Removing this orbit we obtain an open manifold
$\zk\setminus\mathcal Z_{\varnothing}$, which satisfies the
relative version of Poincar\'e duality. Combinatorial invariants
of simplicial complexes $\sK$ therefore can be described in terms
of topological characteristics of the corresponding moment-angle
complexes~$\zk$. In particular, the face numbers of $\sK$, as well
as the more subtle \emph{bigraded Betti numbers} of the face ring
$\Z[\sK]$ can be expressed in terms of the cellular cohomology
groups of~$\zk$. The integral cohomology ring $H^*(\mathcal
Z_{\mathcal K})$ is shown to be isomorphic to the Tor-algebra
$\mathop{\mathrm{Tor}}_{\mathbb Z[v_1,\ldots,v_m]}(\mathbb
Z[\mathcal K],\mathbb Z)$.
%The canonical bigraded structure in the
%Tor groups thereby acquires a geometric interpretation in terms of
%the bigraded cell decomposition of~$\mathcal Z_{\mathcal K}$.
The proof builds upon a construction of a ring model for
\emph{cellular} cochains of~$\zk$ and the corresponding cellular
diagonal approximation, which is functorial with respect to maps
of moment-angle complexes induced by simplicial maps of~$\sK$.
This functorial property of the cellular diagonal approximation
for~$\zk$ is quite special, due to the lack of such a construction
for general cell complexes. Another result of Chapter~4 is a
homotopy equivalence (an equivariant deformation retraction) from
the complement $U(\sK)$ of the arrangement of coordinate subspaces
in~$\C^m$ determined by~$\sK$ to the moment-angle complex~$\zk$.
Particular cases of this result are known in toric geometry and
geometric invariant theory. It opens a new perspective on
moment-angle complexes, linking them to the theory of
configuration spaces and arrangements.

We tried to make the material of the first four chapters of the
book accessible for an undergraduate student, or a reader with a
very basis knowledge of algebra and topology. The general
algebraic and topological constructions required here are
collected in Appendices~\ref{hab} and~\ref{algtop} respectively.
More experienced readers may refer to these appendices purely for
terminology and notation.

Toric varieties are the subject of Chapter~5. This is an extensive
area with vast literature available. We outline the influence of
toric geometry on the emergence  of toric topology and emphasise
combinatorial, topological and symplectic aspects of toric
varieties. The construction of moment-angle manifolds via
nondegenerate intersections of Hermitian quadrics in~$\C^m$,
motivated by symplectic geometry, is also discussed here. Some
basic knowledge of algebraic geometry may be required in
Chapter~5. Appropriate references are given in the introduction to
the chapter.

Geometry of moment-angle manifolds is studied in Chapter~6. The
construction of moment-angle manifolds as the level sets of toric
moment maps is taken as the starting point for the systematic
study of intersections of Hermitian quadrics via Gale duality.
Following a remarkable discovery by Bosio and Meersseman of
complex-analytic structures on moment-angle manifolds
corresponding to simple polytopes, we proceed by showing that
moment-angle manifolds corresponding to a more general class of
complete simplicial fans can also be endowed with complex-analytic
structures. The resulting family of \emph{non-K\"ahler} complex
manifolds includes the classical series of Hopf and
Calabi--Eckmann manifolds. We also describe important invariants
of these complex structures, such as the Hodge numbers and
Dolbeault cohomology rings, study holomorphic torus principal
bundles over toric varieties, and establish collapse results for
the relevant spectral sequences. We conclude by exploring the
construction of A.\,E.~Mironov providing a vast family of
Lagrangian submanifolds with special minimality properties in
complex space, complex projective space and other toric varieties.
Like many other geometric constructions in this chapter, it builds
upon the realisation of the moment-angle manifold as an
intersection of quadrics.

In Chapter 7 we discuss several topological constructions of
even-dimensional manifolds with an effective action of a torus of
half the dimension of the manifold. They can be viewed as
topological analogues and generalisations of compact non-singular
toric varieties (or \emph{toric manifolds}). These include
\emph{quasitoric manifolds} of Davis and Januszkiewicz,
\emph{torus manifolds} of Hattori and Masuda, and
\emph{topological toric manifolds} of Ishida, Fukukawa and Masuda.
For all these classes of toric objects, the equivariant topology
of the action and the combinatorics of the orbit spaces interact
in a nice and peculiar way, leading to a host of results linking
topology with combinatorics. We also discuss the relationship with
GKM-manifolds (named after Goresky, Kottwitz and MacPherson),
another class of toric objects taking its origin in symplectic
topology.

Homotopy-theoretical aspects of toric topology are the subject of
Chapter~8. This is now a very active area. Homotopy techniques
brought to bear on the study of polyhedral products and other
toric spaces include model categories, homotopy limits and
colimits, higher Whitehead and Samelson products. The required
information about categorical methods in topology is collected in
Appendix~C.

The final Chapter~9 we review applications of toric methods in the
classical field of algebraic topology, complex bordism and
cobordism. It is a generalised cohomology theory that combines
both geometric intuition and elaborated algebraic techniques. The
toric viewpoint brings an entirely new perspective on complex
bordism theory in both its non-equivariant and equivariant
versions.

The later chapters require more specific knowledge of algebraic
topology, such as characteristic classes and spectral sequences,
for which we recommend respectively the classical book of Milnor
and Stasheff~\cite{mi-st79} and the excellent guide by
McCleary~\cite{mccl01}. Basic facts and constructions from bordism
and cobordism theory are given in Appendix~D, while the related
techniques of formal group laws and genera are reviewed in
Appendix~E.

\section*{Acknowledgments}
We wish to express our deepest thanks to
\begin{itemize}
\item[$\cdot$] our teacher Sergei Petrovich Novikov for encouragement and support
of our research on toric topology;

\item[$\cdot$] our friends and coathors Mikiya Masuda and Nigel Ray for long and much
fruitful collaboration;

\item[$\cdot$] %Anton Ayzenberg, Nickolai Erochovets, Alexander Gaifullin,
Peter Landweber
%, Ivan Limonchenko, Vadim Volodin
%Volkmar Welker
for most helpful suggestions on improving the text;

\item[$\cdot$] all our colleagues who participated in
conferences on toric topology for the insight gained from
discussions following talks and presentations.
\end{itemize}

The work was supported by the Russian Science Foundation (grant
no.~14-11-00414). We also thank the Russian Foundation for Basic
Research, the President of the Russian Federation Grants Council
and Dmitri Zimin's `Dynasty' foundation for their support of our
research related to this monograph.

%Проверить ссылку [{\cite[Proposition~2.22]{bu-er11}}]
%
%Приложения: разбиение на пермутаэдры (Гайфуллин).
%
%Добавить 2 картинки со связными суммами: в конструкцию и в пример
%
%Map of building sets induces a map of polytopes?
%
%Discussion of the space of polytopes: inequalities and convex hulls approaches, analogous polytopes
%
%Добавить раздел 1.2: проблема классификации и диаграммы Гейла,
%
%Уточнить (убрать) замечание перед теоремой 1.5.18.

\chapter{Geometry and combinatorics of polytopes}\label{combi}
This chapter is an introductory survey of the geometric and
combinatorial theory of convex polytopes, with the emphasis on
those of its aspects related to the topological applications later
in the book. We do not assume any specific knowledge of the reader
here. Algebraic definitions (graded rings and algebras) required
in the end of this chapter are contained in Section~\ref{algmod}
of the Appendix.

Convex polytopes have been studied since ancient times. Nowadays
both combinatorial and geometrical aspects of polytopes are
presented in numerous textbooks and monographs. Among them are the
classical monograph~\cite{grue03} by Gr\"unbaum and Ziegler's more
recent lectures~\cite{zieg95}. Face vectors and other
combinatorial topics are discussed in books by
McMullen--Shephard~\cite{mc-sh71}, Br\o nsted~\cite{bron83}, and
the survey article~\cite{kl-kl95} by Klee and Kleinschmidt; while
Yemelichev--Kovalev--Kravtsov~\cite{y-k-k84} focus on applications
to linear programming and optimisation. All these sources may be
recommended for the subsequent study of the theory of polytopes,
and contain a host of further references.

\section{Convex polytopes}\label{poly}
\subsection*{Definitions and basic constructions}
Let $\R^n$ be $n$-dimensional Euclidean space with the scalar
product $\langle\;,\,\rangle$. There are two constructively
different ways to define a convex polytope in~$\R^n$:

\begin{definition}
\label{pol1} A \emph{convex polytope} is the convex hull
$\mathop{\mathrm{conv}}(\mb v_1,\ldots,\mb v_q)$ of a finite set
of points $\mb v_1,\ldots,\mb v_q\in\R^n$.
\end{definition}

\begin{definition}
\label{pol2} A \emph{convex polyhedron} $P$ is a nonempty
intersection of finitely many half-spaces in some~$\R^n$:
\begin{equation}
\label{ptope}
  P=\bigl\{\mb x\in\R^n\colon\langle\mb a_i,\mb x\rangle+b_i\ge0\quad\text{for }
  i=1,\ldots,m\bigr\},
\end{equation}
where $\mb a_i\in\R^n$ and $b_i\in\R$. A \emph{convex polytope} is
a bounded convex polyhedron.
\end{definition}

All polytopes in this book will be convex. The two definitions
above produce the same geometrical object, i.e. a subset of $\R^n$
is the convex hull of a finite point set if and only if it is a
bounded intersection of finitely many half-spaces. This classical
fact is proved in many textbooks on polytopes and convex geometry,
and it lies at the heart of many applications of polytope theory
to linear programming and optimisation, see
e.g.~\cite[Theorem~1.1]{zieg95}.

The \emph{dimension} of a polyhedron is the dimension of its
affine hull.  We often abbreviate a `polyhedron of dimension $n$'
to \emph{$n$-polyhedron.} A \emph{supporting
hyperplane}\label{suphp} of $P$ is an affine hyperplane $H$ which
has common points with $P$ and for which the polyhedron is
contained in one of the two closed half-spaces determined by the
hyperplane. The intersection $P\cap H$ with a supporting
hyperplane is called a \emph{face}\label{facep} of the polyhedron.
%We also regard the polyhedron $P$ itself as a face; other faces
%are called \emph{proper}.
Denote by $\partial P$ and $\mathop{\mathrm{int}}
P=P\setminus\partial P$ the topological boundary and interior of
$P$ respectively. In the case $\dim P=n$ the boundary $\partial P$
is the union of all faces of~$P$. Each face of an $n$-polyhedron
($n$-polytope) is itself a polyhedron (polytope) of dimension~$\le
n$. Zero-dimensional faces are called \emph{vertices},
one-dimensional faces are \emph{edges}, and faces of codimension
one are \emph{facets}\label{facet}.

Two polytopes $P\subset\R^{n_1}$ and $Q\subset\R^{n_2}$ of the
same dimension are said to be \emph{affinely equivalent} (or
\emph{affinely isomorphic}) if there is an affine map
$\R^{n_1}\to\R^{n_2}$ establishing a bijection between the points
of the two polytopes. Two polytopes are \emph{combinatorially
equivalent}\label{combeq} if there is a bijection between their
faces preserving the inclusion relation. Note that two affinely
isomorphic polytopes are combinatorially equivalent, but the
opposite is not true.

The faces of a given polytope $P$ form a partially ordered set (a
\emph{poset}) with respect to inclusion. It is called the
\emph{face poset}\label{faceposet} of~$P$. Two polytopes are
combinatorially equivalent if and only if their face posets are
isomorphic.

\begin{definition}\label{combinptpe}
A \emph{combinatorial polytope} is a class of combinatorially
equivalent polytopes.
\end{definition}

Many topological constructions later in this book will depend only
on the combinatorial equivalence class of a polytope.
Nevertheless, it is always helpful, and sometimes necessary, to
keep in mind a particular geometric representative~$P$ rather than
thinking in terms of abstract posets. Depending on the context, we
shall denote by $P$, $Q$, etc., geometric polytopes or their
combinatorial equivalent classes (combinatorial polytopes).
Whenever we consider both geometric and combinatorial polytopes,
we shall use the notation $P\approx Q$ for combinatorial
equivalence.

%Most algebraic and topological constructions based on polytopes in
%toric topology will be eventually combinatorial, in the sense that
%the result, be it a ring or a space, will depend only on the
%combinatorial type of the polytope. In this respect the subject
%differs from toric geometry, where the constructions (such as
%projective toric varieties) depend significantly on the
%convex-geometrical structure of the polytope. Nevertheless, even
%in toric topology it is usually important, and sometimes
%necessary, to keep in mind a particular geometric representative
%$P$ within a given combinatorial type.

\medskip

We refer to~\eqref{ptope} as a \emph{presentation}\label{presptpe}
of the polyhedron~$P$ by inequalities. These inequalities contain
more information than the polyhedron~$P$, for the following
reason. It may happen that some of the inequalities $\langle\mb
a_i,\mb x\rangle+b_i\ge0$ can be removed from the presentation
without changing~$P$; we refer to such inequalities as
\emph{redundant}. A presentation without redundant inequalities is
called \emph{irredundant}. An irredundant presentation of a given
polyhedron is unique up to multiplication of pairs $(\mb a_i,b_i)$
by positive numbers.

%We shall assume our polytopes in $\R^n$ to be of full
%dimension~$n$, unless otherwise stated. Under this assumption
%there is no distinction between $\mathop{\mathrm{int}}P$ and the
%interior of $P\subset\R^n$ in the topological sense.
%
%Let $P$ be given by~(\ref{ptope}). We shall further assume, unless
%otherwise stated, that there are no redundant inequalities
%$\langle\mb a_i,\mb x\rangle+b_i\ge0$ in this representation. That
%is, no inequality can be removed from~(\ref{ptope}) without
%changing the polytope $P$. In this case $P$ has exactly $m$
%facets, which are given by
%\[
%  F_i=\bigl\{\mb x\in P\colon\langle\mb a_i,\mb x\rangle+b_i=0\bigr\}
%\]
%for $1\le i\le m$. The vector $\mb a_i$ is orthogonal to the facet
%$F_i$ and points towards the interior of the polytope.

\begin{example}[simplex and cube]\label{simcub}
An $n$-dimensional \emph{simplex} $\varDelta^n$ is the convex hull
of $n+1$ points in $\R^n$ that do not lie on a common affine
hyperplane. All faces of an $n$-simplex are simplices of
dimension~$\le n$. Any two $n$-simplices are affinely equivalent.
Let $\mb e_1,\ldots,\mb e_n$
%\[
%  \{\mb e_k=\mathop{(0,\ldots,0,1,0,\ldots,0)}\limits_k \quad
%  1\le k\le n\}
%\]
be the standard basis in~$\R^n$. The $n$-simplex
$\mathop{\mathrm{conv}}(\mathbf 0,\mb e_1,\ldots,\mb e_n)$ is
called \emph{standard}. Equivalently, the standard $n$-simplex is
specified by the $n+1$ inequalities
\begin{equation}
\label{stsim}
  x_i\ge0\quad\text{for }
  i=1,\ldots,n,\quad\text{and}\quad-x_1-\cdots-x_n+1\ge0.
\end{equation}
The \emph{regular} $n$-simplex is the convex hull of the endpoints
of $\mb e_1,\ldots,\mb e_{n+1}$ in~$\R^{n+1}$.

The \emph{standard $n$-cube} is given by
\begin{equation}
\label{cube}
  \I^n=[0,1]^n=\{(x_1,\ldots,x_n)\in\R^n\colon 0\le x_i\le1
  \quad\text{for }i=1,\ldots,n\}.
\end{equation}
Equivalently, the standard $n$-cube is the convex hull of $2^n$
points $(\varepsilon_1,\ldots,\varepsilon_n)\in\R^n$, where
$\varepsilon_i=0$ or~$1$.% for $1\le i\le n$.
Whenever we work with
combinatorial polytopes, we shall refer to any polytope
combinatorially equivalent to $\I^n$ as a \emph{cube}, and denote
it by~$I^n$.

The cube $\I^n$ has $2n$ facets. We denote by $F_k^0$ the facet
specified by the equation $x_k=0$, and by $F_k^1$ that specified
by the equation $x_k=1$, for $1\le k\le n$.
%Following the notation of the previous Agreement, we have
%$F_k=F_k^0$ and $F_{n+k}=F_k^1$.
\end{example}

\subsection*{Simple and simplicial polytopes. Polarity.}
The notion of a \emph{generic}\label{genpolyt} polytope depends on
the choice of definition. Below we describe the two possibilities.

A set of $q>n$ points in $\R^n$ is in \emph{general position} if
no $n+1$ of them lie on a common affine hyperplane. Now, assuming
Definition~\ref{pol1}, we may say that a polytope is generic if it
is the convex hull of a set of generally positioned points. This
implies that all faces of the polytope are simplices, i.e. every
facet has the minimal number of vertices (namely,~$n$). Polytopes
with this property are called \emph{simplicial}\label{scpolyt}.

Assuming Definition~\ref{pol2}, a presentation~\eqref{ptope} is
said to be \emph{generic} if $P$ is nonempty and the hyperplanes
defined by the equations $\langle\mb a_i,\mb x\rangle+b_i=0$ are
in general position at any point of~$P$ (that is, for any $\mb
x\in P$ the normal vectors $\mb a_i$ of the hyperplanes containing
$\mb x$ are linearly independent). If presentation~\eqref{ptope}
is generic, then $P$ is $n$-dimensional. If $P$ is a polytope,
then the existence of a generic presentation implies that $P$ is
\emph{simple}\label{simpleptpe}, that is, exactly $n$ facets meet
at each vertex of~$P$. Each face of a simple polytope is again a
simple polytope. Every vertex of a simple polytope has a
neighbourhood affinely equivalent to a neighbourhood of $\mathbf
0$ in the \emph{positive orthant}~$\R^n_\ge$. It follows that
every vertex is contained in exactly $n$ edges, and each subset of
$k$ edges with a common vertex spans a $k$-face.

%On the other hand, a set of $m>n$ hyperplanes $\langle\mb a_i,\mb
%x\rangle+b_i=0$, \ $\mb a_i\in\R^n$, \ $b_i\in\R$, \ $1\le i\le
%m$, is in \emph{general position} if no point belongs to more than
%$n$ hyperplanes. From the viewpoint of Definition~\ref{pol2}, a
%polytope $P$ is generic if its bounding hyperplanes are in general
%position. This implies that there are exactly $n$ facets meeting
%at each vertex of~$P$. Such polytopes are called
%\emph{simple}\label{sppolyt}.

A generic presentation may contain redundant inequalities, but,
for any such inequality, the intersection of the corresponding
hyperplane with $P$ is empty (i.e., the inequality is strict at
any $\mb x\in P$). We set
\[
  F_i=\bigl\{\mb x\in P\colon\langle\mb a_i,\mb x\rangle+b_i=0\bigr\}.
\]
If presentation~\eqref{ptope} is generic, then each $F_i$ is
either a facet of~$P$ or is empty.

The \emph{polar set} of a polyhedron $P\subset\R^n$ is defined as
\begin{equation}\label{polarset}
  P^*=\bigl\{\mb u\in\R^n\colon\langle\mb u,\mb x\,\rangle+1\ge0
  \:\text{ for all }\:\mb x\in P\bigr\}.
\end{equation}
The set $P^*$ is a convex polyhedron with $\mathbf 0\in P^*$.

\begin{remark}
$P^*$ is naturally a subset in the dual space $(\R^n)^*$, but we
shall not make this distinction until later on, assuming $\R^n$ to
be Euclidean. Also, in convex geometry the inequality $\langle\mb
u,\mb x\,\rangle\le1$ is usually used in the definition of
polarity, but the definition above is better suited for
applications in toric geometry. These two ways of defining the
polar set are taken into each other by the central symmetry.
\end{remark}

The following properties are well known in convex geometry:

{\samepage
\begin{theorem}[{see~\cite[\S2.9]{bron83} or~\cite[Th.~2.11]{zieg95}}]\label{polarity}\
\begin{itemize}
\item[(a)] $P^*$ is bounded if and only if $\mathbf 0\in\mathop{\mathrm{int}} P$;

\item[(b)] $P\subset (P^*)^*$, and $(P^*)^*=P$ if and only if $\mathbf 0\in
P$;

\item[(c)] if a polytope $Q$ is given as a convex hull,
$Q=\mathop{\mathrm{conv}(\mb a_1,\ldots,\mb a_m)}$, then $Q^*$ is
given by inequalities~\eqref{ptope} with $b_i=1$ for $1\le i\le
m$; in particular, $Q^*$ is a convex polyhedron, but not
necessarily bounded;

\item[(d)] if a polytope $P$ is given by inequalities~\eqref{ptope}
with $b_i=1$, then $P^*=\mathop{\mathrm{conv}(\mb a_1,\ldots,\mb
a_m)}$; furthermore, $\langle\mb a_i,\mb x\,\rangle+1\ge0$ is a
redundant inequality if and only if $\mb a_i\in\conv(\mb a_j\colon
j\ne i)$.
\end{itemize}
\end{theorem}
}

\begin{remark}
A polyhedron $P$ admits a presentation~\eqref{ptope} with $b_i=1$
if and only if $\mathbf 0\in\mathop{\mathrm{int}} P$. In general,
$(P^*)^*=\conv(P,{\bf0})$.
\end{remark}

\begin{example}
The difference between the situations $\mathbf 0\in P$ and
$\mathbf 0\in\mathop{\mathrm{int}} P$ may be illustrated by the
following example. Let $Q=\mathop{\mathrm{conv}}(\mathbf 0,\mb
e_1,\mb e_2)$ be the standard 2-simplex in $\R^2$. By
Theorem~\ref{polarity}~(c), $Q^*$ is specified by the three
inequalities
\[
  \langle\mathbf 0,\mb x\rangle+1\ge0,\quad\langle\mb e_1,\mb x\rangle+1\ge0,\quad
  \langle\mb e_2,\mb x\rangle+1\ge0,
\]
of which the first is satisfied for all $\mb x$, so we obtain an
unbounded polyhedron. Its dual is $\mathop{\mathrm{conv}}(\mathbf
0,\mb e_1,\mb e_2)$ by Theorem~\ref{polarity}~(d), giving back the
standard 2-simplex.
\end{example}

Any combinatorial polytope $P$ has a presentation~\eqref{ptope}
with $b_i=1$ (take the origin to the interior of~$P$ by a parallel
transform, an then divide each of the inequalities
in~\eqref{ptope} by the corresponding~$b_i$). Then $P^*$ is also a
polytope with ${\bf 0}\in P^*$, and $(P^*)^*=P$. We refer to the
combinatorial polytope $P^*$ as the \emph{dual} of the
combinatorial polytope~$P$. (We shall not introduce a new notation
for the dual polytope, keeping in mind that polarity is a
convex-geometric notion, while duality of polytopes is
combinatorial.)

\begin{theorem}[{see~\cite[\S2.10]{bron83}}]\label{dualpol}
If $P$ and $P^*$ are dual polytopes, then the face poset of $P^*$
is obtained from the face poset of $P$ by reversing the inclusion
relation.
\end{theorem}

As a corollary, we obtain that $P$ is simple if and only if its
dual polytope $P^*$ is simplicial. Any polygon is both simple and
simplicial.

\begin{proposition}\label{sssim}
In dimensions $n\ge3$ a simplex is the only polytope which is both
simple and simplicial.
\end{proposition}
\begin{proof}
Let be $P$ be a polytope which is both simple and simplicial.
Choose a vertex $v\in P$. Since $P$ is simple, $v$ is connected by
edges to exactly $n$ other vertices, say $v_1,\ldots,v_n$. We
claim that there are no other vertices in $P$. To prove this it is
enough to show that all vertices $v_1,\ldots,v_n$ are pairwise
connected by edges, as then every vertex from $v,v_1,\ldots,v_n$
will be connected to the remaining $n$ vertices. Indeed, take a
pair $v_i,v_j$. Since $P$ is simple and $v$ is connected to both
$v_i$ and $v_j$, all these three vertices belong to a 2-face.
Since $P$ is simplicial and $n\ge3$, this face is a 2-simplex, in
which $v_i$ and $v_j$ are connected by an edge. We conclude that
$P$ has $n+1$ vertices, so it is an $n$-simplex.
\end{proof}

The proof above also shows that if all 2-faces of a simple
polytope are triangular, then $P$ is a simplex. A similar property
is valid for a cube: if all 2-faces of a simple polytope $P$ are
quadrangular, then $P$ is a cube (an exercise).

\begin{example}\label{cross}
The dual of a simplex is again a simplex.  To describe the dual of
a cube, we consider the cube $[-1,1]^n$ (the standard
cube~\eqref{cube} is not good as 0 is not in its interior). Then
the polar set is the convex hull of the endpoints of $2n$ vectors
$\pm\mb e_k$, \ $1\le k\le n$. It is called the
\emph{cross-polytope}. The 3-dimensional cross-polytope is the
octahedron.
\end{example}

A combinatorial polytope $P$ is called
\emph{self-dual}\label{selfdual} if $P^*$ is combinatorially
equivalent to $P$. There are many examples of self-dual non-simple
polytopes; an infinite family of them is given by $k$-gonal
pyramids for $k\ge4$. Here is a more interesting \emph{regular}
example:

\begin{example}[24-cell]\label{24cell}
Let $Q$ be the 4-polytope obtained by taking the convex hull of
the following 24 points in $\R^4$: endpoints of 8 vectors $\pm\mb
e_i$, \ $1\le i\le 4$, and 16 points of the form
$(\pm\frac12,\pm\frac12,\pm\frac12,\pm\frac12)$. By
Theorem~\ref{polarity}~(c), the polar polytope $Q^*$ is given by
the following 24 inequalities:
\begin{equation}\label{24ineq}
  {}\mathbin\pm x_i+1\ge0\;\text{ for }i=1,\ldots,4,\quad\text{and
  }\;
  \textstyle{\frac12}({}\mathbin\pm x_1\mathbin\pm x_2\mathbin\pm x_3\mathbin\pm x_4)+1\ge0.
\end{equation}
Each of these inequalities turns into equality in exactly one of
the specified 24 points, so it defines a supporting hyperplane
whose intersection with $Q$ is only one point. This implies that
$Q$ has exactly 24 vertices. The vertices of $Q^*$ may be
determined by applying the `elimination process' to~\eqref{24ineq}
(see~\cite[\S1.2]{zieg95}), and as a result we obtain 24 points of
the form ${}\mathbin\pm\mb e_i\mathbin\pm\mb e_j$ for $1\le i<j\le
4$. Each supporting hyperplane defined by~\eqref{24ineq} contains
exactly 6 vertices of~$Q^*$, which form an octahedron. So both $Q$
and $Q^*$ have 24-vertices and 24 octahedral facets. In fact, both
$Q$ and $Q^*$ provide examples of a \emph{regular $4$-polytope}
called a \emph{$24$-cell}. It is the only regular self-dual
polytope different from a simplex. For more details on the 24-cell
and other regular polytopes see~\cite{coxe73}.
\end{example}

\subsection*{Products, hyperplane cuts and connected sums}
\begin{construction}[Product]\label{prodsp}
The \emph{product} $P_1\times P_2$ of two simple polytopes $P_1$
and $P_2$ is again a simple polytope. The dual operation on
simplicial polytopes can be described as follows. Let
$S_1\subset\R^{n_1}$ and $S_2\subset\R^{n_2}$ be two simplicial
polytopes. Assume that both $S_1$ and $S_2$ contain $0$ in their
interiors. Now define
\[
  S_1\circ S_2=\mathop{\mathrm{conv}}
  \bigl(S_1\times0\cup0\times S_2 \bigr)\subset\R^{n_1+n_2}.
\]
Then $S_1\circ S_2$ is again a simplicial polytope. For any two
simple polytopes $P_1$, $P_2$ containing $0$ in their interiors
the following identity holds:
\[
  P_1^*\circ P_2^*=(P_1\times P_2)^*.
\]
Both operations $\times$ and $\circ$ are also defined on
combinatorial polytopes; in this case the above formula holds
without any restrictions.
\end{construction}

\begin{construction}[Hyperplane cuts and face truncations]
\label{hypcut}
Assume given a simple polytope~\eqref{ptope} and a hyperplane
$H=\{\mb x\in\R^n\colon\langle\mb a,\mb x\rangle+b=0\}$ that does
not contain any vertex of~$P$. Then the intersections $P\cap
H_\ge$ and $P\cap H_\le$ of $P$ with either of the halfspaces
determined by $H$ are simple polytopes; we refer to them as
\emph{hyperplane cuts} of~$P$. To see that $P\cap H_\ge$ and
$P\cap H_\le$ are simple we note that their new vertices are
transverse intersections of $H$ with the edges of~$P$. Since $P$
is simple, each of those edges is contained in $n-1$ facets, so
each new vertex of $P\cap H_\ge$ or $P\cap H_\le$ is contained in
exactly~$n$ facets.

If $H$ separates all vertices of a certain $i$-face $G\subset P$
from the other vertices of~$P$ and $G\subset H_\ge$, then $P\cap
H_\ge$ is combinatorially equivalent to $G\times\varDelta^{n-i}$
(an exercise), and we say that the polytope $P\cap H_\le$ is
obtained from $P$ by a \emph{face truncation}.

In particular, if the cut off face $G$ is a vertex, the result is
a \emph{vertex truncation} of~$P$. When the choice of the cut off
vertex is clear or irrelevant we use the notation $\vt(P)$. We
also use the notation $\vt^k(P)$ for a polytope obtained from $P$
by a $k$-fold iteration of the vertex truncation.

We can describe the face poset of the simple polytope $\widetilde
P$ obtained by truncating $P$ at a face $G\subset P$ as follows.
Let $F_1,\ldots,F_m$ be the facets of $P$, and assume that
$G=F_{i_1}\cap\cdots\cap F_{i_k}$. The polytope $\widetilde P$ has
$m$ facets corresponding to $F_1,\ldots,F_m$ (and obtained from
them by truncation), which we denote by the same letters for
simplicity, and a new facet $F=P\cap H$. Then we have
\begin{align*}
  F_{j_1}\cap\cdots\cap F_{j_\ell}&\ne\varnothing
  \text{ in }\widetilde P\\
  &\Longleftrightarrow
  F_{j_1}\cap\cdots\cap F_{j_\ell}\ne\varnothing\text{ in }P,
  \text{ and }F_{j_1}\cap\cdots\cap F_{j_\ell}\not\subset G,\\
  F\cap F_{j_1}\cap\cdots\cap F_{j_\ell}&\ne\varnothing
  \text{ in }\widetilde P\\
  &\Longleftrightarrow
  G\cap F_{j_1}\cap\cdots\cap F_{j_\ell}\ne\varnothing\text{ in }P,
  \text{ and }F_{j_1}\cap\cdots\cap F_{j_\ell}\not\subset G.
\end{align*}
Note that $F_{j_1}\cap\cdots\cap F_{j_\ell}\not\subset G$ if and
only if $\{i_1,\ldots,i_k\}\not\subset\{j_1,\ldots,j_\ell\}$.
\end{construction}

The two previous construction worked for geometric polytopes. Here
is an example of a construction which is more suitable for
combinatorial ones.

\begin{construction}[Connected sum of polytopes]\label{consum}
Suppose we are given two simple polytopes $P$ and $Q$, both of
dimension $n$, with distinguished vertices $v$ and $w$
respectively. An informal way to obtain the \emph{connected sum}
$P\cs_{v,w}Q$ of $P$ at $v$ and $Q$ at $w$ is as follows. We cut
off $v$ from $P$ and $w$ from $Q$; then, after a projective
transformation, we can glue the rest of $P$ to the rest of $Q$
along the new simplex facets to obtain $P\cs_{v,w}Q$.  A more
formal definition is given below, following~\cite[\S6]{bu-ra01}.

First, we introduce an $n$-dimensional polyhedron $\varGamma$,
which will be used as a template for the construction; it arises
by considering the standard $(n-1)$-simplex $\varDelta^{n-1}$ in
the subspace $\{\mb x\colon x_1=0\}\subset\R^n$, and taking its
cartesian product with the first coordinate axis. The facets $G_r$
of $\varGamma$ therefore have the form $\R\times D_r$, where
$D_r$, $1\leq r\leq n$, are the facets of~$\varDelta^{n-1}$. Both
$\varGamma$ and the $G_r$ are divided into positive and negative
halves, determined by the sign of the coordinate~$x_1$.

We order the facets of $P$ meeting in $v$ as $E_1,\ldots,E_n$, and
the facets of $Q$ meeting in $w$ as $F_1,\ldots,F_n$. Denote the
complementary sets of facets by $\mathcal{C}_v$ and
$\mathcal{C}_w$; those in $\mathcal{C}_v$ avoid $v$, and those in
$\mathcal{C}_w$ avoid~$w$.

We now choose projective transformations $\phi_P$ and $\phi_Q$
of~$\R^n$, whose purpose is to map $v$ and $w$ to the infinity of
the $x_1$ axis.  We insist that $\phi_P$ embeds $P$ in $\varGamma$
so as to satisfy two conditions; firstly, that the hyperplane
defining $E_r$ is identified with the hyperplane defining~$G_r$,
for each $1\leq r\leq n$, and secondly, that the images of the
hyperplanes defining $\mathcal{C}_v$ meet $\varGamma$ in its
negative half. Similarly, $\phi_Q$ identifies the hyperplane
defining $F_r$ with that defining~$G_r$, for each $1\leq r\leq n$,
but the images of the hyperplanes defining $\mathcal{C}_w$ meet
$\varGamma$ in its positive half.  We define the \emph{connected
sum} $P\cs_{v,w}Q$ of $P$ at $v$ and $Q$ at $w$ to be the simple
$n$-polytope determined by the images of the hyperplanes defining
$\mathcal{C}_v$ and $\mathcal{C}_w$ and hyperplanes defining
$G_r$, $1\le r\le n$. It is defined only up to combinatorial
equivalence; moreover, different choices for either of $v$
and~$w$, or either of the orderings for $E_r$ and~$F_r$, are
likely to affect the combinatorial type. When the choices are
clear, or their effect on the result irrelevant, we use the
abbreviation $P\cs Q$.

The related construction of connected sum $P\cs S$ of a simple
polytope $P$ and a simplicial polytope $S$ is described
in~\cite[Example~8.41]{zieg95}.
\end{construction}

\begin{example}\label{expcs}\

1. If $P$ is an $m_1$-gon and $Q$ is an $m_2$-gon then $P\cs Q$ is
an $(m_1+m_2-2)$-gon.

2. If $P$ is an $n$-simplex, then $P\cs Q\approx\vt(P)$.

3. If both $P$ and $Q$ are $n$-simplices, then $P\cs
Q\approx\vt(\varDelta^n)\approx\varDelta^{n-1}\times\varDelta^1$.
The combinatorial type of $\vt(\varDelta^n)$ does not depend on
the choice of the cut off vertex. All the vertices of the
resulting polytope $\varDelta^{n-1}\times \varDelta^1$ are
equivalent, therefore, the combinatorial type of
$\vt^2(\varDelta^n)\approx\varDelta^n\cs\varDelta^n\cs\varDelta^n$
is still independent of the choices. The choice of the cut off
vertex becomes significant from the next step, i.e. for
$\vt^3(\varDelta^n)$, see Exercise~\ref{vc3}.
\end{example}

\subsection*{Neighbourly polytopes}
\begin{definition}\label{neighb}
A polytope is called \emph{$k$-neighbourly} if any set of its $k$
or fewer vertices spans a face. According to
Exercise~\ref{kneibex}, the only $n$-polytope which is more than
$\sbr n2$-neighbourly is a simplex. An $n$-polytope which is $\sbr
n2$-neighbourly is called simply \emph{neighbourly}.
%Likewise, a simple polytope $P$ is called
%\emph{$k$-coneighbourly} if any $k$ of its facets have non-empty
%intersection, which is a face of codimension~$k$.
\end{definition}

\begin{example}[neighbourly 4-polytope]
Let $P=\varDelta^2\times\varDelta^2$, the product of two
triangles. Then $P$ is simple, and it is easy to see that any two
facets of $P$ share a common 2-face. Therefore, any two vertices
of $P^*$ are connected by an edge, so $P^*$ is a neighbourly
simplicial 4-polytope.
\end{example}

More generally, if $P_1^*$ is $k_1$-neighbourly and $P_2^*$ is
$k_2$-neighbourly, then $(P_1\times P_2)^*$ is
$\min(k_1,k_2)$-neighbourly. It follows that
$(\varDelta^n\times\varDelta^n)^*$ and
$(\varDelta^n\times\varDelta^{n+1})^*$ are neighbourly $2n$- and
$(2n+1)$-polytopes respectively. The next example gives a
neighbourly polytope with an arbitrary number of vertices.

\begin{example}[cyclic polytopes]\label{cyclic}
The \emph{moment curve} in $\R^n$ is given by
\[
  \mb x\colon\R\longrightarrow\R^n,\qquad t\mapsto\mb x(t)=(t,t^2,\ldots,t^n)
  \in\R^n.
\]
For any $m>n$ define the \emph{cyclic polytope}
$C^n(t_1,\ldots,t_m)$ as the convex hull of $m$ distinct points
$\mb x(t_i)$, $t_1<t_2<\cdots<t_m$, on the moment curve.

\begin{theorem}\
\begin{itemize}
\item[(a)] $C^n(t_1,\ldots,t_m)$ is a simplicial $n$-polytope;

\item[(b)] $C^n(t_1,\ldots,t_m)$
has exactly $m$ vertices $\mb x(t_i)$, $1\le i\le m$;

\item[(c)] the combinatorial type of $C^n(t_1,\ldots,t_m)$ does not depend
on the specific choice of parameters $t_1,\ldots,t_m$;

\item[(d)]
$C^n(t_1,\ldots,t_m)$ is a neighbourly polytope.
\end{itemize}
\end{theorem}
\begin{proof}
This proof is taken from~\cite[Theorem~0.7]{zieg95}. Recall the
well-known Vandermonde determinant identity
%\begin{align*}
\[
  \det\begin{pmatrix}
    1&1&\cdots&1\\
    \mb x(t_0)& \mb x(t_1)&\cdots& \mb x(t_n)
  \end{pmatrix}
  =\left|\begin{array}{llll}
    1&1&\cdots&1\\
    t_0& t_1& \cdots & t_n\\
    \:\vdots &\:\vdots&\ddots&\:\vdots\\
    t_0^{n-1}& t_1^{n-1}& \cdots & t_n^{n-1}\\
    t_0^n& t_1^n& \cdots & t_n^n
  \end{array}\right|=\prod_{0\le i<j\le n}(t_j-t_i).
\]
%\end{align*}
This implies that no $n+1$ points on the moment curve belong to a
common affine hyperplane, proving~(a). Denote
$[m]=\{1,\ldots,m\}$. Properties (b) and (c) follow from the
following statement: an $n$-element subset $\omega\subset[m]$
corresponds to the vertex set of a facet of $C^n(t_1,\ldots,t_m)$
if and only if the following `Gale's evenness condition' is
satisfied:
\begin{quotation}\it
If elements $i<j$ are not in $\omega$, then the number of elements
$k\in\omega$ between $i$ and $j$ is even.
\end{quotation}
To prove this we write $\omega=\{i_1,\ldots,i_n\}$ and consider
the hyperplane $H_\omega$ through the corresponding points $\mb
x(t_{i_s})$, $1\le s\le n$, on the moment curve. We have
\[
  H_\omega=\{\mb x\in\R^n\colon F_\omega(\mb x)=0\},
\]
where
\[
  F_\omega(\mb x)=\det\begin{pmatrix}
    1&1&\ldots&1\\
    \mb x& \mb x(t_{i_1})&\ldots& \mb x(t_{i_n})
  \end{pmatrix}.
\]
(The latter is exactly the linear function vanishing on the
prescribed points.) Now let the point $\mb x(t)$ move on the
moment curve. Then $F_\omega(\mb x(t))$ is a polynomial in $t$ of
degree~$n$. It has $n$ different roots $t_{i_1},\ldots,t_{i_n}$,
and changes sign at each of them. Now $\omega$ corresponds to the
vertex set of a facet if and only if $F_\omega(\mb x(t_i))$ has
the same sign for all the points $\mb x(t_i)$ with
$i\notin\omega$; that is, if $F_\omega(\mb x(t))$ has an even
number of sign changes between $t=t_i$ and $t=t_j$, for  $i<j$ and
$i,j\notin\omega$. This proves Gale's condition, and
statements~(b) and~(c).

It remains to prove~(d). We need to check that any subset
$\tau=\{i_1,\ldots,i_k\}\subset[m]$ of cardinality
$k\le\left[\frac n2\right]$ corresponds to the vertex set of a
face. Choose some $\varepsilon>0$ so that
$t_i<t_i+\varepsilon<t_{i+1}$ for all $i<m$, and some
$N>t_m+\varepsilon$. Define a linear function $F_\tau(\mb x)$ as
\[
  \det\bigl(\mb x, \mb x(t_{i_1}), \mb x(t_{i_1}+\varepsilon),\ldots,
  \mb x(t_{i_k}), \mb x(t_{i_k}+\varepsilon), \mb x(N+1), \ldots,
  \mb x(N+n-2k)\bigr).
\]
It vanishes on $\mb x(t_i)$ with $i\in\tau$. Now $F_\tau(\mb
x(t))$ is a polynomial in $t$ of degree $n$, and it has $n$
different roots
\[
  t_{i_1},t_{i_1}+\varepsilon,\ldots,t_{i_k},t_{i_k}+\varepsilon,N+1,\ldots,
  N+n-2k.
\]
If $i,j\notin\tau$, then there is an even number of roots between
$t=t_i$ and $t=t_j$, because a root $t=t_l$ always come in a pair
with the root $t=t_l+\varepsilon$. Thus, the linear function
$F_\tau(\mb x)$ has the same sign on all the points $\mb x(t_i)$
with $i\notin\tau$. This linear function defines a supporting
hyperplane, so $\tau$ corresponds to the vertex set of a face.
\end{proof}

We shall denote the combinatorial cyclic $n$-polytope with $m$
vertices by~$C^n(m)$.
\end{example}

The vertices and edges of a polytope $P$ determine a graph, which
is called the \emph{graph of polytope}\label{graphptope} and
denoted $\Gamma(P)$. This graph is \emph{simple}, that is, it has
no loops and multiple edges. The following theorem is due to Blind
and Mani, see also~\cite[\S3.4]{zieg95} for a simpler proof given
by Kalai.

\begin{theorem}
The combinatorial type of a simple polytope $P$ is determined by
its graph $\Gamma(P)$. In other words, two simple polytopes are
combinatorially equivalent if their graphs are isomorphic.
\end{theorem}

This theorem fails for general polytopes: the graph of a
neighbourly polytope is isomorphic to that of a simplex with the
same number of vertices. In general, simplicial $n$-polytopes are
determined by their $[\frac n2]$-skeleta. General $n$-polytopes
are determined by their $(n-2)$-skeleta. See~\cite[\S3.4]{zieg95}
for more history and references.

\subsection*{Exercises}
\begin{exercise}\label{ptopesubposet}
Show that if $P$ and $Q$ are $n$-polytopes, and the face poset of
$P$ is a subposet of the face poset of $Q$, then $P$ and $Q$ are
combinatorially equivalent.
\end{exercise}

\begin{exercise}\label{polvertex}
The polyhedron $P$ defined by~\eqref{ptope} has a vertex if and
only if the vectors $\mb a_1,\ldots,\mb a_m$ span the
whole~$\R^n$.
\end{exercise}

\begin{exercise}\label{2facecube}
Show that a simple $n$-polytope all of whose 2-faces are
quadrangular is combinatorially equivalent to an $n$-cube.
\end{exercise}

\begin{exercise}
Show that any hyperplane cut of $\varDelta^n$ is combinatorially
equivalent to a product of two simplices. Conclude that any
combinatorial simple $n$-polytope with $n+2$ facets is
combinatorially equivalent (in fact, projectively equivalent) to a
product of two simplices.
\end{exercise}

\begin{exercise}
Let $P$ be a simple polytope. Show that if a hyperplane $H$
separates all vertices of a certain $i$-face $G\subset P$ from the
other vertices of~$P$ and $G\subset H_\ge$, then $P\cap
H_\ge\approx G\times\varDelta^{n-i}$.
\end{exercise}

\begin{exercise}\label{vc3}
How many combinatorially different polytopes may be obtained
as~$\vt^3(\varDelta^n)$?
\end{exercise}

\begin{exercise}\label{kneibex}
Show that if a polytope is $k$-neighbourly, then every
$(2k-1)$-face is a simplex. Conclude that if an $n$-polytope is
$(\sbr n2+1)$-neighbourly, then it is a simplex. Conclude also
that a neighbourly $2k$-polytope is simplicial. Are there
neighbourly non-simplicial polytopes of odd dimension?
\end{exercise}

\begin{exercise}
Are the polytopes $(\varDelta^n\times\varDelta^n)^*$ and
$(\varDelta^n\times\varDelta^{n+1})^*$ combinatorially equivalent
to cyclic polytopes?
\end{exercise}

\section{Gale duality and Gale diagrams}\label{galed}
The following construction realises any polytope~\eqref{ptope} of
dimension~$n$ by the intersection of the orthant
\begin{equation}
\label{pcone}
  \R^m_\ge=\bigl\{(y_1,\ldots,y_m)\in\R^m\colon y_i\ge0
  \quad\text{for }
  i=1,\ldots,m\bigr\}\subset\R^m
\end{equation}
with an affine $n$-plane.

\begin{construction}\label{dist}
Let~\eqref{ptope} be a presentation of a polyhedron. Consider the
linear map $A\colon\R^m\to\R^n$ sending the $i$th basis vector
$\mb e_i$ to~$\mb a_i$. It is given by the $n\times m$-matrix
(which we also denote by~$A$) whose columns are the vectors $\mb
a_i$ written in the standard basis of~$\R^n$. The dual map
$A^*\colon\R^n\to\R^m$ is given by $\mb x\mapsto(\langle\mb
a_1,\mb x\rangle,\ldots,\langle\mb a_m,\mb x\rangle)$. Note that
$A$ is of rank~$n$ if and only if the polyhedron $P$ has a vertex
(e.g., when $P$ is a polytope, see Exercise~\ref{polvertex}).
Also, let $\mb b=(b_1,\ldots,b_m)^t\in\R^m$ be the column vector
of~$b_i$s. Then we can write~\eqref{ptope} as
\[
  P=P(A,\mb b)=\bigl\{\mb x\in\R^n\colon(A^*\mb x+\mb b)_i\ge 0\quad
  \text{for }i=1,\ldots,m\},
\]
where $\mb x=(x_1,\ldots,x_n)^t$ is the column of coordinates.
Consider the affine map
\[
  i_{A,\mb b}\colon \R^n\to\R^m,\quad i_{A,\mb b}(\mb x)=A^*\mb x+\mb b
  =\bigl(\langle\mb a_1,\mb x\rangle+b_1,\ldots,\langle\mb a_m,\mb x\rangle+b_m\bigr)^t.
\]
If $P$ has a vertex, then the image of $\R^n$ under $i_{A,\mb b}$
is an $n$-dimensional affine plane in $\R^m$, which we can write
by $m-n$ linear equations:
\begin{equation}\label{iabrn}
\begin{aligned}
  i_{A,\mb b}(\R^n)&=\{\mb y\in\R^m\colon\mb y=A^*\mb x+\mb b\quad
  \text{for some }\mb x\in\R^n\}\\
  &=\{\mb y\in\R^m\colon\varGamma\mb y=\varGamma\mb b\},
\end{aligned}
\end{equation}
where $\varGamma=(\gamma_{jk})$ is an $(m-n)\times m$-matrix whose
rows form a basis of linear relations between the vectors~${\mb
a}_i$. That is, $\varGamma$ is of full rank and satisfies the
identity $\varGamma A^*=0$.

The polytopes $P$ and $i_{A,\mb b}(P)=\R^m_\ge\cap i_{A,\mb
b}(\R^n)$ are affinely equivalent.
\end{construction}

\begin{example}\label{simdist}
Consider the standard $n$-simplex $\varDelta^n\subset\R^n$,
see~\eqref{stsim}. It is given by~\eqref{ptope} with $\mb a_i=\mb
e_i$ (the $i$th standard basis vector) for $i=1,\ldots,n$ and $\mb
a_{n+1}=-\mb e_1-\cdots-\mb e_n$; \ $b_1=\cdots=b_n=0$ and
$b_{n+1}=1$. We may take $\varGamma=(1,\ldots,1)$ in
Construction~\ref{dist}. Then $\varGamma\mb y=y_1+\cdots+y_m$, \
$\varGamma\mb b=1$, and we have
\[
  i_{A,\mb b}(\varDelta^n)=\bigl\{\mb y\in\R^{n+1}\colon y_1+\cdots+y_{n+1}=1,\;
  y_i\ge0\quad\text{for }i=1,\ldots,n\bigr\}.
\]
This is the regular $n$-simplex in $\R^{n+1}$.
\end{example}

\begin{construction}[Gale duality]\label{galeduality}
Let $\mb a_1,\ldots,\mb a_m$ be a configuration of vectors that
span the whole~$\R^n$. Form an $(m-n)\times m$-matrix
$\varGamma=(\gamma_{jk})$ whose rows form a basis in the space of
linear relations between the vectors~${\mb a}_i$. The set of
columns $\gamma_1,\ldots,\gamma_m$ of $\varGamma$ is called a
\emph{Gale dual} configuration of $\mb a_1,\ldots,\mb a_m$. The
transition from the configuration of vectors $\mb a_1,\ldots,\mb
a_m$ in $\R^n$ to the configuration of vectors
$\gamma_1,\ldots,\gamma_m$ in $\R^{m-n}$ is called the (linear)
\emph{Gale transform}. Each configuration determines the other
uniquely up to isomorphism of its ambient space. In other words,
each of the matrices $A$ and $\varGamma$ determines the other
uniquely up to multiplication by an invertible matrix from the
left.

Using the coordinate-free notation, we may think of $\mb
a_1,\ldots,\mb a_m$ as a set of linear functions on an
$n$-dimensional space~$W$. Then we have an exact sequence
\[
  0\to W\stackrel
  {A^*}\longrightarrow\R^m\stackrel{\varGamma}\longrightarrow L\to0,
\]
where $A^*$ is given by $\mb x\mapsto \bigl(\langle\mb a_1,\mb
x\rangle,\ldots,\langle\mb a_m,\mb x\rangle\bigr)$, and the map
$\varGamma$ takes $\mb e_i$ to $\gamma_i\in L\cong\R^{m-n}$.
Similarly, in the dual exact sequence
\[
  0\to L^*\stackrel
  {\varGamma^*}\longrightarrow\R^m\stackrel{A}\longrightarrow W^*\to0,
\]
the map $A$ takes $\mb e_i$ to $\mb a_i\in W^*\cong\R^n$.
Therefore, two configurations $\mb a_1,\ldots,\mb a_m$ and
$\gamma_1,\ldots,\gamma_m$ are Gale dual if they are obtained as
the images of the standard basis of $\R^m$ under the maps $A$ and
$\varGamma$ in a pair of dual short exact sequences.
\end{construction}

Here is an important property of Gale dual configurations:

\begin{theorem}\label{galespan}
Let $\mb a_1,\ldots,\mb a_m$ and $\gamma_1,\ldots,\gamma_m$ be
Gale dual configurations of vectors in $\R^n$ and $\R^{m-n}$
respectively, and let $I=\{i_1,\ldots,i_k\}\subset[m]$. Then the
subset $\{\mb a_i\colon i\in I\}$ is linearly independent if and
only if the subset $\{\gamma_j\colon j\notin I\}$ spans the
whole~$\R^{m-n}$.
\end{theorem}

The proof uses an algebraic lemma:

\begin{lemma}\label{2ses}
Let $\k$ be a field or $\Z$, and assume given a diagram
\[
\begin{CD}
  @.@.\begin{array}{c}0\\ \raisebox{2pt}{$\downarrow$}\\
  U\end{array}@.@.\\
  @.@.@VV{i_1}V@.@.\\
  0@>>> R@>{i_2}>> S@>{p_2}>> T @>>>0\\
  @.@.@VV {p_1} V@.@.\\
  @.@.\begin{array}{c} V\\ \downarrow\\ 0\end{array}@.@.\\
\end{CD}
\]
in which both vertical and horizontal lines are short exact
sequences of $\k$-vector spaces or free abelian groups. Then
$p_1i_2$ is surjective (respectively, injective or split
injective) if and only if $p_2i_1$ is surjective (respectively,
injective or split injective).
\end{lemma}
\begin{proof}
This is a simple diagram chase. Assume first that $p_1i_2$ is
surjective. Take $\alpha\in T$; we need to cover it by an element
in~$ U$. There is $\beta\in S$ such that $p_2(\beta)=\alpha$. If
$\beta\in i_1( U)$, then we are done. Otherwise set
$\gamma=p_1(\beta)\ne0$. Since $p_1i_2$ is surjective, we can
choose $\delta\in R$ such that $p_1i_2(\delta)=\gamma$. Set
$\eta=i_2(\delta)\ne0$. Hence, $p_1(\eta)=p_1(\beta)(=\gamma)$ and
there is $\xi\in U$ such that $i_1(\xi)=\beta-\eta$. Then
$p_2i_1(\xi)=p_2(\beta-\eta)=\alpha-p_2i_2(\delta)=\alpha$. Thus,
$p_2i_1$ is surjective.

Now assume that $p_1i_2$ is injective. Suppose $p_2i_1(\alpha)=0$
for a nonzero $\alpha\in U$. Set $\beta=i_1(\alpha)\ne0$. Since
$p_2(\beta)=0$, there is a nonzero $\gamma\in R$ such that
$i_2(\gamma)=\beta$. Then
$p_1i_2(\gamma)=p_1(\beta)=p_1i_1(\alpha)=0$. This contradicts the
assumption that $p_1i_2$ is injective. Thus, $p_2i_1$ is
injective.

Finally, if $p_1i_2$ is split injective, then its dual map
$i_2^*p_1^*\colon V^*\to R^*$ is surjective. Then
$i_1^*p_2^*\colon T^*\to U^*$ is also surjective. Thus, $p_2i_1$
is split injective.
\end{proof}

\begin{proof}[Proof of Theorem~\ref{galespan}]
Let $A$ be the $n\times m$-matrix with column vectors $\mb
a_1,\ldots,\mb a_m$, and let $\varGamma$ be the $(m-n)\times
m$-matrix with columns $\gamma_1,\ldots,\gamma_m$. Denote by $A_I$
the $n\times k$-submatrix formed by the columns $\{\mb a_i\colon
i\in I\}$ and denote by $\varGamma_{\widehat I}$ the
$(m-n)\times(m-k)$-submatrix formed by the columns
$\{\gamma_j\colon j\notin I\}$. We also consider the corresponding
maps $A_I\colon\R^k\to\R^n$ and $\varGamma_{\widehat
I}\colon\R^{m-k}\to\R^{m-n}$.

Let $i\colon\R^k\rightarrow\R^m$ be the inclusion of the
coordinate subspace spanned by the vectors $\mb e_i$, $i\in I$,
and let $p\colon\R^m\rightarrow\R^{m-k}$ the projection sending
every such $\mb e_i$ to zero. Then $A_I=A\cdot i$ and
$\varGamma_{\widehat I}^*=p\cdot\varGamma^*$. The vectors $\{\mb
a_i\colon i\in I\}$ are linearly independent if and only if
$A_I=A\cdot i$ is injective, and the vectors $\{\gamma_j\colon
j\notin I\}$ span~$\R^{m-n}$ if and only if $\varGamma_{\widehat
I}^*=p\cdot\varGamma^*$ is injective. These two conditions are
equivalent by Lemma~\ref{2ses}.
\end{proof}

\begin{construction}[Gale diagram]\label{galediagra}
Let $P$ be a polytope~\eqref{ptope} with $b_i=1$
%(i.e. $\mathbf 0\in\mathop{\mathrm{int}}P$)
and let $P^*=\conv(\mb a_1,\ldots,\mb a_m)$ be the polar polytope.
Let $\widetilde A^*=(A^*\;\mathbf 1)$ be the $m\times(n+1)$-matrix
obtained by appending a column of units to~$A^*$. The matrix
$\widetilde A^*$ has full rank $n+1$ (indeed, otherwise there is
$\mb x\in\R^n$ such that $\langle\mb a_i,\mb x\rangle=1$ for
all~$i$, and then $\lambda\mb x$ is in~$P$ for any $\lambda\ge1$,
so that $P$ is unbounded). A configuration of vectors $G=(\mb
g_1,\ldots,\mb g_m)$ in $\R^{m-n-1}$ which is Gale dual to
$\widetilde A$ is called a \emph{Gale diagram} of~$P^*$. A Gale
diagram $G=(\mb g_1,\ldots,\mb g_m)$ of $P^*$ is therefore
determined by the conditions
\[
  GA^*=0,\quad \rank G=m-n-1,\quad\text{and }
  \sum_{i=1}^m\mb g_i=\mathbf 0.
\]
The rows of the matrix $G$ from a basis of affine dependencies
between the vectors $\mb a_1,\ldots,\mb a_m$, i.e. a basis in the
space of $\mb y=(y_1,\ldots,y_m)^t$ satisfying
\[
  y_1\mb a_1+\cdots+y_m\mb a_m=\mathbf0,\quad
  y_1+\cdots+y_m=0.
\]
\end{construction}

%In the case of polytopes (bounded polyhedra), all equations in the
%system $\varGamma\mb y=\varGamma\mb b$ defining the affine plane
%$i_{A,\mb b}(\R^n)$ in~\eqref{zgamma} except one can be made
%homogeneous:

\begin{proposition}\label{propcf}
The polyhedron $P=P(A,\mb b)$ is bounded if and only if the matrix
$\varGamma=(\gamma_{jk})$ can be chosen so that the affine plane
$i_{A,\mb b}(\R^n)$ is given by
\begin{equation}\label{bplane}
  i_{A,\mb b}(\R^n)= \left\{\begin{array}{ll}
  \mb y\in\R^m\colon&\gamma_{11}y_1+\cdots+\gamma_{1m}y_m=c,\\[1mm]
  &\gamma_{j1}y_1+\cdots+\gamma_{jm}y_m=0,\quad%\text{for }
  2\le j\le m-n.
  \end{array}\right\},
\end{equation}
where $c>0$ and $\gamma_{1k}>0$ for all $k$.

Furthermore, if $b_i=1$ in~\eqref{ptope}, then the vectors $\mb
g_i=(\gamma_{2i},\ldots,\gamma_{m-n,i})^t$, $i=1,\ldots,m$, form a
Gale diagram of the polar polytope $P^*=\conv(\mb a_1,\ldots,\mb
a_m)$.
\end{proposition}
\begin{proof}
The image $i_{A,\mb b}(P)$ is the intersection of the $n$-plane
$L=i_{A,\mb b}(\R^n)$ with $\R^m_\ge$. It is bounded if and only
if $L_0\cap\R_\ge^m=\{\mathbf0\}$, where $L_0$ is the $n$-plane
through~$\mathbf 0$ parallel to~$L$. Choose a hyperplane $H_0$
through $\mathbf 0$ such that $L_0\subset H_0$ and
$H_0\cap\R_\ge^m=\{\mathbf 0\}$. Let $H$ be the affine hyperplane
parallel to $H_0$ and containing~$L$. Since $L\subset H$, we may
take the equation defining $H$ as the first equation in the system
$\varGamma\mb y=\varGamma\mb b$ defining~$L$. The conditions on
$H_0$ imply that $H\cap\R_\ge^m$ is nonempty and bounded, that is,
$c>0$ and $\gamma_{1k}>0$ for all~$k$. Now, subtracting the first
equation from the other equations of the system $\varGamma\mb
y=\varGamma\mb b$ with appropriate coefficients, we achieve that
the right hand sides of the last $m-n-1$ equations become zero.

To prove the last statement, we observe that in our case
\[
  \varGamma=\begin{pmatrix}
  \gamma_{11}&\cdots&\gamma_{1m}\\
  \mb g_1&\cdots&\mb g_m\end{pmatrix}.
\]
The conditions $\varGamma A^t=0$ and $\rank\varGamma=m-n$ imply
that $GA^t=0$ and $\rank G=m-n-1$. Finally, by
comparing~\eqref{iabrn} with~\eqref{bplane} we obtain
$\varGamma\mb b=\begin{pmatrix}c\\\mathbf 0\end{pmatrix}$. Since
$b_i=1$, we get $\sum_{i=1}^m\mb g_i=\mathbf 0$. Thus, $G=(\mb
g_1,\ldots,\mb g_m)$ is a Gale diagram of~$P^*$.
\end{proof}

\begin{corollary}\label{Pbound}
A polyhedron $P=P(A,\mb b)$ is bounded if and only if the vectors
$\mb a_1,\ldots,\mb a_m$ satisfy $\alpha_1\mb
a_1+\cdots+\alpha_m\mb a_m=\mathbf0$ for some positive
numbers~$\alpha_k$.
\end{corollary}
\begin{proof}
If $\mb a_1,\ldots,\mb a_m$ satisfy $\sum_{k=1}^m\alpha_k\mb
a_k=\mathbf0$ with positive~$\alpha_k$, then we can take
$\sum_{k=1}^m\alpha_k y_k=\sum_{k=1}^m\alpha_k b_k$ as the first
equation defining the $n$-plane $i_{A,\mb b}(\R^n)$ in~$\R^m$. It
follows that $i_{A,\mb b}(P)$ is contained in the intersection of
the hyperplane $\sum_{k=1}^m\alpha_k y_k=\sum_{k=1}^m\alpha_k b_k$
with $\R^m_\ge$, which is bounded since all $\alpha_k$ are
positive. Therefore, $P$ is bounded.

Conversely, if $P$ is bounded, then it follows from
Proposition~\ref{propcf} and Gale duality that $\mb a_1,\ldots,\mb
a_m$ satisfy $\gamma_{11}\mb a_1+\cdots+\gamma_{1m}\mb
a_m=\mathbf0$ with $\gamma_{1k}>0$.
\end{proof}
A Gale diagram of $P^*$ encodes its combinatorics (and the
combinatorics of~$P$) completely. We give the corresponding
statement in the generic case only:

\begin{proposition}\label{galecomb}
Assume that~\eqref{ptope} is a generic presentation with $b_i=1$.
Let $P^*=\conv(\mb a_1,\ldots,\mb a_m)$ be the polar simplicial
polytope and $G=(\mb g_1,\ldots,\mb g_m)$ be its Gale diagram.
Then the following conditions are equivalent:
\begin{itemize}
\item[(a)]
$F_{i_1}\cap\cdots\cap F_{i_k}\ne\varnothing$ in~$P=P(A,\mathbf
1)$;

\item[(b)] $\conv(\mb a_{i_1},\ldots,\mb a_{i_k})$ is a face
of~$P^*$;

\item[(c)]
$\mathbf 0\in\conv\bigr(\mb g_j\colon
j\notin\{i_1,\ldots,i_k\}\bigl)$.
\end{itemize}
\end{proposition}
\begin{proof}
The equivalence (a)$\Leftrightarrow\,$(b) follows from
Theorems~\ref{polarity} and~\ref{dualpol}.

(b)$\Rightarrow\,$(c). Let $\conv(\mb a_{i_1},\ldots,\mb a_{i_k})$
be a face of~$P^*$. We first observe that each of $\mb
a_{i_1},\ldots,\mb a_{i_k}$ is a vertex of this face, as otherwise
presentation~\eqref{ptope} is not generic. By definition of a
face, there exists a linear function $\xi$ such that $\xi(\mb
a_j)=0$ for $j\in\{i_1,\ldots,i_k\}$ and $\xi(\mb a_j)>0$
otherwise. The condition $\mathbf 0\in\mathop{\mathrm{int}}P^*$
implies that $\xi({\bf 0})>0$, and we may assume that $\xi$ has
the form $\xi(\mb u)=\langle\mb u,\mb x\rangle+1$ for some $\mb
x\in\R^n$. Set $y_j=\xi(\mb a_j)=\langle\mb a_j,\mb x\rangle+1$,
i.e. $\mb y=A^*\mb x+\mathbf 1$. We have
\[
  \sum_{j\notin\{i_1,\ldots,i_k\}}\mb g_jy_j=
  \sum_{j=1}^m\mb g_j y_j=G\mb y=G(A^*\mb x+\mathbf 1)
  =G\mathbf 1=\sum_{j=1}^m\mb g_j=\mathbf 0,
\]
where $y_j>0$ for $j\notin\{i_1,\ldots,i_k\}$. It follows that
$\mathbf 0\in\conv(\mb g_j\colon j\notin\{i_1,\ldots,i_k\})$.

(c)$\Rightarrow\,$(b). Let $\sum_{j\notin\{i_1,\ldots,i_k\}}\mb
g_jy_j=\mathbf 0$ with $y_j\ge0$ and at least one $y_j$ nonzero.
This is a linear relation between the vectors~$\mb g_j$. The space
of such linear relations has basis formed by the columns of the
matrix $\widetilde A^*=(A^*\;\mathbf 1)$
%, whose rows are $(\mb a_j\;1)$.
Hence, there exists $\mb x\in\R^n$ and $b\in\R$ such that
$y_j=\langle\mb a_j,\mb x\rangle+b$. The linear function $\xi(\mb
u)=\langle\mb u,\mb x\rangle+b$ takes zero values on $\mb a_j$
with $j\in\{i_1,\ldots,i_k\}$ and takes nonnegative values on the
other~$\mb a_j$. Hence, $\mb a_{i_1},\ldots,\mb a_{i_k}$ is a
subset of the vertex set of a face. Since $P^*$ is simplicial,
$\mb a_{i_1},\ldots,\mb a_{i_k}$ is a vertex set of a face.
\end{proof}

\begin{remark}
We allow redundant inequalities in Proposition~\eqref{galecomb}.
In this case we obtain the equivalences
\[
  F_i=\varnothing\quad\Leftrightarrow\quad
  \mb a_i\in\conv(\mb a_j\colon j\ne i)\quad\Leftrightarrow\quad
  \mathbf 0\notin\conv(\mb g_j\colon j\ne i).
\]
\end{remark}

A configuration of vectors $G=(\mb g_1,\ldots,\mb g_m)$ in
$\R^{m-n-1}$ with the property
\[
  \mathbf 0\in\conv\bigr(\mb g_j\colon
  j\notin\{i_1,\ldots,i_k\}\bigl)\quad\Leftrightarrow\quad
  \conv(\mb
  a_{i_1},\ldots,\mb a_{i_k})\text{ is a face of }P^*
\]
is called a \emph{combinatorial Gale diagram}\label{combgalediag}
of~$P^*=\conv(\mb a_1,\ldots,\mb a_m)$. For example, a
configuration obtained by multiplying each vector in a Gale
diagram by a positive number is a combinatorial Gale diagram.
Also, the vectors of a combinatorial Gale diagram can be moved as
long as the origin does not cross the `walls', i.e. affine
hyperplanes spanned by subsets of $\mb g_1,\ldots,\mb g_m$. A
combinatorial Gale diagram of $P^*$ is a Gale diagram of a
polytope which is combinatorially equivalent to~$P^*$.

\begin{example}\label{gdexam}\

1. The Gale diagram of $\varDelta^n$ (when $m=n+1$) consists of
$n+1$ points $\bf 0$ in~$\R^0$, i.e. $\mb g_i=\bf0$,
$i=1,\ldots,m$.

2. Let $P=\varDelta^{p-1}\times\varDelta^{q-1}$, $p+q=m$, i.e.
$m=n+2$. Then $P^*$ is a hyperplane cut of a simplex
$\varDelta^{m-2}$ by a hyperplane that separates some $p-1$ of its
vertices from the other~$q-1$. A combinatorial Gale diagram of
$P^*$ consists of $p$ points $1\in\R^1$ and $q$ points
$-1\in\R^1$. The cases $p=1$ or $q=1$ correspond to a presentation
of $\varDelta^{m-2}$ with one redundant inequality.

3. A combinatorial Gale diagram of a pentagon ($m=n+3$) is shown
in Fig.~\ref{cgd5g}. The property that $\conv(\mb a_1,\mb a_2)$ is
a face translates to ${\bf0}\in\conv(\mb g_3,\mb g_4,\mb g_5)$.
\begin{figure}[h]
\unitlength=0.8mm
  \begin{center}
  \begin{picture}(100,35)
  \put(7,-3){$\mb a_5$}
  \put(-4,18){$\mb a_{\!1}$}
  \put(18,37){$\mb a_2$}
  \put(41,18){$\mb a_3$}
  \put(30,-3){$\mb a_4$}
  \put(0,20){\circle*{1.5}}
  \put(20,35){\circle*{1.5}}
  \put(40,20){\circle*{1.5}}
  \put(30,0){\circle*{1.5}}
  \put(10,0){\circle*{1.5}}
  \put(10,0){\line(1,0){20}}
  \put(10,0){\line(-1,2){10}}
  \put(0,20){\line(4,3){20}}
  \put(20,35){\line(4,-3){20}}
  \put(40,20){\line(-1,-2){10}}
  \put(67.5,-3){$\mb g_3$}
  \put(56,18){$\mb g_{\!1}$}
  \put(79,37){$\mb g_4$}
  \put(101,18){$\mb g_2$}
  \put(90,-3){$\mb g_5$}
  \put(60,20){\circle*{1.5}}
  \put(80,35){\circle*{1.5}}
  \put(100,20){\circle*{1.5}}
  \put(90,0){\circle*{1.5}}
  \put(70,0){\circle*{1.5}}
  \put(70,0){\line(1,0){20}}
  \put(70,0){\line(-1,2){10}}
  \put(60,20){\line(4,3){20}}
  \put(80,35){\line(4,-3){20}}
  \put(100,20){\line(-1,-2){10}}
  \put(70.5,0){\line(1,4){8.7}}
  \put(70,0){\line(3,2){30}}
  \put(60,20){\line(1,0){40}}
  \put(89.5,0){\line(-1,4){8.7}}
  \put(90,0){\line(-3,2){30}}
  \put(80,13){\circle*{1.5}}
  \put(79,14){$\bf 0$}
  \end{picture}
  \end{center}
  \caption{A pentagon and its Gale diagram.}
  \label{cgd5g}
\end{figure}
\end{example}

Gale diagrams provide an efficient tool for studying the
combinatorics of higher-dimensional polytopes with few vertices,
%(or, dually, polytopes with few facets)
as in this case a Gale diagram translates the higher-dimensional
structure to a low-dimensional one. For example, Gale diagrams are
used to classify $n$-polytopes with up to $n+3$ vertices and to
find unusual examples when the number of vertices is~$n+4$,
see~\cite[\S6.5]{zieg95}.

\subsection*{Exercises}
\begin{exercise}
Describe combinatorial Gale diagrams of polytopes shown in
Fig.~\ref{figfv}.
\end{exercise}

\section{Face vectors and Dehn--Sommerville relations}\label{poly2}
The notion of the $f$-vector (or face vector) is a central concept
in the combinatorial theory of polytopes. It has been studied
since the time of Euler.

\begin{definition}\label{fhvect}
Let $P$ be a convex $n$-polytope. Denote by $f_i$ the number of
$i$-dimensional faces of~$P$. The integer sequence $\mb
f(P)=(f_0,f_1,\ldots,f_n)$ is known as the \emph{$f$-vector} (or
the \emph{face vector}) of~$P$. Note that $f_n=1$. The homogeneous
\emph{$F$-polynomial} of $P$ is defined by
\[
  F(P)(s,t)=s^n+f_{n-1}s^{n-1}t+\cdots+f_1st^{n-1}+f_0t^n.
\]
The \emph{$h$-vector} $\mb h(P)=(h_0,h_1,\ldots,h_n)$ and the
\emph{$H$-polynomial} of $P$ are defined by
\begin{equation}
\label{hvector}
  \begin{aligned}
  &h_0s^n+h_1s^{n-1}t+\cdots+h_nt^n=
  (s-t)^n+f_{n-1}(s-t)^{n-1}t+\cdots+f_0t^n,\\
  &H(P)(s,t)=
  h_0s^n+h_1s^{n-1}t+\cdots+h_{n-1}st^{n-1}+h_nt^n=F(P)(s-t,t).
  \end{aligned}
\end{equation}
The \emph{$g$-vector} of a simple polytope $P$ is the vector $\mb
g(P)=(g_0,g_1,\ldots,g_{\sbr n2})$, where $g_0=1$ and
$g_i=h_i-h_{i-1}$ for $i=1,\ldots,[n/2]$.
\end{definition}

\begin{example}\label{Hsim}
We have
\begin{align*}
  F(\varDelta^n)&=s^n+\binom{n+1}1s^{n-1}t+\binom{n+1}2s^{n-2}t^2+
  \cdots+t^n=\frac{(s+t)^{n+1}-t^{n+1}}s,\\
  H(\varDelta^n)&=s^n+s^{n-1}t+s^{n-2}t^2+
  \cdots+t^n=\frac{s^{n+1}-t^{n+1}}{s-t}.
\end{align*}
\end{example}

Obviously, the $f$-vector is a \emph{combinatorial invariant} of a
polytope, that is, it depends only on the face poset. This
invariant is far from being complete, even for simple polytopes:

\begin{example}
Two different combinatorial simple polytopes may have same
$f$-vectors. For instance, let $P$ be the 3-cube and $Q$ a simple
3-polytope with 2 triangular, 2 quadrangular and 2 pentagonal
facets, see Figure~\ref{figfv}. (Note that $Q$ is obtained by
truncating a tetrahedron at two vertices, it is also dual to the
cyclic polytope $C^3(6)$ from Definition~\ref{cyclic}.) Then $\mb
f(P)=\mb f(Q)=(8,12,6)$.
\begin{figure}[h]
\begin{center}
\begin{picture}(120,40)
\put(20,0){\line(-1,1){20}} \put(20,0){\line(3,2){30}}
\put(20,0){\line(0,1){15}} \put(0,20){\line(0,1){10}}
\put(0,30){\line(3,1){30}} \put(50,30){\line(-2,1){20}}
\put(20,15){\line(-4,3){20}} \put(20,15){\line(2,1){30}}
\put(50,20){\line(0,1){10}}
\multiput(0,20)(5.3,1.9){6}{\line(3,1){3.6}}
\multiput(50,20)(-5.6,2.9){4}{\line(-2,1){3.6}}
\multiput(30,31)(0,3.5){3}{\line(0,1){2}}
\put(90,0){\line(-1,1){20}} \put(90,0){\line(3,2){30}}
\put(90,0){\line(0,1){15}} \put(70,20){\line(0,1){10}}
\put(70,30){\line(3,1){15}} \put(85,35){\line(1,0){20}}
\put(105,35){\line(3,-1){15}} \put(90,15){\line(-4,3){20}}
\put(90,15){\line(2,1){30}} \put(120,20){\line(0,1){10}}
\multiput(70,20)(5.7,5.7){3}{\line(1,1){3.6}}
\multiput(120,20)(-5.7,5.7){3}{\line(-1,1){3.6}}
\end{picture}
\caption{Two combinatorially non-equivalent simple polytopes with
the same $f$-vectors.}
\label{figfv}
\end{center}
\end{figure}
\end{example}

The $f$-vector and the $h$-vector contain equivalent combinatorial
information, and determine each other by means of linear
relations, namely
\begin{equation}
\label{hf}
  h_k=\sum_{i=0}^k(-1)^{k-i}\binom{n-i}{n-k}f_{n-i},\quad
  f_k=\sum_{q=k}^n\binom qk h_{n-q},\quad\text{for } 0\le k\le n.
\end{equation}
In particular, $h_0=1$ and $h_n=f_0-f_1+\cdots+(-1)^nf_n$. By the
\emph{Euler formula},
\begin{equation}\label{euler}
  f_0-f_1+\cdots+(-1)^nf_n=1,
\end{equation}
which is equivalent to $h_n=h_0$. This is the first evidence of
the fact that many combinatorial relations for the face numbers
have much simpler form when written in terms of the $h$-vector.
Another example of this phenomenon is given by the following
generalisation of the Euler formula for simple or simplicial
polytopes.

\begin{theorem}[Dehn--Sommerville relations]
\label{ds} The $h$-vector of any simple $n$-polytope is symmetric,
that is,
$$
  H(s,t)=H(t,s),\quad\text{or}\quad h_i=h_{n-i}\quad\text{for }\;0\le i\le n.
$$
\end{theorem}

The Dehn--Sommerville relations can be proved in many different
ways. We present a proof from~\cite{bron83}, which can be viewed
as a combinatorial version of Morse-theoretic arguments. An
alternative proof will be given in Section~\ref{sec:dgrcp}.
%, see Theorem~\ref{FDS}.
%We shall return to this argument in chapter~\ref{tqtm}.

\begin{proof}[Proof of Theorem~{\rm \ref{ds}}.] Let $P\subset\R^n$ be a
simple polytope. Choose a generic linear function
$\varphi\colon\R^n\to\R$ which distinguishes the vertices of~$P$.
Write $\varphi(\mb x)=\langle\bnu,\mb x\rangle$ for some vector
$\bnu$ in~$\R^n$. The assumption on $\varphi$ implies that $\bnu$
is parallel to no edge of~$P$. We can view $\varphi$ as a height
function on~$P$ and turn the 1-skeleton of $P$ into a directed
graph by orienting each edge in such a way that $\varphi$
increases along it, see Figure~\ref{morse}.
\begin{figure}[h]
\begin{center}
\begin{picture}(65,50)
\put(25,0){\vector(-4,1){20}} \put(25,0){\vector(0,1){25}}
\put(25,0){\vector(4,1){20}} \put(5,5){\vector(-1,4){5}}
\put(0,25){\vector(1,4){5}} \put(5,45){\vector(4,1){20}}
\put(45,5){\vector(1,4){5}} \put(50,25){\vector(-1,4){5}}
\put(45,45){\vector(-4,1){20}} \put(25,25){\vector(1,1){5}}
\put(25,25){\vector(-1,1){20}} \put(50,25){\vector(-4,1){20}}
\put(30,30){\vector(-1,4){5}} \put(65,15){\vector(0,1){20}}
\put(62,17){\small$\bnu$} \put(26,3){\scriptsize$\ind=0$}
\put(15,23){\scriptsize$\ind=1$} \put(31,31){\scriptsize$\ind=2$}
\put(16,46){\scriptsize$\ind=3$}
\end{picture}
\caption{Orienting the 1-skeleton of $P$.} \label{morse}
\end{center}
\end{figure}
For each vertex $v$ of $P$ define the index $\ind_\nu(v)$ as the
number of incident edges that point towards~$v$. Denote the number
of vertices of index~$i$ by $I_\nu(i)$.  We claim that
$I_\nu(i)=h_{n-i}$. Indeed, each face of $P$ has a unique top
vertex (the maximum of the height function $\varphi$ restricted to
the face) and a unique bottom vertex (the minimum of $\varphi$).
Let $G$ be a $k$-face of~$P$, and $\va{G}$ its top vertex. Since
$P$ is simple, there are exactly $k$ edges of $G$ meeting at
$\va{G}$, whence $\ind(\va{G})\ge k$. On the other hand, each
vertex of index $q\ge k$ is the top vertex for exactly $\binom qk$
faces of dimension~$k$. It follows that the number of $k$-faces of
$P$ can be calculated as
$$
  f_k=\sum_{q\ge k}\binom qk I_\nu(q).
$$
Now the second identity from (\ref{hf}) shows that
$I_\nu(q)=h_{n-q}$, as claimed. In particular, the number
$I_\nu(q)$ does not depend on~$\bnu$. At the same time, we have
$\ind_\nu(v)=n-\ind_{-\nu}(v)$ for any vertex $v$, which implies
that
\[
  h_{n-q}=I_\nu(q)=I_{-\nu}(n-q)=h_q.
\qedhere\]
\end{proof}

\begin{remark}
The above proof also shows that the numbers $h_k=I_\nu(n-k)$ are
nonnegative, which is not evident from~\eqref{hf}. On the other
hand, the nonnegativity of the $h$-vector translates into certain
conditions on the $f$-vector. This will be important in the
subsequent study of $f$-vectors for combinatorial objects more
general than simple polytopes.
\end{remark}

\begin{theorem}
The $f$-vector of a simple $n$-polytope satisfies
\begin{equation}
\label{DSf}
  f_k=\sum_{i=0}^k(-1)^i\binom{n-i}{n-k}f_i,\quad\text{for }0\le k\le n.
\end{equation}
\end{theorem}
\begin{proof}
By the Dehn--Sommerville relations,
\[
  F(s-t,t)=H(s,t)=H(t,s)=F(t-s,s).
\]
By substituting $u=s-t$ we obtain $F(u,t)=F(-u,t+u)$, or
\begin{multline*}
  u^n+f_{n-1}u^{n-1}t+\cdots+f_1ut^{n-1}+f_0u_n\\=
  (-u)^n+f_{n-1}(-u)^{n-1}(t+u)+\cdots+
  f_1(-u)(t+u)^{n-1}+f_0(t+u)^n.
\end{multline*}
Calculating the coefficient of $u^kt^{n-k}$ in both sides above
yields~\eqref{DSf}.
\end{proof}

By Theorem~\ref{dualpol}, the $f$-vector of the dual $n$-polytope
$P^*$ satisfies
\[
  f_i(P^*)=f_{n-1-i}(P),\quad\text{for }0\le i\le n-1.
\]
Then it follows from~\eqref{DSf} that the $f$-vector of a
simplicial polytope satisfies the relations
$f_{n-1-k}=\sum_{i=0}^k(-1)^i\bin{n-i}{n-k}f_{n-1-i}$ or,
equivalently,
\[
  f_{q-1}=\sum_{j=q}^n(-1)^{n-j}\binom jq f_{j-1},
  \quad\text{for }1\le q\le n+1.
\]

%The Dehn--Sommerville relations were established by
%Dehn~\cite{dehn05} for $n\le5$ in 1905, and by
%Sommerville~\cite{somm27} in the general case in 1927 in a form
%similar to~\eqref{DSf}.

\begin{proposition}\label{multFH}
The $F$- and $H$-polynomial are multiplicative, i.e.
\begin{equation}\label{hvprod}
  F(P_1\times P_2)=F(P_1)F(P_2),  \quad  H(P_1\times P_2)=H(P_1)H(P_2).
\end{equation}
for any convex polytopes $P_1$ and $P_2$.
\end{proposition}
\begin{proof}
Let $\dim P_1=n_1$ and $\dim P_2=n_2$. Each $k$-face of $P_1\times
P_2$ is the product of an $i$-face of $P_1$ and a $(k-i)$-face of
$P_2$ for some $i$, whence
\begin{equation}\label{fprod}
  f_k(P_1\times P_2)=
  \sum_{i=0}^{n_1}f_i(P_1)f_{k-i}(P_2),\quad\text{for }0\le k\le n_1+n_2.
\end{equation}
This implies the first identity, and the second follows
from~\eqref{hvector}.
\end{proof}

\begin{example}\label{Hcube}
We have $F(I^n)=(F(\varDelta^1))^n$ and
$H(I^n)=(H(\varDelta^1))^n$, i.e.
\[
  F(I^n)=(s+2t)^n,\quad H(I^n)=(s+t)^n.
\]
\end{example}

We can also express the $f$-vector and the $h$-vector of the
connected sum $P\cs Q$ in terms of those of $P$ and~$Q$ (the proof
is left as an exercise):

\begin{proposition}\label{fhvcs}
Let $P$ and $Q$ be simple $n$-polytopes. Then
\begin{align*}
  f_0(P\cs Q)&=f_0(P)+f_0(Q)-2;\quad h_0(P\cs Q)=h_n(P\cs Q)=1;\\
  f_i(P\cs Q)&=f_i(P)+f_i(Q)-\binom ni,\quad\text{for }1\le i\le
  n;\\
  h_i(P\cs Q)&=h_i(P)+h_i(Q),\quad\text{for } 1\le i\le n-1.
\end{align*}
%and
%\begin{align*}
%  h_0(P\cs Q)&=h_n(P\cs Q)=1;\\
%  h_i(P\cs Q)&=h_i(P)+h_i(Q)\quad\text{for } 1\le i\le n-1.
%\end{align*}
\end{proposition}

Using the Dehn--Sommerville relations we can show that a
simplicial polytope cannot be `too neighbourly' (see
Definition~\ref{neighb}) if it is not a simplex.

\begin{proposition}\label{hneib}
Let $S$ be a $q$-neighbourly simplicial $n$-polytope, and let
$S\ne\varDelta^n$. Then $q\le\sbr n2$.
\end{proposition}
\begin{proof}
Let $S^*$ be the dual polytope. By Theorem~\ref{dualpol},
$f_{n-i}(S^*)=f_{i-1}(S)=\binom mi$ for $1\le i\le q$, where $m$
is the number of vertices of~$S$. From~(\ref{hf}) we get
\begin{equation}\label{fhneib}
  h_k(S^*)=\sum_{i=0}^k(-1)^{k-i}\binom{n-i}{k-i}\binom mi=
  \binom{m-n+k-1}k,\quad\text{for }k\le q,
\end{equation}
The second equality is obtained by calculating the coefficient of
$t^k$ on both sides of
% the identity
$$
  \frac1{(1+t)^{n-k+1}}(1+t)^m=(1+t)^{m-n+k-1}.
$$
If $S\ne\varDelta^n$, then $m>n+1$, which together
with~(\ref{fhneib}) gives $h_0(S^*)<h_1(S^*)<\cdots<h_q(S^*)$. It
then follows from the Dehn--Sommerville relations that $q\le\sbr
n2$.
\end{proof}

Since the $H$-polynomial of a simple $n$-polytope $P$ satisfies
the identity $H(P)(s,t)=H(P)(t,s)$, we can express it in terms of
elementary symmetric functions as follows:
\begin{equation}\label{Hgamma}
  H(P)(s,t)=\sum_{i=0}^{[n/2]}\gamma_i(s+t)^{n-2i}(st)^i.
\end{equation}
The identity $h_0=h_n=1$ implies that $\gamma_0=1$.

\begin{definition}\label{defgamma}
The integer sequence $\gamma(P)=(\gamma_0,\ldots,\gamma_{[n/2]})$
is called the \emph{$\gamma$-vector} of~$P$. We refer to
\[
  \gamma(P)(\tau)=\gamma_0+\gamma_1\tau+\cdots+
  \gamma_{[n/2]}\tau^{[n/2]}
\]
as the \emph{$\gamma$-polynomial} of~$P$.
\end{definition}

\begin{example}
We have $\gamma(I^1)(\tau)=1$. If $P^2_m$ is an $m$-gon, then
$\gamma(P^2_m)(\tau)=1+(m-4)\tau$.
\end{example}

The components of the $\gamma$-vector can be expressed via the
components of any of the $f$-, $h$- or $g$-vector by means of
linear relations, and vice versa. The explicit transition formulae
between the $g$- and $\gamma$-vectors are given by the next lemma.

\begin{lemma}\label{gammag}
Let $P$ be a simple $n$-polytope. Then
\begin{align*}
g_i &= (n-2i+1)\sum_{j=0}^{i}\frac{1}{n-i-j+1}
\binom{n-2j}{i-j}\gamma_j;\\
\gamma_i &= (-1)^i\sum_{j=0}^{i} (-1)^j \binom{n-i-j}{i-j}g_j, \quad
\quad\quad 0\le i\le\sbr n2.
\end{align*}
\end{lemma}
\begin{proof}
Using the formulae of Examples~\ref{Hsim} and~\ref{Hcube} we
calculate
\begin{align*}
  H(P)& =\sum_{j=0}^{[n/2]}\gamma_j(st)^jH(I^{n-2j}),\\
  H(I^q)&= \sum_{j=0}^q\binom{q}{j}s^jt^{q-j} = \sum_{k=0}^{[q/2]}
  c_{k,q}(st)^k H(\varDelta^{q-2k}),
\end{align*}
where $c_{0,q}=1$ and $c_{k,q}=\binom{q}{k}-\binom{q}{k-1} =
\frac{q-2k+1}{q-k+1} \binom{q}{k} > 0$ for $k=1,\ldots,[\frac
q2]$. Hence
\[
  H(P)(s,t)=\sum_{i=0}^{[n/2]}\Bigl( \sum_{j=0}^{i} c_{i-j,n-2j}\gamma_j
  \Bigr)(st)^iH(\varDelta^{n-2i}).
\]
On the other hand,
\[
  H(P)(s,t) = \sum_{i=0}^{[n/2]} g_i (st)^i H(\varDelta^{n-2i}).
\]
Comparing the last two formulae we obtain $g_i=\sum_{j=0}^{i}
c_{i-j,n-2j}\gamma_j$, which is equivalent to the first required
formula. The second is left as an exercise.
%
%We note that the transition matrix from $\gamma(P)$ to $\mb g(P)$ is
%upper triangular with unit diagonal; the explicit form of its
%inverse matrix given in the second formula of the Lemma is left as
%an exercise.
\end{proof}

\begin{proposition}\label{faceineq}
Let $P,Q$ be two simple $n$-polytopes, and consider the following
four conditions:
\begin{itemize}
\item[(a)] $f_i(P)\ge f_i(Q)$ for $i=0,1,\ldots,n$;
\item[(b)] $h_i(P)\ge h_i(Q)$ for $i=0,1,\ldots,n$;
\item[(c)] $g_i(P)\ge g_i(Q)$ for $i=0,1,\ldots,\sbr n2$;
\item[(d)] $\gamma_i(P)\ge\gamma_i(Q)$ for $i=0,1,\ldots,\sbr n2$.
\end{itemize}
Then {\rm(d)$\Rightarrow$(c)$\Rightarrow$(b)$\Rightarrow$(a)}.
\end{proposition}
\begin{proof}
The first of equations~\eqref{hf} implies that the components of
$\mb f(P)$ are expressed via the components of $\mb h(P)$ with
positive coefficients, which proves the implication
(b)$\Rightarrow$(a). We have $h_k=\sum_{i=0}^kg_i$ for $0\le
k\le\sbr n2$, which gives the implication (c)$\Rightarrow$(b). The
implication (d)$\Rightarrow$(c) follows similarly from the first
formula of Lemma~\ref{gammag}.
\end{proof}

The components of the $f$-vector of any polytope are nonnegative.
The nonnegativity of the components of the $h$-vector of a simple
polytope follows from their geometric interpretation obtained in
the proof of Theorem~\ref{ds}; the $h$-vector of a non-simple
polytope may have negative components (e.g. for the octahedron).
The nonnegativity of the $g$-vector of a simple $n$-polytope (that
is, the inequalities $g_i(P)\ge g_i(\varDelta^n)$) is a much more
subtle property; it follows from the $g$-theorem discussed in the
next section. The $\gamma$-vector of a simple polytope may have
negative components (e.g. for $P=\varDelta^2$); its nonnegativity
for special classes of simple polytopes will be discussed in
Section~\ref{flagpolytopes}. The nonnegativity of the
$\gamma$-vector can be expressed by the inequalities
$\gamma_i(P)\ge\gamma_i(I^n)$, see Exercise~\ref{gammaprodu}.

\subsection*{Exercises.}
\begin{exercise}
Show that any 3-polytope has a face with $\le5$ vertices.
\end{exercise}

\begin{exercise}
Prove Proposition~\ref{fhvcs}.
\end{exercise}

%\begin{problem}
%Describe all integer-valued functions on the set of simple
%polytopes which are linear with respect to the connected sum
%operation. Examples of such functions are the components $h_i$ of
%the $h$-vector for $1\le i\le n-1$.
%\end{problem}

\begin{exercise}
Prove the second transition formula of Lemma~\ref{gammag}.
\end{exercise}

\begin{exercise}\label{gammaprodu}
Show that the $\gamma$-polynomial is multiplicative, that is,
\[
  \gamma(P\times Q)(\tau)=\gamma(P)(\tau)\cdot\gamma(Q)(\tau)
\]
%for any two simple polytopes $P$ and $Q$.
In particular, $\gamma(I^n)(\tau)=1$, i.e.
$\gamma(I^n)=(1,0,\ldots,0)$.
\end{exercise}

\section{Characterising the face vectors of polytopes}
The face numbers are the simplest combinatorial invariants of
polytopes, and they arise in many hard problems of combinatorial
geometry. One of the most natural and basic questions is to
describe all possible face numbers, or, more precisely, determine
which integer vectors arise as the $f$-vectors of polytopes. In
the general case this questions is probably intractable (see the
end of this section), but a particularly nice answer exists in the
case of simple (or, equivalently, simplicial) polytopes.
Obviously, the Dehn--Sommerville relations provide a necessary
condition. As far as only linear equations are concerned, there
are no further restrictions:

\begin{proposition}[Klee~\cite{klee64}]\label{DSgenli}
The Dehn--Sommerville relations are the most general linear
equations satisfied by the face numbers of all simple
%(or simplicial)
polytopes.
\end{proposition}
\begin{proof}
Once again, the use of the $h$-vector simplifies the proof
significantly. Given a simple polytope $P$, we set
\[
  h(P)(t)=H(P)(1,t)=h_0(P)+h_1(P)t+\cdots+h_n(P)t^n.
\]
It is enough to prove that the affine hull of the $h$-vectors
$(h_0,h_1,\ldots,h_n)$ of simple $n$-polytopes is an $\sbr
n2$-dimensional plane. This can be done by presenting $\sbr n2+1$
simple polytopes with affinely independent $h$-vectors. Take
$Q_k=\varDelta^k\times\varDelta^{n-k}$ for $0\le k\le\sbr n2$.
Since $h(\varDelta^k)(t)=1+t+\cdots+t^k$, formula~\eqref{hvprod}
gives
\[
  h(Q_k)(t)=\frac{1-t^{k+1}}{1-t}\cdot\frac{1-t^{n-k+1}}{1-t}.
\]
It follows that the lowest degree term in the polynomial
$h(Q_{k+1})(t)-h(Q_k)(t)$ is $t^{k+1}$, for $0\le k\le\sbr n2-1$.
Therefore, the vectors $\mb h(Q_k)$ are affinely independent for
these values of~$k$.
\end{proof}

\begin{example}
Let $P$ be a simple polytope. Since every vertex is contained in
exactly $n$ edges and each edge connects two vertices, we have a
linear relation
\begin{equation}\label{vered}
  2f_1=nf_0.
\end{equation}
By Proposition~\ref{DSgenli} this must be a consequence of the
Dehn--Sommerville relations (in fact, it is equation~\eqref{DSf}
for $k=1$.)

Equation~\eqref{vered} together with the Euler
identity~\eqref{euler} shows that the $f$-vector of a simple
3-polytope $P^3$ is completely determined by the number of facets,
namely,
\[
  \mb f(P^3)=(2f_2-4,3f_2-6,f_2,1).
\]
Similarly, the $f$-vector of a simplicial 3-polytope $S^3$ is
determined by the number of vertices, namely,
\[
  \mb f(S^3)=(f_0,3f_0-6,2f_0-4,1).
\]
\end{example}

\begin{remark}
Euler's formula~\eqref{euler} is the only linear relation
satisfied by the face vectors of general polytopes. This can be
proved similarly to Proposition~\ref{DSgenli}, by specifying
sufficiently many polytopes with affinely independent face
vectors.
\end{remark}

Apart from the linear equations, the $f$-vectors of polytopes
satisfy certain inequalities. Here are some of the simplest of
them.

\begin{example}\label{simplefaceineq}
There are the following obvious lower bounds for the number of
vertices and the number of facets of an $n$-polytope:
\[
  f_0\ge n+1,\quad f_{n-1}\ge n+1.
\]
Since every pair of vertices is joined by at most one edge, and
every pair of facets intersect at most one face of codimension~2,
we have the upper bounds
\[
  f_1\le\bin{f_0}2,\quad f_{n-2}\le\bin{f_{n-1}}2.
\]
If the polytope is simplicial, then there is also the following
lower bound for~$f_1$:
\[
%f_{n-2}\ge nf_{n-1}-\bin{n+1}2.
  f_1\ge nf_0-\bin{n+1}2.
\]
It is much more difficult to prove though, even for 4-polytopes.
For simplicial 3-polytopes the inequality above turns into
identity.
\end{example}

Historically, the most important inequality-type results preceding
the general characterisation of $f$-vectors were the \emph{Upper
Bound Theorem (UBT)} and the \emph{Lower Bound Theorem (LBT)}.
They give respectively an upper and a lower bound for the number
of faces in a simplicial polytope with the given number of
vertices.

\begin{theorem}[UBT for simplicial polytopes]\label{ubt}
Among all simplicial $n$-polytopes $S$ with $m$ vertices the
cyclic polytope $C^n(m)$ (Example~\ref{cyclic}) has the maximal
number of $i$-faces for $1\le i\le n$. That is, if $f_0(S)=m$,
then
\[
  f_i(S)\le f_i\bigl(C^n(m)\bigr) \quad\text{ for }i=1,\ldots,n.
\]
Equality is achieved for all $i$ if and only if $S$ is a
neighbourly polytope.
%(see Definition~\ref{neighb}).
\end{theorem}

The UBT was conjectured by Motzkin
% in 1957
and proved by McMullen~\cite{mcmu70} in~1970.

Since $C^n(m)$ is neighbourly, we have
\[
  f_i\bigl(C^n(m)\bigr)=\bin m{i+1}\quad\text{ for }0\le i\le\sbr n2-1.
\]
Due to the Dehn--Sommerville relations, this determines the full
$f$-vector of $C^n(m)$. The exact values are given by the
following lemma.

\begin{lemma}\label{cycfn}
The number of $i$-faces of the cyclic polytope $C^n(m)$ (or any
neighbourly $n$-polytope with $m$ vertices) is given by
\[
  f_i=
 %\bigl(C^n(m)\bigr)=
  \sum_{q=0}^{\sbr n2}\binom q{n-1-i}\binom{m-n+q-1}q+
  \sum_{p=0}^{\sbr{n-\!1}2}\binom{n-p}{n-1-i}\binom{m-n+p-1}p
 % \sum_{q=0}^{\sbr n2}\bin
 % q{n-1-i}\bin{m-n+q-1}q+
 % \sum_{p=0}^{\sbr{n-1}2}\bin{n-p}{n-1-i}\bin{m-n+p-1}p
 % \quad\text{ for } 0\le i\le n,
\]
for $0\le i\le n$, where we set $\bin pq=0$ for $p<q$ or $q<0$.
\end{lemma}
\begin{proof}
We set $C=C^n(m)$. Using the second identity from~(\ref{hf}), the
identity $\sbr n2+1=n-\sbr{n-1}2$, the Dehn--Sommerville relations
for~$C^*$, and~(\ref{fhneib}), we calculate
\begin{align*}
f_i(C)&=f_{n-1-i}(C^*)=\sum_{q=0}^n\bin
q{n-1-i}h_{n-q}(C^*)\\
&=\sum_{q=0}^{[n/2]}\bin q{n-1-i}h_q(C^*) +\sum_{q=[n/2]+1}^n\bin q{n-1-i}h_{n-q}(C^*)\\
&=\sum_{q=0}^{[n/2]}\bin q{n-1-i}\bin{m-n+q-1}q+
\sum_{p=0}^{[(n-1)/2]}\bin{n-p}{n-1-i}\bin{m-n+p-1}p.\qedhere
\end{align*}
\end{proof}

\begin{lemma}\label{ubth}
Assume that the inequalities
\[
  h_i(P)\le\bin{m-n+i-1}i,\qquad i=0,\ldots,\sbr n2
\]
hold for the $h$-vector of a simple polytope $P$ with $m$ facets.
Then the dual simplicial polytope $P^*$ satisfies the UBT
inequalities $f_i(P^*)\le f_i\bigl(C^n(m)\bigr)$ for
$i=1,\ldots,n$.
\end{lemma}
\begin{proof}
Since $h_i\bigl((C^n(m))^*\bigr)=\bin{m-n+i-1}i$ for
$i=0,\ldots,\sbr n2$, the statement follows from
Proposition~\ref{faceineq}. Alternatively, replace the last `$=$'
in the calculation from the proof of Lemma~\ref{cycfn} by `$\le$'.
\end{proof}

This was one of the key observations in McMullen's original
proof~\cite{mcmu70} of the UBT for simplicial polytopes (the whole
proof can also be found in~\cite[\S18]{bron83} and
\cite[\S8.4]{zieg95}). R.~Stanley gave an algebraic argument
establishing the inequalities from Lemma~\ref{ubth}, and therefore
the UBT, in a much more general setting of \emph{triangulated
spheres}\label{trisphe}. We shall discuss Stanley's approach and
conclude the proof of the UBT in Section~\ref{cmr}.

\begin{remark}
The UBT holds for all convex polytopes. That is, the cyclic
polytope $C^n(m)$ has the maximal number of $i$-faces among all
convex $n$-polytopes with $m$ vertices. The argument for this
builds on the following observation of Klee and McMullen, which we
reproduce from~\cite[Lemma~8.24]{zieg95}.

\begin{lemma}\label{fperturb}
By a small perturbation of vertices of an $n$-polytope $P$ one can
achieve that the resulting polytope $P'$ is simplicial, and
\[
  f_i(P)\le f_i(P')\quad\text{for }i=1,\ldots,n-1.
\]
\end{lemma}
\end{remark}

\begin{definition}
\label{stacked} A simplicial $n$-polytope $S$ is called
\emph{stacked} if there is a sequence $S_0,S_1,\dots,S_k=S$ of
$n$-polytopes such that $S_0$ is an $n$-simplex and $S_{i+1}$ is
obtained from $S_i$ by adding a pyramid over a facet of~$S_i$ (the
vertex of the added pyramid is chosen close enough to its base, so
that the whole construction remains convex and simplicial). The
polar simple polytopes are those obtained from a simplex by
iterating the vertex cut operation of Example~\ref{expcs}.2. These
are sometimes called \emph{truncation polytopes}.
\end{definition}

The $f$-vector of a stacked polytope is easy to calculate (see
Exercise~\ref{fvsta}).

\begin{theorem}[LBT for simplicial polytopes]
\label{lbt} Among all simplicial $n$-polytopes $S$ with $m$
vertices a stacked polytope has the minimal number of $i$-faces
for $2\le i\le n-1$. That is, if $f_0(S)=m$, then
\begin{align*}
  f_i(S)&\ge \bin ni m-\bin{n+1}{i+1}i \quad\text{ for }
  i=1,\ldots,n-2;\\
  f_{n-1}(S)&\ge(n-1)m-(n+1)(n-2).
\end{align*}
For $n\ne3$ equality is achieved for all $i$ if and only if $S$ is
a stacked polytope.
\end{theorem}

\begin{remark}
For $n=3$ the LBT inequalities $f_1\ge3m-6$ and $f_2\ge2m-4$ turn
into equalities for all simplicial polytopes.
\end{remark}

An inductive argument by McMullen, Perles and
Walkup~\cite{mc-wa71} reduced the LBT to the case $i=1$, namely,
to the inequality $f_1\ge nm-\bin{n+1}2$.
%(or equivalently, the inequality $f_{n-2}\ge nf_{n-1}-\bin{n+1}2$
%for simple polytopes).
It was finally proved by Barnette~\cite{barn71},~\cite{barn73}.
Barnette's proof of the LBT, with some simplifications, can also
be found in~\cite{bron83}. The fact that equality is achieved only
for stacked polytopes (if $n\ne3$) was proved by Billera and
Lee~\cite{bi-le81}.

\begin{remark}
Unlike the UBT, little is know about generalisations of the LBT
to non-simplicial convex polytopes. Some results in this direction
were obtained in~\cite{kala87} along with generalisations of the
LBT to triangulated spheres and manifolds, which we also discuss
later in this book.
\end{remark}

An easy calculation shows that the inequalities $f_{n-1}\ge n+1$
and and $f_{n-2}\ge nf_{n-1}-\bin{n+1}2$ for simple polytopes from
Example~\ref{simplefaceineq} can be written in terms of the
$h$-vector as follows: $h_0\le h_1\le h_2$. Having reduced the
whole LBT to the inequality $h_1\le h_2$, McMullen and
Walkup~\cite{mc-wa71} conjectured that the components of the
$h$-vector `grow up to the middle', that is the inequalities
\begin{equation}
\label{glbt}
  h_0\le h_1\le\cdots\le h_{\sbr n2}
\end{equation}
hold for a simple $n$-polytope. It has since become known as the
\emph{Generalised Lower Bound Conjecture (GLBC)}.

McMullen also suggested a generalisation to the UBT, whose
formulation requires an algebraic digression.

\begin{definition}
For any two positive integers $a$, $i$ there exists a unique
\emph{binomial $i$-expansion} of $a$ of the form
\[
  a=\bin{a_i}i+\bin{a_{i-1}}{i-1}+\cdots+\bin{a_j}j,
\]
where $a_i>a_{i-1}>\cdots>a_j\ge j\ge1$.

The binomial $i$-expansion of $a$ can be constructed by choosing
$a_i$ as the unique number satisfying $\bin{a_i}i\le
a<\bin{a_i+1}i$, then choosing $a_{i-1}$, and so on. One needs
only to check that $a_i>a_{i-1}$, which is straightforward.

%by induction on~$i$. For $i=1$ it has the form $a=\bin a1$. The
%inductive step is achieved by observing that
%\[
%  a=\bin{a_i}i+\Bigl(\text{binomial $(i-1)$-expansion of
%  $a-\bin{a_i}i$}\Bigr),
%\]
%where $a_i\ge i$ is the unique number satisfying $\bin{a_i}i\le
%a<\bin{a_i+1}i$. It remains to check that $a_i>a_{i-1}$, which is
%straightforward.

Now define the $i$th \emph{pseudopower} of~$a$ as
\[
  a^{\langle i\rangle}=
  \bin{a_i+1}{i+1}+\bin{a_{i-1}+1}i+\cdots+\bin{a_j+1}{j+1},\quad
  0^{\langle i\rangle}=0.
\]
\end{definition}

\begin{example}\label{binexpex}\

1. For $a>0$, $a^{\langle\! 1\rangle}=\binom{a+1}2$.

2. If $i\ge a$ then the binomial expansion has the form
$$
  a=\bin ii+\bin{i-1}{i-1}+\cdots+\bin{i-a+1}{i-a+1}=1+\cdots+1,
$$
and therefore $a^{\langle i\rangle}=a$.

3. Let $a=28$, $i=4$. Then
\[
  28=\bin 64+\bin 53+\bin 32,\quad
  28^{\langle 4\rangle}=\bin 75+\bin 64+\bin 43=40.
\]
\end{example}

The importance of the binomial expansion and pseudopowers comes
from the following fundamental result of combinatorial commutative
algebra.

\begin{theorem}[Macaulay, Stanley]
\label{mvect} The following two conditions are equivalent for a
sequence of integers $(k_0,k_1,k_2,\ldots)$:
\begin{enumerate}
\item[(a)] $k_0=1$ and $0\le k_{i+1}\le k_i^{\langle i\rangle}$
for $i\ge1$;

\item[(b)] there exists a connected commutative graded algebra
$A=A^0\oplus A^1\oplus A^2\oplus\cdots$ over a field $\k$ such
that $A$ is generated by its degree-one elements and $\dim_\k
A^{i}=k_i$ for $i\ge0$.
\end{enumerate}
\end{theorem}

Macaulay's original theorem~\cite{maca27} says that (b) above is
equivalent to the existence of a \emph{multicomplex} whose
$h$-vector is given by~$(k_0,k_1,k_2,\ldots)$. The original proof
is long and complicated. The reformulation of Macaulay's condition
in terms of pseudopowers, i.e. condition~(a), is due to
Stanley~\cite[Theorem~2.2]{stan96}. Simpler proofs of
Theorem~\ref{mvect} can be found in~\cite{cl-li69}
and~\cite[\S4.2]{br-he98}.

\begin{definition}\label{defmvect}
A sequence of integers satisfying either of the conditions of
Theorem~\ref{mvect} is called an \emph{M-sequence}. Finite
$M$-sequences are \emph{M-vectors}.
\end{definition}

Now observe that the number $\bin{m-n+i-1}i$ on the right hand
side of the inequality of Lemma~\ref{ubth} equals the number of
degree $i$ monomials in $h_1=m-n$ generators. It follows that the
inequalities of Lemma~\ref{ubth}, and therefore the UBT, hold if
the $h$-vector is an $M$-vector.

In 1970 McMullen~\cite{mcmu71} combined all known and conjectured
information about the $f$-vectors, including the Dehn--Sommerville
relations and the generalisations to the LBT and UBT discussed
above, into a (conjectured) complete characterisation. McMullen's
conjecture is now proved, and remains up to the present time
perhaps the most impressive achievement of the combinatorial
theory of face numbers:

\begin{theorem}[$g$-theorem]\label{gth}
An integer vector $(f_0,f_1,\ldots,f_n)$ is the $f$-vec\-tor of a
simple $n$-polytope if and only if the corresponding sequence
$(h_0,\ldots,h_n)$ determined by~\eqref{hvector} satisfies the
following three conditions:
\begin{enumerate}
\item[(a)] $h_i=h_{n-i}$ for $i=0,1,\ldots,n$ (the Dehn--Sommerville
relations);
\item[(b)] $h_0\le h_1\le\cdots\le h_{\sbr n2}$;
\item[(c)] $h_0=1$, $h_{i+1}-h_i\le(h_i-h_{i-1})^{\langle i\rangle}$
for $i=1,\ldots,\sbr n2-1$.
\end{enumerate}
\end{theorem}

\begin{remark}
Condition (b) says that the components of the $h$-vector `grow up to
the middle', while (c) gives a restriction on the rate of this
growth. Both (b) and (c) can be reformulated by saying that the
$g$-vector $(g_0,\ldots,g_{[n/2]})$ (see Definition~\ref{fhvect}) is
an $M$-vector; this explains the name `$g$-theorem'. The fact that
the $g$-vector is an $M$-vector implies that the $h$-vector is also
an $M$-vector (see Exercise~\ref{ghmvec}).
\end{remark}

%Both $g$-theorem and Theorem~\ref{mvect} are difficult, and we
%shall not provide their proofs in order not to distract the
%reader. We give some more historical remarks and references
%instead.

Both necessity and sufficiency parts of the $g$-theorem were
proved almost simultaneously (around 1980), although by radically
different methods.

The sufficiency part was proved by Billera and
Lee~\cite{bi-le80},~\cite{bi-le81}. The proof is quite elementary
and relies upon a remarkable combinatorial-geometrical
construction combining cyclic polytopes (achieving the upper bound
for the number of faces) with the operation of `adding a pyramid'
(used to produce polytopes achieving the lower bound). As a
result, a simplicial polytope can be produced with any prescribed
$g$-vector between the minimal and the maximal ones. Another
important consequence of the results of~\cite{mc-wa71}
and~\cite{bi-le81} is that the GLBC inequalities~\eqref{glbt} are
the most general linear inequalities satisfied by the $f$-vectors
of simplicial polytopes.

On the other hand, Stanley's proof~\cite{stan80} of the necessity
part of $g$-theorem (i.e. that the $g$-vector of a simple polytope
is an $M$-vector) used deep results from algebraic geometry, in
particular, the \emph{Hard Lefschetz theorem}\label{hardlef} for
the cohomology of \emph{toric varieties}. We shall give this
argument in Section~\ref{cohtm}. After the appearance of Stanley's
paper combinatorists had been looking for a more elementary
combinatorial proof of his theorem, until in 1993 a first such
proof was found by McMullen~\cite{mcmu93}. It builds upon the
notion of the \emph{polytope algebra}\label{polytalg}, which may
be thought of as a combinatorial model for the cohomology algebras
of toric varieties. Despite being elementary, it was a complicated
proof. Later McMullen simplified his approach in~\cite{mcmu96}.
Yet another elementary proof of the $g$-theorem was given by
Timorin~\cite{timo99}. It relies on an interpretation of
McMullen's polytope algebra as the algebra of differential
operators with constant coefficients vanishing on the \emph{volume
polynomial} of the polytope.

By duality, the UBT and the LBT provide upper and lower bounds for
the number of faces of a simple polytope with the given number of
facets. Similarly, the $g$-theorem also provides a
characterisation for the $f$-vectors of simplicial polytopes.
During the last three decades some work was done in extending the
$g$-theorem to objects more general than simplicial (or simple)
polytopes, although the most important conjecture here remains
open since 1971 (see Section~\ref{simsph}). There are basically
two diverging routes for generalisations of the $g$-theorem:
towards non-polytopal objects (like triangulations of spheres or
manifolds), and towards general convex polytopes which are neither
simple nor simplicial. The former requires machinery from
combinatorial topology and commutative algebra; we shall discuss
the corresponding generalisations in more detail in the next
chapters. The generalisations of the $g$-theorem to non-simplicial
convex polytopes are mostly beyond the scope of this book; they
require algebraic geometry techniques such as \emph{intersection
homology}\label{intehomo}, which we only briefly discuss in
Section~\ref{cohtm}.

\subsection*{Exercises.}
\begin{exercise}\label{newub}
Prove the following upper bound for the number of $k$-faces in a
simple $n$-polytope $P$ with $f_0$ vertices:
\[
  f_{k}\le\frac{1}{k+1}\binom{n}{k}f_{0},\quad
  \text{for }k=1,\ldots,n,
\]
where the equality is achieved only for $P=\varDelta^n$. Observe
that this inequality gives a better upper bound for simple
polytopes that the UBT.
\end{exercise}

\begin{exercise}
Prove Lemma~\ref{fperturb}.
\end{exercise}

\begin{exercise}\label{fvsta}
Show that the numbers of faces of a stacked $n$-polytope $S$ with
$m=f_0$ vertices are given by
\begin{align*}
  f_i(S)&= \bin ni m-\bin{n+1}{i+1}i \quad\text{ for }
  1\le i\le n-2;\\
  f_{n-1}(S)&=(n-1)m-(n+1)(n-2).
\end{align*}
\end{exercise}

\begin{exercise}\label{ghmvec}
Let $\mb h=(h_0,h_1,\ldots,h_n)$ be an integer vector with $h_0=1$
and $h_i=h_{n-i}$ for $0\le i\le n$, and let $\mb
g=(g_0,g_1,\ldots,g_{[n/2]})$ where $g_0=1$, $g_i=h_i-h_{i-1}$ for
$i>0$. Show that if $\mb g$ is an $M$-vector, then $\mb h$ is also
an $M$-vector.
\end{exercise}

\begin{exercise}
Prove directly that parts (a) and (b) of the $g$-theorem imply the
LBT, while parts (a) and (c) imply the UBT.
\end{exercise}

\chapter*{\ \ Polytopes: additional topics}

\section{Nestohedra and graph-associahedra}\label{nestohedra}
%For many applications of combinatorial geometry, effective ways
%are needed to construct families of polytopes with special
%properties.
Several constructions of series of simple polytopes with
remarkable properties appeared in the beginning of the 1990s under
the common name of `generalised associahedra'. (The original
associahedron, or \emph{Stasheff polytope}\label{stasheptpe}, was
first introduced in homotopy theory~\cite{stas63}.) Nowadays
generalised associahedra find numerous applications in algebraic
geometry~\cite{g-k-z94}, the theory of knot and link
invariants~\cite{bo-ta94}, representation theory and cluster
algebras~\cite{fo-ze02}, and the theory of operads and geometric
`field theories' originating from quantum physics~\cite{stas97}.

Without attempting to overview all aspects of generalised
associahedra, we describe one particular construction which uses
the Minkowski sum and the combinatorial concept of a
\emph{building set}. Our exposition is much influenced by the
original works of Feichtner--Sturmfels~\cite{fe-st05} and
Postnikov~\cite{post09}. The resulting family of polytopes is
known as \emph{nestohedra}\label{nestohed}. Although it does not
include all possible generalisations of associahedra, this family
is wide enough to contain all classical series, and its
construction is elementary enough so that it requires no specific
knowledge.

The \emph{Minkowski sum} is a classical geometric construction
allowing one to produce new polytopes from known ones, just like
the product, hyperplane cut and connected sum described in
Section~\ref{poly}. However, a Minkowski sum of simple polytopes
usually fails to be simple. Interesting examples of polytopes can
be obtained by taking Minkowski sums of regular simplices.
Simplices in such a Minkowski sum are indexed by a collection
$\mathcal S$ of subsets in a finite set. It was shown
in~\cite{fe-st05} and~\cite{post09} that a Minkowski sum of
regular simplices is a simple polytope if $\mathcal S$ satisfies
certain combinatorial condition, identifying it as a
\emph{building set}\label{buildise}. The resulting family of
simple polytopes was called \emph{nestohedra} in~\cite{p-r-w08}
because of its connection to \emph{nested sets} considered by
De~Concini and Procesi~\cite{de-pr95} in the context of subspace
arrangements. An example of a building set is provided by the
collection of subsets of vertices which span connected subgraphs
in a given graph. The corresponding nestohedra are called
\emph{graph-associahedra}; they were introduced and studied in the
work of Carr and Devadoss~\cite{ca-de04}, and independently in the
work of Toledano Laredo~\cite{tole08} under the name~\emph{De
Concini--Procesi associahedra}. These include the classical series
of \emph{permutahedra} and \emph{associahedra}.
%The latter series, also known as \emph{Stasheff polytopes}, takes
%its origin in the homotopy theory~\cite{stas63}.

\subsection*{Minkowski sums of simplices}
Recall that the \emph{Minkowski sum}\label{minkosu} of two subsets
$A,B\subset\R^n$ is defined as
\[
  A+B=\{\mb x+\mb y\colon\mb x\in A,\;\mb y\in B\}.
\]

\begin{proposition}
The Minkowski sum of two polytopes is a polytope. Moreover, if
$P=\mathop{\mathrm{conv}}(\mb v_1,\ldots,\mb v_k)$ and
$Q=\mathop{\mathrm{conv}}(\mb w_1,\ldots,\mb w_l)$, then
\[
  P+Q=\mathop{\mathrm{conv}}(\mb v_1+\mb w_1,\ldots,\mb v_i+\mb w_j,\ldots,\mb v_k+\mb
  w_l).
\]
\end{proposition}
\begin{proof}
Follows directly from the definition of Minkowski sum.
\end{proof}

For every subset $S\subset[n+1]=\{1,\ldots,n+1\}$ consider the
regular simplex
\[
  \varDelta_S=\mathop{\mathrm{conv}}(\mb e_i\colon i\in
  S)\subset\R^{n+1}.
\]

Let $\mathcal F$ be a collection of nonempty subsets of~$[n+1]$.
We assume that $\mathcal F$ contains all singletons $\{i\}$, \
$1\le i\le n+1$. As usual, denote by $|\mathcal F|$ the number of
elements in $\mathcal F$. Given a subset $N\subset[n+1]$, denote
by $\mathcal F|_N$ the \emph{restriction}\label{restricti} of
$\mathcal F$ to~$N$, i.e.
% the collection
\[
  \mathcal F|_N=\{S\in\mathcal F\colon S\subset N\}
\]
of subsets in~$N$. A collection $\mathcal F$ is \emph{connected}
if $[n+1]$ cannot be represented as a disjoint union of nonempty
subsets $N_1$ and $N_2$ such that for every $S\in\mathcal F$
either $S\subset N_1$ or $S\subset N_2$. Obviously, every
collection $\mathcal F$ splits into the disjoint union of its
\emph{connected components}.

\begin{remark}
%Although it looks too technical,
In some further considerations we shall allow $\mathcal F$ to
contain some subsets of~$[n+1]$ with multiplicities.
\end{remark}

Now consider the convex polytope
\begin{equation}\label{mssim}
  P_{\mathcal F}=\sum_{S\in{\mathcal F}}\varDelta_S\subset\R^{n+1}.
\end{equation}
The following statement gives some of its basic properties.

\begin{proposition}\label{basicms}
Let $\mathcal F=\mathcal F_1\sqcup \cdots\sqcup\mathcal F_q$ be
the decomposition into connected components. Then
\begin{itemize}
\item[\rm(a)]$P_\mathcal F=P_{\mathcal F_1}\times\cdots\times
P_{\mathcal F_q}$;
\item[\rm(b)] $\dim P_\mathcal F=n+1-q$.
\end{itemize}
\end{proposition}
\begin{proof}
Statement~(a) follows from the fact that the polytopes
$P_{\mathcal F_i}$ are contained in complementary subspaces for
$1\le i\le q$, so their Minkowski sum is the product. Because
of~(a), it is enough to verify~(b) for connected $\mathcal F$
only. Then we need to prove that $\dim P_\mathcal F=n$. Observe
that $P_\mathcal F$ is contained in the hyperplane
$
  H_{\mathcal F}=
  \bigl\{\mb x\in\R^{n+1}\colon\sum_{i=1}^{n+1}x_i=
  |\mathcal F|\bigl\},
$
and therefore $\dim P_\mathcal F\le n$. If $\mathcal F$ has a
unique maximal element, then this element is $[n+1]$, because
$\mathcal F$ is connected. Hence, $P_\mathcal F$ has an
$n$-simplex as a Minkowski summand, and therefore $\dim P_\mathcal
F=n$. We shall only need this case in the further considerations,
so we skip the rest of the proof and leave it as an exercise.
% Assume now that
%$\mathcal F$ has $p$ maximal elements $S_1,\ldots,S_p$. By
%reordering the set $[n+1]$ if necessary, we may assume that
%$S_1\cup\ldots\cup S_{p-1}=[k]$ and $S_p=\{l,l+1,\ldots,n+1\}$.
%Since $\mathcal F$ is connected, it is easy to see that $k\ge l$
%and the restriction $\mathcal F|_{[k]}$ is also connected. By the
%induction assumption $\dim P_{\mathcal F|_{[k]}}=k-1$. Now,
%$P_\mathcal F$ has a Minkowski summand of the form $P_{\mathcal
%F|_{[k]}}+\varDelta_{S_p}$. Since $P_{\mathcal F|_{[k]}}$ is a
%$(k-1)$-dimensional polytope in the first $k$ coordinate subspace,
%$\varDelta_{S_p}$ is the regular simplex in the last $n-l+2$
%coordinates, and $k\ge l$, the polytope $P_{\mathcal
%F|_{[k]}}+\varDelta_{S_p}$ is $n$-dimensional. Therefore, $\dim
%P_\mathcal F=n$.
\end{proof}

We now describe two extreme examples of polytopes $P_\mathcal F$,
corresponding to the minimal and the maximal connected
collections.

\begin{example}[Simplex]\label{mssimplex}
Let $\mathcal S$ be the collection consisting of all singletons
and the whole set~$[n+1]$. Then $P_\mathcal S$ is the regular
$n$-simplex $\varDelta_{[n+1]}$ shifted by the vector $\mb
e_1+\cdots+\mb e_{n+1}$.
\end{example}

\subsection*{Permutahedron}\label{mspermut}
Let $\mathcal C$ be the \emph{complete} collection, consisting of
all subsets in~$[n+1]$. The polytope $P_\mathcal C$ is
$n$-dimensional by Proposition~\ref{basicms}.

\begin{theorem}\label{facetsperm}
$P_\mathcal C$ can be described as the intersection of the
hyperplane
\[
  H_{\mathcal C}=
  \bigl\{\mb x\in\R^{n+1}\colon\sum_{i=1}^{n+1}x_i=2^{n+1}-1
  \bigl\}
\]
with the halfspaces
\[
  H_{S,\ge}=\bigl\{\mb x\in\R^{n+1}\colon\sum_{i\in S}x_i\ge2^{|S|}-1
  \bigl\}
\]
for all proper subsets $S\subsetneq[n+1]$. Moreover, every
halfspace above is irredundant, i.e. determines a facet $F_S$ of
$P_\mathcal C$, so there are $|\mathcal C|=2^{n+1}-2$ facets in
total.
\end{theorem}
\begin{proof}
By definition, every point $\mb x=(x_1,\ldots,x_{n+1})\in
P_\mathcal C$ can be written as $\mb x=\sum_{S\in\mathcal C}\mb
x^S$ where $\mb x^S=(x^S_1,\ldots,x^S_{n+1})\in\varDelta_S$. Then
\[
  \sum_{i=1}^{n+1}x_i=\sum_{S\in\mathcal C}
  \sum_{i=1}^{n+1}x^S_i=\sum_{S\in\mathcal C}1=|\mathcal C|=2^{n+1}-1,
\]
which implies that $P_\mathcal C\subset H_\mathcal C$. Similarly,
\begin{equation}\label{minS}
  \sum_{i\in S}x_i=\sum_{T\in\mathcal C}\sum_{i\in S}x^T_i\ge
  \sum_{T\subset S}\sum_{i\in S}x_i^T=\sum_{T\subset S}1=
  \bigl|\mathcal C|_S\bigr|=2^{|S|}-1,
\end{equation}
so $P_\mathcal C$ is contained in all subspaces $H_{S,\ge}$.

It remains to show that any facet of $P_\mathcal C$ has the form
$P_\mathcal C\cap H_S$, where $H_S$ is the bounding hyperplane
for~$H_{S,\ge}$. Since $P_\mathcal C$ is a Minkowski sum of
simplices, each of its faces $G$ is a Minkowski sum of faces of
these simplices. We therefore may write $G$ as $\sum_{S\in\mathcal
C}\varDelta_{T_S}$ where $T_S\subset S$. By
Proposition~\ref{basicms}, if $G$ is a facet (i.e. $\dim G=n-1$),
then the collection $\mathcal T=\{T_S\}$ of subsets in~$[n+1]$ has
exactly two connected components. (Note that $\mathcal T$ may
contain some subsets more than once.) Let $[n+1]=N_1\sqcup N_2$
and $\mathcal T=\mathcal T|_{N_1}\sqcup\mathcal T|_{N_2}$ be the
decomposition into components. Then the hyperplane containing the
facet $G$ is defined by each of the two equations
\begin{equation}\label{n1n2}
  \sum_{i\in N_1}x_i=\bigl|\mathcal T|_{N_1}\bigr|
  \quad\text{ or }\quad
  \sum_{i\in N_2}x_i=\bigl|\mathcal T|_{N_2}\bigr|.
\end{equation}
Since every $T_S$ is contained in the corresponding $S$, we have
\begin{equation}\label{tcn1}
  \bigl|\mathcal T|_{N_1}\bigr|\ge\bigl|\mathcal C|_{N_1}\bigr|
  \quad\text{ and }\quad
  \bigl|\mathcal T|_{N_2}\bigr|\ge\bigl|\mathcal C|_{N_2}\bigr|.
\end{equation}
We claim that at least one of these inequalities turns into
equality. Indeed, assume the converse. By~\eqref{minS} the minimum
of the linear function $\sum_{i\in N_1}x_i$ on $P_\mathcal C$ is
$\bigl|\mathcal C|_{N_1}\bigr|$, so there is a point $\mb x'\in
P_\mathcal C$ with
\begin{equation}\label{sumn1x1}
  \sum_{i\in N_1}x'_i=\bigl|\mathcal C|_{N_1}\bigr|<\bigl|\mathcal
  T|_{N_1}\bigr|.
\end{equation}
Similarly, there is a point $\mb x''\in P_\mathcal C$ with
\[
  \sum_{i\in N_2}x''_i=\bigl|\mathcal C|_{N_2}\bigr|<\bigl|\mathcal
  T|_{N_2}\bigr|.
\]
Since $N_1\sqcup N_2=[n+1]$, the latter inequality is equivalent
to $\sum_{i\in N_1}x''_i>\bigl|\mathcal T|_{N_1}\bigr|$. This
together with~\eqref{sumn1x1} implies that there are points of
$P_\mathcal C$ in both open halfspaces determined by the first of
the equations~\eqref{n1n2}, which contradicts the assumption that
$G$ is a facet. So at least one of~\eqref{tcn1} is an equality,
which implies that the hyperplane~\eqref{n1n2} containing $G$ has
the form~$H_S$ (where $S$ is either $N_1$ or $N_2$).

It follows that every facet is contained in the hyperplane $H_S$
for some~$S$. On the other hand, every subset $S$ can be taken as
$N_1$ in the construction of the preceding paragraph, which shows
that every $H_S$ contains a facet.
\end{proof}

Having identified the facets, we may derive the following
description of the whole face poset of~$P_\mathcal C$.
%Remember
%that a length $k$ \emph{composition} of the set $n+1$ is a
%subdivision of the set $[n+1]$ into ordered sequence of $k$
%disjoint subsets.

\begin{proposition}\label{facesperm}
Faces of~$P_\mathcal C$ of dimension $k$ are in one-to-one
correspondence with ordered partitions of the set $[n+1]$ into
$n+1-k$ nonempty parts. An inclusion of faces $G\subset F$ occurs
whenever the ordered partition corresponding to $G$ can be
obtained by refining the ordered partition corresponding to~$F$.
\end{proposition}
\begin{proof}
This follows from the fact that two facets $F_{S_1}$ and $F_{S_2}$
of $P_\mathcal C$ have nonempty intersection if and only if
$S_1\subset S_2$ or $S_2\subset S_1$. We skip the details to avoid
repetitive arguments; see Theorem~\ref{intfacn} below for a more
general result.
\end{proof}

\begin{corollary}\
\begin{itemize}
\item[(a)] $P_\mathcal C$ is a simple polytope;
\item[(b)] the vertices of $P_\mathcal C$ are obtained by
all permutations of the coordinates of the point
$(1,2,4,\ldots,2^n)\in\mathbb R^{n+1}$.
\end{itemize}
\end{corollary}

The polytope whose vertices are obtained by permuting the
coordinates of a given point is known as the
\emph{permutahedron}\label{permutahe}; it has been studied by
convex geometers since the beginning of the 20th century. More
precisely, given a point $\mb a=(a_1,\ldots,a_{n+1})\in\mathbb
R^{n+1}$ with $a_1<a_2<\cdots<a_{n+1}$, define the corresponding
permutahedron as
\[
  \mathit{P\!e}^n(\mb a)=\mathop{\mathrm{conv}}\bigl(\sigma(a_1),\ldots,\sigma(a_{n+1})
  \colon \sigma\in\Sigma_{n+\!1}\bigr),
\]
where $\Sigma_{n+\!1}$ denotes the group of permutations of $n+1$
elements. In particular, $P_\mathcal
C=\mathit{P\!e}^n(1,2,\ldots,2^n)$. All $n$-dimensional
permutahedra $\mathit{P\!e}^n(\mb a)$ are combinatorially
equivalent; this follows from the following description of their
faces:

\begin{theorem}\label{facesgenperm}\
\begin{itemize}
\item[(a)]
Every facet of $\mathit{P\!e}^n(\mb a)$ is the intersection of
$\mathit{P\!e}^n(\mb a)$ with the hyperplane
\[
  H_S=\bigl\{\mb x\in\R^{n+1}\colon\sum_{i\in
  S}x_i=a_1+a_2+\cdots+a_{|S|}\bigr\}
\]
for a proper subset $S\subsetneq[n+1]$.
\item[(b)]
The faces of $\mathit{P\!e}^n(\mb a)$ are described in the same
way as in Proposition~\ref{facesperm}.
\end{itemize}
\end{theorem}
\begin{proof}
This can be proved by mimicking the proof of
Theorem~\ref{facetsperm}. Another way to proceed is as follows.
Every face of $\mathit{P\!e}^n(\mb a)$ is a set of points where a
certain linear function $\varphi_{\mb b}=(\mb b,\cdot\:)$
restricted to the polytope achieves its minimum. We denote this
face by $G_{\mb b}$. Then $\mb b=(b_1,\ldots,b_{n+1})$ defines an
ordered partition $[n+1]=N_1\sqcup\cdots\sqcup N_k$ according to
the sets of equal coordinates of $\mb b$. Namely, if $b_k=b_l$
then $k$ and $l$ are in the same $N_i$, while if $b_k<b_l$ then
$k\in N_i$ and $l\in N_j$ with $i<j$. It can be shown that (a)
$\dim G_{\mb b}=n+1-k$, and (b) the face $G_{\mb b}$ only depends
on the ordered partition above and does not depend on the
particular values of~$b_k$. In particular, the vertices of
$\mathit{P\!e}^n(\mb a)$ correspond to $\varphi_{\mb b}$ where all
coordinates of $\mb b$ are different, while the facets correspond
to $\varphi_{\mb b}$ where all coordinates of $\mb b$ are either
$0$ or~$1$. We leave the details to the reader.
\end{proof}

We shall denote an $n$-dimensional combinatorial permutahedron by
$\mathit{P\!e}^n$. The classical permutahedron corresponds to $\mb
a=(1,2,\ldots,n+1)$. It can also be obtained as a Minkowski
sum~\eqref{mssim} as follows (we leave the proof as an exercise):

\begin{proposition}\label{clperm}
Let $\mathcal I$ be the collection of all subsets of cardinality
$\le2$ in~$[n+1]$. Then
\[
  P_\mathcal I=\mathit{P\!e}^n(1,2,\ldots,n+1).
\]
\end{proposition}

The polytope $P_\mathcal I$ is the Minkowski sum of segments of
the form $\mathop{\mathrm{conv}}(\!\mb e_i,\mb e_j)$ for $1\le
i<j\le n+1$, shifted by the vector $\mb e_1+\cdots+\mb e_{n+1}$.
Minkowski sums of segments are known as
\emph{zonotopes}\label{zonotop}; this family of polytopes has many
remarkable properties~\cite[\S7.3]{zieg95}. However, general
zonotopes are rarely simple; the permutahedron is one of the few
exceptions.

\subsection*{Building sets and nestohedra}
We now return to general Minkowski sums~\eqref{mssim}. The next
statement gives a description of $P_\mathcal F$ in terms of
inequalities, and generalises the first part of
Theorem~\ref{facetsperm}.

\begin{proposition}[{\cite[Proposition~3.12]{fe-st05}}]\label{pfint}
$P_{\mathcal F}$ can be described as the intersection of the
hyperplane
\[
  H_{\mathcal F}=
  \bigl\{\mb x\in\R^{n+1}\colon\sum_{i=1}^{n+1}x_i=
  |\mathcal F|\bigl\}
\]
with the halfspaces
\[
  H_{T,\ge}=\bigl\{\mb x\in\R^{n+1}\colon\sum_{i\in T}x_i\ge
  \bigl|\mathcal F|_T\bigr|
  \bigl\}
\]
corresponding to all proper subsets $T\subsetneq[n+1]$.
\end{proposition}
\begin{proof}
Since $P_{\mathcal F}$ is a Minkowski summand in the permutahedron
$P_\mathcal C$, it is defined by inequalities of the form
$\sum_{i\in T}x_i\ge b_T$ for some parameters $b_T$. The minimum
value of the linear function $\sum_{i\in T}x_i$ on a simplex
$\varDelta_S$ equals one if $S\subset T$, and equals zero
otherwise. Therefore, $b_T=\bigl|\mathcal F|_T\bigr|$, as needed.
\end{proof}

It follows that $P_\mathcal F$ can be obtained by iteratively
cutting the $n$-simplex
\[
  H_\mathcal F\cap\{\mb x\colon x_i\ge 1\;\text{ for }i=1,\ldots,n+1\}
\]
by the hyperplanes $H_{T,\ge}$ corresponding to subsets $T$ of
cardinality $\ge2$. In the case of the permutahedron, each of
these cuts is nontrivial, that is, the corresponding hyperplane is
not redundant. In general, the description of $P_\mathcal F$ in
Proposition~\ref{pfint} is redundant. The concept of a building
set will allow us to achieve an irredundant description of
$P_\mathcal F$ for certain $\mathcal F$ and therefore describe the
face posets. The resulting polytopes $P_{\mathcal F}$ will be
simple; furthermore they will be obtained from a simplex by a
sequence of face truncations.

\begin{definition}\label{defbs}
A collection $\sB$ of nonempty subsets of $[n+1]$ is called a
\emph{building set} on $[n+1]$ if the following two conditions are
satisfied:
\begin{itemize}
\item[(a)] $S',S'' \in \sB$ with $S'\cap S''\neq \varnothing$ implies $S'\cup S'' \in \sB$;
\item[(b)] $\{i\} \in \sB$ for all $i\in [n+1]$.
\end{itemize}
\end{definition}

\begin{remark}
The terminology comes from a more general notion of a
\emph{building set in a finite lattice}, so that the building set
above corresponds to the case of the Boolean lattice~$2^{[n+1]}$.
We do not give the general definition because it requires too much
poset terminology. See~\cite[\S3]{fe-st05} for the details.
\end{remark}

Note that a building set $\sB$ on $[n+1]$ is \emph{connected} if
and only if $[n+1]\in \sB$. Given $S\subsetneq[n+1]$ define the
\emph{contraction of $S$ from $\sB$}\label{contracti} as
\[
  \sB/S=
  \{T\setminus S\colon T\in\sB,\;T\setminus S\ne\varnothing\}=
  \{ S' \subset [n+1]\!\setminus\! S\colon S' \in \sB\; \text { or }\;S'\cup S \in \sB \}.
\]
The restriction $\sB|_S$ and the contraction $\sB/S$ are building
sets on $S$ and $[n+1]\setminus S$ respectively. Note that
$\sB|_S$ is connected if and only if $S\in\mathcal B$. If $\sB$ is
connected, then $\sB/S$ is also connected for any~$S$.

We now consider polytopes $P_\sB$~\eqref{mssim} corresponding to
building sets~$\sB$. The following specification of
Proposition~\ref{pfint} gives an irredundant description
of~$P_\sB$ as an intersection of halfspaces.

\begin{proposition}[{\cite%
%[Corollary~3.13]
{fe-st05},
\cite%
%[Proposition~7.5]
{post09}}]\label{nestpres} We have
\[
  P_\sB=\Bigl\{\mb x\in\R^{n+1}\colon\sum_{i=1}^{n+1}x_i= |\mathcal
  B|,\quad\sum_{i\in S}x_i\ge
  \bigl|\mathcal B|_S\bigr|\quad\text{for every
  }S\in\sB\Bigr\}.
\]
If $\mathcal B$ is connected, then this representation is
irredundant, that is, every hyperplane $H_S=\bigl\{\mb
x\in\R^{n+1}\colon\sum_{i\in S}x_i= \bigl|\mathcal
B|_S\bigr|\bigl\}$ with $S\ne[n+1]$ defines a facet $F_S$ of
$P_\sB$ (so that the number of facets of $P_\sB$ is $|\sB|-1$).
\end{proposition}
\begin{proof}
The halfspace $H_{T,\ge}$ in the presentation of $P_\sB$ from
Proposition~\ref{pfint} is irredundant if the intersection of
$P_\sB$ with the corresponding hyperplane $H_T$ is a facet. This
intersection is a face of $P_\sB$ given by
\[
  P_{\sB|_T}+P_{\sB/T}
\]
(since $P_{\sB|_T}$ and $P_{\sB/T}$ lie in complementary
subspaces, their Minkowski sum is actually a product). In order
for this face to have codimension one in $P_\sB$, it is necessary,
by Proposition~\ref{basicms}, that the collection $\sB|_T$ is
connected. This condition is equivalent to $T\in\sB$. If $\sB$ is
connected, then this condition is also sufficient, because then
$\sB/T$ is also connected, and $\dim (P_{\sB|_T}+P_{\sB/T})=n-1$
by Proposition~\ref{basicms}.
\end{proof}

\begin{corollary}\label{nestfacet}
If $\sB$ is a connected building set on $[n+1]$, then every facet
of $P_{\sB}$ can be written as
\[
  P_{\sB|_T}\times P_{\sB/T}
\]
for some $T\in\sB\setminus[n+1]$.
\end{corollary}

\begin{theorem}\label{intfacn}
The intersection of facets $F_{S_1}\!\!\cap\cdots\cap F_{S_k}$ is
nonempty (and therefore is a face of $P_\sB$) if and only if the
following two conditions are satisfied:
\begin{itemize}
\item[\rm(a)] for any $i,j$, \ $1\le i<j\le k$, either $S_i\subset S_j$, or $S_j\subset S_i$,
or $S_i\cap S_j=\varnothing$;

\item[\rm(b)] if the sets $S_{i_1},\ldots,S_{i_p}$ are pairwise
nonintersecting, then $S_{i_1}\sqcup\cdots\sqcup
S_{i_p}\notin\sB$.
\end{itemize}
\end{theorem}

\begin{definition}\label{defnestohe}
A subcollection $\{S_1,\ldots,S_k\}\subset\sB$ satisfying
conditions (a) and (b) of Theorem~\ref{intfacn} is called a
\emph{nested set}. Following~\cite{p-r-w08}, we refer to polytopes
$P_\sB$~\eqref{mssim} corresponding to building sets~$\sB$ as
\emph{nestohedra}.
\end{definition}

\begin{proof}[Proof of Theorem~\ref{intfacn}]
Assume $F_{S_1}\!\cap\cdots\cap F_{S_k}\ne\varnothing$.

If (a) fails, then $S_i\cup S_j\in\sB$, and for any $\mb x\in
F_{S_i}\cap F_{S_j}$ we have
\[
  \sum_{k\in S_i}x_k=\bigl|\mathcal B|_{S_i}\bigr|,\quad
  \sum_{k\in S_j}x_k=\bigl|\mathcal B|_{S_j}\bigr|,\quad
  \sum_{k\in S_i\cup S_j}x_k\ge\bigl|\mathcal B|_{S_i\cup
  S_j}\bigr|.
\]
Adding the first two equalities and subtracting the third
inequality we obtain
\[
  \sum_{k\in S_i\cap S_j}x_k \le \bigl|\mathcal B|_{S_i}\bigr|+
  \bigl|\mathcal B|_{S_j}\bigr|-\bigl|\mathcal B|_{S_i\cup
  S_j}\bigr| < \bigl|\mathcal B|_{S_i}\bigr|+
  \bigl|\mathcal B|_{S_j}\bigr|-\bigl|\mathcal B|_{S_i}\cup
  \sB|_{S_j}\bigr|=\bigl|\mathcal B|_{S_i\cap S_j}\bigr|
\]
where the second inequality is strict because $\mathcal
B|_{S_i}\cup\sB|_{S_j}\subsetneq\mathcal B|_{S_i\cup S_j}$. Now
the inequality $\sum_{k\in S_i\cap S_j}x_k < \bigl|\mathcal
B|_{S_i\cap S_j}\bigr|$ contradicts Proposition~\ref{pfint}.

If (b) fails, then $S_{i_1}\sqcup\cdots\sqcup S_{i_p}\in\sB$, and
for any $\mb x\in F_{S_{i_1}}\!\!\cap\cdots\cap F_{S_{i_p}}$ we
have
\[
  \sum_{k\in S_{i_q}}x_k=\bigl|\mathcal B|_{S_{i_q}}\bigr|\quad\text{for
  }1\le q\le p,\quad\text{and }
  \sum_{k\in S_{i_1}\sqcup\cdots\sqcup S_{i_p}}\!\!x_k\;\ge
  \bigl|\mathcal B|_{S_{i_1}\sqcup\cdots\sqcup S_{i_p}}\bigr|.
\]
Subtracting the first $p$ equalities from the last inequality we
obtain
\[
  \bigl|\mathcal B|_{S_{i_1}}\bigr|+\cdots+
  \bigl|\mathcal B|_{S_{i_p}}\bigr|\ge
  \bigl|\mathcal B|_{S_{i_1}\sqcup\cdots\sqcup S_{i_p}}\bigr|.
\]
This leads to a contradiction because $S_{i_1}\sqcup\cdots\sqcup
S_{i_p}\in\sB$.

Now assume that both (a) and (b) are satisfied. We need to show
that $F_{S_1}\cap\cdots\cap F_{S_k}\ne\varnothing$.

We write $\mb x=\sum_{T\in\sB}x^T$ and note that the inequality
$\sum_{i\in S}x_i\ge\bigl|\mathcal B|_{S}\bigr|$ defining the
facet $F_{S}$ turns into equality if and only if $x_i^T=0$ for
every $T\in\sB$, $T\not\subset S$ and $i\in S$ (this follows
from~\eqref{minS}). Hence,
\begin{multline}\label{xinint}
\mb x =\sum_{T\in\sB}x^T\;\in\; F_{S_1}\!\cap\cdots\cap
F_{S_k}\\\Longleftrightarrow\quad x_i^T=0\quad\text{whenever
}T\in\sB,\;T\not\subset S_j,\;i\in S_j,\quad\text{for }1\le j\le
k.
\end{multline}
We therefore need to find $\mb x$ whose coordinates satisfy the
$k$ conditions on the right hand side of~\eqref{xinint}. Given
$T\in\sB$, the $j$th condition is not void only if $T\not\subset
S_j$ and $T\cap S_j\ne\varnothing$. We may assume without the loss
of generality that the first $k'$ conditions in~\eqref{xinint} are
not void, and the rest are void. That is, $T\not\subset S_j$ and
$T\cap S_j\ne\varnothing$ for $1\le j\le k'$, while $T\subset S_j$
or $T\cap S_j=\varnothing$ for $j>k'$. Then we claim that
$T\setminus(S_1\cup\cdots\cup S_{k'})\ne\varnothing$. Indeed,
otherwise choosing among $S_1,\ldots,S_{k'}$ the maximal subsets
$S_{i_1},\ldots,S_{i_p}$ (which are pairwise disjoint by~(a)) we
obtain $S_{i_1}\sqcup\cdots\sqcup S_{i_p}=T\cup
S_{i_1}\cup\cdots\cup S_{i_p}\in\sB$, which contradicts~(b). Now
setting $x^T_i=1$ for only one $i\in T\setminus(S_1\cup\cdots\cup
S_{k'})$ and $x^T_i=0$ for the rest, we obtain the required point
$\mb x$ in the intersection $F_{S_1}\!\cap\cdots\cap F_{S_k}$.
\end{proof}

From the description of the face lattice of nestohedra in
Theorem~\ref{intfacn} it is easy to deduce their following main
property.

\begin{theorem}
Every nestohedron $P_\sB$ is a simple polytope.
\end{theorem}
\begin{proof}
By Proposition~\ref{basicms} we may assume that $\sB$ is
connected. A collection $S_1,\ldots,S_k$ may satisfy both
conditions of Theorem~\ref{intfacn} only if $k\le n$.
\end{proof}

\begin{example}
If $\mathcal B$ is a connected building set on a 2-element set,
then $P_{\mathcal B}$ is an interval~$I^1$. If $\mathcal B$ is a
connected building set on a 3-element set, then $P_\mathcal B$ is
a polygon, and only $m$-gons with $3\le m\le 6$ arise in this way.
\end{example}

More examples will appear in the next subsections.

Proposition~\ref{nestpres} gives a particular way to obtain a
nestohedron $P_\sB$ from a simplex by a sequence of hyperplane
cuts. The next result shows that these hyperplane cuts can be
organised in such a way that we get a sequence of face truncations
(see Construction~\ref{hypcut}).

Let $\sB_0\subset\sB_1$ be building sets on $[n+1]$, and
$S\in\sB_1$. We define a \emph{decomposition of $S$ into elements
of~$\sB_0$} as $S=S_1\sqcup\cdots \sqcup S_k$, where $S_j$ are
pairwise nonintersecting elements of $\sB_0$ and $k$ is minimal
among such disjoint representations of~$S$. It can be easily seen
that this decomposition exists and is unique.

\begin{lemma}\label{nhfacetrun}
Let $\sB_0\subset\sB_1$ be connected building sets on $[n+1]$.
Then $P_{\sB_1}$ is obtained from $P_{\sB_0}$ by a sequence of
truncations at faces $G_i=\bigcap_{j=1}^{k_i} F_{S_j^i}$
corresponding to the decompositions $S^i=S_1^i\sqcup\cdots\sqcup
S_{k_i}^i$ of elements $S^i\in \sB_1\setminus\sB_0$, numbered in
any order that is reverse to inclusion (i.e. $S^i\supseteq
S^{i'}\Rightarrow i\le i'$).
\end{lemma}
\begin{proof}
We use induction on the number $N=|\sB_1|-|\sB_0|$. For $N=1$, we
have $\sB_1=\sB_0\cup \{S^1\}$. We shall show that $P_{\sB_1}$ is
obtained from $P_{\sB_0}$ by a single truncation at the face
$G=F_{S^1_1}\cap\cdots\cap F_{S^1_k}$, where
$S^1=S^1_1\sqcup\cdots\sqcup S^1_k$ is the decomposition of $S^1$
into elements of~$\sB_0$. Let $\widetilde P_{\sB_0}$ denote the
polytope obtained by truncating $P_{\sB_0}$ at~$G$. Since both
$\widetilde P_{\sB_0}$ and $P_{\sB_1}$ are $n$-dimensional
polytopes (here we use the assumption that both $\sB_0$ and
$\sB_1$ are connected), it is enough to verify that the face poset
of $P_{\sB_1}$ is a subposet of the face poset of $\widetilde
P_{\sB_0}$ (see Exercise~\ref{ptopesubposet}).

The facets of $P_{\sB_1}$ are $F_{S^1}$ and $F_{S_j}$ with
$S_j\in\sB_0$. We first consider a nonempty intersection of the
form $F_{S_1}\cap\ldots\cap F_{S_\ell}$ in $P_{\sB_1}$, i.e. a
nested set $\{S_1,\ldots,S_\ell\}$ of~$\sB_1$, with all
$S_j\in\sB_0$. Then obviously, $\{S_1,\ldots,S_\ell\}$ is a nested
set of~$\sB_0$, i.e. $F_{S_1}\cap\ldots\cap
F_{S_\ell}\ne\varnothing$ in $P_{\sB_0}$. Furthermore, since $S^1$
is the only element of $\sB_1\setminus\sB_0$, we have that
\[
  S_{j_1}\sqcup\cdots\sqcup S_{j_p}\ne
  S^1=S^1_1\sqcup\cdots\sqcup S^1_k
\]
for any $\{j_1,\ldots,j_p\}\subset[\ell]$. The latter condition
implies that
$\{S^1_1,\ldots,S^1_k\}\not\subset\{S_1,\ldots,S_\ell\}$, i.e.
$F_{S_1}\cap\ldots\cap F_{S_\ell}\not\subset G$ in the face poset
of $P_{\sB_0}$. By the description of the face poset of
$\widetilde P_{\sB_0}$ given in Construction~\ref{hypcut}, this
implies that $F_{S_1}\cap\ldots\cap F_{S_\ell}\ne\varnothing$ in
$\widetilde P_{\sB_0}$.

Now we consider a nonempty intersection of the form $F_{S^1}\cap
F_{S_1}\cap\ldots\cap F_{S_\ell}$ in $P_{\sB_1}$, i.e. a nested
set $\{S^1,S_1,\ldots,S_\ell\}$ of~$\sB_1$, with $S_j\in\sB_0$ and
$S^1\in\sB_1\setminus\sB_0$. We claim that
$\{S^1_1,\ldots,S^1_k,S_1,\ldots,S_\ell\}$ is a nested set
of~$\sB_0$, i.e. $G\cap F_{S_1}\cap\ldots\cap
F_{S_\ell}\ne\varnothing$ in $P_{\sB_0}$. To do this we need to
verify~(a) and~(b) of Theorem~\ref{intfacn}.

We need to check condition~(a) for pairs of the form $S^1_p,S_q$;
for other pairs it is obvious. That is, we need to check that if
$S^1_p\cap S_q\ne\varnothing$, then one of $S^1_p,S_q$ is
contained in the other. The condition $S^1_p\cap
S_q\ne\varnothing$ implies that $S^1\cap S_q\ne\varnothing$. Since
$\{S^1,S_1,\ldots,S_\ell\}$ is a nested set of~$\sB_1$, we obtain
that $S^1_p\subset S^1\subset S_q$ or $S_q\subset S^1$. By the
minimality of the decomposition $S^1=S^1_1\sqcup\cdots\sqcup
S^1_k$, the inclusion $S_q\subset S^1$ implies that $S_q$ is
contained in some $S^1_r$, which can be only $S^1_p$, since
$S^1_p\cap S_q\ne\varnothing$.

To verify condition (b) of Theorem~\ref{intfacn} for
$\{S^1_1,\ldots,S^1_k,S_1,\ldots,S_\ell\}$, we consider a
subcollection
$\{S^1_{i_1},\ldots,S^1_{i_p},S_{j_1},\ldots,S_{j_q}\}$ consisting
of pairwise nonintersecting subsets. We need to check that its
union is not in~$\sB_0$. For obvious reasons, we may assume that
$p>0$ and $q>0$. Since $\{S^1,S_1,\ldots,S_\ell\}$ is a nested set
of~$\sB_1$, we have that either $S_{j_i}\subset S^1$ or
$S_{j_i}\cap S^1=\varnothing$ for each $i=1,\ldots,q$. Suppose
that $S^1_{i_1}\sqcup\cdots\sqcup S^1_{i_p}\sqcup
S_{j_1}\sqcup\cdots\sqcup S_{j_q}\in\sB_0$. Then $S^1\sqcup
S_{j_1}\sqcup\cdots\sqcup S_{j_q}\in\sB_1$ by the definition of
the building set. If any of $S_{j_i}$ is disjoint with $S^1$, then
we get a contradiction with condition~(b) for the nested set
$\{S^1,S_1,\ldots,S_\ell\}$ of~$\sB_1$. Therefore, $S_{j_i}\subset
S^1$ for $i=1,\ldots,q$, so that $S^1_{i_1}\sqcup\cdots\sqcup
S^1_{i_p}\sqcup S_{j_1}\sqcup\cdots\sqcup S_{j_q}\subset S^1$. By
the argument of the previous paragraph, for each $S_{j_i}$ we have
that $S_{j_i}\subset S^1_r$ or $S^1_r\subset S_{j_i}$ for some
$r=1,\ldots,k$. Then it follows from the minimality of the
decomposition $S^1=S^1_1\sqcup\cdots\sqcup S^1_k$ and the
definition of a building set that $S^1_{i_1}\sqcup\cdots\sqcup
S^1_{i_p}\sqcup S_{j_1}\sqcup\cdots\sqcup S_{j_q}=S^1$, which
contradicts the assumption that $S^1\notin\sB_0$.

Hence, $\{S^1_1,\ldots,S^1_k,S_1,\ldots,S_\ell\}$ is a nested set
of~$\sB_0$. Similarly to the case considered in the previous
paragraph, we also obtain that $F_{S_1}\cap\ldots\cap
F_{S_\ell}\not\subset G$ in the face poset of $P_{\sB_0}$. Once
again, by the description of the face poset of $\widetilde
P_{\sB_0}$ given in Construction~\ref{hypcut}, this implies that
$F_{S^1}\cap F_{S_1}\cap\ldots\cap F_{S_\ell}\ne\varnothing$ in
$\widetilde P_{\sB_0}$.

It follows that the face poset of $P_{\sB_1}$ is indeed contained
as a subposet in the face poset of $\widetilde P_{\sB_0}$, and
thus $P_{\sB_1}=\widetilde P_{\sB_0}$.

It now remains to finish the induction. Assuming the theorem holds
for $M<N$, we shall prove it for $M=N$. Since $S^1$ is not
contained in any other $S^i$, the collection of sets
$\sB'_0=\sB_0\cup \{S^1\}$ is a building set. By the induction
assumption, $P_{\sB'_0}$ is obtained from $P_{\sB_0}$ by
truncation at the face corresponding to the decomposition of
$S^1$, and $P_{\sB_1}$ is obtained from $P_{\sB'_0}$ by a sequence
of truncations corresponding to the decompositions of $S^i$ for
$i=2,\ldots ,N$.
\end{proof}

\begin{remark}
The proof given above only establishes a combinatorial equivalence
between $P_{\sB_1}$ and $\widetilde P_{\sB_0}$. Since
Proposition~\ref{nestpres} gives a geometric presentation of
nestohedra by a sequence of hyperplane cuts, it follows easily
that the face truncations of $P_{\sB_0}$ giving $P_{\sB_1}$ are
also geometric.
%We shall not need this fact though.
\end{remark}

\begin{theorem}\label{nestsimplex}
Every nestohedron $P_{\sB}$ corresponding to a connected building
set $\sB$ can be obtained from a simplex by a sequence of face
truncations.
\end{theorem}
\begin{proof}
Assume that $\sB$ is a connected building set on $[n+1]$. Then we
have $\sS\subset\sB$, where $\mathcal S$ is the connected building
set of Example~\ref{mssimplex}, whose corresponding nestohedron is
an $n$-simplex. Now apply Lemma~\ref{nhfacetrun}.
\end{proof}

The following construction, suggested by N.~Erokhovets, will allow
us to show that every nestohedron can be obtained from a connected
building set, up to combinatorial equivalence.

\begin{construction}[Substitution of building sets]
Let  $\sB_1,\ldots,\sB_{n+1}$ be connected building sets on
$[k_1],\ldots ,[k_{n+1}]$. Then, for every connected building set
$\sB$ on $[n+1]$, we define a connected building set
$\sB(\sB_1,\ldots,\sB_{n+1})$ on
$[k_1]\sqcup\ldots\sqcup[k_{n+1}]=[k_1+\ldots +k_{n+1}]$,
consisting of elements $S^i\in\sB_i$ and $\bigsqcup\limits_{i\in
S}[k_i]$, where $S\in\sB$.

When $\sB_1,\ldots,\sB_n$ are singletons $\{1\},\ldots,\{n\}$, we
shall write $\sB(1,2,\ldots,n,\sB_{n+1})$ instead of
$\sB\bigl(\{1\},\{2\},\ldots,\{n\},\sB_{n+1}\bigr)$.
\end{construction}

\begin{lemma}\label{B_eq1}
Let  $\sB$, $\sB_1,\ldots ,\sB_{n+1}$ be connected building sets
on $[n+1]$, $[k_1],\ldots,[k_{n+1}]$, and let
$\sB'=\sB(\sB_1,\ldots,\sB_{n+1})$. Then $P_{\sB'} \approx P_{\sB}
\times P_{\sB_1}\times\cdots\times P_{\sB_{n+1}}$.
\end{lemma}
\begin{proof}
Set $\sB''=\sB\sqcup\sB_1\sqcup\cdots\sqcup\sB_{n+1}$ and define
the map $\varphi\colon\sB''\to\sB'$ by
\begin{equation*}
  \varphi(S)=\begin{cases} S&\text{ if } S\in\sB_i \\
  \bigsqcup\limits_{i\in S}[k_i]& \text{ if } S\in\sB \; . \end{cases}
\end{equation*}
Then $\varphi$ generates a bijection between $\sB''\setminus
\sB''_{\mathrm{max}}$ and $\sB'\setminus[k_1+\cdots+k_{n+1}]$,
where $\sB''_{\mathrm{max}}=\{[n+1],[k_1],\ldots,[k_{n+1}]\}$
denotes the collection of maximal sets of~$\sB''$. Let
$\mathcal{S}\subset\sB\setminus[n+1]$ and $\mathcal{S}_i\subset
B_i\setminus[k_i]$. Notice that
$\varphi(\mathcal{S})\cup\bigcup_{i=1}^{n+1}\varphi(\mathcal{S}_i)$
is a nested set of $\sB'$ if and only if $\mathcal{S}$ is a nested
set of $\sB$ and $\mathcal{S}_i$ is a nested set of $\sB_i$ for
all~$i$. It follows that $P_{\sB'}\approx P_{\sB''}=P_{\sB}\times
P_{\sB_1}\times\dots\times P_{\sB_{n+1}}$.
\end{proof}

\begin{example}\label{substex}
Assume that each of $\sB$, $\sB_1$, $\sB_2$ is the building set
$\{\{1\},\{2\},\{1,2\}\}$ corresponding to the segment~$\I$. Let
us describe the building set $\sB(\sB_1,\sB_2)$. In the building
set $\{\{a\},\{b\},\{a,b\}\}$, we substitute $a$ by
$\sB_1=\{\{1\},\{2\},\{1,2\}\}$ and $b$ by
$\sB_2=\{\{3\},\{4\},\{3,4\}\}$. As a result, we obtain the
connected building set $\sB'$, consisting of $\{1\},\{2\},\{3\},
\{4\}, \{1,2\}, \{3,4\}, [4]$. Its corresponding nestohedron is
obtained by truncating a 3-simplex at two nonadjacent edges; it is
combinatorially equivalent to a 3-cube.

The facet correspondence $\varphi$ between the combinatorial cubes
$P_{\sB}\times P_{\sB_1}\times P_{\sB_2}$ and $P_{\sB'}$ is given
by
\begin{align*}
\{1\}\in \sB_1&\mapsto\{1\}\in \sB',\qquad
\{2\}\in\sB_1\mapsto\{2\}\in\sB',\\
\{3\}\in \sB_2&\mapsto\{3\}\in \sB',\qquad
\{4\}\in \sB_2\mapsto\{4\}\in \sB',\\
\{a\}\in \sB_2&\mapsto\{1,2\}\in \sB',\quad\{b\}\in
\sB_2\mapsto\{3,4\}\in \sB'.
\end{align*}
\end{example}

\begin{example}\label{B_example}
Let $\sB=\{\{1\},\ldots,\{n+1\},[n+1]\}$ be the building set
corresponding to the simplex $\varDelta^n$ and let $\sB_1,\ldots,
\sB_{n+1}$ be arbitrary connected building sets on
$[k_1],\ldots,[k_{n+1}]$. Then
\[
  \sB'=\sB(\sB_1,\ldots,\sB_{n+1})=(\sB_1\sqcup\dots\sqcup
  \sB_{n+1})\cup[k_1+\dots+k_{n+1}],
\]
and $P_{\sB'}\approx \varDelta^n\times P_{\sB_1}\times\dots\times
P_{\sB_{n+1}}$.
\end{example}

\begin{proposition}\label{B_all_connected}
For each nestohedron $P_{\sB}$ there exists a connected building
set $\sB'$ such that $P_{\sB} \approx P_{\sB'}$.
\end{proposition}
\begin{proof}
Indeed, any building set $\sB$ can be written as
$\sB_1\sqcup\cdots\sqcup\sB_k$, where $\sB_i$ are connected
building sets on $[k_i+1]$. Define a building set
$\widetilde\sB=\sB_1(1,\ldots,k_1,\sB_2)\sqcup\sB_3\sqcup\cdots\sqcup
\sB_k$, giving the same combinatorial polytope. We have that
$\widetilde\sB$ is a product (disjoint union) of $k-1$ connected
building sets. Then we apply again a substitution to
$\widetilde\sB$, and so on. At the end we obtain a connected
building set $\sB'$.
\end{proof}

\subsection*{Graph associahedra}
\begin{definition}\label{graphass}
Let $\Gamma$ be a graph on the vertex set $[n+1]$ without loops
and multiple edges (a \emph{simple} graph). The \emph{graphical
building set} $\mathcal B_\Gamma$ consists of all nonempty subsets
$S\subset[n+1]$ such that the graph $\Gamma|_{S}$ is connected.

The nestohedron $P_\Gamma=P_{\mathcal B_\Gamma}$ corresponding to
a graphical building set is called a
\emph{graph-associahedron}~\cite{ca-de04}.
\end{definition}

\begin{example}[associahedron]\label{defassociah}
Let $\Gamma$ be a `path' with $n$ edges $\{i,i+1\}$ for $1\le i\le
n$. Then $\mathcal B_\Gamma$ consists of all `segments' of the
form $[i,j]=\{i,{i+1},\ldots,j\}$ where $1\le i\le j\le n+1$. To
describe the face poset of the corresponding polytope $P_\Gamma$
it is convenient to use brackets in a string of $n+2$ letters
$a_1a_2\cdots a_{n+2}$. We associate with a segment $[i,j]$ a pair
of brackets before $a_i$ and after $a_{j+1}$. Using
Theorem~\ref{intfacn} it is easy to see that the facets
corresponding to $n$ different segments intersect at a vertex if
and only if the corresponding bracketing of $a_1a_2\cdots a_{n+2}$
with $n$ pairs of brackets is correct. In particular, the vertices
of $\mathcal B_\Gamma$ correspond to all different ways to obtain
a product $a_1a_2\cdots a_{n+2}$ when multiplication is not
associative. The number of vertices is therefore equal to
$\frac1{n+2}\binom{2n+2}{n+1}$, the $(n+1)$th \emph{Catalan
number}. Two vertices are adjacent if and only if the bracketing
corresponding to one vertex can be obtained from the bracketing
corresponding to the other vertex by deleting a pair of brackets
and inserting, in a unique way, another pair of brackets different
from the deleted one. That is, two vertices are adjacent if they
correspond to a single application of the associative law. This
explains the name \emph{associahedron} for the polytope $\mathcal
P_\Gamma$ of this example; we shall denote it $\mathit{As}^n$.
\begin{figure}[h]
\psset{unit=.4cm}
\begin{center}
%Associahedron (Stasheff polytope)
\begin{pspicture}(-4.5,-7)(8.2,6.5)
    \pspolygon(1.3,4.5)(2.3,6)(0.6,6.3)(-0.7,5.5)(-1.7,4)
    \pspolygon[linestyle=dashed](2.3,0)(0.8,2)(-0.7,1)(-1.3,-0.5)(0.3,-1.5)
    \psline(8.2,0)(2.3,6)
    \psline(7.2,-1.4)(1.3,4.5)
    \psline(7.2,-1.4)(8.2,0)
    \psline(-3,-4.5)(7.2,-1.4)
    \psline(-4.4,-3.7)(-3,-4.5)
    \psline(-4.4,-3.7)(-1.7,4)
    \psline[linestyle=dashed](0.8,6.3)(0.8,2)
    \psline[linestyle=dashed](-0.7,5.5)(-0.7,1)
    \psline[linestyle=dashed](-4.4,-3.7)(-1.3,-0.5)
    \psline[linestyle=dashed](-3,-4.5)(0.3,-1.5)
    \psline[linestyle=dashed](2.3,0)(8.2,0)
\end{pspicture}
%Associahedron graph
\begin{pspicture}(-1.8,-1)(2,2)
    \psline[linestyle=dotted](3,3)(5,3)
    \psline[linewidth=1pt](3,3)(3,5)
    \psline[linewidth=1pt](3,3)(1.27,2)
    \psline[linestyle=dotted](1.27,2)(5,3)
    \psline[linestyle=dotted](1.27,2)(3,5)
    \psline[linewidth=1pt](5,3)(3,5)
\end{pspicture}
\vspace{-8mm} \caption{3-dimensional associahedron and the
corresponding graph.}\label{assfig}
\end{center}
\end{figure}

Proposition~\ref{nestpres} describes the associahedron as the
result of consecutive hyperplane cuts of a simplex, see
Figure~\ref{assfig} for the case $n=3$. It can also be obtained by
hyperplane cuts from a cube, as described in the next theorem.

\begin{theorem}\label{asscutcube}
The image of $\mathit{As}^n$ under a certain affine transformation
$\R^{n+1}\to\R^n$ is the intersection of the cube
\[
  \bigl\{\mb y\in\R^n\colon 0\le y_j\le j(n+1-j)\quad
  \text{for }\;1\le j\le n\bigr\}
\]
with the halfspaces
\[
  \bigl\{\mb y\in\R^n\colon y_j-y_k+(j-k)k\ge0 \bigr\}
\]
for $1\le k<j\le n$.
\end{theorem}
\begin{proof}
Proposition~\ref{nestpres} gives the following presentation:
\begin{multline}\label{assx}
  \mathit{As}^n=\Bigl\{\mb x\in\R^{n+1}\colon\sum_{k=1}^{n+1}x_k=\frac{(n+1)(n+2)}2,
  \\\sum_{k=i}^jx_k\ge\frac{(j-i+1)(j-i+2)}2\quad\text{for }1\le i\le j\le
  n+1\Bigr\}.
\end{multline}
Apply the affine transformation $\R^{n+1}\to\R^n$ given by
\[
  (x_1,\ldots,x_{n+1})\mapsto(z_1,\ldots,z_n)\quad\text{where }
  z_l=\sum_{k=1}^l x_k,\quad\text{for } 1\le l\le n.
\]
Now we rewrite the inequalities of~\eqref{assx} in the new
coordinates $(z_1,\ldots,z_n)$. The inequalities with $i=1$
(corresponding to the facets $F_{S}$ with $\{1\}\in S$) become
\begin{equation}\label{ass1}
  z_j\ge\frac{j(j+1)}2, \quad\text{for }1\le j\le
  n.
\end{equation}
Inequalities~\eqref{assx} with $j=n+1$ (corresponding to $F_S$
with $\{n+1\}\in S$) become
\[
  \frac{(n+1)(n+2)}2-z_{i-1}\ge\frac{(n+2-i)(n+3-i)}2,
  \quad\text{for }2\le
  i\le n+1,
\]
or equivalently,
\begin{equation}\label{ass2}
  z_j\le(n+2)j-\frac{j(j+1)}2, \quad\text{for }1\le j\le
  n.
\end{equation}
The remaining inequalities~\eqref{assx} take the form
\begin{equation}\label{ass3}
  z_j-z_{i-1}\ge\frac{(j-i+1)(j-i+2)}2\quad\text{for }1<i\le j<n+1.
\end{equation}
Now the required presentation is obtained from~\eqref{ass1},
\eqref{ass2} and~\eqref{ass3} by applying the shift
$y_j=z_j-\frac{j(j+1)}2$ and setting $k=i-1$.
\end{proof}

\begin{example}\label{2figass} The case $n=3$ of Theorem~\ref{asscutcube} is
shown in Figure~\ref{ascube} (right).
\begin{figure}[h]
\psset{unit=0.5cm}
\begin{center}
\begin{pspicture}(0,0)(24,10)
%left
    \psline[linestyle=dotted](0,0)(1,0)
    \psline(1,0)(7,0)
    \psline[linestyle=dotted](0,6)(7,6)
    \psline(7,0)(7,5)
    \psline[linestyle=dotted](7,5)(7,6)
    \psline[linestyle=dotted](0,0)(0,1)
    \psline(0,1)(0,5)
    \psline[linestyle=dotted](0,5)(0,6)
%   \pspolygon[fillstyle=hlines](1,0)(4,3)(3,4)(0,1)
    \pspolygon[linestyle=dashed](1,0)(4,3)(3,4)(0,1)
    \psline(1,0)(0,1)
    \psline[linestyle=dotted](0,6)(0.5,6.5)
    \psline(0.5,6.5)(3,9)
%   \pspolygon[fillstyle=hlines, hatchangle=70](0,5)(7,5)(7.5,6.5)(0.5,6.5)
    \pspolygon(0,5)(7,5)(7.5,6.5)(0.5,6.5)
    \psline(3,9)(8.5,9)
    \psline[linestyle=dotted](8.5,9)(10,9)
    \psline[linestyle=dotted](7,6)(7.5,6.5)
    \psline(7.5,6.5)(9.5,8.5)
    \psline[linestyle=dotted](9.5,8.5)(10,9)
%   \pspolygon[fillstyle=vlines, hatchangle=60](9.5,2.5)(9.5,8.5)(8.5,9)(8.5,3)
%    \pspolygon[linestyle=dashed](9.5,2.5)(9.5,8.5)(8.5,9)(8.5,3)
%
    \psline[linestyle=dashed](8.5,3)(8.5,9)
    \psline[linestyle=dashed](9.5,2.5)(8.5,3)
    \psline(9.5,2.5)(9.5,8.5)
    \psline(9.5,8.5)(8.5,9)
    \psline[linestyle=dotted](10,9)(10,3)
    \psline(7,0)(9.5,2.5)
    \psline[linestyle=dotted](9.5,2.5)(10,3)
    \psline[linestyle=dotted](0,0)(3,3)
    \psline[linestyle=dotted](3,3)(10,3)
    \psline[linestyle=dotted](3,3)(3,4)
    \psline[linestyle=dashed](4,3)(8.5,3)
    \psline[linestyle=dashed](3,4)(3,9)
%right
  \pspolygon(14,0)(15.9,0)(15.9,6.75)(14,6.75)
  \pspolygon(14,6.75)(18,9)(24,9)(22,7.875)(15.9,6.75)
  \pspolygon(22,2.8125)(23,3.375)(24,6.75)(24,9)(22,7.875)
  \pspolygon(19,0.5625)(20,1.125)(23,3.375)(22,2.8125)
  \psline(15.9,0)(19,0.5625)
  \psline[linestyle=dashed](14,0)(16.1,1.125)
  \psline[linestyle=dashed](16.1,1.125)(20,1.125)
  \psline[linestyle=dashed](16.1,1.125)(18,6.75)
  \psline[linestyle=dashed](18,6.75)(18,9)
  \psline[linestyle=dashed](18,6.75)(24,6.75)
  \put(12.7,4.5){\small$F_{\{1\}}$}
  \put(19,3.75){\small$F_{\{2\}}$}
  \put(23.3,4){\small$F_{\{3\}}$}
  \put(18.5,8){\small$F_{\{4\}}$}
  \put(14.1,3.4){\small$F_{\{1,2\}}$}
  \put(21,1.5){\small$F_{\{2,3\}}$}
  \put(20,9.3){\small$F_{\{3,4\}}$}
  \put(17,-0.3){\small$F_{\{1,2,3\}}$}
  \put(22,6.2){\small$F_{\{\!2,3,4\!\}}$}
\end{pspicture}
%\vspace{-4mm}
\caption{3-dimensional associahedron cut from a
cube.}\label{ascube}
\end{center}
\end{figure}
We start with a 3-cube (more precisely, a parallelepiped) given by
the inequalities
\[
  0\le y_1\le 3,\;0\le y_2\le4,\;0\le y_3\le3,
\]
and cut it by the three hyperplanes
\[
  y_2-y_1+1=0,\;y_3-y_1+2=0,\;y_3-y_2+2=0.
\]
Another way to cut a 3-dimensional combinatorial associahedron
from a 3-cube is shown in Figure~\ref{ascube} (left); this time we
cut three nonadjacent and pairwise orthogonal edges. The two
associahedra in Figure~\ref{ascube} are not affinely equivalent.
\end{example}

The associahedron $\mathit{As}^n$ first appeared (as a
combinatorial object) in the work of Stasheff~\cite{stas63} as the
space of parameters for the higher associativity of the
$(n+2)$-fold product map in an $H$-space.
%Nowadays associahedra
%have become well-known due to the development of the operad theory
%and its applications to symplectic geometry and quantum groups.
For more information about the associahedra we refer
to~\cite{g-k-z94} and~\cite[Lec.~II]{buch08k}, where the reader
may find other geometric and combinatorial realisations
of~$\mathit{As}^n$.
\end{example}

\begin{example}[permutahedron revisited]
Let $\Gamma$ be a complete graph; then $\mathcal B_\Gamma$ is the
complete building set $\mathcal C$ and $P_\Gamma$ is the
permutahedron $\mathit{Pe}^n$, see Figure~\ref{permfig} for the
case $n=3$.
\begin{figure}[h]
\psset{unit=.4cm}
\begin{center}
%Permutahedron
\begin{pspicture}(-4.5,-7)(7.2,6.5)
    \pspolygon(6.3,2.5)(5.3,1)(5.3,-0.5)(5.7,-2.3)(6.6,-1.1)(7.1,0.4)
    \pspolygon(1.3,4.5)(2.3,6)(0.6,6.3)(-0.7,5.5)(-1.7,4)(-.1,3.6)
    \pspolygon(-1.8,-4.4)(-1.8,-6.3)(-3.3,-5.9)(-4.3,-5)(-4.3,-2.6)(-2.8,-3.2)
    \pspolygon[linestyle=dashed](0.8,2)(-0.7,1)(-1.3,-0.5)(-0.3,-2)(1.3,-1.1)(1.8,0.4)
    \psline(-.1,3.6)(-2.8,-3.2)
    \psline(6.3,2.5)(2.3,6)
    \psline(5.3,1)(1.3,4.5)
    \psline(5.3,-.5)(-1.8,-4.4)
    \psline(5.7,-2.3)(-1.8,-6.3)
    \psline(-1.7,4)(-4.3,-2.6)
    \psline[linestyle=dashed](0.8,6.3)(0.8,2)
    \psline[linestyle=dashed](-0.7,5.5)(-0.7,1)
    \psline[linestyle=dashed](-4.3,-5)(-1.3,-0.5)
    \psline[linestyle=dashed](-3.3,-5.9)(-0.3,-2)
    \psline[linestyle=dashed](1.8,0.4)(7.1,0.4)
    \psline[linestyle=dashed](1.3,-1.1)(6.6,-1.1)
\end{pspicture}
%Permutahedron graph
\begin{pspicture}(-1.8,-1)(2,2)
    \psline[linewidth=1pt](3,1)(5,1)
    \psline[linewidth=1pt](3,1)(3,3)
    \psline[linewidth=1pt](3,1)(1.27,0)
    \psline[linewidth=1pt](1.27,0)(5,1)
    \psline[linewidth=1pt](1.27,0)(3,3)
    \psline[linewidth=1pt](5,1)(3,3)
\end{pspicture}
\vspace{-4mm} \caption{3-dimensional permutahedron and the
corresponding graph.}\label{permfig}
\end{center}
\end{figure}
\end{example}

\begin{example}[cyclohedron]\label{defcyclohe}
Let $\Gamma$ be a `cycle' consisting of $n+1$ edges $\{i,i+1\}$
for $1\le i\le n$ and $\{n+1,1\}$. The corresponding $P_\Gamma$ is
known as the \emph{cyclohedron} $\mathit{Cy}^n$, or
\emph{Bott--Taubes polytope}, see Figure~\ref{cyclfig}. It was
first introduced in~\cite{bo-ta94} in connection with the study of
link invariants.
\begin{figure}[h]
\psset{unit=.4cm}
\begin{center}
%Cyclohedron
\begin{pspicture}(-4.5,-7)(7.2,6.5)
    \pspolygon(5.3,1)(5.3,-0.5)(5.7,-2.3)(7.2,0)(6.3,2.5)
    \pspolygon(1.3,4.5)(2.3,6)(0.6,6.3)(-0.7,5.5)(-1.7,4)
    \pspolygon(-4.3,-2.6)(-1.8,-4.4)(-1.8,-6.3)(-3.3,-5.9)(-4.3,-5)
    \pspolygon[linestyle=dashed](2.3,0)(0.8,2)(-0.7,1)(-1.3,-0.5)(-0.3,-2)
    \psline(6.3,2.5)(2.3,6)
    \psline(5.3,1)(1.3,4.5)
    \psline(5.3,-.5)(-1.8,-4.4)
    \psline(5.7,-2.3)(-1.8,-6.3)
    \psline(-1.7,4)(-4.3,-2.6)
    \psline[linestyle=dashed](0.8,6.3)(0.8,2)
    \psline[linestyle=dashed](-0.7,5.5)(-0.7,1)
    \psline[linestyle=dashed](-4.3,-5)(-1.3,-0.5)
    \psline[linestyle=dashed](-3.3,-5.9)(-0.3,-2)
    \psline[linestyle=dashed](2.3,0)(7.2,0)
\end{pspicture}
%Cyclohedron graph
\begin{pspicture}(-1.8,-1)(2,2)
    \psline[linestyle=dotted](3,1)(5,1)
    \psline[linewidth=1pt](3,1)(3,3)
    \psline[linewidth=1pt](3,1)(1.27,0)
    \psline[linewidth=1pt](1.27,0)(5,1)
    \psline[linestyle=dotted](1.27,0)(3,3)
    \psline[linewidth=1pt](5,1)(3,3)
\end{pspicture}
\vspace{-4mm} \caption{3-dimensional cyclohedron and the
corresponding graph.}\label{cyclfig}
\end{center}
\end{figure}
\end{example}

\begin{example}[stellahedron]\label{defstellahe}
Let $\Gamma$ be a `star' consisting of $n$ edges $\{i,{n+1}\}$, \
$1\le i\le n$, emanating from the point~$n+1$. The corresponding
$P_\Gamma$ is known as the \emph{stellahedron} $\mathit{St}^n$,
see Figure~\ref{stelfig}.
\begin{figure}[h]
\psset{unit=.4cm}
\begin{center}
%Stellohedron
\begin{pspicture}(-4.5,-7)(7.2,6.5)
    \pspolygon(6.3,2.5)(5.3,1)(5.3,-0.5)(5.7,-2.3)(7.1,0.2)
    \pspolygon(1.3,4.5)(2.3,6)(0.4,6)(-1.7,4)(-.1,3.6)
    \pspolygon(-1.8,-4.4)(-1.8,-6.3)(-3.8,-5.4)(-4.3,-2.6)(-2.8,-3.2)
    \psline(-.1,3.6)(-2.8,-3.2)
    \psline(6.3,2.5)(2.3,6)
    \psline(5.3,1)(1.3,4.5)
    \psline(5.3,-.5)(-1.8,-4.4)
    \psline(5.7,-2.3)(-1.8,-6.3)
    \psline(-1.7,4)(-4.3,-2.6)
    \psline[linestyle=dashed](0.4,0.2)(7.1,0.2)
    \psline[linestyle=dashed](0.4,0.2)(0.4,6)
    \psline[linestyle=dashed](0.4,0.2)(-3.8,-5.4)
\end{pspicture}
%Stellohedron graph
\begin{pspicture}(-1.8,-1)(2,2)
    \psline[linewidth=1pt](3,1)(5,1)
    \psline[linewidth=1pt](3,1)(3,3)
    \psline[linewidth=1pt](3,1)(1.27,0)
    \psline[linestyle=dotted](1.27,0)(5,1)
    \psline[linestyle=dotted](1.27,0)(3,3)
    \psline[linestyle=dotted](5,1)(3,3)
\end{pspicture}
\vspace{-4mm} \caption{3-dimensional stellahedron and the
corresponding graph.}\label{stelfig}
\end{center}
\end{figure}
\end{example}

\subsection*{Exercises.}\nopagebreak%
\begin{exercise}\label{bscompl}
Every collection $\mathcal F$ of subsets in $[n+1]$ may be
completed in a unique way to a building set by iteratively adding
to $\mathcal F$ the unions $S_1\cup S_2$ of intersecting sets
($S_1\cap S_2\ne\varnothing$) until the process stops. Denote the
resulting building set by $\widehat{\mathcal F}$. Show that
$\widehat{\mathcal F}$ is connected if and only if $\mathcal F$ is
connected, and that $\dim P_{\widehat{\mathcal F}}=\dim
P_{\mathcal F}$ (hint: for a pair of intersecting sets $S_1,S_2$
compare the polytopes $\varDelta_{S_1}+\varDelta_{S_2}$ and
$\varDelta_{S_1}+\varDelta_{S_2}+\varDelta_{S_1\cup S_2}$). Use
this fact to complete the proof of Proposition~\ref{basicms}.
\end{exercise}

\begin{exercise}
Finish the argument in the proof of Proposition~\ref{facesperm}.
\end{exercise}

\begin{exercise}
Complete the details in the proof of Theorem~\ref{facesgenperm}.
\end{exercise}

\begin{exercise}
Prove Proposition~\ref{clperm}.
\end{exercise}

\begin{exercise}
Given two building sets $\sB_0\subset\sB_1$, show that a
decomposition $S=S_1\sqcup\cdots \sqcup S_k$ of $S\in\sB_1$ into
elements $S_i\in\sB_0$ with minimal $k$ is unique.
\end{exercise}

\begin{exercise}\label{trancasim}
A combinatorial polytope obtained from a simplex by a sequence of
face truncations is called a \emph{truncated simplex}. Notice that
a truncated simplex is a simple polytope. By
Theorem~\ref{nestsimplex}, every nestohedron corresponding to a
connected nested set is a truncated simplex. Give an example of
\begin{itemize}
\item[(a)] a simple polytope which is not a truncated simplex;

\item[(b)] a truncated simplex which is not a nestohedron;

\item[(c)] a Minkowski sum of simplices $P_\mathcal F$ given by~\eqref{mssim}
which is not a truncated simplex.
\end{itemize}
\end{exercise}

\begin{exercise}
Consider the following connected building set on~$[4]$:
\[
  \sB=\bigl\{\{1\},\{2\},\{3\},\{4\},\{1,2,3\},\{2,3,4\},[4]\bigr\}.
\]
Then $P_{\sB}$ is combinatorially equivalent to the polytope shown
in Fig.~\ref{figfv}, right.
\end{exercise}

\section{Flagtopes and truncated cubes}\label{flagpolytopes}

\begin{definition}\label{defflp}
A simple polytope $P$ is called a \emph{flagtope} (or \emph{flag
polytope}) if every collection of its pairwise intersecting facets
has a nonempty intersection. (The origins of this terminology will
be explained in Section~\ref{secflc}.)
\end{definition}

Flagtopes and flag simplicial complexes (which are the subject of
Section~\ref{secflc}) feature in the well-known
\emph{Charney--Davis conjecture}~\cite{ch-da95} and its
generalisation due to Gal~\cite{gal05}. According to the Gal
conjecture, the components of the $\gamma$-vector of a flagtope
are nonnegative.

Gal's conjecture is important as it connects the combinatorial
study of polytopes and sphere triangulations to differential
geometry and topology of manifolds. This conjecture has been
proved in special cases.

Although not all nestohedra are flagtopes (a simplex is the
easiest counterexample), flag nestohedra constitute an important
family. In particular, all graph-associahedra (and therefore, the
classical series of associahedra and permutahedra) are flagtopes,
see Proposition~\ref{grassflag}.

As we have seen in Theorem~\ref{asscutcube}, the associahedron can
be obtained by hyperplane cuts from a cube. In this section we
give a proof of a much more general and precise result: a
nestohedron is a flagtope if and only if it can be obtained from a
cube by a sequence of truncations at faces of codimension~2 (we
refer to such polytopes as \emph{$2$-truncated cubes}), see
Theorem~\ref{flnest2cube}. This result was proved by Buchstaber
and Volodin in~\cite[Theorem~6.6]{bu-vo11}. On the other hand, it
can be easily seen that the Gal conjecture is valid for
2-truncated cubes (see Proposition~\ref{cubering}). This
observation (formulated in terms of the dual operation of stellar
subdivision at an edge) was present in the work of Charney--Davis
and later Gal, and used to support their conjectures. As a
corollary we obtain that the Gal conjecture holds for all flag
nestohedra.

\begin{example}\

1. The cube $I^n$ is a flagtope, but the simplex $\varDelta^n$ is
not if $n>1$.

2. The product $P\times Q$ of two flagtopes is a flagtope.

3. The connected sum $P\cs Q$ of two simple $n$-polytopes and the
vertex truncation $\vt(P)$ are not flagtopes if $n>1$.
\end{example}

A polytope $P$ is said to be \emph{triangle-free}\label{trianfre}
if it does not contain a triangular 2-face.

\begin{proposition}\label{flagtrifree}
A flagtope is triangle-free.
\end{proposition}
\begin{proof}
Assume that a flagtope $P$ of dimension~$n$ contains a triangular
2-face~$T$ with vertices $v_1,v_2,v_3$ and edges $e_1,e_2,e_3$,
where $v_i$ is opposite to~$e_i$. Since $P$ is simple, each $e_i$
is an intersection of $n-1$ facets. Hence, for each $i=1,2,3$
there is a unique facet, say~$F_i$, which contains the edge $e_i$
but not the triangle~$T$. Also, $T$ is an intersection of $n-2$
facets, and we may assume that $T=\bigcap_{i=4}^{n+1}F_i$. Now we
observe that $\bigcap_{i=1,i\ne j}^{n+1}F_i=v_j$ for $j=1,2,3$.
This implies that the intersection of any pair of facets among
$F_1,\ldots,F_{n+1}$ is nonempty. On the other hand,
$\bigcap_{i=1}^{n+1}F_i=\varnothing$ because $P$ is simple and
$n$-dimensional. Therefore $P$ cannot be a flagtope.
\end{proof}

The converse to Proposition~\ref{flagtrifree} does not hold in
general (see exercises), but it is valid for polytopes with few
facets:

\begin{theorem}[\cite{bl-bl92}]
If $P$ is a triangle-free convex $n$-polytope then $f_i(P)\ge
f_i(I^n)$ for $i=0,\ldots,n$. In particular, such $P$ has at least
$2n$ facets. Furthermore if $P$ is simple then
\begin{itemize}
\item[(a)] $f_{n-1}(P)=2n$ implies that $P=I^n$;
\item[(b)] $f_{n-1}(P)=2n+1$ implies that $P= P_5\times I^{n-2}$
where $P_5$ is a pentagon;
\item[(c)] $f_{n-1}(P)=2n+2$ implies that $P=P_6\times I^{n-2}$ or $P=Q\times I^{n-3}$ or
$P=P_5\times P_5\times I^{n-4}$ where $P_6$ is a hexagon and $Q$
is a $3$-polytope obtained by truncating a pentagonal prism at one
of its edges forming a pentagonal facet.
\end{itemize}
\end{theorem}

\begin{corollary}
A triangle-free simple $n$-polytope with at most \,$2n+2$\, facets
is a flagtope.
\end{corollary}

Another source of examples of flagtopes is provided by
graph-associahedra (see Definition~\ref{graphass}):

\begin{proposition}\label{grassflag}
Every graph-associahedron $P_\Gamma$ is a flagtope.
\end{proposition}
\begin{proof}
Let $F_{S_1},\ldots,F_{S_k}$ be a set of facets of $P_\Gamma$ with
nonempty pairwise intersections. We need to check that
$F_{S_1}\cap\cdots\cap F_{S_k}\ne\varnothing$, i.e. that
condition~(b) of Theorem~\ref{intfacn} is satisfied (condition~(a)
holds automatically as it depends only on pairwise intersections).
Let $S_{i_1},\ldots,S_{i_p}$ be pairwise nonintersecting sets
among $S_1,\ldots,S_k$; then $S_{i_r}\cup S_{i_s}\notin\mathcal
B_\Gamma$ for $1\le r<s\le p$ because $F_{S_{i_r}}\cap
F_{S_{i_s}}\ne\varnothing$. Therefore, all subgraphs
$\Gamma|_{S_{i_r}\cup S_{i_s}}$ are disconnected, which implies
that $\Gamma|_{S_{i_1}\cup\cdots\cup S_{i_p}}$ is also
disconnected. Thus, $S_{i_1}\cup\cdots\cup S_{i_p}\notin\mathcal
B_\Gamma$, and $F_{S_1}\cap\cdots\cap F_{S_k}\ne\varnothing$ by
Theorem~\ref{intfacn}.
\end{proof}

The main conjecture about the face numbers of flagtopes uses the
notion of the $\gamma$-vector
$\gamma(P)=(\gamma_0,\ldots,\gamma_{[n/2]})$ (see
Definition~\ref{defgamma}):

\begin{conjecture}[Gal~\cite{gal05}]
A flagtope $P$ of dimension $n$ satisfies $\gamma_i(P)\ge0$ for
$i=0,\ldots,\sbr n2$.
\end{conjecture}

We shall verify the Gal conjecture for all flag nestohedra. To do
this we shall show that any flag nestohedron can be obtained by
consecutively truncating a cube at codimension-2 faces.

Our first task is therefore to describe how the $\gamma$-vector
changes under face truncations. For this it is convenient to use
the $H$-polynomial given by~\eqref{hvector}, and the
$\gamma$-polynomial
\[
  \gamma(P)(\tau)=\gamma_0+\gamma_1\tau+\cdots+
  \gamma_{[n/2]}\tau^{[n/2]}.
\]

\begin{proposition}
Let $Q$ be the polytope obtained by truncating a simple
$n$-polytope $P$ at a $k$-dimensional face~$G$. Then
\begin{itemize}
\item[(a)]$H(Q)=H(P)(s,t)+stH(G)H(\varDelta^{n-k-2})$,\\[-7pt]
\item[(b)]$\gamma(Q)=\gamma(P)+\tau\gamma(G)\gamma(\varDelta^{n-k-2})$.
\end{itemize}
\end{proposition}
\begin{proof}
The truncation removes $G$ and creates a face
$G\times\varDelta^{n-k-1}$, so that
\[
  f_i(Q)=f_i(P)+f_i(G\times\varDelta^{n-k-1})-f_i(G),\quad
  \text{for }0\le i\le n.
\]
Hence, $F(Q)=F(P)+tF(G)F(\varDelta^{n-k-1})-t^{n-k}F(G)$, and
\begin{align*}
H(Q) &= H(P)+ t H(G)H(\varDelta^{n-k-1})-t^{n-k} H(G) \\
     &= H(P)+tH(G)\left(\sum_{i=0}^{n-k-1}
       s^i t^{n-k-1-i}-t^{n-k-1}\right) \\
     &= H(P)+s tH(G)\left(\sum_{j=0}^{n-k-2}s^j t^{n-k-2-j}\right)
      = H(P)+s t H(G)H(\varDelta^{n-k-2}) \, ,
\end{align*}
which proves~(a). Furthermore,
\begin{multline*}
\sum_{i=0}^{[\frac{n}{2}]}\gamma_i(Q)(s t)^i(s+t)^{n-2i}=
H(Q)
%\\= H(P)+s t H(G)H(\varDelta^{n-k-2})
= \sum_{i=0}^{[\frac{n}{2}]}\gamma_i(P)(s t)^i(s+t)^{n-2i} \\
  + st\left(\sum_{p=0}^{[\frac{k}{2}]}
  \gamma_p(G)(st)^p(s+t)^{k-2p}\right)
  \left(\sum_{q=0}^{[\frac{n-k-2}{2}]}
  \gamma_q(\varDelta^{n-k-2})(st)^q(s+t)^{n-k-2-2q}\right) \\
= \sum_{i=0}^{[\frac{n}{2}]}\!\gamma_i(P) (s t)^i (s+t)^{n-2i}
  +\sum_{p=0}^{[\frac{k}{2}]}\sum_{p=0}^{[\!\frac{n-k-2}{2}]}
  \!\!\!\!\!\gamma_p(G)\gamma_q(\varDelta^{n-k-2})(st)^{p+q+1}(s+t)^{n-2(p+q+1)} .
\end{multline*}
Hence,
$\gamma_i(Q)=\gamma_i(P)+\sum_{p+q=i-1}\gamma_p(G)\gamma_q(\varDelta^{n-k-2})$,
which proves~(b).
\end{proof}

\begin{definition}\label{truncatedcube}
We refer to a truncation at a codimension-$2$ face as a
\emph{$2$-truncation}. A combinatorial polytope obtained from a
cube by $2$-truncations will be called a \emph{$2$-truncated
cube}.
\end{definition}

The following corollary is the dual of a result from~\cite{gal05}:

\begin{corollary}\label{2truncHgamma}
Let the polytope $Q$ be obtained from a simple polytope $P$ by
$2$-truncation at a face~$G$. Then
\begin{itemize}
\item[(a)]$H(Q)=H(P)+s t H(G)$,
\item[(b)]$\gamma(Q)=\gamma(P)+\tau\gamma(G)$.
\end{itemize}
\end{corollary}

%\begin{lemma}\label{B_flagshave}
%Any $2$-truncation keeps flagness.
%\end{lemma}
%\begin{proof}
%Let $Q$ be obtained from $P$ by truncation of the face $G=F_1\cap
%F_2$. Then $\partial Q^\ast=(v_0, \sigma_G)\partial P^\ast$, where
%$\sigma_G=\{v_1,v_2\}$, and $v_1, v_2$ are the vertices
%corresponding to facets $F_1,F_2$. Let the vertices
%$V\subset\partial Q^\ast$ be pairwise adjacent. Note that one of
%the vertices $v_1, v_2$ is not contained in $V$. The vertices
%$V\setminus \{v_0\}$ are pairwise adjacent in the complex
%$\partial P^\ast$, and then $V\setminus \{v_0\}\in\partial
%P^\ast$. If $v_0\notin V$, then $V\in\partial Q^\ast$ according to
%a) of Remark~\ref{B_fac}. If $v_0\in V$,  then $V\setminus
%\{v_0\}\in\mathrm{lk}(v_0,\partial
%Q^\ast)=\mathrm{lk}(\{v_1,v_2\},\partial P^\ast)$, hence
%$V\in\partial Q^\ast$ according to b) of Remark~\ref{B_fac}.
%\end{proof}

\begin{proposition}\label{truncface}
Each face of a $2$-truncated cube is a $2$-truncated cube.
\end{proposition}
\begin{proof}
It is enough to show that if $P$ is a $2$-truncated cube, then all
the facets of $P$ are $2$-truncated cubes. The proof is by
induction on the number of face truncations. Let the polytope $Q$
be obtained from a $2$-truncated cube $P$ by $2$-truncation at a
face $G$ of codimension~$2$. Then the new facet is $G\times I$,
and it is a $2$-truncated cube by the induction assumption. Every
other facet $F'$ of the polytope $Q$ is either a facet of $P$, or
obtained from a facet $F''$ of $P$ by $2$-truncation at a face
$G'\subset F''$.
\end{proof}

\begin{proposition}\label{cubering}
Any $2$-truncated cube $P$ satisfies $\gamma_i(P)\geq 0$, i.e. the
Gal conjecture holds for $2$-truncated cubes.
\end{proposition}
\begin{proof}
We proceed by induction on the dimension of $P$, using
Proposition~\ref{truncface} and the formula
$\gamma(Q)=\gamma(P)+\tau\gamma(G)$.
\end{proof}

Here is a criterion for a nestohedron to be a flagtope.

\begin{proposition}[\cite{bu-vo11}, \cite{p-r-w08}]
\label{flagnest} Let $\mathcal B$ be a
%connected
building set on $[n+1]$. Then the nestohedron $P_{\mathcal B}$ is
a flagtope if and only if for every element $S\in\sB$ with $|S|>1$
there exist elements $S',S''\in\sB$ such that $S'\sqcup S''=S$.
\end{proposition}
\begin{proof}
Suppose $P_B$ is a flagtope. Consider an element $S\in B$. Then we
may write $S=S_1\sqcup\cdots\sqcup S_k$, where $S_1,\ldots
,S_k\in\sB\setminus \{S\}$ and $k$ is minimal among such
decompositions of~$S$. Then for any subset $J\subset[k]$ with
$1<|J|<k$ we have that $\bigsqcup_{j\in J}S_j\notin\sB$, since
otherwise $k$ can be decreased. If $k>2$ then, by
Theorem~\ref{intfacn}, the facets $F_{S_1},\ldots ,F_{S_k}$ of
$P_{\sB}$ intersect pairwise, but have empty common intersection.
Therefore, $k=2$.

Suppose for each element $S\in\sB$ with $|S|>1$, there exist
elements $S',S''\in\sB$ such that $S'\sqcup S''=S$. Let
$F_{S_1},\ldots, F_{S_k}$, $k\ge3$, be a minimal collection of
facets that intersect pairwise but have empty common intersection.
We shall lead this to a contradiction by finding a nontrivial
subcollection of $F_{S_1},\ldots, F_{S_k}$ with empty common
intersection.

Assume there is a set $\widetilde{S}\in\sB|_S$ intersecting more
than one~$S_i$, but not intersecting every~$S_i$. Then the
subcollection of facets $F_{S_i}$ satisfying $S_i\cap
\widetilde{S}\ne\varnothing$ will have empty common intersection,
since
\[
  \bigsqcup_{S_i\colon S_i\cap \widetilde{S}\ne\varnothing}
  S_i\in\sB
\]
by definition of a building set.

It remains to find $\widetilde{S}\in\sB|_S$ intersecting more than
one~$S_i$, but not intersecting every~$S_i$. By
Theorem~\ref{intfacn}, $S_1\sqcup\dots\sqcup S_k=S\in\sB$.
Therefore, we can write $S=S'\sqcup S''$, where $S',S''\in\sB$.
Let $S^1$ be that of the sets $S'$ and $S''$ which intersects more
elements $S_i$ than the other. Then $S^1$ intersects more than
one~$S_i$. If $S^1$ does not intersect all of $S_i$, then we are
done. Otherwise we write $S^1=S'^1\sqcup S''^1$, where
$S'^1,S''^1\in\sB$, and choose as $S^2$ that of the sets $S'^1$
and $S''^1$ which intersects more elements $S_i$ than the other.
If $S^2$ intersects all the sets $S_i$, choose $S^3$ in the same
way, and so on. Since the cardinality of the sets
$S^1,S^2,S^3,\ldots$ strictly decreases, at some point this
process stops, and we get that one of the sets $S'^i, S''^i$
intersects more than one~$S_i$, but does not intersect
every~$S_i$.
\end{proof}

\begin{proposition}[\cite{p-r-w08}]\label{Bcube}
Let $\sB$ be a building set on $[n+1]$ such that $P_{\sB}$ is a
flagtope. Then there exists a building set $\sB_0\subset\sB$ such
that $P_{\sB_0}$ is a combinatorial cube with $\dim P_{\sB_0}=\dim
P_{\sB}$.
\end{proposition}
\begin{proof}
By Proposition~\ref{basicms}, we need to consider only connected
building sets. For $n=1$, the proposition is true. Assuming that
the assertion holds for $m<n$, we shall prove it for $m=n$. By
Proposition~\ref{flagnest}, we have $[n+1]=S'\sqcup S''$, where
$S',S''\in\sB$. By the induction assumption, the building sets
$\sB|_{S'}$ and $\sB|_{S''}$ have subsets $\sB'_0$ and $\sB''_0$
whose corresponding nestohedra are cubes. The building set
$\sB_0=(\sB'_0\sqcup\sB''_0)\cup[n+1]$ is the desired one (see
Example~\ref{B_example}).
\end{proof}

It now follows from Lemma~\ref{nhfacetrun} that a flag nestohedron
can be obtained from a cube by a sequence of face truncations. The
following lemma shows that there is a sequence consisting only of
codimension-2 face truncations:

\begin{lemma}\label{B12}
Let $\sB_1\subset\sB_2$ be connected building sets on~$[n+1]$
whose corresponding nestohedra $P_{\sB_1}$ and $P_{\sB_2}$ are
flagtopes. Then
\begin{itemize}
\item[(a)]
$P_{\sB_2}$ is obtained from $P_{\sB_1}$ by a sequence of
$2$-truncations;
\item[(b)]
$\gamma_i(P_{\sB_1})\le\gamma_i(P_{\sB_2})$ for
$i=0,1,\ldots,[n/2]$.
\end{itemize}
Furthermore, if $\sB_1\ne\sB_2$, then the inequality above is
strict for some~$i$.
\end{lemma}
\begin{proof}
Let $S$ be a minimal (by inclusion) element of $\sB_2\setminus
\sB_1$. We set $\sB'=\widehat{\sB_1\cup\{S\}}$ to be the minimal
(by inclusion) building set containing $\sB_1\cup \{S\}$. By
Proposition~\ref{flagnest}, there exist $S',S''\in\sB_2$ such that
$S'\sqcup S''=S$. It follows from the choice of $S$ that $S',S''
\in\sB_1$. It is easy to show that $\sB'$ is the collection of
sets $\sB_1\cup\{T=T'\sqcup T''\colon T',T''\in\sB_1, S'\subset
T', S''\subset T''\}$. Hence, the decomposition of any element of
$\sB'\setminus\sB_1$ consists of two elements. Therefore, by
Lemma~\ref{nhfacetrun}, the nestohedron $P_{\sB'}$ is obtained
from $P_{\sB_1}$ by a sequence of 2-truncations.

Since $B_1\subsetneq\sB'\subset\sB_2$, we can finish the proof of
(a) by induction on the number of elements in
$\sB_2\setminus\sB_1$. Statement~(b) follows from
Corollary~\ref{2truncHgamma}.
\end{proof}

Finally we can prove the main result of this section, giving a
characterisation of flag nestohedra:

\begin{theorem}[\cite{bu-vo11}]\label{flnest2cube}
A nestohedron $P_{\sB}$ is a flagtope if and only if it is a
2-truncated cube.

More precisely, if $P_{\sB}$ is a flagtope, then there exists a
sequence of building sets $\sB_0 \subset\sB_1\subset \cdots
\subset \sB_N=\sB$, where $P_{\sB_0}$ is a combinatorial cube,
$\sB_i=\sB_{i-1}\cup \{S_i\}$, and $P_{\sB_i}$ is obtained from
$P_{\sB_{i-1}}$ by $2$-truncation at the face $F_{S_{j_1}}\cap
F_{S_{j_2}}\subset P_{\sB_{i-1}}$, where $S_i=S_{j_1}\sqcup
S_{j_2}$, and $S_{j_1}, S_{j_2}\in\sB_{i-1}$.
\end{theorem}
\begin{proof}
By Proposition~\ref{B_all_connected}, we need to consider only
connected building sets. The `only if' statement follows from
Proposition~\ref{Bcube} and Lemma~\ref{B12}. To prove the `if'
statement we need to check that a polytope obtained from a
flagtope by a 2-truncation is a flagtope. This is left as an
exercise.
\end{proof}

Together with Proposition~\ref{cubering},
Theorem~\ref{flnest2cube} implies

\begin{corollary}\label{galnest}
The Gal conjecture holds for all flag nestohedra~$P_{\mathcal B}$,
i.e. $\gamma_i(P_{\mathcal B})\ge0$.
\end{corollary}

\begin{example}
Let us see how Theorem~\ref{flnest2cube} works in the case of the
$3$-dimensional associahedron. The building set corresponding to
$As^3$ is given by
\[
  \sB=\bigl\{\{1\}, \{2\}, \{3\},
  \{4\}, \{1,2\},\{2,3\},\{3,4\},\{1,2,3\},\{2,3,4\},
  \{1,2,3,4\}\bigr\}
\]
(see Fig.~\ref{ascube}, right). In order to obtain $As^3$ from
$I^3$ by $2$-truncations, we have to specify a building set
$\sB_0\subset\sB$, such that $P_{\sB_0}\approx I^3$, and to order
the elements of $\sB \setminus \sB_0$ in such a way that adding a
new element to the building set corresponds to a $2$-truncation.

First let $\sB_0$ consist of $\{i\},\{1,2\},\{3,4\},[4]$. The
associahedron $P_{\sB}$ is then obtained from $P_{\sB_0}\approx
I^3$ by consecutive truncation at the faces $F_{\{1,2\}} \cap
F_{\{3\}}, F_{\{2\}} \cap F_{\{3,4\}}, F_{\{2\}} \cap F_{\{3\}}$
in this order. (Warning: $P_{\sB_0}$ is \emph{not} the rectangular
parallelepiped, it is a tetrahedron with two opposite edges
truncated. However the three 2-truncations of $P_{\sB_0}$
described here and the three 2-truncations of the rectangular
parallelepiped described in Example~\ref{2figass} give the same
(up to affine equivalence) polytope $As^3$ shown in
Fig.~\ref{ascube}, right.)

Another option is to let $\sB_0$ consist of
$\{i\},\{1,2\},\{1,2,3\},[4]$. Then $P_{\sB_0}$ is obtained from a
tetrahedron by truncating a vertex and then truncating an edge
which contained this vertex; we have that $P_{\sB_0}\approx I^3$.
To obtain the associahedron $P_{\sB}$ from $P_{\sB_0}$ we first
truncate the face $F_{\{2\}}\cap F_{\{3\}}$ of $P_{\sB_0}\approx
I^3$ and get the new facet $F_{\{2,3\}}$. Then we truncate the
faces $F_{\{2,3\}}\cap F_{\{4\}}$ and $F_{\{3\}}\cap F_{\{4\}}$.
\end{example}

\subsection*{Exercises}
\begin{exercise}
In a flagtope, a collection of pairwise intersecting faces (of any
dimension) has a nonempty intersection.
\end{exercise}

\begin{exercise}
A face of a flagtope is a flagtope. This generalises
Proposition~\ref{flagtrifree}.
\end{exercise}

\begin{exercise}
Give an example of a triangle-free simple $n$-polytope with $2n+3$
facets which is not a flagtope.
\end{exercise}

\begin{exercise}
Give an example of a non-flag simple polytope $P$ with
$\gamma_i(P)\ge0$.
\end{exercise}

\begin{exercise}
A polytope obtained from a flagtope by a 2-truncation is a
flagtope.
\end{exercise}

\begin{exercise}[{\cite[Theorem~9.1]{bu-vo11}}]
The face vectors of $n$-dimensional flag nestohedra $P_{\sB}$
satisfy
\begin{itemize}
\item[(a)]$\gamma_i(I^n)\le\gamma_i(P_{\sB})\le\gamma_i(\mathit{Pe}^n)$ for $i=0,1,\ldots,[n/2]$;
\item[(b)]$g_i(I^n)\le g_i(P_{\sB})\le g_i(\mathit{Pe}^n)$ for $i=0,1,\ldots,[n/2]$;
\item[(c)]$h_i(I^n)\le h_i(P_{\sB})\le h_i(\mathit{Pe}^n)$ for $i=0,1,\ldots,n$;
\item[(d)]$f_i(I^n)\le f_i(P_{\sB})\le f_i(\mathit{Pe}^n)$ for
$i=0,1,\ldots,n$.
\end{itemize}
Furthermore, the lower bounds are achieved only for $P_{\sB}=I^n$
and the upper bounds are achieved only for
$P_{\sB}=\mathit{Pe}^n$. (Hint: to prove (a) use Lemma~\ref{B12}.
The other inequalities follow from Proposition~\ref{faceineq}.)
\end{exercise}

\begin{exercise}[{\cite[Theorem~9.2]{bu-vo11}}]
The face vectors of graph associahedra $P_{\Gamma}$ corresponding
to connected graphs $\Gamma$ on~$[n+1]$ satisfy
\begin{itemize}
\item[(a)]$\gamma_i(\mathit{As}^n)\le\gamma_i(P_{\Gamma})\le\gamma_i(\mathit{Pe}^n)$ for $i=0,1,\ldots,[n/2]$;
\item[(b)]$g_i(\mathit{As}^n)\le g_i(P_{\Gamma})\le g_i(\mathit{Pe}^n)$ for $i=0,1,\ldots,[n/2]$;
\item[(c)]$h_i(\mathit{As}^n)\le h_i(P_{\Gamma})\le h_i(\mathit{Pe}^n)$ for $i=0,1,\ldots,n$;
\item[(d)]$f_i(\mathit{As}^n)\le f_i(P_{\Gamma})\le f_i(\mathit{Pe}^n)$ for
$i=0,1,\ldots,n$.
\end{itemize}
Furthermore, the lower bounds are achieved only for
$P_{\sB}=\mathit{As}^n$ and the upper bounds are achieved only for
$P_{\sB}=\mathit{Pe}^n$.
\end{exercise}

\section{Differential algebra of combinatorial
polytopes}\label{sec:dgrcp} Here we develop a
differential-algebraic formalism that will allow us to analyse the
combinatorics of families of polytopes and their face numbers from
the viewpoint of differential equations. All polytopes in this
section are combinatorial.

\subsection*{Ring of polytopes}
Denote by $\mathfrak P^{2n}$ the free abelian group generated by
all combinatorial $n$-polytopes. The group $\mathfrak P^{2n}$ is
infinitely generated for $n>1$, and it splits into a direct sum of
finitely generated components as follows:
\[
  \mathfrak P^{2n}=\bigoplus_{m\ge n+1}\mathfrak P^{2n,2(m-n)}
\]
where $\mathfrak P^{2n,2(m-n)}$ is the group generated by
$n$-polytopes with $m$ facets. (The fact that there are finitely
many combinatorial polytopes with a given number of facets is left
as an exercise.)

\begin{definition}\label{ringofptopes}
The product of polytopes turns the direct sum
\[
  \mathfrak P=\bigoplus_{n\ge0}\mathfrak P^{2n}=
  \mathfrak P^0+\bigoplus_{m\ge\,2}\bigoplus_{n=1}^{m-1}\mathfrak P^{2n,2(m-n)}
\]
into a bigraded commutative ring, the \emph{ring of polytopes}.
The unit is $P^0$, a point.

Simple polytopes generate a bigraded subring of $\mathfrak P$,
which we denote by~$\mathfrak S$.
\end{definition}

A polytope is \emph{indecomposable} if it cannot be represented as
a product of two other polytopes of positive dimension.

\begin{proposition}[{\cite[Proposition~2.22]{bu-er11}}]
$\mathfrak P$ is a polynomial ring generated by indecomposable
combinatorial polytopes.
\end{proposition}
\begin{proof}
We need to show that any polytope $P$ of positive dimension can be
represented as a product of indecomposable polytopes,
$P=P_1\times\cdots\times P_k$, and this representation is unique
up to reordering the factors. The existence of such a
decomposition is clear. Now assume that $P_1\times\cdots\times
P_k\approx Q_1\times\cdots\times Q_s$. Fix a vertex
$v=v_1\times\cdots\times v_k=w_1\times\cdots\times w_s$, where
$v_i\in P$ and $w_j\in Q_j$. The faces of the form
$v_1\times\cdots\times P_i\times\cdots\times v_k$ and
$w_1\times\cdots\times Q_j\times\cdots\times w_s$ are maximal
indecomposable faces of~$P$ containing the vertex~$v$, and
therefore they bijectively correspond to each other under the
combinatorial equivalence. It follows that $k=s$ and $P_i\approx
Q_i$ for $i=1,\ldots,k$ up to a permutation of factors.

Since $\mathfrak P$ is a free abelian group generated by
combinatorial polytopes, and each polytope can be uniquely
represented by a monomial in indecomposable polytopes, it follows
that $\mathfrak P$ is a polynomial ring on indecomposable
polytopes.
\end{proof}

Given $P\in\mathfrak P^{2n}$, denote by $dP\in\mathfrak
P^{2(n-1)}$ the sum of all facets of $P$ in the ring~$\mathfrak
P$. The following lemma is straightforward.

\begin{lemma}
$d\colon\mathfrak P\to\mathfrak P$ is a linear operator of degree
$-2$ satisfying the identity
\[
  d(P_1P_2)=(dP_1)P_2+P_1(dP_2).
\]
\end{lemma}

Therefore, $\mathfrak P$ is a differential ring, and $\mathfrak S$
is its differential subring.

\begin{example}\label{difsimcub}
We have
\[
  dI^n = n(dI)I^{n-1} = 2nI^{n-1},\quad\text{and}\quad
  d\varDelta^n=(n+1)\varDelta^{n-1}.
\]
\end{example}

\subsection*{Face-polynomials revisited}
By Proposition~\ref{multFH}, the $F$-polynomial and the
$H$-polynomial define ring homomorphisms
\[
  F\colon\mathfrak P \longrightarrow \Z[s,t],\quad
  H\colon\mathfrak P \longrightarrow \Z[s,t],
\]
which send $P\in\mathfrak P$ to $F(P)(s,t)$ and $H(P)(s,t)$
respectively.

%\begin{example}
%We have $F(I^1)=s+2t$ and $F(\varDelta^2)=s^2+3st+3t^2$.
%\end{example}

\begin{theorem}\label{diffpol}
For any simple polytope $P$ we have
\begin{equation}\label{fdp}
  F(dP) = \frac{\partial}{\partial t}F(P).
\end{equation}
\end{theorem}
\begin{proof}
Assume that $P$ is a simple $n$-polytopse with facets
$P_1,\ldots,P_m$. Then
\[
  F(dP)=\sum_{i=1}^m
  F(P_i)=\sum_{i=1}^m\sum_{k=0}^{n-1}f_{k}(P_i)s^kt^{n-1-k}.
\]
On the other hand,
\[
  \frac{\partial}{\partial
  t}F(P)=\sum_{k=0}^{n-1}(n-k)f_{k}(P)s^kt^{n-1-k}.
\]
Comparing the coefficients in the two sums above we
reduce~\eqref{fdp} to
\begin{equation}\label{fksimple}
  \sum_{i=1}^mf_{k}(P_i)=(n-k)f_{k}(P).
\end{equation}
Since $P$ is simple, every $k$-face is contained in exactly $n-k$
facets, so it is counted $n-k$ times on the left hand side of the
above identity.
\end{proof}

%\begin{proposition}\label{fpolsimcub}
%We have $F(I^n)=(s+2t)^n$ and
%$F(\varDelta^n)=\frac{(s+t)^{n+1}-t^{n+1}}s$.
%\end{proposition}
%\begin{proof}
%The first identity holds because $F$ is a ring homomorphism. We
%prove the second one by induction. For $n=0$ we have
%$F(\varDelta^0)=1=\frac{(s+t)-t}s$. Assuming the formula is valid in
%dimensions $\le n-1$, we calculate using~\eqref{fdp}
%\[
%  \frac{\partial}{\partial
%  t}F(\varDelta^n)=F(d\varDelta^n)=(n+1)F(\varDelta^{n-1})=
%  (n+1)\frac{(s+t)^n-t^n}s.
%\]
%Integrating over $t$ we obtain
%\[
%  F(\varDelta^n)=\frac{(s+t)^{n+1}-t^{n+1}}s+C(s),
%\]
%where $C(s)$ depends only on~$s$. Since $F(\varDelta^n)(s,0)=s^n$, we
%finally get $C(s)=0$.
%\end{proof}

\begin{corollary}
Let $P_1$ and $P_2$ be two simple $n$-polytopes such that
$dP_1=dP_2$ in $\mathfrak S$. Then $F(P_1)=F(P_2)$.
\end{corollary}
\begin{proof}
Indeed, $F(dP_1)=F(dP_2)$ implies that $\frac{\partial
F(P_1)}{\partial t} =  \frac{\partial F(P_2)}{\partial t}$. Using
that $F(P)(s,0)=s^n$ we obtain $F(P_1)=F(P_2)$.
\end{proof}

\begin{proposition}\label{funiq}
Let $\widetilde F\colon\mathfrak S \to \mathbb{Z}[s,t]$ be a
linear map such that
\[
  \widetilde F(dP) = \frac{\partial}{\partial t}\widetilde F(P)\quad\text{and }\;
  \widetilde F(P)|_{t=0}=s^n.
\]
Then $\widetilde F(P)=F(P)$.
\end{proposition}
\begin{proof}
We have $\widetilde F(P^0)=1=F(P^0)$. Assume by induction the
statement is true in dimensions $\le n-1$, and let $P$ be a simple
$n$-polytope. Then $\widetilde F(dP)(s,t)=F(dP)(s,t)$. It follows
that $\frac{\partial}{\partial t}\widetilde F(P) =
\frac{\partial}{\partial t} F(P)$. Therefore, $\widetilde
F(P)(s,t)=F(P)(s,t)+C(s)$. Setting $t=0$, we obtain
$s^n=s^n+C(s)$, whence $C(s)=0$.
\end{proof}

Theorem~\ref{diffpol} also allows us to reduce the
Dehn--Sommerville relations~\eqref{DSf} to the Euler formula:

\begin{theorem}\label{FDS}
The following identity holds for any simple $n$-polytope $P$:
\begin{equation}\label{newDS}
  F(P)(s,t)=F(P)(-s,s+t).
\end{equation}
\end{theorem}
\begin{proof}
We have $F(P^0)(s,t)=1=F(P^0)(-s,s+t)$. Assume by induction the
identity holds in dimensions $\le n-1$. Then for a given $P$ of
dimension $n$ we have $F(dP)(s,t)=F(dP)(-s,s+t)$. Therefore,
$\frac{\partial}{\partial t}F(P)(s,t)=\frac{\partial}{\partial
t}F(P)(-s,s+t)$ and $F(P)(s,t)-F(P)(-s,s+t)=C(s)$. By Euler's
formula~\eqref{euler},
\[
  F(-s,s)=\bigl(f_0-f_1+\cdots+(-1)^nf_n\bigl)s^n=s^n.
\]
Therefore, $C(s)=F(P)(s,0)-F(P)(-s,s) =0$.
\end{proof}

\begin{example}
For a simple 3-polytope $P^3$ with $m=f_2$ facets we have
$f_1=3(m-2)$, $f_0=2(m-2)$. Assume that all facets are $k$-gons,
so $d P^3=mP_k^2$. Then by Theorem~\ref{diffpol},
\[
  m(s^2+kst+kt^2)=\frac\partial{\partial t}
  \bigl(s^3+ms^2t+3(m-2)st^2+2(m-2)t^3\bigr),
\]
which implies $m(6-k)=12$. It follows that the pair $(k,m)$ may
only take values $(3,4)$, $(4,6)$ and (5,12) corresponding to a
simplex, cube and dodecahedron.
\end{example}

For general polytopes the difference
$\delta=Fd-\frac{\partial}{\partial t}F$ measures the failure of
identity~\eqref{fdp}.

A polytope is said to be \emph{$k$-simple}\label{ksimpleptope}
(for $k\ge0$) if each of its $k$-faces is an intersection of
exactly $n-k$ facets. For example, a simple polytope is
$0$-simple, while 1-simple polytopes are also known as
\emph{simple in edges}. Every $n$-polytope is $(n-1)$- and
$(n-2)$-simple.

A polytope is \emph{$k$-simplicial}\label{ksimplicial} if each of
its $k$-dimensional faces is a simplex. A simplicial $n$-polytope
is $(n-1)$-simplicial, and every polytope is 1-simplicial. By
polarity, if an $n$-polytope $P$ is $k$-simple, then $P^*$ is
$(n-1-k)$-simplicial.

\begin{theorem}\label{deltap}
Let $P\in\mathfrak P$. Then the following identity holds:
\[
  F(dP) = \frac{\partial}{\partial t}F(P) + \delta(P),
\]
where $\delta\colon\mathfrak P\to\Z[s,t]$ is an $F$-derivation,
i.e. it is linear and satisfies the identity
\[
  \delta(P_1P_2)=\delta(P_1)F(P_2)+F(P_1)\delta(P_2).
\]
Moreover, if $P$ is an $n$-polytope, then
$\delta(P)=\delta_2s^{n-3}t^2+\cdots+\delta_{n-1}t^{n-1}$ and
$\delta_i\ge0$ for $2\le i\le n-1$. Also, $P$ is $k$-simple if and
only if $\delta_{n-1-k}=0$; in this case $\delta_i=0$ for $2\le
i\le n-1-k$.
\end{theorem}
\begin{proof}
The fact that $\delta=Fd-\frac{\partial}{\partial t}F$ is an
$F$-derivation is verified by a direct computation.
By~\eqref{fksimple} the coefficient of $s^kt^{n-1-k}$ in
$\delta(P)$ is given by
\[
  \delta_{n-1-k}=\sum_{i=1}^mf_k(P_i)-(n-k)f_k(P).
\]
It vanishes for $k=n-1$ and $k=n-2$ (as each codimension-two face
of $P$ is contained in exactly two facets). Also, $\delta_{n-1-k}$
is non-negative because every $k$-face is contained in at least
$n-k$ facets. Finally, if $P$ is $k$-simple, then every $j$-face
is contained in exactly $n-j$ facets for $j\ge k$, so
$\delta_{n-1-j}$ vanishes for $j\ge k$.
\end{proof}

\begin{example}
Let $P$ be a simple $n$-polytope, and $P^*$ the dual simplicial
polytope. Since the face poset of $P^*$ is the opposite to the
face poset of~$P$,
\[
  F(P^*)(s,t)=s^n+t\,\frac{F(P)(t,s)-t^n}s=s^n+\sum_{k=0}^{n-1}f_{k}(P)s^{n-1-k}t^{k+1}\,.
\]
We have $dP^*=f_{0}(P)\varDelta^{n-1}$, therefore,
\[
  \delta(P^*)=f_{0}(P)\frac{(s+t)^n-t^n}s-
  \frac{\partial}{\partial
  t}\sum_{k=0}^nf_{k}(P)s^{n-1-k}t^{k+1}.
\]
The coefficient of $s^{n-1-k}t^k$ on the right hand side above is
given by
\[
  \delta_k(P^*)=\binom{n}{k}f_{0}(P)-(k+1)f_{k}(P),\quad\text{for }1\le k\le n.
\]
This is non-negative by Theorem~\ref{deltap} or by
Exercise~\ref{newub}.
\end{example}
%
%The inequalities of Theorem~\ref{ubtsimple} may be rewritten as
%\[
%  \frac{f_{k}(P)}{f_{0}(P)} \le\frac{f_{k}(\varDelta^n)}{f_{0}(\varDelta^n)}
%  = \frac{\binom{n+1}{k+1}}{n+1},
%\]
%which implies that from all simple $n$-polytopes $P$ the simplex
%has the maximal ratio $\frac{f_k}{f_0}$. The case $k=1$ is the
%identity $2f_{1}=nf_{0}$, see~\eqref{vered}.
%
%For simple polytopes Theorem~\ref{ubtsimple} gives a better bound
%than the UBT (Theorem~\ref{ubt}), which is valid for all
%polytopes. Indeed, for $k\le\sbr n2-1$ the UBT gives the bound
%$f_{k}\le\bin m{k+1}$, where $m=f_{0}$, and a direct calculations
%shows that $\frac1{k+1}\bin nk m\le\bin m{k+1}$.

The following properties of $H(P)$ follow from the corresponding
properties of $F(P)$ established above and the identity
$H(P)(s,t)= F(P)(s-t,t)$.

\begin{theorem}\
\begin{itemize}
\item[\rm(a)] The ring homomorphism
$H \colon \mathfrak{P}\longrightarrow \Z[s,t]$ satisfies
$H(P)\big|_{t=0}=s^n$. The restriction of $H$ to the subring
$\mathfrak S$ of simple polytopes satisfies the equation
\[
  H(dP)=\partial H(P)
\]
where $\partial=\frac{\partial}{\partial
s}+\frac{\partial}{\partial t}$.

\item[\rm(b)] The image of $H$ is generated by $H(\varDelta^1)=s+t$ and
$H(\varDelta^2)=s^2+st+t^2$.

\item[\rm(c)] If $\widetilde H \colon \mathfrak S\longrightarrow\Z[s,t]$
is a linear map satisfying $\widetilde H(dP) =
\partial\widetilde H(P)$ and $\widetilde H(P)\big|_{t=0}=s^n$ for any simple $n$-polytope $P$,
then $\widetilde H(P)=H(P)$.
\end{itemize}
\end{theorem}

\subsection*{Ring of building sets}
Let $\sB_1$ and $\sB_2$ be building sets on $[n_1]$ and $[n_2]$
respectively. A \emph{map $f\colon \sB_1\to \sB_2$ of building
sets} is a map $f \colon [n_1]\to[n_2]$ satisfying $f^{-1}(S)\in
\sB_1$ for every $S\in \sB_2$. Two building sets $\sB_1$ and
$\sB_2$ are said to be \emph{equivalent}, if there exist maps
$f\colon\sB_1 \to \sB_2$ and $g\colon\sB_2 \to \sB_1$ of building
sets such that the compositions $f\circ g$ and $g\circ f$ are the
identity maps.

The \emph{product}\label{prodbse} of $\sB_1$ and $\sB_2$ is the
building set $\sB_1\cdot \sB_2$ on $[n_1+n_2]$ induced by
appending the interval $[n_2]$ to the interval $[n_1]$.

\begin{example}
The collection $\mathcal S$ of Example~\ref{mssimplex} and the
complete collection $\mathcal C$ are both connected building sets
on~$[n+1]$. They are initial and terminal connected building sets
respectively, since for every connected building set $\mathcal B$
on $[n+1]$ there are maps $\mathcal C\to \mathcal B\to\mathcal S$
of building sets induced by the identity map on~$[n+1]$.
\end{example}

\begin{definition}
Denote by  $\mathfrak{B}^{2n}$ the free abelian group generated by
the equivalence classes of building sets on $[k]$ for $k\le n+1$.
Since $\sB_1\cdot \sB_2$ is equivalent to $\sB_2\cdot \sB_1$, the
product turns $\mathfrak B=\bigoplus_{n\ge0}\mathfrak{B}^{2n}$
into a commutative associative ring.
%The ring $\mathfrak{B}$ is
%multiplicatively generated by connected building sets.

Given a connected building set $\sB$ on $[n+1]$, define
\begin{equation}\label{dbuild}
  d\sB = \sum_{S\in \sB\backslash [n+1]}\sB|_S\cdot \sB/S
\end{equation}
where the sum is taken in the ring $\mathfrak B$. Since $\mathfrak
B$ is generated by connected building sets, we may extend $d$ to a
linear map $d\colon\mathfrak B\to\mathfrak B$ using the Leibniz
formula
\begin{equation}\label{difbuild}
  d(\sB_1\cdot\sB_2) = (d\sB_1)\cdot \sB_2 + \sB_1\cdot (d\sB_2).
\end{equation}
We therefore get a derivation of~$\mathfrak B$.
\end{definition}

\begin{proposition}
The map $\mathcal B\mapsto P_{\mathcal B}$ induces a
%n injective
differential ring homomorphism $\beta\colon\mathfrak B\to\mathfrak
P$. Its image is a graded subring with unit in $\mathfrak P$
multiplicatively generated by nestohedra $P_{\mathcal B}$
corresponding to connected building sets~$\mathcal B$.
\end{proposition}
\begin{proof}
The fact that $\beta$ is a ring homomorphism is obvious. It
follows from Corollary~\ref{nestfacet} and~\eqref{difbuild} that
$\beta$ commutes with the differentials. Every building set
consisting only of singletons is mapped to $P^0$ which represents
a unit in~$\mathfrak P$. Finally, it follows from
Proposition~\ref{basicms}~(a) that the image $\beta(\mathfrak B)$
is generated by $P_{\mathcal B}$ with connected $\mathcal B$.
\end{proof}

\begin{remark}
The ring $\mathfrak B$ does not have unit, and the map $\beta$ is
not injective.
%On the other hand, it is injective on connected
%building sets. That is, if $\mathcal B_1$ and $\mathcal B_2$ are
%connected building sets, then a combinatorial equivalence of
%$P_{\mathcal B_1}$ and $P_{\mathcal B_2}$ implies an equivalence
%of $\mathcal B_1$ and~$\mathcal B_2$. This follows from
%Theorem~\ref{intfacn}.
\end{remark}

Given a subset $S\subset[n+1]$, denote by $\Gamma|_S$ the
restriction of $\Gamma$ to the vertex set $S$, and denote by
$\Gamma/S$ the graph with the vertex set $[n+1]\backslash S$
having an edge between two vertices $i$ and $j$ whenever they are
path connected in $\Gamma_{S\cup \{i,j\}}$.
%The following proposition describes the differential of a
%graph-associahedron $P_\Gamma$ in the ring $\mathfrak B$
%or~$\mathfrak P$.

\begin{proposition} For a connected graph $\Gamma$ on
$n+1$ vertices, we have
\[
  dP_\Gamma=\mathop{\sum_{S\subsetneq[n+1]}}%
  _{\Gamma|_{\!S}\text{\rm\ is connected}}
  P(\Gamma|_S)\times P(\Gamma/S).
\]
\end{proposition}
\begin{proof}
Follows directly from~\eqref{dbuild}.
\end{proof}

\begin{example}\label{difseries}
We have the following formulae for the differentials of the four
graph-associahedra defined in Section~\ref{nestohedra}:
\begin{align*}
&d\,\mathit{As}^n = \sum_{i+j=n-1}(i+2)\mathit{As}^i\times\mathit{As}^j;\\
&d\,\mathit{P\!e}^n = \sum_{i+j=n-1}\binom{n+1}{i+1}\mathit{P\!e}^i\times\mathit{P\!e}^j;\\
&d\,\mathit{Cy}^n = (n+1)\sum_{i+j=n-1}\mathit{As}^i\times Cy^j;\\
&d\,St^n = n\cdot St^{n-1} +
\sum_{i=0}^{n-1}\binom{n}{i}St^i\times\mathit{P\!e}^{n-i-1}.
\end{align*}
For example (see Figs.~\ref{assfig}--\ref{stelfig}),
\begin{align*}
&d\,\mathit{As}^3 = 2As^0\times As^2 + 3As^1\times As^1 +
4As^2\times
As^0\;;\\
&d\,\mathit{P\!e}^3 = 4Pe^0\times Pe^2 + 6Pe^1\times Pe^1 +
4Pe^2\times
Pe^0\;;\\
&d\,\mathit{Cy}^3 = 4(As^0\times Cy^2 + As^1\times Cy^1 +
As^2\times
Cy^0)\;;\\
&d\,\mathit{St}^3 = 3St^2 + St^0\times\mathit{P\!e}^2 +
3St^1\times\mathit{P\!e}^1 + 3St^2\times\mathit{P\!e}^0\;.
\end{align*}
\end{example}

\subsection*{Exercises.}
\begin{exercise}
Show that there are finitely many combinatorial polytopes with a
given number of facets.
\end{exercise}

\begin{exercise}
Show that $\delta=Fd-\frac{\partial}{\partial t}F$ is an
$F$-derivation.
\end{exercise}

\section{Families of polytopes and differential equations}
In this section, using the language of generating series, we
interpret the formulae for the differential of nestohedra and
graph-associahedra as certain partial differential equations.
These differential equations encode the combinatorial information
of the face structure of nestohedra. This exposition builds upon
the results of~\cite{bu-ko07} and~\cite{buch08s}.

Denote by $\mathfrak P[q]$ the polynomial ring in an indeterminate
$q$ with coefficients in~$\mathfrak P$.
\begin{proposition}
Let
\[
  Q\colon\mathfrak P\to\mathfrak P[q],\quad P\mapsto Q(P;q)
\]
be a linear map such that
\[
  Q(dP;q)=\frac\partial{\partial q}Q(P;q)\quad\text{and}\quad
  Q(P;0)=P
\]
for any polytope $P$. Then
\[
  Q(P;q)=\sum_{k=0}^n(d^k\!P)\frac{q^k}{k!}.
\]
\end{proposition}
\begin{proof}
Use induction by dimension, as in the proof of
Proposition~\ref{funiq}.
\end{proof}

Now assume given a sequence $\mathcal P=\{P^n,\;n\ge0 \}$ of
polytopes, one in each dimension. We define its \emph{generating
series}\label{genserie} as the formal power series
\[
  \mathcal P(x) = \sum_{n \ge0}\lambda_nP^nx^{n+n_0}
\]
in $\mathfrak{P}\otimes \mathbb{Q}[[x]]$. The parameter $n_0$ and
the coefficients $\lambda_n$ will be chosen depending on a
particular sequence~$\mathcal P$. Using the transformation $Q$ of
the previous proposition we may define the following 2-parameter
extension of the generating series:
\begin{equation}\label{2pargs}
  \mathcal P(q,x) = \sum_{n \geqslant 0}\lambda_nQ(P^n;q)x^{n+n_0}.
\end{equation}
We have $\mathcal P(0,x) =\mathcal P(x)$.

We consider the following generating series of the six sequences
of nestohedra:
\begin{equation}\label{6gs}
\begin{aligned}
\varDelta(x) &= \sum_{n\geqslant
0}\varDelta^n\frac{x^{n+1}}{(n+1)!}\;;&
I(x) &= \sum_{n\geqslant 0}I^n\frac{x^{n}}{n!}\;;\\[2pt]
\mathit{As}(x) &= \sum_{n\geqslant 0}\mathit{As}^n x^{n+2}\;;&
\mathit{P\!e}(x) &= \sum_{n\geqslant 0}\mathit{P\!e}^n\frac{x^{n+1}}{(n+1)!}\;;\\[2pt]
\mathit{Cy}(x) &= \sum_{n\geqslant
0}\mathit{Cy}^n\frac{x^{n+1}}{n+1}\;;&\mathit{St}(x) &=
\sum_{n\geqslant 0}\mathit{St}^n\frac{x^{n}}{n!}\;.
\end{aligned}
\end{equation}

\begin{lemma}\label{dpx}
The differentials of the generating series above are given by
\begin{align*}
&d\varDelta(x) = x\varDelta(x)\;;&& dI(x) = 2xI(x)\;;\\[3pt]
&d\,\mathit{As}(x) = \mathit{As}(x)\frac{d}{dx}\mathit{As}(x)\;;&
&d\,\mathit{P\!e}(x) = \mathit{P\!e}^2(x)\;;\\[3pt]
&d\,\mathit{Cy}(x) = \mathit{As}(x)\frac{d}{dx}\mathit{Cy}(x)\;;&
&d\,\mathit{St}(x) = \big(x+\mathit{P\!e}(x)\big)\mathit{St}(x)\;.
\end{align*}
\end{lemma}
\begin{proof}
This follows from the formulae of Examples~\ref{difsimcub}
and~\ref{difseries}.
\end{proof}

\begin{theorem}\label{6pde}
The two-parameter extensions of the generating series~\eqref{6gs}
satisfy the following partial differential equations:
\begin{align*}
\frac{\partial}{\partial q}\,\varDelta(q,x) &= x\varDelta(q,x)\;;&
\frac{\partial}{\partial q}\,I(q,x) &= 2xI(q,x)\;;\\[5pt]
\frac{\partial}{\partial q}\,As(q,x) &=
As(q,x)\frac{\partial}{\partial x}\,As(q,x)\;;&
\frac{\partial}{\partial q}\,P\!e(q,x) &= P\!e^2(q,x)\;;\\[5pt]
\frac{\partial}{\partial q}\,C\!y(q,x) &=
As(q,x)\frac{\partial}{\partial x}\,C\!y(q,x)\;;&
\frac{\partial}{\partial q}\,St(q,x) &=
\big(x+P\!e(q,x)\big)St(q,x)\;.
\end{align*}
\end{theorem}
\begin{proof}
A direct calculation using formulae~\eqref{2pargs}
and~\eqref{6gs}.
\end{proof}

\begin{remark}
The role of the parameters $\lambda_n$ in~\eqref{2pargs} can be
illustrated as follows. If we replace the first series
$\varDelta(x)$ of~\eqref{6gs} by
\[
  \widehat \varDelta(x) = \sum_{n\geqslant
  0}\varDelta^n\frac{x^{n+1}}{n+1},
\]
then the first equations of Lemma~\ref{dpx} and Theorem~\ref{6pde}
take the form
\[
  d\widehat\varDelta(x) = x^2\frac d{dx}\widehat\varDelta(x),\quad
  \frac{\partial}{\partial q}\,\widehat\varDelta(q,x) = x^2
  \frac d{dx}\widehat\varDelta(q,x).
\]
\end{remark}

Four of the equations of Theorem~\ref{6pde}, namely those
corresponding to the series $\varDelta$, $I$, $\mathit{P\!e}$ and
$St$, are ordinary differential equations. Their solutions are
completely determined by the initial data $\mathcal P(0,x)
=\mathcal P(x)$ and are given by the explicit formulae
\begin{align*}
  \varDelta(q,x) &= \varDelta(x)e^{qx}\;;&
  I(q,x) &= I(x)e^{2qx}\;;\\[2pt]
  \mathit{P\!e}(q,x) &=\frac{\mathit{P\!e}(x)}{1-q\mathit{P\!e}(x)}\;;&
  St(q,x) &= St(x)\frac{e^{qx}}{1-q\mathit{P\!e}(x)}\;.
\end{align*}

The equation for $U=As(q,x)$ has the form $U_q=UU_x$. This
classical quasilinear partial differential equation was considered
by E.~Hopf, and therefore became known as the \emph{Hopf
equation}\label{Hopfequ}.

\begin{theorem}\
\begin{itemize}
\item[(a)]
The series $As(q,x)$ is given by the solution of the functional
equation (\emph{equation on characteristics})
\begin{equation}\label{eqonchar}
  As(q,x)=As(x+qAs(q,x)),
\end{equation}
where $As(x)=As(0,x)$.

\item[(b)] The series $C\!y(q,x)$ is given by the solution of
\begin{equation}\label{funccy}
  C\!y(q,x)=C\!y(x+qAs(q,x)),
\end{equation}
where $C\!y(x)=C\!y(0,x)$.
\end{itemize}
\end{theorem}
\begin{proof}
Set $U=As(q,x)$ and $As_x=\frac d{dx}As(x)$. If $U$ is a solution
to~\eqref{eqonchar}, then we obtain by differentiating
\[
  U_q=(U+qU_q)As_x,\quad U_x=(1+qU_x)As_x.
\]
Therefore, $(1-qAs_x)U_q=UAs_x$ and $(1-qAs_x)U_x=As_x$, which
implies that $U$ satisfies the Hopf equation $U_q=UU_x$. Its
solution with the initial condition $U(0,x)=As(x)$ is unique by
the general theory of quasilinear equations (in our case the
uniqueness can be also verified using standard arguments with
power series).

Similarly, by differentiating~\eqref{funccy} we obtain for
$V=C\!y(q,x)$:
\[
  V_q=(U+qU_q)C\!y_x,\quad V_x=(1+qU_x)C\!y_x.
\]
Using that $U_q=UU_x$ we rewrite the first equation above as
$V_q=U(1+qU_x)C\!y_x$, which implies $V_q=UV_x$ as claimed. This
is exactly the equation for $V=C\!y(q,x)$ given by
Theorem~\ref{6pde}, and its solution with $V(0,x)=C\!y(x)$ is
unique.
\end{proof}

We can also use Lemma~\ref{dpx} to calculate the face-polynomials
$F(s,t)$ of graph-associahedra. Let $\mathcal P(x)$ be one of the
generating series~\eqref{6gs}, and set
\[
  F_{\mathcal P}=F(\mathcal P(x))=\sum_{n\ge0}\lambda^nx^{n+n_0}
  \sum_{k=0}^nf_{k}(P^n)s^kt^{n-k}.
\]
We refer to $F_{\mathcal P}=F_{\mathcal P}(s,t;x)$ as the
\emph{generating series of face-polynomials}\label{genesefp}; it
is a series in $x$ whose coefficients are polynomials in $s$
and~$t$.

\begin{theorem}\label{6fgs}
The generating series of face-polynomials corresponding
to~\eqref{6gs} satisfy the following differential equations, with
the initial conditions given in the second column:
\begin{align*}
\frac{\partial}{\partial t}\,F_\varDelta &= xF_\varDelta,&
F_\varDelta(s,0;x)&=\frac{e^{sx}-1}s\;;\\[5pt]
\frac{\partial}{\partial t}\,F_I &= 2xF_I,&
F_I(s,0;x)&=e^{sx}\;;\\[5pt]
\frac{\partial}{\partial t}\,F_{As} &=
F_{As}\frac{\partial}{\partial x}\,F_{As},&
F_{As}(s,0;x)&=\frac{x^2}{1-sx}\;;\\[5pt]
\frac{\partial}{\partial t}\,F_{Pe} &= F_{Pe}^2,&
F_{Pe}(s,0;x)&=\frac{e^{sx}-1}s\;;\\[5pt]
\frac{\partial}{\partial t}\,F_{C\!y} &=
F_{As}\frac{\partial}{\partial x}\,F_{C\!y},&
F_{C\!y}(s,0;x)&=-\frac{\ln(1-sx)}s\;;\\[5pt]
\frac{\partial}{\partial t}\,F_{St} &= \big(x+F_{Pe}\big)F_{St},&
F_{St}(s,0;x)&=e^{sx}\;.
\end{align*}
\end{theorem}
\begin{proof}
The differential equations are obtained by applying $F$ to the
equations of Lemma~\ref{dpx} and using the fact that
$F(dP)=\frac\partial{\partial t}F(P)$. The initial conditions
follow by substituting $s^n$ for $P^n$ in~\eqref{6gs} and
calculating the resulting series.
\end{proof}

Again, the four equations for the series $F_\varDelta$, $F_I$,
$F_{\mathit{P\!e}}$ and $F_{\mathit{St}}$ are ordinary
differential equations. Their solutions are completely determined
by the initial data; the explicit formulae are left as exercises.
The remaining two partial differential equations for the series
$F_{\mathit{As}}$ and $F_{\mathit{C\!y}}$ can be explicitly solved
as follows.

\begin{theorem}\
\begin{itemize}
\item[(a)]
The series $U=F_{\mathit{As}}$ satisfies the quadratic equation
\begin{equation}\label{quadu}
  t(s+t)U^2+(2tx+sx-1)U+x^2=0.
\end{equation}
The initial condition $F_{\mathit{As}}(s,0;x)=\frac{x^2}{1-sx}$
determines its solution uniquely.

\item[(b)] The series $F_{\mathit{C\!y}}$ is given by
\[
  F_{\mathit{C\!y}}=-\frac1s\ln\bigl(1-s(x+tF_{\mathit{As}}))\bigr).
\]
\end{itemize}
\end{theorem}
\begin{proof} (a) By analogy with~\eqref{eqonchar} we show that
$U=F_{\mathit{As}}$ satisfies
% the functional equation
\[
  U=\varphi(x+tU),
\]
where $\varphi(x)=F_{\mathit{As}}(s,0;x)=\frac{x^2}{1-sx}$. It is
equivalent to~\eqref{quadu}.

(b) We have that $V=F_{\mathit{C\!y}}$ is given by the solution to
$V_t=UV_x$. By analogy with~\eqref{funccy} we show that it is
given by
\[
  V=\psi(x+tU),
\]
where $U=F_{\mathit{As}}$ and
$\psi(x)=F_{\mathit{C\!y}}(s,0;x)=-\frac{\ln(1-sx)}s$.
\end{proof}

As a corollary, we shall derive a formula for the number of
$k$-faces in $\mathit{As}^n$, which equals the number of
bracketing of a word of $n+2$ letters with $n-k$ pairs of
brackets. These numbers were first calculated by Cayley in 1891:

\begin{theorem}
The number of $k$-dimensional faces in an $n$-dimensional
associahedron is given by
\[
  f_k(\mathit{As}^n)=\frac1{n+2}\binom nk\binom{2n-k+2}{n+1}.
\]
\end{theorem}
\begin{proof}
We use the fact that $F_{As}$ satisfies the Hopf equation, whose
solutions may be obtained using conservation laws. Let
\[
  U(t,x)=\sum_{k\ge0}U_k(x)t^k
\]
be the solution of the Cauchy problem for the Hopf equation:
\begin{equation}\label{hopffi}
  U_t=UU_x,\qquad U(0,x)=\varphi(x).
\end{equation}
This equation has the following conservation laws
\[
  \biggl(\frac{U^{k}}{k}\biggr)_t=
  \biggl(\frac{U^{k+1}}{k+1}\biggr)_x,\quad\text{for } k\ge1.
\]
Hence,
\[
  \frac{d^k}{dt^k}\,U=
  \frac{d^{k-1}}{dt^{k-1}}\biggl(\frac{U^{2}}{2}\biggr)_x=
  \frac{d^{k-2}}{dt^{k-2}}\biggl(\frac{U^{3}}{3}\biggr)_{xx}=\cdots=
  \frac{d^k}{dx^k}\biggl(\frac{U^{k+1}}{k+1}\biggr),
  \quad\text{for } k\ge1.
\]
Evaluating at $t=0$ we obtain
\[
  \frac{d^k}{d t^k}\,U\bigg|_{t=0}=k!\,U_k(x)=
  \frac{d^k}{dx^k}\biggl(\frac{U_0^{k+1}(x)}{k+1}\biggr).
\]
Therefore,
\[
  U_k(x)=\frac{1}{(k+1)!}\,\frac{d^k}{d x^k}\,\varphi^{k+1}(x).
\]

By Theorem~\ref{6fgs} the function
\begin{equation}\label{fasstx}
  U=F_{As}(s,t;x)=\sum_{n\ge0}\sum_{k=0}^nf_{n-k}(\mathit{As}^n)s^{n-k}t^kx^{n+2}
\end{equation}
satisfies the Hopf equation~\eqref{hopffi} with the initial
function $\varphi(s;x)=\frac{x^2}{1-sx}$. We therefore calculate
\begin{align*}
  U_k(s;x)&=\frac{1}{(k+1)!}\,\frac{d^k}{dx^k}\biggl(\frac{x^{2(k+1)}}{(1-sx)^{k+1}}\biggr)\\
  &=\frac{1}{(k+1)!}\,\frac{d^k}{dx^k}\biggl(x^{2(k+1)}\sum_{l\ge0}\binom{l+k}ls^lx^l\biggr)\\
  &=\frac{1}{(k+1)!}\sum_{l\ge0}\binom{l+k}k\frac{(2k+l+2)!}{(k+l+2)!}s^lx^{k+l+2}\\
  &=\sum_{n\ge k}\frac{1}{n+2}\binom nk
  \binom{n+k+2}{n+1}s^{n-k}x^{n+2}.
\end{align*}
On the other hand, it follows from~\eqref{fasstx} that
\[
  U_k(s;x)=\sum_{n\ge k}f_{n-k}(\mathit{As}^n)s^{n-k}x^{n+2}.
\]
Comparing the last two formulae we obtain
\[
  f_{n-k}(\mathit{As}^n)=\frac{1}{n+2}\binom nk\binom{n+k+2}{n+1},
\]
which is equivalent to the required formula.
\end{proof}

\subsection*{Exercises.}
\begin{exercise}
By solving the first two differential equations of
Theorem~\ref{6fgs} show that the generating series for the
$F$-polynomials of simplices and cubes are given by
\begin{align*}
  F_\varDelta(s,t;x)&=\sum_{n\ge0}F(\varDelta^n)\frac{x^{n+1}}{(n+1)!}=
  e^{tx}\frac{e^{sx}-1}s,\\
  F_I(s,t;x)&=\sum_{n\ge0}F(I^n)\frac{x^n}{n!}=
  e^{tx}e^{(s+t)x}.
\end{align*}
Compare this with the formulae for $F(\varDelta^n)$ and~$F(I^n)$.
\end{exercise}

\begin{exercise}
Show by solving the corresponding differential equations from
Theorem~\ref{6fgs} that the generating series for the
face-polynomials of permutahedra and stellahedra are given by
\begin{align*}
  F_{\mathit{P\!e}}(s,t;x)&=\sum_{n\ge0}F(\mathit{P\!e}^n)\frac{x^{n+1}}{(n+1)!}=
  \frac{e^{sx}-1}{s-t(e^{sx}-1)},\\
  F_{\mathit{St}}(s,t;x)&=\sum_{n\ge0}F(\mathit{St}^n)\frac{x^n}{n!}=
  \frac{se^{(s+t)x}}{s-t(e^{sx}-1)}.
\end{align*}
Compute the face-polynomials $F(\mathit{P\!e}^n)$ and
$F(\mathit{St}^n)$ and the face numbers explicitly.
\end{exercise}

\begin{exercise}
Define the generating series $H_{\mathcal P}(s,t;x)$ of
$H$-polynomials by analogy with the generating series of
face-polynomials. Deduce the following formulae for sequences of
nestohedra:
\begin{align*}
  H_\varDelta(s,t;x)&=\sum_{n\ge0}H(\varDelta^n)\frac{x^{n+1}}{(n+1)!}=
  \frac{e^{sx}-e^{tx}}{s-t},\\
  H_I(s,t;x)&=\sum_{n\ge0}H(I^n)\frac{x^n}{n!}=
  e^{(s+t)x},\\
  H_{\mathit{P\!e}}(s,t;x)&=\sum_{n\ge0}H(\mathit{P\!e}^n)\frac{x^{n+1}}{(n+1)!}=
  \frac{e^{sx}-e^{tx}}{se^{tx}-te^{sx}},\\
  H_{\mathit{St}}(s,t;x)&=\sum_{n\ge0}H(\mathit{St}^n)\frac{x^n}{n!}=
  \frac{(s-t)e^{(s+t)x}}{se^{tx}-te^{sx}}.
\end{align*}
\end{exercise}

\begin{exercise}
The series
$Y=H_{\mathit{As}}(s,t;x)=\sum_{n\ge0}H(\mathit{As}^n)x^{n+2}$
satisfies the quadratic equation
\[
  Y=(x+sY)(x+tY).
\]
The initial condition $H_{\mathit{As}}(s,0;x)=\frac{x^2}{1-sx}$
determines its solution uniquely.
\end{exercise}

\begin{exercise}
The components of the $h$-vector of $\mathit{As}^n$ are given by
\[
  h_k=\frac1{n+1}\binom{n+1}k\binom{n+1}{k+1},\qquad 0\le k\le n.
\]
\end{exercise}

%Добавить ссылку на Миллера-Штурмфелься в начале раздела о двойственности Александера
%Комбинаторное пространство Тома нормального вложения (конец раздела о двойств Александера)
%Экзотические (n-1)-сферы не вложимые в R^n?
%Ссылка на Александера в теореме о звёздных преобразованиях
%Ссылка на Гайфуллина в разделе о бизвёздных преобр
%Ссылка на определение копредела в разделе о симпл ч.у.м.
%Добавить дополнительный раздел о минимальных триангуляциях

\setcounter{chapter}1
\chapter{Combinatorial structures}\label{combigt}
The face poset of a convex polytope is a classical example of a
combinatorial structure underlying a decomposition of a geometric
object. With the development of combinatorial topology several new
combinatorial structures emerged, such as simplicial and cubical
complexes, simplicial posets and other types of regular cellular
(or CW) complexes. Many of these structures have eventually become
objects of independent study in geometric combinatorics.

A simplicial complex is the abstract combinatorial structure
behind a simplicial subdivision (or
\emph{triangulation})\label{triangu} of a topological space.
Triangulations were first introduced by Poincar\'e and provide a
rigorous and convenient tool for studying topological invariants
of smooth manifolds by combinatorial methods. The notion of a
\emph{nerve} of a covering of a topological space, introduced by
Alexandroff, provides another source of examples of simplicial
complexes.

The study of triangulations stimulated the development of first
combinatorial and then algebraic topology in the first half of the
XXth century. With the appearance of cellular complexes algebraic
tools gradually replaced the combinatorial ones in mainstream
topology. Simplicial complexes still play a pivotal role in $PL$
(piecewise linear) topology, however nowadays the main source of
interest in them is in discrete and computational geometry. One
reason for that is the emergence of computers, since simplicial
complexes provide the most effective way to translate geometric
and topological structures into machine language.

We therefore may distinguish two different views on the role of
simplicial complexes and triangulations. In topology, simplicial
complexes and their different derivatives such as singular chains
and simplicial sets are used as technical tools to study the
topology of the underlying space. Most combinatorial invariants of
nerves or triangulations (such as the number of faces of a given
dimension) do not have meaning in topology, as they do not reflect
any topological feature of the underlying space. Topologists
therefore tend not to distinguish between simplicial complexes
that have the same underlying topology. For instance, refining a
triangulation (such as passing to the barycentric subdivision)
changes the combinatorics drastically, but does not affect the
underlying topology. On the other hand, in combinatorial geometry
the combinatorics of a simplicial complex is what really matters,
while the underlying topology is often simple or irrelevant.

In toric topology the combinatorist's point of view on
triangulations and similar decompositions is enriched by elaborate
topological techniques. Combinatorial invariants of triangulations
therefore can be analysed by topological methods, and at the same
time combinatorial structures such as simplicial complexes or
posets become a source of examples of topological spaces and
manifolds with nice features and lots of symmetry, e.g. bearing a
torus action. The combinatorial structures are the subject of this
chapter; the associated topological objects will come later.

Here we assume only minimal knowledge of topology. The reader may
also check Appendix~\ref{algtop} for the definition of simplicial
homology groups, etc.

\section{Polyhedral fans}\label{combfan}
Like convex polytopes, polyhedral fans encode both geometrical and
combinatorial information. The geometry of fans is poorer than
that of convex polytopes, but this geometry is still a part of the
structure of a fan, which distinguishes fans from the purely
combinatorial objects considered later in this chapter.

Although fans were considered in convex geometry independently,
the main source of interest to them is in the theory of
\emph{toric varieties}, which are classified by rational fans.
Toric varieties are the subject of Chapter~\ref{toric}, and here
we describe the terminology and constructions related to fans.

\begin{definition}\label{cpcone}
A set of vectors $\mb a_1,\ldots,\mb a_k\in\R^n$ defines a
\emph{convex polyhedral cone}, or simply~\emph{cone},
\[
  \sigma=\R_\ge\langle\mb a_1,\ldots,\mb a_k\rangle
  =\{\mu_1\mb a_1+\cdots+\mu_k\mb a_k\colon\mu_i\in\R_\ge\}.
\]
Here $\mb a_1,\ldots\mb a_k$ are referred to as \emph{generating
vectors} (or \emph{generators}) of~$\sigma$. A \emph{minimal} set
of generators of a cone is defined up to multiplication of vectors
by positive constants. A cone is \emph{rational} if its generators
can be chosen from the integer lattice $\Z^n\subset\R^n$. If
$\sigma$ is a rational cone, then its generators $\mb
a_1,\ldots\mb a_k$ are usually chosen to be \emph{primitive}, i.e.
each $\mb a_i$ is the smallest lattice vector in the ray defined
by it.
%, i.e. each $\mb a_i$ is an integer vector and is
%not an integer multiple of an integer vector.

A cone is \emph{strongly convex} if it does not contain a line. A
cone is \emph{simplicial} if it is generated by a part of basis
of~$\R^n$, and is \emph{regular} if it is generated by a part of
basis of~$\Z^n$. (Regular cones play a special role in the theory
of toric varieties, see Chapter~\ref{toric}.)

A cone is also an (unbounded) intersection of finitely many
halfspaces in $\R^n$, so it is a convex polyhedron in the sense of
Definition~\ref{pol2}. We therefore may define face of a cone in
the same way as we did for polytopes, as intersections of $\sigma$
with supporting hyperplanes. We can only consider supporting
hyperplanes containing~$\bf 0$; such a hyperplane is defined by a
linear function $\mb u$ and will be denoted by~$\mb u^\perp$. A
\emph{face}\label{faceofcone} $\tau$ of a cone $\sigma$ is
therefore an intersection of $\sigma$ with a supporting hyperplane
$\mb u^\perp$, i.e. $\tau=\sigma\cap\mb u^\perp$. Every face of a
cone is itself a cone. If $\sigma\ne\R^n$, then $\sigma$ has the
smallest face $\sigma\cap(-\sigma)$; it is a vertex $\mathbf 0$ if
and only if $\sigma$ strongly convex. A minimal generator set of a
cone consists of nonzero vectors along its edges.

A \emph{fan}\label{deffan} is a finite collection
$\Sigma=\{\sigma_1,\ldots,\sigma_s\}$ of strongly convex cones in
some $\R^n$ such that every face of a cone in $\Sigma$ belongs to
$\Sigma$ and the intersection of any two cones in $\Sigma$ is a
face of each. A fan $\Sigma$ is \emph{rational} (respectively,
\emph{simplicial}, \emph{regular}) if every cone in $\Sigma$ is
rational (respectively, simplicial, regular). A fan
$\Sigma=\{\sigma_1,\ldots,\sigma_s\}$ is called \emph{complete} if
$\sigma_1\cup\cdots\cup\sigma_s=\R^n$.

Given a cone $\sigma\subset\R^n$, define its \emph{dual} as
\begin{equation}\label{dualcone}
  \sigma^{\mathsf{v}}=\bigl\{\mb x\in\R^n\colon \langle\mb u,\mb x\rangle\ge0
  \:\text{ for all }\:\mb u\in\sigma\}.
\end{equation}
(Note the difference with the definition of the polar set of a
polyhedron, see~\eqref{polarset}.) Observe that if $\mb
u\in\sigma^{\mathsf{v}}\ $, then $\mb u^\perp$ is a supporting
hyperplane of~$\sigma$. It can be shown that $\sigma^{\mathsf{v}}\
$ is indeed a cone, $(\sigma^{\mathsf{v}})^\mathsf{v}=\sigma$, and
$\sigma^{\mathsf{v}}\ $ is strongly convex if and only if
$\dim\sigma=n$.
\end{definition}

Cones in a fan can be separated by hyperplanes (or linear
functions):

\begin{lemma}[Separation Lemma]\label{seplemma}
Let $\sigma$ and $\sigma'$ be two cones whose intersection $\tau$
is a face of each. Then there exists $\mb
u\in\sigma^{\mathsf{v}}\,\cap(-\sigma')^\mathsf{v}$ such that
\[
  \tau=\sigma\cap\mb u^\perp=\sigma'\cap\mb u^\perp.
\]
In other words, a supporting hyperplane defining $\tau$ can be
chosen so as to separate $\sigma$ and $\sigma'$.
\end{lemma}
\begin{proof}
We only sketch a proof; the details can be found, e.g.
in~\cite[\S1.2]{fult93}. The fact that $\sigma$ and $\sigma'$
intersect in a face implies that the cone
$\xi=\sigma-\sigma'=\sigma+(-\sigma')$ is not the whole space. Let
$\mb u^\perp$ be a supporting hyperplane defining the smallest
face of~$\xi$, that is,
\[
  \xi\cap\mb u^\perp
  =\xi\cap(-\xi)=(\sigma-\sigma')\cap(\sigma'-\sigma).
\]
We claim that this $\mb u$ has the required properties. Indeed,
$\sigma\subset\xi$ implies $\mb u\in\sigma^{\mathsf{v}}\ $, and
$\tau\subset\xi\cap(-\xi)$ implies $\tau\subset\sigma\cap\mb
u^\perp$. Conversely, if $\mb x\in\sigma\cap\mb u^\perp$, then
$\mb x$ is in $\sigma'-\sigma$, so that $\mb x=\mb y'-\mb y$ for
$\mb y'\in\sigma'$, $\mb y\in\sigma$. Then $\mb x+\mb
y\in\sigma\cap\sigma'=\tau$, which implies that both $\mb x$ and
$\mb y$ are in~$\tau$. Hence $\sigma\cap\mb u^\perp=\tau$, and the
same argument for $-\mb u$ shows that $\mb
u\in(-\sigma')^\mathsf{v}$ and $\sigma'\cap\mb u^\perp=\tau$.
\end{proof}

%In Chapter~\ref{toric} and Chapter~\ref{mamanifolds}, fans will be
%used to construct algebraic varieties or topological spaces,
%either by gluing from open pieces, or by taking quotients by
%non-compact group actions.
Miraculously, the convex-geometrical separation property above
will translate into topological separation (Hausdorffness) of
algebraic varieties and topological spaces constructed from fans
in the latter chapters.

\medskip

The next construction assigns a complete fan to every convex
polytope.

\begin{construction}[Normal fan]\label{nf}
Let $P$ be a polytope~\eqref{ptope} with $m$ facets
$F_1,\ldots,F_m$ and normal vectors $\mb a_1,\ldots,\mb a_m$.
Given a face $Q\subset P$ define the cone
\begin{equation}\label{sqcone}
  \sigma_Q=\{\mb u\in\R^n\colon
  \langle\mb u,\mb x'\rangle\le\langle\mb u,\mb x\rangle
  \text{ for all $\mb x'\in Q$ and $\mb x\in P$}\}.
\end{equation}
The dual cone $\sigma_Q^\mathsf{v}$ is generated by vectors ${\mb
x-\mb x'}$ where $\mb x\in P$ and $\mb x'\in Q$. In other words,
$\sigma_Q^\mathsf{v}$ is the ``polyhedral angle'' at the face $Q$
consisting of all vectors pointing from points of $Q$ to points
of~$P$.

We say that a vector $\mb a_i$ is \emph{normal} to the face $Q$ if
$Q\subset F_i$. The cone $\sigma_Q$ is generated by those $\mb
a_i$ which are normal to~$Q$ (this is an exercise). Then
\[
  \Sigma_P=\{\sigma_Q\colon Q\text{ is a face of }P\}
\]
is a complete fan $\Sigma_P$ in $\R^n$ (this is another exercise),
which is denoted by $\Sigma_P$ and is referred to as the
\emph{normal fan} of the polytope~$P$. If $\mathbf 0$ is contained
in the interior of $P$ then $\Sigma_P$ may be also described as
the set of cones over the faces of the polar polytope~$P^*$ (yet
another exercise).

It is clear from the above descriptions that the normal fan
$\Sigma_P$ is simplicial if and only if $P$ is simple. In this
case the definition of $\Sigma_P$ may be simplified: the cones of
$\Sigma_P$ are generated by those sets $\{\mb a_{i_1},\ldots,\mb
a_{i_k}\}$ for which the intersection $F_{i_1}\cap\cdots\cap
F_{i_k}$ is nonempty.

If all vertices of $P$ are in the lattice $\Z^n$ then the normal
fan $\Sigma_P$ is rational, but the converse is not true.
Polytopes whose normal fans are regular are called
\emph{Delzant}\label{delzantptope} (this name comes from a
symplectic geometry construction discussed in
Section~\ref{symred}). Therefore, $P$ is a Delzant polytope if and
only if for every vertex $\mb x\in P$ the normal vectors to the
facets meeting at $\mb x$ can be chosen to form a basis of~$\Z^n$.
In this definition one may replace `normal vectors to the facets
meeting at~$\mb x$' by `vectors along the edges meeting at~$\mb
x$'. A Delzant polytope is necessarily simple.
\end{construction}

The normal fan $\Sigma_P$ of a polytope $P$ contains the
information about the normals to the facets (the generators $\mb
a_i$ of the edges of~$\Sigma_P$) and the combinatorial structure
of $P$ (which sets of vectors $\mb a_i$ span a cone of $\Sigma_P$
is determined by which facets intersect at a face), however the
scalars $b_i$ in~\eqref{ptope} are lost. Not any complete fan can
be obtained by `forgetting the numbers $b_i$' from a presentation
of a polytope by inequalities, i.e. not any complete fan is a
normal fan. This is fails even for regular fans, as is shown by
the next example, which is taken from~\cite{fult93}.

\begin{example}\label{nonpolytopal}
Consider the complete three-dimensional fan $\Sigma$ with 7 edges
generated by the vectors $\mb a_1=\mb e_1$, $\mb a_2=\mb e_2$,
$\mb a_3=\mb e_3$, $\mb a_4=-\mb e_1-\mb e_2-\mb e_3$, $\mb
a_5=-\mb e_1-\mb e_2$, $\mb a_6=-\mb e_2-\mb e_3$, $\mb a_7=-\mb
e_1-\mb e_3$, and 10 three-dimensional cones with vertex $\mathbf
0$ over the faces of the triangulated boundary of the tetrahedron
shown in Fig.~\ref{npctv}. It is easy to verify that $\Sigma$ is
regular.
\begin{figure}[h]
\begin{center}
\begin{picture}(70,50)
  \put(5,20){\line(2,-1){30}}
  \put(5,20){\line(1,1){30}}
  \put(35,5){\line(0,1){45}}
  \put(35,5){\line(2,1){30}}
  \put(35,5){\line(-1,2){15}}
  \put(35,25){\line(3,2){15}}
  \put(35,25){\line(-3,2){15}}
  \put(65,20){\line(-1,1){30}}
  \qbezier(35,25)(65,20)(65,20)
  \multiput(5,20)(5,0){12}{\line(1,0){3}}
  \multiput(20,35)(5.1,0){6}{\line(1,0){3}}
  \multiput(5,20)(5.1,1.67){9}{\line(3,1){3.6}}
  \put(66,19){1}
  \put(2,19){2}
  \put(34,1){3}
  \put(34,51){4}
  \put(32,22){5}
  \put(52,34){6}
  \put(16,34){7}
\end{picture}%
\caption{Complete regular fan not coming from a polytope.}
\label{npctv}
\end{center}
\end{figure}

Assume that $\Sigma=\Sigma_P$ is the normal fan of a polytope~$P$.
Consider the function $\psi\colon\R^n\to\R$ given by
\[
  \psi(\mb u)=\min_{\mb x\in P}\langle\mb u,\mb x\rangle=
  \min_{\mb v\in V(P)}\langle\mb u,\mb v\rangle,
\]
where $V(P)$ is the set of vertices of~$P$. This function is
continuous and its restriction to every 3-dimensional cone
of~$\Sigma_P$ is linear. Indeed, 3-dimensional cones $\sigma_{\mb
v}$ correspond to vertices $\mb v\in V(P)$, and we have $\psi(\mb
u)=\langle\mb u,\mb v\rangle$ for $\mb u\in\sigma_{\mb v}$ by
definition~\eqref{sqcone} of~$\sigma_{\mb v}$.

Now consider the two 3-dimensional cones of $\Sigma_P$ generated
by the triples $\mb a_1,\mb a_3,\mb a_5$ and $\mb a_1,\mb a_5,\mb
a_6$, and let $\mb v$ and $\mb v'$ be the corresponding vertices
of~$P$. Then $\psi(\mb a_1)=\langle\mb a_1,\mb v\rangle$,
$\psi(\mb a_3)=\langle\mb a_3,\mb v\rangle$, $\psi(\mb
a_5)=\langle\mb a_5,\mb v\rangle$, $\psi(\mb a_6)=\langle\mb
a_6,\mb v'\rangle$, hence,
\[
  \psi(\mb a_1)+\psi(\mb a_5)-\psi(\mb a_3)=\langle\mb a_1+\mb a_5
  -\mb a_3,\mb v\rangle=
  \langle\mb a_6,\mb v\rangle>\langle\mb a_6,\mb v'\rangle=\psi(\mb
  a_6).
\]
Therefore,
\[
  \psi(\mb a_1)+\psi(\mb a_5)>\psi(\mb a_3)+\psi(\mb a_6).
\]
Similarly,
\begin{gather*}
  \psi(\mb a_2)+\psi(\mb a_6)>\psi(\mb a_1)+\psi(\mb a_7),\\
  \psi(\mb a_3)+\psi(\mb a_7)>\psi(\mb a_2)+\psi(\mb a_5).
\end{gather*}
Adding the last three inequalities together we get a
contradiction.
\end{example}

\subsection*{Exercises.}
\begin{exercise}
Let $\sigma$ be a cone in~$\R^n$. Show that $\sigma^{\mathsf{v}}\
$ is also a cone, $(\sigma^{\mathsf{v}})^\mathsf{v}=\sigma$, and
$\sigma^{\mathsf{v}}\ $ is strongly convex if and only if
$\dim\sigma=n$.
\end{exercise}

\begin{exercise}
The cone $\sigma_Q$ given by~\eqref{sqcone} is generated by those
vectors among $\mb a_1,\ldots,\mb a_m$ which are normal to~$Q$.
\end{exercise}

\begin{exercise}
The set $\{\sigma_Q\colon Q\text{ is a face of }P\}$ is a complete
fan in~$\R^n$.
\end{exercise}

\begin{exercise}
If $\mathbf 0$ is contained in the interior of $P$ then $\Sigma_P$
consists of cones over the faces of the polar polytope~$P^*$.
\end{exercise}

\begin{exercise}
Let $P$ be a convex polytope (not necessarily simple). Does the
collection of cones generated by the sets $\{\mb
a_{i_1},\ldots,\mb a_{i_k}\}$ of normal vectors for which
$F_{i_1}\cap\cdots\cap F_{i_k}\ne\varnothing$ form a fan?
\end{exercise}

\section{Simplicial complexes}
A simplex is the convex hull of a set of affinely independent
points in~$\R^n$.

\begin{definition}
\label{polyhed} A \emph{geometric simplicial complex} in $\R^n$ is
a collection $\mathcal P$ of simplices of arbitrary dimension such
that every face of a simplex in $\mathcal P$ belongs to $\mathcal
P$ and the intersection of any two simplices in $\mathcal P$ is
either empty or a face of each.

To make the exposition more streamlined and without creating much
ambiguity we shall not distinguish between the collection
$\mathcal P$ (which is an abstract set of simplices) and the union
of its simplices (which is a subset in $\R^n$). In $PL$ (piecewise
linear) topology the latter union is usually referred to as `the
\emph{polyhedron} of $\mathcal P$'. Although we have already
reserved the term `polyhedron' for a finite intersection of
halfspaces~\eqref{ptope}, we shall also occasionally use it in the
$PL$ topological sense (when it creates no ambiguity).

A \emph{face}\label{faceofsc} of $\mathcal P$ is a face of any of
its simplices. The dimension of $\mathcal P$ is the maximal
dimension of its faces.
\end{definition}

If we know the set of vertices of $\mathcal P$ in $\R^n$ then we
may recover the whole $\mathcal P$ by specifying which subsets of
vertices span simplices. This observation leads to the following
definition.

\begin{definition}\label{absimcom}
An \emph{abstract simplicial complex} on a set $\mathcal V$ is a
collection $\sK$ of subsets $I\subset\mathcal V$ such that if
$I\in\sK$ then any subset of $I$ also belongs to~$\sK$. We always
assume that the empty set $\varnothing$ belongs to~$\sK$. We refer
to $ I\in\sK$ as an (abstract) \emph{simplex} of~$\sK$.

One-element simplices are called \emph{vertices} of~$\sK$. If
$\sK$ contains all one-element subsets of~$\mathcal V$, then we
say that $\sK$ is a simplicial complex \emph{on the vertex
set~$\mathcal V$}.

It is sometimes convenient to consider simplicial complexes $\sK$
whose vertex sets are proper subsets of~$\mathcal V$. In this case
we refer to a one-element subset of $\mathcal V$ which is not a
vertex of~$\sK$ as a \emph{ghost vertex}\label{ghostvertex}.

The \emph{dimension} of a simplex $I\in\sK$ is $\dim I=| I|-1$,
where $|I|$ denotes the number of elements in~$I$. The dimension
of $\sK$ is the maximal dimension of its simplices. A simplicial
complex $\sK$ is \emph{pure} if all its maximal simplices have the
same dimension. A subcollection $\sK'\subset\sK$ which is itself a
simplicial complex is called a
\emph{subcomplex}\label{defsubcomplex} of~$\sK$.

A geometric simplicial complex $\mathcal P$ is said to be a
\emph{geometric realisation} of an abstract simplicial complex
$\sK$ on a set $\mathcal V$ if there is a bijection between
$\mathcal V$ and the vertex set of $\mathcal P$ which maps
abstract simplices of $\sK$ to vertex sets of faces of~$\mathcal
P$.
\end{definition}

Both geometric and abstract simplicial complexes will be assumed
to be \emph{finite}, unless we explicitly specify otherwise. In
most constructions we identify the set $\mathcal V$ with the index
set $[m]=\{1,\ldots,m\}$ and consider abstract simplicial
complexes on~$[m]$. An identification of $\mathcal V$ with $[m]$
fixes an order of vertices, although this order will be irrelevant
in most cases. We drop the parentheses in the notation of
one-element subsets $\{i\}\subset[m]$, so that $i\in\sK$ means
that $\{i\}$ is a vertex of~$\sK$, and $i\notin\sK$ means that
$\{i\}$ is a ghost vertex.

We shall use the common notation $\varDelta^{m-1}$ for any of the
following three objects: an $(m-1)$-simplex (a convex polytope),
the geometric simplicial complex consisting of all faces in an
$(m-1)$-simplex, and the abstract simplicial complex consisting of
all subsets of~$[m]$.

\begin{construction}
Every abstract simplicial complex $\sK$ on~$[m]$ can be realised
geometrically in~$\R^m$ as follows. Let $\mb e_1,\ldots,\mb e_m$
be the standard basis in~$\R^m$, and for each $I\subset[m]$ denote
by $\varDelta^I$ the convex hull of points $\mb e_i$ with $i\in
I$.
%, then $\varDelta^I$ is a (geometric) simplex.
Then
\[
  \bigcup_{I\in\sK}\varDelta^I\subset\R^m
\]
is a geometric realisation of~$\sK$.
\end{construction}

The above construction is just a geometrical interpretation of the
fact that any simplicial complex on~$[m]$ is a subcomplex of the
simplex~$\varDelta^{m-1}$. Also, by a classical
result~\cite{pont52}, a $d$-dimensional abstract simplicial
complex admits a geometric realisation in $\R^{2d+1}$.

\begin{example}
\label{polsph} The boundary of a simplicial $n$-polytope is a
simplicial complex of dimension~$n-1$. For a simple polytope $P$,
we shall denote by~$\sK_P$ the boundary complex $\partial P^*$ of
the dual polytope. It coincides with the nerve of the covering of
$\partial P$ by the facets. That is, the vertices of $\sK_P$ are
the facets of~$P$, and a set of vertices spans a simplex whenever
the intersection of the corresponding facets is nonempty. We refer
to $\sK_P$ as the \emph{nerve complex} of~$P$.
\end{example}

\begin{definition}\label{simfve}
The \emph{$f$-vector} of an $(n-1)$-dimensional simplicial complex
$\sK$ is $\mb f(\sK)=(f_0,f_1,\ldots,f_{n-1})$, where $f_i$ is the
number of $i$-dimensional simplices in~$\sK$. We also set
$f_{-1}=1$; to justify this convention one can assign dimension
$-1$ to the empty simplex. The \emph{$h$-vector} $\mb
h(\sK)=(h_0,h_1,\ldots,h_n)$ is defined by the identity
\begin{equation}\label{hvectors}
  h_0s^n+h_1s^{n-1}+\cdots+h_n=
  (s-1)^n+f_0(s-1)^{n-1}+\cdots+f_{n-1}.
\end{equation}
(Warning: this is not the identity obtained by substituting $t=1$
in~\eqref{hvector}.) The \emph{$g$-vector} $\mb
g(\sK)=(g_0,g_1,\ldots,g_{[n/2]})$ is defined by $g_0=1$ and
$g_i=h_i-h_{i-1}$ for $i=1,\ldots,[n/2]$.
\end{definition}

\begin{remark}
If $\sK=\sK_P$ is the nerve complex of a simple polytope~$P$, then
$\mb f(\sK)=\mb f(P^*)$ and $\mb h(\sK)=\mb h(P)$. This our
notational convention may look artificial, but it seems to be the
best possible way to treat the face vectors of both polytopes and
simplicial complexes consistently.
\end{remark}

\begin{definition}\label{simmap}
Let $\sK_1$, $\sK_2$ be simplicial complexes on the sets $[m_1]$,
$[m_2]$ respectively, and $\mathcal P_1$, $\mathcal P_2$ their
geometric realisations. A map $\phi\colon[m_1]\to[m_2]$ induces a
\emph{simplicial map} between $\sK_1$ and $\sK_2$ if
$\phi(I)\in\sK_2$ for any $I\in\sK_1$. A simplicial map $\phi$ is
said to be \emph{nondegenerate} if $|\phi(I)|=|I|$ for any
$I\in\sK_1$. On the geometric level, a simplicial map extends
linearly on the faces of $\mathcal P_1$ to a map $\mathcal
P_1\to\mathcal P_2$, which we continue to denote by~$\varphi$. A
\emph{simplicial isomorphism} is a simplicial map for which there
exists a simplicial inverse.

There is an obvious isomorphism between any two geometric
realisations of an abstract simplicial complex~$\sK$. We therefore
shall use the common notation $|\sK|$ for any geometric
realisation of~$\sK$. Whenever it is safe, we shall not
distinguish between abstract simplicial complexes and their
geometric realisations. For example, we shall say `simplicial
complex $\sK$ is homeomorphic to~$X$' instead of `the geometric
realisation of $\sK$ is homeomorphic to~$X$'.

We shall refer to a simplicial complex homeomorphic to a
topological space $X$ as a \emph{triangulation of~$X$}, or a
\emph{simplicial subdivision of~$X$}\label{simpsubdi}.

A \emph{subdivision} of a geometric simplicial complex $\mathcal
P$ is a geometric simplicial complex $\mathcal P'$ such that each
simplex of $\mathcal P'$ is contained in a simplex of $\mathcal P$
and each simplex of $\mathcal P$ is a union of finitely many
simplices of~$\mathcal P'$. A $PL$~\emph{map}\label{defplmap}
$\phi\colon\mathcal P_1\to\mathcal P_2$ is a simplicial map from a
subdivision of $\mathcal P_1$ to a subdivision of~$\mathcal P_2$.
A $PL$ \emph{homeomorphism} is a $PL$ map for which there exists a
$PL$ inverse.
%Two $PL$ homeomorphic simplicial complexes are also sometimes
%called \emph{combinatorially equivalent}.
In other words, two
simplicial complexes are $PL$ homeomorphic if there is a
simplicial complex isomorphic to a subdivision of each of them.
\end{definition}

\begin{remark}
%We used a definition of $PL$ maps an homeomorphisms that fits the
%combinatorial flavour of this chapter.
For a topological approach to $PL$ maps (where a $PL$ map is
defined between spaces rather than their triangulations) we refer
to standard sources on $PL$ topology, such as~\cite{huds69}
and~\cite{ro-sa72}.
\end{remark}

\begin{example}\label{simincl}\

1. If $P$ is a simple $n$-polytope then the nerve complex $\sK_P$
(see Example~\ref{polsph}) is a triangulation of an
$(n-1)$-dimensional sphere~$S^{n-1}$.

2. Let $\Sigma$ be a simplicial fan in~$\R^n$ with $m$ edges
generated by vectors $\mb a_1,\ldots,\mb a_m$. Its
\emph{underlying simplicial complex}\label{undecompl} is defined
by
\[
  \sK_\Sigma=\bigl\{\{i_1,\ldots,i_k\}\subset[m]\colon
  \mb a_{i_1},\ldots,\mb a_{i_k}\text{ span a cone of }\Sigma
  \bigr\}.
\]
Informally, $\sK_\Sigma$ may be viewed as the intersection of
$\Sigma$ with a unit sphere. The fan $\Sigma$ is complete if and
only if $\sK_\Sigma$ is a triangulation of~$S^{n-1}$. If $\Sigma$
is a normal fan of a simple $n$-polytope~$P$, then
$\sK_\Sigma=\sK_P$.
\end{example}

\begin{construction}[join]\label{join}
Let $\sK_1$ and $\sK_2$ be simplicial complexes on sets $\mathcal
V_1$ and $\mathcal V_2$ respectively. The \emph{join} of $\sK_1$
and $\sK_2$ is the simplicial complex
\[
  \sK_1*\sK_2=\bigl\{  I\subset\mathcal V_1\sqcup\mathcal V_2\colon
  I= I_1\cup I_2,\; I_1\in\sK_1,  I_2\in\sK_2\bigr\}
\]
on the set $\mathcal V_1\sqcup\mathcal V_2$. The join operation is
associative by inspection.
\end{construction}

\begin{example}\label{consusp}\

1. If $\sK_1=\varDelta^{m_1-1}$, $\sK_2=\varDelta^{m_2-1}$, then
$\sK_1*\sK_2=\varDelta^{m_1+m_2-1}$.

2. The simplicial complex $\varDelta^0*\sK$ (the join of $\sK$ and
a point) is called the \emph{cone} over $\sK$ and denoted
$\cone\sK$.

3. Let $S^0$ be a pair of disjoint points (a $0$-sphere). Then
$S^0*\sK$ is called the \emph{suspension} of $\sK$ and denoted
$\Sigma\sK$. Geometric realisations of $\cone\sK$ and $\Sigma\sK$
are the topological cone and suspension over~$|\sK|$ respectively.

4. Let $P_1$ and $P_2$ be simple polytopes. Then
\[
  \mathcal K_{P_1\times P_2}=
  \sK_{P_1}*\sK_{P_2}.
\]
(see Construction~\ref{prodsp}).
\end{example}

\begin{construction}\label{simprod}
The fact that the product of two simplices is not a simplex makes
triangulations of products of spaces more subtle. There is the
following canonical way to triangulate the product of two
simplicial complexes whose vertices are ordered. Let $\sK_1$ and
$\sK_2$ be simplicial complexes on $[m_1]$ and $[m_2]$
respectively (this is one of the few constructions where the order
of vertices is important: here it is an additional structure). We
construct a simplicial complex $\sK_1\mathbin{\tilde\times}\sK_2$
on $[m_1]\times[m_2]$ as follows. By definition, simplices of
$\sK_1\mathbin{\tilde\times}\sK_2$ are those subsets in products
$I_1\times I_2$ of $I_1\in\sK_1$ and $I_2\in\sK_2$ which establish
non-decreasing relations between $I_1$ and~$I_2$. More formally,
\begin{multline*}
  \sK_1\mathbin{\tilde\times}\sK_2=\bigl\{ I\subset I_1\times I_2\colon
  I_1\in\sK_1,\:I_2\in\sK_2,\\
  \text{and }i\le i'\text{ implies }j\le j'
  \text{ for any two pairs }(i,j),(i',j')\in I \bigr\}.
\end{multline*}
We leave it as an exercise to check that
$|\sK_1\mathbin{\tilde\times}\sK_2|$ defines a triangulation of
$|\sK_1|\times|\sK_2|$. Note that
$\sK_1\mathbin{\tilde\times}\sK_2\ne
\sK_2\mathbin{\tilde\times}\sK_1$ in general. Note also that if
$\sK_1=\sK_2=\sK$, then the diagonal is naturally a subcomplex in
the triangulation of $|\sK|\times|\sK|$.
\end{construction}

\begin{construction}[connected sum of simplicial complexes]\label{simcs}
Let $\sK_1$ and $\sK_2$ be two pure $d$-dimensional simplicial
complexes on sets $\mathcal V_1$ and $\mathcal V_2$ respectively,
where $|\mathcal V_1|=m_1$, $|\mathcal V_2|=m_2$. Choose two
maximal simplices $ I_1\in\sK_1$, $ I_2\in\sK_2$ and fix an
identification of $I_1$ with $I_2$. Let $\mathcal V_1\cup_{
I}\mathcal V_2$ be the union of $\mathcal V_1$ and $\mathcal V_2$
in which $ I_1$ is identified with~$ I_2$, and denote by $I$ the
subset created by the identification. We have $|\mathcal V_1\cup_{
I}\mathcal V_2|=m_1+m_2-d-1$. Both $\sK_1$ and $\sK_2$ now can be
viewed as simplicial complexes on the set $\mathcal V_1\cup_{
I}\mathcal V_2$. We define the \emph{connected sum} of $\sK_1$ and
$\sK_2$ at $I_1$ and $I_2$ as the simplicial complex
\[
  \sK_1\cs_{ I_1, I_2}\sK_2=(\sK_1\cup\sK_2)\setminus\{ I\}
\]
on the set $\mathcal V_1\cup_ I\mathcal V_2$. When the choices are
clear, or their effect on the result irrelevant, we use the
abbreviation $K_1\cs K_2$. Geometrically, the connected sum of
$|\sK_1|$ and $|\sK_2|$ is produced by attaching $|\sK_1|$ to
$|\sK_2|$ along $I_1$ and $I_2$ and then removing the face $I$
obtained by the identification.
\end{construction}

\begin{example}
Connected sum of simple polytopes defined in
Construction~\ref{consum} is dual to the operation described
above. Namely, if $P$ and $Q$ are two simple $n$-polytopes, then
\[
  \sK_P\cs\sK_Q=\sK_{P\cs Q}.
\]
\end{example}

\begin{definition}\label{deflink}
Let $\sK$ be a simplicial complex on a set $\mathcal V$. The
\emph{link} and the \emph{star} of a simplex $I\in\sK$ are the
subcomplexes
\begin{align*}
  \lk_\sK I&=\bigl\{J\in \sK\colon I\cup J\in \sK,\;
  I\cap J=\varnothing\bigr\};\\
  \st_\sK I&=\bigl\{J\in \sK\colon I\cup J\in \sK \bigr\}.
\end{align*}
We also define the subcomplex
\[
  \partial\st_\sK I=\bigl\{J\in \sK\colon I\cup J\in \sK,\;
  I\not\subset J\bigr\}.
\]
Then we have a sequence of inclusions
\[
  \lk_\sK I\subset\partial\st_\sK I\subset\st_\sK I.
\]
For any vertex $v\in\sK$, the subcomplex $\st_\sK v$ is the cone
over $\lk_\sK v=\partial\st_\sK v$. Also, $|\st_\sK v|$ is the
union of all faces of $|\sK|$ that contain~$v$. We omit the
subscript $\sK$ in the notation of link and star whenever the
ambient simplicial complex is clear.
\end{definition}

The links of simplices determine the topological structure of the
space $|\sK|$ near any of its points. In particular, the following
proposition describes the `local cohomology' of~$|\sK|$.

\begin{proposition}\label{lclin}
Let $x$ be an interior point of a simplex $I\in \sK$. Then
\[
  H^i\bigl(|\sK|,|\sK|\setminus x\bigr)\cong\widetilde{H}^{i-|I|}(\lk
  I),
\]
where $H^i(X,A)$ denotes the $i$th relative singular cohomology
group of a pair $A\subset X$, and $\widetilde{H}^i(\sK)$ denotes
the $i$th reduced simplicial cohomology group of~$\sK$.
\end{proposition}
\begin{proof}
We have
\begin{align*}
  H^i\bigl(|\sK|,|\sK|\setminus x\bigr)&\cong
  H^i\bigl(\st I,(\st I)\setminus x\bigr)
  \cong H^i\bigl(\st I,(\partial I)*(\lk I)\bigr)\cong\\
  &\cong\widetilde{H}^{i-1}\bigl((\partial I)*(\lk I)\bigr)
  \cong\widetilde{H}^{i-|I|}(\lk I).
\end{align*}
Here the first isomorphism follows from the excision property, the
second uses the fact that $(\partial I)*(\lk I)$ is a deformation
retract of $(\st I)\setminus x$, the third follows from the
homology sequence of pair, and the fourth is by the suspension
isomorphism.
\end{proof}

Given a subcomplex $\mathcal L\subset\sK$, define the \emph{closed
combinatorial neighbourhood}\label{combnei} of $\mathcal L$ in
$\sK$ by
\[
  U_\sK(\sL)=\bigcup_{I\in\sL}\st_\sK I.
\]
That is, $U_\sK(\sL)$ consists of all simplices of $\sK$, together
with all their faces, having some simplex of $\sL$ as a face.
Define also the \emph{open combinatorial neighbourhood}
${\mathop{\mathit U}\limits^\circ}_\sK(\sL)$ of $|\sL|$ in $|\sK|$
as the union of relative interiors of simplices of $|\sK|$ having
some simplex of $|\sL|$ as a face.

\begin{definition}\label{defcore}
Given a subset $I\subset\mathcal V$, define the corresponding
\emph{full subcomplex} of $\sK$ (or the \emph{restriction} of
$\sK$ to~$I$) as
\begin{equation}\label{ksigma}
  \sK_I=\bigl\{ J\in\sK\colon J\subset I \bigr\}.
\end{equation}
Set $\core\mathcal V=\{v\in\mathcal V\colon\st v\ne\sK\}$. The
\emph{core}\label{core} of $\sK$ is the subcomplex $\core
\sK=\sK_{\core\mathcal V}$. Then we may write
$\sK=\core(\sK)*\varDelta^{s-1}$, where $\varDelta^{s-1}$ is the
simplex on the set $\mathcal V\setminus\core\mathcal V$.
\end{definition}

\begin{example}\

1. $\lk_\sK\varnothing=\sK$.

2. Let $K=\partial\varDelta^3$ be the boundary of the tetrahedron
on four vertices $1,2,3,4$, and $I=\{1,2\}$. Then $\lk I$ consists
of two disjoint points 3 and~4.

3. Let $\sK$ be the cone over $\sL$ with vertex~$v$. Then $\lk
v=\sL$, $\st v=\sK$, and $\core\sK=\core\sL$.
\end{example}

\subsection*{Exercises.}
\begin{exercise}
Assume that $\sK_1$ is realised geometrically in $\R^{n_1}$ and
$\sK_2$ in $\R^{n_2}$. Construct a realisation of the join
$\sK_1*\sK_2$ in $\R^{n_1+n_2+1}$.
\end{exercise}

\begin{exercise}
Show that $|\sK_1\mathbin{\tilde\times}\sK_2|$ is a triangulation
of $|\sK_1|\times|\sK_2|$.
\end{exercise}

\begin{exercise}\label{purelink}
Let $\sK$ be a pure simplicial complex. Then $\lk I$ is pure of
dimension $\dim\sK-|I|$ for any $I\in\sK$.
\end{exercise}

\section{Barycentric subdivision and flag complexes}\label{secflc}
\begin{definition}\label{barsub}
The \emph{barycentric subdivision} of an abstract simplicial
complex $\sK$ is the simplicial complex $\sK'$ defined as follows.
The vertex set of $\sK'$ is the set $\{I\in
\sK,\;I\ne\varnothing\}$ of nonempty simplices of~$\sK$. Simplices
of $\sK'$ are chains of embedded simplices of~$\sK$. That is,
$\{I_1,\ldots,I_r\}\in\sK'$ if and only if $I_1\subset
I_2\subset\cdots\subset I_r$ in $\sK$ (after possible reordering)
and $I_1\ne\varnothing$.

The \emph{barycentre}\label{defbaryce} of a geometric simplex in
$\R^d$ with vertices $\mb v_1,\ldots,\mb v_{d+1}$ is the point
$\frac 1{d+1}(\mb v_1+\cdots+\mb v_{d+1})$. A geometric
realisation of $\sK'$ may be obtained by mapping every vertex of
$\sK'$ to the barycentre of the corresponding simplex of~$|\sK|$;
simplices of $|\sK'|$ are therefore spanned by the sets of
barycentres of chains of embedded simplices of $|\sK|$.
\end{definition}

\begin{example}\label{ndbar}
For any $(n-1)$-dimensional simplicial complex $\sK$ on~$[m]$,
there is a nondegenerate simplicial map $\sK'\to\varDelta^{n-1}$
defined on the vertices by $I\mapsto|I|$ for $I\in\sK$, where
$|I|$ denotes the cardinality of~$I$. Here $I$ is viewed as a
vertex of $\sK'$ and $|I|$ as a vertex of $\varDelta^{n-1}$.
\end{example}

\begin{example}
\label{nabla} Let $\sK$ be a simplicial complex on a set~$\mathcal
V$, and assume we are given a choice function
$c\colon\sK\to\mathcal V$ assigning to each simplex $I\in\sK$ one
of its vertices. For instance, if $\mathcal V=[m]$ we can define
$c(I)$ as the minimal element of~$I$. For every such $f$ there is
a canonical simplicial map $\varphi_c\colon\sK'\to\sK$ constructed
as follows. We define $\varphi_c$ on the vertices of $\sK'$ by
$\varphi_c(I)=c(I)$ and then extend it to simplices of $\sK'$ by
the formula
\[
  \varphi_c(I_1\subset I_2\subset\cdots\subset I_r)=
  \{c(I_1),c(I_2),\ldots,c(I_r)\}.
\]
The right hand side is a subset of $I_r$, and therefore it is a
simplex of~$\sK$.
\end{example}

We shall need explicit formulae for the transformation of the $f$-
and $h$-vectors under the barycentric subdivision. Introduce the
matrix
\[
  B = (b_{ij}),\quad 0\le i,j\le n-1;\quad
  b_{ij} = \sum_{k=0}^i(-1)^k\binom {i+1}k(i-k+1)^{j+1}.
\]
It can be shown that $b_{ij}=0$ for $i>j$ and $b_{ii}=(i+1)!$, so
that $B$ is a nonsingular upper triangular matrix.

\begin{lemma}\label{barfv}
Let $\sK'$ be the barycentric subdivision of an
$(n-1)$-dimensional simplicial complex~$\sK$. Then the $f$-vectors
of $\sK$ and $\sK'$ are related by the identity
\[
  f_i(\sK')= \sum_{j=i}^{n-1}b_{ij}f_j(\sK),\quad 0\le i\le n-1,
\]
that is, $\mb f\:^t(\sK') = B\mb f\:^t(\sK)$ where $\mb
f\:^t(\sK)$ is the column vector with entries~$f_i(\sK)$.
\end{lemma}
\begin{proof}
Consider the barycentric subdivision of a $j$-simplex
$\varDelta^j$, and let $b'_{ij}$ be the number of its
$i$-simplices which lie inside $\varDelta^j$, i.e. not contained
in~$\partial\varDelta^j$. Then we have
$f_i(\sK')=\sum_{j=i}^{n-1}b'_{ij}f_j(\sK)$. It remains to show
that $b_{ij}=b'_{ij}$. Indeed, it is easy to see that the numbers
$b'_{ij}$ satisfy the following recurrence relation:
\[
  b'_{ij}=(j+1)b'_{i-1,j-1}+\bin{j+1}2 b'_{i-1,j-2}+\cdots+
  \bin{j+1}{j-i+1}b'_{i-1,i-1}.
\]
It follows by induction that $b'_{ij}$ is given by the same
formula as~$b_{ij}$.
\end{proof}

Now introduce the matrix
\[
  D = (d_{pq}),\quad 0\le p,q\le n;\quad
  d_{pq} = \sum_{k=0}^p(-1)^k \binom{n+1}k (p-k)^q (p-k+1)^{n-q},
\]
where we set $0^0=1$.

\begin{lemma}\label{barhv}
The $h$-vectors of $\sK$ and $\sK'$ are related by the identity:
\[
  h_p(\sK') = \sum_{q=0}^n d_{pq}h_q(\sK),\quad 0\le p\le n,
\]
that is, $\mb h^t(\sK') = D\mb h^t(\sK)$. Moreover, the matrix $D$
is nonsingular.
\end{lemma}
\begin{proof}
The formula for $h_p(\sK')$ is established by a routine check
using Lemma~\ref{barfv}, relations~\eqref{hvectors} and identities
for the binomial coefficients. If we add the component $f_{-1}=1$
to the $f$-vector and change the matrix $B$ appropriately, then we
obtain $D=C^{-1}BC$, where $C$ is the transition matrix from the
$h$-vector to the $f$-vector (its explicit form can be obtained
easily from relations~\eqref{hvectors}). This implies the
nonsingularity of~$D$.
\end{proof}

%Таким образом, переход к барицентрическому подразбиению индуцирует
%линейные невырожденные преобразования $B$ и $D$ на $f$- и
%$h$-векторах симплициальных комплексов.

\begin{definition}\label{ordercom}
Let $\mathcal P$ be a poset (partially ordered set) with strict
order relation~$<$. Its \emph{order complex} $\ord(\mathcal P)$ is
the collection of all totally ordered chains $x_1< x_2<\cdots<
x_k$ (or \emph{flags}), \ $x_i\in\mathcal P$. Clearly,
$\ord(\mathcal P)$ is a simplicial complex.
\end{definition}

The following proposition is clear from the definition.

\begin{proposition}\label{barorder}
Let $\sK$ be a simplicial complex, viewed as the poset of its
simplices with respect to inclusion. Then $\ord(\mathcal
K\setminus\varnothing)$ is the barycentric subdivision~$\sK'$. The
order complex of $\sK$ (with the empty simplex included) is
$\cone(\sK')$.
\end{proposition}

This observation may be used to define the barycentric subdivision
of other combinatorial objects. For example, let $Q$ be a convex
polytope, and $\mathcal Q$ the poset of its proper faces. Then
$\ord(\mathcal Q)$ is a simplicial complex; moreover, it is the
boundary complex of a simplicial polytope~$Q'$ (an exercise),
called the \emph{barycentric subdivision}\label{bcptope} of~$Q$.
The vertices of $Q'$ correspond to the barycentres of proper faces
of~$Q$.

\begin{proposition}\label{stardual}
Let $P$ be a simple polytope and let $\sK=\sK_P$ be its nerve
complex. Given a facet $F\subset P$, let $v$ be the corresponding
vertex of~$\sK$. Then $\st_{\sK'}v$ is a triangulation of~$F$.
\end{proposition}
\begin{proof}
We identify $\partial P$ with
%(the geometric realisation of)
$\sK'$ by mapping the barycentre of each proper face of~$P$ to the
corresponding vertex of~$\sK'$. Under this identification, $F$ is
mapped to the union of simplices of $\sK'$ corresponding to chains
$G_1\subset\cdots\subset G_k=F$ of faces of $P$ ending at~$F$.
This union is exactly the star of $v$ in~$\sK'$.
\end{proof}

\begin{definition}\label{defflagc}
A simplicial complex $\sK$ is called a \emph{flag complex} if any
set of vertices of $\sK$ which are pairwise connected by edges
spans a simplex.
\end{definition}

A flag complex is therefore completely determined by its
1-skeleton, which is a simple graph. Given such a graph $\Gamma$,
we may reconstruct the flag complex~$\sK_\Gamma$, whose simplices
are the vertex sets of complete subgraphs (or
\emph{cliques})\label{defcliq} of~$\Gamma$.

Flag complexes may be characterised in terms of their missing
faces. A \emph{missing face} of $\sK$ is a subset $I\subset[m]$
such that $I\notin\sK$, but every proper subset of $I$ is a
simplex of~$\sK$. Then $\sK$ is flag if and only if each of its
missing faces has two vertices.

Therefore, $P$ is a flag simple $n$-polytope if and only if any
set of its facets whose pairwise intersections are nonempty
intersects at a face.
%$F_{i_1},\ldots,F_{i_k}$ whose pairwise intersections are nonempty
%intersects at an $(n-k)$-face.

\begin{example}\

1. The order complex of a poset is a flag complex.

2. A simple polytope $P$ is a flag polytope (see
Definition~\ref{defflp}) if and only if its nerve complex $\sK_P$
is a flag complex.

3. The boundary of a 5-gon is a flag complex, but it is not the
order complex of a poset.

4. The boundary of a $d$-simplex is not flag for $d>1$.

5. The join $\sK_1*\sK_2$ of two flag complexes (see
Construction~\ref{join}) is flag (an exercise). Therefore, the
product $P\times Q$ of two flag simple polytopes is flag.

6. The connected sum $\sK_1\cs\sK_2$ of two $d$-dimensional
complexes (see Construction~\ref{simcs}) is not flag if $d>1$.
\end{example}

Since the $h$-vector of a sphere triangulation is symmetric, the
\emph{$\gamma$-vector}
$\gamma(\sK)=(\gamma_0,\gamma_1,\ldots,\gamma_{[n/2]})$ of an
$(n-1)$-dimensional sphere triangulation $\sK$ can be defined by
the equation
\begin{equation}\label{Hgammas}
  \sum_{i=0}^nh_is^it^{n-i}=\sum_{i=0}^{[n/2]}\gamma_i(s+t)^{n-2i}(st)^i.
\end{equation}

\begin{conjecture}[Gal~\cite{gal05}]\label{galconj} Let $\sK$ be a
flag triangulation of an $(n-1)$-sphere. Then $\gamma_i(\sK)\ge0$
for $i=0,\ldots,\sbr n2$.
\end{conjecture}

Substituting $n=2q$, $s=1$ and $t=-1$ into~\eqref{Hgammas} we
obtain
\[
  \sum_{i=0}^{2q}(-1)^ih_i = (-1)^q\gamma_q\,.
\]
The top inequality $\gamma_q\ge0$ from the Gal conjecture
therefore implies the following:

\begin{conjecture}[Charney--Davis~\cite{ch-da95}]\label{cdconj}
The inequality
\[
  (-1)^q\sum_{i=0}^{2q}(-1)^ih_i(\sK)\ge0
\]
holds for any flag triangulation $\sK$ of a $(2q-1)$-sphere.
\end{conjecture}

The Charney--Davis conjecture is a discrete analogue of the
well-known \emph{Hopf conjecture}\label{Hopfconj} of differential
geometry, which states that the Euler characteristic $\chi$ of a
closed nonpositively curved Riemannian manifold $M^{2q}$ satisfies
the inequality $(-1)^q\chi\ge0$. Due to a theorem of
Gromov~\cite{grom87}, a piecewise Euclidean cubical complex
satisfies the discrete analogue of the nonpositive curvature
condition, the so-called CAT(0)~\emph{inequality}, if and only if
the links of all vertices are flag complexes. As was observed
in~\cite{ch-da95}, the Hopf conjecture for piecewise Euclidean
cubical manifolds translates into Conjecture~\ref{cdconj}.

The Hopf and Charney--Davis conjectures are valid for $q=1,2$.
In~\cite{gal05} the Gal conjecture was verified for $n\le5$.

\subsection*{Exercises.}
\begin{exercise}
Show that the matrix $B$ of Lemma~\ref{barfv} satisfies $b_{ij}=0$
for $i>j$ and $b_{ii}=(i+1)!$.
\end{exercise}

\begin{exercise}
Prove the formula for $h_p(\sK')$ of Lemma~\ref{barhv}.
\end{exercise}

\begin{exercise}
Show that the order complex of the poset of proper faces of a
polytope is the boundary complex of a simplicial polytope. (Hint:
use stellar subdivisions, see Section~\ref{stbist} and the
exercises there.)
\end{exercise}

\begin{exercise}
The join of two flag complexes is flag.
\end{exercise}

\begin{exercise}
Links of simplices in a flag complex are flag.
\end{exercise}

\section{Alexander duality}\label{alexd}
For any simplicial complex, a dual complex may be defined on the
same set. This duality has many important combinatorial and
topological consequences.

\begin{definition}[dual complex]
\label{dual} Let $\sK$ be a simplicial complex on~$[m]$ and
$\sK\ne\varDelta^{m-1}$. Define
\[
  \widehat{\sK}=\bigl\{ I\subset[m]\colon\quad
  [m]\setminus I\notin\sK \bigr\}.
\]
Then $\widehat{\sK}$ is also a simplicial complex on~$[m]$, which
we refer to as the \emph{Alexander dual} of~$\sK$. Obviously, the
dual of $\widehat{\sK}$ is~$\sK$.
\end{definition}

\begin{construction}
%The dual simplicial complex $\widehat{K}$ provides the following
%``purely simplicial interpretation" for the {\it Alexander
%duality\/}\label{aldua} (see e.g.~\cite[p.~54]{No1}) between $|K|$
%and $S^{m-1}\setminus|K|$ for any simplicial complex $K$ embedded
%in the $(m-1)$-sphere.
The barycentric subdivisions of both $\sK$ and $\widehat\sK$ can
be realised as subcomplexes in the barycentric subdivision of the
boundary of the standard simplex on the set~$[m]$ in the following
way.

A face of $(\partial\varDelta^{m-1})'$ corresponds to a chain
$I_1\subset\cdots\subset I_r$ of included subsets in~$[m]$ with
$I_1\ne\varnothing$ and $I_r\ne[m]$. We denote this face by
$\varDelta_{I_1\subset\cdots\subset I_r}$. (For example,
$\varDelta_{\{i\}}$ is the $i$th vertex of $\varDelta^{m-1}$
regarded as a vertex of $(\partial\varDelta^{m-1})'$.) Then
\begin{equation}\label{grkprime}
  G(\sK)=\bigcup_{I_1\subset\cdots\subset I_r,\;I_r\in \sK}
  \varDelta_{I_1\subset\cdots\subset I_r}
\end{equation}
is a geometric realisation of $\sK'$. Denote
$\widehat{i}=[m]\setminus\{i\}$ and, more generally,
$\widehat{I}=[m]\setminus{I}$ for any subset $I\subset[m]$. Define
the following subcomplex in $(\partial\varDelta^{m-1})'$:
\[
  D(\sK)=
  \bigcup_{I_1\subset\cdots\subset I_r,\;\widehat{I}_r\notin \sK}
  \varDelta_{\widehat{I}_r\subset\cdots\subset\widehat{I}_1}.
\]
\end{construction}

\begin{proposition}\label{dualbs}
For any simplicial complex $\sK\ne\varDelta^{m-1}$ on the
set~$[m]$, $D(\sK)$ is a geometric realisation of the barycentric
subdivision of the dual complex:
\begin{equation}\label{drkprime}
  D(\sK)=|\widehat{\sK}'|.
\end{equation}
Moreover, the open combinatorial neighbourhood of
complex~\eqref{grkprime} realising $\sK'$ in
$(\partial\varDelta^{m-1})'$ coincides with the complement of the
complex $D(\sK)$ realising $\widehat{\sK}'$:
\[
  {\mathop{\mathit
  U}\limits^\circ}_{(\partial\varDelta^{m-1})'}\bigl(|\sK'|\bigr)=
  \bigl(\partial\varDelta^{m-1}\bigr)'\setminus|\widehat{\sK}'|.
\]
In particular, $|\sK'|$ is a deformation retract of the complement
to $|\widehat{\sK}'|$ in $(\partial\varDelta^{m-1})'$.
\end{proposition}
\begin{proof}
We map a vertex $\{i\}$ of $\widehat \sK$ to the vertex $\widehat
i=\varDelta_{\widehat i}$ of $(\partial\varDelta^{m-1})'$, and map
the barycentre of a face $I\in \widehat{\sK}$ to the vertex
$\varDelta_{\widehat I}$. Then the whole complex $\widehat{\sK}'$
is mapped to the subcomplex
\[
  \bigcup_{I_1\subset\cdots\subset I_r,\;I_r\in{\widehat \sK}}
  \varDelta_{\widehat{I}_r\subset\cdots\subset\widehat{I}_1},
\]
which is the same as $D(\sK)$. The second statement is left as an
exercise.
\end{proof}

\begin{example}
Let $\sK$ be the boundary of a 4-gon with vertices $1,2,3,4$ (see
Fig.~\ref{figaldual}). Then $\widehat{\sK}$ consists of two
disjoint segments. The picture shows both $\sK'$ and
$\widehat{\sK}'$ as subcomplexes in~$(\partial\varDelta^3)'$.
\begin{figure}[h]
\begin{picture}(120,80)
\put(10,30){\circle*{1.6}} \put(7,28){1}
\put(60,80){\circle*{1.6}} \put(56,79){2}
\put(110,30){\circle*{1.6}} \put(112,28){3}
\put(50,0){\circle*{1.6}} \put(46,-2){4}
\multiput(9.6,30)(0.1,0){9}{\line(4,-3){40}}
\multiput(9.6,30)(0.1,0){9}{\line(1,1){50}}
\multiput(59.6,80)(0.1,0){9}{\line(1,-1){50}}
\multiput(49.5,0)(0.1,0){11}{\line(2,1){60}}
\multiput(10,30)(5,0){20}{\line(1,0){3}}
\qbezier(50,0)(60,80)(60,80) \put(35,35){\circle{1.6}}
\put(31,35){$\widehat{3}$} \put(75,35){\circle{1.6}}
\put(77,35){$\widehat{1}$} \put(35,34.2){\line(4,1){20}}
\put(35,35.8){\line(4,1){20}} \put(75,34.2){\line(-4,1){20}}
\put(75,35.8){\line(-4,1){20}} \put(55,20){\circle{1.6}}
\put(54,15){$\widehat{2}$} \put(60,46.7){\circle{1.6}}
\put(60,49){$\widehat{4}$}
\multiput(54.2,20)(2.8,5.6){2}{\line(1,2){1.9}}
\multiput(55.8,20)(2.8,5.6){2}{\line(1,2){1.9}}
\multiput(59.2,46.7)(0,-4.5){4}{\line(0,-3){2.8}}
\multiput(60.8,46.7)(0,-4.5){4}{\line(0,-3){2.8}}
\qbezier[60](35,35)(10,30)(10,30)
\qbezier[50](35,35)(30,50)(30,50)
\qbezier[75](35,35)(60,80)(60,80) \qbezier[65](35,35)(50,0)(50,0)
\qbezier[40](35,35)(26,18)(26,18)
\put(70,42.5){$|\widehat{\sK}'|$} \put(70,41){\line(-1,-2){1.9}}
\put(70,41){\line(-1,0){9}} \put(30,65){\line(1,-2){5}}
\put(25,66){$|\sK'|$}
\end{picture}
\caption{Dual complex and Alexander duality.} \label{figaldual}
\end{figure}
\end{example}

\begin{theorem}[Combinatorial Alexander duality]\label{simaldual}
For any simplicial complex $\sK\ne\varDelta^{m-1}$ on the
set~$[m]$ there is an isomorphism
\[
  \widetilde{H}^j(\sK)\cong
  \widetilde{H}_{m-3-j}(\widehat \sK),\quad\text{for } -1\le j\le m-2,
\]
where $\widetilde{H}_k(\cdot)$ and $\widetilde{H}^k(\cdot)$ denote
the $k$th reduced simplicial homology and cohomology group (with
integer coefficients) respectively. Here we assume that
$\widetilde{H}_{-1}(\varnothing)=\widetilde{H}^{-1}(\varnothing)=\Z$.
\end{theorem}
\begin{proof}
Since $(\partial\varDelta^{m-1})'$ is homeomorphic to $S^{m-2}$,
the Alexander duality theorem~\cite[Theorem~3.44]{hatc02} and
Proposition~\ref{dualbs} imply that
\begin{multline*}
  \widetilde H^j(\sK)=\widetilde H^j
  \bigl({\mathop{\mathit U}\limits^\circ}_{(\partial\varDelta^{m-1})'}(|\sK'|)\bigr)
  =\widetilde
  H^j\bigl((\partial\varDelta^{m-1})'\setminus|\widehat{\sK}'|\bigr)\\
  =\widetilde H^j\bigl(S^{m-2}\setminus|\widehat{\sK}|\bigr)\cong
  \widetilde{H}_{m-3-j}(\widehat \sK).\qedhere
\end{multline*}
%for $-1\le j\le m-2$.
\end{proof}

A more direct topological proof is outlined in
Exercise~\ref{exaldual}. Theorem~\ref{simaldual} can be also
proved in a purely combinatorial way, see~\cite{bj-ta09}. There is
also a proof within `combinatorial commutative algebra' (which is
the subject of Chapter~\ref{facerings}), see
Exercise~\ref{algproofad} or~\cite[Theorem~5.6]{mi-st05}.

The duality between $\sK$ and $\widehat\sK$ extends to a duality
between full subcomplexes of $\sK$ and links of simplices
in~$\widehat\sK$:

\begin{corollary}\label{fulink}
Let $\sK\ne\varDelta^{m-1}$ be a simplicial complex on~$[m]$ and
$I\notin \sK$, that is, $\widehat I\in\widehat{\sK}$. Then there
is an isomorphism
\[
  \widetilde{H}^j(\sK_I)\cong
  \widetilde{H}_{|I|-3-j}\bigl(\lk_{\widehat{\sK}}\widehat
  I\bigr),\quad\text{for } -1\le j\le|I|-2.
\]
\end{corollary}
\begin{proof}
We apply Theorem~\ref{simaldual} to the complex $\sK_I$, viewed as
a simplicial complex on the set~$I$ of $|I|$ elements. It follows
from the definition that the dual complex is
$\lk_{\widehat{\sK}}\widehat I$, which also can be viewed as a
simplicial complex on the set~$I$.
\end{proof}

\begin{example}
Let $\sK$ be the boundary of a pentagon. Then $\widehat{\sK}$ is a
M\"obius band triangulated as shown on Fig.~\ref{5dual}. Note that
this $\widehat{\sK}$ can be realised as a subcomplex in
$\partial\varDelta^4$, and therefore it can be realised in $\R^3$
as a subcomplex in the Schlegel diagram of~$\varDelta^4$, see
Construction~\ref{schldiag}.
\begin{figure}[h]
\begin{picture}(120,45)
\put(15,10){\circle*{1}} \put(10,30){\circle*{1}}
\put(25,40){\circle*{1}} \put(40,30){\circle*{1}}
\put(35,10){\circle*{1}} \put(55,15){\circle*{1}}
\put(75,15){\circle*{1}} \put(95,15){\circle*{1}}
\put(115,15){\circle*{1}} \put(65,35){\circle*{1}}
\put(85,35){\circle*{1}} \put(105,35){\circle*{1}}
\put(15,10){\line(-1,4){5}} \put(10,30){\line(3,2){15}}
\put(25,40){\line(3,-2){15}} \put(40,30){\line(-1,-4){5}}
\put(15,10){\line(1,0){20}}
\multiput(55,15)(20,0){3}{\line(1,2){10}}
\put(55,15){\vector(1,2){5}}
\multiput(65,35)(20,0){3}{\line(1,-2){10}}
\put(105,35){\vector(1,-2){5}} \put(55,15){\line(1,0){60}}
\put(65,35){\line(1,0){40}} \put(13,7){1} \put(7,29){2}
\put(24,41){3} \put(41,29){4} \put(36,7){5}
\put(54,10){$\widehat{1}$} \put(74,10){$\widehat{2}$}
\put(94,10){$\widehat{3}$} \put(114,10){$\widehat{4}$}
\put(64,37){$\widehat{4}$} \put(84,37){$\widehat{5}$}
\put(104,37){$\widehat{1}$} \put(23,1){$\sK$}
\put(83,1){$\widehat{\sK}$}
\end{picture}
\caption{The boundary of a pentagon and its dual complex.}
\label{5dual}
\end{figure}
\end{example}

Adding a ghost vertex to $\sK$ results in suspending
$\widehat{\sK}$, up to homotopy. The precise statement is as
follows (the proof is clear and is omitted).

\begin{proposition}
Let $\sK$ be a simplicial complex on $[m]$ and $\sK^\circ$ be the
complex on $[m+1]$ obtained by adding one ghost vertex $\circ=m+1$
to~$\sK$. Then the maximal simplices of $\widehat{\sK^\circ}$ are
$[m]$ and $I\cup\circ$ for $I\in\widehat{\sK}$, that is,
\[
  \widehat{\sK^\circ}=\varDelta^{m-1}\cup_{\widehat{\sK}}\cone\widehat{\sK}.
\]
In particular, $\widehat{\sK^\circ}$ is homotopy equivalent to the
suspension~$\Sigma\widehat{\sK}$.
\end{proposition}

\subsection*{Exercises.}
\begin{exercise}
Show that
\[
  {\mathop{\mathit
  U}\limits^\circ}_{(\partial\varDelta^{m-1})'}\bigl(|\sK'|\bigr)=
  \bigl(\partial\varDelta^{m-1}\bigr)'\setminus|\widehat{\sK}'|.
\]
\end{exercise}

\begin{exercise}\label{exaldual}
Complete the details in the following direct proof of
combinatorial Alexander duality (Theorem~\ref{simaldual}).

There is a simplicial map
\[
  \varphi\colon (\partial\varDelta^{m-1})'\to\sK'*\widehat{\sK}'
\]
which is constructed as follows. We identify $|\sK'|$ with
$G(\sK)$ and $|\widehat{\sK}'|$ with $D(\sK)$,
see~\eqref{grkprime} and~\eqref{drkprime}. The vertex sets of
these two subcomplexes split the vertex set of
$(\partial\varDelta^{m-1})'$ into two nonintersecting subsets.
Therefore, $\varphi$ is uniquely determined on the vertices. Check
that $\varphi$ is indeed a simplicial map.
%Let
%$\varDelta_{I_1\subset\ldots\subset I_r}\subset(\partial\D^{m-1})'$
%be a face. Then there is such $k$ that $I_k\in \sK$ and
%$I_{k+1}\notin \sK$. In this case we have
%\[
%  \alpha(\varDelta_{I_1\subset\ldots\subset I_r})=
%  \varDelta_{I_1\subset\ldots\subset I_k}*
%  \varDelta_{I_{k+1}\subset\ldots\subset I_r}
%  \subset|\sK'|*|\widehat{\sK}'|.
%\]
It induces a map of simplicial cochains
\[
  \varphi^*\colon C^j(\sK')\otimes C^{m-3-j}(\widehat{\sK}')\to
  C^{m-2}\bigl((\partial\varDelta^{m-1})'\bigr),
\]
or, equivalently,
\[
  C^j(\sK')\to C_{m-3-j}(\widehat{\sK}')\otimes
  C^{m-2}\bigl((\partial\varDelta^{m-1})'\bigr).
\]
By evaluating on the fundamental cycle of
$(\partial\varDelta^{m-1})'$ in $C_{m-2}$ and passing to
simplicial (co)homology, we obtain the required isomorphism
\[
  H^j(\sK)\stackrel\cong\longrightarrow H_{m-3-j}(\widehat{\sK}).
\]
\end{exercise}

\section{Classes of triangulated spheres}\label{simsph}
Boundary complexes of simplicial polytopes form an important
albeit restricted class of triangulated spheres. In this section
we review several other classes of triangulated spheres and
related complexes, and discuss their role in topology, geometry
and combinatorics.

\begin{definition}\label{defts}
A \emph{triangulated sphere} (also known as a \emph{sphere
triangulation} or \emph{simplicial sphere}) of dimension~$d$ is a
simplicial complex $\sK$ homeomorphic to a $d$-sphere~$S^d$. A
\emph{$PL$ sphere} is a triangulated sphere $\sK$ which is $PL$
homeomorphic to the boundary of a simplex (equivalently, there
exists a subdivision of~$\sK$ isomorphic to a subdivision of the
boundary of a simplex).
\end{definition}

\begin{remark}
A $PL$ sphere is not the same as a `$PL$ manifold homeomorphic to
a sphere', but rather a `$PL$ manifold which is $PL$ homeomorphic
to the \emph{standard} sphere', where the standard $PL$ structure
on a sphere is defined by triangulating it as the boundary of a
simplex. Nevertheless, the two notions coincide in dimensions
other than~4, see the discussion in the next section.
\end{remark}

In small dimensions, triangulated spheres can be effectively
visualised using \emph{Schlegel diagrams} and their
generalisations, which we describe below.

\begin{definition}
\label{schldiag} By analogy with the notion of a geometric
simplicial complex, we define a \emph{polyhedral complex} as a
collection $\mathcal C$ of convex polytopes in a space $\R^n$ such
that every face of a polytope in $\mathcal C$ belongs to $\mathcal
C$ and the intersection of any two polytopes in $\mathcal C$ is
either empty or a face of each. Two polyhedral complexes $\mathcal
C_1$ and $\mathcal C_2$ are said to be \emph{combinatorially
equivalent} if there is a one-to-one correspondence between their
polytopes respecting the inclusion of faces.

The boundary $\partial P$ of a convex polytope $P$ is a polyhedral
complex. Another important polyhedral complex associated with $P$
can be constructed as follows. Choose a facet $F$ of $P$ and a
point $p\notin P$ `close enough' to $F$ so that any segment
connecting $p$ to a point in $P\setminus F$ intersects the
relative interior of~$F$ (see Fig.~\ref{schlcube} for the case
$P=I^3$).
\begin{figure}[h]
\begin{picture}(120,70)
  %толстые наклонные линии
  \multiput(80,0)(0.1,0){7}{\line(3,2){30}}
  \multiput(80,50)(0.1,0){7}{\line(3,2){30}}
  \multiput(30,50)(0.1,0){7}{\line(3,2){30}}
  %невидимые линии куба
  \multiput(60,20)(-5,-3.33){6}{\line(-3,-2){3.6}}
  \multiput(60,20)(5,0){10}{\line(1,0){3}}
  \multiput(60,20)(0,5){10}{\line(0,1){3}}
  %проектирующие лучи
  \qbezier[110](110,70)(40,10)(40,10)
  \qbezier[90](110,19.5)(40,10)(40,10)
  \qbezier[80](60,70)(40,10)(40,10)
  \qbezier[30](60,20)(40,10)(40,10)
  %толстые прямые линии
  {\linethickness{0.4mm}
  \put(30,0){\line(1,0){50}}
  \put(30,0){\line(0,1){50}}
  \put(30,50){\line(1,0){50}}
  \put(80,50){\line(0,-1){50}}
  \put(110,20){\line(0,1){50}}
  \put(60,70){\line(1,0){50}}
  \put(110,20){\line(0,1){50}}
  }
  %диаграмма Шлегеля
  \put(63.5,30){\line(-1,0){17}}
  \put(63.5,30){\line(0,-1){17}}
  \put(46.5,30){\line(0,-1){17}}
  \put(46.5,13){\line(1,0){17}}
  \qbezier(46.5,13)(30,0)(30,0)
  \qbezier(46.5,30)(30,50)(30,50)
  \qbezier(63.5,30)(80,50)(80,50)
  \qbezier(63.5,13)(80,0)(80,0)
  \put(37,10){$p$}
\end{picture}
\caption{Schlegel diagram of $I^3$.} \label{schlcube}
\end{figure}
Now project the complex $\partial P$ onto $F$ from the point~$p$.
The projection images of faces of~$P$ different from~$F$ form a
polyhedral complex $\mathcal C$ subdividing~$F$. We refer to
$\mathcal C$, and also to any polyhedral complex combinatorially
equivalent to it, as a~\emph{Schlegel diagram} of the
polytope~$P$.

An \emph{$n$-diagram}\label{ndiag} is a polyhedral
complex~$\mathcal C$ consisting of $n$-polytopes and their faces
and satisfying the following conditions:
\begin{itemize}
\item[(a)] the union of all polytopes in~$\mathcal C$ is an $n$-polytope~$Q$;
\item[(b)] every nonempty intersection of a polytope from $\mathcal C$ with the boundary
of $Q$ belongs to~$\mathcal C$.
\end{itemize}
We refer to $Q$ as the \emph{base}\label{basesch} of the
$n$-diagram~$\mathcal C$. An $n$-diagram is \emph{simplicial} if
it consists of simplices and its base is a simplex.
\end{definition}

By definition, a Schlegel diagram of an $n$-polytope is an
$(n-1)$-diagram.

\begin{proposition}
Let $\mathcal C$ be a simplicial $(n-1)$-diagram with base~$Q$.
Then $\mathcal C\cup_{\partial Q}Q$ is an $(n-1)$-dimensional $PL$
sphere.
\end{proposition}
\begin{proof}
We need to construct a simplicial complex which is a subdivision
of both $\mathcal C\cup_{\partial Q} Q$ and $\partial\varDelta^n$.
This can be done as follows. Replace one of the facets of
$\partial\varDelta^n$ by the diagram $\mathcal C$. The resulting
simplicial complex~$\sK$ is a subdivision
of~$\partial\varDelta^n$. On the other hand, $\sK$ is isomorphic
to the subdivision of $\mathcal C\cup_{\partial Q} Q$ obtained by
replacing $Q$ by a Schlegel diagram of~$\varDelta^n$.
\end{proof}

\begin{corollary}\label{polytpl}
The boundary of a simplicial $n$-polytope is an
$(n-1)$-dimensional $PL$ sphere.
\end{corollary}
%\begin{proof}
%Indeed, a Schlegel diagram of a simplicial $n$-polytope is a
%simplicial $(n-1)$-diagram.
%\end{proof}

In dimensions $n\le3$ every $(n-1)$-diagram is a Schlegel diagram
of an $n$-polytope. Indeed, for $n\le2$ this is obvious, and for
$n=3$ this is one of equivalent formulations of the well-known
\emph{Steinitz Theorem}\label{stein}
(see~\cite[Theorem~5.8]{zieg95}). The first example of a 3-diagram
which is not a Schlegel diagram of a 4-polytope was found by
Gr\"unbaum (\cite[\S11.5]{grue03}, see also~\cite{gr-sr67}) as a
correction of Br\"uckner's result of 1909 on the classification of
simplicial 4-polytopes with 8 vertices. Another example was found
by Barnette~\cite{barn70}. We describe Barnette's example below.
For the original example of Gr\"unbaum, see
Exercise~\ref{brueckner}.

\begin{construction}[Barnette's 3-diagram]\label{barsph}
Here a certain simplicial 3-diagram $\mathcal C$ will be
constructed. Consider the octahedron~$Q$ obtained by twisting the
top face $(abc)$ of a triangular prism (Fig.~\ref{figbarn}~(a))
slightly so that the vertices $a$, $b$, $c$, $d$, $e$ and $f$ are
in general position.
\begin{figure}[h]
\begin{picture}(120,80)
  %октаэдр
  %невидимые линии
  \multiput(0,20)(5.2,1.5){7}{\line(4,1){3.6}}
  \multiput(15,30)(1.5,5.8){7}{\line(1,4){1}}
  \multiput(25,50)(5.1,-6.1){5}{\line(3,-4){2.7}}
  %видимые линии
  \put(0,20){\line(1,0){50}}
  \put(0,20){\line(3,2){15}}
  \qbezier(0,20)(25,50)(25,50)
  \put(0,20){\line(1,2){25}}
  \put(25,70){\line(0,-1){20}}
  \put(25,70){\line(1,-4){10}}
  \put(25,70){\line(1,-2){25}}
  \put(50,20){\line(-3,2){15}}
  \qbezier(50,20)(15,30)(15,30)
  \put(15,30){\line(1,0){20}}
  \put(15,30){\line(1,2){10}}
  \put(25,50){\line(1,-2){10}}
  \put(0,17){$d$}
  \put(49,17){$e$}
  \put(24,71){$f$}
  \put(17,31){\small $a$}
  \put(32,31){\small $b$}
  \put(24,46){\small $c$}
  \put(23,-1){(a)}
  %диаграмма Шлегеля
  \put(60,15){\line(2,-1){30}}
  \multiput(60,15)(5.1,0){12}{\line(1,0){3}}
  \multiput(60,15)(7.3,5.45){6}{\line(4,3){3.6}}
  \put(60,15){\line(2,3){20}}
  \qbezier(60,15)(90,80)(90,80)
  \qbezier(90,80)(80,45)(80,45)
  \put(90,80){\line(0,-1){80}}
  \qbezier(90,80)(95,30)(95,30)
  \qbezier(90,80)(100,45)(100,45)
  \qbezier(90,80)(120,15)(120,15)
  \put(120,15){\line(-2,-1){30}}
  \qbezier(120,15)(95,30)(95,30)
  \put(120,15){\line(-2,3){20}}
  \qbezier(90,0)(80,45)(80,45)
  \qbezier(90,0)(95,30)(95,30)
  \put(95,30){\line(-1,1){15}}
  \put(95,30){\line(1,3){5}}
  \put(80,45){\line(1,0){20}}
  \put(58,12){$d$}
  \put(87,-2){$e$}
  \put(120,15){$f$}
  \put(78,45){\small $a$}
  \put(96,30){\small $b$}
  \put(101,45){\small $c$}
  \put(86,78){$p'$}
  \put(73,-1){(b)}
\end{picture}
\caption{Barnette's 3-diagram.} \label{figbarn}
\end{figure}
Assume that the edges $(bd)$, $(ce)$ and $(af)$ lie inside~$Q$.
The tetrahedra $(abde)$, $(bcef)$ and $(acdf)$ will be included in
the complex~$\mathcal C$. Each of these tetrahedra has two faces
which lie inside~$Q$. These 6 triangles together with $(abc)$ and
$(def)$ form a triangulated 2-sphere inside~$Q$, which we denote
by~$\mathcal S$. Now place a point $p$ inside $\mathcal S$ so that
the line segments from $p$ to the vertices of $\mathcal S$ lie
inside~$\mathcal S$. We add to $\mathcal C$ the eight tetrahedra
obtained by taking cones with vertex $p$ over the faces
of~$\mathcal S$ (namely, the tetrahedra $(pabc)$, $(pdef)$,
$(pabd)$, $(pbed)$, $(pbce)$, $(pcef)$, $(pacf)$ and $(padf)$).
Note that $\mathcal S=\lk p$. Applying a projective transformation
if necessary, we may assume that there is a point $p'$ outside the
octahedron~$Q$ with the property that the segments joining $p'$
with the vertices of~$Q$ lie outside~$Q$ (see
Fig.~\ref{figbarn}~(b)). Finally, we add to $\mathcal C$ the
tetrahedra obtained by taking cones with vertex $p'$ over the
faces of $Q$ other than the face~$(def)$ (there are 7 such
tetrahedra: $(p'abc)$, $(p'abe)$, $(p'ade)$, $(p'acd)$, $(p'cdf)$,
$(p'bcf)$ and $(p'bef)$). The union of all tetrahedra in $\mathcal
C$ is the tetrahedron $(p'def)$; hence, $\mathcal C$ is a
simplicial 3-diagram. It has 8 vertices and 18 tetrahedra.
\end{construction}

\begin{proposition}
The 3-diagram $\mathcal C$ from the previous construction is not a
Schlegel diagram.
\end{proposition}
\begin{proof}
Suppose there is a polytope $P$ whose Schlegel diagram
is~$\mathcal C$. Since $P$ is simplicial, we may assume that its
vertices are in general position. We label the vertices of $P$
with the same letters as the corresponding vertices in~$\mathcal
C$. Consider the complex $\mathcal S=\lk p$. Let $P'$ be the
convex hull of all vertices of~$P$ other than~$p$. Then $P'$ is
still a simplicial polytope, and $\mathcal S$ is a subcomplex
in~$\partial P'$. In the complex $\partial P'$ the sphere
$\mathcal S$ is filled with tetrahedra whose vertices belong
to~$\mathcal S$. Then at least one edge of one of these tetrahedra
lies inside~$\mathcal S$. However, any two vertices of~$\mathcal
S$ which are not joined by an edge on~$\mathcal S$ are joined by
an edge of~$\mathcal C$ lying outside~$\mathcal S$. Since the
polytope $P'$ cannot contain a double edge we have reached a
contradiction.
\end{proof}

Now we can introduce two more classes of triangulated spheres.

\begin{definition}\label{defps}
A \emph{polytopal sphere} is a triangulated sphere isomorphic to
the boundary complex of a simplicial polytope.

A \emph{starshaped sphere} is a triangulated sphere isomorphic to
the underlying complex of a complete simplicial fan. Equivalently,
a triangulated sphere $\sK$ of dimension $n-1$ is starshaped if
there is a geometric realisation of $\sK$ in $\R^n$ and a point
$p\in\R^n$ with the property that each ray emanating from $p$
meets $|\sK|$ in exactly one point. The set of points such
$p\in\R^n$ is called the \emph{kernel} of~$|\sK|$.
\end{definition}

\begin{example}\label{barnstar}
The triangulated 3-sphere coming from Barnette's 3-diagram is
known as the \emph{Barnette sphere}. It is starshaped. Indeed, in
the construction of Barnette's 3-diagram we have the octahedron
$Q\subset\R^3$ and a vertex $p'$ outside $Q$ such that $\lk
p'=\partial Q$. If we choose the vertex $p'$ in
$\R^4\setminus\R^3$, then the Barnette sphere can be realised in
$\R^4$ as the boundary complex of the pyramid with vertex $p'$ and
base $Q$ (subdivided as described in Construction~\ref{barsph}).
This realisation is obviously starshaped.
\end{example}

We therefore have the following hierarchy of triangulations:
\begin{equation}
\label{polplsim}
  \parbox{0.2\textwidth}{\begin{center}polytopal\\spheres\end{center}}\subset
  \parbox{0.2\textwidth}{\begin{center}starshaped\\spheres\end{center}}\subset
  \parbox{0.2\textwidth}{\begin{center}$PL$\\spheres\end{center}}\subset
  \parbox{0.2\textwidth}{\begin{center}triangulated\\spheres\end{center}}
\end{equation}
Here the first inclusion follows from the construction of the
normal fan (Construction~\ref{nf}), and the second is left as an
exercise.

In dimension 2 any triangulated sphere is polytopal (this is
another corollary of the Steinitz theorem). Also, by the result of
Mani~\cite{mani72}, any triangulated $d$-dimensional sphere with
up to $d+4$ vertices is polytopal. However, in general all
inclusions above are strict.

The first inclusion in~\eqref{polplsim} is strict already in
dimension~3, as is seen from Example~\ref{barnstar}. There are 39
combinatorially different triangulations of a 3-sphere with 8
vertices, among which exactly two are nonpolytopal (namely, the
Barnette and Br\"uckner spheres); this classification was
completed by Barnette~\cite{barn73j}.

The second inclusion in~\eqref{polplsim} is also strict in
dimension~3. The first example of a nonstarshaped sphere
triangulation was found by Ewald and Schulz~\cite{ew-sc92}. We
sketch this example below, following~\cite[Theorem~5.5]{ewal96}.

\begin{example}[Nonstarshaped sphere triangulation]
We use the fact, observed by Barnette, that not every tetrahedron
in Barnette's 3-diagram can be chosen as the base of a $3$-diagram
of the Barnette sphere (see Exercise~\ref{barnbase}). For
instance, the tetrahedron $(abcd)$ cannot be chosen as the base.

Now let $\sK$ be a connected sum of two copies of the Barnette
sphere along the tetrahedra $(abcd)$ (the identification of
vertices in the two tetrahedra is irrelevant). Assume that $\sK$
has a starshaped realisation in~$\R^4$. The hyperplane $H$ through
the points $a,b,c,d$ splits $|\sK|$ into two parts $|\sK_1|$ and
$|\sK_2|$. Since the kernel of $|\sK|$ is an open set, it contains
a point $p$ not lying in~$H$. Then by projecting either $|\sK_1|$
or $|\sK_2|$ onto $H$ from~$p$, we obtain a 3-diagram of the
Barnette sphere with base $(abcd)$. This is a contradiction.

The fact that $\sK$ is a $PL$ sphere is an exercise.
\end{example}

The third inclusion in~\eqref{polplsim} is the subtlest one.
%In order to analyse it properly we shall need some deep results
%from geometric topology.
It is known that in dimension~3 any triangulated sphere is~$PL$.
In dimension~4 the corresponding question is open, but starting
from dimension 5 there exist \emph{non-$PL$}\label{nonPLsph}
sphere triangulations. See the discussion in the next section and
Example~\ref{non-PL}.

Many important open problems of combinatorial geometry arise from
analysing the relationships between different classes of sphere
triangulations. We end this section by discussing some of these
problems.

In connection with the condition of realisability of a
triangulated $(n-1)$-sphere in $\R^n$ in the definition of a
starshaped sphere, we note that the existence of such a
realisation is open in general:

\begin{problem}
Does every $PL$ $(n-1)$-sphere admit a geometric realisation in an
$n$-dimensional space?
\end{problem}

The $g$-theorem (Theorem~\ref{gth}) gives a complete
characterisation of integral vectors arising as the $f$-vectors of
polytopal spheres. It is therefore natural to ask whether the
$g$-theorem extends to all sphere triangulations. This question
was posed by McMullen~\cite{mcmu71} as an extension of his
conjecture for simplicial polytopes. Since 1980, when McMullen's
conjecture for simplicial polytopes was proved by Billera, Lee,
and Stanley, its generalisation to spheres has been regarded as
the main open problem in the theory of $f$-vectors:

\begin{problem}[$g$-conjecture for triangulated spheres]
\label{gconj} Does Theorem~{\rm\ref{gth}} hold for triangulated
spheres?
\end{problem}

The $g$-conjecture is open even for starshaped spheres. Note that
only the necessity of the conditions in the $g$-theorem (that is,
the fact that every $g$-vector is an $M$-vector) has to be
verified for triangulated spheres. If correct, the $g$-conjecture
would imply a complete characterisation of $f$-vectors of
triangulated spheres.

The Dehn--Sommerville relations (condition~(a) in
Theorem~\ref{gth}) hold for arbitrary sphere triangulations (see
Corollary~\ref{dsgor} below). The $f$-vectors of triangulated
spheres also satisfy the UBT and LBT inequalities given in
Theorems~\ref{ubt} and~\ref{lbt} respectively. The proof of the
Lower Bound Theorem for simplicial polytopes given by Barnette
in~\cite{barn73} extends to all triangulated spheres (see
also~\cite{kala87}). In particular, this implies the second GLBC
inequality $h_1\le h_2$, see~\eqref{glbt}. The Upper Bound Theorem
for triangulated spheres was proved by Stanley~\cite{stan75} (we
shall give his argument in Section~\ref{cmr}). This implies that
the $g$-conjecture is true for triangulated spheres of
dimension~$\le4$. The third GLBC inequality $h_2\le h_3$ (for
spheres of dimension $\ge 5$) is open.

Many attempts to prove the $g$-conjecture were made after 1980.
Though unsuccessful, these attempts resulted in some very
interesting reformulations of the $g$-conjecture. The results of
Pachner~\cite{pach91} reduce the $g$-conjecture (for $PL$ spheres)
to some properties of \emph{bistellar moves}\label{bist3}; see the
discussion after Theorem~\ref{gconjbs} below.
%We also mention the results of~\cite{TWW} showing that the $g$-conjecture follows from the
%{\it skeletal $r$-rigidity\/} of simplicial $(n-1)$-sphere for $r\le\sbr n2$.
%As was shown independently by Kalai and Stanley
%\cite[Corollary~2.4]{stan93}, the GLBC inequalities (i.e.
%condition~(b) in Theorem~\ref{gth}) hold for the boundary of an
%$n$-dimensional ball which is a subcomplex of the boundary complex
%of a simplicial $(n+1)$-polytope. However, it is not yet clear
%which simplicial complexes arise in this way.

The lack of progress in proving the $g$-conjecture motivated
Bj\"orner and Lutz to launch a computer-aided search for
counterexamples~\cite{bj-lu00}. Though their bistellar flip
algorithm and BISTELLAR software produced many remarkable results
on triangulations of manifolds, no counterexamples to the
$g$-conjecture were found. More information on the $g$-conjecture
and related questions may be found in~\cite{stan96} and
\cite[Lecture~8]{zieg95}.

We may extend the hierarchy~\eqref{polplsim} by considering
polyhedral complexes homeomorphic to spheres (the so-called
\emph{polyhedral spheres}) instead of triangulated spheres. There
are obvious analogues of polytopal and $PL$ spheres in this
generality. However, unlike the case of triangulations, the two
definitions of a starshaped sphere (namely, the one using fans and
the one using the kernel points) no longer produce the same
classes of objects, see~\cite[\S\,III.5]{ewal96} for the
corresponding examples. One of the most notorious and long
standing problems is to find a proper higher dimensional analogue
to the Steinitz theorem. This theorem characterises graphs of
3-dimensional polytopes, and one of its equivalent formulations is
that every polyhedral 2-sphere is polytopal. In higher dimensions,
identification of the class of polytopal spheres inside all
polyhedral spheres is known as the \emph{Steinitz problem}:

\begin{problem}[Steinitz Problem]\label{steinproblem}
Find necessary and sufficient conditions for a polyhedral
decomposition of a sphere to be combinatorially equivalent to the
boundary complex of a convex polytope.
\end{problem}

This is far from being solved even in the case of triangulated
spheres. For more information on the relationships between
different classes of polyhedral spheres and complexes see the
above cited book of Ewald~\cite{ewal96} and the survey article by
Klee and Kleinschmidt~\cite{kl-kl95}.

\subsection*{Exercises.}
\begin{exercise} Show that $\sK$ is the underlying complex of a
complete simplicial fan if and only if there is a geometric
realisation of $\sK$ in $\R^n$ and a point $p\in\R^n$ with the
property that each ray emanating from $p$ meets $|\sK|$ in exactly
one point.
\end{exercise}

\begin{exercise}
Prove that every starshaped sphere is a $PL$ sphere.
\end{exercise}

\begin{exercise}[Br\"uckner sphere]\label{brueckner}
The \emph{Br\"uckner sphere} is obtained by replacing two
tetrahedra $(pabc)$ and $(p'abc)$ in the Barnette sphere by three
tetrahedra $(pp'ab)$, $(pp'ac)$ and $(pp'bc)$ (this is an example
of a bistellar 1-move considered in Section~\ref{stbist}). Show
that the Br\"uckner sphere is starshaped but not polytopal. Note
that the 1-skeleton of the Br\"uckner sphere is a complete graph
(that is, the Br\"uckner sphere is a \emph{neighbourly}
triangulation, see Definition~\ref{neighb}).
%, unlike the Barnette sphere.
\end{exercise}

\begin{exercise}\label{barnbase}
Show that the tetrahedron $(abcd)$ in Barnette's 3-diagram
(Construction~\ref{barsph}) cannot be chosen as the base of a
$3$-diagram of the Barnette sphere. Which tetrahedra can be chosen
as the base?
\end{exercise}

\begin{exercise}
The connected sum of two $PL$ spheres is a $PL$ sphere.
\end{exercise}

\section{Triangulated manifolds}
Piecewise linear topology experienced an intensive development
during the second half of the 20th century, thanks to the efforts
of many topologists. Surgery theory for simply connected manifolds
of dimension $\ge5$ originated from the early work of Milnor,
Kervaire, Browder, Novikov and Wall, culminated in the proof of
the topological invariance of rational Pontryagin classes given by
Novikov in 1965, and was further developed in the work of Lashof,
Rothenberg, Sullivan, Kirby, Siebenmann, and others. It led to a
better understanding of the place of $PL$ manifolds between the
topological and smooth categories. Without attempting to overview
the current state of the subject, which is generally beyond the
scope of this book, we include several important results on the
triangulation of topological manifolds, with a particular emphasis
on various nonexamples. We also provide references for further
reading.

All manifolds here are compact, connected and closed, unless
otherwise stated.

\begin{definition}\label{simman1}
A \emph{triangulated manifold} (or \emph{simplicial manifold}) is
a simplicial complex $\sK$ whose geometric realisation $|\sK|$ is
a topological manifold.

A \emph{$PL$ manifold} is a simplicial complex $\sK$ of
dimension~$d$ such that $\lk I$ is a $PL$ sphere of dimension
$d-|I|$ for every nonempty simplex $I\in\sK$.
\end{definition}

Every $PL$ manifold $|\sK|$ of dimension $d$ is a triangulated
manifold: it has an atlas whose change of coordinates functions
are piecewise linear. Indeed, for each vertex $v\in|\sK|$ the
$(d-1)$-dimensional $PL$ sphere $\lk v$ bounds an open
neighbourhood $U_v$ which is homeomorphic to an open $d$-ball.
Since any point of $|\sK|$ is contained in $U_v$ for some~$v$,
this defines an atlas for~$|\sK|$.

\begin{remark}
The term `$PL$ manifold' is often used for a manifold with a $PL$
atlas, while its particular triangulation with the property above
is referred to as a combinatorial manifold. We shall not
distinguish between these two notions.
\end{remark}

Does every triangulation of a topological manifold yield a
simplicial complex which is a $PL$ manifold?  The answer is `no',
and the question itself ascends to a famous conjecture of the dawn
of topology, known as the \emph{Hauptvermutung}\label{haupt},
which is German for `main conjecture'.  Below we briefly review
the current status of this conjecture, referring to the survey
article~\cite{rani96} by Ranicki for a much more detailed
historical account and more references.

In the early days of topology all of the known topological
invariants were defined in combinatorial terms, and it was very
important to find out whether the topology of a triangulated space
fully determines the combinatorial equivalence class of the
triangulation (in the sense of Definition~\ref{simmap}). In the
modern terminology, the Hauptvermutung states that any two
homeomorphic simplicial complexes are combinatorially equivalent
($PL$ homeomorphic). This is valid in dimensions $\le3$; the
result is due to Rado (1926) for 2-manifolds, Papakyriakopoulos
(1943) for 2-complexes, Moise (1953) for 3-manifolds, and E.~Brown
(1963) for 3-complexes~\cite{brow69}; see~\cite{mois77} for a detailed
exposition. The first examples of complexes disproving the
Hauptvermutung in dimensions $\ge6$ were found by Milnor in the
early 1960s. However, the \emph{manifold Hauptvermutung}, namely
the question of whether two homeomorphic triangulated manifolds
are combinatorially equivalent, remained open until the end of the
1960s. The first counterexamples were found by Siebenmann in 1969, and
relied heavily on the topological surgery theory. The `double
suspension theorem', which we state as Theorem~\ref{cannon} below,
appeared around 1975 and provided much more explicit
counterexamples to the manifold Hauptvermutung.

A \emph{$d$-dimensional homology sphere} (or simply \emph{homology
$d$-sphere}) is a topological $d$-manifold whose integral homology
groups are isomorphic to those of a $d$-sphere~$S^d$.

\begin{theorem}[Edwards, Cannon]\label{cannon}
The double suspension of any homology $d$-sphere is homeomorphic
to $S^{d+2}$.
\end{theorem}

This theorem was proved for most double suspensions and all triple
suspensions by Edwards~\cite{edwa75}; the general case was done by
Cannon~\cite{cann79}. One of its most important consequences is
the existence of non-$PL$ triangulations of 5-spheres, which also
disproves the manifold Hauptvermutung in dimensions~$\ge5$.

\begin{example}[non-$PL$ triangulated 5-sphere]
\label{non-PL} Let $M$ be a triangulated homology 3-sphere which
is not homeomorphic to~$S^3$. An example of such $M$ is provided
by the \emph{Poincar\'e sphere}. It is the homogeneous space
$SO(3)/A_5$, where the alternating group $A_5$ is represented in
$\R^3$ as the group of self-transformations of a dodecahedron. A
particular symmetric triangulation of the Poincar\'e sphere is
given in~\cite{bj-lu00}. By Theorem~\ref{cannon}, the double
suspension $\Sigma^2M$ is homeomorphic to~$S^5$ (and, more
generally, $\Sigma^k M$ is homeomorphic to~$S^{k+3}$ for $k\ge2$).
However, $\Sigma^2M$ cannot be a $PL$ sphere, since $M$ appears as
the link of a 1-simplex in $\Sigma^2M$.

Also, according to a result of Bj\"orner and Lutz~\cite{bj-lu00},
for any $d\ge5$ there is a non-$PL$ triangulation of $S^d$ with
$d+13$ vertices.
\end{example}

%In the positive direction, it is known that two homeomorphic
%simply connected $PL$ manifolds of dimension $\ge5$ with no
%torsion in third homology group are combinatorially equivalent
%($PL$ homeomorphic). This is Sullivan's famous Hauptvermutung
%theorem. The general classification of $PL$ structures on higher
%dimensional topological manifolds was obtained by Kirby and
%Siebenmann, see~\cite{KS}.

Theorem~\ref{cannon} led to progress in the following `manifold
recognition problem': given a simplicial complex, how one can
decide whether its geometric realisation is a topological
manifold? In higher dimensions there is the following result,
which can be viewed as a generalisation of the double suspension
theorem.

\begin{theorem}[Edwards~\cite{edwa80}]
\label{edwards} For $d\ge5$ the realisation of a simplicial
complex $\sK$ is a topological manifold of dimension~$d$ if and
only if $\lk I$ has the homology of a $(d-|I|)$-sphere for each
nonempty simplex $I\in\sK$, and $\lk v$ is simply connected for
each vertex $v$ of~$\sK$.
\end{theorem}

\begin{remark}
From the algorithmic point of view, the homology of links is
easily computable, but their simply connectedness seems to be
undecidable. There is a related result of
Novikov~\cite[Appendix]{v-k-f74} that a triangulated 5-sphere
cannot be algorithmically recognised. On the other hand, the
algorithmic recognition problem for a triangulated 3-sphere has a
positive solution, with the first algorithm provided by
Rubinstein~\cite{rubi95}. See detailed exposition in Matveev's
book~\cite{matv03}, which also contains a proof that all
3-dimensional \emph{Haken manifolds} can be recognised and fully
classified algorithmically.
\end{remark}

With the discovery of exotic smooth structures on 7-spheres by
Milnor and the disproval of the Hauptvermutung it had become
important to understand better the relationship between $PL$ and
smooth structures on topological manifolds. Since a $PL$ structure
implies the existence of a particular sort of triangulation, the
related question of whether a topological manifold admits
\emph{any} triangulation (not necessarily~$PL$) had also become
important.
%Another related question is whether the
%{\it Hauptvermutung\/} is valid in dimension~4. Both questions
%were answered (negatively) by the results of Freedman and
%Donaldson (early 1980s).

Triangulations of 2-manifolds have been known from the early days
of topology. A proof that any 3-manifold can be triangulated was
obtained independently by Moise and Bing in the end of 1950s (the
proof can be found in~\cite{mois77}). Since the link of a vertex
in a triangulated 3-sphere is a 2-sphere, and a 2-sphere is always
$PL$, all topological 3-manifolds are~$PL$.

A smooth manifold of any dimension has a $PL$ triangulation by a
theorem of Whitney (a proof can be found in~\cite{munk66}).
Moreover, in dimensions $\le3$ every topological manifold has a
unique smooth structure, see~\cite{mois77} for a proof. All these
considerations show that in dimensions up to 3 the categories of
topological, $PL$ and smooth manifolds are equivalent.

The situation in dimension 4 is quite different. There exist
topological 4-manifolds that do not admit a $PL$ triangulation. An
example is provided by Freedman's fake $\C
P^2$~\cite[\S8.3,~\S10.1]{fr-qu90}, a topological manifold which
is homeomorphic, but not diffeomorphic to the complex projective
plane~$\C P^2$. This example also shows that the Hauptvermutung is
false in dimension~4. Even worse, some topological 4-manifolds do
not admit any triangulation; an example is the topological
4-manifold with the intersection form~$E_8$, see~\cite{ak-mc90}.

In dimension 4 the categories of $PL$ and smooth manifolds agree,
that is, there is exactly one smooth structure on every $PL$
manifold. However, the classification of $PL$ (or equivalently,
smooth) structures is wide open even for the simplest topological
manifolds. The most notable problem here is the following.

\begin{problem}
Is a $PL$ (or smooth) structure on a 4-sphere unique?
\end{problem}

%Four is the only dimension where the uniqueness of a $PL$
%structure on a topological sphere is open. For dimensions $\le 3$
%the uniqueness was proved by Moise~\cite{Mo}, and for dimensions
%$\ge 5$ it follows from the result of Kirby and
%Siebenmann~\cite{KS}. In dimension~4 the category of $PL$
%manifolds is equivalent to the smooth category, hence, the above
%problem is equivalent to if there exists an exotic (or fake)
%4-sphere.

In dimensions $\ge5$ the $PL$ structure on a topological sphere is
unique (that is, a $PL$ manifold which is homeomorphic to a sphere
is a $PL$ sphere).

%\begin{problem}[Triangulation Conjecture]\label{triconj}
%Is it true that any topological manifold of dimension $\ge5$ can
%be triangulated?
%\end{problem}

For the discussion of the classification of $PL$ structures on
topological manifolds we refer to Ranicki's survey~\cite{rani96},
the original essay~\cite{ki-si77} by Kirby and Siebenmann, and a
more recent survey by Rudyak~\cite{rudy01}.

\section{Stellar subdivisions and bistellar moves}\label{stbist}
By a theorem of Alexander, a common subdivision of two $PL$
homeomorphic $PL$ manifolds can be obtained by iterating
operations from a very simple and explicit list, known as
\emph{stellar subdivisions}. An even more concrete iterative
description of $PL$ homeomorphisms was obtained by
Pachner~\cite{pach87}, who introduced the notion of
\emph{bistellar moves} (in other terminology, \emph{bistellar
flips} or \emph{bistellar operations}), generalising the 2- and
3-dimensional \emph{flips} from low-dimensional topology. These
operations allow us to decompose a $PL$ homeomorphism into a
sequence of simple `moves' and thus provide a very convenient way
to compute and handle topological invariants of $PL$ manifolds.
Starting from a given $PL$ triangulation, bistellar operations may
be used to construct new triangulations with some good properties,
such as ones that are symmetric or have a small number of vertices. On
the other hand, bistellar moves can be used to produce some nasty
triangulations if we start from a non-$PL$ triangulation. Both
approaches were used in the work of Bj\"orner--Lutz~\cite{bj-lu00}
and Lutz~\cite{lutz99} to find many interesting triangulations of
low-dimensional manifolds. In our exposition of bistellar moves we
follow the terminology of~\cite{lutz99}.

Bistellar moves also provide a combinatorial interpretation for
algebraic \emph{flop} operations for projective \emph{toric
varieties} and for certain surgery
operations on moment-angle complexes and torus manifolds. Finally,
bistellar moves may be used to define a metric on the space of
$PL$ triangulations of a given $PL$ manifold, see~\cite{nabu96}.

\begin{definition}[stellar subdivisions and bistellar moves]
\label{bist} Let $I\in\sK$ be a nonempty simplex of a simplicial
complex~$\sK$. The \emph{stellar subdivision} of $\sK$ at $I$ is
obtained by replacing the star of $I$ by the cone over its
boundary:
\[
  \ss_I\sK=
  (\sK\setminus\st_{\sK}I)\cup(\cone\partial\st_{\sK}I).
\]
If $\dim I=0$ then $\ss_I\sK=\sK$. Otherwise the complex
$\ss_I\sK$ acquires an additional vertex (the vertex of the
cone) whose link is $\partial\st_{\sK}I$.
\begin{figure}[h]
\begin{center}
\begin{picture}(90,60)
\put(0,15){\line(1,-1){15}}
\put(15,0){\line(1,1){15}}
\put(30,15){\line(-1,1){15}}
\put(15,30){\line(-1,-1){15}}
\put(15,0){\line(0,1){30}}
\put(75,0){\line(0,1){30}}
\put(75,15){\circle*{1.5}}
\put(35,15){\vector(1,0){20}}
\put(42,17){\small $\ss_I\sK$}
\put(38,11){\small $\dim I=1$}
\put(60,15){\line(1,-1){15}}
\put(75,0){\line(1,1){15}}
\put(90,15){\line(-1,1){15}}
\put(75,30){\line(-1,-1){15}}
\put(60,15){\line(1,0){30}}
\put(0,40){\line(3,4){15}}
\put(15,60){\line(3,-4){15}}
\put(0,40){\line(1,0){30}}
\put(60,40){\line(3,1){15}}
\put(75,60){\line(0,-1){15}}
\put(75,45){\line(3,-1){15}}
\put(75,45){\circle*{1.5}}
\put(35,45){\vector(1,0){20}}
\put(42,47){\small $\ss_I\sK$}
\put(38,41){\small $\dim I=2$}
\put(60,40){\line(3,4){15}}
\put(75,60){\line(3,-4){15}}
\put(60,40){\line(1,0){30}}
\end{picture}
\caption{Stellar subdivisions at a 2-simplex and at an edge.}
\label{stfig}
\end{center}
\end{figure}
Two possible stellar subdivisions of a 2-dimensional complex are
shown in Fig.~\ref{stfig}.

Now let $\sK$ be a triangulated manifold of dimension~$d$. Assume
that $I\in\sK$ is a $(d-j)$-face such that the simplicial complex
$\lk_{\sK}I$ is the boundary of a $j$-simplex $J$ which is not a
face of~$\sK$. Then the operation $\bm_I$ on $\sK$ defined by
$$
  \bm_I\sK=\bigl(\sK\setminus(I*\partial J)\bigr)\cup
  (\partial I*J)
$$
is called a \emph{bistellar $j$-move}. Since $I*\partial
J=\st_{\sK}I$ and $\partial I*J=\st_{\widetilde{\sK}}J$, where
$\widetilde{\sK}=\bm_I\sK$, the bistellar $j$-move is the
composition of a stellar subdivision at $I$ and the inverse
stellar subdivision at~$J$, which explains the term. In
particular, the stellar subdivision $\ss_I\sK$ is a common
subdivision of $\sK$ and $\widetilde{\sK}$, so that $\sK$ and
$\widetilde{\sK}$ are combinatorially equivalent. Note that a
$0$-move is the stellar subdivision at a maximal simplex (we
assume that the boundary of a 0-simplex is empty).

Bistellar $j$-moves with $i\ge\sbr d2$ are also called
\emph{reverse $(d-j)$-moves}. A $0$-move adds a new vertex to the
triangulation, a $d$-move (reverse $0$-move) deletes a vertex, and
all other bistellar moves do not change the number of vertices.
The bistellar moves in dimension 2 and 3 are shown in
Figures~\ref{bm2} and~\ref{bm3}. The bistellar 1-move in dimension
3 replaces two tetrahedra with a common face by 3 tetrahedra with
a common edge.

Two simplicial complexes are said to be \emph{bistellarly
equivalent}\label{bistequ} if one can be transformed to another by
a finite sequence of bistellar moves.
\end{definition}
\begin{figure}[h]
\begin{center}
\begin{picture}(90,60)
\put(0,15){\line(1,-1){15}} \put(15,0){\line(1,1){15}}
\put(30,15){\line(-1,1){15}} \put(15,30){\line(-1,-1){15}}
\put(15,0){\line(0,1){30}} \put(35,15){\vector(1,0){20}}
\put(55,15){\vector(-1,0){20}} \put(40,16){\small 1-move}
\put(60,15){\line(1,-1){15}} \put(75,0){\line(1,1){15}}
\put(90,15){\line(-1,1){15}} \put(75,30){\line(-1,-1){15}}
\put(60,15){\line(1,0){30}}
\put(0,40){\line(3,4){15}} \put(15,60){\line(3,-4){15}}
\put(0,40){\line(1,0){30}} \put(0,40){\line(3,1){15}}
\put(15,60){\line(0,-1){15}} \put(15,45){\line(3,-1){15}}
\put(55,52.5){\vector(-1,0){20}} \put(40,54){\small 0-move}
\put(35,45){\vector(1,0){20}} \put(40,42){\small 2-move}
\put(60,40){\line(3,4){15}} \put(75,60){\line(3,-4){15}}
\put(60,40){\line(1,0){30}}
\end{picture}
\caption{Bistellar moves for $q=2$.} \label{bm2}
\end{center}
\end{figure}
\begin{figure}[h]
\begin{picture}(125,90)
\put(0,25){\line(1,-1){25}} \put(25,0){\line(1,1){25}}
\put(50,25){\line(-1,1){25}} \put(25,50){\line(-1,-1){25}}
%\multiput(25,0)(0,6){8}{\line(0,1){3.5}}
\qbezier[30](25,50)(25,0)(25,0) \qbezier[30](25,0)(25,50)(25,50)
\put(25,48){\line(0,1){2}} \multiput(0,25)(5,0){10}{\line(1,0){2}}
\put(20,15){\line(3,1){30}} \put(20,15){\line(1,-3){5}}
\put(20,15){\line(-2,1){20}} \qbezier(20,15)(25,50)(25,50)
\put(55,20){\vector(1,0){20}} \put(60,17){\small 2-move}
\put(75,30){\vector(-1,0){20}} \put(60,31){\small 1-move}
\put(80,25){\line(1,-1){25}} \put(105,0){\line(1,1){25}}
\put(130,25){\line(-1,1){25}} \put(105,50){\line(-1,-1){25}}
\multiput(80,25)(5,0){10}{\line(1,0){2.5}}
\put(100,15){\line(3,1){30}} \put(100,15){\line(1,-3){5}}
\put(100,15){\line(-2,1){20}} \qbezier(100,15)(105,50)(105,50)
\put(25,60){\line(2,1){20}} \put(25,60){\line(0,1){30}}
\put(25,60){\line(-3,2){15}} \put(25,90){\line(1,-1){20}}
\put(25,90){\line(-3,-4){15}}
\multiput(10,70)(5,0){7}{\line(1,0){2.5}}
\qbezier[30](25,90)(27,73)(27,73)
\qbezier[30](45,70)(27,73)(27,73)
\qbezier[30](10,70)(27,73)(27,73)
\qbezier[25](25,60)(27,73)(27,73) \put(55,70){\vector(1,0){20}}
\put(60,67){\small 3-move} \put(75,80){\vector(-1,0){20}}
\put(60,81){\small 0-move} \put(100,60){\line(2,1){20}}
\put(100,60){\line(0,1){30}} \put(100,60){\line(-3,2){15}}
\put(100,90){\line(1,-1){20}} \put(100,90){\line(-3,-4){15}}
\multiput(85,70)(5,0){7}{\line(1,0){2.5}}
\end{picture}
\caption{Bistellar moves for $q=3$.} \label{bm3}
\end{figure}

%\begin{remark}
%The bistellar 0-move is just the stellar subdivision, or connected
%sum with the boundary of a simplex. In particular, \emph{stacked
%spheres}\label{stackeds} (i.e., the boundaries of stacked
%polytopes, see Definition~\ref{stacked}) are exactly those
%obtained from the boundary of a simplex by applying bistellar
%0-moves.
%\end{remark}

We have seen that bistellar equivalence implies $PL$
homeomorphism. The following result shows that for $PL$ manifolds
the converse is also true.

\begin{theorem}[{Pachner \cite[Theorem~1]{pach87}, \cite[(5.5)]{pach91}}]
\label{pacbs} Two $PL$ manifolds are bistellarly equivalent if and
only if they are $PL$ homeomorphic.
\end{theorem}

The behaviour of the face numbers of a triangulation under
bistellar moves is easily controlled. It can be most effectively
described in terms of the $g$-vector,
$g_i(\sK)=h_i(\sK)-h_{i-1}(\sK)$, \ $0<i\le\sbr d2$:

\begin{theorem}[Pachner~\cite{pach87}]
\label{gconjbs} If a triangulated $d$-manifold $\widetilde{\sK}$
is obtained from $\sK$ by a bistellar $k$-move, $0\le
k\le\sbr{d-1}2$, then
\begin{align*}
  g_{k+1}(\widetilde{\sK})&=g_{k+1}(\sK)+1;\\
  g_i(\widetilde{\sK})&=g_i(\sK)\quad\text{for all }i\ne k+1.
\end{align*}
Furthermore, if $d$ is even and $\widetilde{\sK}$ is obtained from
$\sK$ by a bistellar $\sbr d2$-move, then
\[
  g_i(\widetilde{\sK})=g_i(\sK)\quad\text{for all }i.
\]
\end{theorem}

This theorem allows us to interpret the inequalities from the
$g$-conjecture for $PL$ spheres (see Theorem~\ref{gth}) in terms
of the numbers of bistellar $k$-moves needed to transform a given
$PL$ sphere to the boundary of a simplex. For instance, the
inequality $h_1\le h_2$, \ $n\ge4$, is equivalent to the statement
that the number of 1-moves in the sequence of bistellar moves
taking a given $(n-1)$-dimensional $PL$ sphere to the boundary of
an $n$-simplex is less than or equal to the number of reverse
1-moves. (Note that the $g$-vector of $\partial\varDelta^n$ is
$(1,0,\ldots,0)$.)

\begin{remark}\label{el3}
There is also a generalisation of Theorem~\ref{pacbs} to $PL$
manifolds with boundary, see~\cite[(6.3)]{pach91}.
%, using a larger
%class of operations called \emph{elementary shellings}.
\end{remark}

In the case of polytopal sphere triangulations a stellar
subdivision is related to another familiar operation:

\begin{proposition}\label{stelltrunc}
Assume given a simple polytope $P$ and a proper face $G\subset P$.
Let $P^*$ be the dual simplicial polytope, $\sK_P=\partial P^*$
its nerve complex, and $J\subset P^*$ the face dual to~$G$. Then
the stellar subdivision $\ss_J\sK_P$ is the nerve complex of the
polytope $\widetilde P$ obtained by the face truncation at~$G$.
\end{proposition}
\begin{proof}
This follows directly by comparing the face poset of $\widetilde
P$, described in Construction~\ref{hypcut}, with that of
$\ss_J\sK_P$.
\end{proof}

\subsection*{Exercises.}
\begin{exercise}
The barycentric subdivision of $\sK$ can be obtained as a sequence
of stellar subdivisions at all faces $I\in\sK$, starting from the
maximal ones.
\end{exercise}

%\begin{exercise}
%What is the proper generalisation of Proposition~\ref{stelltrunc}
%to arbitrary convex polytopes?
%\end{exercise}
%
\begin{exercise}
Deduce formulae for the transformation of the $f$-, $h$- and
$g$-vector of $\sK$ under a stellar subdivision. Deduce similar
formulae for a bistellar move (the case of the $g$-vector is
Theorem~\ref{gconjbs}).
\end{exercise}

\section{Simplicial posets and simplicial cell complexes}\label{secsimpos}
Simplicial posets describe the combinatorial structures underlying
`generalised simplicial complexes' whose faces are still
simplices, but two faces are allowed to intersect in any
subcomplex of their boundary, rather than just in a single face.
These are also known as `ideal triangulations' in low-dimensional
topology, or as `simplicial cell complexes'.

\begin{definition}\label{defsp}
A poset (partially ordered set) $\mathcal S$ with order relation
$\le$ is called \emph{simplicial} if it has an initial element
$\hatzero$ and for each $\sigma\in\mathcal S$ the lower segment
\[
  [\hatzero,\sigma]=\{\tau\in\mathcal S\colon\hat0\le \tau\le \sigma\}
\]
is the face poset of a simplex. (The latter poset is also known as
a \emph{Boolean lattice}, and simplicial posets are sometimes
called \emph{Boolean posets}.) We assume all our posets to be
finite. The \emph{rank function} $|\cdot|$ on $\sS$ is defined by
setting $|\sigma|=k$ if $[\hatzero,\sigma]$ is the face poset of a
$(k-1)$-dimensional simplex. The rank of $\sS$ is the maximum of
ranks of its elements, and the \emph{dimension} of $\sS$ is its
rank minus one. A \emph{vertex} of $\sS$ is an element of rank
one. We assume that $\sS$ has $m$ vertices, denote the vertex set
by~$V(\sS)$, and usually identify it with $[m]=\{1,\ldots,m\}$.
Similarly, we denote by $V(\sigma)$ the vertex set of $\sigma$,
that is the set of $i$ with $i\le\sigma$.
\end{definition}

The face poset of a simplicial complex is a simplicial poset, but
there are many simplicial posets that do not arise in this way
(see Example~\ref{exsimpos} below). We identify a simplicial
complex with its face poset, thereby regarding simplicial
complexes as particular cases of simplicial posets.

To each $\sigma\in\mathcal S$ we assign a geometric simplex
$\varDelta^\sigma$ whose face poset is $[\hatzero,\sigma]$, and
glue these geometric simplices together according to the order
relation in~$\mathcal S$. As a result we get a regular cell
complex in which the closure of each cell is identified with a
simplex preserving the face structure, and all attaching and
characteristic maps are inclusions (see~\cite[Appendix]{hatc02}
for the terminology of cell complexes). We call it the
\emph{simplicial cell complex}\label{simpcellcomp} associated with
$\mathcal S$ and denote its underlying space by $|\mathcal S|$.

In the case when $\mathcal S$ is (the face poset of) a simplicial
complex $\mathcal K$ the space $|\mathcal S|$ is the geometric
realisation~$|\sK|$.

\begin{remark}
Using a more formal categorical language, we consider the
\emph{face category}\label{facecategor} $\cat(\sS)$ whose objects
are elements $\sigma\in\sS$ and there is a morphism from $\sigma$
to $\tau$ whenever $\sigma\le\tau$. Define a diagram (covariant
functor) $\varDelta^\sS$ from $\cat(\sS)$ to topological spaces by
sending $\sigma\in\sS$ to the geometric simplex $\varDelta^\sigma$
and sending every morphism $\sigma\le\tau$ to the inclusion
$\varDelta^\sigma\hookrightarrow\varDelta^\tau$. Then we may write
\[
  |\sS|=\colim \varDelta^\sS,
\]
where the colimit (or direct limit) is taken in the category of
topological spaces. This is the first example of colimit
construction over the face category~$\cat(\sS)$. Many other
examples of this sort will appear later.
\end{remark}

In most circumstances we shall not distinguish between simplicial
posets and simplicial cell complexes. We shall also sometimes
refer to elements $\sigma\in\mathcal S$ as \emph{simplices} or
\emph{faces} of~$\sS$\label{faceofsp}.

\begin{example}\label{exsimpos}
Consider the simplicial cell complex obtained by attaching two
$d$-dimensional simplices along their boundaries. Its
corresponding simplicial poset is not the face poset of a
simplicial complex if $d>0$.

Three cellular subdivisions of a circle are shown in
Fig.~\ref{cwcir}. The first is not a simplicial cell complex. The
second is a simplicial cell complex, but not a simplicial complex
(it corresponds to $d=1$ in the previous paragraph). The third one
is a simplicial complex.
\begin{figure}[h]
\begin{center}
\begin{picture}(80,25)
\put(10,15){\circle{14}} \put(3,15){\circle*{1.5}} \put(8,1){(1)}
\put(40,15){\circle{14}} \put(33,15){\circle*{1.5}}
\put(47,15){\circle*{1.5}} \put(38,1){(2)}
\put(70,15){\circle{14}} \put(65,10){\circle*{1.5}}
\put(65,20){\circle*{1.5}} \put(77,15){\circle*{1.5}}
\put(68,1){(3)}
\end{picture}
\caption{Cellular subdivisions of a circle.} \label{cwcir}
\end{center}
\end{figure}
\end{example}

\begin{construction}[folding a simplicial poset
onto a simplicial complex] For every simplicial poset $\sS$ there
is the associated simplicial complex $\sK_\sS$ on the same vertex
set~$V(\sS)$, whose simplices are sets $V(\sigma)$, \
$\sigma\in\sS$. There is a \emph{folding map} of simplicial posets
\begin{equation}\label{fold}
  \sS\longrightarrow\sK_\sS,\quad\sigma\mapsto V(\sigma).
\end{equation}
It is identical on the vertices, and every simplex in $\sK_\sS$
gets covered by a finite number of simplices of~$\sS$.
\end{construction}

For any two elements $\sigma,\tau\in\mathcal S$, denote by
$\sigma\vee\tau$ the set of their least common upper bounds
(\emph{joins}), and denote by $\sigma\wedge\tau$ the set of their
greatest common lower bounds (\emph{meets})\label{joinandmeet}.
Since $\mathcal S$ is a simplicial poset, $\sigma\wedge\tau$
consists of a single element whenever $\sigma\vee\tau$ is
nonempty. It is easy to observe that $\sS$ is a simplicial complex
if and only if for any $\sigma,\tau\in\sS$ the set
$\sigma\vee\tau$ is either empty or consists of a single element.
In this case $\sS$ coincides with $\sK_\sS$.

Applying barycentric subdivision to every simplex
$\sigma\in\sS$ we obtain a new simplicial cell complex~$\mathcal
S'$, called the \emph{barycentric subdivision} of~$\mathcal S$.
From Proposition~\ref{barorder} it is clear that $\mathcal S'$ can
be identified with the (geometric realisation of the) order
complex $\ord(\mathcal S\setminus\hat0)$. We therefore obtain the
following.

\begin{proposition}\label{baryc}
The barycentric subdivision $\mathcal S'$ of a simplicial cell
complex is a simplicial complex.
\end{proposition}

\subsection*{Exercises.}
\begin{exercise}
Show that the following conditions are equivalent for a simplicial
poset~$\sS$:
\begin{itemize}
\item[(a)] $\sS$ is (the face poset of) a simplicial complex;
\item[(b)] for any
$\sigma,\tau\in\sS$ the set $\sigma\wedge\tau$ consists of a
single element;
\item[(c)] for any
$\sigma,\tau\in\sS$ the set $\sigma\vee\tau$ is either empty or
consists of a single element. \end{itemize}
\end{exercise}

\section{Cubical complexes}
At some stage of the development of combinatorial topology,
cubical complexes were considered as an alternative to
triangulations, a new way to study topological invariants
combinatorially. Their nice feature is that the product of cubes
is again a cube, which makes subdivisions of products easier and
leads to a more straightforward definition of the multiplication
in cohomology. Later it turned out, however, that the cubical
(co)homomology itself is not particularly advantageous in
comparison with the simplicial one. Currently combinatorial
geometry is the main field of applications of cubical complexes;
moreover, subdivisions into cubes sometimes are very helpful in
various geometrical and topological considerations. In this
section we collect the necessary definitions and notation, and
then proceed to describe some important cubical decompositions
of simple polytopes and simplical complexes.

\subsection*{Definitions and examples}
As in the case of simplicial complexes, a cubical complex can be
defined either abstractly (as a poset) or geometrically (as a cell
complex).

\begin{definition}\label{cubcom}
An \emph{abstract cubical complex} is a finite poset $(\mathcal
C,\subset)$ containing an initial element $\varnothing$ and
satisfying the following two conditions:
\begin{enumerate}
\item[(a)] for every element $G\in\mathcal C$ the segment
$[\varnothing,G]$ is isomorphic to the face poset of a cube;

\item[(b)] for every two elements $G_1,G_2\in\mathcal C$ there is a
unique meet (greatest lower bound).
%, that is, among all $H\in\mathcal C$ such that $H\subset G_1$
%and $H\subset G_2$ there is a unique maximal element.
\end{enumerate}
Elements $G\in\mathcal C$ are \emph{faces} of the cubical complex.
If $[\varnothing,G]$ is the face poset of the $k$-cube $\I^k$,
then the face $G$ is of \emph{$k$-dimensional}. The dimension of
$\mathcal C$ is the maximal dimension of its faces. The meet of
any two faces $G_1$, $G_2$ is also called their
\emph{intersection} and denoted $G_1\cap G_2$.
\end{definition}

A $d$-dimensional \emph{topological cube} is a $d$-ball with a
face structure defined by a homeomorphism with the standard
cube~$\I^d$. A \emph{face} of a topological $d$-cube is thus the
homeomorphic image of a face of~$\I^d$.

\begin{definition}
\label{topcubcom} A \emph{topological cubical complex} is a set
$\mathcal U$ of topological cubes of arbitrary dimensions which
are all embedded in the same space $\R^n$ and satisfy the
following conditions:
\begin{enumerate}
\item[(a)] every face of a cube in $\mathcal U$ belongs to $\mathcal U$;
\item[(b)] The intersection of any two cubes in $\mathcal U$ is a face of
each.
\end{enumerate}
\end{definition}

Every abstract cubical complex $\mathcal C$ has a \emph{geometric
realisation}\label{geomrelcc}, a topological cubical complex
$\mathcal U$ whose faces form a poset isomorphic to~$\mathcal C$.
Such $\mathcal U$ can be constructed by taking a disjoint union of
topological cubes corresponding to all segments $[\varnothing,
G]\subset\mathcal C$ and identifying faces according to the poset
relation.

From now on we shall not distinguish between abstract cubical
complexes and their geometric realisations.

By analogy with simplicial complexes, define the
\emph{$f$-vector}\label{fvectorcc} of a cubical complex $\mathcal
C$ by $\mb f(\mathcal C)=(f_0,f_1,\ldots)$ where $f_i$ is the
number of $i$-dimensional faces.
%The theory of $f$-vectors of cubical complexes is parallel
%to that of simplicial complexes, but is less developed.
There are also notions of $h$- and $g$-vectors, and cubical
analogues of the UBC, LBC and $g$-conjecture. See~\cite{adin96},
\cite{b-b-c97} and~\cite{stan00} for more details and references.

The difference between Definition \ref{topcubcom} of a geometric
cubical complex and Definition~\ref{polyhed} of a geometric
simplicial complex is that we realise abstract cubes by
topological complexes rather than polytopes. This difference is
substantial: if we replace topological cubes by combinatorial ones
(i.e. by convex polytopes combinatorially equivalent to a cube) in
Definition~\ref{topcubcom}, then we obtain the definition of a
\emph{polyhedral cubical complex}\label{polyhedralcc}. Although
this notion is also important in combinatorial geometry, not every
abstract cubical complex can be realised by a polyhedral complex,
as shown by the next example.

\begin{example}
Consider the decomposition of a M\"obius strip into 3 squares
shown in Fig.~\ref{cubmo}.
\begin{figure}[h]
\begin{center}
\begin{picture}(70,30)
\put(5,5){\line(1,0){60}} \put(5,5){\vector(0,1){20}}
\put(5,25){\line(1,0){60}} \put(25,5){\line(0,1){20}}
\put(45,5){\line(0,1){20}} \put(65,25){\vector(0,-1){20}}
\put(2,1){$A$} \put(23,1){$F$} \put(43,1){$E$} \put(65,1){$B$}
\put(2,26){$B$} \put(23,26){$C$} \put(43,26){$D$} \put(65,26){$A$}
\end{picture}
\caption{Cubical complex which does not admit a polyhedral
realisation.} \label{cubmo}
\end{center}
\end{figure}

\begin{proposition}
The topological cubical complex shown in Fig.~\ref{cubmo} does not
admit a polyhedral cubical realisation.
\end{proposition}
\begin{proof}
Assume that such a realisation exists. Then since $ABED$ is a
convex 4-gon, the points $A$ and $D$ are in the same halfplane
with respect to the line~$BE$, and therefore $A$ and $D$ are in
the same halfspace defined by the plane~$BCE$. Similarly, since
$ABCF$ is a convex 4-gon, the points $A$ and $F$ are in the same
halfspace with respect to~$BCE$. Hence, $D$ and $F$ are also in
the same halfspace with respect to $BCE$. On the other hand, since
$CDEF$ is a convex 4-gon, the points $D$ and $F$ must be in
different subspaces with respect to $BCE$. A contradiction.
\end{proof}
\end{example}

The example above shows that, unlike the case of simplicial
complexes, the theory of abstract cubical complexes cannot be
described by using only convex-geometric considerations. Another
simple manifestation of this is the fact that not every cubical
complex may be realised as a subcomplex in a standard cube, in
contrast to simplicial complexes which are always embeddable in a
standard simplex. The boundary of a triangle is the simplest
example of a cubical complex not embeddable in a cube. It is also
not embeddable into the standard cubical lattice in $\R^n$ (for
any~$n$). On the other hand, every cubical complex admits a
cubical subdivision which is embeddable in a standard cube, as
shown in the next subsection. Without subdivision the question of
embeddability in a standard cube or cubical lattice is nontrivial.
The importance of studying cubical maps (in particular, cubical
embeddings) of 2-dimensional cubical complexes to the cubical
lattice in $\R^3$ was pointed out by S.~Novikov in connection with
the \emph{$3$-dimensional Ising model}.
%A significant advance on this problem has been achieved
In~\cite{d-s-s94} necessary and sufficient conditions were
obtained for a cubical complex to admit a cubical map to the
standard lattice.

\subsection*{Cubical subdivisions of simple polytopes and simplicial complexes}
The particular constructions of cubical complexes given here will
be important in the definition of moment-angle complexes. Neither
of these constructions is particularly new, but they are probably
not well recorded in the literature (see however the references at
the end of the section).

Any face of $\I^m$ has the form
\[
  C_{J\subset I}=\{(y_1,\ldots,y_m)\in \I^m\colon y_j=0
  \text{ for }j\in
  J,\quad y_j=1\text{ for }j\notin I\}
\]
where $J\subset I$ is a pair of embedded (possibly empty) subsets
of~$[m]$. We also set
\begin{equation}\label{Iface}
  C_I=C_{\varnothing\subset I}=
  \{(y_1,\ldots,y_m)\in \I^m\colon y_j=1\text{ for }j\notin I\}
\end{equation}
to simplify the notation.

\begin{construction}[canonical triangulation of~$\I^m$]\label{conbar}
Denote by $\varDelta=\varDelta^{m-1}$ the simplex on~$[m]$. We
assign to a subset $I=\{i_1,\ldots,i_k\}\subset[m]$ the vertex
$v_I=C_{I\subset I}$ of~$\I^m$. That is,
$v_I=(\varepsilon_1,\ldots,\varepsilon_m)$ where $\varepsilon_i=0$
if $i\in I$ and $\varepsilon_i=1$ otherwise. Regarding each $I$ as
a vertex of the barycentric subdivision of~$\varDelta$, we can
extend the correspondence $I\mapsto v_I$ to a piecewise linear
embedding $i_c\colon\varDelta'\to \I^m$. Under this embedding the
vertices of $\varDelta$ are mapped to the vertices of $\I^m$ with
exactly one zero coordinate, and the barycentre of $\varDelta$ is
mapped to $(0,\ldots,0)\in \I^m$ (see Fig.~\ref{barcub}).
\begin{figure}[h]
\begin{picture}(120,145)
  %cube
  \multiput(80,10)(0.1,0){7}{\line(3,2){30}}
  \multiput(80,60)(0.1,0){7}{\line(3,2){30}}
  \multiput(30,60)(0.1,0){7}{\line(3,2){30}}
  \multiput(60,30)(-5,-3.33){6}{\line(-3,-2){3.6}}
  \multiput(60,30)(5,0){10}{\line(1,0){3}}
  \multiput(60,30)(0,5){10}{\line(0,1){3}}
  \qbezier[70](60,30)(30,60)(30,60)
  \qbezier[50](61,30)(80,10)(80,10)
  \qbezier[60](60,30)(80,60)(80,60)
  \qbezier[90](60,30)(110,80)(110,80)
  \put(80,60){\line(-1,-1){50}}
  \put(80,60){\line(-1,1){20}}
  \put(80,59.5){\line(1,-1){30}}
  {\linethickness{0.4mm}
  \put(30,10){\line(1,0){50}}
  \put(30,10){\line(0,1){50}}
  \put(30,60){\line(1,0){50}}
  \put(80,60){\line(0,-1){50}}
  \put(110,30){\line(0,1){50}}
  \put(60,80){\line(1,0){50}}
  \put(110,30){\line(0,1){50}}
  }
  \put(75,66){\footnotesize vertex of}
  \put(78.5,63.5){\footnotesize cone}
  \put(111,81){$i_c(1)$}
  \put(22,56){$i_c(2)$}
  \put(78,5){$i_c(3)$}
  \put(25,5){$i_c(23)$}
  \put(111,28){$i_c(13)$}
  \put(60,81){$i_c(12)$}
  \put(55,24){\small $i_c(123)$}
  %simplex
  \put(0,140){\line(1,0){60}}
  \put(0,140){\line(1,-2){30}}
  \put(0,140){\line(3,-2){45}}
  \put(30,80){\line(0,1){60}}
  \put(30,80){\line(1,2){30}}
  \put(60,140){\line(-3,-2){45}}
  \put(59,141){1}
  \put(0,141){2}
  \put(27,79){3}
  \put(29,141){12}
  \put(11,108){23}
  \put(45.5,108){13}
  \put(22,119){\small 123}
  \put(40,90){\vector(1,-1){10}}
  \put(46,86){$i_c$}
\end{picture}
\caption{Taking cone over the barycentric subdivision of simplex
defines a triangulation of the cube.} \label{barcub}
\end{figure}
The image $i_c(\varDelta')$ is the union of $m$ facets of~$\I^m$
meeting at the vertex $(0,\ldots,0)$. For each pair $I\subset J$,
all simplices of $\varDelta'$ of the form $I=I_1\subset
I_2\subset\cdots\subset I_k= J$ are mapped to the same face $C_{
I\subset J}$ of~$\I^m$. The map $i_c\colon\varDelta'\to \I^m$
extends to $\cone(\varDelta')$ by mapping the cone vertex to
$(1,\ldots,1)\in \I^m$. The image of the resulting map
$\cone(i_c)$ is the whole cube~$\I^m$. Thus, $\cone(i_c)\colon
\cone(\varDelta')\to \I^m$ is a $PL$ homeomorphism which is linear
on simplices of $\cone(\varDelta')$. This defines a canonical
triangulation of $\I^m$, the `triangulation along the main
diagonal'.
\end{construction}

The subdivisions which appear above can be summarised as follows:

\begin{proposition}
The $PL$ map $\cone(i_c)\colon \cone(\varDelta')\to \I^m$ gives
rise to
\begin{itemize}
\item[(a)]
a cubical subdivision of $\varDelta^{m-1}$ isomorphic to ``half of
the boundary of~$\I^m$'' (the union of facets of $\I^m$ containing
the zero vertex);

\item[(b)]
a cubical subdivision of $\cone\varDelta^{m-1}$ (which
is~$\varDelta^m$) isomorphic to~$\I^m$;

\item[(c)]
a simplicial subdivision of $\I^m$ isomorphic to
$\cone((\varDelta^{m-1})')$.
\end{itemize}
\end{proposition}

\begin{construction}[cubical subdivision of a simple polytope]\label{cubpol}
Let $P$ be a simple $n$-polytope with $m$ facets $F_1,\ldots,F_m$.
We shall construct a piecewise linear embedding of $P$ into the
standard cube~$\I^m$, thereby inducing a cubical subdivision
$\mathcal C(P)$ of $P$ by the preimages of faces of~$\I^m$.

Denote by $\mathcal S$ the set of barycentres of faces of~$P$,
including the vertices and the barycentre of the whole polytope.
This will be the vertex set of~$\mathcal C(P)$. Every $(n-k)$-face
$G$ of $P$ is an intersection of $k$ facets:
$G=F_{i_1}\cap\cdots\cap F_{i_k}$. We map the barycentre of $G$ to
the vertex $(\varepsilon_1,\ldots,\varepsilon_m)\in \I^m$, where
$\varepsilon_i=0$ if $i\in\{i_1,\ldots,i_k\}$ and
$\varepsilon_i=1$ otherwise. The resulting map $\mathcal S\to\I^m$
can be extended linearly on the simplices of the barycentric
subdivision of $P$ to an embedding $c_P\colon P\to\I^m$. The case
$n=2$, $m=3$ is shown in Fig.~\ref{ip}.
\begin{figure}[h]
\begin{picture}(120,60)
  \put(10,10){\line(1,2){20}}
  \put(30,50){\line(1,-2){20}}
  \put(10,10){\line(1,0){40}}
  \put(20,30){\line(2,-1){10}}
  \put(40,30){\line(-2,-1){10}}
  \put(30,10){\line(0,1){15}}
  \put(50,30){\vector(1,0){12.5}}
  \put(70,10){\line(1,0){30}}
  \put(70,10){\line(0,1){30}}
  \put(70,40){\line(1,0){30}}
  \put(100,10){\line(0,1){30}}
  \put(70,40){\line(1,1){10}}
  \put(100,40){\line(1,1){10}}
  \put(100,10){\line(1,1){10}}
  \put(80,50){\line(1,0){30}}
  \put(110,20){\line(0,1){30}}
  \put(80,20){\line(-1,-1){3.6}}
  \put(75,15){\line(-1,-1){3.6}}
  \multiput(80,20)(3,0){10}{\line(1,0){2}}
  \multiput(80,20)(0,3){10}{\line(0,1){2}}
  \put(81,21){0}
  \put(7,5){$A$}
  \put(28,5){$F$}
  \put(50,5){$E$}
  \put(32,22){$G$}
  \put(15,29){$B$}
  \put(42,29){$D$}
  \put(29,52){$C$}
  \put(15,52){\Large $P$}
  \put(54,32){$c_P$}
  \put(67,5){$A$}
  \put(98,5){$F$}
  \put(67,41){$B$}
  \put(101,37){$G$}
  \put(80,52){$C$}
  \put(111,52){$D$}
  \put(111,18){$E$}
  \put(66,52){\Large $\I^m$}
\end{picture}
\caption{Embedding $c_P\colon P\to \I^m$ for $n=2$, $m=3$.}
\label{ip}
\end{figure}

The image $c_P(P)\subset\I^m$ is the union of all faces
$C_{J\subset I}$ such that $\bigcap_{i\in I}F_i\ne\varnothing$.
For such $C_{J\subset I}$, the preimage $c_P^{-1}(C_{J\subset I})$
is declared to be a face of the cubical complex~$\mathcal C(P)$.
The vertex set of $c_P^{-1}(C_{J\subset I})$ is the subset of
$\mathcal S$ consisting of barycentres of all faces between the
faces $G$ and $H$ of~$P$, where $G=\bigcap_{j\in J}F_j$ and
$H=\bigcap_{i\in I}F_i$. Therefore, faces of $\mathcal C(P)$
correspond to pairs of embedded faces $G\supset H$ of~$P$, and we
denote them by $C_{G\supset H}$. In particular, maximal
($n$-dimensional) faces of $\mathcal C(P)$ correspond to pairs
$G=P$, $H=v$, where $v$ is a vertex of~$P$. For these maximal
faces we use the abbreviated notation $C_v=C_{P\supset v}$.

For every vertex $v=F_{i_1}\cap\cdots\cap F_{i_n}\in P$ with
$I_v=\{i_1,\ldots,i_n\}$ we have
\begin{equation}
\label{cubpolmap}
  c_P(C_v)=C_{I_v}=\bigl\{(y_1,\ldots,y_m)\in \I^m\colon y_j=1\text{ whenever }
  v\notin F_j\bigr\}.
\end{equation}
\end{construction}

We therefore obtain:

\begin{proposition}
\label{thcubpol} A simple polytope $P$ with $m$ facets admits a
cubical decomposition whose maximal faces $C_v$ correspond to the
vertices $v\in P$. The resulting cubical complex $\mathcal C(P)$
embeds canonically into~$\I^m$, as described
by~{\rm(\ref{cubpolmap})}.
\end{proposition}

\begin{lemma}
The number of $k$-faces of the cubical complex $\mathcal C(P)$ is
given by
\[
  f_k\bigl(\mathcal C(P)\bigr)
  =\sum_{i=0}^{n-k}\binom{n-i}k f_i(P),
%  =\bin nk f_0(P)+\bin{n-1}k f_1(P)+\cdots+f_{n-k}(P),
  \quad\text{for } 0\le k\le n.
\]
\end{lemma}
\begin{proof}
The formula follows from the fact that the $k$-faces of $\mathcal
C(P)$ are in one-to-one correspondence with the pairs
$G^{i+k}\supset H^i$ of faces of~$P$.
\end{proof}

\begin{construction}[cubical subdivision of a simplicial complex]\label{cubk}
Let $\mathcal K$ be a simplicial complex on~$[m]$. Then $\mathcal
K$ is naturally a subcomplex of~$\varDelta^{m-1}$ and its
barycentric subdivision $\mathcal K'$ is a subcomplex
of~$(\varDelta^{m-1})'$. Restricting the $PL$ map from
Construction~\ref{conbar} to $\mathcal K'$, we obtain the
embedding $i_c|_{\mathcal K'}\colon|\mathcal K'|\to \I^m$. Its
image is a cubical subcomplex in~$\I^m$, which we
denote~$\cub(\mathcal K)$. Then $\cub(\mathcal K)$ is the union of
faces $C_{ I\subset J}\subset \I^m$ over all pairs $I\subset J$ of
nonempty simplices of~$\sK$:
\begin{equation}
\label{fcubk}
  \cub(\mathcal K)=\bigcup_{\varnothing\ne I\subset J\in
  \mathcal K}C_{ I \subset
   J}\subset \I^m.
\end{equation}
\end{construction}

\begin{construction}\label{cck}
Since $\cone(\mathcal K')$ is a subcomplex of
$\cone((\varDelta^{m-1})')$, Construction~\ref{conbar} also
provides a $PL$ embedding
$$
  \cone(i_c)|_{\cone(\mathcal K')}\colon|\cone(\mathcal K')|\to \I^m.
$$
The image of this embedding is an $n$-dimensional cubical
subcomplex of~$\I^m$, which we denote $\cc(\mathcal K)$. It is
easy to see that
\begin{equation}
\label{fcck}
  \cc(\mathcal K)=
  \bigcup_{I\subset J\in \mathcal K}C_{ I\subset J}=
  \bigcup_{ J\in \mathcal K}C_ J.
\end{equation}
\end{construction}

\begin{remark}
If $i\in[m]$ is not a vertex of $\mathcal K$ (a ghost vertex),
then $\cc(\mathcal K)$ is contained in the facet $\{y_i=1\}$
of~$\I^m$.
\end{remark}

Here is a summary of the two previous constructions.

\begin{proposition}\label{cubkcck}
For any simplicial complex $\sK$ on the set $[m]$, there is a $PL$
embedding of $|\sK|$ into $\I^m$ linear on the simplices
of~$\sK'$. The  image of this embedding is the cubical
subcomplex~{\rm(\ref{fcubk})}. Moreover, there is a $PL$ embedding
of $|\cone\sK|$ into $\I^m$ linear on the simplices of
$\cone(\sK')$, whose image is the cubical
subcomplex~{\rm(\ref{fcck})}.
\end{proposition}

A cubical complex $\mathcal C'$ is called a \emph{cubical
subdivision}\label{cubicsubd} of a cubical complex $\mathcal C$ if
each face of $\mathcal C'$ is contained in a face of~$\mathcal C$,
and each face of $\mathcal C$ is a union of finitely many faces
of~$\mathcal C'$.

\begin{proposition}
For every cubical complex $\mathcal C$ with $q$ vertices, there
exists a cubical subdivision that is embeddable as a subcomplex
in~$\I^q$.
\end{proposition}
\begin{proof}
We first construct a simplicial complex $\sK_{\mathcal C}$ which
subdivides the cubical complex $\mathcal C$ and has the same
vertices. This can be done by induction on the skeleta
of~$\mathcal C$, by extending the triangulation from the
$k$-dimensional skeleton to the interiors of $(k+1)$-dimensional
faces using a generic convex function $f\colon\R^k\to\R$,
see~\cite[\S5.1]{zieg95} (note that the 1-skeleton of $\mathcal C$
is already a simplicial complex). Then applying
Construction~\ref{cubk} to $\sK_{\mathcal C}$ we get cubical
complex $\cub(\sK_{\mathcal C})$ that subdivides $\sK_{\mathcal
C}$ and therefore~$\mathcal C$. It is embeddable into $\I^q$  by
Proposition~\ref{cubkcck}.
\end{proof}

\begin{example}
\begin{figure}[h]
  \begin{picture}(120,45)
  \put(15,5){\line(1,0){25}}
  \put(15,5){\line(0,1){25}}
  \put(40,5){\line(1,1){10}}
  \put(50,15){\line(0,1){25}}
  \put(50,40){\line(-1,0){25}}
  \put(25,40){\line(-1,-1){10}}
  \put(40,5){\circle*{2}}
  \put(15,30){\circle*{2}}
  \put(50,40){\circle*{2}}
  \put(40,30){\line(1,1){10}}
  \multiput(25,15)(5,0){5}{\line(1,0){3}}
  \multiput(25,15)(0,5){5}{\line(0,1){3}}
  \multiput(25,15)(-5,-5){2}{\line(-1,-1){4}}
  \put(26,16){$0$}
  \put(22,-2){(a)\ $\sK=\text{3 points}$}
  \put(75,5){\line(1,0){25}}
  \put(75,5){\line(0,1){25}}
  \put(100,5){\line(1,1){10}}
  \put(110,15){\line(0,1){25}}
  \put(110,40){\line(-1,0){25}}
  \put(85,40){\line(-1,-1){10}}
  \put(100,5){\circle*{2}}
  \put(110,15){\circle*{2}}
  \put(75,30){\circle*{2}}
  \put(85,40){\circle*{2}}
  \put(75,5){\circle*{2}}
  \put(110,40){\circle*{2}}
  \put(100,30){\line(1,1){10}}
  \multiput(75,29.3)(0,0.1){16}{\line(1,1){10}}
  \multiput(100,4.3)(0,0.1){16}{\line(1,1){10}}
  \multiput(85,15)(5,0){5}{\line(1,0){3}}
  \multiput(85,15)(0,5){5}{\line(0,1){3}}
  \multiput(85,15)(-5,-5){2}{\line(-1,-1){4}}
  \put(86,16){$0$}
  \put(82,-2){(b)\ $\sK=\partial\varDelta^2$}
  \put(15,30){\line(1,0){25}}
  \put(40,5){\line(0,1){25}}
  \put(75,30){\line(1,0){25}}
  \put(100,5){\line(0,1){25}}
  \linethickness{1mm}
  \put(75,5){\line(1,0){25}}
  \put(85,40){\line(1,0){25}}
  \put(75,5){\line(0,1){25}}
  \put(110,15){\line(0,1){25}}
  \end{picture}
  \caption{Cubical complex $\cub(\sK)$.}
  \label{figcubk}
  \end{figure}
  \begin{figure}[h]
  \begin{picture}(120,45)
  \put(15,5){\line(1,0){25}}
  \put(15,5){\line(0,1){25}}
  \put(40,5){\line(1,1){10}}
  \put(50,15){\line(0,1){25}}
  \put(50,40){\line(-1,0){25}}
  \put(25,40){\line(-1,-1){10}}
  \put(40,5){\circle*{2}}
  \put(40,30){\circle*{2}}
  \put(15,30){\circle*{2}}
  \put(50,40){\circle*{2}}
  \multiput(40,29.3)(0,0.1){16}{\line(1,1){10}}
  \multiput(25,15)(5,0){5}{\line(1,0){3}}
  \multiput(25,15)(0,5){5}{\line(0,1){3}}
  \multiput(25,15)(-5,-5){2}{\line(-1,-1){4}}
  \put(26,16){$0$}
  \put(22,-2){(a)\ $\sK=\text{3 points}$}
  \put(75,5){\line(1,0){25}}
  \put(75,5){\line(0,1){25}}
  \put(100,5){\line(1,1){10}}
  \put(110,15){\line(0,1){25}}
  \put(110,40){\line(-1,0){25}}
  \put(85,40){\line(-1,-1){10}}
  \put(100,5){\circle*{2}}
  \put(110,15){\circle*{2}}
  \put(100,30){\circle*{2}}
  \put(75,30){\circle*{2}}
  \put(85,40){\circle*{2}}
  \put(75,5){\circle*{2}}
  \put(110,40){\circle*{2}}
  \multiput(100,29.3)(0,0.1){16}{\line(1,1){10}}
  \multiput(75,29.3)(0,0.1){16}{\line(1,1){10}}
  \multiput(100,4.3)(0,0.1){16}{\line(1,1){10}}
  \multiput(85,15)(5,0){5}{\line(1,0){3}}
  \multiput(85,15)(0,5){5}{\line(0,1){3}}
  \multiput(85,15)(-5,-5){2}{\line(-1,-1){4}}
  \multiput(78,5)(4,0){6}{\line(0,1){25}}
  \multiput(78,30)(4,0){6}{\line(1,1){10}}
  \multiput(100,8)(0,4){6}{\line(1,1){10}}
  \put(86,16){$0$}
  \put(82,-2){(b)\ $\sK=\partial\varDelta^2$}
  \linethickness{1mm}
  \put(15,30){\line(1,0){25}}
  \put(40,5){\line(0,1){25}}
  \put(75,30){\line(1,0){25}}
  \put(75,5){\line(1,0){25}}
  \put(85,40){\line(1,0){25}}
  \put(100,5){\line(0,1){25}}
  \put(75,5){\line(0,1){25}}
  \put(110,15){\line(0,1){25}}
  \end{picture}
  \caption{Cubical complex $\cc(\sK)$.}
  \label{figcck}
  \end{figure}
The cubical complexes $\cub(\sK)$ and $\cc(\sK)$ when $\sK$ is a
disjoint union of 3 vertices or the boundary of a triangle are
shown in Figs.~\ref{figcubk}--\ref{figcck}.
\end{example}

%\begin{remark}
%As a topological space, $\cub(\sK)$ is homeomorphic to $|\sK|$,
%while $\cc(\sK)$ is homeomorphic to $|\cone(\sK)|$. On the other
%hand, cubical complex $\cub(\cone(\sK))$ is homeomorphic to
%$|\cone(\sK)|$. However, as cubical complexes $\cc(K)$ and
%$\cub(\cone(K))$ differ (since $\cone(K')\ne(\cone(K))'$).
%\end{remark}

\begin{remark}
Let $P$ be a simple $n$-polytope, and let $\sK_P$ be its nerve
complex. Then $\cc(\sK_P)=c_P(P)$, i.e. $\cc(\sK_P)$ coincides
with the cubical complex $\mathcal C(P)$ from
Construction~\ref{cubpol}.
%Therefore, Construction~\ref{cubpol} is
%a particular case of Construction~\ref{cck}.
\end{remark}

Different versions of Construction~\ref{cck} can be found
in~\cite[p.~299]{b-b-c97}. A similar construction was also
considered in~\cite[p.~434]{da-ja91}. Finally, a version of
cubical subcomplex $\cub(\sK)\subset\I^m$ appeared
in~\cite{d-s-s94} in connection with the problems described in the
end of the previous subsection.

\subsection*{Exercises.}
\begin{exercise}
Show that the triangulation of $\I^m$ from
Construction~\ref{conbar} coincides with the triangulation of the
product of $m$ one-dimensional simplices from
Construction~\ref{simprod}.
\end{exercise}

%\chapter*{\ \ Combinatorial structures: additional topics}
%\section{Minimal triangulations of manifolds}

f%Проверить, нужно ли условие алгебраической замкнутости k
%в теореме Брунса-Губеладзе.
%
%Доказательство Айзенберга теоремы Райснера - в "Additional topics"

\setcounter{chapter}2
\chapter{Combinatorial algebra of face rings}\label{facerings}
In this chapter we collect the wealth of algebraic notions and
constructions related to face rings. Our choice of material and
notation was guided by the topological applications in the later
chapters of the book. (This explains the unusual for algebraists
even grading in the polynomial rings and their homogeneous
quotients, and also the nonpositive homological grading in free
resolutions and $\Tor$.) In the first sections we review standard
results and constructions of combinatorial commutative algebra,
including the $\Tor$-algebras and algebraic Betti numbers of face
rings, Cohen--Macaulay and Gorenstein complexes. The later
sections contain some more recent developments, including the face
rings of simplicial posets, different characterisations of
Cohen--Macaulay and Gorenstein simplicial posets in terms of their
face rings and $h$-vectors, and generalisations of the
Dehn--Sommerville relations. Although all these algebraic and
combinatorial results have a strong topological flavour and were
indeed originally motivated by topological constructions, we have
tried to keep this chapter mostly algebraic and do not require
much topological knowledge from the reader here.

The preliminary algebraic material of a more general sort, not
directly or exclusively related to the face rings (such as
resolutions and the functor $\Tor$, and Cohen--Macaulay rings) is
collected in Appendix~\ref{hab}.

Alongside the monograph by Stanley~\cite{stan96}, an extensive
survey of Cohen--Macaulay rings by Bruns and Herzog~\cite{br-he98}
and a more recent monograph~\cite{mi-st05} by Miller and Sturmfels
may be recommended for a deeper study of algebraic methods in
combinatorics.

\medskip

We use the common notation $\k$ for the ground ring, which is
always assumed to be the ring~$\Z$ of integers or a field. The
former is preferable for topological applications, but the latter
is more common in the algebraic literature.
%Our algebras $A=\bigoplus_{i\ge0}A^i$ will be graded commutative
%and finitely generated over~$\k$, and all $A$-modules $M$ will be
%graded and finitely generated. All gradings will be nonnegative,
%unless otherwise stated. We always assume $A$ to be connected,
%i.e. $A^0=\k$. Then the augmentation map $A\to\k$, which is
%identical on $A^0$ and takes the positive part $A^+$ to zero,
%determines an $A$-module structure on~$\k$.
We shall often refer to $\k$-modules as `$\k$-vector spaces'; in
the case $\k=\Z$ the latter means `abelian groups'.

We assume graded commutativity instead of commutativity; algebras
commutative in the standard sense will be those whose nontrivial
graded components appear only in even degrees. In particular, the
polynomial algebra $\k[v_1,\ldots,v_m]$, which we often abbreviate
to $\k[m]$, has $\deg v_i=2$. The exterior algebra
$\Lambda[u_1,\ldots,u_m]$ has $\deg u_i=1$. Given a subset
$I=\{i_1,\ldots,i_k\}\subset[m]$ we denote by $v_I$ the
square-free monomial $v_{i_1}\cdots v_{i_k}$ in~$\k[m]$. We also
denote by $u_I$ the exterior monomial $u_{i_1}\cdots u_{i_k}$
where $i_1<\cdots<i_k$.

\section{Face rings of simplicial complexes}
\label{srr}
\begin{definition}
\label{frsim} The \emph{face ring} (or the \emph{Stanley--Reisner
ring}) of a simplicial complex $\mathcal K$ on the set $[m]$ is
the quotient graded ring
$$
  \k[\mathcal K]=\k[v_1,\ldots,v_m]/\mathcal I_\mathcal K,
$$
where $\mathcal I_\mathcal K=(v_I\colon I\notin\sK)$ is the
%homogeneous
ideal generated by those monomials $v_I$ for which $I$
is not a simplex of~$\mathcal K$. The ideal $\mathcal I_\mathcal
K$ is known as the \emph{Stanley--Reisner ideal} of~$\mathcal K$.
\end{definition}

\begin{example}\

1. Let $\mathcal K$ be the 2-dimensional simplicial complex shown
in Fig.~\ref{figfr}. Then
$$
  \mathcal I_\mathcal K=(v_1v_5,v_3v_4,v_1v_2v_3,v_2v_4v_5).
$$
\begin{figure}[h]
\begin{center}
\begin{picture}(60,35)
{\linethickness{0.5mm}
\put(5,5){\line(1,0){50}}
\put(25,15){\line(1,0){10}}
}
\multiput(5,4.9)(0,0.05){10}{\line(1,1){25}}
\multiput(5,4.9)(0,0.05){9}{\line(2,1){20}}
\multiput(25,14.5)(0,0.1){10}{\line(1,3){5}}
\multiput(35,14.5)(0,0.1){10}{\line(-1,3){5}}
\multiput(30,29.9)(0,0.05){10}{\line(1,-1){25}}
\multiput(35,14.9)(0,0.05){9}{\line(2,-1){20}}
\put(7.5,7.5){\line(0,-1){1.25}}
\put(8.75,8.75){\line(0,-1){1.9}}
\put(10,10){\line(0,-1){2.5}}
\put(11.25,11.25){\line(0,-1){3.1}}
\put(12.5,12.5){\line(0,-1){3.75}}
\put(13.75,13.75){\line(0,-1){4.4}}
\put(15,15){\line(0,-1){5}}
\put(16.25,16.25){\line(0,-1){5.7}}
\put(17.5,17.5){\line(0,-1){6.25}}
\put(18.75,18.75){\line(0,-1){6.9}}
\put(20,20){\line(0,-1){7.5}}
\put(21.25,21.25){\line(0,-1){8.1}}
\put(22.5,22.5){\line(0,-1){8.75}}
\put(23.75,23.75){\line(0,-1){9.4}}
\put(25,25){\line(0,-1){10}}
\put(26.25,26.25){\line(0,-1){7.5}}
\put(27.5,27.5){\line(0,-1){5}}
\put(28.75,28.75){\line(0,-1){2.5}}
\put(31.25,28.75){\line(0,-1){2.5}}
\put(32.5,27.5){\line(0,-1){5}}
\put(33.75,26.25){\line(0,-1){7.5}}
\put(35,25){\line(0,-1){10}}
\put(36.25,23.75){\line(0,-1){9.4}}
\put(37.5,22.5){\line(0,-1){8.75}}
\put(38.75,21.25){\line(0,-1){8.1}}
\put(40,20){\line(0,-1){7.5}}
\put(41.25,18.75){\line(0,-1){6.9}}
\put(42.5,17.5){\line(0,-1){6.25}}
\put(43.75,16.25){\line(0,-1){5.7}}
\put(45,15){\line(0,-1){5}}
\put(46.25,13.75){\line(0,-1){4.4}}
\put(47.5,12.5){\line(0,-1){3.75}}
\put(48.75,11.25){\line(0,-1){3.1}}
\put(50,10){\line(0,-1){2.5}}
\put(51.25,8.75){\line(0,-1){1.9}}
\put(52.5,7.5){\line(0,-1){1.25}}
\put(5,5){\circle*{1}} \put(25,15){\circle*{1}}
\put(35,15){\circle*{1}} \put(30,30){\circle*{1}}
\put(55,5){\circle*{1}} \put(3,1){1} \put(56,1){3} \put(29,31){2}
\put(24,11){4} \put(34,11){5}
\end{picture}
\caption{} \label{figfr}
\end{center}
\end{figure}

2. The face ring $\k[\mathcal K]$ is a \emph{quadratic
algebra}\label{quadralg} (that is, the ideal $\mathcal I_\mathcal
K$ is generated by quadratic monomials) if and only if $\mathcal
K$ is a flag complex (an exercise).

3. Let $\mathcal K_1*\mathcal K_2$ be the join of $\mathcal K_1$
and $\mathcal K_2$ (see Construction~\ref{join}). Then
$$
  \k[\mathcal K_1*\mathcal K_2]=\k[\mathcal K_1]\otimes\k[\mathcal K_2].
$$
Here and below $\otimes$ denotes the tensor product over~$\k$.
%
%4. Let $\mathcal K_1\cs\mathcal K_2$ be the connected sum of two
%pure $(n-1)$-dimensional simplicial complexes on sets $\mathcal
%V_1$ и $\mathcal V_2$ respectively, in which two simpleces $ (see
%Construction~\ref{simcs}). Then the ideal $\mathcal I_{\mathcal
%K_1\cs\mathcal K_2}$ is generated by the ideals $\mathcal
%I_{\mathcal K_1}$, $\mathcal I_{\mathcal K_2}$ and the monomials
%$v_{\s}$ и $v_{i_1}v_{i_2}$, где $i_1\in\mathcal
%%M_1\setminus\s_1$ и $i_2\in\mathcal M_2\setminus\s_2$.
\end{example}

We note that $\mathcal I_\mathcal K$ is a \emph{monomial ideal},
and it has a basis consisting of square-free monomials $v_I$
corresponding to the missing faces of~$\mathcal K$.

\begin{proposition}\label{sfmi}
Every square-free monomial ideal $\mathcal I$ in the polynomial
ring is the Stanley--Reisner ideal of a simplicial
complex~$\mathcal K$.
\end{proposition}
\begin{proof}
We set
$$
  \mathcal K=\{I\subset[m]\colon v_I\notin\mathcal I\}.
$$
Then $\mathcal K$ is a simplicial complex and $\mathcal I=\mathcal
I_\mathcal K$.
\end{proof}

Let $P$ be a simple $n$-polytope and let $\mathcal K_P$ be its
nerve complex (see Example~\ref{polsph}). We define the \emph{face
ring}\label{frptope} $\k[P]$ as the face ring of~$\mathcal K_P$.
Explicitly,
$$
  \k[P]=\k[v_1,\ldots,v_m]/\mathcal I_P,
$$
where $\mathcal I_P$ is the ideal generated by those square-free
monomials $v_{i_1}v_{i_2}\cdots v_{i_s}$ whose corresponding
facets intersect trivially, $F_{i_1}\cap\cdots\cap
F_{i_s}=\varnothing$.

\begin{example}\label{frpex}\

1. Let $P$ be an $n$-simplex (viewed as a simple polytope). Then
$$
  \k[P]=\k[v_1,\ldots,v_{n+1}]/(v_1v_2\cdots v_{n+1}).
$$

2. Let $P$ be a 3-cube $I^3$. Then
$$
  \k[P]=\k[v_1,v_2,\ldots,v_6]/(v_1v_4,v_2v_5,v_3v_6).
$$

3. Let $P$  be an $m$-gon, \ $m\ge4$. Then
$$
  \mathcal I_{P}=(v_iv_j\colon i-j\ne0,\pm1\mod m).
$$

4. Given two simple polytopes $P_1$ and $P_2$, we have
$$
  \k[P_1\times P_2]=\k[P_1]\otimes\k[P_2].
$$
\end{example}

\begin{proposition}
\label{frmap} Let $\phi\colon \mathcal K\to \mathcal L$ be a
simplicial map
% (see Definition~\ref{simmap})
between simplicial complexes $\mathcal K$ and $\mathcal L$ on the
vertex sets $[m]$ and $[l]$ respectively. Define the map
$\phi^*\colon \k[w_1,\ldots,w_{l}]\to\k[v_1,\ldots,v_{m}]$ by
$$
  \phi^*(w_j)=\sum_{i\in\varphi^{-1}(j)}v_i.
$$
Then $\phi^*$ descends to a homomorphism $\k[\mathcal
L]\to\k[\mathcal K]$, which we continue to denote~$\phi^*$.
\end{proposition}
\begin{proof}
We only need to check that $\phi^*(\mathcal I_{\mathcal
L})\subset\mathcal I_{\mathcal K}$. Suppose
$J=\{j_1,\ldots,j_s\}\subset[l]$ is not a simplex of~$\mathcal L$.
We have
$$
  \phi^*(w_{j_1}\ldots w_{j_s})=\sum_{i_1\in \phi^{-1}(j_1),\ldots,
  i_s\in \phi^{-1}(j_s)}v_{i_1}\ldots v_{i_s}.
$$
We claim that the right hand side above belongs to $\mathcal
I_{\mathcal K}$, i.e. for any monomial $v_{i_1}\ldots v_{i_s}$ in
the right hand side the set $I=\{i_1,\ldots,i_s\}$ is not a
simplex of~$\mathcal K$. Indeed, otherwise we would have
$\phi(I)=J\in \mathcal L$ by the definition of a simplicial map,
which contradicts the assumption.
\end{proof}

\begin{example}
The face ring of the barycentric subdivision $\sK'$ of $\sK$ is
$$
  \k[\mathcal K']=\k[b_I\colon I\in
  \mathcal K\setminus\varnothing]/\mathcal I_{\mathcal K'},
$$
where $b_I$ is the polynomial generator of degree 2 corresponding
to a nonempty simplex
$I\in \mathcal K$, and $\mathcal I_{\mathcal K'}$ is generated by
quadratic monomials $b_Ib_J$ for which $I\not\subset J$ and
$J\not\subset I$. The simplicial map $\nabla\colon \mathcal K'\to
\mathcal K$ from Example~\ref{nabla} induces a map $\nabla^*$ of
the face ring, given on the generators $v_j\in\k[\mathcal K]$ by
$$
  \nabla^*(v_j)=\sum_{I\in \mathcal K\colon {\min}I=j}b_I.
$$
\end{example}

\begin{example}
The nondegenerate map $\mathcal K'\to\varDelta^{n-1}$ from
Example~\ref{ndbar} induces the following map of the corresponding
face rings:
\begin{align*}
  \k[v_1,\ldots,v_n] &\longrightarrow\k[\mathcal K']\\
  v_i &\longmapsto \sum_{|I|=i} b_I.
\end{align*}
This defines a canonical $\k[v_1,\ldots,v_n]$-module structure on
$\k[\mathcal K']$.
\end{example}

An important tool arising from the functoriality of the face ring
is the restriction homomorphism. For any simplex $I\in \mathcal
K$, the corresponding full subcomplex $\mathcal K_I$ is
$\varDelta^{|I|-1}$ and $\k[\mathcal K_I]$ is the polynomial ring
$\k[v_i\colon i\in I]$ on $|I|$ generators. The inclusion
$\mathcal K_I\subset \mathcal K$ induces the \emph{restriction
homomorphism}
\[
  s_I\colon\k[\mathcal K]\to\k[v_i\colon i\in I],
\]
which maps $v_i$ to zero whenever $i\notin I$.

The following simple proposition will be used in several algebraic
and topo\-lo\-gi\-cal arguments of the later chapters.

\begin{proposition}\label{rmap}
The direct sum
\[
  s={\textstyle\bigoplus\limits_{I\in\sK}} s_I\colon\k[\sK]\longrightarrow
  \bigoplus_{I\in\sK}\k[v_i\colon i\in I]
\]
of all restriction maps is a monomorphism.
\end{proposition}
\begin{proof}
Consider the composite map
$$
\begin{CD}
  \k[v_1,\ldots,v_m] @>p>> \k[\mathcal K] @>s>>
  \bigoplus_{I\in \mathcal K}\k[v_i\colon i\in I]
\end{CD}
$$
where $p$ is the quotient projection. Suppose $s\,p(Q)=0$ where
$Q=Q(v_1,\ldots,v_m)$ is a polynomial. Then for any monomial
$v_{i_1}^{\alpha_1}\cdots v_{i_k}^{\alpha_k}$ which enters $Q$
with a nonzero coefficient we have $I=\{i_1,\ldots,i_k\}\notin
\mathcal K$ (as otherwise the $I$th component of the image under
$s\,p$ is nonzero). Hence $p(Q)=0$ and $s$ is injective.
\end{proof}

\begin{proposition}\label{frbasis}
The face ring $\k[\sK]$ has the $\k$-vector space basis consisting
of monomials $v_{j_1}^{\alpha_1}\cdots v_{j_k}^{\alpha_k}$ where
$\alpha_i>0$ and $\{j_1,\ldots,j_k\}\in\sK$.
\end{proposition}
\begin{proof}
Indeed, the polynomial algebra $\k[m]$ has the $\k$-vector space
basis consisting of all monomials $v_{j_1}^{\alpha_1}\cdots
v_{j_k}^{\alpha_k}$, and such a monomial maps to zero under the
projection $\k[m]\to\k[\sK]$ precisely when
$\{j_1,\ldots,j_k\}\notin\sK$.
\end{proof}

Recall that the Poincar\'e series of a nonnegatively graded
$\k$-vector space $V=\bigoplus_{i=0}^\infty V^i$ is given by
$F(V;\lambda)=\sum_{i=0}^\infty(\dim_\k V^i)\lambda^i$. Since
$\k[\sK]$ is graded by even integers, its Poincar\'e series is
even.

\begin{theorem}[Stanley]
\label{psfr} Let $\sK$ be an $(n-1)$-dimensional simplicial
complex with $f$-vector $(f_0,\ldots,f_{n-1})$ and $h$-vector
$(h_0,\ldots,h_n)$. Then the Poincar\'e series of the face ring
$\k[\mathcal K]$ is
$$
  F\bigl(\k[\mathcal K];\lambda\bigr)=
  \sum_{k=0}^nf_{k-1}\biggl(\frac{\lambda^2}{1-\lambda^2}\biggr)^k=
  \frac{h_0+h_1\lambda^2+\cdots+h_n\lambda^{2n}}{(1-\lambda^2)^n}.
$$
\end{theorem}
\begin{proof}
By Proposition~\ref{frbasis}, a $(k-1)$-dimensional simplex
$\{i_1,\ldots,i_k\}\in\sK$ contributes a summand
$\frac{\lambda^{2k}}{(1-\lambda^2)^k}$ to the Poincar\'e series of
$\sK$ (this summand is just the Poincar\'e series of the subspace
generated by monomials $v_{i_1}^{\alpha_1}\cdots
v_{i_k}^{\alpha_k}$ with positive exponents~$\alpha_i$). This
proves the first identity, and the second follows
from~(\ref{hvectors}).
\end{proof}

\begin{example}\label{simbsim}\

1. Let $\mathcal K=\varDelta^{n-1}$. Then $f_i=\bin{n}{i+1}$ for
$-1\le i\le n-1$, $h_0=1$ and $h_i=0$ for $i>0$. Since every
subset of $[n]$ is a simplex of~$\varDelta^{n-1}$, we have
$\k[\varDelta^{n-1}]=\k[v_1,\ldots,v_n]$ and
$F(\k[\varDelta^{n-1}];\lambda)=(1-\lambda^2)^{-n}$.

2. Let $\mathcal K=\partial\varDelta^n$ be the boundary of an
$n$-simplex. Then $h_i=1$ for $0\le i\le n$, and $\k[\mathcal
\partial\varDelta^n]=\k[v_1,\ldots,v_{n+1}]/(v_1v_2\cdots v_{n+1})$. By
Theorem~\ref{psfr},
$$
  F\bigl(\k[\partial\varDelta^n];\lambda\bigr)=\frac{1+\lambda^2+\cdots+\lambda^{2n}}{(1-\lambda^2)^n}.
$$
\end{example}

The affine algebraic variety corresponding to the commutative
finitely generated $\k$-algebra $\k[\sK]=\k[m]/\mathcal I_\sK$
(i.e. the set of common zeros of elements of $\mathcal I_{\sK}$,
viewed as algebraic functions on~$\k^m$) can be easily identified
as follows.

\begin{proposition}\label{affinefr}
The affine variety corresponding to $\k[\sK]$ is given by
\[
  X(\sK)=\bigcup_{I\in\sK}S_I,
\]
where $S_I=\k\langle\mb e_i\colon i\in I\rangle$ is the coordinate
subspace in $\k^m$ spanned by the set of standard basis vectors
corresponding to~$I$.
\end{proposition}
\begin{proof}
The statement obviously holds in the case $\sK=\varDelta^{m-1}$.
So we assume $\sK\ne\varDelta^{m-1}$. We shall use the following
notation from Section~\ref{alexd}: $\widehat I=[m]\setminus{I}$,
the complement of $I\subset[m]$, and $\widehat{\sK}=\{\widehat
I\in[m]\colon I\notin\sK\}$, the dual complex of~$\sK$. Given a
point $\mb z=(z_1,\ldots,z_m)\in\k^m$, we denote by
\[
  \omega(\!\mb z)=\{i\colon z_i=0\}\subset[m],
\]
the set of zero coordinates of~$\mb z$.

By the definition of the algebraic variety $X(\sK)$ corresponding
to $\k[\sK]$,
\begin{multline*}
  X(\sK)=\bigcap_{J\notin\sK}\bigcup_{j\in J}\{\mb z\colon z_j=0\}=
  \bigcap_{J\notin\sK}\bigl\{\mb z\colon \omega(\!\mb z)\cap J\ne\varnothing\bigr\}\\=
  \bigcap_{\widehat J\in\widehat\sK}\bigl\{\mb z\colon \omega(\!\mb z)\not\subset\widehat J\bigr\}
  =\bigl\{\mb z\colon \omega(\!\mb z)\notin\widehat\sK\bigr\}.
\end{multline*}
On the other hand,
\begin{multline*}
\bigcup_{I\in\sK}S_I=\bigcup_{I\in\sK}\bigcap_{j\in\widehat
I}\{\mb z\colon z_j=0\}=\bigcup_{I\in\sK}\bigl\{\mb
z\colon\widehat I\subset\omega(\!\mb z)\bigr\}
\\=\bigcup_{\widehat I\notin\widehat\sK}\bigl\{\mb
z\colon\omega(\!\mb z)\supset\widehat I\bigr\}=\bigl\{\mb z\colon
\omega(\!\mb z)\notin\widehat\sK\bigr\}.
\end{multline*}
The required identity follows by comparing the two formulae above.
\end{proof}

\begin{remark}\label{arrangem}
The variety $X(\sK)$ is an example of an \emph{arrangement of
coordinate subspaces}, which will be studied further in
Section~\ref{arran}.
\end{remark}

We finish this section with a result showing that the face ring
determines its underlying simplicial complex:

\begin{theorem}[{Bruns--Gubeladze~\cite{br-gu96}}]
Let $\k$ be a field, and $\sK_1$ and $\sK_2$ be two simplicial
complexes on the vertex sets $[m_1]$ and $[m_2]$ respectively.
Suppose $\k[\sK_1]$ and $\k[\sK_2]$ are isomorphic as
$\k$-algebras. Then there exists a bijective map $[m_1]\to[m_2]$
which induces an isomorphism between $\sK_1$ and~$\sK_2$.
\end{theorem}
\begin{proof}
Let $f\colon\k[\sK_1]\to\k[\sK_2]$ be an isomorphism of
$\k$-algebras. An easy argument shows that we can assume that $f$
is a graded isomorphism (an exercise, or
see~\cite[p.~316]{br-gu96}).

Since $f$ is graded, by restriction to the linear components we
observe that $m_1=m_2$ and that $f$ is induced by a linear
isomorphism $F\colon\k[m_1]\to\k[m_2]$. This is described by the
commutative diagram
\[
  \xymatrix{
  \k[v_1,\ldots,v_{m_1}] \ar[r]^F \ar[d] &
  \k[v_1,\ldots,v_{m_2}] \ar[d]\\
  \k[\sK_1] \ar[r]^f & \k[\sK_2]
  }
\]
By passing to the associated affine varieties, we observe that the
isomorphism $f^*\colon X(\sK_2)\to X(\sK_1)$ is the restriction of
the $\k$-linear isomorphism $F^*\colon\k^{m_2}\to\k^{m_1}$. This
is described by the commutative diagram
\[
  \xymatrix{
  \ \ \k^{m_2} \ar[r]^{F^*}  & \ \ \k^{m_1}\\
  X(\sK_2) \ar[u] \ar[r]^{f^*} & X(\sK_1) \ar[u]
  }
\]
The isomorphism $f^*$ establishes a bijective correspondence
\[
  \varPhi\colon\{\text{maximal faces of }\sK_2\}\to
  \{\text{maximal faces of }\sK_1\}
\]
which is defined by the formula $f^*(S_I)=S_{\varPhi(I)}$, where
$I$ is a maximal face of~$\sK_2$. It is also clear that
$|\varPhi(I)|=|I|=\dim S_I$.

We denote by $\mathcal P_1$ the intersection poset of the
subspaces $S_I$, $I\in\sK_1$, with respect to inclusion (i.e. the
elements of $\mathcal P_1$ are nonempty intersections
$S_{I_1}\cap\cdots\cap S_{I_k}$ with $I_j\in\sK_1$). The poset
$\mathcal P_1$ can be also viewed as the intersection poset of the
maximal faces of~$\sK_1$. We define the poset $\mathcal P_2$
corresponding to~$\sK_2$ similarly. The correspondence $\varPhi$
obviously extends to an isomorphism of posets
$\varPhi\colon\mathcal P_2\to\mathcal P_1$, which preserves the
dimension of spaces (or the number of elements in the
intersections of maximal faces).

Now introduce the following equivalence relation on the vertex
sets $[m_1]$ and $[m_2]$: for $i_1,i_2\in [m_1]$ (or
$j_1,j_2\in[m_2]$) we put $i_1\sim i_2$ if and only if the two
sets of maximal faces $\sK_1$ containing $i_1$ and $i_2$
respectively coincide (and similarly for $j_1$ and~$j_2$). The
equivalence classes in $[m_1]$ are the minimal (with respect to
inclusion) nonempty intersections of maximal faces of $\sK_1$, and
similarly for~$[m_2]$. Since $\varPhi$ is an isomorphism of
posets, the two systems of equivalence classes are in natural
bijective correspondence, and the corresponding equivalence
classes  have the same numbers of elements. This gives rise to the
bijective map $\varphi\colon [m_2]\to[m_1]$ which satisfies the
condition that $i\in I$ if and only if $\varphi(i)\in\varPhi(I)$,
where $i\in[m_2]$ and $I\in\sK_2$ is a maximal face. Since any
face of a simplicial complex is contained in a maximal face, we
obtain that $\psi=\varphi^{-1}\colon[m_1]\to[m_2]$ is the required
map.
\end{proof}

\subsection*{Exercises}
\begin{exercise}
Show that the Stanley--Reisner ideal $\mathcal I_\mathcal K$ is
generated by quadratic monomials if and only if $\mathcal K$ is a
flag complex.
\end{exercise}

\begin{exercise}[{see~\cite[(4.7)]{p-r-v04}}]\label{frlimit}
Let $\ca(\sK)$ be the face category of~$\mathcal K$ (objects are
simplices, morphisms are inclusions), $\ca^{op}(\sK)$ the opposite
category (in which the morphisms are reverted), and $\cga$ the
category of commutative graded algebras. (See Appendix~\ref{secmc}
for basics of categories and diagrams.) Consider the diagram
\begin{align*}
  \k[\,\cdot\:]^\sK\colon\ca^{op}(\sK)&\longrightarrow\cga,\\
  I&\longmapsto\k[v_i\colon i\in I]
\end{align*}
whose value on a morphism $I\subset J$ is the surjection
$\k[v_j\colon j\in J]\to\k[v_i\colon i\in I]$ sending each $v_j$
with $j\notin I$ to zero. Show that
\[
  \k[\sK]=\lim\k[\,\cdot\:]^\sK
\]
where the limit is taken in the category $\cga$.
\end{exercise}

\begin{exercise}
If $\k[\sK_1]$ and $\k[\sK_2]$ are isomorphic as $\k$-algebras,
then there is also a graded isomorphism $\k[\sK_1]\to\k[\sK_2]$.
(Hint: Show first that $\k[\sK_1]$ and $\k[\sK_2]$ are isomorphic
as augmented $\k$-algebras, and then pass to the associated graded
algebras with respect to the augmentation ideals.)
\end{exercise}

\section{Tor-algebras and Betti numbers}
\label{hpfr} The algebraic Betti numbers of the face ring
$\k[\sK]$ are the dimensions of the Tor-groups of $\k[\sK]$ viewed
as a module over the polynomial ring. These basic homological
invariants of a simplicial complex~$\sK$ appear to be of great
importance both for combinatorial commutative algebra and toric
topology.

The face ring $\k[\mathcal K]$ acquires a canonical $\k[m]$-module
structure via the quotient projection $\k[m]\to\k[\mathcal K]$. We
therefore may consider the corresponding $\Tor$-modules (see
Appendix, Section~\ref{torapdx}):
$$
  \Tor_{\k[v_1,\ldots,v_m]}\bigl(\k[\mathcal K],\k\bigr)=
  \bigoplus_{i,j\ge0}\Tor^{-i,2j}_{\k[v_1,\ldots,v_m]}\bigl(\k[\mathcal
  K],\k\bigr).
$$
From Lemma~\ref{koscom} we obtain that $\Tor_{\k[m]}(\k[\mathcal
K],\k)$ is a bigraded algebra in a natural way, and there is the
following isomorphism of bigraded algebras:
$$
  \Tor_{\k[v_1,\ldots,v_m]}\bigl(\k[\sK],\k\bigr)\cong
  H\bigl[\Lambda[u_1,\ldots,u_m]\otimes\k[\sK],d\bigr],
$$
where the bigrading and differential on the right hand side are
given by
\begin{equation}\label{diff}
\begin{gathered}
  \bideg u_i=(-1,2),\quad\bideg v_i=(0,2),\\
  du_i=v_i,\quad dv_i=0.
\end{gathered}
\end{equation}

\begin{definition}\label{toralg}
We refer to $\Tor_{\k[v_1,\ldots,v_m]}(\k[\mathcal K],\k)$ as the
\emph{$\Tor$-algebra} of a simplicial complex~$\mathcal K$.

The \emph{bigraded Betti numbers} of $\k[\mathcal K]$ are defined
by
\begin{equation}
\label{bbnfr}
  \beta^{-i,2j}\bigl(\k[\mathcal K]\bigr)=
  \dim_\k\Tor^{-i,2j}_{\k[v_1,\ldots,v_m]}\bigl(\k[\mathcal K],\k\bigr),\qquad
  \text{for }i,j\ge0.
\end{equation}
We also set
$$
  \beta^{-i}\bigl(\k[\mathcal K]\bigr)=
  \dim_\k\Tor^{-i}_{\k[v_1,\ldots,v_m]}\bigl(\k[\mathcal K],\k\bigr)=
  \sum_j\beta^{-i,2j}\bigl(\k[\mathcal K]\bigr).
$$
\end{definition}

The Tor-algebra has the following functorial property:

\begin{proposition}
\label{tamap} A simplicial map $\phi\colon \mathcal K\to \mathcal
L$ between simplicial complexes on the sets $[m]$ and $[l]$
respectively induces a homomorphism
\[
  \phi_{\Tor}^*\colon
  \Tor_{\k[w_1,\ldots,w_{l}]}\bigl(\k[\mathcal L],\k\bigr)\to
  \Tor_{\k[v_1,\ldots,v_{m}]}\bigl(\k[\mathcal K],\k\bigr)
\]
of the corresponding $\Tor$-algebras.
\end{proposition}
\begin{proof}
This follows from Proposition~\ref{frmap} and
Theorem~\ref{torprop}~(b).
\end{proof}

Consider the minimal resolution $(R_{\min},d)$ of the
$\k[m]$-module $\k[\mathcal K]$ (see Construction~\ref{minimal}).
Then $R_{\min}^0\cong1\cdot\k[m]$ is a free module with one
generator of degree~0. The basis of $R_{\min}^{-1}$ is a minimal
generator set for $\mathcal I_\mathcal K$, and these minimal
generators correspond to the missing faces of~$\mathcal K$. Given
a missing face $\{i_1,\ldots,i_k\}\subset[m]$, denote by
$r_{i_1,\ldots,i_k}$ the corresponding generator of
$R_{\min}^{-1}$. Then the map $d\colon R_{\min}^{-1}\to
R_{\min}^0$ takes $r_{i_1,\ldots,i_k}$ to $v_{i_1}\ldots v_{i_k}$.
By Proposition~\ref{mintor}, $\beta^{-1,2j}(\k[\mathcal K])$ is
equal to the number of missing faces with $j$ elements.

\begin{example}
\label{ressqu} Let $\mathcal K=\quad\begin{picture}(5,5)
\put(0,0){\circle*{1}} \put(0,5){\circle*{1}}
\put(5,0){\circle*{1}} \put(5,5){\circle*{1}}
\put(0,0){\line(1,0){5}} \put(0,0){\line(0,1){5}}
\put(5,0){\line(0,1){5}} \put(0,5){\line(1,0){5}}
\put(-2.5,-1){\scriptsize 1} \put(6.3,-1){\scriptsize 2}
\put(-2.5,4){\scriptsize 4} \put(6.3,4){\scriptsize 3}
\end{picture}\quad$,
the boundary of a 4-gon. Then
$$
  \k[\mathcal K]\cong\k[v_1,\ldots,v_4]/(v_1v_3,v_2v_4).
$$
Let us construct a minimal resolution of $\k[\sK]$. The module
$R_{\min}^0$ has one generator 1 (of degree 0). The module
$R_{\min}^{-1}$ has two generators $r_{13}$ and $r_{24}$ of
degree~4, and the differential $d\colon R_{\min}^{-1}\to
R_{\min}^0$ takes $r_{13}$ to $v_1v_3$ and $r_{24}$ to $v_2v_4$.
The kernel $R_{\min}^{-1}\to R_{\min}^0$ is generated by one
element $v_2v_4r_{13}-v_1v_3r_{24}$. Hence, $R_{\min}^{-2}$ has
one generator of degree~8, which we denote by~$a$, and the map
$d\colon R^{-2}_{\min}\to R^{-1}_{\min}$ is injective and takes
$a$ to $v_2v_4r_{13}-v_1v_3r_{24}$. Thus, the minimal resolution
is
$$
\begin{CD}
  0 @>>> R_{\min}^{-2} @>>> R_{\min}^{-1} @>>> R_{\min}^0 @>>> M @>>> 0,
\end{CD}
$$
where $\mathop{\rm rank}R^{0}_{\min}=\beta^{0,0}(\k[\mathcal
K])=1$,\ $\mathop{\rm rank}R^{-1}_{\min}=\beta^{-1,4}(\k[\mathcal
K])=2$,\ $\mathop{\rm rank}R^{-2}_{\min}=\beta^{-2,8}(\k[\mathcal
K])=1$.
\end{example}

The following fundamental result of Hochster reduces the
calculation of the Betti numbers $\beta^{-i,2j}(\k[\mathcal K])$
to the calculation of reduced simplicial cohomology of full
subcomplexes in~$\mathcal K$.
%Having subsequent topological applications in mind, we restrict to
%the case $\k=\Z$ until the end of this Section. This is the most
%basic coefficient ring needed for topologists, but all results
%below are valid with $Z$ replaced by any field.

\begin{theorem}[{Hochster \cite{hoch77}}]\label{hoch}
We have
$$
  \Tor^{-i,2j}_{\k[v_1,\ldots,v_m]}\bigl(\k[\mathcal K],\k\bigr)=\bigoplus_{J\subset[m]\colon|J|=j}
  \widetilde{H}^{j-i-1}(\mathcal K_J;\k),
$$
where $\mathcal K_J$ is the full subcomplex of~$\mathcal K$
obtained by restricting to $J\subset[m]$. We assume
$\widetilde{H}^{-1}(\mathcal K_\varnothing;\k)=\k$ above.
\end{theorem}

We shall give a proof of Hochster's formula
following~\cite{pano08l}. The idea is to first reduce the Koszul
algebra $(\Lambda[u_1,\ldots,u_m]\otimes\k[\sK],d)$ to a certain
finite dimensional quotient~$R^*(\sK)$, without changing the
cohomology, and then identify~$R^*(\sK)$ with the sum of
simplicial cochain complexes of all full subcomplexes in~$\sK$.
The algebra $R^*(\sK)$ will also be used in the cohomological
calculations for moment-angle complexes in Chapter~\ref{macom}.

We use simplified notation $u_Jv_I$ for a monomial $u_J\otimes
v_I$ in the Koszul algebra
$\Lambda[u_1,\ldots,u_m]\otimes\k[\sK]$.

\begin{construction}
\label{astar} We introduce the quotient algebra
$$
  R^*(\sK)=\Lambda[u_1,\ldots,u_m]\otimes\k[\sK]\bigr/(v_i^2=u_iv_i=0,\;
  1\le i\le m).
$$
Since the ideal generated by $v_i^2$ and $u_iv_i$ is homogeneous
and invariant with respect to the differential (since
$d(u_iv_i)=0$ and $d(v_i^2)=0$), we obtain that $R^*(\sK)$ has
differential and bigrading~\eqref{diff}. We also have the quotient
projection
\[
  \varrho\colon\Lambda[u_1,\ldots,u_m]\otimes\k[\sK]\to R^*(\sK).
\]
By definition, the algebra $R^*(\sK)$ has a $\k$-vector space
basis consisting of monomials $u_J v_I$ where $J\subset[m]$,
$I\in\sK$ and $J\cap I=\varnothing$. Therefore,
\begin{equation}\label{rankrpq}
  \dim_{\k}R^{-p,2q}=f_{q-p-1}\bin{m-q+p}p,
\end{equation}
where $(f_0,f_1,\ldots,f_{n-1})$ is the $f$-vector of $\sK$ and
$f_{-1}=1$. We have a $\k$-linear map
$$
  \iota\colon R^*(\sK)\to\Lambda[u_1,\ldots,u_m]\otimes\k[\sK]
$$
sending each $u_J v_I$ identically. The map $\iota$ commutes with
the differentials, and therefore defines a homomorphism of
bigraded differential $\k$-vector spaces satisfying the relation
$\varrho\cdot\iota=\id$. Note that $\iota$ is not a map of
algebras.
\end{construction}

\begin{lemma}
\label{iscoh} The projection homomorphism
$\varrho\colon\Lambda[u_1,\ldots,u_m]\otimes\k[\sK]\to R^*(\sK)$
induces an isomorphism in cohomology.
\end{lemma}
\begin{proof}
The argument is similar to that used for the Koszul resolution
(see Construction~\ref{koszul}). We shall construct a cochain
homotopy between the maps $\id$ and $\iota\cdot\varrho$ from
$\Lambda[u_1,\ldots,u_m]\otimes\k[\sK]$ to itself, that is, a map
$s$ satisfying the identity
\begin{equation}\label{chaineq1}
  ds+sd=\id-\iota\cdot\varrho.
\end{equation}

We first consider the case $\sK=\varDelta^{m-1}$. Then
$\Lambda[u_1,\ldots,u_m]\otimes\k[\varDelta^{m-1}]$ is the Koszul
resolution~\eqref{kosz}, which will be denoted by
\begin{equation}\label{emalg}
  E=E_m=\Lambda[u_1,\ldots,u_m]\otimes\k[v_1,\ldots,v_m],
\end{equation}
and the algebra $R^*(\varDelta^{m-1})$ is isomorphic to
\begin{equation}\label{rmalg}
  \bigl(\Lambda[u]\otimes\k[v]\bigr/(v^2,uv)\bigr)^{\otimes m}.
\end{equation}
For $m=1$, we define the map $s_1\colon E_1^{0,*}=\k[v]\to
E_1^{-1,*}$ by the formula
$$
  s_1(a_0+a_1v+\cdots+a_jv^j)=u(a_2v+a_3v^2+\cdots+a_jv^{j-1}).
$$
We need to check identity~\eqref{chaineq1} for
$x=a_0+a_1v+\cdots+a_jv^j\in E_1^{0,*}$ and for $ux\in
E_1^{-1,*}$, as each element of $E_1$ is the sum of elements of
these two types. In the first case we have
$ds_1x=x-a_0-a_1v=x-\iota\varrho x$, and $s_1dx=0$. In the second
case, i.e. for $ux\in E_1^{-1,*}$, we have $ds_1(ux)=0$, and
$s_1d(ux)=ux-a_0u=ux-\iota\varrho(ux)$. In both
cases~\eqref{chaineq1} holds.

Now we may assume by induction that a cochain homotopy $s_m\colon
E_m\to E_m$ has been already constructed for $m=k-1$. Since
$E_k=E_{k-1}\otimes E_1$, \
$\varrho_k=\varrho_{k-1}\otimes\varrho_1$ and
$\iota_k=\iota_{k-1}\otimes\iota_1$, a direct calculation shows
that the map
\begin{equation}\label{indsk}
  s_k=s_{k-1}\otimes\id+\iota_{k-1}\varrho_{k-1}\otimes s_1
\end{equation}
is a cochain homotopy between $\id$ and $\iota_k\varrho_k$, which
finishes the proof for $\sK=\varDelta^{m-1}$.

In the case of arbitrary $\sK$ the algebras
$\Lambda[u_1,\ldots,u_m]\otimes\k[\sK]$ and $R^*(\sK)$ are
obtained by factorising~\eqref{emalg} and~\eqref{rmalg}
respectively by the ideal $\mathcal I_{\sK}$ in
$\Lambda[u_1,\ldots,u_m]\otimes\k[m]$. Observe that $\mathcal
I_{\sK}$ is generated by $v_I$ with $I\notin\sK$ as an ideal, and
it has a $\k$-vector space basis of monomials
$u_Jv_{i_1}^{\alpha_1}\cdots v_{i_k}^{\alpha_k}$ with
$I=\{i_1,\ldots,i_k\}\notin\sK$ and $\alpha_i>0$. We need to check
that
\[
  d(\mathcal I_{\sK})\subset\mathcal I_{\sK},\quad
  \iota\varrho(\mathcal I_{\sK})\subset\mathcal I_{\sK},\quad
  s(\mathcal I_{\sK})\subset\mathcal I_{\sK}.
\]
The first inclusion is obvious: since $d$ is a derivation, we only
need to check that $dv_I\in\mathcal I_{\sK}$ for $I\notin\sK$, but
$dv_I=0$. The second inclusion is also clear, since
\[
  \iota\varrho(u_Jv_{i_1}^{\alpha_1}\cdots v_{i_k}^{\alpha_k})=
  \begin{cases}
  u_Jv_{i_1}\cdots v_{i_k},&\text{if }\alpha_i=1\text{ and }
  J\cap\{i_1,\ldots,i_k\}=\varnothing;\\
  0,&\text{otherwise.}
  \end{cases}
\]
It remains to check the third inclusion. By expanding the
inductive formula~\eqref{indsk} we obtain
\[
  s_m=s_1\otimes\id\otimes\cdots\otimes\id+
  \iota_1\varrho_1\otimes s_1\otimes\id\otimes\cdots\otimes\id+
  \iota_1\varrho_1\otimes\cdots\otimes
  \iota_1\varrho_1\otimes s_1.
\]
It follows that
\[
  s_m(u_Jv_{i_1}^{\alpha_1}\cdots v_{i_k}^{\alpha_k})=
  \sum_{p\colon\alpha_p>1}\pm
  u_{J}u_{i_p}v_{i_1}^{\alpha_1}\cdots
  v_{i_p}^{\alpha_p-1}\cdots v_{i_k}^{\alpha_k}.
\]
Therefore, $s(\mathcal I_{\sK})\subset\mathcal I_{\sK}$, and
identity~\eqref{chaineq1} holds in
$\Lambda[u_1,\ldots,u_m]\otimes\k[\sK]$.
\end{proof}

As an immediate consequence of Lemma~\ref{iscoh} we obtain

\begin{corollary}
We have that $\beta^{-i,2j}\bigl(\k[\sK]\bigr)=0$ if $i>m$ or
$j>m$.
\end{corollary}
\begin{proof}
Indeed, $R^{-i,2j}(\sK)=0$ if either $i$ or $j$ is greater
than~$m$.
\end{proof}

Now, in order to prove Theorem~\ref{hoch}, we need to show that
the cohomology of $R^*(\sK)$ is isomorphic to the direct sum of
the reduced cohomology of the full subcomplexes on the right hand
side of Hochster's formula. We shall see that this is true even
without passing to cohomology, i.e. $R^*(\sK)$ is isomorphic to
$\bigoplus_{I\subset m}C^*(\sK_I)$, with the appropriate shift in
dimensions, where $C^*$ denotes the simplicial cochain groups. To
do this, it is convenient to refine the grading in $\k[\sK]$ as
follows.

\begin{construction}[multigraded structure in face rings and $\Tor$-algebras]
\label{mgrad} A \emph{multigrading} (more precisely, an
$\N^m$-grading) is defined in $\k[v_1,\ldots,v_m]$ by setting
\[
  \mathop{\rm mdeg}v_1^{i_1}\cdots v_m^{i_m}=(2i_1,\ldots,2i_m).
\]
Since $\k[\mathcal K]$ is the quotient of the polynomial ring by a
monomial ideal, it inherits the multigrading. We may assume that
all free modules in the resolution~(\ref{resol}) are multigraded
and the differentials preserve the multidegree. Then the algebra
$\Tor_{\k[m]}(\k[\mathcal K],\k)$ acquires the canonical
$\Z\oplus\N^m$-grading, i.e.
$$
  \Tor_{\k[v_1,\ldots,v_m]}\bigl(\k[\mathcal K],\k\bigr)
  =\bigoplus_{i\ge0,\;\mb a\in\N^m}
  \Tor^{-i,2\mb a}_{\k[v_1,\ldots,v_m]}\bigl(\k[\mathcal K],\k\bigr).
$$
The differential algebra $R^*(\sK)$ also acquires a
$\Z\oplus\N^m$-grading, and Lemma~\ref{iscoh} implies that
\begin{equation}\label{mgR}
  \Tor^{-i,2\mb a}_{\k[v_1,\ldots,v_m]}\bigl(\k[\mathcal
  K],\k\bigr)\cong H^{-i,2\mb a}(R^*(\sK),d).
\end{equation}
\end{construction}

We may view a subset $J\subset[m]$ as a $(0,1)$-vector in $\N^m$
whose $j$th coordinate is~1 if $j\in J$ and is~0 otherwise. Then
there is the following multigraded version of Hochster's formula:

\begin{theorem}\label{hochmd}
For any subset $J\subset[m]$ we have
$$
  \Tor_{\k[v_1,\ldots,v_m]}^{-i,2J}\bigl(\k[\mathcal
  K],\k\bigr)\cong
  \widetilde{H}^{|J|-i-1}(\mathcal K_J),
$$
and $\Tor_{\k[m]}^{-i,2\mb a}(\k[\mathcal K],\k)=0$ if $\mb a$ is
not a $(0,1)$-vector.
\end{theorem}

\begin{proof}[Proof of Theorem~\ref{hoch} and Theorem~\ref{hochmd}]
Let $C^q(\sK_J)$ denote the $q$th simplicial cochain group with
coefficients in~$\k$. Denote by $\alpha_L\in C^{p-1}(\sK_J)$ the
basis cochain corresponding to an oriented simplex
$L=(l_1,\ldots,l_p)\in\sK_J$; it takes value~$1$ on $L$ and
vanishes on all other simplices. Now we define a $\k$-linear map
\begin{equation}\label{fmapr}
\begin{aligned}
  f\colon C^{p-1}(\sK_J)&\longrightarrow R^{p-|J|,2J}(\sK),\\
  \alpha_L&\longmapsto \varepsilon(L,J)\,u_{J\setminus L}v_L,
\end{aligned}
\end{equation}
where $\varepsilon(L,J)$ is the sign defined by
\[
  \varepsilon(L,J)=\prod_{j\in L}\varepsilon(j,J),
\]
and $\varepsilon(j,J)=(-1)^{r-1}$ if $j$ is the $r$th element of
the set~$J\subset[m]$, written in increasing order. Obviously, $f$
is an isomorphism of $\k$-vector spaces, and a direct check shows
that it commutes with the differentials. Indeed, we have
\begin{align*}
  f(d\alpha_L)&=f\Bigl(\sum_{j\in J\setminus L,\,j\cup L\in\sK_J}
  \varepsilon(j,j\cup L)\,\alpha_{j\cup L}\Bigr)\\
  &=
  \sum_{j\in J\setminus L}\varepsilon(j\cup L,J)
  \varepsilon(j,j\cup L)\,u_{J\setminus(j\cup L)}v_{j\cup L}
\end{align*}
(note that $v_{j\cup L}\in\k[\sK]$, and hence it is zero unless
$j\cup L\in\sK_J$). On the other hand,
\[
  df(\alpha_L)=\sum_{j\in J\setminus L}\varepsilon(L,J)
  \varepsilon(j,J\setminus L)u_{J\setminus(j\cup L)}v_{j\cup L}.
\]
By the definition of $\varepsilon(L,J)$,
\[
  \varepsilon(j\cup L,J)\varepsilon(j,j\cup L)=
  \varepsilon(L,J)\varepsilon(j,J)\varepsilon(j,j\cup L)=
  \varepsilon(L,J)\varepsilon(j,J\setminus L),
\]
which implies that $f(d\alpha_L)=df(\alpha_L)$. Therefore, $f$
together with the map $\k\to R^{-|J|,2J}(\sK)$, $1\mapsto u_J$,
defines an isomorphism of cochain complexes
\[
\begin{array}{ccccccccc}
0 \to\! & \k\ \ & \stackrel{d}\longrightarrow & C^0(\sK_J) &
\stackrel{d}\longrightarrow
% & C^1(\sK_J) & \stackrel{d}\longrightarrow
& \cdots & \stackrel{d}\longrightarrow
& C^{p-1}(\sK_J) & \stackrel{d}\longrightarrow \cdots\\
       & \downarrow\cong
  %& & \downarrow\cong
  & & \downarrow\cong
  && & & \downarrow\cong\\
0 \to\! & R^{-|J|,2J}(\sK) & \stackrel{d}\longrightarrow &
R^{1-|J|,2J}(\sK)
%& \stackrel{d}\longrightarrow & R^{2-|J|,2J}(\sK)
& \stackrel{d}\longrightarrow & \cdots &
\stackrel{d}\longrightarrow & R^{p-|J|,2J}(\sK) &
\stackrel{d}\longrightarrow \cdots
\end{array}
\]
Then it follows from~\eqref{mgR} that
\[
  \widetilde H^{p-1}(K_J)\cong\Tor^{p-|J|,2J}_{\k[v_1,\ldots,v_m]}\bigl(\k[\mathcal
  K],\k\bigr),
\]
which is equivalent to the first isomorphism of
Theorem~\ref{hochmd}. Since $R^{-i,2\mb a}(\sK)=0$ if $\mb a$ is
not a $(0,1)$-vector, $\Tor_{\k[m]}^{-i,2\mb a}(\k[\mathcal
K],\k)$ vanishes for such~$\mb a$.
\end{proof}

Since $\Tor_{\k[m]}(\k[\sK],\k)$ is an algebra, the isomorphisms
of Theorem~\ref{hoch} turn the direct sum
\begin{equation}
\label{dsfullsub}
  \mathop{\bigoplus_{p\ge0}}\limits_{J\subset[m]}\widetilde H^{p-1}(\sK_J)
\end{equation}
into a (multigraded) $\k$-algebra. Consider the product in the
simplicial cochains of full subcomplexes given by
\begin{equation}\label{fullsubcochain}
\begin{aligned}
  \mu\colon C^{p-1}(\sK_I)\otimes C^{q-1}(\sK_J)&\longrightarrow C^{p+q-1}(\sK_{I\cup
  J}),\\
  \alpha_L\otimes \alpha_M\ \ &\longmapsto
    \left\{%
    \begin{array}{ll}
    \alpha_{L\sqcup M}, & \text{if $I\cap J=\varnothing$;} \\
    0, & \hbox{otherwise.} \\
    \end{array}%
    \right.
\end{aligned}
\end{equation}
Here $\alpha_{L\sqcup M}\in C^{p+q-1}(\sK_{I\sqcup J})$ denotes
the basis simplicial cochain corresponding to $L\sqcup M$ if the
latter is a simplex of $\sK_{I\sqcup J}$ and zero otherwise. If
$I\cap J=\varnothing$, then $\sK_{I\sqcup J}$ is a subcomplex in
the join $\sK_I*\sK_J$, and the above product is the restriction
to $\sK_{I\sqcup J}$ of the standard exterior product
\[
  C^{p-1}(\sK_I)\otimes C^{q-1}(\sK_J)\longrightarrow
  C^{p+q-1}(\sK_I*\sK_J).
\]

\begin{proposition}\label{proddirsum}
The product in the direct sum
$\bigoplus_{p\ge0,\;J\subset[m]}\widetilde H^{p-1}(\sK_J)$ induced
by the isomorphisms from Hochster's theorem coincides up to a sign
with the product given by~\eqref{fullsubcochain}.
\end{proposition}
\begin{proof}
This is a direct calculation. We use the isomorphism $f$ given
by~\eqref{fmapr}:
\[
  \alpha_L\cdot\alpha_M=f^{-1}\bigl(f(\alpha_L)\cdot
  f(\alpha_M)\bigr)=
  f^{-1}\bigl(\varepsilon(L,I)\,u_{I\setminus L}v_L\,
  \varepsilon(M,J)\,u_{J\setminus M}v_M\bigr)
\]
If $I\cap J\ne\varnothing$, then the product $u_{I\setminus L}v_L
u_{J\setminus M}v_M$ is zero in~$R^*(\sK)$. Otherwise we have that
$u_{I\setminus L}v_L u_{J\setminus M}v_M=\zeta\,u_{(I\cup
J)\setminus(L\cup M)}v_{L\cup M}$, where $\zeta=\prod_{k\in
I\setminus L}\varepsilon(k,k\cup J\setminus M)$, and we can
continue the above identity as
\[
  \alpha_L\cdot\alpha_M=\varepsilon(L,I)\,\varepsilon(M,J)\,\zeta\,\varepsilon(L\cup
  M,I\cup J)\,\alpha_{L\sqcup M}.
\]
Note that this calculation also gives the explicit value for the
correcting sign, but we shall not need this.
\end{proof}

Let $P$ be a simple polytope. The multigraded components of
$\Tor_{\k[m]}(\k[P],\k)$ can be expressed directly in terms of~$P$
as follows. Let $\{F_1,\ldots,F_m\}$ be the set of facets of~$P$.
Given $I\subset[m]$, we define the following subset of the
boundary of~$P$:
\[
  P_I=\bigcup_{i\in I}F_i.
\]

\begin{proposition}\label{hochpol}
For any subset $J\subset[m]$ we have
$$
  \Tor_{\k[v_1,\ldots,v_m]}^{-i,2J}\bigl(\k[P],\k\bigr)\cong
  \widetilde{H}^{|J|-i-1}(P_J),
$$
and $\Tor_{\k[m]}^{-i,2\mb a}(\k[P],\k)=0$ if $\mb a$ is not a
$(0,1)$-vector.
\end{proposition}
\begin{proof}
Let $\sK=\sK_P$ be the nerve complex of~$P$. Then the statement
follows from Theorem~\ref{hochmd} and the fact that $\sK_J$ is a
deformation retract of~$P_J$. The latter is because $P$ is simple,
and therefore, $P_J=\bigcup_{i\in J}F_J=\bigcup_{i\in
J}\st_{\sK'}\{i\}$ (by Proposition~\ref{stardual}), which is the
combinatorial neighbourhood of $(\sK_J)'$ in~$\sK'$.
\end{proof}

For a description of the multiplication in
$\Tor_{\k[m]}(\k[P],\k)$ in terms of~$P$, see
Exercise~\ref{multppair}.

\begin{example}\label{exbettinum}\ \nopagebreak

1. Let $P$ be a 4-gon, so that $\mathcal
K_P=\quad\begin{picture}(5,5) \put(0,0){\circle*{1}}
\put(0,5){\circle*{1}} \put(5,0){\circle*{1}}
\put(5,5){\circle*{1}} \put(0,0){\line(1,0){5}}
\put(0,0){\line(0,1){5}} \put(5,0){\line(0,1){5}}
\put(0,5){\line(1,0){5}} \put(-2.5,-1){\scriptsize 1}
\put(6.3,-1){\scriptsize 2} \put(-2.5,4){\scriptsize 4}
\put(6.3,4){\scriptsize 3}
\end{picture}\quad$, as in Example~\ref{ressqu}. This time we calculate the
Betti numbers $\beta^{-i,2j}(\k[P])$ using Hochster's formula. We
have
\begin{align*}
  \beta^{0,0}\bigl(\k[P]\bigr)&=
  \dim \widetilde H^{-1}(\varnothing)=1 && 1\\
  \beta^{-1,4}\bigl(\k[P]\bigr)&=
  \dim \widetilde H^0\bigl(P_{\{1,3\}}\bigr)\oplus
  \widetilde H^0\bigl(P_{\{2,4\}}\bigr)=2 && u_1v_3,u_2v_4\\
  \beta^{-2,8}\bigl(\k[P]\bigr)&=
  \dim\widetilde H^1\bigl(P_{\{1,2,3,4\}}\bigr)=1 &&
  u_1u_2v_3v_4
\end{align*}
where in the right column we include cocycles in the Koszul
algebra $\Lambda[u_1,\ldots,u_m]\otimes\k[P]$ representing
generators of the corresponding cohomology groups. All other Betti
numbers are zero. We have a nontrivial product
$[u_1v_3]\cdot[u_2v_4]=[u_1u_2v_3v_4]$; all other products of
positive-dimensional classes are zero. Note that in this example
all $\Tor$-groups have bases represented by monomials in the
Koszul algebra. This is not the case in general, as is shown by
the next example.

2. Now let
$\mathcal K=\quad\begin{picture}(18,2)
\put(0,1){\circle*{1}} \put(5,1){\circle*{1}}
\put(13,1){\circle*{1}} \put(18,1){\circle*{1}}
\put(0,1){\line(1,0){5}} \put(13,1){\line(1,0){5}}
\put(-2.5,0){\scriptsize 1} \put(6.5,0){\scriptsize 2}
\put(10.5,0){\scriptsize 3} \put(19.5,0){\scriptsize 4}
\end{picture}\quad$
be the union of two segments. Then the generator of
\[
  \Tor^{-3,8}_{\k[v_1,\ldots,v_4]}\bigl(\k[\sK],\k\bigr)
  \cong \widetilde H^0\bigl(\sK_{\{1,2,3,4\}};\k\bigr)\cong\k
\]
is represented by the cocycle $u_1u_2u_3v_4-u_1u_2u_4v_3$ in the
Koszul algebra, and it cannot be represented by a monomial.

3. Let us calculate the Betti numbers (both bigraded and
multigraded) of $\k[\sK]$ for the complex shown in
Fig.~\ref{figfr}, using Hochster's formula. We have
\begin{align*}
  \beta^{0,0}&=\dim \widetilde H^0(\varnothing)=1,\\
  \beta^{-1,4}&=\beta^{-1,(2,0,0,0,2)}+\beta^{-1,(0,0,2,2,0)}=
  \dim \widetilde H^0\bigl(\sK_{\{1,5\}}\bigr)\oplus\widetilde H^0\bigl(\sK_{\{3,4\}}\bigr)=2,\\
  \beta^{-1,6}&=\beta^{-1,(2,2,2,0,0)}+\beta^{-1,(0,2,0,2,2)}=
  \dim \widetilde H^1\bigl(\sK_{\{1,2,3\}}\bigr)\oplus\widetilde H^1\bigl(\sK_{\{2,4,5\}}\bigr)=2,\\
  \beta^{-2,8}&=\beta^{-2,(0,2,2,2,2)}+\beta^{-2,(2,0,2,2,2)}
  +\cdots+\beta^{-2,(2,2,2,2,0)}\\
  &\qquad\qquad=\dim \widetilde H^1\bigl(\sK_{\{2,3,4,5\}}\bigr)\oplus\cdots
  \oplus\widetilde H^1\bigl(\sK_{\{1,2,3,4\}}\bigr)=5,\\
  \beta^{-3,10}&=\beta^{-3,(2,2,2,2,2)}=
  \dim \widetilde H^1\bigl(\sK_{\{1,2,3,4,5\}}\bigr)=2.
\end{align*}
All other Betti numbers are zero.

4. Let $\sK$ be a triangulation of the real projective plane $\R
P^2$ with $m$ vertices (the minimal example has $m=6$, see
Fig.~\ref{rptri}, where the vertices with the same labels are
identified, and the boundary edges are identified according to the
orientation shown).
\begin{figure}[h]
  \begin{center}
  \begin{picture}(120,40)
  \put(59,-3){\small 5}
  \put(59,41){\small 5}
  \put(81,29){\small 6}
  \put(81,9){\small 4}
  \put(37.5,29){\small 4}
  \put(37.5,9){\small 6}
  \put(59.5,26){\small 2}
  \put(52,11){\small 3}
  \put(67,11){\small 1}
  \put(60,0){\vector(-2,1){20}}
  \put(40,10){\vector(0,1){20}}
  \put(40,30){\vector(2,1){20}}
  \put(60,40){\vector(2,-1){20}}
  \put(80,30){\vector(0,-1){20}}
  \put(80,10){\vector(-2,-1){20}}
  \put(40,10){\line(1,0){40}}
  \put(40,30){\line(1,0){40}}
  \put(60,0){\line(-1,1){10}}
  \put(60,0){\line(1,1){10}}
  \put(50,10){\line(-1,2){10}}
  \put(50,10){\line(1,2){10}}
  \put(70,10){\line(-1,2){10}}
  \put(70,10){\line(1,2){10}}
  \put(60,30){\line(0,1){10}}
  \end{picture}
  \end{center}
  \caption{6-vertex triangulation of $\R P^2$.}
  \label{rptri}
\end{figure}
Then, by Hochster's formula,
\[
  \Tor^{3-m,2m}_{\k[v_1,\ldots,v_m]}\bigl(\k[\sK],\k\bigr)=\widetilde H^{2}(\sK_{[m]};\k)=\widetilde H^{2}(\R
  P^2;\k)=0
\]
if the characteristic of $\k$ is not~2. On the other hand,
\[
  \Tor^{3-m,2m}_{\Z_2[v_1,\ldots,v_m]}\bigl(\Z_2[\sK],\Z_2\bigr)
  =\widetilde H^{2}(\sK_{[m]};\Z_2)=\widetilde H^{2}(\R
  P^2;\Z_2)=\Z_2.
\]
This example shows that the $\Tor$-groups of $\k[\sK]$, and even
the algebraic Betti numbers, depend on~$\k$. A similar example
shows that $\Tor_{\Z[m]}(\Z[\sK],\Z)$ may have an arbitrary amount
of additive torsion. (This is a well-known fact for the usual
cohomology of spaces, and so we may take $\sK$ to be a
triangulation of a space with the appropriate torsion in
cohomology.)
\end{example}

\subsection*{Exercises.}
\begin{exercise}
Let $P$ be a pentagon. Calculate the bigraded Betti numbers of
$\k[P]$ and the multiplication in $\Tor_{\k[m]}(\k[P],\k)$
\begin{itemize}
\item[(a)] using algebra $R^*(\sK_P)$ and Lemma~\ref{iscoh};
\item[(b)] using Hochster's theorem and Proposition~\ref{proddirsum},
\end{itemize}
and compare the results.
\end{exercise}

\begin{exercise}\label{multppair}
Use Proposition~\ref{hochpol} and the isomorphism
\[
  \widetilde{H}^{|J|-i-1}(P_J)\cong H^{|J|-i}(P,P_J)
\]
to show that the multiplication induced from
$\Tor_{\k[m]}(\k[P],\k)$ in the direct sum
\[
  \mathop{\bigoplus_{p\ge0}}\limits_{J\subset[m]}H^p(P,P_J)
\]
comes from the standard exterior multiplication
\[
  H^p(P,P_I)\otimes H^q(P,P_J)\longrightarrow H^{p+q}(P,P_I\cup P_J)
\]
when $I\cap J=\varnothing$ and is zero otherwise.
\end{exercise}

\begin{exercise}\label{algproofad}
Complete the details in the following algebraic proof of the
Alexander duality (Theorem~\ref{simaldual}); this argument goes
back to the original work of Hochster~\cite{hoch77}:

1. Choose $J\notin\sK$, that is, $\widehat J=[m] \setminus
J\in\widehat\sK$, and show that for any
$L=\{l_1,\ldots,l_q\}\subset J$,
\[
  J\setminus L\notin\sK\;\Longleftrightarrow\;
  L\in\lk_{\widehat\sK}\widehat J.
\]

2. Consider the Koszul algebra
\[
  S(\sK)=\bigl[\Lambda[u_1,\ldots,u_m]\otimes\mathcal
  I_\sK,d\bigr]
\]
of the Stanley--Reisner ideal $\mathcal I_\sK$ (see
Lemma~\ref{koscom} and the remark after it), and show that its
multigraded component $S^{-q,2J}(\sK)$ has a $\k$-basis consisting
of monomials $u_Lv_{J\setminus L}$ where
$L\in\lk_{\widehat\sK}\widehat J$.

3. Consider the $\k$-vector space isomorphism
\[
\begin{aligned}
  g\colon C_{q-1}(\lk_{\widehat\sK}\widehat J)&\longrightarrow S^{-q,2J}(\sK),\\
  [L]&\longmapsto u_Lv_{J\setminus L},
\end{aligned}
\]
where $[L]\in C_{q-1}(\lk_{\widehat\sK}\widehat J)$ is the basis
simplicial chain corresponding to~$L$. Show that $g$ commutes with
the differentials, and therefore defines an isomorphism of chain
complexes (in analogy with~\eqref{fmapr}, but with no correction
sign).

4. Deduce that
\[
  \widetilde H_{q-1}(\lk_{\widehat\sK}\widehat J)\cong
  \Tor^{-q,2J}_{\k[m]}(\mathcal I_\sK,\k)\cong
  \Tor^{-q-1,2J}_{\k[m]}(\k[\sK],\k)\cong
  \widetilde H^{|J|-q-2}(\sK_J),
\]
where the first isomorphism is obtained by passing to homology in
step~3, the second follows from the long exact sequence of
Theorem~\ref{torprop}~(e),
%(applied to the short exact
%sequence $0\to\mathcal I_\sK\to\k[m]\to\k[\sK]\to0$)
and the third is Theorem~\ref{hoch}. It remains to note that the
resulting isomorphism is equivalent to that of
Corollary~\ref{fulink}.
\end{exercise}

\section{Cohen--Macaulay complexes}\label{cmr}
It is usually quite difficult to determine whether a given ring is
Cohen--Macaulay (see Appendix, Section~\ref{cma}). One of the key
results of combinatorial commutative algebra, the \emph{Reisner
Theorem}, gives an effective criterion for the Cohen--Macaulayness
of face rings, in terms of simplicial cohomology of~$\sK$. A
reformulation of Reisner's criterion, due to Munkres and Stanley,
tells us that the Cohen--Macaulayness of the face ring $\k[\sK]$
is a topological property of $\sK$, i.e. it depends only on the
topology of the realisation~$|\sK|$. These results have many
important applications in both combinatorial commutative algebra
and toric topology.

Here we assume that $\k$ is a field, unless otherwise stated. If
$\sK$ is of dimension $n-1$, then the Krull dimension of $\k[\sK]$
is~$n$ (an exercise). We start with the following combinatorial
description of homogeneous systems of parameters (hsop's) in
$\k[\mathcal K]$ in terms of the restriction homomorphisms
$s_I\colon\k[\mathcal K]\to \k[v_i\colon i\in I]$ (see
Proposition~\ref{rmap}).

\begin{lemma}\label{lsopr}
Let $\mathcal K$ be a simplicial complex of dimension $n-1$. A
sequence of homogeneous elements $\mb t=(t_1,\ldots,t_n)$ of the
face ring $\k[\mathcal K]$ is a homogeneous system of parameters
if and only if
\[
  \dim_\k\bigl(\k[v_i\colon i\in I]/s_I(\mb t)\bigr)<\infty
\]
for each simplex $I\in \mathcal K$, where $s_I(\mb t)$ is the
image of the sequence $\mb t$ under the restriction map~$s_I$.
\end{lemma}
\begin{proof}
Assume that $\mb t$ is an hsop. By applying the right exact
functor $\ \otimes_{\k[\mb t]}$ to the epimorphism
$s_I\colon\k[\mathcal K]\to\k[v_i\colon i\in I]$ we obtain that
$\k[\mathcal K]/\mb t\to\k[v_i\colon i\in I]/s_I(\mb t)$ is also
an epimorphism. Hence,
\[
  \dim_\k\bigl(\k[v_i\colon i\in I]/s_I(\mb t)\bigr)\le
  \dim_\k\bigl(\k[\mathcal K]/\mb t\bigr)<\infty.
\]
%Therefore, $s_I(\mb t)$ is an hsop in $\k[v_i\colon i\in I]$ for
%every $I\in\sK$.

For the opposite statement, assume that
$$
  \dim_\k\bigoplus_{I\in \mathcal K}
  \k[v_i\colon i\in I]/s_I(\mb t)<\infty.
$$
Consider the short exact sequence of $\k[\mb t]$-modules
\[
  0\longrightarrow\k[\sK]\stackrel{\oplus s_I}
  \longrightarrow \bigoplus_{I\in \mathcal K}\k[v_i\colon i\in I]
  \longrightarrow Q\longrightarrow0,
\]
where $Q$ is the quotient module, and the following fragment of
the corresponding long exact sequence for~$\Tor$
(Theorem~\ref{torprop}~(e)):
\[
  \cdots\longrightarrow\Tor^{-1}_{\k[\mb t]}(Q,\k)
  \longrightarrow\k[\sK]/\mb t
  \longrightarrow\bigoplus_{I\in \mathcal K}\k[v_i\colon i\in I]/s_I(\mb t)
  \longrightarrow\cdots.
\]
Since $\bigoplus_{I\in \mathcal K}\k[v_i\colon i\in I]$ is a
finitely generated $\k[\mb t]$-module, its quotient $Q$ is also
finitely generated. Hence $\dim_{\k}\Tor^{-1}_{\k[\mb
t]}(Q,\k)<\infty$ (see Proposition~\ref{mintor}), and, by the
exact sequence above, $\dim_{\k}\bigl(\k[\mathcal K]/\mb
t\bigr)<\infty$. Therefore, $\mb t$ is an hsop in~$\k[\sK]$.
\end{proof}

Recall that we refer to a sequence $\mb
t=(t_1,\ldots,t_n)\in\k[\sK]$ as linear if $\deg t_i=2$ for
all~$i$. We may write a linear sequence as
\begin{equation}\label{linearseq}
  t_i=\lambda_{i1}v_1+\cdots+\lambda_{im}v_m,\quad
  \text{for }1\le i\le n.
\end{equation}

Here is a simple characterisation of lsops in a face ring (see
Definition~\ref{krdim}):

\begin{lemma}\label{lsopcrit}
A linear sequence $\mb t=(t_1,\ldots,t_n)\subset\k[\sK]$, \
$\dim\sK=n-1$, is an lsop if and only if the restriction $s_I(\mb
t)$ to each simplex of $I\in\mathcal K$ generates the polynomial
algebra $\k[v_i\colon i\in I]$ (i.e. the rank of the $n\times|I|$
matrix $\Lambda_I=(\lambda_{ij})$, \ $1\le i\le n$, $j\in I$, is
equal to~$|I|$).
\end{lemma}
\begin{proof}
Indeed, if $\mb t$ is linear, then the conditions
$\dim_\k\k[v_i\colon i\in I]/s_I(\mb t)<\infty$ and $\k[v_i\colon
i\in I]/s_I(\mb t)\cong\k$ are equivalent. The latter means that
$s_I(\mb t)$ generates $\k[v_i\colon i\in I]$ as a $\k$-algebra.
\end{proof}

Note that it is enough to verify the conditions of
Lemmata~\ref{lsopr} and~\ref{lsopcrit} only for maximal simplices
$I\in\sK$.

\begin{definition}\label{integrallsop}
A linear sequence $\mb t=(t_1,\ldots,t_n)\subset\Z[\sK]$ is
referred to as a \emph{integral lsop} if its reduction modulo $p$
is an lsop in $\Z_p[\sK]$ for any prime~$p$. Equivalently, $\mb t$
is an integral lsop if $n=\dim\sK+1$ and the restriction $s_I(\mb
t)$ to each simplex $I\in\sK$ generates the polynomial ring
$\Z[v_i\colon i\in I]$ (the equivalence of these two conditions is
an exercise).
\end{definition}

Although the rational face ring $\Q[\sK]$ always admit an lsop by
Theorem~\ref{noether}, an lsop in $\Z_p[\sK]$ for a prime $p$ (or
an integral lsop in $\Z[\sK]$) may fail to exist, as is shown by
the next example.

\begin{example}\label{dj-lsop}\

1. Let $\sK$ be a simplicial complex of dimension $n-1$ on $m\ge
2^n$ vertices whose 1-skeleton is a complete graph. Then the face
ring $\Z_2[\sK]$ does not admit an lsop. Indeed, assume
that~\eqref{linearseq} is an lsop. Then, by
Corollary~\ref{lsopcrit}, each column vector of the $n\times
m$-matrix $(\lambda_{ij})$ is nonzero, and all column vectors are
pairwise different (since each pair of vertices of~$\sK$ spans an
edge). This is a contradiction, since the number of different
nonzero vectors in $\Z_2^n$ is $2^n-1$. By considering the
reduction modulo 2 we obtain that $\Z[\sK]$ also does not admit an
integral lsop.

2. There are also simple polytopes $P$ whose face rings $\k[P]$ do
not admit an lsop over~$\Z_2$ or~$\Z$. Indeed, let $P$ be the dual
of a 2-neighbourly simplicial $n$-polytope (e.g., a cyclic
polytope of dimension $n\ge4$, see Example~\ref{cyclic}) with
$m\ge 2^n$ vertices. Then the 1-skeleton of $\sK_P$ is a complete
graph, and therefore $\Z[P]=\Z[\sK_P]$ does not admit an lsop.
This example is taken from~\cite{da-ja91}.
%[1.22, 5.4]{da-ja91}.
\end{example}

By considering the reduction modulo 2 we observe that the ring
$\Z[\sK]$ for $\sK$ from the previous example also does not admit
an integral lsop. Existence of integral lsop's in the face rings
$\Z[\sK]$ is a very subtle question of great importance for toric
topology; it will be discussed in more detail in
Section~\ref{freetorusactions}.

\begin{definition}\label{CMsim}
$\mathcal K$ is a \emph{Cohen--Macaulay complex over a field~$\k$}
if $\k[\mathcal K]$ is a Cohen--Macaulay algebra. We say that
$\mathcal K$ is a \emph{Cohen--Macaulay complex over~$\Z$}, or
simply a \emph{Cohen--Macaulay complex}, if $\k[\sK]$ is a
Cohen--Macaulay algebra for $\k=\Q$ and any finite field.
\end{definition}

\begin{remark}
We shall often consider $\k[\sK]$ as a $\k[m]$-module rather than
a $\k$-algebra. However, this does not affect regular sequences
and the Cohen--Macaulay property: it is an easy exercise to show
that a sequence $\mb t\subset\k[m]$ is $\k[m]$-regular for
$\k[\sK]$ as a $\k[m]$-module if and only the image of $\mb t$ in
$\k[\sK]$ is $\k[\sK]$-regular. In particular, $\k[\sK]$ is a
Cohen--Macaulay algebra if and only if it is a Cohen--Macaulay
$\k[m]$-module. We shall therefore not distinguish between these
two notions.
\end{remark}

\begin{example}
Let $\mathcal K=\partial\varDelta^2$. Then $\k[\mathcal
K]=\k[v_1,v_2,v_3]/(v_1v_2v_3)$ and the Krull dimension is
$\dim\k[\sK]=2$. The elements $v_1,v_2\in\k[\mathcal K]$ are
algebraically independent, but do not form an lsop, since
$\k[\sK]/(v_1,v_2)\cong\k[v_3]$ and $\dim\bigl(\k[\mathcal
K]/(v_1,v_2)\bigl)=1$. On the other hand, the elements
$t_1=v_1-v_3$, $t_2=v_2-v_3$ form an lsop, since $\k[\mathcal
K]/(t_1,t_2)\cong\k[t]/t^3$. It is easy to see that $\k[\mathcal
K]$ is a free $\k[t_1,t_2]$-module on the basis $\{1,v_1,v_1^2\}$.
Therefore, $\k[\mathcal K]$ is a Cohen--Macaulay ring and
$(t_1,t_2)$ is a regular sequence.
\end{example}

Cohen--Macaulay complexes can be characterised homologically as
follows:

\begin{proposition}\label{cmbet}
The following conditions are equivalent for a simplicial complex
$\mathcal K$ of dimension $n-1$ with $m$ vertices:
\begin{itemize}
\item[(a)] $\sK$ is Cohen--Macaulay over a field $\k$;
\item[(b)] $\beta^{-i}(\k[\mathcal K])=0$ for $i>m-n$ and
$\beta^{-(m-n)}(\k[\mathcal K])\ne0$.
\end{itemize}
\end{proposition}
\begin{proof}
By definition, condition~(a) is that
$\mathop{\mathrm{depth}}\k[\mathcal K]=n$. Condition~(b) is
equivalent to the equality $\pdim\k[\mathcal K]=m-n$, by
Corollary~\ref{hdtor}. Since $\depth\k[m]=m$, the two conditions
are equivalent by Theorem~\ref{abthm}.
\end{proof}

\begin{example}
Let $\sK$ be the 6-vertex triangulation of $\R P^2$, see
Example~\ref{exbettinum}.4 and Figure~\ref{rptri}. Then $m-n=3$,
and, by Theorem~\ref{hoch},
\[
  \beta^{-4}\bigl(\Z_2[\sK]\bigl)=\dim_{\Z_2}\widetilde H^1(\R
  P^2;\Z_2)=1,
\]
so $\sK$ is not Cohen--Macaulay over $\Z_2$. On the other hand, a
similar calculation shows that if the characteristic of $\k$ is
not~2, then $\beta^{-i}(\k[\sK])=0$ for $i>3$ and
$\beta^{-3}(\k[\sK])=6$, i.e. $\sK$ is Cohen--Macaulay over such
fields.
\end{example}

\begin{proposition}[Stanley]\label{stanmv}
If $\mathcal K$ is a Cohen--Macaulay complex of dimension $n-1$,
then $\mb h(\mathcal K)=(h_0,\ldots,h_n)$ is an $M$-vector (see
Definition~\ref{defmvect}).
\end{proposition}
\begin{proof}
Let $\sK$ be Cohen--Macaulay, and let $\mb t=(t_1,\ldots,t_n)$ be
an lsop in $\k[\mathcal K]$, where $\k$ is a field of zero
characteristic. Then, by Proposition~\ref{pscma},
\[
  F\bigl(\k[\sK],\lambda\bigr)=\frac{F\bigl( \k[\sK]/\mb t;\lambda
  \bigr)}{(1-\lambda^2)^n}.
\]
On the other hand, the Poincar\'e series of $\k[\sK]$ is given by
Theorem~\ref{psfr}, which implies that
\[
  F\bigl( \k[\sK]/\mb
  t;\lambda\bigr)=h_0+h_1\lambda^2+\cdots+h_n\lambda^{2n}.
\]
Now, $A=\k[\sK]/\mb t$ is a graded algebra generated by its
degree-two elements and $\dim_\k A^{2i}=h_i$, so
$(h_0,\ldots,h_n)$ is an $M$-vector by definition.
\end{proof}

\begin{remark}
According to a result of Stanley~\cite[Theorem~II.3.3]{stan96}, if
$(h_0,\ldots,h_n)$ is an $M$-vector, then there exists an
$(n-1)$-dimensional Cohen--Macaulay complex $\sK$ such that
$h_i(\sK)=h_i$. Together with Proposition~\ref{stanmv}, this gives
a complete characterisation of face vectors of Cohen--Macaulay
complexes.
\end{remark}

The following fundamental result gives a combinatorial
characterisation of Cohen--Macaulay complexes:

\begin{theorem}[Reisner~\cite{reis76}]
\label{reisner} Let $\k=\Z$ or a field. A simplicial complex
$\mathcal K$ is Cohen--Macaulay over $\k$ if and only if for any
simplex $I\in \mathcal K$ (including $I=\varnothing$) and
$i<\dim(\lk I)$, we have $\widetilde{H}_i(\lk I;\k)=0$.
\end{theorem}

The proof can be found in many texts on combinatorial commutative
algebra, see e.g.~\cite[\S{}II.4]{stan96}
or~\cite[Chapter~13]{mi-st05}. It is not very hard, but uses the
notion of \emph{local cohomology}, which is beyond the scope of
this book.

The following reformulation of Reisner's Theorem shows that the
Cohen--Macaulayness is a topological property of a simplicial
complex.

\begin{proposition}[Munkres, Stanley]\label{munkr}
A simplicial complex $\sK$ is Cohen--Macaulay over $\k$ if and
only if for any point $x\in|\mathcal K|$, we have
\[
  H_i(|\mathcal K|,|\mathcal
  K|\!\setminus x;\k)=\widetilde{H}_i(\mathcal K;\k)=0\quad\text{ for $i<\dim\sK$}.
\]
\end{proposition}
\begin{proof}
Let $I\in\sK$. If $I=\varnothing$, then $\widetilde{H}_i(\mathcal
K;\k)=\widetilde{H}_i(\lk I;\k)$. If $I\ne\varnothing$, then
\begin{equation}\label{Klink}
  H_i(|\mathcal K|,|\mathcal K|\!\setminus x;\k)=\widetilde{H}_{i-|I|}(\lk I;\k)
\end{equation}
for any $x$ in the interior of $I$, by Proposition~\ref{lclin}.

If $\sK$ is Cohen--Macaulay, then it is pure
(Exercise~\ref{cmpure}) and therefore $\lk I$ is pure of dimension
$\dim\sK-|I|$ (Exercise~\ref{purelink}). Hence, $i<\dim\sK$
implies that $i-|I|<\dim\lk I$ and $H_i(|\mathcal K|,|\mathcal
K|\!\setminus x;\k)=0$ by~\eqref{Klink} and Theorem~\ref{reisner}.

On the other hand, if $H_i(|\mathcal K|,|\mathcal
  K|\!\setminus x;\k)=0$ for $i<\dim\sK$, then
$\widetilde{H}_j(\lk I;\k)=H_{j+|I|}(|\mathcal K|,|\mathcal
K|\!\setminus x;\k)=0$ for $j<\dim\lk I$, since $j+|I|<\dim\lk
I+|I|\le\dim\sK$. Thus, $\sK$ is Cohen--Macaulay by
Theorem~\ref{reisner}.
\end{proof}

\begin{corollary}
If a triangulation of a space $X$ is a Cohen--Macaulay complex,
then any other triangulation of $X$ is Cohen--Macaulay as well.
\end{corollary}

\begin{corollary}\label{spherecm}
Any triangulated sphere is a Cohen--Macaulay complex.
\end{corollary}

In particular, the $h$-vector of a triangulated sphere is an
$M$-vector. This fact was used by Stanley in his generalisation of
the UBT (Theorem~\ref{ubt}) to arbitrary sphere triangulations:

\begin{theorem}[UBT for spheres]
\label{ubtss} For any triangulated $(n-1)$-dimensional sphere
$\mathcal K$ with $m$ vertices, the $h$-vector
$(h_0,h_1,\ldots,h_n)$ satisfies the inequalities
$$
  h_i(\mathcal K)\le\bin{m-n+i-1}i.
$$
Therefore, the UBT holds for triangulated spheres, that is,
$$
  f_i(\mathcal K)\le f_i\bigl(C^n(m)\bigr) \quad\text{ for }i=1,\ldots,n-1.
$$
\end{theorem}
\begin{proof}
Let $A=\k[\sK]/\mb t$ be the algebra from the proof of
Proposition~\ref{stanmv}, so that $\dim_\k A^{2i}=h_i$. In
particular, $\dim_\k A^2=h_1=m-n$. Since $A$ is generated
by~$A^2$, the number $h_i$ cannot exceed the number of monomials
of degree $i$ in $m-n$ generators, i.e. $h_i\le\bin{m-n+i-1}i$.
The rest follows from Lemma \ref{ubth}.
\end{proof}

\subsection*{Exercises}
\begin{exercise}
If $\sK$ is a simplicial complex of dimension $n-1$, then
$\dim\k[\sK]=n$.
\end{exercise}

\begin{exercise}
Let $\mb t$ be an lsop in $\k[\sK]$. Then the $\k$-vector space
$\k[\sK]/\mb t$ is generated by monomials $v_I$ for $I\in\sK$.
(Hint: prove that $v_iv_I=0$ in $\k[\sK]/\mb t$ for any $i\in[m]$
and for any \emph{maximal} simplex $I\in\sK$, and then use
Proposition~\ref{frbasis}).
\end{exercise}

\begin{exercise}
A sequence $\mb t\subset\k[m]$ is $\k[m]$-regular for $\k[\sK]$ as
a $\k[m]$-module if and only the image of $\mb t$ in $\k[\sK]$ is
$\k[\sK]$-regular.
\end{exercise}

\begin{exercise}
Let $\mb t=(t_1,\ldots,t_n)\subset\Z[\sK]$ be a linear sequence,
$\dim\sK=n-1$. The following conditions are equivalent:
\begin{itemize}
\item[(a)] the reduction of $\mb t$ modulo $p$ is an lsop in
$\Z_p[\sK]$ for any prime~$p$;
\item[(b)] the restriction $s_I(\mb t)$
to each simplex $I\in\sK$ generates the polynomial ring
$\Z[v_i\colon i\in I]$;
\item[(c)] for each $I\in\sK$ the columns of the $n\times|I|$
matrix $(\lambda_{ij})$, \ $1\le i\le n$, $j\in I$, generate the
integer lattice~$\Z^{|I|}$.
\end{itemize}
\end{exercise}

\begin{exercise}\label{frcomplint}
A finitely generated commutative $\k$-algebra is called a
\emph{complete intersection algebra} if it is the quotient of a
polynomial algebra by a regular sequence. Observe that a complete
intersection algebra is Cohen--Macaulay. Show that a face ring
$\k[\sK]$ is a complete intersection algebra if and only if it is
isomorphic to the quotient of the form
\[
  \k[v_1,\ldots,v_m]/(v_1v_2\!\cdots v_{k_1},
  v_{k_1+1}v_{k_1+2}\!\cdots v_{k_1+k_2},\ldots,
  v_{k_1+\cdots+k_{p-1}+1}\!\cdots v_{k_1+\cdots+k_p}).
\]
This is equivalent to $\sK$ being decomposable into a join of the
form
\[
  \partial\varDelta^{k_1-1}*\partial\varDelta^{k_2-1}*\cdots*
  \partial\varDelta^{k_p-1}*\varDelta^{m-s-1},
\]
where $s=k_1+\cdots+k_p$ and the join factor $\varDelta^{m-s-1}$
is void if $s=m$.
\end{exercise}

\begin{exercise}\label{cmpure}
A Cohen--Macaulay complex is pure. (Hint: given a maximal simplex
$J\in\sK$, consider the ideal $\mathcal I_J=(v_i\colon i\notin J)$
in~$\k[m]$,
%(the \emph{associated prime ideal} of $\k[\sK]$ corresponding to~$J$),
and show that
\[
  \depth\k[\sK]\le\dim(\k[m]/\mathcal I_J)=\dim\k[v_j\colon j\in
  J]=|J|.)
\]
\end{exercise}

\section[Gorenstein complexes]{Gorenstein complexes and
Dehn--Sommerville relations}\label{gc} Gorenstein rings are a
class of Cohen--Macaulay rings with a special duality property. As
in the case of Cohen--Macaulayness, simplicial complexes whose
face rings are Gorenstein play an important role in combinatorial
commutative algebra. In a sense, non-acyclic Gorenstein complexes
provide a `best possible algebraic approximation' to sphere
triangulations. We review here the most important aspects of
Gorenstein complexes. The proofs of the main results of this
section, in particular Theorems~\ref{gorencom} and~\ref{av-go},
require considerably more commutative algebraic techniques than
those from the previous sections. We refer the reader
to~\cite[Chapter~3]{br-he98} for the general theory of Gorenstein
rings and the missing proofs.

We recall from Proposition~\ref{cmbet} that nonzero Betti numbers
of a Cohen--Macaulay complex $\sK$ of dimension $n-1$ with $m$
vertices appear up to homological degree $-(m-n)$, and
$\beta^{-(m-n)}(\k[\sK])\ne0$.

\begin{definition}
\label{goren} A Cohen--Macaulay complex $\mathcal K$ of dimension
$n-1$ with $m$ vertices is called \emph{Gorenstein} (over a field
$\k$) if $\beta^{-(m-n)}(\k[\mathcal K])=1$, that is, if
$\Tor^{-(m-n)}_{\k[m]}(\k[\mathcal K],\k)\cong\k$. Furthermore,
$\sK$ is \emph{Gorenstein*} if $\sK$ is Gorenstein and $\mathcal
K=\core \mathcal K$ (see Definition~\ref{defcore}).
\end{definition}

Since $\mathcal K=\core(\mathcal K)*\varDelta^{s-1}$ for some~$s$,
we have $\k[\mathcal K]=\k[\core(\mathcal K)]\otimes\k[s]$. Then
Lemma~\ref{tortor} implies that
$$
  \Tor^{-i}_{\k[m]}\bigl(\k[\mathcal K],\k\bigr)\cong
  \Tor^{-i}_{\k[m-s]}\bigl(\k[\core \mathcal K],\k\bigr).
$$
Therefore, $\mathcal K$ is Gorenstein if and only if $\core
\mathcal K$ is Gorenstein*.

As in the case of Cohen--Macaulay complexes, Gorenstein* complexes
can be characterised topologically as follows.

\begin{theorem}[{\cite[\S II.5]{stan96} or \cite[Theorem~5.6.1]{br-he98}}]
\label{gorencom} The following conditions are equivalent:
\begin{itemize}
\item[(a)] $\mathcal K$ is a Gorenstein* complex over~$\k$;
\item[(b)] for any simplex $I\in \mathcal K$ (including
$I=\varnothing$) the subcomplex $\lk I$ has homology of a sphere
of dimension $\dim(\lk I)$;
\item[(c)] for any $x\in|\mathcal K|$,
\[
  H^i\bigl(|\mathcal K|,|\mathcal
  K|\!\setminus x;\k\bigr)=\widetilde{H}^i(\mathcal K;\k)=\begin{cases}\k&\text{if \ $i=\dim\sK$;}\\
  0&\text{otherwise.}\end{cases}
\]
\end{itemize}
\end{theorem}

In topology, polyhedra $|\sK|$ satisfying the conditions of the
previous theorem are sometimes called \emph{generalised homology
spheres} (`generalised' because a homology sphere is usually
assumed to be a manifold). In particular, triangulated spheres are
Gorenstein* complexes. Triangulated manifolds are not Gorenstein*
or even Cohen--Macaulay in general (\emph{Buchsbaum
complexes}\label{buchsbaumcplx} provide a proper algebraic
approximation to triangulated manifolds,
see~\cite[\S~II.8]{stan96}). Nevertheless, the $\Tor$-algebra of a
Gorenstein* complex behaves like the cohomology algebra of a
manifold: it satisfies \emph{Poincar\'e duality}. This fundamental
result was proved by Avramov and Golod for Noetherian local rings;
here we state the graded version of their theorem in the case of
face rings.

\begin{definition}\label{pdalgebra}
A graded commutative connected $\k$-algebra $A$ is called a
\emph{Poincar\'e algebra} if it is finite dimensional over $\k$,
i.e. $A=\bigoplus_{i=0}^dA^i$, and the $\k$-linear maps
\begin{align*}
  A^i&\to\Hom_\k(A^{d-i},A^d),\\
  a&\mapsto\phi_a,\quad\text{where } \phi_a(b)=ab
\end{align*}
are isomorphisms for $0\le i\le d$. The classical example of a
Poincar\'e algebra is the cohomology algebra of a manifold.
\end{definition}

\begin{theorem}[Avramov--Golod, {\cite[Theorem~3.4.5]{br-he98}}]
\label{av-go} A simplicial complex $\sK$ is Gorenstein* if and
only if the algebra $T=\bigoplus_{i=0}^d T^i$, where
$T^i=\Tor^{-i}_{\k[m]}(\k[\sK],\k)$ and $d=\max\{ j\colon
\Tor^{-j}_{\k[m]}(\k[\sK],\k)\ne0\}$, is a Poincar\'e algebra.
\end{theorem}

\begin{corollary}\label{tordual}
Let $\mathcal K$ be a Gorensten* complex of dimension $n-1$ on the
set~$[m]$. Then the Betti numbers and the Poincar\'e series of the
$\Tor$ groups satisfy
\begin{gather*}
  \beta^{-i,\,2j}\bigl(\k[\mathcal K]\bigr)
  =\beta^{-(m-n)+i,\,2(m-j)}\bigl(\k[\mathcal K]\bigr),\quad
  0\le i\le m-n,\; 0\le j\le m,
  \\
  F\bigl(\Tor^{-i}_{\k[m]}(\k[\mathcal K],\k);\;\lambda\bigr)=
  \lambda^{2m}F\bigl(\Tor^{-(m-n)+i}_{\k[m]}
  (\k[\mathcal K],\k);\;{\textstyle\frac 1\lambda}\bigr).
\end{gather*}
\end{corollary}
\begin{proof}
Theorems \ref{hoch} and \ref{gorencom} imply that
\[
%\begin{equation}\label{gstar}
  \beta^{-(m-n)}\bigl(\k[\mathcal K]\bigr)=\beta^{-(m-n),2m}\bigl(\k[\mathcal K]\bigr)=1.
%\end{equation}
\]
We therefore have $d=m-n$ and $T^d=
\Tor^{-(m-n),\,2m}_{\k[m]}(\k[\sK],\k)\cong\k$ in the notation of
Theorem~\ref{av-go}. Since the multiplication in the
$\Tor$-algebra preserves the bigrading, the isomorphisms
$T^i\stackrel\cong\to\Hom_\k(T^{m-n-i},T^{m-n})$ from the
definition of a Poincar\'e algebra can be refined to isomorphisms
\[
  T^{i,\,2j}\stackrel\cong\longrightarrow
  \Hom_\k(T^{m-n-i,\,2(m-j)},T^{m-n,\,2m}),
\]
where $T^{m-n,\,2m}\cong\k$. This implies the first identity, and
the second is a direct corollary.
\end{proof}

As a further corollary we obtain the following symmetry property
for the Poincar\'e series of the face ring:

\begin{corollary}
\label{frdual} If $\mathcal K$ is Gorenstein* of dimension $n-1$,
then
$$
  F\bigl(\k[\mathcal K],\lambda\bigr)=(-1)^nF\bigl(\k[\mathcal K],{\textstyle\frac1{\lambda}}\bigr).
$$
\end{corollary}
\begin{proof}
We apply Proposition~\ref{psresol} to the minimal resolution of
the $\k[m]$-module $\k[\mathcal K]$. Note that
$F(\k[m];\lambda)=(1-\lambda^2)^{-m}$. It follows from the formula
for the Poincar\'e series from Proposition~\ref{psresol} and
Proposition~\ref{mintor} that
$$
  F\bigl(\k[\mathcal K];\lambda\bigr)=(1-\lambda^2)^{-m}\sum_{i=0}^{m-n}(-1)^i
  F\Bigl(\Tor^{-i}_{\k[m]}\bigl(\k[\mathcal K],\k\bigr);\lambda\Bigr).
$$

Using Corollary~\ref{tordual}, we calculate
\begin{multline*}
  F\bigl(\k[\mathcal K];\lambda\bigr)=(1-\lambda^2)^{-m}\sum_{i=0}^{m-n}(-1)^i \lambda^{2m}
  F\Bigl(\Tor^{-(m-n)+i}_{\k[m]}\bigl(\k[\mathcal K],\k\bigr);
  {\textstyle\frac1{\lambda}}\Bigr)=\\
  =\bigl(1-({\textstyle\frac1{\lambda}})^2\bigr)^{-m}
  (-1)^m\sum_{j=0}^{m-n}(-1)^{m-n-j}
  F\Bigl(\Tor^{-j}_{\k[m]}\bigl(\k[\mathcal K],\k\bigr);
  {\textstyle\frac1\lambda}\Bigr)=\\
  =(-1)^nF\bigl(\k[\mathcal
  K];{\textstyle\frac1{\lambda}}\bigr).\qedhere
\end{multline*}
\end{proof}

\begin{corollary}
\label{dsgor} The Dehn--Sommerville relations $h_i=h_{n-i}$ hold
for any Gorenstein* complex of dimension $n-1$ (in particular, for
any triangulation of an $(n-1)$-sphere).
\end{corollary}
\begin{proof}
This follows from the explicit form of the Poincar\'e series for
$\k[\mathcal K]$ (Theorem~\ref{psfr}) and the previous corollary.
\end{proof}

The Dehn--Sommerville relations may be further generalised to
wider classes of complexes and posets; we give some of these
generalisations in Section~\ref{gendsr} below.

Unlike the situation with Cohen--Macaulay complexes, a
characterisation of $h$-vectors (or, equivalently, $f$-vectors) of
Gorenstein complexes is not known:

\begin{problem}[Stanley]
Characterise the $h$-vectors of Gorenstein complexes.
\end{problem}

If the \emph{$g$-conjecture}\label{gconje} (i.e. the inequalities
of Theorem~\ref{gth}~(b) and~(c)) holds for Gorenstein complexes,
then this would imply a solution to the above problem.

\subsection*{Exercises}
\begin{exercise}
A Gorenstein complex $\sK$ is Gorenstein* if and only if it is
\emph{non-acyclic} (i.e. $\widetilde H^*(\sK;\k)\ne0$).
\end{exercise}

\begin{exercise}[{\cite[Theorem~II.5.1]{stan96}}]\label{gorchi} Show that $\sK$ is Gorenstein* if and only if
the following three conditions are satisfied:
\begin{itemize}
\item[(a)] $\sK$ is Cohen--Macaulay;
\item[(b)] every $(n-2)$-dimensional simplex is contained in exactly two
$(n-1)$-\-di\-men\-si\-o\-nal simplices;
\item[(c)] $\chi(\sK)=\chi(S^{n-1})$, where $\chi(\:\cdot\:)$ denotes the Euler characteristic.
\end{itemize}
\end{exercise}

\begin{exercise}
Let $\sK$ be a Gorenstein* complex of dimension $n-1$ on~$[m]$.
Show that
\[
  \widetilde H^k(\sK_J)\cong\widetilde H^{n-2-k}(\sK_{\widehat J})
\]
for any $J\subset[m]$, where $\widehat J=[m]\setminus J$. This is
known as \emph{Alexander duality for non-acyclic Gorenstein
complexes}\label{alexadua}.
\end{exercise}

\section{Face rings of simplicial posets}\label{frsimpo}
The whole theory of face rings may be extended to simplicial
posets (defined in Section~\ref{secsimpos}), thereby leading to
new important classes of rings in combinatorial commutative
algebra and applications in toric topology.

The \emph{face ring} $\k[\sS]$ of a simplicial poset $\sS$ was
introduced by Stanley~\cite{stan91} as a quotient of a certain
graded polynomial ring by a homogeneous ideal determined by the
poset relation in~$\sS$.
% (see Definition~\ref{frspo} below).
The rings~$\k[\sS]$ have remarkable algebraic and homological
properties, albeit they are much more complicated than the
Stanley--Reisner face rings~$\k[\sK]$. Unlike~$\k[\sK]$, the
ring~$\k[\sS]$ is not generated in the lowest positive degree.
Face rings of simplicial posets were further studied by
Duval~\cite{duva97} and Maeda--Masuda--Panov \cite{ma-pa06},
\cite{m-m-p07}, among others. \emph{Cohen--Macaulay} and
\emph{Gorenstein*} face rings are particularly important; both
properties are topological, that is, depend only on the
topological type of the geometric realisation~$|\sS|$.

As usual, we shall not distinguish between simplicial posets
$\mathcal S$ and their geometric realisations (simplicial cell
complexes) $|\mathcal S|$. Given two elements
$\sigma,\tau\in\mathcal S$, we denote by $\sigma\vee\tau$ the set
of their joins, and denote by $\sigma\wedge\tau$ the set of their
meets\label{joinmee}. Whenever either of these sets consists of a
single element, we use the same notation for this particular
element of~$\mathcal S$.

To make clear the idea behind the definition of the face ring of a
simplicial poset, we first consider the case when $\mathcal S$ is
a simplicial complex~$\mathcal K$. Then $\sigma\wedge\tau$
consists of a single element (possibly $\varnothing$), and
$\sigma\vee\tau$ is either empty or consists of a single element.
We consider the graded polynomial ring $\k[v_\sigma\colon
\sigma\in \mathcal K]$ with one generator $v_\sigma$ of degree
$\deg v_\sigma=2|\sigma|$ for each simplex $\sigma\in\sK$. The
following proposition provides an alternative presentation of the
face ring $\k[\mathcal K]$, with a larger set of generators:

\begin{proposition}\label{frext}
There is a canonical isomorphism of graded rings
$$
  \k[\mathcal K]\cong
  \k[v_\sigma\colon\sigma\in \mathcal K]/\mathcal I'_{\sK},
$$
where $\mathcal I'_{\sK}$ is the ideal generated by the element
$v_\varnothing-1$ and all elements of the form
$$
  v_\sigma v_\tau-v_{\sigma\wedge\tau}v_{\sigma\vee\tau}.
$$
Here we set $v_{\sigma\vee\tau}=0$ whenever $\sigma\vee\tau$ is
empty.
\end{proposition}
\begin{proof}
The isomorphism is established by the map taking $v_\sigma$ to
$\prod_{i\in\sigma}v_i$. The rest is left as an exercise.
\end{proof}

Now let $\mathcal S$ be an arbitrary simplicial poset with the
vertex set $V(\sS)=[m]$. In this case both $\sigma\vee\tau$ and
$\sigma\wedge\tau$ may consist of more than one element, but
$\sigma\wedge\tau$ consists of a single element whenever
$\sigma\vee\tau$ is nonempty.

We consider the graded polynomial ring $\k[v_\sigma\colon
\sigma\in \mathcal S]$ with one generator $v_\sigma$ of degree
$\deg v_\sigma=2|\sigma|$ for every element $\sigma\in\sS$.

\begin{definition}[\cite{stan91}]\label{frspo}
The \emph{face ring} of a simplicial poset $\mathcal S$ is the
quotient
\[
  \k[\mathcal S]=
  \k[v_\sigma\colon \sigma\in\mathcal S]\;/\;
  \mathcal I_{\mathcal S},
\]
where $\mathcal I_{\mathcal S}$ is the ideal generated by the
elements $v_{\hatzero}-1$ and
\begin{equation}\label{estre}
  v_\sigma v_\tau-v_{\sigma\wedge\tau}\cdot\!\!\!\sum_{\eta\in\sigma\vee\tau}\!\!\!v_\eta.
\end{equation}
The sum over the empty set is zero, so we have $v_\sigma v_\tau=0$
in $\k[\sS]$ if $\sigma\vee\tau$ is empty.

The grading may be refined to an $\N^m$-grading by setting
$\mathop{\mathrm{mdeg}}v_\sigma=2V(\sigma)$. Here
$V(\sigma)\subset[m]$ is the vertex set of $\sigma$, and we
identify subsets of $[m]$ with $(0,1)$-vectors in
$\{0,1\}^m\subset\N^m$ as usual. In particular,
$\mathop{\mathrm{mdeg}} v_i=2\mb e_i$.
%(two times the $i$th basis vector).
\end{definition}

\begin{example}\label{exsp}\
\begin{figure}
\begin{picture}(0,0)%
\includegraphics{sposet1.pstex}%
\end{picture}%
\setlength{\unitlength}{1973sp}%
\begingroup\makeatletter\ifx\SetFigFont\undefined%
\gdef\SetFigFont#1#2#3#4#5{%
  \reset@font\fontsize{#1}{#2pt}%
  \fontfamily{#3}\fontseries{#4}\fontshape{#5}%
  \selectfont}%
\fi\endgroup%
\begin{picture}(11769,2907)(676,-5149)
\put(10801,-4036){\makebox(0,0)[lb]{\smash{{\SetFigFont{8}{9.6}{\rmdefault}{\mddefault}{\updefault}
$\sigma$}}}}
\put(4276,-3436){\makebox(0,0)[lb]{\smash{{\SetFigFont{8}{9.6}{\rmdefault}{\mddefault}{\updefault}
$2$}}}}
\put(676,-3436){\makebox(0,0)[lb]{\smash{{\SetFigFont{8}{9.6}{\rmdefault}{\mddefault}{\updefault}
$1$}}}}
\put(2551,-2461){\makebox(0,0)[lb]{\smash{{\SetFigFont{8}{9.6}{\rmdefault}{\mddefault}{\updefault}
$\tau$}}}}
\put(2476,-4336){\makebox(0,0)[lb]{\smash{{\SetFigFont{8}{9.6}{\rmdefault}{\mddefault}{\updefault}
$\sigma$}}}}
\put(2251,-5086){\makebox(0,0)[lb]{\smash{{\SetFigFont{8}{9.6}{\rmdefault}{\mddefault}{\updefault}
(a)   \ \ \ $ r=2$}}}}
\put(9301,-5011){\makebox(0,0)[lb]{\smash{{\SetFigFont{8}{9.6}{\rmdefault}{\mddefault}{\updefault}
(b) \ \ \ $ r=3$}}}}
\put(6751,-3436){\makebox(0,0)[lb]{\smash{{\SetFigFont{8}{9.6}{\rmdefault}{\mddefault}{\updefault}
$1$}}}}
\put(12076,-2461){\makebox(0,0)[lb]{\smash{{\SetFigFont{8}{9.6}{\rmdefault}{\mddefault}{\updefault}
$3$}}}}
\put(9976,-4036){\makebox(0,0)[lb]{\smash{{\SetFigFont{8}{9.6}{\rmdefault}{\mddefault}{\updefault}
$2$}}}}
\put(8626,-2386){\makebox(0,0)[lb]{\smash{{\SetFigFont{8}{9.6}{\rmdefault}{\mddefault}{\updefault}
$\tau$}}}}
\put(8476,-3586){\makebox(0,0)[lb]{\smash{{\SetFigFont{8}{9.6}{\rmdefault}{\mddefault}{\updefault}
$e$}}}}
\end{picture}%
\centering
   \caption{Simplicial cell complexes.} \label{sccfig}
        \end{figure}

1. The simplicial cell complex shown in Fig.~\ref{sccfig}~(a) is
obtained by gluing two segments along their boundaries and has
rank~$2$. The vertices are $1,2$ and we denote the 1-dimensional
simplices by $\sigma$ and $\tau$. Then the face ring $\k[\sS]$ is
the quotient of the graded polynomial ring
\[
  \k[v_1,v_2,v_\sigma,v_\tau],\quad\deg v_1= \deg v_2=2,\quad
  \deg v_\sigma=\deg v_\tau=4
\]
by the two relations
\[
  v_1v_2=v_\sigma+v_\tau,\quad v_\sigma v_\tau=0.
\]

2. The simplicial cell complex in Fig.~\ref{sccfig}~(b) is
obtained by gluing two triangles along their boundaries and has
rank~$3$. The vertices are $1,2,3$ and we denote the 1-dimensional
simplices (edges) by $e$, $f$ and $g$, and the 2-dimensional
simplices by $\sigma$ and~$\tau$. The face ring $\k[\sS]$ is
isomorphic to the quotient of the polynomial ring
\[
  \k[v_1,v_2,v_3,v_\sigma,v_\tau],\quad\deg v_1= \deg v_2=\deg v_3=2,
  \quad \deg v_\sigma=\deg v_\tau=6
\]
by the two relations
\[
  v_1v_2v_3=v_\sigma+v_\tau,\quad v_\sigma v_\tau=0.
\]
The generators corresponding to the edges can be excluded because
of the relations $v_e=v_1v_2$, $v_f=v_2v_3$ and $v_g=v_1v_3$.
\end{example}

\begin{remark}\label{straighrel}
The ideal $\mathcal I_{\mathcal S}$ is generated by
\emph{straightening relations}~\eqref{estre}; these relations
allow us to express the product of any pair of generators via
products of generators corresponding to pairs of ordered elements
of the poset. This can be restated by saying that $\k[\sS]$ is an
example of an \emph{algebra with straightening law} (ASL for
short, also known as a \emph{Hodge algebra}). Lemma~\ref{order}
and Theorem~\ref{chdec} below reflect algebraic properties of
ASL's, and may be restated in this generality. For more on the
theory of ASL's see~\cite[\S~III.6]{stan96}
and~\cite[Chapter~7]{br-he98}.
\end{remark}

A monomial $v_{\sigma_1}^{i_1}v_{\sigma_2}^{i_2}\cdots
v_{\sigma_k}^{i_k}\in\k[v_\sigma\colon\sigma\in\sS]$ is
\emph{standard} if $\sigma_1<\sigma_2<\cdots<\sigma_k$.

\begin{lemma}\label{order}
Any element of $\k[\sS]$ can be written as a linear combination of
standard monomials.
%$v_{\sigma_1}^{i_1}v_{\sigma_2}^{i_2}\cdots
%v_{\sigma_k}^{i_k}$ corresponding to chains of totally ordered
%elements $\sigma_1<\sigma_2<\cdots<\sigma_k$ of~$\sS$.
%of~$\sS\!\setminus\!\hatzero$.
\end{lemma}
\begin{proof}
It is enough to prove the statement for elements of $\k[\sS]$
represented by monomials in generators~$v_\sigma$. We write such a
monomial as $a=v_{\tau_1}v_{\tau_2}\cdots v_{\tau_k}$ where some
of the $\tau_i$ may coincide. We need to show that any such
monomial can be expressed as a sum of monomials $\sum
v_{\sigma_1}\cdots v_{\sigma_l}$ with
$\sigma_1\le\cdots\le\sigma_l$. We may assume by induction that
$\tau_2\le\cdots\le\tau_k$. Using relation~\eqref{estre} we can
replace $a$ by
$$
  v_{\tau_1\wedge\tau_2}\Bigl(\sum_{\rho\in\tau_1\vee
  \tau_2}\!\!\!v_\rho\Bigr) v_{\tau_3}\cdots v_{\tau_k}.
$$
Now the first two factors in each summand above correspond to
ordered elements of~$\sS$. We proceed by replacing the products
$v_\rho v_{\tau_3}$ by
$v_{\rho\wedge\tau_3}(\sum_{\pi\in\rho\vee\tau_3}v_\pi)$. Since
$\tau_1\wedge\tau_2\le\rho\wedge\tau_3$, now the first three
factors in each monomial are in order. Continuing this process, we
obtain in the end a sum of monomials corresponding to totally
ordered sets of elements of~$\sS$.
\end{proof}

We refer to the presentation from Lemma~\ref{order} as a
\emph{standard representation} of an element $a\in\k[\mathcal S]$.

Given $\sigma\in\mathcal S$, we define the corresponding
\emph{restriction homomorphism} as
$$
  s_\sigma\colon\k[\mathcal S]\to\k[\mathcal S]/(v_\tau\colon
  \tau\not\le\sigma).
$$
The following result is straightforward.

\begin{proposition}\label{restr}
Let $|\sigma|=k$ with $V(\sigma)=\{i_1,\ldots,i_k\}$. Then the
image of the homomorphism $s_\sigma$ is the polynomial ring
$\k[v_{i_1},\ldots,v_{i_k}]$.
\end{proposition}

The next result generalises Proposition~\ref{rmap} to simplicial
posets.

\begin{theorem}\label{alres}
The direct sum
\[
  s={\textstyle\bigoplus_{\sigma\in\sS}} s_\sigma\colon\k[\sS]\longrightarrow
  \bigoplus_{\sigma\in\sS}\k[v_i\colon i\in V(\sigma)]
\]
of all restriction maps is a monomorphism.
\end{theorem}
\begin{proof}
Take a nonzero element $a\in\k[\mathcal S]$ and write its standard
representation. Fix a standard monomial $v_{\sigma_1}^{i_1}\cdots
v_{\sigma_k}^{i_k}$ which enters into this decomposition with a
nonzero coefficient. Allowing some of the exponents $i_j$ to be
zero, we may assume that $\sigma_k$ is a maximal element in
$\mathcal S$ and $|\sigma_j|=j$ for $1\le j\le k$. We shall prove
that $s_{\sigma_k}(a)\ne0$. Identify $s_{\sigma_k}(\k[\mathcal
S])$ with the polynomial ring $\k[t_1,\ldots,t_k]$ (so that
$t_j=v_{i_j}$ in the notation of Proposition~\ref{restr}). Then
$s_{\sigma_k}(v_{\sigma_k})=t_1\cdots t_k$, and we may assume
without loss of generality that
$s_{\sigma_k}(v_{\sigma_j})=t_1\cdots t_j$ for $1\le j\le k$.
Hence,
$$
  s_{\sigma_k}\bigl( v_{\sigma_1}^{i_1}\cdots v_{\sigma_k}^{i_k} \bigr)=
  t_1^{i_1}(t_1t_2)^{i_2}\cdots(t_1\cdots t_k)^{i_k}.
$$
If we prove that no other monomial $v_{\tau_1}^{j_1}\cdots
v_{\tau_m}^{j_m}$ is mapped by $s_{\sigma_k}$ to the same element
of $\k[t_1,\ldots,t_k]$, then this would imply that
$s_{\sigma_k}(a)\ne0$. Note that
$$
  s_{\sigma_k}(v_{\tau_1}^{j_1}\cdots v_{\tau_m}^{j_m})=0\qquad
  \text{if } \tau_i\not\le\sigma_k\text{ for some }i
  \text{ with } j_i\ne0,
$$
so that we may assume that $m=k$. Now suppose that
\begin{equation}\label{rpvgh}
  s_{\sigma_k}\bigl( v_{\sigma_1}^{i_1}\cdots v_{\sigma_k}^{i_k} \bigr)=
  s_{\sigma_k}\bigl( v_{\tau_1}^{j_1}\cdots v_{\tau_k}^{j_k} \bigr).
\end{equation}
We shall prove that $v_{\sigma_1}^{i_1}\cdots v_{\sigma_k}^{i_k}=
v_{\tau_1}^{j_1}\ldots v_{\tau_k}^{j_k}$. We may assume by
induction that the `tails' of these monomial coincide, that is,
there is some $q$, $1\le q\le k$, such that $i_p=j_p$ and
$\sigma_p=\tau_p$ for $i_p\ne0$ whenever $p>q$. We shall prove
that $i_q=j_q$ and $\sigma_q=\tau_q$ if $i_q\ne0$. We obtain
from~\eqref{rpvgh} that
\begin{multline*}
  s_{\sigma_k}\bigl( v_{\sigma_1}^{i_1}\cdots v_{\sigma_q}^{i_q} \bigr)
  (t_1\cdots t_{q+1})^{i_{q+1}}\cdots(t_1\cdots t_k)^{i_k}=\\
  =s_{\sigma_k}\bigl( v_{\tau_1}^{j_1}\cdots v_{\tau_q}^{j_q} \bigr)
  (t_1\cdots t_{q+1})^{i_{q+1}}\cdots(t_1\cdots
  t_k)^{j_k},
\end{multline*}
hence, $s_{\sigma_k}(v_{\sigma_1}^{i_1}\cdots v_{\sigma_q}^{i_q})=
s_{\sigma_k}(v_{\tau_1}^{j_1}\cdots v_{\tau_q}^{j_q})$. Let $j_l$
be the last nonzero exponent in $v_{\tau_1}^{j_1}\cdots
v_{\tau_q}^{j_q}$ (i.e. $j_{l+1}=\cdots=j_q=0$). Then we also have
$i_{l+1}=\cdots=i_q=0$, as otherwise
$s_{\sigma_k}(v_{\sigma_1}^{i_1}\cdots v_{\sigma_q}^{i_q})$ is
divisible by $t_1\cdots t_{l+1}$, while
$s_{\sigma_k}(v_{\tau_1}^{j_1}\cdots v_{\tau_q}^{j_q})$ is not. We
also have $i_l=j_l$ and $\sigma_l=\tau_l$, since $i_l$ is the
maximal power of the monomial $t_1\cdots t_l$ which divides
$s_{\sigma_k}(v_{\sigma_1}^{i_1}\cdots v_{\sigma_q}^{i_q})$. We
conclude by induction that $v_{\sigma_1}^{i_1}\cdots
v_{\sigma_k}^{i_k}= v_{\tau_1}^{j_1}\cdots v_{\tau_k}^{j_k}$, and
$s_{\sigma_k}(a)\ne0$.
\end{proof}

\begin{remark}
The proof above also shows that the map $s=\bigoplus_\sigma
s_\sigma$ in Theorem~\ref{alres} can be defined as the sum over
only the maximal elements $\sigma\in\mathcal S$.
\end{remark}

\begin{theorem}\label{chdec}
The standard representation of an element $a\in\k[\mathcal S]$ is
unique. In other words, the standard monomials
$v_{\sigma_1}^{i_1}\cdots v_{\sigma_k}^{i_k}$ form a $\k$-basis
of~$\k[\mathcal S]$.
\end{theorem}
\begin{proof}
This follows directly from Lemma~\ref{order} and
Theorem~\ref{alres}.
\end{proof}

Lemma~\ref{lsopr}, which describes hsop's in the face rings of
simplicial complexes, can be readily extended to simplicial posets
(the same proof based on the properties of the restriction map $s$
works):

\begin{lemma}\label{lsopscc}
Let $\mathcal S$ be a simplicial poset of rank \ $n$. A sequence
of homogeneous elements $\mb t=(t_1,\ldots,t_n)$ of $\k[\mathcal
S]$ is a homogeneous system of parameters if and only
\[
  \dim_\k\bigl(\k[v_i\colon i\in V(\sigma)]/s_\sigma(\mb t)\bigr)<\infty
\]
for each element $\sigma\in \mathcal S$.
\end{lemma}

The \emph{$f$-vector}\label{fvsimpose} of a simplicial poset
$\mathcal S$ is $\mb f(\mathcal S)=(f_0,\ldots,f_{n-1})$, where
$n-1=\dim\mathcal S$ and $f_i$ is the number of elements of rank
$i+1$ (i.e.~the number of faces of dimension~$i$ in the simplicial
cell complex). We also set $f_{-1}=1$. The \emph{$h$-vector} $\mb
h(\mathcal S)=(h_0,\ldots,h_n)$ is then defined
by~\eqref{hvectors}.

The Poincar\'e series of the face ring $\k[\sS]$ has exactly the
same form as in the case of simplicial complexes:

\begin{theorem}\label{psgfr}
We have
$$
  F\bigl(\k[\mathcal S];\lambda\bigr)
  =\sum_{k=0}^{n}\frac{f_{k-1}\lambda^{2k}}{(1-\lambda^2)^k}
  =\frac{h_0+h_1\lambda^2+\cdots+h_n\lambda^{2n}}{(1-\lambda^2)^n}.
$$
\end{theorem}
\begin{proof}
By Theorem \ref{chdec}, we need to calculate the Poincar\'e series
of the $\k$-vector space generated by the monomials
$v_{\sigma_1}^{i_1}\cdots v_{\sigma_k}^{i_k}$ with
$\sigma_1<\cdots<\sigma_k$. For every $\sigma\in\mathcal S$ denote
by $\mathcal M_\sigma$ the set of such monomials with
$\sigma_k=\sigma$ and $i_k>0$. Let $|\sigma|=k$; consider the
restriction homomorphism $s_{\sigma}$ to the polynomial ring
$\k[t_1,\ldots,t_k]$. Then $s_{\sigma}(\mathcal M_{\sigma})$ is
the set of monomials in $\k[t_1,\ldots,t_k]$ which are divisible
by $t_1\cdots t_k$. Therefore, the Poincar\'e series of the
subspace generated by the set $\mathcal M_{\sigma}$ is
$\frac{\lambda^{2k}}{(1-\lambda^2)^k}$. Now, to finish the proof
of the first identity we note that $\mathcal S$ is the union
$\bigcup_{\sigma\in\mathcal S}\mathcal M_\sigma$ of the
nonintersecting subsets~$\mathcal M_\sigma$. The second identity
follows from~\eqref{hvectors}.
\end{proof}

As we have seen in Exercise~\ref{frlimit}, the face ring $\k[\sK]$
of a simplicial complex can be realised as the limit of a diagram
of polynomial algebras over $\ca^{op}(\sK)$. A similar description
exists for the face ring $\k[\sS]$:

\begin{construction}[{$\k[\sS]$ as a limit}]\label{kslimit}
We consider the diagram $\k[\,\cdot\:]^\mathcal S$ similar to that
of Exercise~\ref{frlimit}:
\begin{align*}
  \k[\,\cdot\:]^\sS\colon\ca^{op}(\sS)&\longrightarrow\cga,\\
  \sigma&\longmapsto\k[v_i\colon i\in V(\sigma)],
\end{align*}
whose value on a morphism $\sigma\le\tau$ is the surjection
\[
  \k[v_i\colon i\in V(\tau)]\to\k[v_i\colon i\in V(\sigma)]
\]
sending each $v_i$ with $i\notin V(\sigma)$ to zero.

\begin{lemma}\label{zslim}
We have
\[
  \k[\sS]=\lim\k[\,\cdot\:]^\mathcal S
\]
where the limit is taken in the category $\cga$.
\end{lemma}
\begin{proof}
We set up a total order on the elements of~$\sS$ so that the rank
function does not decrease, and proceed by induction. We therefore
may assume the statement is proved for a simplicial
poset~$\mathcal T$, and need to prove it for $\sS$ which is
obtained from $\mathcal T$ by adding one element~$\sigma$. Then
$\sS_{<\sigma}=\{\tau\in\sS\colon\tau<\sigma\}$ is the face poset
of the boundary of the simplex~$\varDelta^\sigma$. Geometrically,
we may think of $|\sS|$ as obtained from $|\mathcal T|$ by
attaching one simplex $\varDelta^\sigma$ along its boundary (if
$|\sigma|=1$, then $\varDelta^\sigma$ is a single point, so
$|\sS|$ is a disjoint union of $|\mathcal T|$ and a point). We
therefore need to prove that the following is a pullback diagram:
\begin{equation}\label{stpb}
\begin{CD}
  \k[\sS] @>>> \k[\sS_{\le\sigma}]\\
  @VVV @VVV\\
  \k[\mathcal T] @>>> \k[\sS_{<\sigma}].
\end{CD}
\end{equation}
Here the vertical arrows map $v_\sigma$ to 0, while the horizontal
ones map $v_\tau$ to 0 for $\tau\not\le\sigma$. Denote by $A$ the
pullback of~\eqref{stpb} with $\k[\sS]$ dropped. We need to show
that the natural map $\k[\sS]\to A$ is an isomorphism.

Since the limits in $\cga$ are created in the underlying category
of graded $\k$-vector spaces, the space of $A$ is the direct sum
of $\k[\mathcal T]$ and $\k[\sS_{\le\sigma}]$ with the pieces
$\k[\sS_{<\sigma}]$ identified in both spaces. In other words,
\begin{equation}\label{3dec}
  A=T\oplus\k[\sS_{<\sigma}]\oplus S,
\end{equation}
where $T$ is the complement to $\k[\sS_{<\sigma}]$ in $\k[\mathcal
T]$, and $S$ is the complement to $\k[\sS_{<\sigma}]$ in
$\k[\sS_{\le\sigma}]$. By Theorem~\ref{chdec}, the space
$\k[\sS_{<\sigma}]$ has basis of standard monomials
$v_{\tau_1}^{j_1}v_{\tau_2}^{j_2}\cdots v_{\tau_k}^{j_k}$ with
$\tau_k<\sigma$. Similarly, $S$ has basis of those monomials with
$\tau_k=\sigma$ and $j_k>0$, while $T$ has basis of those
monomials with $\tau_k\not\le\sigma$ and $j_k>0$. Yet another
application of Theorem~\ref{chdec} gives a decomposition
of~$\k[\sS]$ identical to~\eqref{3dec}: a standard basis monomial
$v_{\tau_1}^{j_1}v_{\tau_2}^{j_2}\cdots v_{\tau_k}^{j_k}$ with
$j_k>0$ has either $\tau_k\not\le\sigma$, or $\tau_k<\sigma$, or
$\tau_k=\sigma$. These three possibilities map to $T$,
$\k[\sS_{<\sigma}]$ and $S$ respectively. It follows that
$\k[\sS]\to A$ is an isomorphism of $\k$-vector spaces. Since it
is an algebra map, it is also an isomorphism of algebras, thus
finishing the proof.
\end{proof}
\end{construction}

The description of $\k[\sS]$ as a limit has the following
important corollary, describing the functorial properties of the
face rings and generalising Proposition~\ref{frmap}.

\begin{proposition}\label{funct}
Let $f\colon\sS\to\mathcal T$ be a rank-preserving map of
simplicial posets. Define a homomorphism
\[
  f^*\colon\k[w_\tau\colon
  \tau\in\mathcal T]\to\k[v_\sigma\colon
  \sigma\in\mathcal S],\quad
  f^*(w_\tau)=\sum_{\sigma\in f^{-1}(\tau),\;|\sigma|=|\tau|}
  v_\sigma.
\]
Then $f^*$ descends to a ring homomorphism $\k[\mathcal
T]\to\k[\sS]$, which we continue to denote by~$f^*$.
\end{proposition}
\begin{proof}
The poset map $f$ gives rise to a functor
$f\colon\ca^{op}(\sS)\to\ca^{op}(\mathcal T)$ and therefore to a
natural transformation
\[
  f^*\colon[\ca^{op}(\mathcal T),\cga]\to[\ca^{op}(\sS),\cga],
\]
where $[\ca^{op}(\sS),\cga]$ denotes the functors from
$\ca^{op}(\sS)$ to $\cga$. It is easy to see that
$f^*\,\k[\,\cdot\,]_\mathcal T=\k[\,\cdot\,]_\mathcal S$ in the
notation of Construction~\ref{kslimit}, so we have the induced map
of limits $f^*\colon\k[\mathcal T]\to\k[\sS]$. We also have that
$f^*(w_\tau)=\sum_{\sigma\in f^{-1}(\tau)}v_\sigma$ by the
construction of $\lim$ in $\cga$.
\end{proof}

\begin{example}
The folding map~\eqref{fold} induces a monomorphism
$\k[\sK_\sS]\to\k[\sS]$, which embeds $\k[\sK_\sS]$ in $\k[\sS]$
as the subring generated by the elements~$v_i$.
\end{example}

\begin{remark}
An attempt to prove Proposition~\ref{funct} directly from the
definition, by showing that $f^*(\mathcal I_{\mathcal
T})\subset\mathcal I_\sS$, runs into a complicated combinatorial
analysis of the poset structure. This is an example of a situation
where the use of an abstract categorical description of $\k[\sS]$
proves to be beneficial.
\end{remark}

Let $\k[m]=\k[v_1,\ldots,v_m]$ be the polynomial algebra on $m$
generators of degree~2 corresponding to the vertices of~$\sS$. The
face ring $\k[\sS]$ acquires a $\k[m]$-algebra structure via the
map $\k[m]\to\k[\sS]$ sending each $v_i$ identically. (Unlike the
case of simplicial complexes, this map is generally not
surjective.) We thereby obtain a $\Z\oplus\N^m$-graded
$\Tor$-algebra of $\k[\sS]$:
\[
  \Tor_{\k[v_1,\ldots,v_m]}\bigl(\k[\sS],\k\bigr)=
  \bigoplus_{i\ge0,\mb a\in\N^m}
  \Tor^{-i,\,2\mb a}_{\k[v_1,\ldots,v_m]}\bigl(\k[\sS],\k\bigr),
\]
by analogy with Construction~\ref{mgrad} for simplicial complexes.

We finish this section by stating a generalisation of Hochster's
theorem to simplicial posets, and deriving some of its
corollaries.

\begin{theorem}[{Duval~\cite{duva97}, see also~\cite{lu-pa11}}]\label{hochsp}
For any subset $J\subset[m]$ we have
$$
  \Tor_{\k[v_1,\ldots,v_m]}^{-i,2J}\bigl(\k[\mathcal
  S],\k\bigr)\cong
  \widetilde{H}^{|J|-i-1}(|\mathcal S_J|;\k),
$$
where $\sS_J$ the subposet of $\sS$ consisting of those $\sigma$
for which $V(\sigma)\subset J$. Also, $\Tor_{\k[m]}^{-i,2\mb
a}(\k[\mathcal S],\k)=0$ if $\mb a$ is not a $(0,1)$-vector.
\end{theorem}
\begin{proof}
The argument follows the lines of the proof of
Theorem~\ref{hochmd}. We define the quotient differential graded
algebra
\[
  R^*(\sS)=\Lambda[u_1,\ldots,u_m]\otimes\k[\sS]/\mathcal I_R
\]
where $\mathcal I_R$ is the ideal generated by the elements
\[
  u_iv_\sigma\quad\text{with }i\in V(\sigma),\quad\text{and}\quad
  v_\sigma v_\tau\;\;\;\text{with }\sigma\wedge\tau\ne\hatzero.
\]
Note that the latter condition is equivalent to $V(\sigma)\cap
V(\tau)\ne\varnothing$.

Then we need to prove the analogue of Lemma~\ref{iscoh}, that is,
to show that the quotient projection
\[
  \varrho\colon\Lambda[u_1,\ldots,u_m]\otimes\k[\sS]\to R^*(\sS)
\]
induces an isomorphism in cohomology. This can be done by
providing the appropriate chain homotopy, as in the proof of
Lemma~\ref{iscoh}, but the formulae will be more complicated.
Alternatively, we can use a topological argument, see the proof of
Theorem~\ref{hzs2} below.

Let $\sC^{p-1}(|\mathcal S_J|)$ denote the $(p-1)$th cellular
cochain group of $|\mathcal S_J|$ with coefficients in~$\k$. It
has a basis of cochains $\alpha_\sigma$ corresponding to elements
$\sigma\in\mathcal S_J$ with $|\sigma|=p$. We define a $\k$-linear
map
\[
  f\colon\sC^{p-1}(|\mathcal S_J|)\longrightarrow
  R^{p-|J|,2J}(\sS),\quad
  \alpha_\sigma \longmapsto \varepsilon(V(\sigma),J)\,u_{J\setminus V(\sigma)}v_\sigma,
\]
where $\varepsilon(V(\sigma),J)$ is the sign from the proof of
Theorem~\ref{hochmd}. This map is an isomorphism of cochain
complexes; the details are left to the reader. Therefore,
\[
  \widetilde H^{p-1}(|\mathcal S_J|)\cong\Tor^{p-|J|,2J}_{\k[v_1,\ldots,v_m]}\bigl(\k[\mathcal
  S],\k\bigr),
\]
which is equivalent to the first required isomorphism. Since
$R^{-i,2\mb a}(\sS)=0$ if $\mb a$ is not a $(0,1)$-vector,
$\Tor_{\k[m]}^{-i,2\mb a}(\k[\mathcal K],\k)$ vanishes for
such~$\mb a$.
\end{proof}

We define the \emph{multigraded algebraic Betti
numbers}\label{aBettisp} of $\k[\sS]$ as
\[
  \beta^{-i,\,2\mb a}\bigl(\k[\sS]\bigr)=
  \dim_\k\Tor^{-i,\,2\mb
  a}_{\k[v_1,\ldots,v_m]}\bigl(\k[\sS],\k\bigr),
\]
for $0\le i\le m$, $\mb a\in\N^m$. We also set
\[
  \beta^{-i}\bigl(\k[\sS]\bigr)=\dim_\k\Tor^{-i}_{\k[v_1,\ldots,v_m]}
  \bigl(\k[\sS],\k\bigr)=\sum_{\mb a\in\N^m}\beta^{-i,\,2\mb a}\bigl(\k[\sS]\bigr).
\]

\begin{example}
Let $\sS$ be the simplicial poset of Example~\ref{exsp}.1. By
Theorem~\ref{hochsp},
$\beta^{0,(0,0)}(\k[\sS])=\beta^{0,(2,2)}(\k[\sS])=1$, and the
other Betti numbers are zero. This implies that $\k[\sS]$ is a
free $\k[v_1,v_2]$-module with two generators, $1$ and~$v_\sigma$,
of degree 0 and 4 respectively.
\end{example}

Note that unlike the case of simplicial complexes,
$\beta^0(\k[\sS])$ may be larger than~1. In fact, the following
proposition follows easily from Theorem~\ref{hochsp}.

\begin{proposition}
The number of generators of $\k[\sS]$ as a $\k[m]$-module equals
\[
  \beta^0(\k[\sS])=\sum_{J\subset[m]}
  \dim\widetilde H^{|J|-1}\bigl(|\sS_J|\bigr).
\]
\end{proposition}

\subsection*{Exercises}

\begin{exercise}
Finish the proof of Proposition~\ref{frext}.
\end{exercise}

\begin{exercise}
Fill in the details in the proof of Theorem~\ref{hochsp}.
\end{exercise}

\begin{exercise}
Calculate the multigraded Betti numbers for the simplicial poset
of Example~\ref{exsp}.2.
\end{exercise}

\chapter*{\ \ Face rings: additional topics}

\section{Cohen--Macaulay simplicial posets}\label{scccm}

Assume given a property $A$ of simplicial complexes. Then we can
extend this property to posets by postulating that a poset
$\mathcal P$ has the property $A$ if the order complex
$\ord(\mathcal P)$ (see Definition~\ref{ordercom}) has the
property~$A$. In particular, Cohen--Macaulay and Gorenstein posets
can be defined in this way. Simplicial posets $\sS$ are of
particular interest to us; in this case the order complex is
identified with the barycentric subdivision $\sS'$ (to be precise,
with the cone over the barycentric subdivision, as we include the
empty simplex, but this difference is inessential for the
definitions to follow).

\begin{definition}\label{CMsimpose}
A simplicial poset $\mathcal S$ is \emph{Cohen--Macaulay}
(over~$\k$) if its barycentric subdivision $\mathcal S'$ is a
Cohen--Macaulay simplicial complex.
\end{definition}

By definition, $\mathcal S$ is a Cohen--Macaulay simplicial poset
if and only if the face ring $\k[\mathcal S']$ is Cohen--Macaulay.
Since the face ring is also defined for the face poset $\mathcal
S$ itself (and not only for its barycentric subdivision), it is
perfectly natural to ask whether the class of Cohen--Macaulay
simplicial posets admits an intrinsic description in terms of
their face rings~$\k[\mathcal S]$. One would achieve such a
description by proving that the ring $\k[\sS']$ is Cohen--Macaulay
if and only if the ring $\k[\sS]$ is Cohen--Macaulay. The `if'
part follows from the general theory of ASL's,
see~\cite[Corollary~3.7]{stan91}. The `only if' part was proved
in~\cite{m-m-p07}; the proof uses the decomposition of the
barycentric subdivision into a sequence of stellar subdivisions
and then goes on to show that the Cohen--Macaulay property is
preserved under stellar subdivisions. We include this
characterisation of Cohen--Macaulay simplicial posets in terms of
their face rings in Theorem~\ref{theo:cmpos} below.

Since many of the constructions in this section are geometric, we
often talk about simplicial cell complexes rather than simplicial
posets. We say that a simplicial subdivision of a simplicial cell
complex $\mathcal S$ is \emph{regular}\label{regsubsp} if it is a
simplicial complex. For instance, the barycentric subdivision is
regular. Since the Cohen--Macaulayness of a simplicial complex is
a topological property (see Proposition~\ref{munkr}), we have the
following statement.

\begin{proposition}\label{coro:cmcha}
The following conditions are equivalent:
\begin{itemize}
\item[(a)] the barycentric subdivision of a simplicial cell
complex $\mathcal S$ is a Cohen--Macaulay complex;

\item[(b)] any regular subdivision of $\mathcal S$ is a
Cohen--Macaulay complex;

\item[(c)] a regular subdivision of $\mathcal S$ is a
Cohen--Macaulay complex.
\end{itemize}
\end{proposition}

As a further corollary we obtain that Proposition~\ref{munkr}
itself extends to simplicial cell complexes, i.e. the property of
a simplicial cell complex to be Cohen--Macaulay is also
topological.

By analogy with Definition~\ref{deflink}, we define the
\emph{star}\label{starlinksp} and the \emph{link} of
$\sigma\in\sS$ as the following subcomplexes:
\begin{align*}
  \st_{\mathcal S}\sigma&=
  \{\tau\in\mathcal S\colon\sigma\vee\tau\text{ is nonempty}\};\\
  %\partial\st_{\mathcal S}\sigma&=
  %\{\tau\in\mathcal S\colon
  %\sigma\vee\tau\text{ is nonempty, and }\sigma\not\le\tau\};\\
  \lk_{\mathcal S}\sigma&=\{\tau\in\mathcal S\colon
  \sigma\vee\tau\text{ is nonempty, and }
  \tau\wedge\sigma=\hat 0\}.
\end{align*}

\begin{remark}
%In the simplicial cell complex $\sS$, the star $\st_{\mathcal
%S}\sigma$ is the `closed combinatorial neighbourhood'
%of~$\varDelta^\sigma$.
If $\mathcal S$ is a simplicial complex, then the poset
$\lk_{\sS}\sigma$ is isomorphic to the open semi-interval
\[
  \mathcal S_{>\sigma}=\{\rho\in\mathcal S\colon\rho>\sigma\},
\]
and $|\st_{\sS}\sigma|\cong\varDelta^\sigma\!*|\lk_{\sS}\sigma|$,
where $*$ denotes the join. However, none of these isomorphisms
holds for general~$\sS$, see Example~\ref{2tria} below.
\end{remark}

Because of this remark, we cannot simply extend the definition of
stellar subdivisions (Definition~\ref{bist}) to simplicial cell
complexes. Instead, we define the \emph{stellar subdivision}
$\ss_\sigma\sS$ of $\sS$ at $\sigma$ as the simplicial cell
complex obtained by stellarly subdividing each face
containing~$\sigma$ in a compatible way.

\begin{proposition}\label{prop:barst}
The barycentric subdivision $\sS'$ can be obtained as a sequence
of stellar subdivisions, one at each face $\sigma\in\sS$, starting
from the maximal faces. Moreover, each stellar subdivision in the
sequence is applied to a face whose star is a simplicial complex.
\end{proposition}
\begin{proof}
Assume $\dim\mathcal S=n-1$. We start by applying to $\sS$ stellar
subdivisions at all $(n-1)$-dimensional faces. Denote the
resulting complex by $\sS_1$. The $(n-2)$-faces of $\mathcal S_1$
are of two types: the ``old'' ones, remaining from~$\mathcal S$,
and the ``new'' ones, appearing as the result of the stellar
subdivisions. Then we take stellar subdivisions of $\sS_1$ at all
``old'' $(n-2)$-faces, and denote the result by $\sS_2$. Next we
apply to $\sS_2$ stellar subdivisions at all $(n-3)$-faces
remaining from~$\sS$. Proceeding in this way, at the end we get
$\sS_{n-1}=\sS'$. To prove the second statement, consider two
subsequent complexes $\mathcal R$ and $\widetilde{\mathcal R}$ in
the sequence, so that $\widetilde{\mathcal R}$ is obtained from
$\mathcal R$ by a single stellar subdivision at
some~$\sigma\in\sS$. Then $\st_{\mathcal R}\sigma$ is isomorphic
to $\varDelta^\sigma\mathbin{*}(\sS_{>\sigma})'$ and therefore it
is a simplicial complex.
\end{proof}

We proceed with two lemmata necessary to prove our main result.

\begin{lemma}\label{lemm:frppt}
Let $\sS$ be a simplicial poset of rank $n$ with vertex set
$V(\sS)=[m]$, and assume that the first $k$ vertices span a
face~$\sigma$. Assume further that $\st_{\sS}\sigma$ is a
simplicial complex, and let $\widetilde\sS$ be the stellar
subdivision of $\sS$ at~$\sigma$. Let $v$ denote the degree-two
generator of $\k[\widetilde\sS]$ corresponding to the added
vertex. Then there exists a unique homomorphism
$\beta\colon\k[\sS]\to\k[\widetilde\sS]$ such that
\begin{align*}
  v_\tau & \mapsto v_\tau &&\text{for \ }\tau\notin\st_{\sS}\sigma;\\
  v_i & \mapsto v+v_i, &&\text{for \ } i=1,\ldots,k;\\
  v_i & \mapsto v_i, &&\text{for \ } i=k+1\ldots,m.
\end{align*}
Moreover, $\beta$ is injective, and if $\mb t$ is an hsop in
$\k[\sS]$, then $\beta(\mb t)$ is an hsop in~$\k[\widetilde\sS]$.
\end{lemma}
\begin{proof}
In order to define the map $\beta$ we first need to specify the
images of $v_\tau$ for all $\tau\in\st_{\sS}\sigma$. Choose such a
$v_\tau$ and let $V(\tau)=\{i_1,\ldots,i_\ell\}$ be its vertex
set. Then we have the following identity in the ring
$\k[\sS]=\k[v_\tau\colon\tau\in\sS]/\mathcal I_{\sS}$:
\begin{equation}\label{expan}
  v_{i_1}\cdots v_{i_\ell}=
  v_\tau+\sum_{\eta\colon V(\eta)=V(\tau),\,\eta\ne\tau}v_\eta.
\end{equation}
For any $v_\eta$ in the latter sum we have
$\eta\notin\st_{\sS}\sigma$, since $\st_{\sS}\sigma$ is a
simplicial complex, in which any set of vertices spans at most one
face. Since $\beta$ is already defined on the product on the left
hand side and on the sum on the right hand side above, this
determines $\beta(v_\tau)$ uniquely.

We therefore obtain a map of polynomial algebras
$\k[v_\tau\colon\tau\in\sS]\to\k[v_\tau\colon\tau\in\widetilde{\sS}]$
(which we denote by the same letter $\beta$ for a moment), and
need to check that it descends to a map of face rings,
$\k[\sS]\to\k[\widetilde{\sS}]$. In other words, we need to verify
that $\beta(\mathcal I_{\sS})\subset\mathcal I_{\widetilde{\sS}}$.

It is clear from the definition of $\beta$ that we have the
commutative diagram
\[
  \xymatrix{
  \k[v_\tau\colon\tau\in\sS] \ar[r]^-p \ar[d]^\beta &
  \k[\sS] \ar[r]^-s \ar@{-->}[d]^\beta &
  \bigoplus_{\tau\in\sS}\k[v_i\colon i\in V(\tau)] \ar[d]^{s(\beta)}\\
  \k[v_\tau\colon\tau\in\widetilde{\sS}] \ar[r]^-{\widetilde p} &
  \k[\widetilde{\sS}] \ar[r]^-{\widetilde s} &
  \bigoplus_{\tau\in\widetilde{\sS}}\k[v_i\colon i\in V(\tau)],
  }
\]
in which the middle vertical map is not yet defined. Here by $s$
and $\widetilde s$ we denote the restriction maps from
Theorem~\ref{alres}, and $s(\beta)$ is the map induced by $\beta$
on the direct sum of polynomial algebras.
%Note that $s(\beta)$ is defined on the generators as described in
%the Lemma (i.e. it takes $v_i$ to $v+v_i$ or $v_i$ depending
%on~$i$).
Now let $x\in \mathcal I_{\sS}$, i.e. $p(x)=0$. Then, by the
commutativity of the diagram, $\widetilde s\widetilde
p\beta(x)=0$. Since $\widetilde s$ is injective, we have
$\widetilde p\beta(x)=0$. Hence, $\beta(x)\in\mathcal
I_{\widetilde{\sS}}$, which implies that the middle vertical map
is well defined.

The last statement also follows from the commutative diagram
above. The map $s(\beta)$ sends each direct summand of its domain
isomorphically to at least one summand of its range, and therefore
it is injective. Thus, $\beta\colon\k[\sS]\to\k[\sS']$ is also
injective. The statement about hsop's then follows
Lemma~\ref{lsopscc}.
\end{proof}

\begin{remark} If we defined the map $\beta$ by sending each $v_i$ identically,
then it would still give rise to a ring homomorphism $\k[\mathcal
S]\to\k[\widetilde{\mathcal S}]$, but the latter would not be
injective (for example, it would map $v_\sigma\in\k[\sS]$ to
zero).
\end{remark}

\begin{example}\label{2tria}
The assumption on $\st_{\sS}\sigma$ in Lemma~\ref{lemm:frppt} is
not always satisfied. For example, if $\sS$ is obtained by
identifying two 2-simplices along their boundaries, and $\sigma$
is any edge, then $\st_{\sS}\sigma=\sS$, which is not a simplicial
complex.

Note also that if $\st_{\sS}\sigma$ is not a simplicial complex,
then the map $\beta\colon\k[\sS]\to\k[\sS']$ is not determined
uniquely by the conditions specified in Lemma~\ref{lemm:frppt} (we
cannot determine the images of $v_\tau$ with
$\tau\in\st_{\sS}\sigma$). Nevertheless, it is still possible to
define the map $\beta\colon\k[\sS]\to\k[\sS']$ for an arbitrary
simplicial poset~$\sS$, see Section~\ref{weightgraphs}.
\end{example}

\begin{lemma}\label{lemm:cmtil}
Assume that $\k[\sS]$ is a Cohen--Macaulay ring, and let
$\widetilde\sS$ be a stellar subdivision of $\sS$ at $\sigma$ such
that $\st_{\sS}\sigma$ is a simplicial complex. Then
%\begin{itemize}
%\item[(a)] $\st_{\sS}\sigma$ is a Cohen--Macaulay complex;
%
%\item[(b)]
$\k[\widetilde{\mathcal S}]$ is a Cohen--Macaulay ring.
%\end{itemize}
\end{lemma}
\begin{proof}
We first prove that $\st_{\sS}\sigma$ is a Cohen--Macaulay
complex. Since
$\st_{\sS}\sigma=\varDelta^\sigma\mathbin{*}\,\lk_{\sS}\sigma$, it
is enough to verify that $\lk_{\sS}\sigma$ is Cohen--Macaulay.
This follows from Reisner's Theorem (Theorem~\ref{reisner}) and
the fact that simplicial cohomology of $\lk_{\sS}\sigma$ is a
direct summand in local cohomology of~$\k[\sS]$
(see~\cite[Theorem~II.4.1]{stan96}
or~\cite[Theorem~5.3.8]{br-he98}).

Choose an hsop $\mb t=(t_1,\ldots,t_n)$ in $\k[\sS]$ and set
$\widetilde{\mb t}=\beta(\mb t)$. Let
\[
  p\colon \k[\sS]\to\k[\sS]/(v_\tau\colon\tau\notin\st_{\sS}\sigma)
  =\k[\st_{\sS}\sigma]
\]
be the quotient projection. Set $R=\ker p$. Similarly, set
\[
  \widetilde R=\ker\bigl(\widetilde
  p\colon\k[\widetilde\sS]\to\k[\st_{\widetilde\sS}v]\bigr),
\]
where $v$ is the new vertex added in the process of stellar
subdivision. Since the simplicial cell complexes $\sS$ and
$\widetilde\sS$ do not differ on the complements of
$\st_{\sS}\sigma$ and $\st_{\widetilde\sS}v$ respectively, the map
$\beta$ restricts to the identity isomorphism $R\to\widetilde R$.
We therefore have the following commutative diagram with exact
rows:
\[
\begin{CD}
  0 @>>> R @>>> \k[\sS] @>p>>
  \k[\st_{\sS}\sigma] @>>> 0\\
  @. @VV\cong V @VV\beta V @VVV\\
  0 @>>> \widetilde R @>>> \k[\widetilde\sS] @>\widetilde p>>
  \k[\st_{\widetilde\sS}v] @>>>
  0,
\end{CD}
\]
Applying the functors $\ \otimes_{\k[\mb t]}\k$ and $\
\otimes_{\k[\widetilde{\mb t}]}\k$ to the diagram above, we get a
map between the long exact sequences for $\Tor$. Consider the
following fragment:
\[
\begin{array}{cccc}
  \Tor^{-2}_{\k[\mb t]}(\k[\st\sigma],\k)
  \stackrel{f}{\to} &\!\Tor^{-1}_{\k[\mb t]}(R,\k)
  \to &\! \Tor^{-1}_{\k[\mb t]}\k[\sS],\k)  \to &\!
  \Tor^{-1}_{\k[\mb t]}(\k[\st\sigma],\k)\\
  \downarrow & \downarrow{\scriptstyle\cong} & \downarrow & \downarrow\\
  \Tor^{-2}_{\k[\widetilde{\mb t}]}(\k[\st v],\k)
  \stackrel{\widetilde f}{\to} &\!
  \Tor^{-1}_{\k[\widetilde{\mb t}]}(\widetilde R,\k)  \to &\!
  \Tor^{-1}_{\k[\widetilde{\mb t}]}(\k[\widetilde\sS],\k)  \to &\!
  \Tor^{-1}_{\k[\widetilde{\mb t}]}(\k[\st v],\k).
\end{array}
\]
Since $\k[\sS]$ is Cohen--Macaulay, $\Tor^{-1}_{\k[\mb
t]}(\k[\sS],\k)=0$ and the map $f$ is surjective. Then $\widetilde
f$ is also surjective. Since $\st_{\sS}\sigma$ is a
Cohen--Macaulay simplicial complex and
$|\st_{\sS}\sigma|\cong|\st_{\widetilde\sS}v|$,
Proposition~\ref{munkr} implies that $\k[\st v]$ is
Cohen--Macaulay. Therefore, $\Tor^{-1}_{\k[\widetilde{\mb
t}]}(\k[\st v],\k)=0$. Since $\widetilde f$ is surjective, we also
have $\Tor^{-1}_{\k[\widetilde{\mb t}]}(\k[\widetilde\sS],\k)=0$.
Hence $\k[\widetilde\sS]$ is free as a ${\k[\widetilde{\mb
t}]}$-module (see~\cite[Lemma~VII.6.2]{macl63}) and thereby is
Cohen--Macaulay.
\end{proof}

Now we can prove the main result of this section:

\begin{theorem}\label{theo:cmpos}
A simplicial poset $\sS$ is Cohen--Macaulay if and only if the
face ring $\k[\sS]$ is Cohen--Macaulay.
\end{theorem}
\begin{proof}
The fact that the face ring of a Cohen--Macaulay simplicial
poset~$\sS$ is Cohen--Macaulay is proved
in~\cite[Corollary~3.7]{stan91} (see also~\cite[\S~III.6]{stan96})
using the theory of ASL's.
%For an alternative proof, see Exercise~\ref{exsstil}.

Assume now that $\k[\mathcal S]$ is a Cohen--Macaulay ring. Since
the barycentric subdivision $\sS'$ is obtained by a sequence of
stellar subdivisions, subsequent application of
Lemma~\ref{lemm:cmtil} shows that $\k[\sS']$ is also
Cohen--Macaulay. Thus, $\sS'$ is a Cohen--Macaulay poset.
\end{proof}

We end this section by giving Stanley's characterisation of
$h$-vectors of Cohen--Macaulay simplicial posets.

\begin{theorem}[Stanley]\label{hvcmsp}
The integer vector $\mb h=(h_0,h_1,\ldots,h_n)$ is the $h$-vector
of a Cohen--Macaulay simplicial poset if and only if $h_0=1$ and
$h_i\ge0$.
\end{theorem}
\begin{proof}
Let $\mb h=\mb h(\mathcal S)$ for a Cohen--Macaulay simplicial
poset~$\mathcal S$. The condition $h_0=1$ follows from the
definition of the $h$-vector, see~\eqref{hvectors}. Let $\k$ be a
field of zero characteristic, and $\mb t=(t_1,\ldots,t_n)$ an lsop
in $\k[\mathcal S]$ (since $\k[\mathcal S]$ is not generated by
linear elements, the existence of an lsop is not automatic and is
left as an exercise; alternatively, see~\cite[Lemma~3.9]{stan91}).
Comparing the formula for the Poincar\'e series from
Proposition~\ref{pscma} with that of Theorem~\ref{psgfr}, we
obtain
$$
  F\bigl(\k[\mathcal S]/\mb t;\lambda\bigr)=h_0+h_1\lambda^2+\cdots+h_n\lambda^{2n}.
$$
Hence, $h_i\ge0$, as needed.

Now we construct a Cohen--Macaulay simplicial cell complex
$\mathcal S$ with any given $h$-vector such that $h_0=1$ and
$h_i\ge0$. First note that $\mb h(\varDelta^{n-1})=(1,0,\ldots,0)$
and $\varDelta^{n-1}$ is a Cohen--Macaulay simplicial (cell)
complex. Now, given an $(n-1)$-dimensional Cohen--Macaulay
simplicial cell complex $\mathcal S$ with the $h$-vector
$(h_0,\ldots,h_n)$, it suffices to construct, for any
$k=1,\ldots,n$, a new Cohen--Macaulay simplicial cell complex
$\mathcal S_k$ with the $h$-vector given by
\begin{equation}\label{hs1}
  \mb h(\mathcal
  S_k)=(h_0,\ldots,h_{k-1},h_k+1,h_{k+1},\ldots,h_n).
\end{equation}
To do this, we choose an $(n-1)$-face of $\mathcal S$, and in this
face choose some $k$ faces of dimension~$n-2$. Then add to
$\mathcal S$ a new $(n-1)$-simplex by attaching it along some $k$
faces of dimension $n-2$ to the chosen $k$ faces of~$\mathcal S$.
A direct check shows that the $h$-vector of the resulting
simplicial cell complex $\sS_k$ is given by~\eqref{hs1}. The fact
that $\mathcal S_k$ is Cohen--Macaulay follows directly from
Proposition~\ref{munkr}.
\end{proof}

Note that this characterisation is substantially simpler than that
for simplicial complexes (see Propositions~\ref{stanmv} and the
remark after it).

\subsection*{Exercises}
\begin{exercise}
The  map of face rings $\k[\sS]\to\k[\widetilde{\sS}]$ of
Lemma~\ref{lemm:frppt} is not induced by any poset map
$\widetilde{\sS}\to\sS$.
\end{exercise}

\begin{exercise}
Let $\widetilde\sS$ be a stellar subdivision of $\sS$ at $\sigma$
such that $\st_{\sS}\sigma$ is a simplicial complex. Show that the
ring $\k[\sS]$ is Cohen--Macaulay if and only if
$\k[\widetilde{\mathcal S}]$ is Cohen--Macaulay, i.e. the converse
of Lemma~\ref{lemm:cmtil} holds.
\end{exercise}

\begin{exercise}
If $\k$ is of characteristic zero, then $\k[\sS]$ admits an lsop.
\end{exercise}

\section{Gorenstein simplicial posets}\label{gscc}
Gorenstein simplicial posets arise in toric topology as the
combinatorial structures associated to the orbit quotients of
\emph{torus manifolds}, which are the subject of
Chapter~\ref{torus}. It was exactly this particular feature of
Gorenstein simplicial posets which allowed Masuda~\cite{masu05} to
complete the characterisation of their $h$-vectors, conjectured by
Stanley in~\cite{stan91}. We include Masuda's result here as
Theorem~\ref{heven}.

\begin{definition}\label{goresp}
A simplicial poset~$\mathcal S$ is \emph{Gorenstein}
(respectively, \emph{Gorenstein*}) if its barycentric subdivision
$\mathcal S'$ is a Gorenstein (respectively, Gorenstein*)
simplicial complex.
\end{definition}

Like Cohen--Macaulayness, the property of a simplicial poset $\sS$
being Gorenstein* depends only on the topology of the realisation
$|\mathcal S|$ (this follows from Theorem~\ref{gorencom}). In
particular, simplicial cell subdivisions of spheres are
Gorenstein*.

The problem of characterisation of $h$-vectors of Gorenstein*
simplicial posets is more subtle than the corresponding question
in the Cohen--Macaulay case. (Although this problem is much easier
for simplicial posets than for simplicial complexes, see the
discussion of the $g$-conjecture in Sections~\ref{simsph}
and~\ref{gc}.)

\begin{theorem}\label{dsgsp}
Let $\mb h(\mathcal S)=(h_0,h_1,\ldots,h_n)$ be the $h$-vector of
a Gorenstein* simplicial poset of rank~$n$. Then $h_0=1$,
$h_i\ge0$ and $h_i=h_{n-i}$ for any~$i$.
\end{theorem}
\begin{proof}
The inequalities $h_i\ge0$ follow from the fact that $\mathcal S$
is Cohen--Macaulay (Theorem~\ref{hvcmsp}). The identities
$h_i=h_{n-i}$ will follow from the expression of the $h$-vector of
the barycentric subdivision $\mathcal S'$ via $\mb h(\mathcal S)$
and from the Dehn--Sommerville relations for the Gorenstein*
simplicial complex~$\mathcal S'$. Indeed, repeating the argument
from Lemmata~\ref{barfv} and~\ref{barhv} we  obtain the identity
$\mb h(\mathcal S')=D\mb h(\mathcal S)$, in which the vector $\mb
h(\mathcal S')$ is symmetric, i.e. satisfies the Dehn--Sommerville
relations. It can be checked directly using some identities for
binomial coefficients that the operator $D$ (and its inverse)
takes symmetric vectors to symmetric ones (which is equivalent to
the identity $d_{pq}=d_{n+1-p,n+1-q}$). This calculation can be
avoided by using the following argument. The Dehn--Sommerville
relations specify a linear subspace $W$ of dimension $k=\sbr n2+1$
in the space $\R^{n+1}$ with coordinates $h_0,\ldots,h_n$. We need
to check that this subspace is $D$-invariant. To do this it
suffices to choose a basis $\mb e_1,\ldots,\mb e_k$ in $W$ and
check that $D\mb e_i\in W$ for all~$i$. There is a basis in $W$
consisting of $h$-vectors of simplicial spheres (and even
simplicial polytopes, see the proof of Proposition~\ref{DSgenli}).
Since the barycentric subdivision of a simplicial sphere is a
simplicial sphere, the vectors $D\mb e_i$, \ $1\le i\le k$, are
also symmetric, and $W$ is a $D$-invariant subspace. Thus, the
vector $\mb h(\mathcal S)=D^{-1}\mb h(\mathcal S')$ satisfies the
Dehn--Sommerville relations.
\end{proof}

%В~\cite[Theorem~4.3]{stan91} были получены следующие достаточные
%условия.

\begin{theorem}[{\cite[Theorem~4.3]{stan91}}]\label{hvsc}
Let $\mb h=(h_0,h_1,\ldots,h_n)$ be an integer vector with
$h_0=1$, $h_i\ge0$ and $h_i=h_{n-i}$. Any of the following
(mutually exclusive) conditions are sufficient for the existence
of a Gorenstein* simplicial poset of rank $n$ and $h$-vector
$h(\sS)=\mb h$:
\begin{itemize}
\item[(a)] $n$ is odd;

\item[(b)] $n$ is even and $h_{n/2}$ is even;
\item[(c)] $n$ is even, $h_{n/2}$ is odd, and $h_i>0$ for all~$i$.
\end{itemize}
\end{theorem}
\begin{proof}
We start with the following two basic examples of
$(n-1)$-dimensional simplicial cell complexes of dimension:
$\partial\varDelta^n$, with $h$-vector $\mb
h(\partial\varDelta^n)=(1,1,\ldots,1)$; and $\mathcal S_n$, the
simplicial cell complex obtained by identifying two
$(n-1)$-simplices along their boundaries, with $\mb h(\mathcal
S_n)=(1,0,\ldots,0,1)$. By applying the standard operations of
join and connected sum (Constructions~\ref{join} and~\ref{simcs})
to these two complexes we can obtain a simplicial cell complex
with any prescribed $h$-vector satisfying the conditions of the
theorem. Indeed, for $k\ne n-k$ we have
$$
  \mb h(\mathcal S_k*\mathcal S_{n-k})=
  (1,0,\ldots,0,1,0,\ldots,0,1,0,\ldots,0,1),
$$
where $h_k=h_{n-k}=1$, and the other entries are zero. Also, for
$n=2k$ we have
$$
  \mb h(\mathcal S_k*\mathcal S_k)=
  (1,0,\ldots,0,2,0,\ldots,0,1),
$$
where $h_k=2$. Now, by taking connected sum of the appropriate
number of complexes $\partial \varDelta^n$, $\mathcal S_n$ and
$\mathcal S_k*\mathcal S_{n-k}$ and using the identity
$$
  h_i(\mathcal S\cs\widetilde{\mathcal S})=h_i(\mathcal S)+h_i(\widetilde{\mathcal S})\quad
  \text{for }1\le i\le n-1,
$$
(see Example~\ref{fhvcs}, which is valid for any two pure
$(n-1)$-dimensional simplicial cell complexes), we obtain any
required $h$-vector.
\end{proof}

The subtlest part of the characterisation of $h$-vectors of
Gorenstein* simplicial posets is the following result, which was
proved by Masuda:

\begin{theorem}[\cite{masu05}]\label{heven}
Let $\mb h(\mathcal S)=(h_0,h_1,\ldots,h_n)$ be the $h$-vector of
a Gorenstein* simplicial poset $\mathcal S$ of even rank~$n$, and
let $h_i=0$ for some~$i$. Then the number $h_{n/2}$ is even.
\end{theorem}

Note that the evenness of $h_{n/2}$ is equivalent to the evenness
of the number of facets $f_{n-1}=\sum_{i=0}^nh_i$. The idea behind
Masuda's proof of Theorem~\ref{heven} lies within the topological
theory of torus manifolds, which is the subject of
Section~\ref{torusman}.

We combine the results of Theorems~\ref{dsgsp}, \ref{hvsc}
and~\ref{heven} in the following characterisation result for the
$h$-vectors of Gorenstein* simplicial posets.

\begin{theorem}
An integer vector $\mb h=(h_0,h_1,\ldots,h_n)$ is the $h$-vector
of a Gorenstein* simplicial poset of rank~$n$ (or a simplicial
cell subdivision of an $(n-1)$-sphere) if and only if the
following conditions are satisfied:
\begin{itemize}
\item[(a)] $h_0=1$ and $h_i\ge0$;
\item[(b)] $h_i=h_{n-i}$ for all~$i$;
\item[(c)] either $h_i>0$ for all $i$ or $\sum_{i=0}^nh_i$ is
even.
\end{itemize}
\end{theorem}
%
%The same conditions also characterise the $h$-vectors of
%simplicial cell subdivisions of spheres.

\section{Generalised Dehn--Sommerville relations}\label{gendsr}
In this section we obtain some further generalisations of the
Dehn--Sommerville relations, in particular, to arbitrary
triangulated manifolds.

Let $\mathcal S$ be a simplicial poset of rank~$n$. Given
$\sigma\in\sS$, consider the closed upper semi-interval $\mathcal
S_{\ge\sigma}=\{\tau\in\mathcal S\colon\tau\ge\sigma\}$ with the
induced rank function, and set
\begin{equation} \label{eqn:tchi}
\chi(\mathcal S_{\ge\sigma})=\sum_{\tau\ge \sigma}(-1)^{|\tau|-1}.
\end{equation}

A simplicial poset $\mathcal S$ of rank~$n$ satisfying
$\chi(\mathcal S_{\ge\sigma})=(-1)^{n-1}$ for all
$\sigma\in\mathcal S$ is called \emph{Eulerian}. According to a
result of~\cite[(3.40)]{stan86}, the Dehn--Sommerville relations
$h_i=h_{n-i}$ hold for Eulerian posets. This can be generalised as
follows.

\begin{theorem}[{see~\cite[Theorem~9.1]{m-m-p07}}]\label{theo:DS}
The following identity holds for the $h$-vector $\mb
h(\sS)=(h_0,\ldots,h_n)$ of a simplicial poset $S$ of rank~$n$:
\[
\sum_{i=0}^n(h_{n-i}-h_{i})t^i=\sum_{\sigma\in\mathcal S}
\Bigl(1+(-1)^n\chi(\mathcal
S_{\ge\sigma})\Bigr)(t-1)^{n-|\sigma|}.
\]
In particular, if $\sS$ is Eulerian, then $h_i=h_{n-i}$.
\end{theorem}
\begin{proof}We have
\begin{equation} \label{eqn:DS}
\begin{split}
\sum_{i=0}^n h_it^i&=t^n \sum_{i=0}^n
h_i({\textstyle\frac1t})^{n-i}
=t^n\sum_{i=0}^n f_{i-1}\bigl({\textstyle\frac{1-t}t}\bigr)^{n-i}\\
&=\sum_{i=0}^n f_{i-1}t^i(1-t)^{n-i}
=\sum_{\tau\in\mathcal S} t^{|\tau|}(1-t)^{n-|\tau|}\\
&=\sum_{\tau\in\mathcal S}\sum_{\sigma\le
\tau}(t-1)^{|\tau|-|\sigma|} (1-t)^{n-|\tau|}
=\sum_{\tau\in\mathcal S}\sum_{\sigma\le \tau}(-1)^{n-|\tau|}(t-1)^{n-|\sigma|}\\
&=\sum_{\sigma\in\mathcal S}(t-1)^{n-|\sigma|}\sum_{\tau\ge
\sigma} (-1)^{n-|\tau|}=\sum_{\sigma\in\mathcal
S}(t-1)^{n-|\sigma|}(-1)^{n-1}\chi(\mathcal S_{\ge\sigma}),
\end{split}
\end{equation}
where the fifth identity follows from the binomial expansion of
the right hand side of the identity
$t^{|\tau|}=((t-1)+1)^{|\tau|}$ and the fact that
$[\hatzero,\tau]=\{\sigma\in\sS\colon\sigma\le\tau\}$ is a Boolean
lattice of rank~$|\tau|$.

On the other hand, we have
\begin{equation} \label{eqn:DS2}
\sum_{i=0}^nh_{n-i}t^i=\sum_{i=0}^n
h_it^{n-i}=\sum_{i=0}^nf_{i-1}(t-1)^{n-i} =\sum_{\sigma\in\mathcal
S}(t-1)^{n-|\sigma|}.
\end{equation}

Subtracting (\ref{eqn:DS}) from (\ref{eqn:DS2}) we obtain the
required identity.
\end{proof}

As a corollary we obtain a generalisation of the Dehn--Sommerville
relations to triangulated manifolds. This formula appeared
in~\cite[p.~74]{stan96} (the orientability assumption there can be
removed by passing to the orientation double cover, see
also~\cite[Corollary~4.5.4]{bu-pa00-2}):

\begin{theorem}\label{dstm}
Let $\mathcal K$ be a triangulation of a closed
$(n-1)$-dimensional mani\-fold. Then the $h$-vector $\mb
h(\sK)=(h_0,\ldots,h_n)$ satisfies the identities
\[
  h_{n-i}-h_i=(-1)^i\binom ni\bigl(\chi(\mathcal
  K)-\chi(S^{n-1})\bigr),\quad 0\le i\le n.
\]
Here $\chi(\sK)=f_0-f_1+\cdots+(-1)^{n-1}f_{n-1}=1+(-1)^{n-1}h_n$
is the Euler characteristic of~$\sK$ and
$\chi(S^{n-1})=1+(-1)^{n-1}$.
\end{theorem}
\begin{proof}
Viewing $\sK$ as a simplicial poset, we calculate
\begin{align*}
  \chi(\mathcal K_{\ge\sigma})&=\sum_{\tau>\sigma}(-1)^{|\tau|-1}+(-1)^{|\sigma|-1}
  =(-1)^{|\sigma|}\Bigl(
  \sum_{\tau>\sigma}(-1)^{|\tau|-|\sigma|-1}-1\Bigr)\\
  &=(-1)^{|\sigma|}\Bigl(\sum_{\varnothing\ne\rho\in\lk_\mathcal K\sigma}(-1)^{|\rho|-1}-1\Bigr)
  =(-1)^{|\sigma|}\bigl(\chi(\lk_\mathcal K\sigma)-1\bigr).
\end{align*}
Here we have used the fact that the poset of nonempty faces of
$\lk_\mathcal K\sigma$ is isomorphic to $\mathcal S_{>\sigma}$,
with the rank function shifted by~$|\sigma|$. Now since $\mathcal
K$ is a triangulated $(n-1)$-dimensional manifold, the link of a
nonempty face $\sigma\in\sK$ has the homology of a sphere of
dimension $(n-|\sigma|-1)$. Hence, $\chi(\lk_\mathcal
K\sigma)=1+(-1)^{n-|\sigma|-1}$, and therefore $\chi(\mathcal
K_{\ge\sigma})=(-1)^{n-1}$ for $\sigma\ne\varnothing$. Also,
$\lk_\mathcal K\varnothing=\mathcal K$. Now using the identity of
Theorem~\ref{theo:DS} we calculate
\[
%\begin{split}
\sum_{i=0}^n(h_{n-i}-h_{i})t^i =\bigl(1+(-1)^{n}(\chi(\mathcal
K)-1)\bigr)(t-1)^n
%=\\[-4pt]
=(-1)^n\bigl(\chi(\mathcal K)-\chi(S^{n-1})\bigr)(t-1)^n.
%\end{split}
\]
The required identity follows by comparing the coefficients
of~$t^i$.
\end{proof}

For other generalisations of Dehn--Sommerville relations
see~\cite{ba-bi85} and~\cite{g-s-j}.

\subsection*{Exercises}
\begin{exercise}
The identity of Theorem~\ref{dstm} holds for simplicial posets.
\end{exercise}

%обобщение полиэдр. произведения на случай произведений над Y
%(X_i,A_i)\to Y
%
%R_K для многоугольников: Coxeter'1938
%
%R_K - ациклическое многообразие для флаговых K
%
%DJ(K) - формально. Дать док-во используя, что кольцо граней -
%предел колец многочленов, а A_{PL} переводит копределы в пределы.
%
%Z_K как конфигурационное пространство шарнирных механизмов
%(статьи Изместьева, японцев)
%
%Дать ссылку на общую теорему о связных суммах в упражнении \ref{exzk5gon}
%
%Кольцо когомологий R_K: описание, ссылки
%
%В раздел 4.8 о произведениях Масси вставить ссылки на Берглунда и нас о св-ве Голода

\setcounter{chapter}3
\chapter{Moment-angle complexes}\label{macom}

This is the first genuinely `toric' chapter of this book; it links
the combinatorial and algebraic constructions of the previous
chapters to the world of toric spaces.

The term `moment-angle complex' refers to a decomposition of a
certain toric space $\zk$ into products of polydiscs and tori
parametrised by simplices in a given simplicial complex~$\sK$. The
underlying space $\zk$ features in several important algebraic,
symplectic and topological constructions of torus actions, which
are the subject of the next three chapters. The decomposition of
$\zk$ as the `moment-angle complex', which first appeared
in~\cite{bu-pa98},~\cite{bu-pa99}, provided an effective
topological instrument for studying these spaces.

The basic building block in the `moment-angle' decomposition of
$\zk$ is the pair $(D^2,S^1)$ of a unit disc and circle, and the
whole construction can be extended naturally to arbitrary pairs of
spaces $(X,A)$. The resulting complex $(X,A)^\sK$ is now known as
the `polyhedral product space' over a simplicial complex~$\sK$;
this terminology was suggested by William Browder,
cf.~\cite{b-b-c-g10}. Many spaces important for toric topology
admit polyhedral product decompositions.

It has soon become clear that the construction of the moment-angle
complex $\zk$ and its generalisation $(X,A)^\sK$ is of truly
universal nature, and has remarkable functorial properties. The
most basic of these is that the construction of $\zk$ establishes
a functor from simplicial complexes and simplicial maps to spaces
with torus actions and equivariant maps. If $\sK$ is a simplicial
subdivision of a sphere (a triangulated sphere), then $\zk$ is a
manifold, and most important geometric examples of $\zk$ arise in
this way. In the case when $A$ is a point, the polyhedral product
$(X,\pt)^\sK$ interpolates between the $m$-fold wedge, or bouquet,
of $X$ (corresponding to $m$ discrete points as~$\sK$) and the
$m$-fold product of $X$ (corresponding to the full simplex as
$\sK$). Parallel to the topological and geometric study of
moment-angle complexes $\zk$ and related toric spaces, a
homotopy-theoretic study of polyhedral products $(X,A)^\sK $ has
now gained its own momentum. Basic homotopy properties of
moment-angle complexes are given in Section~\ref{basichomotopy},
while more advanced homotopy-theoretic aspects of toric topology
are the subject of Chapter~\ref{homotopy}.

The key result of this chapter is the calculation of the integral
cohomology ring of~$\zk$, carried out in Section~\ref{cohma}. The
ring $H^*(\mathcal Z_{\mathcal K})$ is shown to be isomorphic to
the Tor-algebra $\mathop{\mathrm{Tor}}_{\mathbb
Z[v_1,\ldots,v_m]}(\mathbb Z[\mathcal K],\mathbb Z)$, where
$\mathbb Z[\mathcal K]$ is the face ring of~$\mathcal K$. The
canonical bigraded structure in the Tor groups thereby acquires a
geometric interpretation in terms of the bigraded cell
decomposition of~$\mathcal Z_{\mathcal K}$. The calculation of
$H^*(\zk)$ builds upon a construction of a ring model for
\emph{cellular} cochains of~$\zk$ and the corresponding cellular
diagonal approximation, which is functorial with respect to maps
of moment-angle complexes induced by simplicial maps of~$\sK$.
This functorial property of the cellular diagonal approximation
for~$\zk$ is quite special, due to the lack of such a construction
for general cell complexes.

The construction of the moment-angle complex therefore brings the
methods of equivariant topology to bear on the study of
combinatorics of simplicial complexes, and gives a new geometric
dimension to combinatorial commutative algebra. In particular,
homological invariants of face rings, such as Tor-algebras or
algebraic Betti numbers, can be now interpreted geometrically in
terms of cohomology of moment-angle complexes. This link is
explored further in Sections~\ref{bettizk} and~\ref{macsp}.

Another important aspect of the theory of moment-angle complexes
is their connection to coordinate subspace arrangements and their
complements. As we have already seen in
Proposition~\ref{affinefr}, coordinate subspace arrangements arise
as affine varieties corresponding to face rings. Their complements
have played an important role in toric geometry and singularity
theory, and, more recently, in the theory of linkages and robotic
motion planning. Arrangements of coordinate subspaces in $\C^m$
correspond bijectively to simplicial complexes $\mathcal K$ on the
set~$[m]$, and the complement of such an arrangement deformation
retracts onto the corresponding moment-angle complex~$\zk$. In
particular, the moment-angle complex and the complement of the
arrangement have the same homotopy type. We therefore may use the
results on moment-angle complexes to obtain a description of the
cohomology groups and cup product structure of a coordinate
subspace arrangement complement. The formula obtained for the
cohomology groups is related to the general formula of
Goresky--MacPherson~\cite{go-ma88} by means of Alexander duality.

The material of this chapter, with the exception of `Additional
topics', is mainly a freshened and modernised exposition of the
results obtained by the authors and their collaborators in
\cite{bu-pa99}, \cite{bu-pa02}, \cite{b-b-p04}, \cite{pano08l}.

\medskip

Spaces with torus actions, or \emph{toric
spaces}\label{toricspac}, will be the main players throughout the
rest of this book. (See Appendix~\ref{gractions} for the key
concepts of the theory of group actions on topological spaces.)
The most basic example of a toric space is the complex
$m$-dimensional space~$\C^m$, on which the \emph{standard torus}
\[
  \T^m=\bigl\{\mb t=(t_1,\ldots,t_m)\in\C^m\colon|t_i|=1
  \quad\text{for }i=1,\ldots,m\bigr\}
\]
acts coordinatewise. That is, the action is given by
\begin{align*}
  &\T^m\times\C^m\longrightarrow\C^m,\\
  &(t_1,\ldots,t_m)\cdot(z_1,\ldots,z_m)=(t_1z_1,\ldots,t_mz_m).
\end{align*}
The quotient $\C^m/\T^m$ of this action is the \emph{positive
orthant}\label{posiorth}
\[
  \R^m_\ge=\bigl\{(y_1,\ldots,y_m)\in\R^m\colon y_i\ge0
  \quad\text{for }
  i=1,\ldots,m\bigr\},
\]
with the quotient projection given by
\begin{align*}
  \mu\colon
  \C^m&\longrightarrow\R^m_\ge,\\
  (z_1,\ldots,z_m)&\longmapsto(|z_1|^2,\ldots,|z_m|^2)
\end{align*}
(or by $(z_1,\ldots,z_m)\longmapsto(|z_1|,\ldots,|z_m|)$, but the
former is usually more preferable).

We shall use the blackboard bold capital $\T$ in the notation for
the standard torus $\T^m$ only, and use italic $T^m$ to denote an
abstract $m$-torus, i.e. a compact abelian Lie group isomorphic to
a product of $m$ circles. We shall also denote the standard unit
circle by $\mathbb S$ or $\mathbb T$ occasionally, to distinguish
it from an abstract circle~$S^1$.

All homology and cohomology groups in this chapter are with
integer coefficients, unless another coefficient group is
explicitly specified.

\section{Basic definitions}\label{defzk}

\subsection*{Moment-angle complex $\zk$} We consider the unit polydisc in the $m$-dimensional
complex space~$\C^m$:
\[
  \D^m=\bigl\{ (z_1,\ldots,z_m)\in\C^m\colon |z_i|^2\le1
  \quad\text{for } i=1,\ldots,m
  \bigr\}.
\]
The polydisc $\D^m$ is a $\T^m$-invariant subspace of $\C^m$, and
the quotient $\D^m/\T^m$ is identified with the standard unit cube
$\I^m\subset\R_\ge^m$.

\begin{construction}[moment-angle complex]\label{constrmac}
Let $\mathcal K$ be a simplicial complex on the set~$[m]$. We
recall the cubical subcomplex $\cc(\mathcal K)$ in~$\I^m$ from
Construction~\ref{cck}, which subdivides $\cone\sK$. The
\emph{moment-angle complex} $\zk$ corresponding to $\sK$ is
defined from the pullback square
\[
  \xymatrix{
  \zk \ar@{^{(}->}[r] \ar[d] & \D^m\ar[d]^\mu\\
  \cc(\sK)\ar@{^{(}->}[r] & \I^m
  }
\]
Explicitly, $\zk=\mu^{-1}(\cc(\sK))$. By construction, $\zk$ is a
$\T^m$-invariant subspace in the polydisc~$\D^m$, and the quotient
$\zk/\T^m$ is homeomorphic to $|\cone \mathcal K|$.

Using the decomposition $\cc(\sK)=\bigcup_{I\in\sK}C_I$ into
faces, see~\eqref{Iface}, it follows that
\begin{equation}\label{zkbj}
  \zk=\bigcup_{I\in\sK}B_I,
\end{equation}
where
\[
  B_I=\mu^{-1}(C_I)=\bigl\{(z_1,\ldots,z_m)\in
  \D^m\colon |z_j|^2=1\text{ for }j\notin I\bigl\},
\]
and the union in~\eqref{zkbj} is understood as the union of
subsets inside the polydisc~$\D^m$. Note that $B_I$ is a product
of $|I|$ discs and $m-|I|$ circles. Following our notational
tradition, we denote a topological 2-disc (the underlying space
of~$\D$) by $D^2$. Then we may rewrite~\eqref{zkbj} as the
following decomposition of $\zk$ into products of discs and
circles:
\begin{equation}\label{zkd2s1}
  \zk=\bigcup_{I\in\mathcal K}
  \Bigl(\prod_{i\in I}D^2\times\prod_{i\notin I}S^1 \Bigl),
\end{equation}
From now on we shall denote the space $B_I$ by~$(D^2,S^1)^I$.
Obviously, the union in~\eqref{zkbj} or~\eqref{zkd2s1} can be
taken over the maximal simplices $I\in\sK$ only.

Using the categorical language, we may consider the face category
$\mathop{\text{\sc cat}}(\mathcal K)$, and define the functor (or
diagram, see Appendix~\ref{secmc})
\begin{equation}\label{ddsdiag}
\begin{aligned}
  \mathcal D_\sK(D^2,S^1)\colon \ca(\sK)&\longrightarrow \top,\\
  I&\longmapsto B_I=(D^2,S^1)^I,
\end{aligned}
\end{equation}
which maps the morphism $I\subset J$ of $\ca(\sK)$ to the
inclusion of spaces $(D^2,S^1)^I\subset(D^2,S^1)^J$. Then we have
\[
  \zk=\mathop{\mathrm{colim}}
  \mathcal D_\sK(D^2,S^1)=\mathop{\mathrm{colim}}_{I\in\sK}
  (D^2,S^1)^I.
\]
\end{construction}

\begin{example}\

1. Let $\sK=\varDelta^{m-1}$ be the full simplex. Then
$\cc(\sK)=\I^m$ and $\zk=\D^m$.

2. Let $\sK$ be a simplicial complex on $[m]$, and let $\sK^\circ$
be the complex on $[m+1]$ obtained by adding one ghost vertex
$\circ=\{m+1\}$ to~$\sK$. Then the cubical complex
$\cc(\sK^\circ)$ is contained in the facet $\{y_{m+1}=1\}$ of the
cube $\I^{m+1}$, and
\[
  \mathcal Z_{\sK^\circ}=\zk\times S^1.
\]
In particular, if $\sK$ is the `empty' simplicial complex on
$[m]$, consisting of the empty simplex $\varnothing$ only, then
$\cc(\sK)$ is the vertex $(1,\ldots,1)\in\I^m$ and
$\zk=\mu^{-1}(1,\ldots,1)=\T^m$ is the standard $m$-torus.

For arbitrary $\sK$ on $[m]$, the moment-angle complex $\zk$
contains the $m$-torus $\T^m$ (corresponding to $\sK=\varnothing$)
and is contained in the polydisc $\D^m$ (corresponding to
$\sK=\varDelta^{m-1}$).

3. Let $\sK$ be the complex consisting of two disjoint points.
Then
\[
  \zk=(D^2\times S^1)\cup(S^1\times D^2)=\partial(D^2\times
  D^2)\cong S^3,
\]
the standard decomposition of a 3-sphere into the union of two
solid tori.

4. More generally, if $\sK=\partial\varDelta^{m-1}$ (the boundary
of a simplex), then
\begin{align*}
\zk&=(D^2\times\cdots\times D^2\times S^1)\cup
     (D^2\times\cdots\times S^1\times D^2)\cup\cdots\cup
     (S^1\times\cdots\times D^2\times
     D^2)\\
   &=\partial\bigl((D^2)^m\bigr)\cong S^{2m-1}.
\end{align*}

5. Let $\mathcal K=\quad\begin{picture}(5,5)
\put(0,0){\circle*{1}} \put(0,5){\circle*{1}}
\put(5,0){\circle*{1}} \put(5,5){\circle*{1}}
\put(0,0){\line(1,0){5}} \put(0,0){\line(0,1){5}}
\put(5,0){\line(0,1){5}} \put(0,5){\line(1,0){5}}
\put(-2.5,-1){\scriptsize 1} \put(6.3,-1){\scriptsize 3}
\put(-2.5,4){\scriptsize 4} \put(6.3,4){\scriptsize 2}
\end{picture}\quad$,
the boundary of a 4-gon. Then we have four maximal simplices
$\{1,3\}$, $\{2,3\}$, $\{1,4\}$ and $\{2,4\}$, and
\begin{align*}
  \zk&=(D^2\times S^1\times D^2\times S^1)\cup
  (S^1\times D^2\times D^2\times S^1)\\
  &\qquad\cup(D^2\times S^1\times S^1\times D^2)\cup
  (S^1\times D^2\times S^1\times D^2)\\
  &=\bigl((D^2\times S^1)\cup(S^1\times D^2)\bigr)
  \times D^2\times S^1
  \cup\bigl((D^2\times S^1)\cup(S^1\times D^2)\bigr)\times
  S^1\times D^2\\
  &=\bigl((D^2\times S^1)\cup(S^1\times D^2)\bigr)\times
  \bigl((D^2\times S^1)\cup(S^1\times D^2)\bigr)
   \cong S^3\times S^3.
\end{align*}
\end{example}

In the last example, $\sK$ is the join of $\{1,2\}$ and $\{3,4\}$.
More generally,

\begin{proposition}\label{zkjoin}
Let $\sK=\sK_1*\sK_2$; then
\[
  \zk=\mathcal Z_{\sK_1}\times\mathcal Z_{\sK_2}.
\]
\end{proposition}
\begin{proof}
By the definition of join (Construction~\ref{join}), we have
\begin{multline*}
  \mathcal Z_{\sK_1*\sK_2}=\bigcup_{I_1\in\sK_1,\,I_2\in\sK_2}
  (D^2,S^1)^{I_1\sqcup I_2}=\bigcup_{I_1\in\sK_1,\,I_2\in\sK_2}
  (D^2,S^1)^{I_1}\times(D^2,S^1)^{I_2}\\
  =\Bigl(\bigcup_{I_1\in\sK_1}(D^2,S^1)^{I_1}\Bigr)\times
  \Bigl(\bigcup_{I_2\in\sK_2}(D^2,S^1)^{I_2}\Bigr)=
  \mathcal Z_{\sK_1}\times\mathcal Z_{\sK_2}.\qedhere
\end{multline*}
\end{proof}

The topological structure of~$\zk$ is quite complicated even for
simplicial complexes~$\sK$ with few vertices. Several different
techniques will be developed to describe the topology of~$\zk$;
this is one of the main subjects of the book. To get an idea on
how $\zk$ may look like we included Exercises~\ref{exzk3p}
and~\ref{exzk5gon}, in which the topological structure of $\zk$ is
more complicated than in the examples above, but which are still
accessible by relatively elementary topological methods. As we
shall see below, the cohomology of $\zk$ may have an arbitrary
torsion (Corollary~\ref{arbitrarytorsion}), as well as Massey
products (Section~\ref{Masseymac}).

Moment-angle complexes corresponding to triangulated spheres and
manifolds are of particular interest:

\begin{theorem}\label{zkman}
Let $\mathcal K$ be a triangulation of an $(n-1)$-dimensional
sphere with $m$ vertices. Then $\zk$ is a (closed) topological
manifold of dimension~$m+n$.

If $\mathcal K$ is a triangulated manifold, then
$\zk\setminus\mu^{-1}(1,\ldots,1)$ is an (open) non-compact
manifold. Here $(1,\ldots,1)\in\I^m$ is the cone vertex, and
$\mu^{-1}(1,\ldots,1)=\T^m$.
\end{theorem}
\begin{proof}
We first construct a decomposition of the polyhedron
$|\cone(\mathcal K')|$ into `faces' similar to the faces of a
simple polytope (in fact, in the case of the nerve complex
$\mathcal K=\mathcal K_P$, our faces will be exactly the faces
of~$P$). The vertices $i\in\mathcal K$ are also vertices of the
barycentric subdivision $\mathcal K'$, and we set
\begin{equation}\label{kfacet}
  F_i=\st_{\mathcal K'}\{i\}, \quad 1\le i\le m,
\end{equation}
(i.e. $F_i$ is the star of the $i$th vertex of~$\sK$ in the
barycentric subdivision~$\sK'$). We refer to $F_i$ as
\emph{facets} of our face decomposition, and define a face of
codimension~$k$ as a nonempty intersection of a set of $k$ facets.
In particular, the vertices of our face decomposition are the
barycentres of $(n-1)$-dimensional simplices of~$\mathcal K$. Such
a barycentre $b$ corresponds to a maximal simplex
$I=\{i_1,\ldots,i_n\}$ of~$\sK$, and we denote by $U_I$ the open
subset in $|\cone(\mathcal K')|$ obtained by removing all faces
not containing~$b$. We observe that any point of $|\cone\mathcal
K'|$ is contained in $U_I$ for some~$I\in\sK$.

Under the map $\cone(i_c)\colon|\cone(\sK')|\to\I^m$ of
Construction~\ref{cck} the facet $F_i$ is mapped to the
intersection of $\cc(\sK)$ with the $i$th coordinate plane
$y_i=0$. Therefore, the image of $U_I$ under the map $\cone(i_c)$
is given by
\begin{equation}\label{coneui}
  W_I=\cone(i_c)(U_I)=\{(y_1,\ldots,y_m)\in\cc(\sK)\colon y_i\ne0\quad\text{for }i\notin
  I\}.
\end{equation}

Assume now that $\mathcal K$ is a triangulated sphere. Then
$|\cone\mathcal K'|$ is homeomorphic to an $n$-dimensional
disc~$D^n$, and each $U_I$ is homeomorphic to an open subset in
$\R^n_\ge$ preserving the dimension of faces. (This means that
$|\cone\mathcal K'|$ is a \emph{manifold with
corners}\label{maniwithcor}, see Section~\ref{locst}). By
identifying $|\cone\mathcal K'|$ with $\cc(\mathcal K)$ and
further identifying $\cc(\mathcal K)$ with the quotient
$\zk/\T^m$, we obtain that each point of $\zk$ has a neighbourhood
of the form $\mu^{-1}(W_I)$. It follows from~\eqref{coneui} that
the latter is homeomorphic to an open subset in
$\mu_n^{-1}(\R^n_\ge)\times T^{m-n}=\C^n\times T^{m-n}$, where
$\R^n$ is the coordinate $n$-plane corresponding to
$i_1,\ldots,i_n$, the map $\mu_n\colon\C^n\to\R^n_\ge$ is the
restriction of $\mu\colon\C^m\to\R^m_\ge$ to the corresponding
coordinate plane in~$\C^m$, and the torus $T^{m-n}$ sits in the
complementary coordinate $(m-n)$-plane in~$\C^m$. An open subset
in $\C^n\times T^{m-n}$ with $n\ge1$ can be regarded as an open
subset in~$\R^{m+n}$, and therefore $\zk$ is an
$(m+n)$-dimensional manifold.

If $\mathcal K$ is a triangulated manifold, then $|\cone\mathcal
K'|$ is not a manifold because of the singularity at the cone
vertex~$v$. However, by removing this vertex we obtain a
(non-compact) manifold whose boundary is~$|\mathcal K|$. Using the
face decomposition defined by~\eqref{kfacet} we obtain that
$|\cone\mathcal K'|\setminus v$ is locally homeomorphic to
$\R^n_\ge$ preserving the dimension of faces (i.e. $|\cone\mathcal
K'|\setminus v$ is a non-compact manifold with corners). Under the
identification of $|\cone\mathcal K'|$ with $\cc(\mathcal K)$ the
vertex of the cone is mapped to the vertex $(1,\ldots,1)\in\I^m$,
and $\mu^{-1}(1,\ldots,1)=\T^m$. Therefore,
\[
  \mu^{-1}\bigl(\cc(\mathcal K)\backslash(1,\ldots,1)\bigr)=
  \zk\setminus\mu^{-1}(1,\ldots,1)
\]
is an $(m+n)$-dimensional non-compact manifold.
\end{proof}

\begin{remark}
A pair of spaces $(X,A)$ where $A$ is a compact subset in $X$ is
called a \emph{Lefschetz pair}\label{lefpair} if $X\setminus A$ is
an open (non-compact) manifold. We therefore obtain that
$(\zk,\mu^{-1}(1,\ldots,1))$ is a Lefschetz pair whenever $\sK$ is
a triangulated manifold.
\end{remark}

We therefore refer to moment-angle complexes $\zk$ corresponding
to triangulated spheres as \emph{moment-angle
manifolds}\label{mamanifo}. \emph{Polytopal} moment-angle
manifolds $\mathcal Z_{\mathcal K_P}$, corresponding to the nerve
complexes $\sK_P$ of simple polytopes (see Example~\ref{polsph}),
are particularly important. As we shall see in
Chapter~\ref{mamanifolds}, polytopal moment-angle manifolds are
smooth. A smooth structure also exists on moment-angle manifolds
$\zk$ corresponding to starshaped spheres $\sK$ (i.e. underlying
complexes of complete simplicial fans, see Section~\ref{simsph}).
In general, the smoothness of $\zk$ is open.

%\begin{problem}
%Describe the class of triangulated spheres $\sK$ for which the
%corresponding moment-angle manifolds $\zk$ are smooth.
%\end{problem}

The geometry of moment-angle manifolds is nice and rich; it is the
subject of Chapter~\ref{mamanifolds}.

\subsection*{Real moment-angle complex $\rk$}
The construction of the moment-angle complex $\zk$ has a real
analogue, in which the complex space $\C^m$ is replaced by the
real space $\R^m$, the complex polydisc $\D^m$ is replaced by the
`big' cube
\[
  [-1,1]^m=\bigl\{ (u_1,\ldots,u_m)\in\R^m\colon |u_i|^2\le1
  \quad\text{for } i=1,\ldots,m
  \bigr\},
\]
the standard torus $\T^m$ is replaced by the `real torus' $\Z_2^m$
(the product of $m$ copies of the group $\Z_2=\{-1,1\}$), and the
pair $(D^2,S^1)$ is replaced by $(D^1,S^0)$, where $S^0$ (a pair
of points) is the boundary of the segment~$D^1$. The group
$(\Z_2)^m$ acts on the big cube $[-1,1]^m$ coordinatewise, with
quotient the standard `small' cube~$\I^m$. The quotient projection
$[-1,1]^m\to\I^m$ may be described by the map
\[
  \rho\colon(u_1,\ldots,u_m)\mapsto(u_1^2,\ldots,u_m^2).
\]

\begin{construction}[real moment-angle complex]\label{realmac}
Given a simplicial complex $\mathcal K$ on~$[m]$, define the
\emph{real moment-angle complex} $\rk$ from the pullback square
\[
  \xymatrix{
  \rk \ar@{^{(}->}[r] \ar[d] & [-1,1]^m\ar[d]^\rho\\
  \cc(\sK)\ar@{^{(}->}[r] & \I^m
  }
\]
Explicitly, $\rk=\rho^{-1}(\cc(\sK))$. By construction, $\rk$ is a
$\Z_2^m$-invariant subspace in the `big' cube $[-1,1]^m$, and the
quotient $\rk/\Z_2^m$ is homeomorphic to $|\cone \mathcal K|$.

$\rk$ is a cubical subcomplex in $[-1,1]^m$ obtained by reflecting
the subcomplex $\cc(\sK)\subset\I^m=[0,1]^m$ at all $m$ coordinate
hyperplanes of~$\R^m$. If $\sK=\sK_P$ is the nerve complex of a
simple polytope~$P$, then $\cc(\sK_P)$ can be viewed as a cubical
subdivision of $P$ embedded piecewise linearly into $\R^m_\ge$
(see Construction~\ref{cubpol}). In this case $\mathcal
R_{\mathcal K_P}$ is obtained by reflecting the image of $P$ at
all coordinate planes.

By analogy with \eqref{zkd2s1}, we have
\[
  \rk=\bigcup_{I\in\mathcal K}
  \Bigl(\prod_{i\in I}D^1\times\prod_{i\notin I}S^0 \Bigl).
\]
\end{construction}

\begin{example}\label{rkexam}\

1. Let $\sK=\partial\varDelta^{m-1}$ be the boundary of the
standard simplex. Then $\rk\cong S^{m-1}$ is the boundary of the
cube $[-1,1]^m$. For $m=3$ this complex is obtained by reflecting
the complex shown in Fig.~\ref{figcck}~(b) at all 3 coordinate
planes of~$\R^3$.

2. Let $\sK$ consist of $m$ disjoint points. Then $\rk$ is the
1-dimensional skeleton (graph) of the cube $[-1,1]^m$. For $m=3$
this complex is obtained by reflecting the complex shown in
Fig.~\ref{figcck}~(a) at all 3 coordinate planes of~$\R^3$.

3. More generally, let $\sK=\sk^i\varDelta^{m-1}$ be the
\emph{$i$-dimensional skeleton} of $\varDelta^{m-1}$ (i.e. the set
of all faces of $\varDelta^{m-1}$ of dimension $\le i$). Then
$\rk$ is the $(i+1)$-dimensional skeleton of the cube~$[-1,1]^m$.
\end{example}

The following analogue of Theorem~\ref{zkman} holds, and is proved
similarly:

\begin{theorem}
Let $\mathcal K$ be a triangulation of an $(n-1)$-dimensional
sphere with $m$ vertices. Then $\rk$ is a (closed) topological
manifold of dimension~$n$.

If $\mathcal K$ is a triangulated manifold, then
$\rk\setminus\rho^{-1}(1,\ldots,1)$ is an (open) non-compact
manifold, where $\rho^{-1}(1,\ldots,1)=\{-1,1\}^m$.
\end{theorem}

%Unlike the case of moment-angle complexes $\zk$, t
The real moment-angle complexes corresponding to polygons can be
identified easily (compare Exercise~\ref{exzk5gon}):

\begin{proposition}\label{rpolygon}
Let $\sK$ be the boundary of an $m$-gon. Then $\rk$ is
homeomorphic to an oriented surface $S_g$ of genus
$g=1+(m-4)2^{m-3}$.
\end{proposition}
\begin{proof}
We observe that the manifold $\rk$ is orientable (an exercise).
Since it is 2-dimensional, its topological type is determined by
the Euler characteristic. Now $\rk$ is obtained by reflecting the
$m$-gon embedded into $\R^m_\ge$ at $m$ coordinate hyperplanes, so
that $\rk$ is patched from $2^m$ polygons, meeting by 4 at each
vertex. Therefore, the number of vertices is $m2^{m-2}$, the
number of edges is $m2^{m-1}$, and the Euler characteristic is
\[
  \chi(\rk)=2^{m-2}(4-m)=2-2g.
\]
The result follows.
\end{proof}

\section{Polyhedral products}
Decomposition~\eqref{zkd2s1} of $\zk$
which uses the disc and circle $(D^2,S^1)$ is readily generalised
to arbitrary pairs of spaces:

\begin{construction}[polyhedral product]\label{nsc}
Let $\sK$ be a simplicial complex on~$[m]$ and
\[
  (\mb X,\mb A)=\{(X_1,A_1),\ldots,(X_m,A_m)\}
\]
be a collection of $m$ pairs of spaces, $A_i\subset X_i$. For each
simplex $I\in[m]$ we set
\begin{equation}\label{XAI}
  (\mb X,\mb A)^I=\bigl\{(x_1,\ldots,x_m)\in
  \prod_{j=1}^m X_j\colon\; x_j\in A_j\quad\text{for }j\notin I\bigl\}
\end{equation}
and define the \emph{polyhedral product} of $(\mb X,\mb A)$
corresponding to $\sK$ by
\[
  (\mb X,\mb A)^{\sK}=\bigcup_{I\in\mathcal K}(\mb X,\mb A)^I=
  \bigcup_{I\in\mathcal K}
  \Bigl(\prod_{i\in I}X_i\times\prod_{i\notin I}A_i\Bigl).
\]
Using the categorical language, we can define the
$\ca{(\sK)}$-diagram
\begin{equation}\label{dxadiag}
\begin{aligned}
  \mathcal D_\sK(\mb X,\mb A)\colon \ca(\sK)&\longrightarrow \top,\\
  I&\longmapsto (\mb X,\mb A)^I,
\end{aligned}
\end{equation}
which maps the morphism $I\subset J$ of $\ca(\sK)$ to the
inclusion of spaces $(\mb X,\mb A)^I\subset(\mb X,\mb A)^J$. Then
we have
\[
  (\mb X,\mb A)^\sK=\mathop{\mathrm{colim}}
  \mathcal D_\sK(\mb X,\mb A)=\mathop{\mathrm{colim}}_{I\in\sK}
  (\mb X,\mb A)^I.
\]

In the case when all the pairs $(X_i,A_i)$ are the same, i.e.
$X_i=X$ and $A_i=A$ for $i=1,\ldots,m$, we use the notation
$(X,A)^\sK$ for $(\mb X,\mb A)^\sK$. Also, if each $X_i$ is a
pointed space and $A_i=\pt$, then we use the abbreviated notation
$\mb X^\sK$ for $(\mb X,\pt)^\sK$, and $X^\sK$ for $(X,\pt)^\sK$.
\end{construction}

\begin{remark}
The decomposition of $\zk$ into a union of products of discs and
circles first appeared in~\cite{bu-pa98} (in the polytopal case)
and in~\cite{bu-pa00-2} (in general). The term `moment-angle
complex' for $\zk=(D^2,S^1)^\sK$ was also introduced
in~\cite{bu-pa00-2}, where several other examples of polyhedral
products $(X,A)^\sK$ were considered. The definition of
$(X,A)^\sK$ for an arbitrary pair of spaces $(X,A)$ was suggested
to the authors by N.~Strickland (in a private communication, and
also in an unpublished note) as a general framework for the
constructions of~\cite{bu-pa00-2}; it was also included in the
final version of~\cite{bu-pa00-2} and in~\cite{bu-pa02}. Further
generalisations of $(X,A)^\sK$ to a set of pairs of spaces $(\mb
X,\mb A)$ were studied in the work of Grbi\'c and
Theriault~\cite{gr-th07}, as well as Bahri, Bendersky, Cohen and
Gitler~\cite{b-b-c-g10}, where the term `polyhedral product' was
introduced (following a suggestion of W.~Browder). Since 2000, the
terms `generalised moment-angle complex', `$\sK$-product' and
`partial product space' have been also used to refer to the
spaces~$(X,A)^\sK$.
\end{remark}

Recall that a \emph{map of pairs} $(X,A)\to(X',A')$ is a
commutative diagram
\begin{equation}\label{mapofpairs}
\begin{array}{ccc}
  A&\hookrightarrow&X\\
  \downarrow&&\downarrow\\
  A'&\hookrightarrow&X'.
\end{array}
\end{equation}
We refer to $(X,A)$ as a \emph{monoid pair}\label{monoidpair} if
$X$ is a topological monoid (a space with a continuous associative
multiplication and unit), and $A$ is a submonoid. A \emph{map of
monoid pairs} is a map of pairs in which all maps
in~\eqref{mapofpairs} are homomorphisms.

If $(X,A)$ is a monoid pair, then a set map
$\varphi\colon[l]\to[m]$ induces a map
\begin{equation}\label{monprod}
  \psi\colon \prod_{i=1}^l X\to\prod_{i=1}^m X,\quad
  (x_1,\ldots,x_l) \mapsto (y_1,\ldots,y_m),
\end{equation}
where
\[
  y_j=\prod_{i\in\varphi^{-1}(j)}x_i,\qquad\text{for } j=1,\ldots,m,
\]
and we set $y_j=1$ if $\varphi^{-1}(j)=\varnothing$.

\begin{proposition}
If $(X,A)$ is a monoid pair, then $(X,A)^\sK$ is an invariant
subspace of $\prod_{i=1}^m X$  with respect to the coordinatewise
action of $\prod_{i=1}^m A$ on $\prod_{i=1}^m X$.
\end{proposition}
\begin{proof}
Indeed, $(X,A)^I\subset\prod_{i=1}^m X$ is an invariant subset for
each $I\in\sK$.
\end{proof}

The following proposition describes the functorial properties of
the polyhedral product in its $\sK$ and $(\mb X,\mb A)$ arguments.

\begin{proposition}\label{mafunc}\ \nopagebreak
\begin{itemize}
\item[(a)] A set of maps of pairs $(\mb X,\mb A)\to(\mb X',\mb
A')$ induces a map of polyhedral products $(\mb X,\mb A)^\sK
\to(\mb X',\mb A')^\sK$. If two sets of maps $(\mb X,\mb A)\to(\mb
X',\mb A')$ are componentwise homotopic, then the induced maps
$(\mb X,\mb A)^\sK \to(\mb X',\mb A')^\sK$ are also homotopic.

\item[(b)] An inclusion of a simplicial subcomplex
$\mathcal L\hookrightarrow\sK$ induces an inclusion of polyhedral
products $(\mb X,\mb A)^{\mathcal L}\hookrightarrow(\mb X,\mb
A)^\sK $.

\item[(c)] If $(X,A)$ is a monoid pair, then for any
simplicial map $\varphi\colon\mathcal L\to\sK$ of simplicial
complexes on the sets $[l]$ and $[m]$ respectively, the
map~\eqref{monprod} restricts to a map of polyhedral products
$\varphi_{\mathcal Z}\colon(X,A)^{\mathcal L}\to(X,A)^\sK $.

\item[(d)] If $(X,A)$ is a commutative monoid pair, then the
restriction
\[
  \psi|_A\colon \prod_{i=1}^l A\to\prod_{i=1}^m A
\]
of~\eqref{monprod} is a homomorphism, and the induced map
$\varphi_{\mathcal Z}\colon(X,A)^{\mathcal L}\to(X,A)^\sK$ is
$\psi|_A$-equivariant, i.e.
\[
  \varphi_{\mathcal Z}(\mb a\cdot\mb x)=
  \psi|_A(\mb a)\cdot\varphi_{\mathcal Z}(\mb x)
\]
for all $\mb a=(a_1,\ldots,a_l)\in\prod_{i=1}^l A$ and $\mb
x=(x_1,\ldots,x_l)\in(X,A)^{\mathcal L}$.
\end{itemize}
\end{proposition}
\begin{proof}
For (a), we observe that a set of maps $(\mb X,\mb A)\to(\mb
X',\mb A')$ induces a map $(\mb X,\mb A)^I\to(\mb X',\mb A')^I$
for each $I\in\sK$, and these maps corresponding to different
$I,J\in\sK$ are compatible on the intersection $(\mb X,\mb
A)^I\cap(\mb X,\mb A)^J=(\mb X,\mb A)^{I\cap J}$. We therefore
obtain a map $(\mb X,\mb A)^\sK \to(\mb X',\mb A')^\sK$. A
componentwise homotopy between two maps $(\mb X,\mb A)\to(\mb
X',\mb A')$ can be thought of as a map of pairs $(\mb
X\times\I,\mb A\times\I)\to(\mb X',\mb A')$, where $\mb X\times\I$
consists of spaces $X_i\times\I$. It therefore induces a map of
polyhedral products
\[
  (\mb X\times\I,\mb A\times\I)^\sK\to(\mb X',\mb A')^\sK
\]
where $(\mb X\times\I,\mb A\times\I)^\sK\cong(\mb X,\mb A)^\sK
\times(\I,\I)^\sK=(\mb X,\mb A)^\sK \times\I^m$. By restricting
the resulting map $(\mb X,\mb A)^\sK \times\I^m\to(\mb X',\mb
A')^\sK$ to the diagonal of the cube $\I^m$ we obtain a homotopy
between the two induced maps $(\mb X,\mb A)^\sK \to(\mb X',\mb
A')^\sK$.

To prove~(b) we just observe that if $\mathcal L$ is a subcomplex
of~$\sK$, then for each $I\in\mathcal L$ we have $(\mb X,\mb
A)^I\subset\zk$.

To prove (c) we observe that for any subset $I\subset[m]$ we have
$\psi\bigl((X,A)^I\bigr)\subset(X,A)^{\varphi(I)}$. Let
$I\in\mathcal L$, so that $(X,A)^I\subset(X,A)^{\mathcal L}$.
Since $\varphi$ is a simplicial map, we have $\varphi(I)\in\sK$
and $(X,A)^{\varphi(I)}\subset(X,A)^\sK $. Therefore, the map
$\psi$ restricts to a map of polyhedral products $(X,A)^{\mathcal
L}\to(X,A)^\sK $.

Statement (d) is proved by a direct calculation:
\begin{align*}
  \varphi_{\mathcal Z}(\mb a\cdot\mb x)&=
  \varphi_{\mathcal Z}(a_1x_1,\ldots,a_lx_l)=
  \Bigr(\prod_{i\in\varphi^{-1}(1)}a_ix_i,\ldots,
  \prod_{i\in\varphi^{-1}(m)}a_ix_i\Bigr)\\
  &=
  \Bigr(\prod_{i\in\varphi^{-1}(1)}a_i\prod_{i\in\varphi^{-1}(1)}x_i,\ldots,
  \prod_{i\in\varphi^{-1}(m)}a_i\prod_{i\in\varphi^{-1}(m)}x_i\Bigr)
  \\
  &=  \Bigr(\prod_{i\in\varphi^{-1}(1)}a_i,\ldots,
  \prod_{i\in\varphi^{-1}(m)}a_i\Bigr)\cdot
  \Bigr(\prod_{i\in\varphi^{-1}(1)}x_i,\ldots,
  \prod_{i\in\varphi^{-1}(m)}x_i\Bigr)\\
  &=\psi|_A(\mb a)\cdot\varphi_{\mathcal Z}(\mb x).
\end{align*}
Note that we have used the commutativity of $X$ in the third
identity.
\end{proof}

We state the most important particular case of
Proposition~\ref{mafunc} separately:

\begin{proposition}\label{functmac1}
A simplicial map $\mathcal L\to\sK$ of simplicial complexes on the
sets $[l]$ and $[m]$ gives rise to a map of moment-angle complexes
$\mathcal Z_{\mathcal L}\to\zk$ which is equivariant with respect
to the induced homomorphism of tori $\T^l\to\T^m$.
\end{proposition}

The polyhedral product construction behaves nicely with respect to
the join operation (the proof is exactly the same as that of
Proposition~\ref{zkjoin}):

\begin{proposition}
We have that
\[
  (\mb X,\mb A)^{\sK_1*\sK_2}=
  (\mb X,\mb A)^{\sK_1}\times(\mb X,\mb A)^{\sK_2}.
\]
\end{proposition}

\begin{example}\label{exkpo}\

1. The moment-angle complex $\zk$ is the polyhedral product
$(D^2,S^1)^\sK$ (when considered abstractly) or $(\D,\mathbb
S)^\sK$ (when viewed as a subcomplex in~$\D^m$).

2. For the cubical complex of~\eqref{fcck} we have
\[
  \cc(\mathcal K)=(I^1,1)^\sK,
\]
where $I^1=[0,1]$ is the unit interval and $1$ is its edge. The
quotient map $\zk\to|\cone\sK|$ is the map of polyhedral products
$(D^2,S^1)^\sK\to(I^1,1)^\sK$ induced by the map of pairs
$(D^2,S^1)\to(I^1,1)$, which is the quotient map by the
$S^1$-action.

3. For the real moment-angle complex we have
\[
  \rk=(D^1,S^0)^\sK,
\]
Where $D^1=[-1,1]$ is a 1-disc, and $S^0=\{-1,1\}$ is its
boundary.

4. If $\sK$ consists of $m$ disjoint points and $A_i=\pt$, then
\[
  (\mb X,pt)^\sK=\mb X^\sK=X_1\vee X_2\vee\cdots\vee X_m
\]
is the \emph{wedge}\label{wedgedef} (or \emph{bouquet}) of the
$X_i$'s.

5. More generally, consider the sequence of inclusions of
skeleta
\[
  \sk^0\varDelta^{m-1}%\subset\sk^1\varDelta^{m-1}
  \subset\cdots
  \subset\sk^{m-2}\varDelta^{m-1}\subset\sk^{m-1}\varDelta^{m-1}=\varDelta^{m-1}.
\]
It gives rise to a filtration in the product $X_1\times
X_2\times\cdots\times X_m$:
\[
\begin{array}{ccccc}
  \mb X^{{\sk^0\varDelta^{m-1}}}\;\;\subset&\cdots\;\;\subset
  &\mb X^{\sk^{m-2}\varDelta^{m-1}}&\subset
  &\mb X^{\sk^{m-1}\varDelta^{m-1}}\\
  \|&&\|&&\|\\
  X_1\vee X_2\vee\cdots\vee X_m&&\mb X^{\partial\varDelta^{m-1}}&&
  X_1\times X_2\times\cdots\times X_m
\end{array}
\]
Its second-to-last term, $\mb X^{\partial\varDelta^{m-1}}$ is
known to topologists as the \emph{fat wedge}\label{dfatwedge} of
the $X_i$'s. Explicitly, the fat wedge of a sequence of pointed
spaces $X_1,\ldots,X_m$ is
\[
  (X_1\times\cdots\times X_{m-1}\times\pt)\cup
     (X_1\times\cdots\times\pt\times X_m)\cup\cdots\cup
     (\pt\times\cdots\times X_{m-1}\times
     X_m),
\]
where the union is taken inside the product $X_1\times
X_2\times\cdots\times X_m$.

The filtration above was considered by G.~Porter~\cite{port66},
who obtained a decomposition of its loop spaces into a wedge in
the case when each $X_i$ is a suspension, generalising the
\emph{Hilton--Milnor Theorem}\label{HMTheorem}. We shall consider
this decomposition in more detail in Section~\ref{stabledec}.
\end{example}

\subsection*{Exercises.}
\begin{exercise}
Show that if $\sK$ is a triangulated sphere, then the manifold
$\rk$ is orientable. (Hint: use the fact that $\rk$ is obtained by
reflecting the $n$-ball $|\cone\sK|$ at all $m$ coordinate
hyperplanes of~$\R^m$ to extend the orientation from $|\cone\sK|$
to the whole of~$\rk$.)
\end{exercise}

\begin{exercise}\label{realskel}
Show that if $\sK=\sk^i\varDelta^{m-1}$, then $\rk$ is homotopy
equivalent to a wedge of $(i+1)$-dimensional spheres (see
Example~\ref{rkexam}.3). The number of spheres is given by
$\sum_{k=i+2}^m\bin mk\bin{k-1}{i+1}$.
%$i=0$ $2^{m-1}(m-2)+1$ spheres
\end{exercise}

\begin{exercise}\label{exzk3p}
Let $\sK$ be the complex consisting of three disjoint points. Show
that $\zk$ is homotopy equivalent to the following wedge (bouquet)
of spheres:
\[
  \zk\cong S^3\vee S^3\vee S^3\vee S^4\vee S^4.
\]
(Hint: compare the case $m=3,i=0$ of the previous exercise; it may
also help to look at the realisation of $\zk$ as the complement of
a coordinate subspace arrangement (up to homotopy), see
Section~\ref{arran}.)
\end{exercise}

\begin{exercise}\label{exzk5gon}
Let $\sK$ be the boundary of a 5-gon. Show that the manifold $\zk$
is homeomorphic to $(S^3\times S^4)^{\cs5}$, a connected sum of 5
copies of $S^3\times S^4$. (This may be a difficult one; the
general statement is given in Theorem~\ref{zpstacked} below. The
reader may return to this exercise after reading
Section~\ref{mampol}, see Exercise~\ref{pentsurg}.)
\end{exercise}

\begin{exercise}
If $\sK$ is a triangulation of an $(n-1)$-sphere with $m$
vertices, then for any $k>0$ the polyhedral product
$(D^k,S^{k-1})^\sK$ is a manifold of dimension $m(k-1)+n$.
\end{exercise}

\begin{exercise}
More generally, if $(M,\partial M)$ is a manifold with boundary
and $\sK$ is a triangulated sphere, then $(M,\partial M)^\sK$ is a
manifold (without boundary). If $(M,\partial M)$ is a $PL$
manifold, and $\sK$ is a $PL$ sphere, then the manifold
$(M,\partial M)^\sK$ is also $PL$.
\end{exercise}

\begin{exercise}
Let $\sK_I$ be the full subcomplex corresponding to $I\subset[m]$.
Then $\mathcal Z_{\mathcal K_I}$ is a retract of~$\zk$.
\end{exercise}

\section{Homotopical properties}\label{basichomotopy}
Two key observations of this section constitute the basis for the
subsequent applications of the commutative algebra apparatus of
Chapter~\ref{facerings} to toric topology. First, the cohomology
of the polyhedral product space of the form $(\C P^\infty)^\sK=(\C
P^\infty,\pt)^\sK$ is isomorphic to the face ring $\Z[\sK]$
(Proposition~\ref{homsrs}). Second, the moment-angle complex
$\zk=(D^2,S^1)^\sK$ is the homotopy fibre of the canonical
inclusion $(\C P^\infty)^\sK\hookrightarrow(\C P^\infty,\C
P^\infty)^\sK=(\C P^\infty)^m$ (Theorem~\ref{zkhofib}).

The classifying space $BS^1$ of the circle $S^1$ is the
infinite-dimensional complex projective space $\C P^\infty$. The
classifying space $BT^m$ of the $m$-torus $T^m$ is the product
$(\C P^\infty)^m$ of $m$ copies of $\C P^\infty$. The universal
principal $S^1$-bundle is the infinite Hopf bundle $S^\infty\to\C
P^\infty$ (the direct limit of Hopf bundles $S^{2k+1}\to\C P^k$),
so the total space $ET^m$ of the universal principal $T^m$-bundle
over $BT^m$ can be identified with the $m$-fold product of the
infinite-dimensional sphere~$S^\infty$.

The integral cohomology ring of $BT^m$ is isomorphic to the
polynomial ring $\Z[v_1,\ldots,v_m]$, \ $\deg v_i=2$ (this
explains our choice of grading). The space $BT^m$ has the
canonical cell decomposition, in which each factor $\C P^\infty$
has one cell in every even dimension. The polyhedral product
\[
  (\C P^\infty)^\sK=\bigcup_{I\in\sK}(\C P^\infty,\pt)^I
\]
is a cellular subcomplex in $BT^m=(\C P^\infty)^m$.

\begin{proposition}
\label{homsrs} The cohomology ring of $(\C P^\infty)^\sK$ is
isomorphic to the face ring~$\Z[\sK]$. The inclusion of a cellular
subcomplex
\[
  i\colon(\C P^\infty)^\sK\hookrightarrow(\C P^\infty)^m
\]
induces the quotient projection in cohomology:
\[
  i^*\colon \Z[v_1,\ldots,v_m]\to
  \Z[v_1,\ldots,v_m]/\mathcal I_{\mathcal K}=\Z[\sK].
\]
\end{proposition}
\begin{proof}
Since $(\C P^\infty)^m$ has only even-dimensional cells and $(\C
P^\infty)^\sK$ is a cellular subcomplex, the cohomology of both
spaces coincides with their cellular cochains. Let $D^{2k}_j$
denote the $2k$-dimensional cell in the $j$th factor of $(\C
P^\infty)^m$. The cellular cochain group $\sC^*((\C P^\infty)^m)$
has basis of cochains $\bigl(D^{2k_1}_{j_1}\cdots
D^{2k_p}_{j_p}\bigr)^*$ dual to the products of cells
$D^{2k_1}_{j_1}\times\cdots\times D^{2k_p}_{j_p}$. The cochain map
\[
  \sC^*\bigl((\C P^\infty)^m\bigr)\to\sC^*\bigl((\C P^\infty)^\sK\bigr)
\]
induced by the inclusion $(\C P^\infty)^\sK\hookrightarrow (\C
P^\infty)^m$ is an epimorphism with kernel generated by those
cochains $\bigl(D^{2k_1}_{j_1}\cdots D^{2k_p}_{j_p}\bigr)^*$ for
which $\{j_1,\ldots,j_p\}\notin\sK$. Under the identification of
$\sC^*((\C P^\infty)^m)$ with $\Z[v_1,\ldots,v_m]$, a cochain
$\bigl(D^{2k_1}_{j_1}\cdots D^{2k_p}_{j_p}\bigr)^*$ is mapped to
the monomial $v_{j_1}^{k_1}\cdots v_{j_p}^{k_p}$. Therefore,
$\sC^*((\C P^\infty)^\sK)$ is identified with the quotient of
$\Z[v_1,\ldots,v_m]$ by the subgroup generated by all monomials
$v_{j_1}^{k_1}\cdots v_{j_p}^{k_p}$ with
$\{j_1,\ldots,j_p\}\notin\sK$. By Proposition~\ref{frbasis}, this
quotient is exactly~$\Z[\sK]$.
\end{proof}

We consider the Borel construction $ET^m\times_{T^m}\zk$ for the
$T^m$-space $\zk$ (see Appendix~\ref{gractions}).

\begin{theorem}\label{zkhofib}
The inclusion $i\colon(\C P^\infty)^\sK\hookrightarrow(\C
P^\infty)^m$ is decomposed into a composition of a homotopy
equivalence
\[
  h\colon(\C P^\infty)^\sK\stackrel{\simeq}\longrightarrow ET^m\times_{T^m}\zk
\]
and the fibration $p\colon ET^m\times_{T^m}\zk\to BT^m=(\C
P^\infty)^m$ with fibre~$\zk$.

In particular,  the moment-angle complex $\zk$ is the homotopy
fibre of the canonical inclusion $i\colon(\C
P^\infty)^\sK\hookrightarrow(\C P^\infty)^m$.
\end{theorem}
\begin{proof}
We use functoriality and homotopy invariance of the polyhedral
product construction (Proposition~\ref{mafunc}). We have
$\zk=\bigcup_{I\in\sK}(D^2,S^1)^I$, see~\eqref{zkbj}, and each
subset $(D^2,S^1)^I$ is $T^m$-invariant. We therefore obtain the
following decomposition of the Borel construction as a polyhedral
product:
\begin{multline*}
  ET^m\times_{T^m}\zk=
  \bigcup_{I\in\sK}\bigl(ET^m\times_{T^m}(D^2,S^1)^I\bigr)=
  \bigcup_{I\in\sK}\bigl(S^\infty\times_{S^1}D^2,
  S^\infty\times_{S^1}S^1\bigr)^I\\=
  (S^\infty\times_{S^1}D^2,
  S^\infty\times_{S^1}S^1\bigr)^\sK,
\end{multline*}
where $S^\infty=ET^1=ES^1$.

Now consider the commutative diagram
\begin{equation}\label{diwdr}
\begin{CD}
  \pt @>>> S^\infty\times_{S^1}S^1 @>>> \pt\\
  @VVV @VVV @VVV\\
  \C P^\infty @>j>> S^\infty\times_{S^1}D^2 @>f>> \C P^\infty,
\end{CD}
\end{equation}
where $j$ is the inclusion of the zero section in a disc bundle
(note that $\C P^\infty=S^\infty/S^1$), and $f$ is the projection
map from the bundle to its Thom space,
\[
  f\colon S^\infty\times_{S^1}D^2\to(S^\infty\times_{S^1}D^2)/(S^\infty\times_{S^1}S^1)
  \cong\C P^\infty.
\]
Since $S^\infty\times_{S^1}S^1=S^\infty$ and $D^2$ are
contractible, the composite maps $f\circ j$ and $j\circ f$ are
homotopic to the identity. It follows that we have a homotopy
equivalence of pairs $(\C
P^\infty,\pt)\to(S^\infty\times_{S^1}D^2,S^\infty\times_{S^1}S^1)$
which induces a homotopy equivalence of polyhedral products
\[
  h\colon(\C P^\infty)^\sK\to
  (S^\infty\times_{S^1}D^2,S^\infty\times_{S^1}S^1)^\sK=ET^m\times_{T^m}\zk.
\]

In order to establish the factorisation $i=p\circ h$ we consider
the diagram
\[
\begin{CD}
  \pt @>>> S^\infty\times_{S^1}S^1 @>>> \C P^\infty\\
  @VVV @VVV @|\\
  \C P^\infty @>j>> S^\infty\times_{S^1}D^2 @>g>> \C P^\infty,
\end{CD}
\]
where $g$ is the projection of the disc bundle onto its base (note
that his map is different from the map $f$ above). By passing to
the induced maps of polyhedral products we obtain the
factorisation of
\[
  i\colon(\C P^\infty)^\sK\stackrel h\longrightarrow
  (S^\infty\times_{S^1}D^2,S^\infty\times_{S^1}S^1)^\sK
  \stackrel p\longrightarrow(\C P^\infty,\C P^\infty)^\sK
\]
into the composition of $h$ and $p$.
\end{proof}

%\begin{remark}
%If we replace the map $T$ in diagram~\eqref{diwdr} by the
%projection map of the disc bundle, then the lower row would define
%a deformation retraction $S^\infty\times_{S^1}D^2 \to \C
%P^\infty$. However, diagram~\eqref{diwdr} would no longer be
%commutative (only commutative up to homotopy), and the required
%retraction $ET^m\times_{T^m}\zk\to(\C P^\infty,\pt)^\sK$ would not
%follow directly from the functoriality of the polyhedral product.
%\end{remark}

The following statement is originally due to Davis and
Januszkiewicz~\cite[Theorem~4.8]{da-ja91} (they used a different
model for~$\zk$, which will be discussed in Section~\ref{mampol}).

\begin{corollary}
\label{cohombk} The equivariant cohomology ring of the
moment-angle complex $\zk$ is isomorphic to the face ring
of~$\sK$:
\[
  H^*_{T^m}(\zk)\cong\Z[\sK].
\]
In equivariant cohomology, the projection $p\colon
ET^m\times_{T^m}\zk\to BT^m$ induces the quotient projection
$$
  p^*\colon\Z[v_1,\ldots,v_m]\to\Z[\mathcal K]=
  \Z[v_1,\ldots,v_m]/\mathcal I_\mathcal K.
$$
\end{corollary}
\begin{proof}
This follows from Theorem~\ref{zkhofib} and
Proposition~\ref{homsrs}.
\end{proof}

In view of this result, the Borel construction
$ET^m\times_{T^m}\zk$ is often called the
\emph{Davis--Januszkiewicz space}\label{ddjs} and denoted
by~$\djs(\sK)$. According to Theorem~\ref{zkhofib} it is modelled
(up to homotopy) on the polyhedral product~$(\C P^\infty)^{\sK}$.

\begin{example}\label{zk2wedge}
Let $\sK$ be the complex consisting of two disjoint points. Then
$\zk\cong S^3$ and $(\C P^\infty)^{\sK}=\C P^\infty\vee\C
P^\infty$ (a wedge of two copies of $\C P^\infty$). The Borel
construction $ET^2\times_{T^2}\zk$ can be identified with the
total space of the sphere bundle $S(\eta\times\eta)$ associated
with the product of two universal (Hopf) complex line bundles
$\eta$ over~$BT^1=\C P^\infty$. By Theorem~\ref{zkhofib}, the
space $ET^2\times_{T^2}\zk$ is homotopy equivalent to $\C
P^\infty\vee\C P^\infty$, and the bundle projection
$S(\eta\times\eta)\to\C P^\infty\times\C P^\infty$ induces the
quotient projection $\Z[v_1,v_2]\to\Z[v_1,v_2]/(v_1v_2)$ in
cohomology.
\end{example}

We have the following basic information about the homotopy groups
of $\zk$:

\begin{proposition}
\label{homgr}\
\begin{itemize}
\item[(a)] If $\sK$ is a simplicial complex on the vertex
set~$[m]$ (i.e. there are no ghost vertices), then the
moment-angle complex $\zk$ is 2-connected (i.e.
$\pi_1(\zk)=\pi_2(\zk)=0$), and
\[
  \pi_i(\zk)=\pi_i\bigl((\C P^\infty)^\sK\bigr)\quad
  \text{for }i\ge 3.
\]

\item[(b)] If $\mathcal K$ is a $q$-neighbourly simplicial
complex,
%(i.e. each $q$-element subset of $[m]$ spans a simplex),
then
$\pi_i(\zk)=0$ for $i<2q+1$. Furthermore, $\pi_{2q+1}(\zk)$ is a
free abelian group of rank equal to the number of $(q+1)$-element
missing faces of~$\mathcal K$.
\end{itemize}
\end{proposition}
\begin{proof}
We observe that $(\C P^\infty)^m$ is the Eilenberg--Mac Lane space
$K(\Z^m,2)$, and the 3-dimensional skeleton of $(\C P^\infty)^\sK$
coincides with the 3-skeleton of $(\C P^\infty)^m$. If $\mathcal
K$ is $q$-neighbourly, then the $(2q+1)$-skeletons of $(\C
P^\infty)^\sK$ and $(\C P^\infty)^m$ agree. Now both statements
follow from the exact homotopy sequence of the map $(\C
P^\infty)^\sK\to(\C P^\infty)^m$ with homotopy fibre~$\zk$.
\end{proof}

\begin{example}
Let $\sK=\sk_1\varDelta^3$ (a complete graph on 4 vertices). Then
$\sK$ is 2-neighbourly and has 4 missing triangles, so $\zk$ is
4-connected and $\pi_5(\zk)=\Z^4$.
\end{example}

\subsection*{Exercises}
\begin{exercise}
By analogy with Proposition~\ref{homsrs}, show that
\[
  H^*\bigl((\R
  P^\infty)^\sK,\Z_2\bigr)\cong\Z_2[v_1,\ldots,v_m]\big/
  \bigl(v_{i_1}\cdots
  v_{i_k}\colon\{i_1,\ldots,i_k\notin\sK\}\big),\qquad\deg v_i=1.
\]
\end{exercise}

\begin{exercise}
Consider a sequence of pointed odd-dimensional spheres
\[
  \mb S=(S^{2p_1-1},\ldots,S^{2p_m-1}).
\]
Show that there is an isomorphism of rings
\[
  H^*\bigl((\mb S)^\sK\bigr)\cong\Lambda[u_1,\ldots,u_m]\big/
  \bigl(u_{i_1}\!\cdots
  u_{i_k}\colon\{i_1,\ldots,i_k\notin\sK\}\bigr),\qquad\deg
  u_i=2p_i-1.
\]
This ring is known as the \emph{exterior face ring}\label{extefr}
of~$\sK$. In the case $p_1=\cdots=p_m=1$ we obtain $(\mb
S)^\sK=(S^1)^\sK$, which is a cell subcomplex in the torus~$T^m$.

For a more general statement describing the cohomology of the
polyhedral product $\mb X^\sK$, see Theorem~\ref{xkcohomo}.
\end{exercise}

\begin{exercise}[{\cite[Lemma~2.3.1]{de-su07}}]
Assume given commutative diagrams
\[
\diagram
  F'_i\rto\dto & E'_i\rto\dto & B_i\ddouble\\
  F_i\rto & E_i\rto & B_i
\enddiagram
\]
of cell complexes where the horizontal arrows are fibrations and
the vertical arrows are inclusions of cell subcomplexes, for
$i=1,\ldots,m$. Denote by $(\mb F,\mb F')$, $(\mb E,\mb E')$ and
$(\mb B,\mb B)$ the corresponding sequences of pairs. Show that
there is a fibration of polyhedral products
\[
  (\mb F,\mb F')^\sK\to(\mb E,\mb E')^\sK\to(\mb B,\mb B)^\sK
\]
where $(\mb B,\mb B)^\sK=B_1\times\cdots\times B_m$.
\end{exercise}

\begin{exercise}\label{plooppp}
Use the previous exercise, the path-loop fibration $\varOmega X\to
PX\to X$ and homotopy invariance of the polyhedral product
(Proposition~\ref{mafunc}) to show that the homotopy fibre of the
inclusion $\mb X^\sK\to\mb X^m$ is $(P\mb X,\varOmega\mb X)^\sK$
or, equivalently, $(\cone\varOmega\mb X,\varOmega\mb X)^\sK$. That
is, construct a homotopy fibration
\[
  (P\mb X,\varOmega\mb X)^\sK\to(\mb X,pt)^\sK\to(\mb X,\mb X)^\sK.
\]
When $X_i=\C P^\infty$ we obtain the homotopy fibration $\zk\to(\C
P^\infty)^\sK\to(\C P^\infty)^m$ of Theorem~\ref{zkhofib}.
\end{exercise}

\begin{exercise}
When $\sK$ is a pair of points and $\mb X=(X_1,X_2)$, show that
$(P\mb X,\varOmega\mb X)^\sK$ is homotopy equivalent to
$\varSigma\varOmega X_1\wedge\varOmega X_2$. Deduce \emph{Ganea's
Theorem}\label{GaneaT} identifying the homotopy fibre of the
inclusion $X_1\vee X_2\to X_1\times X_2$ with $\varSigma\varOmega
X_1\wedge\varOmega X_2$.
\end{exercise}

\begin{exercise}
In the setting of Example~\ref{zk2wedge}, consider the diagonal
circle $S^1_d\subset T^2$. Show that it acts freely on $\zk\cong
S^3$. Deduce that the Borel construction $ET^2\times_{T^2}\zk$ is
homotopy equivalent to $ES^1\times_{S^1}\C P^1$, where $\C
P^1=S^3/S^1_d$ with $S^1$-action given by
$t\cdot[z_0:z_1]=[z_0:tz_1]$. It follows that $ES^1\times_{S^1}\C
P^1\simeq\C P^\infty\vee\C P^\infty$. Show that
$ES^1\times_{S^1}\C P^1$ can be identified with the complex
projectivisation $\C P(\eta\oplus\bar\eta)$, where $\eta$ is the
tautological line bundle over $\C P^\infty$ and $\bar\eta$ is its
complex conjugate. What can be said about the complex
projectivisation $\C P(\eta\oplus\eta)$?
\end{exercise}

\section{Cell decomposition}\label{celld}
{\noindent \hangindent=27mm \hangafter=5 We consider the following
decomposition of the disc $\D$ into 3 cells: the point $1\in\D$ is
the 0-cell; the complement to $1$ in the boundary circle is the
1-cell, which we denote by~$T$; and the interior of $\D$ is the
2-cell, which we denote by~$D$. By taking product we obtain a
cellular decomposition of~$\D^m$ whose cells are parametrised by
pairs of subsets $J,I\subset [m]$ with $J\cap I=\varnothing$: the
set $J$ parametrises the $T$-cells in the product and $I$
parametrises the $D$-cells. We denote the cell of $\D^m$
corresponding to a pair $J,I$ by $\varkappa(J,I)$; it is a product
of $|J|$ cells of $T$-type and $|I|$ cells of $D$-type (the
positions in $[m]\setminus I\cup J$ are filled by 0-cells). Then
$\zk$ embeds as a cellular subcomplex in~$\D^m$; we have
$\varkappa(J,I)\subset\zk$ whenever $I\in\sK$.

} \noindent\raisebox{0.7\baselineskip}[0pt]
{%
\begin{picture}(20,20)
  \put(12,10){\circle{20}}
  \put(19,10){\circle*{1.3}}
  \put(20,9){$1$}
  \put(11,9){$D$}
  \put(16,16.5){$T$}
\end{picture}%
}

\vspace{-0.9\baselineskip}

Let $\sC^*(\zk)$ be the cellular cochains of~$\zk$. It has a basis
of cochains $\varkappa(J,I)^*$ dual to the corresponding cells. We
introduce the bigrading by setting
\[
  \bideg\varkappa(J,I)^*=(-|J|,2|I|+2|J|),
\]
so that $\bideg D=(0,2)$, $\bideg T=(-1,2)$ and $\bideg 1=(0,0)$.
Since the cellular differential preserves the second grading, the
complex $\sC^*(\zk)$ splits into the sum of its components with
fixed second degree:
\begin{equation}\label{bicel}
  \sC^*(\zk)=\bigoplus_{q=0}^m \sC^{*,2q}(\zk).
\end{equation}
The cohomology of the moment-angle complex therefore acquires an
additional grading, and we define the \emph{bigraded Betti
numbers} of~$\zk$ by
\begin{equation}\label{defbn}
  b^{-p,2q}(\zk)=\rank H^{-p,2q}(\zk),
  \quad\text{for } 1\le p,q\le m.
\end{equation}
The ordinary Betti numbers of $\zk$ therefore satisfy
\begin{equation}\label{ordb}
  b^k(\zk)=\sum_{-p+2q=k}b^{-p,2q}(\zk).
\end{equation}

The map of the moment-angle complexes $\zk\to\mathcal Z_{\mathcal
L}$ induced by a simplicial map $\sK\to\mathcal L$ (see
Proposition~\ref{functmac1}) is clearly a cellular map. We
therefore obtain

\begin{proposition}\label{functmac2}
The correspondence $\mathcal K\mapsto\zk$ gives rise to a functor
from the category of simplicial complexes and simplicial maps to
the category of cell complexes with torus actions and equivariant
maps. It also induces a natural transformation between the functor
of simplicial cochains of~$\mathcal K$ and the functor of cellular
cochains of~$\zk$.
\end{proposition}

The map $\zk\to\mathcal Z_{\mathcal L}$ induced by a simplicial
map $\sK\to\mathcal L$ preserves the cellular bigrading, so that
the bigraded cohomology groups are also functorial.

\section{Cohomology ring}\label{cohma}
The main result of this section, Theorem~\ref{zkcoh}, establishes
an isomorphism between the integral cohomology ring of the
moment-angle complex $\zk$ and the Tor-algebra of the simplicial
complex~$\sK$. This result was first proved in~\cite{bu-pa99} for
field coefficients using the Eilenberg--Moore spectral sequence of
the fibration $\zk\to(\C P^\infty)^\sK\to(\C P^\infty)^m$.
%(see Section~\ref{apemss} of the Appendix).
The proof given here is taken from~\cite{b-b-p04}
and~\cite{bu-pa04-2}; it establishes the isomorphism over the
integers and works with the cellular cochains.

One of the key steps in the proof is the construction of a
cellular approximation of the diagonal map
$\varDelta\colon\zk\to\zk\times\zk$ which is functorial with
respect to maps of moment-angle complexes induced by simplicial
maps. The resulting cellular cochain algebra is isomorphic to the
algebra $R^*(\sK)$ from Construction~\ref{astar} (obtained by
factorising the Koszul algebra of the face ring $\Z[\mathcal K]$
by an acyclic ideal); its cohomology is isomorphic to the
Tor-algebra of~$\sK$.

Another proof of Theorem~\ref{zkcoh} was given by
Franz~\cite{fran06}.

\subsection*{Algebraic model for cellular cochains}\label{algmodsc}
We recall the algebra $R^*(\sK)$ from Construction~\ref{astar}:
$$
  R^*(\sK)=\Lambda[u_1,\ldots,u_m]\otimes\Z[\sK]\bigr/(v_i^2=u_iv_i=0,\;
  1\le i\le m),
$$
with the bigrading and differential given by
\[
  \bideg u_i=(-1,2),\quad\bideg v_i=(0,2),\quad
  du_i=v_i,\quad dv_i=0.
\]
The algebra $R^*(\sK)$ has finite rank as an abelian group, with a
basis of monomials $u_Jv_I$ where $J\subset[m]$, $I\in\sK$ and
$J\cap I=\varnothing$.

Comparing the differential graded module structures in $R^*(\sK)$
and $\sC^*(\zk)$ we observe that they coincide, as described in
the following statement:

\begin{lemma}\label{cellcom}
The map
\[
  g\colon R^*(\sK) \to \sC^*(\zk),\quad
  u_J v_I \mapsto\mathcal \varkappa(J,I)^*,
\]
is an isomorphism of cochain complexes. Hence, there is an
additive isomorphism
\[
  H[R^*(\sK)]\cong H^*(\zk).
\]
\end{lemma}
\begin{proof}
Since $g$ arises from a bijective correspondence between bases of
$R^*(\sK)$ and $\sC^*(\zk)$, it is an isomorphism of bigraded
modules (or groups). It also clearly commutes with the
differentials:
\[
  \delta g(u_i)=\delta(T_i^\ast)=D_i^*=g(v_i)=g(du_i),\quad
  \delta g(v_i)=\delta(D_i^*)=0=g(dv_i),
\]
where $T_i$ denotes the cell $\varkappa(\{i\}, \varnothing)$, and
$D_i=\varkappa(\varnothing,\{i\})$.
\end{proof}

Having identified the algebra $R^*(\mathcal K)$ with the cellular
cochains of the moment-angle complex, we can give a topological
interpretation to the quasi-isomorphism of Lemma~\ref{iscoh}. To
do this we shall identify the Koszul algebra
$\Lambda[u_1,\ldots,u_m]\otimes\Z[\mathcal K]$ with the cellular
cochains of a space homotopy equivalent to~$\zk$.

The infinite-dimensional sphere $S^\infty$ is the direct limit
(union) of standardly embedded odd-dimensional spheres. Each odd
sphere $S^{2k+1}$ can be obtained from $S^{2k-1}$ by attaching two
cells of dimensions $2k$ and $2k+1$:
\[
  S^{2k+1}\cong (S^{2k-1}\cup_{f}D^{2k})\cup_g
  D^{2k+1}.
\]
Here the map $f\colon\partial D^{2k}\to S^{2k-1}$ is the identity
(and has degree~1), and the map $g\colon\partial
D^{2k+1}=S^{2k}\to D^{2k}$ is the projection of the standard
sphere onto its equatorial plane (and has degree~0). This implies
that $S^{\infty}$ is contractible and has a cell decomposition
with one cell in each dimension; the boundary of an even cell is
the closure of an odd cell, and the boundary of an odd cell is
zero. The 2-dimensional skeleton of this cell decomposition is the
disc~$D^2$ decomposed into three cells as described in
Section~\ref{celld}. The cellular cochain complex of $S^\infty$
can be identified with the Koszul algebra
$$
  \Lambda[u]\otimes\Z[v],\quad\deg u=1,\,\deg v=2,\quad
  du=v,\,dv=0.
$$

The functoriality of the polyhedral product
(Proposition~\ref{mafunc}~(a)) implies that there is the following
deformation retraction onto a cellular subcomplex:
$$
  \zk=(D^2,S^1)^\sK\hookrightarrow
  (S^\infty,S^1)^\sK
  \longrightarrow(D^2,S^1)^\sK.
$$
Furthermore, the cellular cochains of the polyhedral product
$(S^\infty,S^1)^\sK$ are identified with the Koszul algebra
$\Lambda[u_1,\ldots,u_m]\otimes\Z[\sK]$, in the same way as
$\sC^*(\zk)$ is identified with $R^*(\sK)$. Since
$\zk\subset(S^\infty,S^1)^\sK$ is a deformation retract, the
corresponding cellular cochain map
$$
  \Lambda[u_1,\ldots,u_m]\otimes\Z[\sK]=
  \sC^*\bigl((S^\infty,S^1)^\sK\bigr)\longrightarrow
  \sC^*(\zk)=R^*(\sK)
$$
induces an isomorphism in cohomology. The above map is a
homomorphism of algebras, so it is a quasi-isomorphism. In fact,
the cochain homotopy map constructed in the proof of
Lemma~\ref{iscoh} is nothing but the cellular cochain map induced
by the homotopy above.

\subsection*{Cellular diagonal approximation}
Here we establish the cohomology ring isomorphism in
Lemma~\ref{cellcom}. The difficulty of working with cellular
cochains is that they do not admit a functorial associative
multiplication. The diagonal map used in the definition of the
cohomology product is not cellular, and a cellular approximation
cannot be made functorial with respect to arbitrary cellular maps.
Here we construct a canonical cellular diagonal approximation
$\widetilde\varDelta\colon\zk\to\zk\times\zk$ which is functorial
with respect to maps of $\zk$ induced by simplicial maps of~$\sK$,
and show that the resulting product in the cellular cochains of
$\zk$ coincides with the product in~$R^*(\sK)$.

The product in the cohomology of a cell complex $X$ is defined as
follows (see~\cite{novi96}, \cite{hatc02}). Consider the composite
map of cellular cochain complexes
\begin{equation}\label{cemul}
\begin{CD}
  \sC^*(X)\otimes \sC^*(X) @>\times>> \sC^*(X\times X)
  @>\widetilde{\varDelta}^*>> \sC^*(X).
\end{CD}
\end{equation}
Here the map $\times$ sends a cellular cochain $c_1\otimes c_2\in
\sC^{q_1}(X)\otimes \sC^{q_2}(X)$ to the cochain $c_1\times c_2\in
\sC^{q_1+q_2}(X\times X)$, whose value on a cell $e_1\times e_2\in
X\times X$ is $(-1)^{q_1q_2}c_1(e_1)c_2(e_2)$. The map
$\widetilde{\varDelta}^*$ is induced by a cellular map
$\widetilde{\varDelta}$ (a cellular \emph{diagonal
approximation})\label{diagappro} homotopic to the diagonal
$\varDelta\colon X\to X\times X$. In cohomology, map~\eqref{cemul}
induces a multiplication $H^*(X)\otimes H^*(X)\to H^*(X)$ which
does not depend on a choice of cellular approximation and is
functorial. However, map~\eqref{cemul} itself is not functorial
because the choice of a cellular approximation is not canonical.

Nevertheless, in the case $X=\zk$ we can use the following
construction.

\begin{construction}[cellular approximation for $\varDelta\colon
\zk\to\zk\times\zk$] Consider the map $\widetilde{\varDelta}\colon
\D\to\D\times\D$ given in the polar coordinates $z=\rho
e^{i\varphi}\in\D$, $0\le\rho\le1$, $0\le\varphi<2\pi$, by the
formula
\begin{equation}\label{Ddiagmap}
  \rho e^{i\varphi}\mapsto\left\{
  \begin{array}{ll}
    (1-\rho+\rho e^{2i\varphi},1)&\text{ for }0\le\varphi\le\pi,\\
    (1,1-\rho+\rho e^{2i\varphi})&\text{ for }\pi\le\varphi<2\pi.
  \end{array}
  \right.
\end{equation}
It is easy to see that this is a cellular map homotopic to the
diagonal $\varDelta\colon\D\to\D\times\D$, and its restriction to
the boundary circle $\mathbb S$ is a diagonal approximation for
$\mathbb S$, as described by the following diagram:
\begin{equation}\label{DSdiag}
\begin{CD}
  \mathbb S  @>>> \D\\
  @V{\widetilde{\varDelta}}VV @VV{\widetilde{\varDelta}}V\\
  \mathbb S\times\mathbb S @>>> \mathbb D\times\mathbb D
\end{CD}
\end{equation}
(explicit formulae for the homotopies involved can be found in
Exercise~\ref{exdiag}). Taking the $m$-fold product we obtain a
cellular approximation
$\widetilde{\varDelta}\colon\D^m\to\D^m\times\D^m$. Applying
Proposition~\ref{mafunc}~(a) to the map of pairs
$\widetilde\varDelta\colon(\D,\mathbb S)\to(\D\times\D,\mathbb
S\times\mathbb S)$ and observing that $(\D\times\D,\mathbb
S\times\mathbb S)^\sK\cong\zk\times\zk$, we obtain that
$\widetilde{\varDelta}\colon\D^m\to\D^m\times\D^m$ restricts to a
cellular approximation of the diagonal map of~$\zk$, as described
in the following diagram:
$$
\begin{CD}
  \zk @>>> \D^m\\
  @V\widetilde{\varDelta}VV @VV\widetilde{\varDelta}V\\
  \zk\times\zk @>>> \D^m\times \D^m.
\end{CD}
$$
Finally, applying Proposition~\ref{mafunc}~(c) to
diagram~\eqref{DSdiag} we obtain that the approximation
$\widetilde{\varDelta}$ is functorial with respect to the maps of
moment-angle-complexes $\mathcal Z_{\mathcal K}\to\mathcal
Z_{\mathcal L}$ induced by simplicial maps $\mathcal K\to\mathcal
L$.
\end{construction}

\begin{lemma}\label{cellappr}
The cellular cochain algebra $\sC^*(\zk)$ with the product defined
via the diagonal approximation
$\widetilde\varDelta\colon\zk\to\zk\times\zk$ and
map~\eqref{cemul} is isomorphic to the algebra~$R^*(\mathcal K)$.
We therefore have an isomorphism of cohomology rings
\[
  H[R^*(\mathcal K)]\cong H^*(\zk).
\]
\end{lemma}
\begin{proof}
We first consider the case $\mathcal K=\varDelta^0$, i.e.
$\zk=\D$. The cellular cochain complex has basis of cochains $1\in
\sC^0(\D)$, $T^*\in \sC^1(\D)$ and $D^*\in \sC^2(\D)$ dual to the
cells introduced in Section~\ref{celld}. The multiplication
defined by~\eqref{cemul} in $\sC^*(\D)$ is trivial, so we have a
ring isomorphism
$$
  R^*(\varDelta^0)=\Lambda[u]\otimes\Z[v]\big/(v^2=uv=0)\longrightarrow \sC^*(\D).
$$

Taking an $m$-fold tensor product we obtain a ring isomorphism for
$\mathcal K=\varDelta^{m-1}$:
$$
  f\colon R^*(\varDelta^{m-1})=
  \Lambda[u_1,\ldots,u_m]\otimes\Z[v_1,\ldots,v_m]\big/
  (v_i^2=u_iv_i=0)
  \longrightarrow \sC^*(\D^m).
$$

Now for arbitrary $\mathcal K$ we have an inclusion
$\zk\subset\D^m=\mathcal Z_{\varDelta^{m-1}}$ of a cellular
subcomplex and the corresponding ring homomorphism $q\colon
\sC^*(\D^m)\to \sC^*(\zk)$. Consider the commutative diagram
$$
\begin{CD}
  R^*(\varDelta^{m-1}) @>f>> \sC^*(\D^m)\\
  @V{p}VV       @VVq V\\
  R^*(\mathcal K)        @>g>> \sC^*(\zk).
\end{CD}
$$
Here the maps $p$, $f$ and $q$ are ring homomorphisms, and $g$ is
an isomorphism of groups by Lemma~\ref{cellcom}. We claim that $g$
is also a ring isomorphism. Indeed, take $\alpha,\beta\in
R^*(\mathcal K)$. Since $p$ is onto, we have $\alpha=p(\alpha')$
and $\beta=p(\beta')$. Then
$$
  g(\alpha\beta)=gp(\alpha'\beta')=qf(\alpha'\beta')=
  qf(\alpha')qf(\beta')=gp(\alpha')gp(\beta')=g(\alpha)g(\beta),
$$
as claimed. Thus, $g$ is a ring isomorphism.
\end{proof}

\subsection*{Main result}
By combining the results of Lemmata~\ref{koscom}, \ref{iscoh}
and~\ref{cellappr} we obtain the main result of this section:

\begin{theorem}\label{zkcoh}
There are isomorphisms, functorial in~$\mathcal K$, of bigraded
algebras
\begin{align*}
  H^{*,*}(\zk)&\cong
  \Tor_{\Z[v_1,\ldots,v_m]}\bigl(\Z[\mathcal K],\Z\bigr)\\
  &\cong
  H\bigl(\Lambda[u_1,\ldots,u_m]\otimes\Z[\mathcal K],d\bigr),
\end{align*}
where the bigrading and the differential on the right hand side
are defined by
\[
  \bideg u_i=(-1,2),\quad\bideg v_i=(0,2),\qquad
  du_i=v_i,\quad dv_i=0.
\]
\end{theorem}

The algebraic Betti numbers~\eqref{bbnfr} of the face
ring~$\Z[\sK]$ therefore acquire a topological interpretation as
the bigraded Betti numbers~\eqref{defbn} of the moment-angle
complex~$\zk$.

Now we combine results of Propositions~\ref{frmap}, \ref{tamap}
and~\ref{functmac2}, Corollary~\ref{cohombk} and
Theorem~\ref{zkcoh} in the following statement describing the
functorial properties of the correspondence $\mathcal
K\mapsto\zk$.

\begin{proposition}\label{functors}
Consider the following functors:
\begin{itemize}
\item[(a)] $\mathcal Z$, the covariant functor $\mathcal
K\mapsto\zk$ from the category of finite simplicial complexes and
simplicial maps to the category of spaces with torus actions and
equivariant maps (the moment-angle complex functor);

\item[(b)] $\k[\cdot]$, the contravariant functor $\mathcal
K\mapsto\k[\mathcal K]$ from simplicial complexes to graded
$\k$-algebras (the face ring functor);

\item[(c)] $\mbox{\rm Tor-alg}$, the contravariant functor
$$
  \mathcal K\mapsto\Tor_{\k[v_1,\ldots,v_m]}\bigl(\k[\mathcal K],\k\bigr)
$$
from simplicial complexes to bigraded $\k$-algebras (the
$\Tor$-algebra functor; it is the composition of $\k[\cdot]$ and
$\Tor_{\k[v_1,\ldots,v_m]}(\cdot\,,\k)$);

\item[(d)] $H^*_T$, the contravariant functor
$X\mapsto H^*_T(X;\k)$ from spaces with torus actions to
$\k$-algebras (the equivariant cohomology functor);

\item[(e)] $H^*$, the contravariant functor $X\mapsto H^*(X;\k)$
from spaces to $\k$-algebras (the ordinary cohomology functor).
\end{itemize}
Then we have the following identities:
$$
  H^*_T\circ\mathcal Z=\k[\cdot],\qquad
  H^*\circ\mathcal Z=\mbox{\rm Tor-alg}.
$$
The second identity implies that for any simplicial map
$\phi\colon\sK\to\mathcal L$ the corresponding cohomology map
\[
  \phi^*_{\mathcal Z}\colon H^*(\mathcal Z_{\mathcal L})\to
  H^*(\mathcal Z_{\sK})
\]
coincides with the induced homomorphism of $\Tor$-algebras
$\phi_{\Tor}^*$ from Proposition~\ref{tamap}. In particular, the
map $\phi$ gives rise to a map
\[
  H^{-q,2p}(\mathcal Z_{\mathcal L})\to
  H^{-q,2p}(\mathcal Z_{\sK})
\]
of bigraded cohomology.
\end{proposition}

In the case of Cohen--Macaulay complexes $\sK$ (see
Section~\ref{cmr}) we have the following version of
Theorem~\ref{zkcoh}.

\begin{proposition}\label{reduc}
Let $\mathcal K$ be an $(n-1)$-dimensional Cohen--Macaulay
complex, and let $\mb t$ be an hsop in $\k[\mathcal K]$. Then we
have the following isomorphism of algebras:
$$
  H^*(\zk;\k)\cong\Tor_{\k[v_1,\ldots,v_m]/\mb t}
  \bigl(\k[\mathcal K]/\mb t,\k\bigr).
$$
\end{proposition}
\begin{proof}
This follows from Theorem~\ref{zkcoh} and Lemma~\ref{tortor}.
\end{proof}

Note that the algebra $\k[\mathcal K]/\mb t$ is finite-dimensional
as a $\k$-vector space, unlike $\k[\mathcal K]$. In some
circumstances this observation allows us to calculate the
cohomology of $\zk$ more effectively.

\subsection*{Description of the product in terms of full subcomplexes}
The Hochster formula (Theorem~\ref{hoch}) for the components of
the $\Tor$-algebra can be used to obtain an alternative
description of the product structure in~$H^*(\zk)$.

We recall from Section~\ref{hpfr} that the bigraded structure in
the Tor-algebra is refined to a multigrading\label{multigrator},
and the multigraded components of $\Tor$ can be calculated in
terms of the full subcomplexes of~$\sK$:
$$
  \Tor_{\Z[v_1,\ldots,v_m]}^{-i,2J}\bigl(\Z[\mathcal
  K],\Z\bigr)\cong
  \widetilde{H}^{|J|-i-1}(\mathcal K_J),
$$
where $J\subset[m]$, see Theorem~\ref{hochmd}. Furthermore, the
product in the $\Tor$-algebra defines a product in the direct sum
$\bigoplus_{p\ge0,\;J\subset[m]}\widetilde H^{p-1}(\sK_J)$ given
by~\eqref{fullsubcochain}.

The bigraded structure in the cellular cochain complex of $\zk$
defined in Section~\ref{celld} can be also refined to a
multigrading (a $\Z\oplus\Z^m$-grading):
\[
  \sC^{*}(\zk)=\bigoplus_{J\subset[m]}\sC^{*,\,2J}(\zk),
\]
where $\sC^{*,\,2J}(\zk)$ is the subcomplex spanned by the
cochains $\varkappa(J\backslash I,I)^*$ with $I\subset J$ and
$I\in\sK$. The bigraded cohomology groups are decomposed as
follows:
\[
  H^{-i,\,2j}(\zk)=
  \bigoplus_{J\subset[m]\colon|J|=j}H^{-i,\,2J}(\zk),
\]
where $H^{-i,\,2J}(\zk)=H^{-i}[\sC^{*,\,2J}(\zk)]$.

\begin{theorem}[{Baskakov~\cite{bask02}}]\label{zkhoch}
There are isomorphisms
\[
  \widetilde H^{p-1}(\sK_J)\stackrel{\cong}{\longrightarrow}
  H^{p-|J|,2J}(\zk),
\]
which are functorial with respect to simplicial maps and induce a
ring isomorphism
\[
  h\colon%\mathop{\bigoplus_{p\ge0}}\limits
  \sum_{J\subset[m]}\widetilde{H}^*(\sK_J)\stackrel{\cong}{\longrightarrow}
  H^{*}(\zk).
\]
\end{theorem}
\begin{proof}
The statement about the additive isomorphisms follows from
Theorems~\ref{hochmd} and~\ref{zkcoh}. In explicit terms, the
cohomology isomorphisms are induced by the cochain isomorphisms
given by
\[
\begin{aligned}
  C^{p-1}(\sK_J)&\longrightarrow \sC^{p-|J|,2J}(\zk),\\
  \alpha_L&\longmapsto\varepsilon(L,J)\varkappa(J\setminus L,L)^*
\end{aligned}
\]
similar to~\eqref{fmapr}, where $\alpha_L\in C^{p-1}(\sK_J)$ is
the cochain dual to a simplex $L\in\sK_J$.

The ring isomorphism follows from Proposition~\ref{proddirsum} and
Theorem~\ref{zkcoh}.
\end{proof}

We summarise the results above in the following description of the
cohomology groups and the product structure of $H^*(\zk)$ in terms
of full subcomplexes of~$\sK$:

\begin{theorem}\label{zkadd}
There are isomorphisms of groups
\[
  H^{-i,2j}(\zk)\cong\bigoplus_{J\subset[m]\colon|J|=j}
  \widetilde H^{j-i-1}(\sK_J),\qquad
  H^\ell(\zk)\cong\bigoplus_{J\subset[m]}
  \widetilde H^{\ell-|J|-1}(\sK_J).
\]
These isomorphisms sum up into a ring isomorphism
\[
  H^*(\zk)\cong\bigoplus_{J\subset[m]} \widetilde H^*(\sK_J),
\]
where the ring structure on the right hand side is given by the
canonical maps
\[
  H^{k-|I|-1}(\sK_{I})\otimes H^{\ell-|J|-1}(\sK_{J})\to
  H^{k+\ell-|I|-|J|-1}(\sK_{I\cup J})
\]
which are induced by simplicial maps $\sK_{I\cup
J}\to\sK_I\mathbin{*}\sK_J$ for $I\cap J=\varnothing$ and zero
otherwise.
\end{theorem}

It follows that the cohomology of $\zk$ may have arbitrary
torsion:

\begin{corollary}\label{arbitrarytorsion}
Any finite abelian group can appear as a summand in a cohomology
group of $H^*(\zk)$ for some~$\sK$.
\end{corollary}
\begin{proof}
It follows from the Theorem~\ref{zkadd} that $\widetilde H^*(\sK)$
is a direct summand in~$H^*(\zk)$ (with appropriate shifts in
dimension). Therefore, we can take $\sK$ whose simplicial
cohomology contains the appropriate torsion.
\end{proof}

\subsection*{Exercises}
\begin{exercise}
Let $\mathbb S$ be the standard unit circle decomposed into two
cells, where the 0-cell is the unit. The map
\[
  \widetilde\varDelta\colon\mathbb S\to\mathbb S\times\mathbb S,
  \qquad
  e^{i\varphi}\mapsto\left\{
  \begin{array}{ll}
    (e^{2i\varphi},1)&\text{ for }0\le\varphi\le\pi,\\
    (1,e^{2i\varphi})&\text{ for }\pi\le\varphi<2\pi
  \end{array}
  \right.
\]
is a cellular diagonal approximation. It is obtained by
restricting map~\eqref{Ddiagmap} to the boundary circle
($\rho=1$). A homotopy $F_t$ between the diagonal
$\varDelta\colon\mathbb S\to\mathbb S\times\mathbb S$ ($t=0$) and
its cellular approximation $\widetilde\varDelta$ ($t=1$) is given
by
\[
F_t\colon\mathbb S\to\mathbb S\times\mathbb S, \qquad
e^{i\varphi}\mapsto\left\{
  \begin{array}{ll}
    \bigl(e^{i(1+t)\varphi},e^{i(1-t)\varphi}\bigr)
      &\text{ for }0\le\varphi\le\pi,\\[2pt]
    \bigl(e^{i(1-t)\varphi+2\pi it},e^{i(1+t)\varphi-2\pi it}\bigr)
      &\text{ for }\pi\le\varphi<2\pi.
  \end{array}
  \right.
\]
\end{exercise}

\begin{exercise}\label{exdiag}
Show that the formula
\[
\rho e^{i\varphi}\mapsto\left\{
  \begin{array}{ll}
    \bigl((1-\rho)t+\rho e^{i(1+t)\varphi},
      (1-\rho)t+\rho e^{i(1-t)\varphi}\bigr)
      &\text{ for }0\le\varphi\le\pi,\\[2pt]
      %\qquad\\[2pt]\hfill\text{ for }0\le\varphi\le\pi,\\[4pt]
    \bigl((1-\rho)t+\rho e^{i(1-t)\varphi+2\pi it},
    (1-\rho)t+\rho e^{i(1+t)\varphi-2\pi
    it}\bigr)
    %\qquad\\[2pt]\hfill\text{ for }\pi\le\varphi<2\pi
    &\text{ for }\pi\le\varphi<2\pi.
  \end{array}
  \right.
\]
defines a homotopy $G_t\colon\mathbb D\to\mathbb D\times\mathbb D$
between the diagonal $\varDelta\colon\D\to\D\times\D$ ($t=0$) and
its approximation $\widetilde\varDelta$~\eqref{DSdiag} ($t=1$).
For $\rho=1$ the homotopy $G_t$ restricts to the homotopy $F_t$ of
the previous exercise.
\end{exercise}

\section{Bigraded Betti numbers}\label{bettizk}
Here we describe the main properties of the bigraded Betti
numbers~\eqref{defbn} of moment-angle complexes and give some
examples of explicit calculations.

\begin{lemma}\label{bbgen}
Let $\sK$ be a simplicial complex of dimension $n-1$ with $f_0=m$
vertices and $f_1$ edges, so that $\dim\zk=m+n$. We have
\begin{itemize}
\item[(a)] $b^{0,0}(\zk)=b^0(\zk)=1$ and $b^{0,\,2q}(\zk)=0$ for
$q\ne0$;

\item[(b)] $b^{-p,\,2q}=0$ for $q>m$ or $p>q$;

\item[(c)] $b^1(\zk)=b^2(\zk)=0$;

\item[(d)] $b^3(\zk)=b^{-1,4}(\zk)=\bin m2-f_1$;

\item[(e)] $b^{-p,\,2q}(\zk)=0$ for $p\ge q>0$ or $q-p>n$;

\item[(f)] $b^{m+n}(\zk)=b^{-(m-n),\,2m}(\zk)=\rank\widetilde H^{n-1}(\mathcal K)$.
\end{itemize}
\end{lemma}
\begin{proof}
We consider the algebra $R^{*}(\mathcal K)$ whose cohomology is
$H^*(\zk)$. Recall that $R^*(\sK)$ has additive basis of monomials
$u_Jv_I$ with $I\in\sK$ and $I\cap J=\varnothing$. Since $\bideg
v_i=(0,2)$, $\bideg u_j=(-1,2)$, the bigraded component
$R^{-p,2q}(\sK)$ has basis of monomials $u_J v_I$ with $|I|=q-p$
and $|J|=p$. In particular, $R^{-p,2q}(\sK)=0$ for $q>m$ or $p>q$,
which implies~(b). To prove~(a) we observe that $R^{0,0}(\sK)=\k$
and each $v_I\in R^{0,2q}(\sK)$ with $q>0$ is a coboundary, hence,
$H^{0,2q}(\zk)=0$ for $q>0$.

Now we prove (e). Let $u_J v_I\in R^{-p,2q}(\sK)$; then $|I|=q-p$
and $I\in \sK$. Since a simplex of $\sK$ has at most $n$ vertices,
$R^{-p,2q}(\sK)=0$ for $q-p>n$. We have $b^{-p,2q}(\zk)=0$ for
$p>q$ by~(b) so we need only to check that $b^{-q,2q}(\zk)=0$ for
$q>0$. The group $R^{-q,2q}(\sK)$ has basis of monomials $u_J$
with $|J|=q$. Since $d(u_i)=v_i$, there are no nonzero cocycles in
$R^{-q,2q}(\sK)$ for $q>0$, hence, $H^{-q,2q}(\zk)=0$.

Statement~(c) follows from (e) and~\eqref{ordb}.

We also have $H^{3}(\zk)=H^{-1,4}(\zk)$, by~(e). There is a basis
in $R^{-1,4}(\sK)$ consisting of monomials $u_jv_i$ with $i\ne j$.
We have $d(u_jv_i)=v_iv_j$ and $d(u_iu_j)=u_jv_i-u_iv_j$. Hence,
$u_jv_i$ is a cocycle if and only if $\{i,j\}\notin\sK$; in this
case the two cocycles $u_jv_i$ and $u_iv_j$ represent the same
cohomology class. This proves~(d).

It remains to prove~(f). The total degree of a monomial $u_J
v_I\in R^*(\sK)$ is $2|I|+|J|$, and there are constraints
$|I|+|J|\le m$ and $|I|\le n$. Therefore, the maximum of the total
degree is achieved for $|I|=n$ and $|J|=m-n$. This proves the
first identity of~(f), and the second follows from
Theorem~\ref{zkadd}.
\end{proof}

Lemma~\ref{bbgen} shows that nonzero bigraded Betti numbers
$b^{r,2q}(\zk)$ with $r\ne0$ appear only in the strip bounded by
the lines $r=-1$, \ $q=m$, \ $r+q=1$ and $r+q=n$ in the second
quadrant, see Fig.~\ref{bnloc}~(a).
\begin{figure}[h] \begin{picture}(115,60)
\multiput(45,10)(0,5){10}{\line(-1,0){47}}
  \multiput(45,10)(-5,0){10}{\line(0,1){47}}
  \put(46,11){\small 0}
  \put(46,16){\small 2}
  \put(46,21){\small 4}
  \put(46,31){\vdots}
  \put(46,51){\small $2m$}
  \put(42,7){\footnotesize 0}
  \put(35.5,7){\footnotesize $-1$}
  \put(30,7){\small $\ldots$}
  \put(22.5,10){\line(0,-1){2}}
  \put(15,6){\footnotesize $-(m-n)$}
  \put(-1,7){\footnotesize $-m$}
  \put(42,11){\small $1$}
  \multiput(36,21)(-5,5){7}{\Large $*$}
  \multiput(36,26)(-5,5){6}{\Large $*$}
  \multiput(36,31)(-5,5){5}{\Large $*$}
  \multiput(36,36)(-5,5){4}{\Large $*$}
  \put(7,0){(a)\ arbitrary $\sK$}
  \multiput(110,10)(0,5){10}{\line(-1,0){47}}
  \multiput(110,10)(-5,0){10}{\line(0,1){47}}
  \put(111,11){\small 0}
  \put(111,16){\small 2}
  \put(111,21){\small 4}
  \put(111,31){\vdots}
  \put(111,51){\small $2m$}
  \put(107,7){\footnotesize 0}
  \put(100.5,7){\footnotesize $-1$}
  \put(95,7){\small $\ldots$}
  \put(87.5,10){\line(0,-1){2}}
  \put(80,6){\footnotesize $-(m-n)$}
  \put(64,7){\footnotesize $-m$}
  \put(107,11){\small $1$}
  \multiput(101,21)(-5,5){3}{\Large $*$}
  \multiput(101,26)(-5,5){3}{\Large $*$}
  \multiput(101,31)(-5,5){3}{\Large $*$}
  \put(87,51){\small $1$}
  \put(77,0){(b)\ $|\sK|\cong S^{n-1}$}
  \end{picture}
  \caption{Possible locations on nonzero $b^{-p,2q}(\zk)$
  (marked by $*$).}
  \label{bnloc}
\end{figure}

The next result allows us to express the numbers of faces of~$\sK$
(i.e. its $f$- and $h$-vectors) via the bigraded Betti numbers. We
consider the Euler characteristics of the complexes
$\sC^{*,2q}(\zk)$ (see~\eqref{bicel}),
\begin{equation}\label{chip}
  \chi_q(\zk)=\sum_{p=0}^m(-1)^p\rank \sC^{-p,2q}(\zk)
  =\sum_{p=0}^m(-1)^pb^{-p,2q}(\zk)
\end{equation}
and define the generating series
$$
  \chi(\zk;t)=\sum_{q=0}^m\chi_q(\zk)t^{2q}.
$$

\begin{theorem}
\label{gpz} The following identity holds for an
$(n-1)$-dimensional simplicial complex $\sK$ with $m$ vertices:
\[
  \chi(\zk;t)=(1-t^2)^{m-n}(h_0+h_1t^2+\cdots+h_nt^{2n}).
\]
Here $(h_0,h_1,\ldots,h_n)$ is the $h$-vector of~$\sK$.
\end{theorem}
\begin{proof}
The bigraded component $\sC^{-p,2q}(\zk)$ has basis of cellular
cochains $\varkappa(J,I)^*$ with $I\in\sK$, \ $|I|=q-p$ and
$|J|=p$. Therefore, $\rank
\sC^{-p,2q}(\zk)=f_{q-p-1}\bin{m-q+p}p$, where
$(f_0,f_1,\ldots,f_{n-1})$ is the $f$-vector of $\sK$ and
$f_{-1}=1$. By substituting this into~\eqref{chip} we obtain
\[
  \chi_q(\zk)=\sum_{j=0}^m(-1)^{q-j}f_{j-1}\bin{m-j}{q-j}.
\]
Then
\begin{multline}
\label{chidir}
  \chi(\zk;t)=
  %\sum_{q=0}^m\chi_q(\zk)t^{2q}=
  \sum_{q=0}^m\sum_{j=0}^mt^{2j}t^{2(q-j)}(-1)^{q-j}f_{j-1}
  \bin{m-j}{q-j}\\
  =\sum_{j=0}^mf_{j-1}t^{2j}(1-t^2)^{m-j}=
  (1-t^2)^m\sum_{j=0}^nf_{j-1}(t^{-2}-1)^{-j}.
\end{multline}
Set $h(t)=h_0+h_1t+\cdots+h_nt^n$. Then it follows
from~\eqref{hvectors} that
$$
  t^nh(t^{-1})=(t-1)^n\sum_{j=0}^nf_{j-1}(t-1)^{-j}.
$$
By substituting $t^{-2}$ for $t$ in the identity above we finally
rewrite~\eqref{chidir} as
$$
  \frac{\chi(\zk;t)}{(1-t^2)^m}=\frac{t^{-2n}h(t^2)}{(t^{-2}-1)^n}=
  \frac{h(t^2)}{(1-t^2)^n},
$$
which is equivalent to the required identity.
\end{proof}

\begin{corollary}\label{chizk}
If $\sK\ne\varDelta^{m-1}$, then the Euler characteristic of $\zk$
is zero.
\end{corollary}
\begin{proof}
We have
$$
  \chi(\zk)=\sum_{p,q=0}^m(-1)^{-p+2q}b^{-p,2q}(\zk)=
  \sum_{q=0}^m\chi_q(\zk)=\chi(\zk;1)=0
$$
by Theorem~\ref{gpz} (note that $\sK\ne\varDelta^{m-1}$ implies
that $m>n$).
\end{proof}

We proceed by describing the properties of bigraded Betti numbers
for particular classes of simplicial complexes.

\begin{definition}\label{defpm}
A finite simplicial complex $\sK$ is called a $d$-dimensional
\emph{pseudomanifold} if the following three conditions are
satisfied:
\begin{itemize}
\item[(a)] all maximal simplices of $\sK$ have dimension $d$ (i.e. $\sK$ is pure $d$-dimensional);

\item[(b)] each $(d-1)$-simplex of $\sK$ is the face of exactly
two $d$-simplices of~$\sK$.

\item[(c)] if $I$ and $I'$ are $d$-simplices of $\sK$, then there
is a sequence $I=I_1,I_2,\ldots,I_k=I'$ of $d$-simplices of $\sK$
such that $I_j$ and $I_{j+1}$ have a common $(d-1)$-face for $1\le
i\le k-1$.
\end{itemize}
\end{definition}

If $\sK$ is a $d$-dimensional pseudomanifold, then either
$H_d(\sK)\cong\Z$ or~$0$ (an exercise). In the former case the
pseudomanifold $\sK$ is called \emph{orientable}.

\begin{lemma}\label{fc}
Let $\sK$ be an orientable pseudomanifold of dimension $n-1$ with
$m$ vertices. Then
\[
  H^{m+n}(\zk)=\widetilde H^{n-1}(\sK)\cong\Z.
\]
Under the isomorphism $H^*(\zk)\cong H[R^*(\sK)]$, the group above
is generated by the class of any monomial $u_Jv_I\in R^*(\sK)$ of
bidegree $(-(m-n),2m)$ such that $I\in\sK$ and $J=[m]\setminus I$.
\end{lemma}
\begin{proof}
The isomorphism of groups follows from Theorem~\ref{zkadd} and the
fact that $\sK$ is orientable. We have
$H^{m+n}(\zk)=H^{-(m-n),2m}(\zk)$. The group $R^{-(m-n),2m}(\sK)$
has basis of monomials $u_Jv_I$ with $I\in\sK$, \ $|I|=n$ and
$J=[m]\setminus I$. Each of these monomials is a cocycle. Let
$I,I'$ be two $(n-1)$-simplices of $\sK$ having a common
$(n-2)$-face. Consider the corresponding cocycles $u_Jv_I$ and
$u_{J'}v_{I'}$ (where $J=[m]\setminus I$, \ $J'=[m]\setminus I'$):
\begin{align*}
  u_Jv_I
  &=u_{j_1}u_{j_2}\cdots u_{j_{m-n}}v_{i_1}\cdots v_{i_{n-1}}v_{i_n},\\
  u_{J'}v_{I'}
  &=u_{i_n}u_{j_2}\cdots u_{j_{m-n}}v_{i_1}\cdots v_{i_{n-1}}v_{j_1}.
\end{align*}
Since $\sK$ is a pseudomanifold, the $(n-2)$-face
$\{i_1,\ldots,i_{n-1}\}$ is contained in exactly two
$(n-1)$-faces, namely $I=\{i_1,\ldots,i_{n-1},i_n\}$ and
$I'=\{i_1,\ldots,i_{n-1},j_1\}$. Therefore we have the following
identity in $R^*(\sK)$
\begin{multline*}
  d(u_{i_n}u_{j_1}u_{j_2}\cdots u_{j_{m-n}}
  v_{i_1}\cdots v_{i_{n-1}})=\\
  =u_{j_1}u_{j_2}\cdots u_{j_{m-n}}v_{i_1}\cdots v_{i_{n-1}}v_{i_n}-
  u_{i_n}u_{j_2}\cdots u_{j_{m-n}}v_{i_1}\cdots v_{i_{n-1}}v_{j_1}
\end{multline*}
Hence, $[u_Jv_I]=[u_{J'}v_{I'}]$ (as cohomology classes).
Property~(c) from the definition of a pseudomanifold implies that
all monomials $u_Jv_I\in R^{-(m-n),2m}(\sK)$ represent the same
cohomology class up to sign. The isomorphism~\eqref{fmapr} takes
$u_Jv_I$ to $\pm\alpha_I\in C^{n-1}(\sK)$, which represents a
generator of $\widetilde H^{n-1}(\sK)\cong\Z$ (see
Exercise~\ref{fundpseudo}).
\end{proof}

\begin{remark}
If $\sK$ is a non-orientable pseudomanifold, then the same
argument shows that any monomial $u_Jv_I\in R^*(\sK)$ as above
represents the generator of $H^{m+n}(\zk)=H^{n-1}(\sK)\cong\Z_2$.
\end{remark}

\begin{proposition}\label{bpd}
Let $\sK$ be a triangulated sphere of dimension $n-1$. Then
Poincar\'e duality for the moment-angle manifold $\zk$ respects
the bigrading in cohomology. In particular,
\[
  b^{-p,\,2q}(\zk)=b^{-(m-n)+p,\,2(m-q)}(\zk)\quad\text{for } 0\le p\le m-n,\,
  0\le q\le m.
\]
\end{proposition}
\begin{proof}
The Poincar\'e duality maps (see Definition~\ref{pdalgebra}) are
defined via the cohomology multiplication in~$H^*(\zk)$, which
respects the bigrading. We have $\dim\zk=m+n$, and
\[
  H^{m+n}(\zk)=
  \Tor_{\Z[v_1,\ldots,v_m]}^{-(m-n),\,2m}\bigl(\Z[\sK],\Z\bigr)\cong\Z,
\]
by Lemma~\ref{fc}. This implies the required identity for the
Betti numbers.
\end{proof}

\begin{corollary}\label{bbss}
Let $\sK$ be a triangulated $(n-1)$-sphere and $\zk$ the
corresponding moment-angle manifold, $\dim\zk=m+n$. Then
\begin{itemize}
\item[(a)] $b^{-p,2q}(\zk)=0$ for $p\ge m-n$,
with the only exception $b^{-(m-n),2m}=1$;

\item[(b)] $b^{-p,2q}(\zk)=0$ for $q-p\ge n$,
with the only exception $b^{-(m-n),2m}=1$.
\end{itemize}
\end{corollary}

It follows that if $|\sK|\cong S^{n-1}$, then nonzero bigraded
Betti numbers $b^{r,2q}(\zk)$, except $b^{0,0}(\zk)$ and
$b^{-(m-n),2m}(\zk)$, appear only in the strip bounded by the
lines $r=-(m-n-1)$, \ $r=-1$, \ $r+q=1$ and $r+q=n-1$ in the
second quadrant, see Fig.~\ref{bnloc}~(b).

A space $X$ is called a \emph{Poincar\'e duality space} (over
$\k$)\label{defPds} if $H^*(X;\k)$ is a Poincar\'e algebra (see
Definition~\ref{pdalgebra}). We have the following
characterisation of moment-angle complexes with Poincar\'e
duality, extending the result of Corollary~\ref{bpd}.

\begin{theorem}\label{pdcom}
$\zk$ is a Poincar\'e duality space over a field $\k$ if and only
$\sK$ is a Gorenstein complex over~$\k$.
\end{theorem}
\begin{proof}
Assume that $\sK$ is a Gorenstein complex. Consider the algebra
$T$ defined in Theorem~\ref{av-go}, i.e. $T=\bigoplus_{i=0}^d
T^i$, where $T^i=\Tor^{-i}_{\k[m]}(\k[\sK],\k)$ and $d=\max\{
j\colon \Tor^{-j}_{\k[m]}(\k[\sK],\k)\ne0\}$. Since $T$ is
Poincar\'e algebra, $\k\cong T^0\cong\Hom_\k(T^d,T^d)$, which
implies that $T^d\cong\k$. Since $T$ has a bigrading, we obtain
$T^d=T^{d,2q}$ for some~$q\ge0$. Since the multiplication in $T$
respects the bigrading, the isomorphisms
$T^i\stackrel\cong\to\Hom_\k(T^{d-i},T^d)$ from the definition of
a Poincar\'e algebra split into isomorphisms
\[
  T^{i,\,2j}\stackrel\cong\longrightarrow
  \Hom_\k(T^{d-i,\,2(q-j)},T^{d,\,2q}).
\]
Let $H^k=H^k(\zk;\k)$ and $H=\bigoplus_{k=0}^rH^k$; then
$H^k=\bigoplus_{i+2j=k}T^{i,\,2j}$ and $r=d+2q$. Therefore, we
have the isomorphisms
\[
  H^k=\!\!\bigoplus_{i+2j=k}T^{i,2j}\stackrel\cong\longrightarrow
  \!\!\bigoplus_{i+2j=k}\Hom_\k(T^{d-i,2(q-j)},T^{d,2q})=
  \Hom_\k(H^{r-k},H^r),
\]
which imply that $H$ is a Poincar\'e algebra.

Now assume that $H=\bigoplus_{k=0}^rH^k$ is a Poincar\'e algebra.
Then
\[
  \k\cong H^r=\Tor^{-d,\,2q}_{\k[v_1,\ldots,v_m]}
  \bigl(\k[\sK],\k\bigr)=T^{d,\,2q}
\]
for some $d,q\ge0$. Since the multiplication in the cohomology of
$\zk$ respects the bigrading, the isomorphisms
$H^k\stackrel\cong\to \Hom_\k(H^{r-k},H^r)$ split into
isomorphisms
\[
  H^{-i,\,2j}=T^{i,\,2j}\stackrel\cong\longrightarrow
  \Hom_\k(T^{d-i,\,2(q-j)},T^{d,\,2q}),
\]
which in their turn define the isomorphisms
\[
  T^i=\bigoplus_j T^{i,\,2j}\stackrel\cong\longrightarrow
  \bigoplus_j\Hom_\k(T^{d-i,\,2(q-j)},T^d)=\Hom_\k(T^{d-i},T^d).
\]
Thus, $T$ is a Poincar\'e algebra.
\end{proof}

\begin{remark}
We do not assume that $r=\max\{k\colon H^k(\zk;\k)\ne0\}$ is equal
to $\dim\zk=m+n$ in Theorem~\ref{pdcom}. It follows from the proof
above that $\zk$ is a Poincar\'e duality space with $r=\dim\zk$ if
and only if $\sK$ is a Gorenstein* complex.
\end{remark}

Here are some explicit examples of calculations of $H^*(\zk)$
using Theorem~\ref{zkcoh}.

\begin{example}
Let $\sK=\partial\varDelta^{m-1}$. Then
\[
  \Z[\sK]=\Z[v_1,\ldots,v_m]/(v_1\cdots v_m).
\]
The cocycle $u_1v_2v_3\cdots
v_m\in\Lambda[u_1,\ldots,u_m]\otimes\Z[\sK]$ of bidegree $(-1,2m)$
represents a generator of the top degree cohomology group of
$\zk\cong S^{2m-1}$.
\end{example}

\begin{example}\label{b5-gon}
Let $\sK$ be the boundary of 5-gon. We have $\dim\zk=7$. We
enumerate the vertices of $\sK$ clockwise. The face ring of $\sK$
is given in Example~\ref{frpex}.3. The group $H^3(\zk)$ has 5
generators corresponding to the diagonals of the 5-gon; these
generators are represented by the cocycles
$u_iv_{i+2}\in\Z[\sK]\otimes\Lambda[u_1,\ldots,u_5]$, \ $1\le
i\le5$ (the summation of indices is modulo~5). A direct
calculation shows that $H^4(\zk)$ also has 5 generators,
represented by the cocycles $u_ju_{j+1}v_{j+3}$, \ $1\le j\le5$.
The Betti vector of $\zk$ is therefore given by
\[
  \bigl(b^0(\zk),b^1(\zk),\ldots,b^7(\zk)\bigr)=
  (1,0,0,5,5,0,0,1).
\]
By Lemma~\ref{fc}, the product of cocycles $u_iv_{i+2}$ and
$u_ju_{j+1}v_{j+3}$ represents a generator of $H^7(\zk)$ if and
only if all the indices $i,i+2,j,j+1,j+3$ are different. Hence,
for each cohomology class $[u_iv_{i+2}]\in H^3(\zk)$ there is a
unique class $[u_ju_{j+1}v_{j+3}]\in H^4(\zk)$ such that the
product $[u_iv_{i+2}]\cdot[u_ju_{j+1}v_{j+3}]$ is nonzero. These
calculations are summarised by the cohomology ring isomorphism
\[
  H^*(\zk)\cong H^*\bigl((S^3\times S^4)^{\cs5}\bigr),
\]
in accordance with Exercise~\ref{exzk5gon}.
\end{example}

%Напомним, что \emph{когомологической длиной} пространства~$X$
%называется максимальное число $l=l(X)$, для которого существует
%набор классов когомологий $\{\alpha_i\in H^{>0}(X),\,1\le i\le
%l\}$ с ненулевым произведением $\alpha_1\cdots\alpha_l$.

\begin{example}\label{zkgon}
Now we calculate the Betti numbers and the cohomology product for
$\zk$ in the case when $\sK$ is a boundary of an $m$-gon with
$m\ge4$.

It follows from Corollary~\ref{bbss} that the only nonzero
bigraded Betti numbers of $\zk$ are  $b^{-(p-1),2p}$ for $2\le
p\le m-2$ and $b^{0,0}(\zk)=b^{-(m-2),2m}(\zk)=1$. The ordinary
Betti numbers are therefore given by
\[
  b^0(\zk)=b^{m+2}(\zk)=1,\quad
  b^k(\zk)=b^{-(k-2),2(k-1)}(\zk)\quad\text{for }3\le k\le m-1.
\]

To calculate $b^{-(k-2),2(k-1)}(\zk)$ for $3\le k\le m-1$ we use
the algebra $R^*(\sK)$ of Construction~\ref{astar}. We have
\begin{multline}\label{bmgon}
  b^{-(k-2),2(k-1)}(\zk)=\rank H^{-(k-2),2(k-1)}(R^*,d)=
  \\\rank\ker\bigl[d\colon R^{-(k-2),2(k-1)}\to
  R^{-(k-3),2(k-1)}\bigr]-\rank d R^{-(k-1),2(k-1)}.
\end{multline}
Since $H^{-(k-1),2(k-1)}(\zk)=0$ for $k>1$, the differential $d$
from $R^{-(k-1),2(k-1)}$ is monomorphic, and
\[
  \rank d R^{-(k-1),2(k-1)}=\rank R^{-(k-1),2(k-1)}.
\]
Similarly, since $H^{-(k-3),2(k-1)}(\zk)=0$ for $k\le m$, the
differential from $R^{-(k-2),2(k-1)}$ is epimorphic, and
\begin{multline*}
  \rank\ker\bigl[d\colon R^{-(k-2),2(k-1)}\to
  R^{-(k-3),2(k-1)}\bigr]\\=\rank R^{-(k-2),2(k-1)}-\rank
  R^{-(k-3),2(k-1)}.
\end{multline*}
Substituting the above two expressions into~\eqref{bmgon} and
using formula~\eqref{rankrpq} for the dimensions of $R^{-p,2q}$,
we calculate
\begin{align*}
  b^k(\zk)&=b^{-(k-2),2(k-1)}(\zk)\\
  &=\rank R^{-(k-2),2(k-1)}-\rank
  R^{-(k-3),2(k-1)}-\rank R^{-(k-1),2(k-1)}\\
  &=m\bin{m-1}{k-2}-m\bin{m-2}{k-3}-\bin m{k-1}\\
  &=(k-2)\bin{m-2}{k-1}+(m-k)\bin{m-2}{m-k+1}
  \quad\text{for }3\le k\le m-1.
\end{align*}

Note that the cohomology of $\zk$ does not have torsion in this
example (this follows from Theorem~\ref{zkadd}).

The product of any three classes of positive degree in $H^*(\zk)$
is zero (i.e. the \emph{cohomology product
length}\label{cohoprole} of $\zk$ is~2). Indeed, if $\alpha_i\in
H^{-(p_i-1),2p_i}(\zk)$ for $i=1,2,3$, then
\[
  \alpha_1\alpha_2\alpha_3\in
  H^{-(p_1+p_2+p_3-3),2(p_1+p_2+p_3)}(\zk),
\]
which is zero by Lemma~\ref{bbgen}~(e) (note that $n=2$ in this
example). Hence all nontrivial products in $H^*(\zk)$ arise from
Poincar\'e duality.

The above calculations of the Betti numbers together with the
observations about the cohomology product can be summarised by
saying that $\zk$ cohomologically looks like a connected sum of
sphere products, namely,
\begin{equation}\label{zpcs}
  H^*(\zk)\cong H^*\Bigl(\mathop{\#}_{k=3}^{m-1}
  \bigl(S^k\times S^{m+2-k}\bigr)^{\#(k-2)\bin{m-2}{k-1}}\Bigr),
\end{equation}
as rings.
\end{example}

According to a result of McGavran~\cite{mcga79}, the cohomology
isomorphism of~\eqref{zpcs} is induced by a homeomorphism of
manifolds. This is true even in a more general situation, when
$\sK$ is the boundary of a stacked polytope (see
Definition~\ref{stacked}):

\begin{theorem}[{\cite{mcga79}, \cite[Theorem~6.3]{bo-me06}}]\label{zpstacked}
Let $\sK$ be the boundary of a stacked polytope of dimension $n$
with $m>n+1$ vertices. Then the corresponding moment-angle
manifold is homeomorphic to a connected sum of sphere products,
\[
  \zk\cong\mathop{\#}_{k=3}^{m-n+1}
  \bigl(S^k\times S^{m+n-k}\bigr)^{\#(k-2)\bin{m-n}{k-1}}.
\]
\end{theorem}

For $n=2$ we obtain a homeomorphism of manifolds underlying the
isomorphism~\eqref{zpcs} (note that any $2$-polytope is stacked).

The bigraded Betti numbers for stacked polytopes can also be
calculated:
%effectively:

\begin{theorem}[{\cite{te-hi97}, \cite{ch-ki10}, \cite{limo12}}]
Let $\sK$ be as in Theorem~\ref{zpstacked} and $n\ge3$. Then the
nonzero bigraded Betti numbers of $\zk$ are given by
\begin{align*}
b^{0,0}(\zk)&=b^{-(m-n),2m}(\zk)=1,\\
b^{-i,2(i+1)}(\zk)&=i\bin{m-n}{i+1}\quad\text{for }1\le i\le m-n-1,\\
b^{-i,2(i+n-1)}(\zk)&=(m-n-i)\bin{m-n}{m-n+1-i}\quad\text{for
}1\le i\le m-n-1.
\end{align*}
\end{theorem}

\subsection*{Exercises}

\begin{exercise}
A triangulated manifold $\sK$ is a pseudomanifold.
\end{exercise}

\begin{exercise}
If $\sK$ is a $d$-dimensional pseudomanifold, then either
$H_d(\sK)\cong\Z$ or~$H_d(\sK)\cong0$. What happens for homology
and cohomology with coefficients in a field~$\k$?
\end{exercise}

\begin{exercise}
If $\sK$ is an orientable $d$-dimensional pseudomanifold (i.e.
$H_d(\sK)\cong\Z$) then $H^d(\sK)\cong\Z$, and if $\sK$ is
non-orientable then $H^d(\sK)\cong\Z_2$. What happens for
cohomology with coefficients in a field~$\k$?
\end{exercise}

\begin{exercise}\label{fundpseudo}
Let $\sK$ be an orientable $d$-dimensional pseudomanifold. The
homology group $H_d(\sK)\cong\Z$ is generated by the class of
simplicial chain $\langle\sK\rangle=\sum_{I\in\sK,\,\dim I=d}I$
where the $d$-simplices $I\in\sK$ are oriented properly (the
\emph{fundamental homology class} of~$\sK$). The cohomology group
$H^d(\sK)\cong\Z$ is generated by the class of any cochain
$\alpha_I$ taking value 1 on an oriented $d$-simplex $I\in\sK$ and
vanishing on all other simplices.
%(the \emph{fundamental cohomology class}).
\end{exercise}

\begin{exercise}
Calculate $H^*(\zk)$ where $\sK$ is the boundary of a pentagon
using Theorem~\ref{zkadd} and the description of the cohomology
product given in Proposition~\ref{proddirsum} (or
Exercise~\ref{multppair}).
\end{exercise}

\section{Coordinate subspace arrangements}\label{arran}
Here we establish a homotopy equivalence between the moment-angle
complex $\zk$ and the complement of the arrangement of coordinate
subspaces in $\C^m$ corresponding to a simplicial complex~$\sK$.
As a corollary we obtain an explicit description of the cohomology
ring of a coordinate subspace arrangement complement. In some
cases, knowing the cohomology ring allows us to identify the
homotopy type of arrangement complements.

Coordinate subspace arrangements already appeared in
Section~\ref{srr} as affine algebraic varieties corresponding to
face rings (see Proposition~\ref{affinefr}). Here we consider
these arrangements from the general point of view.

A \emph{coordinate subspace} in $\C^m$ can be given as
\begin{equation}\label{li}
  L_I=\{(z_1,\ldots,z_m)
  \in\C^m\colon z_{i_1}=\cdots=z_{i_k}=0\},
\end{equation}
where $I=\{i_1,\ldots,i_k\}$ is a subset of~$[m]$.

\begin{construction}\label{casim}
We assign to a simplicial complex $\sK$ the arrangement of complex
coordinate subspaces (or \emph{coordinate subspace arrangement})
given by
$$
  \mathcal A(\sK)=\{L_I\colon I\notin\sK\}.
$$
We denote by $U(\sK)$ the complement to $\mathcal A(\sK)$ in
$\C^m$, that is,
\begin{equation}
\label{compl}
  U(\sK)=\C^m\setminus\bigcup_{I\notin\sK}L_I.
\end{equation}
Observe that if $\sK'\subset\sK$ is a subcomplex, then
$U(\sK')\subset U(\sK)$.
\end{construction}

\begin{proposition}
\label{cacor} The assignment $\sK\mapsto U(\sK)$ defines a
bijective inclusion-preserving correspondence between simplicial
complexes on the set $[m]$ and complements of coordinate subspace
arrangements in~$\C^m$.
\end{proposition}
\begin{proof}
We need to reconstruct a simplicial complex from the complement
and check that it indeed defines the inverse correspondence. Let
$\mathcal A$ be a coordinate subspace arrangement in $\C^m$, and
let $U$ be its complement. Set
\[
  \sK(U)=
  \{I\subset[m]\colon
  L_I\cap U\ne\varnothing\}.
\]
It is easy to see that $\sK(U)$ is a simplicial complex satisfying
$U(\sK(U))=U$ and $\sK(U(\sK))=\sK$.
\end{proof}

If $\{i\}$ is a ghost vertex of $\sK$, then the coordinate
subspace arrangement $\mathcal{A}(\sK)$ contains the hyperplane
$\{z_i=0\}$. The arrangement $\mathcal A(\sK)$ does not contain
hyperplanes if and only the vertex set of $\sK$ is the
whole~$[m]$.

The complement $U(\sK)$ is an example of a polyhedral product
space (see Construction~\ref{nsc}), as is shown by the next
proposition.

\begin{proposition}\label{ukppd}
$U(\sK)=(\C,\C^\times)^\sK$.
\end{proposition}
\begin{proof}
Given a point $\mb z=(z_1,\ldots,z_m)\in\C^m$, consider its zero
set $\omega(\mb z)=\{i\in[m]\colon z_i=0\}\subset[m]$. We have
\begin{multline*}
U(\sK)=\C^m\setminus\bigcup_{I\notin\sK}L_I=
\C^m\setminus\bigcup_{I\notin\sK}\{\mb z\colon\omega(\mb z)\supset
I\}=\C^m\setminus\bigcup_{I\notin\sK}\{\mb z\colon\omega(\mb
z)=I\}\\=\bigcup_{I\in\sK}\{\mb z\colon\omega(\mb z)=I\}=
\bigcup_{I\in\sK}\{\mb z\colon\omega(\mb z)\subset I\}=
\bigcup_{I\in\sK}(\C,\C^\times)^I=(\C,\C^\times)^\sK.\qedhere
\end{multline*}
\end{proof}

\begin{example}
\label{uk}\

1. If $\sK=\varDelta^{m-1}$, then $U(\sK)=\C^m$.

2. If $\sK=\partial\varDelta^{m-1}$, then
$U(\sK)=\C^m\setminus\{0\}$.

3. Let $\sK$ be the discrete complex with $m$ vertices. Then
\[
  U(\sK)=\C^m\setminus\bigcup_{1\le i<j\le m}\{z_i=z_j=0\}
\]
is the complement to all coordinate subspaces of codimension two.

4. More generally, if $\sK$ is the $i$-dimensional skeleton of
$\varDelta^{m-1}$, then $U(\sK)$ is the complement to all
coordinate subspaces of codimension $i+2$.
\end{example}

Since each coordinate subspace is invariant under the standard
action of $\T^m$ on $\C^m$, the complement $U(\sK)$ is also a
$\T^m$-invariant subset in~$\C^m$.

A \emph{deformation retraction} of a space $X$ onto a subspace $A$
is a homotopy $F_t\colon X\to X$, $t\in\I$, such that $F_0=\id$
(the identity map), $F_1(X)=A$ and $F_t|_A=\id$ for all~$t$. The
term `deformation retraction' is often used only for the last map
$f=F_1\colon X\to A$; this map is a homotopy equivalence.

\begin{theorem}\label{deret}
The moment-angle complex $\zk$ is a $\T^m$-invariant subspace
of~$U(\sK)$, and there is a $\T^m$-equivariant deformation
retraction
\[
  \zk\hookrightarrow U(\sK)\stackrel{\simeq}\longrightarrow\zk.
\]
\end{theorem}
\begin{proof}
We have $\zk=(\D,\mathbb S)^\sK\subset(\C,\C^\times)^\sK=U(\sK)$
by the functoriality of the polyhedral product, so the
moment-angle complex $\zk$ is indeed contained in the complement
$U(\sK)$ as a $\T^m$-invariant subset.

The deformation retraction $U(\sK)\to\zk$ will be constructed by
induction. We remove simplices from $\varDelta^{m-1}$ until we
obtain $\sK$, in such a way that we get a simplicial complex at
each intermediate step.

The base of induction is clear: if $\sK=\varDelta^{m-1}$, then
$U(\sK)=\C^m$, $\zk=\D^m$, and the retraction $\C^m\to\D^m$ is
evident.

The orbit space $\zk/\T^m$ is the cubical complex
$\cc(\sK)=(\mathbb I,1)^\sK$ (see Construction~\ref{cck}). The
orbit space $U(\sK)/\T^m$ can be identified with
\[
  U(\sK)_\ge=U(\sK)\cap\R^m_\ge=(\R_\ge,\R_>)^\sK
\]
where $\R^m_\ge$ is viewed as a subset in~$\C^m$.

We shall first construct a deformation retraction $r\colon
U(\sK)_\ge\to\cc(\sK)$ of orbit spaces, and then cover it by a
deformation retraction $\widetilde r\colon U(\sK)\to\zk$.

Now assume that $\sK$ is obtained from a simplicial complex $\sK'$
by removing one maximal simplex $J=\{j_1,\ldots,j_k\}$, i.e.
$\sK\cup J=\sK'$. Then the cubical complex $\cc(\sK')$ is obtained
from $\cc(\sK)$ by adding a single $k$-dimensional face
$C_J=(\mathbb I ,1)^J$. We also have $U(\sK)=U(\sK')\setminus
L_J$, so that
\[
  U(\sK)_\ge=U(\sK')_\ge\setminus\{\mb y\colon y_{j_1}=\cdots=y_{j_k}=0\}.
\]
We may assume by induction that there is a deformation retraction
$r'\colon U(\sK')_\ge\to\cc(\sK')$ such that $\omega(r'(\mb
y))=\omega(\mb y)$, where $\omega(\mb y)$ is the set of zero
coordinates of~$\mb y$. In particular, $r'$ restricts to a
deformation retraction
\[
  r'\colon U(\sK')_\ge\setminus\{\mb y\colon y_{j_1}=\cdots=y_{j_k}=0\}
  \longrightarrow\cc(\sK')\backslash\,\mb y_J
\]
where $\mb y_J$ is the point with coordinates
$y_{j_1}=\cdots=y_{j_k}=0$ and $y_j=1$ for $j\notin J$.

Since $J\notin\sK$, we have $\mb y_J\notin\cc(\sK)$. On the other
hand, $\mb y_J$ belongs to the extra face $C_J=(\mathbb I ,1)^J$
of $\cc(\sK')$. We therefore may apply the deformation retraction
$r_J$ shown in Fig.~\ref{retr} on the face $C_J$, with centre
at~$\mb y_J$.
\begin{figure}[h]
  \begin{picture}(120,45)
  \put(45,5){\circle{2}}
  \put(80,5){\circle*{2}}
  \put(45,40){\circle*{2}}
  \put(80,40){\circle*{2}}
  \put(46,5){\vector(1,0){24}}
  \put(46,5){\line(1,0){34}}
  \put(45.8,5.2){\vector(4,1){24}}
  \put(45.8,5.2){\line(4,1){34}}
  \put(45.5,5.5){\vector(2,1){24}}
  \put(45.5,5.5){\line(2,1){34}}
  \put(45.8,5.8){\vector(4,3){24}}
  \put(45.8,5.8){\line(4,3){34}}
  \put(46,6){\vector(1,1){24}}
  \put(46,6){\line(1,1){34}}
  \put(45,6){\vector(0,1){24}}
  \put(45,6){\line(0,1){34}}
  \put(45.2,5.8){\vector(1,4){6}}
  \put(45.2,5.8){\line(1,4){8.5}}
  \put(45.5,5.5){\vector(1,2){12}}
  \put(45.5,5.5){\line(1,2){17}}
  \put(45.5,5.5){\vector(3,4){18}}
  \put(45.5,5.5){\line(3,4){26}}
  \linethickness{1mm}
  \put(80,5){\line(0,1){35}}
  \put(45,40){\line(1,0){35}}
  \put(39,4){$\mb y_J$}
  \end{picture}
  \caption{Retraction $r_J\colon \cc(\sK')\backslash\,\mb y_J
  \to\cc(\sK)$.}
  \label{retr}
\end{figure}
In coordinates, a homotopy $F_t$ between the identity map
$\cc(\sK')\backslash\,\mb y_J\to\cc(\sK')\backslash\,\mb y_J$ (for
$t=0$) and the retraction $r_J\colon\cc(\sK')\backslash\,\mb
y_J\to\cc(\sK)$ (for $t=1$) is given by
\begin{align*}
  F_t\colon\cc(\sK')\backslash\,\mb y_J&\longrightarrow
  \cc(\sK')\backslash\,\mb y_J,\\
  (y_1,\ldots,y_m,t)&\longmapsto(y_1+t\alpha_1y_1,\ldots,
  y_m+t\alpha_my_m)
\end{align*}
where
\[
  \alpha_i=\begin{cases}
  \frac{1-\max_{j\in J}y_j}{\max_{j\in J}y_j},&
  \text{if $i\in J$,}\\
  0,&\text{if $i\notin J$,}\end{cases}\qquad\text{for }1\le i\le m.
\]
We observe that $\omega(F_t(\mb y))=\omega(\mb y)$ for any $t$ and
$\mb y\in\cc(\sK')$. Now, the composition
\begin{equation}\label{indretr}
  r\colon U(\sK)_\ge=U(\sK')_\ge\!\backslash\,
  \{\mb y\colon y_{j_1}=\cdots=y_{j_k}=0\}
  \stackrel{r'}\longrightarrow
  \cc(\sK')\backslash\,\mb y_J
  \stackrel{r_J}\longrightarrow\cc(\sK)
\end{equation}
is a deformation retraction, and it satisfies $\omega(r(\mb
y))=\omega(\mb y)$ as this is true for $r_J$ and~$r'$. The
inductive step is now complete. The required retraction
$\widetilde r\colon U(\sK)\to\zk$ covers $r$ as shown in the
following commutative diagram:
\[
\xymatrix{
  \zk \ar@{^{(}->}[r]\ar[d]^{\mu} &
  U(\sK) \ar[r]^{\widetilde r}\ar[d]^{\mu} & \zk \ar[d]^{\mu}\\
  \cc(\sK) \ar@{^{(}->}[r] & U_\ge(\sK) \ar[r]^r & \cc(\sK)
}
\]
Explicitly, $\widetilde r$ is decomposed inductively in a way
similar to~\eqref{indretr},
\[
  \widetilde r\colon U(\sK)=
  U(\sK')\backslash\,L_J
  \stackrel{\widetilde r'}\longrightarrow
  \mathcal Z_{\sK'}\backslash\,\mu^{-1}(\mb y_J)
  \stackrel{\widetilde r_J}\longrightarrow\zk,
\]
where $\mu^{-1}(\mb y_J)=\prod_{j\in J}\{0\}\times\prod_{j\notin
J}\mathbb S$, and $\widetilde r_J$ is given in coordinates
$(z_1,\ldots,z_m)=\bigl(\sqrt{y_1}e^{i\varphi_1},\ldots,\sqrt{y_m}e^{i\varphi_m}\bigr)$
by
\[
  \bigl(\sqrt{y_1}e^{i\varphi_1},\ldots,
  \sqrt{y_m}e^{i\varphi_m}\bigr)\mapsto
  \bigl(\sqrt{y_1+\alpha_1y_1}e^{i\varphi_1},\ldots,
  \sqrt{y_m+\alpha_my_m}e^{i\varphi_m}\bigr)
\]
with $\alpha_i$ as above.
\end{proof}

Since $U(\sK)$ and $\zk$ are homotopy equivalent, we can use the
results on the cohomology of $\zk$ (such as Theorems~\ref{zkcoh}
and~\ref{zkadd}) to describe the cohomology rings of coordinate
subspace arrangement complements. The additive isomorphism
$H^k(U(\sK))\cong\bigoplus_{-i+2j=k}
\Tor^{-i,2j}_{\k[v_1,\ldots,v_m]}(\k[\sK],\k)$ has been also
proved in~\cite{g-p-w01}.

\begin{example}\label{codim2ar}
Let $\sK$ be the set of $m$ disjoint points. Then $\zk$ is
homotopy equivalent to the complement $U(\sK)$ of
Example~\ref{uk}.3, and
\[
  \Z[\sK]=\Z[v_1,\ldots,v_m]/(v_iv_j, \ i\ne j).
\]
The subspace of cocycles in $R^*(\sK)$ has basis of monomials
\[
  u_{i_1}u_{i_2}\cdots u_{i_{k-1}} v_{i_k}\quad\text{with }
  i_p\ne i_q\text{ for }p\ne q.
\]
Since the total degree of $u_{i_1}u_{i_2}\cdots u_{i_{k-1}}
v_{i_k}$ is $k+1$, the space of cocycles of degree ${k+1}$ has
dimension $m\bin{m-1}{k-1}$. The subspace of coboundaries of
degree $k+1$ is spanned by the elements of the form
\[
  d(u_{i_1}\cdots u_{i_k})
\]
and has dimension $\bin mk$. Therefore,
\[
\begin{array}{l}
  \rank H^{0}(U(\sK))=1,\\[2mm]
  \rank H^{1}(U(\sK))=H^{2}(U(\sK))=0,\\[2mm]
  \rank H^{k+1}(U(\sK))=
  m\bin{m-1}{k-1}-\bin mk=(k-1)\bin mk,
  \quad\text{for } 2\le k\le m,
\end{array}
\]
and the multiplication in the cohomology of $U(\sK)$ is trivial.
\end{example}

The calculation of the previous example shows that if $\sK$ is the
set of $m$ points, then there is a cohomology ring isomorphism
\begin{equation}\label{ukwed}
  H^*\bigl(U(\sK)\bigr)\cong
  H^*\Bigl(\bigvee_{k=2}^m\bigl(S^{k+1}\bigr)^{\vee(k-1)\bin mk}\Bigr),
\end{equation}
where $X^{\vee k}$ denotes the $k$-fold wedge of a space~$X$. This
cohomology isomorphism is induced by a homotopy equivalence, as is
shown by the following result.

\begin{theorem}[{\cite{gr-th04},
\cite[Corollary~9.5]{gr-th07}}]\label{gr-th} Let $\sK$ be the
$i$-dimensional skeleton of the simplex $\varDelta^{m-1}$, so that
$U(\sK)$ is the complement to all coordinate planes of codimension
$i+2$ in~$\C^m$. Then $U(\sK)$ is homotopy equivalent to a wedge
of spheres:
\[
  U(\sK)\simeq\bigvee_{k=i+2}^{m}\bigl(S^{i+k+1}\bigr)^{\vee\binom mk
  \binom{k-1}{i+1}}.
\]
\end{theorem}
The proof uses the homotopy fibration $\zk\to(\C
P^\infty)^\sK\to(\C P^\infty)^m$ (see Theorem~\ref{zkhofib}); it
will be further discussed in Section~\ref{homtypes}. For $i=0$ we
obtain the homotopy equivalence behind cohomology
isomorphism~\eqref{ukwed}.

Real coordinate subspace arrangements
\[
  U_{\R}(\sK)=(\R,\R^\times)^\sK
\]
in $\R^m$ are defined and treated similarly; there are real
analogues of all results and constructions of this section. (The
real version of Theorem~\ref{gr-th} is much simpler: see
Exercise~\ref{realskel}).

\medskip

A coordinate subspace can be given either by setting some
coordinates to zero, as in~\eqref{li}, or as the linear span of a
subset of the standard basis $\mb e_1,\ldots,\mb e_m$. The latter
approach leads to an alternative way of parametrising complements
of coordinate subspace arrangements by simplicial complexes, which
is related to the former one by Alexander duality.

Given a subset $I\subset [m]$ we set $S_I=\C\langle\mb e_i\colon
i\in I\rangle$ (the $\C$-span of the basis vectors corresponding
to~$I$), and use the notation $\widehat I=[m]\setminus I$ and
$\widehat{\sK}=\{\widehat I\subset[m]\colon I\notin\sK\}$ from
Construction~\ref{dual}. Then the coordinate subspace arrangement
corresponding to a simplicial complex~$\sK$ can be written in the
following two ways:
\[
  \mathcal A(\sK)=\{L_I\colon I\notin\sK\}=
  \{S_{\widehat I}\colon
  \widehat I\in \widehat{\sK} \}.
\]

Using Alexander duality we can reformulate the description of the
cohomology of $U(\sK)$ in terms of full subcomplexes of~$\sK$
(Theorem~\ref{zkadd}) as follows.

\begin{proposition}\label{gmcsa}
There are isomorphisms
\[
  \widetilde{H}^q\bigl(U(\sK)\bigr)\cong
  \bigoplus_{\widehat I\in\widehat{\sK}}
  \widetilde{H}_{2m-2|\widehat I|-q-2}
  (\lk_{\widehat{\sK}}\widehat I).
\]
\end{proposition}
\begin{proof}
By Theorems~\ref{deret} and \ref{zkadd},
\[
  H^q\bigl(U(\sK)\bigr)\cong
  \bigoplus_{I\subset[m]}\widetilde H^{q-|I|-1}(\sK_I).
\]
Nonempty simplices $I\in\sK$ do not contribute to the sum above,
since the corresponding subcomplexes $\sK_I$ are contractible.
Since $\widetilde{H}^{-1}(\varnothing)=\Z$, the empty simplex
contributes $\Z$ to $H^0(U(\sK))$. Therefore, we can rewrite the
isomorphism above as
\[
  \widetilde{H}^q\bigr(U(\sK)\bigl)
  \cong\bigoplus_{I\notin\sK}\widetilde{H}^{q-|I|-1}(\sK_I).
\]
Using Alexander duality (Proposition~\ref{fulink}) we calculate
\[
  \widetilde{H}^{q-|I|-1}(\sK_I)\cong
  \widetilde H_{|I|-3-q+|I|+1}(\lk_{\widehat
  \sK}\widehat I)=
  \widetilde H_{2m-2|\widehat I|-q-2}(\lk_{\widehat
  \sK}\widehat I),
\]
where $\widehat I=[m]\backslash I$ is a simplex of~$\widehat\sK$.
\end{proof}

Proposition~\ref{gmcsa} is a particular case of the well-known
\emph{Goresky--MacPherson
formula}~\cite[Chapter~III]{go-ma88}\label{gomaformu}, which
calculates the (co)homology groups of the complement of an
arrangement of affine subspaces in terms of its \emph{intersection
poset}. In the case of coordinate subspace arrangements $\mathcal
A(\sK)$ the intersection poset is the face poset of the dual
complex~$\widehat\sK$. For more on the relationships between
general affine subspace arrangements and moment-angle complexes
see~\cite[Chapter~8]{bu-pa02}.

\subsection*{Exercises}
\begin{exercise}
The affine algebraic variety $X(\sK)$ corresponding to the face
ring $\C[\sK]$ (see Proposition~\ref{affinefr}) and the coordinate
subspace arrangement $\mathcal A(\sK)$ of Construction~\ref{casim}
are related by the identity $X(\widehat\sK)=\mathcal A(\sK)$,
where $\widehat\sK=\{I\subset[m]\colon[m]\setminus I\notin\sK\}$
is the Alexander dual complex.
\end{exercise}

\begin{exercise}\label{ex3clines}
Show directly that the complement to the 3 coordinate lines in
$\C^3$ is homotopy equivalent to the wedge of spheres $S^3\vee
S^3\vee S^3\vee S^4\vee S^4$; this corresponds to $m=3$
in~\eqref{ukwed}.
\end{exercise}

\begin{exercise}\label{u2con}
Show directly (without referring to Theorem~\ref{deret} and
Proposition~\ref{homgr}) that the complement $U(\sK)$ is
2-connected if $\sK$ does not have ghost vertices, and that
$U(\sK)$ is $2q$-connected if $\sK$ is $q$-neighbourly.
\end{exercise}

\chapter*{\ \ Moment-angle complexes: additional topics}

\section{Free and almost free torus actions on moment-angle complexes}
\label{freetorusactions} Here consider free and almost free
actions of toric subgroups $T^k\subset\T^m$ on~$\zk$. As usual,
$\sK$ is an $(n-1)$-dimensional simplicial complex on~$[m]$, and
$\zk$ is the corresponding moment-angle complex.

We start with a simple characterisation of the isotropy subgroups
of the standard $\T^m$-action on~$\zk$. For each $I\subset[m]$ we
consider the coordinate subtorus
\[
  \T^I=\{(t_1,\ldots,t_m)\in\T^m\colon t_j=1\text{ for
  }j\notin I\}=\prod_{i\in I}\T\;\subset\;\T^m
\]
(note that $\T^I=(\T,1)^I$ in the notation of
Construction~\ref{nsc}).

\begin{proposition}\label{isotrzk}
Let $\mb z\in\zk$, and set $\omega(\mb z)=\{i\in[m]\colon
z_i=0\}\in\sK$. Then the isotropy subgroup of $\mb z$ with respect
to the $\T^m$-action is $\T^{\omega(\mb z)}$. Furthermore, each
coordinate subtorus $\T^I$ for $I\in\sK$ is the isotropy subgroup
for a point $\mb z\in\zk$.
\end{proposition}
\begin{proof}
An element $\mb t=(t_1,\ldots,t_m)\in\T^m$ fixes $\mb z$ if and
only if $t_i=1$ whenever $z_i\ne0$, which is equivalent to that
$\mb t\in\T^{\omega(\mb z)}$. The last statement is also clear:
$\T^I$ is the isotropy subgroup for any $\mb z\in(\D,\mathbb
S)^I\subset\zk$ with $\omega(\mb z)=I$.
\end{proof}

Recall that an action of a group on a topological space is
\emph{almost free} if all isotropy subgroups are finite.

\begin{definition}\label{defftra}
We define the \emph{free toral rank} of $\zk$, denoted $\ftr\zk$,
as the maximal dimension of toric subgroups $T^k\subset\T^m$
acting on $\zk$ freely. Similarly, the \emph{almost free toral
rank} of $\zk$, denoted $\atr\zk$, is the maximal dimension of
toric subgroups $T^k\subset\T^m$ acting on $\zk$ almost freely.
\end{definition}

\begin{proposition}
Let $\sK$ be a simplicial complex of dimension $n-1$ on $m$
vertices and $\sK\ne\varDelta^{m-1}$. The toral ranks of $\zk$
satisfy the following inequalities:
\[
  1\le\ftr\zk\le\atr\zk\le m-n.
\]
\end{proposition}
\begin{proof}
By Proposition~\ref{isotrzk}, isotropy subgroups for the
$\T^m$-action on~$\zk$ are coordinate subgroups of the form
$\T^I$. The diagonal circle in $\T^m$ intersects each of these
coordinate subgroups trivially (since $I\ne[m]$), and therefore
acts freely on~$\zk$. This proves the first inequality. The second
is obvious. To prove the third one, assume that $T^k\subset\T^m$
acts almost freely on~$\zk$. Then the intersection of $T^k$ with
every $\T^m$-isotropy subgroup $\T^I$ is a finite group. Choose a
maximal simplex $I\in\sK$, \ $|I|=n$. Then $\T^I\cap T^k$ can be
finite only if $k\le m-n$.
\end{proof}

The map  $\R^m\to\T^m$, $(\phi_1,\ldots,\phi_m)\mapsto
  (e^{2\pi i\phi_1},\ldots,e^{2\pi i\phi_m})$,
identifies
%the standard torus
$\T^m$ with the quotient
$\R^m/\Z^m$. Subtori $T^k\subset\T^m$ of dimension $k$ bijectively
correspond to unimodular sublattices $L\subset\Z^m$ of rank~$k$ (a
sublattice is \emph{unimodular} if it is a direct summand
in~$\Z^m$). The inclusion $T^k\subset\T^m$ can be viewed as
$L_\R/L\subset\R^m/\Z^m$, where $L_\R$ is the $k$-dimensional
subspace in $\R^m$ spanned by~$L$.

Choosing a basis in $L$ we obtain an integer $m\times k$-matrix
$S=(s_{ij})$, \ $1\le i\le m$, $1\le j\le k$, so that $L$ is
identified with the image of $S\colon\Z^k\to\Z^m$. The $k$-torus
$T^k$ is the image of the corresponding monomorphism of tori
$\T^k\to\T^m$, namely,
\begin{equation}\label{k-subtorus}
  T^k=\bigl\{\bigl(e^{2\pi i(s_{11}\psi_1+\cdots+s_{1k}\psi_k)},
  \ldots,e^{2\pi i(s_{m1}\psi_1+\cdots+s_{mk}\psi_k)} \bigr)\bigr\}
  \subset\T^m,
\end{equation}
where $(\psi_1,\ldots,\psi_k)\in\R^k$. Since $L$ is unimodular,
the columns of $S$ form a part of basis of the lattice~$\Z^m$.
%
%For each subset $I\subset[m]$ we get the corresponding unimodular
%coordinate sublattice $\Z^I\subset\Z^m$ spanned by the basis
%vectors $\mb e_i$ with $i\in I$.

\begin{lemma}\label{fafcrit}
Let $T^k$ be a $k$-dimensional subtorus in $\T^m$ and let $L$ be
the corresponding unimodular sublattice of rank $k$ in~$\Z^m$. Let
$S=(s_{ij})$, \ $1\le i\le m$, $1\le j\le k$, be a matrix defining
$L$, so that $T^k$ is given by~\eqref{k-subtorus}.

\begin{itemize}
\item[(a)] The torus $T^k$ acts on $\zk$ almost freely if and
only if for each $I\in\sK$ the intersection of subspaces $L_\R$
and $\R^I$ in $\R^m$ is zero. Equivalently, the $(m-|I|)\times
k$-matrix $S_{\widehat I}$ obtained by deleting from $S$ the rows
with numbers $i\in I$ has rank~$k$.

\item[(b)] The torus $T^k$ acts on $\zk$ freely if and
only if for each $I\in\sK$ the sublattice spanned by $L$ and
$\Z^I$ in $\Z^m$ is unimodular of rank~$k+|I|$. Equivalently, the
columns of the $(m-|I|)\times k$-matrix $S_{\widehat I}$ form a
part of basis of~$\Z^{m-|I|}$.
\end{itemize}
\end{lemma}

\begin{proof}
We prove (a) first. By Proposition~\ref{isotrzk}, the $T^k$-action
on $\zk$ is almost free if and only if the intersection
$T^k\cap\T^I\subset\T^m$ is finite for each $I\in\sK$. This
intersection can be identified with the kernel of the map $f\colon
T^k\times\T^I\to\T^m$ (the product of the inclusion maps
$T^k\to\T^m$ and $\T^I\to\T^m$). This kernel is finite if and only
if the corresponding map of real spaces $L_\R\times\R^I\to\R^m$ is
injective, which is equivalent to that $L_\R\cap\R^I=\{0\}$. Let
$I=\{i_1,\ldots,i_p\}$, then the matrix of $f$ has the form $(S\,|
\mb e_{i_1}|\cdots|\mb e_{i_p})$, where $\mb e_i$ is the $i$th
standard basis column vector. Clearly, this matrix has rank
$k+|I|$ if and only if the matrix $S_{\widehat I}$ has rank~$k$.

Now we prove~(b). The $T^k$-action on $\zk$ is free if and only if
the the kernel of $f\colon T^k\times\T^I\to\T^m$ is trivial for
each $I\in\sK$, i.e. $T^k\times\T^I$ embeds as a subtorus. This is
equivalent to the conditions stated in~(b).
\end{proof}

Let $\mb t=(t_1,\ldots,t_n)\in\Z[\sK]$ be a linear sequence given
by
\begin{equation}\label{linearseq1}
  t_i=\lambda_{i1}v_1+\cdots+\lambda_{im}v_m,\quad
  \text{for }1\le i\le n.
\end{equation}
We consider the integer $n\times m$-matrix
$\Lambda=(\lambda_{ij})$, \ $1\le i\le n$, $1\le j\le m$. It
defines a homomorphism of lattices $\Lambda\colon\Z^m\to\Z^n$ and
a homomorphism of tori $\Lambda\colon\T^m\to\T^n$.

\begin{theorem}\label{afreeeq}
The following conditions are equivalent:
\begin{itemize}
\item[(a)] the sequence $(t_1,\ldots,t_n)$ given by~\eqref{linearseq1} is an lsop in the
rational face ring $\Q[\sK]$;
\item[(b)] the kernel $T_\Lambda=\Ker(\Lambda\colon\T^m\to\T^n)$
is a product of an $(m-n)$-torus and a finite group, and
$T_\Lambda$ acts almost freely on~$\zk$.
\end{itemize}
\end{theorem}
\begin{proof}
We first observe that under any of conditions~(a) or~(b) the
rational map $\Lambda\colon\Q^m\to\Q^n$ is surjective. For each
simplex $I\in\sK$ we consider the diagram
\[
\begin{CD}
  @.@.\begin{array}{c}0\\ \raisebox{2pt}{$\downarrow$}\\
  \Q^{m-n}\end{array}@.@.\\
  @.@.@VV{i_1}V@.@.\\
  0@>>>\Q^{|I|}@>{i_2}>>\Q^m@>{p_2}>>\Q^{m-|I|} @>>>0\\
  @.@.@VV {\Lambda} V@.@.\\
  @.@.\begin{array}{c}\Q^{n}\\ \downarrow\\ 0\end{array}@.@.\\
\end{CD}
\]
where $i_2$ is the inclusion of the coordinate subspace
$\Q^I\to\Q^m$. The map $\Lambda i_2$ is given by the $n\times|I|$
matrix $\Lambda_I=(\lambda_{ij})$, \ $1\le i\le n$, $j\in I$. By
Lemma~\ref{lsopcrit}, the sequence $(t_1,\ldots,t_n)$ is a
rational lsop if and only if the rank of $\Lambda_I$ is~$|I|$ for
each $I\in\sK$. Hence condition~(a) of the theorem is equivalent
to the injectivity of the map $\Lambda i_2$ for each $I\in\sK$.

On the other hand, by Lemma~\ref{fafcrit}~(a), the action of the
$(m-n)$-torus $T_\Lambda$ on $\zk$ is almost free if and only if
$i_1(\Q^{m-n})\cap\Q^I=\{0\}$ for each $I\in\sK$. The latter
condition is equivalent to that $i_1(\Q^{m-n})\cap\Ker p_2=\{0\}$,
i.e. that $p_2i_1$ is injective. Hence condition~(b) of the
theorem is equivalent to the injectivity of the map $p_2i_1$ for
each $I\in\sK$. Now the theorem follows from Lemma~\ref{2ses}.
\end{proof}

The almost free toral rank of $\zk$ can now be easily determined:

\begin{corollary}\label{detatr}
We have $\atr\zk=m-n$, i.e. for each simplicial complex $\sK$ of
dimension $(n-1)$ on $[m]$ there is an $(m-n)$-dimensional
subtorus in $\T^m$ acting on $\zk$ almost freely.
\end{corollary}
\begin{proof}
Choose a rational lsop in $\Q[\sK]$ by Theorem~\ref{noether}, and
multiply it by a common denominator to get an integral
sequence~\eqref{linearseq1}. It is still an lsop in $\Q[\sK]$ (but
it may fail to be an lsop in $\Z[\sK]$), and therefore the
$(m-n)$-torus $T_\Lambda$ acts on $\zk$ almost freely by
Theorem~\ref{afreeeq}.
\end{proof}

There is an analogue of Theorem~\ref{afreeeq} for free torus
actions:

\begin{theorem}
The following conditions are equivalent:
\begin{itemize}
\item[(a)] the sequence $(t_1,\ldots,t_n)$ given by~\eqref{linearseq1}
is an lsop in $\Z[\sK]$;
\item[(b)] $T_\Lambda=\Ker(\Lambda\colon\T^m\to\T^n)$
is an $(m-n)$-torus acting freely on~$\zk$.
\end{itemize}
\end{theorem}
\begin{proof}
The argument is the same as in the proof of Theorem~\ref{afreeeq}:
Lemma~\ref{2ses} is now applied to the diagram of integral
lattices instead of rational vector spaces.
\end{proof}

Nevertheless, there is no analogue of Corollary~\ref{detatr} for
free torus actions, as integral lsop's may fail to exist. Indeed,
the free toral rank of $\zk$ where $\sK$ is the boundary of a
cyclic $n$-polytope with $m\ge2^n$ vertices is strictly less than
$m-n$, as shown by Example~\ref{dj-lsop}. The free toral rank of
$\zk$ is a combinatorial characteristic of~$\sK$, also known as
the \emph{Buchstaber invariant}\label{buchinvari}. Its
determination is a very subtle problem; see~\cite{erok08}
and~\cite{fu-ma11} for some partial results in this direction.

\medskip

There is the following important conjecture of equivariant
topology and rational homotopy theory concerning almost free torus
action.

\begin{conjecture}[{Toral Rank Conjecture, Halperin~\cite{halp85}}]
Assume that a torus $T^k$ acts almost freely on a
finite-dimensional topological space $X$. Then
\[
  \dim H^*(X;\Q)\ge2^k,
\]
i.e. the total dimension of the cohomology of $X$ is at least that
of the torus~$T^k$.
\end{conjecture}

The Toral Rank Conjecture is valid for $k\le3$ and is open in
general. See~\cite{pupp09} and~\cite{usti12} for the discussion of
the current status of this conjecture.

In the case of moment-angle complexes we have the following
result:

\begin{theorem}[{\cite{ca-lu12}, \cite{usti11}}]\label{trzk}
Let $\sK$ be a simplicial complex of dimension~$n-1$ with $m$
vertices, and let $\zk$ be the corresponding moment-angle complex.
Then
\[
  \rank H^*(\zk)=\sum_{k=0}^{m+n}\rank H^k(\zk)
  %=\sum_{J\subset[m],\,k\ge0}\widetilde H^{k-|J|-1}(\sK_J)
  \ge 2^{m-n}.
\]
\end{theorem}

The proof of this theorem given in~\cite{usti11} uses a
construction of independent interest and a couple of technical
lemmata. We include this proof below.

\begin{corollary}
The Toral Rank Conjecture is valid for subtori $T^k\subset\T^m$
acting almost freely on~$\zk$.
\end{corollary}

The following is a particular case of the so-called
\emph{simplicial wedge
construction}\label{simplicialwedge}~\cite{pr-bi80}. It has been
brought into toric topology by the work~\cite{b-b-c-g}.

\begin{construction}[Simplicial doubling]
Let $\sK$ be a simplicial complex on the vertex set $[m]$. The
\emph{double} of $\sK$ is the simplicial complex $\mathcal D(\sK)$
on the vertex set~$[2m]=\{1,1',2,2',\ldots,m,m'\}$ whose missing
faces (minimal non-faces) are $\{i_1,i'_1,\ldots,i_k,i'_k\}$ where
$\{i_1,\ldots,i_k\}$ is a missing face of~$\sK$. In other words,
$\mathcal D(\sK)$ is determined by its face ring given by
\[
  \k\bigl[\mathcal D(\sK)\bigr]=\k[v_1,v_{1'},\ldots,v_m,v_{m'}]\big/
  \bigl(\{v_{i_1}v_{i'_1}\cdots v_{i_k}v_{i'_k}\}\colon
  \{i_1,\ldots,i_k\}\notin\sK\bigr).
\]
\end{construction}

\begin{example}\

1. If $\sK=\varDelta^{m-1}$ (the full simplex on $m$ vertices),
then $\mathcal D(\sK)=\varDelta^{2m-1}$.

2. If $\sK=\partial\varDelta^{m-1}$, then $\mathcal
D(\sK)=\partial\varDelta^{2m-1}$.
\end{example}

The doubling construction interacts nicely with the polyhedral
product:

\begin{theorem}\label{doublezk}
Let $(X,A)$ be a pair of spaces, let $\sK$ be a simplicial complex
on~$[m]$ and $\mathcal D(\sK)$ its double. Then
\[
  (X,A)^{\mathcal D(\sK)}=
  (X\times X,X\times A\cup A\times X)^\sK.
\]
\end{theorem}
\begin{proof}
Set $(Y,B)=(X\times X,X\times A\cup A\times X)$. Given a point
$\mb y=(y_1,\ldots,y_m)\in Y^m$, we set
\[
  \omega_Y(\mb y)=\{i\in[m]\colon y_i\notin B\}\subset[m].
\]
Similarly, given $\mb x=(x_1,x_{1'},\ldots,x_m,x_{m'})\in X^{2m}$,
we set
\[
  \omega_X(\mb x)=\bigl\{j\in\{1,1',\ldots,m,m'\}\colon x_j\notin A\bigr\}\subset
  \{1,1',\ldots,m,m'\}.
\]
We identify $\mb y$ with $\mb x$ by the formula
$(y_1,\ldots,y_m)=((x_1,x_{1'}),\ldots,(x_m,x_{m'}))\in
Y^m=X^{2m}$. It follows from the definition of the polyhedral
product that $\mb y\notin(Y,B)^\sK$ if and only if $\omega_Y(\mb
y)\notin\sK$. The latter is equivalent to the condition
$\omega_X(\mb x)\notin \mathcal D(\sK)$, since if $\omega_Y(\mb
y)=\{i_1,\ldots,i_k\}$ then $\omega_X(\mb
x)\supset\{i_1,i'_1,\ldots,i_k,i'_k\}$. Therefore,
\[
  \mb y\notin(Y,B)^\sK\;\Longleftrightarrow\;\mb x\notin
  (X,A)^{\mathcal D(\sK)},
\]
which implies that $(X,A)^{\mathcal D(\sK)}=(Y,B)^\sK$.
\end{proof}

\begin{remark}
The simplicial wedge~\cite{pr-bi80}, \cite{b-b-c-g} is a
generalisation of the doubling construction, in which the $i$th
vertex of $\sK$ is replaced by a $j_i$-tuple of vertices, for some
vector $\mb j=(j_1,\ldots,j_m)$ of natural numbers. The double
corresponds to $\mb j=(2,\ldots,2)$. There is an analogue of
Theorem~\ref{doublezk} in this setting, see~\cite[\S7]{b-b-c-g}.
\end{remark}

As an important consequence of Theorem~\ref{doublezk} we obtain
the following relationship between the moment-angle complex $\zk$
and its real analogue~$\rk$:

\begin{corollary}
We have $\zk\cong\mathcal R_{\mathcal D(\sK)}$.
\end{corollary}
\begin{proof}
Apply Theorem~\ref{doublezk} to the pair $(X,A)=(D^1,S^0)$,
observing that $(D^1\times D^1,D^1\times S^0\cup S^0\times
D^1)\cong(D^2,S^1)$.
\end{proof}

\begin{lemma}\label{mvdouble}
Let $(X,A)$ be a pair of cell complexes such that $A$ has a collar
neighbourhood $U(A)$ in $X$ (i.e. there is a homeomorphism of
pairs $(U(A),A)\cong(A\times[0,1),A\times\{0\})$). Let
$Y=X_1\cup_A X_2$ be the space obtained by attaching two copies of
$X$ along~$A$. Then
%the total amount of cohomology of $Y$ is
%bounded from below by the total amount of cohomology of~$A$:
\[
  \rank H^*(Y)\ge\rank H^*(A).
\]
\end{lemma}
\begin{proof}
The assumption on $(X,A)$ implies that we can apply the
Mayer--Vietoris sequence to the decomposition $Y=X_1\cup_A X_2$:
\[
  \cdots\stackrel{\beta_{k-1}}\longrightarrow
  H^{k-1}(A)\stackrel{\delta_{k-1}}{\longrightarrow}
  H^k(Y)\stackrel{\alpha_k}\longrightarrow
  H^k(X_1)\oplus H^k(X_2)\stackrel{\beta_k}\longrightarrow
  H^k(A)\to\cdots.
\]
The map $\beta_k$ is $i_1^*\oplus(-i_2^*)$, where $i_1\colon A\to
X_1$ and $i_2\colon A\to X_2$ are the inclusions. Since
$X_1=X_2=X$, and the inclusions $i_1$ and $i_2$ coincide, we have
$\rank\Ker\beta_k\ge\rank H^k(X)$ and $\rank\Im\beta_k\le\rank
H^k(X)$. Using these inequalities we calculate
\begin{align*}
  \rank H^k(Y)&=\rank\Ker\alpha_k+\rank\Im\alpha_k=
  \rank\Im\delta_{k-1}+\rank\Ker\beta_k\\&\ge
  \rank H^{k-1}(A)-\rank\Im\beta_{k-1}+\rank H^k(X)\\
  &\ge\rank H^{k-1}(A)-\rank H^{k-1}(X)+\rank H^k(X).
\end{align*}
The required inequality is obtained by summing up over $k$.
\end{proof}

\begin{theorem}\label{trrk}
Let $\sK$ be a simplicial complex of dimension~$n-1$ with $m$
vertices, and let $\rk$ be the corresponding real moment-angle
complex. Then
\[
  \rank H^*(\rk)
  \ge 2^{m-n'}\ge 2^{m-n},
\]
where $n'$ is the minimum of the cardinality of maximal simplices
of~$\sK$ (so that $n'=n=\dim\sK+1$ if and only if $\sK$ is pure).
\end{theorem}
\begin{proof}[Proofs of Theorems~\ref{trzk} and~\ref{trrk}]
We first prove the inequality for $\rk$, by induction on the
number of vertices~$m$. For $m=1$ the statement is clear. We embed
$\rk$ as a subcomplex in the `big' cube $[-1,1]^m$ (see
Construction~\ref{realmac}) with coordinates $\mb
u=(u_1,\ldots,u_m)$, $-1\le u_i\le 1$. Assume that the first
vertex of $\sK$ belongs to an $(n'-1)$-dimensional maximal
simplex of~$\sK$, and consider the following subspaces of~$\rk$:
\begin{gather*}
X_+=\{\mb u\in\rk\colon u_1\ge0\},\quad X_-=\{\mb u\in\rk\colon
u_1\le0\},\\
A=X_+\cap X_-=\{\mb u\in\rk\colon u_1=0\}.
\end{gather*}
Applying Lemma~\ref{mvdouble} to the decomposition $\rk=X_+\cup_A
X_-$ we obtain
\[
  \rank H^*(\rk)\ge\rank H^*(A).
\]
On the other hand, $A$ is the disjoint union of $2^{m-m_1-1}$
copies of $\mathcal R_{\lk_{\sK}\{1\}}$, where $m_1$ is the number
of vertices of $\lk_{\sK}\{1\}$. Since $\{1\}$ is a vertex of a
maximal simplex of~$\sK$ of minimal cardinality~$n'$, the minimal
cardinality of maximal simplices in $\lk_{\sK}\{1\}$ is
$n'_1=n'-1$. Now using the inductive hypothesis we obtain
\[
  %\dim H^*(\rk;\Q)\ge
  \rank H^*(A)=
  2^{m-m_1-1}\rank H^*(\mathcal R_{\lk_{\sK}\{1\}})\ge
  2^{m-m_1-1}2^{m_1-n'_1}=2^{m-n'}.
\]
Theorem~\ref{trrk} is therefore proved.

To prove Theorem~\ref{trzk} we use the fact that $\zk\cong\mathcal
R_{\mathcal D(\sK)}$, and observe that the numbers $m-n$ (and
$m-n'$) for $\sK$ and $\mathcal D(\sK)$ coincide.
\end{proof}

Using Theorems~\ref{zkcoh} and~\ref{zkadd} we may reformulate
Theorem~\ref{trzk} in both algebraic and combinatorial terms:

\begin{theorem}
Let $\sK$ be a simplicial complex of dimension~$n-1$ with $m$
vertices, and let $\k$ be a field. Then
\[
  \dim\Tor_{\k[v_1,\ldots,v_m]}\bigl(\k[\sK],\k\bigr)
  =\sum_{J\subset[m],\,k\ge0}\dim\widetilde H^{k-|J|-1}(\sK_J;\k)
  \ge 2^{m-n}.
\]
\end{theorem}

As a corollary we obtain that the weak Horrocks Conjecture
(Conjecture~\ref{weakhorrocks}) holds for a particular class of
rings:

\begin{corollary}
Let $\sK$ be a Cohen--Macaulay simplicial complex (over a
field~$\k$) of dimension~$n-1$ with $m$ vertices. Let $\mb
t=(t_1,\ldots,t_n)$ be an lsop in $\k[\sK]$, so that $\k[m]/\!\mb
t\cong\k[w_1,\ldots,w_{m-n}]$ and $\dim\k[\sK]/\!\mb t<\infty$.
Then
\[
  \dim\Tor_{\k[w_1,\ldots,w_{m-n}]}\bigl(\k[\sK]/\!\mb t\,,\k\bigr)
  \ge 2^{m-n},
\]
i.e. the weak Horrocks Conjecture holds for the rings
$\k[\sK]/\!\mb t$.
\end{corollary}
\begin{proof}
This follows from the previous theorem and
Proposition~\ref{reduc}.
\end{proof}

\subsection*{Exercises}
\begin{exercise}
Show that $\ftr\zk=1$ if and only if
$\sK=\partial\varDelta^{m-1}$.
\end{exercise}

\begin{exercise}
The Toral Rank Conjecture fails if $\dim X=\infty$.
\end{exercise}

\begin{exercise}
Show that the doubling operation respects the join, that is,
$\mathcal D(\sK*\mathcal L)=\mathcal D(\sK)*\mathcal D(\mathcal
L)$.
\end{exercise}

\begin{exercise}
Assume that $\sK$ is the boundary complex of a simplicial
$n$-polytope~$Q\subset\R^n$ with~$m$ vertices $v_1,\ldots,v_m$.
Then $\mathcal D(\sK)$ is the boundary of a simplicial polytope
$\mathcal D(Q)$ of dimension~$m+n$ with $2m$ vertices, which can
be obtained in the following way. We embed $\R^n$ as the
coordinate subspace in~$\R^{m+n}$ on the last $n$ coordinates. For
each vertex $v_i\in Q\subset\R^n$ take the line
$l_i\subset\R^{m+n}$ through $v_i$ parallel to the $i$th
coordinate line of~$\R^{m+n}$, for $1\le i\le m$. Then replace
each $v_i$ by a pair of points $v'_i,v''_i\in l_i$ such that $v_i$
the centre of the segment with the vertices $v'_i,v''_i$. Then the
boundary of
\[
  \mathcal D(Q)=\conv(v'_1,v''_1,\ldots,v'_m,v''_m)\subset\R^{m+n}
\]
is~$\mathcal D(\sK)$.
\end{exercise}

\begin{exercise}
There is the following generalisation of Theorem~\ref{doublezk}.
Let
\[
  (\mb X,\mb A)=\bigl\{(X_1,A_1),(X_{1'},A_{1'}),\ldots,
  (X_m,A_m),(X_{m'},A_{m'})\bigr\}
\]
be a set of $2m$ pairs of spaces. Define a new set $(\mb Y,\mb
B)=\{(Y_1,B_1),\ldots,(Y_m,B_m)\}$ of $m$ pairs, where
\[
  (Y_i,B_i)=\bigl(X_{i}\times X_{i'},
  X_{i}\times A_{i'}\cup A_i\times X_{i'}\bigr).
\]
Then
\[
  (\mb X,\mb A)^{\mathcal D(\sK)}\cong
  (\mb Y,\mb B)^\sK.
\]
Further generalisations can be found in~\cite[\S7]{b-b-c-g}.
\end{exercise}

\begin{exercise}
The inequality $\rank H^*(\rk)\ge 2^{m-n}$ of Theorem~\ref{trrk}
(or the inequality $\rank H^*(\zk)\ge 2^{m-n}$) turns into
identity if and only if
\[
  \sK=\partial\varDelta^{k_1-1}*\partial\varDelta^{k_2-1}*\cdots*
  \partial\varDelta^{k_p-1}*\varDelta^{m-s-1},
\]
where $s=k_1+\cdots+k_p$ and the join factor $\varDelta^{m-s-1}$
is void if $s=m$ (compare Exercise~\ref{frcomplint}). In this case
both $\rk$ and $\zk$ are products of spheres and a disc.
\end{exercise}
%
%\begin{exercise}\label{hortrc}
%Assume that a torus $T^k$ acts freely on a finite-dimensional
%space $Z$, and let $M=Z/T^k$ be the quotient. Consider the
%Leray--Serre spectral sequence of the principal $T^k$-bundle $Z\to
%M$; its $E_2$-term is the Koszul algebra of the form
%$(\Lambda[u_1,\ldots,u_k]\otimes H^*(M),d)$. Show that the
%spectral sequence collapses at the $E_3$-term, so that we have an
%additive isomorphism
%\[
%  H^*(Z;\k)\cong\Tor_{\k[w_1,\ldots,w_k]}\bigr(H^*(M;\k),\k\bigr),
%\]
%where $\dim H^*(M;\k)<\infty$. Deduce that the weak Horrocks
%Conjecture implies the Toral Rank Conjecture for free torus
%actions. For almost free actions,
%see~\cite[Proposition~2.1]{usti12}.
%\end{exercise}

\section{Massey products in the cohomology of moment-angle
complexes}\label{Masseymac} Here we address the question of
existence of nontrivial triple Massey products in the Koszul
complex
\[
  (\Lambda[u_1,\ldots,u_m]\otimes\Z[\sK],d)
\]
of the face ring, and therefore in the cohomology of~$\zk$. The
general definition of Massey products in the cohomology of a
differential graded algebra is reviewed in Section~\ref{dgaap} of
the Appendix. A geometrical approach to constructing nontrivial
triple Massey products in the Koszul complex of the face ring was
developed by Baskakov in~\cite{bask03} as an extension of the
cohomology calculation in Theorem~\ref{zkhoch}. It is well-known
that non-trivial higher Massey products obstruct the
\emph{formality} of a differential graded algebra, which in our
case leads to a family of nonformal moment-angle manifolds $\zk$
(see Section~\ref{aprht} for the background material).

\begin{construction}[Baskakov~\cite{bask03}]\label{baskconstr}
Let $\sK_i$ be a triangulation of a sphere $S^{n_i-1}$ with
$|V_i|=m_i$ vertices, $i=1,2,3$. Set $m=m_1+m_2+m_3$,
$n=n_1+n_2+n_3$,
\[
  \sK=\sK_1*\sK_2*\sK_3,\quad
  \zk=\mathcal Z_{\sK_1}\times\mathcal Z_{\sK_2}\times\mathcal
  Z_{\sK_3}.
\]
Then $\sK$ is a triangulation of $S^{n-1}$ and therefore $\zk$ is
an $(m+n)$-manifold.

We choose maximal simplices $I_1\in\sK_1$, \ $I'_2,I''_2\in\sK_2$
such that $I_2'\cap I_2''=\varnothing$, and $I_3\in\sK_3$. Set
\[
  \widetilde\sK=\ss_{I_1\cup I_2'}(\ss_{I_2''\cup I_3}\sK),
\]
where $\ss_I$ denotes the stellar subdivision at~$I$, see
Definition~\ref{bist}. Then $\widetilde\sK$ is a triangulation of
$S^{n-1}$ with $m+2$ vertices. Take generators
\[
  \beta_i\in\widetilde H^{n_i-1}(\widetilde\sK_{V_i})
  \cong\widetilde H^{n_i-1}(S^{n_i-1}),\quad\text{for } i=1,2,3,
\]
where $\widetilde\sK_{V_i}$ is the restriction of $\widetilde\sK$
to the vertex set of $\sK_i$, and set
\[
  \alpha_i=h(\beta_i)\in H^{n_i-m_i,2m_i}
  (\mathcal Z_{\widetilde\sK})\subset H^{m_i+n_i}
  (\mathcal Z_{\widetilde\sK}),
\]
where $h$ is the isomorphism of Theorem~\ref{zkhoch}. Then
\[
  \beta_1\beta_2\in\widetilde H^{n_1+n_2-1}
  (\widetilde\sK_{V_1\sqcup V_2})\cong\widetilde
  H^{n_1+n_2-1}(S^{n_1+n_2-1}\setminus\pt=0,
\]
hence $\alpha_1\alpha_2=h(\beta_1\beta_2)=0$, and similarly
$\alpha_2\alpha_3=0$. Therefore, the triple Massey product
$\langle \alpha_1,\alpha_2,\alpha_3\rangle\subset
H^{m+n-1}(\mathcal Z_{\widetilde\sK})$ is defined. By definition,
it is the set of cohomology classes represented by the cocycles
$(-1)^{\deg a_1+1}a_1f+ea_3$ where $a_i$ is a cocycle representing
$\alpha_i$, and $e$, $f$ are cochains satisfying $de=a_1a_2$,
$df=a_2a_3$.

For the simplest example of this series, take $\sK_i=S^0$ (two
points), so that $\sK$ is the boundary of an octahedron, and
$\widetilde\sK$ is obtained by applying stellar subdivisions at
two skew edges. We shall consider this example in more detail
below.
\end{construction}

Recall that a Massey product is \emph{trivial} if it contains
zero.

\begin{theorem}\label{ntmas}
The above defined triple Massey product
\[
  \langle \alpha_1,\alpha_2,\alpha_3\rangle\subset
  H^{m+n-1}(\mathcal Z_{\widetilde K})
\]
%in the cohomology of $(m+n+2)$-manifold $\mathcal Z_{\widetilde K}$
is nontrivial.
\end{theorem}
\begin{proof}
Consider the subcomplex of $\widetilde\sK$ consisting of the two
new vertices added to $\sK$ in the process of stellar subdivision.
By Proposition~\ref{mafunc}~(b), the inclusion of this subcomplex
induces an embedding of a 3-dimensional sphere
$S^3\hookrightarrow\mathcal Z_{\widetilde K}$. Since the two new
vertices are not joined by an edge in $\mathcal
Z_{\widetilde\sK}$, the embedded 3-sphere defines a non-trivial
class $x\in H_3(\mathcal Z_{\widetilde\sK})$. Its Poincar\'e dual
cohomology class in $H^{m+n-1}(\mathcal Z_{\widetilde\sK})$ is
contained in the Massey product
$\langle\alpha_1,\alpha_2,\alpha_3\rangle$. We need only to check
that this element cannot be turned into zero by adding elements
from the indeterminacy of the Massey product, i.e. from the
subspace
\[
  \alpha_1\cdot H^{m_2+m_3+n_2+n_3-1}(\mathcal Z_{\widetilde\sK})+
  \alpha_3\cdot H^{m_1+m_2+n_1+n_2-1}(\mathcal Z_{\widetilde\sK}).
\]
To do this we use the multigraded structure in~$H^*(\mathcal
Z_{\widetilde\sK})$. The multigraded components of the group
$H^{m_2+m_3+n_2+n_3-1}(\mathcal Z_{\widetilde\sK})$ different from
the component determined by the full subcomplex
$\widetilde\sK_{V_2\sqcup V_3}$ do not affect the nontriviality of
the Massey product, while the multigraded component corresponding
to $\widetilde\sK_{V_2\sqcup V_3}$ is zero since
$\widetilde\sK_{V_2\sqcup V_3}\cong S^{n_2+n_3-1}\setminus\pt$ is
contractible. The group $H^{m_1+m_2+n_1+n_2-1}(\mathcal
Z_{\widetilde\sK})$ is considered similarly. It follows that the
Massey product contains a unique nonzero element in its
multigraded component and therefore it is nontrivial.
\end{proof}

\begin{corollary}\label{nonformalmac}
Let $\widetilde\sK$ be a triangulated sphere obtained from another
sphere $\sK$ by applying two stellar subdivisions as described
above. Then the corresponding 2-connected moment-angle manifold
$\mathcal Z_{\widetilde\sK}$ is non-formal.
\end{corollary}

In the proof of Theorem~\ref{ntmas} the nontriviality of the
Massey product is established geometrically. A parallel argument
may be carried out algebraically using the Koszul complex or its
quotient algebra~$R^*(\sK)$, as illustrated in the next example.
To be precise, the nonformality of a manifold $\zk$ is equivalent
to the nonformality of its \emph{singular} cochain algebra
$C^*(\zk;\Q)$ (or Sullivan's algebra $\APL(\zk)$), while the
rational Koszul complex or the algebra $R^*(\sK)\otimes\Q$ are
quasi-isomorphic to the \emph{cellular} cochain algebra
$\sC^*(\zk;\Q)$. However, this difference is irrelevant: it can be
easily seen that the existence of a nontrivial triple Massey
product for the cellular cochain algebra $\sC^*(\zk;\Q)$ implies
that $\zk$ is nonformal (an exercise, or
see~\cite[Proposition~5.6.1]{de-su07}).

\begin{example}
We consider the simplest case of Construction~\ref{baskconstr},
when $\widetilde\sK$ is obtained from the boundary of an
octahedron by applying stellar subdivisions at two skew edges.
Then $\widetilde\sK=\sK_P$ is the nerve complex of a simple
3-polytope $P$ obtained by truncating a cube at two edges as shown
in Figure~\ref{masseyfig}. The face ring is given by
\[
  \Z[\sK_P]=\Z[v_1,\ldots,v_6,w_1,w_2]/\mathcal I_{\sK_P},
\]
where $v_i$, $i=1,\ldots,6$, are the generators coming from the
facets of the cube and $w_1,w_2$ are the generators corresponding
to the two new facets, and
\[
  \mathcal I_{\sK_P}=(v_1v_2,v_3v_4,v_5v_6,w_1w_2,v_1v_3,v_4v_5,
  w_1v_3,w_1v_6,w_2v_2,w_2v_4).
\]
\begin{figure}[h]
\psset{unit=0.8mm}
\begin{center}
%\vspace{4mm}
\begin{pspicture}(0,0)(60,55)
  \psline[linestyle=dashed](0,2.5)(15,10)
  \psline[linestyle=dashed](15,10)(15,55)
  \psline[linestyle=dashed](15,10)(60,10)
  \pspolygon(0,2.5)(0,47.5)(10,45)(10,0)
  \pspolygon(0,47.5)(15,55)(50,55)(30,45)(10,45)
  \pspolygon(50,55)(30,45)(40,35)(60,45)
  \pspolygon(40,35)(60,45)(60,10)(40,0)
  \psline(10,0)(40,0)
  \put(22,49){\small$v_1$}
  \put(25,-5){\vector(0,1){3}}
  \put(24,-7.5){\small$v_2$}
  \put(49,24){\small$v_3$}
  \put(-5,30){\vector(1,0){3}}
  \put(-9,29){\small$v_4$}
  \put(22,21){\small$v_5$}
  \put(65,30){\vector(-1,0){3}}
  \put(65.5,29){\small$v_6$}
  \put(3,24){\small$w_1$}
  \put(44,44){\small$w_2$}
\end{pspicture}
\vspace{4mm} \caption{}\label{masseyfig}
\end{center}
\end{figure}%
We denote the corresponding exterior generators of $R^*(\sK_P)$ by
$u_1,\ldots,u_6,t_1,t_2$; they satisfy $du_i=v_i$ and $dt_i=w_i$.
Consider the cocycles
\[
  a=v_1u_2,\quad b=v_3u_4,\quad c=v_5u_6
\]
and the corresponding cohomology classes $\alpha,\beta,\gamma\in
H^{-1,4}[R^*(\sK)]$. The equations
\[
  ab=de,\quad bc=df
\]
have a solution $e=0$, $f=v_5u_3u_4u_6$, so the triple Massey
product $\langle\alpha,\beta,\gamma\rangle\in H^{-4,12}[R^*(\sK)]$
is defined. This Massey product is represented by the cocycle
\[
  af+ec=v_1v_5u_2u_3u_4u_6
\]
and is nontrivial. The differential graded algebra $R^*(\sK_P)$
and the 11-dimensional manifold $\mathcal Z_{\sK_P}$ are not
formal.
\end{example}

\begin{remark}
We can truncate the polytope $P$ from the previous example at
another edge to obtain a 3-dimensional associahedron
$\mathit{As}^3$, shown in Figure~\ref{ascube} (left). By
considering similar nontrivial Massey products (now there will be
three of them, corresponding to each pair of cut off edges) we
deduce that the 12-dimensional moment-angle manifold corresponding
to $\mathit{As}^3$ is also nonformal.
\end{remark}

In view of Theorem~\ref{ntmas},
the question arises %%%%%a question
of describing the class of simplicial complexes $\sK$ for which
the algebra $R^*(\sK)$ (equivalently, the Koszul algebra
$(\Lambda[u_1,\ldots,u_m]\otimes\Z[\sK],d)$)
% or the space $\zk$)
is formal.
%(in particular, does not carry nontrivial Massey products).
For example, this is the case if $\sK$ is the boundary of a
polygon or, more generally, if $\sK$ is of the form described in
Theorem~\ref{zpstacked}.

Triple Massey products in the cohomology of $\zk$ were further
studied in the work of Denham and Suciu~\cite{de-su07}. According
to~\cite[Theorem~6.1.1]{de-su07}, there exists a nontrivial triple
Massey product of 3-dimensional cohomology classes
$\alpha,\beta,\gamma\in H^3(\zk)$ if and only if the 1-skeleton of
$\sK$ contains an induced subgraph isomorphic to one of the five
explicitly described `obstruction' graphs.
In~\cite[Example~8.5.1]{de-su07} there is also constructed an
example of~$\sK$ for which the corresponding $\zk$ has an
\emph{indecomposable}\label{indec3massey} triple Massey product in
the cohomology (a triple Massey product is indecomposable if it
does not contain a cohomology class that can be written as a
product of two cohomology classes of positive dimension).

\smallskip

To conclude this section, we mention that the algebraic study of
Massey products in the cohomology of Koszul complexes has a long
history. It goes back to the work of Golod~\cite{golo62}, who
studied the Poincar\'e series of $\Tor_R(\k,\k)$ for a N\"otherian
local ring~$R$. The main result of Golod is a calculation of the
Poincar\'e series for the class of rings whose Koszul complexes
have all Massey products vanishing (including the cohomology
multiplication). Such rings were called \emph{Golod} in the
monograph~\cite{gu-le69} of Gulliksen and Levin, where the reader
can find a detailed exposition of Golod's theorem together with
several further applications.

\begin{definition}\label{defgolod}
We refer to a simplicial complex $\sK$ as \emph{Golod} (over a
ring~$\k$) if its face ring $\k[\sK]$ has the Golod property, i.e.
if the multiplication and all higher Massey products in
$\Tor_{\k[v_1,\ldots,v_m]}(\k[\sK],\k)=H(\Lambda[u_1,\ldots,u_m]\otimes\k[\sK],d)$
are trivial.
\end{definition}

Golod complexes were studied in~\cite{h-r-w99}, where several
combinatorial criteria for Golodness were given. The appearance of
moment-angle complexes added a topological dimension to the whole
study. In particular, Theorem~\ref{zkcoh} implies that $\sK$ is
Golod whenever $\zk$ is homotopy equivalent to a wedge of spheres.
This observation was used in~\cite[Theorems~9.1,~11.2]{gr-th07} to
produce new classes of Golod simplicial complexes, including
skeleta of simplices considered in Theorem~\ref{gr-th}, and, more
generally, all \emph{shifted}\label{defshiftedco} complexes. For
all such $\sK$ the corresponding moment-angle complex $\zk$ is
homotopy equivalent to a wedge of spheres. There are examples of
Golod complexes $\sK$ for which $\zk$ is \emph{not} homotopy
equivalent to a wedge of spheres (see Exercise~\ref{rp2ng}
and~\cite[Example~3.3]{g-p-t-w}). More explicit series of Golod
complexes were constructed by Seyed Fakhari and Welker
in~\cite{se-we12}.

By a result of Berglund and
J\"ollenbek~\cite[Theorem~5.1]{be-jo07}, the face ring $\k[\sK]$
is Golod if and only if the multiplication in
$\Tor_{\k[v_1,\ldots,v_m]}(\k[\sK],\k)$ is trivial (i.e.
triviality of the cup-product implies that all higher Massey
products are also trivial).

More details on the relationship between the Golod property for
$\sK$ and the homotopy theory of moment-angle complexes $\zk$ can
be found in~\cite{gr-th07}, \cite{g-p-t-w}, as well as in
Section~\ref{homotflag}.

%
%The face rings of sphere triangulations $\widetilde\sK$ considered
%earlier in this section do not qualify for Golodness, as the
%corresponding moment-angle complexes $\mathcal Z_{\widetilde\sK}$
%are manifolds, and therefore, the cohomology of the Koszul complex
%of $\Q[\widetilde\sK]$ has many nontrivial products.
%The question of formality of these manifolds is related to
%vanishing of \emph{higher order} Massey products, while the
%Golodness is vanishing of \emph{all} Massey products, including
%the cohomology multiplication.

\subsection*{Exercises}

\begin{exercise}
If the cellular cochain algebra $\sC^*(\zk;\Q)$ carries a
nontrivial triple Massey product, then $\zk$ is nonformal.
\end{exercise}

\begin{exercise}\label{rp2ng}
Let $\sK$ be the triangulation of $\R P^2$ from
Example~\ref{exbettinum}.4. Then $\sK$ is a Golod complex, but
$\zk$ is not homotopy equivalent to a wedge of spheres.
\end{exercise}

\section{Moment-angle complexes from simplicial posets}\label{macsp}
Simplicial posets $\sS$ generalise naturally abstract simplicial
complexes (see Section~\ref{secsimpos}). Algebraic properties of
their face rings $\k[\sS]$ were discussed in
Section~\ref{frsimpo}. Following the categorical description of
the moment-angle complex $\zk$ outlined in the end of
Construction~\ref{constrmac}, it is easy to extend the definition
of $\zk$ to simplicial posets. The resulting space $\zs$ carries a
torus action, and its equivariant and ordinary cohomology is
expressed in terms of the face ring $\Z[\sS]$ in the same way as
for the standard moment-angle complexes~$\zk$. Simplicial posets
and associated moment-angle complexes therefore provide a broader
context for studying the link between torus actions and
combinatorial commutative algebra. These developments are
originally due to L\"u and Panov~\cite{lu-pa11}.

Let $\sS$ be a finite simplicial poset with the vertex set
$V(\sS)=[m]$.

\begin{construction}[moment-angle complex]\label{macsimpos}
We consider the face category $\mathop{\text{\sc cat}}(\mathcal
S)$ whose objects are elements $\sigma\in\sS$ and there is a
morphism from $\sigma$ to $\tau$ whenever $\sigma\le\tau$. For
each element $\sigma\in\sS$ we define the following subset in the
standard unit polydisc $\D^m\subset\C^m$:
\[
  (D^2,S^1)^\sigma=\bigl\{(z_1,\ldots,z_m)\in
  \D^m\colon |z_j|^2=1\text{ for }j\not\le\sigma\bigl\}.
\]
Then $(D^2,S^1)^\sigma$ is homeomorphic to a product of $|\sigma|$
discs and $m-|\sigma|$ circles. We have an inclusion
$(D^2,S^1)^\tau\subset(D^2,S^1)^\sigma$ whenever $\tau\le\sigma$.
Now define a diagram
\begin{align*}
  \mathcal D_\sS(D^2,S^1)\colon \ca(\sS)&\longrightarrow \top,\\
  \sigma&\longmapsto (D^2,S^1)^\sigma,
\end{align*}
which maps a morphism $\sigma\le\tau$ of $\ca(\sS)$ to the
inclusion $(D^2,S^1)^\sigma\subset(D^2,S^1)^\tau$ (see
Appendix~\ref{secmc} for the definition of diagrams and their
colimits).

We define the \emph{moment-angle complex}\label{colimitmac}
corresponding to $\sS$ by
\[
  \zs=\mathop{\mathrm{colim}}
  \mathcal D_\sS(D^2,S^1)=\mathop{\mathrm{colim}}_{\sigma\in\sS}
  (D^2,S^1)^\sigma.
\]
The space $\zs$ is therefore glued from the blocks
$(D^2,S^1)^\sigma$ according to the poset relation in~$\sS$. When
$\sS$ is (the face poset of) a simplicial complex $\sK$ it becomes
the standard moment-angle complex~$\zk$.

Since every subset $(D^2,S^1)^\sigma\subset\D^m$ is invariant with
respect to the coordinatewise action of the $m$-torus $\T^m$, the
moment-angle complex $\zs$ acquires a $\T^m$-action.

This definition extends to a set $(\mb X,\mb A)$ of $m$ pairs of
spaces (see Construction~\ref{nsc}), so we define the
\emph{polyhedral product}\label{polyprosp} of $(\mb X,\mb A)$
corresponding to $\sS$ by
\[
  (\mb X,\mb A)^\sS=\mathop{\mathrm{colim}}
  \mathcal D_\sS(\mb X,\mb A)=\mathop{\mathrm{colim}}_{\sigma\in\sS}
  (\mb X,\mb A)^\sigma.
\]
\end{construction}

The construction of the polyhedral product $(\mb X,\mb A)^\sS$ is
functorial in all arguments: there are straightforward analogues
of Propositions~\ref{mafunc} and~\ref{functmac1}, which are proved
in a similar way.

\begin{example}
Let $\sS$ be the simplicial poset of Fig.~\ref{sccfig}~(a). Then
$\zs$ is obtained by gluing two copies of $D^2\times D^2$ along
their boundary $S^3=D^2\times S^1\cup S^1\times D^2$. Therefore,
$\zs\cong S^4$. Here, $\sK_\sS=\varDelta^1$ (a segment), and the
moment-angle complex map induced by the map
$\sS\to\sK_{\sS}$~\eqref{fold} folds $S^4$ onto~$D^4$. Similarly,
if $\sS$ is of Fig.~\ref{sccfig}~(b), then $\zs\cong S^6$. Note
that even-dimensional spheres do not appear as moment-angle
complexes~$\zk$ for simplicial complexes~$\sK$.
\end{example}

The \emph{join}\label{joinsimpose} of simplicial posets $\sS_1$
and $\sS_2$ is the simplicial poset $\sS_1\mathbin{*}\sS_2$ whose
elements are pairs $(\sigma_1,\sigma_2)$, with
$(\sigma_1,\sigma_2)\le(\tau_1,\tau_2)$ whenever
$\sigma_1\le\tau_1$ in $\sS_1$ and $\sigma_2\le\tau_2$ in $\sS_2$.
The following properties of $\zs$ are similar to those of~$\zk$.

{\samepage
\begin{theorem}\label{zsman}\
\begin{itemize}
\item[(a)] $\mathcal Z_{\sS_1\mathbin{*}\sS_2}\cong
\mathcal Z_{\sS_1}\times\mathcal Z_{\sS_2}$;

\item[(b)] the quotient $\zs/\T^m$ is homeomorphic to the cone over
$|\sS|$;

\item[(c)] if $|\sS|\cong S^{n-1}$, then $\zs$ is a manifold of dimension~$m+n$.
\end{itemize}
\end{theorem}
}
\begin{proof}
Statements (a) and (b) are proved in the same way as the
corresponding statements for simplicial complexes, see
Section~\ref{defzk}. To prove~(c) we use the `dual' decomposition
of the boundary of the $n$-ball $\cone|\sS|$ into faces, in the
same way as in the proof of Theorem~\ref{zkman}.
\end{proof}

\begin{construction}[cell decomposition]
We proceed by analogy with the construction of
Section~\ref{celld}. The disc $\D$ is decomposed into 3 cells: the
point $1\in\D$ is the 0-cell; the complement to $1$ in the
boundary circle is the 1-cell, which we denote~$T$; and the
interior of $\D$ is the 2-cell, which we denote~$D$. The polydisc
$\D^m$ then acquires the product cell decomposition, with each
$(D^2,S^1)^\sigma\subset(D^2,S^1)^\tau$ being an inclusion of
cellular subcomplexes for $\sigma\le\tau$. We therefore obtain a
cell decomposition of $\zs$ (although this time it is not a
subcomplex in~$\D^m$ in general!). Each cell in $\zs$ is
determined by an element $\sigma\in\sS$ and a subset $\omega\in
V(\sS)$ with $V(\sigma)\cap\omega=\varnothing$. Such a cell is a
product of $|\sigma|$ cells of $D$-type, $|\omega|$ cells of
$T$-type and the rest of $1$-type. We denote this cell by
$\kappa(\omega,\sigma)$.

The resulting cellular cochain complex $\sC^*(\zs)$ has an
additive basis consisting of cochains $\kappa(\omega,\sigma)^*$
dual to the corresponding cells. We introduce a
$\Z\oplus\Z^m$-grading on the cochains by setting
\[
  \mathop{\mathrm{mdeg}}\kappa(\omega,\sigma)^*=(-|\omega|,2V(\sigma)+2\omega),
\]
where we think of both $V(\sigma)$ and $\omega$ as vectors
in~$\{0,1\}^m\subset\Z^m$. The cellular differential preserves the
$\Z^m$-part of the multigrading, so we obtain a decomposition
\[
  \sC^*(\zs)=\bigoplus_{\mb a\in\Z^m}\sC^{*,2\mb a}(\zk)
\]
into a sum of subcomplexes. The only nontrivial subcomplexes are
those for which $\mb a$ is in~$\{0,1\}^m$. The cellular cohomology
of $\zs$ thereby acquires a multigrading, and we define the
\emph{multigraded Betti numbers}\label{multibnsp} $b^{-i,2\mb
a}(\zs)$ by
\[
  b^{-i,2\mb a}(\zs)=\mathop{\mathrm{rank}} H^{-i,2\mb a}(\zs),
  \quad\text{for } 1\le i\le m,\;\mb a\in\Z^m.
\]
For the ordinary Betti numbers we have $b^k(\zs)=\sum_{-i+2|\mb
a|=k}b^{-i,2\mb a}(\zs)$.

The map of moment-angle complexes $\mathcal Z_{\sS_1}\to\mathcal
Z_{\mathcal S_2}$ induced by a simplicial poset map
$\sS_1\to\sS_2$ is clearly a cellular map, and therefore the
cellular cohomology is functorial with respect to such maps.
\end{construction}

We now recall from Section~\ref{frsimpo} that the face ring
$\Z[\sS]$ is a $\Z^m$-graded $\Z[v_1,\ldots,v_m]$-module via the
map sending each $v_i$ identically, and we have the
$\Z\oplus\Z^m$-graded $\Tor$-algebra of $\Z[\sS]$:
\[
  \Tor_{\Z[v_1,\ldots,v_m]}\bigl(\Z[\sS],\Z\bigr)=
  \bigoplus_{i\ge0,\mb a\in\N^m}
  \Tor^{-i,\,2\mb a}_{\Z[v_1,\ldots,v_m]}\bigl(\Z[\sS],\Z\bigr).
\]

\begin{theorem}\label{hzs1}
There is a graded ring isomorphism
\[
  H^*(\zs)\cong\Tor_{\Z[v_1,\ldots,v_m]}(\Z[\sS],\Z)
\]
whose graded components are given by the group isomorphisms
\begin{equation}\label{additive}
  H^p(\zs)\cong\bigoplus_{-i+2|\mb a|=p}
  \Tor^{-i,\,2\mb a}_{\Z[v_1,\ldots,v_m]}(\Z[\sS],\Z)
\end{equation}
in each degree $p$. Here $|\mb a|=j_1+\cdots+j_m$ for $\mb
a=(j_1,\ldots,j_m)$.
\end{theorem}

Using the Koszul complex we restate the above theorem as follows:

\begin{theorem}\label{hzs2}
There is a graded ring isomorphism
\[
  H^*(\zs)\cong
  H\bigl(\Lambda[u_1,\ldots,u_m]\otimes\Z[\sS],d\bigr),
\]
where the $\Z\oplus\Z^m$-grading and the differential on the right
hand side are defined by
\[
  \mathop{\mathrm{mdeg}}u_i=(-1,2\mb e_i),\quad\mathop{\mathrm{mdeg}}v_\sigma=(0,2V(\sigma)),
  \qquad du_i=v_i, \quad dv_\sigma=0,
\]
and $\mb e_i\in\Z^m$ is the $i$th basis vector, for
$i=1,\ldots,m$.
\end{theorem}
\begin{proof}
The proof follows the lines of the proof of Theorem~\ref{zkcoh},
but the analogues of Lemmata~\ref{iscoh} and~\ref{cellappr} are
proved in a different way.

We first set up the quotient differential graded ring
\[
  R^*(\sS)=\Lambda[u_1,\ldots,u_m]\otimes\Z[\sS]/\mathcal I_R
\]
where $\mathcal I_R$ is the ideal generated by the elements
\[
  u_iv_\sigma\quad\text{with }i\in V(\sigma),\quad\text{and}\quad
  v_\sigma v_\tau\;\;\;\text{with }\sigma\wedge\tau\ne\hatzero.
\]
Note that the latter condition is equivalent to $V(\sigma)\cap
V(\tau)\ne\varnothing$. The ring $R^*(\sS)$ will serve as an
algebraic model for the cellular cochains of~$\zs$.

We need to prove an analogue of Lemma~\ref{iscoh}, i.e. show that
the quotient map
\[
  \varrho\colon\Lambda[u_1,\ldots,u_m]\otimes\Z[\sS]\to R^*(\sS)
\]
induces an isomorphism in cohomology. Instead of constructing a
chain homotopy directly, we shall identify both $R^*(\sS)$ and
$\Lambda[u_1,\ldots,u_m]\otimes\Z[\sS]$ with the cellular cochains
of homotopy equivalent spaces.

Theorem~\ref{chdec} implies that $R^*(\sS)$ has basis of monomials
$u_\omega v_\sigma$ where $\omega\subset V(\sS)$, \
$\sigma\in\sS$, \ $\omega\cap V(\sigma)=\varnothing$, and
$u_\omega=u_{i_1}\ldots u_{i_k}$ for $\omega=\{i_1,\ldots,i_k\}$.
In particular, $R^*(\sS)$ is a free abelian group of finite rank.
The map
\begin{equation}\label{gmap}
  g\colon R^*(\sS) \to \sC^*(\zs),\quad
  u_\omega v_\sigma \mapsto\mathcal \kappa(\omega,\sigma)^*
\end{equation}
is an isomorphism of cochain complexes. Indeed, the additive bases
of the two groups are in one-to-one correspondence, and the
differential in $R^*(\sS)$ acts (in the case $|\omega|=1$ and
$i\notin V(\sigma)$) as
\[
  d(u_iv_\sigma)=v_iv_\sigma=\sum_{\eta\in i\vee\sigma}v_\eta.
\]
This is exactly how the cellular differential in $\sC^*(\zs)$ acts
on $\kappa(i,\sigma)^*$. The case of arbitrary $\omega$ is treated
similarly. It follows that we have an isomorphism of cohomology
groups $H^j[R^*(\sS)]\cong H^j(\zs)$ for all~$j$.

The differential algebra
$(\Lambda[u_1,\ldots,u_m]\otimes\Z[\sS],d)$ also may be identified
with the cellular cochains of a certain space. Namely, consider
the polyhedral product $(S^\infty,S^1)^\sS$. Then we show in the
same way as in Subsection~\ref{algmodsc} that there is an
isomorphism of cochain complexes
\[
  g'\colon\Lambda[u_1,\ldots,u_m]\otimes\Z[\sS]\to\sC^*\bigl((S^\infty,S^1)^\sS\bigr).
\]
Furthermore, the standard functoriality arguments give a
deformation retraction
\[
  \zs=(D^2,S^1)^{\mathcal S}\hookrightarrow(S^\infty,S^1)^{\mathcal S}
  \longrightarrow(D^2,S^1)^\sS
\]
onto a cellular subcomplex. Hence the cochain map
$\sC^*((S^\infty,S^1)^{\mathcal S})\to\sC^*(\zs)$ corresponding to
the inclusion $\zs\hookrightarrow(S^\infty,S^1)^{\mathcal S}$
induces an isomorphism in cohomology.

Summarising the above observations we obtain the commutative
square
\begin{equation}\label{comsq}
\begin{CD}
  \Lambda[u_1,\ldots,u_m]\otimes\Z[\sS] @>g'>>
  \sC^*((S^\infty,S^1)^\sS)\\
  @V\varrho VV @VVV\\
  R^*(\sS) @>g>> \sC^*(\zs)
\end{CD}
\end{equation}
in which the horizontal arrows are isomorphisms of cochain
complexes, and the right vertical arrow induces an isomorphism in
cohomology. It follows that the left arrow also induces an
isomorphism in cohomology, as claimed.

The additive isomorphism of~\eqref{additive} now follows
from~\eqref{comsq}. To establish the ring isomorphism we need to
analyse the multiplication of cellular cochains.
% in~$\sC^*(\zs)$.

We consider the diagonal approximation map
$\widetilde{\Delta}\colon\D^m\to\D^m\times\D^m$ given on each
coordinate by~\eqref{Ddiagmap}. It restricts to a map
$(D^2,S^1)^\sigma\to(D^2,S^1)^\sigma\times(D^2,S^1)^\sigma$ for
every $\sigma\in\sS$ and gives rise to a map of diagrams
\[
  \mathcal D_\sS(D^2,S^1)\to\mathcal D_\sS(D^2,S^1)\times\mathcal D_\sS(D^2,S^1).
\]
By definition, the colimit of the latter is $\mathcal
Z_{\sS\mathop{*}\sS}$, which is identified with $\zs\times\zs$. We
therefore obtain a cellular approximation
$\widetilde\Delta\colon\zs\to\zs\times\zs$ for the diagonal map
of~$\zs$. It induces a ring structure on the cellular cochains via
the composition
\[
\begin{CD}
  \sC^*(\zs)\otimes\sC^*(\zs) @>\times>>\sC^*(\zs\times\zs)
  @>\widetilde{\Delta}^*>>\sC^*(\zs).
\end{CD}
\]
We claim that, with this multiplication in $\sC^*(\zs)$, the
map~\eqref{gmap} becomes a ring isomorphism. To see this we first
observe that since~\eqref{gmap} is a linear map, it is enough to
consider the product of two generators $u_\omega v_\sigma$ and
$u_\psi v_\tau$. If any two of the subsets $\omega$, $V(\sigma)$,
$\psi$ and $V(\tau)$ have nonempty intersection, then $u_\omega
v_\sigma\cdot u_\psi v_\tau=0$ in $R^*(\sS)$. Otherwise (if all
four subsets are disjoint) we have
\begin{equation}\label{gprod}
  g(u_\omega v_\sigma\cdot u_\psi v_\tau)=
  g\bigl(u_{\omega\sqcup\psi}\cdot
  \sum_{\eta\in\,\sigma\vee\tau}v_\eta\bigr)=
  \sum_{\eta\in\,\sigma\vee\tau}\kappa(\omega\sqcup\psi,\eta)^*.
\end{equation}
We also observe that for any cell $\kappa(\chi,\eta)$ of $\zs$
(with $\chi\cap V(\eta)=\varnothing$) we have
\[
  \widetilde\Delta\kappa(\chi,\eta)=
  \mathop{\sum_{\omega\sqcup\psi=\chi}}\limits_{\sigma\vee\tau\,\ni\;\eta}
  \kappa(\omega,\sigma)\times\kappa(\psi,\tau).
\]
Therefore,
\begin{multline*}
  g(u_\omega v_\sigma)\cdot g(u_\psi v_\tau)=
  \kappa(\omega,\sigma)^*\cdot\kappa(\psi,\tau)^*\\=
  \widetilde\Delta^*
  \bigl(\kappa(\omega,\sigma)\times\kappa(\psi,\tau)\bigr)^*=
  \sum_{\eta\,\in\,\sigma\vee\tau}\kappa(\omega\sqcup\psi,\eta)^*.
\end{multline*}
Comparing this with~\eqref{gprod} we deduce that~\eqref{gmap} is a
ring map.
\end{proof}

\begin{remark}
Using the monoid structure on $\D$ as in
Proposition~\ref{functmac1} one easily sees that the construction
of $\zs$ is functorial with respect to maps of simplicial posets.
This together with Proposition~\ref{funct} makes the isomorphism
of Theorem~\ref{hzs1} functorial.
\end{remark}

Using Hochster's formula for simplicial posets
(Theorem~\ref{hochsp}) we can calculate the cohomology of $\zs$
via the cohomology of full subposets $\sS_J\subset\sS$. Here is an
example of calculation using this method.

\begin{example}\label{exmwc}
Let $\sS$ be the simplicial poset shown on Fig.~\ref{mwc}~(right).
It has $m=5$ vertices, 9 edges and 6 triangular faces. The dual
3-dimensional `ball with corners' $Q$ (see the proof of
Theorem~\ref{zsman}) is shown in Figure~\ref{mwc}~(left). We
denote its facets $F_1,\ldots,F_5$, edges $e,f,g$ and the vertex
$\sigma=F_3\cap f$ as shown. The corresponding moment-angle
complex $\zs$ is an 8-dimensional manifold.
\begin{figure}
\begin{picture}(0,0)%
\includegraphics{sposet2.pstex}%
\end{picture}%
\setlength{\unitlength}{1776sp}%
\begingroup\makeatletter\ifx\SetFigFont\undefined%
\gdef\SetFigFont#1#2#3#4#5{%
  \reset@font\fontsize{#1}{#2pt}%
  \fontfamily{#3}\fontseries{#4}\fontshape{#5}%
  \selectfont}%
\fi\endgroup%
\begin{picture}(12648,5550)(8018,-6769)
\put(18151,-6661){\makebox(0,0)[lb]{\smash{{\SetFigFont{11}{13.2}{\rmdefault}{\mddefault}{\updefault}$\mathcal{S}$}}}}
\put(10201,-5986){\makebox(0,0)[lb]{\smash{{\SetFigFont{8}{9.6}{\rmdefault}{\mddefault}{\updefault}$F_2$}}}}
\put(12526,-2761){\makebox(0,0)[lb]{\smash{{\SetFigFont{8}{9.6}{\rmdefault}{\mddefault}{\updefault}$e$}}}}
\put(10501,-2311){\makebox(0,0)[lb]{\smash{{\SetFigFont{8}{9.6}{\rmdefault}{\mddefault}{\updefault}$F_4$}}}}
\put(8626,-4561){\makebox(0,0)[lb]{\smash{{\SetFigFont{8}{9.6}{\rmdefault}{\mddefault}{\updefault}$F_3$}}}}
\put(12376,-4411){\makebox(0,0)[lb]{\smash{{\SetFigFont{8}{9.6}{\rmdefault}{\mddefault}{\updefault}$F_5$}}}}
\put(20251,-3436){\makebox(0,0)[lb]{\smash{{\SetFigFont{8}{9.6}{\rmdefault}{\mddefault}{\updefault}$4$}}}}
\put(18226,-5536){\makebox(0,0)[lb]{\smash{{\SetFigFont{8}{9.6}{\rmdefault}{\mddefault}{\updefault}$2$}}}}
\put(18451,-3436){\makebox(0,0)[lb]{\smash{{\SetFigFont{8}{9.6}{\rmdefault}{\mddefault}{\updefault}$5$}}}}
\put(18226,-1411){\makebox(0,0)[lb]{\smash{{\SetFigFont{8}{9.6}{\rmdefault}{\mddefault}{\updefault}$1$}}}}
\put(8551,-2761){\makebox(0,0)[lb]{\smash{{\SetFigFont{8}{9.6}{\rmdefault}{\mddefault}{\updefault}$g$}}}}
\put(16126,-3436){\makebox(0,0)[lb]{\smash{{\SetFigFont{8}{9.6}{\rmdefault}{\mddefault}{\updefault}$3$}}}}
\put(10426,-3586){\makebox(0,0)[lb]{\smash{{\SetFigFont{8}{9.6}{\rmdefault}{\mddefault}{\updefault}$F_1$}}}}
\put(9301,-5461){\makebox(0,0)[lb]{\smash{{\SetFigFont{8}{9.6}{\rmdefault}{\mddefault}{\updefault}$\sigma$}}}}
\put(10501,-5011){\makebox(0,0)[lb]{\smash{{\SetFigFont{8}{9.6}{\rmdefault}{\mddefault}{\updefault}$f$}}}}
\put(10351,-6661){\makebox(0,0)[lb]{\smash{{\SetFigFont{11}{13.2}{\rmdefault}{\mddefault}{\updefault}$Q$}}}}
\end{picture}%
\centering
  \caption{`Ball with corners' $Q$ dual to the simplicial poset~$\sS$.}
\label{mwc}
\end{figure}

The face ring $\Z[\sS]$ is the quotient of the polynomial ring
\[
  \Z[\sS]=\Z[v_1,\ldots,v_5,v_e,v_f,v_g],\quad\deg v_i=2,\quad
  \deg v_e=\deg v_f=\deg v_e=4
\]
by the relations
\begin{align*}
  v_1v_2&=v_e+v_f+v_g,\\
  v_3v_4&=v_3v_5=v_4v_5=
  v_3v_e=v_4v_f=v_5v_g=v_ev_f=v_ev_g=v_ev_f=0.
\end{align*}
The other generators and relations in the original presentation of
$\Z[\sS]$ can be derived from the above; e.g., $v_\sigma=v_3v_f$.

Given $J\subset[m]$ we define the following subset in the boundary
of~$Q$:
\[
  Q_J=\bigcup_{j\in J}F_j\subset Q.
\]
By analogy with Proposition~\ref{hochpol}
%(using the fact that $|\sS_J|$ is a deformation retract of~$Q_J$)
we prove that
\begin{equation}\label{faces}
  H^{-i,\,2J}(\zs)\cong
  \widetilde H^{|J|-i-1}(Q_J).
\end{equation}
Using this formula we calculate the nontrivial cohomology groups
of $\zs$ as follows:
\begin{align*}
H^{0,(0,0,0,0,0)}(\zs)&=\widetilde H^{-1}(\varnothing)=\Z
 &&1\\
H^{-1,(0,0,2,2,0)}(\zs)&=\widetilde H^{0}(F_3\cup F_4)=\Z
 &&u_3v_4\\
H^{-1,(0,0,2,0,2)}(\zs)&=\widetilde H^{0}(F_3\cup F_5)=\Z
 &&u_5v_3\\
H^{-1,(0,0,0,2,2)}(\zs)&=\widetilde H^{0}(F_4\cup F_5)=\Z
 &&u_4v_5\\
H^{-2,(0,0,2,2,2)}(\zs)&=\widetilde H^{0}(F_3\cup F_4\cup
F_5)=\Z\oplus\Z
 &&u_5u_3v_4,\;u_5u_4v_3\\
H^{0,(2,2,0,0,0)}(\zs)&=\widetilde H^{1}(F_1\cup F_2)=\Z\oplus\Z
 &&v_e,\;v_f\\
H^{-1,(2,2,2,0,0)}(\zs)&=\widetilde H^{1}(F_1\cup F_2\cup F_3)=\Z
 &&u_3v_e\\
H^{-1,(2,2,0,2,0)}(\zs)&=\widetilde H^{1}(F_1\cup F_2\cup F_4)=\Z
 &&u_4v_f\\
H^{-1,(2,2,0,0,2)}(\zs)&=\widetilde H^{1}(F_1\cup F_2\cup F_5)=\Z
 &&u_5v_g\\
H^{-2,(2,2,2,2,2)}(\zs)&=\widetilde H^{2}(F_1\cup\cdots\cup
F_5)=\Z && u_5u_4v_3v_f=u_5u_4v_\sigma
\end{align*}
It follows that the ordinary (1-graded) Betti numbers of $\zs$ are
given by the sequence $(1,0,0,3,4,3,0,0,1)$. In the right column
of the table above we include the cocycles in the differential
graded ring $\Lambda[u_1,\ldots,u_5]\otimes\Z[\sS]$ representing
generators of the corresponding cohomology group. This allows us
to determine the ring structure in $H^*(\zs)$. For example,
\[
  [u_5u_3v_4]\cdot[v_f]=[u_5u_3v_4v_f]=0=[u_5u_4v_3]\cdot[v_e].
\]
On the other hand,
\begin{multline*}
  [u_5u_3v_4]\cdot[v_e]=-[u_3u_5v_4v_e]=-[u_3u_4v_5v_e]
  \\=[u_3u_4v_5v_f]=[u_5u_4v_3v_f]=[u_5u_4v_3]\cdot[v_f].
\end{multline*}
Here we have used the relations
$d(u_3u_4u_5v_e)=u_3u_4v_5v_e-u_3u_5v_4v_e$ and
$d(u_1u_3u_4v_2v_5)=u_3u_4v_5v_e+u_3u_4v_5v_f$. All nontrivial
products come from Poincar\'e duality. These calculations are
summarised by the cohomology ring isomorphism
\[
  H^*(\zs)\cong H^*\bigl((S^3\times S^5)^{\#3}\cs
  (S^4\times S^4)^{\#2}\bigr)
\]
where the manifold on the right hand side is the connected sum of
three copies of $S^3\times S^5$ and two copies of $S^4\times S^4$.
We expect that the isomorphism above is induced by a
homeomorphism.
%; one might be able to prove this by using the
%surgery techniques of~\cite{gi-lo13}.
\end{example}

\subsection*{Exercises}
\begin{exercise}
Generalise Proposition~\ref{homsrs} to simplicial posets, i.e.
establish a ring isomorphism
%\[
$H^*\bigl((\C P^\infty,pt)^\sS\bigr)\cong\Z[\sS]$.
%\]
\end{exercise}

\begin{exercise}\label{evcohzs}
Construct a homotopy equivalence
\[
  h\colon(\C P^\infty,pt)^\sS\stackrel{\simeq}\longrightarrow
  E\T^m\times_{\T^m}\zs
\]
by extending the argument of Theorem~\ref{zkhofib}, and deduce
that $H^*_{\T^m}(\zs)\cong\Z[\sS]$.
\end{exercise}

%Сослаться на обзор Винберга-Попова про алгебраические факторпространства
%Сослаться на "Topological toric manifolds" в связи с
%   хаусдорфовостью факторпространства в конструкции Кокса
%В Example 5.4.3 добавить ссылку на компактный пример S^3/S^1 для двух разных действий
%Симплектоморфность алгебраической и гамильтоновой симплектических структур
%Формулировка и доказательство теоремы Дельзанта
%Раздутия торических многообразий и срезки граней
%Сослаться на книгу Аржанцева

\setcounter{chapter}4
\chapter{Toric varieties and manifolds}\label{toric}
A toric variety is an algebraic variety on which an
\emph{algebraic torus}~$(\C^\times)^n$ acts with a dense (Zariski
open) orbit. An algebraic torus contains a (compact) torus $T^n$,
so toric varieties are toric spaces in our usual sense. Toric
varieties are described by combinatorial-geometric objects,
rational fans (see Section~\ref{combfan}), and the combinatorics
of the fan determines the orbit structure of the torus action.

Toric varieties were introduced in 1970 in the pioneering work of
Demazure on \emph{Cremona group}. The geometry of toric varieties,
or \emph{toric geometry}, very quickly became one of the most
fascinating topics in algebraic geometry and found applications in
many other mathematical sciences, sometimes distant from each
other. We have already mentioned the proof for the necessity part
of the $g$-theorem for simplicial polytopes given by Stanley.
Other remarkable applications include counting lattice points and
volumes of lattice polytopes; relations with Newton polytopes and
the number of solutions of a system of algebraic equations (after
Khovanskii and Kushnirenko); discriminants, resultants and
hypergeometric functions (after Gelfand, Kapranov and Zelevinsky);
reflexive polytopes and mirror symmetry for Calabi--Yau toric
hypersurfaces and complete intersections (after Batyrev). Standard
references on toric geometry include Danilov's
survey~\cite{dani78} and books by Oda~\cite{oda88},
Fulton~\cite{fult93} and Ewald~\cite{ewal96}. The most recent
exhaustive account by Cox, Little and Schenck~\cite{c-l-s11}
covers many new applications, including those mentioned above.
Without attempting to give another review of toric geometry, in
this chapter we collect the basic definitions and constructions,
and emphasise topological and combinatorial aspects of toric
varieties.

We review the three main approaches to toric varieties in the
appropriate sections: the `classical' construction via fans, the
`algebraic quotient' construction, and the `symplectic reduction'
construction. The intersection of Hermitian quadrics appearing in
the symplectic construction of toric varieties links toric
geometry to moment-angle complexes. This link will be developed
further in the next chapter.
%In the last section we consider the topological properties which
%will be taken as the base for different topological
%generalisations of toric varieties in the later chapters.

A basic knowledge of algebraic geometry would much help the reader
of this chapter, although it is not absolutely necessary.

\section{Classical construction from rational
fans}\label{toricfan}

An \emph{algebraic torus} is a commutative complex algebraic group
isomorphic to a product $(\C^\times)^n$ of copies of the
multiplicative group $\C^\times=\C\setminus\{0\}$. It contains a
compact torus $T^n$ as a Lie (but not algebraic) subgroup.

We shall often identify an algebraic torus with the standard
model~$(\C^\times)^n$.

\begin{definition}\label{toric2}
A \emph{toric variety} is a normal complex algebraic variety~$V$
containing an algebraic torus $(\C^\times)^n$ as a Zariski open
subset in such a way that the natural action of $(\C^\times)^n$ on
itself extends to an action on~$V$.
\end{definition}

It follows that $(\C^\times)^n$ acts on $V$ with a dense orbit.
Toric varieties originally appeared as equivariant
compactifications of an algebraic torus, although non-compact
(e.g., affine) examples are now of equal importance.

\begin{example}\label{tvcpn}
The algebraic torus $(\C^\times)^n$ and the affine space $\C^n$
are the simplest examples of toric varieties. A compact example is
given by the projective space $\C P^n$ on which the torus acts in
homogeneous coordinates as follows:
\[
  (t_1,\ldots,t_n)\cdot(z_0:z_1:\ldots:z_n)=
  (z_0:t_1z_1:\ldots:t_nz_n).
\]
\end{example}

Algebraic geometry of toric varieties is translated completely
into the language of combinatorial and convex geometry. Namely,
there is a bijective correspondence between rational fans in
$n$-dimensional space (see Section~\ref{combfan}) and complex
$n$-dimensional toric varieties. Under this correspondence,
\begin{align*}
\text{cones }&\longleftrightarrow\text{ affine varieties}\\
\text{complete fans }&\longleftrightarrow\text{ compact (complete) varieties}\\
\text{normal fans of polytopes}&\longleftrightarrow\text{ projective varieties}\\
\text{regular fans }&\longleftrightarrow\text{ nonsingular varieties}\\
\text{simplicial fans }&\longleftrightarrow\text{ orbifolds}
\end{align*}
We review this construction below; the details can be found in the
sources mentioned above. Following the algebraic tradition, we use
the coordinate-free notation.

We fix a lattice $N$ of rank~$n$ (isomorphic to~$\Z^n$), and
denote by $N_\R$ its ambient $n$-dimensional real vector space
$N\otimes_\Z\R\cong\R^n$. Define the algebraic torus
$\C^\times_N=N\otimes_\Z\C^\times\cong(\C^\times)^n$. All cones
and fans in this chapter are rational.

\begin{construction}\label{tvff}
We first describe how to assign an affine toric variety to a cone
$\sigma\subset N_\R$. Consider the dual cone
$\sigma^\mathsf{v}\subset N_\R^*$ (see~\eqref{dualcone}) and
denote by
\[
  S_\sigma=\sigma^\mathsf{v}\cap N^*
\]
the set of its lattice points. Then $S_\sigma$ is a finitely
generated semigroup (with respect to addition). Let
$A_\sigma=\C[S_\sigma]$ be the semigroup ring of~$S_\sigma$. It is
a commutative finitely generated $\C$-algebra, with a $\C$-vector
space basis $\{\chi^{\mb u}\colon \mb u\in S_\sigma\}$. The
multiplication in $A_\sigma$ is defined via the addition
in~$S_\sigma$:
$$
  \chi^{\mb u}\cdot\chi^{\mb u'}=\chi^{\mb u+\mb u'},
$$
so $\chi^0$ is the unit. The \emph{affine toric
variety}\label{affinetv} $V_\sigma$ corresponding to~$\sigma$ is
the affine algebraic variety corresponding to $A_\sigma$:
$$
  V_\sigma=\mathop{\mathrm{Spec}}(A_\sigma),\quad A_\sigma=\C[V_\sigma].
$$
By choosing a multiplicative generator set in $A_\sigma$ we
represent it as a quotient
$$
  A_\sigma=\C[x_1,\ldots,x_r]/\mathcal I;
$$
then the variety $V_\sigma$ is the common zero set of polynomials
from the ideal~$\mathcal I$. Each point of $V_\sigma$ corresponds
to a semigroup homomorphism $\Hom_\mathrm{sg}(S_\sigma,\C_\mathrm
m)$, where $\C_\mathrm m=\C^\times\cup\{0\}$ is the multiplicative
semigroup of complex numbers.

Now if $\tau$ is a face of~$\sigma$, then
$\sigma^\mathsf{v}\subset\tau^\mathsf{v}$, and the inclusion of
semigroup algebras $\C[S_\sigma]\to\C[S_\tau]$ induces a morphism
$V_\tau\to V_\sigma$, which is an inclusion of a Zariski open
subset. This allows us to glue the affine varieties $V_\sigma$
corresponding to all cones $\sigma$ in a fan $\Sigma$ into an
algebraic variety~$V_\Sigma$, referred to as the \emph{toric
variety} corresponding to the fan~$\Sigma$. More formally,
$V_\Sigma$ may be defined as the \emph{colimit}\label{colimittv}
of algebraic varieties $V_\sigma$ over the partially ordered set
of cones of~$\Sigma$:
$$
  V_\Sigma=\mathop{\mathrm{colim}}\limits_{\sigma\in\Sigma}V_\sigma.
$$
Here is the crucial point: the fact that the cones $\sigma$ patch
into a fan $\Sigma$ guarantees that the variety $V_\Sigma$
obtained by gluing the pieces $V_\sigma$ is Hausdorff in the usual
topology. In algebraic geometry, the Hausdorffness is replaced by
the related notion of separatedness: a variety $V$ is
\emph{separated} if the image of the diagonal map $\varDelta\colon
V\to V\times V$ is Zariski closed. A separated variety is
Hausdorff in the usual topology.

\begin{lemma}\label{sepvar}
If a collection of cones $\{\sigma\}$ forms a fan $\Sigma$, then
the variety
$V_\Sigma=\mathop{\mathrm{colim}}_{\sigma\in\Sigma}V_\sigma$ is
separated.
\end{lemma}
\begin{proof}
Using the separatedness criterion of~\cite[Ch.~V,~\S4.3]{shaf94}
(see also~\cite[Proposition~5.4]{dani78}), it is enough to verify
the following: if cones $\sigma$ and $\sigma'$ intersect in a
common face $\tau$, then the diagonal map $V_\tau\to
V_\sigma\times V_{\sigma'}$ is a closed embedding. This is
equivalent to the assertion that the natural homomorphism
$A_\sigma\otimes A_{\sigma'}\to A_\tau$ is surjective. To prove
this, we use Separation Lemma (Lemma~\ref{seplemma}). According to
it, there is a linear function $\mb u$ which is nonnegative on
$\sigma$, nonpositive on $\sigma'$ and the intersection of the
hyperplane $\mb u^\perp$ with $\sigma$ is $\tau$. Now take $\mb
u'\in\tau^\mathsf{v}$, i.e. $\mb u'$ is nonnegative on~$\tau$.
Then there is an integer $k\ge0$ such that $\mb u'+k\mb u$ is
nonnegative on~$\sigma$, i.e. $\mb u'+k\mb u\in\sigma^\mathsf{v}$.
Then $\mb u'=(\mb u'+k\mb u)+(-k\mb
u)\in\sigma^\mathsf{v}+{\sigma'}{}^\mathsf{v}$. It follows that
$S_\sigma\oplus S_{\sigma'}\to S_\tau$ is surjective map of vector
spaces (or semigroups), hence $A_\sigma\otimes A_{\sigma'}\to
A_\tau$ is a surjective homomorphism.
\end{proof}

We shall consider only separated varieties in what follows.

The variety $V_\sigma$ carries an algebraic action of the torus
$\C^\times_N=N\otimes_\Z\C^\times$
\begin{equation}\label{atact}
  \C^\times_N\times V_\sigma\to V_\sigma,\quad (\mb t,x)\mapsto \mb t\cdot x
\end{equation}
which is defined as follows. A point $\mb t\in\C^\times_N$ is
determined by a group homomorphism $N^*\to\C^\times$. In
coordinates, the homomorphism $\Z^n\cong N^*\to\C^\times$
corresponding to $\mb t=(t_1,\ldots,t_n)$ is given by
\[
  \mb u=(u_1,\ldots,u_n)\mapsto \mb t(\mb u)=(t_1^{u_1}\cdots
  t_n^{u_n}).
\]
A point $x\in V_\sigma$ corresponds to a semigroup homomorphism
$S_\sigma\to\C_\mathrm m$. Then we define $\mb t\cdot x$ as the
point in $V_\sigma$ corresponding to the semigroup homomorphism
$S_\sigma\to\C_\mathrm m$ given by
\[
  \mb u\mapsto t(\mb u)x(\mb u).
\]
The homomorphism of algebras $A_\sigma\to A_\sigma\otimes\C[N^*]$
dual to the action~\eqref{atact} maps $\chi^\mb u$ to $\chi^\mb
u\otimes\chi^\mb u$ for $\mb u\in S_\sigma$. If $\sigma=\{\mathbf
0\}$, then we obtain the multiplication in the algebraic group
$\C^\times_N$. The actions on the varieties $V_\sigma$ are
compatible with the inclusions of open sets $V_{\tau}\to
V_{\sigma}$ corresponding to the inclusions of faces
$\tau\subset\sigma$. Therefore, for each fan $\Sigma$ we obtain a
$\C^\times_N$-action on the variety~$V_\Sigma$, which extends the
$\C^\times_N$-action on itself.
\end{construction}

\begin{example}\label{affns}
Let $N=\Z^n$ and let $\sigma$ be the cone spanned by the first $k$
basis vectors $\mb e_1,\ldots,\mb e_k$, where $0\le k\le n$. The
semigroup $S_\sigma=\sigma^\mathsf{v}\cap N^*$ is generated by the
dual elements $\mb e_1^*,\ldots,\mb e_k^*$ and $\pm \mb
e_{k+1}^*,\ldots,\pm \mb e_n^*$. Therefore,
$$
  A_\sigma\cong\C[x_1,\ldots,x_k,x_{k+1},x_{k+1}^{-1},\ldots,x_n,x_n^{-1}],
$$
where we set $x_i=\chi^{\mb e_i^*}$. It follows that the
corresponding affine variety is
$$
  V_\sigma\cong\C\times\cdots\times\C\times\C^\times\times\cdots\times\C^\times=
  \C^k\times(\C^\times)^{n-k}.
$$
In particular, for $k=n$ we obtain an $n$-dimensional affine
space, and for $k=0$ (i.e. $\sigma=\{0\}$) we obtain the algebraic
torus~$(\C^\times)^n$.
\end{example}

\begin{example}\label{econe}
Let $\sigma\subset\R^2$ be the cone generated by the vectors $\mb
e_2$ and $2\mb e_1-\mb e_2$ (note that these two vectors do not
span $\Z^2$, so this cone is not regular). The dual cone
$\sigma^\mathsf{v}$ is generated by $\mb e_1^*$ and $\mb
e_1^*+2\mb e_2^*$. The semigroup $S_\sigma$ is generated by $\mb
e_1^*$, $\mb e_1^*+\mb e_2^*$ and $\mb e_1^*+2\mb e_2^*$, with one
relation among them. Therefore,
$$
  A_\sigma=\C[x,xy,xy^2]\cong\C[u,v,w]/(v^2-uw)
$$
and $V_\sigma$ is a quadratic cone (a singular variety).
\end{example}

\begin{figure}[h]
\begin{center}
\begin{picture}(100,35)
  \put(20,15){\line(1,0){20}}
  \put(20,15){\line(0,1){20}}
  \put(20,15){\line(-1,-1){15}}
  \put(20,15){\vector(1,0){10}}
  \put(20,15){\vector(0,1){10}}
  \put(20,15){\vector(-1,-1){10}}
  \put(28,17){$\mb e_1$}
  \put(21,24){$\mb e_2$}
  \put(12,5){$-\mb e_1-\mb e_2$}
  \put(21,-2){(a)}
  \put(80,15){\line(1,0){20}}
  \put(80,15){\line(0,1){20}}
  \put(80,15){\line(0,-1){15}}
  \put(80,15){\vector(1,0){10}}
  \put(80,15){\vector(0,1){10}}
  \put(80,15){\vector(-1,2){10}}
  \put(80,15){\vector(0,-1){10}}
  \put(88,17){$\mb e_1$}
  \put(81,24){$\mb e_2$}
  \put(81,5){$-\mb e_2$}
  \put(53,33){$-\mb e_1+k\mb e_2$}
  \put(83,-2){(b)}
\end{picture}%
\caption{Complete fans in~$\R^2$} \label{compfans}
\end{center}
\end{figure}

\begin{example}\label{cp2fa}
Let $\Sigma$ be the complete fan in $\R^2$ with the following
three maximal cones: the cone $\sigma_0$ generated by $\mb e_1$
and $\mb e_2$, the cone $\sigma_1$ generated by $\mb e_2$ and
$-\mb e_1-\mb e_2$, and the cone $\sigma_2$ generated by $-\mb
e_1-\mb e_2$ and $\mb e_1$, see Fig.~\ref{compfans}~(a). Then each
affine variety $V_{\sigma_i}$ is isomorphic to~$\C^2$, with
coordinates $(x,y)$ for $\sigma_0$, \ $(x^{-1},x^{-1}y)$ for
$\sigma_1$, and $(y^{-1},xy^{-1})$ for $\sigma_2$. These three
affine charts glue together into the complex projective plane
$V_\Sigma=\C P^2$ in the standard way: if $(z_0:z_1:z_2)$ are the
homogeneous coordinats in~$\C P^2$, then we have $x=z_1/z_0$ and
$y=z_2/z_0$.
\end{example}

\begin{example}\label{hirzebruch}
Fix $k\in\Z$ and consider the complete fan in~$\R^2$ with the four
two-dimensional cones generated by the pairs of vectors $(\mb
e_1,\mb e_2)$, $(\mb e_1,-\mb e_2)$, $(-\mb e_1+k\mb e_2,-\mb
e_2)$ and $(-\mb e_1+k\mb e_2,\mb e_2)$, see
Fig~\ref{compfans}~(b). It can be shown that the corresponding
toric variety $F_k$ is the projectivisation $\C
P(\underline\C\oplus\mathcal O(k))$ of the sum of a trivial line
bundle $\underline{\C}$ and the $k$th power $\mathcal
O(k)=\gamma^{\otimes k}$ of the canonical line bundle~$\gamma$
over~$\C P^1$ (an exercise). These 2-dimensional complex varieties
$F_k$ are known as \emph{Hirzebruch surfaces}.
\end{example}

\begin{example}\label{uktoric}
Let $\sK$ be a simplicial complex on~$[m]$. The complement
$U(\sK)$ of the coordinate subspace arrangement corresponding
to~$\sK$ (see~\eqref{compl}) is a $(\C^\times)^m$-invariant subset
in~$\C^m$, and therefore it is a nonsingular toric variety. This
variety is not affine in general; it is quasiaffine (the
complement to a Zariski closed subset in an affine variety).

The fan $\Sigma_{\sK}$ corresponding to $U(\sK)$ consists of the
cones $\sigma_I\subset\R^m$ generated by the basis vectors $\mb
e_{i_1},\ldots,\mb e_{i_k}$, for all simplices
$I=\{i_1,\ldots,i_k\}\in\sK$. The affine toric variety
corresponding to $\sigma_I$ is $(\C,\C^\times)^I$, and the affine
cover of $U(\sK)$ is its polyhedral product decomposition
$U(\sK)=\bigcup_{I\in\sK}(\C,\C^\times)^I$ given by
Proposition~\ref{ukppd}.
\end{example}

The inclusion poset of closures of $\C^\times_N$-orbits of
$V_\Sigma$ is isomorphic to the reversed inclusion poset of faces
of~$\Sigma$. That is, $k$-dimensional cones of $\Sigma$ correspond
to codimension-$k$ orbits of the algebraic torus action
on~$V_\Sigma$. In particular, $n$-dimensional cones correspond to
fixed points, and the apex (the zero cone) corresponds to the
dense orbit. Furthermore, if a subcollection of cones of $\Sigma$
forms a fan $\Sigma'$, then the toric variety $V_{\Sigma'}$ is
embedded into $V_{\Sigma}$ as a Zariski open subset.

We recall from Section~\ref{combfan} that a fan is called
simplicial (respectively, regular) if each of its cones is
generated by a part of basis of the space~$N_\R$ (respectively, of
the lattice~$N$), and a fan is called complete if the union of its
cones is the whole~$N_\R$.

A toric variety $V_\Sigma$ is compact (in usual topology) if and
only if the fan $\Sigma$ is complete. If $\Sigma$ is a simplicial
fan, then $V_\Sigma$ is an \emph{orbifold}\label{dorbifold}, that
is, it is locally isomorphic to a quotient of $\C^n$ by a finite
group action. A toric variety $V_\Sigma$ is nonsingular (smooth)
if and only if the fan $\Sigma$ is regular.

\subsection*{Exercises.}
%\begin{exercise}[{\cite[\S1.4]{fult93}
%or~\cite[Theorem~3.1.5]{c-l-s11}}]\label{septoric} Show that for
%any pair of cones $\sigma_1,\sigma_2$ that intersect in common
%face~$\tau$, the diagonal map of toric varieties $V_{\tau}\to
%V_{\sigma_1}\times V_{\sigma_2}$ is a Zariski closed embedding.
%Deduce that the toric variety $V_\Sigma$ obtained by gluing the
%affine toric varieties $V_\sigma$ corresponding to
%$\sigma\in\Sigma$ is separated.
%\end{exercise}
%
\begin{exercise}\label{nstoric}
Let $\Sigma$ be the `multifan' in $\R^1$ consisting of two
identical 1-dimensional cones generated by $\mb e_1$ and a
0-dimensional cone $\bf 0$. Describe the algebraic variety
$V_{\Sigma}$ obtained by gluing the affine varieties corresponding
to this `multifan' and show that $V_\Sigma$ is not separated (or
non-Hausdorff in the usual topology).
\end{exercise}

\begin{exercise}\label{cp2without3lines}
Describe the toric variety corresponding to the fan with 3
one-dimensional cones generated by the vectors $\mb e_1$, $\mb
e_2$ and $-\mb e_1-\mb e_2$.
\end{exercise}

\begin{exercise}
Show that the toric variety of Example~\ref{hirzebruch} is
isomorphic to the Hirzebruch surface $F_k=\C
P(\underline\C\oplus\mathcal O(k))$.
\end{exercise}

\begin{exercise}\label{hstop}
Show that the Hirzebruch surface $F_k$ is homeomorphic to
$S^2\times S^2$ for even~$k$ and is homeomorphic to $\C
P^2\mathbin{\#}\overline{\C P}{}^2$ for odd~$k$, where $\#$
denotes the connected sum, and $\overline{\C P}{}^2$ is $\C P^2$
with the orientation reversed.
\end{exercise}

\section{Projective toric varieties and
polytopes}\label{projtvpol}
\begin{construction}[projective toric varieties]\label{constrproj}
Let $P$ be a convex polytope with vertices in the dual lattice
$N^*$ (a \emph{lattice polytope}), and let $\Sigma_P$ be the
normal fan of $P$ (see Construction~\ref{nf}). Since $P\subset
N_{\R}^*$, the fan $\Sigma_P$ belongs to the space~$N_\R$. It has
a maximal cone $\sigma_v$ for each vertex $v\in P$. The dual cone
$\sigma^*_v$ is the `vertex cone' at~$v$, generated by all vectors
pointing from $v$ to other points of~$P$.

Define the toric variety $V_P=V_{\Sigma_P}$. Since the normal fan
$\Sigma_P$ does not depend on the linear size of the polytope, we
may assume that for each vertex $v$ the semigroup $S_{\sigma_v}$
is generated by the lattice points of the polytope (this can
always be achieved by replacing $P$ by $kP$ with sufficiently
large~$k$). Since $N^*$ is the lattice of characters of the
algebraic torus $\C^\times_N=N\otimes_\Z\C^\times$ (that is,
$N^*=\Hom_\Z(\C^\times_N,\C^\times)$), the lattice points of the
polytope $P\subset N^*$ define an embedding
\[
  i_P\colon\C^\times_N\to(\C^\times)^{|N^*\cap P|},
\]
where $|N^*\cap P|$ is the number of lattice points in~$P$.

\begin{proposition}[{see~\cite[\S2.2]{c-l-s11} or~\cite[\S3.4]{fult93}}]
The toric variety $V_P$ is identified with the projective closure
$\overline{i_P(\C^\times_N)}\subset\C P^{|N^*\cap P|}$.
\end{proposition}

It follows that toric varieties arising from polytopes are
projective, i.e. can be defined by a set of homogeneous equations
in a projective space. The converse is also true: the fan
corresponding to a projective toric variety is the normal fan of a
lattice polytope.

The polytope $P$ carries more geometric information than the
normal fan $\Sigma_P$: different lattice polytopes with the same
normal fan $\Sigma$ correspond to different projective embeddings
of the toric variety~$V_\Sigma$.

Nonsingular projective toric varieties correspond to lattice
polytopes~$P$ which are simple and Delzant\label{delzaptope} (that
is, for each vertex $v$, the normal vectors of facets meeting
at~$v$ form a basis of the lattice~$N$).
\end{construction}

Constructions of Chapter~\ref{combi} provide explicit series of
Delzant polytopes and therefore nonsingular projective toric
varieties. Basic examples include simplices and cubes in all
dimensions. The product of two Delzant polytopes is obviously
Delzant. If $P$ is a Delzant polytope, then its face truncation
$P\cap H_\ge$ (Construction~\ref{hypcut}) by an appropriately
chosen hyperplane~$H$ is also Delzant (an exercise). All
nestohedra (in particular, permutahedra and associahedra) admit
Delzant realisations. This fact, first observed
in~\cite[Proposition~7.10]{post09}, can be proved either by using
the sequence of face truncations described in
Lemma~\ref{nhfacetrun} or directly from the presentation of
nestohedra given in Proposition~\ref{nestpres}.

\begin{example}\label{nonproj}
The fan $\Sigma$ described in Example~\ref{nonpolytopal} is
regular, but cannot be obtained as the normal fan of a simple
polytope. The corresponding 3-dimensional toric variety $V_\Sigma$
is compact and nonsingular, but not projective.
\end{example}

We note that although the fan $\Sigma$ from the previous example
cannot be realised \emph{geometrically} as the fan over the faces
of a simplicial polytope (or, equivalently, as the normal fan of a
simple polytope), its underlying simplicial complex $\sK_\Sigma$
is nevertheless \emph{combinatorially} equivalent to the boundary
complex of a simplicial polytope (namely, an octahedron with a
pyramid over one of its facets, see Fig.~\ref{npctv}). In other
words, using the terminology of Section~\ref{simsph}, the
starshaped sphere triangulation $\sK_\Sigma$ is polytopal (in the
combinatorial sense).

There are, of course, nonpolytopal starshaped spheres, such as the
Barnette sphere\label{barnesph} $\sK$ (see
Example~\ref{barnstar}). It is easy to see that a simplicial fan
realising the Barnette sphere can be chosen rational. However, to
realise a nonpolytopal sphere by a \emph{regular} fan turned out
to be a more difficult task. This question has been finally
settled by Suyama~\cite{suya}; his example is obtained by
subdividing a simplicial fan realising the Barnette sphere.
% In fact, not such examples are known:
%
%\begin{problem}\label{regnpt}
%Does there exist a complete regular fan $\Sigma$ whose underlying
%complex $\sK_\Sigma$ is a nonpolytopal sphere triangulation?
%%(in the combinatorial sense)?
%\end{problem}
This is also important for the study of \emph{quasitoric
manifolds} and other topological generalisations of toric
varieties discussed in Chapter~\ref{torus}.
%
%There are also several other topological and combinatorial
%questions related to projective toric varieties, some of which we
%discuss below.

We observe that each combinatorial simple polytope admits a convex
realisation as a lattice polytope. Indeed by a small perturbation
of the defining inequalities in~(\ref{ptope}) we can make all of
them rational (that is, with rational $\mb a_i$ and $b_i$). Such a
perturbation does not change the combinatorial type, as the
half-spaces defined by the inequalities are in general position.
As a result, we obtain a simple polytope $P'$ of the same
combinatorial type with rational vertex coordinates. To get a
lattice polytope (say, with vertices in the standard
lattice~$\Z^n$) we just take the magnified polytope $kP'$ for
appropriate $k\in\Z$. Similarly, by perturbing the vertices
instead of the hyperplanes, we can obtain a lattice realisation
for an arbitrary simplicial polytope (and we can obtain a rational
fan realisation of any starshaped sphere triangulation). However,
this argument does not work for convex polytopes which are neither
simple nor simplicial. In fact, there are exist
\emph{nonrational}\label{nrptope} combinatorial polytopes, which
cannot be realised with rational vertex coordinates,
see~\cite[Example~6.21]{zieg95} and the discussion there.

Toric geometry, even in its topological part, does not translate
to a purely combinatorial study of fans and polytopes: the
underlying convex geometry is what really matters. This is
illustrated by the simple observation that different realisations
of a combinatorial polytope by lattice polytopes often produce
different (even topologically) toric varieties:

\begin{example}
The complete regular fan $\Sigma_k$ corresponding to the
Hirzebruch surface $F_k$ (see Fig.~\ref{compfans}~(b)) is the
normal fan of a lattice quadrilateral (trapezoid), e.g. given by
\[
  P_k=\{(x_1,x_2)\in\R^2\colon
  x_1\ge0,\;x_2\ge0,\;-x_1+kx_2\ge-1,\;-x_2\ge-1\}.
\]
The polytopes $P_k$ corresponding to different $k$ are
combinatorially equivalent (as they are all quadrilaterals), but
the topology of the corresponding toric varieties $F_k$ is
different for even and odd $k$ (see Exercise~\ref{hstop}).
\end{example}

Furthermore, there exist combinatorial simple polytopes that do
not admit any lattice realisation $P$ with smooth~$V_P$:

\begin{example}[{\cite[1.22]{da-ja91}}]
Let $P$ be the dual of a 2-neighbourly simplicial $n$-polytope
(e.g., a cyclic polytope of dimension $n\ge4$, see
Example~\ref{cyclic}) with $m\ge 2^n$ vertices. Then for any
lattice realisation of $P$ the corresponding normal fan $\Sigma_P$
is not regular, and the toric variety $V_P$ is singular. Indeed,
assume that the normal fan $\Sigma_P$ is regular. Since the
1-skeleton of $\sK_P$ is a complete graph, each pair of primitive
generators $\mb a_i,\mb a_j$ of one-dimensional cones of
$\Sigma_P$ is a part of basis of $\Z^n$, and therefore $\mb a_i$
and $\mb a_j$ must be different modulo~2. This is a contradiction,
since the number of different nonzero vectors in $\Z_2^n$ is
$2^n-1$.
\end{example}

\subsection*{Exercises.}

\begin{exercise}
Let $P$ be a Delzant polytope and $G\subset P$ a face. Show that a
hyperplane $H$ truncating $G$ from $P$ in
Construction~\ref{hypcut} can be chosen so that the truncated
polytope $P\cap H_\ge$ is Delzant.
\end{exercise}

\begin{exercise}
The presentation of a nestohedron $P_{\mathcal B}$ from
Proposition~\ref{nestpres} is Delzant.
\end{exercise}

\begin{exercise}
Describe explicitly a rational fan realising the Barnette sphere
by writing down its primitive integral generator vectors.
\end{exercise}

\begin{exercise}
Write down a system of homogeneous equations defining each
Hirzebruch surface in a projective space. (Hint: use
Construction~\ref{constrproj}.)
\end{exercise}

\section{Cohomology of toric manifolds}\label{cohtm}
A \emph{toric manifold} is a smooth compact toric variety.
(Compactness will be always understood in the sense of usual
topology; it corresponds to algebraic geometer's notion of
\emph{completeness.}) Toric manifolds $V_\Sigma$ correspond to
complete regular fans~$\Sigma$. \emph{Projective toric manifolds}
$V_P$ correspond to lattice polytopes $P$ whose normal fans are
regular.

The cohomology of a toric manifold $V_\Sigma$ can be calculated
effectively from the fan~$\Sigma$. The Betti numbers are
determined by the combinatorics of~$\Sigma$ only, while the ring
structure of $H^*(V_\Sigma)$ depends on the geometric data. The
required combinatorial ingredients are the $h$-vector $\mb
h(\sK_\Sigma)=(h_0,h_1,\ldots,h_n)$ (see Definition~\ref{simfve})
of the underlying simplicial complex~$\sK_\Sigma$ and its face
ring $\Z[\sK_\Sigma]$ (Definition~\ref{frsim}). The geometric data
consists of the primitive generators $\mb a_1,\ldots,\mb a_m$ of
one-dimensional cones (edges) of~$\Sigma$.

\begin{theorem}[Danilov--Jurkiewicz]
\label{danjur} Let $V_\Sigma$ be the toric manifold corresponding
to a complete regular fan~$\Sigma$ in~$N_\R$. The cohomology ring
of $V_\Sigma$ is given by
\[
  H^*(V_\Sigma)\cong \Z[v_1,\ldots,v_m]/\mathcal I,
\]
where $v_1,\ldots,v_m\in H^2(V_\Sigma)$ are the cohomology classes
dual the invariant divisors corresponding to the one-dimensional
cones of~$\Sigma$, and $\mathcal I$ is the ideal generated by
elements of the following two types:
\begin{itemize}
\item[(a)] $v_{i_1}\cdots v_{i_k}$ with
$\{i_1,\ldots,i_k\}\notin\sK_\Sigma$ (the Stanley--Reisner
relations);

\item[(b)] $\displaystyle\sum_{i=1}^m\langle\mb a_j,\mb u\rangle v_i$, for
any $\mb u\in N^*$.
\end{itemize}

The homology groups of $V_\Sigma$ vanish in odd dimensions, and
are free abelian in even dimensions, with ranks given by
$$
  b_{2i}(V_\Sigma)=h_i(\sK_\Sigma),
  %\quad\text{for } i=0,1,\ldots,n,
$$
where $h_i(\sK_\Sigma)$, $i=0,1,\ldots,n$, are the components of
the $h$-vector of~$\sK_\Sigma$.
\end{theorem}

This theorem was proved by Jurkiewicz for projective toric
manifolds and by Danilov~\cite[Theorem~10.8]{dani78} in the
general case. We shall give a topological proof of a more general
result in Section~\ref{torusman} (see Theorem~\ref{theo:stcoh}).

To obtain an explicit presentation of the ring $H^*(V_\Sigma)$ we
choose a basis of $N$ and write the vectors $\mb a_j$ in
coordinates: $\mb a_j=(a_{j1},\ldots,a_{jn})^t$, $1\le j\le m$.
Then the ideal $J_\Sigma$ is generated by the $n$ linear forms
\[
  t_i=a_{1i}v_1+\cdots+a_{mi}v_m\in\Z[v_1,\ldots,v_m],\quad 1\le i\le
  n.
\]
By Lemma~\ref{lsopcrit}, the sequence $t_1,\ldots,t_n$ is an lsop
in the Cohen--Macaulay ring~$\Z[\sK_\Sigma]$, so it is a regular
sequence. Hence, $\Z[\sK_\Sigma]$ is a free
$\Z[t_1,\ldots,t_n]$-module, and statement~(a) of
Theorem~\ref{danjur} follows from~(b) and Theorem~\ref{psfr}.

\begin{remark}
Theorem~\ref{danjur} remains valid for complete simplicial fans
and corresponding toric orbifolds if the integer coefficients are
replaced by the rationals~\cite{dani78}. The integral cohomology
of toric orbifolds often has torsion, and the ring structure is
subtle even in the simplest case of weighted projective
spaces~\cite{kawa73},~\cite{b-f-n-r12}.
\end{remark}

It follows from Theorem~\ref{danjur} that the cohomology ring of
$V_\Sigma$ is generated by two-dimensional classes. This is the
first property to check if one wishes to determine whether a given
algebraic variety or smooth manifold has a structure of a toric
manifold. For instance, this rules out flag varieties and
Grassmanians different from projective spaces. Another important
property of toric manifolds is that the \emph{Chow
ring}\label{chowring} of $V_\Sigma$ coincides with its integer
cohomology ring~\cite[\S~5.1]{fult93}.

Assume now that $\Sigma=\Sigma_P$ is the normal fan of a lattice
polytope $P=P(A,\mb b)$ given by~\eqref{ptope}. Let $V_P$ be the
corresponding projective toric variety, see
Construction~\ref{constrproj}. In the notation of
Theorem~\ref{danjur},
%it is well known in toric geometry that
the linear combination $D_P=b_1D_1+\dots+b_mD_m$ is an \emph{ample
divisor}\label{amplediviso} on~$V_P$ (see, e.g.,
\cite[Proposition~6.1.10]{c-l-s11}). This means that, when $k$ is
sufficiently large, $kD_P$ is a hyperplane section divisor for a
projective embedding $V_P\subset\C P^r$. In fact, the space of
sections $H^0(V_P,kD_P)$ of (the line bundle corresponding
to)~$kD_P$ has basis corresponding to the lattice points in~$kP$.
One may take $k$ so that $kP$ has `enough lattice points' to get
an embedding of $V_P$ into the projectivisation of
$H^0(V_P,kD_P)$; this is exactly the embedding described in
Construction~\ref{constrproj}. Let $\omega=b_1v_1+\dots+b_mv_m\in
H^2(V_P;\C)$ be the complex cohomology class of~$D_P$.

\begin{theorem}[Hard Lefschetz Theorem for toric orbifolds]
\label{hlth} Let $P$ be a lattice simple polytope~\eqref{ptope},
let $V_P$ be the corresponding projective toric variety, and let
$\omega=b_1v_1+\dots+b_mv_m\in H^2(V_P;\C)$ be the class defined
above. Then the maps
\[
  H^{n-i}(V_P;\C)\stackrel{\omega^i}\longrightarrow H^{n+i}(V_P;\C)
\]
are isomorphisms for all $i=1,\ldots,n$.
\end{theorem}
If $V_P$ is smooth, then it is K\"ahler, and $\omega$ is the class
of the K\"ahler 2-form.

The proof of the Hard Lefschetz Theorem is well beyond the scope
of this book. In fact, it is a corollary of a more general version
of Hard Lefschetz Theorem for the (middle perversity)
\emph{intersection cohomology}\label{intecohomol}, which is valid
for all projective varieties (not necessary orbifolds). See the
discussion in~\cite[\S5.2]{fult93} or~\cite[\S12.6]{c-l-s11}.

Now we are ready to give Stanley's argument for the `only if' part
of the $g$-theorem for simple polytopes:

\begin{proof}[Proof of the necessity part of
Theorem~{\rm\ref{gth}}]\label{proofgth} We need to establish
conditions~(a)--(c) for a combinatorial simple polytope. Realise
it by a lattice polytope $P\subset\R^n$ as described in the
previous section. Let $V_P$ be the corresponding toric variety.
Part~(a) is already proved (Theorem~\ref{ds}). It follows from
Theorem~\ref{hlth} that the multiplication by $\omega\in
H^2(V_P;\Q)$ is a monomorphism $H^{2i-2}(V_P;\Q)\to
H^{2i}(V_P;\Q)$ for $i\le \sbr n2$. This together with part~(a) of
Theorem~\ref{danjur} gives that $h_{i-1}\le h_i$ for $0\le i\le
\sbr n2$, thus proving~(b). To prove~(c), define the graded
commutative $\Q$-algebra $A=H^*(V_P;\Q)/(\omega)$, where
$(\omega)$ is the ideal generated by~$\omega$. Then $A^0=\Q$,
$A^{2i}=H^{2i}(V_P;\Q)/(\omega\cdot H^{2i-2}(V_P;\Q))$ for $1\le
i\le\sbr n2$, and $A$ is generated by degree-two elements (since
so is $H^*(V_P;\Q)$). It follows from Theorem~\ref{mvect} that the
numbers $\dim A^{2i}=h_i-h_{i-1}$, \ $0\le i\le\sbr n2$, are the
components of an $M$-vector, thus proving~(c) and the whole
theorem.
\end{proof}

\begin{remark}
The Dehn--Sommerville equations\label{ds2} now can be interpreted
as Poincar\'e duality for~$V_P$. (Even though $V_P$ may be not
smooth, the rational cohomology algebra of a toric orbifold
satisfies Poincar\'e duality.)
\end{remark}

The Hard Lefschetz Theorem holds for projective varieties only.
Therefore, Stanley's argument cannot be generalised to
nonpolytopal spheres.
%So far this case is the only generality in which methods involving
%the Hard Lefschetz Theorem have been efficient for proving the
%$g$-theorem.
However, the cohomology of toric varieties can be used to prove
statements generalising the $g$-theorem in a different direction,
namely, to the case of general (not necessarily simple or
simplicial) convex polytopes. So suppose $P$ is a convex lattice
$n$-polytope, and $V_P$ is the corresponding projective toric
variety. If $P$ is not simple then $V_P$ has worse than
orbifold-type singularities and its ordinary cohomology behaves
badly. The Betti numbers of $V_P$ are not determined by the
combinatorial type of~$P$ and do not satisfy Poincar\'e duality.
On the other hand, the dimensions $\widehat
h_i=\dim\textit{IH}_{2i}(V_P)$ of the intersection homology groups
of $V_P$ are combinatorial invariants of~$P$ (see the description
in~\cite{stan87} or~\cite[\S12.5]{c-l-s11}). The vector
$$
  \widehat{\mb h}(P)=(\widehat h_0,\widehat h_1,\dots,\widehat h_n)
$$
is called the \emph{intersection $h$-vector}, or the \emph{toric
$h$-vector} of~$P$. If $P$ is simple, then the toric $h$-vector
coincides with the standard $h$-vector, but in general
$\widehat{\mb h}(P)$ is not determined by the face numbers of~$P$.
The toric $h$-vector satisfies the `Dehn--Sommerville equations'
$\widehat h_i=\widehat h_{n-i}$\label{glbctor}, and the Hard
Lefschetz Theorem for intersection cohomology shows that it also
satisfies the GLBC inequalities:
$$
  \widehat h_0\le\widehat h_1\le\cdots\le\widehat h_{\sbr n2}.
$$
In the case when $P$ cannot be realised by a lattice polytope
(i.e., when $P$ is nonrational), the toric $h$-vector can still be
defined combinatorially, but the GLBC inequalities require a
separate proof. Partial results in this direction were obtained by
several people, before the Hard Lefschetz Theorem for nonrational
polytopes was eventually proved in the work of Karu~\cite{karu04}.
This result also gives a purely combinatorial proof of the Hard
Lefschetz Theorem for projective toric varieties.

\subsection*{Exercises.}
\begin{exercise}
Show that the complex Grassmanian $\mathop{\mathrm{Gr}}_k(\C^n)$
($k$-planes in~$\C^n$) with $1<k<n-1$ does not support an
algebraic torus action turning it into a toric variety.
\end{exercise}

\section{Algebraic quotient construction}\label{algtq}
Along with the classical construction of toric varieties from
fans, described in Section~\ref{toricfan}, there is an alternative
way to define a toric variety: as the quotient of a Zariski open
subset in $\C^m$ (more precisely, the complement of a coordinate
subspace arrangement) by an action of an abelian algebraic group
(a product of an algebraic torus and a finite group). Different
versions of this construction, which we refer to as simply the
`quotient construction', have appeared in the work of several
authors since the early 1990s. In our exposition we mainly follow
the work of Cox~\cite{cox95} (and also its modernised exposition
in~\cite[Chapter~5]{c-l-s11}); more historical remarks can be also
found in these sources.

\subsection*{Quotients in algebraic geometry} Taking quotients of algebraic varieties by algebraic
group actions is tricky for both topological and algebraic
reasons. First, as algebraic groups are non-compact  (as algebraic
tori), their orbits may be not closed, and the quotients may be
non-Hausdorff. Second, even if the quotient is Hausdorff as a
topological space, it may fail to be an algebraic variety. This
may be remedied to some extent by the notion of the categorical
quotient.

Let $X$ be an algebraic variety with an action of an affine
algebraic group~$G$. An algebraic variety $Y$ is called a
\emph{categorical quotient}\label{catequot} of $X$ by the action
of~$G$ if there exists a morphism $\pi\colon X\to Y$ which is
constant on $G$-orbits of~$X$ and has the following universal
property: for any morphism $\varphi\colon X\to Z$ which is
constant on $G$-orbits, there is a unique morphism
$\widehat\varphi\colon Y\to Z$ such that
$\widehat\varphi\circ\pi=\varphi$. This is described by the
diagram
\[
\xymatrix{
  X \ar[rr]^{\varphi} \ar[dr]^{\pi} && Z\\
  & Y \ar@{-->}[ur]^{\widehat\varphi}
}
\]
A categorical quotient $Y$ is unique up to isomorphism, and we
denote it by $X\dbs G$.
%(although sometimes this notation is
%reserved for categorical quotients with extra good properties).

Assume first that $X=\Spec A$ is an affine variety, where
$A=\C[X]$ is the algebra of regular functions on~$X$. Let
$\C[X]^{G}$ be the subalgebra of $G$-invariant functions (i.e.
such functions $f$ that $f(gx)=f(x)$ for any $g\in G$ and $x\in
X$). If $G$ is an algebraic torus (or any \emph{reductive} affine
algebraic group), then $\C[X]^{G}$ is finitely generated. The
corresponding affine variety $\Spec \C[X]^{G}$ is the categorical
quotient $X\dbs G$. The quotient morphism $\pi\colon X\to X\dbs G$
is dual to the inclusion of algebras $\C[X]^{G}\to \C[X]$. The
morphism $\pi$ is surjective and induces a one-to-one
correspondence between points of $X\dbs G$ and \emph{closed}
$G$-orbits of~$X$ (i.e. $\pi^{-1}(x)$ contains a unique closed
$G$-orbit for any $x\in X\dbs G$,
see~\cite[Proposition~5.0.7]{c-l-s11}).

Therefore, if all $G$-orbits of an affine variety $X$ are closed,
then the categorical quotient $X\dbs G$ is identified as a
topological space with the ordinary `topological' quotient $X/G$.
Quotients of this type are called
\emph{geometric}\label{geomequot} and also denoted by~$X/G$.

\begin{example}
Let $\C^\times$ act on $\C=\Spec(\C[z])$ by scalar multiplication.
There are two orbits: the closed orbit $0$ and the open orbit
$\C^\times$. The topological quotient $\C/\C^\times$ is a
non-Hausdorff two-point space.

On the other hand, the categorical quotient
$\C\dbs\C^\times=\Spec(\C[z]^{\C^\times})$ is a point, since any
$\C^\times$-invariant polynomial is constant (and there is only
one closed orbit).

Similarly, if $\C^\times$ acts on $\C^n=\Spec(\C[z_1,\ldots,z_n])$
diagonally, then an invariant polynomial satisfies $f(\lambda
z_1,\ldots,\lambda z_n)=f(z_1,\ldots,z_n)$ for all
$\lambda\in\C^\times$. Such a polynomial must be constant, hence
$\C^n\dbs\C^{\times}$ is a point.
%, and the only closed orbit is the origin.
\end{example}

In good cases categorical quotients of general (non-affine)
varieties $X$ may be constructed by `gluing from affine pieces' as
follows. Assume that $G$ acts on $X$ and $\pi\colon X\to Y$ is a
morphism of varieties that is constant on $G$-orbits. If $Y$ has
an open affine cover $Y=\bigcup_\alpha V_\alpha$ such that
$\pi^{-1}(V_\alpha)$ is affine and $V_\alpha$ is the categorical
quotient (that is,
$\pi|_{\pi^{-1}(V_\alpha)}\colon\pi^{-1}(V_\alpha)\to V_\alpha$ is
the morphism dual to the inclusion of algebras
$\C[\pi^{-1}(V_\alpha)]^G\to\C[\pi^{-1}(V_\alpha)]$), then $Y$ is
the categorical quotient $X\dbs G$.

\begin{example}
Let $\C^\times$ act on $\C^2\setminus\{\mathbf0\}$ diagonally,
where $\C^2=\Spec(\C[z_0,z_1])$. We have an open affine cover
$\C^2\setminus\{\mathbf0\}=U_0\cup U_1$, where
\begin{align*}
  U_0&=\C^2\setminus\{z_0=0\}=\C^\times\times\C=\Spec(\C[z_0^{\pm1},z_1]),\\
  U_1&=\C^2\setminus\{z_1=0\}=\C\times\C^\times=\Spec(\C[z_0,z_1^{\pm1}]),\\
  U_0\cap
  U_1&=\C^2\setminus\{z_0z_1=0\}=\C^\times\times\C^\times=\Spec(\C[z_0^{\pm1},z_1^{\pm1}]).
\end{align*}
The algebras of $\C^\times$-invariant functions are
\[
  \C[z_0^{\pm1},z_1]^{\C^\times}\!\!=\C[z_1/z_0],\quad
  \C[z_0,z_1^{\pm1}]^{\C^\times}\!\!=\C[z_0/z_1],\quad
  \C[z_0^{\pm1},z_1^{\pm1}]^{\C^\times}\!\!=\C[(z_1/z_0)^{\pm1}].
\]
It follows that $V_0=U_0\dbs\C^\times\cong\C$ and
$V_1=U_1\dbs\C^\times\cong\C$ glue together along $V_0\cap
V_1=(U_0\cap U_1)\dbs\C^\times\cong\C^\times$ in the standard way
to produce~$\C P^1$. All $\C^\times$-orbits are closed in
$\C^2\setminus\{\mathbf0\}$, hence $\C
P^1=(\C^2\setminus\{\mathbf0\})/\C^\times$ is the geometric
quotient.

Similarly, $\C P^n=(\C^{n+1}\setminus\{\mathbf0\})/\C^\times$ for
the diagonal action of~$\C^\times$.
\end{example}

\begin{example}
Now we let $\C^\times$ act on $\C^2\setminus\{\mathbf0\}$ by
$\lambda\cdot(z_0,z_1)=(\lambda z_0,\lambda^{-1}z_1)$. Using the
same affine cover of $\C^2\setminus\{\mathbf0\}$ as in the
previous example, we obtain the following algebras of
$\C^\times$-invariant functions:
\[
  \C[z_0^{\pm1},z_1]^{\C^\times}\!\!=\C[z_0z_1],\quad
  \C[z_0,z_1^{\pm1}]^{\C^\times}\!\!=\C[z_0z_1],\quad
  \C[z_0^{\pm1},z_1^{\pm1}]^{\C^\times}\!\!=\C[(z_0z_1)^{\pm1}].
\]
This times gluing together $V_0$ and $V_1$ along $V_0\cap
V_1\cong\C^\times$ gives the variety obtained from two copies of
$\C$ by identifying all nonzero points. This variety is
nonseparated (the two zeros do not have nonintersecting
neighbourhoods in the usual topology). Since we only consider
separated varieties, this is not a categorical quotient.
\end{example}

A toric variety $V_\Sigma$ will be described as the
%categorical (or, in good cases, geometric)
quotient of the `total space' $U(\Sigma)$ by an action of a
commutative algebraic group~$G$.
%(which is an algebraic torus in good cases, and is a product of a
%torus and a finite group in general).
We now proceed to describe $G$ and~$U(\Sigma)$.

\subsection*{The total space $U(\Sigma)$ and the acting group~$G$}
Let $\Sigma$ be a rational fan in the $n$-dimensional space $N_\R$
with $m$ one-dimensional cones generated by primitive vectors $\mb
a_1,\ldots,\mb a_m$. We shall assume that the linear span of $\mb
a_1,\ldots,\mb a_m$ is the whole~$N_\R$. (Equivalently, the toric
variety~$V_\Sigma$ does not have torus factors, i.e. cannot be
written as $V_\Sigma=V_{\Sigma'}\times\C^\times$. For the general
case see~\cite[\S5.1]{c-l-s11}.)

We consider the map of lattices $\Amap\colon\Z^m\to N$ sending the
$i$th basis vector of $\Z^m$ to $\mb a_i\in N$. Our assumption
implies that the corresponding map of algebraic tori,
\[
  \Amap\otimes_\Z\C^\times\colon(\C^\times)^m\to \C^\times_N
\]
is surjective.

Define the group $G=G(\Sigma)$ as the kernel of the map
$\Amap\otimes_\Z\C^\times$, which we denote by $\exp\Amap$.
%(which we continue to denote by $\Amap$ for simplicity).
We therefore have an exact sequence of groups
\begin{equation}\label{ggrou}
  1\longrightarrow G\longrightarrow (\C^\times)^m
  \stackrel{\exp\Amap}\longrightarrow \C^\times_N\longrightarrow1.
\end{equation}
Explicitly, $G$ is given by
\begin{equation}\label{gexpl}
  G=\bigl\{(z_1,\ldots,z_m)\in(\C^\times)^m\colon
  \prod_{i=1}^m z_i^{\langle\mb a_i,\mb u\rangle}=1
  \quad\text{for any }\mb u\in N^*\bigr\}.
\end{equation}
The group $G$ is isomorphic to a product of $(\C^\times)^{m-n}$
and a finite abelian group. If $\Sigma$ is a regular fan with at
least one $n$-dimensional cone, then $G\cong(\C^\times)^{m-n}$.

Given a cone $\sigma\in\Sigma$, we set
$g(\sigma)=\{i_1,\ldots,i_k\}\subset[m]$ if $\sigma$ is generated
by $\mb a_{i_1}\ldots,\mb a_{i_k}$, and consider the monomial
$z^{\hat\sigma}=\prod_{j\notin g(\sigma)}z_j$. The quasiaffine
variety
\[
  U(\Sigma)=\C^m\setminus\{\mb z\in\C^m\colon z^{\hat\sigma}=0\text{ for all
  }\;\sigma\in\Sigma\}
\]
has the affine cover
\begin{equation}\label{affinecover}
  U(\Sigma)=\bigcup_{\sigma\in\Sigma}U(\sigma).
\end{equation}
by affine varieties
\[
  U(\sigma)=\{\mb z\in\C^m\colon z^{\hat\sigma}\ne0\}=
  \{\mb z\in\C^m\colon z_j\ne0\;\text{ for }\;j\notin
  g(\sigma)\}=(\C,\C^\times)^{g(\sigma)}.
\]
Here we used the notation of Construction~\ref{nsc}, so that
$U(\sigma)\cong\C^k\times(\C^\times)^{m-k}$. Each subset
$U(\sigma)\subset\C^m$ is invariant under the coordinatewise
action of $(\C^\times)^m$ on $\C^m$, so that $U(\Sigma)$ is also
invariant.

By definition, $U(\Sigma)$ is the complement of a union of
coordinate subspaces, so we know from Proposition~\ref{cacor} that
it has the form
\begin{equation}\label{UK}
  U(\sK)=\C^m\big\backslash\bigcup_{\{i_1,\ldots,i_k\}\notin\sK}\bigl\{\mb z
  \in\C^m\colon z_{i_1}=\cdots=z_{i_k}=0\bigr\}=\zk(\C,\C^\times)
\end{equation}
for some simplicial complex $\sK$ on $[m]$. What is this
simplicial complex?

The answer is suggested by decomposition~\eqref{affinecover}. We
define the simplicial complex $\sK_\Sigma$ generated by all
subsets $g(\sigma)\subset[m]$:
\[
  \sK_\Sigma=\{I\colon I\subset g(\sigma)\quad\text{for some
  }\sigma\in\Sigma\}.
\]
If $\Sigma$ is a simplicial fan, then each $I\subset g(\sigma)$ is
$g(\tau)$ for some $\tau\in\Sigma$, and we obtain the `underlying
complex' of~$\Sigma$ defined in Example~\ref{simincl}. In
particular, if $\Sigma$ is simplicial and complete, then
$\sK_\Sigma$ is a triangulation of~$S^{n-1}$; and if $\Sigma$ is a
normal fan of a simple polytope, then $\sK_\Sigma$ is the boundary
complex of the dual simplicial polytope. If $\Sigma$ is the normal
fan of a non-simple polytope~$P$ (i.e. the fan over the faces of
the polar polytope~$P^*$), then $\sK_\Sigma$ is obtained by
replacing each face of $P^*$ by a simplex with the same set of
vertices; such a simplicial complex is not pure in general.

\begin{proposition}
We have $U(\Sigma)=U(\sK_\Sigma)$.
\end{proposition}
\begin{proof}
We have $(\C,\C^\times)^I\subset(\C,\C^\times)^{g(\sigma)}$
whenever $I\subset g(\sigma)$, hence,
\[
  U(\Sigma)=\bigcup_{\sigma\in\Sigma}U(\sigma)=
  \bigcup_{\sigma\in\Sigma}(\C,\C^\times)^{g(\sigma)}=\!\!
  \bigcup_{\sigma\in\Sigma,I\subset g(\sigma)}\!\!(\C,\C^\times)^I=
  \bigcup_{I\in\sK_\Sigma}(\C,\C^\times)^I=U(\sK_\Sigma).
\]
\end{proof}

We observe that the subset $U(\Sigma)\subset\C^m$ depends only on
the combinatorial structure of the fan~$\Sigma$, while the
subgroup $G\subset(\C^\times)^m$ depends on the geometric data,
namely, the primitive generators of one-dimensional cones.

\subsection*{Toric variety as a quotient}
Since $U(\Sigma)\subset\C^m$ is invariant under the coordinatewise
action of $(\C^\times)^m$, we obtain a $G$-action on $U(\Sigma)$
by restriction.

\begin{theorem}[{Cox~\cite[Theorem~2.1]{cox95}}]\label{coxth}
Assume that the linear span of one-dimensional cones of $\Sigma$
is the whole space~$N_\R$.\\
{\rm(a)} The toric variety $V_\Sigma$ is naturally isomorphic
to the categorical quotient $U(\Sigma)\dbs G$.\\
{\rm(b)} $V_\Sigma$ is the geometric quotient $U(\Sigma)/G$ if and
only if the fan $\Sigma$ is simplicial.
\end{theorem}
\begin{proof}
We first prove that the affine variety $V_\sigma$ corresponding to
a cone $\sigma\in\Sigma$ is the categorical quotient
$U(\sigma)\dbs G$. The algebra of regular functions
$\C[U(\sigma)]$ is isomorphic to $\C[z_i,z_j^{-1}\colon 1\le i\le
m,\;j\notin g(\sigma)]$ and is generated by Laurent monomials
$\prod_{i=1}^mz_i^{k_i}$ with $k_i\ge0$ for $i\in g(\sigma)$.

It follows easily from~\eqref{gexpl} that a monomial
$\prod_{i=1}^mz_i^{k_i}$ is invariant under the $G$-action on
$U(\sigma)$ if and only if it has the form
$\prod_{i=1}^mz_i^{\langle\mb u,\mb a_i\rangle}$ for some $\mb
u\in N^*$.

Conditions $\langle\mb u,\mb a_i\rangle\ge0$ for $i\in g(\sigma)$
specify the dual cone $\sigma^\mathsf{v}\subset N_\R^*$,
see~\eqref{dualcone}. Hence the invariant subalgebra
$\C[U(\sigma)]^{G}$ is isomorphic to $\C[\sigma^\mathsf{v}\cap
N^*]=A_\sigma=\C[V_\sigma]$ (the isomorphism is given by
$\prod_{i=1}^my_i^{\langle\mb u,\mb a_i\rangle}\mapsto\chi^{\mb
u}$). Thus, $U(\sigma)\dbs G\cong V_\sigma$.

The next step is to glue the isomorphisms $U(\sigma)\dbs G\cong
V_\sigma$ together into an isomorphism $U(\Sigma)\dbs G\cong
V_\Sigma$. To do this we need to check that the isomorphisms
$\C[U(\sigma)]^{G}\to \C[V_\sigma]$ are compatible when we pass to
the faces of~$\sigma$. In other words, for each face
$\tau\subset\sigma$ we need to establish the commutativity of the
diagram
\begin{equation}\label{compatible}
\begin{CD}
  \C[U(\sigma)]^{G} @>>> \C[U(\tau)]^{G}\\
  @V\cong VV @VV\cong V\\
  \C[V_\sigma] @>>> \C[V_\tau].
\end{CD}
\end{equation}
By the definition of a face, we have $\tau=\sigma\cap\mb u^\bot$
for some $\mb u\in\sigma^\mathsf{v}\cap N^*$, where $\mb u^\bot$
denotes the hyperplane in $N_\R$ normal to~$\mb u$. Consider the
monomial $\mb z(\mb u)=\prod_{i=1}^mz_i^{\langle\mb u,\mb
a_i\rangle}$. Since $\tau=\sigma\cap\mb u^\bot$, the monomial $\mb
z(\mb u)$ has positive exponent of $z_i$ for $i\in
g(\sigma)\setminus g(\tau)$ and zero exponent of $z_j$ for $j\in
g(\tau)$. It follows that the algebra $\C[U(\tau)]$ is the
localisation of $\C[U(\sigma)]$ by the ideal generated by $\mb
z(\mb u)$, i.e. $\C[U(\tau)]=\C[U(\sigma)]_{\mb z(\mb u)}$. Since
$\mb z(\mb u)$ is a $G$-invariant monomial, the localisation
commutes with passing to invariant subalgebras, i.e.
$\C[U(\tau)]^{G}=\C[U(\sigma)]^{G}_{\mb z(\mb u)}$. Similarly,
$\C[V_\tau]=A_\tau=(A_\sigma)_{\chi^{\mb u}}$.
Diagram~\eqref{compatible} then takes the form
\[
\begin{CD}
  \C[U(\sigma)]^{G} @>>> \C[U(\sigma)]^{G}_{\mb z(\mb u)}\\
  @V\cong VV @VV\cong V\\
  \C[V_\sigma] @>>> \C[V_\sigma]_{\chi^{\mb u}},
\end{CD}
\]
where the vertical arrows are localisation maps. It is obviously
commutative.

Now using the affine cover~\eqref{affinecover} and the
compatibility of the isomorphisms on affine varieties we obtain
the isomorphism $U(\Sigma)\dbs G\cong
V_\Sigma=\bigcup_{\sigma\in\Sigma}V_\sigma$. Statement~(a) is
therefore proved.

\smallskip

To verify (b) we need to check that all orbits of the $G$-action
on $U(\Sigma)$ are closed if and only if the fan $\Sigma$ is
simplicial.

Assume $\Sigma$ is simplicial, and consider any $G$-orbit $G\mb
z$, \ $\mb z\in U(\Sigma)$. We shall prove that $G\mb z$ is closed
in the usual topology, which is sufficient since the closures of
orbits in the usual and Zariski topologies coincide. We need to
check that whenever a sequence $\{\mb w^{(k)}\colon
k=1,2,\ldots\}$ of points of $G\mb z$ has a limit $\mb w\in
U(\Sigma)$, this limit is in~$G\mb z$. Write $\mb w^{(k)}=\mb
g^{(k)}\mb z$ with $\mb g^{(k)}\in G$. Then it is enough to show
that a subsequence of $\{\mb g^{(k)}\}$ converges to a point $\mb
g\in G$, as in this case $\lim_{k\to\infty}\mb w^{(k)}=\mb g\mb
z\in G\mb z$. We write
\[
  \mb g^{(k)}=\bigl(g_1^{(k)},\ldots,g_m^{(k)}\bigr)=
  \bigl(e^{\alpha^{(k)}_1+i\beta^{(k)}_1},\ldots,
  e^{\alpha^{(k)}_m+i\beta^{(k)}_m}\bigr)\in G\subset (\C^\times)^m
\]
where $g_j^{(k)}\in\C^\times$ and
$\alpha_j^{(k)},\beta_j^{(k)}\in\R$. Since
$e^{i\beta_j^{(k)}}\in\mathbb S$ and the circle is compact, we may
assume by passing to a subsequence that the sequence
$\{e^{i\beta_j^{(k)}}\}$ has a limit $e^{i\beta_j}$ as
$k\to\infty$ for each $j=1,\ldots,m$. It remains to consider the
sequences $\{e^{\alpha^{(k)}_j}\}$.

By passing to a subsequence we may assume that each sequence
$\{\alpha^{(k)}_j\}$, \ $j=1,\ldots,m$, has a finite or infinite
limit (including $\pm\infty$). Let
\[
  I_+=\{j\colon\alpha^{(k)}_j\to +\infty\}\subset[m],\quad
  I_-=\{j\colon\alpha^{(k)}_j\to -\infty\}\subset[m].
\]
Since the sequence $\{\mb w^{(k)}=\mb g^{(k)}\mb z\}$ is
converging to $\mb w=(w_1,\ldots,w_m)\in U(\Sigma)$,
%and is therefore bounded,
we have $z_j=0$ for $j\in I_+$ and $w_j=0$ for $j\in I_-$. Then it
follows from the decomposition
$U(\Sigma)=\bigcup_{I\in\sK_\Sigma}(\C,\C^\times)^I$ that $I_+$
and $I_-$ are simplices of~$\sK_\Sigma$. Let $\sigma_+,\sigma_-$
be the corresponding cones of~$\Sigma$ (here we use the fact that
$\Sigma$ is simplicial). Then $\sigma_+\cap\sigma_-=\{\mathbf 0\}$
by definition of a fan. By Lemma~\ref{seplemma}, there is a linear
function $\mb u\in N^*$ such that $\langle\mb u,\mb a\rangle>0$
for any nonzero $\mb a\in \sigma_+$, and $\langle\mb u,\mb
a\rangle<0$ for any nonzero $\mb a\in \sigma_-$. Now, since $\mb
g^{(k)}\in G$, it follows from~\eqref{gexpl} that
\begin{equation}\label{alphak}
%  0=\Bigl\langle\mb u,\sum_{j=1}^m \alpha^{(k)}_j\mb a_j\Bigr\rangle=
  \sum_{j=1}^m \alpha^{(k)}_j\langle\mb u,\mb a_j\rangle=0.
\end{equation}
This implies that both $I_+$ and $I_-$ are empty, as otherwise the
sum above tends to infinity. Thus, each sequence
$\{\alpha^{(k)}_j\}$ has a finite limit $\alpha_j$, and a
subsequence of $\{\mb g^{(k)}\}$ converges to
$(e^{\alpha_1+i\beta_1},\ldots,e^{\alpha_m+i\beta_m})$. Passing to
the limit in~\eqref{alphak} and in the similar equation for
$\beta^{(k)}_j$ as $k\to\infty$ we obtain that
$(e^{\alpha_1+i\beta_1},\ldots,e^{\alpha_m+i\beta_m})\in G$.

The fact that the $G$-action on $U(\Sigma)$ with non-simplicial
$\Sigma$ has non-closed orbits is left as an exercise
(alternatively, see~\cite[\S2]{cox95}
or~\cite[Theorem~5.1.11]{c-l-s11}.
\end{proof}

The quotient torus $\C^\times_N=(\C^\times)^m/G$ acts on
$V_\Sigma=U(\Sigma)\dbs G$ with a dense orbit.

\begin{remark}
Observe that $U(\Sigma)$ is itself a toric variety by
Example~\ref{uktoric}. The map $\Amap\colon\R^m\to N_\R$, $\mb
e_i\mapsto\mb a_i$, projects the fan $\Sigma_{\sK}$ corresponding
to $U(\Sigma)$ to the fan $\Sigma$ corresponding to~$V_\Sigma$.
This projection defines a morphism of toric varieties
$U(\Sigma)\to V_\Sigma$, which is exactly the quotient described
above. Both fans $\Sigma_{\sK}$ and $\Sigma$ have the same
underlying simplicial complex~$\sK$.
\end{remark}

Another way to see that the orbits of the $G$-action on
$U(\Sigma)$ are closed is to use the almost freeness of this
action (see Exercise~\ref{afclosed}):

{\samepage
\begin{proposition}\label{freeaction}\
\begin{itemize}
\item[(a)] If $\Sigma$ is a simplicial fan, then the $G$-action on $U(\Sigma)$
is almost free (i.e., all stabiliser subgroups are finite);

\item[(b)] If $\Sigma$ is regular, then the $G$-action on $U(\Sigma)$ is free.
\end{itemize}
\end{proposition}
}
\begin{proof}
The stabiliser of a point $\mb z\in\C^m$ under the action of
$(\C^\times)^m$ is
\[
  (\C^\times)^{\omega(\mb z)}=\{(t_1,\ldots,t_m)\in(\C^\times)^m\colon
  t_i=1\text{ if }z_i\ne0\},
\]
where $\omega(\mb z)$ be the set of zero coordinates of~$\mb z$.
The stabiliser of $\mb z$ under the $G$-action is $G_{\mb
z}=(\C^\times)^{\omega(\mb z)}\cap G$. Since $G$ is the kernel of
the map $\exp\Amap\colon(\C^\times)^m\to \C^\times_N$ induced by
the map of lattices $A\colon\Z^m\to N$, the subgroup $G_{\mb z}$
is the kernel of the composite map
\begin{equation}\label{finitegroup}
  (\C^\times)^{\omega(\mb z)}\hookrightarrow(\C^\times)^m
  \stackrel{\exp\Amap}\longrightarrow \C^\times_N.
\end{equation}
This homomorphism of tori is induced by the map of lattices
$\Z^{\omega(\mb z)}\to\Z^m\to N$, where $\Z^{\omega(\mb
z)}\to\Z^m$ is the inclusion of a coordinate sublattice.

Now let $\Sigma$ be a simplicial fan and $\mb z\in U(\Sigma)$.
Then $\omega(\mb z)=g(\sigma)$ for a cone $\sigma\in\Sigma$.
Therefore, the set of primitive generators $\{\mb a_i\colon
i\in\omega(\mb z)\}$ is linearly independent. Hence, the map
$\Z^{\omega(\mb z)}\to\Z^m\to N$ taking $\mb e_i$ to $\mb a_i$ is
a monomorphism, which implies that the kernel
of~\eqref{finitegroup} is a finite group.

If the fan $\Sigma$ is regular, then $\{\mb a_i\colon
i\in\omega(\mb z)\}$ is a part of basis of~$N$. In this
case~\eqref{finitegroup} is a monomorphism and $G_{\mb z}=\{1\}$.
\end{proof}

\begin{remark}
The closedness of orbits is a necessary condition for the
topological quotient $U(\Sigma)/G$ to be Hausdorff
(Exercise~\ref{hausclosed}). The proof of Theorem~\ref{coxth}~(b)
above uses Separation Lemma (Lemma~\ref{seplemma}); this is
another example of situation when convex-geometric separation
translates into Hausdorffness.

We may consider the following more general setup. Let $\mb
a_1,\ldots,\mb a_m$ be a set of primitive vectors in $N$, and let
$\sK$ be a simplicial complex on~$[m]$. Assume further that for
any $I\in\sK$ the set of vectors $\{\mb a_i\colon i\in I\}$ is
linearly independent. The latter set spans a simplicial cone,
which we denote by $\sigma_I$. The data $\{\sK,\mb a_1,\ldots,\mb
a_m\}$ defines the coordinate subspace arrangement complement
$U(\sK)\subset\C^m$ and the group $G$~\eqref{gexpl}. Furthermore,
the action of $G$ on $U(\sK)$ is almost free (and it is free if
all cones $\sigma_I$ are regular; this is proved in the same way
as Proposition~\ref{freeaction}). However, the quotient $U(\sK)/G$
is Hausdorff precisely when the cones $\{\sigma_I\colon I\in\sK\}$
form a fan~$\Sigma$ (see~\cite[Proposition~II.3.1.6]{a-d-h-l} and
use the previous remark). In this case $\sK=\sK_\Sigma$ and
$U(\sK)=U(\Sigma)$.

To see that the quotient $U(\Sigma)/G$ is Hausdorff in the case
when $\Sigma$ is a simplicial fan, it is not necessary to use the
algebraic criterion of separatedness as in the proof of
Lemma~\ref{sepvar}. Instead, we may modify the argument for the
closedness of orbits in the proof of Theorem~\ref{coxth} and show
that action of $G$ on $U(\Sigma)$ is
\emph{proper}\label{propeacti} (Exercise~\ref{Gprop}), which
guarantees that $U(\Sigma)/G$ is Hausdorff
(Exercise~\ref{hausquot}).
\end{remark}

\begin{example}
Let $V_\sigma$ be the affine toric variety corresponding to an
$n$-dimensional simplicial cone~$\sigma$. We may write
$V_\sigma=V_\Sigma$ where $\Sigma$ is the simplicial fan
consisting of all faces of~$\sigma$. Then $m=n$, $U(\Sigma)=\C^n$,
and $\Amap\colon\Z^n\to N$ is the monomorphism onto the full rank
sublattice generated by $\mb a_1,\ldots,\mb a_n$. Therefore, $G$
is a finite group and $V_\sigma=\C^n/G=\Spec\C[z_1,\ldots,z_n]^G$.

In particular, if we consider the cone $\sigma$ generated by $2\mb
e_1-\mb e_2$ and $\mb e_2$ in $\R^2$ (see Example~\ref{econe}),
then $G$ is $\Z_2$ embedded as $\{(1,1),(-1,-1)\}$
in~$(\C^\times)^2$. The quotient construction realises the
quadratic cone
$V_\sigma=\Spec\C[z_1,z_2]^G=\Spec\C[z_1^2,z_1z_2,z_2^2]$ as a
quotient of $\C^2$ by~$\Z_2$.
\end{example}

\begin{example}
Consider the complete fan of Example~\ref{cp2fa}. Then
\[
  U(\Sigma)=\C^3\setminus\{z_1=z_2=z_3=0\}=\C^3\setminus\{\mathbf0\}
\]
The subgroup $G$ defined by~\eqref{ggrou} is the diagonal
$\C^\times$ in $(\C^\times)^3$. We therefore obtain
$V_\Sigma=U(\Sigma)/G=\C P^2$.
\end{example}

\begin{example}\label{C3minus3lines1}
Consider the fan $\Sigma$ in~$\R^2$ with three one-dimensional
cones generated by the vectors $\mb e_1$, $\mb e_2$ and $-\mb
e_1-\mb e_2$. This fan is not complete, but its one-dimensional
cones generate $\R^2$, so we may apply Theorem~\ref{coxth}. The
simplicial complex $\mathcal K_\Sigma$ consists of 3 disjoint
points. The space $U(\Sigma)=U(\sK_\Sigma)$ is therefore the
complement to 3 coordinate lines in~$\C^3$:
\[
  U(\Sigma)=\C^3\big\backslash\bigl(\{z_1=z_2=0\}\cup
  \{z_1=z_3=0\}\cup\{z_2=z_3=0\}\bigr)
\]
The group $G$ is the diagonal $\C^\times$ in $(\C^\times)^3$.
Hence $V_\Sigma=U(\Sigma)/G$ is a quasiprojective variety obtained
by removing three points from~$\C P^2$.
\end{example}

\subsection*{Exercises.}
\begin{exercise}
The one-dimensional cones of a fan $\Sigma\subset N_\R$
span~$N_\R$ if and only if the toric variety~$V_\Sigma$ does not
have $\C^\times$-factors.
\end{exercise}

\begin{exercise}
If $\Sigma$ is a regular fan of full dimension, then
$G\cong(\C^\times)^{m-n}$.
\end{exercise}

\begin{exercise}
A $G$-invariant %Laurent
monomial has the form
$\prod_{i=1}^mz_i^{\langle\mb u,\mb a_i\rangle}$ for some $\mb
u\in N^*$.
\end{exercise}

\begin{exercise}
If the fan $\Sigma$ is non-simplicial, then there exists a
non-closed orbit of the $G$-action on~$U(\Sigma)$.
\end{exercise}

\begin{exercise}\label{afclosed}
Assume that the action of an algebraic group $G$ on an algebraic
variety $X$ is almost free. Show that then all $G$-orbits are
Zariski closed.
\end{exercise}

\begin{exercise}\label{Gprop}
A $G$-action on $X$ is called \emph{proper} if the map $G\times
X\to X\times X$, $(g,x)\mapsto(gx,x)$ is proper, i.e. the preimage
of a compact subset is compact. Modify the argument in the proof
of Theorem~\ref{coxth}~(b) to show that the $G$-action on
$U(\Sigma)$ is proper whenever $\Sigma$ is simplicial. (Hint: show
that if sequences $\{\mb z^{(k)}\}$ and $\{\mb g^{(k)}\mb
z^{(k)}\}$ have limits in $U(\Sigma)$, then a subsequence of
$\{\mb g^{(k)}\}$ has a limit in~$G$.)
\end{exercise}

\begin{exercise}\label{hausquot}
Show that the quotient $X/G$ of a locally compact Hausdorff space
(e.g. a manifold) $X$ by a proper $G$-action is Hausdorff. Deduce
that the quotient $U(\Sigma)/G$ corresponding to a simplicial fan
$\Sigma$ is Hausdorff. This gives an alternative topological proof
of Lemma~\ref{sepvar} (separatedness of toric varieties) in the
simplicial case.
\end{exercise}

\begin{exercise}
Let $\Sigma$ be a regular fan whose one-dimensional cones
span~$N_\R$. Observe that $U(\Sigma)$ is 2-connected (see
Exercise~\ref{u2con}). Show by considering the exact homotopy
sequence of the principal $G$-bundle $U(\Sigma)\to V_\Sigma$ that
the nonsingular toric variety $V_\Sigma$ is simply connected, and
$H_2(V_\Sigma;\Z)$ is naturally identified with the kernel of the
map $\Amap\colon\Z^m\to N$. Hence $H^2(V_\Sigma;\Z)=\Z^m/N^*$,
which coincides with the \emph{Picard group}\label{picardgro}
of~$V_\Sigma$, see~\cite[\S3.4]{fult93}.
\end{exercise}

\section{Hamiltonian actions and symplectic reduction}\label{symred}
Here we describe projective toric manifolds as \emph{symplectic
quotients} of Hamiltonian torus actions on~$\C^m$. This approach
may be viewed as a symplectic geometry version of the algebraic
quotient construction from the previous section, although
historically the symplectic construction preceded the algebraic
one~\cite{audi91},~\cite{guil94}.

\subsection*{Symplectic reduction}
We briefly review the background material on symplectic geometry,
referring the reader to monographs by Audin~\cite{audi91},
Guillemin~\cite{guil94} and
Guillemin--Ginzburg--Karshon~\cite{g-g-k02} for further details.

A \emph{symplectic manifold}\label{sympmani} is a pair
$(W,\omega)$ consisting of a smooth manifold $W$ and a closed
differential 2-form $\omega$ which is nondegenerate at each point.
The dimension of a symplectic manifold $W$ is necessarily even.

Assume now that a (compact) torus $T$ acts on $W$ preserving the
symplectic form~$\omega$. We denote the Lie algebra of the
torus~$T$ by $\mathfrak t$ (since $T$ is commutative, its Lie
algebra is trivial, but the construction can be generalised to
noncommutative Lie groups). Given an element $\mb v\in\mathfrak
t$, we denote by $X_{\mb v}$ the corresponding $T$-invariant
vector field on~$W$. The torus action is called
\emph{Hamiltonian}\label{hamiltac} if the 1-form $\omega(X_{\mb
v},\,\cdot\,)$ is exact for any $\mb v\in\mathfrak t$. In other
words, an action is Hamiltonian if for any $\mb v\in\mathfrak t$
there exist a function $H_{\mb v}$ on $W$ (called a
\emph{Hamiltonian}) satisfying the condition
\[
  \omega(X_{\mb v},Y)=dH_{\mb v}(Y)
\]
for any vector field $Y$ on~$W$. The function $H_{\mb v}$ is
defined up to addition of a constant.  Choose a basis $\{\mb
e_i\}$ in $\mathfrak t$ and the corresponding Hamiltonians
$\{H_{\mb e_i}\}$. Then the \emph{moment map}\label{dmomentmap}
\[
  \mu\colon W\to\mathfrak t^*,\qquad (x,\mb e_i)\mapsto
  H_{\mb e_i}(x)
\]
(where $x\in W$) is defined. Observe that changing the
Hamiltonians $H_{\mb e_i}$ by constants results in shifting the
image of $\mu$ by a vector in~$\mathfrak t^*$. According to a
theorem of Atiyah~\cite{atiy82} and
Guillemin--Sternberg~\cite{gu-st82}, the image $\mu(W)$ of the
moment map is convex, and if $W$ is compact then $\mu(W)$ is a
convex polytope in~$\mathfrak t^*$.

\begin{example}\label{simcm}
The most basic example is $W=\C^m$ with symplectic form
\[
  \omega=i\sum_{k=1}^m dz_k\wedge d\overline{z}_k=
  2\sum_{k=1}^mdx_k\wedge dy_k,
\]
where $z_k=x_k+iy_k$. The coordinatewise action of the torus
$\T^m$ on $\C^m$ is Hamiltonian. The moment map
$\mu\colon\C^m\to\R^m$ is given by
$\mu(z_1,\ldots,z_m)=(|z_1|^2,\ldots,|z_m|^2)$ (an exercise).
%; the scaling factor 2 in the symplectic form is chosen so as to
%avoid the factor $\frac12$ in the moment map).
The image of the moment map $\mu$ is the positive
orthant~$\R^m_\ge$.
\end{example}

\begin{construction}[symplectic reduction]\label{constsymred}
Assume given a Hamiltonian action of a torus $T$ on a symplectic
manifold~$W$. Assume further that the moment map $\mu\colon
W\to\mathfrak t^*$ is \emph{proper}, i.e. $\mu^{-1}(V)$ is compact
for any compact subset $V\subset\mathfrak t^*$ (this is always the
case if $W$ itself is compact). Let $\mb u\in\mathfrak t^*$ be a
\emph{regular value} of the moment map, i.e. the differential
$\mathcal T_x W\to\mathfrak t^*$ is surjective for all
$x\in\mu^{-1}(\mb u)$. Then the level set $\mu^{-1}(\mb u)$ is a
smooth compact $T$-invariant submanifold in~$W$. Furthermore, the
$T$-action on $\mu^{-1}(\mb u)$ is almost free (an exercise).

Assume now that the $T$-action on $\mu^{-1}(\mb u)$ is free. The
restriction of the symplectic form $\omega$ to $\mu^{-1}(\mb u)$
may be degenerate. However, the quotient manifold $\mu^{-1}(\mb
u)/T$ is endowed with is a unique symplectic form $\omega'$ such
that
\[
  p^*\omega'=i^*\omega,
\]
where $i\colon\mu^{-1}(\mb u)\to W$ is the inclusion and $p\colon
\mu^{-1}(\mb u)\to \mu^{-1}(\mb u)/T$ the projection.

We therefore obtain a new symplectic manifold $(\mu^{-1}(\mb
u)/T,\omega')$ which is referred to as the \emph{symplectic
reduction}, or the \emph{symplectic quotient}\label{sympquoti} of
$(W,\omega)$ by~$T$.

The construction of symplectic reduction works also under milder
assumptions on the action (see~\cite{du-he82} and more references
there), but the generality described here will be enough for our
purposes.
\end{construction}

\subsection*{The toric case}
The algebraic quotient construction describes a toric manifold
$V_\Sigma$ as a quotient of a noncompact set $U(\Sigma)$ by a
noncompact group~$G$. Using symplectic reduction, the projective
toric manifold $V_P$ corresponding to a simple lattice polytope
$P$ can be obtained as the quotient of a compact submanifold
$\zp\subset U(\Sigma_P)$ by a free action of a compact torus.

Let $\Sigma$ be a complete regular fan in $N_\R\cong\R^n$ with $m$
one-dimensional cones generated by $\mb a_1,\ldots,\mb a_m$.
Consider the exact sequence of maximal compact subgroups (tori)
corresponding to exact sequence of algebraic tori~\eqref{ggrou}:
\begin{equation}\label{kgrou}
  1\longrightarrow K \longrightarrow\T^m\stackrel{\exp\Amap}\longrightarrow
  T_N
  \longrightarrow1,
\end{equation}
where $T_N=N\otimes_\Z\mathbb S\cong\T^n$, $\exp\Amap\colon\T^m\to
T_N$ is the map of tori corresponding to the map of lattices
$\Amap\colon\Z^m\to N$, \ $\mb e_i\mapsto\mb a_i$, and
$K=\Ker\Amap$. The group $K$ is isomorphic to $\T^{m-n}$ because
$\Sigma$ is complete and regular.

We let $K\subset\T^m$ act on $\C^m$ by restriction of the
coordinatewise action of $\T^m$. This $K$-action on $\C^m$ is also
Hamiltonian, and the corresponding moment map is given by the
composition
\begin{equation}\label{tmoma}
  \mu_\Sigma\colon\C^m\stackrel\mu\longrightarrow\R^m\longrightarrow \mathfrak
  k^*,
\end{equation}
where $\R^m\to\mathfrak k^*$ is the map of dual Lie algebras
corresponding to the inclusion $K\to\T^m$. By choosing a basis in
the weight lattice $L \subset \mathfrak k^*$ of the
$(m-n)$-torus~$K$ we write the linear map $\R^m\to\mathfrak k^*$
by an integer $(m-n)\times m$-matrix $\varGamma=(\gamma_{jk})$.
Moment map~\eqref{tmoma} is then given by
\[
  (z_1,\ldots,z_m)\longmapsto
  \Bigl(\sum_{k=1}^m\gamma_{1k}|z_k|^2,\ldots,\sum_{k=1}^m
  \gamma_{m-n,k}|z_k|^2\Bigr),
\]
and its level set $\mu_\Sigma^{-1}(\delta)$ corresponding to a
value $\delta=(\delta_1,\ldots,\delta_{m-n})\in\mathfrak k^*$ is
the intersection of $m-n$ Hermitian quadrics in~$\C^m$:
\begin{equation}\label{hquadrics}
  \sum_{k=1}^m\gamma_{jk}|z_k|^2=\delta_j\qquad\text{for }j=1,\ldots,m-n.
\end{equation}

To apply the symplectic reduction we need to identify the regular
values of the moment map $\mu_\Sigma$. We recall from
Section~\ref{combfan} that a polytope~\eqref{ptope} is called
Delzant if its normal fan is regular. If $\Sigma=\Sigma_P$ is the
normal fan of a Delzant polytope~$P$, then the $(m-n)\times
m$-matrix $\varGamma$ above is the one considered in
Construction~\ref{dist}. Namely, the rows of $\varGamma$ form a
basis of linear dependencies between the vectors~$\mb a_i$.

Given a polytope~\eqref{ptope}, we denote by $\zp$ the
intersection of quadrics~\eqref{hquadrics} corresponding
to~$\delta=\varGamma\mb b_P$, where $\mb b_P=(b_1,\ldots,b_m)^t$:
\begin{equation}\label{levelset}
  \zp=\mu_\Sigma^{-1}(\varGamma\mb b_P)=\Bigl\{\mb z\in\C^m\colon\!\!
  \sum_{k=1}^m \gamma_{jk}\bigl(|z_k|^2-b_k\bigr)=0\;\text{ for
  }j=1,\ldots,m-n\Bigr\}.
\end{equation}

\begin{proposition}\label{regvalue}
Assume that $\Sigma=\Sigma_P$ is the normal fan of a Delzant
polytope~$P$ given by~\eqref{ptope}. Then $\delta=\varGamma\mb
b_P$ is a regular value of the moment map~$\mu_\Sigma$.
\end{proposition}
\begin{proof}
We only sketch the proof here; a more general statement will be
proved as Theorem~\ref{propmmap} in the next chapter. We need to
check that this intersection is nondegenerate at each point $\mb
z\in\zp$, i.e. that $\zp$ is a smooth submanifold in~$\C^m$. This
means that the $m-n$ gradient vectors of the left hand sides of
quadratic equations~\eqref{hquadrics} are linearly independent at
each $\mb z\in\zp$. This can be shown to be equivalent to that the
polytope $P$ defined by~\eqref{ptope} is simple.
\end{proof}

As we shall see in Section~\ref{mampol}, the manifold $\zp$ is
$\T^m$-equivariantly homeomorphic to the moment-angle
manifold\label{momeanma} $\mathcal Z_{\sK_P}$ (where $\sK_P$ is
the nerve complex of~$P$, or the underlying simplicial complex of
the normal fan~$\Sigma_P$). Furthermore, we have $\zp\subset
U(\Sigma_P)$. (The reader may either wait until
Chapter~\ref{mamanifolds}, or view these two statements as
exercises.)

We therefore may consider the symplectic quotient of $\C^m$
by~$K\cong\T^{m-n}$, the $2n$-dimensional symplectic manifold
$\mu_\Sigma^{-1}(\varGamma\mb b_P)/K=\zp/K$.

\begin{theorem}\label{symplectictoric}
Let $P$ be a lattice Delzant polytope with the normal fan
$\Sigma_P$, and let $V_P$ be the corresponding projective toric
manifold. The inclusion $\zp\subset U(\Sigma_P)$ induces a
diffeomorphism
\[
  \zp/K
  \stackrel{\cong}\longrightarrow U(\Sigma_P)/G=V_P.
\]
Therefore, any projective toric manifold $V_P$ is obtained as the
symplectic quotient of $\C^m$ by an action of a torus
$K\cong\T^{m-n}$.
\end{theorem}
\begin{proof}
We sketch the proof given in~\cite[Prop.~VI.3.1.1]{audi91}; a
different proof of a more general statement will be given in
Section~\ref{camam}.

Consider the function
\[
  f\colon\C^m\to\R,\quad f(\mb z)=\|\mu_\Sigma(\mb z)-\varGamma\mb b_P\|^2.
\]
It is nonnegative and minimal on the
set~$\zp=\mu_\Sigma^{-1}(\varGamma\mb b_P)$. The only critical
points of $f$ in $U(\Sigma_P)$ are $\mb z\in\zp$. Hence, for any
$\mb z\in U(\Sigma_P)$, the gradient trajectory descending from
$\mb z$ will reach a point in~$\zp$. Furthermore, any gradient
trajectory is contained in a $G$-orbit. We therefore obtain that
each $G$-orbit of $U(\Sigma_P)$ intersects~$\zp$. Finally, it can
be shown that each $G$-orbit intersects $\zp$ at a unique
$K$-orbit, i.e. for each $\mb z\in\zp$ we have that
\[
  G\cdot \mb z\cap\zp=K\cdot \mb z.
\]
The statement follows.
\end{proof}

The toric manifold $V_P$ therefore acquires a symplectic structure
as the symplectic quotient $\mu_\Sigma^{-1}(\varGamma\mb b_P)/K$.
On the other hand, the projective embedding of $V_P$ defined by
the lattice polytope $P$ provides a symplectic form on $V_P$ by
restriction of the standard symplectic form on the complex
projective space. It can be shown~\cite[Appendix~2]{guil94} that
the diffeomorphism from Theorem~\ref{symplectictoric} preserves
the cohomology class of the symplectic form, or equivalently, the
two symplectic structures on $V_P$ are $T_N$-equivariantly
symplectomorphic.

The symplectic quotient $\mu_\Sigma^{-1}(\varGamma\mb b_P)/K$ has
a residual action of the quotient $n$-torus $T_N=\T^m/K$, which is
obviously Hamiltonian. This action is identified, via
Theorem~\ref{symplectictoric}, with the action of the maximal
compact subgroup $T_N\subset\C_N^\times$ on the toric
variety~$V_P$. We denote by $\mu_V\colon V_P\to\mathfrak t_N^*$
the moment map for the Hamiltonian action of $T_N$ on~$V_P$, where
$\mathfrak t_N\cong\R^n$ is the Lie algebra of~$T_N$. It follows
from~\eqref{kgrou} that $\mathfrak t_N^*$ embeds in~$\R^m$ by the
map~$A^*$.

\begin{proposition}\label{tvmomentmap}
The image of the moment map $\mu_V\colon V_P\to\mathfrak t_N^*$ is
the polytope $P$, up to translation.
\end{proposition}
\begin{proof}
Let $\omega$ be the standard symplectic form on $\C^m$ and
$\mu\colon\C^m\to\R^m$ the moment map for the standard action
of~$\T^m$ (see Example~\ref{simcm}). Let $p\colon\zp\to V_P$ be
the quotient projection by the action of~$K$, and let
$i\colon\zp\to\C^m$ be the inclusion, so that the symplectic form
$\omega'$ on $V_P$ satisfies $p^*\omega'=i^*\omega$. Let $H_{\mb
e_i}\colon\C^m\to\R$ be the Hamiltonian of the $\T^m$-action on
$\C^m$ corresponding to the $i$th basis vector $\mb e_i$
(explicitly, $H_{\mb e_i}(\!\mb z)=|z_i|^2$), and let $H_{\mb
a_i}\colon V_P\to\R$ be the Hamiltonian of the $T_N$-action on
$V_P$ corresponding to $\mb a_i\in\mathfrak t_N$. Denote by
$X_{\mb e_i}$ the vector field on $\zp$ generated by $\mb e_i$,
and denote by $Y_{\mb a_i}$ the vector field on $V_P$ generated
by~$\mb a_i$. Observe that $p_*X_{\mb e_i}=Y_{\mb a_i}$. For any
vector field $Z$ on $\zp$ we have
\begin{multline*}
  dH_{\mb e_i}(Z)=i^*\omega(X_{\mb e_i},Z)=
  p^*\omega'(X_{\mb e_i},Z)\\=\omega'(Y_{\mb a_i},p_*Z)=
  dH_{\mb a_i}(p_*Z)=d(p^*H_{\mb a_i})(Z),
\end{multline*}
hence $H_{\mb e_i}=p^*H_{\mb a_i}$ or $H_{\mb e_i}(\!\mb
z)=H_{\!\mb a_i}(p\,(\!\mb z))$ up to constant. By definition of
the moment map this implies that $\mu_V(V_P)\subset\mathfrak
t^*_N\subset\R^m$ is identified with $\mu(\zp)\subset\R^m$ up to
shift by a vector in~$\R^m$. The inclusion $\mathfrak
t^*_N\subset\R^m$ is the map~$A^*$, and $\mu(\zp)=i_{A,\mb
b}(P)=A^*(P)+\mb b$ by definition of $\zp$, see~\eqref{iabrn}
and~\eqref{levelset}. We therefore obtain that there exists $\mb
c\in\R^m$ such that
\[
  A^*(\mu_V(V_P))+\mb c=A^*(P)+\mb b,
\]
i.e. $A^*(\mu_V(V_P))$ and $A^*(P)$ differ by $\mb b-\mb c\in
A^*(\mathfrak t^*_N)$. Since $A^*$ is monomorphic, the result
follows.
\end{proof}

The symplectic quotient $V_P=\mu_\Sigma^{-1}(\varGamma\mb b_P)/K$
with the Hamiltonian action of the $n$-torus $\T^m/K$ is called
the \emph{Hamiltonian toric manifold}\label{hamitoma}
corresponding to a Delzant polytope~$P$.

According to the theorem of Delzant~\cite{delz88}, any
$2n$-dimensional compact connected symplectic manifold $W$ with an
effective Hamiltonian action of an $n$-torus~$T$ is equivariantly
symplectomorphic to a Hamiltonian toric manifold $V_P$, where $P$
is the image of the moment map $\mu\colon W\to\mathfrak t^*$
(whence the name `Delzant polytope').

\begin{remark}
There is a canonical lattice in the dual Lie algebra $\mathfrak
t^*$ (the weight lattice of the torus~$T$), and the moment
polytope $P=\mu(W)\subset\mathfrak t^*$ satisfies the Delzant
condition with respect to this lattice.
\end{remark}

\begin{example}
Let $P=\varDelta^n$ be the standard simplex (see
Example~\ref{simdist}). The cones of the normal fan $\Sigma_P$ are
spanned by the proper subsets of the set of $n+1$ vectors $\{\mb
e_1,\ldots,\mb e_n,-\mb e_1-\cdots-\mb e_n\}$. The groups
$G\cong\C^\times$ and $K\cong\mathbb S^1$ are the diagonal
subgroups in $(\C^\times)^{n+1}$ and $\T^{n+1}$ respectively, and
$U(\Sigma_P)=\C^{n+1}\setminus\{0\}$. The matrix $\varGamma$ is a
row of $n+1$ units. The moment map~\eqref{tmoma} is given by
$\mu_\Sigma(z_1,\cdots,z_{n+1})=|z_1|^2+\ldots+|z_{n+1}|^2$. Since
$\varGamma b_P=1$, the manifold $\zp=\mu^{-1}_P(1)$ is the unit
sphere $\mathbb S^{2n+1}\subset\C^{n+1}$, and $\mathbb
S^{2n+1}/K\cong(\C^{n+1}\setminus\{0\})/G=V_P$ is the complex
projective space~$\C P^n$.
\end{example}

\subsection*{Exercises.}
\begin{exercise}
The moment map $\mu\colon\C^m\to\R^m$ for the coordinatewise
action of $\T^m$ on $\C^m$ is given by
$\mu(z_1,\ldots,z_m)=(|z_1|^2,\ldots,|z_m|^2)$.
\end{exercise}

\begin{exercise}
The $T$-action on $\mu^{-1}(\mb u)$ is almost free.
\end{exercise}

\begin{exercise}
The level set $\zp=\mu_\Sigma^{-1}(\varGamma\mb b_P)$ is
$\T^m$-equivariantly homeomorphic to the moment-angle manifold
$\mathcal Z_{\sK_P}$.
\end{exercise}

\begin{exercise}
Show that $\zp\subset U(\Sigma_P)$.
\end{exercise}

%\chapter*{\ \ Toric varieties: additional topics}
%\section{Face truncations and blowups}

%В разделе про когомологии Дольбо, где используется формальность
% торических многообразий дать ссылку на док-во этого рез-та в главе 8.
%
%Описать нетривиальность расслоений в случае двух квадрик
%
%В последнем разделе появляются small covers
%
%В замечании со ссылкой на {ay-bu11} дать также ссылку на главу 8

\setcounter{chapter}5
\chapter{Geometric structures on moment-angle manifolds}\label{mamanifolds}

In this chapter we study the geometry of moment-angle manifolds,
in its convex, complex-analytic, symplectic and Lagrangian
aspects.

As we have seen in Theorem~\ref{zkman}, the moment-angle complex
$\zk$ corresponding to a triangulated sphere~$\sK$ is a
topological manifold. Moment-angle manifolds corresponding to
simplicial polytopes or, more generally, complete simplicial fans,
are smooth. In the polytopal case a smooth structure arises from
the realisation of $\zk$ by a nonsingular (transverse)
intersection of Hermitian quadrics in $\C^m$, similar to a level
set of the moment map in the construction of symplectic quotients
(see Section~\ref{symred}). The relationship between polytopes and
systems of quadrics is described in terms of Gale duality (see
Sections~\ref{intquad} and~\ref{mampol}).

Another way to give $\zk$ a smooth structure is to realise it as
the quotient of an open subset in $\C^m$ (the complement $U(\sK)$
of the coordinate subspace arrangement defined by~$\sK$) by an
action of the multiplicative group~$\R^{m-n}_>$. As in the case of
the quotient construction of toric varieties
(Section~\ref{algtq}), the quotient of a non-compact manifold
$U(\sK)$ by the action of a non-compact group $\R^{m-n}_>$ is
Hausdorff precisely when $\sK$ is the underlying complex of a
simplicial fan.

If $m-n=2\ell$ is even, then the action of the real group
$\R^{m-n}_>$ on $U(\sK)$ can be turned into a holomorphic action
of a complex (but not algebraic) group isomorphic to~$\C^{\ell}$.
In this way the moment-angle manifold $\zk\cong U(\sK)/\C^\ell$
acquires a complex-analytic structure. The resulting family of
non-K\"ahler complex manifolds generalises the well-known series
as of Hopf and Calabi--Eckmann manifolds~\cite{cher56}, as well as
LVM-manifolds (a class of complex structures on nonsingular
intersections of Hermitian quadrics arising in holomorphic
dynamics as transverse sets to certain complex
foliations~\cite{lo-ve97},~\cite{meer00}).

The intersections of Hermitian quadrics defining polytopal
moment-angle manifolds can be also used to construct new families
of Lagrangian submanifolds in~$\C^m$, $\C P^m$ and toric
varieties.

Particularly interesting examples of geometric structures
(nonsingular intersections of quadrics, polytopal moment-angle
manifolds, non-K\"ahler LVM-manifolds, Hamiltonian-minimal
Lagrangian submanifolds) arise from Delzant polytopes. These
polytopes are abundant in toric topology (see
Section~\ref{projtvpol}).
%providing series of explicit examples of
%moment-angle manifolds with interesting geometric structures.

\section{Intersections of quadrics}\label{intquad}
Here we describe the correspondence between convex polyhedra and
intersections of quadrics. It will be used in the next section to
define a smooth structure on moment-angle manifolds coming from
polytopes.
%We continue assuming that $P$ has a vertex.

\subsection*{From polyhedra to quadrics}
The following construction originally appeared
in~\cite[Construction~3.1.8]{bu-pa02} (see
also~\cite[\S3]{b-p-r07}):
\begin{construction}
Consider a presentation of a convex polyhedron
\begin{equation}
\label{ptope1}
  P=P(A,\mb b)=\bigl\{\mb x\in\R^n\colon\langle\mb a_i,\mb x\rangle+b_i\ge0\quad\text{for }
  i=1,\ldots,m\bigr\}.
\end{equation}
Assume that $\mb a_1,\ldots,\mb a_m$ span $\R^n$ (i.e., $P$ has a
vertex) and recall the map
\[
  i_{A,\mb b}\colon \R^n\to\R^m,\quad i_{A,\mb b}(\mb x)=A^*\mb x+\mb b
  =\bigl(\langle\mb a_1,\mb x\rangle+b_1,\ldots,\langle\mb a_m,\mb x\rangle+b_m\bigr)^t
\]
(see Construction~\ref{dist}). It embeds $P$ into $\R^m_\ge$. We
define the space $\mathcal Z_{A,\mb b}$ from the commutative
diagram
\begin{equation}\label{cdiz}
\begin{CD}
  \mathcal Z_{A,\mb b} @>i_{\mathcal Z}>>\C^m\\
  @VVV\hspace{-0.2em} @VV\mu V @.\\
  P @>i_{A,\mb b}>> \R^m_\ge
\end{CD}
\end{equation}
where $\mu(z_1,\ldots,z_m)=(|z_1|^2,\ldots,|z_m|^2)$. The torus
$\T^m$ acts on $\mathcal Z_{A,\mb b}$ with quotient $P$, and
$i_{\mathcal Z}$ is a $\T^m$-equivariant embedding.

By replacing $y_k$ by $|z_k|^2$ in the equations defining the
affine plane $i_{A,\mb b}(\R^n)$ (see~\eqref{iabrn}) we obtain
that $\mathcal Z_{A,\mb b}$ embeds into $\C^m$ as the set of
common zeros of $m-n$ quadratic equations (\emph{Hermitian
quadrics}):
\begin{equation}\label{zpqua}
  i_{\mathcal Z}(\mathcal Z_{A,\mb b})=\Bigl\{\mb z\in\C^m\colon\sum_{k=1}^m\gamma_{jk}|z_k|^2=
  \sum_{k=1}^m\gamma_{jk}b_k,\;\text{ for }1\le j\le m-n\Bigr\}.
\end{equation}
%We shall often not distinguish between $\mathcal Z_{A,\mb b}$ and
%its embedding $i_{\mathcal Z}(\mathcal Z_{A,\mb b})\subset\C^m$.
\end{construction}

The following property of the space $\mathcal Z_{A,\mb b}$ follows
easily from its construction.

\begin{proposition}\label{easyzp}
%\begin{itemize}
%\item[(a)]
Given a point $z\in\mathcal Z_{A,\mb b}$, the $j$th coordinate of
$i_{\mathcal Z}(z)\in\C^m$ vanishes if and only if $z$ projects
onto a point $\mb x\in P$ such that $\mb x\in F_j$.
%
%\item[(b)] Adding a redundant inequality to~\eqref{ptope1} results in multiplying $\mathcal Z_{A,\mb b}$ by a circle.
%\end{itemize}
\end{proposition}

\begin{theorem}%[{\cite[Cor.~3.9]{bu-pa02}}]
\label{zpsmooth} The following conditions are equivalent:
\begin{itemize}
\item[(a)]
presentation~\eqref{ptope1} determined by the data $(A,\mb b)$ is
generic;
\item[(b)]
the intersection of quadrics in~\eqref{zpqua} is nonempty and
nonsingular, so that $\mathcal Z_{A,\mb b}$ is a smooth manifold
of dimension $m+n$.
\end{itemize}
Furthermore, under these conditions the embedding $i_{\mathcal
Z}\colon\mathcal Z_{A,\mb b}\to\C^m$ has $\T^m$-equivariantly
framed normal bundle; a trivialisation is defined by a choice of
matrix $\varGamma=(\gamma_{jk})$ in~\eqref{iabrn}.
\end{theorem}
\begin{proof}
For simplicity we identify the space $\mathcal Z_{A,\mb b}$ with
its embedded image $i_{\mathcal Z}(\mathcal Z_{A,\mb
b})\subset\C^m$. We calculate the gradients of the $m-n$ quadrics
in~\eqref{zpqua} at a point $\mb
z=(x_1,y_1,\ldots,x_m,y_m)\in\mathcal Z_{A,\mb b}$, where
$z_k=x_k+iy_k$:
\begin{equation}\label{grve}
  2\left(\gamma_{j1}x_1,\,\gamma_{j1}y_1,\,\dots,\,\gamma_{jm}x_m,\,\gamma_{jm}y_m\right),\quad
  1\le j\le m-n.
\end{equation}
These gradients form the rows of the $(m-n)\times 2m$ matrix
$2\varGamma\varDelta$, where
\[
\varDelta\;=\;
\begin{pmatrix}
  x_1&y_1& \ldots & 0  &0 \\
  \vdots & \vdots & \ddots & \vdots & \vdots\\
  0  & 0 & \ldots &x_m &y_m
\end{pmatrix}.
\]
Let $I=\{i_1,\ldots,i_k\}=\{i\colon z_i=0\}$ be the set of zero
coordinates of~$\mb z$. Then the rank of the gradient matrix
$2\varGamma\varDelta$ at $\mb z$ is equal to the rank of the
$(m-n)\times(m-k)$ matrix $\varGamma_{\widehat I}$ obtained by
deleting the columns $i_1,\ldots,i_k$ from~$\varGamma$.

Now let~\eqref{ptope1} be a generic presentation. By
Proposition~\ref{easyzp}, a point $\mb z$ with
$z_{i_1}=\cdots=z_{i_k}=0$ projects to a point in
$F_{i_1}\cap\cdots\cap F_{i_k}\ne\varnothing$. Hence the vectors
$\mb a_{i_1},\ldots,\mb a_{i_k}$ are linearly independent. By
Theorem~\ref{galespan}, the rank of $\varGamma_{\widehat I}$
is~$m-n$. Therefore, the intersection of quadrics~\eqref{zpqua} is
nonsingular.

On the other hand, if~\eqref{ptope1} is not generic, then there is
a point $\mb z\in\mathcal Z_{A,\mb b}$ such that the vectors
$\{\mb a_{i_1},\ldots,\mb a_{i_k}\colon
z_{i_1}=\cdots=z_{i_k}=0\}$ are linearly dependent. By
Theorem~\ref{galespan}, the columns of the corresponding matrix
$\varGamma_{\widehat I}$ do not span~$\R^{m-n}$, so its rank is
less than~$m-n$ and the intersection of quadrics~\eqref{zpqua} is
degenerate at~$\mb z$.

The last statement follows from the fact that $\mathcal Z_{A,\mb
b}$ is a nonsingular intersection of quadratic surfaces, each of
which is $\T^m$-invariant.
\end{proof}

\subsection*{From quadrics to polyhedra}
This time we start with an intersection of $m-n$ Hermitian
quadrics in~$\C^m$:
\begin{equation}\label{zgamma}
  \mathcal Z_{\varGamma,\delta}=\Bigl\{\mb z=(z_1,\ldots,z_m)\in\C^m\colon
  \sum_{k=1}^m\gamma_{jk}|z_k|^2=\delta_j,\quad\text{for }
  1\le j\le m-n\Bigr\}.
\end{equation}
The coefficients of quadrics form an $(m-n)\times m$-matrix
$\varGamma=(\gamma_{jk})$, and we denote its column vectors by
$\gamma_1,\ldots,\gamma_m$. We also consider the column vector of
the right hand sides,
$\delta=(\delta_1,\ldots,\delta_{m-n})^t\in\R^{m-n}$.

These intersections of quadrics are considered up to \emph{linear
equivalence}, which corresponds to applying a nondegenerate linear
transformation of $\R^{m-n}$ to $\varGamma$ and~$\delta$.
Obviously, such a linear equivalence does not change the
set~$\mathcal Z_{\varGamma,\delta}$.

We denote by $\R_\ge\langle\gamma_1,\ldots,\gamma_m\rangle$ the
cone generated by the vectors $\gamma_1,\ldots,\gamma_m$ (i.e.,
the set of linear combinations of these vectors with nonnegative
coefficients).

A version of the following proposition appeared in~\cite{lope89},
and the proof below is a modification of the argument
in~\cite[Lemma~0.3]{bo-me06}. It allows us to determine the
nondegeneracy of an intersection of quadrics directly from the
data $(\varGamma,\delta)$:

\begin{proposition}\label{zgsmooth}
The intersection of quadrics $\mathcal Z_{\varGamma,\delta}$ given
by~\eqref{zgamma} is nonempty and nonsingular if and only if the
following two conditions are satisfied:
\begin{itemize}
\item[(a)] $\delta\in
\R_\ge\langle\gamma_1,\ldots,\gamma_m\rangle$;\\[-0.7\baselineskip]

\item[(b)] if $\delta\in\R_\ge\langle
\gamma_{i_1},\ldots\gamma_{i_k}\rangle$, then $k\ge m-n$.
\end{itemize}
Under these conditions, $\mathcal Z_{\varGamma,\delta}$ is a
smooth submanifold in $\C^m$ of dimension $m+n$, and the vectors
$\gamma_1,\ldots,\gamma_m$ span~$\R^{m-n}$.
\end{proposition}
\begin{proof}
First assume that (a) and (b) are satisfied. Then (a) implies that
$\mathcal Z_{\varGamma,\delta}\ne\varnothing$. Let $\mb
z\in\mathcal Z_{\varGamma,\delta}$. Then the rank of the matrix of
gradients of~\eqref{zgamma} at $\mb z$ is equal to
$\mathop{\mathrm{rk}}\{\gamma_k\colon z_k\ne0\}$. Since $\mb
z\in\mathcal Z_{\varGamma,\delta}$, the vector $\delta$ is in the
cone generated by those $\gamma_k$ for which $z_k\ne0$. By the
Carath\'eodory Theorem~(see~\cite[\S1.6]{zieg95}), $\delta$
belongs to the cone generated by some $m-n$ of these vectors, that
is,
$\delta\in\R_\ge\langle\gamma_{k_1},\ldots,\gamma_{k_{m-n}}\rangle$,
where $z_{k_i}\ne0$ for $i=1,\ldots,m-n$. Moreover, the vectors
$\gamma_{k_1},\ldots,\gamma_{k_{m-n}}$ are linearly independent
(otherwise, again by the Carath\'eodory Theorem, we obtain a
contradiction with~(b)). This implies that the $m-n$ gradients of
quadrics in~\eqref{zgamma} are linearly independent at~$\mb z$,
and therefore $\mathcal Z_{\varGamma,\delta}$ is smooth and
$(m+n)$-dimensional.

To prove the other implication we observe that if (b) fails, that
is, $\delta$ is in the cone generated by some $m-n-1$ vectors of
$\gamma_1,\ldots,\gamma_m$, then there is a point $\mb
z\in\mathcal Z_{\varGamma,\delta}$ with at least $n+1$ zero
coordinates. The gradients of quadrics in~\eqref{zgamma} cannot be
linearly independent at such~$\mb z$.
\end{proof}

The torus $\T^m$ acts on $\mathcal Z_{\varGamma,\delta}$, and the
quotient $\mathcal Z_{\varGamma,\delta}/\T^m$ is identified with
the set of nonnegative solutions of the system of $m-n$ linear
equations
\begin{equation}\label{linsys}
  \sum_{k=1}^m\gamma_ky_k=\delta.
\end{equation}
This set can be described as a convex polyhedron $P(A,\mb b)$
given by~\eqref{ptope1}, where $(b_1,\ldots,b_m)$ is any solution
of~\eqref{linsys} and the vectors $\mb a_1,\ldots,\mb a_m\in\R^n$
form the transpose of a basis of solutions of the homogeneous
system $\sum_{k=1}^m\gamma_ky_k=\mathbf 0$.  We refer to $P(A,\mb
b)$ as the \emph{associated polyhedron} of the intersection of
quadrics~$\mathcal Z_{\varGamma,\delta}$. If the vectors
$\gamma_1,\ldots,\gamma_m$ span $\R^{m-n}$, then $\mb
a_1,\ldots,\mb a_m$ span~$\R^n$. In this case the two vector
configurations are Gale dual.

We summarise the results and constructions of this section as
follows:

\begin{theorem}\label{polquad}
A presentation of a polyhedron
\[
  P=P(A,\mb b)=\bigl\{\mb x\in\R^n\colon\langle\mb a_i,\mb
  x\rangle+b_i\ge0\quad\text{for }
  i=1,\ldots,m\bigr\}
\]
(with $\mb a_1,\ldots,\mb a_m$ spanning~$\R^n$) defines an
intersection of Hermitian quadrics
\[
  \mathcal Z_{\varGamma,\delta}=\Bigl\{\mb
  z=(z_1,\ldots,z_m)\in\C^m\colon
  \sum_{k=1}^m\gamma_{jk}|z_k|^2=\delta_j\quad\text{for }
  j=1,\ldots,m-n\Bigr\}.
\]
(with $\gamma_1,\ldots,\gamma_m$ spanning~$\R^{m-n}$) uniquely up
to a linear isomorphism of~$\R^{m-n}$, and an intersection of
quadrics~$\mathcal Z_{\varGamma,\delta}$ defines a
presentation~$P(A,\mb b)$ uniquely up to an isomorphism of~$\R^n$.

The systems of vectors $\mb a_1,\ldots,\mb a_m\in\R^n$ and
$\gamma_1,\ldots,\gamma_m\in\R^{m-n}$ are Gale dual, and the
vectors $\mb b\in\R^{m}$ and $\delta\in\R^{m-n}$ are related by
the identity $\delta=\varGamma\mb b$.

The intersection of quadrics $\mathcal Z_{\varGamma,\delta}$ is
nonempty and nonsingular if and only if the presentation $P(A,\mb
b)$ is generic.
\end{theorem}

\begin{example}[$m=n+1$: one quadric]\label{mamsimplex}
If presentation~\eqref{ptope1} is generic and $P$ is bounded, then
$m\ge n+1$. The case $m=n+1$ corresponds to a simplex. If
$P=P(A,\mb b)$ is the standard simplex (see Example~\ref{simdist})
we obtain
\[
  \mathcal Z_{A,\mb b}=\mathbb S^{2n+1}=\{\mb z\in\C^{n+1}\colon
  |z_1|^2+\cdots+|z_{n+1}|^2=1\}.
\]

More generally, a presentation~\eqref{ptope1} with $m=n+1$ and
$\mb a_1,\ldots,\mb a_n$ spanning $\R^n$ can be taken by an
isomorphism of $\R^n$ to the form
\[
  P=\bigl\{\mb x\in\R^n\colon x_i+b_i\ge0\quad\text{for }
  i=1,\ldots,n,\quad\text{and }
  -c_1x_1-\cdots-c_nx_n+b_{n+1}\ge0\bigr\}.
\]
We therefore have $\varGamma=(c_1\cdots c_n\, 1)$, and $\mathcal
Z_{A,\mb b}$ is given by the single equation
\[
  c_1|z_1|^2+\cdots+c_n|z_n|^2+|z_{n+1}|^2=c_1b_1+\cdots+c_nb_n+b_{n+1}.
\]
If the presentation is generic and bounded, then $\mathcal
Z_{A,\mb b}$ is nonempty, nonsingular and bounded by
Theorem~\ref{zpsmooth}. This implies that all $c_i$ and the right
hand side above are positive, and $\mathcal Z_{A,\mb b}$ is an
ellipsoid.
\end{example}

\section{Moment-angle manifolds from polytopes}\label{mampol}
In this section we identify the polytopal moment-angle manifold
$\mathcal Z_{\mathcal K_P}$ (the moment-angle complex
corresponding to the nerve complex $\sK_P$ of a simple
polytope~$P$) with the intersections of quadrics $\mathcal
Z_{A,\mb b}$~\eqref{zpqua}.
% defined by~$P=P(A,\mb b)$.

A $\T^m$-equivariant homeomorphism $\mathcal Z_{\mathcal
K_P}\cong\mathcal Z_{A,\mb b}$ will be established using the
following construction of an identification space, which goes back
to the work of Vinberg~\cite{vinb71} on Coxeter groups and was
presented in the form described below in the work of Davis and
Januszkiewicz~\cite{da-ja91}. This was the first construction of
what later became known as the moment-angle manifold.

\begin{construction}
Let $[m]=\{1,\ldots,m\}$ be the standard $m$-element set. For each
$I\subset[m]$ we consider the coordinate subtorus
\[
  \T^I=\{(t_1,\ldots,t_m)\in\T^m\colon t_j=1\text{ for
  }j\notin I\}\;\subset\;\T^m.
\]
In particular, $\T^\varnothing$ is the trivial
subgroup~$\{1\}\subset\T^m$.

Define the map $\R_\ge\times\T\to\C$ by $(y,t)\mapsto yt$. Taking
product we obtain a map $\R^m_\ge\times\T^m\to\C^m$. The preimage
of a point $\mb z\in\C^m$ under this map is $\mb
y\times\T^{\omega(\mb z)}$, where $y_i=|z_i|$ for $1\le i\le m$
and $\omega(\mb z)=\{i\colon z_i=0\}\subset[m]$ is the set of zero
coordinates of~$\mb z$. Therefore, $\C^m$ can be identified with
the quotient space
\begin{equation}\label{Cmids}
  \R^m_\ge\times\T^m/{\sim\:}\quad\text{where }(\mb y,\mb t_1)\sim(\mb y,\mb
  t_2)\text{ if }\mb t_1^{-1}\mb t_2\in\T^{\omega(\mb y)}.
\end{equation}

Given $\mb x\in P$, set $I_{\mb x}=\{i\in[m]\colon{\mb x}\in
F_i\}$ (the set of facets containing~$\mb x$).
\end{construction}

\begin{proposition}\label{zpids}
The space $\mathcal Z_{A,\mb b}$ given by intersection of
quadrics~\eqref{zpqua} corresponding to a presentation $P=P(A,\mb
b)$ is $\T^m$-equivariantly homeomorphic to the quotient
\[
  P\times\T^m/{\sim\:}\quad\text{where }
  (\mb x,\mb t_1)\sim(\mb x,\mb t_2)\:\text{ if }\:\mb t_1^{-1}\mb t_2\in
  \T^{I_{\mb x}}.
\]
\end{proposition}
\begin{proof}
Using~\eqref{cdiz}, we identify the space $\mathcal Z_{A,\mb b}$
with $i_{A,\mb b}(P)\times\T^m/{\sim\:}$, where $\sim$ is the
equivalence relation from~\eqref{Cmids}. A point $\mb x\in P$ is
mapped by $i_{A,\mb b}$ to $\mb y\in\R^m_\ge$ with $I_{\mb
x}=\omega(\mb y)$.
\end{proof}

An important corollary of this construction is that the
topological type of the intersection of quadrics $\mathcal
Z_{A,\mb b}$ depends only on the combinatorics of~$P$:

\begin{proposition}\label{zpcomb}
Assume given two generic presentations:
\[
  P=\bigl\{\mb x\in\R^n\colon(A^*\mb x+\mb b)_i\ge0\bigr\}\quad\text{and}
  \quad P'=\bigl\{\mb x\in\R^n\colon(A'^*\mb x+\mb b')_i\ge 0\bigr\}
\]
such that $P$ and $P'$ are combinatorially equivalent simple
polytopes.
\begin{itemize}
\item[(a)] If both presentations are irredundant, then the corresponding manifolds
$\mathcal Z_{A,\mb b}$ and $\mathcal Z_{A',\mb b'}$ are
$\T^m$-equivariantly homeomorphic.

\item[(b)] If the second presentation is obtained from the first one by
adding $k$ redundant inequalities, then $\mathcal Z_{A',\mb b'}$
is homeomorphic to a product of $\mathcal Z_{A,\mb b}$ and a
$k$-torus~$T^k$.
\end{itemize}
\end{proposition}
\begin{proof} (a) By Proposition~\ref{zpids}, we have
$\mathcal Z_{A,\mb b}\cong P\times\T^m/{\sim\:}$ and $\mathcal
Z_{A',\mb b'}\cong P'\times\T^m/{\sim\:}$. If both presentations
are irredundant, then any $F_i$ is a facet of $P$, and the
equivalence relation $\sim$ depends on the face structure of $P$
only. Therefore, any homeomorphism $P\to P'$ preserving the face
structure extends to a $\T^m$-homeomorphism
$P\times\T^m/{\sim\:}\to P'\times\T^m/{\sim\:}$.

(b) Suppose the first presentation has $m$ inequalities, and the
second has~$m'$ inequalities, so that $m'-m=k$. Let $J\subset[m']$
be the subset corresponding to the added redundant inequalities;
we may assume that $J=\{m+1,\ldots,m'\}$. Since $F_j=\varnothing$
for any $j\in J$, we have $I_{\mb x}\cap J=\varnothing$ for any
$\mb x\in P'$. Therefore, the equivalence relation $\sim$ does not
affect the factor $\T^J\subset\T^{m'}$, and we have
\[
  \mathcal Z_{A',\mb b'}\cong P'\times\T^{m'}/{\sim\:}\cong
  (P\times\T^m/{\sim\:})\times\T^J\cong\mathcal Z_{A,\mb b}\times
  T^k.\qedhere
\]
\end{proof}

\begin{remark}
A $\T^m$-homeomorphism in Proposition~\ref{zpcomb}~(a) can be
replaced by a $\T^m$-diffeomorphism (with respect to the smooth
structures of Theorem~\ref{zpsmooth}), but the proof is more
technical. It follows from the fact that two simple polytopes are
combinatorially equivalent if and only if they are diffeomorphic
as `smooth manifolds with corners'. For an alternative argument,
see~\cite[Corollary~4.7]{bo-me06}.

Statement (a) remains valid without assuming that the presentation
is generic or bounded, although $\mathcal Z_{A,\mb b}$ is not a
manifold in this case.
\end{remark}

Now we recall the moment-angle manifold $\mathcal Z_{\mathcal
K_P}$ corresponding to the nerve complex $\sK_P$ of a simple
polytope~$P$ (see Section~\ref{defzk} and Example~\ref{polsph}).
Observe that in our formalism, redundant inequalities in a
presentation of $P$ correspond to ghost vertices of~$\sK_P$.

\begin{theorem}
Let~\eqref{ptope1} be a generic bounded presentation, so that
$P=P(A,\mb b)$ is a simple $n$-polytope. The moment-angle manifold
$\mathcal Z_{\sK_P}$ is $\T^m$-\-equi\-va\-ri\-ant\-ly
homeomorphic to the intersection of quadrics $\mathcal Z_{A,\mb
b}$ given by~\eqref{zpqua}.
\end{theorem}
\begin{proof}
Recall from Construction~\ref{constrmac} that the moment-angle
complex $\mathcal Z_{\sK_P}$ is defined from  a diagram similar
to~\eqref{cdiz}, in which the bottom map is replaced by the
piecewise linear embedding $c_P\colon P\to\I^m$ from
Construction~\ref{cubpol}:
\begin{equation}\label{cdwiz}
\begin{CD}
  \mathcal Z_{\sK_P} @>>>\D^m\\
  @VVV\hspace{-0.2em} @VV\mu V @.\\
  P @>c_P>> \I^m
\end{CD}
\end{equation}

As we have seen in Proposition~\ref{zpids}, the intersection of
quadrics $\mathcal Z_{A,\mb b}$ is $\T^m$-homeomorphic to the
identification space
\[
  P\times\T^m/{\sim\:}\quad\text{where }
  (\mb x,\mb t_1)\sim(\mb x,\mb t_2)\:\text{ if }\:\mb t_1^{-1}\mb t_2\in
  \T^{I_{\mb x}}.
\]
By restricting the equivalence relation~\eqref{Cmids} to
$\D^m\subset\C^m$ we obtain that
\[
  \D^m\cong   \I^m\times\T^m/{\sim\:}\quad\text{where }(\mb y,\mb t_1)\sim(\mb y,\mb
  t_2)\text{ if }\mb t_1^{-1}\mb t_2\in\T^{\omega(\mb y)}.
\]
As in the proof of Proposition~\ref{zpids}, the space $\mathcal
Z_{\sK_P}$ is identified with $c_P(P)\times\T^m/{\sim\:}$. A point
$\mb x\in P$ is mapped by $c_P$ to $\mb y\in\I^m$ with $I_{\mb
x}=\omega(\mb y)=\{i\in[m]\colon\mb x\in F_i\}$. We therefore
obtain that both $\mathcal Z_{\sK_P}$ and $\mathcal Z_{A,\mb b}$
are $\T^m$-homeomorphic to~$P\times\T^m/{\sim\:}$.
\end{proof}

\begin{corollary}\label{zpss}
The moment-angle manifold $\mathcal Z_{\sK_P}$ corresponding to
the nerve complex $\sK_P$ of a simple polytope $P$ has a smooth
structure in which the $\T^m$-action is smooth.
\end{corollary}

\begin{definition}\label{polymamanif}
Given a bounded generic presentation~\eqref{ptope1} defining a
simple polytope~$P$, we shall use the common notation $\zp$ for
both the moment-angle manifold $\mathcal Z_{\sK_P}$ and the
intersection of quadrics~\eqref{zpqua}. We refer to $\zp$ as a
\emph{polytopal moment-angle manifold}. We endow $\zp$ with the
smooth structure coming from the nonsingular intersection of
quadrics.
\end{definition}

\begin{remark}
As it is shown in~\cite[Corollary~4.7]{bo-me06} a $\T^m$-invariant
smooth structure on $\zp$ is unique.
%provided that it satisfies compatibility assumptions with respect
%to the $\T^m$-action and the combinatorics of~$P$ (namely, the
%action is smooth and the normal bundles of the submanifolds
%corresponding to the faces of $P$ are equivariantly trivial).
If the condition of the invariance under the torus action is
dropped, then a smooth structure on $\zp$ is not unique, as is
shown by examples of odd-dimensional spheres and products of
spheres. It would be interesting to relate exotic smooth
structures with the construction of moment-angle complexes.
\end{remark}

\begin{remark}
If the polytope $P=P(A,\mb b)$ is not simple, then moment-angle
complex $\mathcal Z_{\sK_P}$ corresponding to the nerve complex
$\sK_P$ is not homeomorphic to the (singular) intersection of
quadrics $\mathcal Z_{A,\mb b}$. However, the two spaces are
homotopy equivalent (see~\cite{ay-bu11}).
\end{remark}

An intersection of quadrics representing $\zp$ can be chosen more
canonically:

\begin{proposition}\label{malink}
The moment-angle manifold $\zp$ is $\T^m$-equivariantly
diffeomorphic to a nonsingular intersection of quadrics of the
following form:
\begin{equation}\label{link}
  \left\{\begin{array}{lrcl}
  \mb z\in\C^m\colon&\sum_{k=1}^m|z_k|^2&=&1,\\[1mm]
  &\sum_{k=1}^m\mb g_k|z_k|^2&=&\mathbf0,
  \end{array}\right\}
\end{equation}
where $(\mb g_1,\ldots,\mb g_m)\subset\R^{m-n-1}$ is a
combinatorial Gale diagram of~$P^*$.
\end{proposition}
\begin{proof}
It follows from Proposition~\ref{propcf} that $\zp$ is given by
\[
  \left\{\begin{array}{lc}
  \mb z\in\C^m\colon&\gamma_{11}|z_1|^2+\cdots+\gamma_{1m}|z_m|^2=c,\\[1mm]
  &\mb g_1|z_1|^2+\cdots+\mb g_m|z_m|^2=\mathbf0,
  \end{array}\right\}
\]
where $\gamma_{1k}$ and $c$ are positive. Divide the first
equation by~$c$, and then replace each $z_k$ by $\sqrt\frac
c{\gamma_{1k}}z_k$. As a result, each $\mb g_k$ is multiplied by a
positive number, so that $(\mb g_1,\ldots,\mb g_m)$ remains to be
a combinatorial Gale diagram for~$P^*$.
\end{proof}

By adapting Proposition~\ref{zgsmooth} to the special case of
quadrics~\eqref{link}, we obtain

\begin{proposition}\label{nondeglink}
The intersection of quadrics given by~\eqref{link} is nonempty
nonsingular if and only if the following two conditions are
satisfied:
\begin{itemize}
\item[(a)] $\mathbf0\in
\conv(\mb g_1,\ldots,\mb g_m)$;

\item[(b)] if ${\mathbf 0}\in\conv(\mb g_{i_1},\ldots,\mb g_{i_k})$,
then $k\ge m-n$.
\end{itemize}
\end{proposition}

Following~\cite{bo-me06}, we refer to a nonsingular
intersection~\eqref{link} of $m-n-1$ homogeneous quadrics with a
unit sphere in~$\C^m$ as a \emph{link}\label{linkintq}. We
therefore obtain that the class of links coincides with the class
of polytopal moment-angle manifolds.

As we have seen in Example~\ref{mamsimplex}, the moment-angle
manifold corresponding to an $n$-simplex is diffeomorphic to a
sphere~$S^{2n+1}$. This is also the link of an empty system of
homogeneous quadrics, corresponding to the case $m=n+1$.

\begin{example}[$m=n+2$: two quadrics]\label{prodsimex}
A polytope $P$ defined by $m=n+2$ inequalities either is
combinatorially equivalent to a product of two simplices (when
there are no redundant inequalities), or is a simplex (when one
inequality is redundant). In the case $m=n+2$ the
link~\eqref{link} has the form
\[
  \left\{\begin{array}{lc}
  \mb z\in\C^m\colon&|z_1|^2+\cdots+|z_m|^2=1,\\[1mm]
  &g_1|z_1|^2+\cdots+g_m|z_m|^2=0,
  \end{array}\right\}
\]
where $g_k\in\R$. Condition~(b) of Proposition~\ref{nondeglink}
implies that all $g_i$ are nonzero; assume that there are $p$
positive and $q=m-p$ negative numbers among them. Then
condition~(a) implies that $p>0$ and $q>0$. Therefore, the link is
the intersection of the cone over a product of two ellipsoids of
dimensions $2p-1$ and $2q-1$ (given by the second quadric) with a
unit sphere of dimension $2m-1$ (given by the first quadric). Such
a link is diffeomorphic to $S^{2p-1}\times S^{2q-1}$. The case
$p=1$ or $q=1$ corresponds to one redundant inequality. In the
irredundant case ($P$ is a product
$\varDelta^{p-1}\times\varDelta^{q-1}$, $p,q>1$) we obtain that
$\zp\cong S^{2p-1}\times S^{2q-1}$.
\end{example}

The case of three quadrics was resolved by Lopez de
Medrano~\cite{lope89} in~1989. Here is a reformulation of his
result in terms of moment-angle manifolds:

\begin{theorem}\label{zp3quad}
Let $\zp$ be the moment-angle manifold given by a nonempty and
nodegenerate intersection of three quadrics~\eqref{link}, i.e.
$m=n+3$. Then $\zp$ is diffeomorphic to a product of three
odd-dimensional spheres or to a connected sum of products of
spheres with two spheres in each product.
\end{theorem}

The proof of this theorem uses Gale duality, surgery theory and
the $h$-cobordism Theorem. The original statement of~\cite{lope89}
contained some restrictions on the types of quadrics, which later
were lifted in~\cite{gi-lo13}. The moment-angle manifold $\zp$ is
diffeomorphic to a product of three odd-dimensional spheres
precisely when $P$ is combinatorially equivalent to a product of
three simplices. In all other cases, the manifold $\zp$ given by
three quadrics is diffeomorphic to a connected sum of the form
$\mathop{\#}_{k=3}^{m-2}(S^k\times S^{2m-3-k})^{\#q_k}$. The
numbers $q_k$ of products $S^k\times S^{2m-3-k}$ in the connected
sum can be described explicitly in terms of the planar Gale
diagram of the associated $n$-polytope $P$ with $m=n+3$ facets.

Theorem~\ref{zpstacked} together with Theorem~\ref{zp3quad} and
the previous examples gives a description of the topology of
moment-angle manifolds $\zp$ corresponding to the following
classes of simple $n$-polytopes~$P$: dual stacked polytopes
(including polygons), polytopes with $m\le n+3$ facets, and
products of them. More examples of polytopes $P$ whose
corresponding manifolds $\zp$ are diffeomorphic to a connected sum
of sphere products were described in~\cite{gi-lo13} (these include
some dual cyclic polytopes). In general, the topology of
moment-angle manifolds is much more complicated than in these
series of examples (see Section~\ref{Masseymac} where examples of
$\zp$ with nontrivial Massey products were constructed). On the
other hand, no other explicit topological types of moment-angle
manifolds $\zp$ are known. Furthermore, the following question
remains open:

\begin{problem}
Does there exist a moment-angle manifold $\zp$ decomposable into
nontrivial connected sum where one of the summands is
diffeomorphic to a product of more than two spheres.
\end{problem}

Here is an example illustrating how the structure of the
moment-angle manifold $\zp$ changes when one truncates~$P$ at a
vertex:

\begin{example}\label{ex3quadr}
Consider the following presentation of a polygon:
\begin{multline*}
  P=\bigl\{(x_1,x_2)\in\R^2\colon x_1\ge0,\;x_2\ge0,\\-x_1+1+\delta\ge0,\;-x_2+1+\varepsilon\ge0,\;
  -x_1-x_2+1\ge0\bigr\},
\end{multline*}
where $\delta,\varepsilon$ are parameters. First we fix some small
positive $\delta$ and vary $\varepsilon$. If $\varepsilon$ is
positive, then we get a presentation of a triangle with two
redundant inequalities. The corresponding moment-angle manifold is
diffeomorphic to $S^5\times S^1\times S^1$. If $\varepsilon$ is
negative, then we get a presentation of a quadrangle with one
redundant inequality. The corresponding manifold is diffeomorphic
to $S^3\times S^3\times S^1$. We see that when the hyperplane
$-x_2+1+\varepsilon=0$ crosses the vertex of the triangle, the
corresponding moment-angle manifold $\zp$ undergoes a surgery
turning $S^5\times S^1\times S^1$ to $S^3\times S^3\times S^1$.
Now if we fix some small negative $\varepsilon$ and decrease
$\delta$ so that the hyperplane $-x_1+1+\delta=0$ crosses the
vertex of the quadrangle, then the manifold $\zp$ undergoes a more
complicated surgery turning the product $S^3\times S^3\times S^1$
to the connected sum $(S^3\times S^4)^{\#5}$ (corresponding to the
case $n=2,m=5$ in Theorem~\ref{zpstacked}).

On the dual language of Gale diagrams, the surgeries described
above happen when the origin crosses the ``walls'' inside the
pentagonal Gale diagram corresponding to a presentation of a
2-polytope with 5 inequalities (see Figure~\ref{cgd5g}). For more
information about surgeries of moment-angle manifolds and their
relation to ``wall crossing'', see~\cite[\S6.4]{bu-pa02},
\cite{bo-me06} and~\cite{gi-lo13}.
\end{example}

\subsection*{Exercises}
\begin{exercise}\label{facesubmzp}
Let $G\subset P$ be a face of codimension~$k$ in a simple
$n$-polytope~$P$, let $\zp$ be the corresponding moment-angle
manifold with the quotient projection $p\colon\zp\to P$. Show that
$p^{-1}(G)$ is a smooth submanifold of $\zp$ of codimension $2k$.
Furthermore, $p^{-1}(G)$ is diffeomorphic to $\mathcal Z_G\times
T^\ell$, where $\mathcal Z_G$ is the moment-angle manifold
corresponding to~$G$ and $\ell$ is the number of facets of $P$ not
containing~$G$.
\end{exercise}

\begin{exercise}
Let $P=\vt(I^3)$ be the polytope obtained by truncating a 3-cube
at a vertex (see Construction~\ref{hypcut}). Write down
intersection of quadrics~\eqref{link} defining the corresponding
10-dimensional moment-angle manifold~$\zp$. Use
Theorem~\ref{zkcoh} or Theorem~\ref{zkhoch} to describe the
cohomology ring~$H^*(\zp)$. Deduce that $\zp$ cannot be
diffeomorphic to a connected sum of sphere products
(cf.~\cite[Example~11.5]{bo-me06}).
\end{exercise}

\begin{exercise}\label{pentsurg}
Write down the quadratic equations defining the moment-angle
manifolds $S^5\times S^1\times S^1$, $S^3\times S^3\times S^1$ and
$(S^3\times S^4)^{\#5}$ from Example~\ref{ex3quadr}, and describe
explicitly the surgeries between them.
\end{exercise}

\section{Symplectic reduction and moment maps revisited}\label{symred1}
As we have seen in Section~\ref{symred}, particular examples of
polytopal moment-angle manifolds $\zp$ appear as level sets for
the moment maps used in the construction of Hamiltonian toric
manifolds. In this case, the left hand sides of the equations
in~\eqref{zpqua} are quadratic Hamiltonians of a torus action
on~$\C^m$. Here we investigate the relationship between symplectic
quotients of $\C^m$ and intersections of quadrics more thoroughly.
As a corollary, we obtain that any symplectic quotient of $\C^m$
by a torus action is a Hamiltonian toric manifold.

We want to study symplectic quotients of $\C^m$ by torus subgroups
$T\subset\T^m$. Such a subgroup of dimension $m-n$ has the form
\begin{equation}\label{Tgamma}
  T_\varGamma=
  \bigl\{\bigr(e^{2\pi i\langle\gamma_1,\varphi\rangle},
  \ldots,e^{2\pi i\langle\gamma_m,\varphi\rangle}\bigl)
  \in\T^m\bigr\},
\end{equation}
where $\varphi\in\R^{m-n}$ is an $(m-n)$-dimensional parameter,
and $\varGamma=(\gamma_1,\ldots,\gamma_m)$ is a set of $m$ vectors
in~$\R^{m-n}$. In order for $T_\varGamma$ to be a torus, the
configuration of vectors $\gamma_1,\ldots,\gamma_m$ must be
\emph{rational}, i.e. their integer span
$L=\Z\langle\gamma_1,\ldots,\gamma_m\rangle$ must be a full rank
discrete subgroup (a \emph{lattice})\label{defnlattice}
in~$\R^{m-n}$. The lattice $L$ is identified canonically with
$\Hom(T_\varGamma,\mathbb S^1)$ and is called the \emph{weight
lattice} of torus~\eqref{Tgamma}. Let
\[
  L^*=\{\lambda^*\in{\mathbb R}^{m-n}\colon
  \langle\lambda^*,\lambda\rangle\in{\mathbb Z}
  \text{ for all }\lambda\in L\}
\]
be the dual lattice. We shall represent elements of $T_\varGamma$
by $\varphi\in\R^{m-n}$ occasionally, so that $T_\varGamma$ is
identified with the quotient $\R^{m-n}/\displaystyle L^*$.

The restricted action of $T_\varGamma\subset\T^m$ on $\C^m$ is
obviously Hamiltonian, and the corresponding moment map is the
composition
\begin{equation}\label{mugamma}
  \mu_\varGamma\colon\C^m\stackrel\mu\longrightarrow\R^m\longrightarrow
  \mathfrak t_\varGamma^*,
\end{equation}
where $\R^m\to\mathfrak t_\varGamma^*$ is the map of the dual Lie
algebras corresponding to $T_\varGamma\to\T^m$. The map
$\R^m\to\mathfrak t_\varGamma^*$ takes the $k$th basis vector $\mb
e_k\in\R^m$ to $\gamma_k\in\mathfrak t_\varGamma^*$. By
identifying $t_\varGamma^*$ with $\R^{m-n}$ we write the map
$\R^m\to\mathfrak t_\varGamma^*$ by the
matrix~$\varGamma=(\gamma_{jk})$, where $\gamma_{jk}$ is the $j$th
coordinate of~$\gamma_k$. The moment map~\eqref{mugamma} is then
given by
\[
  (z_1,\ldots,z_m)\longmapsto
  \Bigl(\sum_{k=1}^m\gamma_{1k}|z_k|^2,\ldots,\sum_{k=1}^m
  \gamma_{m-n,k}|z_k|^2\Bigr).
\]
Its level set $\mu_\varGamma^{-1}(\delta)$ corresponding to a
value $\delta=(\delta_1,\ldots,\delta_{m-n})^t\in\mathfrak
t_\varGamma^*$ is exactly the intersection of quadrics $\mathcal
Z_{\varGamma,\delta}$ given by~\eqref{zgamma}.

To apply symplectic reduction we need to identify when the moment
map~$\mu_\varGamma$ is proper, find its regular values~$\delta$,
and finally identify when the action of $T_\varGamma$ on
$\mu_\varGamma^{-1}(\delta)=\mathcal Z_{\varGamma,\delta}$ is
free. In Theorem~\ref{propmmap} below, all these conditions are
expressed in terms of the polyhedron $P$ associated with $\mathcal
Z_{\varGamma,\delta}$ as described in Section~\ref{intquad}.

It follows from Gale duality that $\gamma_1,\ldots,\gamma_m$ span
a lattice $L$ in $\R^{m-n}$ if and only if the dual configuration
$\mb a_1,\ldots,\mb a_m$ spans a lattice $N=\Z\langle\mb
a_1,\ldots,\mb a_m\rangle$ in~$\R^n$. We refer to a
presentation~\eqref{ptope1} as~\emph{rational}\label{rationprese}
if $\Z\langle\mb a_1,\ldots,\mb a_m\rangle$ is a lattice.

%\begin{remark}
%Any lattice polytope has a rational presentation, but the vertices
%of a polytope with a rational presentation do not necessarily belong
%to any lattice.
%\end{remakr}

%Recall that for each $\mb x\in P$ we defined
%\[
%  I_{\mb x}=\{i\in[m]\colon\langle\mb a_i,\mb x\rangle+b_i=0\}
%  =\{i\in[m]\colon{\mb x}\in F_i\}
%\]
%(the set of facets containing~$\mb x$).
A polyhedron $P$ is called \emph{Delzant}\label{delzpolyh} if it
has a vertex and there is a rational presentation~\eqref{ptope1}
such that for any vertex $\mb x\in P$ the vectors
%$\{\mb a_i\colon i\in I_{\mb x}\}$
%constitute a part of a basis of $N=\Z\langle\mb a_1,\ldots,\mb
%a_m\rangle$. It is enough to verify this condition at vertices
%$\mb x\in P$ only; in which case the vectors
$\mb a_i$ normal to the facets meeting at~$\mb x$ form a basis of
the lattice~$N=\Z\langle\mb a_1,\ldots,\mb a_m\rangle$. In the
case when $P$ is bounded and irredundant we obtain the definition
used before: a polytope $P$ is Delzant when its normal fan is
regular.

Now let $\delta\in\mathfrak t_\varGamma$ be a value of the moment
map $\mu_\varGamma\colon\C^m\to\mathfrak t_\varGamma^*$, and
$\mu_\varGamma^{-1}(\delta)=\mathcal Z_{\varGamma,\delta}$ the
corresponding level set, which is an intersection of
quadrics~\eqref{zgamma}. We associate with $\mathcal
Z_{\varGamma,\delta}$ a presentation~\eqref{ptope1} as described
in Section~\ref{intquad} (see Theorem~\ref{polquad}).

\begin{theorem}\label{propmmap}
Let $T_\varGamma\subset\T^m$ be a torus subgroup~\eqref{Tgamma},
determined by a rational configuration of vectors
$\gamma_1,\ldots,\gamma_m$.
\begin{itemize}
\item[(a)] The moment map
$\mu_\varGamma\colon\C^m\to\mathfrak t_\varGamma^*$ is proper if
and only if its level set $\mu_\varGamma^{-1}(\delta)$ is bounded
for some (and then for any) value $\delta\in\mathfrak
t_\varGamma^*$. Equivalently, the map $\mu_\varGamma$ is proper if
and only if the Gale dual configuration $\mb a_1,\ldots,\mb a_m$
satisfies $\alpha_1\mb a_1+\cdots+\alpha_m\mb a_m=\mathbf0$ for
some positive numbers~$\alpha_k$.

\item[(b)] $\delta\in\mathfrak t^*_\varGamma$ is a regular value of
$\mu_\varGamma$ if and only if the intersection of quadrics
$\mu_\varGamma^{-1}(\delta)=\mathcal Z_{\varGamma,\delta}$ is
nonempty and nonsingular. Equivalently, $\delta$ is a regular
value if and only if the associated presentation $P=P(A,\mb b)$ is
generic.

\item[(c)] The action of $T_\varGamma$ on
$\mu_\varGamma^{-1}(\delta)=\mathcal Z_{\varGamma,\delta}$ is free
if and only if the associated polyhedron $P$ is Delzant.
\end{itemize}
\end{theorem}
\begin{proof}
(a) If $\mu_\varGamma$ is proper then
$\mu_\varGamma^{-1}(\delta)\subset\mathfrak t^*_\varGamma$ is
compact, so it is bounded.

Now assume that $\mu_\varGamma^{-1}(\delta)=\mathcal
Z_{\varGamma,\delta}$ is bounded for some~$\delta$. Then the
associated polyhedron~$P$ is also bounded. By
Corollary~\ref{Pbound}, this is equivalent to vanishing of a
positive linear combination of $\mb a_1,\ldots,\mb a_m$. This
condition is independent of~$\delta$, and we conclude that
$\mu_\varGamma^{-1}(\delta)$ is bounded for any~$\delta$. Let
$X\subset\mathfrak t^*_\varGamma$ be a compact subset. Since
$\mu_\varGamma^{-1}(X)$ is closed, it is compact whenever it is
bounded. By Proposition~\ref{propcf} we may assume that, for any
$\delta\in X$, the first quadric defining
$\mu_\varGamma^{-1}(\delta)=\mathcal Z_{\varGamma,\delta}$ is
given by $\gamma_{11}|z_1|^2+\cdots+\gamma_{1m}|z_m|^2=\delta_1$
with $\gamma_{1k}>0$. Let $c=\max_{\delta\in X}\delta_1$. Then
$\mu_\varGamma^{-1}(X)$ is contained in the bounded set
\[
  \{\mb z\in\C^m\colon
  \gamma_{11}|z_1|^2+\cdots+\gamma_{1m}|z_m|^2\le c\}
\]
and is therefore bounded. Hence $\mu_\varGamma^{-1}(X)$ is
compact, and $\mu_\varGamma$ is proper.

(b) The first statement is the definition of a regular value. The
equivalent statement is already proved as Theorem~\ref{zpsmooth}.

(c) We first need to identify the stabilisers of the
$T_\varGamma$-action on $\mu_\varGamma^{-1}(\delta)$. Although the
fact that these stabilisers are finite for a regular value
$\delta$ follows from the general construction of symplectic
reduction, we can prove this directly.

%For any $I\subset[m]$, we define the sublattice generated by
%$\gamma_i$ which are not in~$I$:
%\[
%  L_{\widehat I}=\Z\langle\gamma_i\colon i\notin I\rangle\subset L=
%  \Z\langle\gamma_1,\ldots,\gamma_m\rangle.
%\]
%Recall that for any point $\mb z\in\mathcal Z_{\varGamma,\delta}$
%we have its zero set $\omega(\mb z)=\{i\colon z_i=0\}\subset[m]$.

Given a point $\mb z=(z_1,\ldots,z_m)\in\mathcal
Z_{\varGamma,\delta}$, we define the sublattice
\[
  L_{\mb z}=\Z\langle\gamma_i\colon z_i\ne0\rangle\subset L=
  \Z\langle\gamma_1,\ldots,\gamma_m\rangle.
\]

\begin{lemma}\label{afree}
The stabiliser subgroup of $\mb z\in\mathcal Z_{\varGamma,\delta}$
under the action of $T_\varGamma$ is given by $\displaystyle
L^*_{\mb z}/L^*$. Furthermore, if $\mathcal Z_{\varGamma,\delta}$
is nonsingular, then all these stabilisers are finite, i.e. the
action of $T_\varGamma$ on $\mathcal Z_{\varGamma,\delta}$ is
almost free.
\end{lemma}
\begin{proof}
An element $(e^{2\pi
i\langle\gamma_1,\varphi\rangle},\ldots,e^{2\pi
i\langle\gamma_m,\varphi\rangle})\in T_\varGamma$ fixes a point
$\mb z\in\mathcal Z_\varGamma$ if and only if $e^{2\pi
i\langle\gamma_k,\varphi\rangle}=1$ whenever $z_k\ne0$. In other
words, $\varphi\in T_\varGamma$ fixes $\mb z$ if and only if
$\langle\gamma_k,\varphi\rangle\in\Z$ whenever $z_k\ne0$. The
latter means that $\displaystyle\varphi\in\displaystyle L^*_{\mb
z}$. Since $\displaystyle\varphi\in L^*$ maps to $1\in
T_\varGamma$, the stabiliser of $\mb z$ is $\displaystyle L^*_{\mb
z}/L^*$.

Assume now that $\mathcal Z_{\varGamma,\delta}$ is nonsingular. In
order to see that $\displaystyle L^*_{\mb z}/L^*$ is finite we
need to check that the sublattice $L_{\mb
z}=\Z\langle\gamma_i\colon z_i\ne0\rangle\subset L$ has full
rank~$m-n$. Indeed, $\mathop{\mathrm{rk}}\{\gamma_i\colon
z_i\ne0\}$ is the rank of the matrix of gradients of quadrics
in~\eqref{zgamma} at~$\mb z$. Since $\mathcal
Z_{\varGamma,\delta}$ is nonsingular, this rank is $m-n$, as
needed.
\end{proof}

Now we can finish the proof of Theorem~\ref{propmmap}~(c). Assume
that $P$ is Delzant. By Lemma~\ref{afree}, the
$T_\varGamma$-action on $\mathcal Z_{\varGamma,\delta}$ is free if
and only if $L_{\mb z}=L$ for any $\mb z\in\mathcal
Z_{\varGamma,\delta}$. Let $i\colon\Z^k\rightarrow\Z^m$ be the
inclusion of the coordinate sublattice spanned by those $\mb e_i$
for which $z_i=0$, and let $p\colon\Z^m\rightarrow\Z^{m-k}$ be the
projection sending every such $\mb e_i$ to zero. We also have the
maps of lattices
\[
  \varGamma^*\colon L^*\to\Z^m,\;
  \mb l\mapsto\bigl(\langle\gamma_1,\mb l\rangle,\ldots,
  \langle\gamma_m,\mb l\rangle\bigr),\quad\text{and}\quad
  A\colon\Z^m\to N,\;\mb e_k\mapsto\mb a_k.
\]
%\[
%  A^*\colon N^*\to\Z^m,\;
%  \mb x\mapsto\bigl(\langle\mb a_1,\mb x\rangle,\ldots,
%  \langle\mb a_m,\mb x\rangle\bigr)\quad\text{and}\quad
%  \varGamma\colon\Z^m\to L,\;\mb e_k\mapsto\gamma_k.
%\]
Consider the diagram
\begin{equation}\label{2Zses}
\begin{CD}
  @.@.\begin{array}{c}0\\ \raisebox{2pt}{$\downarrow$}\\
  L^*\end{array}@.@.\\
  @.@.@VV{\varGamma^*}V@.@.\\
  0@>>>\Z^{k}@>i>>\Z^m@>p>>\Z^{m-k} @>>>0\\
  @.@.@VV A V@.@.\\
  @.@.\begin{array}{c}N\\ \downarrow\\ 0\end{array}@.@.\\
\end{CD}
\end{equation}
in which the vertical and horizontal sequences are exact. Then the
Delzant condition is equivalent to that the composition $A\cdot i$
is split injective. The condition $L_{\mb z}=L$ is equivalent to
that $\varGamma\cdot p^*$ is surjective, or $p\cdot\varGamma^*$ is
split injective. These two conditions are equivalent by
Lemma~\ref{2ses}.
\end{proof}

In the case of polytopes we obtain the following version of
Proposition~\ref{regvalue}:

\begin{corollary}
Let $P=P(A,\mb b)$ be a Delzant polytope,
$\varGamma=(\gamma_1,\ldots,\gamma_m)$ the Gale
dual\label{galeduadef} configuration, and $\zp$ the corresponding
moment-angle manifold. Then
\begin{itemize}
\item[(a)]
$\delta=\varGamma\mb b$ is a regular value of the moment map
$\mu_\varGamma\colon\C^m\to\mathfrak t^*_\varGamma$ for the
Hamiltonian action of $T_\varGamma\subset\T^m$ on~$\C^m$;
\item[(b)] $\zp$ is the regular level set $\mu_\varGamma^{-1}(\varGamma\mb
b)$;
\item[(c)] the action of $T_\varGamma$ on $\zp$ is free.
\end{itemize}
\end{corollary}

In Section~\ref{symred} we defined the Hamiltonian toric manifold
$V_P$\label{hamitomamdef} corresponding to a Delzant polytope~$P$
as the symplectic quotient of $\C^m$ by the torus subgroup
$K\subset\T^m$ determined by the normal fan of~$P$. By comparing
the vertical exact sequence in~\eqref{2Zses} with~\eqref{kgrou} we
obtain that $K=T_\varGamma$, and the quotient $n$-torus
$\T^m/T_\varGamma$ acting on $V_P=\zp/T_\varGamma$ is
$T_N=N\otimes_\Z\mathbb S=\R^n/N$.

\begin{corollary}
Any symplectic quotient of $\C^m$ by a torus subgroup
$T\subset\T^m$ is a Hamiltonian toric manifold.
\end{corollary}

\begin{example}\label{cp2sq}
Consider the case $m-n=1$, i.e. $T_\varGamma$ is 1-dimensional,
and $\gamma_k\in\R$. By Theorem~\ref{propmmap}~(a), the moment map
$\mu_\varGamma$ is proper whenever
\[
  \mu^{-1}_\varGamma(\delta)=\{\mb z\in\C^m\colon
  \gamma_1|z_1|^2+\cdots+\gamma_m|z_m|^2=\delta\}
\]
is bounded for any $\delta\in\R$. By Theorem~\ref{propmmap}~(b),
$\delta$ is a regular value whenever the quadratic hypersurface
$\gamma_1|z_1|^2+\cdots+\gamma_m|z_m|^2=\delta$ is nonempty and
nonsingular. These two conditions together imply that the
hypersurface is an ellipsoid, and the associated polyhedron is an
$n$-simplex (see Example~\ref{mamsimplex}). By Lemma~\ref{afree},
the $T_\varGamma$-action on $\mu^{-1}_\varGamma(\delta)$ is free
if and only if $L_{\mb z}=L$ for any $\mb
z\in\mu^{-1}_\varGamma(\delta)$. This means that each $\gamma_k$
generates the same lattice as the whole set
$\gamma_1,\ldots,\gamma_m$, which implies that
$\gamma_1=\cdots=\gamma_m$. The Gale dual configuration satisfies
$\mb a_1+\cdots+\mb a_m=\mathbf0$. Then $T_\varGamma$ is the
diagonal circle in~$\T^m$, the hypersurface
$\mu^{-1}_\varGamma(\delta)=\zp$ is a sphere, and the associated
polytope~$P$ is the standard simplex up to shift and magnification
by a positive factor~$\delta$. The Hamiltonian toric manifold
$V_P=\zp/T_\varGamma$ is the complex projective space~$\C P^n$.
\end{example}

\section{Complex structures on intersections of quadrics}\label{lvmma}
Bosio and Meersseman~\cite{bo-me06} identified polytopal
moment-angle manifolds $\zp$ with a class of non-K\"ahler
complex-analytic manifolds introduced in the works of Lopez de
Medrano, Verjovsky and Meersseman (\emph{LVM-manifolds}). This was
the starting point in the subsequent study of the complex geometry
of moment-angle manifolds. We review the construction of
LVM-manifolds and its connection to polytopal moment-angle
manifolds here.

The initial data of the construction of an LVM-manifold is a link
of a homogeneous system of quadrics similar to~\eqref{link}, but
with \emph{complex} coefficients:
\begin{equation}\label{clink}
  \mathcal L=\left\{\begin{array}{lrcl}
  \mb z\in\C^m\colon&\sum_{k=1}^m|z_k|^2&=&1,\\[1mm]
  &\sum_{k=1}^m \zeta_k|z_k|^2&=&\mathbf0
  \end{array}\right\},
\end{equation}
where $\zeta_k\in\C^s$. We can obviously turn this link into the
form~\eqref{link} by identifying $\C^s$ with $\R^{2s}$ in the
standard way, so that each $\zeta_k$ becomes $\mb
g_k\in\R^{m-n-1}$ with $n=m-2s-1$. We assume that the link is
nonsingular, i.e. the system of complex vectors
$(\zeta_1,\ldots,\zeta_m)$ (or the corresponding system of real
vectors $(\mb g_1,\ldots,\mb g_m)$) satisfies the conditions (a)
and~(b) of Proposition~\ref{nondeglink}.

Now define the manifold $\mathcal N$ as the projectivisation of
the intersection of homogeneous quadrics in~\eqref{clink}:
\begin{equation}\label{ndef}
\mathcal N=\bigl\{\mb z\in\C
P^{m-1}\colon\zeta_1|z_1|^2+\cdots+\zeta_m|z_m|^2=\textbf
0\bigr\},\quad\zeta_k\in\C^s.
\end{equation}

We therefore have a principal $S^1$-bundle $\mathcal L\to\mathcal
N$.

\begin{theorem}[Meersseman~\cite{meer00}]\label{thlvm}
The manifold $\mathcal N$ has a holomorphic atlas describing it as
a compact complex manifold of complex dimension $m-1-s$.
\end{theorem}
\begin{proof}[Sketch of proof]
Consider a holomorphic action of $\C^s$ on $\C^m$ given by
\begin{equation}\label{csact}
\begin{aligned}
  \C^s\times\C^m&\longrightarrow\C^m\\
  (\mb w,\mb z)&\mapsto \bigl(z_1e^{\langle\zeta_1,\mb
  w\rangle},\ldots,z_me^{\langle\zeta_m,\mb w\rangle}\bigr),
\end{aligned}
\end{equation}
where $\mb w=(w_1,\ldots,w_s)\in\C^s$, and $\langle\zeta_k,\mb
w\rangle=\zeta_{1k}w_1+\cdots+\zeta_{sk}w_s$.

Let $\sK$ be the simplicial complex consisting of zero-sets of
points of the link~$\mathcal L$:
\[
  \sK=\{\omega(\mb z)\colon\mb z\in \mathcal L\}.
\]
Observe that $\sK=\sK_P$, where $P$ is the simple polytope
associated with the link~$\mathcal L$. Let $U=U(\sK)$ be the
corresponding subspace arrangement complement given by~\eqref{uk}.
Note that Proposition~\ref{galecomb} implies that $U$ can be also
defined as
\[
  U=\bigl\{(z_1,\ldots,z_m)\in\C^m\colon\mathbf 0\in\conv(\zeta_j\colon z_j\ne0)\bigr\}.
\]

An argument similar to the proof of Lemma~\ref{afree} shows that
the restriction of the action~\eqref{csact} to $U\subset\C^m$ is
free. Also, this restricted action of $\C^s$ on $U$ is proper (we
shall prove this in more general context in
Theorem~\ref{zkcomplex} below), so the quotient $U/\C^s$ is
Hausdorff. Using a holomorphic atlas transverse to the orbits of
the free action of $\C^s$ on the complex manifold $U$ we obtain
that the quotient $U/\C^s$ has a structure of a complex manifold.

On the other hand, it can be shown that the function
$|z_1|^2+\cdots+|z_m|^2$ (the square of the distance to the origin
in~$\C^m$) has a unique minimum when restricted to an orbit of the
free action of $\C^s$ on $U$. The set of these minima (i.e. the
set of points closest to the origin in each orbit) can be
described as
\[
  \mathcal T=\bigl\{\mb z\in\C^m\setminus\!\{\mathbf 0\}
  \colon\;\zeta_1|z_1|^2+\cdots+\zeta_m|z_m|^2=\textbf 0\bigr\}.
\]
It follows that the quotient $U/\C^s$ can be identified
with~$\mathcal T$, and therefore $\mathcal T$ acquires a structure
of a complex manifold of dimension $m-s$.

By projectivising the construction we identify $\mathcal N$ with
the quotient of a complement of coordinate subspace arrangement in
$\C P^{m-1}$ (the projectivisation of~$U$) by a holomorphic action
of~$\C^s$. In this way $\mathcal N$ becomes a compact complex
manifold.
\end{proof}

The manifold $\mathcal N$ with the complex structure of
Theorem~\ref{thlvm} is referred to as an
\emph{LVM-manifold}\label{lvmmanid}. These manifolds were
described by Meersseman~\cite{meer00} as a generalisation of the
construction of Lopez de Medrano and Verjovsky~\cite{lo-ve97}.

\begin{remark}
The embedding of $\mathcal T$ in $\C^m$ and of $\mathcal N$ in $\C
P^{m-1}$ given by~\eqref{ndef} is not holomorphic.
\end{remark}

A polytopal moment-angle manifold $\zp$ is diffeomorphic to a
link~\eqref{link}, which can be turned into a complex
link~\eqref{clink} whenever $m+n$ is odd. It follows that the
quotient $\zp/S^1$ of an odd-dimensional moment-angle manifold has
a complex-analytic structure as an LVM-manifold. By adding
redundant inequalities and using the $S^1$-bundle $\mathcal
L\to\mathcal N$, Bosio--Meersseman observed that $\zp$ or
$\zp\times S^1$ has a structure of an LVM-manifold, depending on
whether $m+n$ is even or odd.

We first summarise the effects that a redundant inequality
in~\eqref{ptope1} has on different spaces appeared above:

\begin{proposition}
Assume that~\eqref{ptope1} is a generic presentation. The
following conditions are equivalent:
\begin{itemize}
\item[(a)] $\langle\mb a_i,\mb x\:\rangle+b_i\ge0$ is a redundant
inequality in~\eqref{ptope1} (i.e. $F_i=\varnothing$);
\item[(b)] $\zp\subset\{\mb z\in\C^m\colon z_i\ne0\}$;
\item[(c)] $\{i\}$ is a ghost vertex of~$\sK_P$;
\item[(d)] $U(\sK_P)$ has a factor $\C^\times$ on the $i$th
coordinate;
\item[(e)] $\mathbf 0\notin\conv(\mb g_k\colon k\ne i)$.
\end{itemize}
\end{proposition}
\begin{proof}
The equivalence of the first four conditions follows directly from
the definitions. The equivalence (a)$\Leftrightarrow$(e) follows
from Proposition~\ref{galecomb}.
\end{proof}

\begin{theorem}[\cite{bo-me06}]\label{bometh}
Let $\zp$ be the moment angle manifold corresponding to an
$n$-dimensional simple polytope~\eqref{ptope1} defined by $m$
inequalities.
\begin{itemize}
\item[(a)] If $m+n$ is even then $\zp$ has a complex structure as
an LVM-manifold.
\item[(b)] If $m+n$ is odd then $\zp\times S^1$ has a complex structure as
an LVM-manifold.
\end{itemize}
\end{theorem}
\begin{proof}
(a) We add one redundant inequality of the form $1\ge0$
to~\eqref{ptope1}, and denote the resulting manifold
of~\eqref{cdiz} by $\zp'$. We have $\zp'\cong\zp\times S^1$. By
Proposition~\ref{malink}, $\zp$ is diffeomorphic to a link given
by~\eqref{link}. Then $\zp'$ is given by the intersection of
quadrics
\[
  \left\{\begin{array}{lrcccrcr}
  \mb z\in\C^{m+1}\colon&|z_1|^2&+&\cdots&+&|z_m|^2&&=1,\\[1mm]
                 &\mb g_1|z_1|^2&+&\cdots&+&\mb g_m|z_m|^2&&=\mathbf 0,\\[1mm]
  &&&&&&&|z_{m+1}|^2=1,
  \end{array}\right\}
\]
which is diffeomorphic to the link given by
\[
  \left\{\begin{array}{lrcccrcr}
  \mb z\in\C^{m+1}\colon  &|z_1|^2&+&\cdots&+&|z_m|^2&+&|z_{m+1}|^2=1,\\[1mm]
  &\mb g_1|z_1|^2&+&\cdots&+&\mb g_m|z_m|^2&&=\mathbf 0,\\[1mm]
  &|z_1|^2&+&\cdots&+&|z_m|^2&-&|z_{m+1}|^2=0.
  \end{array}\right\}
\]
If we denote by $\varGamma^\star=(\mb g_1\ldots\:\mb g_m)$ the
$(m-n-1)\times m$-matrix of coefficients of the homogeneous
quadrics for~$\zp$, then the corresponding matrix for $\zp'$ is
\[
  {\varGamma^\star}^\prime=\begin{pmatrix}\mb g_1&\cdots&\mb
  g_m&0\\1&\cdots&1&-1
  \end{pmatrix}.
\]
Its height $m-n$ is even, so that we may think of its $k$th column
as a complex vector $\zeta_k$ (by identifying $\R^{m-n}$ with
$\C^{\frac{m-n}2}$), for $k=1,\ldots,m+1$. Now define
\begin{equation}\label{nprimeeven}
  \mathcal N'=\bigl\{\mb z\in\C P^m\colon
  \zeta_1|z_1|^2+\cdots+\zeta_{m+1}|z_{m+1}|^2=\mathbf
  0\bigr\}.
\end{equation}
Then $\mathcal N'$ has a complex structure as an LVM-manifold by
Theorem~\ref{thlvm}. On the other hand,
\[
  \mathcal N'\cong\zp'/S^1=(\zp\times S^1)/S^1\cong\zp,
\]
so that $\zp$ also acquires a complex structure.

\smallskip

(b) The proof here is similar, but we have to add two redundant
inequalities $1\ge0$ to~\eqref{ptope1}. Then $\zp'\cong\zp\times
S^1\times S^1$ is given by
\[
  \left\{\begin{array}{lrcrclcr}
  \mb
  z\in\C^{m+2}\colon&|z_1|^2&+\;\cdots\;+&|z_m|^2&+&|z_{m+1}|^2&+&|z_{m+2}|^2=1,\\[1mm]
  &\mb g_1|z_1|^2&+\;\cdots\;+&\mb g_m|z_m|^2&&&&=\mathbf 0,\\[1mm]
  &|z_1|^2&+\;\cdots\;+&|z_m|^2&-&|z_{m+1}|^2&&=0,\\[1mm]
  &|z_1|^2&+\;\cdots\;+&|z_m|^2&&&-&|z_{m+2}|^2=0.
  \end{array}\right\}
\]
The matrix of coefficients of the homogeneous quadrics is
therefore
\[
  {\varGamma^\star}^\prime=\begin{pmatrix}\mb g_1&\cdots&\mb
  g_m&0&0\\1&\cdots&1&-1&0\\1&\cdots&1&0&-1
  \end{pmatrix}.
\]
We think of its columns as a set of $m+2$ complex vectors
$\zeta_1,\ldots,\zeta_{m+2}$, and define
\begin{equation}\label{nprimeodd}
  \mathcal N'=\bigl\{\mb z\in\C P^{m+1}\colon
  \zeta_1|z_1|^2+\cdots+\zeta_{m+2}|z_{m+2}|^2=\mathbf
  0\bigr\}.
\end{equation}
Then $\mathcal N'$ has a complex structure as an LVM-manifold. On
the other hand,
\[
  \mathcal N'\cong\zp'/S^1=(\zp\times S^1\times S^1)/S^1\cong\zp\times S^1,
\]
and therefore $\zp\times S^1$ has a complex structure.
\end{proof}

In the next two sections we describe a more direct method of
endowing $\zp$ with a complex structure, without referring to
projectivised quadrics and LVM-manifolds. This approach, developed
in~\cite{pa-us12}, works not only in the polytopal case, but also
for the moment-angle manifolds $\zk$ corresponding to underlying
complexes $\sK$ of complete simplicial fans.

\section{Moment-angle manifolds from simplicial
fans}\label{mamsf}

Let $\sK=\sK_\Sigma$ be the underlying complex of a complete
simplicial fan~$\Sigma$, and $U(\sK)$ the complement of the
coordinate subspace arrangement~\eqref{uk} defined by~$\sK$. Here
we shall identify the moment-angle manifold $\zk$ with the
quotient of $U(\sK)$ by a smooth action of a non-compact group
isomorphic to~$\R^{m-n}$, thereby defining a smooth structure
on~$\zk$. A modification of this construction will be used in the
next section to endow $\zk$ with a complex structure.
%These results were obtained in the work~\cite{pa-us12}.

Let $\Sigma$ be a simplicial fan in $N_\R\cong\R^n$ with
generators $\mb a_1,\ldots,\mb a_m$. Recall that the underlying
simplicial complex $\sK=\sK_\Sigma$ is the collection of subsets
$I\subset[m]$ such that $\{\mb a_i\colon i\in I\}$ spans a cone
of~$\Sigma$.

A simplicial fan $\Sigma$ is therefore determined by two pieces of
data:
\begin{itemize}
\item[--] a simplicial complex $\sK$ on $[m]$;
\item[--] a configuration of vectors $\mb a_1,\ldots,\mb a_m$ in
$\R^n$ such that  for any simplex $I\in\sK$ the subset $\{\mb
a_i\colon i\in I\}$ is linearly independent.
\end{itemize}
Then for each $I\in\sK$ we can define the simplicial cone
$\sigma_I$ spanned by $\mb a_i$ with $i\in I$. The `bunch of
cones' $\{\sigma_I\colon I\in\sK\}$ patches into a fan $\Sigma$
whenever any two cones $\sigma_I$ and $\sigma_J$ intersect in a
common face (which has to be $\sigma_{I\cap J}$). Equivalently,
the relative interiors of cones $\sigma_I$ are pairwise
non-intersecting. Under this condition, we say that the data
$\{\sK;\mb a_1,\ldots,\mb a_m\}$ \emph{define a fan}~$\Sigma$.

We do allow ghost vertices in $\sK$; they do not affect the
fan~$\Sigma$. The vector $\mb a_i$ corresponding to a ghost vertex
$\{i\}$ can be zero as it does not correspond to a one-dimensional
cone of~$\Sigma$. This formalism will be important for the
construction of a complex structure on~$\zk$; it was also used
in~\cite{ba-za} under the name \emph{triangulated vector
configurations}.

\begin{construction}
For a set of vectors $\mb a_1,\ldots,\mb a_m$, consider the linear
map
\begin{equation}\label{lambdar}
  A\colon\R^m\to N_\R,\quad\mb e_i\mapsto\mb a_i,
\end{equation}
where $\mb e_1,\ldots,\mb e_m$ is the standard basis of~$\R^m$.
Let
\[
  \R^m_>=\{(y_1,\ldots,y_m)\in\R^m\colon y_i>0\}
\]
be the multiplicative group of $m$-tuples of positive real
numbers, and define
\begin{equation}\label{rsigma}
\begin{aligned}
  R&=\exp(\Ker A)=\bigl\{\bigl(e^{y_1},\ldots,e^{y_m}\bigr)
  \colon(y_1,\ldots,y_m)\in\Ker A\bigr\}\\
  &=\bigl\{(t_1,\ldots,t_m)\in\R^m_>\colon
  \prod_{i=1}^mt_i^{\langle\mb a_i,\mb u\rangle}=1
  \text{ for all }\mb u\in N^*_\R\bigr\}.
\end{aligned}
\end{equation}
\end{construction}

%If $\mb a_1,\ldots,\mb a_m$ span $N_\R\cong\R^n$ (e.g. if $\Sigma$
%is a complete fan), then $R\cong\R^{m-n}_>$.

We let $\R^m_>$ act on the complement $U(\sK)\subset\C^m$ by
coordinatewise multiplications and consider the restricted action
of the subgroup $R\subset\R^m_>$. Recall that a $G$-action on a
space $X$ is \emph{proper}\label{propergad} if the \emph{group
action map} $h\colon G\times X\to X\times X$, \ $(g,x)\mapsto
(gx,x)$ is proper (the preimage of a compact subset is compact).

\begin{theorem}[\cite{pa-us12}]\label{zksmooth}
Assume given data $\{\sK;\mb a_1,\ldots,\mb a_m\}$ satisfying the
conditions above. Then
\begin{itemize}
\item[(a)]
the group $R\cong\R^{m-n}$ given by \eqref{rsigma} acts on
$U(\sK)$ freely;

\item[(b)] if the data $\{\sK;\mb a_1,\ldots,\mb a_m\}$ define a simplicial fan~$\Sigma$,
then $R$ acts on $U(\sK)$ properly, so the quotient $U(\sK)/R$ is
a smooth Hausdorff $(m+n)$-dimensional manifold;

\item[(c)] if the fan $\Sigma$ is complete, then $U(\sK)/R$ is
homeomorphic to the moment-angle manifold~$\zk$.
\end{itemize}
Therefore, $\zk$ can be smoothed whenever $\sK=\sK_\Sigma$ for a
complete simplicial fan~$\Sigma$.
\end{theorem}

\begin{proof}
Statement~(a) is proved in the same way as
Proposition~\ref{freeaction}. Indeed, a point $\mb z\in U(\sK)$
has a nontrivial stabiliser with respect to the action of $\R^m_>$
only if some of its coordinates vanish. These $\R^m_>$-stabilisers
are of the form $(\R_>,1)^I$, see~\eqref{XAI}, for some $I\in\sK$.
The restriction of $\exp A$ to any such $(\R_>,1)^I$ is an
injection. Therefore, $R=\exp(\Ker A)$ intersects any
$\R^m_>$-stabilisers only at the unit, which implies that the
$R$-action on $U(\sK)$ is free.

Let us prove~(b) (compare the proof of Theorem~\ref{coxth}~(b) and
Exercise~\ref{Gprop}). Consider the map
\[
  h\colon R\times U(\sK)\to U(\sK)\times U(\sK), \quad (\mb
  g,\mb z)\mapsto (\mb r\mb z,\mb z),
\]
for $\mb r\in R$, $\mb z\in U(\sK)$. Let $V\subset U(\sK)\times
U(\sK)$ be a compact subset; we need to show that $h^{-1}(V)$ is
compact. Since $R\times U(\sK)$ is metrisable, it suffices to
check that any infinite sequence $\{(\mb r^{(k)},\mb
z^{(k)})\colon k=1,2,\ldots\}$ of points in $h^{-1}(V)$ contains a
converging subsequence. Since $V\subset U(\sK)\times U(\sK)$ is
compact, by passing to a subsequence we may assume that the
sequence
\[
  \{h(\mb r^{(k)},\mb z^{(k)})\}=\{(\mb r^{(k)}\mb z^{(k)},\mb z^{(k)})\}
\]
has a limit in $U(\sK)\times U(\sK)$. We set $\mb w^{(k)}=\mb
r^{(k)}\mb z^{(k)}$, and assume that
\[
  \{\mb w^{(k)}\}\to {\mb w}=(w_1,\dots,w_m),\quad
  \{\mb z^{(k)}\}\to {\mb z}=(z_1,\dots,z_m)
\]
for some $\mb w,\mb z\in U(\sK)$. We need to show that a
subsequence of $\{\mb r^{(k)}\}$ has limit in~$R$. We write
\[
  \mb r^{(k)}=\bigl(g_1^{(k)},\ldots,g_m^{(k)}\bigr)=
  \bigl(e^{\alpha^{(k)}_1},\ldots,
  e^{\alpha^{(k)}_m}\bigr)\in R\subset \R^m_>,
\]
$\alpha_j^{(k)}\in\R$. By passing to a subsequence we may assume
that each sequence $\{\alpha^{(k)}_j\}$, \ $j=1,\ldots,m$, has a
finite or infinite limit (including $\pm\infty$). Let
\[
  I_+=\{j\colon\alpha^{(k)}_j\to +\infty\}\subset[m],\quad
  I_-=\{j\colon\alpha^{(k)}_j\to -\infty\}\subset[m].
\]
Since the sequences $\{\mb z^{(k)}\}$, $\{\mb w^{(k)}=\mb
r^{(k)}\mb z^{(k)}\}$ converge to $\mb z,\mb w\in U(\sK)$
respectively, we have $z_j=0$ for $j\in I_+$ and $w_j=0$ for $j\in
I_-$. Then it follows from the decomposition
$U(\sK)=\bigcup_{I\in\sK}(\C,\C^\times)^I$ that $I_+$ and $I_-$
are simplices of~$\sK$. Let $\sigma_+,\sigma_-$ be the
corresponding cones of the simplicial fan~$\Sigma$. Then
$\sigma_+\cap\sigma_-=\{\mathbf 0\}$ by definition of a fan. By
Lemma~\ref{seplemma}, there exists a linear function $\mb u\in
N^*_\R$ such that $\langle\mb u,\mb a\rangle>0$ for any nonzero
$\mb a\in \sigma_+$, and $\langle\mb u,\mb a\rangle<0$ for any
nonzero $\mb a\in \sigma_-$. Since $\mb r^{(k)}\in R$, it follows
from~\eqref{rsigma} that
\begin{equation}\label{alphak1}
%  0=\Bigl\langle\mb u,\sum_{j=1}^m \alpha^{(k)}_j\mb a_j\Bigr\rangle=
  \sum_{j=1}^m \alpha^{(k)}_j\langle\mb u,\mb a_j\rangle=0.
\end{equation}
This implies that both $I_+$ and $I_-$ are empty, as otherwise the
latter sum tends to infinity. Thus, each sequence
$\{\alpha^{(k)}_j\}$ has a finite limit $\alpha_j$, and a
subsequence of $\{\mb r^{(k)}\}$ converges to
$(e^{\alpha_1},\ldots,e^{\alpha_m})$. Passing to the limit
in~\eqref{alphak1} we obtain that
$(e^{\alpha_1},\ldots,e^{\alpha_m})\in R$. This proves the
properness of the action. Since the Lie group $R$ acts smoothly,
freely and properly on the smooth manifold $U(\sK)$, the orbit
space $U(\sK)/R$ is Hausdorff and smooth by the standard
result~\cite[Theorem~9.16]{lee00}.

In the case of complete fan it is possible to construct a smooth
atlas on $U(\sK)/R$ explicitly. To do this, it is convenient to
pre-factorise everything by the action of $\T^m$, as in the proof
of Theorem~\ref{deret}. We  have
\[
  U(\sK)/\T^m=(\R_\ge,\R_>)^\sK=\bigcup_{I\in\sK}(\R_\ge,\R_>)^I.
\]
Since the fan $\Sigma$ is complete, we may take the union above
only over $n$-element simplices $I=\{i_1,\ldots,i_n\}\in\sK$.
Consider one such simplex~$I$; the generators of the corresponding
$n$-dimensional cone $\sigma\in\Sigma$ are $\mb a_{i_1},\ldots,\mb
a_{i_n}$. Let $\mb u_1,\ldots,\mb u_n$ be the dual basis of
$N_\R^*$ (which is a generator set of the dual cone
$\sigma^\mathsf{v}$). Then we have $\langle\mb a_{i_k},\mb
u_j\rangle=\delta_{kj}$. Now consider the map
\begin{align*}
  p_I\colon(\R_\ge,\R_>)^I&\to\R^n_\ge\\
  (y_1,\ldots,y_m)&\mapsto
  \Bigr(\prod_{i=1}^my_i^{\langle\mb a_i,\mb u_1\rangle},\ldots,\:
  \prod_{i=1}^my_i^{\langle\mb a_i,\mb u_n\rangle}\Bigl),
\end{align*}
where we set $0^0=1$. Note that zero does not occur with a
negative exponent in the right hand side, hence $p_I$ is well
defined as a continuous map. Each $(\R_\ge,\R_>)^I$ is
$R$-invariant, and it follows from~\eqref{rsigma} that $p_I$
induces an injective map
\[
  q_I\colon(\R_\ge,\R_>)^I/R\to\R^n_\ge.
\]
This map is also surjective since every
$(x_1,\ldots,x_n)\in\R^n_\ge$ is covered by $(y_1,\ldots,y_m)$
where $y_{i_j}=x_j$ for $1\le j\le n$ and $y_k=1$ for
$k\notin\{i_1,\ldots,i_n\}$. Hence, $q_I$ is a homeomorphism. It
is covered by a $\T^m$-equivariant homeomorphism
\[
  \overline q_I\colon
  (\C,\C^\times)^I/R\to\C^n\times\T^{m-n},
\]
where $\C^n$ is identified with the quotient
$\R_\ge^n\times\T^n/\!\sim\,$, see~\eqref{Cmids}. Since $U(\sK)/R$
is covered by open subsets $(\C,\C^\times)^I/R$, and
$\C^n\times\T^{m-n}$ embeds as an open subset in $\R^{m+n}$, the
set of homeomorphisms $\{\overline q_I\colon I\in\sK\}$ provides
an atlas for $U(\sK)/R$.
%(compare~\cite[proof of Lemma~6.2]{bu-pa02}).
The change of coordinates transformations $\overline q_J\overline
q_I^{-1}\colon\C^n\times\T^{m-n}\to\C^n\times\T^{m-n}$ are smooth
by inspection; thus $U(\sK)/R$ is a smooth manifold.

\begin{remark}
The set of homeomorphisms
$\{q_I\colon(\R_\ge,\R_>)^I/R\to\R^n_\ge\}$ defines an atlas for
the smooth manifold with corners $\zk/\T^m$. If $\sK=\sK_P$ for a
simple polytope $P$, then this smooth structure with corners
coincides with that of~$P$.
\end{remark}

It remains to prove statement~(c), that is, identify $U(\sK)/R$
with $\zk$. If $X$ is a Hausdorff locally compact space with a
proper $G$-action, and $Y\subset X$ a compact subspace which
intersects every $G$-orbit at a single point, then $Y$ is
homeomorphic to the orbit space $X/G$. Therefore, we need to
verify that each $R$-orbit intersects $\zk\subset U(\sK)$ at a
single point. We first prove that the $R$-orbit of any $\mb y\in
U(\sK)/\T^m=(\R_\ge,\R_>)^\sK$ intersects $\zk/\T^m$ at a single
point. For this we use the cubical decomposition
$\cc(\sK)=(\I,1)^\sK$ of $\zk/\T^m$, see Example~\ref{exkpo}.2.

Assume first that $\mb y\in\R^m_>$. The $R$-action on $\R^m_>$ is
obtained by exponentiating the linear action of $\Ker A$ on
$\R^m$. Consider the subset $(\R_\le,0)^\sK\subset\R^m$, where
$\R_\le$ denotes the set of nonpositive reals. It is taken by the
exponential map $\exp\colon\R^m\to\R^m_>$ homeomorphically onto
$\cc^\circ(\sK)=((0,1],1)^\sK\subset\R^m_>$, where $(0,1]$ is
denotes the semi-interval $\{y\in\R\colon 0<y\le1\}$.
%, where the latter is the relative interior of $\cc(\sK)$.
The map
\begin{equation}\label{oto}
  A\colon(\R_\le,0)^\sK\to N_\R
\end{equation}
takes every $(\R_\le,0)^I$ to $-\sigma$, where $\sigma\in\Sigma$
is the cone corresponding to $I\in\sK$. Since $\Sigma$ is
complete, map~\eqref{oto} is one-to-one.

The orbit of $\mb y$ under the action of $R$ consists of points
$\mb w\in\R^m_>$ such that $\exp A\mb w=\exp A\mb y$. Since $A\mb
y\in N_\R$ and map~\eqref{oto} is one-to-one, there is a unique
point $\mb y'\in(\R_\le,0)^\sK$ such that $A\mb y'=A\mb y$. Since
$\exp A\mb y'\subset\cc^\circ(\sK)$, the $R$-orbit of $\mb y$
intersects $\cc^\circ(\sK)$ and therefore $\cc(\sK)$ at a unique
point.

Now let $\mb y\in(\R_\ge,\R_>)^\sK$ be an arbitrary point. Let
$\omega(\mb y)\in\sK$ be the set of zero coordinates of~$\mb y$,
and let $\sigma\in\Sigma$ be the cone corresponding to $\omega(\mb
y)$. The cones containing $\sigma$ constitute a fan $\st\sigma$
(called the \emph{star}\label{starinfan} of~$\sigma$) in the
quotient space $N_\R/\R\langle\mb a_i\colon i\in\omega(\mb
y)\rangle$. The underlying simplicial complex of $\st\sigma$ is
the link $\lk_{\sK}\omega(\mb y)$ of $\omega(\mb y)$ in~$\sK$. Now
observe that the action of $R$ on the set
\[
  \{(y_1,\ldots,y_m)\in(\R_\ge,\R_>)^\sK\colon y_i=0\text{ for }i\in \omega(\mb
  y)\}\cong(\R_\ge,\R_>)^{\lk\omega(\mb y)}
\]
coincides with the action of the group $R_{\st\sigma}$ (defined by
the fan~$\st\sigma$). Now we can repeat the above arguments for
the complete fan $\st\sigma$ and the action of $R_{\st\sigma}$ on
$(\R_\ge,\R_>)^{\lk\omega(\mb y)}$. As a result, we obtain that
every $R$-orbit intersects $\cc(\sK)$ at a unique point.

To finish the proof of~(c) we consider the commutative diagram
\[
\begin{CD}
\zk @>>> U(\sK)\\
@VVV @VV\pi V\\
\cc(\sK) @>>> (\R_\ge,\R_>)^\sK
\end{CD}
\]
where the horizontal arrows are embeddings and the vertical ones
are projections onto the quotients of $\T^m$-actions. Note that
the projection $\pi$ commutes with the $R$-actions on $U(\sK)$ and
$(\R_\ge,\R_>)^\sK$, and the subgroups $R$ and $\T^m$ of
$(\C^\times)^m$ intersect trivially. It follows that every
$R$-orbit intersects the full preimage $\pi^{-1}(\cc(\sK))=\zk$ at
a unique point. Indeed, assume that $\mb z$ and $\mb r\mb z$ are
in $\zk$ for some $\mb z\in U(\sK)$ and $\mb r\in R$. Then
$\pi(\mb z)$ and $\pi(\mb r\mb z)=\mb r\pi(\mb z)$ are in
$\cc(\sK)$, which implies that $\pi(\mb z)=\pi(\mb r\mb z)$.
Hence, $\mb z=\mb t\mb r\mb z$ for some $\mb t\in\T^m$. We may
assume that $\mb z\in(\C^\times)^m$, so that the action of both
$R$ and $\T^m$ is free (otherwise consider the action on
$U(\lk\omega(\mb z))$).
%, like in the argument above).
It follows that $\mb t\mb r=\mathbf 1$, which implies that $\mb
r=\mathbf 1$, since $R\cap\T^m=\{\mathbf 1\}$.
\end{proof}

%\begin{remark} Our construction of a smooth structure on
%$\mathcal Z_{\sK_\Sigma}$ depends on the geometry of~$\Sigma$.
%However, we expect that the smooth structures coming from fans
%$\Sigma$ and $\Sigma'$ are the same whenever the underlying
%simplicial complexes $\sK_\Sigma$ and $\sK_{\Sigma'}$ are
%isomorphic. Equivalently, the quotients $\mathcal
%Z_{\sK_\Sigma}/\T^m$ and $\mathcal Z_{\sK_{\Sigma'}}/\T^m$ are
%diffeomorphic as manifolds with corners
%%(see the remark in the
%%proof of Theorem~\ref{zksmooth})
%whenever $\sK_\Sigma=\sK_{\Sigma'}$. It is true in the polytopal
%case (see also discussion in Section~4), and also for those fans
%$\Sigma$ which are \emph{shellable}. (A shelling order allows us
%to use an inductive argument, at each step extending a
%diffeomorphism between two $(k-1)$-balls in the boundaries of
%$k$-balls to the whole $k$-balls.)
%\end{remark}

We do not know if Theorem~\ref{zksmooth} generalises to other
sphere triangulations:

\begin{problem}
Describe the class of sphere triangulations $\sK$ for which the
moment-angle manifold $\zk$ admits a smooth structure.
\end{problem}

\begin{remark}
Even if $\zk$ admits a smooth structure for some simplicial
complexes $\sK$ not arising from fans, such a structure does not
come from a quotient $U(\sK)/R$ determined by data $\{\sK;\mb
a_1,\ldots,\mb a_m\}$. Like in the toric case (see
Section~\ref{algtq}), the $R$-action on $U(\sK)$ is proper and the
quotient $U(\sK)/R$ is Hausdorff \emph{precisely when} $\{\sK;\mb
a_1,\ldots,\mb a_m\}$ defines a fan, i.e. the simplicial cones
generated by any two subsets $\{\mb a_i\colon i\in I\}$ and $\{\mb
a_j\colon j\in J\}$ with $I,J\in\sK$ can be separated by a
hyperplane. This observation is originally due to
Bosio~\cite{bosi01}, see also~\cite[\S II.3]{a-d-h-l}
and~\cite{ba-za}.
\end{remark}

\section{Complex structures on moment-angle manifolds}\label{camam}
Let $\zk$ be the moment-angle manifold corresponding to a complete
simplicial fan $\Sigma$ defined by data $\{\sK;\mb a_1,\ldots,\mb
a_m\}$. We assume that the dimension $m+n$ of $\zk$ is even, and
set $m-n=2\ell$. This can always be achieved by adding a ghost
vertex with any corresponding vector to our data $\{\sK;\mb
a_1,\ldots,\mb a_m\}$; topologically this results in multiplying
$\zk$ by a circle. Here we show that $\zk$ admits a structure of a
complex manifold. The idea is to replace the action of
$R\cong\R^{m-n}$ on $U(\sK)$ (whose quotient is $\zk$) by a
holomorphic action of $\C^{\frac{m-n}2}$ on the same space.

We identify $\C^m$ (as a real vector space) with $\R^{2m}$ using
the map
\[
  (z_1,\ldots,z_m)\mapsto(x_1,y_1,\ldots,x_m,y_m),
\]
where $z_k=x_k+iy_k$, and consider the $\R$-linear map
\[
  \Re\colon\C^m\to\R^m,\qquad (z_1,\ldots,z_m)\mapsto(x_1,\ldots,x_m).
\]

In order to obtain a complex structure on the quotient $\zk\cong
U(\sK)/R$ we replace the action of $R$ by the action of a
holomorphic subgroup $C\subset(\C^\times)^m$ by means of the
following construction.

\begin{construction}\label{psi}
Let $\mb a_1,\ldots,\mb a_m$ be a configuration of vectors that
span $N_\R\cong\R^n$, and assume that $m-n=2\ell$. Some of the
$\mb a_i$'s may be zero. Recall the map $A\colon\R^m\to N_\R$, \
$\mb e_i\mapsto\mb a_i$.

We choose a complex $\ell$-dimensional subspace in $\C^m$ which
projects isomorphically onto the real $(m-n)$-dimensional subspace
$\Ker A\subset\R^m$. More precisely, let $\mathfrak
c\cong\C^\ell$, and choose a linear map $\varPsi\colon \mathfrak
c\to\C^m$ satisfying the two conditions:
\begin{itemize}
\item[(a)] the composite map
$\mathfrak c\stackrel{\varPsi}\longrightarrow\C^m
\stackrel{\Re}\longrightarrow\R^m$ is a monomorphism;
%$\Re\circ\varPsi\colon \mathfrak c\to\R^m$ is a monomorphism.

\item[(b)] the composite map
$\mathfrak c\stackrel{\varPsi}\longrightarrow\C^m
\stackrel{\Re}\longrightarrow\R^m \stackrel{A}\longrightarrow
N_\R$ is zero.
%$A\circ\Re\circ\varPsi=0$.
\end{itemize}
These two conditions are equivalent to the following:
\begin{itemize}
\item[(a')] $\varPsi(\mathfrak c)\cap\overline{\varPsi(\mathfrak c)\!}=\{\mathbf0\}$;
\item[(b')] $\varPsi(\mathfrak c)\subset\Ker(A_\C\colon\C^m\to N_\C)$,
\end{itemize}
where $\overline{\varPsi(\mathfrak c)\!}$ is the complex conjugate
space and $A_\C\colon\C^m\to N_\C$ is the complexification of the
real map $A\colon\R^m\to N_\R$. Consider the following commutative
diagram:
\begin{equation}\label{cdiag}
\begin{CD}
  \mathfrak c @>\varPsi>> \C^m @>\Re>> \R^m @>A>> N_\R\\
  @. @VV\exp V @VV\exp V\\% @VV\exp V\\
  \ @. (\C^\times )^m@>|\,\cdot|>> \R^m_>% @>\exp\Lambda_\R>> \R^n_>
\end{CD}
\end{equation}
where the vertical arrows are the componentwise exponential maps,
and $|\cdot|$ denotes the map
$(z_1,\ldots,z_m)\mapsto(|z_1|,\ldots,|z_m|)$. Now set
\begin{equation}\label{csigma}
  C_\varPsi=\exp\varPsi(\mathfrak c)
  =\bigl\{\bigl(e^{\langle\psi_1,\mb w\rangle},\ldots,
  e^{\langle\psi_m,\mb w\rangle}\bigr)\in(\C^\times )^m\bigr\}
\end{equation}
where $\mb w\in \mathfrak c$ and $\psi_k\in \mathfrak c^*$ is
given by the $k$th coordinate projection $\mathfrak
c\stackrel\varPsi\longrightarrow\C^m\to\C$.
%where $\mb w=(w_1,\ldots,w_\ell)\in \mathfrak c$, \ $\psi_i$ denotes the
%$i$th row of the $m\times\ell$-matrix $\varPsi=(\psi_{ij})$, and
%$\langle\psi_i,\mb
%w\rangle=\psi_{i1}w_1+\ldots+\psi_{i\ell}w_\ell$.
Then $C_\varPsi\cong\C^\ell$ is a complex-analytic (but not
algebraic) subgroup in~$(\C^\times)^m$, and therefore there is a
holomorphic action of $C_\varPsi$ on $\C^m$ and $U(\sK)$ by
restriction.
\end{construction}

\begin{example}\label{2torus}
Let $\mb a_1,\ldots,\mb a_m$ be the configuration of $m=2\ell$
zero vectors. We supplement it by the empty simplicial complex
$\sK$ on $[m]$ (with $m$ ghost vertices), so that the data
$\{\sK;\mb a_1,\ldots,\mb a_m\}$ define a complete fan in
0-dimensional space. Then $A\colon\R^m\to\R^0$ is a zero map, and
condition~(b) of Construction~\ref{psi} is void. Condition~(a)
means that $\mathfrak c\stackrel{\varPsi}\longrightarrow\C^{2\ell}
\stackrel{\mathrm{Re}}\longrightarrow\R^{2\ell}$ is an isomorphism
of real spaces.

Consider the quotient $(\C^\times)^m/C_\varPsi$ (note that
$U(\sK)=(\C^\times)^m$ in our case). The exponential map
$\C^m\to(\C^\times)^m$ identifies $(\C^\times)^m$ with the
quotient of $\C^m$ by the imaginary lattice $\Gamma=\Z\langle2\pi
i\mb e_1,\ldots,2\pi i\mb e_m\rangle$. Condition~(a) implies that
the projection $p\colon\C^m\to\C^m/\varPsi(\mathfrak c)$ is
nondegenerate on the imaginary subspace of~$\C^m$. In particular,
$p\,(\Gamma)$ is a lattice of rank $m=2\ell$ in
$\C^m/\varPsi(\mathfrak c)\cong\C^\ell$. Therefore,
\[
  (\C^\times)^m/C_\varPsi\cong\bigl(\C^m/\Gamma\bigr)/\varPsi(\mathfrak c)
  =\bigl(\C^m/\varPsi(\mathfrak c)\bigr)
  \big/p\,(\Gamma)\cong\C^\ell/p\,(\Gamma)
\]
is a complex compact $\ell$-dimensional torus.

Any complex torus can be obtained in this way. Indeed, let
$\varPsi\colon\mathfrak c\to\C^m$ be given by an
$2\ell\times\ell$-matrix $\begin{pmatrix}-B\\I\end{pmatrix}$ where
$I$ is the unit matrix and $B$ is a square matrix of size~$\ell$.
Then $p\colon\C^m\to\C^m/\varPsi(\mathfrak c)$ is given by the
matrix $(I\,B)$ in appropriate bases, and
$(\C^\times)^m/C_\varPsi$ is isomorphic to the quotient of
$\C^\ell$ by the lattice $\Z\langle\mb e_1,\ldots,\mb e_\ell,\mb
b_1,\ldots,\mb b_\ell\rangle$, where $\mb b_k$ is the $k$th column
of~$B$. (Condition~(b) implies that the imaginary part of $B$ is
nondegenerate.)

For example, if $\ell=1$, then $\varPsi\colon\C\to\C^2$ is given
by $w\mapsto(\beta w,w)$ for some $\beta\in\C$, so that
subgroup~\eqref{csigma} is
\[
  C_\varPsi=\{(e^{\beta w},e^w)\}\subset(\C^\times )^2.
\]
Condition~(a) implies that $\beta\notin\R$. Then
$\exp\varPsi\colon\C\to(\C^\times )^2$ is an embedding, and
\[
  (\C^\times )^2/C_\varPsi\cong\C/(\Z\oplus\beta\Z)=T^1_\C(\beta)
\]
is a complex 1-dimensional torus with lattice parameter
$\beta\in\C$.
\end{example}

\begin{theorem}[\cite{pa-us12}]\label{zkcomplex}
Assume that data $\{\sK;\mb a_1,\ldots,\mb a_m\}$ define a
complete fan~$\Sigma$ in~$N_\R\cong\R^n$, and $m-n=2\ell$. Let
$C_\varPsi\cong\C^\ell$ be the group given by~\eqref{csigma}. Then
\begin{itemize}
\item[(a)]
the holomorphic action of $C_\varPsi$ on $U(\sK)$ is free and
proper, and the quotient $U(\sK)/C_\varPsi$ has a structure of a
compact complex manifold;

\item[(b)] $U(\sK)/C_\varPsi$ is diffeomorphic to the moment-angle manifold~$\zk$.
\end{itemize}
Therefore, $\zk$ has a complex structure, in which each element of
$\T^m$ acts by a holomorphic transformation.
\end{theorem}

\begin{remark}
A result similar to Theorem~\ref{zkcomplex} was obtained by
Tambour~\cite{tamb12}. The approach of Tambour was somewhat
different; he constructed complex structures on manifolds $\zk$
arising from \emph{rationally}\label{starshsphd} starshaped
spheres~$\sK$ (underlying complexes of complete rational
simplicial fans) by relating them to a class of generalised
LVM-manifolds described by Bosio in~\cite{bosi01}.
\end{remark}

\begin{proof}[Proof of Theorem~\ref{zkcomplex}]
We first prove statement (a). The stabilisers of the
$(\C^\times)^m$-action on $U(\sK)$ are of the form
$(\C^\times,1)^I$ for $I\in\sK$. In order to show that
$C_\varPsi\subset(\C^\times)^m$ acts freely we need to check that
$C_\varPsi$ has trivial intersection with any stabiliser of
the~$(\C^\times)^m$-action. Since $C_\varPsi$ embeds into $\R^m_>$
by~\eqref{cdiag}, it enough to check that the image of $C_\varPsi$
in $\R^m_>$ intersects the image of $(\C^\times,1)^I$ in $\R^m_>$
trivially. The former image is $R$ and the latter image is
$(\R_>,1)^I$; the triviality of their intersection follows from
Theorem~\ref{zksmooth}~(a).

Now we prove the properness of this action. Consider the
projection $\pi\colon U(\sK)\to(\R_\ge,\R_>)^\sK$ onto the
quotient of the $\T^m$-action, and the commutative square
\[
\begin{CD}
  C_{\varPsi}\times U(\sK) @>h_\C>> U(\sK)\times U(\sK)\\
  @VVf\times\pi V @VV\pi\times\pi V\\
  R\times (\R_\ge,\R_>)^\sK@>h_\R>> (\R_\ge,\R_>)^\sK\times(\R_\ge,\R_>)^\sK
\end{CD}
\]
where $h_\C$ and $h_\R$ denote the group action maps, and $f\colon
C_{\varPsi}\to R$ is the isomorphism given by the restriction of
$|\cdot|\colon (\C^\times )^m\to \R_>^m$. The preimage
$h^{-1}_\C(V)$ of a compact subset $V\in U(\sK)\times U(\sK)$ is a
closed subset in $W=(f\times\pi)^{-1}\circ
h_\R^{-1}\circ(\pi\times\pi)(V)$. The image $(\pi\times\pi)(V)$ is
compact, the action of $R$ on $(\R_\ge,\R_>)^\sK$ is proper by
Theorem~\ref{zksmooth}~(b), and the map $f\times \pi$ is proper as
the quotient projection for a compact group action. Hence, $W$ is
a compact subset in $C_{\varPsi}\times U(\sK)$, and $h^{-1}_\C(V)$
is compact as a closed subset in~$W$.

The group $C_{\varPsi}\cong\C^l$ acts holomorphically, freely and
properly on the complex manifold $U(\sK)$, therefore the quotient
manifold $U(\sK)/C_\varPsi$ has a complex structure.

As in the proof of Theorem~\ref{zksmooth}, it is possible to
describe a holomorphic atlas of $U(\sK)/C_{\varPsi}$. Since the
action of $C_{\varPsi}$ on the quotient
$U(\sK)/\T^m=(\R_\ge,\R_>)^\sK$ coincides with the action of $R$
on the same space, the quotient of $U(\sK)/C_{\varPsi}$ by the
action of $\T^m$ has exactly the same structure of a smooth
manifold with corners as the quotient of $U(\sK)/R$ by $\T^m$ (see
the proof of Theorem~\ref{zksmooth}). This structure is determined
by the atlas $\{q_I\colon(\R_\ge,\R_>)^I/R\to\R^n_\ge\}$, which
lifts to a covering of $U(\sK)/C_{\varPsi}$ by the open subsets
$(\C,\C^\times)^I/C_{\varPsi}$. For any $I\in\sK$, the subset
$(\C,\T)^I\subset(\C,\C^\times)^I$ intersects each orbit of the
$C_{\varPsi}$-action on $(\C,\C^\times)^I$ transversely at a
single point. Therefore, every
$(\C,\C^\times)^I/C_{\varPsi}\cong(\C,\T)^I$ acquires a structure
of a complex manifold. Since
$(\C,\C^\times)^I\cong\C^n\times(\C^\times)^{m-n}$, and the action
of $C_{\varPsi}$ on the $(\C^\times)^{m-n}$ factor is free, the
complex manifold $(\C,\C^\times)^I/C_{\varPsi}$ is the total space
of a holomorphic $\C^n$-bundle over the complex torus
$(\C^\times)^{m-n}/C_{\varPsi}$ (see Example~\ref{2torus}).
Writing trivialisations of these $\C^n$-bundles for every~$I$, we
obtain a holomorphic atlas for $U(\sK)/C_{\varPsi}$.

The proof of statement (b) follows the lines of the proof of
Theorem~\ref{zksmooth}~(b). We need to show that each
$C_{\varPsi}$-orbit intersects $\zk\subset U(\sK)$ at a single
point. First we show that the $C_{\varPsi}$-orbit of any point in
$U(\sK)/\T^m$ intersects $\zk/\T^m=\cc(\sK)$ at a single point;
this follows from the fact that the actions of $C_{\varPsi}$ and
$R$ coincide on $U(\sK)/\T^m$. Then we show that each
$C_{\varPsi}$-orbit intersects the preimage $\pi^{-1}(\cc(\sK))$
at a single point, using the fact that $C_{\varPsi}$ and $\T^m$
have trivial intersection in $(\C^\times)^m$.
\end{proof}

\begin{example}[Hopf manifold]\label{hopf}
Let $\mb a_1,\ldots,\mb a_{n+1}$ be a set of vectors which span
$N_\R\cong\R^n$ and satisfy a linear relation $\lambda_1\mb
a_1+\cdots+\lambda_{n+1}\mb a_{n+1}=\mathbf0$ with all
$\lambda_k>0$. Let $\Sigma$ be the complete simplicial fan in
$N_\R$ whose cones are generated by all proper subsets of $\mb
a_1,\ldots,\mb a_{n+1}$. To make $m-n$ even we add one more ghost
vector $\mb a_{n+2}$. Hence $m=n+2$, $\ell=1$, and we have one
more linear relation $\mu_1\mb a_1+\cdots+\mu_{n+1}\mb a_{n+1}+\mb
a_{n+2}=\mathbf0$ with $\mu_k\in\R$. The subspace $\Ker
A\subset\R^{n+2}$ is spanned by
$(\lambda_1,\ldots,\lambda_{n+1},0)$ and
$(\mu_1,\ldots,\mu_{n+1},1)$.

Then $\sK=\sK_\Sigma$ is the boundary of an $n$-dimensional
simplex with $n+1$ vertices and one ghost vertex, $\zk\cong
S^{2n+1}\times S^1$, and
$U(\sK)=(\C^{n+1}\setminus\{{\bf0}\})\times\C^\times$.

Conditions~(a) and~(b) of Construction~\ref{psi} imply that
$C_\varPsi$ is a 1-dimensional subgroup in $(\C^\times )^m$ given
in appropriate coordinates by
\[
  C_\varPsi=\bigl\{(e^{\zeta_1 w},\ldots,e^{\zeta_{n+1}w},e^w)
  \colon w\in\C\bigl\}\subset(\C^\times )^m,
\]
where $\zeta_k=\mu_k+\alpha\lambda_k$ for some
$\alpha\in\C\setminus\R$. By changing the basis of $\Ker A$ if
necessary, we may assume that $\alpha=i$. The moment-angle
manifold $\zk\cong S^{2n+1}\times S^1$ acquires a complex
structure as the quotient $U(\sK)/C_\varPsi$:
\begin{multline*}
  \bigl(\C^{n+1}\setminus\{\mathbf 0\}\bigr)\times\C^\times\bigl/\;
  \bigl\{(z_1,\ldots,z_{n+1},t)\!\sim
  (e^{\zeta_1w}z_1,\ldots,e^{\zeta_{n+1}w}z_{n+1},
  e^w t)\bigr\}
  \\
  \cong\bigl(\C^{n+1}\setminus\{\mathbf0\}\bigr)\bigl/\;
  \bigl\{(z_1,\ldots,z_{n+1})\!\sim
  (e^{2\pi i\zeta_1}z_1,\ldots,
  e^{2\pi i\zeta_{n+1}}z_{n+1})\bigr\},
\end{multline*}
where $\mb z\in\C^{n+1}\setminus\{\mathbf 0\}$, $t\in\C^\times$.
The latter is the quotient of $\C^{n+1}\setminus\{\mathbf 0\}$ by
a diagonalisable action of $\Z$. It is known as a \emph{Hopf
manifold}. For $n=0$ we obtain the complex torus of
Example~\ref{2torus}.
\end{example}

Theorem~\ref{zkcomplex} can be generalised to the quotients of
$\zk$ by freely acting subgroups $H\subset\T^m$, or \emph{partial
quotients} of~$\zk$ in the sense of~\cite[\S7.5]{bu-pa02}. These
include both toric manifolds and LVM-manifolds:

\begin{construction}\label{psi_part}
Let $\Sigma$ be a complete simplicial fan in $N_\R$ defined by
data $\{\sK;\mb a_1,\ldots,\mb a_m\}$, and let $H\subset\T^m$ be a
subgroup which acts freely on the corresponding moment-angle
manifold~$\zk$. Then $H$ is a product of a torus and a finite
group, and $h=\dim H\le m-n$ by Proposition~\ref{freeaction} ($H$
must intersect trivially with an $n$-dimensional coordinate
subtorus in~$\T^m$). Under an additional assumption on~$H$, we
shall define a holomorphic subgroup $D$ in $(\C^\times)^m$ and
introduce a complex structure on $\zk/H$ by identifying it with
the quotient $U(\sK)/D$.

The additional assumption is the compatibility with the fan data.
Recall the map $A\colon \R^m\to N_\R$, \ $\mb e_i\mapsto\mb a_i$,
and let $\mathfrak h\subset\R^m$ be the Lie algebra of
$H\subset\T^m$. We assume that $\mathfrak h\subset\Ker A$. We also
assume that $2\ell=m-n-h$ is even (this can be satisfied by adding
a zero vector to $\mb a_1,\ldots,\mb a_m$). Let $T=\T^m/H$ be the
quotient torus, $\mathfrak t$ its Lie algebra, and
$\rho\colon\R^m\to\mathfrak t$ the map of Lie algebras
corresponding to the quotient projection $\T^m\to T$.

Let $\mathfrak c\cong\C^\ell$, and choose a linear map
$\varOmega\colon \mathfrak c\to\C^m$ satisfying the two
conditions:
\begin{itemize}
\item[(a)] the composite map
$\mathfrak c\stackrel{\varOmega}\longrightarrow\C^m
\stackrel{\Re}\longrightarrow\R^m \stackrel{\rho}\longrightarrow
\mathfrak t$ is a monomorphism;
\item[(b)] the composite map
$\mathfrak c\stackrel{\varOmega}\longrightarrow\C^m
\stackrel{\Re}\longrightarrow\R^m \stackrel{A}\longrightarrow
N_\R$ is zero.
%$A\circ\Re\circ\,\varOmega=0$.
\end{itemize}
Equivalently, choose a complex subspace $\mathfrak
c\subset\mathfrak t_\C$ such that the composite map $\mathfrak
c\to\mathfrak t_\C\stackrel{\Re}\longrightarrow\mathfrak t$ is a
monomorphism.

As in Construction~\ref{psi}, $\exp\varOmega(\mathfrak
c)\subset(\C^\times)^m$ is a holomorphic subgroup isomorphic
to~$\C^\ell$. Let $H_\C\subset(\C^\times)^m$ be the
complexification of $H$ (it is a product of $(\C^\times)^h$ and a
finite group). It follows from (a) that the subgroups $H_\C$ and
$\exp\varOmega(\mathfrak c)$ intersect trivially
in~$(\C^\times)^m$. We can define a complex $(h+\ell)$-dimensional
subgroup
\begin{equation} \label{csigma_part}
  D_{H,\varOmega}=H_\C\times\exp\varOmega(\mathfrak c)
  \subset(\C^\times)^m.
\end{equation}
\end{construction}

\begin{theorem}[{\cite[Theorem~3.7]{pa-us12}}]\label{zkcomplex_part}
Let $\Sigma$, $\sK$ and $D_{H,\varOmega}$ be as above. Then
\begin{itemize}
\item[(a)] the holomorphic action of the group
$D_{H,\varOmega}$ on $U(\sK)$ is free and proper, and the quotient
$U(\sK)/D_{H,\varOmega}$ has a structure of a compact complex
manifold of complex dimension $m-h-\ell$;

\item[(b)] there is a diffeomorphism between
$U(\sK)/D_{H,\varOmega}$ and $\zk/H$ defining a complex structure
on the quotient $\zk/H$, in which each element of $T=\T^m/H$ acts
by a holomorphic transformation.
\end{itemize}
\end{theorem}

The proof is similar to that of Theorem~\ref{zkcomplex} and is
omitted.

\begin{example}\

1. If $H$ is trivial ($h=0$) then we obtain
Theorem~\ref{zkcomplex}.

2. Let $H$ be the diagonal circle in~$\T^m$. The condition
$\mathfrak h\subset\Ker A_\R$ implies that the vectors $\mb
a_1,\ldots,\mb a_m$ sum up to zero, which can always be achieved
by rescaling them (as $\Sigma$ is a complete fan). As the result,
we obtain a complex structure on the quotient $\zk/S^1$ by the
diagonal circle in~$\T^m$, provided that $m-n$ is odd. In the
polytopal case $\sK=\sK_P$, the quotient $\zk/S^1$ embeds into
$\C^m\setminus\{\mathbf 0\}/\C^\times=\C P^{m-1}$ as an
intersection of homogeneous quadrics~\eqref{ndef}, and the complex
structure on $\zk/S^1$ coincides with that of an
\emph{LVM-manifold}\label{lvmmanicomp}, see Section~\ref{lvmma}.

3. Let $h=\dim H=m-n$. Then $\mathfrak h=\Ker A$. Since $\mathfrak
h$ is the Lie algebra of a torus, the $(m-n)$-dimensional subspace
$\Ker A\subset\R^m$ is rational. By Gale duality, this implies
that the fan $\Sigma$ is also rational. We have $\ell=0$,
$D_{H,\varOmega}=H_\C\cong(\C^\times)^{m-n}$ and
$U(\sK)/H_\C=\zk/H$ is the toric variety corresponding
to~$\Sigma$.
\end{example}

An effective action of $T^k$ on an $m$-dimensional manifold $M$ is
called \emph{maximal}\label{maxitorusact} if there exists a point
$x\in M$ whose stabiliser has dimension $m-k$; the two extreme
cases are the free action of a torus on itself and the
half-dimensional torus action on a toric manifold. As it is shown
by Ishida~\cite{isid}, any compact complex manifold with a maximal
effective holomorphic action of a torus is biholomorphic to a
quotient $\zk/H$ of a moment-angle manifold with a complex
structure described by Theorem~\ref{zkcomplex_part}. The argument
of~\cite{isid} recovering a fan $\Sigma$ from a maximal
holomorphic torus action builds up on the works~\cite{i-f-m13}
and~\cite{is-ka}, where the result was proved in particular cases.
The main result of~\cite{is-ka} provides a purely complex-analytic
description of toric manifolds~$V_\Sigma$:

\begin{theorem}[{\cite[Theorem~1]{is-ka}}]\label{iskath}
Let $M$ be a compact connected complex manifold of complex
dimension $n$, equipped with an effective action of $T^n$ by
holomorphic transformations. If the action has fixed points, then
there exists a complete regular fan $\Sigma$ and a
$T^n$-equivariant biholomorphism of $V_\Sigma$ with~$M$.
\end{theorem}

\section{Holomorphic principal bundles and Dolbeault
cohomology} In the case of rational simplicial normal fans
$\Sigma_P$ a construction of Meersseman--Verjovsky~\cite{me-ve04}
identifies the corresponding projective toric variety $V_P$ as the
base of a holomorphic principal \emph{Seifert fibration}, whose
total space is the moment-angle manifold~$\zp$ equipped with a
complex structure of an LVM-manifold, and fibre is a compact
complex torus of complex dimension $\ell=\frac{m-n}2$. (Seifert
fibrations are generalisations of holomorphic fibre bundles to the
case when the base is an orbifold.) If $V_P$ is a projective toric
manifold, then there is a holomorphic free action of a complex
$\ell$-dimensional torus $T^{\ell}_\C$ on $\zp$ with
quotient~$V_P$.

Using the construction of a complex structure on $\zk$ described
in the previous section, in~\cite{pa-us12} holomorphic (Seifert)
fibrations with total space $\zk$ were defined for arbitrary
complete rational simplicial fans~$\Sigma$. By an application of
the Borel spectral sequence to the holomorphic fibration $\zk\to
V_\Sigma$, the Dolbeault cohomology of $\zk$ can be described and
some Hodge numbers can be calculated explicitly.
%We review these and related results in this section.

Here we make additional assumption that the set of integral linear
combinations of the vectors $\mb a_1,\ldots,\mb a_m$ is a
full-rank lattice (a discrete subgroup isomorphic to~$\Z^n$)
in~$N_\R\cong\R^n$. We denote this lattice by $N_\Z$ or
simply~$N$. This assumption implies that the complete simplicial
fan $\Sigma$ defined by the data $\{\sK;\mb a_1,\ldots,\mb a_m\}$
is rational. We also continue assuming that $m-n$ is even and
setting $\ell=\frac{m-n}2$.

Because of our rationality assumption, the algebraic group $G$ is
defined by~\eqref{gexpl}. Furthermore, since we defined $N$ as the
lattice generated by $\mb a_1,\ldots,\mb a_m$, the group $G$ is
isomorphic to~$(\C^\times)^{2\ell}$ (i.e. there are no finite
factors). We also observe that $C_{\varPsi}$ lies in~$G$ as an
$\ell$-dimensional complex subgroup. This follows from
condition~(b') of Construction~\ref{psi}.

The quotient construction (Section~\ref{algtq}) identifies the
toric variety $V_\Sigma$ with $U(\sK)/G$, provided that $\mb
a_1,\ldots,\mb a_m$ are \emph{primitive} generators of the edges
of~$\Sigma$. In our data $\{\sK;\mb a_1,\ldots,\mb a_m\}$, the
vectors $\mb a_1,\ldots,\mb a_m$ are not necessarily primitive in
the lattice $N$ generated by them. Nevertheless, the quotient
$U(\sK)/G$ is still isomorphic to~$V_\Sigma$, see~\cite[Ch.~II,
Proposition~3.1.7]{a-d-h-l}. Indeed, let $\mb a'_i\in N$ be the
primitive generator along $\mb a_i$, so that $\mb a_i=r_i\mb a'_i$
for some positive integer~$r_i$. Then we have a finite branched
covering
\[
  U(\sK)\to U(\sK),\quad(z_1,\ldots,z_m)\mapsto
  (z_1^{r_1},\ldots,z_m^{r_m}),
\]
which maps the group $G$ defined by $\mb a_1,\ldots,\mb a_m$ to
the group $G'$ defined by $\mb a'_1,\ldots,\mb a'_m$,
see~\eqref{gexpl}. We therefore obtain a covering $U(\sK)/G\to
U(\sK)/G'$ of the toric variety $V_\Sigma\cong U(\sK)/G\cong
U(\sK)/G'$ over itself. Having this in mind, we can relate the
quotients $V_\Sigma\cong U(\sK)/G$ and $\zk\cong U(\sK)/C_\varPsi$
as follows:

\begin{proposition}\label{toricfib}Assume that data
$\{\sK;\mb a_1,\ldots,\mb a_m\}$ define a complete simplicial
rational fan~$\Sigma$, and let $G$ and $C_\varPsi$ be the groups
defined by~\eqref{gexpl} and~\eqref{csigma}.
\begin{itemize}
\item[(a)]The toric variety $V_\Sigma$ is identified, as a
topological space, with the quotient of $\zk$ by the holomorphic
action of the complex compact torus~$G/C_{\varPsi}$.

\item[(b)]If the fan $\Sigma$ is regular, then $V_\Sigma$ is the base of a
holomorphic principal bundle with total space $\zk$ and fibre the
complex compact torus $G/C_{\varPsi}$.
\end{itemize}
\end{proposition}
\begin{proof}
To prove~(a) we just observe that
\[
  V_\Sigma=U(\sK)/G=
  \bigl(U(\sK)/C_{\varPsi}\bigr)\big/(G/C_{\varPsi})\cong
  \mathcal Z_{\sK}\big/(G/C_{\varPsi}),
\]
where we used Theorem~\ref{zkcomplex}. The quotient
$G/C_{\varPsi}$ is a compact complex $\ell$-torus by
Example~\ref{2torus}. To prove~(b) we observe that the holomorphic
action of $G$ on $U(\sK)$ is free by Proposition~\ref{freeaction},
and the same is true for the action of $G/C_{\varPsi}$ on~$\zk$. A
holomorphic free action of the torus $G/C_{\varPsi}$ gives rise to
a principal bundle.
\end{proof}

\begin{remark}
As in the projective situation of~\cite{me-ve04}, if the fan
$\Sigma$ is not regular, then the quotient projection $\zk\to
V_\Sigma$ of Proposition~\ref{toricfib}~(a) is a holomorphic
principal \emph{Seifert fibration}\label{seifertfib} for an
appropriate orbifold structure on~$V_\Sigma$.
\end{remark}

Let $M$ be a complex $n$-dimensional manifold. The space
$\varOmega_\C^*(M)$ of complex differential forms on $M$
decomposes into a direct sum of the subspaces of
\emph{$(p,q)$-forms}, $\varOmega_\C^*(M)=\bigoplus_{0\le p,q\le
n}\varOmega^{p,q}(M)$, and there is the \emph{Dolbeault
differential}\label{dolbcohom}
$\bar\partial\colon\varOmega^{p,q}(M)\to \varOmega^{p,q+1}(M)$.
The dimensions $h^{p,q}(M)$ of the Dolbeault cohomology groups
$H_{\bar\partial}^{p,q}(M)$, $0\le p,q\le n$, are known as the
\emph{Hodge numbers} of~$M$. They are important invariants of the
complex structure of~$M$.

The Dolbeault cohomology of a compact complex $\ell$-torus
$T_\C^{\ell}$ is isomorphic to an exterior algebra on $2\ell$
generators:
\begin{equation}\label{dolbtorus}
  H_{\bar\partial}^{*,*}(T_\C^{\ell})\cong
  \Lambda[\xi_1,\ldots,\xi_\ell,\eta_1,\ldots,\eta_\ell],
\end{equation}
where $\xi_1,\ldots,\xi_\ell\in
H_{\bar\partial}^{1,0}(T_\C^{\ell})$ are the classes of basis
holomorphic 1-forms, and $\eta_1,\ldots,\eta_\ell\in
H_{\bar\partial}^{0,1}(T_\C^{\ell})$ are the classes of basis
antiholomorphic 1-forms. In particular, the Hodge numbers are
given by $h^{p,q}(T_\C^{\ell})=\binom\ell p\binom\ell q$.

The de Rham cohomology of a toric manifold $V_\Sigma$ admits a
Hodge decomposition with only nontrivial components of
bidegree~$(p,p)$, $0\le p\le n$~{\cite[\S12]{dani78}}. This
together with Theorem~\ref{danjur} gives the following description
of the Dolbeault cohomology:
\begin{equation}\label{dolbtoric}
  H_{\bar\partial}^{*,*}(V_\Sigma)\cong
  \C[v_1,\ldots,v_m]/(\mathcal I_{\sK}+\mathcal J_{\Sigma}),
\end{equation}
where $v_i\in H_{\bar\partial}^{1,1}(V_\Sigma)$ are the cohomology
classes corresponding to torus-invariant divisors (one for each
one-dimensional cone of~$\Sigma$), the ideal $\mathcal I_{\sK}$ is
generated by the monomials $v_{i_1}\cdots v_{i_k}$ for which $\mb
a_{i_1},\ldots,\mb a_{i_k}$ do not span a cone of $\Sigma$ (the
Stanley--Reisner ideal of~$\sK$), and $\mathcal J_{\Sigma}$ is
generated by the linear forms $\sum_{j=1}^m\langle\mb a_j,\mb
u\rangle v_j$, $\mb u\in N^*$. For the Hodge numbers,
$h^{p,p}(V_\Sigma)=h_p$, where $(h_0,h_1,\ldots,h_n)$ is the
$h$-vector of~$\sK$, and $h^{p,q}(V_\Sigma)=0$ for $p\ne q$.

\begin{theorem}[{\cite{pa-us12}}]\label{dolbzp}
Assume that data $\{\sK;\mb a_1,\ldots,\mb a_m\}$ define a
complete rational regular fan~$\Sigma$ in~$N_\R\cong\R^n$,
$m-n=2\ell$, and let $\zk$ be the corresponding moment-angle
manifold with a complex structure defined by
Theorem~\ref{zkcomplex}. Then the Dolbeault cohomology algebra
$H_{\bar\partial}^{*,*}(\zk)$ is isomorphic to the cohomology of
the differential bigraded algebra
\begin{equation}\label{zkmult}
\bigl[\Lambda[\xi_1,\ldots,\xi_\ell,\eta_1,\ldots,\eta_\ell]\otimes
  H_{\bar\partial}^{*,*}(V_\Sigma),d\bigr]
\end{equation}
%whose bigrading is defined by~\eqref{dolbtorus}
%and~\eqref{dolbtoric},
with differential $d$ of bidegree $(0,1)$ defined on the
generators as follows:
\[
  dv_i=d\eta_j=0,\quad d\xi_j=c(\xi_j),\quad
  1\le i\le m,\;1\le j\le\ell,
\]
where $c\colon H^{1,0}_{\bar\partial}(T_\C^{\ell})\to
H^2(V_\Sigma, \C)=  H_{\bar\partial}^{1,1}(V_\Sigma)$ is the first
Chern class map of the principal $T^\ell_\C$-bundle $\zk\to
V_\Sigma$.
\end{theorem}
\begin{proof}
We use the notion of a minimal Dolbeault model of a complex
manifold~\cite[\S4.3]{f-o-t08}. Let $[B,d_B]$ be such a model
for~$V_\Sigma$, i.e. $[B,d_B]$ is a minimal commutative bigraded
differential algebra together with a quasi-isomorphism $f\colon
B^{*,*}\to \varOmega^{*,*}(V_\Sigma)$.
%(i.e. $f$ commutes with the differentials $d_B$ and
%$\bar\partial$, and induces an isomorphism in cohomology).
Consider the differential bigraded algebra
\begin{equation}\label{dolbmod}
\begin{gathered}
  \bigl[\Lambda[\xi_1,\ldots,\xi_\ell,\eta_1,\ldots,\eta_\ell]\otimes
  B, d\bigr],\qquad\text{where}\\
  d|_B=d_B,\quad  d(\xi_i)=c(\xi_i)\in
  B^{1,1}= H^{1,1}_{\bar\partial}(V_\Sigma),\quad d(\eta_i)=0.
\end{gathered}
\end{equation}
By \cite[Corollary~4.66]{f-o-t08}, this is a model for the
Dolbeault cohomology algebra of the total space $\zk$ of the
principal $T^{\ell}_\C$-bundle $\zk\to V_\Sigma$, provided that
$V_\Sigma$ is strictly formal. Recall
from~\cite[Definition~4.58]{f-o-t08} that a complex manifold $M$
is \emph{strictly formal} if there exists a differential bigraded
algebra $[Z,\delta]$ together with quasi-isomorphisms
\[
\xymatrix{
  [\varOmega^{*,*},\bar\partial]  &
  [Z,\delta] \ar[l]_{\simeq} \ar[r]^{\simeq} \ar[d]^{\simeq}&
  [\varOmega^{*},d_{\mathrm{DR}}]\\
  & [H^{*,*}_{\bar\partial}(M),0]
}
\]
linking together the de Rham algebra, the Dolbeault algebra and
the Dolbeault cohomology.

According to~\cite[Corollary~7.2]{pa-ra08}, the toric manifold
$V_\Sigma$ is formal in the usual (de Rham) sense. Also, the above
mentioned Hodge decomposition of~\cite[\S12]{dani78} implies that
$V_\Sigma$ satisfies the
$\partial\bar\partial$-lemma~\cite[Lemma~4.24]{f-o-t08}. Therefore
$V_\Sigma$ is strictly formal by the same argument
as~\cite[Theorem~4.59]{f-o-t08}, and~\eqref{dolbmod} is a model
for its Dolbeault cohomology.

The usual formality of $V_\Sigma$ implies the existence of a
quasi-isomorphism $\phi_B\colon B\to
H_{\bar\partial}^{*,*}(V_\Sigma)$, which extends to a
quasi-isomorphism
\[
  \mbox{id}\otimes\phi_B\colon
  \bigl[\Lambda[\xi_1,\ldots,\xi_\ell,\eta_1,\ldots,\eta_\ell]
  \otimes B, d\bigr]\to
  \bigl[\Lambda[\xi_1,\ldots,\xi_\ell,\eta_1,\ldots,\eta_\ell]
  \otimes H_{\bar\partial}^{*,*}(V_\Sigma), d\bigr]
\]
by~\cite[Lemma~14.2]{f-h-t01}. Thus, the differential algebra
$\bigl[\Lambda[\xi_1,\ldots,\xi_\ell,\eta_1,\ldots,\eta_\ell]\otimes
H^{*,*}_{\bar\partial}(V_\Sigma), d\bigr]$ provides a model for
the Dolbeault cohomology of $\zk$, as claimed.
\end{proof}

\begin{remark}
If $V_\Sigma$ is projective, then it is K\"ahler; in this case the
model of Theorem~\ref{dolbzp} coincides with the model for the
Dolbeault cohomology of the total space of a holomorphic torus
principal bundle over a K\"ahler
manifold~\cite[Theorem~4.65]{f-o-t08}.
\end{remark}

The first Chern class map $c$ from Theorem~\ref{dolbzp} can be
described explicitly in terms of the map $\varPsi$ defining the
complex structure on $\zk$. We recall the map $A_\C\colon\C^m\to
N_\C$, $\mb e_i\mapsto\mb a_i$ and the Gale dual $(m-n)\times
m$-matrix $\varGamma=(\gamma_{jk})$ whose rows form a basis of
linear relations between $\mb a_1,\ldots,\mb a_m$. By
Construction~\ref{psi}, $\Im\varPsi\subset\Ker A_\C$. Denote by
$\Ann U$ the annihilator of a linear subspace $U\subset\C^m$, i.e.
the subspace of linear functions on $\C^m$ vanishing on~$U$.

\begin{lemma}\label{mumatrix}
Let $k$ be the number of zero vectors among $\mb a_1,\ldots,\mb
a_m$. The first Chern class map
\[
  c\colon H^{1,0}_{\bar\partial}(T_\C^{\ell})\to H^2(V_\Sigma, \C)=
  H_{\bar\partial}^{1,1}(V_\Sigma)
\]
of the principal $T^\ell_\C$-bundle $\zk\to V_\Sigma$ is given by
the composition
\[
\begin{CD}
  \Ann\Im\varPsi/\Ann\Ker A_\C @>i>> \C^m/\Ann\Ker A_\C
  @>p>> \C^{m-k}/\Ann\Ker A_\C
\end{CD}
\]
where $i$ is the inclusion and $p$ is the projection forgetting
the coordinates in $\C^m$ corresponding to zero vectors.

Explicitly, the map $c$ is given on the generators of
$H^{1,0}_{\bar\partial}(T_\C^{\ell})$ by
\[
  c(\xi_j)=\mu_{j1}v_1+\cdots+\mu_{jm}v_m,\quad 1\le j\le \ell,
\]
where $M=(\mu_{jk})$ is an $\ell\times m$-matrix satisfying the
two conditions:
\begin{itemize}
\item[(a)] $\varGamma M^t\colon\C^\ell\to\C^{2\ell}$ is a
monomorphism;

\item[(b)] $M\varPsi=0$.
\end{itemize}
\end{lemma}
\begin{proof}
Let $A^*_\C\colon N_\C^*\to\C^m$, $\mb u\mapsto(\langle\mb a_1,\mb
u\rangle,\ldots,\langle\mb a_m,\mb u\rangle)$, be the dual map. We
have $H^1(T_\C^{\ell};\C)=(\Ker A_\C)^*=\C^m/\Im A^*_\C$ and
$H^2(V_\Sigma;\C)=\C^{m-k}/\Im A^*_\C$. The first Chern class map
$c\colon H^1(T_\C^{\ell};\C)\to H^2(V_\Sigma;\C)$ (the
transgression) is then given by $p\colon \C^m/\Im
A^*_\C\to\C^{m-k}/\Im A^*_\C$. In order to separate the
holomorphic part of $c$ we need to identify the subspace of
holomorphic differentials
$H^{1,0}_{\bar\partial}(T_\C^{\ell})\cong\C^\ell$ inside the space
of all 1-forms $H^1(T_\C^{\ell};\C)\cong\C^{2\ell}$. Since
\[
  T_\C^{\ell}=G/C_{\varPsi}=(\Ker\exp A_\C)/(\exp\Im\varPsi),
\]
holomorphic differentials on $T_\C^{\ell}$ correspond to
$\C$-linear functions on $\Ker A_\C$ which vanish on $\Im\varPsi$.
The space of functions on $\Ker A_\C$ is $\C^m/\Im
A^*_\C=\C^m/\Ann\Ker A_\C$, and the functions vanishing on
$\Im\varPsi$ form the subspace $\Ann\Im\varPsi/\Ann\Ker A_\C$.
Condition~(b) says exactly that the linear functions on $\C^m$
corresponding to the rows of $M$ vanish on $\Im\varPsi$.
Condition~(a) says that the rows of $M$ constitute a basis in the
complement of $\Ann\Ker A_\C$ in $\Ann\Im\varPsi$.
\end{proof}

It is interesting to compare Theorem~\ref{dolbzp} with the
following description of the de Rham cohomology of~$\zk$:

\begin{theorem}%[{\cite[Th.~7.36]{bu-pa02}}]
\label{cohomzpred}
Let $\zk$ and $V_\Sigma$ be as in
Theorem~\ref{dolbzp}. The de Rham cohomology $H^*(\zk)$ is
isomorphic to the cohomology of the differential graded algebra
\[
  \bigl[\Lambda[u_1,\ldots,u_{m-n}]\otimes
  H^*(V_\Sigma),d\bigr],
\]
with $\deg u_j=1$, $\deg v_i=2$, and differential $d$ defined on
the generators as
\[
  dv_i=0,\quad du_j=\gamma_{j1}v_1+\cdots+\gamma_{jm}v_m,\quad
  %1\le i\le m,\quad
  1\le j\le m-n.
\]
\end{theorem}
\begin{proof}
The de Rham cohomology of the manifold $\zk$ is isomorphic to its
cellular cohomology (with coefficients in~$\R$). By
Theorem~\ref{zkcoh},
\[
  H^*(\zk)\cong
  \Tor_{\R[v_1,\ldots,v_m]}\bigl(\R[\mathcal K],\R\bigr)
\]
Since $\Sigma$ is a complete fan, $\sK=\sK_\Sigma$ is a sphere
triangulation, and therefore the face ring
$\R[\sK]=\R[v_1,\ldots,v_m]/\mathcal I_\sK$ is Cohen--Macaulay by
Corollary~\ref{spherecm}. The ideal $\mathcal J_{\Sigma}$ is
generated by a regular sequence, so we obtain by
Lemma~\ref{tortor},
\[
  \Tor_{\R[v_1,\ldots,v_m]}\bigl(\R[\mathcal K],\R\bigr)\cong
  \Tor_{\R[v_1,\ldots,v_m]/\mathcal J_\Sigma}\bigl
  (\R[\mathcal K]/\mathcal J_\Sigma,\R\bigr).
\]
Since $\R[v_1,\ldots,v_m]/\mathcal J_\Sigma$ is a polynomial ring
in $m-n$ variables, Lemma~\ref{koscom} implies that the
$\Tor$-algebra above is isomorphic to
\[
  H\bigl[\Lambda[u_1,\ldots,u_{m-n}]\otimes
  \R[\mathcal K]/\mathcal J_\Sigma,d\bigr]\cong
  H\bigl[\Lambda[u_1,\ldots,u_{m-n}]\otimes
  H^*(V_\Sigma),d\bigr]
\]
where the explicit form of the differential $d$ follows from the
definition of the Gale dual configuration
$\varGamma=(\gamma_1,\ldots,\gamma_m)$.
\end{proof}

There are two classical spectral sequences for the Dolbeault
cohomology. First, the \emph{Borel spectral
sequence}~\cite{bore66}\label{borelss} of a holomorphic bundle
$E\to B$ with a compact K\"ahler fibre~$F$, which has
$E_2=H_{\bar\partial}(B)\otimes H_{\bar\partial}(F)$ and converges
to $H_{\bar\partial}(E)$. Second, the \emph{Fr\"olicher spectral
sequence}~\cite[\S3.5]{gr-ha78}, whose $E_1$-term is the Dolbeault
cohomology of a complex manifold $M$ and which converges to the de
Rham cohomology of~$M$. Theorem~\ref{dolbzp} implies a collapse
result for these spectral sequences:

{\samepage
\begin{corollary}\
\begin{itemize}
\item[(a)]
The Borel spectral sequence of the holomorphic principal bundle
$\zk\to V_\Sigma$ collapses at the $E_3$-term, i.e.
$E_3=E_\infty$;

\item[(b)]
the Fr\"olicher spectral sequence of $\zk$ collapses at the
$E_2$-term.
\end{itemize}
\end{corollary}
}
\begin{proof}
To prove~(a) we just observe that the differential
algebra~\eqref{zkmult} is the $E_2$-term of the Borel spectral
sequence, and its cohomology is the $E_3$-term.

By comparing the Dolbeault and de Rham cohomology algebras of
$\zk$ given by Theorems~\ref{dolbzp} and~\ref{cohomzpred} we
observe that the elements $\eta_1,\ldots,\eta_\ell\in E_1^{0,1}$
cannot survive in the~$E_\infty$-term of the Fr\"olicher spectral
sequence. The only possible nontrivial differential on them is
$d_1\colon E_1^{0,1}\to E_1^{1,1}$. By Theorem~\ref{cohomzpred},
the cohomology algebra of $[E_1,d_1]$ is exactly the de Rham
cohomology of~$\zk$, proving~(b).
\end{proof}

Theorem~\ref{cohomzpred} can also be interpreted as a collapse
result for the Leray--Serre spectral sequence of the principal
$T^{m-n}$-bundle $\zk\to V_\Sigma$.

%There is also a bigrading in the ordinary cohomology of $\zk$,
%which is different from the bigrading in the Dolbeault cohomology.
%These two bigradings may be merged in the Dolbeault cohomology,
%providing a four-graded structure.

In order to proceed with calculation of Hodge numbers, we need the
following bounds for the dimension of $\Ker c$ in
Lemma~\ref{mumatrix}:

\begin{lemma}\label{cbounds}
Let $k$ be the number of zero vectors among $\mb a_1,\ldots,\mb
a_m$. Then
\[
  k-\ell\le\dim_\C\Ker\bigl(c\colon
  H^{1,0}_{\bar\partial}(T_\C^{\ell})\to
  H_{\bar\partial}^{1,1}(V_\Sigma)\bigr)\le{\textstyle\frac k2}.
\]
In particular, if $k\le1$ then $c$ is monomorphism.
\end{lemma}
\begin{proof}
Consider the commutative diagram
\[
\begin{CD}
  \Ann\Im\varPsi/\Ann\Ker A_\C @> i>> \C^m/\Ann\Ker A_\C
  @>p>> \C^{m-k}/\Ann\Ker A_\C\\
  @VV{\cong}V @VV\mathrm{Re}V @VV\mathrm{Re}V\\
  \R^{m-n}@=\R^{m-n} @>p'>> \R^{m-n-k}.
\end{CD}
\]
The composition $\Re\cdot i$ is an $\R$-linear isomorphism, as it
has the form $H^{1,0}_{\bar\partial}(T_\C^{\ell})\to
H^1(T_\C^{\ell},\C)\to H^1(T_\C^{\ell},\R)$, and any real-valued
function on the lattice $\Gamma$ defining the torus
$T_\C^{\ell}=\C^\ell/\Gamma$ is the real part of the restriction
to $\Gamma$ of a $\C$-linear function on~$\C^\ell$.

Since the diagram above is commutative, the kernel of $c=p\circ i$
has real dimension at most~$k$, which implies the upper bound on
its complex dimension. For the lower bound, $\dim_\C\Ker c\ge\dim
H^{1,0}_{\bar\partial}(T_\C^{\ell})-\dim
H_{\bar\partial}^{1,1}(V_\Sigma) =\ell-(2\ell-k)=k-\ell$.
\end{proof}

%Now we can obtain the following information about the Hodge
%numbers of~$\zk$.
%
\begin{theorem}\label{hodge}
Let $\zk$ be as in Theorem~\ref{dolbzp}, and let $k$ be the number
of zero vectors among $\mb a_1,\ldots,\mb a_m$. Then the Hodge
numbers $h^{p,q}=h^{p,q}(\zk)$ satisfy
\begin{itemize}
\item[(a)]
$\binom{k-\ell}p \le h^{p,0}\le\binom {[k/2]}p$ for $p\ge0$; in
particular, $h^{p,0}=0$ for $p>0$ if $k\le1$;
\item[(b)] $h^{0,q}=\binom\ell q$ for $q\ge0$;
\item[(c)] $h^{1,q}
  =(\ell-k)\binom\ell{q-1}+h^{1,0}\binom{\ell+1}q$ for $q\ge1$;
\item[(d)] $\frac{\ell(3\ell+1)}2-h_2(\sK)-\ell k+(\ell+1)h^{2,0}
 \le h^{2,1}\le\frac{\ell(3\ell+1)}2-\ell k+(\ell+1)h^{2,0}$.
\end{itemize}
%where $h_2(\sK)=\binom n2-(n-1)f_0(\sK)+f_1(\sK)$.
\end{theorem}
\begin{proof}
Let $A^{p,q}$ denote the bidegree $(p,q)$ component of the
differential algebra from Theorem~\ref{dolbzp}, and let
$Z^{p,q}\subset A^{p,q}$ denote the subspace of cocycles. Then
$d^{1,0}\colon A^{1,0}\to Z^{1,1}$ coincides with the map~$c$, and
the required bounds for $h^{1,0}=\Ker d^{1,0}$ are already
established in Lemma~\ref{cbounds}. Since $h^{p,0}=\dim\Ker
d^{p,0}$, and $\Ker d^{p,0}$ is the $p$th exterior power of the
space $\Ker d^{1,0}$, statement~(a) follows.

The differential is trivial on $A^{0,q}$, hence $h^{0,q}=\dim
A^{0,q}$, proving~(b).

The space $Z^{1,1}$ is spanned by the cocycles $v_i$ and
$\xi_i\eta_j$ with $\xi_i\in\Ker d^{1,0}$. Hence $\dim
Z^{1,1}=2\ell-k+h^{1,0}\ell$. Also, $\dim
d(A^{1,0})=\ell-h^{1,0}$, therefore,
$h^{1,1}=\ell-k+h^{1,0}(\ell+1)$. Similarly, $\dim
Z^{1,q}=(2\ell-k)\binom\ell{q-1}+h^{1,0}\binom\ell q$ (with basis
of $v_i\eta_{j_1}\cdots\eta_{j_{q-1}}$ and
$\xi_i\eta_{j_1}\cdots\eta_{j_q}$ where $\xi_i\in\Ker d^{1,0}$,
$j_1<\cdots<j_q$), and $d\colon A^{1,q-1}\to Z^{1,q}$ hits a
subspace of dimension $(\ell-h^{1,0})\binom\ell{q-1}$. This
proves~(c).

We have $A^{2,1}=U\oplus W$, where $U$ has basis of monomials
$\xi_iv_j$ and $W$ has basis of monomials $\xi_i\xi_j\eta_k$.
Therefore,
\begin{equation}\label{h21}
  h^{2,1}=\dim U-\dim dU+\dim W-\dim dW-\dim dA^{2,0}.
\end{equation}
Now $\dim U=\ell(2\ell-k)$, $0\le\dim dU\le h_2(\sK)$ (since
$dU\subset H_{\bar\partial}^{2,2}(V_\Sigma))$, $\dim W-\dim
dW=\dim\Ker d|_W=\ell h^{2,0}$, and $\dim
dA^{2,0}=\binom\ell2-h^{2,0}$. By substituting all this
into~\eqref{h21} we obtain the inequalities of~(d).
\end{proof}

\begin{remark}
At most one ghost vertex needs to be added to $\sK$ to make
$\dim\zk=m+n$ even. Since $h^{p,0}(\zk)=0$ when $k\le1$, the
manifold $\zk$ does not have holomorphic forms of any degree in
this case.

If $\zk$ is a torus (so that $\sK$ is empty), then $m=k=2\ell$,
and $h^{1,0}(\zk)=h^{0,1}(\zk)=\ell$. Otherwise
Theorem~\ref{hodge} implies that $h^{1,0}(\zk)<h^{0,1}(\zk)$, and
therefore $\zk$ is not K\"ahler.
%(in the polytopal case this was observed in~\cite[Th.~3]{meer00}).
\end{remark}

\begin{example}Let $\zk\cong S^1\times S^{2n+1}$ be a Hopf manifold of
Example~\ref{hopf}. Our rationality assumption is that $\mb
a_1\ldots,\mb a_{n+2}$ span an $n$-dimensional lattice $N$
in~$N_\R\cong\R^n$; in particular, the fan $\Sigma$ defined by the
proper subsets of $\mb a_1,\ldots,\mb a_{n+1}$ is rational. We
assume further that $\Sigma$ is regular (this is equivalent to the
condition $\mb a_1+\cdots+\mb a_{n+1}=\bf0$), so that $\Sigma$ is
a the normal fan of a Delzant $n$-dimensional
simplex~$\varDelta^n$. We have $V_\Sigma=\C P^n$,
and~\eqref{dolbtoric} describes its cohomology as the quotient of
$\C[v_1,\ldots,v_{n+2}]$ by the two ideals: $\mathcal I$ generated
by $v_1\cdots v_{n+1}$ and $v_{n+2}$, and $\mathcal J$ generated
by $v_1-v_{n+1},\ldots,v_{n}-v_{n+1}$. The differential algebra of
Theorem~\ref{dolbzp} is therefore given by
$\bigl[\Lambda[\xi,\eta]\otimes\C[t]/t^{n+1},d\bigr]$, %yura
with $dt=d\eta=0$ and $d\xi=t$ for a proper choice of~$t$. The
nontrivial cohomology classes are represented by the cocycles $1$,
$\eta$, $\xi t^n$ and $\xi\eta t^n$, which gives the following
nonzero Hodge numbers of~$\zk$:
$h^{0,0}=h^{0,1}=h^{n+1,n}=h^{n+1,n+1}=1$. Observe that the
Dolbeault cohomology and Hodge numbers do not depend on a choice
of complex structure (the map~$\varPsi$).
\end{example}

\begin{example}[Calabi--Eckmann manifold]\label{calabieck}
Let $\{\sK;\mb a_1,\ldots,\mb a_{n+2}\}$ be the data defining the
normal fan of the product $P=\varDelta^p\times\varDelta^q$ of two
Delzant simplices with $p+q=n$, $1\le p\le q\le n-1$. That is,
$\mb a_1,\ldots,\mb a_p,\mb a_{p+2},\ldots,\mb a_{n+1}$ is a basis
of lattice~$N$ and there are two relations $\mb a_1+\cdots+\mb
a_{p+1}=\bf0$ and $\mb a_{p+2}+\cdots+\mb a_{n+2}=\bf0$. The
corresponding toric variety $V_\Sigma$ is $\C P^p\times \C P^q$
and its cohomology ring is isomorphic to $\C[x,y]/(x^{p+1},
y^{q+1})$. Consider the map
\[
  \varPsi\colon\C\to\C^{n+2},\quad w\mapsto
  (w,\ldots,w,\alpha w,\ldots,\alpha w),
\]
where $\alpha\in\C\setminus\R$ and $\alpha w$ appears $q+1$ times.
The map $\varPsi$ satisfies the conditions of
Construction~\ref{psi}. The resulting complex structure on
$\zp\cong S^{2p+1}\times S^{2q+1}$ is that of a
\emph{Calabi--Eckmann manifold}. We denote complex manifolds
obtained in this way by $\mbox{\textit{C\!E}}(p,q)$ (the complex
structure depends on the choice of $\varPsi$, but we do not
reflect this in the notation). Each manifold
$\mbox{\textit{C\!E}}(p,q)$ is the total space of a holomorphic
principal bundle over $\C P^p\times \C P^q$ with fibre the complex
1-torus~$\C/(\Z\oplus\alpha\Z)$.

Theorem~\ref{dolbzp} and Lemma~\ref{mumatrix} provide the
following description of the Dolbeault cohomology of
$\mbox{\textit{C\!E}}(p,q)$:
\[
  H^{*,*}_{\bar\partial}\bigl(\mbox{\textit{C\!E}}(p,q)\bigr)\cong
  H\bigl[\Lambda[\xi,\eta]\otimes\C[x,y]/(x^{p+1},y^{q+1}),d\bigr],
\]
where $dx=dy=d\eta=0$ and $d\xi=x-y$ for an appropriate choice of
$x,y$. We therefore obtain
\begin{equation}\label{dolbce}
  H^{*,*}_{\bar\partial}\bigl(\mbox{\textit{C\!E}}(p,q)\bigr)\cong
  \Lambda[\omega,\eta]\otimes\C[x]/(x^{p+1}),
\end{equation}
where $\omega\in
H^{q+1,q}_{\bar\partial}\bigl(\mbox{\textit{C\!E}}(p,q)\bigr)$ is
the cohomology class of the cocycle
$\xi\frac{x^{q+1}-y^{q+1}}{x-y}$. This calculation is originally
due to~\cite[\S9]{bore66}. We note that the Dolbeault cohomology
of a Calabi--Eckmann manifold depends only on $p,q$ and does not
depend on the complex parameter~$\alpha$ (or the map~$\varPsi$).
\end{example}

\begin{example}
Now let
$P=\varDelta^1\times\varDelta^1\times\varDelta^2\times\varDelta^2$.
Then the moment-angle manifold $\zp$ has two structures of a
product of Calabi--Eckmann manifolds, namely,
$\mbox{\textit{C\!E}}(1,1)\times\mbox{\textit{C\!E}}(2,2)$ and
$\mbox{\textit{C\!E}}(1,2)\times\mbox{\textit{C\!E}}(1,2)$. Using
isomorphism~\eqref{dolbce} we observe that these two complex
manifolds have different Hodge numbers: $h^{2,1}=1$ in the first
case and $h^{2,1}=0$ in the second. This shows that the choice of
the map $\varPsi$ affects not only the complex structure of~$\zk$,
but also its Hodge numbers, unlike the previous examples of
complex tori, Hopf and Calabi--Eckmann manifolds. Certainly it is
not highly surprising from the complex-analytic point of view.
\end{example}

\section{Hamiltonian-minimal Lagrangian submanifolds}
In this last section we apply the accumulated knowledge on
topology of moment-angle manifolds in a somewhat different area,
Lagrangian geometry. Systems of real quadrics, which we used in
Sections~\ref{intquad} and~\ref{mampol} to define moment-angle
manifolds, also give rise to a family of Hamiltonian-minimal
Lagrangian submanifolds in a complex space or more general toric
varieties.

Hamiltonian minimality ($H$-minimality for short) for Lagrangian
submanifolds is a symplectic analogue of minimality in Riemannian
geometry. A Lagrangian immersion is called $H$-minimal if the
variations of its volume along all Hamiltonian vector fields are
zero. This notion was introduced in the work of
Y.-G.~Oh~\cite{oh93} in connection with the celebrated
\emph{Arnold conjecture} on the number of fixed points of a
Hamiltonian symplectomorphism. The simplest example of an
$H$-minimal Lagrangian submanifold is the coordinate
torus~\cite{oh93} $S^1_{r_1}\times\dots \times S^1_{r_m}\subset
{\mathbb C}^m$, where $S^1_{r_k}$ denotes the circle of radius
$r_k>0$ in the $k$th coordinate subspace of~$\C^m$. More examples
of $H$-minimal Lagrangian submanifolds in a complex space were
constructed in the
works~\cite{ca-ur98},~\cite{he-ro02},~\cite{an-ca11}, among
others.

In~\cite{miro04} Mironov suggested a general construction of
$H$-minimal Lagrangian immersions $N\looparrowright\C^m$ from
intersections of real quadrics. These systems of quadrics are the
same as those we used to define moment-angle manifolds, and
therefore one can apply toric methods for analysing the
topological structure of~$N$. In~\cite{mi-pa13f} an effective
criterion was obtained for $N\looparrowright\C^m$ to be an
embedding: the polytope corresponding to the intersection of
quadrics must be Delzant. As a consequence, any Delzant polytope
gives rise to an $H$-minimal Lagrangian submanifold
$N\subset\C^m$. As in the case of moment-angle manifolds, the
topology of $N$ is quite complicated even for low-dimensional
polytopes: for example, a Delzant 5-gon gives rise to a manifold
$N$ which is the total space of a bundle over a 3-torus with fibre
a surface of genus~5. Furthermore, by combining Mironov's
construction with symplectic reduction, a new family of
$H$-minimal Lagrangian submanifolds in of toric varieties was
defined in~\cite{mi-pa13u}. This family includes many previously
constructed explicit examples in~$\C^m$ and~$\C P^{m-1}$.

\subsection*{Preliminaries}
Let $(M,\omega)$ be a symplectic manifold of dimension $2n$. An
immersion $i\colon N\looparrowright M$ of an $n$-dimensional
manifold $N$ is called \emph{Lagrangian} if $i^*(\omega)=0$. If
$i$ is an embedding, then $i(N)$ is a \emph{Lagrangian
submanifold} of~$M$. A vector field $X$ on $M$ is
\emph{Hamiltonian}\label{lagrsubm} if the 1-form
$\omega(X,\,\cdot\,)$ is exact.

Now assume that $M$ is K\"ahler, so that it has compatible
Riemannian metric and symplectic structure. A Lagrangian immersion
$i\colon N\looparrowright M$ is called \emph{Hamiltonian minimal}
(\emph{$H$-minimal})\label{hminimalsubd} if the variations of the
volume of $i(N)$ along all Hamiltonian vector fields with compact
support are zero, that is,
\[
  \frac d{dt}\mathop{\mathrm{vol}}\bigl(i_t(N)\bigr)\big|_{t=0}=0,
\]
where $i_t(N)$ is a deformation of $i(N)$ along a Hamiltonian
vector field, $i_0(N)=i(N)$, and $\mathop{\mathrm{vol}}(i_t(N))$
is the volume of the deformed part of $i_t(N)$. An immersion $i$
is \emph{minimal} if the variations of the volume of $i(N)$ along
\emph{all} vector fields are zero.

Our basic example is $M=\C^m$ with the Hermitian metric
$2\sum_{k=1}^m d\overline{z}_k\otimes dz_k$. Its imaginary part is
the symplectic form of Example~\ref{simcm}. In the end we consider
a more general case when $M$ is a toric manifold.

\subsection*{The construction}
We consider an intersection of quadrics similar to~\eqref{zgamma},
but in the real space:
\begin{equation}\label{rgamma}
  \mathcal R=\Bigl\{\mb u=(u_1,\ldots,u_m)\in\R^m\colon
  \sum_{k=1}^m\gamma_{jk}u_k^2=\delta_j,\quad\text{for }
  1\le j\le m-n\Bigr\}.
\end{equation}

We assume the nondegeneracy and rationality conditions on the
coefficient vectors
$\gamma_i=(\gamma_{1i},\ldots,\gamma_{m-n,i})^t\in\R^{m-n}$,
$i=1,\ldots,m$:
\begin{itemize}
\item[(a)] $\delta\in
\R_\ge\langle\gamma_1,\ldots,\gamma_m\rangle$;

\item[(b)] if $\delta\in\R_\ge\langle
\gamma_{i_1},\ldots\gamma_{i_k}\rangle$, then $k\ge m-n$;

\item[(c)]
the vectors $\gamma_1,\ldots,\gamma_m$ generate a lattice
$L\cong\Z^{m-n}$ in~$\R^{m-n}$.
\end{itemize}

These conditions guarantee that $\mathcal R$ is a smooth
$n$-dimensional submanifold in $\R^m$ (by the argument of
Proposition~\ref{zgsmooth}) and that
\[
  T_\varGamma=
  \bigl\{\bigr(e^{2\pi i\langle\gamma_1,\varphi\rangle},
  \ldots,e^{2\pi i\langle\gamma_m,\varphi\rangle}\bigl)
  \in\T^m\bigr\}
\]
is an $(m-n)$-dimensional torus subgroup in~$\T^m$. We identify
the torus $T_\varGamma$ with $\R^{m-n}/\displaystyle L^*$ and
represent its elements by $\varphi\in\R^{m-n}$. We also define
\[
  D_\varGamma=\frac12 L^*/L^*\cong(\Z_2)^{m-n}.
\]
Note that $D_\varGamma$ embeds canonically as a subgroup in
$T_\varGamma$.

Now we view the intersection $\mathcal R$ as a subset in the
intersection $\mathcal Z$ or Hermitian quadrics given
by~\eqref{zgamma}, or as a subset in the whole space $\C^m$. Then
we `spread' $\mathcal R$ by the action of $T_\varGamma$, that is,
consider the set of $T_\varGamma$-orbits through $\mathcal R$.
More precisely, we consider the map
\begin{align*}
  j\colon\mathcal R\times T_\varGamma &\longrightarrow \C^m,\\
  (\mb u,\varphi) &\mapsto \mb u\cdot\varphi=\bigl(u_1e^{2\pi
i\langle\gamma_1,\varphi\rangle},\ldots,u_me^{2\pi
i\langle\gamma_m,\varphi\rangle}\bigr)
\end{align*}
and observe that $j(\mathcal R\times T_\varGamma)\subset\mathcal
Z$. We let $D_\varGamma$ act on $\mathcal R_\varGamma\times
T_\varGamma$ diagonally; this action is free, since it is free on
the second factor. The quotient
\[
  N=\mathcal R\times_{D_\varGamma} T_\varGamma
\]
is an $m$-dimensional manifold.

For any $\mb u=(u_1,\ldots,u_m)\in\mathcal R$, we have the
sublattice
\[
  L_{\mb u}=\Z\langle\gamma_k\colon u_k\ne0\rangle\subset
  L=\Z\langle\gamma_1,\ldots,\gamma_m\rangle.
\]
The set of $T_\varGamma$-orbits through $\mathcal R$ is an
immersion of~$N$:

\begin{lemma}\label{immer}\
\begin{itemize}
\item[(a)] The map $j\colon\mathcal R\times
T_\varGamma \to\C^m$ induces an immersion $i\colon
N\looparrowright\C^m$.

\item[(b)] The immersion $i$ is an embedding if and only
if $L_{\mb u}=L$ for any $\mb u\in\mathcal R$.
\end{itemize}
\end{lemma}
\begin{proof}
Take $\mb u\in\mathcal R$, $\varphi\in T_\varGamma$ and $g\in
D_\varGamma$. We have $\mb u\cdot g\in\mathcal R$, and $j(\mb
u\cdot g,g\varphi)=\mb u\cdot g^2\varphi=\mb u\cdot\varphi=j(\mb
u,\varphi)$. Hence, the map $j$ is constant on
$D_\varGamma$-orbits, and therefore induces a map of the quotient
$N=(\mathcal R\times T_\varGamma)/D_\varGamma$, which we denote
by~$i$.

%Since $D_\varGamma$ is a finite group, $i$ is an immersion if and
%only if $j$ is an immersion.
Assume that $j(\mb u,\varphi)=j(\mb u',\varphi')$. Then $L_{\mb
u}=L_{\mb u'}$ and
\begin{equation}\label{uu'}
  u_ke^{2\pi i\langle\gamma_k,\varphi\rangle}=u'_ke^{2\pi
  i\langle\gamma_k,\varphi'\rangle}\quad
  \text{for }k=1,\ldots,m.
\end{equation}
Since both $u_k$ and $u'_k$ are real, this implies that $e^{2\pi
i\langle\gamma_k,\varphi-\varphi'\rangle}=\pm1$ whenever
$u_k\ne0$, or, equivalently,
$\varphi-\varphi'\in\frac12\displaystyle L^*_{\mb u}/L^*$. In
other words,~\eqref{uu'} implies that $\mb u'=\mb u\cdot g$ and
$\varphi'=g\varphi$ for some $g\in\frac12{\displaystyle L^*_{\mb
u}/L^*}$. The latter is a finite group by Lemma~\ref{afree}; hence
the preimage of any point of $\C^m$ under $j$ consists of a finite
number of points. If $L_{\mb u}=L$, then $\frac12{\displaystyle
L^*_{\mb u}/L^*}=\frac12\displaystyle L^*/L^*=D_\varGamma$; hence
$(\mb u,\varphi)$ and $(\mb u',\varphi')$ represent the same point
in~$N$. Statement~(b) follows; to prove~(a), it remains to observe
that we have $L_{\mb u}=L$ for generic $\mb u$ (with all
coordinates nonzero).
\end{proof}

\begin{theorem}[{\cite[Theorem~1]{miro04}}]\label{hmin}
The immersion $i\colon N\looparrowright\C^m$ is $H$-minimal
Lagrangian. Moreover, if $\sum_{k=1}^m\gamma_k=0$, then $i$ is a
minimal Lagrangian immersion.
\end{theorem}
\begin{proof}
We only prove that $i$ is a Lagrangian immersion here. Let
\[
  (\mb x,\varphi)\mapsto
  \mb z(\mb x,\varphi)=\Bigl(u_1(\mb x)e^{2\pi
  i\langle\gamma_1,\varphi\rangle},\ldots,u_m(\mb x)e^{2\pi
  i\langle\gamma_m,\varphi\rangle}\Bigr)
\]
be a local coordinate system on $N=\mathcal R\times_{D_\varGamma}
T_\varGamma$, where $\mb x=(x_1,\ldots,x_n)\in\R^n$ and
$\varphi=(\varphi_1,\ldots,\varphi_{m-n})\in\R^{m-n}$. Let
$\langle\xi,\eta\rangle_\C=\sum_{i=1}^m\overline\xi_i\eta_i
=\langle\xi,\eta\rangle+i\omega(\xi,\eta)$ be the Hermitian scalar
product of $\xi,\eta\in\C^m$. Then
\[
  \Bigl\langle\frac{\partial\mb z}{\partial x_k},
  \frac{\partial\mb z}{\partial\varphi_j}\Bigr\rangle_\C=
  2\pi i\Bigr(\gamma_{j1}u_1\frac{\partial u_1}{\partial x_k}+\cdots+
  \gamma_{jm}u_m\frac{\partial u_m}{\partial x_k}\Bigl)=0
\]
where the second identity follows by differentiating the quadrics
equations~\eqref{rgamma}. Also, $\bigl\langle\frac{\partial\mb
z}{\partial x_k},\frac{\partial\mb z}{\partial
x_j}\bigr\rangle_\C\in\R$ and $\bigl\langle\frac{\partial\mb
z}{\partial\varphi_k},\frac{\partial\mb
z}{\partial\varphi_j}\bigr\rangle_\C\in\R$. It follows that
\[
  \omega\Bigl(\frac{\partial\mb z}{\partial x_k},
  \frac{\partial\mb z}{\partial\varphi_j}\Bigr)=
  \omega\Bigl(\frac{\partial\mb z}{\partial x_k},
  \frac{\partial\mb z}{\partial x_j}\Bigr)=
  \omega\Bigl(\frac{\partial\mb z}{\partial \varphi_k},
  \frac{\partial\mb z}{\partial\varphi_j}\Bigr)=0,
\]
i.e. the restriction of the symplectic form to the tangent space
of~$N$ is zero.
\end{proof}

\begin{remark}
The identity $\sum_{k=1}^m\gamma_k=0$ can not hold for a compact
$\mathcal R$ (or~$N$).
\end{remark}

We recall from Theorem~\ref{polquad} that a nonsingular
intersection of quadrics~\eqref{zgamma} or~\eqref{rgamma} defines
a simple polyhedron~\eqref{ptope}, and $\mathcal Z$ is identified
with the moment-angle manifold~$\zp$. Now we can summarise the
results of the previous sections in the following criterion for
$i\colon N\to\C^m$ to be an embedding:

\begin{theorem}\label{nembed}
Let $\mathcal Z$ and $\mathcal R$ be the intersections of
Hermitian and real quadrics defined by~\eqref{zgamma}
and~\eqref{rgamma}, respectively, satisfying
conditions~{\rm(a)--\,(c)} above. Let $P$ be the associated simple
polyhedron, and $N=\mathcal R\times_{D_\varGamma} T_\varGamma$.
The following conditions are equivalent:
\begin{itemize}
\item[(a)] $i\colon N\to\C^m$ is an embedding of an
$H$-minimal Lagrangian submanifold;

\item[(b)] $L_{\mb u}=L$ for any $\mb u\in\mathcal R$;

\item[(c)] the torus $T_\varGamma$ acts freely on the moment-angle manifold $\mathcal
Z=\zp$;

\item[(d)] $P$ is a Delzant polyhedron.
\end{itemize}
\end{theorem}
\begin{proof}
Equivalence (a)$\,\Leftrightarrow\,$(b) follows from
Lemma~\ref{immer} and Theorem~\ref{hmin}. Equivalence
(b)$\,\Leftrightarrow\,$(c) is Lemma~\ref{afree}, and
(c)$\,\Leftrightarrow\,$(d) is Theorem~\ref{propmmap}~(c).
\end{proof}

\subsection*{Topology of Lagrangian submanifolds~$N\subset\C^m$}
We start by reviewing three simple properties linking the
topological structure of $N$ to that of the intersections of
quadrics $\mathcal Z$ and~$\mathcal R$.

\begin{proposition}\label{nprop}\
\begin{itemize}
\item[(a)] The immersion of $N$ in $\C^m$ factors as
$N\looparrowright \mathcal Z\hookrightarrow\C^m$;
\item[(b)] $N$ is the total space of a bundle over the torus
$T^{m-n}$ with fibre $\mathcal R$;
\item[(c)] if $N\to\C^m$ is an embedding, then $N$ is the total space of a principal
$T^{m-n}$-bundle over the $n$-dimensional manifold $\mathcal
R/D_\varGamma$.
\end{itemize}
\end{proposition}
\begin{proof}
Statement (a) is clear. Since $D_\varGamma$ acts freely on
$T_\varGamma$, the projection $N=\mathcal
R\times_{D_\varGamma}T_\varGamma\to T_\varGamma/D_\varGamma$ onto
the second factor is a fibre bundle with fibre~$\mathcal R$.
Then~(b) follows from the fact that $T_\varGamma/D_\varGamma\cong
T^{m-n}$.

If $N\to\C^m$ is an embedding, then $T_\varGamma$ acts freely
on~$\mathcal Z$ by Theorem~\ref{nembed} and the action of
$D_\varGamma$ on $\mathcal R$ is also free. Therefore, the
projection $N=\mathcal R\times_{D_\varGamma}T_\varGamma\to
\mathcal R/D_\varGamma$ onto the first factor is a principal
$T_\varGamma$-bundle, which proves~(c).
\end{proof}

\begin{remark}\label{smalcove}
The quotient $\mathcal R/D_\varGamma$ is a \emph{real toric
variety}, or a \emph{small cover}, over the corresponding
polytope~$P$, see~\cite{da-ja91} and~\cite{bu-pa02}.
\end{remark}

\begin{example}[one quadric]\label{1quad}
Suppose that $\mathcal R$ is given by a single equation
\begin{equation}\label{1q}
  \gamma_1u_1^2+\cdots+\gamma_mu_m^2=\delta
\end{equation}
in $\R^m$. We assume that $\mathcal R$ is compact, so that
$\gamma_i$ and $\delta$ are positive reals, $\mathcal R\cong
S^{m-1}$, and the corresponding polytope $P$ is an
$n$-simplex~$\varDelta^n$. Then $N\cong S^{m-1}\times_{\Z_2}S^1$,
where the generator of $\Z_2$ acts by the standard free involution
on $S^1$ and by a certain involution $\tau$ on~$S^{m-1}$. The
topological type of $N$ depends on~$\tau$. Namely,
\[
  N\cong\begin{cases}S^{m-1}\times S^1&\text{if $\tau$ preserves the orientation of }S^{m-1},\\
  \mathcal K^{m}&\text{if $\tau$ reverses the orientation of }S^{m-1},\end{cases}
\]
where $\mathcal K^m$ is known as the \emph{$m$-dimensional Klein
bottle}\label{Kleinbot}.

\begin{proposition}\label{1qemb}
Let $m-n=1$ (one quadric). We obtain an $H$-minimal Lagrangian
embedding of $N\cong S^{m-1}\times_{\Z_2}S^1$ in $\C^m$ if and
only if $\gamma_1=\cdots=\gamma_m$ in~\eqref{1q}. In this case,
the topological type of $N=N(m)$ depends only on the parity of~$m$
and is given by
\begin{align*}\label{m=1}
  N(m)&\cong S^{m-1}\times S^1&&\text{if $m$ is even},\\
  N(m)&\cong\mathcal K^{m}&&\text{if $m$ is odd}.
\end{align*}
\end{proposition}
\begin{proof}
Since there exists $\mb u\in\mathcal R$ with only one nonzero
coordinate, Theorem~\ref{nembed} implies that $N$ embeds in $\C^m$
if only if $\gamma_i$ generates the same lattice as the whole set
$\gamma_1,\ldots,\gamma_m$ for each~$i$. Therefore,
$\gamma_1=\cdots=\gamma_m$. In this case $D_\varGamma\cong\Z_2$
acts by the standard antipodal involution on $S^{m-1}$, which
preserves orientation if $m$ is even and reverses orientation
otherwise.
\end{proof}

Both examples of $H$-minimal Lagrangian embeddings given by
Proposition~\ref{1qemb} are well known. The Klein bottle $\mathcal
K^m$ with even $m$ does not admit Lagrangian embeddings in~$\C^m$
(see~\cite{nemi09} and~\cite{shev09}).
\end{example}

\begin{example}[two quadrics]
In the case $m-n=2$, the topology of $\mathcal R$ and $N$ can be
described completely by analysing the action of the two commuting
involutions on the intersection of quadrics. We consider the
compact case here.

Using Proposition~\ref{propcf}, we write $\mathcal R$ in the form
\begin{equation}\label{2q}
\begin{aligned}
  \gamma_{11}u_1^2+\cdots+\gamma_{1m}u_m^2&=c,\\
  \gamma_{21}u_1^2+\cdots+\gamma_{2m}u_m^2&=0,
\end{aligned}
\end{equation}
where $c>0$ and $\gamma_{1i}>0$ for all $i$.

\begin{proposition}
There is a number $p$, \ $0<p<m$, such that $\gamma_{2i}>0$ for
$i=1,\ldots,p$ and $\gamma_{2i}<0$ for $i=p+1,\ldots,m$
in~\eqref{2q}, possibly after reordering the coordinates
$u_1,\ldots,u_m$. The corresponding manifold $\mathcal R=\mathcal
R(p,q)$, where $q=m-p$, is diffeomorphic to $S^{p-1}\times
S^{q-1}$. Its associated polytope $P$ either coincides with
$\varDelta^{m-2}$ (if one of the inequalities in~\eqref{ptope} is
redundant) or is combinatorially equivalent to the product
$\varDelta^{p-1}\times\varDelta^{q-1}$ (if there are no redundant
inequalities).
\end{proposition}
\begin{proof}
We observe that $\gamma_{2i}\ne0$ for all~$i$ in~\eqref{2q}, as
$\gamma_{2i}=0$ implies that the vector
$\delta=\Bigl(\!\!\begin{array}{l}c\\0\end{array}\!\!\Bigr)$ is in
the cone generated by one $\gamma_i$, which contradicts
Proposition~\ref{zgsmooth}~(b). By reordering the coordinates, we
can achieve that the first $p$ of $\gamma_{2i}$ are positive and
the rest are negative. Then $1<p<m$, because otherwise~\eqref{2q}
is empty. Now,~\eqref{2q} is the intersection of the cone over the
product of two ellipsoids of dimensions $p-1$ and $q-1$ (given by
the second quadric) with an $(m-1)$-dimensional ellipsoid (given
by the first quadric). Therefore, $\mathcal R(p,q)\cong
S^{p-1}\times S^{p-1}$. The statement about the polytope follows
from the combinatorial fact that a simple $n$-polytope with up to
$n+2$ facets is combinatorially equivalent to a product of
simplices;
%(see, e.g.~\cite[Ex.~I.8]{pano10});
the case of one redundant inequality corresponds to $p=1$ or
$q=1$.
\end{proof}

An element $\varphi\in
D_\varGamma=\frac12L^*/L^*\cong\Z_2\times\Z_2$ acts on $\mathcal
R(p,q)$ by
\[
  (u_1,\ldots,u_m)\mapsto
  (\varepsilon_1(\varphi)u_1,\ldots,\varepsilon_m(\varphi)u_m),
\]
where $\varepsilon_k(\varphi)=e^{2\pi
i\langle\gamma_k,\varphi\rangle}=\pm1$ for $1\le k\le m$.

\begin{lemma}\label{free1}
Suppose that $D_\varGamma$ acts on $\mathcal R(p,q)$ freely and
$\varepsilon_i(\varphi)=1$ for some~$i$, $1\le i\le p$, and
$\varphi\in D_\varGamma$. Then $\varepsilon_l(\varphi)=-1$ for all
$l$ with $p+1\le l\le m$.
\end{lemma}
\begin{proof}
Assume the opposite, that is, $\varepsilon_i(\varphi)=1$ for some
$1\le i\le p$ and $\varepsilon_j(\varphi)=1$ for some $p+1\le j\le
m$. Then $\gamma_{2i}>0$ and $\gamma_{2j}<0$ in~\eqref{2q}, so we
can choose $\mb u\in\mathcal R(p,q)$ whose only nonzero
coordinates are $u_i$ and~$u_j$. The element $\varphi\in
D_\varGamma$ fixes this $\mb u$, leading to contradiction.
\end{proof}

\begin{lemma}\label{free2}
Suppose $D_\varGamma$ acts on $\mathcal R(p,q)$ freely. Then there
exist two generating involutions $\varphi_1,\varphi_2\in
D_\varGamma\cong\Z_2\times\Z_2$ whose action on $\mathcal R(p,q)$
is described by either~{\rm(a)} or~{\rm(b)} below, possibly after
reordering the coordinates:
\begin{itemize}
\item[(a)]
$\begin{aligned}
  \varphi_1\colon(u_1,\ldots,u_m)&\mapsto
  (u_1,\ldots,u_k,-u_{k+1},\ldots,-u_p,-u_{p+1},\ldots,-u_m),\\[-2pt]
  \varphi_2\colon(u_1,\ldots,u_m)&\mapsto
  (-u_1,\ldots,-u_k,u_{k+1},\ldots,u_p,-u_{p+1},\ldots,-u_m);
\end{aligned}$\\
\item[(b)]
$\begin{aligned}
  \varphi_1\colon(u_1,\ldots,u_m)&\mapsto
  (-u_1,\ldots,-u_p,u_{p+1},\ldots,u_{p+l},-u_{p+l+1},\ldots,-u_m),\\[-2pt]
  \varphi_2\colon(u_1,\ldots,u_m)&\mapsto
  (-u_1,\ldots,-u_p,-u_{p+1},\ldots,-u_{p+l},u_{p+l+1},\ldots,u_m);
\end{aligned}$
\end{itemize}
here $0\le k\le p$ and $0\le l\le q$.
\end{lemma}
\begin{proof}
By Lemma~\ref{free1}, for each of the three nonzero elements
$\varphi\in D_\varGamma$, we have either
$\varepsilon_i(\varphi)=-1$ for $1\le i\le p$ or
$\varepsilon_i(\varphi)=-1$ for $p+1\le i\le m$. Therefore, we may
choose two different nonzero elements $\varphi_1,\varphi_2\in
D_\varGamma$ such that either $\varepsilon_i(\varphi_j)=-1$ for
$j=1,2$ and $p+1\le i\le m$, or $\varepsilon_i(\varphi_j)=-1$ for
$j=1,2$ and $1\le i\le p$. This corresponds to the cases (a) and
(b) above, respectively. In the former case, after reordering the
coordinates, we may assume that $\varphi_1$ acts as in~(a). Then
$\varphi_2$ also acts as in~(a), since otherwise the composition
$\varphi_1\cdot \varphi_2$ cannot act freely by Lemma~\ref{free1}.
The second case is treated similarly.
\end{proof}

Each of the actions of $D_\varGamma$ described in
Lemma~\ref{free2} can be realised by a particular intersection of
quadrics~\eqref{2q}. For example, the system of quadrics
\begin{equation}\label{2qex}
\begin{aligned}
  2u_1^2+\cdots+2u_k^2+u_{k+1}^2+\cdots+u_p^2+
  u_{p+1}^2+\cdots+u_m^2&=3,\\
  u_1^2+\cdots+u_k^2+2u_{k+1}^2+\cdots+2u_p^2-
  u_{p+1}^2-\cdots-u_m^2&=0
\end{aligned}
\end{equation}
gives the first action of Lemma~\ref{free2}; the second action is
realised similarly. Note that the lattice $L$ corresponding
to~\eqref{2qex} is a sublattice of index 3 in~$\Z^2$. We can
rewrite~\eqref{2qex} as
\begin{equation}\label{2qex1}
\begin{aligned}
  u_1^2+\cdots+u_k^2&+u_{k+1}^2+\cdots+u_p^2&&=1,\\
  u_1^2+\cdots+u_k^2&&+u_{p+1}^2+\cdots+u_m^2&=2,
\end{aligned}
\end{equation}
in which case $L=\Z^2$. The action of the two involutions
$\psi_1,\psi_2\in D_\varGamma=\frac12\Z^2/\Z^2$ corresponding to
the standard basis vectors of $\frac12\Z^2$ is given by
\begin{equation}\label{2inv}
\begin{aligned}
  \psi_1\colon(u_1,\ldots,u_m)&\mapsto
  (-u_1,\ldots,-u_k,-u_{k+1},\ldots,-u_p,u_{p+1},\ldots,u_m),\\
  \psi_2\colon(u_1,\ldots,u_m)&\mapsto
  (-u_1,\ldots,-u_k,u_{k+1},\ldots,u_p,-u_{p+1},\ldots,-u_m).
\end{aligned}
\end{equation}

We denote the manifold $N$ corresponding to~\eqref{2qex1} by
$N_k(p,q)$. We have
\begin{equation}\label{nkpq}
  N_k(p,q)\cong\mathcal (S^{p-1}\times
  S^{q-1})\times_{\Z_2\times\Z_2}(S^1\times S^1),
\end{equation}
where the action of the two involutions on $S^{p-1}\times S^{q-1}$
is given by $\psi_1,\psi_2$ above. Note that $\psi_1$ acts
trivially on $S^{q-1}$ and acts antipodally on $S^{p-1}$.
Therefore,
\[
  N_k(p,q)\cong N(p)\times_{\Z_2}(S^{q-1}\times S^1),
\]
where $N(p)$ is the manifold from Proposition~\ref{1qemb}. If
$k=0$ then the second involution $\psi_2$ acts trivially on
$N(p)$, and $N_0(p,q)$ coincides with the product $N(p)\times
N(q)$ of the two manifolds from Example~\ref{1quad}. In general,
the projection
\[
  N_k(p,q)\to S^{q-1}\times_{\Z_2}S^1=N(q)
\]
describes $N_k(p,q)$ as the total space of a fibration over $N(q)$
with fibre~$N(p)$.

We summarise the above facts and observations in the following
topological classification result for compact $H$-minimal
Lagrangian submanifolds $N\subset\C^m$ obtained from intersections
of two quadrics.

\begin{theorem}\label{2qemb}
Let $N\to\C^m$ be the embedding of the $H$-minimal Lagrangian
submanifold corresponding to a compact intersection of two
quadrics. Then $N$ is diffeomorphic to some $N_k(p,q)$ given
by~\eqref{nkpq}, where $p+q=m$, $0<p<m$ and $0\le k\le p$.
Moreover, every such triple $(k,p,q)$ can be realised by~$N$.
\end{theorem}
\end{example}

In the case of up to two quadrics considered above, the topology
of $\mathcal R$ is relatively simple, and in order to analyse the
topology of $N$, one only needs to describe the action of
involutions on~$\mathcal R$. When the number of quadrics is more
than two, the topology of $\mathcal R$ becomes an issue as well.

\begin{example}[three quadrics]\label{3quad}
In the case $m-n=3$, the topology of compact manifolds $\mathcal
R$ and $\mathcal Z$ was fully described
in~\cite[Theorem~2]{lope89}. Each of these manifolds is
diffeomorphic to a product of three spheres or to a connected sum
of products of spheres with two spheres in each product.

Note that, for $m-n=3$, the manifolds $\mathcal R$ (or~$\mathcal
Z$) can be distinguished topologically by looking at the planar
Gale diagrams of the corresponding simple polytopes~$P$ (see
Section~\ref{galed}). This chimes with the classification of
simple $n$-polytopes with $n+3$ facets, well-known in
combinatorial geometry~\cite[\S6.5]{zieg95}.

The smallest polytope with $m-n=3$ is a pentagon. It has many
Delzant realisations, for instance,
\[
  P=\bigl\{(x_1,x_2)\in\R^2\colon x_1\ge0,\;x_2\ge0,\;-x_1+2\ge0,\;-x_2+2\ge0,\;
  -x_1-x_2+3\ge0\bigr\}.
\]
In this case, $\mathcal R$ is an oriented surface of genus~5 (see
Proposition~\ref{rpolygon}), and the moment-angle manifold
$\mathcal Z$ is diffeomorphic to a connected sum of 5 copies of
$S^3\times S^4$.

We therefore obtain an $H$-minimal Lagrangian submanifold
$N\subset\C^5$, which is the total space of a bundle over $T^3$
with fibre a surface of genus~5.
\end{example}

Now assume that the polytope $P$ associated with intersection of
quadrics~\eqref{rgamma} is a polygon (i.e., $n=2$). If there are
no redundant inequalities, then $P$ is an $m$-gon and $\mathcal R$
is an orientable surface $S_g$ of genus $g=1+2^{m-3}(m-4)$ by
Proposition~\ref{rpolygon}. If there are $k$ redundant
inequalities, then $P$ is an $(m-k)$-gon. In this case $\mathcal
R\cong\mathcal R'\times(S^0)^k$, where $\mathcal R'$ corresponds
to an $(m-k)$-gon without redundant inequalities. That is,
$\mathcal R$ is a disjoint union of $2^k$ surfaces of genus
$1+2^{m-k-3}(m-k-4)$.

The corresponding $H$-minimal submanifold $N\subset\C^m$ is the
total space of a bundle over $T^{m-2}$ with fibre~$S_g$. This is
an aspherical manifold\label{asphHm} for $m\ge4$.

\subsection*{Generalisation to toric manifolds}
Consider two sets of quadrics:
\begin{align*}
  \mathcal Z_{\varGamma}&=\Bigl\{\mb z\in\C^m\colon
  \sum\nolimits_{k=1}^m\gamma_k|z_k|^2=\mb c\Bigr\},\quad \gamma_k,\mb c\in\R^{m-n};\\
  \mathcal Z_{\varDelta}&=\Bigl\{\mb z\in\C^m\colon
  \sum\nolimits_{k=1}^m\delta_k|z_k|^2=\mb d\Bigr\},\quad \delta_k,\mb
  d\in\R^{m-\ell};
\end{align*}
such that $\mathcal Z_{\varGamma}$, $\mathcal Z_{\varDelta}$
\emph{and} $\mathcal Z_{\varGamma}\cap\mathcal Z_{\varDelta}$
satisfy the nondegeneracy and rationality conditions (a)--(c) from
the beginning of this section. Assume also that the polyhedra
corresponding to $\mathcal Z_{\varGamma}$, $\mathcal
Z_{\varDelta}$ and $\mathcal Z_{\varGamma}\cap\mathcal
Z_{\varDelta}$ are Delzant.

The idea is to use the first set of quadrics to produce a toric
manifold $V$ via symplectic reduction (as described in
Section~\ref{symred}), and then use the second set of quadrics to
define an $H$-minimal Lagrangian submanifold in~$V$.

\begin{construction}\label{Hminimaltoricd}
Define the real intersections of quadrics $\mathcal
R_{\varGamma}$, $\mathcal R_{\varDelta}$, the tori
$T_{\varGamma}\cong\T^{m-n}$, $T_{\varDelta}\cong\T^{m-\ell}$, and
the groups $D_{\varGamma}\cong\Z_2^{m-n}$,
$D_{\varDelta}\cong\Z_2^{m-\ell}$ as before.

We consider the toric variety $V$ obtained as the symplectic
quotient of $\C^m$ by the torus corresponding to the first set of
quadrics: $V=\mathcal Z_{\varGamma}/T_{\varGamma}$. It is a
K\"ahler manifold of real dimension~$2n$. The quotient $\mathcal
R_{\varGamma}/D_{\varGamma}$ is the set of real points of $V$ (the
fixed point set of the complex conjugation, or the real toric
manifold); it has dimension~$n$. Consider the subset of $\mathcal
R_{\varGamma}/D_{\varGamma}$ defined by the second set of
quadrics:
\[
  \mathcal S=(\mathcal R_{\varGamma}\cap\mathcal
  R_{\varDelta})/D_{\varGamma},
\]
we have $\dim\mathcal S=n+\ell-m$. Finally, define the
$n$-dimensional submanifold of $V$:
\[
  N=\mathcal S\times_{D_{\varDelta}}T_{\varDelta}.
\]
\end{construction}

\begin{theorem}
$N$ is an $H$-minimal Lagrangian submanifold in the toric
manifold~$V$.
\end{theorem}
\begin{proof}
Let $\widehat V$ be the symplectic quotient of $V$ by the torus
corresponding to the second set of quadrics, that is, $\widehat
V=(V\cap\mathcal Z_\varDelta)/T_\varDelta=(\mathcal
Z_{\varGamma}\cap\mathcal Z_{\varDelta})/(T_{\varGamma}\times
T_{\varDelta})$. It is a toric manifold of real dimension
$2(n+\ell-m)$. The submanifold of real points
\[
  \widehat N=N/T_{\varDelta}=(\mathcal R_{\varGamma}\cap\mathcal
  R_{\varDelta})/(D_{\varGamma}\times D_{\varDelta})\hookrightarrow(\mathcal
  Z_{\varGamma}\cap\mathcal Z_{\varDelta})/(T_{\varGamma}\times
  T_{\varDelta})=\widehat V
\]
is the fixed point set of the complex conjugation, hence it is a
totally geodesic submanifold. In particular, $\widehat N$ is a
minimal submanifold in~$\widehat V$. According
to~\cite[Corollary~2.7]{dong07}, $N$ is an $H$-minimal submanifold
in~$V$.
\end{proof}

\begin{example}\

1. If $m-\ell=0$, i.e. $\mathcal Z_{\varDelta}=\varnothing$, then
$V=\C^m$ and we obtain the original construction of $H$-minimal
Lagrangian submanifolds $N$ in~$\C^m$.

2. If $m-n=0$, i.e. $\mathcal Z_{\varGamma}=\varnothing$, then $N$
is set of real points of~$V$. It is minimal (totally geodesic).

3. If $m-\ell=1$, i.e. $\mathcal Z_{\varDelta}\cong S^{2m-1}$,
then we get $H$-minimal Lagrangian submanifolds in $V=\C P^{m-1}$.
This includes the families of projective examples constructed
in~\cite{miro03},~\cite{ma05} and~\cite{mi-zu08}.
\end{example}

%Дописать раздел про локально стандартные действия
%
%В связи с результатом Кустарёва этим также упомянуть
%H. Geiges. Chern numbers of almost complex manifolds.
%Proceedings of the AMS 129 (2001), no. 12, 3749--3752.
%Any system of numbers that can be realised as the system of Chern numbers o
%a stably complex manifold of dim 2n, n\ge 2, can also be realised
%in this way by a connected almost complex manifold.
%Problem: which cobordism classes are realised by almost complex q-t manifolds.
%
%Добавить обсуждение задач: когда имеется (эквивариантное)
%отображение квазиторических многообразий M_1\to M_2 степени 1?
%Существование такого отображения - ключевой шаг в проблеме Новикова
%классификации гладких односвязных гомотопически эквивалентных многообразий с точностью до диффеоморфизма.
%Связь с результатами Новикова и результатами о когомологической жёсткости.
%
%Доказательство теоремы локализации (Сян)

\chapter{Half-dimensional torus actions}\label{torus}
In this chapter we study several topological generalisations of
toric varieties. All of them are smooth manifolds with an action
of a compact torus, the so-called \emph{$T$-manifolds}. Most of
the results here concern the case when the dimension of the acting
torus is half the dimension of the manifold.

A compact nonsingular toric variety (or toric manifold) $V_\Sigma$
of real dimension~$2n$ has an action of $n$-dimensional torus
$T^n$ obtained by restricting the action of the algebraic
torus~$(\C^\times)^n$. The action of $T^n$ on $V_\Sigma$ is
\emph{locally standard}, that is, it locally looks like the
standard coordinatewise action of $ T^n$ on~$\C^n$ (more precise
definitions are given below). If the variety $V_\Sigma$ is
projective, then the quotient $V_\Sigma/T^n$ can be identified
with a simple convex $n$-polytope $P$ via the moment map. In
general, the quotient $V_\Sigma/T^n$ may be not a polytope, even
combinatorially, but it still has a face structure of a manifold
with corners, which is dual to the face structure of the
simplicial fan~$\Sigma$. In particular, the fixed points of the
$T^n$-action on $V_\Sigma$ correspond to the vertices of the
quotient~$V_\Sigma/T^n$. These basic properties of the torus
action can be taken as the starting point for several different
topological generalisations of toric manifolds.

A \emph{quasitoric manifold} is a $2n$-dimensional manifold $M$
with a locally standard action of $T^n$ such that the quotient
$M/T^n$ can be identified with a simple $n$-polytope~$P$. This
class of manifolds was introduced in the seminal
paper~\cite{da-ja91} of Davis and Januszkiewicz. They showed,
among other things, that the cohomology ring structure of a
quasitoric manifold is exactly the same as that of a toric
manifold (see Theorem~\ref{danjur}). Quasitoric manifolds have
been studied intensively since the second half of 1990s, and the
work~\cite{bu-pa99} summarised these early developments and
emphasised the role of moment-angle manifolds~$\zp$.

Around the same time, an alternative way to generalise toric
varieties was developed in the works of Masuda~\cite{masu99} and
Hattori--Masuda~\cite{ha-ma03}, which led to a wider class of
\emph{torus manifolds}. Along with the usual conditions on the
$T^n$-action such as smoothness and effectiveness, the crucial
point in the definition of a torus manifold is the non-emptiness
of the fixed point set. Torus manifolds also admit a combinatorial
treatment similar to that of classic toric varieties in terms of
fans and polytopes. Namely, torus manifolds may be described by
\emph{multi-fans}\label{multifan} and \emph{multi-polytopes}; a
multi-fan is a collection of cones parametrised by a simplicial
complex, where some cones may overlap unlike in the usual fan. The
cohomology ring of a torus manifold has a more complicated
structure than that of a (quasi)toric manifold; in particular, it
is no longer generated by two-dimensional classes. Face rings of
simplicial posets (as described in Section~\ref{frsimpo}) arise
naturally in this context.

Another interesting generalisation of toric manifolds was
suggested in the work~\cite{i-f-m13} of Ishida, Fukukawa and
Masuda under the name \emph{topological toric manifolds}. The idea
is to consider not only the actions of a compact torus, but rather
two commuting actions of $T^n$ and $\mathbb R_>^n$, which patch
together to a \emph{smooth} $(\mathbb C^\times)^n$-action on a
manifold. This smooth $(\mathbb C^\times)^n$-action is what
replaces the algebraic $(\mathbb C^\times)^n$-action on a toric
manifold. The resulting class of manifolds, when properly defined,
turns out to be both wide and tractable, and admits a
combinatorial description in terms of (generalised) fans in a way
similar to toric varieties. The topological characteristics of
topological toric manifolds, including their integral cohomology
rings and characteristic classes, also have much similarity with
those of toric manifolds.

In Section~\ref{locst} we discuss the notion of a locally standard
torus action and related combinatorics of orbits; it features in
many subsequent generalisations of toric manifolds.
Section~\ref{tos} is a brief account of topological properties of
the (compact) torus action on toric manifolds, these properties
are taken as the base for subsequent topological generalisations.
Sections~\ref{qtman}, \ref{torusman} and~\ref{ttman} describe the
classes of quasitoric manifolds, torus manifolds and topological
toric manifolds, respectively. Section~\ref{torusman} can be also
viewed as an account on topological properties of locally standard
$T^n$-manifolds, because the existence of a fixed point (required
in the definition of a torus manifold) often comes as a
consequence of other topological restrictions on the locally
standard $T^n$-manifold~$M$; for example, a fixed point
automatically exists when the odd-degree cohomology of $M$
vanishes. Similarly, a torus manifold $M$ with $H^{odd}(M)=0$ is
necessarily locally standard. In Section~\ref{reltm} we describe
the relationship between the different classes of half-dimensional
torus actions. In Section~\ref{bfman} we discuss an important
class of examples of projective toric manifolds obtained as spaces
of bounded flags in a complex space. These manifolds illustrate
nicely many previous constructions with toric and quasitoric
manifolds, they will also feature in the last chapter on toric
cobordism. Another class of examples, Bott towers, is the subject
of Section~\ref{bott}. The study of Bott towers has become an
important part of toric topology, and many interesting open
questions arise here. In the last Section~\ref{weightgraphs} we
explore connections with another active area, the theory of
GKM-manifolds and GKM-graphs, and also study blow-ups of
$T$-manifolds and their related combinatorial objects. As usual,
more specific introductory remarks are available at the beginning
of each section.

\section[Locally standard actions]{Locally standard actions and manifolds with corners}\label{locst}
We have collected here background material on locally standard
$\T^n$- and $\Z_2^n$-actions and combinatorial structures on their
orbit spaces. The latter include the notions of \emph{manifolds
with faces} and \emph{manifolds with corners}, which have been
studied in differential topology since 1960s in the works of
J\"anich~\cite{jaen66}, Bredon~\cite{bred72}, Davis~\cite{davi83},
Davis--Januszkiewicz~\cite{da-ja91}, Izmestiev~\cite{izme00},
among others.

As usual, we denote by $\Z_2$ (respectively, by $\mathbb S$ or
$\T$) the multiplicative group of real (respectively, complex)
numbers of absolute value one. In this section we denote by $G$
one of the groups $\Z_2$ or $\mathbb S=\T$, and denote by $\mathbb
F$ its ambient field $\R$ or $\C$ respectively. We refer to the
coordinatewise action of $G^n$ on $\mathbb F^n$ given by
\[
  (g_1,\ldots,g_n)\cdot(z_1,\ldots,z_n)=(g_1z_1,\ldots,g_nz_n)
\]
as the \emph{standard action} or \emph{standard representation}.

\begin{definition}\label{defls}
Let $M$ be a manifold with an action of~$G^n$. A \emph{standard
chart} on $M$ is a triple $(U,f,\psi)$, where $U\subset M$ is a
$G^n$-invariant open subset, $\psi$ is an automorphism of~$G^n$,
and $f$ is a $\psi$-equivariant homeomorphism $f\colon U\to W$
onto a $G^n$-invariant open subset $W\subset\mathbb F^n$. (Recall
from Appendix~\ref{gractions} that the latter means that $f(\mb
t\cdot y)=\psi(\mb t)f(y)$ for all $\mb t\in G^n$, \ $y\in U$.) A
$G^n$-action on $M$ is said to be \emph{locally standard} if $M$
has a standard atlas, i.e. if any point of $M$ belongs to a
standard chart.
\end{definition}

The dimension of a manifold with a locally standard $G^n$-action
is $n$ if $G=\Z_2$ and $2n$ if $G=\T$. The orbit space of the
standard $G^n$-action on $\mathbb F^n$ is the orthant
\[
  \R^n_\ge=\bigl\{\mb x=(x_1,\ldots,x_n)\in\R^n\colon x_i\ge0
  \quad\text{for }
  i=1,\ldots,n\bigr\}.
\]
Therefore, the orbit space of a locally standard action is locally
modelled by~$\R^n_\ge$. There are two slightly different ways to
formalise this property, depending on whether we work in the
topological or smooth category.

\begin{definition}\label{mawfa}
A \emph{manifold with faces} (of dimension $n$) is a topological
manifold $Q$ with boundary $\partial Q$ together with a covering
of $\partial Q$ by closed connected subsets $\{F_i\}_{i\in\mathcal
S}$, called \emph{facets}, satisfying the following properties:
\begin{itemize}
\item[(a)] any facet $F_i$ is an $(n-1)$-dimensional submanifold
(with boundary) of~$\partial Q$;

\item[(b)] for any finite subset $I\subset\sS$, the intersection
$\bigcap_{i\in I}F_i$ is either empty or a disjoint union of
submanifolds of codimension~$|I|$;
% (i.e. of dimension $n-|I|$);
in the latter case we refer to a connected component of the
intersection $\bigcap_{i\in I}F_i$ as a \emph{face} of~$Q$;

\item[(c)] for any point $q\in Q$, there exist an open
neighbourhood $U\ni q$ and a homeomorphism $\phi\colon U\to W$
onto an open subset $W\subset\R^n_\ge$ such that
\[
  \varphi^{-1}(W\cap\{x_j=0\})=U\cap F_i
\]
for some~$i\in\mathcal S$.
\end{itemize}
%A \emph{facet} (respectively, \emph{face}) of $Q$ is a connected
%component of a pre-facet (respectively, pre-face).
\end{definition}

Observe that the set $\mathcal S$ has to be countable, hence the
set of faces in a manifold with faces $Q$ is also countable. If
$Q$ is compact, then both these sets are finite.

The orthant $\R^n_\ge$ itself has the canonical structure of a
manifold with faces; each face has the form $\R_\ge^I$ for
$I\subset[n]$, where $\R_\ge^I=\{\mb x\in\R^n_\ge\colon
x_j=0\quad\text{for }j\notin I\}$. The \emph{codimension} $c(\mb
x)$ of point $\mb x\in\R^n_\ge$ is the number of zero coordinates
of~$\mb x$. A point of~$Q$ has \emph{codimension~$k$} if it
belongs to a face of codimension $k$ and does not belong to a face
of codimension~$k+1$.

\begin{definition}\label{mawco}
A \emph{manifold with corners} (of dimension~$n$) is a topological
manifold $Q$ with boundary together with an atlas $\{U_i,\phi_i\}$
consisting of homeomorphisms $\phi_i\colon U_i\to W_i$ onto open
subsets $W_i\subset\R^n_\ge$ such that
$\phi_i\phi^{-1}_j\colon\phi_j(U_i\cap U_j)\to\phi_i(U_i\cap U_j)$
is a diffeomorphism for all $i,j$. (A homeomorphism between open
subsets in $\R^n_\ge$ is called a \emph{diffeomorphism} if it can
be obtained by restriction of a diffeomorphism of open subsets in
$\R^n$.)

For any $q\in Q$, its codimension $c(q)$ is well defined. An
\emph{open face} of $Q$ of codimension $k$ is a connected
component of $c^{-1}(k)$. A \emph{closed face} (or simply
\emph{face})\label{facemwco} of $Q$ is the closure of an open
face. A \emph{facet} is a face of codimension~1.

A manifold with corners $Q$ is said to be
\emph{nice}\label{nicemwco} if the covering of $Q$ by its facets
satisfies condition (c) of Definition~\ref{mawfa} (conditions (a)
and (b) are satisfied automatically). Equivalently, a manifold
with corners is nice if and only if each of its faces of
codimension~$2$ is contained in exactly $2$ facets (an exercise).
\end{definition}

\begin{remark}
Our definitions of faces differ from those of~\cite{davi83}
and~\cite{izme00}; the reason is that we want our faces to be
connected.
\end{remark}

\begin{example}\label{examplemwc}\ \nopagebreak

1. The 2-disc with a single `corner point' on its boundary is a
manifold with corners which is not nice. All other examples in
this list will be nice.

2. A smooth manifold $Q$ with boundary is a manifold with corners,
whose facets are connected components of~$\partial Q$, and there
are no other faces.

3. A direct product of manifolds with corners is a manifold with
corners. In particular, a product of smooth manifolds with
boundary is a manifold with corners.

4. Let $P$ be a simple polytope. For each vertex $v\in P$ we
denote by $U_v$ the open subset in~$P$ obtained by removing all
faces of $P$ that do not contain~$v$. The subset $U_v$ is affinely
isomorphic to a neighbourhood of zero in~$\R^n_\ge$. Therefore $P$
is a compact manifold with corners, with atlas~$\{U_v\}$.
\end{example}

The orbit space $Q=M/G^n$ of a locally standard action is a
manifold with faces. As we shall see in Proposition~\ref{oslsa},
if the action is smooth, then $Q$ is a nice manifold with corners.

\subsection*{Exercises}
\begin{exercise}
A manifold with corners $Q$ is nice if an only if any face of
codimension two is contained in exactly two facets.
\end{exercise}

\begin{exercise}
Two simple polytopes are diffeomorphic as manifolds with corners
if and only if they are combinatorially equivalent
(see~\cite{davi13} for a more general statement).
\end{exercise}

\section{Toric manifolds and their quotients}\label{tos}
By way of motivation, here we take a closer look at the action of
the (compact) torus $T_N\cong T^n$ on a toric manifold $V_\Sigma$.
The topological properties of the quotient projection $\pi\colon
V_\Sigma\to V_\Sigma/T_N$ will be taken as the starting point for
subsequent topological generalisations of toric manifolds.

We first recall from Sections~\ref{symred} and~\ref{symred1} that
a projective (or Hamiltonian) toric manifold $V_P$ can be
identified with the symplectic quotient of $\C^m$ by an action of
$K\cong T^{m-n}$, i.e. with the quotient manifold $\zp/K$ where
$\zp$ is the moment-angle manifold corresponding to~$P$. Therefore
the quotient of $V_P$ by the action of the $n$-torus $T_N=\T^m/K$
coincides with the quotient of $\zp$ by~$\T^m$. Both quotients are
identified with the Delzant polytope~$P$; in fact, the moment map
$\mu_V\colon V_P\to P$ is the quotient projection (see
Proposition~\ref{tvmomentmap}). In the non-projective smooth case,
there is no moment map, and there is no canonical way to identify
the quotient $V_\Sigma/T_N$ with a convex polytope. However, there
is a face decomposition (stratification) of $V_\Sigma/T_N$
according to orbit types, and this face structure is very similar
to that of a simple polytope.

Following Davis--Januszkiewicz~\cite{da-ja91}, we can describe a
projective toric manifold $V_P$ as an identification space similar
to~\eqref{Cmids}.

As usual, for each $I\subset[m]$ we denote by $\T^I=\prod_{i\in
I}\T$ the corresponding coordinate subgroup in~$\T^m$. Given $\mb
x\in P$, set $I_{\mb x}=\{i\in[m]\colon{\mb x}\in F_i\}$ (the set
of facets containing~$\mb x$). We recall the map $A\colon\R^m\to
N_\R$, \ $\mb e_i\mapsto\mb a_i$, and its exponential
$\exp\Amap\colon\T^m\to T_N$. For each $\mb x\in P$ define the
subtorus $T_{\mb x}=(\exp\Amap)(\T^{I_{\mb x}})\subset T_N$. If
$\mb x$ is a vertex then $T_{\mb x}=T_N$, and if $\mb x$ is an
interior point of $P$ then $T_{\mb x}=\{1\}$.

\begin{proposition}
A projective toric manifold $V_P$ is $T_N$-equivariantly
homeomorphic to the quotient
\begin{equation}\label{identificationspace}
  P\times T_N/{\approx}\quad\text{where }
  (\mb x,t_1)\approx(\mb x,t_2)\:\text{ if }\:
  t_1^{-1}t_2\in T_{\mb x}.
\end{equation}
\end{proposition}
\begin{proof}
Using Proposition~\ref{zpids}, we obtain
\[
  V_P=\zp/K=\bigl(P\times(\T^m/K)\bigr)/{\sim\:}
  =\bigl(P\times(\exp\Amap)(\T^m)\bigr)/{\sim\:}
  =P\times T_N/{\approx}.\qedhere
\]
\end{proof}

The projection $\pi\colon V_P=P\times T_N/{\approx}\to P$ is the
quotient map for the $T_N$-action, and its fibre $\pi^{-1}(\mb
x)=T_{\mb x}$ is the stabiliser of the $T_N$-orbit corresponding
to~$\mb x$. The $T_N$-action on $V_P$ is therefore free over the
interior of the polytope, vertices of the polytope correspond to
fixed points, and points in the relative interior of a
codimension-$k$ face correspond to orbits with the same
$k$-dimensional stabiliser.

\begin{construction}\label{nsorb}
Let $\Sigma$ be a simplicial fan and $V_\Sigma$ the corresponding
toric variety. Consider the affine cover
$\{V_\sigma\colon\sigma\in\Sigma\}$ (see Construction~\ref{tvff}).
The quotient $(V_\sigma)_\ge=V_\sigma/T_N$ can be identified with
the set of semigroup homomorphisms from $S_\sigma$ to the
semigroup $\R_\ge$ of nonnegative real numbers:
\[
  V_\sigma/T_N=\Hom_{\mathrm{sg}}(S_\sigma,\R_\ge),
  \quad\sigma\in\Sigma
\]
(the details of this construction can be found
in~\cite[\S4.1]{fult93}).

If the fan $\Sigma$ is regular, then
$V_\sigma\cong\C^k\times(\C^\times)^{n-k}$ where $k=\dim\sigma$,
see Example~\ref{affns}. It follows easily that the $T_N$-action
on a nonsingular toric variety $V_\Sigma$ is locally standard, and
the cover $\{V_\sigma\colon\sigma\in\Sigma\}$ provides an atlas of
standard charts (see Section~\ref{locst}). Furthermore,
$(V_\sigma)_\ge\cong\R_\ge^k\times\R^{n-k}$ and the orbit space
$Q=V_\Sigma/T_N$ is a manifold with corners, with atlas
$\{(V_\sigma)_\ge\colon\sigma\in\Sigma\}$.

In the singular case the varieties $V_\sigma$ may be not
isomorphic to $\C^k\times(\C^\times)^{n-k}$ and the $T_N$-action
on $V_\Sigma$ may fail to be locally standard. Nevertheless, the
cover $\{(V_\sigma)_\ge\colon\sigma\in\Sigma\}$ defines a
structure of a manifold with faces on $Q=V_\Sigma/T_N$.

If $\Sigma$ is a complete simplicial fan, then the face
decomposition of the manifold with corners $Q=V_\Sigma/T_N$ is
Poincar\'e dual to the sphere triangulation defined by the
simplicial complex~$\mathcal K_\Sigma$. There is also a projection
$Q\times T_N\to V_\Sigma$ defining a homeomorphism $V_\Sigma\cong
Q\times T_N/{\approx}$, by analogy
with~\eqref{identificationspace}.
\end{construction}

%As it follows from the analysis in the beginning of this section,
%the quotient $V_P/T_N$ of a projective toric manifold $V_P$ is
%homeomorphic to $P$ as a manifold with corners. Furthermore, if
%$\pi\colon V_P\to P$ is the quotient projection, then
%$\pi^{-1}(U_v)\subset V_P$ is the affine subvariety
%$V_{\sigma(v)}$ where $\sigma(v)\in\Sigma_P$ is the maximal cone
%in the normal fan corresponding to the vertex $v\in P$.

We finish by summarising the observations of this section as
follows:

\begin{proposition}\label{2prop}\
\begin{itemize}
\item[(а)] The action of the torus $T_N$ on a nonsingular toric variety $V_\Sigma$
is locally standard.

\item[(b)] If $V_P$ is a projective toric manifold then the quotient $V_P/T_N$ is diffeomorphic to the simple polytope $P$
as a manifold with corners.
\end{itemize}
\end{proposition}

If a toric manifold $V_\Sigma$ is not projective, then the
manifold with corners $Q=V_\Sigma/T_N$ may be not homeomorphic to
a simple polytope (see Section~\ref{reltm}).

\section{Quasitoric manifolds}\label{qtman}
Quasitoric manifolds were introduced by Davis and Januszkiewicz as
a topological alternative to (nonsingular projective) toric
varieties. Originally, the term `toric manifolds' was used
in~\cite{da-ja91} to describe this class of manifolds, but later
is was replaced by `quasitoric', as `toric manifold' is often used
by algebraic and symplectic geometers as a synonym for
`non-singular complete toric variety'.

Any quasitoric manifold over $P$ can be obtained as a quotient of
the moment-angle manifold $\mathcal Z_P$ by a freely acting
subtorus. This can be viewed as a topological version of the
symplectic quotient construction of projective toric manifolds
(see Section~\ref{symred}). As a result we obtain a canonical
smooth structure on a quasitoric manifold, in which the torus
action is smooth.

Unlike toric manifolds, quasitoric manifolds are not complex
varieties in general, and they may even not admit an almost
complex structure. However, an equivariant stably complex
structure always exists on a quasitoric manifold $M$ and is
defined canonically by the underlying combinatorial data. These
structures will be used in Chapter~\ref{cobtoric} to define
quasitoric representatives in complex cobordism classes.

\subsection*{Definition and basic constructions}
\begin{definition}
\label{qtm} Let $P$ be a combinatorial simple polytope of
dimension~$n$. A \emph{quasitoric manifold} over $P$ is a smooth
$2n$-dimensional manifold $M$ with a smooth action of the torus
$T^n$ satisfying the two conditions:
\begin{enumerate}
  \item[(a)] the action is locally standard (see Definition~\ref{defls});
  \item[(b)] there is a continuous projection $\pi\colon M\to P$ whose fibres are $T^n$-orbits.
\end{enumerate}
Property (a) implies that the quotient $M/P$ is a manifold with
corners, and property (b) implies that the quotient is
homeomorphic, as a manifold with corners, to the simple
polytope~$P$.
\end{definition}

%\begin{remark}
%Smoothness of $M$ and of the $T$-action was not assumed in the
%original definition of~\cite{da-ja91}.
%
%The manifold $M$ is not required to be smooth in the definition
%above. It can be shown that any quasitoric manifold has a smooth
%structure in which the torus action is smooth (see~\cite{davi78}).
%An alternative construction providing a quasitoric manifold with a
%canonical smooth structure will be given below.
%\end{remark}
%
It follows from the definition that the projection $\pi\colon M\to
P$ maps a $k$-dimensional orbit of the $T^n$-action to a point in
the relative interior of a $k$-dimensional face of~$P$. %, for $k=0,\ldots,n$.
In particular, the action is free over the interior of the
polytope, while vertices of $P$ correspond to fixed points of the
torus action on~$M$.

\begin{proposition}\label{prqua}
A nonsingular projective toric variety $V_P$ is a quasitoric
manifold over~$P$.
\end{proposition}
\begin{proof}
This follows from Proposition~\ref{2prop}.
\end{proof}

\begin{example}
The complex projective space $\C P^n$ with the action of
$T^n\subset(\C^\times)^n$ described in Example~\ref{tvcpn} is a
quasitoric manifold over the simplex~$\varDelta^n$. The projection
$\pi\colon\C P^n\to\varDelta^n$ is given by
\[
  (z_0:z_1:\cdots:z_n) \mapsto
  \frac{1}{\sum_{i=0}^n |z_i|^2}(|z_1|^2,\ldots,|z_n|^2).
\]
\end{example}

\begin{remark}
For any 2- or 3-dimensional polytope $P$ there exists a quasitoric
manifold over~$P$. For $n\ge4$, there exist $n$-dimensional
polytopes which do not arise as quotients of quasitoric manifolds.
(See Exercises~\ref{23dqt} and~\ref{4dqt}.)
\end{remark}

Let $\mathcal F=\{F_1,\ldots,F_m\}$ be the set of facets of~$P$.
Consider the preimages
\[
  M_{j}=\pi^{-1}(F_j), \quad 1\le j\le m.
\]
Points in the relative interior of a facet $F_j$ correspond to
orbits with the same one-dimensional stabiliser subgroup, which we
denote by~$T_{F_j}$. It follows that $M_{j}$ is a connected
component of the fixed point set of the circle subgroup
$T_{F_j}\subset T^n$. This implies that $M_{j}$ is a
$T^n$-invariant submanifold of codimension 2 in~$M$, and $M_{j}$
is a quasitoric manifold over~$F_j$ with the action of the
quotient torus $T^n/T_{F_j}\cong T^{n-1}$.
Following~\cite{da-ja91}, we refer to $M_{j}$ as the
\emph{characteristic submanifold} corresponding to the $j$th face
$F_j\subset P$. The mapping
\begin{equation}
\label{charf}
  \lambda\colon F_j\mapsto T_{F_j},\quad 1\le j\le m,
\end{equation}
is called the \emph{characteristic function} of the quasitoric
manifold~$M$.

Now let $G$ be a codimension-$k$ face of~$P$. We can write it as
an intersection of $k$ facets: $G=F_{j_1}\cap\cdots\cap F_{j_k}$.
Then $M_G=\pi^{-1}(G)$ is a $T^n$-invariant submanifold of
codimension $2k$ in~$M$, and $M_G$ is fixed under each circle
subgroup $T(F_{j_p})$, $1\le p\le k$. By considering any vertex
$v\in G$ and using the local standardness of the $T^n$-action on
$M$ near $v$, we observe that the characteristic submanifolds
$M_{j_1},\ldots,M_{j_k}$ intersect transversely at the submanifold
$M_G$, and the map
\[
  T_{F_{j_1}}\times\cdots\times T_{F_{j_k}}\to T^n
\]
is a monomorphism onto the $k$-dimensional stabiliser of~$M_G$.
The mapping
\[
  G\mapsto\text{ the stabiliser of }M_G
\]
extends the characteristic function~\eqref{charf} to a map from
the face poset of $P$ to the poset of torus subgroups in~$T^n$.

\begin{definition}
\label{charpair} Let $P$ be a combinatorial $n$-dimensional simple
polytope and let $\lambda$ be a map from the set of facets of~$P$
to the set of circle subgroups of the torus~$T^n$. We refer to
$(P,\lambda)$ as a \emph{characteristic pair} if the map
$\lambda(F_{j_1})\times\cdots\times\lambda(F_{j_k})\to T^n$ is a
monomorphism whenever $F_{j_1}\cap\cdots\cap
F_{j_k}\ne\varnothing$.
\end{definition}

If $(P,\lambda)$ is a characteristic pair, then the map $\lambda$
extends to the face poset of~$P$, and we have a torus subgroup
$T_G=\lambda(G)\subset T^n$ for each face $G\subset P$.

As we shall see below, a quasitoric manifold can be reconstructed
from its characteristic pair $(P,\lambda)$ up to a weakly
$T$-equivariant homeomorphism.

\begin{construction}[canonical model $M(P,\lambda)$]
\label{der} Assume given a characteristic pair $(P,\lambda)$. For
any point $x\in P$, we denote by $G(x)$ the smallest face
containing~$x$. By analogy with~\eqref{identificationspace}, we
define the identification space
$$
  M(P,\lambda)=P\times T^n/{\sim}\quad\text{where }
  (x,t_1)\sim(x,t_2)\:\text{ if }\:
  t_1^{-1}t_2\in\lambda(G(x)).
$$
The free action of $T^n$ on $P\times  T^n$ descends to an action
on $P\times  T^n/{\sim}$. This action is free over the interior of
the polytope (since no identifications are made
over~$\mathop{\mathrm{int}} P$), and the fixed points correspond
to the vertices. The space $P\times T^n/{\sim}$ is covered by the
open subsets $U_v\times T^n/{\sim}$, indexed  by the vertices
$v\in P$ (see Example~\ref{examplemwc}.4), and each $U_v\times
T^n/{\sim}$ is equivariantly homeomorphic to $\C^n=\R^n_\ge\times
T^n/{\sim}$. This implies that the canonical model
$M(P,\varLambda)$ is a (topological) manifold with a locally
standard $T^n$-action and quotient~$P$. As we shall see from
Definition~\ref{cansmstqt}, $M(P,\varLambda)$ has a canonical
smooth structure, so it is a quasitoric manifold over~$P$.
\end{construction}

\begin{proposition}[{\cite[Lemma~1.4]{da-ja91}}]\label{equivar}
There exists a weakly $T^n$-equivariant homeomorphism
\[
  M(P,\lambda)=P\times T^n/{\sim}\to M
\]
covering the identity map on~$P$.
\end{proposition}
\begin{proof}
One first constructs a weakly $T^n$-equivariant map $f\colon
P\times T^n\to M$ such that $f$ maps $x\times T^n$ onto
$\pi^{-1}(x)$ for any point $x\in P$. Such a map $f$ induces a
weakly $T^n$-equivariant map $\widehat f\colon
M(P,\lambda)=P\times T^n/{\sim}\to M$ covering the identity
on~$P$. Furthermore, the map $\widehat f$ is one-to-one on
$T^n$-orbits, so it is a homeomorphism.

It remains to construct a map $f\colon P\times T^n\to M$. The
argument of Davis--Januszkiewicz which we present here actually
works for a more general class of locally standard
$T^n$-manifolds~$M$. There is the manifold with boundary
$\widetilde M$ obtained by consecutive blowing up the singular
strata of~$M$ consisting of non-principal $T^n$-orbits. The
$T^n$-action on $\widetilde M$ is free and $\widetilde M$ is
equivariantly diffeomorphic to the complement in $M$ of the union
of tubular neighbourhoods of the singular strata (the latter are
characteristic submanifolds~$M_j$ in our case). There is the
following canonical inductive procedure for constructing
$\widetilde M$ from~$M$. One begins by removing a minimal stratum
(a fixed point in our case) from $M$ and replacing it by the
sphere bundle of its normal bundle. One continues in this fashion,
blowing up minimal strata, until only top stratum is left. There
is the canonical projection from the union of sphere bundles to
the union of their base spaces (i.e. to the union of singular
strata of~$M$). Using the construction of~\cite[p.~344]{davi78}
one extends this projection to a map $\widetilde M\to M$ which is
the identity over top stratum. Now if $M$ is locally standard,
then the quotient $\widetilde M/T^n$ is canonically identified
with~$M/T^n$. In our particular case the latter is a simple
polytope $P$. Since $P$ is acyclic together with all its faces,
the resulting principal $T^n$-bundle
$\widetilde\pi\colon\widetilde M\to P$ is trivial. Therefore,
there is an equivariant diffeomorphism $\widetilde f\colon P\times
T^n\to\widetilde M$ inducing the identity on~$P$.
Composing~$\widetilde f$ with the collapse map $\widetilde M\to
M$, we obtain the map~$f$.
\end{proof}

\begin{remark}
Note that existence a map $f\colon P\times T^n\to M$ in the proof
above is equivalent to existence of a section $s\colon P\to M$ of
the quotient projection $\pi\colon M\to P$. Indeed, given a
section $s$, one defines $f(x,t)=t\cdot s(x)$ for $x\in P$ and
$t\in T^n$. Conversely, given a map~$f$, one defines a section by
$s(x)=f(x,1)$.

It may seem that constructing a section $s\colon P\to M$ is an
easier task than constructing a map~$f$, taking into account that
$P$ is contractible. However, this comes out to be subtle; it
would be interesting to find a more explicit way to construct a
section~$s$.
%
%A section $s\colon P\to M$ can be constructed easily using the
%atlas $\{U_v\}$ of $P$ described in Example~\ref{examplemwc}.4.
%Let $\widetilde U_v=\pi^{-1}(U_v)$. Then $\widetilde U_v\cong W$,
%where $W$ is a $T^n$-invariant open subset in~$\C^n$, and
%$\{\widetilde U_v\}$ is an atlas of~$M$. We pick any initial
%vertex $v$ and choose a section $U_v\to\widetilde U_v$ arbitrarily
%(e.g. by identifying $U_v$ with a subset of $W$ consisting of
%points with real coordinates). Then for any other vertex $v'$ the
%section $U_{v'}\to\widetilde U_{v'}$ is already defined on a `big'
%open subset $U_v\cap U_{v'}\subset U_{v'}$ whose closure is the
%whole~$U_{v'}$. Therefore, the section $U_v\cap
%U_{v'}\to\widetilde U_{v'}$ extends uniquely to
%$U_{v'}\to\widetilde U_{v'}$. By doing such an extension for any
%vertex we obtain a section $s\colon P\to M$.
\end{remark}

%\begin{remark} The above sketch of proof shows that the only condition
%on the quotient $M/T^n$ which guarantees the equivalence between
%$M$ and the canonical model is that the principal $T^n$-bundle
%over $M/T^n$ obtained after all blow-ups is trivial. Since any
%principal $T^n$-bundle over $M/T^n$ with $H^2(M/T^n;\Z)=0$ is
%trivial, the proposition above can be extended to a wider class of
%torus actions. We shall use this observation in
%Section~\ref{torusman}.
%\end{remark}

\begin{definition}[equivalences]\label{psieq}
Quasitoric manifolds $M_1$ and $M_2$ over the same polytope~$P$
are said to be \emph{equivalent over~$P$} if there exists a weak
$T^n$-equivariant homeomorphism $f\colon M_1\to M_2$ covering the
identity map on~$P$.

Two characteristic pairs $(P,\lambda_1)$ and $(P,\lambda_2)$ are
said to be \emph{equivalent} if there exists an automorphism
$\psi\colon  T^n\to T^n$ such that $\lambda_2=\psi\cdot\lambda_1$.
\end{definition}

\begin{proposition}[{\cite[Proposition~1.8]{da-ja91}}]\label{equivar1}
There is a one-two-one correspondence between equivalence classes
of quasitoric manifolds and characteristic pairs. In particular,
for any quasitoric manifold $M$ over $P$ with characteristic
function~$\lambda$, there is a homeomorphism $M\cong
M(P,\lambda)$.
\end{proposition}
\begin{proof}
Obviously, if two quasitoric manifolds are equivalent over~$P$,
then their characteristic pairs are also equivalent. To establish
the other implication, it is enough to show that a quasitoric
manifold $M$ is equivalent to the canonical model $M(P,\lambda)$.
This follows from Proposition~\ref{equivar}.
\end{proof}

\subsection*{Omniorientations and combinatorial quasitoric
data} Here we elaborate on the combinatorial description of
quasitoric manifolds~$M$. Characteristic pairs $(P,\lambda)$ are
replaced by more naturally defined \emph{combinatorial quasitoric
pairs} $(P,\varLambda)$ consisting of an oriented simple polytope
and an integer matrix of special type. Compared with the
characteristic pair, the pair $(P,\varLambda)$ carries some
additional information, which is equivalent to a choice of
orientation for the manifold $M$ and its characteristic
submanifolds. The terminology and constructions described here
were introduced in~\cite{bu-ra01} and~\cite{b-p-r07}.

\begin{definition}\label{defomni}
An \emph{omniorientation} of a quasitoric manifold $M$ consists of
a choice of orientation for $M$ and each characteristic
submanifold $M_{j}$, $1\le j\le m$.
\end{definition}

In general, an omniorientation on~$M$ cannot be chosen
canonically. However, if $M$ admits a \emph{$T^n$-invariant almost
complex structure}\label{omniacs} (see Definition~\ref{iscs}),
then a choice of such structure provides canonical orientations
for~$M$ and the invariant submanifolds $M_{j}$. We therefore
obtain an omniorientation \emph{associated} with the invariant
almost complex structure. In the case when $M$ has an invariant
almost complex structure (for example, when $M$ is a toric
manifold), we always choose the associated omniorientation.
Otherwise we choose an omniorientation arbitrarily.

The stabiliser $T_{F_j}$ of a characteristic submanifold
$M_{j}\subset M$ can be written as
\begin{equation}
\label{fisotr}
  T_{F_j}=\bigl\{\bigl(e^{2\pi i\lambda_{1j}\varphi},\ldots,e^{2\pi
  i\lambda_{nj}\varphi}\bigr)\in\T^n\bigr\},
\end{equation}
where $\varphi\in\R$ и
$\lambda_j=(\lambda_{1j},\ldots,\lambda_{nj})^t\in\Z^n$ is a
primitive vector. (In the coordinate-free notation used in
Chapter~\ref{toric}, the vector $\lambda_j$ belongs to the lattice
$N$ of one-parameters subgroups of the torus.) This vector is
determined by the subgroup $T_{F_j}$ up to sign. A choice of this
sign (and therefore an unambiguous choice of the vector) defines a
parametrisation of the circle subgroup~$T_{F_j}$.

An omniorientation of $M$ provides a canonical way to choose the
vectors $\lambda_j$. Indeed, the action of a parametrised circle
$T_{M_{j}}\subset\T^n$ defines an orientation in the normal bundle
$\nu_j$ of the embedding $M_{j}\subset M$. An omniorientation also
defines an orientation on $\nu_j$ by means of the following
decomposition of the tangent bundle:
\[
  {\mathcal T}\!M|_{M_{j}}={\mathcal T}\!M_{j}\oplus\nu_j.
\]
Now we choose the direction of the primitive vector $\lambda_j$ so
that these two orientations coincide.

Having fixed an omniorientation, we can extend
correspondence~\eqref{charf} to a map of lattices
\[
  \varLambda\colon\Z^m\to\Z^n,\quad\mb e_j\mapsto\lambda_j,
\]
which we refer to as a \emph{directed}\label{dircharf}
characteristic function. A characteristic function is assumed to
be directed whenever an omniorientation is chosen.

We can think of a directed characteristic function as an integer
$n\times m$-matrix~$\varLambda$ with the following property: if
the intersection of facets $F_{j_1},\ldots,F_{j_k}$ is nonempty,
then the vectors $\lambda_{j_1},\ldots,\lambda_{j_k}$ form a part
of basis of the lattice~$\Z^n$. We refer to such matrices as
\emph{characteristic}\label{charmat}. In particular, we can write
any vertex $v\in P$ as an intersection of $n$ facets:
$v=F_{j_1}\cap\cdots\cap F_{j_n}$, and consider the maximal minor
$\varLambda_{v}=\varLambda_{j_1,\ldots,j_n}$ formed by the columns
$j_1,\ldots,j_n$ of matrix~$\varLambda$. Then
\begin{equation}
\label{detLv}
  \det\varLambda_{v}=\pm1.
\end{equation}

Since $M\cong M(P,\lambda)=P\times T^n/\!\sim$, and no
identifications are made over the interior of~$P$, a choice of an
orientation for $M$ is equivalent to a choice of an orientation
for the polytope $P$ (once we assume that the torus $T^n$ is
oriented canonically).
%we have a decomposition
%\begin{equation}\label{orientM}
%  \mathcal T_{(p,t)}M\cong\mathcal T_p P\oplus\mathcal T_t T^n.
%\end{equation}
%of the tangent space at a point $(p,t)\in M(P,\lambda)$, where
%$p\in P^\circ$ and $t\in T^n$. Now we orient the polytope $P$ by
%the following condition: $(\xi_1,\ldots,\xi_n)$ is a positive
%basis of $\mathcal T_p P$ whenever
%$(\xi_1,\ldots,\xi_n,\eta_1,\ldots,\eta_n)$ is a positive basis of
%$\mathcal T_{(p,t)}M(P,\lambda)$ for any positive basis
%$(\eta_1,\ldots,\eta_n)$ of $\mathcal T_t T^n$.

\begin{definition}\label{cqtp}
Let $P$ be an oriented combinatorial simple $n$-polytope with $m$
facets,
% ordered in such a way the the first $n$ of them intersect
%at a vertex,
and let $\varLambda$ be an integer $n\times m$-matrix satisfying
condition~\eqref{detLv} for any vertex $v\in P$. Then
$(P,\varLambda)$ is called a \emph{combinatorial quasitoric pair}.
\end{definition}

An \emph{equivalence} of omnioriented quasitoric manifolds is
assumed to preserve the omniorientation; in this case the
automorphism $\psi\colon  T^n\to T^n$ in Definition~\ref{psieq} is
orientation-preserving. Similarly, two combinatorial quasitoric
pairs $(P,\varLambda_1)$ and $(P,\varLambda_2)$ are said to be
\emph{equivalent} if the orientation of $P$ coincide and there
exists a square integer matrix $\Psi$ with determinant~1 such that
$\varLambda_2=\Psi\cdot\varLambda_1$.

We can summarise the observations above in the following refined
version of Proposition~\ref{equivar1}:

\begin{proposition}\label{cqtpe}
There is a one-to-one correspondence between equivalence classes
of omnioriented quasitoric manifolds and combinatorial quasitoric
pairs.
\end{proposition}

We shall denote the omnioriented quasitoric manifold corresponding
to a combinatorial quasitoric pair $(P,\varLambda)$ by
$M(P,\varLambda)$.

Let $\mathop{\mathrm{chf}}(P)$ denote the set of directed
characteristic functions $\lambda$ for~$P$. The group $GL(n,\Z)$
of automorphisms of the torus $T^n$ acts on the set
$\mathop{\mathrm{chf}}(P)$ from the left. Proposition~\ref{cqtpe}
establishes a one-two-one correspondence
\begin{equation}\label{lcset}
  GL(n,\Z)\backslash\mathop{\mathrm{chf}}(P)\longleftrightarrow
  \{\text{equivalence classes of omnioriented $M$ over $P$}\}.
\end{equation}

If the facets of $P$ are ordered in such a way that the first $n$
of them meet at a vertex~$v$, i.e. $F_1\cap\cdots\cap F_n=v$, then
each coset from $GL(n,\Z)\backslash\mathop{\mathrm{chf}}(P)$
contains a unique directed characteristic function given by a
matrix of the form
\begin{equation}\label{lamat}
\varLambda=\left(I\;|\;\varLambda_\star\right)\;=\;\begin{pmatrix}
  1&0&\ldots&0&\lambda_{1,n+1}&\ldots&\lambda_{1,m}\\
  0&1&\ldots&0&\lambda_{2,n+1}&\ldots&\lambda_{2,m}\\
  \vdots&\vdots&\ddots&\vdots&\vdots&\ddots&\vdots\\
  0&0&\ldots&1&\lambda_{n,n+1}&\ldots&\lambda_{n,m}
\end{pmatrix}
\end{equation}
where where $I$ is the unit matrix and $\varLambda_\star$ is an
$n\times(m-n)$-matrix. We refer to~\eqref{lamat} as a
\emph{refined} characteristic matrix, and to $\varLambda_*$ as its
\emph{refined submatrix}. If a characteristic matrix is given in a
non-refined form $\varLambda=(A\;|\;B)$, where $A$ has size
$n\times n$, then its refined representative is given by
$(I\;|\;A^{-1}B)$.

\subsection*{Smooth and stably complex structures}
Let $(P,\varLambda)$ be a combinatorial quasitoric pair. The
matrix $\varLambda$ defines an epimorphism of tori
$\exp\varLambda\colon\mathbb T^m\to \mathbb T^n$, whose kernel we
denote by $K=K(\varLambda)$. Condition~\eqref{detLv} implies that
$K\cong\mathbb T^{m-n}$. We therefore have an exact sequence of
tori similar to~\eqref{kgrou}. Also, there is the moment-angle
manifold $\zp$ corresponding to the polytope~$P$ (see
Section~\ref{mampol}).

\begin{proposition}\label{kact}
The group $K(\varLambda)\cong T^{m-n}$ acts freely and smoothly
on~$\zp$. There is a $T^n$-equivariant homeomorphism
\[
  \zp/K(\varLambda)\stackrel{\cong}{\longrightarrow}M(P,\varLambda)
\]
between the quotient $\zp/K(\varLambda)$ and the canonical model
$M(P,\varLambda)$.
\end{proposition}
\begin{proof}
The fact that $K$ acts freely on $\zp$ is proved in the same way
as Proposition~\ref{freeaction}~(a) and
Theorem~\ref{zksmooth}~(a). The stabiliser of a point $z\in\zp$
with respect to the $\T^m$-action is a coordinate subtorus $\T^I$
for some $I\in\sK_P$. Namely, if $z\in\zp$ projects to $x\in P$,
then $I=I_{\mb x}=\{i\in[m]\colon x\in F_i\}$.
Condition~\eqref{detLv} implies that the restriction of the
homomorphism $\exp\varLambda\colon\mathbb T^m\to \mathbb T^n$ to
any such subtorus $\T^I$ is injective, and therefore the kernel
$K$ intersects each $\T^I$, $I\in\sK_P$, trivially. Hence $K$ acts
freely on $\zp$, and the quotient $\zp/K$ is a $2n$-dimensional
manifold with an action of the torus $\mathbb T^m/K\cong\T^n$.

To prove the second statement, we identify $\zp$ with the quotient
$P\times\T^m/{\sim}$, see Section~\ref{mampol}. Then the
projection $\zp\to\zp/K$ is identified with the projection
\[
  \zp=P\times\T^m/{\sim}\to P\times\T^n/{\sim},
\]
induced by the homomorphism $\exp\varLambda\colon\T^m\to\T^n$. By
Construction~\ref{der}, the quotient $P\times\T^n/{\sim}$ is the
canonical quasitoric manifold $M(P,\varLambda)$.
\end{proof}

\begin{definition}\label{cansmstqt}
Using Proposition~\ref{kact} we obtain a smooth structure on the
canonical model $M(P,\varLambda)$ as the quotient of the smooth
manifold $\zp$ by the smooth action of $K(\varLambda)$, with the
induced smooth action of the $n$-torus $\T^m/K(\varLambda)$. We
refer to this smooth structure on $M(P,\varLambda)$ and on any
quasitoric manifold equivalent to it as \emph{canonical}.
\end{definition}

Now let $p\colon \zp\to P$ be the projection, and consider the
submanifolds $p^{-1}(F_i)\subset\zp$ corresponding to facets $F_i$
of~$P$ (see Exercise~\ref{facesubmzp}). The submanifold
$p^{-1}(F_i)$ is fixed by the $i$th coordinate subcircle
in~$\T^m$. Denote by $\C_i$ the space of the 1-dimensional complex
representation of the torus~$\T^m$ induced from the standard
representation in $\C^m$ by the projection $\C^m\to\C_i$ onto the
$i$th coordinate. Let $\zp\times\C_i\to\zp$ be the trivial complex
line bundle; we view it as an equivariant $\T^m$-bundle with
diagonal action of~$\T^m$. Then the restriction of the bundle
$\zp\times\C_i\to\zp$ to the invariant submanifold $p^{-1}(F_i)$
is $\T^m$-isomorphic to the normal bundle of the embedding
$p^{-1}(F_i)\subset\zp$. By taking quotient with respect to the
diagonal action of $K=K(\varLambda)$ we obtain a $T^n$-equivariant
complex line bundle
\begin{equation}\label{rhoi}
  \rho_i\colon \zp\times_K\C_i\to\zp/K=M(P,\varLambda)
\end{equation}
over the quasitoric manifold $M=M(P,\varLambda)$. The restriction
of the bundle $\rho_i$ to the characteristic submanifold
$p^{-1}(F_i)/K=M_{i}$ is isomorphic to the normal bundle of
$M_{i}\subset M$ (an exercise). The resulting complex structure on
this normal bundle is the one defined by the omniorientation of
$M(P,\varLambda)$.

\begin{theorem}[{\cite[Theorem 6.6]{da-ja91}}]\label{taum}
There is the following isomorphism of real $T^n$-bundles over
$M=M(P,\varLambda)$:
\begin{equation}\label{stabsplitqt}
  {\mathcal
  T}\!M\oplus\underline{\R}^{2(m-n)}\cong\rho_1\oplus\cdots\oplus\rho_m;
\end{equation}
here $\underline{\R}^{2(m-n)}$ denotes the trivial real
$2(m-n)$-dimensional $T^n$-bundle over~$M$.
\end{theorem}
\begin{proof}
The proof given here differs from that of~\cite{da-ja91}: we use
the equivariant framing of $\zp$ coming from its realisation by an
intersection of quadrics, see Section~\ref{mampol}. Let
$i_{\mathcal Z}\colon\zp\to\C^m$ be the $\T^m$-equivariant
embedding, see~\eqref{cdiz}. We have a $\T^m$-equivariant
decomposition
\begin{equation}\label{tauzp1}
  \mathcal T\zp\oplus\nu(i_{\mathcal Z})=\zp\times\C^m
\end{equation}
obtained by restricting the tangent bundle $\mathcal T\C^m$
to~$\zp$. The normal bundle $\nu(i_{\mathcal Z})$ is
$\T^m$-equivariantly trivial by Theorem~\ref{zpsmooth}, and
$\zp\times\C^m$ is isomorphic, as a $\T^m$-bundle, to the sum of
line bundles $\zp\times\C_i$, \ $1\le i\le m$. Further, we have
\begin{equation}\label{tauzp2}
  \mathcal T\zp=q^*({\mathcal T}\!M)\oplus\xi,
\end{equation}
where $\xi$ is the tangent bundle along the fibres of the
principal $K$-bundle $q\colon\zp\to\zp/K=M$,
see~\cite[Corollary~6.2]{szcz64}. The bundle $\xi$ is induced by
the projection $q$ from the vector bundle over $M$ associated with
the principal bundle $\zp\to M$ through the adjoint representation
of the Lie group~$K$; since $K$ is abelian, this bundle is
trivial. Taking quotient of identity~\eqref{tauzp1} by the action
of~$K$ and using~\eqref{tauzp2}, we obtain a decomposition
\begin{equation}\label{xinu}
  {\mathcal T}\!M\oplus(\xi/K)\oplus(\nu(i_{\mathcal Z})/K)\cong
  \zp\times_{K}\C^m.
\end{equation}
As we have seen above, both $\xi$ and $\nu(i_{\mathcal Z})$ are
trivial $\T^m$-bundles, so that $(\xi/K)\oplus(\nu(i_{\mathcal
Z})/K)\cong\underline{\R}^{2(m-n)}$. Also,
$\zp\times_{K}\C^m\cong\rho_1\oplus\cdots\oplus\rho_m$ as
$T^n$-bundles.
\end{proof}

The isomorphism of Theorem~\ref{taum} gives us an isomorphism of
the stable tangent bundle of $M$ with a complex vector bundle:
this is the setup for a stably complex structure (see
Definition~\ref{iscs}).

\begin{corollary}\label{qtscs}
An omnioriented quasitoric manifold $M$ has a canonical
$T^n$-invariant stably complex structure $c_{\mathcal T}$ defined
by the isomorphism of~\eqref{stabsplitqt}.
\end{corollary}

The corresponding bordism classes $[M]\in\varOmega^U_{2n}$ will be
studied in Section~\ref{qtgenera}.

\begin{example}\label{2cp1qt}
Let us see which stably complex structures we can obtain from
different omniorientations on the simplest quasitoric manifold
$S^2$ over the segment~$I^1$ using Theorem~\ref{taum}. The
standard complex structure of $\C P^1$ has the characteristic
matrix $(1\,-\!\!1)$. The group $K$ is the diagonal circle
$\{(t,t)\}\subset\T^2$, and~\eqref{xinu} becomes the standard
decomposition
\[
  \mathcal T\C P^1\oplus\underline\C\cong
  S^3\times_K\C^2=\bar\eta\oplus\bar\eta,
\]
where $\eta$ is the tautological and $\bar\eta$ is the canonical
line bundle, see Exercise~\ref{hopfexercise}.

On the other hand, there is an omniorientation of $S^2$
corresponding to the characteristic matrix $(1\;1)$. Then
$K=\{(t^{-1},t)\}\subset\T^2$, and~\eqref{xinu} becomes
\[
  \mathcal T S^2\oplus\underline{\R}^2\cong
  S^3\times_K\C^2=\eta\oplus\bar\eta,
\]
which is the trivial stably complex structure on~$S^2$ (see
Example~\ref{2cp1}).
\end{example}

\subsection*{Weights and signs at fixed points}
Any fixed point $v$ of the $T^n$-action on a quasitoric manifold
$M$ is isolated. It can be obtained as an intersection
$M_{{j_1}}\cap\cdots\cap M_{{j_n}}$ of $n$ characteristic
submanifolds and corresponds to a vertex
${F_{j_1}}\cap\cdots\cap{F_{j_n}}$ of the polytope~$P$. Therefore,
the tangent space to $M$ at $v$ decomposes into the sum of normal
spaces to $M_{{j_k}}$ for $1\le k\le n$:
\begin{equation}\label{2orie}
  \mathcal T_v M=(\rho_{j_1}\oplus\cdots\oplus\rho_{j_n})|_v.
\end{equation}
We use this decomposition to identify $\mathcal T_v M$ with
$\mathbb C^n$; then the tangent space to $M_{{j_k}}$ is given in
the corresponding coordinates $(z_1,\ldots,z_n)$ by the equation
$z_k=0$. The representation of the torus $T^n$ in the tangent
space $\mathcal T_v M\cong\C^n$ is determined by its set of
\emph{weights}\label{weighttoru} $\mb w_k(v)\in\Z^n$, \ $1\le k\le
n$. Namely, for $\mb t=(e^{2\pi i\varphi_1},\ldots,e^{2\pi
i\varphi_n})\in T^n$ and $\mb z=(z_1,\ldots,z_n)\in\mathcal T_v
M$, we have
\[
  \mb t\cdot\!\mb z=(e^{2\pi i\langle\mb w_1(v),\,\varphi\rangle}z_1,
  \ldots,e^{2\pi i\langle\mb w_n(v),\,\varphi\rangle}z_n),
\]
where $\varphi=(\varphi_1,\ldots,\varphi_n)\in\R^n$. The weights
can be found from the combinatorial quasitoric pair
$(P,\varLambda)$ using the following statement:

\begin{proposition}\label{qtw}
Let $M=M(P,\varLambda)$ be a quasitoric manifold. The weights $\mb
w_1(v),\ldots,\mb w_n(v)$ of the tangent representation of $T^n$
at a fixed point $v=M_{{j_1}}\cap\cdots\cap M_{{j_n}}$ are given
by the columns of the square matrix $W_v$ satisfying the identity
\[
  W^t_v\,\varLambda_v^{\phantom{t}}=I.
\]
In other words, $\{\mb w_1(v),\ldots,\mb w_n(v)\}$ is the lattice
basis conjugate to $\{\lambda_{j_1},\ldots,\lambda_{j_n}\}$.
\end{proposition}
\begin{proof}
First, note that the local standardness of the action implies that
$\{\mb w_1(v),\ldots,\mb w_n(v)\}$ is a lattice basis. (The fact
that $\{\lambda_{j_1},\ldots,\lambda_{j_n}\}$ is a lattice basis
is expressed by identity~\eqref{detLv}.)

Since the one-parameter subgroup $T_{F_{j_k}}\subset T^n$
(see~\eqref{fisotr}) fixes the hyperplane $z_k=0$ tangent to
$M_{{j_k}}$, we obtain that $\langle\mb
w_i(v),\lambda_{j_k}\rangle=0$ for $i\ne k$. Therefore,
$W_v^t\varLambda_v^{\phantom{t}}$ is a diagonal matrix. Now, the
columns of both $W_v$ and $\varLambda_v$ are lattice bases, which
implies $\langle\mb w_k(v),\lambda_{j_k}\rangle=\pm1$ for $1\le
k\le n$.

On the other hand, the complex structure on the line bundle
$\rho_{j_k}$ comes from the orientation induced by the action of
the one-parameter subgroup of $T^n$ corresponding to the vector
$\lambda_{j_k}$. Hence $\langle\mb w_k(v),\lambda_{j_k}\rangle>0$
for $1\le k\le n$.
\end{proof}

The signs of the fixed points defined by the $T^n$-invariant
stably complex structure on $M$ (see Definition~\ref{defsign}) can
be calculated in terms of the combinatorial quasitoric pair
$(P,\varLambda)$ as follows:

\begin{lemma}\label{qts}
Let $v=M_{j_1}\cap\cdots\cap M_{j_n}$ be a fixed point.
\begin{itemize}
\item[(a)] In terms of decomposition~\eqref{2orie}, we have $\sigma(v)=1$ if the orientation of the space $\mathcal T_v
M$ determined by the orientation of~$M$ coincides with the
orientation of the space
$(\rho_{j_1}\oplus\cdots\oplus\rho_{j_n})|_v$ determined by the
orientation of the line bundles $\rho_{j_k}$, \ $1\le k\le n$.
Otherwise, $\sigma(v)=-1$.

\item[(b)] In terms of the combinatorial quasitoric pair $(P,\varLambda)$, we have
\[
  \sigma(v)=\mathop{\mathrm{sign}}\bigl(
  \det({\lambda_{j_1},\ldots,\lambda_{j_n}})
  \det(\mb a_{j_1},\ldots,\mb a_{j_n})\bigr),
\]
where $\mb a_{j_1},\ldots,\mb a_{j_n}$ are inward-pointing normals
to the facets $F_{j_1},\ldots,F_{j_n}$.
\end{itemize}
\end{lemma}
\begin{proof}
To prove (a), we note that the complex line bundle $\rho_i$ is
trivial over the complement to $M_i$ in~$M$. Therefore, the
nontrivial part of the $T^n$-representation
$(\rho_1\oplus\cdots\oplus\rho_m)|_v$ is exactly
$(\rho_{j_1}\oplus\cdots\oplus\rho_{j_n})|_v$. This implies that
the composite map in Definition~\ref{defsign} is given by
\begin{equation}\label{ormap}
  \mathcal T_v M\to (\rho_{j_1}\oplus\cdots\oplus\rho_{j_n})|_v.
\end{equation}

To prove (b), we write map~\eqref{ormap} in coordinates. To do
this, we identify $\C^m$ with $\R^{2m}$ by mapping a point
$(z_1,\ldots,z_m)\in\C^m$ to
$(x_1,\ldots,x_m,y_1,\ldots,y_m)\in\R^{2m}$, where $z_k=x_k+iy_k$.
Using decomposition~~\eqref{xinu}, we obtain that the map
\[
  \mathcal T_v M\to\mathcal T_v M\oplus\R^{2(m-n)}
  \stackrel{c_{\mathcal T}}\longrightarrow
  (\rho_1\oplus\cdots\oplus\rho_m)|_v\cong\R^{2m}
\]
from Definition~\ref{defsign} is given by the $2m\times2n$-matrix
\[
\begin{pmatrix}
  A^t&\mathbf0\\
  \mathbf0&\varLambda^t
\end{pmatrix}
\]
where $A$ denotes the $n\times m$-matrix with columns vectors $\mb
a_i$ defined by presentation~\eqref{ptope} of the polytope~$P$.
Map~\eqref{ormap} is obtained by restricting to the submatrices
of~$A^t$ and $\varLambda^t$ formed by rows with numbers
$j_1,\ldots j_n$, which implies the required formula for the sign.
\end{proof}

\begin{example}\label{torvarsig}
Let $V_P$ be the projective toric manifold corresponding to a
simple lattice polytope $P$ given by~\eqref{ptope}. Then
$\lambda_i=\mb a_i$ for $1\le i\le m$. Proposition~\ref{qtw}
implies that the weights $\mb w_1(v),\ldots,\mb w_n(v)$ of the
tangent representation of the torus at a fixed point $v\in V_P$
are the primitive vectors along the edges of $P$ pointing out
of~$v$. Furthermore, Lemma~\ref{qts} implies that $\sigma(v)=1$
for all~$v$.
\end{example}

For a general quasitoric manifold, the weights $\mb
w_1(v),\ldots,\mb w_n(v)$ are not vectors along edges. However, by
Proposition~\ref{qtw}, there is a natural one-to-one
correspondence
\begin{multline}\label{qtaf}
  \{\text{oriented edges of }P\}\\ \leftrightarrow
  \{\text{weights of the tangential $T^n$-representations at fixed points}\}.
\end{multline}
Under this correspondence, an edge $e$ coming out of a vertex
$v=F_{j_1}\cap\cdots\cap F_{j_n}\in P$ maps to $\mb w_k(v)$, where
$F_{j_k}$ is the unique facet containing $v$ and not
containing~$e$ and $\mb w_k(v)$ is the $k$th vector of the
conjugate basis of $\lambda_{j_1},\ldots,\lambda_{j_n}$.

\begin{lemma}\label{ws}
Let $\mb e_1,\ldots,\mb e_n$ be vectors along the edges coming out
of a vertex~$v\in P$, and let $\mb w_1(v),\ldots,\mb w_n(v)$ be
the corresponding weights. Then there is the following formula for
the sign of the vertex:
\[
  \sigma(v)=\mathop{\mathrm{sign}}\bigl(
  \det(\mb w_1(v),\ldots,\mb w_n(v))
  \det(\mb e_1,\ldots,\mb e_n)\bigr).
\]
\end{lemma}
\begin{proof}
This follows from Lemma~\ref{qts} and the fact that $\{\mb
w_1(v),\ldots,\mb w_n(v)\}$ is the conjugate basis
of~$\{\lambda_{j_1},\ldots,\lambda_{j_n}\}$, while the vectors
$\{\mb e_1,\ldots,\mb e_n\}$ can be scaled so that they form the
conjugate basis of $\{\mb a_{j_1},\ldots,\mb a_{j_n}\}$.
\end{proof}

\begin{remark}
Formulae of Lemmata~\ref{qts} and~\ref{ws} can be rephrased as the
following practical rule for calculation of signs, which will be
used below. Write the vectors $\lambda_{j_1},\ldots,\lambda_{j_n}$
(respectively, $\mb w_1(v),\ldots,\mb w_n(v)$) into a square
matrix in the order given by the orientation of~$P$, that is, in
such a way that inward-pointing normals of the corresponding
facets (respectively, vectors along the corresponding edges) form
a positive basis of~$\R^n$. Then the determinant of this matrix
is~$\sigma(v)$. In other words, assuming that the weights $\mb
w_1(v),\ldots,\mb w_n(v)$ are ordered so that vectors along their
corresponding edges form a positive basis, we can rewrite the
formula of Lemma~\ref{ws} as
\begin{equation}\label{wsimpl}
  \sigma(v)=\det(\mb w_1(v),\ldots,\mb w_n(v)).
\end{equation}
\end{remark}

\begin{example}
\label{cpnoo} The complex projective plane $\C P^2$ has the
standard stably complex complex structure coming from the bundle
isomorphism
\[
  \mathcal T(\C P^2)\oplus\underline{\C}\cong\bar\eta\oplus\bar\eta\oplus\bar\eta,
\]
where $\eta$ is the tautological line bundle. The orientation is
determined by the complex structure. The toric manifold $\C P^2$
corresponds to the lattice 2-simplex $\varDelta^2$ with vertices
$(0,0)$, $(1,0)$ and $(0,1)$. The column vectors $\lambda_1$,
$\lambda_2$, $\lambda_3$ of the matrix $\varLambda$ are the
primitive inward-pointing normals to the facets (i.e. they
coincide with the vectors $\mb a_1,\mb a_2,\mb a_3$ for the
standard presentation of~$\varDelta^2$). The weights of the
tangential $T^2$-representation at a fixed point are the primitive
vectors along the edges coming out of the corresponding vertex.
This is shown in Fig.~\ref{qtm1}. We have
$\sigma(v_1)=\sigma(v_2)=\sigma(v_3)=1$.
\begin{figure}[h]
\begin{center}
\begin{picture}(100,60)
  \put(20,10){\line(0,1){45}}
  \put(20,10){\line(1,0){45}}
  \put(20,55){\line(1,-1){45}}
  \put(33,23){\oval(13,13)[b]}
  \put(33,23){\oval(13,13)[tr]}
  \put(34,29.5){\vector(-1,0){2}}
  \put(16,6){$v_1$}
  \put(22,6){\small $(1,0)$}
  \put(54,6){\small $(-1,0)$}
  \put(67,6){$v_2$}
  \put(33,0){$\lambda_2=(0,1)$}
  \put(11,13){\small $(0,1)$}
  \put(0,30){$\lambda_1=(1,0)$}
  \put(8,50){\small $(0,-1)$}
  \put(16,57){$v_3$}
  \put(25,52){\small $(1,-1)$}
  \put(62,15){\small $(-1,1)$}
  \put(47,32){$\lambda_3=(-1,-1)$}
\end{picture}%
\caption{$\mathcal T(\C
P^2)\oplus\C\cong\bar\eta\oplus\bar\eta\oplus\bar\eta$}
\label{qtm1}
\end{center}
\end{figure}
\end{example}

\begin{example}\label{cpquas}
Now consider $\C P^2$ with the omniorientation defined by the
three vectors $\lambda_1,\lambda_2,\lambda_3$ in Fig.~\ref{qtm2}.
This omniorientation differs from the previous one by the sign
of~$\lambda_3$. The stably complex structure is defined by the
isomorphism
\[
  \mathcal T(\C
  P^2)\oplus\underline{\R}^2\cong\bar\eta\oplus\bar\eta\oplus\eta.
\]
Using formula~\eqref{wsimpl}, we calculate
$$
  \sigma(v_1)=\begin{vmatrix} 1&0\\0&1 \end{vmatrix}=1,\quad
  \sigma(v_2)=\begin{vmatrix} -1&1\\1&0 \end{vmatrix}=-1,\quad
  \sigma(v_3)=\begin{vmatrix} 0&1\\1&-1 \end{vmatrix}=-1.
$$
\begin{figure}[h]
\begin{center}
\begin{picture}(100,60)
  \put(20,10){\line(0,1){45}}
  \put(20,10){\line(1,0){45}}
  \put(20,55){\line(1,-1){45}}
  \put(33,23){\oval(13,13)[b]}
  \put(33,23){\oval(13,13)[tr]}
  \put(34,29.5){\vector(-1,0){2}}
  \put(16,6){$v_1$}
  \put(22,6){\small $(1,0)$}
  \put(54,6){\small $(1,0)$}
  \put(67,6){$v_2$}
  \put(32,0){$\lambda_2=(0,1)$}
  \put(11,13){\small $(0,1)$}
  \put(0,30){$\lambda_1=(1,0)$}
  \put(11,50){\small $(0,1)$}
  \put(16,57){$v_3$}
  \put(25,52){\small $(1,-1)$}
  \put(62,15){\small $(-1,1)$}
  \put(47,32){$\lambda_3=(1,1)$}
\end{picture}%
\caption{$\mathcal T(\C
P^2)\oplus\C\cong\bar\eta\oplus\bar\eta\oplus\eta$}\label{qtm2}
\end{center}
\end{figure}
\end{example}

If the omniorientation of $M$ comes from a $T^n$-invariant almost
complex structure, then all signs of fixed points are positive
because the two orientation in~\eqref{2orie} coincide. As we shall
see below in this section, the top Chern number of the stably
complex structure of an omnioriented quasitoric manifold $M$ is
equal to the sum of signs of fixed points: $c_n[M]=\sum_{v\in
M}\sigma(v)$. On the other hand, the Euler characteristic
$\chi(M)$\label{eulerqtori} is equal to the number of fixed
points. Therefore, the positivity of all signs is equivalent to
the condition $c_n[M]=\chi(M)$. According to the classical result
of Thomas~\cite{thom67}, this condition is sufficient for a stably
complex manifold $M$ to be almost complex. The following result of
Kustarev shows that this almost complex structure can be made
$T^n$-invariant:

\begin{theorem}[\cite{kust09}]\label{kusth} A quasitoric manifold $M$ admits a
$T^n$-invariant almost complex structure if and only if it admits
an omniorientation in which all signs of fixed points are
positive.
\end{theorem}

Theorem~\ref{kusth} allows one to construct an invariantly almost
complex quasitoric manifold which is not toric (see
Exercise~\ref{qtnt} below). Such an almost complex structure
cannot be integrable. Indeed, according to the result of
Ishida--Karshon~\cite{is-ka} (Theorem~\ref{iskath}, see
also~\cite{is-ma12}), a quasitoric manifold with an invariant
complex structure is biholomorphic to a compact toric variety.

\subsection*{Cohomology ring and characteristic classes}
The cohomology ring $H^*(M)$ of a quasitoric manifold $M$ has the
same structure as the cohomology ring of a nonsingular compact
toric variety (see Theorem~\ref{danjur}). In particular, the ring
$H^*(M)$ is generated by two-dimensional classes $v_i$ dual to the
characteristic submanifolds (or equivalently, by the first Chern
classes of line bundles $\rho_i$ from Theorem~\ref{taum}). The
elements $v_i$ satisfy two types of relations: monomial relations
coming from the face ring of the polytope~$P$ and linear relations
coming from the characteristic matrix~$\varLambda$.

Let $M=M(P,\varLambda)$ be an omnioriented quasitoric manifold,
and let $\pi\colon M\to P$ be the projection onto the orbit space.
We start by describing a canonical cell decomposition of $M$ with
only even-dimensional cells. It was first constructed by
Khovanskii~\cite{khov86} for toric manifolds.

\begin{construction}\label{morse1}
Recall the `Morse-theoretic' arguments used in the proof of the
Dehn--Sommerville relations (Theorem~\ref{ds}). There we turned
the 1-skeleton of $P$ into an oriented graph and defined the index
$\ind(v)$ of a vertex $v\in P$ as the number of incoming edges.
The incoming edges of $v$ span a face $G_v$ of dimension
$\ind(v)$. Denote by $\widehat{G}_v$ the subset obtained by
removing from $G_v$ all faces not containing~$v$. Then
$\widehat{G}_v$ is homeomorphic to $\R^{\ind(v)}_\ge$ and is
contained in the open set $U_v\subset P$ from
Example~\ref{examplemwc}.4. The preimage
$e_v=\pi^{-1}\widehat{G}_v$ is homeomorphic to $\C^{\ind(v)}$ and
the union of subsets $e_v\subset M$ over all vertices $v\in P$
defines a cell decomposition of~$M$. Observe that all cells have
even dimension, and the closure of the cell $e_v$ is the
quasitoric submanifold $\pi^{-1}(G_v)\subset M$.
\end{construction}

\begin{proposition}
\label{qtbn} The homology groups of a quasitoric manifold
$M=M(P,\varLambda)$ vanish in odd dimensions, and therefore are
free abelian in even dimensions. Their ranks (Betti numbers) are
given by
\[
  b_{2i}(M)=h_i(P),
\]
where $h_i(P)$, $i=0,1,\ldots,n$, are the components of the
$h$-vector of~$P$.
\end{proposition}
\begin{proof}
The rank of $H_{2i}(M;\Z)$ is equal to the number of
$2i$-dimensional cells in the above cell decomposition. This
number is equal to the number of vertices of index~$i$, which is
$h_i(P)$ by the argument from the proof of Theorem~\ref{ds}.
\end{proof}

Now consider the face ring $\Z[P]$ (see Section~\ref{srr}) and
define its elements
\begin{equation}
\label{theta}
  t_i=\lambda_{i1}v_1+\cdots+\lambda_{im}v_m\in\Z[P],\quad 1\le i\le n,
\end{equation}
corresponding to the rows of the characteristic
matrix~$\varLambda$.

\begin{lemma}
\label{crs} The elements $t_1,\ldots,t_n$ form a linear regular
sequence (an lsop) in the ring~$\Z[P]$. Conversely, any
lsop~\eqref{theta} in the ring $\Z[P]$ defines a combinatorial
quasitoric pair $(P,\varLambda)$.
\end{lemma}
\begin{proof}
Since $\Z[P]$ is a Cohen--Macaulay ring
(Corollary~\ref{spherecm}), any lsop is a regular sequence by
Proposition~\ref{rlsop}. So it is enough to show that
$t_1,\ldots,t_n$ is an lsop. Condition~\eqref{detLv} implies that
for any vertex $v=F_{j_1}\cap\cdots\cap F_{j_n}$ the restrictions
of the elements $t_1,\ldots,t_n$ form a basis in the linear part
of the polynomial ring $\Z[v_{j_1},\ldots,v_{j_n}]$. By
Lemma~\ref{lsopr}, this condition specifies lsop's in the
ring~$\Z[P]$.
\end{proof}

\begin{theorem}[\cite{da-ja91}]\label{qtcoh}
Let $M=M(P,\varLambda)$ be a quasitoric manifold with
$\varLambda=(\lambda_{ij})$, $1\le i\le n$, $1\le j\le m$. The
cohomology of $M$ is given by
\[
  H^*(M)=\Z[v_1,\ldots,v_m]/\mathcal I,
\]
where $v_i\in H^2(M)$ is the class dual to the characteristic
submanifold~$M_i$, and $\mathcal I$ is the ideal generated by
elements of the following two types:
\begin{itemize}
\item[(a)] $v_{i_1}\cdots v_{i_k}$ whenever $M_{i_1}\cap\cdots\cap M_{i_k}=\varnothing$ (the Stanley--Reisner relations);

\item[(b)] the linear forms $t_i=\lambda_{i1}v_1+\cdots+\lambda_{im}v_m$, \ $1\le i\le n$.
\end{itemize}
\end{theorem}

In other words, $H^*(M)$ is the quotient of the face ring $\Z[P]$
by the ideal generated by linear forms~\eqref{theta}. We shall
prove a more general result, covering both Theorem~\ref{danjur}
and Theorem~\ref{qtcoh}, in the next section (see
Theorem~\ref{theo:stcoh}).

If the matrix $\varLambda$ has refined form~\eqref{lamat}, then
the linear relations between the cohomology classes can be written
as
\begin{equation}\label{viref}
  v_i\;=\;-\lambda_{i,n+1}v_{n+1}-\cdots-\lambda_{i,m}v_m,
  \quad 1\le i\le n.
\end{equation}
It follows that the classes $v_{n+1},\ldots,v_m$ multiplicatively
generate the ring~$H^*(M)$.

If the cohomology ring of a manifold is not generated by
two-dimensional classes, then the manifold does not support a
torus action turning it into a quasitoric manifold. For example,
complex Grassmanians (except complex projective spaces) are not
quasitoric.

\begin{remark}
The structure of the cohomology ring of $M$ given in
Theorem~\ref{qtcoh} allows one to describe $H^*(M)$ as a module
over the Steenrod algebra. Furthermore, the ring $E^*(M)$ can be
described easily for any complex-oriented cohomology theory~$E^*$.
We shall give such a description for the complex cobordism ring
$\varOmega^*_U(M)$ and the action of the Landweber--Novikov
algebra in Chapter~\ref{cobtoric}.
\end{remark}

The Chern classes of the stably complex structure on $M$ (see
Corollary~\ref{qtscs}) can be also described easily:

\begin{theorem}\label{qtchern}
Let $(M,c_{\mathcal T})$ be a quasitoric manifold with the
canonical stably complex structure defined by an omniorientation.
Then, in the notation of Theorem~{\rm\ref{qtcoh}}, we have the
following expression for the total Chern class:
\[
  c(M)=1+c_1(M)+\cdots+c_n(M)=(1+v_1)\cdots(1+v_m)\in H^*(M).
\]
The homology class dual to $c_k(M)\in H^{2k}(M)$ is represented by
the sum of the submanifolds $\pi^{-1}(G)\subset M$ corresponding
to all $(n-k)$-dimensional faces $G\subset P$.
\end{theorem}
\begin{proof}
The first statement holds since the stably complex structure on
$M$ is defined by the isomorphism with the complex bundle
$\rho_1\oplus\cdots\oplus\rho_m$, and $c(\rho_i)=1+v_i$, $1\le
i\le m$. To prove the second statement, note that
\[
  c_k(M)=\sum_{1\le j_1<\cdots< j_k\le m}
  v_{j_1}\!\cdots v_{j_k}\in H^{2k}(M).
\]
Here summands $v_{j_1}\!\cdots v_{j_k}$ for which
$F_{j_1}\cap\cdots\cap F_{j_k}=\varnothing$ are zero by
Theorem~\ref{qtcoh}, and the remaining summands are dual to the
submanifolds $\pi^{-1}(G)$, where $G=F_{j_1}\cap\cdots\cap
F_{j_k}$ is a $(n-k)$-dimensional face.
\end{proof}

\subsection*{Exercises.}

\begin{exercise}\label{23dqt}
For any 2- or 3-dimensional simple polytope $P$, there exists a
quasitoric manifold over~$P$. (Hint: use the model $M(P,\lambda)$
and the 4-colour Theorem, cf.~\cite{da-ja91}.)
\end{exercise}

\begin{exercise}\label{4dqt}
Let $P$ be the dual of a 2-neighbourly (e.g., cyclic) simplicial
polytope of dimension $n\ge4$ with $m\ge 2^n$ vertices. Then there
is no quasitoric manifold over~$P$. (Hint: use the argument from
Example~\ref{dj-lsop}.)
\end{exercise}

\begin{exercise}
Let $M=M(P,\varLambda)$ be a quasitoric manifold and $\rho_i$ the
complex line bundle over $M$ defined by~\eqref{rhoi}. Show that
the restriction of $\rho_i$ to the characteristic submanifold
$M_i$ is isomorphic to the normal bundle of~$M_i$, and the
restriction of $\rho_i$ to the complement $M\setminus M_i$ is
trivial.
\end{exercise}

\begin{exercise}
Calculate the signs of the fixed points for the two
omniorientations on $S^2$ in Example~\ref{2cp1qt} and compare this
with the sign calculation from Example~\ref{2cp1eq}.
\end{exercise}

\begin{exercise}
Let $M$ be a quasitoric manifold over an $n$-polytope~$P$ and let
$\zp$ be the corresponding moment-angle manifold. Show that the
Borel constructions $ET^n\times_{T^n}M$ and $ET^m\times_{T^m}\zp$
are homotopy equivalent. Deduce that the equivariant cohomology
$H^*_{T^n}(M)$ is isomorphic to the face ring~$\Z[P]$.
%(Hint: use Corollary~\ref{cohombk}.)
\end{exercise}

\begin{exercise}\label{qtnt}
Let $(P,\varLambda)$ be the quasitoric pair shown in
Fig.~\ref{hexag}. Show that the corresponding quasitoric manifold
$M(P,\varLambda)$ admits a $T^2$-invariant almost complex
structure, but is not homeomorphic to a toric variety. (Hint: read
Section~\ref{qtgenera}.)
\begin{figure}[h]
\begin{picture}(115,33)
  %lines:
  \put(45,10){\line(3,-2){12.5}}
  \put(45,10){\line(0,1){15}}
  \put(45,25){\line(3,2){12.5}}
  \put(70,10){\line(-3,-2){12.5}}
  \put(70,10){\line(0,1){15}}
  \put(70,25){\line(-3,2){12.5}}
  %orientations:
  \put(57.5,17.5){\oval(8,8)[b]}
  \put(57.5,17.5){\oval(8,8)[tr]}
  \put(58.5,21.5){\vector(-1,0){2}}
  %tuples:
  \put(35.5,17){$\scriptscriptstyle(-1,-1)$}
  \put(46.5,30.2){$\scriptscriptstyle(0,1)$}
  \put(63,30.2){$\scriptscriptstyle(1,0)$}
  \put(46.5,3.7){$\scriptscriptstyle(1,0)$}
  \put(63,3.7){$\scriptscriptstyle(0,1)$}
  \put(71,17){$\scriptscriptstyle(-1,-1)$}
\end{picture}
\vspace{-5mm} \caption{A quasitoric manifold over a
hexagon.}\label{hexag}
\end{figure}
\end{exercise}

\begin{exercise}
The complex Grassmanian $\mbox{\it Gr\/}_k(\C^n)$ of $k$-planes in
$\C^n$ with $2\le k\le n-2$ does not support a torus action
turning it into a (quasi)toric manifold.
\end{exercise}

\begin{exercise}
Let $M$ be an omnioriented quasitoric manifold over a
polytope~$P$. Let $x=F_{i_1}\cap\cdots\cap F_{i_n}$ be a vertex of
$P$ written as an intersection of $n$ facets; it corresponds to
fixed point $x$ of~$M$. Show that
\[
  \bigl\langle v_{i_1}v_{i_2}\!\cdots v_{i_n},[M]\bigr\rangle=\sigma(x),
\]
where $v_i\in H^2(M)$ is the generator corresponding to $F_i$ (see
Theorem~\ref{qtcoh}), $[M]\in H_{2n}(M)$ is the fundamental
homology class of the oriented manifold~$M$, and $\sigma(x)$ is
the sign of~$x$ (see Lemma~\ref{qts}).
\end{exercise}

\begin{exercise}\label{cnsign}
The top Chern number of $M$ is given by $\langle
c_n(M),[M]\rangle=\sum_{v\in P}\sigma(v)$ where the sum is taken
over all vertices $v$ of~$P$. (Hint: use Theorem~\ref{qtchern}.)
\end{exercise}

\begin{exercise}
Show that if $M$ is a toric manifold, then~$c_1(M)\ne0$. Given an
example of an omnioriented quasitoric manifold $M$ with
$c_1(M)=0$.
\end{exercise}

\section{Locally standard $T$-manifolds and torus manifolds}\label{torusman}
In this section we consider two closely related classes of
half-dimensional torus actions on even-dimensional manifolds.

\begin{definition}\label{locastatm}
A \emph{locally standard $T$-manifold} is a smooth connected
closed orientable $2n$-dimensional manifold $M$ with a locally
standard action of an $n$-dimensional torus~$T=T^n$.
\end{definition}

The orbit space $Q=M/T$ of a locally standard $T$-manifold is a
manifold with corners, but, unlike the case of quasitoric
manifolds, it may fail to be a simple polytope. For example, a
free smooth $T$-action with $\dim M=2\dim T$ is locally standard,
but the quotient $Q$ does not have faces at all. The richer is the
combinatorics of $Q$, the more information about the topology of $M$ can be retrieved from this combinatorics.
%Nevertheless, the topology of locally standard $T$-manifolds to a
%large extent can be described via the combinatorial-algebraic
%constructions of Chapters~\ref{combigt} and~\ref{facerings}.
The easiest way to make the combinatorics of the orbit space $Q$
`rich enough' is to require the existence of $T$-fixed points, or
0-faces of~$Q$; then the local standardness condition would also
imply the existence of faces of all dimensions between 0 and~$n$.
We therefore come to the following definition:

\begin{definition}\label{deftoma}
A \emph{torus manifold} is a smooth connected closed orientable
$2n$-dimensional manifold $M$ with an effective smooth action of
an $n$-torus~$T$ such that the fixed point set $M^T$ is nonempty.
\end{definition}

Since $M$ is connected, the $T$-action is effective and  $\dim
M=2\dim T$, it follows that the fixed point set $M^T$ is isolated.
Indeed, the tangential $T^n$-representation at a fixed point is
faithful (Exercise~\ref{faisre}), which implies that this normal
space has dimension at least~$2n$. Furthermore, $M^T$ is finite
because $M$ is compact. So the last condition in the definition
above only excludes the possibility when the number of fixed
points is zero.

Torus manifolds were introduced and studied in the
works~\cite{masu99}, \cite{ha-ma03},~\cite{ha-ma05} of Hattori and
Masuda. Homological aspects of this study which we present here
were developed in~\cite{ma-pa06} and~\cite{m-m-p07}.
%The topology of torus manifolds can be described to some extent
%combinatorially, using the notion of a \emph{multifan}. A multi-fan
%is a bunch of simplicial cones, which may overlap, unlike in the
%usual fan.

The classes of locally standard $T$-manifolds and torus manifolds
contain both toric and quasitoric manifolds in their intersection
(see Section~\ref{reltm}). Quasitoric manifolds are examples with
the most regular structure of the orbit space: a simple polytope
is contractible as a topological space, and the topology of the
manifold is described fully in terms of the combinatorics of
faces. In general, the topology of a locally standard $T$-manifold
depends not only on the combinatorics of the orbit space, but also
on its topology. However, even orbit spaces with trivial topology
may have combinatorics more complicated than that of a polytope.
Examples are simplicial posets reviewed in
Sections~\ref{secsimpos} and~\ref{frsimpo}.

Even in the case of toric manifolds, we have encountered with
combinatorial structures more general than simple polytopes. As we
have seen in Section~\ref{tos}, the orbit space $Q=V/T$ of a
nonsingular compact toric variety $V$ is a manifold with corners
in which all faces, including~$Q$ itself, are acyclic, and all
nonempty intersections of faces are connected. We refer to such a
manifold with corners as a \emph{homology polytope}. It is a
genuine polytope when the toric manifold is projective, but in
general $Q$ may be not combinatorially equivalent to a convex
polytope (see the discussion in Section~\ref{reltm}). As a result,
the class of quasitoric manifolds does not include all nonsingular
compact toric varieties (toric manifolds), which is not very
convenient. At the same time, one could expect that all
topological properties characterising quasitoric manifolds also
hold in a more general situation when the orbit space is a
homology (rather than combinatorial) polytope. This is indeed the
case, as is seen by the results of this section.

The cohomology ring of a locally standard $T$-manifold $M$ is
generated by its degree-two elements if and only if the orbit
space $Q$ is a homology polytope (Theorem~\ref{hptpe}). In this
case, the cohomology ring has the structure familiar from toric
geometry: $H^*(M)$ is isomorphic to the quotient of the face ring
of the orbit space $Q$ by an ideal generated by certain linear
forms.

More generally, we consider $T$-manifolds $M$ whose cohomology
vanishes in odd dimensions. Such a $T$-manifold necessarily has a
fixed point (Lemma~\ref{oddfp}), and is locally standard whenever
$\dim M=2\dim T$ (Theorem~\ref{theo:local standardness}).
Therefore, under the condition $H^{odd}(M)=0$, the classes of
locally standard $T$-manifolds and torus manifolds coincide. In
this case, the equivariant cohomology of $M$ is a free finitely
generated module over the $T$-equivariant cohomology of point,
i.e. over the polynomial ring $H^*(BT)\cong\Z[t_1,\ldots,t_n]$. In
other words, $H^*_T(M)$ is a Cohen--Macaulay ring. The orbit space
of a $T$-manifold with vanishing odd-degree cohomology is not
necessarily a homology polytope, as is seen from a simple example
of a half-dimensional torus action on an even-dimensional sphere
(Example~\ref{2nsphere}). There is a more general notion of a
\emph{face-acyclic} manifold with corners~$Q$, in which all faces
are acyclic, but the intersections of faces are not necessarily
connected. It turns out that the odd-degree cohomology of a
$T$-manifold $M$ vanishes if and only if the orbit space $Q$ is
face-acyclic (Theorem~\ref{cohfa}). The equivariant cohomology of
such $M$ is isomorphic to the face ring of the simplicial poset of
faces of~$Q$ (this ring may no longer be generated in degree two,
see Section~\ref{frsimpo}).

The proofs use several results from the theory of
\emph{GKM-manifolds} and related \emph{GKM-graphs}, whose
foundations were laid in the work~\cite{g-k-m98} of Goresky,
Kottwitz and MacPherson. Relationship between this theory and
torus manifolds is explored further in Section~\ref{weightgraphs}.

\subsection*{Preliminaries: cohomology and fixed points}
Here we obtain some preliminary results about torus actions,
without assuming that the action is locally standard. In this
subsection $M$ is a closed connected smooth orientable manifold
equipped with an effective smooth action of a torus $T$ of
arbitrary dimension.

We denote by $M^T$ the set of $T$-fixed points of~$M$, which is a
disjoint union of finitely many connected submanifolds.

\begin{lemma}\label{sa}
There exists a circle subgroup $S\subset T$ such that $M^S=M^T$.
\end{lemma}
\begin{proof}
According to the standard result~\cite[Theorem~IV.10.5]{bred72},
there is only a finite number of orbit types of the $T$-action
on~$M$, i.e. only finitely many subgroups of $T$ appear as the
stabilisers of the action. We have $M^S=M^T$ whenever $S$ is not
contained in any proper stabiliser subgroup $G\subset T$. This
condition is obviously satisfied for a generic circle $S\subset
T$.
\end{proof}

\begin{lemma}\label{oddfp}
If $H^{odd}(M;\Q)=0$, then $M$ has a $T$-fixed point.
\end{lemma}
\begin{proof}
Choose a generic circle $S\subset T$ satisfying $M^S=M^T$ (see
Lemma~\ref{sa}). If $M$ does not have $T$-fixed points, then it
also does not have $S$-fixed points, which implies that the Euler
characteristic $\chi(M)$ is zero. On the other hand, if
$H^{odd}(M;\Q)=0$, then $\chi(M)>0$: a contradiction.
%The condition $H^{odd}(M;\Q)=0$ implies
%$\chi(M^T)=\chi(M^S)=\chi(M)>0$ (here $\chi(\cdot)$ denote the
%Euler characteristic), hence the set $M^T$ is nonempty.
\end{proof}

The inclusion $M^T\to M$ induces the \emph{restriction map} in
equivariant cohomology:
\begin{equation}\label{rmevc}
  r\colon H^*_T(M)\to H^*_T(M^T)=H^*(BT)\otimes H^*(M^T).
\end{equation}
The equivariant cohomology $H^*_T(M)$ is a $H^*(BT)$-module. We
denote by $H^+(BT)$ the positive-degree part of $H^*(BT)$. We
shall need the following version of the `localisation theorem':

\begin{theorem}\label{locth}
The restriction map $r\colon H^*_T(M)\to H^*_T(M^T)$ becomes an
isomorphism when localised at~$H^+(BT)$.
\end{theorem}

For the proof, see~\cite[p.~40]{hsia75}
or~\cite[Theorem~11.44]{gu-st99}.

\begin{corollary}\label{kerres}
The kernel of the restriction map $r\colon H^*_T(M)\to H^*_T(M^T)$
is a $H^*(BT)$-torsion module.
\end{corollary}

Since $H^*(BT)\cong\Z[t_1,\ldots,t_n]$, we obtain that $H^*_T(M)$
is a free $H^*(BT)$-module if and only if $H^*_T(M)$ is a
Cohen--Macaulay ring (note that the $H^*(BT)$-module $H^*_T(M)$ is
finitely generated, as $M$ is compact). The next statement gives a
topological characterisation of $T$-manifolds with this property:

\begin{lemma}\label{lemm:freeness}
Assume that $M^T$ is finite. Then $H^*_T(M)$ is a free
$H^*(BT)$-module if and only if $H^{odd}(M)=0$. In this case,
$H^*_T(M)\cong H^*(BT)\otimes H^*(M)$ as $H^*(BT)$-modules.
\end{lemma}
\begin{proof}
Assume that $H^{odd}(M)=0$. Then $H^*(M)$ is torsion-free (an
exercise) and the Serre spectral sequence of the bundle
$\rho\colon ET\times_T M\to BT$ collapses at~$E_2$. It follows
that $H^*_T(M)=H^*(ET\times_T M)$ is isomorphic (as a
$H^*(BT)$-module) to the tensor product $H^*(BT)\otimes H^*(M)$,
and therefore it is a free $H^*(BT)$-module.

Assume now that $H^*_T(M)$ is a free $H^*(BT)$-module. Consider
the Eilenberg--Moore spectral sequence of the bundle $\rho\colon
ET\times_T M\to BT$ (see Corollary~\ref{onefib}). It converges to
$H^*(M)$ and has
$$
  E_2^{*,*}=\Tor^{*,*}_{H^*(BT)}\bigl(H^*_T(M),\Z\bigr).
$$
Since $H_T^*(M)$ is a free $H^*(BT)$-module,
\begin{multline*}
  \Tor^{*,*}_{H^*(BT)}\bigl(H^*_T(M),\Z\bigr)=
  \Tor^{0,*}_{H^*(BT)}\bigl(H^*_T(M),\Z\bigr)
  =H^*_T(M)\otimes_{H^*(BT)}\Z\\
  =  H^*_T(M)\big/\bigl(\rho^*(H^+(BT))\bigr).
\end{multline*}
Hence $E_2^{0,*}=H^*_T(M)\bigr/(\rho^*(H^+(BT)))$ and
$E_2^{-p,*}=0$ for $p>0$. It follows that the spectral sequence
collapses at $E_2$, and
\begin{equation} \label{eqn:H(M)}
  H^*(M)=H^*_T(M)\big/\bigl(\rho^*(H^+(BT))\bigr).
\end{equation}
Since $H^*_T(M)$ is a free $H^*(BT)$-module, the restriction
map~\eqref{rmevc} is a monomorphism (see Corollary~\ref{kerres}).
As $M^T$ is finite, $H^*_T(M^T)$ is a sum of polynomial rings,
hence $H^{odd}_T(M)=0$. This together with~(\ref{eqn:H(M)})
implies $H^{odd}(M)=0$.
\end{proof}

We shall be interested in the two special classes of $T$-manifolds
$M$: those with vanishing odd-degree cohomology, and those with
cohomology generated by the degree-two classes. Both these
cohomological properties are inherited by the fixed point set
$M^H$ with respect to the action of any torus subgroup $H\subset
T$. This fact is proved in the next two lemmata; it enables us to
use inductive arguments.

\begin{lemma}\label{lemm:odd=0}
Let $H$ be a torus subgroup of~$T$, and let $N$ be a connected
component of $M^H$. If $H^{odd}(M)=0$, then $H^{odd}(N)=0$.
\end{lemma}
\begin{proof}
Choose a generic circle subgroup $S\subset H$ with $M^S=M^H$, as
in Lemma~\ref{sa}. Let $G\subset S$ be a subgroup of prime
order~$p$. The action of $G$ on $H^*(M)$ is trivial, because $G$
is contained in a connected group~$S$.
By~\cite[Theorem~VII.2.2]{bred72}, $\dim H^{odd}(M^G;\Z_p)\le\dim
H^{odd}(M;\Z_p)$. Hence $H^{odd}(M^G;\Z_p)=0$. Repeating the same
argument for the set $M^G$ with the induced action of the quotient
group $S/G$ (which is again a circle), we conclude that
$H^{odd}(M^K;\Z_p)=0$ for any $p$-subgroup $K$ of~$S$. However,
$M^K=M^S=M^H$ if the order of $K$ is sufficiently large, so we
obtain $H^{odd}(M^H;\Z_p)=0$. Since $p$ is an arbitrary prime, we
conclude that $H^{odd}(M^H)=0$.
\end{proof}

\begin{lemma}\label{sconn}
Let $M,H,N$ be as in Lemma {\rm\ref{lemm:odd=0}}.  If the ring
$H^*(M)$ is generated by its degree-two part, then the restriction
map $H^*(M)\to H^*(N)$ is surjective; in particular, the ring
$H^*(N)$ is also generated by its degree-two part.
\end{lemma}
\begin{proof}
It suffices to prove that the restriction map $H^*(M;\Z_p)\to
H^*(N;\Z_p)$ is surjective for any prime~$p$, because
$H^{odd}(N)=0$ by Lemma~\ref{lemm:odd=0}.

The argument below is similar to that used
in~\cite[Theorem~VII.3.1]{bred72}. As in the proof of
Lemma~\ref{lemm:odd=0}, let $S\subset H$ be a generic circle with
$M^S=M^H$ and let $G\subset S$ be a subgroup of prime order $p$.
By~\cite[Theorem~VII.1.5]{bred72}, the restriction map
$H^k_G(M;\Z_p) \to H^k_G(M^G;\Z_p)$ is an isomorphism for
sufficiently large~$k$. Hence, for any connected component $N'$
of~$M^G$, the restriction $r\colon H^k_G(M;\Z_p) \to
H^k_G(N';\Z_p)$ is surjective when $k$ is large. Now consider the
commutative diagram
\[
\begin{CD}
  H^*_{G}(M;\Z_p) @>r>> H^*_{G}(N';\Z_p)@.\;\cong H^*(BG;\Z_p)
  \otimes H^*(N';\Z_p)\\
  @VVV @VVV @.\\
  H^*(M;\Z_p) @>s>> H^*(N';\Z_p)@.
\end{CD}.
\]
Choose a basis $v_1,\ldots,v_d\in H^2(M;\Z_p)$; then these
elements multiplicatively generate $H^*(M;\Z_p)$. Since
$H^{odd}(M;\Z_p)=H^{odd}(M^G;\Z_p)=0$ and
$\chi(M)=\chi(M^T)=\chi(M^G)$, it follows that $\sum_i\dim
H^i(M;\Z_p)=\sum_i\dim H^i(M^G;\Z_p)$.
By~\cite[Theorem~VII.1.6]{bred72}, the Serre spectral sequence of
the bundle $EG\times_G M\to BG$ collapses at~$E_2$. Therefore, the
vertical map $H^*_G(M;\Z_p)\to H^*(M;\Z_p)$ in the above diagram
is surjective. Let $\xi_j\in H^2_{G}(M;\Z_p)$ be a lift of $v_j$,
and $w_j=s(v_j)$. Let $t$ be a generator of
$H^2(BG;\Z_p)\cong\Z_p$.  We have $H^1(N';\Z_p)=0$ by
Lemma~\ref{lemm:odd=0}, which together with the commutativity of
the diagram above implies $r(\xi_j)=\alpha_jt+w_j$ for some
$\alpha_j\in \Z_p$. Now let $a\in H^*(N';\Z_p)$ be an arbitrary
element. Then $t^{\ell}a$ is in the image of $r$ when $\ell$ is
large, i.e. there exists a polynomial $P(\xi_1,\dots,\xi_d)$ such
that
\[
  r\bigl(P(\xi_1,\dots,\xi_d)\bigr)=t^{\ell}a.
\]
On the other hand,
\[
  r\bigl(P(\xi_1,\dots,\xi_d)\bigr)=P(\alpha_1t+w_1,\dots,\alpha_dt+w_d)=
  \sum_{k\ge0}t^kQ_k(w_1,\dots,w_d)
\]
for some polynomials $Q_k$. Therefore, $a=Q_\ell(w_1,\dots,w_d)$,
the restriction map $H^*(M;\Z_p)\to H^*(N';\Z_p)$ is surjective,
and $H^*(N';\Z_p)$ is generated by the degree-two elements
$w_1,\dots,w_d$.

Now we can repeat the same argument for $N'$ with the induced
action of $S/G$, which is again a circle. We conclude that the
restriction map $H^*(M;\Z_p)\to H^*(N';\Z_p)$ is surjective for
any $p$-subgroup $K$ of~$S$ and any connected component $N'$
of~$M^K$. Now, if the order of $K$ is sufficiently large, then
$M^K=M^S=M^H$ and hence $N'=N$. It follows that the restriction
map $H^*(M;\Z_p)\to H^*(N;\Z_p)$ is surjective for any connected
component $N$ of~$M^H$ and for arbitrary prime~$p$.
\end{proof}

\subsection*{Characteristic submanifolds}
From now on we assume that $M$ is a torus manifold, i.e. $M$ is
closed, connected, smooth and orientatble, $\dim M=2\dim T=2n$,
the $T$-action is smooth and effective and $M^T\ne\varnothing$.

A closed connected codimension-two submanifold of $M$ is called
\emph{characteristic}\label{charm} if it is fixed pointwise by a
circle subgroup $S\subset T$. Since $M^T\ne\varnothing$, there
exists at least one characteristic submanifold (an exercise).
Furthermore, any fixed point is contained in an intersection of
some $n$ characteristic submanifolds.

There are finitely many characteristic submanifolds in $M$, and we
shall denote them by $M_i$, \ $1\le i\le m$. The intersection of
any $k\le n$ characteristic submanifolds is either empty or a
(possibly disconnected) submanifold of codimension $2k$ fixed
pointwise by a $k$-dimensional subtorus of~$T$. In particular, the
intersection of any $n$ characteristic submanifolds consists of
finitely many $T$-fixed points.

Since $M$ is orientable, each $M_i$ is also orientable (since it
is fixed by a circle action). We say that $M$ is
\emph{omnioriented}\label{omnitorma} if an orientation is
specified for $M$ and for each characteristic submanifold $M_i$.
There are $2^{m+1}$ choices of omniorientations on~$M$. In what
follows we shall always assume $M$ to be omnioriented. This allows
us to view the circle fixing $M_i$ as an element in the integer
lattice $\Hom(\mathbb S,T)\cong\Z^n$.

Since both $M$ and $M_i$ are oriented, the equivariant Gysin
homomorphism\label{gysin1} $H^*_T(M_i)\to H^{*+2}_T(M)$ is defined
(see Section~\ref{gractions}). Denote by $\tau_i\in H^2_T(M)$ the
image of $1\in H^0_T(M_i)$ under this homomorphism. The
restriction of $\tau_i$ to $H^2_T(M_i)$ is the equivariant Euler
class $e^T(\nu_i)$ of the normal bundle $\nu_i=\nu(M_i\subset M)$.

\begin{proposition}\label{prop:a_i}
Let $\lambda_i\in H_2(BT)$ be the element corresponding to the
circle subgroup fixing $M_i$ via the identification
$H_2(BT)=\Hom(\mathbb S,T)$.

\begin{itemize}
\item[(a)] Let $v$ be a $T$-fixed point of $M$, and assume that $v$
is contained in the intersection $M_{i_1}\cap\cdots\cap M_{i_n}$.
Then the corresponding elements
$\lambda_{i_1},\ldots,\lambda_{i_n}$ form a basis of
$H_2(BT)\cong\Z^n$. This basis is conjugate to the set of weights
$\mb w_1(v),\ldots,\mb w_n(v)$ of the tangential
$T$-representation at~$v$.

%\item[(a)] If $M_{i_1}\cap\cdots\cap M_{i_k}\ne\varnothing$, then the elements
%$\lambda_{i_1},\ldots,\lambda_{i_k}$ form a part of basis of
%$H_2(BT)\cong\Z^n$.

\item[(b)] Let $\rho\colon ET\times_T M\to BT$ be the projection. For any $t\in H^2(BT)$,
\[
  \rho^*(t)=\sum_{i=1}^m\langle t,\lambda_i\rangle\tau_i\quad
  \text{\rm modulo $H^*(BT)$-torsion.}
\]
\end{itemize}
\end{proposition}
\begin{proof}
Since the action is effective and $M$ is connected, the tangential
$T$-representation at $v$ is faithful, so the set of weights
$\{\mb w_1(v),\ldots,\mb w_n(v)\}$ is a basis of $\Hom(T,\mathbb
S)=H^2(BT)$. The proof of~(a) is the same as that of
Proposition~\ref{qtw}.

By Theorem~\ref{locth}, the restriction map~\eqref{rmevc} is an
isomorphism after localisation at~$H^+(BT)$. Therefore, to
prove~(b), we can apply the map $r$ to the both sides of the
identity in question and verify the resulting identity in
$H^*(BT)\otimes H^*(M^T)$. We may write $r=\bigoplus_v r_v$, where
$r_v\colon H^*_T(M)\to H^*(BT)$ is the map induced by the
inclusion of the fixed point $v\to M$. We have $r_v(\rho^*(t))=t$
and
\[
  r_v\Bigl(\sum_{i=1}^m\langle
  t,\lambda_i\rangle\tau_i\Bigr)=\sum_{k=1}^n\langle
  t,\lambda_{i_k}\rangle e^T(\nu_{i_k})|_v=\sum_{k=1}^n\langle
  t,\lambda_{i_k}\rangle\mb w_k(v)=t,
\]
where the last identity holds because
$\lambda_{i_1},\ldots,\lambda_{i_n}$ and $\mb w_1(v),\ldots,\mb
w_n(v)$ are conjugate bases.
\end{proof}

Here is an example of a locally standard torus manifold, which is
not quasitoric:

\begin{example}\label{2nsphere}
Consider the unit $2n$-sphere in~$\C^n\times\R$:
$$
  S^{2n}=\bigl\{ (z_1,\ldots,z_n,y)\in\C^n\times\R\colon|z_1|^2+\cdots+|z_n|^2+y^2=1
  \bigr\}.
$$
Define a $T$-action by the formula
$$
  (t_1,\ldots,t_n)\cdot(z_1,\ldots,z_n,y)=(t_1z_1,\ldots,t_nz_n,y).
$$
There are two fixed points $(0,\ldots,0,\pm 1)$ and $n$
characteristic submanifolds given by
$\{z_1=0\},\ldots,\{{z_n=0}\}$. The intersection of any $k$
characteristic submanifolds is connected if $k\le n-1$, and
consists of two disjoint fixed points if $k=n>1$.
\end{example}

Unlike quasitoric manifolds, some intersections of characteristic
submanifolds are disconnected in the example above. The cohomology
ring a quasitoric manifold is generated in degree two
(Theorem~\ref{qtcoh}). The next lemma shows that this is exactly
the condition that guarantees the connectedness of intersections
of characteristic submanifolds.

\begin{lemma}\label{cmcon}
Suppose that $H^*(M)$ is generated in degree two. Then all
nonempty multiple intersections of characteristic submanifolds are
connected and have cohomology generated in degree two.
\end{lemma}
\begin{proof}
By Lemma~\ref{sconn}, the cohomology $H^*(M_i)$ is generated by
the degree-two part and the restriction map $H^*(M)\to H^*(M_i)$
is surjective for any characteristic submanifold~$M_i$. Then it
follows from Lemma~\ref{lemm:freeness} that the restriction map
$H^*_T(M)\to H^*_T(M_i)$ in equivariant cohomology is also
surjective.

Now we prove that multiple intersections are connected.  Suppose
that $M_{i_1}\cap\dots\cap M_{i_k}\ne\varnothing$, \ $1<k\le n$,
and let $N$ be a connected component of this intersection. Since
$N$ is fixed by a subtorus, it contains a $T$-fixed point by
Lemmata~\ref{oddfp} and~\ref{lemm:odd=0}. For each
$i\in\{i_1,\dots,i_k\}$, there are embeddings $\varphi_i\colon
N\to M_i$, $\psi_i\colon M_i\to M$, and the corresponding
equivariant Gysin homomorphisms:
\[
\begin{CD}
  H^0_T(N) @>\varphi_{i_!}>> H^{2k-2}_T(M_i) @>\psi_{i_!}>> H^{2k}_T(M).
\end{CD}
\]
The map $\psi_i^*\colon H^*_T(M)\to H^*_T(M_i)$ is surjective, so
we obtain $\varphi_{i_!}(1)=\psi^*_i(u)$ for some $u\in
H^{2k-2}_T(M)$. Now we calculate
\[
  (\psi_i\circ\varphi_i)_!(1)=\psi_{i_!}(\varphi_{i_!}(1))=
  \psi_{i_!}\bigl(\psi_i^*(u)\bigr)=\psi_{i_!}(1)u=\tau_i u.
\]
Hence $(\psi_i\circ\varphi_i)_!(1)$ is divisible by $\tau_i$, for
each $i\in\{i_1,\dots,i_k\}$. By \cite[Proposition~3.4]{masu99}
(see also Theorem~\ref{th:eqcoh} below), the degree-$2k$ part of
$H^*_T(M)$ is additively generated by monomials
$\tau^{k_1}_{j_1}\cdots \tau^{k_p}_{j_p}$ such that
$M_{j_1}\cap\dots\cap M_{j_p}\not=\varnothing$ and
$k_1+\dots+k_p=k$. It follows that $(\psi_i\circ\varphi_i)_!(1)$
is a nonzero integral multiple of $\tau_{i_1}\cdots \tau_{i_k}\in
H^{2k}_T(M)$. By the definition of Gysin homomorphism,
$(\psi_i\circ\varphi_i)_!(1)$ maps to zero under the restriction
map $H_T^*(M)\to H_T^*(x)$ for any point $x\in (M\backslash N)^T$.
On the other hand, the image of $\tau_{i_1}\dots \tau_{i_k}$ under
the map $H_T^*(M)\to H_T^*(x)$ is nonzero for any $T$-fixed point
$x\in M_{i_1}\cap\dots\cap M_{i_k}$. Thus, $N$ is the only
connected component of the latter intersection. The fact that
$H^*(N)$ is generated by its degree-two part follows from
Lemma~\ref{sconn}.
\end{proof}

\subsection*{Orbit quotients and manifolds with corners}
Let $Q=M/T$ be the orbit space of a locally standard $T$-manifold
$M$, and let $\pi\colon M\to Q$ be the quotient projection. Then
$Q$ is a manifold with corners (see Definition~\ref{mawco}).
%; with atlas consisting of quotients of locally standard charts.
The facets of $Q$ are the projections of the characteristic
submanifolds: $\F_i=\pi(M_i)$, \ $1\le i\le m$. Faces of $Q$ are
connected components of intersections of facets (note that the
these intersections may be disconnected, as in
Example~\ref{2nsphere})\label{facesmwc}. For convenience, we
regard $Q$ itself as a face; all other faces are
called~\emph{proper}.

If $H^{odd}(M)=0$, then each face has a vertex by
Lemmata~\ref{lemm:odd=0} and~\ref{oddfp}. Moreover, if $H^*(M)$ is
generated in degree two, then all intersections of facets are
connected by Lemma~\ref{cmcon}.

\begin{proposition}\label{oslsa}
The orbit space $Q$ of a locally standard $T$-manifold is a nice
manifold with corners.
\end{proposition}
\begin{proof}
We need to show that any face $G$ of codimension~$k$ in~$Q$ is an
intersection of exactly~$k$ facets. Let $q$ be a point in the
interior of~$G$, and let $x\in\pi^{-1}(q)$. Let $z_1,\ldots,z_n$
be the coordinates in a locally standard chart containing~$x$. The
point $x$ has exactly $k$ of these coordinates vanishing, assume
that these are $z_1,\ldots,z_k$. Then, for any $i=1,\ldots,k$, the
equation $z_i=0$ specifies the tangent space at~$x$ to a
characteristic submanifold of~$M$. Therefore, $x$ is contained in
exactly $k$ characteristic submanifolds. Each characteristic
submanifold is projected onto a facet of~$Q$, so that $G$ is
contained in exactly $k$ facets.
\end{proof}

\begin{theorem}%[\cite{ma-pa06}]
\label{theo:local standardness}
A torus manifold $M$ with $H^{odd}(M)=0$ is locally standard.
\end{theorem}
\begin{proof}
We first show that there are no nontrivial finite stabilisers for
the $T$-action on~$M$. Assume the opposite, i.e. there is a point
$x\in M$ with finite nontrivial stabiliser~$T_x$. Then $T_x$
contains a nontrivial cyclic subgroup $G$ of prime order~$p$. Let
$N$ be the connected component of $M^G$ containing $x$.  Since $N$
contains $x$ and $T_x$ is finite, the principal (i.e. the
smallest) stabiliser of the induced $T$-action on $N$ is finite.
As in the proof of Lemma~\ref{lemm:odd=0}, it follows
from~\cite[Theorem~VII.2.2]{bred72} that $H^{odd}(N;\Z_p)=0$. In
particular, the Euler characteristic of $N$ is non-zero, hence $N$
has a $T$-fixed point, say $y$. The tangential $T$-representation
$\mathcal T_y M$ at $y$ is faithful, $\dim M=2\dim T$ and
$\mathcal T_y N$ is a proper $T$-subrepresentation of $\mathcal
T_yM$. It follows that there is a nontrivial subtorus $T'\subset
T$ which fixes $\mathcal T_yN$ and does not fix the complement of
$\mathcal T_yN$ in $\mathcal T_yM$. Then $T'$ is the principal
stabiliser of $N$, which contradicts the above observation that
the principal stabiliser of $N$ is finite.

If the stabiliser $T_x$ is trivial, $M$ is obviously locally
standard near~$x$. Suppose that $T_x$ is non-trivial. Then it
cannot be finite, i.e. $\dim T_x>0$.  Let $H$ be the identity
component of~$T_x$, and $N$ the connected component of $M^{H}$
containing~$x$. By Lemmata~\ref{lemm:odd=0} and~\ref{oddfp}, $N$
has a $T$-fixed point, say $y$. Looking at the tangential
representation at $y$, we observe that the induced action of $T/H$
on $N$ is effective. By the previous argument, no point of $N$ has
a nontrivial finite stabiliser for the induced action of $T/H$,
which implies that $T_x=H$. Now $x$ and $y$ are both in the same
connected submanifold $N$ fixed pointwise by $T_x$, hence the
$T_x$-representation $\mathcal T_x M$ agrees with the restriction
of the tangential $T$-representation $\mathcal T_yM$ to $T_x$.
This implies that $M$ is locally standard near~$x$.
\end{proof}

%From now on we assume that $M$ is a locally standard torus manifold.

Recall that a space $X$ is \emph{acyclic}\label{acysp} if
$\widetilde H_i(X)=0$ for any~$i$.

\begin{definition}\label{deffaceacy}
We say that a manifold with corners $Q$ is \emph{face-acyclic} if
$Q$ and all its faces are acyclic.  We call $Q$ a \emph{homology
polytope} if it is face-acyclic and all nonempty multiple
intersections of facets are connected.
\end{definition}

A simple polytope is a homology polytope. Here is an example that
does not arise in this way:

\begin{example}\label{S2n orbit}
The torus manifold $S^{2n}$ with the $T$-action from
Example~\ref{2nsphere} is locally standard, and the map
\[
  (z_1,\dots,z_n,y)\to (|z_1|,\dots,|z_n|,y)
\]
induces a face preserving homeomorphism from the orbit space
$S^{2n}/T$ to the space
\[
  \{ (x_1,\dots,x_n,y)\in \R^{n+1}\colon x_1^2+\dots+x_n^2+y^2=1,\ x_1\ge 0,
  \dots,x_n\ge 0\}.
\]
This manifold with corners is face-acyclic, but is not a homology
polytope if $n>1$.
\end{example}

Proposition~\ref{prop:a_i} allows us to define a characteristic
map for locally standard torus manifolds:
\begin{equation}\label{eqn:characteristic map}
\begin{aligned}
  \lambda\colon \{F_1,\ldots,F_m\} &\to H_2(BT)=\Hom(\mathbb S,T)\cong
  \Z^n,\\
  F_i &\mapsto \lambda_i.
\end{aligned}
\end{equation}
%The map $\lambda$ удовлетворяет условию неособости, аналогичному
%соответствующему условию для квазиторических многообразий:
%\begin{quote}
%если $F_{i_1}\cap\ldots\cap F_{i_k}\ne\varnothing$, то
%$\lambda(F_{i_1}),\ldots,\lambda(F_{i_k})$ задают часть базиса
%решётки $\Hom(\mathbb S^1,T)\cong\Z^n$.
%\end{quote}
Given a point $q\in Q$, consider the smallest face $G(q)$
containing~$q$. This face is a connected component of an
intersection of facets $F_{i_1}\cap\cdots\cap F_{i_k}$. We define
the subtorus $T(q)\subset T$ generated by the circle subgroups
corresponding to $\lambda(F_{i_1}),\ldots,\lambda(F_{i_k})$, and
the identification space
\begin{equation}\label{idspa}
  M(Q,\lambda)=Q\times T/{\sim}\quad\text{where }
  (x,t_1)\sim(x,t_2)\:\text{ if }\:
  t_1^{-1}t_2\in T(q).
\end{equation}
It is easy to see that $M(Q,\lambda)$ is a closed manifold with a
$T$-action. Here is a straightforward generalisation
of~\cite[Lemma~1.4]{da-ja91} (see Proposition~\ref{equivar}):

\begin{proposition}\label{acydj}
Let $M$ be a locally standard torus manifold with orbit space $Q$,
and let $\lambda$ be the map defined by~\eqref{eqn:characteristic
map}. If $Q$ is face-acyclic, then there is a weakly
$T$-equivariant homeomorphism
\[
  M(Q,\lambda)\to M
\]
covering the identity on~$Q$.
\end{proposition}

\begin{remark}
Instead of face-acyclicity, one can only require that the second
cohomology group of each face of $Q$ vanishes, as this is the
condition implying the triviality of torus principal bundles
obtained after blowing up the singular strata of lower dimension.
\end{remark}

\subsection*{Face rings of manifolds with corners}
Let $Q$ be a nice manifold with corners. The set of faces of $Q$
containing a given face is isomorphic to the poset of faces of a
simplex. In other words, the faces of~$Q$ form a simplicial poset
$\mathcal S$ (see Definition~\ref{defsp}) with respect to the
reverse inclusion. The initial element $\hatzero$ of this poset
is~$Q$. We refer to this simplicial poset $\mathcal S$ as the
\emph{dual}\label{dualsp} of~$Q$; this duality extends the
combinatorial duality between simple polytopes and their boundary
sphere triangulations. The dual poset $\mathcal S$ is the poset of
faces of a simplicial complex $\mathcal K$ if and only if all
nonempty multiple intersections of facets of $Q$ are connected. In
this case, $\mathcal K$ is the nerve of the covering of $\partial
Q$ by facets.

\begin{example}\label{3mwc}
Consider the three structures of a manifold with corners on a disc
$D^2$, shown in Fig.~\ref{mwc2d}. The manifold with corners shown
on the left is not nice. The middle one is nice and is
face-acyclic, but is not a homology polytope. The right one is
homeomorphic to a 2-simplex, so it is a homology polytope. Compare
this with simplicial posets from Example~\ref{exsimpos}.
\begin{figure}[h]
\begin{center}
\begin{picture}(80,25)
\put(10,15){\circle{14}} \put(17,15){\circle*{1.5}} \put(9,1){(1)}
\put(40,15){\circle{14}} \put(40,8){\circle*{1.5}}
\put(40,22){\circle*{1.5}} \put(39,1){(2)}
\put(70,15){\circle{14}} \put(75,10){\circle*{1.5}}
\put(75,20){\circle*{1.5}} \put(63,15){\circle*{1.5}}
\put(69,1){(3)}
\end{picture}
\caption{2-disc as a manifold with corners.} \label{mwc2d}
\end{center}
\end{figure}
\end{example}

\begin{example}
Let $Q=S^{2n}/T$ be the orbit space of the torus manifold
$S^{2n}$, see Examples~\ref{2nsphere} and~\ref{S2n orbit} (the
case $n=2$ is shown in Fig.~\ref{mwc2d}~(2)). Here we have $n$
facets, the intersection of any $k$ facets is connected if $k\le
n-1$, but the intersection of $n$ facets consists of two points.
The dual simplicial cell complex is obtained by gluing two
$(n-1)$-simplices along their boundaries.
\end{example}

We can define the face ring of the orbit space $Q$ as the face
ring of the dual simplicial poset (see Definition~\ref{frspo}).
However, for reader's convenience, we give the definition and
state the main properties directly in terms of the combinatorics
of faces of~$Q$. The proofs of the statements in this subsection
are obtained by obvious dualisation of the corresponding
statements in Section~\ref{frsimpo}.

The intersection of two faces $G,H$ in a manifold with corners can
be disconnected. We consider $G\cap H$ as a set of its connected
components and use the notation $E\in G\cap H$ for connected
components $E$ of this intersection. If $G\cap H\ne\varnothing$,
then there exists a unique minimal face $G\vee H$ containing both
$G$ and~$H$.

\begin{definition}\label{qfrgen}
The \emph{face ring} of a nice manifold with corners $Q$ is the
quotient
$$
  \Z[Q]=\Z[\va{G}\colon G\text{ a face}]/\mathcal I_Q,
$$
where $\mathcal I_Q$ is the ideal generated by $\va{Q}-1$ and all
elements of the form
$$
  \va{G}\va{H}-
    \va{G\vee H}\cdot\!\!\!\sum_{E\in{G\cap H}}\!\!\!\va{E}.
$$
In particular, if $G\cap H=\varnothing$, then $\va{G}\va{H}=0$ in
$\Z[Q]$.

The grading is given by $\deg \va{G}=2\codim G$.
\end{definition}

\begin{example}\label{4sphere1}
Consider the torus action on $S^{2n}$ from Examples~\ref{2nsphere}
and~\ref{S2n orbit} and let $n=2$. Then $Q$ is a 2-disc with two
0-faces, say $p$ and $q$, and two 1-faces, say $G$ and~$H$. Then
\[
  \Z[Q]=\Z[\va G,\va H,v_p,v_q]/(\va G\va H=v_p+v_q,\ v_pv_q=0),
\]
where $\deg\va G=\deg\va H=2$, $\deg v_p=\deg v_q=4$. This is the
same ring as the one described in Example~\ref{exsp}.1, but
written in the dual notation.
\end{example}

Here is a dualisation of Theorem~\ref{chdec}:

\begin{theorem}\label{qorder}
Any element $a\in\Z[Q]$ can be written uniquely as an integral
linear combination of monomials
$\va{G_1}^{i_1}\va{G_2}^{i_2}\cdots\va{G_n}^{i_n}$ corresponding
to chains of faces $G_1\supset G_2\supset\cdots\supset G_n$
of~$Q$.
\end{theorem}

Given a vertex $v\in Q$, define the \emph{restriction
map}\label{rmapsp}
$$
  s_v\colon\Z[Q]\to\Z[Q]/(\va{G}\colon G\not\ni v).
$$
The ring $\Z[Q]/(\va{G}\colon G\not\ni v)$ is identified with the
polynomial ring $\Z[\va{\F_{i_1}},\ldots,\va{\F_{i_n}}]$ on $n$
degree-two generators corresponding to the facets
$\F_{i_1},\ldots,\F_{i_n}$ containing~$v$.

\begin{theorem}\label{qalres}
Assume that each face of $Q$ has a vertex. Then the sum
$s=\bigoplus_vs_v$ of the restriction maps over all vertices $v\in
Q$ is a monomorphism from $\Z[Q]$ to a direct sum of polynomial
rings.
\end{theorem}

\subsection*{Equivariant cohomology}
Here we construct a natural homomorphism from the face ring
$\Z[Q]$ to the equivariant cohomology ring $H^*_T(M)$ of a locally
standard torus manifold modulo $H^*(BT)$-torsion. Then we obtain
conditions under which this homomorphism is monic and epic; in
particular, we show that $\Z[Q]\to H^*_T(M)$ is an isomorphism
when $H^{odd}(M)=0$.

Since the fixed point set $M^T$ is finite, the restriction
map~\eqref{rmevc} defines a map
\begin{equation}\label{eqn:restriction map}
  r=\bigoplus\limits_{v\in M^T} r_v\colon H^*_T(M)\to H^*_T(M^T)=
  \bigoplus\limits_{v\in M^T}H^*(BT)
\end{equation}
to a direct sum of polynomial rings. Its kernel is a
$H^*(BT)$-torsion by Corollary~\ref{kerres}, and $r$ is a
monomorphism when $H^{odd}(M)=0$.

The $1$-skeleton of $Q$ is an $n$-valent graph. We identify $M^T$
with the vertices of $Q$ and denote by $E(Q)$ the set of oriented
edges.  Given an element $e\in E(Q)$, denote the initial point and
the terminal point of $e$ by $i(e)$ and $t(e)$, respectively. Then
$M_e=\pi^{-1}(e)$ is a 2-sphere fixed pointwise by a
codimension-one subtorus in $T$, and it contains two $T$-fixed
points $i(e)$ and~$t(e)$. The 2-dimensional subspace $\mathcal
T_{i(e)}M_e\subset \mathcal T_{i(e)}M$ is an irreducible component
of the tangential $T$-representation $\mathcal T_{i(e)}M$. The
same is true for the other point $t(e)$, and the
$T$-representations $\mathcal T_{i(e)}M$ and $\mathcal T_{t(e)}M$
are isomorphic. There is a unique characteristic submanifold, say
$M_i$, intersecting $M_e$ at $i(e)$ transversely. The
omniorientation defines an orientation for the normal bundle
$\nu_i$ of $M_i$ and, therefore, an orientation of~$\mathcal
T_{i(e)}M_e$. We therefore can view the representation $\mathcal
T_{i(e)}M_e$ as an element of $\Hom(T,\mathbb S)=H^2(BT)$, and
denote this element by~$\alpha(e)$.

Let $e^T(\nu_i)\in H^2_T(M_i)$ be the Euler class of the normal
bundle, and denote its restriction to a fixed point $v\in M_i^T$
by $e^T(\nu_i)|_v\in H^2_T(v)=H^2(BT)$. Then
\begin{equation} \label{eqn:euler}
  e^T(\nu_i)|_v=\alpha(e),
\end{equation}
where $e$ is the unique edge such that $i(e)=v$ and $e\notin
\F_i=\pi(M_i)$. Using the terminology of~\cite{gu-za01}, we refer
to the map
\[
  \alpha\colon E(Q)\to H^2(BT),\quad e\mapsto\alpha(e),
\]
as an \emph{axial function}.

\begin{lemma}\label{torusaf}
The axial function $\alpha$ has the following properties:
\begin{itemize}
\item[(a)] $\alpha(\bar e)=\pm \alpha(e)$ for any $e\in E(Q)$, where $\bar e$
denotes $e$ with the opposite orientation;
\item[(b)] for any vertex $v$, the set $\alpha_v=\{
\alpha(e)\colon i(e)=v\}$ is a basis of $H^2(BT)$;
\item[(c)] for $e\in E(Q)$, we have
$\alpha_{i(e)}= \alpha_{t(e)}\mod \alpha(e)$.
\end{itemize}
\end{lemma}
\begin{proof}
Property (a) follows from the fact that $\mathcal T_{i(e)}M_e$ and
$\mathcal T_{t(e)}M_e$ are isomorphic as real $T$-representations,
and~(b) holds since the $T$-representation $\mathcal T_{i(e)}M$ is
faithful. Let $T_e$ be the codimension-one subtorus fixing~$M_e$.
Then the $T_e$-representations $\mathcal T_{i(e)}M$ and $\mathcal
T_{t(e)}M$ are isomorphic, since the points $i(e)$ and $t(e)$ are
contained in the same connected component $M_e$ of~$M^{T_e}$. This
implies~(c).
\end{proof}

\begin{remark}
The original definition of an axial function in~\cite{gu-za01}
requires the property $\alpha(\bar e)=-\alpha(e)$, but we allow
$\alpha(\bar e)=\alpha(e)$. For example, $\alpha(\bar
e)=\alpha(e)$ for the $T^2$-action on $S^{4}$ from
Example~\ref{2nsphere}.
\end{remark}

\begin{lemma} \label{lemm:necessity}
Given $\eta\in H^*_T(M)$ and $e\in E(Q)$, the difference
$r_{i(e)}(\eta)-r_{t(e)}(\eta)$ is divisible by $\alpha(e)$.
\end{lemma}
\begin{proof}
Consider the commutative diagram of restrictions
\[
\begin{CD}
  H^*_T(M) @>>> H^*_T(i(e))\oplus H^*_T(t(e))= @. H^*(BT)\oplus H^*(BT)\\
  @VVV          @VVV \\
  H^*_{T_e}(M_e) @>>> H^*_{T_e}(i(e))\oplus H^*_{T_e}(t(e))= @.\:
  H^*(BT_e)\oplus H^*(BT_e)
\end{CD}
\]
Since $H^*_{T_e}(M_e)=H^*(BT_e)\otimes H^*(M_e)$, the two
components of the image of $\eta$ in $H^*(BT_e)\oplus H^*(BT_e)$
above coincide. Then it follows from the commutativity of the
diagram that the restrictions of $r_{i(e)}(\eta)$ and
$r_{t(e)}(\eta)$ to $H^*(BT_e)$ coincide. Now the result follows
from the fact that the kernel of the restriction map $H^*(BT)\to
H^*(BT_e)$ is the ideal generated by~$\alpha(e)$.
\end{proof}

The preimage $M_G=\pi^{-1}(G)$ of a codimension-$k$ face $G\subset
Q$ is a closed $T$-invariant submanifold of~$M$. It is a connected
component of an intersection of $k$ characteristic submanifolds.
The omniorientation defines an orientation of~$M_G$ and the
equivariant Gysin homomorphism $H_T^0(M_G)\to H_T^{2k}(M)$. Let
$\ta{G}$ denote the image of~1 under this homomorphism; it is
called the \emph{Thom class}\label{thomclasstm} of~$M_G$. The
restriction of $\ta{G}\in H_T^{2k}(M)$ to $H_T^{2k}(M_G)$ is the
equivariant Euler class of the normal bundle $\nu(M_G\subset M)$,
and $r_v(\ta{G})=0$ for $v\notin(M_G)^T$. It follows from
(\ref{eqn:euler}) that
\begin{equation} \label{eqn:tau_G}
  r_v(\ta{G})=\left\{%
  \begin{array}{ll}
    \prod\limits_{i(e)=v,\ e\not\subset G}\alpha(e)
       & \hbox{ \ if \ $v\in (M_G)^T$;}\\[2mm]
    \qquad 0 & \hbox{ \ otherwise.}\\
  \end{array}%
  \right.
\end{equation}

Define the quotient ring
\[
  \widehat{H}^*_T(M)=H^*_T(M)/H^*(BT)\text{-torsion}.
\]
The restriction map $r$ from \eqref{eqn:restriction map} induces a
monomorphism $\widehat{H}^*_T(M)\to H^*_T(M^T)$, which we continue
to denote by~$r$. In particular, $\ta{G}=0$ in $\widehat H^*_T(M)$
if $M_G$ has no $T$-fixed points. The next lemma shows that the
face ring relations from Definition~\ref{qfrgen} hold in
$\widehat{H}^*_T(M)$ with $\va{G}$ replaced by $\ta{G}$.

\begin{lemma}\label{taurel}
For any two faces $G$ and $H$ of $Q$, the relation
$$
  \ta{G}\ta{H}=
    \ta{G\vee H}\cdot\!\!\!\sum_{E\in{G\cap H}}\!\!\!\ta{E}
$$
holds in $\widehat{H}^*_T(M)$.
\end{lemma}
\begin{proof}
Since the map $r\colon\!\widehat{H}^*_T(M)\to H^*_T(M^T)$ is
injective, it suffices to show that $r_v$ maps both sides of the
identity to the same element, for any $v\in M^T$.

Let $v\in M^T$. Given a face $G\ni v$, set
\[
  N_v(G)= \{ e\in E(Q) \colon i(e)=v,\ e\not\subset G\},
\]
which may be thought of as the set of directions transverse to $G$
at~$v$. Then we can rewrite~\eqref{eqn:tau_G} as follows:
\begin{equation} \label{eqn:tau_G2}
  r_v(\ta{G})=\prod_{e\in N_v(G)}\alpha(e),
\end{equation}
where the right hand side is understood to be 1 if
$N_v(G)=\varnothing$ and to be 0 if $v\notin G$. Assume that
$v\notin G\cap H$; then either $v\notin G$ or $v\notin H$, and
$v\notin E$ for any $E\in G\cap H$. Hence both sides of the
identity from the lemma map to zero by~$r_v$. Assume that $v\in
G\cap H$; then
\[
  N_v(G)\cup N_v(H)=N_v(G\vee H)\cup N_v(E),
\]
where $E$ is the connected component of the intersection $G\cap H$
containing~$v$, and $v\notin E'$ for any other $E'\in G\cap H$.
This together with (\ref{eqn:tau_G2}) implies that both sides of
the identity map to the same element by~$r_v$.
\end{proof}

Lemma~\ref{taurel} implies that the map
\begin{align*}
  \Z[\va{G}\colon G \text{ a face}]&\to H^*_T(M),\\
  \va{G}&\mapsto\ta{G}
\end{align*}
induces a homomorphism
\begin{equation}\label{eqn:varphi2}
  \varphi\colon \Z[Q]\to \widehat{H}_T^*(M).
\end{equation}

\begin{lemma}\label{lemm:monof}
The homomorphism $\varphi$ is injective if any face of $Q$ has a
vertex.
\end{lemma}
\begin{proof}
We have $s=r\circ\varphi$, where $s$ is the algebraic restriction
map from Lemma~\ref{qalres}. Since $s$ is injective, $\varphi$ is
also injective.
\end{proof}

The next theorem, which is a particular case of one of the main
results of~\cite{g-k-m98}, says that when $H^{odd}(M)=0$, the
condition from Lemma~\ref{lemm:necessity} specifies precisely the
image of the equivariant cohomology under the restriction map:

\begin{theorem}[{\cite{g-k-m98}, see also \cite[Chapter~11]{gu-st99}}]
\label{theo:GKM} Let $M$ be a torus manifold with $H^{odd}(M)=0$.
Assume given an element $\eta_v\in H^*(BT)$ for each $v\in M^T$.
Then $\{\eta_v\}\in \bigoplus_{v\in M^T}H^*(BT)$ belongs to the
image of the restriction map $r$ in~\eqref{eqn:restriction map} if
and only if $\eta_{i(e)}-\eta_{t(e)}$ is divisible by $\alpha(e)$
for any $e\in E(Q)$.
\end{theorem}

\begin{corollary}\label{1skel}
If $H^{odd}(M)=0$, then the $1$-skeleton of any face of $Q$
(including $Q$ itself) is connected.
\end{corollary}
\begin{proof}
By Theorem~\ref{theo:GKM}, $\{\eta_v\}\in \bigoplus_{v\in
M^T}H^0(BT)$ belongs to $r(H^0_T(M))$ if $\eta_v$ is a locally
constant function on the 1-skeleton of~$Q$. On the other hand,
since $M$ is connected, the image $r(H^0_T(M))$ is isomorphic
to~$\Z$. Hence the 1-skeleton of $Q$ is connected. Similarly, the
1-skeleton of any face $G$ of $Q$ is connected, because
$M_G=\pi^{-1}(G)$ is also a torus manifold with $H^{odd}(M_G)=0$
(see Lemma~\ref{lemm:odd=0}).
\end{proof}

\begin{remark}
Connectedness of $1$-skeletons of faces of $Q$ can be proven
without referring to Theorem~\ref{theo:GKM}, see the remark after
Theorem~\ref{cohfa}.
\end{remark}

For a face $G\subset Q$, we denote by $I(G)$ the ideal in
$H^*(BT)$ generated by all elements $\alpha(e)$ with $e\in G$.

\begin{lemma} \label{lemm:I(G)}
Suppose that the $1$-skeleton of a face $G$ is connected. Given
$\eta\in H^*_T(M)$, if $r_v(\eta)\notin I(G)$ for some vertex
$v\in G$, then $r_{w}(\eta)\notin I(G)$ for any vertex $w\in G$.
\end{lemma}
\begin{proof}
Suppose $r_{w}(\eta)\in I(G)$ for some vertex $w\in G$. Then
$r_u(\eta)\in I(G)$ for any vertex $u\in G$ joined with $w$ by an
edge $f\subset G$, because $r_w(\eta)-r_u(\eta)$ is divisible by
$\alpha(f)$ by Lemma~\ref{lemm:necessity}. Since the $1$-skeleton
of $G$ is connected, $r_w(\eta)\in I(G)$ for any vertex $w\in G$,
which contradicts the assumption.
\end{proof}

\begin{lemma}\label{modgn}
If each face of $Q$ has connected 1-skeleton, then $\widehat
H^*_T(M)$ is generated by the elements $\ta{G}$ as an
$H^*(BT)$-module.
\end{lemma}
\begin{proof}
Let $\eta\in H^+_T(M)$ be a nonzero element. Set
$$
  Z(\eta)=\{v\in M^T\colon r_v(\eta)=0\}.
$$
Take $v\notin Z(\eta)$. Then $r_v(\eta)\in H^*(BT)$ is nonzero and
we can express it as a polynomial in $\{\alpha(e)\colon i(e)=v\}$,
as the latter is a basis of $H^2(BT)$. Let
\begin{equation}
\label{monom}
  \prod_{i(e)=v}\alpha(e)^{n_e},\quad n_e\ge0
\end{equation}
be a monomial entering $r_v(\eta)$ with a nonzero coefficient. Let
$G$ be the face spanned by the edges $e$ with $n_e=0$. Then
$r_v(\eta)\notin I(G)$ since $r_v(\eta)$ contains
monomial~\eqref{monom}. Hence $r_w(\eta)\notin I(G)$ (in
particular, $r_w(\eta)\ne0$) for any vertex $w\in G$, by
Lemma~\ref{lemm:I(G)}.

On the other hand, it follows from (\ref{eqn:tau_G}) that
monomial~\eqref{monom} can be written as $r_v(\u{G}\ta{G})$ with
some $\u{G}\in H^*(BT)$. Set $\eta'=\eta-\u{G}\ta{G}\in H^*_T(M)$.
We have $r_w(\ta{G})=0$ for any $w\notin G$, which implies
$r_w(\eta')=r_w(\eta)$ for $w\notin G$.  At the same time,
$r_u(\eta)\not=0$ for $u\in G$ (see above). It follows that
$Z(\eta')\supset Z(\eta)$. The number of monomials in $r_v(\eta')$
is less than that in $r_v(\eta)$. Therefore, by subtracting from
$\eta$ a linear combination of elements $\ta{G}$ with coefficients
in $H^*(BT)$, we obtain an element $\lambda$ such that
$Z(\lambda)$ contains $Z(\eta)$ as a proper subset. By iterating
this procedure, we end up at an element whose restriction to every
vertex is zero. Since the restriction map $r\colon\widehat
H^*_T(M) \to H^*_T(M^T)$ is injective, the result follows.
\end{proof}

\begin{theorem}\label{theo:eqcoh}
Let $M$ be a locally standard torus manifold with orbit space~$Q$.
If each face of $Q$ has connected $1$-skeleton and contains a
vertex, then the mo\-no\-mor\-phism $\varphi\colon \Z[Q]\to
\widehat H^*_T(M)$ of \eqref{eqn:varphi2} is an isomorphism.
\end{theorem}
\begin{proof}
The homomorphism $\varphi$ is injective by Lemma~\ref{lemm:monof}.
To prove that $\varphi$ is surjective it suffices to show that
$\widehat H^*_T(M)$ is generated by the elements $\ta{G}$ as a
ring. By Proposition~\ref{prop:a_i}, the group $\widehat H^2_T(M)$
is generated by the elements $\ta{\F_i}$ corresponding to the
facets~$F_i$. (Note: the notation $\tau_i$ is used for $\ta{\F_i}$
in Proposition~\ref{prop:a_i}.) In particular, any element in
$H^2(BT)\subset \widehat H^*_T(M)$ can be written as a linear
combination of the elements $\ta{F_i}$. Hence any element in
$H^*(BT)$ is a polynomial in $\ta{F_i}$. The rest follows from
Lemma~\ref{modgn}.
\end{proof}

As a corollary we obtain a complete description of the equivariant
cohomology in the case $H^{odd}(M)=0$:

\begin{theorem}[{\cite[Corollary~7.6]{ma-pa06}}]\label{th:eqcoh}
Let $M$ be a locally standard $T$-manifold with $H^{odd}(M)=0$.
Then the equivariant cohomology $H^*_T(M)$ is isomorphic to the
face ring $\Z[Q]$ of the manifold with corners $Q=M/T$.
\end{theorem}
\begin{proof}
Indeed, if $H^{odd}(M)=0$, then $M$ is a torus manifold by
Lemma~\ref{oddfp}. Furthermore, $H^*_T(M)$ is a free
$H^*(BT)$-module by Lemma~\ref{lemm:freeness}, i.e. $\widehat
H^*_T(M)=H^*_T(M)$. The result follows from
Theorem~\ref{theo:eqcoh}.
\end{proof}

The condition $H^{odd}(M)=0$ can be interpreted in terms of the
simplicial poset $\mathcal S$ dual to~$Q$ as follows:

\begin{lemma} \label{theo:face CM}
Let $M$ be a torus manifold with quotient $Q$, and $\mathcal S$ be
the face poset of $Q$. Then $H^{odd}(M)=0$ if and only if the
following conditions are satisfied:
\begin{itemize}
\item[(a)] the ring $H^*_T(M)$ is isomorphic to $\Z[\mathcal S](=\Z[Q])$;
\item[(b)] $\Z[\mathcal S]$ is a Cohen--Macaulay ring.
\end{itemize}
Furthermore, the ring $H^*(M)$ is generated in degree two if and
only if $\mathcal S$ is (the face poset of) a simplicial complex
in addition to the above two conditions.
\end{lemma}
\begin{proof}
If $H^{odd}(M)=0$, then $H^*_T(M)\cong \Z[Q]$ by
Theorem~\ref{th:eqcoh}, and $\Z[\mathcal S]$ is a Cohen--Macaulay
ring by Lemma~\ref{lemm:freeness}.

Now we prove that $H^{odd}(M)=0$ under conditions (a) and~(b). The
composite
\[
  H^*(BT)\stackrel{\rho^*}\longrightarrow  H^*_T(M)
  \stackrel{r}\longrightarrow \bigoplus_{v\in M^T}H^*(BT).
\]
is the diagonal map. By Lemma~\ref{lsopscc}, this implies that
$\rho^*(t_1),\ldots,\rho^*(t_n)$ is an lsop. Since $H^*_T(M)$ is a
Cohen--Macaulay ring, any lsop is a regular sequence
(Proposition~\ref{rlsop}). It follows that $H^*_T(M)$ is a free
$H^*(BT)$-module and hence $H^{odd}(M)=0$, by
Lemma~\ref{lemm:freeness}.

It remains to prove the last statement. Assume that $H^*(M)$ is
generated in degree two. By Lemma~\ref{cmcon}, all non-empty
multiple intersections of facets are connected. Then $\mathcal S$
is the nerve of the covering of $\partial Q$ by facets.

Assume that $\mathcal S$ is a simplicial complex. Then
$\Z[\mathcal S]$ is generated in degree two. Furthermore,
$H^*_T(M)\cong \Z[\mathcal S]$ is a free $H^*(BT)$-module by the
first part of the theorem, whence $H^*(M)$ is a quotient ring of
$H^*_T(M)$. It follows that $H^*(M)$ is also generated by its
degree-two part.
\end{proof}

\subsection*{Ordinary cohomology}
We can now describe the ordinary cohomology of a locally standard
$T$-manifold (or torus manifold) with $H^{odd}(M)=0$. This result
generalises the corresponding statements for toric and quasitoric
manifolds (Theorems~\ref{danjur} and~\ref{qtcoh}):

\begin{theorem}\label{theo:stcoh}
Let $M$ be a locally standard $T$-manifold with $H^{odd}(M)=0$,
and let $Q=M/T$ be the orbit space. Then there is a ring
isomorphism
$$
  H^*(M)\cong \Z[\va{G}\colon G\text{ a face of $Q$}]/\mathcal I,
$$
where $\mathcal I$ is the ideal generated by elements of the
following two types:
\begin{itemize}
\item[(a)]  $\displaystyle{\va{G}\va{H}-\va{G\vee
H}\sum_{E\in{G\cap H}}\va{E}}$;

\item[(b)] $\displaystyle{\sum_{i=1}^m\< t,\lambda_i\> \va{F_i}}$, for
$t\in H^2(BT)$.
\end{itemize}
Here $F_i$ are the facets of $Q$, $i=1,\ldots,m$, and the element
$\lambda_i\in H_2(BT)$ corresponds to the circle subgroup fixing
the characteristic submanifold~$M_i=\pi^{-1}(F_i)$.

The Betti numbers are given by the formula
\[
  \rank H^{2i}(M)=h_i, \qquad 0\le i\le n,
\]
where $h_i$ denote the components of the $h$-vector of the dual
simplicial poset~$\mathcal S$.
\end{theorem}
\begin{proof}
Since the Serre spectral sequence of the bundle $\rho\colon
ET\times_T M\to BT$ collapses at~$E_2$, the map $H^*_T(M)\to
H^*(M)$ is surjective and its kernel is the ideal generated by all
elements $\rho^*(t)$, \ $t\in H^2(BT)$. Therefore, the statement
about the cohomology ring follows from Proposition~\ref{prop:a_i}
and Theorem~\ref{th:eqcoh}.

By Lemma~\ref{lemm:freeness},  $H^*_T(M)\cong H^*(BT)\otimes
H^*(M)$ as $H^*(BT)$-modules, so we have the following formula for
the Poincar\'e series:
\[
  F\bigl(H^*_T(M);\lambda\bigr)=\frac{\sum_{i=0}^n \rank
  H^{2i}(M)\lambda^{2i}}{(1-\lambda^2)^n}.
\]
On the other hand, the Poincar\'e series of the face ring $\Z[Q]$
is given by Theorem~\ref{psgfr}, and the two series coincide by
Theorem~\ref{th:eqcoh}. This implies the statement about the Betti
numbers.
\end{proof}

\begin{example}
The equivariant cohomology ring of the torus manifold $S^4$ from
Example~\ref{4sphere1} is isomorphic to the ring $\Z[Q]$ described
there. The ordinary cohomology ring is obtained by taking quotient
by the ideal generated by $\va{G}$ and~$\va{H}$.
\end{example}

\subsection*{Torus manifolds over homology polytopes}
Using the previous results on the equivariant and ordinary
cohomology of torus manifolds with $H^{odd}(M)=0$, we can now
proceed to describing the relationship between the cohomology of
$M$ and the cohomology of its orbit space~$Q$. Here we prove
Theorem~\ref{hptpe}, which gives a cohomological characterisation
of $T$-manifolds whose orbit spaces are homology polytopes. In the
next subsection we prove that $Q$ is face-acyclic if and only if
$H^{odd}(M)=0$.

\begin{lemma}\label{lemm:H1=0}
If $H^{odd}(M)=0$, then $H^1(Q)=0$.
\end{lemma}
\begin{proof}
We use the Leray spectral sequence of the projection map
$ET\times_T M\to M/T=Q$ onto the second factor.  It has
$E_2^{p,q}=H^p(M/T;\mathcal H^q)$ where $\mathcal H^q$ is the
sheaf with stalk $H^q(BT_x)$ over a point $x\in M/T$, and the
spectral sequence converges to $H^*_T(M)$. Since the $T$-action on
$M$ is locally standard, the isotropy group $T_x$ at $x\in M$ is a
subtorus; so $H^{odd}(BT_x)=0$. Hence $\mathcal H^{odd}=0$, in
particular, $\mathcal H^1=0$. Moreover, $\mathcal H^0=\Z$ (the
constant sheaf). Therefore, we have $E_2^{0,1}=0$ and
$E_2^{1,0}=H^1(M/T)$, whence $H^1(M/T)\cong H^1_T(M)$. On the
other hand, since $H^{odd}(M)=0$ by assumption, $H^*_T(M)$ is a
free $H^*(BT)$-module. Therefore, $H^{odd}_T(M)=0$ by the
universal coefficient theorem. In particular, $H^1_T(M)=0$.
\end{proof}

\begin{lemma}\label{kgore}
If either
\begin{itemize}
\item[(1)] $Q$ is a homology polytope, or
\item[(2)] $H^*(M)$ is generated by its degree-two part,
\end{itemize}
then the dual poset $\mathcal S$ of $Q$ is (the face poset of) a
Gorenstein* simplicial complex.
%In particular, $\Z[\mathcal S]$ is Cohen--Macaulay and the
%geometric realisation $|\mathcal S|$ has homology of an
%$(n-1)$-sphere.
\end{lemma}
\begin{proof}
Under any of the assumptions (1) or (2), all nonempty multiple
intersections of facets of $Q$ are connected, so $\mathcal S$ is
the face poset of the nerve simplicial complex $\mathcal K$ of the
covering of~$\partial Q$. For simplicity, we identify $\mathcal S$
with~$\mathcal K$.

We first prove that $\mathcal S$ is Gorenstein* under
assumption~(1). According to Theorem~\ref{gorencom}, it is enough
to show that the link $\lk\sigma$ of any simplex
$\sigma\in\mathcal K$, has homology of a sphere of dimension $\dim
\lk\sigma=n-2-\dim\sigma$. If $\sigma=\varnothing$ then
$\lk\sigma$ is $\mathcal K$ itself, and it has homology of an
$(n-1)$-sphere, since $Q$ is a homology polytope. If
$\sigma\not=\varnothing$ then $\lk\sigma$ is the nerve of a face
of~$Q$.  Since any face of $Q$ is again a homology polytope,
$\lk\sigma$ has homology of a sphere of dimension $\dim\lk\sigma$.

Now we prove that $\mathcal K$ is Gorenstein* under
assumption~(2). By Exercise~\ref{gorchi}, it is enough to show
that
\begin{itemize}
\item[(a)] $\sK$ is Cohen--Macaulay;
\item[(b)] every $(n-2)$-dimensional simplex is contained in exactly two
$(n-1)$-\-di\-men\-si\-o\-nal simplices;
\item[(c)] $\chi(\sK)=\chi(S^{n-1})$.
\end{itemize}
Condition (a) follows from Lemma~\ref{theo:face CM}. By
definition, every $k$-dimensional simplex of $\sK$ corresponds to
a set of $k+1$ characteristic submanifolds with nonempty
intersection. By Lemma~\ref{cmcon}, the intersection of any $n$
characteristic submanifolds is either empty or consists of a
single $T$-fixed point. This means that $(n-1)$-simplices of $\sK$
are in one-to-one correspondence with $T$-fixed points of~$M$.
Now, each $(n-2)$-simplex of $\sK$ corresponds to a non-empty
intersection of $n-1$ characteristic submanifolds of~$M$. The
latter intersection is connected by Lemma~\ref{cmcon}, so it is a
$2$-sphere. Every $2$-sphere contains exactly two $T$-fixed
points, which implies~(b). Finally,~(c) is just the equation
$h_0=h_n$, which is valid as $h_n=\rank H^{2n}(M)=1$.
\end{proof}

Consider the order complex $\ord(\mathcal S)$ (see
Definition~\ref{ordercom}) and denote by $C$ its geometric
realisation. Then $C$ is the cone over~$|\mathcal S|$. The space
$C$ has a face structure, as in the proof of Theorem~\ref{zkman}.
Namely, for each simplex $\sigma\in \ord(\mathcal S)$ we denote by
$C_\sigma$ the geometric realisation of the simplicial complex
$\st\sigma=\{\tau\in\ord(\mathcal S) \colon\sigma\subset\tau\}$.
If $\sigma$ has dimension $(k-1)$, then we say that $C_\sigma$ is
a codimension-$k$ face of~$C$. Each facet $C_i$ (a face of
codimension one) is the star of a vertex of~$\ord(\mathcal S)$, as
in~\eqref{kfacet}. A face of codimension~$k$ is a connected
component of an intersection of $k$ facets. Since any face is a
cone, it is acyclic.

Although the face posets of $C$ and~$Q$ coincide, the spaces
themselves are different: faces $C_\sigma$ are defined abstractly
and they are contractible (being cones), but faces of $Q$ may be
not contractible even when $Q$ is a homology polytope.
Nevertheless, we can define the characteristic map $\lambda$ for
the face structure of $C$ by~(\ref{eqn:characteristic map}), and
define a $T$-space
\[
  M(C,\lambda)=C\times T/{\sim}
\]
by analogy with~\eqref{idspa}. Since $C$ may be not a manifold
with corners, the space $M(C,\lambda)$ is not a manifold in
general. By a straightforward generalisation of
Proposition~\ref{kact}, the space $M(C,\lambda)$ can be identified
with the quotient $\zs/K$ of the moment-angle complex
corresponding to~$\mathcal S$ (see Section~\ref{macsp}) by a
freely acting torus of dimension~$(m-n)$.

\begin{proposition}\label{prop:srfr}
We have $H^*_T(M(C,\lambda))\cong\Z[\mathcal S]$.
\end{proposition}
\begin{proof}
We have $H^*_{\T^m}(\zs)\cong\Z[\sS]$ by Exercise~\ref{evcohzs}.
Now, $H^*_{\T^m}(\zs)\cong H^*_T(M(C,\lambda))$, because
$M(C,\lambda)\cong\zs/K$ and $T=\T^m/K$.
\end{proof}

\begin{proposition}\label{eqn:Phi}
There is a face-preserving map $Q\to C$, which is covered by a
$T$-equivariant map
\[
  \Phi \colon M(Q,\lambda)\to M(C,\lambda).
\]
\end{proposition}
\begin{proof}
The map $Q\to C$ is constructed inductively; we start with a
bijection between vertices, and the extend the map to faces of
higher dimension. Each face of $C$ is a cone, there are no
obstructions to such extensions. Since the map $Q\to C$ preserves
the face structure, it is covered by a $T$-equivariant map
\[
  M(Q,\lambda)=T\times Q/\!\sim\;\longrightarrow T\times C/\!\sim\: =
  M(C,\lambda).\qedhere
\]
\end{proof}

Now we can prove the main result of this subsection:

\begin{theorem}[{\cite%[Theorem~8.3]
{ma-pa06}}]\label{hptpe}
The cohomology of a locally standard $T$-manifold $M$ is generated
in degree two if and only if the orbit space $Q$ is a homology
polytope.
\end{theorem}
\begin{proof}
Assume that $Q$ is a homology polytope. Then $M$ is homeomorphic
to the canonical model $M(Q,\lambda)$ (Lemma~\ref{acydj}), and we
can view the map $\Phi$ from Proposition~\ref{eqn:Phi} as a map
$M\to M(C,\lambda)$. For simplicity, we denote $M(C,\lambda)$ by
$M_C$ in this proof. Let $M_{C,i}=\pi^{-1}(C_i)$, \ $1\le i\le m$,
be the `characteristic' subspaces of~$M_C$. We also denote by
$\partial C$ the union of all facets $C_i$ of~$C$; topologically,
$\partial C$ is the simplicial cell complex~$|\mathcal S|$. The
$T$-action is free on $M_C\backslash \cup_i M_{C,i}$ and on
$M\backslash \cup_i M_{i}$, so we have
\begin{equation*} \label{eqn:relative}
  H_T^*(M_C,\cup_i M_{C,i})\cong H^*(C,\partial C), \quad
  H^*_T(M,\cup_i M_{i})\cong H^*(Q,\partial Q).
\end{equation*}
Therefore, the map $\Phi$ induces a map between exact sequences
\begin{equation}\label{seqma}
\begin{CD}
  \longrightarrow\; H^*(C,\partial C) @>>> H_T^*(M_C) @>>>
  H_T^*(\cup_i M_{C,i})\;\longrightarrow\\
  @VVV @VV\Phi^*V @VVV\\
  \longrightarrow\; H^*(Q,\partial Q) @>>> H_T^*(M) @>>>
  H_T^*(\cup_i M_{i})\;\longrightarrow
\end{CD}
\end{equation}
Each $M_{i}$ itself is a torus manifold with quotient homology
polytope~$F_i$. Using induction and the Mayer--Vietoris sequence,
we may assume that the map $H_T^*(\cup_i M_{C,i})\to H_T^*(\cup_i
M_{i})$ above is an isomorphism. By Lemma~\ref{kgore}, $\partial
C\cong|\sS|$ has homology of an $(n-1)$-sphere. Hence
$H^*(C,\partial C)\cong H^*(D^n,S^{n-1})$, because $C$ is the cone
over~$\partial C$. We also have $H^*(Q,\partial Q)\cong
H^*(D^n,S^{n-1})$, because $Q$ is a homology polytope. Using these
isomorphisms, we see from the construction of the map $\Phi$ that
the induced map $H^*(C,\partial C)\to H^*(Q,\partial Q)$ is the
identity on $H^*(D^n,S^{n-1})$. By applying the 5-lemma
to~\eqref{seqma} we obtain that $\Phi^*\colon H_T^*(M_C)\to
H_T^*(M)$ is an isomorphism. This together with
Proposition~\ref{prop:srfr} implies $H_T^*(M)\cong \Z[\mathcal
S]$. The ring $\Z[\sS]$ is Cohen--Macaulay by Lemma~\ref{kgore}.
Therefore, all conditions of Lemma~\ref{theo:face CM} are
satisfied, and $H^*(M)$ is generated by its degree-two part.

\smallskip

Assume now that that $H^*(M)$ is generated in degree two. Since
all nonempty intersections of characteristic submanifolds are
connected and their cohomology rings are generated in degree two
(Lemma~\ref{cmcon}), we may assume by induction that all proper
faces of $Q$ are homology polytopes. In particular, the proper
faces are acyclic, whence $H^*(\partial Q)\cong H^*(\partial C)$,
because both $\partial Q$ and $\partial C$ have acyclic coverings
with the same nerve~$\mathcal S$. This together with
Lemma~\ref{kgore} implies
\begin{equation}\label{eqn:boundary Q}
  H^*(\partial Q)\cong H^*(S^{n-1}).
\end{equation}
We need to show that $Q$ itself is acyclic. We first prove the
following:

\begin{claim}
$H^2(Q)=0$.
\end{claim}
\begin{proof}
The claim is trivial for $n=1$. If $n=2$ then $Q$ is a surface
with boundary, hence $H^2(Q)=0$. Now assume $n\ge 3$. We consider
the equivariant cohomology exact sequence of pair $(M,\cup_i M_i)$
(the bottom row of~\eqref{seqma}). All the maps in the exact
sequence are $H^*(BT)$-module maps. By Lemma~\ref{lemm:freeness},
$H^*_T(M)$ is a free $H^*(BT)$-module. On the other hand,
$H^*(Q,\partial Q)$ is finitely generated over~$\Z$, so it is a
torsion $H^*(BT)$-module. It follows that the whole sequence
splits into short exact sequences:
\begin{equation}\label{short exact}
  0 \to H^k_T(M) \to H^k_T(\cup_i M_i) \to H^{k+1}(Q,\partial Q) \to 0
\end{equation}
By setting $k=1$ we obtain
\[
  H^1_T(\cup_iM_i)\cong H^2(Q,\partial Q).
\]
The same argument as in Lemma~\ref{lemm:H1=0} shows that
\[
  H^1_T(\cup_iM_i)=H^1((\cup_iM_i)/T)=H^1(\partial Q).
\]
By considering the projection $(ET\times M)/T\to M/T=Q$ we
conclude that the coboundary map $H^1(\partial Q)\to
H^2(Q,\partial Q)$ is an isomorphism. Therefore, we get the
following fragment of the exact sequence of pair:
\[
  0\to H^2(Q) \to H^2(\partial Q) \to H^3(Q,\partial Q).
\]
By (\ref{eqn:boundary Q}), $H^2(\partial Q)\cong H^2(S^{n-1})$,
whence $H^2(Q)=0$ for $n\ge 4$. If $n=3$, the coboundary map
$H^2(\partial Q) \to H^3(Q,\partial Q)$ above is an isomorphism
because $Q$ is orientable by Lemma~\ref{lemm:H1=0}, whence
$H^2(Q)=0$ again.
\end{proof}

Now we resume the proof of the theorem. We obtain a
$T$-homeomorphism $M\to M(Q,\lambda)$ (because $H^2(Q)=0$ and all
proper faces are acyclic by the inductive assumption, see the
remark after Proposition~\ref{acydj}), and therefore a $T$-map
$\Phi\colon M\to M_C$, as in the proof of the `if' part of the
theorem. We consider diagram~\eqref{seqma} again. Using induction
and a Mayer--Vietoris argument, we may assume that $H^*_T(\cup_i
M_{C,i})\to H^*_T(\cup_i M_{i})$ is an isomorphism. By
Lemma~\ref{kgore}, $H^*(C,\partial C)\cong H^*(D^n,S^{n-1})$. The
map $Q\to C$ used in the construction of $\Phi$ induces a map
\begin{equation} \label{eqn:QP}
  H^*(D^n,S^{n-1})\cong H^*(C,\partial C)\to H^*(Q,\partial Q),
\end{equation}
which is an isomorphism in dimension~$n$.
%, and hence it is injective in each degree.
Applying the 5-lemma (Exercise~\ref{5lemma}) to~\eqref{seqma}, we
obtain that $\Phi^*\colon H^*_T(M_C)\to H^*_T(M)$ is injective.
Theorem~\ref{th:eqcoh} and Proposition~\ref{prop:srfr} imply that
$H^*_T(M)\cong \Z[Q]\cong H^*_T(M_C)$, and all graded components
of these rings are finitely generated. The same argument works
with any field coefficients, so that $\Phi^*\colon H^*_T(M_C)\to
H^*_T(M)$ is actually an isomorphism. By applying the 5-lemma
again to diagram~\eqref{seqma}, we obtain that~(\ref{eqn:QP}) is
an isomorphism, i.e. $H^*(Q,\partial Q)\cong H^*(D^n,S^{n-1})$.
This together with~(\ref{eqn:boundary Q}) implies that $Q$ is
acyclic.
\end{proof}

Theorem~\ref{hptpe} shows that if the cohomology ring of a locally
standard $T$-manifold is generated in degree two, then the
combinatorics of the orbit space $Q$ is fully determined by its
nerve simplicial complex. The following result gives a
characterisation of simplicial complexes arising in this way.

\begin{proposition}
A simplicial complex $\sK$ can be the nerve of a locally standard
$T$-manifold with cohomology generated in degree two if and only
if $\sK$ is Gorenstein* and $\Z[\sK]$ admits an lsop.
\end{proposition}
\begin{proof}
If $H^*(M)$ is generated in degree two, then $\mathcal K$ is
Gorenstein* by Lemma~\ref{kgore}. In particular, $\Z[\mathcal K]$
is a Cohen--Macaulay ring. Furthermore, $H^*_T(M)\cong \Z[\mathcal
K]$ by Theorem~\ref{th:eqcoh}. Since $H^*_T(M)\cong H^*(BT)\otimes
H^*(M)$ as a $H^*(BT)$-module, the ring $\Z[\mathcal K]$ admits an
lsop.

Now assume that $\Z[\mathcal K]$ is Gorenstein* and admits an
lsop. By~\cite[Theorem~12.2]{davi83}, there exists a homology
polytope $Q$ with nerve~$\mathcal K$. Since $\Z[\mathcal K]$
admits an lsop, any element $t\in H^2(BT)$ can be written as
\[
  t=\sum_{i=1}^m\lambda_i(t)v_i
\]
with $\lambda_i(t)\in\Z$. Clearly, $\lambda_i(t)$ is linear
in~$t$, so that we can view $\lambda_i$ as an element of the dual
lattice $H_2(BT)$ (see Proposition~\ref{prop:a_i}). Now define a
map $\lambda$~\eqref{eqn:characteristic map} which sends $F_i$
to~$\lambda_i$. Then $M=M(Q,\lambda)$ (see~\eqref{idspa}) is a
locally standard $T$-manifold, and its cohomology is generated in
degree two by Theorem~\ref{hptpe}.
\end{proof}

\subsection*{Torus manifolds over face-acyclic manifolds with corners}
Here we prove the second main result on the cohomology of
$T$-manifolds, Theorem~\ref{cohfa}. It states that the orbit space
$Q$ of a locally standard $T$-manifold $M$ is face-acyclic if and
only if $H^{odd}(M)=0$. The proof is by reduction to
Theorem~\ref{hptpe} on $T$-manifolds over homology polytopes; it
relies upon the operation of blow-up and the algebraic results of
Section~\ref{scccm}.

As before, $M$ is a locally standard $T$-manifold with orbit
projection $\pi\colon M\to Q$.

\begin{construction}[Blow-up of a $T$-manifold]\label{lstbu}
Let $M_G=\pi^{-1}(G)$ be the submanifold corresponding to a face
$G\subset Q$, and $\n{G}=\nu(M_G\subset M)$ the normal bundle.
Since $M_G$ is a transverse intersection of characteristic
submanifolds, $\n{G}$ is the Whitney sum of their normal bundles.
The omniorientation on $M$ makes $\n{G}$ into a complex
$T$-bundle.

Consider the $T$-bundle $\n{G}\oplus\underline{\C}$, where the
$T$-action on the trivial summand $\underline{\C}$ is trivial. The
projectivisation $\C P(\n{G}\oplus\underline{\C})$ is a locally
standard $T$-manifold containing~$M_G$, and there are invariant
neighbourhoods of $M_G$ in $M$ and of $M_G$ in $\C
P(\n{G}\oplus\underline{\C})$ which are $T$-diffeomorphic. After
removing these invariant neighbourhoods of $M_G$ from $M$ and $\C
P(\n{G}\oplus\underline{\C})$ and reversing orientation on the
latter, we can identify the resulting $T$-manifolds along their
boundaries. As a result, we obtain a locally standard
$T$-manifold~$\widetilde M$, which is called the
\emph{blow-up}\label{dblowup} of $M$ at~$M_G$.

If $G$ is a vertex, then $\widetilde M$ is diffeomorphic to the
connected sum $M\cs\overline{\C P^n\!}$.

There is a blow-down map $\widetilde M\to M$, which collapses the
total space $\C P(\n{G}\oplus\underline{\C})$ onto $M_G$ and is
the identity on the remaining part of~$\widetilde M$.

The orbit space $\widetilde Q$ of $\widetilde M$ is obtained by
truncating $Q$ at the face~$G$. As a result, $\widetilde Q$
acquires a new facet, which we denote by~$\widetilde G$. The
simplicial cell complex dual to $\widetilde Q$ is obtained from
the dual of~$Q$ by applying a stellar subdivision at the face dual
to~$G$ (see Definition~\ref{bist}).
\end{construction}

\begin{lemma}\label{qtild}
$\widetilde Q$ is face-acyclic if and only if $Q$ is face-acyclic.
\end{lemma}
\begin{proof}
All new faces appearing as the result of truncating $Q$ at $G$ are
contained in the facet $\widetilde G\subset\widetilde Q$. The
blow-down map $\widetilde M\to M$ induces the projection
$\widetilde Q\to Q$ collapsing $\widetilde G$ onto~$G$. The face
$G$ is a deformation retract of~$\widetilde G$ (combinatorially,
$\widetilde G$ is a product of $G$ and a simplex). Hence $G$ is
acyclic if and only if $\widetilde G$ is acyclic. Similarly, any
other new face $\widetilde Q$ deformation retracts onto a face
of~$Q$. Furthermore, the map $\widetilde Q\to Q$ is also a
deformation retraction.
\end{proof}

\begin{lemma}\label{mtild}
$H^{odd}(\widetilde M)=0$ if and only if $H^{odd}(M)=0$.
\end{lemma}
\begin{proof}
Assume that $H^{odd}(M)=0$. By Lemma~\ref{lemm:odd=0},
$H^{odd}(M_G)=0$. The facial submanifold $M_G\subset M$ is blown
up to a codimension-two submanifold $\widetilde M_{\widetilde
G}\cong \C P(\n{G})$. The cohomology of the projectivisation of a
complex vector bundle over $M_G$ is a free $H^*(M_G)$-module on
even-dimensional generators (see, e.g.~\cite[Chapter~V]{ston68}).
Therefore, $H^{odd}(\widetilde M_{\widetilde G})=0$.

The blow-down map $\widetilde M\to M$ induces a map between exact
sequences of pairs
$$
\begin{CD}
  H^{k-1}(M_G) @>>> H^k(M,M_G) @>>> H^k(M) @>>> H^k(M_G) \\
    @VVV              @VV\cong V            @VVV        @VVV     \\
  H^{k-1}(\widetilde M_{\widetilde G}) @>>> H^k(\widetilde M,\widetilde M_{\widetilde G}) @>>>
   H^k(\widetilde M) @>>> H^k(\widetilde M_{\widetilde G})
\end{CD}
$$
where the second vertical arrow is an isomorphism by excision.
Assume that $k$ is odd. Since $H^k(M)=0$, the map $H^{k-1}(M_G)\to
H^k(M,M_G)$ is onto. Therefore, $H^{k-1}(\widetilde M_{\widetilde
G})\to H^k(\widetilde M,\widetilde M_{\widetilde G})$ is onto.
Since $H^k(\widetilde M_{\widetilde G})=0$, this implies
$H^k(\widetilde M)=0$.

\smallskip

To prove the opposite statement, we use the algebraic results from
Section~\ref{scccm}. Assume $H^{odd}(\widetilde M)=0$. Let
$\mathcal S$ be the dual simplicial poset of~$Q$, and let
$\widetilde{\mathcal S}$ be the dual poset of~$\widetilde Q$. Then
$\widetilde{\mathcal S}$ is obtained from $\mathcal S$ by stellar
subdivision at the face dual to~$G$. By Lemma~\ref{theo:face CM},
$\Z[\widetilde{\mathcal S}]$ is a Cohen--Macaulay ring. We claim
that $\Z[\mathcal S]$ is also Cohen--Macaulay (i.e. the converse
of Lemma~\ref{lemm:cmtil} holds). Indeed, Theorem~\ref{theo:cmpos}
implies that $\widetilde{\mathcal S}$ is a Cohen--Macaulay
simplicial poset. Let $\mathcal K$ be a simplicial complex which
is a common subdivision of simplicial cell complexes
$\widetilde{\mathcal S}$ and $\mathcal S$ (for example, we may
take $\mathcal K$ to be the barycentric subdivision of
$\widetilde{\mathcal S}$). Then $\mathcal K$ is a Cohen--Macaulay
complex by Corollary~\ref{coro:cmcha}, hence $\mathcal S$ is a
Cohen-Macaulay simplicial poset. Another application of
Theorem~\ref{theo:cmpos} gives that $\Z[\mathcal S]$ is a
Cohen--Macaulay ring. Finally, Lemma~\ref{theo:face CM} implies
that $H^{odd}(M)=0$.
\end{proof}

Now we can prove our final result:

\begin{theorem}[{\cite%[Theorem~9.3]
{ma-pa06}}]\label{cohfa}
The odd-degree cohomology of a locally standard $T$-manifold $M$
vanishes if and only if the orbit space $Q$ is face-acyclic.
\end{theorem}
\begin{proof}
The idea is to reduce to Theorem~\ref{hptpe} by blowing up
sufficiently many facial submanifolds. If the orbit space of $M$
is face-acyclic, then it becomes a homology polytope after
sufficiently many blow-ups.

Let $\mathcal S$ be the simplicial poset dual to~$Q$. Since the
barycentric subdivision is a sequence of stellar subdivisions
(Proposition~\ref{prop:barst}), by applying appropriate blow-ups
we get a torus manifold $M'$ with orbit space $Q'$ such that the
face poset of $Q'$ is the barycentric subdivision of the face
poset of $Q$. The collapse map $M'\to M$ is a composition of
blow-down maps:
\begin{equation}\label{blseq}
  M=M_0 \longleftarrow M_1 \longleftarrow \cdots \longleftarrow M_k=M'.
\end{equation}

Assume that $H^{odd}(M)=0$. Then $M$ is locally standard by
Theorem~\ref{theo:local standardness}.  By applying
Lemma~\ref{mtild} successively, we get $H^{odd}(M')=0$. By
construction, all intersections of faces of $Q'$ are connected, so
$H^*(M')$ is generated in degree two by Lemma~\ref{theo:face CM}
and $Q'$ is a homology polytope by Theorem~\ref{hptpe}. In
particular, $Q'$ is face-acyclic. Finally, by applying
Lemma~\ref{qtild} successively, we conclude that $Q$ is also
face-acyclic.

Assume now that $Q$ is face-acyclic. By applying Lemma~\ref{qtild}
successively, we obtain that $Q'$ is also face-acyclic. On the
other hand, $\mathcal S'$ is a simplicial complex, hence $Q'$ is a
homology polytope. By Theorem~\ref{hptpe}, $H^{odd}(M')=0$. By
applying lemma~\ref{mtild} successively, we finally conclude that
$H^{odd}(M)=0$.
\end{proof}

\begin{remark}
As one can easily observe, the argument in the ``only if'' part of
the above theorem is independent of Theorem~\ref{theo:GKM} and
Theorem~\ref{th:eqcoh}. Now, given that $Q$ is face-acyclic, one
readily deduces that the 1-skeleton of $Q$ is connected. Indeed,
otherwise the smallest face containing vertices from two different
connected components of the 1-skeleton would be a manifold with at
least two boundary components and thereby non-acyclic. Thus, our
reference to Theorem~\ref{theo:GKM} was actually irrelevant,
although it made the argument more straightforward.
\end{remark}

\subsection*{Exercises}
\begin{exercise}
If $M$ is orientable and $H^{odd}(M)=0$, then $H^*(M)$ is
torsion-free.
\end{exercise}

\begin{exercise}
Any torus manifold $M$ has at least two characteristic
submanifolds.
\end{exercise}

\begin{exercise}
Any $T$-fixed point of a torus $2n$-manifold $M$ is contained in
an intersection of $n$ characteristic submanifolds.
\end{exercise}

\begin{exercise}
Each face of a face-acyclic manifold with corners has a vertex.
\end{exercise}

\begin{exercise}
The equivariant Chern class of torus manifold $M$ with an
invariant stably complex structure is given by
\[
  c^T(M)=\prod_{i=1}^m(1+\tau_i)\quad\text{ modulo $H^*(BT)$-torsion}
\]
where $\tau_i\in H_T^2(M)$ is the Thom class defined before
Proposition~\ref{prop:a_i}. (Hint: use Corollary~\ref{kerres},
see~\cite[Theorem~3.1]{masu99}.)
\end{exercise}

\begin{exercise}
Let $M$ be a torus manifold of dimension~$2n$ with
$H^{odd}(M;\Z_2)=0$ and let $G\subset T$ denote the discrete
subgroup isomorphic to $\Z_2^n$. Show that the $G$-equivariant
Stiefel--Whitney class of $M$ is given by
\[
  w^G(M)=\prod_{i=1}^m(1+\tau_i),
\]
where $\tau_i\in H_G^2(M;\Z_2)$ is the mod-2 Thom class.
\end{exercise}

\begin{exercise}
Prove the following particular case of Theorem~\ref{heven}. Let
$\mathcal S$ be a Gorenstein* simplicial poset such that there
exists a torus manifold $M$ with $H^{odd}(M)=0$ and orbit space
$Q$ whose dual poset is~$\sS$. Let $\mb h(\mathcal
S)=(h_0,h_1,\ldots,h_n)$ be the $h$-vector. Assume that $n$ is
even and $h_i=0$ for some~$i$. Then $h_{n/2}$ is even. (Hint: use
the previous exercise, see~\cite[Theorem~10.1]{ma-pa06}. An
algebraic version of this argument was used in~\cite{masu05} to
prove Theorem~\ref{heven} completely.)
\end{exercise}

\section{Topological toric manifolds}\label{ttman}
Recall that a toric manifold is a smooth complete (compact)
algebraic variety with an effective algebraic action of an
algebraic torus $(\C^\times)^n$ having an open dense orbit. More
constructively, toric varieties can be defined via the fan-variety
correspondence; this gives a covering of the variety by invariant
affine charts, which in the smooth complete case are algebraic
representation spaces of~$(\C^\times)^n$.

The idea behind Ishida, Fukukawa and Masuda's generalisation of
toric manifolds is to combine topological versions of these two
definitions of toric varieties:

\begin{definition}[{\cite{i-f-m13}}]\label{defittm}
A \emph{topological toric manifold} is a closed smooth manifold
$X$ of dimension $2n$ with an effective \emph{smooth} action of
$(\C^\times)^n$ having an open dense orbit and covered by finitely
many invariant open subsets each equivariantly diffeomorphic to a
smooth representation space of~$(\C^\times)^n$.
\end{definition}

As is pointed out in~\cite{i-f-m13}, keeping only the first part
of the definition (i.e. a smooth $(\C^\times)^n$-action with a
dense orbit) leads to a vast and untractable class of objects;
therefore it is important to include the covering by equivariant
charts.

In this section we review the main properties of topological toric
manifolds, and outline the construction of the correspondence
between topological toric manifolds and generalised fans, called
\emph{topological fans}. We mainly follow the notation and
terminology of~\cite{i-f-m13}. The details of proofs can be found
in the original paper.

The effectiveness of the $(\C^\times)^n$-action on a topological
toric manifold $X$ implies that the smooth representation of
$(\C^\times)^n$ modelling each invariant chart of $X$ is faithful.
A faithful smooth real $2n$-dimensional representation of
$(\C^\times)^n$ is isomorphic to a direct sum of complex
one-dimensional representations.

To get a hold on smooth representations of~$(\C^\times)^n$, we
consider the case $n=1$ first. Since $GL(1,\C)=\C^\times$, a
smooth representation of $\C^\times$ in~$\C$ can be viewed as a
smooth endomorphism of~$\C^\times$. Such an endomorphism has the
form
\[
  z\mapsto z^\mu=|z|^{b+ic}\Bigl(\frac z{|z|}\Bigr)^a\qquad
  \text{with }\mu=(b+ic,a)\in\C\times\Z.
\]
The representation given by $z\to z^\mu$ is algebraic if and only
if $c=0$ and $b=a$.

Composition of smooth endomorphisms of $\C^\times$ defines a
(noncommutative) product on $\C\times\Z$, given by
\[
  (z^{\mu_1})^{\mu_2}=z^{\mu_2\mu_1},\qquad
  \mu_2\mu_1=(b_1b_2 + i(b_1c_2+c_1a_2),a_1a_2).
\]
This product becomes the matrix product if we represent $\mu$ by
$2\times 2$-matrices:
\[
  \begin{pmatrix}b_2&0\\c_2&a_2\end{pmatrix}
  \begin{pmatrix}b_1&0\\c_1&a_1\end{pmatrix}=
  \begin{pmatrix}b_2b_1&0\\c_2b_1+a_2c_1&a_2a_1\end{pmatrix}.
\]

Let $\mathcal R$ be the ring consisting of elements of
$\C\times\Z$ with componentwise addition and multiplication
defined above. The ring $\mathcal R$ is therefore isomorphic to
the ring $\Hom_{\mathrm{sm}}(\C^\times,\C^\times)$ of smooth
endomorphisms of~$\C^\times$.

Given $\alpha=(\alpha^1,\dots,\alpha^n)\in \mathcal R^n$ and
$\beta=(\beta^1,\dots,\beta^n)\in \mathcal R^n$, define smooth
homomorphisms $\chi^\alpha\in\Hom((\C^\times)^n,\C^\times)$ and
$\lambda_\beta\in \Hom(\C^\times,(\C^\times)^n)$ by
\[
%\begin{equation}\label{chilambda}
  \chi^\alpha(z_1,\dots,z_n)=\prod_{k=1}^n z_k^{\alpha^k},
  \qquad\lambda_\beta(z)=(z^{\beta^1},\ldots,z^{\beta^n}),
%\end{equation}
\]
and also define
\[
  \langle \alpha,\beta\rangle=\sum_{k=1}^n \alpha^k\beta^k\in \mathcal R.
\]

The following properties are checked easily:
\begin{itemize}
%\item[(a)] $(z_1,\ldots,z_n)\in (\C^\times)^n$ is the identity element if and only if
%$\chi^{\alpha}(z_1,\ldots,z_n)=1$ for any $\alpha\in \mathcal
%R^n$.
\item[(a)] $\chi^\alpha(\lambda_{\beta}(z))=z^{\langle
\alpha,\beta\rangle}$;
\item[(b)] ${\displaystyle\lambda_{\beta}(\chi^{\alpha}(z_1,\ldots,z_n))
=\Bigl(\prod\limits_{k=1}^nz_k^{\beta^1\alpha^k},\ldots,
\prod\limits_{k=1}^nz_k^{\beta^n\alpha^k}\Bigr)}$.
\end{itemize}

Given $\alpha_1,\ldots,\alpha_n\in\mathcal R^n$, define the
endomorphism $\bigoplus_{i=1}^n\chi^{\alpha_i}$ of $(\C^\times)^n$
by
\[
  \Bigl(\bigoplus_{i=1}^n\chi^{\alpha_i}\Bigr)(z_1,\ldots,z_n)
  =(\chi^{\alpha_1}(z_1,\ldots,z_n),\ldots,\chi^{\alpha_n}(z_1,\dots,z_n)).
\]

\begin{proposition}
Any smooth representation of $(\C^\times)^n$ in $\C^n$ has the
form  $\bigoplus_{i=1}^n\chi^{\alpha_i}$ with $\alpha_i\in\mathcal
R^n$. This representation is faithful if and only if the $n\times
n$-matrix formed by the coordinates of $\alpha_i$ has an inverse
in $\mathop{\mathrm{Mat}}_n(\mathcal R)$.
\end{proposition}

A faithful representation of $(\C^\times)^n$ in $\C^n$ has a
unique fixed point~$\mathbf0$, hence the fixed point set
$X^{(\C^\times)^n}$ of a topological toric manifold is finite.

A closed connected submanifold of real codimension two in a
topological toric manifold $X$ is called
\emph{characteristic}\label{charmttm} if it is fixed pointwise by
a subgroup isomorphic to~$\C^\times$. There are finitely many
characteristic submanifolds in~$X$, and we denote them by
$X_1,\ldots,X_m$.

It can be easily seen that a topological toric manifold $X$ is
simply connected (\cite[Proposition~3.2]{i-f-m13}). In particular,
$X$ is orientable. For each characteristic submanifold~$X_j$, the
normal bundle $\nu_j=\nu(X_j\subset X)$ is orientable as a
$(\C^\times)$-equivariant bundle. Therefore, $X_j$ itself is also
orientable. A choice of an orientation for each $X_j$ together
with an orientation of~$X$ is called an
\emph{omniorientation}\label{omniorittm} on~$X$. Topological toric
manifolds are assumed to be omnioriented below.

\begin{lemma}[{\cite[Lemma~3.3]{i-f-m13}}]\label{betai}
For each characteristic submanifold~$X_j$, there is a unique
$\beta_j(X)\in \mathcal R^n$ such that the subgroup
$\lambda_{\beta_j(X)}(\C^\times)\subset(\C^\times)^n$ fixes $X_j$
pointwise and $\lambda_{\beta_j(X)}(z)_*\xi=z\xi$ for any
$z\in\C^\times, \xi\in\nu_i$, where $\lambda_{\beta_j(X)}(z)_*$
denotes the differential of~$\lambda_{\beta_j(X)}(z)$.
\end{lemma}

Characteristic submanifolds of $X$ intersect transversely.
Furthermore, multiple intersections of characteristic submanifolds
are all connected (as in the case of quasitoric manifolds, and
unlike general torus manifolds), see~\cite[Lemma~3.6]{i-f-m13}. In
particular, any fixed point $v\subset X^{(\C^\times)^n}$ is an
intersection of an $n$-tuple $X_{j_1},\ldots,X_{j_n}$ of
characteristic submanifolds, so we have an isomorphism of real
$(\C^\times)^n$-representation spaces
\[
  \mathcal T_v X\cong(\nu_{j_1}\oplus\cdots\oplus\nu_{j_n})|_v.
\]
The omniorientation of $X$ defines orientations for the left and
right hand side of the identity above. These two orientations may
be different, so the sign of~$v$ is defined (compare
Lemma~\ref{qts}~(a)).

The element $\beta_j(X)$ defined in Lemma~\ref{betai} can be
written as
\begin{equation}\label{betaj}
  \beta_j(X)=(\mb b_j(X)+i\mb c_j(X),\mb a_j(X))\in \C^n\times\Z^n.
\end{equation}
Here is an analogue of Proposition~\ref{qtw} for topological toric
manifolds:

\begin{lemma}[{\cite[Lemma~3.4]{i-f-m13}}]\label{ttmdb}
Let $v=X_{j_1}\cap\cdots\cap X_{j_n}$ be a fixed point of~$X$.
Then $\{\mb b_{j_1}(X),\ldots,\mb b_{j_n}(X)\}$ and $\{\mb
a_{j_1}(X),\ldots,\mb a_{j_n}(X)\}$ are bases of $\R^n$ and $\Z^n$
respectively.

The complex $\C^\times$-representation space
$(\nu_{j_1}\oplus\cdots\oplus\nu_{j_n})|_v$ is isomorphic to
$\bigoplus_{k=1}^n\chi^{\alpha^{(v)}_k}$, where
$\{\alpha^{(v)}_1,\ldots,\alpha^{(v)}_n\}$ is the dual set of
$\{\beta_{j_1}(X),\ldots,\beta_{j_n}(X)\}$, defined uniquely by
the condition
\[
  \bigl\langle\alpha^{(v)}_k,\beta_{j_\ell}(X)\bigr\rangle=\delta_{k\ell}
\]
(here $\delta_{k\ell}$ denotes the Kronecker delta).
\end{lemma}

Define the simplicial complex $\sK(X)$ on $[m]$ whose simplices
correspond to nonempty intersections of characteristic
submanifolds:
\[
  \sK(X)=\bigl\{I=\{i_1,\ldots,i_k\}\in[m]\colon
  X_{i_1}\cap\cdots\cap X_{i_k}\ne\varnothing\bigr\}.
\]

When $X$ is a toric manifold, we have $\mb b_j(X)=\mb a_j(X)$ and
$\mb c_j(X)=\bf0$ in~\eqref{betaj}, and the primitive vector $\mb
a_j(X)$ corresponds to the 1-parameter algebraic subgroup of
$(\C^\times)^n$ fixing the divisor~$X_j$. Furthermore, the data
$\{\sK(X);\mb a_1(X),\ldots,\mb a_m(X)\}$ define a complete
simplicial regular fan (see Section~\ref{mamsf}). Now let us see
what kind of combinatorial structure replaces a fan in the case of
topological toric manifolds.

Given $I\in\sK(X)$, let $\sigma_I=\R_\ge\langle\mb b_i\colon i\in
I\rangle$ denote the cone spanned by the vectors $\mb b_i\in\R^n$
with $i\in I$. By Lemma~\ref{ttmdb}, $\sigma_I$ is a is a
simplicial cone of dimension~$|I|$.

\begin{lemma}[{\cite[Lemma~3.7]{i-f-m13}}]\label{DeltaX}
$\bigcup_{I\in\sK(X)}\sigma_I=\R^n$ and
$\sigma_I\cap\sigma_J=\sigma_{I\cap J}$. In other words, the data
$\{\sK(X);\mb b_1(X),\ldots,\mb b_m(X)\}$ define a complete
simplicial fan.
\end{lemma}

\begin{definition}\label{deftopfan}
Let $\sK$ be a simplicial complex on~$[m]$, and let
\[
  \beta_j=(\mb b_j+i\mb c_j,\mb a_j)\in \C^n\times\Z^n,\quad
  j=1\ldots,m
\]
be a collection of $m$ elements of $\C^n\times\Z^n$. The data
$\{\sK;\beta_1,\ldots,\beta_m\}$ is said to define a (regular)
\emph{topological fan}~$\Delta$ if the following two conditions
are satisfied:
\begin{itemize}
\item[(a)]
the data $\{\sK;\mb b_1,\ldots,\mb b_m\}$ define a simplicial fan
in~$\R^n$;
\item[(b)] for each $I\in\sK$, the set $\{\mb a_i\colon i\in I\}$
is a part of basis of~$\Z^n$.
\end{itemize}
A topological fan $\Delta$ is said to be \emph{complete} if the
ordinary fan from~(a) is complete.
\end{definition}

Note that the fan of~(a) is not required to be rational or
regular, but if $\mb a_j=\mb b_j$ for all~$j$, then $\Delta$
becomes a regular ordinary fan.

\begin{theorem}[{\cite[Theorem~8.1]{i-f-m13}}]\label{ttmtf} There is a
bijective correspondence between omnioriented topological toric
manifolds of dimension~$2n$ and complete topological fans of
dimension~$n$.
\end{theorem}
\begin{proof}[Sketch of proof]
Let $X$ be a topological toric manifold. By Lemma~\ref{DeltaX},
the data $(\sK(X);\beta_1(X),\ldots,\beta_m(X))$ define a complete
topological fan~$\Delta(X)$.

Now let $\Delta$ be a complete topological fan, defined by data
$(\sK;\beta_1,\ldots,\beta_m)$. For each maximal simplex
$I=\{i_1,\ldots,i_n\}\in\sK$, let
$\{\alpha^I_1,\ldots,\alpha^I_n\}$ be the dual set of
$\{\beta_{i_1},\ldots,\beta_{i_n}\}$ (compare Lemma~\ref{ttmdb}).
Condition~(b) of Definition~\ref{deftopfan} guarantees that the
complex $n$-dimensional representation
$\bigoplus_{k=1}^n\chi^{\alpha_k^I}$ of $(\C^\times)^n$ is
faithful. These representation spaces corresponding to all maximal
$I\in\sK$ patch together into a topological space $X(\Delta)$
locally homeomorphic to~$\C^n$, and $(\C^\times)^n$ acts on $X$
smoothly with an open dense orbit. As in the case of ordinary
fans, condition~(b) of Definition~\ref{deftopfan} guarantees that
the space $X(\Delta)$ is Hausdorff, so it is a smooth manifold.
(However, the algebraic criterion for separatedness used in the
proof of Lemma~\ref{sepvar} cannot be used here; a topological
argument is needed.) Finally the condition that $\Delta$ is
complete gives that $X(\Delta)$ is compact, i.e. closed.

An alternative way to proceed is to use an analogue of the
quotient construction of toric varieties, described in
Section~\ref{algtq}. To do this, define the coordinate subspace
arrangement complement $U(\sK)$ by~\eqref{UK}, and define the
homomorphism
\[
  \lambda\colon(\C^\times)^m\to(\C^\times)^n,\qquad
  \lambda(z_1,\ldots,z_m)=\prod_{k=1}^m\lambda_{\beta_k}(z_k).
\]
Then $\lambda$ is surjective and its kernel is given by
\[
  \Ker\lambda=\bigl\{(z_1,\ldots,z_m)\in(\C^\times)^m\colon
  \prod_{i=1}^m z_i^{\langle\alpha,\beta_i\rangle}=1
  \quad\text{for any }\alpha\in\mathcal R^n\bigr\},
\]
by analogy with~\eqref{gexpl}. Then define
\[
  X(\Delta)=U(\sK)/\Ker\lambda=\bigcup_{I\in\sK}(\C,\C^\times)^I/\Ker\lambda.
\]
The space $X(\Delta)$ has a smooth action of
$(\C^\times)^m/\Ker\lambda\cong(\C^\times)^n$ with an open dense
orbit. Furthermore, for each maximal $I\in\sK$ there is an
equivariant diffeomorphism
\[
  \varphi_I\colon(\C,\C^\times)^I/\Ker\lambda\to
  \bigoplus_{k=1}^n\chi^{\alpha_k^I},
\]
where the latter is the faithful smooth
$(\C^\times)^n$-representation space defined above. Condition~(b)
of Definition~\ref{deftopfan} translates into the condition
of~$X(\Delta)$ being Hausdorff. The fact that only the `real part'
$(\sK;\mb b_1,\ldots,\mb b_m)$ of the topological fan data matters
when deciding whether the quotient is Hausdorff should be clear
from the similar argument in the proof of
Theorem~\ref{zksmooth}~(b). Finally, $X(\Delta)$ is compact
because $\Delta$ is complete. Thus, $X(\Delta)$ with the local
charts $\{(\C,\C^\times)^I/\Ker\lambda,\varphi_I\}$ is a
topological toric manifold.
\end{proof}

For the classification of topological toric manifolds up to
equivariant diffeomorphism or homeomorphism,
see~\cite[Corollary~8.2]{i-f-m13}

One can restrict the $(\C^\times)^n$-action on $X$ to the compact
$n$-torus~$\T^n$. The resulting $\T^n$-manifold $X$ is obviously a
locally standard torus manifold, so the quotient $X/\T^n$ is a
manifold with corners.

\begin{lemma}[{\cite[Lemma~7.1]{i-f-m13}}]
All faces of $X/\T^n$ are contractible and the face poset of
$X/\T^n$ coincides with the inverse poset of~$\sK(X)$.
\end{lemma}

It follows that $X/\T^n$ is a homology polytope. As a corollary of
Theorem~\ref{theo:stcoh} and Theorem~\ref{hptpe}, we obtain the
following description of the cohomology of~$X$, similar to toric
or quasitoric manifolds:

\begin{proposition}\label{cohttm}
Let $X$ be a topological toric manifold, whose associated
topological fan is defined by the data
$(\sK(X);\beta_1(X),\ldots,\beta_m(X))$. Then the cohomology ring
of $X$ is given by
\[
  H^*(X)\cong \Z[v_1,\ldots,v_m]/\mathcal I,
\]
where $v_i\in H^2(X)$ is the class dual to the characteristic
submanifold~$X_i$, and $\mathcal I$ is the ideal generated by
elements of the following two types:
\begin{itemize}
\item[(a)] $v_{i_1}\cdots v_{i_k}$ with $\{i_1,\ldots,i_k\}\notin\sK(X)$;

\item[(b)] $\displaystyle\sum_{i=1}^m\bigl\langle\mb u,\mb a_i(X)\bigr\rangle v_i$, for
any $\mb u\in\Z^n$.
\end{itemize}
The element $\mb a_i(X)\in\Z^n$ here is the second coordinate
of~$\beta_i(X)$, see~\eqref{betaj}.
\end{proposition}

\section{Relationship between different classes of
$T$-manifolds}\label{reltm}

The relationship is described schematically in
Fig.~\ref{comparison}. Each class shown in an oval is contained as
a proper subclass in the next larger oval, except for one case
(topological toric manifolds and quasitoric manifolds), where the
relation is slightly more subtle. Different examples are discussed
below.
\begin{figure}[h]
  \begin{picture}(120,55)
  \put(62.5,22.5){\oval(20,15)}
  \put(72.5,22.5){\oval(45,20)}
  \put(52.5,22.5){\oval(45,25)}
  \put(62.5,25){\oval(70,35)}
  \put(62.5,27.5){\oval(75,45)}
  \put(51.25,27.5){\oval(102.5,50)}
  \put(71.25,27.5){\oval(97.5,55)}
  \put(55.5,25){\small projective}
  \put(58.5,22){\small toric}
  \put(55.5,19){\small manifolds}
  \put(32.75,23.5){\small quasitoric}
  \put(33,20){\small manifolds}
  \put(81.5,23.5){\small toric}
  \put(78.5,20){\small manifolds}
  \put(44.5,37.5){\small topological toric manifolds}
  \put(41.5,45){\small torus manifolds with $H^{odd}=0$}
  \put(6.8,29){\small locally}
  \put(5.6,26){\small standard}
  \put(3.6,22.5){\small $T$-manifolds}
  \put(107.5,27){\small torus}
  \put(104.5,24){\small manifolds}
  \end{picture}
  \caption{Classes of $T$-manifolds.}
  \label{comparison}
\end{figure}

Projective toric manifolds are also Hamiltonian toric manifolds
(see Sections~\ref{symred} and~\ref{symred1}). However, when
viewed as symplectic manifolds, projective toric manifolds form a
smaller class: their symplectic forms represent integral
cohomology classes and their moment polytopes are lattice Delzant,
while arbitrary Delzant polytope can be realised as the moment
polytope of a Hamiltonian toric manifold.

A toric manifold (nonsingular compact toric variety) which is not
projective is described in Example~\ref{nonproj}.

A projective toric manifold is quasitoric by
Proposition~\ref{prqua}.

Many examples of quasitoric manifolds which are not toric can be
constructed using the equivariant connected sum operation
(Construction~\ref{equcs}). The simplest example is $\C P^2\cs\C
P^2$. It can be easily seen to be a quasitoric manifold over a
4-gon, but it does not admit an almost complex structure, and
therefore cannot be a complex algebraic variety. A non-toric
example with an invariant almost complex structure is given in
Exercise~\ref{qtnt}.

Examples of toric manifolds which is not quasitoric are
constructed by Suyama~\cite{suya}. The basic example is of real
dimension~4; its corresponding regular simplicial fan is obtained
by subdividing a singular fan whose underlying simplicial complex
is the Barnette sphere. More examples in arbitrary dimension can
be constructed by subsequent subdivision and suspension.

%There are no known examples of toric manifolds which are not
%quasitoric, although such examples are likely to exist: there is
%no reason why the orbit space of a (non-projective) toric manifold
%has to be combinatorially equivalent to a simple polytope. This
%question is equivalent to Problem~\ref{regnpt}.

Any quasitoric manifold $M$ is $T$-equivariantly homeomorphic to a
$T$-manifold obtained by restricting the $(\C^\times)^n$-action on
a topological toric manifold to the compact torus
$\T^n\subset(\C^\times)^n$ (see~\cite[Theorem~10.2]{i-f-m13}). The
easiest way to see this is to use the classification results
(Proposition~\ref{cqtpe} and Theorem~\ref{ttmtf}), and construct a
topological fan from the combinatorial quasitoric pair
$(P,\varLambda)$ corresponding to~$M$. To do this, consider any
convex realisation~\eqref{ptope1} of the polytope~$P$, and define
\[
  \beta_j=(\mb a_j+i\mb c_j,\lambda_j)\in\C^n\times\Z^n,\quad
  j=1,\ldots,m,
\]
where $\mb a_j$ are the normal vectors to the facets of~$P$, and
$\lambda_j$ are the columns of the characteristic
matrix~$\varLambda$. The vectors $\mb c_j\in\R^n$ can be chosen
arbitrarily. Then the data $(\sK_P;\beta_1,\ldots,\beta_m)$ define
a topological fan~$\Delta$. Indeed, condition~(a) from
Definition~\ref{deftopfan} is satisfied because $(\sK_P;\mb
a_1,\ldots,\mb a_m)$ define the normal fan~$\Sigma_P$, and
condition~(b) is equivalent to~\eqref{detLv}. Then the topological
toric manifold $X(\Delta)$ is $T$-homeomorphic to~$M$ and the
restriction of the $(\C^\times)^n$-action to the compact torus $T$
gives the $T$-action on~$M$ (this follows by comparing the
construction of $X(\Delta)$ with Proposition~\ref{kact}).

\begin{remark}
One would expect that the $T$-action on a quasitoric manifold $M$
can be extended to a $(\C^\times)^n$-action which gives $M$ a
structure of a topologically toric manifold. This stronger
statement would hold if one can replace an equivariant
homeomorphism by an equivariant diffeomorphism in
Proposition~\ref{equivar}.
\end{remark}

In \cite[\S9]{i-f-m13} there is constructed a topological toric
manifold $X$ whose associated simplicial complex $\sK(X)$ is the
Barnette sphere (see Construction~\ref{barsph}). Since the
Barnette sphere is not polytopal, this $X$ is not a quasitoric
manifold.

A toric manifold is topologically toric by definition. An example
of a topological toric manifold which is not toric can be
constructed from the quasitoric manifold $\C P^2\cs\C P^2$ as
described above. An explicit topological toric atlas on $\C
P^2\cs\C P^2$ and the corresponding topological fan are described
in~\cite[\S5,~Example]{i-f-m13}.

Any topological toric manifolds is a torus manifold with
$H^{odd}=0$ by Proposition~\ref{cohttm}. An even-dimensional
sphere is a torus manifold with $H^{odd}=0$
(Example~\ref{2nsphere}), but it is not a topological toric
manifold.

An example of a torus manifold with $H^{odd}\ne0$ can be
constructed as follows. Take any torus manifold $M$ whose quotient
$Q$ is face-acyclic. Let $R$ be any closed manifold which is not a
homology sphere. Consider the connected sum $\widehat Q=Q\cs R$
taken near an interior point of~$Q$. Then $\widehat Q$ is a
manifold with corners with $\partial\widehat Q=\partial Q$. In
particular, all proper faces of $\widehat Q$ are acyclic, but
$\widehat Q$ itself is not (we have $H^*(\widehat Q)\cong
H^*(R\setminus\pt)$). Consider the manifold $\widehat M=M(\widehat
Q,\lambda)$ constructed using the characteristic map of~$M$,
see~\eqref{idspa}. The singular $T$-orbits of $\widehat M$ are the
same as those of~$M$, but the free orbits are different. Now the
quotient of $\widehat M$ is $\widehat Q$, which is not
face-acyclic. Hence $H^{odd}(\widehat M)\ne0$ by
Theorem~\ref{cohfa} (this can be also easily seen directly). The
simplest example is obtained when $M=\C P^2$ and $R$ is a 2-torus.

Any torus manifold with $H^{odd}=0$ is locally standard by
Theorem~\ref{theo:local standardness}.

An example of a torus manifold which is not locally standard is
given in~\cite[\S11]{i-f-m13}. A free action of $T^n$ on the first
factor of a product manifold $T^n\times N^n$ gives an example of a
locally standard $T$-manifold which is not a torus manifold.

We conclude this section by mentioning that there are
\emph{real}\label{realtoriv} analogues of all classes of
$T$-manifolds considered here, in which the torus $\T^n$ is
replaced by the `real torus' $(\Z_2)^n$ and the algebraic torus
$(\C^\times)^n$ is replaced by $(\R^\times)^n$, where
$\R^\times\cong\R_>\times\Z_2$ is the multiplicative group of real
numbers. The `real' versions of the results of this chapter are
often simpler; the reader may recover the details of the proofs
himself. \emph{Real toric varieties} feature in \emph{tropical
geometry}~\cite{i-m-s09}. Real quasitoric manifolds are known as
\emph{small covers} of simple polytopes; they were introduced by
Davis and Januszkiewicz in~\cite{da-ja91} along with quasitoric
manifolds.

\section{Bounded flag manifolds}\label{bfman}
Bounded flag manifolds $\BF_n$ introduced by Buchstaber and Ray
in~\cite{bu-ra98e} and subsequently studied in~\cite{bu-ra98r}
and~\cite{bu-ra01}. Each $\BF_n$ is a projective toric manifold
whose moment polytope is combinatorially equivalent to an
$n$-cube, so that $\BF_n$ is also a quasitoric manifold over a
cube. Bounded flag manifolds are the examples of iterated
projective bundles, or \emph{Bott towers}, which are studied in
the next section. The manifolds $\BF_n$ find numerous application
in cobordism theory; they are implicitly present in the work of
Conner--Floyd~\cite{co-fl64} and in the construction of
\emph{Ray's basis}~\cite{ray86} in complex bordism of $\C
P^\infty$ (see details in Section~\ref{utg}), as well as used in
the construction of toric representatives in complex bordism
classes, described in Section~\ref{toricrep}. Bounded flag
manifolds also illustrate nicely many constructions and results
related to toric and quasitoric manifolds.

\begin{construction}[Bounded flag manifold]\label{bfm}
A \emph{bounded flag} in $\C^{n+1}$ is a complete flag
\[
  \mathcal U=\{U_1\subset U_2\subset\cdots\subset
  U_{n+1}=\C^{n+1},\quad \dim U_i=i\}
\]
for which $U_k$, \ $2\le k\le n$, contains the coordinate subspace
$\C^{k-1}=\<\mb e_1,\ldots,\mb e_{k-1}\>$ spanned by the first
$k-1$ standard basis vectors. Denote by $\BF_n$ the set of all
bounded flags in~$\C^{n+1}$.

Every bounded flag $\mathcal U$ in $\C^{n+1}$ is uniquely
determined by the set of $n$ lines
\begin{equation}\label{flaglines}
  \mathcal L=\{l_1,\ldots,l_n\colon
  \quad l_k\subset\C_k\oplus l_{k+1}
  \text{ for }1\le k\le n,\quad l_{n+1}=\C_{n+1}\},
\end{equation}
where $\C_k=\<\mb e_k\>$ is the $k$th coordinate line
in~$\C^{n+1}$. Indeed, given a set of lines $\mathcal L$ as above,
we can construct a bounded flag $\mathcal U$ by setting
$U_k=\C^{k-1}\oplus l_k$ for $1\le k\le n+1$. Conversely, the
conditions $l_k\subset\C_k\oplus l_{k+1}$ and $U_k=\C^{k-1}\oplus
l_k$ allows us to reconstruct the set of lines $\mathcal L$ from a
flag $\mathcal U$ in the reverse order $l_{n+1},l_n,\ldots,l_1$.
\end{construction}

\begin{theorem}\label{btoric}
The action of the algebraic torus $(\C^\times)^n$ on $\C^{n+1}$
given by
\[
  (t_1,\ldots,t_n)\cdot(w_1,\ldots,w_n,w_{n+1})=(t_1w_1,\ldots,t_nw_n,w_{n+1}),
\]
where $(t_1,\ldots,t_n)\in(\C^\times)^n$ and
$(w_1,\ldots,w_n,w_{n+1})\in\C^{n+1}$, induces an action on
bounded flags, and therefore makes $\BF_n$ into a smooth toric
variety.
\end{theorem}
\begin{proof}
We first construct a covering of $\BF_n$ by smooth affine charts
with regular change of coordinates functions, thereby giving
$\BF_n$ a structure of a smooth affine variety. We parametrise
bounded flags by sets of lines~\eqref{flaglines}. Let $\mb v_k$ be
a nonzero vector in~$l_k$, for $1\le k\le n$, and set $\mb
v_{n+1}=\mb e_{n+1}$ for the last line. Consider two collections
of $n$ opens subsets in~$\BF_n$:
\[
  V_k^0=\{\mathcal U\in \BF_n\colon l_k\ne \C_k\},\qquad
  V_k^1=\{\mathcal U\in \BF_n\colon \langle\mb v_k,\mb e_k\rangle\ne0\},
  \qquad 1\le k\le n.
\]
Now define $2^n$ open subsets
\[
  V^{\varepsilon_1,\ldots,\varepsilon_n}=V_1^{\varepsilon_1}\cap\cdots\cap V_n^{\varepsilon_n},\qquad
  \text{where }\varepsilon_k=0,1,
\]
Then $\{V^{\varepsilon_1,\ldots,\varepsilon_n}\}$ is a covering
of~$\BF_n$, because $V_k^0\cup V_k^1=\BF_n$ for any~$k$. The
condition $l_k\subset\C_k\oplus l_{k+1}$ implies
\begin{equation}\label{vvector}
  \mb v_k=z_k\mb e_k+z_{k+n}\mb v_{k+1}, \qquad 1\le k\le n,
\end{equation}
for some $z_i\in\C$, \ $1\le i\le 2n$. We have $z_k\ne0$ if
$\mathcal U\in V_k^1$, and $z_{k+n}\ne0$ if $\mathcal U\in V_k^0$.
Let $\mathcal U\in V^{\varepsilon_1,\ldots,\varepsilon_n}$; then
we can choose the vectors~\eqref{vvector} in the form $\mb
v_k=x^0_k\mb e_k+\mb v_{k+1}$ if $\varepsilon_k=0$, and $\mb
v_k=\mb e_k+x^1_k\mb v_{k+1}$ if $\varepsilon_k=1$, for $1\le k\le
n$. Then we can identify $V^{\varepsilon_1,\ldots,\varepsilon_n}$
with $\C^n$ using the affine coordinates
$(x^{\varepsilon_1}_1,\ldots,x^{\varepsilon_n}_n)$. The change of
coordinate functions are regular on intersections of charts by
inspection, so that $\BF_n$ is a smooth algebraic variety, with
affine atlas $\{V^{\varepsilon_1,\ldots,\varepsilon_n}\}$.

Furthermore, the change of coordinate functions are Laurent
monomials, which implies that $\BF_n$ is a toric variety. This can
also be seen directly, as the torus action defined in the theorem
is standard in the affine chart $V^{0,\ldots,0}$, that is,
\[
  (t_1,\ldots,t_n)\cdot(x^0_1,\ldots,x^0_n)=
  (t_1x^0_1,\ldots,t_nx^0_n).\qedhere
\]
\end{proof}

\begin{proposition}\label{bfan}
The complete fan $\Sigma$ corresponding to the toric variety
$\BF_n$ has $2n$ one-dimensional cones generated by the vectors
\[
  \mb a_k^0=\mb e_k,\quad
  \mb a_k^1=-\mb e_1-\cdots-\mb e_k,\quad 1\le k\le
  n,
\]
and $2^n$ maximal cones generated by the sets of vectors $\mb
a_1^{\varepsilon_1},\ldots,\mb a_n^{\varepsilon_n}$, where
$\varepsilon_k=0,1$.
\end{proposition}
\begin{proof}
Each affine chart $V^{\varepsilon_1,\ldots,\varepsilon_n}\subset
\BF_n$ constructed in the proof of Theorem~\ref{btoric}
corresponds to an $n$-dimensional cone
$\sigma^{\varepsilon_1,\ldots,\varepsilon_n}$ of the fan~$\Sigma$,
so there are $2^n$ maximal cones in total. One-dimensional cones
of $\Sigma$ correspond to $(\C^\times)^n$-invariant submanifolds
of complex codimension~1 in~$\BF_n$. Each of these submanifolds is
defined by vanishing of one of the affine coordinates, i.e. by an
equation $x^{\varepsilon_k}_k=0$, so there are $2n$ such
submanifolds.

In order to find the generators of the cone
$\sigma^{\varepsilon_1,\ldots,\varepsilon_n}$, we note that the
primitive generators of the dual cone
$(\sigma^{\varepsilon_1,\ldots,\varepsilon_n})^{\mathsf{v}}$ are
the weights of the $(\C^\times)^n$-representation in the affine
space $\C^n$ corresponding to the chart
$V^{\varepsilon_1,\ldots,\varepsilon_n}$. The
$(\C^\times)^n$-representation in the chart $V^{0,\ldots,0}$ is
standard, so we have $\mb a_k^0=\mb e_k$ for $1\le k\le n$. In
order to find the remaining vectors it is enough to calculate the
weights of the torus representation in the chart~$V^{1,\ldots,1}$.

The coordinates $(x_1^1,\ldots,x_n^1)$ in the chart
$V^{1,\ldots,1}$ are defined from the relations $\mb v_k=\mb
e_k+x^1_k\mb v_{k+1}$, $1\le k\le n$, and $\mb v_{n+1}=\mb
e_{n+1}$. An element $(t_1,\ldots,t_n)\in(\C^\times)^n$ acts on
$\mb e_k$ by multiplication by~$t_k$ for  $1\le k\le n$ and acts
on $\mb e_{n+1}$ identically (see Theorem~\ref{btoric}). Then it
is easy to see that the torus representation is written in the
coordinates $(x_1^1,\ldots,x_n^1)$ as follows:
\[
  (t_1,\ldots,t_n)\cdot (x_1^1,\ldots,x_n^1)=
  (t_1^{-1}t^{}_2x_1^1,\ldots,t_{n-1}^{-1}t^{}_nx^1_{n-1},t_n^{-1}x^1_n).
\]
In other words, the weights of this representation are the columns
of the matrix
\[
  W=\begin{pmatrix}
  -1& 0&\ldots&0&0\\
   1&-1&\ldots&0&0\\
  \vdots&\vdots&\ddots&\vdots&\vdots\\
   0& 0&\ldots&-1&0\\
   0& 0&\ldots&1&-1\\
  \end{pmatrix}.
\]
The generators $\mb a^1_1,\ldots,\mb a_n^1$ of the cone
$\sigma^{1,\ldots,1}$ form the dual basis, i.e. they are columns
of the matrix~$(W^{-1})^t$. These are the vectors listed in the
lemma.
\end{proof}

\begin{proposition}\label{bproj}
The bounded flag manifold $\BF_n$ is the projective toric variety
corresponding to the polytope
\[
  P=\bigl\{\mb x\in\R^n\colon x_k\ge0,\quad
  x_1+\cdots+x_k\le k, \quad\text{for }  1\le k\le n\bigr\}.
\]
\end{proposition}
\begin{proof}
We need to check that the normal fan of this $P$ is the fan from
Proposition~\ref{bfan}. Indeed, the $2n$ inequalities specifying
$P$ can be written as $\langle\mb a^0_k,\mb x\,\rangle\ge0$ and
$\langle\mb a^1_k,\mb x\,\rangle+k\ge0$, for $1\le k\le n$. Set
\begin{equation}\label{facetsbf}
  F^0_k=\{\mb x\in P\colon \langle\mb a^0_k,\mb x\,\rangle=0\}
  \quad\text{and}\quad
  F^1_k=\{\mb x\in P\colon \langle\mb a^1_k,\mb x\,\rangle+k=0\}.
\end{equation}
Then we need to check that
\begin{itemize}
\item[(a)] each $F^\varepsilon_k$ ($\varepsilon=0,1$) is a facet of~$P$;
\item[(b)] $F^0_k\cap F^1_k=\varnothing$,  for $1\le k\le n$;
\item[(c)] the intersection of any $n$-tuple
$F_1^{\varepsilon_1},\ldots,F_n^{\varepsilon_n}$ is a vertex
of~$P$.
\end{itemize}
This is left as an exercise.
\end{proof}

\begin{proposition}[{\cite%[Example~2.8]
{bu-ra01}}]\label{bqt}
The bounded flag manifold $\BF_n$ is a quasitoric manifold over a
combinatorial $n$-cube $I^n$, with characteristic matrix
\begin{equation}\label{bmatr}
  \varLambda=\left(\begin{matrix}
  1&0&\ldots&0\\
  0&1&\ldots&0\\
  \vdots&\vdots&\ddots&\vdots\\
  0&0&\ldots&1\\
  \end{matrix}\ \right|\left.
  \begin{matrix}
  -1&-1&\ldots&-1\\
   0&-1&\ldots&-1\\
  \vdots&\vdots&\ddots&\vdots\\
   0& 0&\ldots&-1\\
  \end{matrix}\right).
\end{equation}
\end{proposition}
\begin{proof}
Indeed, the polytope from Proposition~\ref{bproj} is
combinatorially equivalent to a cube.
\end{proof}

\begin{example}
The manifold $\BF_2$ is isomorphic to the Hirzebruch surface
$F_{1}$ (or $F_{-1}$) from Example~\ref{hirzebruch}.
\end{example}

We can also describe $\BF_n$ as a toric manifold using the
quotient construction (Section~\ref{algtq}) or symplectic
reduction (Section~\ref{symred}), as follows. The moment-angle
manifold corresponding to a cube~$I^n$ is a product of $n$
three-dimensional spheres:
\[
  \mathcal Z_{I^n}=
  \{(z_1,\ldots,z_{2n})\in\C^{2n}\colon|z_k|^2+|z_{k+n}|^2=1,\ 1\le
  k\le n\}.
\]
The manifold $\BF_n$ is obtained by taking quotient of $\mathcal
Z_{I^n}$ by the kernel $K$ of the map $\T^{2n}\to\T^n$ given by
matrix~\eqref{bmatr}. We have $K\cong\T^n$, and the inclusion
$K(\varLambda)\subset\T^{2n}$ is given by
\[
  (t_1,\ldots,t_n)\mapsto(t_1t_2\!\cdots t_{n-1}t_n,t_2\!\cdots
  t_{n-1}t_n,\ldots,t_{n-1}t_n,t_n,t_1,t_2,\ldots,t_n).
\]
Geometrically, the projection $\mathcal Z_{I^n}\to \BF_n$ maps
$\mb z=(z_1,\ldots,z_{2n})$ to the bounded flag defined by the set
of lines $l_1,\ldots,l_{n+1}$, where $l_k=\langle\mb v_k\rangle$
and $\mb v_k$ is given by~\eqref{vvector} (an exercise).

The algebraic quotient description of $\BF_n$ is very much
similar. Instead of the moment-angle manifold $\mathcal
Z_{I^n}\cong (S^3)^{n}$ we have the space
$U(\Sigma)=(\C^2\setminus\{\mathbf0\})^{n}$. The manifold $\BF_n$
is obtained by taking quotient of $U(\Sigma)$ by the kernel $G$ of
the map of algebraic tori $(\C^\times)^{2n}\to(\C^\times)^n$ given
by matrix~\eqref{bmatr}.

Now we describe the characteristic submanifolds and their
corresponding line bundles~\eqref{rhoi}. Let $\pi\colon \BF_n\to
P$ be the quotient projection for the torus action, and let
$\rho_k^{\varepsilon}$ denote the line bundle corresponding to the
characteristic submanifold (or $(\C^\times)^n$-invariant divisor)
$\pi^{-1}(F^\varepsilon_k)$, for $1\le k\le n$, \
$\varepsilon=0,1$, see~\eqref{facetsbf}.

\begin{proposition}[{\cite%[Example~2.8]
{bu-ra01}}]\label{linebunbt}\

\begin{itemize}
\item[(a)] The characteristic submanifold $\pi^{-1}(F_k^0)$ is
isomorphic to $\BF_{n-1}$, and $\pi^{-1}(F_k^1)$ is isomorphic to
$\BF_{k-1}\times \BF_{n-k}$, for $1\le k\le n$.

\item[(b)] The line bundle $\rho_k^0$ is isomorphic to the bundle
whose fibre over $\mathcal U\in \BF_n$ is the
line~$l_k=U_k/\C^{k-1}$. The line bundle $\rho_k^1$ is isomorphic
to the bundle whose fibre over $\mathcal U\in \BF_n$ is the
quotient $(\C_k\oplus l_{k+1})/l_k=U_{k+1}/U_k$.
%orthogonal complement to $l_k$ in $\C_k\oplus l_{k+1}$.
\end{itemize}
\end{proposition}
\begin{proof}
The submanifold $\pi^{-1}(F_k^0)\subset \BF_n$ is obtained by
projecting the submanifold of $\mathcal Z_{I^n}$ given by the
equation $z_k=0$ onto~$\BF_n$. If $z_k=0$, then the vectors $\mb
v_1,\ldots,\mb v_n$ defined by~\eqref{vvector} all belong to the
subspace $\C^{\{1,\ldots,\,n+1\}\setminus k}$. It follows that
$\pi^{-1}(F_k^0)$ can be identified with the set of bounded flags
in $\C^{\{1,\ldots,\,n+1\}\setminus k}$, that is,
with~$\BF_{n-1}$.

Similarly, the submanifold $\pi^{-1}(F_k^1)$ is the projection of
the submanifold of $\mathcal Z_{I^n}$ given by the equation
$z_k=1$. Then~\eqref{vvector} implies that the vectors $\mb
v_1,\ldots,\mb v_k$ belong to the subspace $\C^k\subset\C^{n+1}$,
and the vectors $\mb v_{k+1},\ldots,\mb v_n$ belong to the
subspace $\C^{\{k+1,\ldots,\,n+1\}}$. Therefore, $\pi^{-1}(F_k^1)$
can be identified with $\BF_{k-1}\times \BF_{n-k}$.

Statement (b) also follows from~\eqref{vvector}, because the line
bundle $\rho_k^0$ is isomorphic to $\mathcal Z_{I^n}\times_K\C_k$,
and $\rho_k^1$ is isomorphic to $\mathcal
Z_{I^n}\times_K\C_{n+k}$.
\end{proof}

\begin{proposition}\label{bbott}
The manifold $\BF_n$ is the complex projectivisation of the
complex plane bundle $\underline{\C}\oplus\rho_1^0$
over~$\BF_{n-1}$.
\end{proposition}
\begin{proof}
Consider the projection $\BF_n\to \BF_{n-1}$ taking a bounded flag
$\mathcal U=\{U_1\subset U_2\subset\cdots\subset
U_n\subset\C^{n+1}\}$ to the flag $\mathcal U'=\mathcal U/\C_1$ in
$\C^{2,\ldots,\,n+1}\cong\C^n$. (More precisely, $\mathcal
U'=\{U'_1\subset U'_2\subset\cdots\subset U'_{n-1}\}$, where
$U'_{k}=U_{k+1}/\C_1$, \ $1\le k\le n-1$.) The set of
lines~\eqref{flaglines} corresponding to $\mathcal U'$ is obtained
from the set of lines corresponding to $\mathcal U$ by forgetting
the first line. In order to recover the flag $\mathcal U$ from the
flag $\mathcal U'$, one needs to choose a line $l_1$ in the plane
$\C_1\oplus l_2$. Since $l_2$ is the first line in the set
corresponding to the flag $\mathcal U'\in \BF_{n-1}$, we obtain
$\BF_n\cong\C P(\underline{\C}\oplus\rho_1^0)$, as needed.
\end{proof}

We therefore obtain a tower of fibrations $\BF_n\to
\BF_{n-1}\to\cdots\to \BF_1=\C P^1$, where each $\BF_k$ is the
projectivisation of a complex 2-plane bundle over~$\BF_{k-1}$. In
particular, the fibre of each bundle in the tower is~$\C P^1$.
Towers of fibrations arising in this way are called \emph{Bott
towers}; they are the subject of the next section.

\subsection*{Exercises}

\begin{exercise}
The fan described in Proposition~\ref{bfan} is the normal fan of
the polytope from Proposition~\ref{bproj}.
\end{exercise}

\begin{exercise}
Given $\mb z=(z_1,\ldots,z_{2n})$, define the vectors $\mb
v_{n+1}=\mb e_{n+1},\mb v_n,\ldots,\mb v_1$ by~\eqref{vvector},
and set $l_k=\langle\mb v_k\rangle$, \ $1\le k\le n+1$. Then the
projection $\mathcal Z_{I^n}\to \BF_n$ maps $\mb z\in\mathcal
Z_{I^n}$ to the bounded flag in~$\C^{n+1}$ defined by the set of
lines $l_1,\ldots,l_{n+1}$.
\end{exercise}

\section{Bott towers}\label{bott}
In their study of symmetric spaces, Bott and
Samelson~\cite{bo-sa58} introduced a family of complex manifolds
obtained as the total spaces of iterated bundles over $\mathbb C
P^1$ with fibre~$\mathbb C P^1$. Grossberg and
Karshon~\cite{gr-ka94} showed that these manifolds carry an
algebraic torus action, and therefore constitute an important
family of smooth projective toric varieties, and called them Bott
towers. Civan and Ray~\cite{ci-ra05} developed significantly the
study of Bott towers by enumerating the invariant stably complex
structures and calculating their complex and real $K$-theory
rings, and cobordism.

Each Bott tower is a projective toric manifold whose corresponding
simple polytope is combinatorially equivalent to a cube (a
\emph{toric manifold over cube}\label{tmancube} for short). We
have the following hierarchy of classes of $T$-manifolds:
\[
  \text{Bott towers }
  \subset\text{ toric manifolds over cubes }
  \subset\text{ quasitoric manifolds over cubes}
\]
By the result of Dobrinskaya~\cite{dobr01}, the first inclusion
above is in fact an identity (we explain this in
Corollary~\ref{torbt}).

Two results were obtained in~\cite{ma-pa08} relating circle
actions on Bott towers, their topological structure, and
cohomology rings. First (Theorem~\ref{btpro}), if a Bott tower
admits a semifree $\mathbb S^1$-action with isolated fixed points,
then it is $\mathbb S^1$-equivariantly diffeomorphic to a product
of 2-spheres. Second (Theorem~\ref{trivc}), a Bott tower whose
cohomology ring is isomorphic to that of a product of spheres is
actually diffeomorphic to this product. Both theorems can be
extended to quasitoric manifolds over cubes, but only in the
topological category (Theorems~\ref{qtmsf} and~\ref{qtcpn}). These
results have been further extended by several authors, and led to
the study of the so-called \emph{cohomological rigidity} property
for different classes of manifolds with torus actions; we discuss
this circle of problems in the end of this section.

\subsection*{Definition and main properties}
\begin{definition}\label{defbotttower}
A \emph{Bott tower} of height $n$ is a tower of fibre bundles
\[
  B_n\stackrel{p_n}\longrightarrow B_{n-1}
  \stackrel{p_{n-1}}\longrightarrow\cdots
  \longrightarrow B_2\stackrel{p_2}\longrightarrow B_1\longrightarrow\pt,
\]
of complex manifolds, where $B_1=\C P^1$ and $B_k=\C
P(\underline{\C}\oplus\xi_{k-1})$ for $2\le k\le n$. Here $\C
P(\cdot)$ denotes complex projectivisation, $\xi_{k-1}$ is a
complex line bundle over~$B_{k-1}$ and $\underline{\C}$ is a
trivial line bundle. The fibre of the bundle $p_k\colon B_k\to
B_{k-1}$ is~$\C P^1$.

A Bott tower $B_n$ is said to be \emph{topologically
trivial}\label{btowetri} if each $p_k\colon B_k\to B_{k-1}$ is
trivial as a smooth fibre bundle; in particular, such $B_n$ is
diffeomorphic to a product of 2-dimensional spheres.

We shall refer to the last stage $B_n$ in a Bott tower as a
\emph{Bott manifold}\label{defbottmani} (although it is also often
called by the same name `Bott tower').
\end{definition}

In order to describe the cohomology ring of a Bott manifold, we
need the following general result:

\begin{theorem}[{see \cite[Chapter~V]{ston68}}]\label{cohomproj}
Let $\xi$ be a complex $n$-dimensional vector bundle over a finite
cell complex~$X$ with complex projectivisation $\C P(\xi)$, and
let $u\in H^2(\C P(\xi))$ be the first Chern class of the
tautological line bundle over $\C P(\xi)$. The integral cohomology
ring of $\C P(\xi)$ is the quotient of the polynomial ring
$H^*(X)[u]$ on a generator $u$ with coefficients in $H^*(X)$ by
the single relation
\[
  u^n-c_1(\xi)u^{n-1}+\cdots+(-1)^nc_n(\xi)=0.
\]
\end{theorem}

\begin{corollary}\label{cohbt}
$H^*(B_k)$ is a free module over $H^*(B_{k-1})$ on generators $1$
and $u_k$, where $u_k$ is the first Chern class of the
tautological line bundle over $B_k=\C
P(\underline{\C}\oplus\xi_{k-1})$. The ring structure is
determined by the single relation
\[
  u_k^2=c_1(\xi_{k-1})u_k.
\]
\end{corollary}

For simplicity, we denote the element $p_k^*(u_{k-1})\in H^2(B_k)$
by~$u_{k-1}$; similarly, we denote by $u_i$ each of the elements
in $H^*(B_k)$, \ $k\ge i$, which map to each other by the
homomorphisms~$p_k^*$. Each line bundle $\xi_{k-1}$ is determined
by its first Chern class, which can be written as a linear
combination
\[
  c_1(\xi_{k-1})=a_{1k}u_1+a_{2k}u_2+\cdots+a_{k-1,k}u_{k-1}\in
  H^2(B_{k-1}).
\]
It follows that a Bott tower of height $n$ is uniquely determined
by the list of integers $\{a_{ij}\colon 1\le i<j\le n\}$, where
\begin{equation}\label{crobt}
  u_k^2=\sum_{i=1}^{k-1}a_{ik}u_iu_k, \qquad 1\le k\le n.
\end{equation}
The cohomology ring of $B_n$ is the quotient of
$\Z[u_1,\ldots,u_n]$ by relations~\eqref{crobt}.

It is convenient to organise the integers $a_{ij}$ into an upper
triangular matrix,
\begin{equation}\label{amatr}
  A=\begin{pmatrix}
    -1 & a_{12} & \cdots & a_{1n}\\
    0  & -1     & \cdots & a_{2n}\\
    \vdots & \vdots & \ddots & \vdots\\
    0 & 0 & \cdots & -1
  \end{pmatrix}.
\end{equation}

\begin{example}\label{2dbto}
Let $n=2$. Then a Bott tower $B_2\to B_1$ is determined by a line
bundle $\xi_1$ over $B_1=\C P^1$, i.e. $B_2$ is a Hirzebruch
surface (see Example~\ref{hirzebruch}). We have $\xi_1=\eta^k$ for
some $k\in\Z$ where $\eta^k$ denotes the $k$th tensor power of the
tautological line bundle over $\C P^1$. The cohomology ring is
given by
\[
  H^*(B_2)=\Z[u_1,u_2]/(u_1^2,\, u_2^2-ku_1u_2).
\]
We have
\[
  \C P(\underline{\C}\oplus\eta^k)\cong
  \C P(\underline{\C}\oplus\eta^{k'})
  \quad\Leftrightarrow\quad k=k'
  \mod 2,
\]
where $\cong$ denotes a diffeomorphism. This is proved by the
following observation. First, note that $\C P(\xi)\cong \C
P(\xi\otimes\eta)$ for any complex vector bundle $\xi$ and line
bundle~$\eta$. Let $k'-k=2\ell$ for some $\ell\in \Z$, then
\[
  \C P(\underline{\C}\oplus\eta^k)\cong
  \C P((\underline{\C}\oplus\eta^k)\otimes\eta^{\ell})=
  \C P(\eta^\ell\oplus\eta^{k+\ell})\cong
  \C P(\underline{\C}\oplus\eta^{k'}),
\]
where the last diffeomorphism is induced by the vector bundle
isomorphism
$\eta^\ell\oplus\eta^{k+\ell}\cong\underline{\C}\oplus\eta^{k'}$,
as both are plane bundles over $\C P^1$ with equal Chern classes.

On the other hand, a cohomology ring isomorphism $H^*(\C
P(\underline{\C}\oplus\eta^k)) \cong H^*(\C
P(\underline{\C}\oplus\eta^{k'}))$ implies that $k=k'\mod2$ (an
exercise).
\end{example}

This example shows that the cohomology ring determines the
diffeomorphism type of a Bott manifold $B_n$ for $n=2$. We may ask
if this is true for arbitrary~$n$; the questions of this sort are
discussed in the last subsection.

\begin{example}\label{btbfn}
The bounded flag manifold $\BF_n$ is a Bott manifold. By
Proposition~\ref{bbott}, $\BF_n=\C
P(\underline{\C}\oplus\rho_1^0)$, where $\rho_1^0$ is the line
bundle over~$\BF_{n-1}$ whose fiber over a bounded flag $\mathcal
U$ is its first space~$U_1$. By Corollary~\ref{cohbt}, the ring
structure of $H^*(\BF_n)$ is determined by the relation
$u_n^2=c_1(\rho_1^0)u_n$. As it is clear from the proof of
Proposition~\ref{bbott}, $\rho_1^0$ is the tautological line
bundle over~$\BF_{n-1}$ (considered as a complex
projectivisation), so $c_1(\rho_1^0)=u_{n-1}$. (Warning: the line
bundle $\rho_1^0$ over $\BF_n$ is \emph{not} the pullback of the
line bundle $\rho_1^0$ over $\BF_{n-1}$ by the projection
$p_n\colon \BF_n\to \BF_{n-1}$, because
$p_n^*(\rho^0_1)=\rho^0_2$.)

We therefore obtain the identity $u_n^2=u_{n-1}u_n$ in
$H^*(\BF_n)$, and matrix~\eqref{amatr} for the structure of the
Bott tower on the bounded flag manifold has the form
\[
  A=\begin{pmatrix}
  -1 &  1 & 0  & \cdots & 0\\
  0  & -1 & 1  & \cdots & 0\\
  \vdots & \vdots & \ddots & \ddots & \vdots\\
  0 & 0 & 0 & \cdots & 1\\
  0 & 0 & 0 & \cdots & -1
  \end{pmatrix}.
\]
\end{example}

As a corollary of this example we obtain the following
characterisation of bounded flag manifolds:

\begin{proposition}\label{charbfmbt}
The bounded flag manifold $\BF_n$ is the Bott manifold whose tower
structure is defined as follows: in each $B_k=\C
P(\underline\C\oplus\xi_{k-1})$ the bundle $\xi_{k-1}$ is the
tautological line bundle over $B_{k-1}=\C
P(\underline\C\oplus\xi_{k-2})$, for $2\le k\le n$.
\end{proposition}

\subsection*{Bott towers as toric manifolds}
\begin{theorem}\label{almat}
The Bott manifold $B_n$ corresponding to a matrix $A$ given
by~\eqref{amatr} is isomorphic to the toric manifold corresponding
to the complete fan $\Sigma$ with $2n$ one-dimensional cones
generated by the vectors
\[
  \mb a_k^0=\mb e_k,\quad
  \mb a_k^1=-\mb e_k+a_{k,k+1}\mb e_{k+1}+\cdots+a_{kn}\mb e_n\quad (1\le k\le
  n),
\]
and $2^n$ maximal cones generated by the sets of vectors $\mb
a_1^{\varepsilon_1},\ldots,\mb a_n^{\varepsilon_n}$, where
$\varepsilon_k=0,1$.
\end{theorem}
\begin{proof}
Let $X_n$ denote the toric manifold corresponding to the fan
described in the theorem. We may assume by induction that
$X_{n-1}=B_{n-1}$ (the base of the induction is clear, as
$X_1=B_1=\C P^1$). By the construction of Section~\ref{algtq}, the
manifold $X_n$ can be obtained as the quotient of
\[
  U_n=
  \{(z_1,\ldots,z_{2n})\in\C^{2n}\colon|z_k|^2+|z_{k+n}|^2\ne0,\ 1\le
  k\le n\}\cong(\C^2\setminus\{\mathbf 0\})^n
\]
by the action of the group $G_n\cong(\C^\times)^n$ given
by~\eqref{gexpl} (we have $m=2n$ here). Explicitly, the inclusion
$G_n\to(\C^\times)^{2n}$ is given by
\[
  (t_1,\ldots,t_n)\mapsto
  (t_1,\:t_1^{-a_{12}}t_2,\:\ldots,\:t_1^{-a_{1n}}t_2^{-a_{2n}}\!\cdots
  t_{n-1}^{-a_{n-1,n}}t_n,\:t_1,\:t_2,\:\ldots,\:t_n).
\]
Observe that $U_n=U_{n-1}\times(\C^2\setminus\{\mathbf0\})$,
$G_n=G_{n-1}\times\C^\times$, and the last factor $\C^\times$
(corresponding to~$t_n$) acts trivially on $U_{n-1}$. Therefore,
we have
\begin{multline*}
  X_n=U_n/G_n
  %=\bigl(U_{n-1}\times(\C^2\setminus\{\mathbf0\})\bigr)\big/(G_{n-1}\times\C^\times)
  =\bigl(U_{n-1}\times(\C^2\setminus\{\mathbf0\})/\C^\times\bigr)/G_{n-1}\\
  =U_{n-1}\times_{G_{n-1}}\C P^1=
  \C P\bigl(U_{n-1}\times_{G_{n-1}}(\C\oplus\C)\bigr),
\end{multline*}
where $U_{n-1}\times_{G_{n-1}}(\C\oplus\C)$ is a complex 2-plane
bundle over $U_{n-1}/G_{n-1}=X_{n-1}$ defined by the
representation of the algebraic torus
$G_{n-1}\cong(\C^\times)^{n-1}$ in $\C\oplus\C$ which is given by
the character $(t_1,\ldots,t_{n-1})\mapsto
t_1^{-a_{1n}}t_2^{-a_{2n}}\!\cdots t_{n-1}^{-a_{n-1,n}}$ on the
first summand and is trivial on the second summand.

By the inductive assumption, $X_{n-1}=B_{n-1}$. The bundle
$U_{n-1}\times_{G_{n-1}}(\C\oplus\C)$ is
$\xi_{n-1}\oplus\underline{\C}$, where $\xi_{n-1}$ is the line
bundle over $B_{n-1}$ with first Chern class
$\sum_{i=1}^{n-1}a_{in}u_i$. Thus, $X_n=\C
P(\xi_{n-1}\oplus\underline{\C})=B_n$, and the inductive step is
complete.
\end{proof}

Note that the associated simplicial complex of the fan described
in Theorem~\ref{almat} is the boundary of a cross-polytope, so the
toric manifold is also a quasitoric manifold over a cube. We
therefore obtain:

\begin{corollary}\label{bottquasi}
A Bott tower of height $n$ determined by matrix $A$ has a natural
action of the torus $T^n$ making it into a quasitoric manifold
over a cube with~refined characteristic matrix
$\varLambda=(I\;|\;A^t)$, see~\eqref{lamat}.
%
%Conversely, a quasitoric manifold over a cube with refined
%characteristic matrix $(I\;|\;\varLambda_\star)$, where
%$\varLambda_\star$ is a lower triangular square matrix, has a
%structure of a Bott tower.
\end{corollary}

\begin{remark}
The bounded flag manifold $\BF_n$ has a toric structure described
in Proposition~\ref{bfan}, and another toric structure coming from
its Bott tower structure via Theorem~\ref{almat} (its matrix $A$
is given in Example~\ref{btbfn}). The corresponding fans are not
the same, but isomorphic (see Exercise \ref{bfm2toric}).
\end{remark}

\begin{remark}
The relations in the face ring of an $n$-cube are $v_iv_{i+n}=0$,
\ $1\le i\le n$. These relations together with~\eqref{viref} give
the relations~\eqref{crobt} after substitution
$\varLambda=(I\;|\;A^t)$ and $u_i=-v_{i+n}$. In fact, one has
$u_i=c_1(\bar\rho_{i+n})$, where $\rho_{i+n}$ is the line
bundle~\eqref{rhoi} over the toric manifold~$B_n$ (an exercise).
It follows that the description of the cohomology ring of a Bott
manifold from Corollary~\ref{cohbt} agrees with the description of
the cohomology of a toric manifold from Theorem~\ref{qtcoh}.
\end{remark}

Given a permutation $\sigma$ of $n$ elements, denote by
$P(\sigma)$ the corresponding \emph{permutation matrix}, the
square matrix of size $n$ with ones at the positions
$(i,\sigma(i))$ for $1\le i\le n$, and zeros elsewhere. There is
an action of the symmetric group $S_n$ on square $n$-matrices by
conjugations, $A\mapsto P(\sigma)^{-1}AP(\sigma)$, or,
equivalently, by permutations of the rows and columns of~$A$.

\begin{proposition}\label{qtbtc}
A quasitoric manifold $M$ over a cube with refined characteristic
matrix $\varLambda=(I\;|\;\varLambda_\star)$ is equivalent to a
Bott manifold if and only if $\varLambda_\star$ is conjugate by
means of a permutation matrix to an upper triangular matrix.
\end{proposition}
\begin{proof}
Assume that $\varLambda_\star$ is conjugate by means of a
permutation matrix to an upper triangular matrix. Clearly, this
condition is equivalent to the conjugacy of $\varLambda_\star$ to
a lower triangular matrix. Consider the action of $S_n$ on the set
of facets of the cube $\mathbb I^n$ by permuting pairs of opposite
facets. A rearrangement of facets corresponds to a rearrangement
of columns in the characteristic $n\times2n$-matrix $\varLambda$,
so an element $\sigma\in S_n$ acts as
\[
  \varLambda\mapsto\varLambda\cdot
  \begin{pmatrix}P(\sigma)&0\\0&P(\sigma)\end{pmatrix}.
\]
This action does not preserve the refined form of $\varLambda$, as
$(I\;|\;\varLambda_\star)$ becomes
$(P(\sigma)\;|\;\varLambda_\star P(\sigma))$. The refined
representative in the left coset~\eqref{lcset} of the latter
matrix is given by $(I\;|\;P(\sigma)^{-1}\varLambda_\star
P(\sigma))$. (In other words, we must compensate for the
permutation of pairs of facets by an automorphism of the torus
$\mathbb T^n$ permuting the coordinate subcircles to keep the
characteristic matrix in the refined form.) This implies that the
action by permutations on pairs of opposite facets induces an
action by conjugations on refined submatrices~$\varLambda_\star$.
Hence we may assume, up to an equivalence, that the refined
characteristic submatrix~$\varLambda_\star$ of $M$ is lower
triangular. The non-singularity condition~\eqref{detLv} guarantees
that the diagonal entries of~$\varLambda_\star$ are equal to
$\pm1$, and we can set all of them equal to $-1$ by changing the
omniorientation of $M$ if necessary. Now, $M$ has the same
characteristic matrix as the Bott manifold corresponding to the
matrix $A=\varLambda_\star^t$ (see Corollary~\ref{bottquasi}).
Therefore, $M$ and the Bott manifold are equivalent by
Proposition~\ref{cqtpe}.

The converse statement follows from Corollary~\ref{bottquasi}.
\end{proof}

Our next goal is to characterise Bott manifolds within the class
of quasitoric manifolds over cubes more explicitly. Given a subset
$\{i_1,\ldots,i_k\}\subset[n]$, the \emph{principal
minor}\label{priminor} of a square $n$-matrix $A$ is the
determinant of the submatrix formed by the elements in columns and
rows with numbers $i_1,\ldots,i_k$. In the case of Bott manifolds,
according to Corollary~\ref{bottquasi}, all principal minors of
the matrix $-\varLambda_\star$ are equal to~1; for an arbitrary
quasitoric manifold the non-singularity condition~\eqref{detLv}
only guarantees that all principal minors of $\varLambda_\star$
are equal to~$\pm1$.

Recall that an upper triangular matrix is
\emph{unipotent}\label{dunipot} if all its diagonal entries are
ones. The following key technical lemma can be retrieved from the
proof of Dobrinskaya’s general result~\cite[Theorem~6]{dobr01}
characterising quasitoric manifolds over products of simplices
which can be decomposed into towers of fibrations.

\begin{lemma}%[\cite{dobr01}]
\label{minor} Let $R$ be a commutative integral domain with
identity element~$1$, and let $A$ be a square $n$-matrix ($n\ge2$)
with entries in~$R$. Suppose that every proper principal minor
of~$A$ is equal to~$1$. If $\det A=1$, then A is conjugate by
means of a permutation matrix to a unipotent upper triangular
matrix, otherwise it is permutation-conjugate to a matrix of the
following form:
\begin{equation}\label{exce}
\begin{pmatrix}
1 & b_1 & 0 & \dots &  0\\
0 & 1 & b_2 & \dots &  0\\
\vdots& \vdots  &\ddots &\ddots & \vdots \\
0 & 0 & \dots & 1 &  b_{n-1}\\
b_n & 0 & \dots & 0  & 1
\end{pmatrix}
\end{equation}
where $b_i\not=0$ for all $i$.
\end{lemma}
\begin{proof}
By assumption the diagonal entries of $A$ must be ones. We say
that the $i$th row is \emph{elementary} if its $i$th entry is~$1$
and the other entries are~$0$. Assuming by induction that the
theorem holds for matrices of size $(n-1)$ we deduce that $A$ is
itself permutation-conjugate to a unipotent upper triangular
matrix if and only if it contains an elementary row. We denote by
$A_i$ the square $(n-1)$-matrix obtained by removing from $A$ the
$i$th column and the $i$th row.

We may assume by induction that $A_n$ is a unipotent upper
triangular matrix. Next we apply the induction assumption
to~$A_1$. The permutation of rows and columns transforming $A_1$
into a unipotent upper triangular matrix turns $A$ into an
`almost' unipotent upper triangular matrix; the latter may have
only one non-zero entry below the diagonal, which must be in the
first column. If $a_{n1}=0$, then the $n$th row of $A$ is
elementary and $A$ is permutation-conjugate to a unipotent upper
triangular matrix. Otherwise we have
\begin{equation*}
  A=
  \begin{pmatrix}
    1 & * & * & \cdots &  *\\
    0 & 1 & * & \cdots &  *\\
    \vdots& \vdots &\ddots &\ddots & \vdots \\
    0 & 0 & \cdots & 1 &  b_{n-1}\\
    b_n & 0 & \cdots & 0  & 1
  \end{pmatrix},
\end{equation*}
where $b_{n-1}\ne0$ and $b_n\ne0$ (otherwise $A$ contains an
elementary row). Now let $a_{1j_1}$ be the last non-zero entry in
the first row of~$A$. If $A$ does not contain an elementary row,
then we may define by induction $a_{j_ij_{i+1}}$ as the last
non-zero non-diagonal entry in the $j_i$th row of~$A$. Clearly, we
have
\[
  1<j_1<\cdots<j_i<j_{i+1}<\cdots<j_k=n
\]
for some $k<n$. Now, if $j_i=i+1$ for $1\le i\le n-1$, then $A$ is
the matrix~\eqref{exce} with $b_i=a_{j_{i-1}j_i}$, $1\le i\le
n-1$. Otherwise, the submatrix
\begin{equation*}
  S=
  \begin{pmatrix}
    1 & a_{1j_1} & 0 & \dots &  0\\
    0 & 1 & a_{j_1j_2} & \dots &  0\\
    \vdots& \vdots &\ddots &\ddots & \vdots \\
    0 & 0 & \dots & 1 & a_{j_{k-1}n} \\
    b_n & 0 & \dots & 0  & 1
  \end{pmatrix},
\end{equation*}
of $A$ formed by the columns and rows with indices
$1,j_1,\ldots,j_k$ is proper and has determinant $1\pm b_n\prod
a_{j_ij_{i+1}}\ne1$. This contradiction finishes the proof.
\end{proof}

\begin{theorem}\label{qtcub}
Let $M=M(I^n,\varLambda)$ be a quasitoric manifold over a cube
with canonical $T^n$-invariant smooth structure, and
$\varLambda_\star$ the corresponding refined submatrix. Then the
following conditions are equivalent:
\begin{itemize}
\item[(a)] $M$ is equivalent to a Bott manifold;

\item[(b)] all the principal minors of $-\varLambda_\star$ are equal to~$1$;

\item[(c)] $M$ admits a $T^n$-invariant almost complex structure (with the associated omniorientation).
\end{itemize}
\end{theorem}
\begin{proof}
The implication $\text{(b)}\Rightarrow\text{(a)}$ follows from
Lemma~\ref{minor} and Proposition~\ref{qtbtc}. The implication
$\text{(a)}\Rightarrow\text{(c)}$ is obvious. Let us prove
$\text{(c)}\Rightarrow\text{(b)}$. Recall the definition of the
sign $\sigma(v)$ and the formula from Lemma~\ref{qts}~(b)
expressing this sign in terms of the combinatorial data. Denote
the facets of the cube $I^n$ by
$F^\varepsilon_1,\ldots,F^\varepsilon_n$ ($\varepsilon=0,1$),
assuming that $F^0_k\cap F^1_k=\varnothing$, for $1\le k\le n$.
The normal vectors of facets are $\mb
a^{\varepsilon}_k=(-1)^\varepsilon\mb e_k$. A vertex of $I^n$ is
given by
\[
  v=F_1^{\varepsilon_1}\cap\cdots\cap F_n^{\varepsilon_n}.
\]
Therefore, the expression for the sign $\sigma(v)$ on the right
hand side of the formula from Lemma~\ref{qts} is equal to a
principal minor of the matrix $-\varLambda_\star$ (namely, the
minor formed by the columns and rows with numbers $i$ such that
$v\in F^1_i$). It remains to note that in the almost complex case
the sign of every vertex is~$1$.
\end{proof}

\begin{remark}
The equivalence (a)$\Leftrightarrow$(b) is a particular case
of~\cite[Theorem~6]{dobr01}.
\end{remark}

\begin{corollary}\label{torbt}
Let $V$ be a toric manifold whose associated fan is
combinatorially equivalent to the fan consisting of cones over the
faces of a cross-polytope. Then $V$ is a Bott manifold.
\end{corollary}
\begin{proof}
If we view $V$ as a quasitoric manifold (over a cube), then all
the principal minors of the corresponding matrix
$\varLambda_\star$ are equal to~$1$ by the same reason as in the
proof of Theorem~\ref{qtcub}. By Lemma~\ref{minor}, the matrix
$\varLambda_\star$ is permutation-conjugate to a unipotent upper
triangular matrix, so the full characteristic matrix $\varLambda$
has the same form as the characteristic matrix of a Bott manifold.
The columns of $\varLambda$ are the primitive vectors along edges
of the fan corresponding to~$V$, so the combinatorial type of the
fan and the matrix $\varLambda$ determine the fan completely. It
follows that the fan of $V$ is the same as the fan of some Bott
manifold, which implies that $V$ has the structure of a Bott
tower.
\end{proof}

\begin{example}\label{qtovercube}
Corollary~\ref{torbt} shows that the class of Bott manifolds
coincides with the class of toric manifolds over cubes, i.e. the
first inclusion in the hierarchy described in the beginning of
this section is an identity. This is not the case for the second
inclusion. For example the quasitoric manifold $M$ over a square
with refined characteristic submatrix
$\varLambda_\star=\begin{pmatrix}1&2\\1&1\end{pmatrix}$ is not a
Bott manifold, because $\varLambda_\star$ is not
permutation-conjugate to an upper triangular matrix. This $M$ is
diffeomorphic to $\C P^2\mathbin{\#}\C P^2$.
\end{example}

\subsection*{Semifree circle actions}
Recall that an action of a group is called
\emph{semifree}\label{semifreeHam} if it is free on the complement
to fixed points. A particularly interesting class of
\emph{Hamiltonian} semifree circle actions was studied by Hattori,
who proved in~\cite{hatt92} that a compact symplectic manifold $M$
carrying a semifree Hamiltonian $\mathbb S^1$-action with nonempty
isolated fixed point set has the same cohomology ring and the same
Chern classes as $\C P^1\times\cdots\times\C P^1$, thus imposing a
severe restriction on the topological structure of the manifold.
Hattori's results were further extended by Tolman and Weitsman,
who showed in~\cite{to-we00} that a semifree symplectic $\mathbb
S^1$-action with nonempty isolated fixed point set is
automatically Hamiltonian, and the \emph{equivariant} cohomology
ring and Chern classes of $M$ also agree with those of $\C
P^1\times\cdots\times\C P^1$. In dimensions up to 6 it is known
that a symplectic manifold with a $\mathbb S^1$-action satisfying
the properties above is diffeomorphic to a product of 2-spheres,
but in higher dimensions this remains open.

Ilinskii considered in~\cite{ilin06} an algebraic version of
Hattori's question on semifree symplectic $\mathbb S^1$-actions.
Namely, he conjectured that a smooth compact complex algebraic
variety $V$ carrying a semifree action of the algebraic 1-torus
$\C^\times$ with positive number of isolated fixed points is
homeomorphic to $S^2\times\cdots\times S^2$. The algebraic and
symplectic versions of the conjecture are related via the common
subclass of projective varieties; a smooth projective variety is a
symplectic manifold. Ilinskii proved the \emph{toric} version of
his algebraic conjecture, namely, when $V$ is a toric manifold and
the semifree 1-torus is a subgroup of the acting torus (of
dimension $\dim_\C V$). The first step of Ilinskii's argument was
to show that if $V$ admits a semifree action of a subcircle with
isolated fixed points, then the corresponding fan is
combinatorially equivalent to the fan over the faces of a
cross-polytope. By Corollary~\ref{torbt}, such a toric manifold
$V$ admits a structure of a Bott tower.

The above described classification of Bott towers can therefore be
applied to obtain results on semifree circle actions. According to
a result of~\cite{ma-pa08} (Theorem~\ref{qtbt} below), a
quasitoric manifold over a cube with a semifree circle action is a
Bott tower. By another result of~\cite{ma-pa08}, all such Bott
towers are topologically trivial, i.e. diffeomorphic to a product
of 2-dimensional spheres (Theorem~\ref{btpro}).

A complex $n$-dimensional representation of the circle $\mathbb
S^1$ is determined by the set of weights $k_j\in\Z$, \ $1\le j\le
n$. In appropriate coordinates an element $s=e^{2\pi i\varphi}\in
\mathbb S^1$ acts as follows:
\begin{equation}\label{cirac}
  s\cdot(z_1,\ldots,z_n)=(e^{2\pi ik_1\varphi}z_1,\ldots,e^{2\pi
  ik_n\varphi}z_n).
\end{equation}
The following result is straightforward.

\begin{proposition}\label{sfrep}
A representation of $\mathbb S^1$ in $\C^n$ is semifree if and
only if $k_j=\pm1$ for $1\le j\le n$.
\end{proposition}

Let $M=M(P,\varLambda)$ be a quasitoric manifold. A circle
subgroup in $\mathbb T^n$ is determined by a primitive integer
vector $\nu=(\nu_1,\ldots,\nu_n)$:
\begin{equation}\label{cirnu}
  S(\nu)=\{(e^{2\pi i\nu_1\varphi},\ldots,e^{2\pi i\nu_n\varphi})\in\mathbb T^n\colon
  \varphi\in\R\}.
\end{equation}
Given a vertex $v=F_{j_1}\cap\cdots\cap F_{j_n}$ of~$P$, we can
decompose $\nu$ in terms of the basis
$\lambda_{j_1},\ldots,\lambda_{j_n}$:
\begin{equation}\label{nucoe}
  \nu=k_1(\nu,v)\lambda_{j_1}+\cdots+k_n(\nu,v)\lambda_{j_n}.
\end{equation}

\begin{proposition}\label{sfcoe}
A circle $S(\nu)\subset\mathbb T^n$ acts on a quasitoric manifold
$M=M(P,\varLambda)$ semifreely and with isolated fixed points if
and only if, for every vertex $v=F_{j_1}\cap\cdots\cap F_{j_n}$,
the coefficients in~\eqref{nucoe} satisfy $k_i(\nu,v)=\pm1$ for
$1\le i\le n$.
\end{proposition}
\begin{proof}
It follows from Proposition~\ref{qtw} that the coefficients
$k_i(\nu,v)$ are the weights of the representation of the circle
$S(\nu)$ in the tangent space to $M$ at~$v$. The statement follows
from Proposition~\ref{sfrep}.
\end{proof}

\begin{theorem}[\cite{ma-pa08}]\label{qtbt}
Let $M$ be a quasitoric manifold over a cube~$I^n$. Assume that
the torus acting on $M$ has a circle subgroup acting semifreely
and with isolated fixed points. Then $M$ is equivalent to a Bott
tower.
\end{theorem}
\begin{proof}
Let $\varLambda_\star$ be the refined characteristic submatrix
of~$M$. We may assume by induction that every characteristic
submanifold of~$M$ is a Bott manifold, so that every proper
principal minor of the matrix $-\varLambda_\star$ is~$1$.
Therefore, we are in the situation of Lemma~\ref{minor}, and
$-\varLambda_\star$ is a matrix of one of the two types described
there. The second type is ruled out because of the semifreeness
assumption. Indeed, let $\varLambda=\left(I\;|\;-B\right)$, where
$B$ is the matrix~\eqref{exce} and assume that $S(\nu)\subset
\mathbb T^n$ acts semifreely with isolated fixed points. Applying
the criterion from Proposition~\ref{sfcoe} to the vertex
$v=F^0_1\cap\cdots\cap F^0_n$ we obtain $\nu_i=\pm 1$ for $1\le
i\le n$. Now we apply the same criterion to the vertex
$v'=F^1_1\cap\cdots\cap F^1_n$. Since the submatrix formed by the
corresponding columns of~$\varLambda$ is precisely~$-B$, it
follows that $\det B=\pm1$. This implies that $b_i=\pm1$
in~\eqref{exce} for some~$i$. Therefore, if all the coefficients
$k_j(\nu,v')$ in the expression
$\nu=k_1(\nu,v')\lambda_{n+1}+\cdots+k_n(\nu,v')\lambda_{2n}$ are
equal to~$\pm1$, then the $i$th coordinate of~$\nu$ is
$\nu_i=\pm1\pm b_i\ne\pm1$: a contradiction.
\end{proof}

Our next result shows that a Bott tower with a semifree circle
subgroup and isolated fixed points is topologically trivial.
%, that is, diffeomorphic to a product of 2-spheres.
Let $t$ (respectively,~$\C$) be the standard (respectively, the
trivial) complex one-dimensional representation of the
circle~$\mathbb S^1$, and let $\underline{V}$ denote the trivial
bundle with fibre~$V$ over a given base. We say that an action of
a group~$G$ on a Bott manifold~$B_n$ \emph{preserves the tower
structure} if for each stage $B_k=\C
P(\underline{\C}\oplus\xi_{k-1})$ the line bundle $\xi_{k-1}$ is
$G$-equivariant. The intrinsic $\mathbb T^n$-action on $B_n$ (see
Corollary~\ref{bottquasi}) obviously preserves the tower
structure.

\begin{theorem}[\cite{ma-pa08}]\label{btpro}
Assume that a Bott manifold $B_n$ admits a semifree $\mathbb
S^1$-action with isolated fixed points preserving the tower
structure. Then the Bott tower is topologically trivial;
furthermore, $B_n$ is $\mathbb S^1$-equivariantly diffeomorphic to
the product $(\C P(\C\oplus t))^n$.
\end{theorem}
\begin{proof}
We may assume by induction that the $(n-1)$th stage of the Bott
tower is diffeomorphic to $(\C P(\C\oplus t))^{n-1}$ and $B_n=\C
P(\underline{\C}\oplus \xi)$ for some $\mathbb S^1$-equivariant
line bundle $\xi$ over ${(\C P(\C\oplus t))^{n-1}}$.

Let $\gamma$ be the canonical line bundle over $\C P(\C\oplus
t)\cong\C P^1$. It carries a unique structure of an $\mathbb
S^1$-equivariant bundle such that
\begin{equation} \label{gamma}
  \gamma|_{(1:0)}=\C \qquad\text{and}\qquad \gamma|_{(0:1)}=t.
\end{equation}
We denote by $x\in{H^2(\C P(\C\oplus t))}$ the first Chern class
of~$\gamma$, and let $x_i=\pi_i^*(x)\in H^2(\C P(\C\oplus
t)^{n-1})$ be the pullback of~$x$ by the projection $\pi_i$ onto
the $i$th factor. Then $c_1(\xi)=\sum_{i=1}^{n-1}a_i x_i$ for some
$a_i\in \Z$. The $\mathbb S^1$-equivariant line bundles $\xi$ and
$\otimes_{i=1}^{n-1}\pi_i^*(\gamma^{a_i})$ have the same
underlying bundles, so there is an integer $k$ such that
\begin{equation} \label{xi}
  \xi= \underline t^k\otimes\bigotimes_{i=1}^{n-1}\pi_i^*(\gamma^{a_i})
\end{equation}
as $\mathbb S^1$-equivariant line bundles
(see~\cite[Corollary~4.2]{ha-yo76}).

We encode $\mathbb S^1$-fixed points in $\C P(\C\oplus t)^{n-1}$
by sequences
$(p_1^{\varepsilon_1},\dots,p_{n-1}^{\varepsilon_{n-1}})$, where
$\varepsilon_i=0$ or $1$, and $p_i^{\varepsilon_i}$ denotes
$(1:0)$ if $\varepsilon_i=0$ and $(0:1)$ if $\varepsilon_i=1$.
Then it follows from~\eqref{gamma} and~\eqref{xi} that
\[
  \xi|_{(p_1^{\varepsilon_1},\dots, p_{n-1}^{\varepsilon_{n-1}})}
  =t^{k+\sum_{i=1}^{n-1} \varepsilon_i a_i}.
\]
The $\mathbb S^1$-action on $B_n=\C P(\underline{\C}\oplus\xi)$ is
semifree if and only if $|k+\sum_{i=1}^{n-1} \varepsilon_i a_i|=1$
for all possible values of~$\varepsilon_i$. Setting
$\varepsilon_i=0$ for all~$i$ we obtain $|k|=1$. Let $k=1$ (the
case $k=-1$ is treated similarly). Then
$(a_1,\dots,a_{n-1})=(0,\dots,0)$ or $(0,\dots,0,-2,0,\dots,0)$.
In the former case, $\xi=\underline{t}$ and $B_n=\C
P(\underline{\C}\oplus\xi)\cong\C P(\C\oplus t)^n$. In the latter
case we have $\xi=t\pi_i^*(\gamma^{-2})$ for some~$i$, so that
$B_n=\pi_i^*\C P(\underline{\C}\oplus t\gamma^{-2})$. Since for
any $\mathbb S^1$-vector bundle $E$ and $\mathbb S^1$-line bundle
$\eta$ the projectivisations $\C P(E)$ and $\C P(E\otimes\eta)$
are $\mathbb S^1$-diffeomorphic, it follows that $\C
P(\underline{\C}\oplus t\gamma^{-2})\cong\C P(\gamma\oplus
t\gamma^{-1})$. The first Chern class of $\gamma\oplus
t\gamma^{-1}$ is zero, so its underlying bundle is trivial. The
$\mathbb S^1$-representation in the fibre of $\gamma\oplus
t\gamma^{-1}$ over a fixed point is isomorphic to $\C\oplus t$
by~\eqref{gamma}. Theorefore, $\gamma\oplus
t\gamma^{-1}=\underline{\mathbb C}\oplus\underline t$ as $\mathbb
S^1$-bundles. It follows that $\C P(\underline{\C}\oplus
t\gamma^{-2})\cong\C P(\underline{\mathbb C}\oplus\underline t)$
and $B_n\cong(\C P(\C\oplus t))^n$.
\end{proof}

\begin{remark}
The diffeomorphism of Theorem~\ref{btpro} is not $\mathbb
T^n$-equivariant.
\end{remark}

The next example shows that Theorem~\ref{btpro} cannot be
generalised to quasitoric manifolds. However, as we shall see, it
holds under the additional assumption that the quotient polytope
of the quasitoric manifold is a cube.

\begin{example}\label{nonexamp}
Let $M$ be a quasitoric manifold over a $2k$-gon with
\[
  \varLambda=\begin{pmatrix}
  1&0&1&0&\cdots&1&0\\
  0&1&0&1&\cdots&0&1
  \end{pmatrix}.
\]
By Corollary~\ref{sfcoe}, the circle subgroup determined by the
vector $\nu=(1,1)$ acts semifreely on~$M$. However, the quotient
of $M$ is not a 2-cube if $k>2$, so $M$ cannot be homeomorphic to
a product of spheres (it can be shown that $M$ is a connected sum
of $k-1$ copies of $S^2\times S^2$).
\end{example}

\begin{theorem}[\cite{ma-pa08}]\label{qtmsf}
Let $M$ a quasitoric manifold $M$ over a cube~$I^n$. Assume that
the torus acting on $M$ has a circle subgroup acting semifreely.
Then $M$ is $\mathbb S^1$-equivariantly homeomorphic to a product
of 2-spheres.
\end{theorem}
\begin{proof}
By Theorem~\ref{qtbt}, $M$ is equivalent to a Bott tower. By
Theorem~\ref{btpro}, it is $\mathbb S^1$-homeomorphic to a product
of spheres.
\end{proof}

We can also derive Ilinskii's result on semifree actions on toric
manifolds:

\begin{theorem}[\cite{ilin06}]\label{ilins}
A toric manifold $V$ carrying a semifree action of a circle
subgroup with isolated fixed points is diffeomorphic to a product
of 2-spheres.
\end{theorem}
\begin{proof}
By Theorem~\ref{btpro}, it is sufficient to show that $V$ is a
Bott manifold. A semifree circle subgroup acting on $V$ also acts
semifreely and with isolated fixed points on every characteristic
submanifold $V_j$ of~$V$. We use induction on the dimension. The
base of induction is the case $\dim_\C V=2$; we consider it below.
By the inductive hypothesis, each $V_j$ is a Bott manifold, i.e.
its quotient polytope is a combinatorial cube. On the other hand,
the quotient polytope of $V_j$ is the facet $F_j$ of the quotient
polytope $P$ of~$V$; since each $F_j$ is a cube, $P$ is also a
cube by Exercise~\ref{2facecube}. By Corollary~\ref{torbt}, $V$ is
a Bott manifold, and we are done.

It remains to consider the case $\dim_\C V=2$. We need to show
that the quotient polytope of a complex 2-dimensional toric
manifold $V$ with semifree circle subgroup action and isolated
fixed points is a 4-gon.

Let $\Sigma$ be the fan corresponding to~$V$. One-dimensional
cones of $\Sigma$ correspond to facets (or edges) of the quotient
polygon~$P^2$. We must show that there are precisely 4
one-dimensional cones. The values of the characteristic function
on the facets of $P^2$ are given by the primitive vectors
generating the corresponding one-dimensional cones of~$\Sigma$.
Let $\nu$ be the vector generating the semifree circle subgroup.
We may choose an initial vertex $v$ of $P^2$ so that $\nu$ belongs
to the 2-dimensional cone of $\Sigma$ corresponding to~$v$. Then
we index the primitive generators $\mb a_i$, $1\le i\le m$, of
1-cones so that $\nu$ is in the cone generated by $\mb a_1$ and
$\mb a_2$, and any two consecutive vectors span a two-dimensional
cone (see Fig.~\ref{figur}).
\begin{figure}[h]
\begin{center}
\begin{picture}(60,40)
  \put(40,20){\vector(1,0){20}}
  \put(40,20){\vector(1,1){20}}
  \put(40,20){\vector(0,1){20}}
  \put(40,20){\vector(-1,0){20}}
  \put(40,20){\vector(-2,-1){40}}
  \put(59.5,17){$\mb a_1$}
  \put(59.5,37){$\bnu$}
  \put(35.5,37){$\mb a_2$}
  \put(19,21.5){$\mb a_3$}
  \put(0,3.5){$\mb a_m$}
  \multiput(9,9)(1.5,1.5){3}{$\cdot$}
\end{picture}
\caption{ } \label{figur}
\end{center}
\end{figure}
This provides us with a refined characteristic matrix $\varLambda$
of size $2\times m$. We have $\mb a_1=(1,0)$ and $\mb a_2=(0,1)$,
and applying the criterion from Proposition~\ref{sfcoe} to the
first cone $\R_\ge\langle\mb a_1,\mb a_2\rangle$ (that is, to the
initial vertex of the polygon) we obtain $\nu=(1,1)$.

The rest of the proof is a case by case analysis using the
non-singularity condition~\eqref{detLv} and
Proposition~\ref{sfcoe}. The reader may be willing to do this
himself as an exercise rather than following the argument below.

Consider now the second cone. The non-singularity
conditions~\eqref{detLv} gives us $\det(\mb a_2,\mb a_3)=1$, hence
$\mb a_3=(-1,*)$. Writing $\bnu=k_1\mb a_2+k_2\mb a_3$ and
applying Proposition~\ref{sfcoe} to the second cone
$\R_\ge\langle\mb a_2,\mb a_3\rangle$ we obtain
\[
  (1,1)=\pm(0,1)\pm(-1,*).
\]
Therefore, $\mb a_3=(-1,0)$ or $\mb a_3=(-1,-2)$. Similarly,
considering the last cone $\R_\ge\langle\mb a_m,\mb a_1\rangle$ we
obtain $\mb a_m=(*,-1)$, and then, applying
Proposition~\ref{sfcoe}, we see that $\mb a_m=(0,-1)$ or $\mb
a_m=(-2,-1)$. The case when $\mb a_3=(-1,-2)$ and $\mb
a_m=(-2,-1)$ is impossible since then the second and the last
cones overlap.

Let $\mb a_3=(-1,-2)$. Then $\mb a_m=(0,-1)$. Considering the cone
$\R_\ge\langle\mb a_{m-1},\mb a_m\rangle$ we obtain $\mb
a_{m-1}=(-1,0)$ or $\mb a_{m-1}=(-1,-2)$. In the former case cones
overlap, and in the latter case we get $\mb a_{m-1}=\mb a_3$. This
implies $m=4$ and we are done.

Let $\mb a_3=(-1,0)$. Then considering the third cone
$\R_\ge\langle\mb a_3,\mb a_4\rangle$ we obtain $\mb a_4=(0,-1)$
or $\mb a_4=(-2,-1)$. If $\mb a_4=(0,-1)$, then $\mb a_4=\mb a_m$
(otherwise cones overlap), and we are done. Let $\mb a_4=(-2,-1)$.
Then either $\mb a_m=(-2,-1)=\mb a_4$ and we are done, or $\mb
a_m=(0,-1)$. In the latter case, we get $\mb a_{m-1}=(-1,-2)$ (see
the previous paragraph).

We are left with the case $\mb a_4=(-2,-1)$ and $\mb
a_{m-1}=(-1,-2)$. The only way to satisfy both~\eqref{detLv} and
the condition of Proposition~\ref{sfcoe} without overlapping cones
is to set $\mb a_5=(-3,-2)$ and $\mb a_{m-2}=(-2,-3)$. Continuing
this process, we obtain $\mb a_k=(-k+2,-k+3)$, $k\ge2$, and $\mb
a_{m-l}=(-l,-l-1)$, $l\ge0$. This process never stops, as we never
get a complete regular fan. So this case is impossible.
\end{proof}

The proof above leaves three possibilities for the vectors $\mb
a_3$ and $\mb a_4$ of the 2-dimensional fan: $(-1,0)$ and
$(0,-1)$, or $(-1,0)$ and $(-2,-1)$, or $(-1,-2)$ and $(0,-1)$.
The last two pairs correspond to isomorphic fans. The refined
characteristic submatrices corresponding to the first two pairs
are
\[
  \begin{pmatrix}-1&0\\0&-1\end{pmatrix}\quad\text{and}\quad
  \begin{pmatrix}-1&-2\\0&-1\end{pmatrix}.
\]
The first corresponds to $\C P^1\times\C P^1$, and the second to a
Bott manifold (Hirzebruch surface) with $a_{12}=-2$.

The next result gives an explicit description of
matrices~\eqref{amatr} corresponding to our specific class of Bott
towers:

\begin{theorem}[\cite{ma-pa08}]\label{sfdec}
A Bott manifold $B_n$ admits a semifree circle subgroup with
isolated fixed points if and only if its matrix~\eqref{amatr}
satisfies the identity
\[
  \frac12(E-A)=C_1C_2\cdots C_n,
\]
where $C_k$ is either the identity matrix or a unipotent upper
triangular matrix with only one nonzero element above the
diagonal; this element is $1$ in the $k$th column.
\end{theorem}
\begin{proof}
Assume first that $B_n$ admits a semifree circle subgroup with
isolated fixed points. We have two sets of multiplicative
generators for the ring $H^*(B_n)$: the set $\{u_1,\ldots,u_n\}$
from Corollary~\ref{cohbt} satisfying identities~\eqref{crobt},
and the set $\{x_1,\ldots,x_n\}$ satisfying $x_i^2=0$ (the latter
set exists since $B_n$ is diffeomorphic to a product of
2-spheres). The reduced sets with $i\le k$ can be regarded as the
corresponding sets of generators for the $k$th stage~$B_k$. As is
clear from the proof of Theorem~\ref{btpro}, we have
$c_1(\xi_{k-1})=-2c_{i_kk}x_{i_k}$ for some $i_k<k$, where
$c_{i_kk}=1$ or~$0$. From $u_k^2+2c_{i_kk}x_{i_k}u_k=0$ we obtain
$x_k=u_k+c_{i_kk}x_{i_k}$. In other words, the transition matrix
$C_k$ from the basis $x_1,\ldots,x_{k-1},u_k,\ldots,u_n$ of
$H^2(B_n)$ to $x_1,\ldots,x_k,u_{k+1},\ldots,u_n$ may have only
one nonzero entry off the diagonal, which is~$c_{i_kk}$. The
transition matrix from $u_1,\ldots,u_n$ to $x_1,\ldots,x_n$ is the
product $D=C_1C_2\cdots C_n$.
%(here $C_1$ is the identity matrix since $x_1=u_1$).
Then $D=(d_{jk})$ is a unipotent upper
triangular matrix consisting of zeros and ones,
$x_k=\sum_{j=1}^nd_{jk}u_j$ and
\[\textstyle
  0=x_k^2=\bigl(u_k+\sum_{j=1}^{k-1}d_{jk}u_j\bigr)^2=u_k^2+
  2\sum_{j=1}^{k-1}d_{jk}u_ju_k+\cdots, \qquad
  1\le k\le n.
\]
On the other hand, $0=u_k^2-\sum_{j=1}^{k-1}a_{jk}u_ju_k$
by~\eqref{crobt}. Comparing the coefficients of $u_ju_k$ for $1\le
j\le k-1$ in the last two equations and observing that these
elements are linearly independent in $H^4(B^{2k})$ we obtain
$2d_{jk}=-a_{jk}$ for $1\le j<k\le n$. As both $D$ and $-A$ are
unipotent upper triangular matrices, this implies $2D=E-A$.

Assume now that the matrix $A$ satisfies $E-A=2C_1C_2\cdots C_n$.
Then for the corresponding Bott tower we have
$\xi_{k-1}=\pi_{i_k}^*(\gamma^{-2c_{i_kk}})$. Therefore, we may
choose a circle subgroup such that $\xi_{k-1}$ becomes
$t\pi_{i_k}^*(\gamma^{-2c_{i_kk}})$ (as an $\mathbb
S^1$-equivariant bundle), for $1<k\le n$. This circle subgroup
acts semifreely and with isolated fixed points as seen from the
same argument as in the proof of Theorem~\ref{btpro}.
\end{proof}

\begin{example}
The condition of Theorem~\ref{sfdec} implies in particular that
the matrix~\eqref{amatr} may have only entries equal to $0$ or
$-2$ above the diagonal. However, the hypothesis of
Theorem~\ref{sfdec} is stronger. For instance, if
\[
  A=\begin{pmatrix} -1&0&-2\\ 0&-1&-2\\ 0&0&-1 \end{pmatrix}
\]
then the matrix $(E-A)/2$ cannot be factored as $C_1C_2C_3$.
Consequently, the corresponding 3-stage Bott tower does not admit
a subcircle acting semifreely and with isolated fixed points. On
the other hand, if
\[
  A=\begin{pmatrix} -1&-2&-2\\ 0&-1&0\\ 0&0&-1 \end{pmatrix}
\]
then we have
\[
  \frac12(E-A)=\begin{pmatrix} 1&1&1\\ 0&1&0\\ 0&0&1 \end{pmatrix}=
  \begin{pmatrix} 1&0&0\\ 0&1&0\\ 0&0&1 \end{pmatrix}
  \begin{pmatrix} 1&1&0\\ 0&1&0\\ 0&0&1 \end{pmatrix}
  \begin{pmatrix} 1&0&1\\ 0&1&0\\ 0&0&1 \end{pmatrix}.
\]
\end{example}

It is clear that not every topologically trivial Bott manifold
admits a semifree subcircle action with isolated fixed points (the
latter condition is stronger even for $n=2$). We now consider the
former class in more detail.

\subsection*{Topologically trivial Bott towers}\label{tocla}
The topological triviality of a Bott tower can be detected by its
cohomology ring:

\begin{theorem}[{\cite{ma-pa08}}]\label{trivc}
A Bott tower $B_n$ is topologically trivial if and only if there
is an isomorphism $H^*(B_n)\cong H^*((S^2)^n)$ of graded rings.
\end{theorem}
\begin{proof}
Corollary~\ref{cohbt} implies
\[
  H^*(B_n)=H^*\bigl(B_{n-1}\bigr)[u_n]\bigl/
  \bigl(u_n^2-c_1(\xi_{n-1})u_n\bigr).
\]
We may therefore write any element of $H^2(B_n)$ as $x+bu_n$,
where $x\in H^2(B_{n-1})$ and $b\in\Z$. We have
\[
  (x+bu_n)^2=x^2+2bxu_n+b^2u_n^2=x^2+b(2x+bc_1(\xi_{n-1}))u_n,
\]
so that the square of $x+bu_n$ with $b\not=0$ is zero if and only
if $x^2=0$ and $2x+bc_1(\xi_{n-1})=0$. This shows that elements of
the form $x+bu_n$ with $b\ne0$ whose squares are zero generate a
rank-one free subgroup of $H^2(B_n)$.

Assume that $H^*(B_n)\cong H^*((\C P^1)^n)$. Then there is a basis
$\{x_1,\dots,x_n\}$ in $H^2(B_n)$ such that $x_i^2=0$ for all~$i$.
By the observation from the previous paragraph, we may assume that
the elements $x_1,\dots,x_{n-1}$ lie in $H^2(B_{n-1})$, and $x_n$
is not in $H^2(B_{n-1})$. Then we can have
$x_n=\sum_{i=1}^{n-1}b_ix_i+u_n$ for some $b_i\in\Z$. A product
$\prod_{i\in I}x_i$ with $I\subset\{1,\dots,n\}$ lies in
$H^*(B_{n-1})$ if and only if $n\notin I$. This implies that the
ring $H^*(B_{n-1})$ is generated by the elements
$x_1,\dots,x_{n-1}$ and is isomorphic to the cohomology ring of
$(\C P^1)^{n-1}$. Therefore, we may assume by induction that
$B_{n-1}\cong (\C P^1)^{n-1}$.

Writing $c_1(\xi_{n-1})=\sum_{i=1}^{n-1}a_ix_i$, we obtain
\[\textstyle
  0=x_n^2=\bigl(u_n+\sum_{i=1}^{n-1} b_ix_i\bigr)^2
  =\sum_{i=1}^{n-1}(a_i+2b_i)x_iu_n+
  \bigl(\sum_{i=1}^{n-1} b_ix_i\bigr)^2.
\]
This may hold only if at most one of the $a_i$ is nonzero (and
equal to~$-2b_i$) because the elements $x_ix_j$ and $x_iu_n$ with
$i<j<n$ form a basis of~$H^4(B_n)$. Therefore, $\xi_{n-1}$ is the
pullback of the bundle $\gamma^{-2b_i}$ over $\C P^1$ by the $i$th
projection map $B_{n-1}=(\C P^1)^{n-1}\to \C P^1$. Since $\C
P(\underline{\C}\oplus\gamma^{-2b_i})$ is a topologically trivial
bundle (see Example~\ref{2dbto}), the bundle $B_n=\C
P(\underline{\C}\oplus\xi_{n-1})$ is also trivial.
\end{proof}

We can now also describe effectively the class of
matrices~\eqref{amatr} corresponding to topologically trivial Bott
towers:

\begin{theorem}
A Bott tower is topologically trivial if and only if its
corresponding matrix~\eqref{amatr} satisfies the identity
\[
  \frac12(E-A)=C_1C_2\cdots C_n,
\]
where $C_k$ is either the identity matrix or a unipotent upper
triangular matrix with only one nonzero element above the
diagonal; this element lies in the $k$th column.
\end{theorem}
\begin{proof}
The argument is the same as in the proof of Theorem~\ref{sfdec}.
The only difference is that the number $c_{i_kk}$ in the formula
$c_1(\xi_{k-1})=-2c_{i_kk}x_{i_k}$ is now an arbitrary integer.
\end{proof}

Theorem~\ref{trivc} can be generalised to quasitoric manifolds,
but only in the topological category:

\begin{theorem}[{\cite[Theorem~5.7]{ma-pa08}}]\label{qtcpn}
A quasitoric manifold $M$ is homeomorphic to a product $(S^2)^n$
if and only if there is an isomorphism $H^*(M)\cong H^*((S^2)^n)$
of graded rings.
\end{theorem}

\subsection*{Generalisations and cohomological rigidity}
\begin{definition}\label{genbott}
A \emph{generalised Bott tower} of height $n$ is a tower of
bundles
\[
  B_n\stackrel{p_n}\longrightarrow B_{n-1}
  \stackrel{p_{n-1}}\longrightarrow\cdots
  \longrightarrow B_2\stackrel{p_2}\longrightarrow B_1\longrightarrow\pt,
\]
of complex manifolds, where $B_1=\C P^{j_1}$ and each $B_k$ is the
complex projectivisation of a sum of $j_k+1$ complex line bundles
over $B_{k-1}$. The fibre of the bundle $p_k\colon B_k\to B_{k-1}$
is~$\C P^{j_k}$. A generalised Bott tower is \emph{topologically
trivial} if each $p_k$ is trivial as a smooth bundle.

The last stage $B_n$ in a generalised Bott tower is a
\emph{generalised Bott manifold}.
\end{definition}

\begin{remark}
A generalised Bott manifold $B_n$ is a projective toric manifold
whose corresponding polytope is combinatorially equivalent to a
product of simplices
$\varDelta^{j_1}\times\cdots\times\varDelta^{j_n}$ (an exercise).
In particular, $B_n$ is a quasitoric manifold over a product of
simplices. If we replace `a sum of complex line bundles' by `a
complex vector bundle' in the definition, then the resulting tower
will not be a toric manifold in general: the torus action on
$B_{k-1}$ lifts to the projectivisation of a sum of line bundles,
but not to the projectivisation of an arbitrary vector bundle
over~$B_{k-1}$.
\end{remark}

Generalised Bott towers were considered by
Dobrinskaya~\cite{dobr01}, who proved the following result (see
Lemma~\ref{qts} for the information about the signs of vertices):

\begin{theorem}[{\cite[Corollary~7]{dobr01}}]
A quasitoric manifold over a product of simplices
$P=\varDelta^{j_1}\times\cdots\times\varDelta^{j_n}$ is a
generalised Bott manifold if and only if the sign of each vertex
of~$P$ is the product of the signs of its corresponding vertices
of the simplices $\varDelta^{j_k}$ according to the decomposition
of~$P$. In particular, a toric manifold over a product of
simplices is always a generalised Bott manifold.
\end{theorem}

\begin{remark}
According to the general result~\cite[Theorem~6]{dobr01}, a
quasitoric manifold over any product polytope
$P=P_{j_1}\times\cdots\times P_{j_n}$ decomposes into a tower of
quasitoric fibre bundles if and only if the sign of each vertex
of~$P$ decomposes into the corresponding product.
\end{remark}

Theorem~\ref{trivc} can be extended to generalised Bott towers:

\begin{theorem}[{\cite[Theorem~1.1]{c-m-s10T}}]\label{gentrivc}
If the integral cohomology ring of a generalised Bott manifold
$B_n$ is isomorphic to $H^*(\C P^{j_1}\times\cdots\times\C
P^{j_n})$, then the generalised Bott tower is topologically
trivial; in particular, $B_n$ is diffeomorphic to $\C
P^{j_1}\times\cdots\times\C P^{j_n}$.
\end{theorem}

According to another result of Choi, Masuda and
Suh~\cite[Theorem~8.1]{c-m-s10Q}, a quasitoric manifold $M$ over a
product of simplices
$P=\varDelta^{j_1}\times\cdots\times\varDelta^{j_n}$ is
homeomorphic to a generalised Bott manifold if $H^*(M)\cong H^*(\C
P^{j_1}\times\cdots\times\C P^{j_n})$. Therefore, such a
quasitoric manifold $M$ is homeomorphic to $\C
P^{j_1}\times\cdots\times\C P^{j_n}$.

Since trivial (generalised) Bott towers can be detected by their
cohomology rings, a natural question arises of whether \emph{any}
two (generalised) Bott manifolds can be distinguished, either
smoothly or topologically, by their cohomology rings. This leads
to the notion of cohomological rigidity:

\begin{definition}\label{defcohorig}
Fix a commutative ring~$\k$ with unit. We say that a family of
closed manifolds is \emph{cohomologically rigid} over $\mathbf k$
if manifolds in the family are distinguished up to homeomorphism
by their cohomology rings with coefficients in~$\mathbf k$. That
is, a family is cohomologically rigid if a graded ring isomorphism
$H^*(M_1;\mathbf k)\cong H^*(M_2;\mathbf k)$ implies a
homeomorphism $M_1\cong M_2$ whenever $M_1$ and $M_2$ are in the
family.

A manifold $M$ in the given family is said to be
\emph{cohomologically rigid} if for any other manifold $M'$ in the
family a ring isomorphism $H^*(M;\mathbf k)\cong H^*(M';\mathbf
k)$ implies a homeomorphism $M\cong M'$. Obviously a family is
cohomologically rigid whenever every its element is rigid.

There is a smooth version of cohomological rigidity for families
of smooth manifolds, with homeomorphisms replaced by
diffeomorphisms.
\end{definition}

\begin{problem}\label{crbott}
Is the family of Bott manifolds cohomologically rigid (over~$\Z$)?
Namely, is it true that any two Bott manifolds $M_1$ and $M_2$
with isomorphic integral cohomology rings are homeomorphic (or
even diffeomorphic)?
\end{problem}

This question is open even for the much larger families of toric
and quasitoric manifolds. Namely, there are no known examples of
non-homeomorphic (quasi)toric manifolds with isomorphic integral
cohomology rings. By the result of Choi, Park and
Suh~\cite{c-p-s12}, quasitoric manifolds with second Betti
number~2 (i.e. over a product of two simplices) are homeomorphic
when their cohomology rings are isomorphic.

In the positive direction, Theorem~\ref{trivc} shows that a
topologically trivial Bott manifold is cohomologically rigid in
the family of Bott manifolds, in the smooth category. By
Theorem~\ref{qtcpn}, a topologically trivial Bott manifold is
cohomologically rigid in the wider family of quasitoric manifolds,
but only in the topological category. Smooth cohomological
rigidity was established in Choi~\cite{choi11} for Bott manifolds
up to dimension~8 (i.e. up to height~4).

There is an $\mathbb R$-version of this circle of questions, with
quasitoric manifolds replaced by \emph{small
covers}\label{bottsmalcov}, (generalised) Bott towers replaced by
\emph{real (generalised) Bott towers}, and the cohomology rings
taken with coefficients in~$\mathbb Z_2$. A real generalised Bott
tower is defined similarly to a complex tower, with the complex
projectivisation replaced by the real one.

The family of real Bott manifolds is cohomologically rigid over
$\mathbb Z_2$ by the result of Kamishima--Masuda~\cite{ka-ma09}
(see also~\cite{c-m-o}), but the family of generalised Bott
manifolds is not~\cite{masu10}. Also, every small cover over a
product of simplices is a generalised real Bott
tower~\cite{c-m-s10Q} (this is not true for quasitoric manifolds,
see Example~\ref{qtovercube}). For more results on the
cohomological rigidity of (generalised) Bott towers, see the
survey articles~\cite{ma-su08} and~\cite{c-m-s11}.

\subsection*{Exercises}
\begin{exercise}
Let $\eta^k$ denote the $k$th tensor power of the tautological
line bundle over~$\C P^1$, and $\underline{\C}$ denote the trivial
line bundle. Show that there is a cohomology ring isomorphism
$H^*(\C P(\underline{\C}\oplus\eta^k)) \cong H^*(\C
P(\underline{\C}\oplus\eta^{k'}))$ if and only if $k=k'\mod2$.
\end{exercise}

\begin{exercise}
Show that the toric manifold described in Theorem~\ref{almat} is
projective by providing explicitly a polytope $P$ whose normal fan
is~$\Sigma$. Observe that $P$ is combinatorially equivalent to a
cube (compare Proposition~\ref{bproj}).
\end{exercise}

\begin{exercise}\label{bfm2toric}
The bounded flag manifold $\BF_n$ is a Bott tower, as described in
Example~\ref{btbfn}. By Theorem~\ref{almat}, it is isomorphic to
the toric manifold whose corresponding fan has generators
\[
  \mb a_k^0=\mb e_k\quad (1\le k\le n),\quad
  \mb a_k^1=-\mb e_k+\mb e_{k+1}\quad (1\le k\le
  n-1),\quad\mb a_n^1=\mb e_n.
\]
This fan is not the same as the one described in
Proposition~\ref{bfan}. However, there is a linear isomorphism
of~$\R^n$ taking one fan to another, and therefore the
corresponding toric manifolds are isomorphic.
\end{exercise}

\begin{exercise}\label{reluivi}
Let $B_n$ be the Bott manifold and $u_i\in H^2(B_n)$ the canonical
cohomology ring generator (obtained by pulling back the first
Chern class of the tautological line bundle over~$B_i$ to the top
stage~$B_n$). Show that $u_i=c_1(\bar\rho_{i+n})$, where
$\rho_{i+n}$ is the line bundle~\eqref{rhoi} over the toric
manifold~$B_n$. (Hint: use induction.)
\end{exercise}

\begin{exercise}
A generalised Bott manifold is a projective toric manifold over a
product of simplices.
\end{exercise}

\begin{exercise}
The projectivisation of a complex $k$-plane bundle ($k>1$) over
$\C P^n$ is not necessarily a toric manifold.
\end{exercise}

\section{Weight graphs}\label{weightgraphs}
The quotient space $M^{2n}/T^n$ of a half-dimensional torus action
has the orbit stratification (the face structure), therefore
providing a natural combinatorial object associated with the
action and allowing us to translate equivariant topology into
combinatorics. For the classes of half-dimensional torus actions
considered in the previous sections, the quotient $M^{2n}/T^n$ is
contractible or acyclic, so many topological invariants of the
action can be expressed in purely combinatorial terms.

In this section we consider another type of combinatorial objects
associated to $T$-manifolds, which catch some important
information about the $T$-action and its orbit structure: the
so-called \emph{weight graphs} (for particular classes of
$T$-manifolds these are also known as \emph{GKM-graphs}). A weight
graph is an oriented graph with a special labelling of edges,
which can be assigned to an effective $T^k$-action on $M^{2n}$
($k\le n$) with isolated fixed points under very mild assumption
on the action. When $k<n$ the topology of the quotient
$M^{2n}/T^k$ is often quite complicated, and the weight graph can
be viewed as its combinatorial approximation.

In the study of quasitoric manifolds we constructed the
correspondence~\eqref{qtaf}, which assigns to each oriented edge
of the quotient simple polytope $P$ a weight of the
$T$-representation at the fixed point corresponding to the origin
of the edge. This correspondence can be viewed as a graph $\Gamma$
with special labels on its oriented edges, and we refer to such an
object as a \emph{weight graph}. As it follows from
Proposition~\ref{qtw}, defining the correspondence~\eqref{qtaf} is
equivalent to defining the characteristic matrix~$\varLambda$. At
the same time, it is well-known that the 1-skeleton of a simple
polytope determines its entire combinatorial structure (see
e.g.~\cite[Theorem~3.12]{zieg95}). It follows that the weight
graph contains the same information as the combinatorial
quasitoric pair $(P,\varLambda)$, and therefore it completely
determines the torus action on the quasitoric manifold~$M$. For
more general classes of torus actions considered in this chapter,
one cannot expect that the weight graph determines the action, but
it still contains an important piece of information.

Graphs whose oriented edges are labelled by the weights of a torus
action were considered in the works of Musin~\cite{musi80},
Hattori and other authors since the 1970s. A renewed interest to
these graphs was stimulated by the works of
Goresky--Kottwitz--MacPherson~\cite{g-k-m98} and
Guillemin--Zara~\cite{gu-za99} in connection with the study of
symplectic manifolds with Hamiltonian torus actions. A related
more general class of $T$-manifolds has become known as
\emph{GKM-manifolds}\label{dGKMmani}, and their weight graphs are
often referred to as \emph{GKM-graphs}.

A closed $2n$-dimensional manifold $M$ with an effective smooth
action of torus $T^k$ ($k\le n$) is called a \emph{GKM-manifold}
if the fixed point set is finite and nonempty, a $T^k$-invariant
almost complex structure is given on~$M$, and the weights of the
tangential $T^k$-representation at any fixed point are pairwise
linearly independent. As in the case of quasitoric manifolds,
there is a weight graph associated with each GKM-manifold~$M$
(vertices of the graph correspond to fixed points of~$M$, and
edges correspond to connected components of the set of points of
$M$ with codimension-one stabilisers). As it was shown
in~\cite{g-k-m98}, many important topological characteristics of a
GKM-manifold~$M$ (such as the Betti numbers or the equivariant
cohomology) can be described in terms of the weight graph.
Axiomatisation of the properties of the weight graph of a
GKM-manifold led to the notion of a GKM-graph~\cite{gu-za99}.

Weight graphs arising from locally standard $T$-manifolds (or
torus manifolds) were studied in~\cite{m-m-p07}. As in the case of
GKM-manifolds, axiomatisation of the properties of weight graphs
leads to an interesting combinatorial object, known as a
\emph{$T$-graph} (or \emph{torus graph}).

A \emph{$T$-graph} is a finite $n$-valent graph $\Gamma$ (without
loops, but with multiple edges allowed) with an \emph{axial
function} on the set $E(\Gamma)$ of oriented edges taking values
in $\Hom(T^n,\mathbb S^1)=H^2(BT^n)$ and satisfying certain
compatibility conditions. These conditions (described below) are
similar to those for GKM-graphs, but not exactly the same. The
weight graph of a torus manifold is an example of a $T$-graph; in
this case the values of the axial function are the weights of the
tangential representations of $T^n$ at fixed points.

The equivariant cohomology ring $H^*_T(\Gamma)$ of a $T$-graph
$\Gamma$ can be defined in the same way as for GKM-graphs; when
the $T$-graph arises from a locally standard torus manifold $M$ we
have $H^*_T(\Gamma)=H^*_T(M)$. Furthermore, unlike the case of
GKM-graphs, the equivariant cohomology ring of a $T$-graph can be
described in terms of generators and relations (see
Theorem~\ref{tequi}). Such a description is obtained by defining
the simplicial poset $\mathcal S(\Gamma)$ associated with a
$T$-graph~$\Gamma$; then $H^*_T(\Gamma)$ is shown to be isomorphic
to the face ring $\Z[\mathcal S(\Gamma)]$. This theorem continues
the series of results identifying the equivariant cohomology of a
(quasi)toric manifold~\cite{da-ja91} and a locally standard
$T$-manifold (Theorem~\ref{th:eqcoh}) with the face ring of the
associated polytope, simplicial complex, or simplicial poset.

Although the classes of GKM- and $T$-graphs diverge in general,
they contain an important subclass of $n$-independent GKM-graphs
in their intersection.

\subsection*{Definition of a $T$-graph}\label{stogr}
This definition is a natural adaptation of the notion of
GKM-graph~\cite{gu-za99} to torus manifolds.

Let $\Gamma$ be a connected $n$-valent graph without loops but
possibly with multiple edges. Denote by $V(\Gamma)$ the set of
vertices and by $E(\Gamma)$ the set of oriented edges (so that
each edge enters $E(\Gamma)$ twice with the opposite
orientations). We further denote by $i(e)$ and $t(e)$ the initial
and terminal points of $e\in E(\Gamma)$, respectively, and denote
by $\bar e$ the edge $e$ with the reversed orientation. For $v\in
V(\Gamma)$ we set
\[
  E(\Gamma)_v=\{ e\in E(\Gamma)\colon i(e)=v\}.
\]
A collection $\theta=\{\theta_e\}$ of bijections
\[
  \theta_e \colon E(\Gamma)_{i(e)}\to E(\Gamma)_{t(e)},\qquad e\in E(\Gamma),
\]
is called a \emph{connection}\label{defconn} on $\Gamma$ if
\begin{itemize}
\item[(a)] $\theta_{\bar e}$ is the inverse of
$\theta_e$;

\item[(b)] $\theta_e(e)=\bar e$.
\end{itemize}

An $n$-valent graph $\Gamma$ with $g$ edges admits $((n-1)!)^{g}$
different connections. Let $T=T^n$ be an $n$-torus. A map
\[
  \alpha\colon E(\Gamma)\to \Hom(T,S^1)=H^2(BT)
\]
is called an \emph{axial function}\label{defiaxif} (associated
with the connection $\theta$) if it satisfies the following three
conditions:
\begin{itemize}
\item[(a)] $\alpha(\bar e)=\pm \alpha(e)$;

\item[(b)] elements of $\alpha(E(\Gamma)_v)$ are pairwise linearly independent
(2-\emph{independent}) for each vertex $v\in V(\Gamma)$;

\item[(c)] $\alpha(\theta_e(e'))\equiv \alpha(e') \mod \alpha(e)$
for any $e\in E(\Gamma)$ and $e'\in E(\Gamma)_{i(e)}$.
\end{itemize}

We also denote by $T_e=\ker\alpha(e)$ the codimension-one subtorus
in $T$ determined by $\alpha$ and~$e$. Then we may reformulate the
condition~(c) above as follows: the restrictions of
$\alpha(\theta_e(e'))$ and $\alpha(e')$ to $H^*(BT_e)$ coincide.

\begin{remark} Guillemin and Zara required $\alpha(\bar
e)=-\alpha(e)$ in their definition of axial function. A connection
$\theta$ satisfying condition~(c) above is unique if elements of
$\alpha(E(\Gamma)_v)$ are 3-independent for each vertex~$v$ (an
exercise, see~\cite{gu-za99}).
\end{remark}

\begin{definition}\label{defiTgraph}
We call $\alpha$ a \emph{$T$-axial function} if it is
$n$-independent, i.e. if $\alpha(E(\Gamma)_v)$ is a basis of
$H^2(BT)$ for each $v\in V(\Gamma)$. A triple
$(\Gamma,\theta,\alpha)$ consisting of a graph $\Gamma$, a
connection $\theta$ and a $T$-axial function~$\alpha$ is called a
\emph{$T$-graph}. Since a connection $\theta$ is uniquely
determined by~$\alpha$, we often suppress it in the notation.
\end{definition}

\begin{remark}
Compared with GKM-graphs, the definition of a $T$-graph has weaker
condition~(a) (we only require $\alpha(\bar e)=\pm\alpha(e)$
instead of $\alpha(\bar e)=-\alpha(e)$), but stronger
condition~(b) ($\alpha$ is required to be $n$-independent rather
than 2-independent).
\end{remark}

\begin{example}\label{etmtg}
Let $M$ be a locally standard torus manifold. Denote by $\Gamma_M$
the 1-skeleton of the orbit space $Q=M/T$ (it is easy to see that
the 1-skeleton can be defined without the local standardness
assumption), and let $\alpha_M$ be the axial function of
Lemma~\ref{torusaf}. Then $(\Gamma_M,\alpha_M)$ is a $T$-graph.
\end{example}

\begin{example}\label{ebana}
Two $T$-graphs are shown in Fig.~\ref{ftogr}. The first is
2-valent and the second is 3-valent. The axial function $\alpha$
takes the edges, regardless of their orientation, to the
generators $t_1,t_2\in H^2(BT^2)$ (respectively, $t_1,t_2,t_3\in
H^2(BT^3)$). These $T$-graphs are not GKM-graphs, as the condition
$\alpha(\bar e)=-\alpha(e)$ is not satisfied. Both come from torus
manifolds, $S^4$ and $S^6$, respectively (see Example~\ref{S2n
orbit}).
\begin{figure}[h]
\begin{picture}(120,30)
\put(30,0){\circle*{1}} \put(30,30){\circle*{1}}
\put(90,0){\circle*{1}} \put(90,30){\circle*{1}}
\qbezier(30,0)(15,15)(30,30) \qbezier(30,0)(45,15)(30,30)
\qbezier(90,0)(70,15)(90,30) \qbezier(90,0)(110,15)(90,30)
\qbezier(90,0)(95,15)(90,30) \put(19,14){$t_1$}
\put(38.5,14){$t_2$} \put(76.5,14){$t_1$} \put(89,14){$t_2$}
\put(101,14){$t_3$} \put(21,-5){(a) \ $n=2$} \put(81,-5){(b) \
$n=3$}
\end{picture}
\medskip
\caption{$T$-graphs.} \label{ftogr}
\end{figure}
\end{example}

\begin{definition}\label{eqcohTgr}
The \emph{equivariant cohomology} $H^*_T(\Gamma)$ of a $T$-graph
$\Gamma$ is a set of maps
\[
  f\colon V(\Gamma)\to H^*(BT)
\]
such that for every $e\in E(\Gamma)$ the restrictions of $f(i(e))$
and $f(t(e))$ to $H^*(BT_e)$ coincide. Since $H^*(BT)$ is a ring,
the set of maps from $V(\Gamma)$ to $H^*(BT)$, denoted by
$H^*(BT)^{V(\Gamma)}$, is also a ring with respect to the
vertex-wise multiplication. Its subspace $H_T^*(\Gamma)$ is a
subring because the restriction map $H^*(BT)\to H^*(BT_e)$ is
multiplicative. Furthermore, $H_T^*(\Gamma)$ is a
$H^*(BT)$-algebra.
\end{definition}

If $M$ is a torus manifold with $H^{odd}(M)=0$, then for the
corresponding $T$-graph $\Gamma_M$ we have $H^*_T(\Gamma_M)\cong
H^*_T(M)$ by Theorem~\ref{theo:GKM}.

\subsection*{Calculation of equivariant cohomology}
Here we interpret the results of Section~\ref{torusman} on
equivariant cohomology of torus manifolds in terms of their
associated $T$-graphs, thereby providing a purely combinatorial
model for this calculation, which is applicable to a wider class
of objects.

\begin{definition}\label{faceTgr}
Let $(\Gamma,\theta,\alpha)$ be a $T$-graph and $\Gamma'$ a
connected $k$-valent subgraph of~$\Gamma$, where $0\le k\le n$. If
$\Gamma'$ is invariant under the connection $\theta$, then we say
that $(\Gamma',\alpha|E(\Gamma'))$ is a \emph{$k$-dimensional
face} of $\Gamma$. As usual, $(n-1)$-dimensional faces as called
\emph{facets}.
\end{definition}

An intersection of faces is invariant under the connection, but
can be disconnected. In other words, such an intersection is a
union of faces.

The \emph{Thom class}\label{ThomTgr} of a $k$-dimensional face
$G=(\Gamma',\alpha|E(G'))$ is the map $\ta{G}\colon V(\Gamma)\to
H^{2(n-k)}(BT)$ defined by
\begin{equation}\label{etaf1}
  \ta{G}(v)=\begin{cases}\displaystyle{\prod_{i(e)=v,\
     e\notin\Gamma'}
  \alpha(e)}\quad&\text{if $v\in V(\Gamma')$,}\\
  \qquad 0 \quad&\text{otherwise.}
  \end{cases}
\end{equation}

\begin{lemma}
The Thom class $\ta{G}$ is an element of $H^*_T(\Gamma)$.
\end{lemma}
\begin{proof}
Let $e\in E(\Gamma)$. If neither of the vertices of $e$ is
contained in~$G$, then the values of $\ta{G}$ on both vertices of
$e$ are zero. If only one vertex of $e$, say $i(e)$, is contained
in~$G$, then $\ta{G}(t(e))=0$, while
$\ta{G}(i(e))=0\mod\alpha(e)$, so that the restriction of
$\ta{G}(i(e))$ to $H^*(BT_e)$ is also zero. Finally, assume that
the whole $e$ is contained in~$G$. Let $e'$ be an edge such that
$i(e')=i(e)$ and $e'\notin G$, so that $\alpha(e')$ is a factor in
$\ta{G}(i(e))$. Since $G$ is invariant under the connection, it
follows that $\theta_{e}(e')\notin G$. Therefore,
$\alpha(\theta_e(e'))$ is a factor in $\ta{G}(t(e))$. Now we have
$\alpha(\theta_e(e'))\equiv\alpha(e')\mod\alpha(e)$ by the
definition of axial function. The same holds for any other factor
in $\ta{G}(i(e))$, whence the restrictions of $\ta{G}(i(e))$ and
$\ta{G}(t(e))$ to $H^*(BT_e)$ coincide.
\end{proof}

\begin{lemma} \label{lemm:face}
If $\Gamma$ is a $T$-graph, then there is a unique $k$-face
containing any given $k$ elements of $E(\Gamma)_v$.
\end{lemma}
\begin{proof}
Let $S\subset E(\Gamma)_v$ be the set of given oriented $k$ edges
with the common origin~$v$. Consider the graph $\Gamma'$ obtained
by `spreading' $S$ using the connection~$\theta$. In more detail,
at the first step we add to $S$ all oriented edges of the form
$\theta_e(e')$ where $e,e'\in S$. Denote the resulting set by
$S_1$. At the second step we add to $S_1$ all oriented edges of
the form $\theta_e(e')$ where $e,e'\in S_1$ and $i(e)=i(e')$, an
so on. Since $\Gamma$ is a finite graph, this process stabilises
after a finite number of steps, and we obtain a subgraph~$\Gamma'$
of~$\Gamma$. This subgraph $\Gamma'$ is obviously
$\theta$-invariant. We claim that $\Gamma'$ is $k$-valent. To see
this, define for any vertex $w\in V(\Gamma')$ the subgroup
\[
  N_w=\Z\langle\alpha(e)\colon e\in E(\Gamma')_w\rangle\subset H^2(BT).
\]
Condition~(c) from the definition of the axial function implies
that $N_w=N_{w'}$ for any vertices $w,w'\in V(\Gamma')$. Since
$N_v\cong\Z^k$ for the initial vertex~$v$, it follows that
$N_w\cong\Z^k$ for any $w\in V(\Gamma')$. Now, then
$n$-independence of the axial function implies that there are
exactly $k$ edges in the set $\{e\in E(\Gamma')_w\}$, for any
vertex $w$ of~$\Gamma'$. In other words, $\Gamma'$ is $k$-valent
and therefore it defines a $k$-face of $\Gamma$.
\end{proof}

\begin{corollary}
Faces of a $T$-graph $\Gamma$ form a simplicial poset $\mathcal
S(\Gamma)$ of rank $n$ with respect to reversed inclusion.
\end{corollary}

Denote by $G\vee H$ a minimal face containing both $G$ and $H$. In
general such a least upper bound may fail to exist or be
non-unique; however it exists and is unique provided that the
intersection $G\cap H$ is non-empty.

\begin{lemma}\label{ltcre}
For any two faces $G$ and $H$ of $\Gamma$ the corresponding Thom
classes satisfy the relation
\begin{equation}\label{etocl}
  \ta{G}\ta{H}=
    \ta{G\vee H}\cdot\!\!\!\sum_{E\in{G\cap H}}\!\!\!\ta{E},
\end{equation}
where we formally set $\ta{\Gamma}=1$ and $\tau_{\varnothing}=0$,
and the sum in the right hand side is taken over connected
components $E$ of $G\cap H$.
\end{lemma}
\begin{proof} We need to check that both sides of the identity take the same values on any vertex $v$.
The argument is the same as in the proof of Lemma~\ref{taurel}.
\end{proof}

\begin{lemma}\label{ltcmg}
The Thom classes $\ta{G}$ corresponding to all proper faces of
$\Gamma$ constitute a set of ring generators for $H^*_T(\Gamma)$.
\end{lemma}
\begin{proof}
This is proved in the same way as Lemma~\ref{modgn}.
\end{proof}

Consider the face ring $\Z[\mathcal S(\Gamma)]$, obtained by
taking quotient of the polynomial ring on generators $\va{G}$
corresponding to non-empty faces of~$\Gamma$ by
relations~\eqref{etocl}. The grading is given by
$\deg\va{G}=2(n-\dim G)$.

\begin{example}
Let $\Gamma$ be the torus graph shown in Fig.~\ref{ftogr}~(b).
Denote its vertices by $p$ and $q$, the edges by $e$, $g$, $h$,
and their opposite 2-faces by $E$, $G$, $H$, respectively. The
simplicial cell complex $\mathcal S(\Gamma)$ is obtained by gluing
two triangles along their boundaries. The face ring $\Z[\mathcal
S(\Gamma)]$ is the quotient of the graded polynomial ring
\[
  \Z[\va{E},\va{G},\va{H},v_p,v_q],\quad\deg\va{E}=\deg\va{G}=\deg\va{H}=2,\quad
  \deg v_p=\deg v_q=6
\]
by the two relations
\[
  \va{E}\va{G}\va{H}=v_p+v_q,\quad v_pv_q=0.
\]
(The generators $v_e,v_g,v_h$ can be excluded using relations like
$v_e=\va{G}\va{H}$.)
\end{example}

By definition, the equivariant cohomology of a $T$-graph comes
together with a monomorphism into the sum of polynomial rings:
\[
  r\colon H^*_T(\Gamma)\longrightarrow\bigoplus_{V(\Gamma)}H^*(BT).
\]
A similar map for the face ring $\Z[\mathcal S(\Gamma)]$ is given
by Theorem~\ref{alres} (or by Theorem~\ref{qalres}). The latter
map can be written in our case as
\[
  s\colon\Z[\mathcal S(\Gamma)]\longrightarrow
  \bigoplus_{v\in V(\Gamma)}
  \Z[\mathcal S(\Gamma)]/(\va{G}\colon G\not\ni v).
\]

\begin{theorem}[\cite{m-m-p07}]\label{tequi}
The equivariant cohomology ring $H^*_T(\Gamma)$ of a $T$-graph
$\Gamma$ is isomorphic to the face ring $\Z[\mathcal S(\Gamma)]$.
In other words, $H^*_T(\Gamma)$ is isomorphic to the quotient of
the polynomial ring on the Thom classes $\ta{G}$ by
relations~\eqref{etocl}.
\end{theorem}
\begin{proof}
We define a map
\[
  \Z[\va{G}\colon G\text{ is a face}]\longrightarrow
  H^*_T(\Gamma)
\]
by sending $\va{G}$ в $\ta{G}$. By Lemma~\ref{ltcre}, it factors
through a map $\varphi\colon\Z[\mathcal S(\Gamma)]\to
H^*_T(\Gamma)$. This map is surjective by Lemma~\ref{ltcmg}. Also,
$\varphi$ is injective, because we have $s=r\circ\varphi$ and $s$
is injective by Theorem~\ref{alres}.
\end{proof}

\subsection*{Pseudomanifolds and orientations}\label{spseu}
Which simplicial posets arise as the posets of faces of
$T$-graphs? Here we obtain a partial answer to this question.

The definition of pseudomanifold (Definition~\ref{defpm}) can be
extended easily to simplicial posets:

\begin{definition}\label{pseudosimpos}
A simplicial poset $\mathcal S$ of rank~$n$ is called an
\emph{$(n-1)$-dimensional pseudomanifold} (without boundary) if
\begin{itemize}
\item[(a)] for any element $\sigma\in\mathcal S$, there is an element
$\tau$ of rank $n$ such that $\sigma\le\tau$ (in other words,
$\mathcal S$ is \emph{pure $(n-1)$-dimensional});

\item[(b)] for any element $\sigma\in\mathcal S$ of rank $(n-1)$ there
are exactly two elements $\tau$ of rank $n$ such that
$\sigma<\tau$;

\item[(c)] for any two elements $\tau$ and $\tau'$ of rank $n$
there is a sequence of elements
$\tau=\tau_1,\tau_2,\ldots,\tau_k=\tau'$ such that $\rank\tau_i=n$
and $\tau_i\wedge\tau_{i+1}$ contains an element of rank $(n-1)$
for $i=1,\ldots,k-1$.
\end{itemize}
\end{definition}

Simplicial cell decompositions of topological manifolds are
pseudomanifolds, but there are pseudomanifolds that do not arise
in this way, see Example~\ref{expsm}.

\begin{theorem}[\cite{m-m-p07}]\label{theo:eqcal}\
\begin{itemize}
\item[(a)] Let $\Gamma$ be a torus graph; then $\mathcal S(\Gamma)$ is a
pseudomanifold, and the face ring $\Z[\mathcal S(\Gamma)]$ admits
an lsop;

\item[(b)] Given an arbitrary pseudomanifold $\mathcal S$ and an lsop in
$\Z[\mathcal S]$, one can canonically construct a torus graph
$\Gamma_{\mathcal S}$.
\end{itemize}
Furthermore, $\Gamma_{\mathcal S(\Gamma)}=\Gamma$.
\end{theorem}
\begin{proof}
(a) Vertices of $\mathcal S(\Gamma)$ correspond to $(n-1)$-faces
of~$\Gamma$. Since any face of $\Gamma$ contains a vertex and
$\Gamma$ is $n$-valent, $\mathcal S(\Gamma)$ is pure
$(n-1)$-dimensional. Condition~(b) from the definition of a
pseudomanifold follows from the fact that an edge of $\Gamma$ has
exactly two vertices, and~(c) follows from the connectivity
of~$\Gamma$. In order to find an lsop, we identify $\Z[\mathcal
S(\Gamma)]$ with a subset of $H^*(BT)^{V(\Gamma)}$ (see
Theorem~\ref{tequi}) and consider the constant map $H^*(BT)\to
H^*(BT)^{V(\Gamma)}$. It factors through a monomorphism
$H^*(BT)\to \Z[\mathcal S(\Gamma)]$, and Lemma~\ref{lsopscc}
implies that the image of a basis in $H^2(BT)$ is an lsop.

\smallskip

(b) Let $\mathcal S$ be a pseudomanifold of dimension $(n-1)$.
Define a graph $\Gamma_{\mathcal S}$ whose vertices correspond to
$(n-1)$-dimensional simplices $\sigma\in\mathcal S$, and in which
the number of edges between two vertices $\sigma$ and $\sigma'$ is
equal to the number of $(n-2)$-dimensional simplices in
$\sigma\wedge\sigma'$. Then $\Gamma_{\mathcal S}$ is a connected
$n$-valent graph, and we need to define an axial function.

We can regard an lsop as a map $\lambda\colon
H^*(BT)\to\Z[\mathcal S]$. Assume that $\mathcal S$ has $m$
elements of rank~1 (we do not call them vertices to avoid
confusion with the vertices of~$\Gamma$) and let $v_1,\ldots,v_m$
be the corresponding degree-two generators of $\Z[\mathcal S]$.
Then for $t\in H^2(BT)$ we can write
\[
  \lambda(t)=\sum_{i=1}^m\lambda_i(t)u_i,
\]
where $\lambda_i$ is a linear function on $H^2(BT)$, that is, an
element of $H_2(BT)$. Let $e$ be an oriented edge of $\Gamma$ with
initial vertex $v=i(e)$. Then $v$ corresponds to an
$(n-1)$-simplex of $\mathcal S$, and we denote by
$I(v)\subset\{1,\ldots,m\}$ the corresponding set of rank~1
elements of~$\mathcal S$; note that $|I(v)|=n$. Since $\lambda$ is
an lsop, the set $\{\lambda_i\colon i\in I(v)\}$ is a basis of
$H_2(BT)$. Now we define the axial function $\alpha\colon
E(\Gamma)\to H^2(BT)$ by requiring that its value on $E(\Gamma)_v$
is the dual basis of $\{\lambda_i\colon i\in I(v)\}$. In more
detail, the edge $e$ corresponds to an $(n-2)$-simplex of
$\mathcal S$ and let $\ell\in I(v)$ be the unique element which is
not in this $(n-2)$-simplex. Then we define $\alpha(e)$ by
requiring that
\begin{equation}\label{eaxfu}
  \langle\alpha(e),\lambda_i\rangle=\delta_{i\ell},\quad i\in
  I(v),
\end{equation}
where $\delta_{i\ell}$ is the Kronecker delta. We need to check
the three conditions from the definition of axial function. Let
$v'=t(e)=i(\bar e)$. Note that the intersection of $I(v)$ and
$I(v')$ consists of at least $(n-1)$ elements. If $I(v)=I(v')$
then $\Gamma$ has only two vertices, like in Example~\ref{ebana},
while $\mathcal S$ is obtained by gluing together two
$(n-1)$-simplices along their boundaries, see Example~\ref{2tria}.
Otherwise, $|I(v)\cap I(v')|=n-1$ and we have $\ell\notin I(v')$.
Let $\ell'$ be an element such that $\ell'\in I(v')$, but
$\ell'\notin I(v)$. Then~\eqref{eaxfu} implies that
$\langle\alpha(e),\lambda_i\rangle= \langle\alpha(\bar
e),\lambda_i\rangle=0$ for $i\in I(v)\cap I(v')$. As we work with
integral bases, this implies $\alpha(\bar e)=\pm\alpha(e)$. It
also follows that $\alpha(E(\Gamma)_v\setminus e)$ and
$\alpha(E(\Gamma)_{v'}\setminus\bar e)$ give the same bases in the
quotient space $H^2(BT)/\alpha(e)$. Identifying these bases, we
obtain a connection $\theta_e\colon E(\Gamma)_v\to E(\Gamma)_{v'}$
satisfying $\alpha(\theta_e(e'))\equiv\alpha(e')\mod\alpha(e)$ for
any $e'\in E(\Gamma)_v$, as needed.

The identity $\Gamma_{\mathcal S(\Gamma)}=\Gamma$ is obvious.
\end{proof}

Theorem~\ref{theo:eqcal} would have provided a complete
characterisation of simplicial posets arising from $T$-graphs if
one had $\mathcal S(\Gamma_{\mathcal S})=\mathcal S$. However,
this is not the case in general, as is shown by the next example:

\begin{example}\label{expsm}
Let $\mathcal K$ be a triangulation of a 2-dimensional sphere
different from the boundary of a simplex. Choose two vertices that
are not joined by an edge. Let $\widehat{\mathcal K}$ be the
complex obtained by identifying these two vertices. Then
$\widehat{\mathcal K}$ is a pseudomanifold. If $\Z[\mathcal K]$
admits an lsop, then $\Z[\widehat{\mathcal K}]$ also admits an
lsop (this follows easily from Lemma~\ref{lsopscc}). However,
$\mathcal S(\Gamma_{\widehat{\mathcal K}})\ne\widehat{\mathcal K}$
(in fact, $\mathcal S(\Gamma_{\widehat{\mathcal K}})=\mathcal K$).
It follows that $\widehat{\mathcal K}$ does not arise from any
$T$-graph.
\end{example}

\begin{definition}
A map $o\colon V(\Gamma)\to \{\pm 1\}$ is called an
\emph{orientation} of a $T$-graph $\Gamma$ if
$o(i(e))\alpha(e)=-o(i(\bar e))\alpha(\bar e)$ for any $e\in
E(\Gamma)$.
\end{definition}

\begin{example}
Let $M$ be a torus manifold which admits a $T$-invariant almost
complex structure. The associated axial function $\alpha_M$
satisfies $\alpha_M(\bar e)=-\alpha_M(e)$ for any oriented
edge~$e$. In this case we can take $o(v)=1$ for every $v\in
V(\Gamma_M)$.
\end{example}

\begin{proposition}\label{omniorieorie}
An omniorientation of a torus manifold $M$ induces an orientation
of the associated $T$-graph~$\Gamma_M$.
\end{proposition}
\begin{proof}
For any vertex $v\in M^T=V(\Gamma_M)$ we set $o(v)=\sigma(v)$,
where $\sigma(v)$ is the sign of~$v$ (see Lemma~\ref{qts}).
\end{proof}

\begin{example}
Let $\Gamma$ be a complete graph on four vertices
$v_1,v_2,v_3,v_4$. Choose a basis $t_1,t_2,t_3\in H^2(BT^3)$ and
define an axial function by setting
\[
  \alpha(v_1v_2)=\alpha(v_3v_4)=t_1,\quad
  \alpha(v_1v_3)=\alpha(v_2v_4)=t_2,\quad
  \alpha(v_1v_4)=\alpha(v_2v_3)=t_3
\]
and $\alpha(\bar e)=\alpha(e)$ for any oriented edge~$e$. A direct
check shows that this $T$-graph is not orientable. This graph is
associated with the pseudomanifold (simplicial cell complex) shown
in Fig.~\ref{frp2c} via the construction of
Theorem~\ref{theo:eqcal}~(b). This pseudomanifold $\mathcal S$ is
homeomorphic to $\R P^2$ (the opposite outer edges are identified
according to the arrows shown), the ring $\Z[\mathcal S]$ has
three two-dimensional generators $v_p,v_q,v_r$, which constitute
an lsop. Note that $\R P^2$ itself is non-orientable.
\begin{figure}[h]
\begin{picture}(120,30)
\put(45,0){\vector(1,0){29.5}} \put(75,0){\vector(0,1){29.5}}
\put(75,30){\vector(-1,0){29.5}} \put(45,30){\vector(0,-1){29.5}}
\put(45,0){\line(1,1){30}} \put(45,30){\line(1,-1){30}}
\put(45,0){\circle*{1}} \put(45,30){\circle*{1}}
\put(75,0){\circle*{1}} \put(75,30){\circle*{1}}
\put(60,15){\circle*{1}} \put(42,-1){$q$} \put(76,-1){$p$}
\put(42,29){$p$} \put(76,29){$q$} \put(60,16){$r$}
\end{picture}
\caption{Simplicial cell decomposition of $\mathbb RP^2$ with 3
vertices.} \label{frp2c}
\end{figure}

It follows that this $T$-graph does not arise from a torus
manifold.
\end{example}

\begin{proposition} A $T$-graph is $\Gamma$ is orientable if and only if the associated pseudomanifold
$\mathcal S(\Gamma)$ is orientable.
\end{proposition}
\begin{proof}
Let $v\in V(\Gamma)$ and let $\sigma$ be the corresponding
$(n-1)$-simplex of $\mathcal S(\Gamma)$. The oriented edges in
$E(\Gamma)_v$ canonically correspond to the vertices of~$\sigma$.
Choose a basis of $H^2(BT)$. Assume first that $\mathcal
S(\Gamma)$ is oriented. Choose a `positive' (i.e. compatible with
the orientation) order of vertices of~$\sigma$; this allows us to
regard $\alpha(E(\Gamma)_v)$ as a basis of $H^2(BT)$. We set
$o(v)=1$ if it is a positively oriented basis, and $o(v)=-1$
otherwise. This defines an orientation on~$\Gamma$. To prove the
opposite statement we just reverse this procedure.
\end{proof}

\subsection*{Blowing up $T$-manifolds and $T$-graphs}\label{blowu}
Here we elaborate on relating the following three geometric
constructions:
\begin{itemize}
\item[(a)] blowing up a torus manifold at a facial
submanifold (Construction~\ref{lstbu});

\item[(b)] truncating a simple polytope at a face (Construction~\ref{hypcut}) or, more
generally, blowing up a GKM graph or a $T$-graph;

\item[(c)] stellar subdivision of a simplicial complex or simplicial poset
(see Definition~\ref{bist} and Section~\ref{scccm}).
\end{itemize}

The construction of blow-up of a GKM-graph is described
in~\cite[\S2.2.1]{gu-za99}. It also applies to $T$-graphs and
agrees with the topological picture for graphs arising from torus
manifolds.

Let $(\Gamma,\alpha,\theta)$ be a $T$-graph and $G$ its $k$-face.
The \emph{blow-up}\label{blowupTgr} of $\Gamma$ at~$G$ is a
$T$-graph $\widetilde\Gamma$ which is defined as follows. Its
vertex set is $V(\widetilde\Gamma)=(V(\Gamma)\setminus V(G))\cup
V(G)^{n-k}$, that is, each vertex $p\in V(G)$ is replaced by $n-k$
new vertices $\widetilde p_1,\ldots,\widetilde p_{n-k}$. It is
convenient to choose these new vertices close to $p$ on the edges
from the set $E_p(\Gamma)\setminus E_p(G)$, and we denote by
$p'_i$ the endpoint of the edge of $\Gamma$ containing both $p$
and~$\widetilde p_i$. Furthermore, for any two vertices $p,q\in G$
joined by an edge $pq$ we index the corresponding new vertices of
$\widetilde G$ in the way compatible with the connection, i.e. so
that $\theta_{pq}(pp'_i)=qq'_i$. Now we need to define the edges
of the new graph $\widetilde\Gamma$ and the axial function
$\widetilde\alpha\colon E(\widetilde\Gamma)\to H^*(BT)$. We have
four types of edges in $\widetilde\Gamma$, which are given in the
following list together with the values of the axial function:
\begin{itemize}
\item[(a)] $\widetilde p_i\widetilde p_j$ for every $p\in V(G)$; \
$\widetilde\alpha(\widetilde p_i\widetilde
p_j)=\alpha(pp'_j)-\alpha(pp'_i)$;

\item[(b)] $\widetilde p_i\widetilde q_i$ if $p$ and $q$ where joined by an edge in~$G$;
\ $\widetilde\alpha(\widetilde p_i\widetilde q_i)=\alpha(pq)$;

\item[(c)] $\widetilde p_ip'_i$ for every $p\in V(G)$; \
$\widetilde\alpha(\widetilde p_ip'_i)=\alpha(pp'_i)$;

\item[(d)] edges `left over from $\Gamma$', i.e. $e\in
E(\Gamma)$ with $i(e)\notin V(G)$ and $t(e)\notin V(G)$; \
$\widetilde\alpha(e)=\alpha(e)$,
\end{itemize}
see Fig.~\ref{fblo1} ($n=3$, $k=1$) and Fig.~\ref{fblo2} ($n=3$,
$k=0$).

\begin{figure}[h]
\begin{picture}(120,40)
\put(15,30){\vector(-1,-2){15}} \put(15,30){\vector(1,-2){15}}
\put(15,30){\line(5,2){25}} \put(40,40){\vector(1,-2){15}}
\multiput(40,40)(-3,-6){4}{\line(-1,-2){1.8}}
\put(15,30){\circle*{1}} \put(40,40){\circle*{1}}
\put(70,20){\vector(-1,-2){7.5}} \put(70,20){\vector(4,-1){20}}
\put(90,15){\vector(1,-2){5}} \put(70,20){\line(5,2){25}}
\put(90,15){\line(5,2){25}} \put(95,30){\vector(4,-1){20}}
\put(115,25){\vector(1,-2){5}}
\multiput(95,30)(-3,-6){3}{\line(-1,-2){1.8}}
\put(70,20){\circle*{1}} \put(90,15){\circle*{1}}
\put(95,30){\circle*{1}} \put(115,25){\circle*{1}}
\put(0,14){$\scriptstyle\alpha(e_1)$}
\put(15,14){$\scriptstyle\alpha(e_2)$}
\put(22,36){$\scriptstyle\alpha(e_3)$}
\put(35,28){$\scriptstyle\alpha(e_4)$}
\put(45,14){$\scriptstyle\alpha(e_5)$}
\put(65,7){$\scriptstyle\alpha(e_1)$}
\put(95,7){$\scriptstyle\alpha(e_2)$}
\put(70,15){$\scriptscriptstyle\alpha(e_2)-\alpha(e_1)$}
\put(97,16){$\scriptstyle\alpha(e_3)$}
\put(77,26){$\scriptstyle\alpha(e_3)$}
\put(111,17){$\scriptstyle\alpha(e_5)$}
\put(100,29){$\scriptscriptstyle\alpha(e_5)-\alpha(e_4)$}
\put(65,25){\vector(-1,0){10}} \put(59,26){$b$}
\put(20,-5){\large$\Gamma$} \put(95,-5){\large$\widetilde\Gamma$}
\end{picture}
\medskip
\caption{Blow-up at an edge} \label{fblo1}
\end{figure}

\begin{figure}[h]
\begin{picture}(120,40)
\put(35,40){\vector(-4,-3){35}} \put(35,40){\vector(-1,-4){10}}
\put(35,40){\vector(2,-3){20}} \put(35,40){\circle*{1}}
\put(75,30){\vector(-4,-3){15}} \put(75,30){\vector(1,-1){12.6}}
\put(75,30){\vector(1,0){29.5}} \put(88,17){\vector(-1,-4){4}}
\put(87.5,16.8){\vector(4,3){17.2}}
\put(105,30){\vector(3,-4){15}} \put(75,30){\circle*{1}}
\put(88,17){\circle*{1}} \put(105,30){\circle*{1}}
\put(2,21){$\scriptstyle\alpha(e_1)$}
\put(27.5,7){$\scriptstyle\alpha(e_2)$}
\put(46,11){$\scriptstyle\alpha(e_3)$}
\put(55,21){$\scriptstyle\alpha(e_1)$}
\put(86.5,7){$\scriptstyle\alpha(e_2)$}
\put(70,19){$\scriptscriptstyle\alpha(e_2)-\alpha(e_1)$}
\put(111,11){$\scriptstyle\alpha(e_3)$}
\put(82,31){$\scriptscriptstyle\alpha(e_3)-\alpha(e_1)$}
\put(92,19){$\scriptscriptstyle\alpha(e_3)-\alpha(e_2)$}
\put(62,28){\vector(-1,0){10}} \put(56,29){$b$}
\put(20,-5){\large$\Gamma$} \put(95,-5){\large$\widetilde\Gamma$}
\end{picture}
\medskip
\caption{Blow-up at a vertex} \label{fblo2}
\end{figure}

There is a \emph{blow-down map}\label{deblowdown}
$b\colon\widetilde\Gamma\to\Gamma$ taking faces to faces. The face
$G$ is blown up to a new facet $\widetilde
G\subset\widetilde\Gamma$ (unless $G$ is a facet itself, in which
case $\widetilde\Gamma=\Gamma$). For any face $H\subset\Gamma$ not
contained in~$G$, there is a unique face $\widetilde
H\subset\widetilde\Gamma$ that is blown down onto~$H$. The
blow-down map induces a homomorphism in equivariant cohomology
$b^*\colon H^*_T(\Gamma)\to H^*_T(\widetilde\Gamma)$, which is
defined by the following commutative diagram
\begin{equation}
\label{egrre}
\begin{CD}
  H^*_T(\Gamma) @>b^*>> H^*_T(\widetilde\Gamma)\\
  @VrVV @VV{\widetilde r}V\\
  H^*(BT)^{V(\Gamma)} @>V\!b^*>> H^*(BT)^{V(\widetilde\Gamma)},
\end{CD}
\end{equation}
where $r$ and $\widetilde r$ are the monomorphisms from the
definition of equivariant cohomology of a $T$-graph, and $V(b)^*$
is the homomorphism induced by the set map $V(b)\colon
V(\widetilde\Gamma)\to V(\Gamma)$. The next lemma describes the
images of the two-dimensional generators $\ta{F}\in H^*_T(\Gamma)$
corresponding to facets $F\subset\Gamma$.

\begin{lemma}
For a given facet $F\subset\Gamma$, we have
$b^*(\ta{F})=\ta{\widetilde G}+\ta{\widetilde F}$, if $G\subset F$
and $b^*(\ta{F})=\ta{\widetilde F}$ otherwise.
\end{lemma}
\begin{proof}
We consider diagram~\eqref{egrre} and check that the images of
$b^*(\ta{F})$ and $\ta{\widetilde G}+\ta{\widetilde F}$ (or
$\ta{\widetilde F}$) under the map $\widetilde r$ agree. Take a
vertex $p\in V(\Gamma)$. If $p\notin G$, then $b^{-1}(p)=p$ and
$r(\ta{F})(p)=\widetilde r(\ta{\widetilde F})(p)$, \ $\widetilde
r(\ta{\widetilde G})(p)=0$. Now let $p\in G$; then $b(\widetilde
p_i)=p$ for $1\le i\le n-k$.

First assume $G\not\subset F$ (see Fig.~\ref{fres1}). If $p\notin
F$, then $r(\ta{F})(p)=\widetilde r(\ta{\widetilde F})(\widetilde
p_i)=0$. Otherwise $p\in F\cap G$. Let $e$ be the unique edge such
that $e\in E_p(\Gamma)$ and $e\notin F$. Then $e=pq$ for some
$q\in V(G)$ (because $G\not\subset F$). From~\eqref{etaf1} we
obtain
\[
  r(\ta{F})(p)=\alpha(pq)=\widetilde\alpha(\widetilde p_i\widetilde q_i)=
  \widetilde r(\ta{\widetilde F})(\widetilde p_i),
  \quad 1\le i\le n-k.
\]
Then $V\!b^* r(\ta{F})=\widetilde r(\ta{\widetilde F})$, and
therefore, $b^*(\ta{F})=\ta{\widetilde F}$.
\begin{figure}[h]
\begin{picture}(120,40)
\put(15,30){\line(-1,-2){15}} \put(15,30){\line(1,-2){15}}
\put(15,30){\line(5,2){25}} \put(40,40){\line(1,-2){15}}
\put(15,30){\circle*{1}} \put(40,40){\circle*{1}}
\put(70,20){\line(-1,-2){7.5}} \put(70,20){\line(4,-1){20}}
\put(90,15){\line(1,-2){5}} \put(70,20){\line(5,2){25}}
\put(90,15){\line(5,2){25}} \put(95,30){\line(4,-1){20}}
\put(115,25){\line(1,-2){5}} \put(70,20){\circle*{1}}
\put(90,15){\circle*{1}} \put(95,30){\circle*{1}}
\put(115,25){\circle*{1}}
\put(26,35.5){$G$} \put(13,31){$p$} \put(38.5,37){$q$}
\put(14,10){$F$}
\put(92,21){$\widetilde G$} \put(88,11){$\widetilde p_i$}
\put(113,21){$\widetilde q_i$} \put(77,10){$\widetilde F$}
\put(65,25){\vector(-1,0){10}} \put(59,26){$b$}
\end{picture}
\caption{} \label{fres1}
\end{figure}

Now assume $G\subset F$ (see Fig.~\ref{fres2}). In this case the
unique edge~$e$ for which $e\in E_p(\Gamma)$ and $e\notin F$ is of
type $pp'_j$. Using~\eqref{etaf1} we calculate
\begin{align*}
  r(\ta{F})(p)&=\alpha(pp'_j),\\
  \widetilde r(\ta{\widetilde F})(\widetilde p_i)&=\widetilde\alpha(\widetilde p_i\widetilde p_j)
  =\alpha(pp'_j)-\alpha(pp'_i),\\
  \widetilde r(\ta{\widetilde G})(\widetilde p_i)&=\widetilde\alpha(\widetilde
  p_ip'_i)=\alpha(pp'_i),\quad 1\le i\le n-k.
\end{align*}
Then $V\!b^* r(\ta{F})=\widetilde r(\ta{\widetilde G})+\widetilde
r(\ta{\widetilde F})$, and therefore, $b^*(\ta{F})=\ta{\widetilde
G}+\ta{\widetilde F}$.
\begin{figure}[h]
\begin{picture}(120,40)
\put(15,30){\line(-1,-2){15}} \put(15,30){\line(1,-2){15}}
\put(15,30){\line(5,2){25}} \put(40,40){\line(1,-2){15}}
\put(15,30){\circle*{1}} \put(40,40){\circle*{1}}
\put(0,0){\circle*{1}} \put(30,0){\circle*{1}}
\put(70,20){\vector(-1,-2){7.3}} \put(70,20){\line(4,-1){20}}
\put(90,15){\vector(1,-2){4.8}} \put(70,20){\line(5,2){25}}
\put(90,15){\line(5,2){25}} \put(90,15){\vector(-4,1){19.5}}
\put(115,25){\vector(-4,1){19.5}} \put(115,25){\line(1,-2){5}}
\put(70,20){\circle*{1}} \put(90,15){\circle*{1}}
\put(95,30){\circle*{1}} \put(115,25){\circle*{1}}
\put(62.5,5){\circle*{1}} \put(95,5){\circle*{1}}
\put(26,35.5){$G$} \put(13,31){$p$} \put(1.5,0){$p'_j$}
\put(25.5,0){$p'_i$} \put(35,15){$F$}
\put(92,21){$\widetilde G$} \put(88,11){$\widetilde p_i$}
\put(70,16){$\widetilde p_j$} \put(91,4){$p'_i$}
\put(64,4){$p'_j$} \put(105,12){$\widetilde F$}
\put(65,25){\vector(-1,0){10}} \put(59,26){$b$}
\end{picture}
\caption{} \label{fres2}
\end{figure}
\end{proof}

\begin{corollary}
After the identifications $H^*_T(\Gamma)\cong\Z[\mathcal
S(\Gamma)]$ and $H^*_T(\widetilde\Gamma)\cong\Z[\mathcal
S(\widetilde\Gamma)]$, the equivariant cohomology homomorphism
$b^*$ induced by the blow-down map
$b\colon\widetilde\Gamma\to\Gamma$ coincides with the homomorphism
$\beta$ from Lemma~{\rm\ref{lemm:frppt}}.
\end{corollary}
\begin{proof}
Recall from Theorem~\ref{tequi} that the poset $\mathcal
S(\Gamma)$ is formed by faces of $\Gamma$ with the reversed
inclusion relation, and the isomorphism
$H^*_T(\Gamma)\cong\Z[\mathcal S(\Gamma)]$ is established by
identifying Thom classes $\ta{H}$ with generators $\va{H}$
corresponding to faces $H\subset\Gamma$. Let $\sigma\in\mathcal
S(\Gamma)$ be the element corresponding to the blown up face~$G$.
Then an element $\tau\in\mathcal S(\Gamma)$ belongs to
$\st_{\mathcal S(\Gamma)}\sigma$ if and only if its corresponding
face $H\subset\Gamma$ satisfies $G\cap H\ne\varnothing$. The
degree-two generators $v_1,\ldots,v_m$ of $\Z[\mathcal S(\Gamma)]$
and of $\Z[\mathcal S(\widetilde\Gamma)]$ correspond to the
generators $\ta{F_1},\ldots,\ta{F_m}$ of $H^*_T(\Gamma)$ and the
generators $\ta{{\widetilde F}_1},\ldots,\ta{{\widetilde F}_m}$ of
$H^*_T(\widetilde\Gamma)$, respectively. Making the appropriate
identifications, we see that the homomorphism from
Lemma~\ref{lemm:frppt} is determined uniquely by the conditions
\begin{align*}
  \ta{H} & \mapsto\ta{H} &&\text{if \ }G\cap H=\varnothing,\\
  \ta{F_i} & \mapsto\ta{\widetilde G}+\ta{\widetilde F_i} &&\text{if \ }G\subset F_i,\\
  \ta{F_i} & \mapsto\ta{\widetilde F_i} &&\text{if \ }G\not\subset F_i.
\end{align*}
The blow-down map $b^*$ satisfies these conditions, thus finishing
the proof.
\end{proof}

\subsection*{Exercises}
\begin{exercise}
A connection $\theta$ satisfying condition~(c) from the definition
of axial function is unique if elements of $\alpha(E(\Gamma)_v)$
are 3-independent for any vertex~$v$.
\end{exercise}

%Единственность пр-ва DJ(K) - результат Notbohm-Ray
%
%момент-угол-комплексы для произвольных (не простых) многогранников
%(Айзенберг-Бухштабер)
%изучить кольцо когомологий конструкции Бореля на таком м-у-комплексе
%
%Доказательство результата Фрёберга
%
%Discuss minimal non-Golodness

\chapter{Homotopy theory of polyhedral products}\label{homotopy}
The homotopy-theoretical study of toric spaces, such as
moment-angle complexes or general polyhedral products, has
recently evolved into a separate branch linking toric topology to
unstable homotopy theory. Like toric varieties in algebraic
geometry, polyhedral product spaces provide an effective testing
ground for many important homotopy-theoretical techniques.

Basic homotopical properties of polyhedral products were described
in our 2002 text~\cite{bu-pa02}; these are included in
Section~\ref{basichomotopy} of this book. Two important
developments have followed shortly. First, Grbi\'c and
Theriault~\cite{gr-th04},~\cite{gr-th07} described a wide class of
simplicial complexes $\sK$ whose corresponding moment-angle
complex~$\zk$ is homotopy equivalent to a wedge of spheres. (This
class includes, for example, $i$-dimensional skeleta of a simplex
$\varDelta^{m-1}$ for all $i$ and~$m$.) Second, formality of the
Davis--Januszkiewicz space $\djs(\sK)$ (or, equivalently, the
polyhedral product~$(\C P^\infty)^{\sK}$) was established by
Notbohm and Ray in~\cite{no-ra05}. (An alternative proof of this
result was also given in~\cite[Lemma~7.35]{bu-pa04-2}.) The
importance of the homotopy-theoretical viewpoint on toric spaces
has been emphasised in two papers~\cite{bu-ra08}
and~\cite{pa-ra08}, both coauthored with Ray and appeared in the
proceedings of Osaka~2006 conference on toric topology. Following
the earlier work of Panov, Ray and Vogt~\cite{p-r-v04},
in~\cite{pa-ra08} categorical methods have been brought to bear on
toric topology. We include the main results of~\cite{pa-ra08}
and~\cite{p-r-v04} in Sections~\ref{models} and~\ref{loops}.

The idea is to exhibit a toric space as the homotopy colimit of a
diagram of spaces over the small category {\sc cat}$(\sK)$, whose
objects are the faces of a finite simplicial complex $\sK$ and
morphisms are their inclusions. The corresponding {\sc
cat}$(K)$-diagrams can also be studied in various algebraic
Quillen model categories, and their homotopy (co)limits can be
interpreted as algebraic models for toric spaces. Such models
encode many standard algebraic invariants, and their existence is
assured by the Quillen structure. Several illustrative
calculations will be provided. In particular, we it is proved that
toric and quasitoric manifolds (and various generalisations) are
rationally formal, and that the polyhedral product $(\C
P^\infty)^{\sK}$ (or the Davis--Januszkiewicz space) is coformal
precisely when $\sK$ is flag.

A number of papers on the homotopy-theoretical aspects of toric
spaces has appeared since 2008. One of the most important
contributions was the 2010 work of Bahri, Bendersky, Cohen and
Gitler~\cite{b-b-c-g10} establishing a decomposition of the
suspension of a polyhedral product $(X,A)^{\sK}$ into a wedge of
suspensions corresponding to subsets $I\subset[m]$. In particular,
the moment-angle complex $\zk$ breaks up into a wedge of
suspensions of full subcomplexes $\sK_I$ after one suspension,
providing a homotopy-theoretical interpretation of the cohomology
calculation of Theorem~\ref{zkadd}. We include the results
of~\cite{b-b-c-g10} and related results on stable decompositions
of polyhedral products in Section~\ref{stabledec}.

In Section~\ref{loops} we study the loop spaces on moment-angle
complexes and polyhedral products, by applying both categorical
decompositions and the classical homotopy-theoretical approach via
the (higher) Whitehead and Samelson products.

In the last section we restrict our attention to the case of flag
complexes~$\sK$. We describe the Pontryagin algebra
$H_*(\varOmega(\C P^\infty)^\sK)$ explicitly by generators and
relations in Theorem~\ref{hldj}, and also exhibit it as a colimit
(or \emph{graph product}) of exterior algebras. For the commutator
subalgebra $H_*(\varOmega\zk)$ a minimal set of generators was
constructed in~\cite{g-p-t-w}; it is included as
Theorem~\ref{multgen}. Another result of~\cite{g-p-t-w}
(Theorem~\ref{flws}) gives a complete characterisation of the
class of flag complexes~$\sK$ for which the moment-angle complex
$\zk$ is homotopy equivalent to a wedge of spheres.

Background material on model categories and homotopy (co)limits is
given in Appendix~\ref{catconstr}, further references can be also
found there.

\section{Rational homotopy theory of polyhedral products}\label{models}
Here we construct several decompositions of polyhedral products
into colimits and homotopy colimits of diagrams over $\ca(\sK)$.
We also establish formality of polyhedral powers $X^\sK$ with
formal $X$, as well as formality of (quasi)toric manifolds and
some torus manifolds. This contrasts the situation with
moment-angle complexes $\zk=(D^2,S^1)^\sK$, which are not formal
in general (see Section~\ref{Masseymac}).

The indexing category for all diagrams in this section is
$\ca(\sK)$\label{facecate} (simplices of a finite simplicial
complex $\sK$ and their inclusions) or its opposite
$\ca(\sK)^{op}$. We recall the following diagrams and their
(co)limits which appeared earlier in this book:
\begin{itemize}
\item[$\cdot$] $\ca(\sK)^{op}$-diagram $\k[\,\cdot\:]^\sK$ in $\cat{cga}$, see
Exercise~\ref{frlimit} and Lemma~\ref{zslim}; its limit is the
face ring $\k[\sK]$;

\item[$\cdot$] $\ca(\sK)$-diagram $\mathcal D_\sK(D^2,S^1)$ in $\cat{top}$, see~\eqref{ddsdiag};
its colimit is the moment-angle complex $\zk=(D^2,S^1)^\sK$;

\item[$\cdot$] $\ca(\sK)$-diagram $\mathcal D_\sK(\mb X,\mb A)$ in $\cat{top}$, see~\eqref{dxadiag};
its colimit is the polyhedral product $(\mb X,\mb A)^\sK$;
\end{itemize}

We can generalise the first diagram above as follows. Given a
sequence $\mb C=(C_1,\ldots,C_m)$ of commutative dg-algebras
over~$\Q$, define the diagram
\begin{equation}\label{DsupK}
  \mathcal D^\sK(\mb
  C)\colon\ca(\sK)^{op}\to\cat{cdga},\qquad I\mapsto \bigotimes_{i\in
  I}C_i,
\end{equation}
by mapping a morphism $I\subset J$ to the surjection
$\bigotimes_{i\in J}C_i\to\bigotimes_{i\in I}C_i$ sending each
$C_i$ with $i\notin I$ to~$1$. Then the diagram
$\k[\,\cdot\:]^\sK$ corresponds to the case $C_i=\k[v]$, the
polynomial algebra on one generator of degree~2.

\begin{proposition}\label{Dcofib}\
\begin{itemize}
\item[(a)] The diagram $\mathcal D^\sK(\mb C)$ is Reedy fibrant. Therefore, there is a weak equivalence
$\lim\mathcal D^\sK(\mb
C)\stackrel\simeq\longrightarrow\holim\mathcal D^\sK(\mb C)$. In
particular, there is a weak equivalence
$\Q[\sK]=\lim\Q[\,\cdot\:]^\sK\stackrel\simeq\longrightarrow\holim\Q[\,\cdot\:]^\sK$.

\item[(b)] The diagram $\mathcal D_\sK(\mb X,\mb A)$ is Reedy cofibrant
whenever each $A_i\to X_i$ is a cofibration (e.g. when $(X_i,A_i)$
is a cellular pair). Under this condition, there is a weak
equivalence $\hocolim\mathcal D_\sK(\mb X,\mb
A)\stackrel\simeq\longrightarrow(\mb X,\mb A)^{\sK}$.
\end{itemize}
\end{proposition}
\begin{proof}
(a) Recall from Section~\ref{secmc} that a
$\cat{cat}^{op}(\sK)$-diagram $\mathcal C$ is Reedy fibrant when
the canonical map $\mathcal C(I)\rightarrow\lim \mathcal
C|_{\scat{cat}^{op}(\partial\varDelta(I))}$ is a fibration for
each $I\in\sK$. In our case,
\[
  \mathcal D^\sK(\mb C)(I)=\bigotimes_{i\in
  I}C_i,\qquad
  \lim\mathcal D^\sK(\mb C)|_{\scat{cat}^{op}(\partial\varDelta(I))}
  =\bigotimes_{i\in I}C_i/\mathcal I,
\]
where $\mathcal I$ is the ideal generated by all products
$\prod_{i\in I}c_i$ with $c_i\in C_i^+$. Hence the Reedy fibrance
condition is satisfied. Proposition~\ref{hotolim} implies that the
canonical map $\lim\mathcal D^\sK(\mb
C)\stackrel\simeq\longrightarrow\holim\mathcal D^\sK(\mb C)$ is a
is a weak equivalence.

(b) A $\ca{(\sK)}$-diagram $\mathcal D$ in $\cat{top}$ is Reedy
cofibrant whenever each map $\colim \mathcal
D|_{\scat{cat}(\partial\varDelta(I))}\rightarrow \mathcal D(I)$ is
a cofibration. In our case,
\[
  \colim \mathcal D_\sK(\mb X,\mb
  A)|_{\scat{cat}(\partial\varDelta(I))}=(\mb X,\mb
  A)^{\partial\varDelta(I)}\times\mb A^{[m]\setminus I},\quad
  \mathcal D_\sK(\mb X,\mb A)(I)=(\mb X,\mb A)^I,
\]
so the Reedy cofibrance condition is satisfied.
\end{proof}

We now consider more specific models for particular polyhedral
products and quasitoric manifolds. All cohomology in this section
is with rational coefficients.

\subsection*{Formality of $\mb X^{\sK}$ and $(\C P^\infty)^\sK$}
Recall that we use the notation $\mb X^{\sK}$ for the polyhedral
product $(\mb X,\pt)^{\sK}$, and a space $X$ is formal if the
commutative dg-algebra $\APL(X)=A^*(S_\bullet X)$ is weakly
equivalent to its cohomology $H^*(X)$.

\begin{theorem}\label{xkform}
If each space $X_i$ in $\mb X=(X_1,\ldots,X_m)$ is formal, then
the polyhedral product $\mb X^\sK$ is also formal.
\end{theorem}

\begin{proof}
For notational clarity, we denote the colimit of the diagram
$\mathcal D_{\mathcal K}(\mb X,\pt)$ defining $\mb X^\sK$ by
$\colim_I\mb X^I$, where $I\in\sK$. We shall prove that
$\APL=A^*S_\bullet$ maps this colimit to the limit of dg-algebras
$\APL(\mb X^I)\cong\bigotimes_{i\in I}\APL(X_i)$. For this it is
convenient to work with simplicial sets, as the definition of the
polynomial de Rham functor $A^*$ implies that it takes colimits in
$\cat{sset}$ to limits in $\cat{cdga}$,
see~\cite[\S13.5]{bo-gu76}.

We have a natural weak equivalence $|S_\bullet(\mb X^I)|\to\mb
X^I$. Since the total singular complex functor $S_\bullet$ is
right adjoint, it preserves products, so we have an equivalence
$|(S_\bullet\mb X)^I|\to\mb X^I$. The $\ca{(\sK)}$-diagram
$\mathcal D_{\mathcal K}(\mb X,\pt)$ given by $I\mapsto\mb X^I$ is
cofibrant by Proposition~\ref{Dcofib}, and the diagram
$I\mapsto|(S_\bullet\mb X)^I|$ is cofibrant by the same reason. We
therefore have a weak equivalence $\colim_I|(S_\bullet\mb
X)^I|\to\colim_I\mb X^I$ by Proposition~\ref{wecolim}. Applying
$\APL$, we obtain the zigzag
\[
  \APL(\mb X^\sK)=\APL\colim_I\mb X^I\stackrel\simeq\longrightarrow
  \APL\colim_I|(S_\bullet\mb X)^I|\cong
  \APL|\colim_I(S_\bullet\mb X)^I|,
\]
where in the last identity we used the fact that the realisation
functor is left adjoint and therefore preserves colimits. Given an
arbitrary simplicial set $Y_\bullet$, there is a quasi-isomorphism
$\APL(|Y_\bullet|)=A^*(S_\bullet|Y_\bullet|)\stackrel\simeq\longrightarrow
A^*(Y_\bullet)$ induced by the equivalence $Y_\bullet\to
S_\bullet|Y_\bullet|$. We therefore can continue the zigzag above
as
\[
  \APL|\colim_I(S_\bullet\mb X)^I|\stackrel\simeq\longrightarrow A^*(\colim_I(S_\bullet\mb
  X)^I)\cong\lim_I A^*(S_\bullet\mb X)^I=\lim_I\APL(\mb X^I).
\]

Now, since each $X_i$ is formal, there is a zigzag of
quasi-isomorphisms $\APL(X_i)\gets\cdots\to H^*(X_i)$. Applying
Proposition~\ref{Dcofib}~(a) for the case $C_i=\APL(X_i)$ and
$C_i=H^*(X_i)$ we obtain that both the corresponding diagrams
$\mathcal D^\sK(\mb C)$ are fibrant, so their limits are weakly
equivalent by Proposition~\ref{wecolim}:
\[
  \lim_I\APL(\mb X^I)\stackrel\simeq\longleftarrow\cdots
  \stackrel\simeq\longrightarrow\lim_I H^*(\mb X^I)
\]
(here we also use the fact that $H^*(\mb X^I)\cong\bigotimes_{i\in
I}H^*(X_i)$, as we work with rational coefficients). The proof is
finished by appealing to the isomorphism
\[
  \lim_I H^*(\mb X^I)\cong H^*(\mb X^\sK).
\]
In the case $X_i=\C P^\infty$ this is proved in
Proposition~\ref{homsrs} (we only need this case in this section).
For the general case, the proof will be given in
Section~\ref{stabledec}.
\end{proof}

\begin{corollary}\label{djsformalsp}
The Davis--Januszkiewicz space $ET^m\times_{T^m}\zk$ is formal.
\end{corollary}
\begin{proof}
This follows from the homotopy equivalence $(\C
P^\infty)^\sK\simeq ET^m\times_{T^m}\zk$ of Theorem~\ref{zkhofib}.
\end{proof}

\begin{remark}
The case $X_i=\C P^\infty$ of Theorem~\ref{xkform} (formality of
the Davis--Januszkiewicz space) was proved
in~\cite[Theorem~5.5]{no-ra05} and~\cite[Lemma~7.35]{bu-pa04-2}.
According to~\cite[Theorem~4.8]{no-ra05}, the space $(\C
P^\infty)^\sK$ is \emph{integrally}\label{intformality} formal,
i.e. the singular cochain algebra $C^*((\C P^\infty)^\sK;\Z)$ is
formal as a non-commutative dg-algebra.

The result of Theorem~\ref{xkform} cannot be extended to
polyhedral products of the form $(\mb X,\mb A)^\sK$. Although
$\lim_I\APL((\mb X,\mb A)^I)$ is still a model for $\APL(\mb X,\mb
A)^\sK$ (see the next subsection), the $\ca(\sK)^{op}$-diagram
$I\mapsto H^*((\mb X,\mb A)^I)$ is \emph{not} fibrant in general,
and therefore its limit is neither isomorphic to $\lim_I\APL((\mb
X,\mb A)^I)$, nor to $H^*((\mb X,\mb A)^\sK)$. Indeed, as we have
seen in Section~\ref{Masseymac}, the moment-angle complex
$\zk=(D^2,S^1)^\sK$ is not formal in general.

The argument in the proof of Theorem~\ref{xkform} shows that
rational cohomology, as a functor from $\cat{top}$ to
$\cat{cdga}$, maps homotopy colimit of diagram $\mathcal D_\sK(\mb
X,\pt)$ %with formal $\mb X$
to homotopy limit of diagram $\mathcal
D^{\sK}(H^*(\mb X))$.

The coformality of $(\C P^\infty)^\sK$ is explored in
Theorem~\ref{djcoform}.
\end{remark}

\subsection*{Models for $(\mb X,\mb A)^\sK$ and $\zk$}
Given a polyhedral product $(\mb X,\mb A)^\sK=\colim_I(\mb X,\mb
A)^I$, we can consider the $\ca({\sK})^{op}$-diagram of
commutative dg-algebras defined by $I\mapsto\APL((\mb X,\mb
A)^I)$, and denote its limit by $\lim_I\APL((\mb X,\mb A)^I)$.

\begin{proposition}\label{aplXA}
There is a quasi-isomorphism of commutative dg-algebras
\[
  \APL((\mb X,\mb A)^\sK)\stackrel\simeq\longrightarrow
  \lim_I\APL((\mb X,\mb A)^I).
\]
\end{proposition}
\begin{proof}
Repeat literally the argument in the first part of proof of
Theorem~\ref{xkform} (before appealing to the formality of~$X_i$).
\end{proof}

More specific models can be obtained in the particular case
$\zk=(D^2,S^1)^{\sK}$.

Alongside with the diagram $\mathcal D_\sK(D^2,S^1)$ we consider
the $\ca{(\sK)}$-diagram $\mathcal D_\sK(\pt,S^1)$ in $\cat{top}$
whose value on $I\subset J$ is the quotient map of tori
\begin{equation}\label{dkpts}
  T^m/T^I=(S^1)^{[m]\setminus I}\to(S^1)^{[m]\setminus J}=T^m/T^J.
\end{equation}
This diagram is not cofibrant; we denote its homotopy colimit by
$\hocolim_I T^m/T^I$.

\begin{proposition}[{\cite{p-r-v04}}]\label{zkhocolim}
There is a weak equivalence $\zk\!\simeq\hocolim_I T^m\!/T^I$.
\end{proposition}
\begin{proof}
Objectwise projections $(D^2,S^1)^I\to(S^1)^{[m]\setminus I}$
induce a weak equivalence of diagrams $\mathcal D_\sK(D^2,S^1)\to
\mathcal D_\sK(\pt,S^1)$, whose source is Reedy cofibrant but
whose target is not. Proposition \ref{hclhlpreswe} therefore
determines a weak equivalence
\[
  \zk=\colim_I(D^2,S^1)^I\stackrel{\simeq}\longrightarrow\hocolim_I
  T^m/T^I.\qedhere
\]
\end{proof}

In order to obtain a rational model of $\zk$ from this homotopy
limit decomposition, we consider $\cat{cat}^{op}(\sK)$-diagrams
$\Lambda[m]\otimes\Q[\,\cdot\:]^\sK$ and
$\Lambda[m]/\Lambda[\,\cdot\:]^\sK$ in $\cat{cdga}$. The first is
obtained by objectwise tensoring the diagram $\Q[\,\cdot\:]^\sK$
with the exterior algebra $\Lambda[m]=\Lambda[u_1,\ldots,u_m]$,
$\deg u_i=1$, and imposing the standard Koszul differential. So
the value of $\Lambda[m]\otimes\Q[\,\cdot\:]^\sK$ on $I\subset J$
is the quotient map
\[
  \bigl(\Lambda[m]\otimes\Q[v_i\colon i\in J],d\bigr)
  \to \bigl(\Lambda[m]\otimes\Q[v_i\colon i\in I],d\bigr),
\]
where $d$ is defined on $\Lambda[m]\otimes\Q[v_i\colon i\in I]$ by
$du_i=v_i$ for $i\in I$ and $du_i=0$ otherwise. The value of the
second diagram $\Lambda[m]/\Lambda[\,\cdot\:]^\sK$ on $I\subset J$
is the monomorphism
\[
  \Lambda[u_i\colon i\notin J]\to\Lambda[u_i\colon i\notin I]
\]
of algebras with zero differential. Objectwise projections induce
a weak equivalence
\begin{equation}\label{uskuk}
  \Lambda[m]\otimes\Q[\,\cdot\:]^\sK\to\Lambda[m]/\Lambda[\,\cdot\:]^\sK
\end{equation}
in $\fcat{cat$^{op}(\sK)$}{cdga}$, whose source is Reedy fibrant
but whose target is not.

The following result describes the first algebraic model of $\zk$,
and also recovers the main result of Section~\ref{cohma} with
rational coefficients.

\begin{theorem}[\cite{pa-ra08}]\label{zkmodel1}
The commutative differential graded algebras $\APL(\zk)$ and
$\holim\Lambda[m]/\Lambda[\,\cdot\:]^\sK=
\holim_I\Lambda[u_i\colon i\notin I]$ are weakly equivalent in
$\cat{cdga}$.
\end{theorem}
\begin{proof}
We first construct models for the products of disks and circles
$(D^2,S^1)^I$ which are compatible with the inclusions
$(D^2,S^1)^I\subset(D^2,S^1)^J$ forming the diagram $\mathcal
D_\sK(D^2,S^1)$. The model of a single disk $D^2$ is the Koszul
algebra $(\Lambda[u]\otimes\Q[v],d)$. The map
$\Lambda[u]\otimes\Q[v]\to\APL(D^2)$ takes $u$ to the form
$\omega=xdy-ydx$ (where $(x,y)\in D^2$ are the standard cartesian
coordinates), and takes $v$ to its differential $2dx\wedge dy$. It
is important that the map $\Lambda[u]\otimes\Q[v]\to\APL(D^2)$
restricts to the quasi-isomorphism $\Lambda[u]\to\APL(S^1)$ taking
$u$ to the form $d\varphi$ representing a generator of $H^1(S^1)$.
(Here $\varphi$ is the polar angle; note that the restriction of
$xdy-ydx=r^2d\varphi$ to $S^1$ is closed, but not exact, because
$\varphi$ is not a globally defined function.) Taking product
yields compatible quasi-isomorphisms
\[
  \bigl(\Lambda[u_1,\ldots,u_m]\otimes\Q[v_i\colon i\in I],d\bigr)
  \stackrel\simeq\longrightarrow
  \APL((D^2,S^1)^I),
\]
and therefore a weak equivalence of Reedy fibrant diagrams in
$\cat{cdga}$. Their limits are therefore quasi-isomorphic, which
together with the quasi-isomorphisms of Proposition~\ref{aplXA},
Proposition~\ref{hotolim} and~\eqref{uskuk} gives the required
zigzag
\begin{multline*}
  \APL(\zk)=\APL\bigl((D^2,S^1)^\sK\bigr)
  \stackrel\simeq\longrightarrow
  \lim_I\APL\bigr((D^2,S^1)^I\bigl)\\
  \stackrel{\simeq}{\longleftarrow}
  \lim_I\bigl(\Lambda[u_1,\ldots,u_m]\otimes\Q[v_i\colon i\in I],d\bigr)
  \stackrel{\simeq}{\longrightarrow}
  \holim_I\Lambda[u_i\colon i\notin I].\qedhere
\end{multline*}
\end{proof}

On the other hand, the zigzag above together with
Exercise~\ref{frlimit} (or Lemma~\ref{zslim}) gives a weak
equivalence
\begin{equation}\label{aplzkkoszul}
  \APL(\zk)\simeq\lim_I\bigl(\Lambda[m]\otimes\Q[v_i\colon i\in
  I],d\bigr)=\bigl(\Lambda[m]\otimes\Q[\sK],d\bigr),
\end{equation}
which implies the isomorphism of Theorem~\ref{zkcoh} in the case
of rational coefficients:
\[
  H^*(\zk;\Q)\cong
  H\bigl(\Lambda[u_1,\ldots,u_m]\otimes\Q[\sK],d\bigr).
\]

To obtain the second algebraic model for $\zk$ we regard it as the
homotopy fibre of the inclusion $(\C
P^\infty)^\sK\hookrightarrow(\C P^\infty)^m=BT^m$ (see
Theorem~\ref{zkhofib}). We therefore can identify $\zk$ with the
limit of the following diagram in $\cat{top}$:
\begin{equation}\label{zkpul}
\begin{CD}
   @. ET^m\\
  @. @VVV\\
  (\C P^\infty)^\sK @>>> BT^m
\end{CD}
\end{equation}
This diagram is fibrant for the appropriate Reedy structure, as
$ET^m\to BT^m$ is a fibration (see Proposition~\ref{hopu}~(b)), so
its limit is weakly equivalent to the homotopy limit.

Now we apply $\APL$ to~\eqref{zkpul}, on the understanding that it
does not generally convert pullbacks to
pushouts~\cite[\S3]{bo-gu76}. We obtain the diagram
\begin{equation}\label{aplzkpul}
\begin{CD}
  \APL(BT^m) @>>> \APL(ET^m)\\
  @VVV @.\\
  \APL((\C P^\infty)^\sK) @.
\end{CD}
\end{equation}
in $\cat{cdga}$, which is not Reedy cofibrant. Here is the second
model for~$\zk$:

\begin{theorem}[\cite{pa-ra08}]\label{zkmodel2}
The algebra $\APL(\zk)$ is weakly equivalent the homotopy colimit
of~\eqref{aplzkpul} in $\cat{cdga}$.
\end{theorem}
\begin{proof}
There is an objectwise weak equivalence mapping the diagram
\[
\begin{CD}
  \Q[v_1,\ldots,v_m] @>>> \bigl(\Lambda[u_1,\ldots,u_m]\otimes\Q[v_1,\ldots,v_m],d\bigr)\\
  @VVV @.\\
  \Q[\sK] @.
\end{CD}
\]
to~\eqref{aplzkpul}; here the upper arrow
$\Q[m]\to(\Lambda[m]\otimes\Q[m],d)$ is the standard model for the
fibration $ET^m\to BT^m$, and $\Q[\sK]$ is a model for $\APL((\C
P^\infty)^\sK)$ by Theorem~\ref{xkform}. The above diagram is
cofibrant, because the upper arrow is a cofibration
in~$\cat{cdga}$, and its colimit is
$(\Lambda[m]\otimes\Q[\sK],d\bigr)$. So the result follows
from~\eqref{aplzkkoszul}.
\end{proof}

Theorem~\ref{zkmodel2} chimes with the rational models of
fibrations from \cite[\S15(c)]{f-h-t01}.

\begin{remark}
We may summarise the results of Theorems \ref{zkmodel1} and
\ref{zkmodel2} as follows. As functors $\cat{top}\to\cat{cdga}$,
both rational cohomology and $\APL$ map homotopy colimits to
homotopy limits on diagrams $\mathcal D_\sK(\pt,S^1)$, $I\mapsto
T^m/T^I$, and map homotopy limits to homotopy colimits on
diagrams~\eqref{zkpul}.
\end{remark}

\subsection*{Models for toric and quasitoric manifolds}
Let $M=M(P,\varLambda)$ be a quasitoric manifold\label{qtmanmode}
over a simple $n$-polytope $P$ with characteristic map
$\varLambda\colon\Z^m\to\Z^n$. Here we denote by $\sK$ the dual
sphere triangulation $\sK_P$ of~$P$. We recall from
Proposition~\ref{kact} that $M$ can be identified with the
quotient of the moment-angle manifold $\zp=\zk$ by the free action
of the $(m-n)$-torus $K(\varLambda)=\Ker(\varLambda\colon T^m\to
T^n)$; we denote the map of tori defined by $\varLambda$ by the
same letter for simplicity. All results below are equally
applicable to toric manifolds $M$, in which case $\sK$ is the
underlying complex of the corresponding complete regular
simplicial fan.

Our topological and algebraic models of $M$ are obtained from
those of $\zk$ by the appropriate factorisation by the action
of~$K(\varLambda)$.

We consider the $\ca{(\sK)}$-diagram $\mathcal
D_\sK(\pt,S^1)/K(\varLambda)$ obtained by factorisation
of~\eqref{dkpts}; its value on $I\subset J$ is the quotient map of
tori
\[
  T^n/\varLambda(T^I)\to T^n/\varLambda(T^J),
\]
where we identified $T^m/K(\varLambda)$ with $T^n$. This diagram
is not cofibrant; we denote its homotopy colimit by $\hocolim_I
T^n/\varLambda(T^I)$.

The following result was first proved by Welker, Ziegler and \v
Zivaljevi\'c in~\cite{w-z-z99}. It appears to be the earliest
mention of homotopy colimits in the toric context, and refers to a
diagram that is clearly not Reedy cofibrant:

\begin{proposition}\label{qthocolim}
There is a weak equivalence $M\simeq\hocolim_I
T^n/\varLambda(T^I)$.
\end{proposition}
\begin{proof}
As in the proof of Proposition~\ref{zkhocolim}, we consider
objectwise projections $(D^2,S^1)^I\to T^m/T^I=(S^1)^{[m]\setminus
I}$, and take quotients by the action of $K(\varLambda)$. As a
result, we obtain a weak equivalence
\[
  M=\zk/K(\varLambda)=\colim_I\bigl((D^2,S^1)^I/K(\varLambda)\bigr)
  \stackrel{\simeq}\longrightarrow\hocolim_I
  T^n/\varLambda(T^I).\qedhere
\]
\end{proof}

Next we construct an analogue of the algebraic
model~\eqref{aplzkkoszul} for quasitoric manifolds. For this we
consider the elements
\[
  t_i=\lambda_{i1}v_1+\cdots+\lambda_{im}v_m,\quad 1\le i\le n,
\]
in the face ring $\Q[\sK]=\Q[v_1,\ldots,v_m]/\mathcal I_\sK$
corresponding to the rows of matrix~$\varLambda$.

\begin{lemma}\label{qtmodel}
For a toric or quasitoric manifold $M$, the algebra $\APL(M)$ is
weakly equivalent to the commutative dg-algebra
\[
  \bigl(\Lambda[x_1,\ldots,x_n]\otimes\Q[\sK],d\bigr),\quad
  \text{with}\quad dx_i=t_i,\; dv_i=0.
\]
\end{lemma}
\begin{proof}
The argument is similar to that for Theorem~\ref{zkmodel1}. We
consider a $\cat{cat}^{op}(\sK)$-diagram
$\Lambda[n]\otimes\Q[\,\cdot\:]^\sK$, whose value on $I\subset J$
is the quotient map
\[
  \bigl(\Lambda[x_1,\ldots,x_n]\otimes\Q[v_i\colon i\in J],d\bigr)
  \to \bigl(\Lambda[x_1,\ldots,x_n]\otimes\Q[v_i\colon i\in I],d\bigr),
\]
where $d$ is defined on $\Lambda[n]\otimes\Q[v_i\colon i\in I]$ by
$dx_i=t_i$ for $i\in I$ and $dx_i=0$ otherwise. Then define
quasi-isomorphisms
\[
  \bigl(\Lambda[x_1,\ldots,x_n]\otimes\Q[v_i\colon i\in I],d\bigr)
  \stackrel\simeq\longrightarrow
  \APL\bigl((D^2,S^1)^I/K(\varLambda)\bigr)
\]
by sending  each $x_i$ to the $K(\varLambda)$-invariant 1-form
$\lambda_{i1}r_1^2d\varphi_1+\cdots+\lambda_{im}r_m^2d\varphi_m$,
where $(r_i,\varphi_i)$ are polar coordinates on the $i$th disk or
circle. These quasi-isomorphisms are compatible with the maps
corresponding to inclusions of simplices $I\subset J$ and
therefore provide a weak equivalence of Reedy fibrant diagrams in
$\cat{cdga}$. Their limits are therefore quasi-isomorphic, and we
obtain the required zigzag
\begin{multline*}
  \APL(M)=\APL\bigl((D^2,S^1)^\sK/K(\varLambda)\bigr)
  \stackrel\simeq\longrightarrow
  \lim_I\APL\bigr((D^2,S^1)^I/K(\varLambda)\bigl)\\
  \stackrel{\simeq}{\longleftarrow}
  \lim_I\bigl(\Lambda[x_1,\ldots,x_n]\otimes\Q[v_i\colon i\in
  I],d\bigr)=
  \bigl(\Lambda[x_1,\ldots,x_n]\otimes\Q[\sK],d\bigr).\qedhere
\end{multline*}
\end{proof}

\begin{remark}
The (quasi)toric manifold $M=M(P,\varLambda)$ is the homotopy
fibre of the composition of the inclusion $(\C P^\infty)^\sK\to
BT^m$ and the map $B\varLambda\colon BT^m\to BT^n$ (an exercise).
In other words, there is a homotopy pullback diagram
\begin{equation}\label{qtfibre}
  \xymatrix{
  M \ar[rr]\ar[d] & & ET^n \ar[d]\\
  (\C P^\infty)^\sK \ar[r]^i & BT^m \ar[r]^{B\varLambda} & BT^n
  }
\end{equation}
The model of Lemma~\ref{qtmodel} can be obtained by applying the
results of~\cite[\S15(c)]{f-h-t01} to the fibration $M\to(\C
P^\infty)^\sK$ above.
\end{remark}

Now we can prove the main result of this subsection:

\begin{theorem}[\cite{pa-ra08}]\label{formtvman}
Every toric or quasitoric manifold is formal.
\end{theorem}
\begin{proof}
We use the model of Lemma~\ref{qtmodel} and utilise the fact that
$\Q[\sK]$ is Cohen--Macaulay (Corollary~\ref{spherecm}), i.e.
$\Q[\sK]$ is free as module over $\Q[t_1,\ldots,t_n]$. Hence
$\otimes_{\Q[t_1,\ldots,t_n]}\Q[\sK]$ is a right exact functor,
and applying it to the quasi-isomorphism
$(\Lambda[u_1,\ldots,u_n]\otimes\Q[t_1,\ldots,t_n],d)\to\Q$ yields
a quasi-isomorphism
\[
  (\Lambda[u_1,\ldots,u_n]\otimes\Q[\sK],d)\to\Q[\sK]/(t_1,\ldots,t_n),
\]
which is given by the projection onto the second factor. Since
$\Q[\sK]/(t_1,\ldots,t_n)\cong H^*(M)$ by Theorem~\ref{qtcoh}, the
result follows from Lemma~\ref{qtmodel}.
\end{proof}

Similar arguments apply more generally to torus manifolds over
homology polytopes (see Section~\ref{torusman}), and even to
arbitrary torus manifolds with zero odd dimensional cohomology. In
the latter case, $\Q[\sK]$ is replaced by the face ring
$\Q[\mathcal S]$ of an appropriate simplicial poset $\mathcal S$
(see exercises below).

Note also that the formality of projective toric manifolds follows
immediately from the fact that they are K\"ahler.

\subsection*{Exercises}

\begin{exercise}
Construct the homotopy pullback diagram~\eqref{qtfibre}.
\end{exercise}

\begin{exercise}
Extend the argument of Theorem~\ref{xkform} to show that the
polyhedral power $(\C P^\infty,\pt)^\sS$ (see
Construction~\ref{macsimpos}) is formal for any simplicial
poset~$\sS$.
\end{exercise}

\begin{exercise}
Show that a torus manifold $M$ with $H^{odd}(M;\Z)=0$ is formal.
(Hint use Theorem~\ref{theo:stcoh} establishing an isomorphism
$H^*(M)\cong\Q[\sS]$, where $\Q[\sS]$ is the face ring of the
simplicial poset dual to the quotient $M/T=Q$, and use
Lemma~\ref{theo:face CM} to show that $\Q[\sS]$ is free over
$\Q[t_1,\ldots,t_n]$.)
\end{exercise}

\section{Wedges of spheres and connected sums of sphere products}\label{homtypes}

There are two situations when the homotopy type of the
moment-angle complex $\zk$ can be described explicitly. The first
concerns a family of examples of polytopal sphere triangulations
$\sK$ for which $\zk$ is homeomorphic to a connected sum of sphere
products, with two spheres in each product. The proofs use
differential topology and surgery theory; we included sample
results as Theorem~\ref{zpstacked} and Theorem~\ref{zp3quad} and
refer to a more detailed account in the work of Gitler and L\'opez
de Medrano~\cite{gi-lo13}. The second situation is when $\zk$ is
homotopy equivalent to a wedge of spheres; the corresponding
families of examples of $\sK$ were constructed by Grbi\'c and
Theriault in~\cite{gr-th04},~\cite{gr-th07} and are reviewed here.

\begin{theorem}[\cite{gr-th07}]\label{zkglue}
Let $\sK=\sK_1\cup_I\sK_2$ be a simplicial complex obtained by
gluing $\sK_1$ and $\sK_2$ along a common face, which may be
empty. If $\mathcal Z_{\sK_1}$ and $\mathcal Z_{\sK_2}$ are
homotopy equivalent to wedges of spheres, then $\zk$ is also
homotopy equivalent to a wedge of spheres.
\end{theorem}

We reproduce the proof from~\cite{gr-th07}, which uses two
lemmata.

\begin{lemma}[Cube Lemma]\label{Lcube}
Suppose there is a homotopy commutative diagram of spaces
\[
\diagram
  E\rrto\drto\ddto & & F\dline\drto & \\
  & G\rrto\ddto & \dto & H\ddto \\
  A\rline\drto & \rto & B\drto & \\
  & C\rrto & & D.
\enddiagram
\]
Suppose the bottom face is a homotopy pushout and the four
sides are homotopy pullbacks. Then the top face is a homotopy
pushout.
\end{lemma}
\begin{proof}
See~\cite[Theorem~25]{math76}. Note that the statement is fairly
straightforward in the case when $D=B\cup_A C$ and the vertical
arrows are locally trivial fibre bundles obtained by pulling back
the bundle $H\to D$ along the arrows of the bottom face. This will
be enough for our purposes.
\end{proof}

The \emph{join}\label{joinspace} of spaces $A$, $B$ is defined as
the identification space
\[
  A\mathop{*}B=A\times B\times \mathbb
  I\,/\;(a,b_1,0)\sim(a,b_2,0),\;(a_1,b,1)\sim(a_2,b,1).
\]
The product $A\times B$ embeds into the join as $A\times
B\times\frac12$. Furthermore, the lower half
\[
  A\mathop{*}B_{\le\frac12}=\{(a,b,t)\in A\mathop{*}B\colon t\le\textstyle\frac12\},
\]
of the join is the mapping cylinder of the first projection
$\pi_1\colon A\times B\to A$, and the upper half
$A\mathop{*}B_{\ge\frac12}$ is the mapping cylinder of the second
projection $\pi_2\colon A\times B\to B$. It follows that there is
a homotopy pushout diagram
\begin{equation}\label{diagjoin}
\diagram
  A\times B\rto^-{\pi_2}\dto_{\pi_{1}} & B\dto \\
  A\rto & A\mathop{*}B,
\enddiagram
\end{equation}
which can be viewed as the homotopy-theoretic definition of the
join.

For pointed spaces $A$, $B$, there are canonical homotopy
equivalences $\varSigma A\wedge B\simeq\varSigma(A\wedge B)\simeq
A\mathop{*}B$, where $A\wedge B=(A\times B)/(A\times
\pt\cup\pt\times B)$ is the smash product\label{halfsmash}. The
\emph{left half-smash product} is $A\ltimes B=A\times B/(A\times
\pt)$ and the \emph{right half-smash product} is $A\rtimes
B=A\times B/(\pt\times B)$. Let $\epsilon_A$ denote the map
collapsing $A$ to a point.

\begin{lemma}\label{gluingpo}
Let $A$, $B$, $C$ and $D$ be spaces. Define $Q$ as the homotopy
pushout
\[
  \diagram
  A\times B\rrto^-{\epsilon_A\times\id_B}\dto_{\id_A\times\epsilon_B}
  && C\times B\dto \\
  A\times D\rrto && Q.
  \enddiagram
\]
Then $Q\simeq (A\ast B)\vee (C\rtimes B)\vee (A\ltimes D)$.
\end{lemma}
\begin{proof}
We can decompose the pushout square above as
\begin{equation}\label{podeco}
  \diagram
  A\times B\rto^{\pi_2}\dto_{\pi_1}& B\rto^{i_2}\dto & C\times B\dto\\
  A\rto\dto_{i_1} & A\mathop{*}B \rto^{j_2}\dto_{j_1} & E\dto\\
  A\times D \rto & F\rto & Q,
  \enddiagram
\end{equation}
where $i_1$ and $i_2$ denote the inclusions into the first and
second factor, and each small square is a homotopy pushout.

Since the map $A\to A\mathop{*}B$ is null homotopic, we can pinch
out $A$ in the left bottom square above to obtain a homotopy
pushout
\[
  \diagram
  \pt \rto\dto & A\mathop{*}B\dto_{j_1}\\
  A\ltimes D \rto & F.
  \enddiagram
\]
Hence $F\simeq(A\mathop{*}B)\vee(A\ltimes D)$ and $j_1$ is
homotopic to the inclusion into the first wedge summand.
Similarly, $E\simeq(A\mathop{*}B)\vee(C\rtimes B)$. The required
decomposition of $Q$ now follows by considering the right bottom
square in~\eqref{podeco}.
\end{proof}

\begin{proof}[Proof of Theorem~\ref{zkglue}]
Recall from Theorem~\ref{zkhofib} that $\zk$ is the homotopy fibre
of the canonical inclusion $i\colon(\C P^\infty)^\sK\to(\C
P^\infty)^m=BT^m$. Assume that $\sK_1$ is a simplicial complex
on~$[m_1]$, $\sK_2$ is a simplicial complex on~$[m_2]$ and $\sK$
is a simplicial complex on~$[m]$, so that $m=m_1+m_2-|I|$. By some
abuse of notation we assume that the set $[m_1]$ is included as
the first $m_1$ elements in $[m]$, and $[m_2]$ is included as the
last $m_2$ elements. This defines the inclusions like $(\C
P^\infty)^{m_1}\to(\C P^\infty)^m$, $(\C P^\infty)^{\sK_2}\to(\C
P^\infty)^{\sK}$, etc. We have a pushout square
\[
  \diagram
  (\C P^\infty)^I \rto\dto & (\C P^\infty)^{\sK_1}\dto\\
  (\C P^\infty)^{\sK_2} \rto & (\C P^\infty)^\sK.
  \enddiagram
\]
Map each of the four corners of this pushout into $(\C
P^\infty)^m$ and take homotopy fibres. This gives homotopy
fibrations
\[
\begin{aligned}
  &T^{m-m_2}\times T^{m-m_1}=T^{m-|I|}\to(\C P^\infty)^{I}\to(\C
  P^\infty)^m,\\
  &\mathcal Z_{\sK_1}\times T^{m-m_1}\to(\C P^\infty)^{\sK_1}\to(\C P^\infty)^m,\\
  &T^{m-m_2}\times \mathcal Z_{\sK_2}\to(\C P^\infty)^{\sK_2}\to(\C P^\infty)^m,\\
  &\zk\to(\C P^\infty)^\sK\to(\C P^\infty)^m.\\
\end{aligned}
\]
%where $\mathcal Z_{\sK_1}$ is the homotopy fibre of $(\C
%P^\infty)^{\sK_1}\to(\C P^\infty)^{m_1}$ and $\mathcal Z_{\sK_2}$
%is the homotopy fibre of $(\C P^\infty)^{\sK_2}\to(\C
%P^\infty)^{m_2}$.

Including $(\C P^\infty)^I$ into $(\C P^\infty)^{\sK_1}$ gives a
homotopy pullback diagram
\[
\diagram
  \varOmega BT^m\rto\ddouble & T^{m-|I|}\rto\dto^{\theta}
  & (\C P^\infty)^I\rto\dto & BT^m\ddouble \\
  \varOmega BT^m\rto & \mathcal Z_{\sK_1}\times T^{m-m_1}\rto
  & (\C P^\infty)^{\sK_1}\rto & BT^m
\enddiagram
\]
for some map $\theta$ of fibres. We now identify $\theta$. With
$BT^{m}=\prod_{i=1}^{m}\mathbb{C}P^{\infty}$, the pullback just
described is the product of the homotopy pullback
\[
\diagram
  \varOmega BT^{m_1}\rto\ddouble & T^{m-m_2}\rto\dto^{\theta^{\prime}}
  & (\C P^\infty)^I\rto\dto & BT^{m_1}\ddouble \\
  \varOmega BT^{m_1}\rto & \mathcal Z_{\sK_1}\rto
  & (\C P^\infty)^{\sK_1}\rto & BT^{m_1}
\enddiagram
\]
and the path-loop fibration $T^{m-m_1}\to\pt\to BT^{m-m_1}$. So
$\theta=\theta^{\prime}\times\id_{T^{m-m_1}}$. Further,
$T^{m-m_2}$ is a retract of $\varOmega BT^{m_1}\simeq T^{m_1}$ and
$\varOmega BT^{m_1}\to\mathcal Z_{\sK_1}$ is null homotopic since
$\varOmega BT^{m_1}$ is a retract of $\varOmega((\C
P^\infty)^{\sK_1})$. Hence
$\theta^{\prime}\simeq\epsilon_{T^{m-m_2}}$ and therefore
$\theta\simeq\epsilon_{T^{m-m_2}}\times\id_{T^{m-m_1}}$. A similar
argument for the inclusion of $(\C P^\infty)^I$ into $(\C
P^\infty)^{\sK_2}$ shows that the map of fibres $T^{m-m_2}\times
T^{m-m_1}\to T^{m-m_2}\times\mathcal Z_{\sK_2}$ is homotopic to
$\id_{T^{m-m_2}}\times\epsilon_{T^{m-m_1}}$.

Collecting all this information about homotopy fibres,
Lemma~\ref{Lcube} shows that there is a homotopy pushout
\[
\diagram
  T^{m-m_2}\times T^{m-m_1}\rto^-{\epsilon\times\,\id}
  \dto_{\id\times\epsilon}
  & \mathcal Z_{\sK_1}\times T^{m-m_1}\dto \\
  T^{m-m_2}\times \mathcal Z_{\sK_2}\rto & \zk.
\enddiagram
\]
Lemma~\ref{gluingpo} then gives a homotopy decomposition
\begin{equation}\label{zkwedgedec}
  \zk\simeq(T^{m-m_2}\mathop{*}T^{m-m_1})\vee
  (\mathcal Z_{\sK_1}\rtimes T^{m-m_1})\vee(T^{m-m_2}\ltimes\mathcal Z_{\sK_2}).
\end{equation}

To show $\zk$ is homotopy equivalent to a wedge of spheres, we
show that each of $T^{m-m_2}\mathop{*}T^{m-m_1}$, $\mathcal
Z_{\sK_1}\rtimes T^{m-m_1}$ and $T^{m-m_2}\ltimes\mathcal
Z_{\sK_2}$ is homotopy equivalent to a wedge of spheres. First,
observe that the suspension of a product of spheres is homotopy
equivalent to a wedge of spheres, so
$T^{m-m_2}\mathop{*}T^{m-m_1}$ is homotopy equivalent to a wedge
of spheres. Second, as $\mathcal Z_{\sK_1}$ is homotopy equivalent
to a wedge of spheres we can write $\mathcal
Z_{\sK_1}\simeq\varSigma W$, where $W$ is a wedge of spheres (note
that $\mathcal Z_{\sK_1}$ is 2-connected by
Proposition~\ref{homgr}~(a)). We then have $\mathcal
Z_{\sK_1}\rtimes T^{m-m_1}\simeq\varSigma W\rtimes
T^{m-m_1}\simeq\varSigma W\vee(\varSigma T^{m-m_1}\wedge W)$. Now
$\varSigma T^{m-m_1}$ is homotopy equivalent to a wedge of
spheres. Therefore, as $W$ is homotopy equivalent to a wedge of
spheres so is $\varSigma T^{m-m_1}\wedge W$. Hence $\mathcal
Z_{\sK_1}\rtimes T^{m-m_1}$ is homotopy equivalent to a wedge of
spheres. The decomposition of the summand
$T^{m-m_2}\ltimes\mathcal Z_{\sK_2}$ into a wedge of spheres is
exactly as for $\mathcal Z_{\sK_1}\rtimes T^{m-m_1}$.
\end{proof}

\begin{corollary}\label{orderws}
Assume that there is an order $I_1,\ldots,I_s$ of the maximal
faces of~$\sK$ such that $\bigr(\bigcup_{j<k}I_j\bigl)\cap I_k$ is
a single face for each $k=1,\ldots,s$. Then $\zk$ has homotopy
type of a wedge of spheres.
\end{corollary}

As an application, we describe the homotopy type of $\zk$ for two
particular series of~$\sK$: 0-dimensional complexes and trees
(connected graphs without cycles).

\begin{proposition}\label{zkmpoints}
Let $\sK$ be $m$ disjoint points, $m\ge2$. Then
\begin{equation}\label{zkmwedge}
  \mathcal Z_{\sK}\simeq
  \bigvee_{k=2}^m\bigl(S^{k+1}\bigr)^{\vee(k-1)\bin mk}.
\end{equation}
\end{proposition}

\begin{proof}
Let $\sK_m$ denote the complex consisting of $m$ disjoint points.
Applying~\eqref{zkwedgedec} to the decomposition
$\sK_m=\sK_{m-1}\sqcup\sK_1$ we obtain
\begin{equation}\label{zkmdec}
  \mathcal Z_{\sK_m}\simeq(T^{m-1}\mathop{*}T^1)\vee
  (\mathcal Z_{\sK_{m-1}}\rtimes T^1)
\end{equation}
(the third wedge summand vanishes because $\mathcal
Z_{\sK_1}\simeq\pt$). An inductive argument using the
decomposition $\varSigma(A\times B)\simeq\varSigma A\vee\varSigma
B\vee(\varSigma A\wedge B)$ shows that
\[
  T^{m-1}\mathop{*}T^1\simeq\varSigma\varSigma T^{m-1}\simeq
  \bigvee_{k=2}^m\bigl(S^{k+1}\bigr)^{\vee\bin {m-1}{k-1}}.
\]
Assuming by induction that~\eqref{zkmwedge} holds for $\sK_{m-1}$,
we obtain
\[
  \mathcal Z_{\sK_{m-1}}\rtimes T^1\simeq\mathcal
  Z_{\sK_{m-1}}\vee\varSigma\mathcal Z_{\sK_{m-1}}\simeq
  \bigvee_{k=2}^{m-1}\bigl(S^{k+1}\bigr)^{\vee(k-1)\bin{m-1}k}\vee
  \bigvee_{k=3}^m\bigl(S^{k+1}\bigr)^{\vee(k-2)\bin{m-1}{k-1}}.
\]
Substituting the last two formulae into~\eqref{zkmdec} we finally
obtain~\eqref{zkmwedge}.
\end{proof}

Observe that $\zk$ corresponding to $m$ disjoint points is the
homotopy fibre of the inclusion of the $m$-fold wedge $(\C
P^\infty)^{\vee m}$ into the $m$-fold product $(\C P^\infty)^m$.
As we have seen in Example~\ref{codim2ar}, $\mathcal Z_{\sK_m}$ is
homotopy equivalent to the complement of the union of all
coordinate planes of codimension two in~$\C^m$. In this context
the result of Proposition~\ref{zkmpoints} was obtained
in~\cite{gr-th04}.

The homotopy type of $\zk$ corresponding to a tree with $m+1$
vertices depends only on the number of vertices and does not
depend of the form of the tree; also, the homotopy type of $\zk$
corresponding to a tree with $m+1$ vertices is the same as that of
$\zk$ corresponding to $m$ disjoint points:

\begin{proposition}\label{zktree}
Let $\sK$ be a tree with $m+1$ vertices, $m\ge2$. Then
\[
  \mathcal Z_{\sK}\simeq
  \bigvee_{k=2}^m\bigl(S^{k+1}\bigr)^{\vee(k-1)\bin mk}.
\]
\end{proposition}
\begin{proof}
This time we use the decomposition $\sK_m=\sK_{m-1}\cup_v\sK_1$,
where $\sK_m$ denotes a tree with $m+1$ vertices (so that $\sK_1$
is a segment), and the union in taken along a common vertex~$v$.
The rest of the proof is as for Proposition~\ref{zkmpoints}.
\end{proof}

One can notice the similarity between the wedge decomposition of
Proposition~\ref{zktree} and the connected sum decomposition of
Theorem~\ref{zpstacked}. The nature of this similarity is
explained in the work of Theriault~\cite{ther}.

A simplicial complex $\sK$ is \emph{shifted} if there is an
ordering of its vertices such that whenever $I\in\sK$, $i\in I$
and $i<j$, then $(I\setminus\{i\})\cup\{j\}\in\sK$.

\begin{theorem}[{\cite[Theorem~9.4]{gr-th07}}]\label{zkshifted}
If $\sK$ is a shifted complex, then $\zk$ is homotopy equivalent
to a wedge of spheres.
\end{theorem}
\begin{proof}[Idea of proof]
For any simplicial complex $\sK$ on $[m]$ there is a pushout
square
\[
  \diagram
  \lk_{\{m\}}\sK \rto\dto & \sK_{\{1,\ldots,m-1\}}\dto\\
  \st_{\{m\}}\sK \rto & \sK,
  \enddiagram
\]
where $\sK_{\{1,\ldots,m-1\}}$ denotes the restriction of $\sK$ to
the first $m-1$ vertices. It gives rise to a pushout square of the
corresponding polyhedral products $(\C P^\infty)^\sK$ and, by
application of Lemma~\ref{Lcube}, to a pushout square of the
moment-angle complexes~$\zk$. The key observation is that if $\sK$
is shifted,  then all three subcomplexes $\lk_{\{m\}}\sK$,
$\st_{\{m\}}\sK$ and $\sK_{\{1,\ldots,m-1\}}$ are also shifted,
with respect to the induced ordering of vertices. This allows us
to use induction, in a way similar to the proof of
Theorem~\ref{zkglue}.
\end{proof}

The $i$-dimensional skeleton of a simplex $\varDelta^{m-1}$ is a
shifted complex. In this case the dimensions of spheres in the
wedge decomposition of $\zk$ can be described explicitly; this
result was given as Theorem~\ref{gr-th}. It also follows from a
result of Porter~\cite{port66} for general polyhedral products,
which we state in the next section.

Not all complexes $\sK$ obtained by iterative gluing along a
common face are shifted (see Exercise~\ref{2trees}), and not all
shifted complexes can be obtained by iterative gluing along a
common face. So one can obtain even a wider class of simplicial
complexes $\sK$ whose corresponding $\zk$ are wedges of spheres by
combining the results of Theorem~\ref{zkglue} and
Theorem~\ref{zkshifted}.

%Let $\mathrm{WS}$ denote the class of simplicial complexes $\sK$
%whose corresponding moment-angle complex $\zk$ is homotopy
%equivalent to a wedge of spheres, and let $\mathrm{CS}$ denote the
%class of simple polytopes $P$ for which the corresponding
%moment-angle manifold $\zp$ is diffeomorphic to a connected sum of
%products $S^k\times S^{m-k}$ ($m$ is fixed, $k$ may vary). We have
%the following results about the classes $\mathrm{WS}$ and
%$\mathrm{CS}$:
%
%\begin{itemize}
%\item[$\cdot$] If $\sK$ is a tree, then $\sK\in\mathrm{WS}$ (Proposition~\ref{zktree});
%
%\item[$\cdot$] If $\sK$ is a shifted complex (in particular, if $\sK$ is the $i$-dimensional skeleton skeleton of
%a simplex $\varDelta^{m-1}$ for some~$i$), then
%$\sK\in\mathrm{WS}$ (Theorem~\ref{zkshifted});
%
%\item[$\cdot$] If $\sK=\sK_1\cup_I\sK_2$ where
%$\sK_1,\sK_2\in\mathrm{WS}$, then $\sK\in\mathrm{WS}$
%(Theorem~\ref{zkglue});
%
%\item[$\cdot$] If $P$ is obtained from a simplex by iterative
%vertex truncations
%\end{itemize}

\subsection*{Exercises}

\begin{exercise}\label{2trees}
The tree $\quad\begin{picture}(10,2) \put(0,1){\circle*{1}}
\put(5,1){\circle*{1}} \put(10,1){\circle*{1}}
\put(0,1){\line(1,0){10}}
\end{picture}\quad$
is a shifted complex, but $\quad\begin{picture}(15,2)
\put(0,1){\circle*{1}} \put(5,1){\circle*{1}}
\put(10,1){\circle*{1}} \put(15,1){\circle*{1}}
\put(0,1){\line(1,0){15}}
\end{picture}\quad$ is not.
\end{exercise}

\begin{exercise}
Let $\sK$ be the graph $\quad\begin{picture}(10,5)
\put(0,2.5){\circle*{1}} \put(5,0){\circle*{1}}
\put(5,5){\circle*{1}} \put(10,2.5){\circle*{1}}
\put(5,0){\line(0,1){5}} \put(0,2.5){\line(2,1){5}}
\put(0,2.5){\line(2,-1){5}} \put(10,2.5){\line(-2,1){5}}
\put(10,2.5){\line(-2,-1){5}}
\end{picture}\quad$. Describe the homotopy type of~$\zk$.

\begin{exercise}
Let $\sK$ be a complex obtained by iteration of the operation of
attaching a $k$-simplex along a common $(k-1)$-face, starting from
a $k$-simplex (so $\sK$ is a tree when $k=1$). Describe the
homotopy type of~$\zk$.
\end{exercise}
\end{exercise}

\section{Stable decompositions of polyhedral products}\label{stabledec}
Several important results on wedge decomposition of polyhedral
products after one suspension were obtained in the work of Bahri,
Bendersky, Cohen and Gitler~\cite{b-b-c-g10}. These can be seen as
far-reaching generalisations of the classical decomposition
$\varSigma(A\times B)\simeq\varSigma A\vee\varSigma
B\vee(\varSigma A\wedge B)$. The proofs given below are reproduced
from~\cite{b-b-c-g10} with few or no modifications. Homotopy
theory of polyhedral products has become quite an active area, and
we also review some recent results on stable and unstable
decompositions in the end of this section.

We start by defining the smash version of the polyhedral product.

\begin{construction}[polyhedral smash product]\label{phsp}
The initial setup is again a simplicial complex $\sK$ on~$[m]$ and
a sequence of $m$ pairs of pointed cell complexes
\[
  (\mb X,\mb A)=\{(X_1,A_1),\ldots,(X_m,A_m)\}.
\]
We denote the $m$-fold smash product of the $X_i$ by
\[
  \mb X^{\wedge m}=X_1\wedge X_2\wedge\cdots\wedge X_m.
\]
Then the \emph{polyhedral smash product} $(\mb X,\mb A)^{\wedge
\sK}$ is defined as the image of $(\mb X,\mb A)^\sK$ under the
projection $\mb X^m\to\mb X^{\wedge m}$. More specifically, for
each $I\subset[m]$ we set
\[
  (\mb X,\mb A)^{\wedge I}=\bigl\{(x_1,\ldots,x_m)\in
  X_1\wedge X_2\wedge\cdots\wedge X_m\colon\; x_j\in A_j\quad\text{for }j\notin
  I\bigl\},
\]
then
\[
  (\mb X,\mb A)^{\wedge\sK}=\bigcup_{I\in\mathcal K}(\mb X,\mb A)^{\wedge I}=
  \bigcup_{I\in\mathcal K}
  \Bigl(\bigwedge_{i\in I}X_i\wedge\bigwedge_{i\notin I}A_i\Bigl).
\]
Using the categorical language, define the $\ca{(\sK)}$-diagram
\begin{equation}\label{dhatdiag}
\begin{aligned}
  \widehat{\mathcal D}_\sK(\mb X,\mb A)\colon \ca(\sK)&\longrightarrow \top,\\
  I&\longmapsto (\mb X,\mb A)^{\wedge I},
\end{aligned}
\end{equation}
which maps the morphism $I\subset J$ to the inclusion $(\mb
X,\mb A)^{\wedge I}\subset(\mb X,\mb A)^{\wedge J}$. Then
\[
  (\mb X,\mb A)^{\wedge\sK}=\mathop{\mathrm{colim}}
  \widehat{\mathcal D}_\sK(\mb X,\mb A)=\mathop{\mathrm{colim}}_{I\in\sK}
  (\mb X,\mb A)^{\wedge I}.
\]

In the case when all the pairs $(X_i,A_i)$ are the same, i.e.
$X_i=X$ and $A_i=A$, we use the notation $(X,A)^{\wedge\sK}$ for
$(\mb X,\mb A)^{\wedge\sK}$. Also, if each $A_i=\pt$, then we use
the abbreviated notation $\mb X^{\wedge\sK}$ for $(\mb
X,\pt)^{\wedge\sK}$, and $X^{\wedge\sK}$ for
$(X,\pt)^{\wedge\sK}$.
\end{construction}

An inductive argument using the decomposition $\varSigma(A\times
B)\simeq\varSigma A\vee\varSigma B\vee(\varSigma A\wedge B)$ shows
that there is a natural pointed homotopy equivalence
\begin{equation}\label{suspprod}
  \varSigma(X_1\times\cdots\times\ X_m)
  \stackrel\simeq\longrightarrow
  \varSigma\Bigl(\bigvee_{J\subset[m]}\mb X^{\wedge J}\Bigr).
\end{equation}

For each $J\subset[m]$, define the subfamily
\[
  (\mb X_J,\mb A_J)=\{(X_j,A_j)\colon j\in J\}.
\]
The first result shows that the polyhedral product splits after
one suspension into a wedge of polyhedral smash products
corresponding to all full subcomplexes of~$\sK$:

\begin{theorem}[\cite{b-b-c-g10}]\label{wedgethm}
For any sequence $(\mb X,\mb A)$ of pairs of pointed cell
complexes, homotopy equivalence~\eqref{suspprod} induces a natural
pointed homotopy equivalence
\[
  \varSigma(\mb X,\mb A)^\sK\stackrel\simeq\longrightarrow
  \varSigma\Bigl(\bigvee_{J\subset[m]}(\mb X_J,\mb A_J)^{\wedge\mathcal K_J}\Bigl).
\]
\end{theorem}
\begin{proof}
We have $(\mb X,\mb A)^\sK=\mathop{\mathrm{colim}} \mathcal
D_\sK(\mb X,\mb A)=\mathop{\mathrm{colim}}_{I\in\sK} (\mb X,\mb
A)^I$, where $\mathcal D_\sK(\mb X,\mb A)$ is
diagram~\eqref{dxadiag}. Define another diagram
\[
\begin{aligned}
  \mathcal E_\sK(\mb X,\mb A)\colon \ca(\sK)&\longrightarrow \top,\\
  I&\longmapsto\bigvee_{J\subset[m]}(\mb X_J,\mb A_J)^{\wedge(I\cap J)}.
\end{aligned}
\]
By~\eqref{suspprod}, there is a natural pointed homotopy
equivalence
\[
  \varSigma(\mb X,\mb A)^I\stackrel\simeq\longrightarrow
  \varSigma\Bigl(\bigvee_{J\subset[m]}(\mb X_J,\mb A_J)^{\wedge(I\cap J)}\Bigr)
\]
The diagrams $\mathcal D_{\sK}(\mb X,\mb A)$, $\mathcal
E_{\sK}(\mb X,\mb A)$, as well as $\varSigma\mathcal D_{\sK}(\mb
X,\mb A)$, $\varSigma\mathcal E_{\sK}(\mb X,\mb A)$, are obviously
cofibrant, so the objectwise homotopy equivalence above induces a
homotopy equivalence of their colimits (see
Proposition~\ref{wecolim}). It remains to note that
\begin{align*}
  \colim\varSigma\mathcal D_{\sK}(\mb X,\mb A)&
  =\varSigma\colim\mathcal D_{\sK}(\mb X,\mb A)=\varSigma(\mb X,\mb
  A)^\sK,\\
  \colim\varSigma\mathcal E_{\sK}(\mb X,\mb A)&=
  \colim_{I\in\sK} \varSigma\Bigl(\bigvee_{J\subset[m]}(\mb
  X_J,\mb A_J)^{\wedge(I\cap J)}\Bigr)\\
  &=\varSigma\Bigl(\bigvee_{J\subset[m]}\colim_{(I\cap J)\in\sK_J} (\mb
  X_J,\mb A_J)^{\wedge(I\cap J)}\Bigr)\\
  &=\varSigma\Bigl(\bigvee_{J\subset[m]}(\mb
  X_J,\mb A_J)^{\wedge\sK_J}\Bigr).\qedhere
\end{align*}
\end{proof}

\vspace{-\baselineskip}

The homotopy type of the wedge summands $(\mb X_J,\mb
A_J)^{\wedge\mathcal K_J}$ can be described explicitly in the case
when the inclusions $A_k\hookrightarrow X_k$ are null-homotopic:

\begin{theorem}[\cite{b-b-c-g10}]\label{wedgenh}
Let $\sK$ be a simplicial complex on~$[m]$, and let $(\mb X,\mb
A)$ be a sequence of pairs of cell complexes with the property
that the inclusion $A_k\hookrightarrow X_k$ is null-homotopic for
all~$k$. Then there is a homotopy equivalence
\[
  (\mb X,\mb A)^{\wedge\mathcal K}
  \stackrel\simeq\longrightarrow
  \bigvee_{I\in\sK}|\lk_{\sK}I|\mathbin{*}(\mb X,\mb A)^{\wedge
  I},
\]
where $|\lk_{\sK}I|$ is the geometric realisation of the link of
$I$ in~$\sK$.
\end{theorem}
\begin{proof}
By hypothesis, there is a homotopy $F_k\colon A_k\times\I\to X_k$
such that $F_k(a,0)=i_k(a)$ and $F_k(a,1)=\pt$, where $i_k\colon
A_k\to X_k$ is the inclusion. By the homotopy extension property,
there exists $\widehat F_k\colon X_k\times\I\to X_k$ with
$\widehat F_k(x,0)=x$, $\widehat F_k(x,1)=g_k(x)$ where $g_k\colon
X_k\to X_k$ is a map such that $g_k(a)=\pt$ for all $a\in A_k$.
Hence there is a commutative diagram
\[
  \diagram
  A_k\rto^{\id}\dto_{i_k} & A_k\dto^\epsilon\\
  X_k\rto^{g_k} & X_k
  \enddiagram
\]
where $\epsilon\colon A_k\to X_k$ is the constant map to the
basepoint. Along with the diagram $\widehat{\mathcal
D}_\sK=\widehat{\mathcal D}_\sK(\mb X,\mb A)$ given
by~\eqref{dhatdiag}, define a new diagram
\[
  \widehat{\mathcal E}_\sK\colon \ca(\sK)\longrightarrow
  \top,\quad
  I\longmapsto (\mb X,\mb A)^{\wedge I},
\]
which maps the non-identity morphism $I\subset J$ to the constant
map $(\mb X,\mb A)^{\wedge I}\to(\mb X,\mb A)^{\wedge J}$ to the
basepoint. For every $I\in\sK$, define
\[
  \alpha(I)\colon\widehat{\mathcal D}_\sK(I)\to\widehat{\mathcal E}_\sK(I)
\]
by $\alpha(I)=\alpha_1(I)\wedge\cdots\wedge\alpha_m(I)$ where
\[
  \alpha_k(I)=\left\{%
  \begin{array}{ll}
    g_k\colon X_k\to X_k&\quad\text{if}\;\; k\in I, \\[2pt]
    \id\colon A_k\to A_k&\quad\text{if}\;\; k\notin I.\\
  \end{array}%
  \right.
\]
Since the $g_k$ are homotopy equivalences, so is $\alpha(I)$ for
all $I\in\sK$. Furthermore, if $I\subset J$, the following diagram
commutes:
\[
  \diagram
  \widehat{\mathcal D}_\sK(I)\rto^{\alpha(I)}\dto & \widehat{\mathcal E}_\sK(I)\dto\\
  \widehat{\mathcal D}_\sK(J)\rto^{\alpha(J)} & \widehat{\mathcal E}_\sK(J).
  \enddiagram
\]
Hence the maps $\alpha(I)$ give a weak equivalence of diagrams
$\widehat{\mathcal D}_\sK\to\widehat{\mathcal E}_\sK$, which gives
a homotopy equivalence
\[
  \hocolim\mathcal{\widehat D}_\sK\stackrel\simeq\longrightarrow
  \hocolim\mathcal{\widehat E}_\sK.
\]
Finally, the diagram $\mathcal{\widehat E}_\sK$ satisfies the
conditions of Lemma~\ref{wedgelemma}, so we get a homotopy
equivalence
\[
  \hocolim\mathcal{\widehat E}_\sK\stackrel\simeq\longrightarrow
  \bigvee_{I\in\sK}|\lk_{\sK}I|\mathbin{*}(\mb X,\mb A)^{\wedge
  I}
\]
(upper semi-intervals in the face poset of~$\sK$ are links). The
result follows since $\hocolim{\widehat D}_\sK=(\mb X,\mb
A)^{\wedge\mathcal K}$.
\end{proof}

Two special cases of Theorem~\ref{wedgenh} are presented next
where either $A_i$ are contractible for all~$i$ or $X_i$ are
contractible for all~$i$.

\begin{theorem}[\cite{b-b-c-g10}]\label{wedgethmA}
Let $\sK$ be a simplicial complex on~$[m]$, and let $(\mb X,\mb
A)$ be a sequence of pairs of cell complexes with the property
that all the $A_i$ are contractible. Then there is a homotopy
equivalence
\[
  \varSigma(\mb X,\mb A)^\sK\stackrel\simeq\longrightarrow
  \varSigma\Bigl(\bigvee_{I\in\sK}\mb X^{\wedge I}\Bigr).
\]
\end{theorem}
\begin{proof}
When all the $A_i$ are contractible, the space $(\mb X,\mb
A)^{\wedge I}$ is also contractible unless $I=[m]$. By
Theorem~\ref{wedgenh},
\[
  (\mb X_J,\mb A_J)^{\wedge\mathcal K_J}\simeq
  \bigvee_{I\in\sK_J}|\lk_{\sK_J}I|\mathbin{*}(\mb X_J,\mb A_J)^{\wedge
  I},
\]
which is contractible unless $J\in\sK_J$, i.e. $J\in\sK$. In the
latter case we have $(\mb X_J,\mb A_J)^{\wedge\mathcal K_J}=\mb
X^{\wedge J}$. By Theorem~\ref{wedgethm},
\[
  \varSigma(\mb X,\mb A)^\sK\simeq
  \varSigma\Bigl(\bigvee_{J\subset[m]}(\mb X_J,\mb A_J)^{\wedge\mathcal K_J}\Bigl)
  =\varSigma\Bigl(\bigvee_{J\in\sK}\mb X^{\wedge J}\Bigl).\qedhere
\]
\end{proof}

\begin{remark}
An interesting corollary of Theorem~\ref{wedgethmA} is that the
polyhedral products $X^\sK=(X,\pt)^\sK$ corresponding to
simplicial complexes $\sK$ with the same $f$-vectors become
homotopy equivalent after one suspension.
\end{remark}

\begin{theorem}[\cite{b-b-c-g10}]\label{wedgethmX}
Let $\sK$ be a simplicial complex on~$[m]$, and let $(\mb X,\mb
A)$ be a sequence of pairs of cell complexes with the property
that all the $X_i$ are contractible. Then there is a homotopy
equivalence
\[
  \varSigma(\mb X,\mb A)^\sK\stackrel\simeq\longrightarrow
  \varSigma\Bigl(\bigvee_{J\notin\sK}|\sK_J|\mathbin{*}\mb A^{\wedge J}\Bigr).
\]
\end{theorem}
\begin{proof}
Since all the $X_i$ are contractible, all of the spaces $(\mb
X,\mb A)^{\wedge I}$ are also contractible with the possible
exception of $(\mb X,\mb A)^{\wedge\varnothing}=\mb A^{\wedge m}$.
By Theorem~\ref{wedgenh},
\[
  (\mb X_J,\mb A_J)^{\wedge\mathcal K_J}\simeq
  \bigvee_{I\in\sK_J}|\lk_{\sK_J}I|\mathbin{*}(\mb X_J,\mb A_J)^{\wedge
  I}=|\lk_{\sK_J}\varnothing|\mathbin{*}(\mb X_J,\mb A_J)^{\wedge\varnothing}=
  |\sK_J|\mathbin{*}\mb A^{\wedge J}
\]
which is contractible if $J\in\sK$. By Theorem~\ref{wedgethm},
\[
  \varSigma(\mb X,\mb A)^\sK\simeq
  \varSigma\Bigl(\bigvee_{J\subset[m]}(\mb X_J,\mb A_J)^{\wedge\mathcal K_J}\Bigl)
  =\varSigma\Bigl(\bigvee_{J\notin\sK}|\sK_J|\mathbin{*}\mb A^{\wedge J}\Bigr).\qedhere
\]
\end{proof}

{\samepage
\begin{corollary}\
\begin{itemize}
\item[(a)] Let $(\mb X,\mb A)=(D^1,S^0)$, so that $(\mb X,\mb
A)^\sK$ is the real moment-angle complex~$\rk$. Then there is a
homotopy equivalence
\[
  \varSigma\rk\stackrel\simeq\longrightarrow
  \bigvee_{J\notin\sK}\varSigma^2|\sK_J|.
\]
\item[(b)] Let $(\mb X,\mb A)=(D^2,S^1)$, so that $(\mb X,\mb
A)^\sK$ is the moment-angle complex~$\zk$. Then there is a
homotopy equivalence
\[
  \varSigma\zk\stackrel\simeq\longrightarrow
  \bigvee_{J\notin\sK}\varSigma^{2+|J|}|\sK_J|.
\]
\end{itemize}
\end{corollary}
}

The above decomposition of $\varSigma\zk$ implies the additive
isomorphism
\[
  H^k(\zk;\Z)\cong\bigoplus_{J\subset[m]}\widetilde H^{k-|J|-1}(\sK_J;\Z)
\]
of Theorem~\ref{zkadd}. Similarly, the decomposition of
$\varSigma\rk$ implies the isomorphism
\[
  H^k(\rk;\Z)\cong\bigoplus_{J\subset[m]}\widetilde
  H^{k-1}(\sK_J;\Z).
\]

Another result of~\cite{b-b-c-g10} describes the cohomology ring
of a polyhedral product $\mb X^\sK$ and generalises the
isomorphism $H^*((\C P^\infty)^\sK;\Z)\cong\Z[\sK]$ of
Proposition~\ref{homsrs}:

\begin{theorem}[\cite{b-b-c-g10}]\label{xkcohomo}
Let $\mb X=(X_1,\ldots,X_m)$ be a sequence of pointed cell
complexes, and let $\k$ be a ring such that the natural map
\[
  H^*(X_{j_1};\k)\otimes\cdots\otimes H^*(X_{j_k};\k)\to
  H^*(X_{j_1}\times\cdots\times X_{j_k};\k)
\]
is an isomorphism for any $\{j_1,\ldots,j_k\}\subset[m]$. There is
an isomorphism of algebras
\[
  H^*(\mb X^{\sK};\k)\cong\bigl(H^*(X_1;\k)\otimes\cdots\otimes
  H^*(X_m;\k)\bigr)/\mathcal I,
\]
where $\mathcal I$ is the generalised Stanley--Reisner ideal,
generated by elements $x_{j_1}\otimes\cdots\otimes x_{j_k}$ for
which $x_{j_i}\in\widetilde H^*(X_{j_i};\k)$ and
$\{j_1,\ldots,j_k\}\notin\sK$.
%\[
%  \mathcal I=\bigl(x_{i_1}\!\cdots\, x_{i_k}\colon x_{i_j}\in\widetilde H^*(X_{i_j};\k),\quad
%  \{i_1,\ldots,i_k\}\notin\sK\bigr).
%\]
Furthermore, the inclusion $\mb X^\sK\to X_1\times\cdots\times
X_m$ induces the quotient projection in cohomology.
\end{theorem}
\begin{proof}
By Theorem~\ref{wedgethmA}, there are homotopy equivalences
\[
  \varSigma(X_1\times\cdots\times\ X_m)
  \stackrel\simeq\longrightarrow
  \varSigma\Bigl(\bigvee_{J\subset[m]}\mb X^{\wedge J}\Bigr),\quad
  \varSigma\mb X^\sK\stackrel\simeq\longrightarrow
  \varSigma\Bigl(\bigvee_{J\in\sK}\mb X^{\wedge J}\Bigr).
\]
Naturality implies that the map $\varSigma\mb
X^\sK\to\varSigma(X_1\times\cdots\times X_m)$ is split with
cofibre $\varSigma\bigl(\bigvee_{J\notin\sK}\mb X^{\wedge
J}\bigr)$. Hence there is a split cofibration
\[
  \varSigma\Bigl(\bigvee_{J\in\sK}\mb X^{\wedge J}\Bigr)\to
  \varSigma\Bigl(\bigvee_{J\subset[m]}\mb X^{\wedge J}\Bigr)\to
  \varSigma\Bigl(\bigvee_{J\notin\sK}\mb X^{\wedge J}\Bigr).
\]
Under the given condition on~$\k$ there is a ring isomorphism
\[
  \widetilde H^*(X_1\times\cdots\times X_m;\k)\cong
  \bigoplus_{J\subset[m]}
  \widetilde H^*(X_{j_1};\k)\otimes\cdots\otimes\widetilde
  H^*(X_{j_k};\k).
\]
The natural inclusion map $\varSigma\mb
X^\sK\to\varSigma(X_1\times\cdots\times X_m)$ induces a map
\[
  \widetilde H^*(X_1\times\cdots\times X_m;\k)\to
  \widetilde H^*(\mb X^\sK;\k),
\]
which corresponds to the projection map
\[
  \bigoplus_{J\subset[m]}
  \widetilde H^*(X_{j_1};\k)\otimes\cdots\otimes\widetilde
  H^*(X_{j_k};\k)\to\bigoplus_{J\in\sK}
  \widetilde H^*(X_{j_1};\k)\otimes\cdots\otimes\widetilde
  H^*(X_{j_k};\k).
\]
Its kernel is exactly
\[
  \bigoplus_{J\notin\sK}
  \widetilde H^*(X_{j_1};\k)\otimes\cdots\otimes\widetilde
  H^*(X_{j_k};\k)
\]
which is the generalised Stanley--Reisner ideal $\mathcal I$ by
inspection.
\end{proof}

As in the case of the face ring $\k[\sK]=H^*((\C
P^\infty)^\sK;\k)$, which can be decomposed as the limit of the
$\ca{(\sK)}^{op}$-diagram of polynomial algebras $\k[v_i\colon
i\in I]=H^*((\C P^\infty)^I;\k)$, the isomorphism of
Theorem~\ref{xkcohomo} can be interpreted as
\[
  H^*(\mb X^\sK;\k)\cong\lim_{I\in\sK}H^*(\mb X^I;\k).
\]
This isomorphism was used in the proof of Theorem~\ref{xkform}.

There are several important situations when the isomorphism of
Theorem~\ref{wedgethmX} can be desuspended. As we have seen in the
previous section, this is the case for the pairs $(\mb X,\mb
A)=(D^2,S^1)$ when $\sK$ is obtained by iterative gluing along a
common face or when $\sK$ is shifted complex. More generally, the
following result, conjectured in~\cite{b-b-c-g10}, was proved
independently by Grbi\'c--Theriault and Iriye--Kishimoto:

\begin{theorem}[{\cite[Theorem~1.1]{gr-th13},
\cite[Theorem~1.7]{ir-ki13}}]\label{grteirki} Let $\sK$ be a
shifted complex. Let $\mb A=(A_1,\ldots,A_m)$ be a sequence of
pointed cell complexes, and let $\cone A$ denote the cone on~$A$.
Then there is a homotopy equivalence
\[
  (\cone\mb A,\mb A)^\sK\stackrel\simeq\longrightarrow
  \bigvee_{J\notin\sK}|\sK_J|\mathbin{*}\mb A^{\wedge J}.
\]
\end{theorem}

Here is a case where the wedge summands can be described very
explicitly:

\begin{corollary}\label{xskele}
Let $\sK_i$ be the $i$-dimensional skeleton of the simplex
$\varDelta^{m-1}$. Then there is a homotopy equivalence
\[
  (\cone\mb A,\mb A)^{\sK_i}\stackrel\simeq\longrightarrow
  \bigvee_{k=i+2}^{m}\Bigl(\bigvee_{1\le j_1<\cdots<j_k\le m}
  \bigl(\varSigma^{i+1}A_{j_1}\wedge\cdots\wedge A_{j_k}
  \bigr)^{\vee\binom{k-1}{i+1}}\Bigr).
\]
\end{corollary}

The proof of this corollary is left as an exercise. In the case
when each $A_i$ is a loop space, $A_i=\varOmega B_i$, the space
$(\cone\varOmega\mb B,\varOmega\mb B)^{\sK_i}$ is the homotopy
fibre of the inclusion $\mb B^{\sK_i}\to\mb B^m$ (see
Example~\ref{exkpo}.5 and Exercise~\ref{plooppp}) and
decomposition above was obtained by Porter~\cite{port66} . In the
case $A_i=S^1$, Corollary~\ref{xskele} turns into
Theorem~\ref{gr-th}.

The wedge decomposition of Theorem~\ref{wedgethmX} can be used to
describe the ring structure for the cohomology of $(\mb X,\mb
A)^\sK$, see~\cite{b-b-c-g12}.

\subsection*{Exercises}
\begin{exercise}
Deduce Corollary~\ref{xskele} from Theorem~\ref{grteirki}.
\end{exercise}

\section{Loop spaces, Whitehead and Samelson products}\label{loops}
We now turn our attention to topological and algebraic models for
the loop spaces $\varOmega(\C P^\infty)^\sK$ and $\varOmega\zk$.
We can view the latter as objects in the category $\cat{tmon}$ of
topological monoids by considering Moore loops\label{mooreloop}
(of arbitrary length), whose composition is strictly associative.
%The properties of $\varOmega(\C P^\infty)^\sK$ and $\varOmega\zk$
%are considerably simplified when $\sK$ is a flag complex, but we
%postpone discussion of this situation until the following section.

\subsection*{Pontryagin algebras, Whitehead and Samelson products}
We loop the fibration $\zk\to(\C P^\infty)^\sK\to(\C P^\infty)^m$
to obtain a fibration
\begin{equation}\label{omegazkfib}
  \varOmega\zk\longrightarrow\varOmega(\C P^\infty)^\sK\longrightarrow T^m.
\end{equation}
It admits a section, defined by the $m$ generators of $\pi_2((\C
P^\infty)^\sK)\cong\Z^m$, and therefore splits in~$\cat{top}$. So
we have a homotopy equivalence
\[
  \varOmega(\C P^\infty)^\sK\stackrel{\simeq}{\longrightarrow}\varOmega\zk\times
  T^m
\]
which does \emph{not} preserve monoid structures.

\begin{proposition}
There is an exact sequence of homotopy Lie algebras
\[
  0\longrightarrow \pi_*(\varOmega\zk)\otimes\Q\longrightarrow
  \pi_*(\varOmega(\C P^\infty)^\sK)\otimes\Q\longrightarrow
  \cl(u_1,\ldots,u_m)\longrightarrow 0,
\]
where $\cl(u_1,\ldots,u_m)$ denotes the commutative Lie algebra
with generators $u_i$, $\deg u_i=1$, and an exact sequence of
Pontryagin algebras
\begin{equation}\label{paseq}
  0\longrightarrow H_*(\varOmega\zk;\k)\longrightarrow
  H_*\bigl(\varOmega(\C P^\infty)^\sK;\k\bigr)\longrightarrow
  \Lambda[u_1,\ldots,u_m]\longrightarrow 0,
\end{equation}
for any commutative ring $\k$ with unit.
\end{proposition}
\begin{proof}
The first exact sequence follows by considering the homotopy exact
sequence of the fibration~\eqref{omegazkfib}, whose connecting
homomorphism is zero because the fibration is trivial.

Since $H_*(T^m;\k)=\Lambda[u_1,\ldots,u_m]$ is a finitely
generated free $\k$-module, the K\"unneth formula gives an
isomorphism of $\k$-modules
\[
  H_*\bigl(\varOmega(\C P^\infty)^\sK;\k\bigr)\cong
  H_*(\varOmega\zk;\k)\otimes\Lambda[u_1,\ldots,u_m]
\]
and therefore an exact sequence of $\k$-algebras~\eqref{paseq}.
\end{proof}

The homotopy group $\pi_2((\C P^\infty)^\sK)\cong\Z^m$ has $m$
canonical generators represented by the maps
\[
  \widehat\mu_{i}\colon S^{2}\longrightarrow
  \C P^{\infty}\longrightarrow
  (\mathbb{C}P^{\infty})^{\vee m}\longrightarrow(\C P^\infty)^\sK
\]
for $1\le i\le m$, where the left map is the inclusion of the
bottom cell, the middle map is the inclusion of the $i$th wedge
summand, and the right map is the canonical inclusion of
polyhedral powers corresponding to the inclusion of the discrete
$m$-point complex into~$\sK$. Let
\[
  \mu_i\colon S^1\longrightarrow\varOmega(\C P^\infty)^\sK
\]
be the adjoint of $\widehat\mu_i$, and let $u_i$ denote the
Hurewicz image of $\mu_{i}$ in $H_1(\varOmega(\C P^\infty)^\sK)$.

We shall be interested in elements of $\pi_*(\varOmega(\C
P^\infty)^\sK)$ represented by Samelson products\label{samelprodu}
of the~$\mu_i$ (see Section~\ref{homothomol} for the definition).

\begin{proposition}
The Samelson products of the canonical generators
$\mu_i\in\pi_1(\varOmega(\C P^\infty)^\sK)$ satisfy the identities
\[
  [\mu_i,\mu_i]_s=0,\quad[\mu_i,\mu_j]_s=0\quad\text{if and only
  if }\quad\{i,j\}\in\sK.
\]
\end{proposition}
\begin{proof}
By adjunction, we can work with the Whitehead products instead.
The Whitehead square $[\widehat\mu_i,\widehat\mu_i]_w$ is zero in
$\pi_3((\C P^\infty)^\sK)$, because it is zero in $\pi_3(\C
P^\infty)=0$. Furthermore, the map
$\widehat\mu_i\vee\widehat\mu_j\colon S^2\vee S^2\to(\C
P^\infty)^\sK$ with $i\ne j$ extends to a map $S^2\times S^2\to(\C
P^\infty)^\sK$ whenever $\{i,j\}$ is an edge of $\sK$, which
implies that $[\widehat\mu_i,\widehat\mu_j]_w=0$ whenever
$\{i,j\}\in\sK$.
\end{proof}

\begin{corollary}\label{cortensgr}
The algebra $H_*(\varOmega(\C P^\infty)^\sK;\k)$ contains the
subalgebra
\begin{equation}\label{tensorgraph}
  T\langle u_1,\ldots,u_m\rangle/(u_i^2=0,\quad u_iu_j+u_ju_i=0\text{ if
  }\{i,j\}\in\sK),
\end{equation}
where $u_i\in H_1(\varOmega(\C P^\infty)^\sK;\k)$ is the Hurewitz
image of $\mu_i\in\pi_1(\varOmega(\C P^\infty)^\sK)$.
\end{corollary}

The subalgebra above maps onto the `fully commutative' algebra
$\Lambda[u_1,\ldots,u_m]$ under the projection map
of~\eqref{paseq}.

In the homotopy fibration~\eqref{omegazkfib}, since $\pi_k(T^m)=0$
for $k>1$, any iterated Samelson product of the form
$[\mu_{i_1},[\mu_{i_2},\cdots[\mu_{i_{k-1}},\mu_{i_k}]\cdots]]$
with $k>1$ composes trivially into $T^m$ and so lifts to
$\varOmega\zk$.

The Whitehead product $[\widehat\mu_i,\widehat\mu_j]_w\colon
S^3\to(\C P^\infty)^\sK$ is nontrivial whenever $\{i,j\}$ is a
missing edge of~$\sK$. We may generalise this construction by
considering missing faces $I=\{i_1,\ldots,i_k\}$ of~$\sK$ (recall
that this means that $I\notin\sK$, but any proper subset of $I$ is
in~$\sK$). Geometrically a missing face defines a subcomplex
$\partial\varDelta(I)\subset\sK$. Define the \emph{$k$-fold higher
Whitehead product} $[\widehat\mu_{i_1},\ldots,\widehat
\mu_{i_k}]_w$ as the composite
\begin{equation}\label{higherWP}
  [\widehat\mu_{i_1},\ldots,\widehat\mu_{i_k}]_w\colon S^{2k-1}
  \stackrel w\longrightarrow
  (S^2)^{\partial\varDelta(I)}\longrightarrow
  (\C P^\infty)^{\partial\varDelta(I)}\longrightarrow(\C P^\infty)^{\sK}
\end{equation}
where $(S^2)^{\partial\varDelta(I)}$ is the fat wedge of $k$
spheres, $w$ is the attaching map of the $2k$-cell in the product
$(S^2)^I$, and the last two maps of the polyhedral products are
induced by the inclusions $S^2\to\C P^\infty$ and
$\partial\varDelta(I)\to\sK$. The \emph{$k$-fold higher Samelson
product}\label{hsamelprodu} $[\mu_{i_1},\ldots,\mu_{i_k}]_s$ is
defined as the adjoint of $[\widehat\mu_{i_1},\ldots,\widehat
\mu_{i_k}]_w$:
\[
  [\mu_{i_1},\ldots,\mu_{i_k}]_s\colon S^{2k-2}\longrightarrow
  \varOmega(\C P^\infty)^{\sK}.
\]

\begin{remark}
As it is standard with higher operations, the higher product
$[\mu_{i_1},\ldots,\mu_{i_k}]$ is defined only when all shorter
higher products of the $\mu_{i_1},\ldots,\mu_{i_k}$ (corresponding
to proper subsets of $I$) are trivial. The general definition of
higher Whitehead and Samelson products (see~\cite{will72})requires
treatment of the indeterminacy, which we avoided in the case of
the polyhedral product $(\C P^\infty)^\sK$ by the canonical choice
of map~\eqref{higherWP}.
\end{remark}

As in the case of standard (2-fold) products, higher Whitehead and
Samelson products of the $\mu_i$ can be iterated and lifted to
$\varOmega\zk$. We summarise this observation as follows:

\begin{proposition}
Any iterated higher Whitehead product $\widehat\nu\colon S^p\to(\C
P^\infty)^\sK$ of the canonical maps $\widehat\mu_i\colon
S^2\to(\C P^\infty)^\sK$ lifts to a map $S^p\to\zk$.

Similarly, any iterated higher Samelson product $\nu\colon
S^{p-1}\to\varOmega(\C P^\infty)^\sK$ of the $\mu_i\colon
S^1\to\varOmega(\C P^\infty)^\sK$ lifts to a map
$S^{p-1}\to\varOmega\zk$.
\end{proposition}

Lifts $S^p\to\zk$ of higher iterated Whitehead products of the
$\widehat\mu_i$ provide important family of spherical classes in
$H_*(\zk)$. We may ask the following question:

\begin{problem}\label{WPrep}
Assume that $\zk$ is homotopy equivalent to a wedge of spheres. Is
it true that all wedge summands are represented by lifts
$S^p\to\zk$ of higher iterated Whitehead products of the canonical
maps $\mu_i\colon S^2\to(\C P^\infty)^\sK$?
\end{problem}

For all known classes of examples when $\zk$ is a wedge of
spheres, the answer to the above question is positive. We shall
give some evidence below.

\begin{example}\label{PAex}\

1. Let $\sK$ be two points. The fibration~\eqref{omegazkfib}
becomes
\[
  \varOmega S^3\to\varOmega(\C P^\infty\vee\C P^\infty)\to S^1\times
  S^1,
\]
and the corresponding sequence of Pontryagin
algebras~\eqref{paseq} is
\[
  0\longrightarrow \k[w]\stackrel i\longrightarrow
  T\langle u_1,u_2\rangle/(u_1^2,u_2^2)\stackrel j\longrightarrow
  \Lambda[u_1,u_2]\longrightarrow 0
\]
where $\k[w]=H_*(\varOmega S^3;\k)$, $\deg w=2$, the map $i$ takes
$w$ to the commutator $u_1u_2+u_2u_1$, and $j$ is the projection
to the quotient by the ideal generated by $u_1u_2+u_2u_1$ (an
exercise). So
\[
  H_*(\varOmega(\C P^\infty\vee\C P^\infty);\k)=T\langle
  u_1,u_2\rangle/(u_1^2,u_2^2)=\Lambda[u_1]\mathbin\star\Lambda[u_2]
\]
is the free product of two exterior algebras and $i(\k[w])$ is its
commutator subalgebra. In particular, exact sequence~\eqref{paseq}
does not split multiplicatively in this example.

Here $u_1,u_2$ are the Hurewitz images of $\mu_1,\mu_2$, the
commutator $u_1u_2+u_2u_1$ is the Hurewitz image of the Samelson
product $[\mu_1,\mu_2]_s\colon S^2\to\varOmega(\C P^\infty\vee\C
P^\infty)$, and $w\in H_2(\varOmega S^3)$ is the Hurewitz image of
the lift of $[\mu_1,\mu_2]_s$ to $\varOmega S^3$.

2. Now let $\sK=\partial\varDelta^2$, the boundary of a triangle.
The fibration~\eqref{omegazkfib} becomes
\[
  \varOmega S^5\to\varOmega(\C P^\infty)^{\partial\varDelta^2}\to
  T^3,
\]
where $\varOmega(\C P^\infty)^{\partial\varDelta^2}$ is the fat
wedge of 3 copies of $\C P^\infty$. We have $H_*(\varOmega
S^5;\k)=\k[w]$, $\deg w=4$. Algebra~\eqref{tensorgraph} is
isomorphic to $\Lambda[u_1,u_2,u_3]$, so the sequence of
Pontryagin algebras~\eqref{paseq} splits multiplicatively in this
example:
\[
  0\longrightarrow \k[w]\longrightarrow
  \k[w]\otimes\Lambda[u_1,u_2,u_3]\longrightarrow
  \Lambda[u_1,u_2,u_3]\longrightarrow 0.
\]
Here $w\in H_4(\varOmega(\C P^\infty)^{\partial\varDelta^2};\k)$
is the Hurewicz image of the higher Samelson product
$[\mu_1,\mu_2,\mu_3]_s\in\pi_4(\varOmega(\C
P^\infty)^{\partial\varDelta^2})$, which lifts to $\varOmega S^5$.
The fact that $[\mu_1,\mu_2,\mu_3]_s$ is a nontrivial higher
Samelson product (and its Hurewicz image $w$ is the `higher
commutator product' of $u_1,u_2,u_3$) constitutes the additional
information necessary to distinguish between the topological
monoids $\varOmega(\C P^\infty)^{\partial\varDelta^2}$ and
$\varOmega S^5\times T^3$.

This calculation generalises easily to the case
$\sK=\partial\varDelta^{m-1}$, showing that
\[
  H_*(\varOmega(\C
  P^\infty)^\sK;\k)\cong\k[w]\otimes\Lambda[u_1,\ldots,u_m]
\]
where $\deg w=2m-2$.
\end{example}

\subsection*{Topological models for loop spaces}
We consider the diagram
\begin{equation}\label{DKS1}
  \mathcal D_\sK(S^1)\colon\ca{(\sK)}\to\cat{tmon}
\end{equation}
whose value on the morphism $I\subset J$ is the monomorphism of
tori $T^I\to T^J$. The classifying space diagram $B\mathcal
D_\sK(S^1)$ is $\mathcal D_\sK(\C
P^\infty,\pt)\colon\ca{(\sK)}\to\cat{top}$ with colimit $(\C
P^\infty)^\sK$. We denote the colimit and homotopy colimit of
$\mathcal D_\sK(S^1)$ by $\colim^{\scat{tmon}}_{I\in\sK}T^I$ and
$\hocolim^{\scat{tmon}}_{I\in\sK}T^I$, respectively.

\begin{theorem}[\cite{p-r-v04}]\label{gomegahocolim}
There is a commutative diagram
\begin{equation}\label{omegahocolimdiag}
\begin{CD}
  @.\varOmega\hocolim^{\scat{top}}_{I\in\sK} BT^I@>g>\simeq>\hocolim^{\scat{tmon}}_{I\in\sK}T^I\\
  @.@V\varOmega p^{\scat{top}}V\simeq V@VV p^{\scat{tmon}}V\\
  \varOmega(\C P^\infty)^\sK@=\varOmega\colim^{\scat{top}}_{I\in\sK} BT^I@>>>\colim^{\scat{tmon}}_{I\in\sK}T^I
\end{CD}
\end{equation}\\
in $Ho(\cat{tmon})$, where the top and left homomorphisms are
homotopy equivalences.
\end{theorem}
\begin{proof}
We apply Corollary~\ref{loopshocolim} with $\mathcal D=\mathcal
D_\sK(S^1)$. The left projection $\varOmega p^{\scat{top}}\colon
\varOmega\hocolim^{\scat{top}}_{I\in\sK} BT^I\to
\varOmega\colim^{\scat{top}}_{I\in\sK} BT^I$ is a weak equivalence
because $B\mathcal D_\sK(S^1)=\mathcal D_\sK(\C P^\infty,\pt)$ is
a cofibrant diagram in~$\cat{top}$.
\end{proof}

\begin{corollary}
There is a weak equivalence
\[
  \varOmega(\C
  P^\infty)^\sK\simeq\hocolim^{\scat{tmon}}_{I\in\sK}T^I
\]
in~$\cat{tmon}$.
\end{corollary}

The right projection $p^{\scat{tmon}}$ (and therefore the bottom
homomorphism in~\eqref{omegahocolimdiag}) is not a weak
equivalence in general, because $\mathcal D_\sK(S^1)$ is
\emph{not} a cofibrant diagram in~$\cat{tmon}$. The appropriate
examples are discussed below.

\begin{example}\label{tmonex}\

1. Let $\sK$ be two points. Then
\[
  (\C P^\infty)^\sK=\C P^\infty\vee\C P^\infty,\qquad
  \colim^{\scat{tmon}}\mathcal D_\sK(S^1)=S^1\mathbin{\star}S^1
\]
where $\star$ denotes the free product of topological monoids,
i.e. the coproduct in~\cat{tmon}. The bottom homomorphism
in~\eqref{omegahocolimdiag} is $\varOmega(\C P^\infty\vee\C
P^\infty)\to S^1\mathbin{\star}S^1$, it is a weak equivalence
in~$\cat{tmon}$.

2. Now let $\sK=\partial\varDelta^2$. The loop space $\varOmega(\C
P^\infty)^{\partial\varDelta^2}$ is described in
Example~\ref{PAex}.2. On the other hand, the colimit of $\mathcal
D_\sK(S^1)$ is obtained by taking quotient of
$T^{\{1\}}\mathbin{\star}T^{\{2\}}\mathbin{\star}T^{\{3\}}$ by the
commutativity relations
\[
  t_1\mathbin{\star}t_2=t_2\mathbin{\star}t_1,\quad
  t_2\mathbin{\star}t_3=t_3\mathbin{\star}t_2,\quad
  t_3\mathbin{\star}t_1=t_1\mathbin{\star}t_3
\]
for $t_i\in T^{\{i\}}$, so that
$\colim^{\scat{tmon}}_{I\in\sK}T^I=T^3$. It follows that the
bottom map $\varOmega(\C P^\infty)^{\partial\varDelta^2}\to T^3$
in~\eqref{omegahocolimdiag} is not a weak equivalence; it has
kernel $\varOmega S^5$.

The diagram $\mathcal D=\mathcal
D_\sK(S^1)\colon\ca{(\sK)}\to\cat{tmon}$ is not cofibrant in this
example. Indeed, if we take $I=\{1,2\}$, then the induced diagram
over the overcategory $\ca{(K)}\under I=\cat{cat}(\varDelta(I))$
has the form
\[
\begin{CD}
  \pt @>>> T^{\{1\}}\\
  @VVV @VVV\\
  T^{\{2\}} @>>> T^{\{1\}}\times T^{\{2\}}
\end{CD}
\]
The map $\colim \mathcal
D|_{\scat{cat}(\partial\varDelta(I))}\to\mathcal D(I)$ is the
projection $T^{\{1\}}\mathbin{\star}T^{\{2\}}\to T^{\{1\}}\times
T^{\{2\}}$ from the free product to the cartesian product, which
is not a cofibration in~\cat{tmon}.
\end{example}

\subsection*{Algebraic models for loop spaces}
Our next aim is to obtain an algebraic analogue of
Theorem~\ref{gomegahocolim}. We work over over a commutative
ring~$\k$.

We define the \emph{face coalgebra}\label{facecoalgebra}
$\k\langle\sK\rangle$ as the graded dual of the face ring
$\k[\sK]$. As a $\k$-module, $\k\langle\sK\rangle$ is free on
generators $v_\sigma$ corresponding to multisets of $m$ elements
of the form
\[
  \sigma=\{\underbrace{1,\ldots,1}_{k_1},\underbrace{2,\ldots,2}_{k_2},
  \ldots,\underbrace{m,\ldots,m}_{k_m}\}
\]
such that \emph{support} of $\sigma$ (i.e. the set
$I_\sigma=\{i\in[m]\colon k_i\ne0\}$) is a simplex of~$\sK$. The
element $v_\sigma$ is dual to the monomial
$v_1^{k_1}v_2^{k_2}\cdots v_m^{k_m}\in\k[\sK]$. The
comultiplication takes the form
\[
  \Delta v_\sigma=
  \sum_{\sigma=\tau\,\sqcup\,\tau'}v_\tau\otimes v_{\tau'},
\]
where the sum ranges over all partitions of $\sigma$ into
submultisets $\tau$ and $\tau'$.

We recall the Adams cobar construction
$\varOmega_*\colon\cat{dgc}\to\cat{dga}$, see~\eqref{barcobar},
and the Quillen functor $L_*\colon\cat{dgc}\to\cat{dgl}$,
see~\eqref{lcadj}. The loop algebra
$\varOmega_*\k\langle\sK\rangle$ is our first algebraic model for
$\varOmega(\C P^\infty)^\sK$:

\begin{proposition}\label{adamsSR}
There is an isomorphism of graded algebras
\[
  H_*(\varOmega(\C
  P^\infty)^\sK;\k)\cong H(\varOmega_*\k\langle\sK\rangle)=
  \mathop{\mathrm{Cotor}}\nolimits_{\k\langle\sK\rangle}(\k,\k).
\]
\end{proposition}
\begin{proof}
By dualising the integral formality results of \cite[Theorem
4.8]{no-ra05}, we obtain a zigzag of quasi-isomorphisms
\begin{equation}\label{cdgam}
  C_*((\C P^\infty)^\sK;\k)\stackrel\simeq\longleftarrow\cdots
  \stackrel\simeq\longrightarrow \k\langle\sK\rangle.
\end{equation}
in $\cat{dgc}$ (when~$\k=\Q$ this follows from
Theorem~\ref{xkform}). Since $\varOmega_*$ preserves
quasi-isomorphisms, the zigzag above combines with Adams' result
(Theorem~\ref{adamscobar}) to obtain the required isomorphism of
algebras.
\end{proof}

\begin{remark}
When $\k$ is a field, there are isomorphisms
\[
  H_*(\varOmega(\C
  P^\infty)^\sK;\k)\cong
  \mathop{\mathrm{Cotor}}\nolimits_{\k\langle\sK\rangle}(\k,\k)
  \cong\Ext_{\k[\sK]}(\k,\k).
\]
\end{remark}

The graded algebra underlying the cobar construction
$\varOmega_*\k\langle\sK\rangle$ is the tensor algebra
$T(s^{-1}\,\overline{\k\langle\sK\rangle\!\!}\,)$ on the
desuspended $\k$-module
$\overline{\k\langle\sK\rangle\!\!}\,=\Ker(\varepsilon\colon\k\langle\sK\rangle
\to\k)$; the differential is defined on generators by
\[
  d(s^{-1}v_\sigma)=
  \sum_{\sigma=\tau\sqcup\tau';\;\tau,\,\tau'\ne\varnothing}
  s^{-1}v_\tau\otimes s^{-1}v_{\tau'},
\]
because $d=0$ on $\k\langle\sK\rangle$. For future purposes it is
convenient to write $s^{-1}v_\sigma$ as $\chi_\sigma$ for any
multiset~$\sigma$.

We define now some algebraic diagrams over $\ca{(\sK)}$. Our
previous algebraic diagrams such as~\eqref{DsupK} were
commutative, contravariant and cohomological, but to investigate
the loop space $\varOmega(\C P^\infty)^\sK$ we introduce models
that are covariant and homological. We consider the diagram
\[
  \k[\,\cdot\:]_\sK\colon\ca{(\sK)}\to\cat{dga}, \quad
  I\mapsto\k[v_i\colon i\in I]
\]
which maps a morphism $I\subset J$ to the monomorphism of
polynomial algebras $\k[v_i\colon i\in I]\to \k[v_i\colon i\in J]$
with $\deg v_i=2$ and zero differential. Similarly, we define the
diagrams
\begin{equation}\label{3diagrams}
\begin{aligned}
  \Lambda[\,\cdot\:]_\sK\colon\ca{(\sK)}\to\cat{dga},&\quad
  I\mapsto\Lambda[u_i\colon i\in I],& \deg u_i=1,\\
  \k\langle\,\cdot\:\rangle_\sK\colon\ca{(\sK)}\to\cat{dgc},&\quad
  I\mapsto\k\langle v_i\colon i\in I\rangle,& \deg v_i=2,\\
  \cl(\,\cdot\:)_\sK\colon\ca{(\sK)}\to\cat{dgl},&\quad
  I\mapsto\cl(u_i\colon i\in I),& \deg u_i=1,
\end{aligned}
\end{equation}
where $\k\langle v_i\colon i\in I\rangle$ denotes the free
commutative coalgebra and $\cl(u_i\colon i\in I)$ denotes the
commutative Lie algebra on $|I|$ generators.

The individual algebras and coalgebras in these diagrams are all
commutative, but the context demands they be interpreted in the
non-commutative categories; this is especially important when
forming limits and colimits.

Note that
$\colim^{\scat{dgc}}\k\langle\,\cdot\:\rangle_\sK=\k\langle\sK\rangle$,
while $\colim^{\scat{dga}}\Lambda[\,\cdot\:]_\sK$ is the
non-commutative algebra~\eqref{tensorgraph} (an exercise).

\begin{proposition}\label{symwe}
There are acyclic fibrations
\[
  \varOmega_*\k\langle v_i\colon i\in I\rangle\stackrel\simeq\longrightarrow\Lambda[u_i\colon
  i\in I]\quad\text{ and }\quad L_*\k\langle v_i\colon i\in I\rangle
  \stackrel\simeq\longrightarrow\cl(u_i\colon i\in I)
\]
in $\cat{dga}$ and $\cat{dgl}$ respectively, for any set
$I\subset[m]$.
\end{proposition}
\begin{proof}
We define the first map by $\chi_i\mapsto u_i$ for $1\le i\le m$.
Because
\begin{equation}\label{cobarrel}
  d\chi_{ii}=\chi_i\otimes\chi_i,\qquad
  d\chi_{ij}=\chi_i\otimes\chi_j+\chi_j\otimes\chi_i\quad
  \text{for}\quad i\ne j
\end{equation}
hold in $\varOmega_*\k\langle v_i\colon i\in I\rangle$, the map is
consistent with the exterior relations in its target. So it is an
epimorphism and quasi-isomorphism in $\cat{dga}$, and hence an
acyclic fibration. The corresponding result for $\cat{dgl}$
follows by restriction to primitives.
\end{proof}

Observe that the diagram $\Lambda[\,\cdot\:]_\sK$ in $\cat{dga}$
can be thought of as the diagram of homology algebras of
topological monoids in the diagram $\mathcal D_\sK(S^1)$,
see~\eqref{DKS1}, and the diagram $\k\langle\,\cdot\:\rangle_\sK$
in $\cat{dgc}$ is the diagram of homology coalgebras of spaces in
the classifying diagram~$B\mathcal D_\sK(S^1)$. This relationship
extends to the following algebraic analogue of
Theorem~\ref{gomegahocolim}:

\begin{theorem}[\cite{pa-ra08}]\label{homegahocolim}
There is a commutative diagram
\begin{equation}\label{aomegahocolimdiag}
\begin{CD}
  \varOmega_*\hocolim^{\scat{dgc}}\k\langle\,\cdot\:\rangle_\sK
  @>h>\simeq>\hocolim^{\scat{dga}}\Lambda[\,\cdot\:]_\sK\\
  @V\varOmega_* p^{\scat{dgc}}V\simeq V@VV p^{\scat{dga}}V\\
  \varOmega_*\k\langle\sK\rangle@>>>\colim^{\scat{dga}}\Lambda[\,\cdot\:]_\sK
\end{CD}
\end{equation}\\
in $Ho(\cat{dga})$, where the top and left arrows are
isomorphisms.
\end{theorem}
\begin{proof}
This follows by considering the diagram
\begin{align*}
&\begin{CD}
  \varOmega_*\hocolim^{\scat{dgc}}\k\langle\,\cdot\:\rangle_\sK
  @.\hocolim^{\scat{dga}}\varOmega_*\k\langle\,\cdot\:\rangle_\sK
  @>\simeq>>\hocolim^{\scat{dga}}\Lambda[\,\cdot\:]_\sK\\
  @V\varOmega_*{p^{\scat{dgc}}}V{\simeq}V
  @V{p^{\scat{dga}}}V{\simeq}V
  @V{p^{\scat{dga}}}VV\\
  \varOmega_*\colim^{\scat{dgc}}\k\langle\,\cdot\:\rangle_\sK
  @<\cong<<\colim^{\scat{dga}}\varOmega_*\k\langle\,\cdot\:\rangle_\sK
  @>>>\colim^{\scat{dga}}\Lambda[\,\cdot\:]_\sK\\
\end{CD}\\
&\qquad\qquad\;\|\\
&\qquad\quad\;\varOmega_*\k\langle\sK\rangle
\end{align*}
Here the top right horizontal map is induced by the map of
diagrams
$\varOmega_*\k\langle\,\cdot\:\rangle_\sK\to\Lambda[\,\cdot\:]_\sK$
whose objectwise maps are acyclic fibrations from
Proposition~\ref{symwe}; the map of homotopy colimits is a weak
equivalence because the map of diagrams is an acyclic Reedy
fibration. The right square is commutative. The central vertical
map is a weak equivalence because the diagram
$\varOmega_*\k\langle\,\cdot\:\rangle_\sK$ is Reedy cofibrant. The
bottom left horizontal map is an isomorphism because $\varOmega_*$
is left adjoint. The left vertical map is a weak equivalence
because the diagram $\k\langle\,\cdot\:\rangle_\sK$ is Reedy
cofibrant and $\varOmega_*$ preserves weak equivalences. The
resulting zigzag of quasi-isomorphisms
\[
  \varOmega_*\hocolim^{\scat{dgc}}\k\langle\,\cdot\:\rangle_\sK
  \stackrel\simeq\longrightarrow\cdots
  \stackrel\simeq\longrightarrow
  \hocolim^{\scat{dga}}\Lambda[\,\cdot\:]_\sK
\]
induces an isomorphism in the homotopy category $Ho(\cat{dga})$;
we denote it by~$h$.
\end{proof}

As in the case of diagram~\eqref{omegahocolimdiag}, the right and
bottom maps in~\eqref{aomegahocolimdiag} are not weak equivalences
in general, because $\Lambda[\,\cdot\:]_\sK$ is \emph{not} a
cofibrant diagram in~$\cat{dga}$.

The following statement is proved similarly.

\begin{theorem}[\cite{pa-ra08}]
There is a homotopy commutative diagram
\begin{equation}\label{lomegahocolimdiag}
\begin{CD}
  L_*\hocolim^{\scat{dgc}}\Q\langle\,\cdot\:\rangle_\sK
  @>>\simeq>\hocolim^{\scat{dgl}}\cl(\,\cdot\:)_\sK\\
  @V L_* p^{\scat{dgc}}V\simeq V@VV p^{\scat{dgl}}V\\
  L_*\Q\langle\sK\rangle@>>>\colim^{\scat{dgl}}\cl(\,\cdot\:)_\sK
\end{CD}
\end{equation}\\
in $Ho(\cat{dgl})$, where the top and left arrows are
isomorphisms.
\end{theorem}

The homotopy colimit decomposition above defines our second
algebraic model for $\varOmega(\C P^\infty)^\sK$:

\begin{corollary}\label{hghocolim}
For any simplicial complex $\sK$ and commutative ring~$\k$, there
are isomorphisms
\begin{align*}
H_*\bigl(\varOmega(\C P^\infty)^\sK;\k\bigr)&\cong
H\bigl(\hocolim^{\scat{dga}}\Lambda[\,\cdot\:]_\sK\bigr)\\
\pi_*\bigl(\varOmega(\C P^\infty)^\sK)\otimes_{\Z}\Q&\cong
H\bigl(\hocolim^{\scat{dgl}}\cl(\,\cdot\:)_\sK\bigr)
\end{align*}
of graded algebras and Lie algebras respectively.
\end{corollary}

\begin{example}\label{exmpzk}\

1. Let $\sK$ be a discrete complex on $m$ vertices, so that $(\C
P^\infty)^\sK$ is a wedge of $m$ copies of~$\C P^\infty$. The
cobar construction $\varOmega_*\k\langle\sK\rangle$ on the
corresponding face coalgebra is generated as an algebra by the
elements of the form $\chi_{i\ldots i}$ with $i\in[m]$. The first
identity of~\eqref{cobarrel} still holds, but
$\chi_i\otimes\chi_j+\chi_j\otimes\chi_i$ is no longer a boundary
for $i\ne j$ since there is no element $\chi_{ij}$ in
$\varOmega_*\k\langle\sK\rangle$. We obtain a quasi-isomorphism
\[
  \varOmega_*\k\langle\sK\rangle\stackrel\simeq\longrightarrow
  T_\k(u_1,\ldots,u_m)/(u_i^2=0,\; 1\le i\le m)
\]
that maps $\chi_i$ to $u_i$. The right hand side is isomorphic to
$H_*(\varOmega(\C P^\infty)^\sK;\k)$; it is the coproduct (the
free product) of $m$ algebras~$\Lambda[u_i]$. Therefore, the
bottom map
$\varOmega_*\k\langle\sK\rangle\to\colim^{\scat{dga}}\Lambda[\,\cdot\:]_\sK$
in~\eqref{aomegahocolimdiag} is a quasi-isomorphism in this
example.

The algebra $H_*(\varOmega\zk;\k)$ is the commutator subalgebra of
$H_*(\varOmega(\C P^\infty)^\sK;\k)$ according to~\eqref{paseq}.
In contains iterated commutators
$[u_{i_1},[u_{i_2},\cdots[u_{i_{k-1}},u_{i_k}]\cdots]]$ with $k\ge
2$ corresponding to iterated Samelson products in
$\pi_*(\varOmega(\C P^\infty)^\sK)$. On the other hand, $\zk$ is a
wedge of spheres given by~\eqref{zkmwedge}. Therefore,
$H_*(\varOmega\zk;\k)$ is a free (tensor) algebra on the
generators corresponding to the wedge summands. So the number of
independent iterated commutators of length~$k$ is $(k-1)\bin mk$.
This fact can be proved purely algebraically, see
Corollary~\ref{CNzk} below. There are no higher Samelson products
(and higher iterated commutators) in this example, as $\sK$ does
not have missing faces with $>2$ vertices.

\smallskip

2. Let $\sK=\partial\varDelta^2$. As we have seen in
Example~\ref{PAex}.2,
\[
  H_*(\varOmega(\C
  P^\infty)^\sK;\k)\cong\k[w]\otimes\Lambda[u_1,u_2,u_3].
\]
This can also be seen algebraically using the cobar model
$\varOmega_*\k\langle\sK\rangle$. Here $u_i$ is the homology
classes of the element $\chi_i$, and $w\in H_4(\varOmega(\C
P^\infty)^\sK;\k)$ is the homology class of the 4-dimensional
cycle
\[
  \psi=\chi_1\chi_{23}+\chi_2\chi_{13}+\chi_3\chi_{12}+
  \chi_{12}\chi_3+\chi_{13}\chi_2+\chi_{23}\chi_1,
\]
whose failure to bound is due to the non-existence
of~$\chi_{123}$. Relations~\eqref{cobarrel} hold, and give rise to
the exterior relations between $u_1,u_2,u_3$. Furthermore, a
direct check shows that $\chi_i\psi-\psi\chi_i$ is a boundary,
which implies that $u_i$ commutes with $w$ in $H_*(\varOmega(\C
P^\infty)^\sK;\k)$ for $1\le i\le 3$.

Here the colimit of  $\Lambda[\,\cdot\:]_\sK$ in $\cat{dga}$ is
obtained by taking quotient of $T_\k(u_1,u_2,u_3)$ by all exterior
relations. Therefore,
$\colim^{\scat{dga}}\Lambda[\,\cdot\:]_\sK=\Lambda[u_1,u_2,u_3]$
and the bottom map
$\varOmega_*\k\langle\sK\rangle\to\colim^{\scat{dga}}\Lambda[\,\cdot\:]_\sK$
in~\eqref{aomegahocolimdiag} is \emph{not} a quasi-isomorphism in
this example.

\smallskip

3. Let $\sK=\sk^1\partial\varDelta^3$, the 1-skeleton of a
3-simplex, or a complete graph on 4 vertices. Arguments similar to
those of the previous example show that the Pontryagin algebra
$H_*(\varOmega(\C P^\infty)^\sK;\k)$ contains $1$-dimensional
classes $u_1,\ldots,u_4$ and $4$-dimensional classes $w_{123}$,
$w_{124}$, $w_{134}$, $w_{234}$, corresponding to the four missing
faces with three vertices each. For example, $w_{123}$ is the
homology class of the cycle $\psi_{123}$ which may be thought of
as the `boundary of the non-existing element $\chi_{123}$'.
Identities~\eqref{cobarrel} give rise to the exterior relations
between $u_1,\ldots,u_4$. We may easily check that $u_i$ commutes
with $w_{jkl}$ if $i\in\{j,k,l\}$. There are four remaining
non-trivial commutators of the form $[u_i,w_{jkl}]$ with all
$i,j,k,l$ different.

It follows that the commutator subalgebra $H_*(\varOmega\zk;\k)$
contains four higher commutators $w_{jkl}=[u_j,u_k,u_l]$ (the
Hurewitz images of the higher Samelson products
$[\mu_j,\mu_k,\mu_l]_s$) and four iterated commutators
$[u_i,w_{jkl}]$. On the other hand, Theorem~\ref{gr-th} gives a
homotopy equivalence
\[
  \zk\simeq(S^5)^{\vee 4}\vee(S^6)^{\vee 3},
\]
which implies that $H_*(\varOmega\zk;\k)$ is a free algebra on
four 4-dimensional and \emph{three} 5-dimensional generators. The
point is that the commutators $[u_i,w_{jkl}]$ are subject to one
extra relation, which can be derived as follows. Consider the
relation
\begin{equation}\label{reln}
  d\chi_{1234}\;=\;(\chi_1\chi_{234}+\chi_{234}\chi_1)
  +\cdots+(\chi_4\chi_{123}+\chi_{123}\chi_4)+\beta
\end{equation}
in $\varOmega_*\k\langle v_1,v_2,v_3,v_4\rangle$, where $\beta$
consists of terms $\chi_{\sigma}\chi_{\tau}$ such that
$|\sigma|=|\tau|=2$. Denote the first four summands on the right
hand side of \eqref{reln} by $\alpha_1$, $\alpha_2$, $\alpha_3$,
$\alpha_4$ respectively, and apply the differential to both sides.
Observing that
$d\alpha_1=-\chi_1\psi_{234}+\psi_{234}\chi_1=-[\chi_1,\psi_{234}]$,
and similarly for $d\alpha_2$, $d\alpha_3$ and $d\alpha_4$, we
obtain
\[
  [\chi_1,\psi_{234}]+[\chi_2,\psi_{134}]+
  [\chi_3,\psi_{124}]+[\chi_4,\psi_{123}]\;=\;d\beta.
\]

The outcome is an isomorphism
\[
  H_*(\varOmega(\C P^\infty)^\sK;\k)\;\cong\;
  T_\k(u_1,u_2,u_3,u_4,w_{123},w_{124},w_{134},w_{123})/\mathcal I,
\]
where $\deg w_{ijk}=4$ and $\mathcal I$ is generated by three
types of relation:
\begin{enumerate}
\item[$\cdot$] exterior algebra relations for $u_1,u_2,u_3,u_4$;
\item[$\cdot$] $[u_i,w_{jkl}]=0$ for $i\in\{j,k,l\}$;
\item[$\cdot$]
  $[u_1,w_{234}]+[u_2,w_{134}]+[u_3,w_{124}]+[u_4,w_{123}]=0$.
\end{enumerate}
As $w_{ijk}$ is the higher commutator of $u_i$, $u_j$ and $u_k$,
the third relation may be considered as a higher analogue of the
Jacobi identity.
\end{example}

It is a challenging task to construct explicit algebraic models
for $\varOmega(\C P^\infty)^\sK$ and $\varOmega\zk$ which would
include a description of the Pontryagin algebra structure, as well
as higher Samelson and commutator products. The situation is
considerably simpler when $\sK$ is a flag complex, as there are no
higher products; this is the subject of the next section.

\subsection*{Exercises}
\begin{exercise}
For any $2n$-dimensional (quasi)toric manifold $M$, show that
there is a fibration
\[
  \varOmega M\longrightarrow\varOmega(\C P^\infty)^\sK\longrightarrow
  T^n,
\]
which splits in~$\cat{top}$.
\end{exercise}

\begin{exercise}
For the classes of simplicial complexes $\sK$ described in
Propositions~\ref{zkmpoints} and~\ref{zktree} (discrete complexes
and trees), show that each wedge summand of $\zk$ is represented
by a lift $S^p\to\zk$ of iterated Whitehead products of the
$\mu_i\colon S^2\to(\C P^\infty)^\sK$ (no higher products appear
here). Describe the corresponding iterated brackets explicitly. In
particular, the answer to Problem~\ref{WPrep} is positive for
these two classes of~$\sK$.
\end{exercise}

\begin{exercise}
Show that $H^*(\varOmega(\C P^\infty\vee\C P^\infty);\k)=T\langle
u_1,u_2\rangle/(u_1^2,u_2^2)$ and describe the sequence of
Pontryagin algebras corresponding to the fibration $\varOmega
S^3\to\varOmega(\C P^\infty\vee\C P^\infty)\to S^1\times S^1$.
\end{exercise}

\begin{exercise}
Show that $\colim^{\scat{dga}}\Lambda[\,\cdot\:]_\sK$ is the
algebra given by~\eqref{tensorgraph}.
\end{exercise}

\section{The case of flag complexes}\label{homotflag}
In this section we study the loop spaces associated with
\emph{flag complexes} $\sK$. Such complexes have significantly
simpler combinatorial properties, which are reflected in the
homotopy theory of the toric spaces. We modify results of the
previous section in this context, and focus on applications to the
Pontryagin rings and homotopy Lie algebras of $\varOmega(\C
P^\infty)^\sK$ and~$\varOmega\zk$. We also describe completely the
class of flag complexes $\sK$ for which $\zk$ is homotopy
equivalent to a wedge of spheres.

For any simplicial complex $\sK$ on~$[m]$, recall that a subset
$I\subset[m]$ is called a missing face when every proper subset
lies in $\sK$, but $I$ itself does not. If every missing face of
$\sK$ has $2$ vertices, then $\sK$ is a flag complex;
equivalently, $\sK$ is flag when every set of vertices that is
pairwise connected spans a simplex. A flag complex is therefore
determined by its 1-skeleton, which is a graph. When $\sK$ is
flag, we may express the face ring as
\[
%\begin{equation}\label{srflag}
  \k[\sK]=T_\k(v_1,\ldots,v_m)\bigr/
  (v_iv_j-v_jv_i=0\text{ for }\{i,j\}\in\sK,\;
  v_iv_j=0\text{ for }\{i,j\}\notin\sK).
%\end{equation}
\]
It is therefore \emph{quadratic}\label{quadralge}, in the sense
that it is the quotient of a free algebra by quadratic relations.

The following result of Fr\"oberg allows us to calculate the
Yoneda algebras $\Ext_A(\k,\k)$ explicitly for a class of
quadratic algebras $A$ that includes face rings of flag complexes.

\begin{proposition}[{\cite[\S3]{frob75}}]\label{srkos}
When $\k$ is a field and $\sK$ is a flag complex, there is an
isomorphism of graded algebras
\begin{equation}\label{srext}
  \Ext_{\k[\sK]}(\k,\k)\cong
  T_\k(u_1,\ldots,u_m)\bigr/(u_i^2=0,\; u_iu_j+u_ju_i=0\text{ for
}\{i,j\}\in\sK).
\end{equation}
\end{proposition}

\begin{remark}
The algebra on the right hand side of~\eqref{srext} is the
\emph{quadratic dual} of~$\k[\sK]$. A quadratic algebra $A$ is
called \emph{Koszul}\label{koszulalge} if its quadratic dual
coincides with $\Ext_A(\k,\k)$, so Proposition~\ref{srkos} asserts
that $\k[\sK]$ is Koszul when $\sK$ is flag.
\end{remark}

When $\sK$ is flag, \eqref{tensorgraph} is the whole Pontryagin
algebra $H_*(\varOmega(\C P^\infty)^\sK;\k)$:

\begin{theorem}[\cite{pa-ra08}]\label{hldj}
For any flag complex $\sK$, there are isomorphisms
\begin{align*}
H_*(\varOmega(\C P^\infty)^\sK;\k)&\;\cong\;
T_\k(u_1,\ldots,u_m)\bigr/
  (u_i^2=0,\; u_iu_j+u_ju_i=0\text{ for }\{i,j\}\in\sK)\\
\pi_*(\varOmega(\C P^\infty)^\sK)\otimes_\Z\Q&\;\cong\;
\fl(u_1,\ldots,u_m)\bigr/
  \bigl([u_i,u_i]=0,\; [u_i,u_j]=0\text{ for }\{i,j\}\in\sK\bigr),
\end{align*}
where $\k$ is $\Z$ or a field, $\fl(\,\cdot\,)$ denotes a free Lie
algebra and $\deg u_i=1$.
\end{theorem}
\begin{proof}
By Proposition~\ref{adamsSR}, $H_*(\varOmega(\C
P^\infty)^\sK;\k)\cong\mathop{\mathrm{Cotor}}\nolimits_{\k\langle\sK\rangle}(\k,\k)$.
When $\k$ is a field,
$\mathop{\mathrm{Cotor}}\nolimits_{\k\langle\sK\rangle}(\k,\k)\cong\Ext_{\k[\sK](\k,\k)}$
by~\eqref{cotorext}, and the first required isomorphism follows
from Proposition~\ref{srkos}.

Now let $\k=\Z$. Denote by $A$ algebra~\eqref{tensorgraph} with
$\k=\Z$. Then $A$ includes as a subalgebra in $H_*(\varOmega(\C
P^\infty)^\sK;\Z)$ by Corollary~\ref{cortensgr}. Consider the
exact sequence
\[
  0\to A\stackrel i\longrightarrow H_*(\varOmega(\C
  P^\infty)^\sK;\Z)\longrightarrow R\to 0
\]
where $R$ is the cokernel of~$i$. By tensoring with a field $\k$
we obtain an exact sequence
\[
  A\otimes\k\stackrel{i\otimes\k}\longrightarrow H_*(\varOmega(\C
  P^\infty)^\sK;\Z)\otimes\k\longrightarrow R\otimes\k\to 0
\]
The composite map
\[
  A\otimes\k\stackrel{i\otimes\k}\longrightarrow H_*(\varOmega(\C
  P^\infty)^\sK;\Z)\otimes\k\stackrel{j}\longrightarrow H_*(\varOmega(\C
  P^\infty)^\sK;\k)
\]
is an isomorphism by the argument in the previous paragraph, and
$j$ is a monomorphism by the universal coefficient theorem.
Therefore, $i\otimes\k$ is also an isomorphism, which implies that
$R\otimes\k=0$ for any field~$\k$. Thus, $R=0$ and $i$ is an
isomorphism.

The isomorphism of Lie algebras follows by restriction to
primitives.
\end{proof}

The authors are grateful to Kouyemon Iriye and Jie Wu for pointing
out that the original argument of~\cite{pa-ra08} for the theorem
above can be extended to the case~$\k=\Z$.

The associative algebra and the Lie algebra from
Theorem~\ref{hldj} are examples of \emph{graph
products}\label{graphprodu}, by which one usually means algebraic
objects described by generators corresponding to the vertices in a
simple graph, with each edge giving rise to a commutativity
relation between the generators corresponding to its two ends.
These two graph products can be also described as the colimits of
the diagrams $\Lambda[\,\cdot\:]_\sK$ and $\cl(\,\cdot\:)_\sK$ in
$\cat{dga}$ and $\cat{dgl}$ respectively, see~\eqref{3diagrams}.
It follows that the bottom maps in~\eqref{aomegahocolimdiag}
and~\eqref{lomegahocolimdiag} are quasi-isomorphisms when $\sK$ is
flag, and homotopy colimit in the models of
Corollary~\ref{hghocolim} can be replaced by colimit:

\begin{corollary}\label{flagcolim}
For any flag complex $\sK$, there are isomorphisms
\begin{align*}
H_*\bigl(\varOmega(\C P^\infty)^\sK;\k\bigr)&\cong
\colim^{\scat{dga}}\Lambda[\,\cdot\:]_\sK\\
\pi_*\bigl(\varOmega(\C P^\infty)^\sK\bigr)\otimes_{\Z}\Q&\cong
\colim^{\scat{dgl}}\cl(\,\cdot\:)_\sK
\end{align*}
of graded algebras and Lie algebras respectively, where $\k=\Z$ or
a field.
\end{corollary}

\begin{remark}
The bottom homomorphism
\[
  \varOmega(\C P^\infty)^\sK\to \colim^{\scat{tmon}}_{I\in\sK}T^I
\]
in the topological model~\eqref{omegahocolimdiag} is also a
homotopy equivalence when $\sK$ is flag,
by~\cite[Proposition~6.3]{p-r-v04}. Furthermore, there are
analogues of this result for other polyhedral powers. Interesting
cases are $(\R P^\infty)^\sK$ and $(S^1)^\sK$, for which there are
homotopy equivalence homomorphisms
\begin{equation}\label{raCAcolim}
\begin{aligned}
  \varOmega(\R P^\infty)^\sK & \stackrel\simeq\longrightarrow
  \colim^{\scat{tmon}}_{I\in\sK}\Z_2^I,\\
  \varOmega(S^1)^\sK & \stackrel\simeq\longrightarrow
  \colim^{\scat{tmon}}_{I\in\sK}\Z^I.
\end{aligned}
\end{equation}
Since $\Z_2$ and $\Z$ are discrete groups, the colimit in
\cat{tmon} is the ordinary colimit of groups, and the two colimits
above have the following graph product presentations:
\begin{align*}
\colim^{\scat{grp}}_{I\in\sK}\Z_2^I&= F(g_1,\ldots,g_m)\big/
(g_i^2=0,\; g_ig_j=g_jg_i\text{ for }\{i,j\}\in\sK),\\
\colim^{\scat{grp}}_{I\in\sK}\Z^I&= F(g_1,\ldots,g_m)\big/
(g_ig_j=g_jg_i\text{ for }\{i,j\}\in\sK)
\end{align*}
where $F(g_1,\ldots,g_m)$ denotes a free group on $m$ generators.
The two groups above are known as the \emph{right-angled Coxeter
group}\label{raCoxeter} and the \emph{right-angled Artin group}
corresponding to the 1-skeleton (graph) of the flag complex~$\sK$.

Homotopy equivalences~\eqref{raCAcolim} imply that the polyhedral
power $(\R P^\infty)^\sK$ is the classifying space for the
right-angled Coxeter group $\colim^{\scat{grp}}_{I\in\sK}\Z_2^I$
and $(S^1)^\sK$ is the classifying space for the right-angled
Artin group $\colim^{\scat{grp}}_{I\in\sK}\Z^I$. The former result
is implicit in the work of Davis and Januszkiewicz~\cite{da-ja91},
and the latter is due to Kim and Roush~\cite{ki-ro80}. Note that
$(S^1)^\sK$ is a finite cell complex.
%(a subcomplex in the $m$-torus~$T^m$).
\end{remark}

\begin{proposition}\label{loopsps}
For a flag complex $\sK$, the Poincar\'e series of the Pontryagin
algebras $H_*(\varOmega(\C P^\infty)^\sK;\k)$ and
$H_*(\varOmega\zk;\k)$ are given by
\begin{align*}
  F\bigl( H_*(\varOmega(\C P^\infty)^\sK;\k);\lambda \bigr)&=\;
  \frac{(1+\lambda)^n}{1-h_1\lambda+\cdots+(-1)^nh_n\lambda^n}\;,\\
  F\bigl( H_*(\varOmega\zk;\k);\lambda \bigr)&=\;
  \frac1{(1+\lambda)^{m-n}(1-h_1\lambda+\cdots+(-1)^nh_n\lambda^n)}\;,
\end{align*}
where $(h_0,h_1,\ldots,h_n)$ is the $h$-vector of~$\sK$.
\end{proposition}
\begin{proof}
Since $H_*(\varOmega(\C P^\infty)^\sK;\k)$ is the quadratic dual
of $\k[\sK]$, the identity
\[
  F\bigl(\k[\sK];-\lambda\bigr)\cdot F\bigl(H_*(\varOmega(\C P^\infty)^\sK);\lambda\bigr)\;=\;1
\]
follows from Fr\"oberg~\cite[\S4]{frob75} (in the identity above
it is assumed that the generators of $\k[\sK]$ have degree one).
The Poincar\'e series of the face ring is given by
Theorem~\ref{psfr}, whence the first formula follows. The second
formula follows from exact sequence~\eqref{paseq}.
\end{proof}

The Poincar\'e series of $\pi_*(\varOmega(\C
P^\infty)^\sK)\otimes_\Z\Q$ (and therefore the ranks of homotopy
groups of $(\C P^\infty)^\sK$ and $\zk$) can also be calculated in
the flag case, although less explicitly.
See~\cite[\S4.2]{de-su07}.

\begin{example}
Let $\sK=\partial\varDelta^n$, so that $h_0=\cdots=h_n=1$. Then
$\varOmega(\C P^\infty)^\sK\simeq\varOmega S^{2n+1}\times
T^{n+1}$, and
\[
  F\bigl( H_*(\varOmega(\C P^\infty)^\sK);\lambda \bigr)\;=\;
  \frac{(1+\lambda)^{n+1}}{1-\lambda^{2n}}.
\]
On the other hand, Proposition~\ref{loopsps} gives
\[
  \frac{(1+\lambda)^n}{1-\lambda+\lambda^2+\cdots+(-1)^n\lambda^n}\;=\;
  \frac{(1+\lambda)^{n+1}}{1+(-1)^n\lambda^{n+1}}.
\]
The formulae agree if $n=1$, in which case $\sK$ is flag, but
differ otherwise.
\end{example}

We recall from Section~\ref{secmc} that a space $X$ is
\emph{coformal} when its Quillen model $Q(X)$ is weakly equivalent
to the rational homotopy Lie algebra $\pi_*(\varOmega
X)\otimes_{\Z}\Q$ as objects in $\cat{dgl}$. The space $(\C
P^\infty)^\sK$ is always formal by Theorem~\ref{xkform}, while its
coformality depends on~$\sK$:

\begin{theorem}\label{djcoform}
The space $(\C P^\infty)^\sK$ is coformal if and only if $\sK$ is
flag.
\end{theorem}
\begin{proof}
If $\sK$ is flag, Theorem~\ref{homegahocolim} together with
Corollary~\ref{flagcolim} provide an acyclic fibration
\[
  \varOmega_*\Q\langle\sK\rangle\cong
  \varOmega_*\colim^{\scat{dgc}}\Q\langle\,\cdot\:\rangle_\sK
  \stackrel\simeq\longrightarrow\colim^{\scat{dga}}\Lambda[\,\cdot\:]_\sK
\]
in $\cat{dga}_\Q$. Restricting to primitives yields a
quasi-isomorphism $e\colon
L_*\Q\langle\sK\rangle\stackrel\simeq\longrightarrow\pi_*(\varOmega(\C
P^\infty)^\sK)\otimes_\Z\Q$ in $\cat{dgl}$.

Now choose a minimal model $M_\sK\to\Q[\sK]$ for the face ring in
$\cat{cdga}_\Q$. Its graded dual $\Q\langle\sK\rangle\to C_\sK$ is
a minimal model for $\Q\langle\sK\rangle$ in $\cat{cdgc}$ (see~
\cite[\S 5]{ne-mi78}), so
$\varOmega_*\Q\langle\sK\rangle\to\varOmega_*C_\sK$ is a weak
equivalence in $\cat{dga}$. Restricting to primitives provides the
central map in the zigzag
\begin{equation}\label{cform}
  L_\sK\stackrel{\simeq}{\longrightarrow}L_*C_\sK
  \stackrel{\simeq\!}{\longleftarrow}L_*\Q\langle\sK\rangle
  \stackrel{e}{\longrightarrow}\pi_*(\varOmega(\C
  P^\infty)^\sK)\otimes_\Z\Q
\end{equation}
of quasi-isomorphisms in $\cat{dgl}$, where $L_\sK$ is a minimal
model for $(\C P^\infty)^\sK$ in $\cat{dgl}$ \cite[\S 8]{ne-mi78}.
Hence $(\C P^\infty)^\sK$ is coformal.

On the other hand, every missing face of $\sK$ with $>2$ vertices
determines a nontrivial higher Samelson bracket in
$\pi_*(\varOmega(\C P^\infty)^\sK)\otimes_\Z\Q$. The existence of
such brackets in $\pi_*(\varOmega X)\otimes_\Z\Q$ ensures that $X$
cannot be coformal, just as higher Massey products in $H^*(X;\Q)$
obstruct formality.
\end{proof}

Unlike the situation with the Pontryagin algebra $H_*(\varOmega(\C
P^\infty)^\sK)$, we are unable to describe the structure of its
commutator subalgebra $H_*( \varOmega\zk)$ completely even in the
flag case. However, the following result identifies a minimal set
of multiplicative generators as a specific set of iterated
commutators of the~$u_i$:

\begin{theorem}[{\cite[Theorem~4.3]{g-p-t-w}}]\label{multgen}
Assume that $\sK$ is flag and $\k$ is a field. The algebra $H_*(
\varOmega\zk;\k)$, viewed as the commutator subalgebra of
\eqref{tensorgraph} via exact sequence~\eqref{paseq}, is
multiplicatively generated by $\sum_{I\subset[m]}\dim\widetilde
H^0(\sK_I)$ iterated commutators of the form
\[
  [u_j,u_i],\quad [u_{k_1},[u_j,u_i]],\quad\ldots,\quad
  [u_{k_1},[u_{k_2},\cdots[u_{k_{m-2}},[u_j,u_i]]\cdots]]
\]
where $k_1<k_2<\cdots<k_p<j>i$, $k_s\ne i$ for any~$s$, and $i$ is
the smallest vertex in a connected component not containing~$j$ of
the subcomplex $\sK_{\{k_1,\ldots,k_p,j,i\}}$. Furthermore, this
multiplicative generating set is minimal, that is, the commutators
above form a basis in the submodule of indecomposables in~$H_*(
\varOmega\zk;\k)$.
\end{theorem}

\begin{remark}
To help clarify the statement of Theorem~\ref{multgen}, it is
useful to consider which brackets $[u_{j},u_{i}]$ are in the list
of multiplicative generators for $H_*(\varOmega\zk;\k)$. If
$\{j,i\}\in\sK$ then $i$ and $j$ are in the same connected
component of the subcomplex $\sK_{\{j,i\}}$, so $[u_{j},u_{i}]$ is
not a multiplicative generator. On the other hand, if
$\{j,i\}\notin\sK$ then the subcomplex $\sK_{\{j,i\}}$ consists of
the two distinct points $i$ and $j$, and $i$ is the smallest
vertex in its connected component of $\sK_{\{j,i\}}$ which does
not contain~$j$, so $[u_{j},u_{i}]$ is a multiplicative generator.

For a given $I=\{k_1,\ldots,k_p,j,i\}$, the number of the
commutators containing all $u_{k_1},\ldots,u_{k_p},u_j,u_i$ in the
set above is equal to $\dim\widetilde H^0(\sK_I)$ (one less the
number of connected components in $\sK_I$), so there are indeed
$\sum_{I\subset[m]}\dim\widetilde H^0(\sK_I)$ commutators in
total. More details are given in examples below.
\end{remark}

An important particular case of Theorem~\ref{multgen} corresponds
to $\sK$ consisting of $m$ disjoint points. This result may be of
independent algebraic interest, as it is an analogue of the
description of a basis in the commutator subalgebra of a free
algebra, given by Cohen and Neisendorfer~\cite{co-ne84}:

\begin{corollary}\sloppy
\label{CNzk} Let $A$ be the commutator subalgebra of the algebra
$T\langle u_1,\ldots,u_m\rangle/(u_i^2=0)$, that is, $A$ is the
algebra defined by the exact sequence
\[
  1\longrightarrow A\longrightarrow
  T\langle u_1,\ldots,u_m\rangle/
  (u_i^2=0)\longrightarrow\Lambda[u_1,\ldots,u_m]
  \longrightarrow1
\]
where $\deg u_i=1$. Then $A$ is a free associative algebra
minimally generated by the iterated commutators of the form
\[
  [u_j,u_i],\quad [u_{k_1},[u_j,u_i]],\quad\ldots,\quad
  [u_{k_1},[u_{k_2},\cdots[u_{k_{m-2}},[u_j,u_i]]\cdots]]
\]
where $k_1<k_2<\cdots<k_p<j>i$ and $k_s\ne i$ for any~$s$. Here,
the number of commutators of length $\ell$ is equal to
$(\ell-1)\bin m\ell$.
\end{corollary}

\begin{example}\label{exPApent}
Let $\sK$ be the boundary of pentagon, shown in Fig.~\ref{fbpent}.
\begin{figure}[h]
\unitlength=0.8mm
  \begin{center}
  \begin{picture}(40,37)
  \put(8.5,-3){\small 5}
  \put(-2,18){\small 1}
  \put(19.5,36){\small 2}
  \put(41,18){\small 3}
  \put(30,-3){\small 4}
  \put(0,20){\circle*{1}}
  \put(20,35){\circle*{1}}
  \put(40,20){\circle*{1}}
  \put(30,0){\circle*{1}}
  \put(10,0){\circle*{1}}
  \put(10,0){\line(1,0){20}}
  \put(10,0){\line(-1,2){10}}
  \put(0,20){\line(4,3){20}}
  \put(20,35){\line(4,-3){20}}
  \put(40,20){\line(-1,-2){10}}
  \end{picture}
  \end{center}
  \caption{Boundary of pentagon.}
  \label{fbpent}
\end{figure}
Theorem~\ref{multgen} gives the following 10 generators for the
algebra $H_*(\varOmega\zk)$:
\begin{gather*}
  a_1=[u_3,u_1],\quad a_2=[u_4,u_1],\quad a_3=[u_4,u_2],\quad
  a_4=[u_5,u_2],\quad a_5=[u_5,u_3],\\
  b_1=[u_4,[u_5,u_2]],\quad b_2=[u_3,[u_5,u_2]],\quad
  b_3=[u_1 ,[u_5,u_3]],\\
  b_4=[u_3,[u_4,u_1]],\quad b_5=[u_2,[u_4,u_1]],
\end{gather*}
where $\deg a_i=2$ and $\deg b_i=3$. In the notation of the
beginning of the previous section, $a_1$ is the Hurewicz image of
the Samelson product $[\mu_3,\mu_1]\colon S^2\to\varOmega(\C
P^\infty)^\sK$ lifted to $\varOmega\zk$, and $b_1$ is the Hurewicz
image of the iterated Samelson product
$[\mu_4,[\mu_5,\mu_2]]\colon S^3\to\varOmega(\C P^\infty)^\sK$
lifted to $\varOmega\zk$; the other $a_i$ and $b_i$ are described
similarly. We therefore have adjoint maps
\[
  \iota\colon(S^2\vee
  S^3)^{\vee5}\to\varOmega\zk\quad\text{and}\quad
  j\colon(S^3\vee S^4)^{\vee5}\to\zk
\]
corresponding to the wedge of all $a_i$ and~$b_i$. Now a
calculation using relations from Theorem~\ref{hldj} and the Jacobi
identity shows that $a_i$ and $b_i$ satisfy the relation
\begin{equation}\label{onerel}
  -[a_1,b_1]+[a_2,b_2]+[a_3,b_3]-[a_4,b_4]+[a_5,b_5]=0,
\end{equation}
where $[a_i,b_i]=a_ib_i-b_ia_i$. (One can make all the commutators
to enter the sum with positive signs by changing the order the
elements in the commutators defining $a_i,b_i$.) This relation has
a topological meaning. In general, suppose that $M$ and $N$ are
$d$-dimensional manifolds. Let $\overline{M}$ be the
$(d-1)$-skeleton of $M$, or equivalently, $\overline{M}$ is
obtained from $M$ by removing a disc in the interior of the
$d$-cell of $M$. Define $\overline{N}$ similarly. Suppose that
$f\colon S^{d-1}\to\overline{M}$ and $g\colon
S^{d-1}\to\overline{N}$ are the attaching maps for the top cells
in $M$ and $N$. Then the attaching map for the top cell in the
connected sum $M\cs N$\label{connectesu} is
$S^{d-1}\stackrel{f+g}\longrightarrow
\overline{M}\vee\overline{N}$. In our case, $S^{3}\times S^{4}$ is
a manifold and the attaching map $S^6\to S^{3}\vee S^{4}$ for its
top cell is the Whitehead product $[s_{1},s_{2}]_w$, where $s_{1}$
and $s_{2}$ respectively are the inclusions of $S^{3}$ and $S^{4}$
into $S^{3}\vee S^{4}$. The attaching map for the top cell of the
$5$-fold connected sum $(S^{3}\times S^{4})^{\cs 5}$ is therefore
the sum of five such Whitehead products. Composing it with $j$
into~$\zk$ and passing to the adjoint map we obtain
$\sum_{i=1}^{5}\pm[a_{i},b_{i}]$ (the signs depend on the
orientation chosen, see Construction~\ref{cobcs}).
By~(\ref{onerel}), this sum is null homotopic. Thus the inclusion
$j\colon(S^3\vee S^4)^{\vee5}\to\zk$ extends to a map
\[
  \widetilde j\colon(S^3\times S^4)^{\cs5}\to\zk.
\]
Furthermore, a calculation using Theorem~\ref{zkcoh} shows that
$\widetilde j$ induces an isomorphism in cohomology (see
Example~\ref{b5-gon}), that is, $\widetilde j$ is a homotopy
equivalence. Since both $(S^3\times S^4)^{\cs5}$ and $\zk$ are
manifolds, the complement of $(S^3\vee S^4)^{\vee5}$ in
$(S^3\times S^4)^{\cs5}$ and $\zk$ is a 7-disc, so that the
extension map $\widetilde j$ can be chosen to be one-to-one, which
implies that $\widetilde j$ is a homeomorphism.

We also obtain that $H_*(\varOmega\zk)$ is the quotient of a free
algebra on ten generators $a_i,b_i$ by relation~\eqref{onerel}.
Its Poicar\'e series is given by Proposition~\ref{loopsps}:
\[
  P\bigl( H_*( \varOmega\zk);\lambda \bigr)\;=\;
  \frac1{1-5\lambda^2-5\lambda^3+\lambda^5}\;.
\]
The summand $t^5$ in the denominator is what differs the
Poincar\'e series of the one-relator algebra $H_*(\varOmega\zk)$
from that of the free algebra $H_*(\varOmega(S^3\vee
S^4)^{\vee5})$.

A similar argument can be used to show that $\zk$ is homeomorphic
to a connected sum of sphere products when $\sK$ is a boundary of
a $m$-gon with $m\ge4$, therefore giving a homotopical proof of a
particular case of Theorem~\ref{zpstacked}. It would be
interesting to give a homotopical proof of this theorem in
general.
\end{example}

Another result of~\cite{g-p-t-w} identifies the class of flag
complexes $\sK$ for which $\zk$ has homotopy type of a wedge of
spheres.

We recall from Definition~\ref{defgolod} that $\k[\sK]$ is a
\emph{Golod ring}\label{Golodrin} and $\sK$ is a \emph{Golod
complex} when the multiplication and all higher Massey products in
$\Tor_{\k[v_1,\ldots,v_m]}(\k[\sK],\k)=H(\Lambda[u_1,\ldots,u_m]\otimes\k[\sK],d)$
are trivial.

We also need some terminology from graph theory. Let $\Gamma$ be a
graph on the vertex set~$[m]$. A \emph{clique} of $\Gamma$ is a
subset $I$ of vertices such that every two vertices in $I$ are
connected by an edge. Each flag complex $\sK$ is the \emph{clique
complex}\label{dclique} of its one-skeleton $\Gamma=\sK^1$, that
is, the simplicial complex formed by filling in each clique of
$\Gamma$ by a face.

A graph $\Gamma$ is called \emph{chordal}\label{dchordal} if each
of its cycles with $\ge 4$ vertices has a chord (an edge joining
two vertices that are not adjacent in the cycle). Equivalently, a
chordal graph is a graph with no induced cycles of length more
than three.

The following result gives an alternative characterisation of
chordal graphs.

\begin{theorem}[Fulkerson--Gross~\cite{fu-gr65}]
A graph is chordal if and only if its vertices can be ordered in
such a way that, for each vertex~$i$, the lesser neighbours of~$i$
form a clique.
\end{theorem}

Such an order of vertices is called a \emph{perfect elimination
order}\label{peorder}.

\begin{theorem}[\cite{g-p-t-w}]\label{flws}
Let $\sK$ be a flag complex and $\k$ a field. The following
conditions are equivalent:
\begin{itemize}
\item[(a)] $\k[\sK]$ is a Golod ring;
\item[(b)] the multiplication in $H^*(\zk;\k)$ is trivial;
\item[(c)] $\Gamma=\sK^1$ is a chordal graph;
\item[(d)] $\zk$ has homotopy type of a wedge of spheres.
\end{itemize}
\end{theorem}
\begin{proof}
(a)$\Rightarrow$(b) This is by definition of Golodness and
Theorem~\ref{zkcoh}.

(b)$\Rightarrow$(c) Assume that $\sK^1$ is not chordal, and choose
an induced chordless cycle $I$ with $|I|\ge4$. Then the full
subcomplex $\sK_I$ is the same cycle (the boundary of an
$|I|$-gon), and therefore $\mathcal Z_{\sK_I}$ is a connected sum
of sphere products. Hence, $H^*(\mathcal Z_{\sK_I})$ has
nontrivial products (this can be also seen directly by using
Theorem~\ref{zkadd}). Then, by Theorem~\ref{zkadd}, the same
nontrivial products appear in~$H^*(\zk)$.

(c)$\Rightarrow$(d) Assume that the vertices of $\sK$ are in
perfect elimination order. We assign to each vertex $i$ the clique
$I_i$ consisting of $i$ and the lesser neighbours of~$i$. Each
maximal face of $\sK$ (that is, each maximal clique of~$\sK^1$) is
obtained in this way, so we get an induced order on the maximal
faces: $I_{i_1},\ldots,I_{i_s}$. Then, for each~$k=1,\ldots,s$,
the simplicial complex $\bigcup_{j<k}I_{i_j}$ is flag (since it is
the full subcomplex $\sK_{\{1,2,\ldots,i_{k-1}\}}$ in a flag
complex). The intersection $\bigr(\bigcup_{j<k}I_{i_j}\bigl)\cap
I_{i_k}$ is a clique, so it is a face of $\bigcup_{j<k}I_{i_j}$.
Therefore, $\zk$ has homotopy type of a wedge of spheres by
Corollary~\ref{orderws}.

(d)$\Rightarrow$(a) This is by definition of the Golod property
and the fact that the cohomology of the wedge of spheres contains
only trivial Massey products.
\end{proof}

\begin{corollary}
Assume that $\sK$ is flag with $m$ vertices, and $\zk$ has
homotopy type of a wedge of spheres. Then
\begin{itemize}
\item[(a)]
the maximal dimension of spheres in the wedge is~$m+1$;
\item[(b)]
the number of spheres of dimension $\ell+1$ in the wedge is given
by $\sum_{|I|=\ell}\dim\widetilde H^0(\sK_I)$, for $2\le\ell\le
m$;
\item[(c)] $H^i(\sK_I)=0$ for $i>0$ and all~$I$.
\end{itemize}
\end{corollary}
\begin{proof}
If $\zk$ is a wedge of spheres, then $H_*(\varOmega\zk)$ is a free
algebra on generators described by Theorem~\ref{multgen}, which
implies (a) and (b). It also follows that
$H^*(\zk)\cong\bigoplus_{J\subset[m]}\widetilde H^0(\sK_J)$. On
the other hand, $H^*(\zk)\cong\bigoplus_{J\subset[m]}\widetilde
H^*(\sK_J)$ by Theorem~\ref{zkadd}, whence (c) follows.
\end{proof}

\begin{remark}
The equivalence of (a), (b) and (c) in Theorem~\ref{flws} was
proved in~\cite{be-jo07}. All the implications in the proof of
Theorem~\ref{flws} except (c)$\Rightarrow$(d) are valid for
arbitrary $\sK$, with the same arguments. However,
(c)$\Rightarrow$(d) fails in the non-flag case. Indeed, if $\sK$
be the triangulation of $\R P^2$ from Example~\ref{exbettinum}.4,
then $\sK^1$ is a complete graph, so it is chordal. However, $\zk$
is not homotopy equivalent to a wedge of spheres, because it has
2-torsion in homology. Furthermore, this $\sK$ is a Golod complex
by Exercise~\ref{rp2ng}. The following question is still open:
\end{remark}

\begin{problem}[\cite{g-p-t-w}]\label{GcoHs}
Assume that $H^*(\zk)$ has trivial multiplication, so that $\sK$
is Golod, over any field. Is it true that $\zk$ is a co-$H$-space,
or even a suspension, as in all known examples?
\end{problem}

\subsection*{Exercises}
\begin{exercise}
Shown that the generators $a_1,\ldots,a_5,b_1,\ldots,b_5$ of the
Pontryagin algebra $H^*(\varOmega\zk)$ from Example~\ref{exPApent}
satisfy relation~\eqref{onerel}.
\end{exercise}

\begin{exercise}[{\cite[Example~3.3]{g-p-t-w}}]
Describe explicitly the homotopy type of $\zk$ when $\sK$ is the
triangulation of $\R P^2$ from Example~\ref{exbettinum}.4.
\end{exercise}

%Уточнить ссылку из главы 7 на материал об эквивар. св. суммах (пока стоит ссылка на раздел~\ref{toricrep}
%Уточнить ссылку на главу 9 из упражнения об эквивар. почти компл. структурах в главе 7
%Ссылка на Гусейн-Заде
%Соотношение из док-ва Proposition 9.3.2. Спросить В.М.

\chapter{Torus actions and complex cobordism}\label{cobtoric}

Here we consider applications of toric methods in the theory of
complex cobordism. In particular, we describe new families of
toric generators of complex bordism ring and quasitoric
representatives in bordism classes. We also develop the theory of
torus-equivariant genera with applications to rigidity and fibre
multiplicativity problems, and provide explicit formulae for
bordism classes and genera of quasitoric manifolds and their
generalisations via localisation techniques.

We refer to Appendices~\ref{cobor} and~\ref{genera} for the
background material on complex (co)bordism and Hirzebruch genera.

As usual, when working with cobordism we assume all manifolds to
be smooth and compact. We denote by $[M]$ the bordism class in
$\varOmega^U_{2n}$ of $2n$-dimensional stably complex
manifold~$M$, and denote by $\langle M\rangle\in H_{2n}(M)$ its
fundamental homology class defined by the orientation arising from
the stably complex structure.

\section[Toric representatives in complex bordism
classes]{Toric and quasitoric representatives in
complex bordism classes}\label{toricrep}
Describing multiplicative generators for the complex bordism ring
$\varOmega^U$ and representing bordism classes by manifolds with
specific nice properties are well-known questions in cobordism
theory. For the application of toric methods, it is important to
represent complex bordism classes by manifolds with nicely
behaving torus actions preserving the stably complex structure. In
the context of oriented bordism, this question goes back to the
fundamental work of Conner--Floyd~\cite{co-fl64}.

The most `nicely behaving' torus actions that we have at our
disposal are the torus actions on toric and quasitoric manifolds.
Such an action does not exist on Milnor hypersurfaces $H_{ij}$,
which constitute the most well-known multiplicative generator set
for~$\varOmega^U$ (see Theorem~\ref{hijnt}). An alternative
multiplicative generator set for $\varOmega^U$ consisting of
projective toric manifolds~$B_{ij}$ was constructed by Buchstaber
and Ray in~\cite{bu-ra98r}. Each $B_{ij}$ is a complex
projectivisation of a sum of line bundles over the bounded flag
manifold $B_i$; in particular, $B_{ij}$ is a generalised Bott
manifold.

It seems likely that a \emph{minimal} set of ring generators of
$\varOmega^U$ can be found among toric manifolds, i.e. there exist
projective toric manifolds $X_i$ whose bordism classes $a_i=[X_i]$
are polynomial generators of the bordism ring:
$\varOmega^U\cong\Z[a_1,a_2,\ldots]$. A partial result in this
direction was obtained by Wilfong~\cite{wilf}. Nevertheless, not
any complex bordism class can be represented by a toric manifold.
The reason is that toric manifolds are very special algebraic
varieties, and there are many restrictions on their characteristic
numbers. For example, the Todd genus of a toric manifold is equal
to~1, which implies that the bordism class of a disjoint union of
toric manifolds cannot be represented by a toric manifold.

The main result of this section is Theorem~\ref{6.11} (originally
proved in~\cite{b-p-r07}), which shows that any complex bordism
class (in dimensions~$>2$) contains a quasitoric manifold. A
canonical torus-invariant stably complex structure is induced by
an omniorientation of a quasitoric manifold, see
Corollary~\ref{qtscs}. The above mentioned result
of~\cite{bu-ra98r} provides an additive basis for each bordism
group $\varOmega_{2n}^U$ represented by toric manifolds; it
implies that any complex bordism class can be represented by a
disjoint union of toric manifolds. The next step is to replace
disjoint unions by connected sums. There is the standard
construction of connected sum of $T$-manifolds at their fixed
points. (The connected sum of two stably complex manifolds $M_1$
and $M_2$ always admits a stably complex structure representing
the bordism class $[M_1]+[M_2]$.)
%, but this structure is not necessarily compatible with the torus action.)
Davis--Januzkiewicz~\cite{da-ja91} proposed to use this
construction to make the connected sum $M_1\cs M_2$ of two
quasitoric manifolds $M_1$ and $M_2$
%over polytopes $P_1$ and $P_2$, respectively,
into a quasitoric manifold over the connected sum $P_1\cs P_2$ of
quotient polytopes. However, the main difficulty here is that one
needs to keep track of both the torus action and the stably
complex structure on the connected sum of manifolds. It turns out
that the connected sum $M_1\cs M_2$ does not always admit an
omniorientation such that the bordism class $[M_1\cs
M_2]\in\varOmega^U$ of the induced \emph{$T^k$-invariant} stably
complex structure represents the sum $[M_1]+[M_2]$; this depends
on the sign pattern of fixed points of the manifolds.

In order to overcome the difficulty describe above, we replace
$M_2$ by a bordant quasitoric manifold $M_2'$ whose quotient
polytope is $I^n\cs P_2$ (a connected sum of $P_2$ with an
$n$-cube). Then we show that the connected sum $M_1\cs M_2'$
admits a stably complex structure which is invariant under the
torus action and represents the bordism class
$[M_1]+[M_2']=[M_1]+[M_2]$; the quotient polytope of $M_1\cs M_2'$
is $P_1\cs I^n\cs P_2$. This allows us to finish the proof of the
main result.

This result on quasitoric representatives can be viewed as an
answer to a toric version of the famous Hirzebruch question (see
Problem~\ref{hirz}) on bordism classes representable by connected
nonsingular algebraic varieties. Note that quasitoric manifolds
are connected by definition.

Using the constructions of Chapter~\ref{mamanifolds} we can
interpret this result as follows: each complex bordism class can
be represented by the quotient of a nonsingular complete
intersection of real quadrics by a free torus action.

\subsection*{Milnor hypersurfaces $H_{ij}$ are not quasitoric}
Milnor hypersurfaces
\[
   H_{ij}=\{
  (z_0:\cdots :z_i)\times (w_0:\cdots :w_j)\in
  \mathbb{C}P^i\times \mathbb{C}P^j\colon z_0w_0+\cdots +z_iw_i=0\}
\]
(corresponding to pairs of integers $j\ge i\ge0$) form the most
well-known set of multiplicative generators of the complex bordism
ring $\varOmega^U$, see Theorem~\ref{hijgen}. However, as it was
shown in~\cite{bu-ra98r}, the manifold $H_{ij}$ is not
(quasi)toric when $i>1$. We give the argument below.

\begin{construction}\label{msbun}
Let $\C^{i+1}\subset\C^{j+1}$ be the subspace generated by the
first $i+1$ vectors of the standard basis of~$\C^{j+1}$. We
identify $\C P^i$ with the set of lines $l\subset\C^{i+1}$. To
each line $l$ we assign the set of hyperplanes $W\subset\C^{j+1}$
containing~$l$. This set can be identified with~$\C P^{j-1}$.
Consider the set of pairs
\[
  E=\bigl\{(l,W)\colon \ l\subset W,
  \;l\subset\C^{i+1},\;W\subset\C^{j+1}\bigr\}.
\]
The projection $(l,W)\mapsto l$ defines a bundle $E\to\C P^i$ with
fibre~$\C P^{j-1}$.
\end{construction}

\begin{lemma}\label{hije}
Milnor hypersurface $H_{ij}$ is identified with the space $E$
above.
\end{lemma}
\begin{proof}
Indeed, a line $l\subset\C^{i+1}$ can be given by its generating
vector with homogeneous coordinates $(z_0:z_1:\cdots:z_i)$. A
hyperplane $W\subset\C^{j+1}$ is given by a linear form with
coefficients $w_0,w_1,\ldots,w_j$. The condition $l\subset W$ is
equivalent to the equation in the definition of~$H_{ij}$.
\end{proof}

\begin{theorem}\label{hijcoh}
The cohomology ring of $H_{ij}$ is given by
\[
  H^*(H_{ij})\cong\Z[u,v]\bigl/
  \bigl(u^{i+1},\;(u^i+u^{i-1}v+\cdots+uv^{i-1}+v^i)v^{j-i}\bigr),
\]
where $\deg u=\deg v=2$.
\end{theorem}
\begin{proof}
We use the notation from Construction~\ref{msbun}. Let $\zeta$
denotes the vector bundle over~$\C P^i$ whose fibre over $l\in\C
P^i$ is the $j$-plane $l^\bot\subset\C^{j+1}$. Then $H_{ij}$ is
identified with the projectivisation $\C P(\zeta)$. Indeed, for
any line $l'\subset l^\bot$ representing a point in the fibre of
the bundle $\C P(\zeta)$ over $l\in\C P^i$, the hyperplane
$W=(l')^\bot\subset\C^{j+1}$ contains~$l$, so that the pair
$(l,W)$ defines a point in $H_{ij}$ by Lemma~\ref{hije}. The rest
of the proof reproduces the general argument for the description
of the cohomology of a complex projectivisation (see
Theorem~\ref{cohomproj}).

Denote by $\eta$ the tautological line bundle over $\C P^i$ (its
fibre over $l\in\C P^i$ is the line~$l$). Then $\eta\oplus\zeta$
is a trivial $(j+1)$-plane bundle. Set $w=c_1(\bar\eta)\in H^2(\C
P^i)$ and consider the total Chern class
$c(\eta)=1+c_1(\eta)+c_2(\eta)+\cdots$. Since $c(\eta)c(\zeta)=1$
and $c(\eta)=1-w$, it follows that
\begin{equation}\label{cxi}
  c(\zeta)=1+w+\cdots+w^i.
\end{equation}

Consider the projection $p\colon \C P(\zeta)\to\C P^i$. Denote by
$\gamma$ the tautological line bundle over $\C P(\zeta)$, whose
fibre over $l'\in\C P(\zeta)$ is the line~$l'$. Let $\gamma^\bot$
denote the $(j-1)$-plane bundle over $\C P(\zeta)$ whose fibre
over $l'\subset l^\bot$ is the orthogonal complement to $l'$
in~$l^\bot$ (by definition, a point of $\C P(\zeta)$ is
represented by a line $l'$ in a fibre $l^\bot$ of bundle~$\zeta$).
It is easy to see that $p^*(\zeta)=\gamma\oplus\gamma^\bot$. We
set $v=c_1(\bar\gamma)\in H^2(\C P(\zeta))$ and $u=p^*(w)\in
H^2(\C P(\zeta))$. Then $u^{i+1}=0$. We have $c(\gamma)=1-v$ and
$c(p^*(\zeta))=c(\gamma)c(\gamma^\bot)$, hence
$$
  c(\gamma^\bot)=p^*\bigl(c(\zeta)\bigr)(1-v)^{-1}=
  (1+u+\cdots+u^i)(1+v+v^2+\cdots),
$$
see~\eqref{cxi}. Since $\gamma^\bot$ is a $(j-1)$-plane bundle, it
follows that $c_j(\gamma^\bot)=0$. Calculating the homogeneous
component of degree~$j$ in the identity above, we obtain the
second relation $v^{j-i}\sum_{k=0}^iu^kv^{i-k}=0$. It follows that
there is a homomorphism $\Z[u,v]\to H^*(\C P(\zeta))$ which
factors through a homomorphism $R\to H^*(\C P(\zeta))$, where $R$
is the quotient ring of $\Z[u,v]$ given in the theorem.

It remains to observe that $R\to H^*(\C P(\zeta))$ is actually an
isomorphism. This follows by considering the Serre spectral
sequence of the bundle~$p\colon \C P(\zeta)\to\C P^i$. Since both
$\C P^i$ and $\C P^{j-1}$ have only even-dimensional cells, the
spectral sequence collapses at~$E_2$. It follows that $\Z[u,v]\to
H^*(\C P(\zeta))$ is an epimorphism, and therefore so is $R\to
H^*(\C P(\zeta))$. Furthermore, the cohomology groups of $\C
P(\zeta)$ are the same as those of $\C P^i\times\C P^{j-1}$, which
implies that $R\to H^*(\C P(\zeta))$ is a monomorphism.
\end{proof}

\begin{theorem}\label{hijnt}
There is no torus action on the Milnor hypersurface $H_{ij}$ with
$i>1$ making it into a quasitoric manifold.
\end{theorem}
\begin{proof}
The cohomology of a quasitoric manifold has the form
$\Z[v_1,\ldots,v_m]/\mathcal I$, where $\mathcal I$ is the sum of
two ideals, one generated by square-free monomials and another
generated by linear forms (see Theorem~\ref{qtcoh}). We may assume
that the characteristic matrix $\varLambda$ has reduced
form~\eqref{lamat} and express the first $n$ generators
$v_1,\ldots,v_n$ via the last $m-n$ ones by means of linear
relations with integer coefficients. Therefore, we have
$$
  \Z[v_1,\ldots,v_m]/\mathcal I\cong
  \Z[w_1,\ldots,w_{m-n}]/\mathcal I',
$$
where $\mathcal I'$ is an ideal with basis consisting of products
of $\ge2$ integers linear forms.

Assume now that $H_{ij}$ is a quasitoric manifold. Then we
have
$$
  \Z[w_1,\ldots,w_{m-n}]/\mathcal I'\cong\Z[u,v]/\mathcal I'',
$$
where $\mathcal I''$ is the ideal from Theorem~\ref{hijcoh}.
Comparing the dimensions of linear (degree-two) components above,
we obtain $m-n=2$, so that $w_1,w_2$ can be identified with $u,v$
after a linear change of variables. Thus, the ideal $\mathcal I''$
has a basis consisting of polynomials which are decomposable into
linear factors over~$\Z$, which is impossible for~$i>1$.
\end{proof}

\begin{remark}
$H_{ij}$ is a projectivisation of a complex $j$-plane bundle over
$\C P^i$, but this bundle does not split into a sum of line
bundles, preventing $H_{ij}$ from carrying an effective
$T^{i+j-1}$-action. See the remark after Definition~\ref{genbott}
and Exercise~\ref{hijnoef}.
\end{remark}

\subsection*{Toric generators for the bordism ring $\varOmega^U$}
Here we describe, following~\cite{bu-ra98r} and~\cite{b-p-r07}, a
family of toric manifolds $\{B_{ij},\:0\le i\le j\}$ satisfying
the condition $s_{i+j-1}[B_{ij}]=s_{i+j-1}[H_{ij}]$, where
$s_n[M]$ denotes the characteristic number defined
by~\eqref{defsn}. This implies that the family $\{B_{ij}\}$
multiplicatively generates the complex bordism ring, by the same
argument as Theorem~\ref{hijgen}.

\begin{construction}
Given a pair of integers $0\le i\le j$, we introduce the manifold
$B_{ij}$ consisting of pairs $(\mathcal U,W)$, where
\[
  \mathcal U=\{U_1\subset U_2\subset\cdots\subset
  U_{i+1}=\C^{i+1},\quad \dim U_k=k\}
\]
is a bounded flag in $\C^{i+1}$ (that is, $U_k\supset\C^{k-1}$,
see Construction~\ref{bfm}) and $W$ is a hyperplane in $\C^{j+1}$
containing~$U_1$.
%line in $U_1^\bot\oplus\C^{j-i}$. (Here $U_1^\bot$ denotes the
%orthogonal complement to $U_1$ in $\C^{i+1}$, so that
%$U_1^\bot\oplus\C^{j-i}$ is the orthogonal complement to $U_1$
%in~$\C^{j+1}$.)
The projection $(\mathcal U,W)\mapsto\mathcal U$ describes
$B_{ij}$ as the projectivisation of a $j$-plane bundle over the
bounded flag manifold~$\BF_i$. This bundle splits into a sum of
line bundles:
\[
  B_{ij}=\C
  P(\rho^1_1\oplus\cdots\oplus\rho^1_i\oplus\underline{\C}^{j-i}),
\]
where $\rho^1_1,\ldots,\rho^1_i$ are the line bundles over $\BF_i$
described in Proposition~\ref{linebunbt} (the splitting follows
from the fact that $W$ can be identified with a line in
$U_1^\bot\oplus\C^{j-i}$ using the Hermitian scalar product
in~$\C^{j+1}$). Therefore, $B_{ij}$ is a generalised Bott manifold
(see Definition~\ref{genbott}). It follows that $B_{ij}$ is a
%complex $(i+j-1)$-dimensional
toric manifold, and also a quasitoric manifold over the
combinatorial polytope $I^i\times\varDelta^{j-1}$. The description
of the corresponding characteristic matrix and characteristic
submanifolds can be found in~\cite[Examples~2.9, 4.5]{bu-ra01}
or~\cite[Example~3.13]{b-p-r07}.
\end{construction}

\begin{proposition}\label{bicpi}
Let $f\colon\BF_i\to\C P^i$ be the map sending a bounded flag
$\mathcal U$ to its first line $U_1\subset\C^{i+1}$. Then the
bundle $B_{ij}\to\BF_i$ is induced from the bundle $H_{ij}\to\C
P^i$ by means of the map~$f$:
\[
\begin{CD}
  B_{ij} @>>> H_{ij}\\
  @VVV @VVV\\
  \BF_i @>f>> \C P^i
\end{CD}.
\]
\end{proposition}
\begin{proof}
%By viewing $\C^{j+1}$ as a Hermitian space, we can identify lines
%$W$ in $U_1^\bot\oplus\C^{j-i}$ with hyperplanes in $\C^{j+1}$
%containing~$U_1$. Now the statement
This follows from Lemma~\ref{hije}.
\end{proof}

\begin{theorem}\label{shbij}
We have $s_{i+j-1}[B_{ij}]=s_{i+j-1}[H_{ij}]$, where the
characteristic number $s_{i+j-1}$ of Milnor hypersurface~$H_{ij}$
is given by Lemma~{\rm\ref{shij}}.
\end{theorem}
\begin{proof}
We first prove a lemma:

\begin{lemma}\label{sidegmap}
Let $f\colon M\to N$ be a degree $d$ map of $2i$-dimensional
stably complex manifolds, and let $\xi$ be a complex $j$-plane
bundle over $N$, \ $j>1$. Then
\[
  s_{i+j-1}[\C P(f^*\xi)]=d\cdot s_{i+j-1}[\C P(\xi)].
\]
\end{lemma}
\begin{proof}
Let $p\colon \C P(\xi)\to N$ be the projection, $\gamma$ the
tautological bundle over $\C P(\xi)$, and $\gamma^\bot$ the
complementary bundle, so that $\gamma\oplus\gamma^\bot=p^*(\xi)$.
Then we have
\[
  \mathcal T(\C P(\xi))=p^*{\mathcal T}\!N\oplus\mathcal T_F(\C P(\xi)),
\]
where $\mathcal T_F(\C P(\xi))$ is the tangent bundle along the
fibres of the projection~$p$. Since $\mathcal T_F(\C
P(\xi))=\Hom(\gamma,\gamma^\bot)$ and
$\Hom(\gamma,\gamma)=\underline{\C}$, it follows that
\[
  \mathcal T_F(\C P(\xi))\oplus\underline{\C}=
  \Hom(\gamma,\gamma\oplus\gamma^\bot).
\]
Therefore,
\begin{multline}\label{projd}
\mathcal T(\C P(\xi))\oplus\underline{\C}
   =p^*{\mathcal T}\!N\oplus\Hom(\gamma,\gamma\oplus\gamma^\bot)=\\
   =p^*{\mathcal T}\!N\oplus\Hom(\gamma,p^*\xi)=
   p^*{\mathcal T}\!N\oplus(\bar\gamma\otimes p^*\xi),
\end{multline}
where $\bar\gamma=\Hom(\gamma,\underline{\C})$.

The map $f\colon M\to N$ induces a map $\widetilde f\colon\C
P(f^*\xi)\to\C P(\xi)$ such that
\begin{itemize}
\item[$\cdot$] $p\widetilde f=fp_1$, where $p_1\colon\C P(f^*\xi)\to M$ is the projection;

\item[$\cdot$] $\deg\widetilde f=\deg f$;

\item[$\cdot$] $\widetilde f^*\gamma$ is the tautological bundle over $\C P(f^*\xi)$.
\end{itemize}
Using~\eqref{projd}, we calculate
\[
  s_{i+j-1}\bigl(\mathcal T(\C P(\xi))\bigr)=
  p^*s_{i+j-1}({\mathcal T}\!N)+s_{i+j-1}(\bar\gamma\otimes
  p^*\xi)=s_{i+j-1}(\bar\gamma\otimes p^*\xi)
\]
(since $i+j-1>i$), and similarly for $\mathcal T(\C P(f^*\xi))$.
Thus,
\begin{multline*}
s_{i+j-1}[\C P(f^*\xi)]= s_{i+j-1}\bigl(\mathcal T(\C
P(f^*\xi))\bigr) \bigl\langle\C P(f^*\xi)\bigr\rangle\\ =
s_{i+j-1}\bigl((\widetilde f^*\bar\gamma)\otimes p_1^*f^*\xi\bigr)
\bigl\langle\C P(f^*\xi)\bigr\rangle = s_{i+j-1}\bigl(\widetilde
f^*(\bar\gamma\otimes p^*\xi)\bigr) \bigl\langle\C
P(f^*\xi)\bigr\rangle\\ = s_{i+j-1}(\bar\gamma\otimes p^*\xi)
\widetilde f_*\langle\C P(f^*\xi)\rangle = d\cdot
s_{i+j-1}(\bar\gamma\otimes p^*\xi) \langle\C P(\xi)\rangle
\\=d\cdot s_{i+j-1}[\C P(\xi)].\qedhere
\end{multline*}
\end{proof}

To finish the proof of Theorem~\ref{shbij} we note that the map
$f\colon \BF_i\to\C P^i$ from Proposition~\ref{bicpi} has
degree~1. (The map $f\colon \BF_i\to\C P^i$ is birational: it is
an isomorphism on the affine chart $V^{0,\ldots,0}=\{\mathcal U\in
\BF_i\colon U_1\not\subset \C^i\}$, because a bounded flag in
$V^{0,\ldots,0}\subset \BF_i$ is uniquely determined by its first
line~$U_1$.)
\end{proof}

\begin{theorem}[{\cite{bu-ra98r}}]\label{ducob}
The bordism classes of toric varieties $B_{ij}$, \ $0\le i\le j$,
multiplicatively generate the complex bordism ring
$\varOmega_*^U$. Therefore, every complex bordism class contains a
disjoint union of toric manifolds.
\end{theorem}
\begin{proof}
The first statement follows from Theorems~\ref{shbij}
and~\ref{hijgen}. A product of toric manifolds is toric, but a
disjoint union of toric manifolds is not, since toric manifolds
are connected by definition.
\end{proof}

\begin{remark}
The manifolds $H_{ij}$ and $B_{ij}$ are not bordant in general,
although $H_{0j}=B_{0j}=\C P^{j-1}$ and $H_{1j}=B_{1j}$. The proof
of Lemma~\ref{sidegmap} uses specific properties of the
number~$s_n$, and it does not work for arbitrary characteristic
numbers.
\end{remark}

\subsection*{Connected sums}\label{ssecs}
Our next goal is to replace the disjoint union of toric manifolds
by a version of connected sum, which will be a quasitoric
manifold.

\begin{construction}[Equivariant connected sum at fixed points]\label{equcs}
We give the construction for quasitoric manifolds only, although
it can be generalised easily to locally standard $T$-manifolds.
Let $M'=M(P',\varLambda'')$ and $M''=M(P'',\varLambda'')$ be two
quasitoric manifolds over $n$-polytopes $P'$ and $P''$,
respectively (see Section~\ref{qtman}). We assume that both
characteristic matrices $\varLambda'$ and $\varLambda''$ are in
the refined form~\eqref{lamat}, and denote by $x'$ and $x''$ the
initial vertices (given by the intersection of the first $n$
facets) of $P'$ and $P''$, respectively.

Consider the connected sum of polytopes $P'\cs
P''=P'\cs_{x',x''}P''$ (see Construction~\ref{consum}). By
definition, the \emph{equivariant connected sum} $M'\scs M''=
M'\scs_{x',x''}M''$ is the quasitoric manifold over $P'\cs P''$
with characteristic matrix
\begin{equation}\label{cmcs}
\varLambda_\#\;=\;\begin{pmatrix}
1&0&\cdots&0&\lambda'_{1,n+1}&\cdots&\lambda'_{1,m'}&\lambda''_{1,n+1}&
\cdots&\lambda''_{1,m''}\\
  0&1&\cdots&0&\lambda'_{2,n+1}&\cdots&\lambda'_{2,m'}&\lambda''_{2,n+1}&
\cdots&\lambda''_{2,m''}\\
  \vdots&\vdots&\ddots&\vdots&\vdots&\ddots&\vdots&\vdots&\ddots&\vdots\\
  0&0&\cdots&1&\lambda'_{n,n+1}&\cdots&\lambda'_{n,m'}&\lambda''_{n,n+1}&
\cdots&\lambda''_{n,m''}
\end{pmatrix}.
\end{equation}
Note that the matrix $\varLambda_\#$ is not refined, because the
first $n$ facets of $P'\cs P''$ do not intersect.

The manifold $M'\scs M''$ is $T^n$-equivariantly diffeomorphic to
the manifold obtained by removing from $M'$ and $M''$ invariant
neighbourhoods of the fixed points corresponding to $x'$ and $x''$
with subsequent $T^n$-equivariant identification of the boundaries
of these neighbourhoods. The latter manifold becomes the standard
connected sum $M'\cs M''$ (see Construction~\ref{cobcs}) if we
forget the action.
\end{construction}

In order to define an omniorientation (and therefore an invariant
stably complex structure) on $M'\scs M''$ we need to specify an
orientation of~$M'\scs M''$ along with matrix~\eqref{cmcs}.

Since both $M'$ and $M''$ are oriented, the quasitoric manifold
$M'\scs M''$ can be oriented so as to be oriented diffeomorphic
either to the oriented connected sum $M'\cs M''$, or to
$M'\cs\overline{M''}$ (see Construction~\ref{cobcs}). In the first
case we say that the orientation of $M'\scs M''$ is
\emph{compatible} with the orientations of $M'$ and~$M''$.

The existence of a compatible orientation on $M'\scs M''$ can be
detected from the combinatorial quasitoric pairs
$(P',\varLambda')$ and $(P'',\varLambda'')$. We recall the notion
of sign of a fixed point of a quasitoric manifold~$M$ (or a vertex
of the quotient polytope~$P$). By Lemma~\ref{qts}, the sign
$\sigma(x)$ of a vertex $x=F_{j_1}\cap\cdots\cap F_{j_n}$ measures
the difference between the orientations of $\mathcal T_x M$ and
$(\rho_{j_1}\oplus\cdots\oplus\rho_{j_n})|_x$. Therefore,
\[
  \sigma(x)=v_{j_1}\cdots v_{j_n}\langle M\rangle,
\]
where $v_i=c_1(\rho_i)\in H^2(M)$, $1\le i\le m$, are the ring
generators of~$H^*(M)$ and $\langle M\rangle\in H_{2n}(M)$ is the
fundamental homology class.

\begin{lemma}\label{cssign}
The equivariant connected sum $M'\scs_{x',x''}M''$ of omnioriented
quasitoric manifolds admits an orientation compatible with the
orientations of $M'$ and $M''$ if and only if
$\sigma(x')=-\sigma(x'')$.
\end{lemma}
\begin{proof}
Denote by $\rho'_j$, \ $1\le j\le m'$, the complex line
bundles~\eqref{rhoi} corresponding to the characteristic
submanifolds of~$M'$ (or to the facets of~$P'$), and similarly for
$\rho''_k$, \ $1\le k\le m''$, and $M''$. We also denote
\[
  c_1(\rho'_j)=v'_j,\quad c_1(\rho''_k)=v''_k,\quad 1\le j\le m',
  \quad 1\le k\le m''.
\]
The facets of the polytope $P'\cs P''$ are of three types: $n$
facets arising from the identifications of facets meeting at $x'$
and $x''$, $(m'-n)$ facets coming from $P'$, and $(m''-n)$ facets
coming from~$P''$. We denote the corresponding line bundles over
$M'\scs M''$ by $\xi_i$, $\xi'_j$ and $\xi''_k$, respectively
(they correspond to the columns of the characteristic
matrix~\eqref{cmcs}). Consider their first Chern classes in
$H^2(M'\scs M'')$:
\begin{align*}
  &w_i=c_1(\xi_i),&& w'_j=c_1(\xi'_j),&&w''_k=c_1(\xi''_k),\\
  &1\le i\le n,&&n+1\le j\le m',&&n+1\le k\le m''.
\end{align*}
Now consider the maps $p'\colon M'\scs M''\to M'$ and $p''\colon
M'\scs M''\to M''$ pinching one of the connected summands to a
point. We have $p'^{\,*}(\rho'_j)=\xi'_j$ for $n+1\le j\le m'$ and
$p''^{\,*}(\rho''_k)=\xi''_k$ for $n+1\le k\le m''$.
Relations~\eqref{viref} in the cohomology ring of $M'\scs M''$
take the form
\[
  w_i=-\lambda'_{i,n+1}w'_{n+1}-\cdots-\lambda'_{i,m'}w'_{m'}
  -\lambda''_{i,n+1}w''_{n+1}-\cdots-\lambda''_{i,m''}w''_{m''}.
\]
It follows that
\begin{equation}\label{wi}
  w_i=p'^{\,*}v'_i+p''^{\,*}v''_i,\quad 1\le i\le n.
\end{equation}
Since the first $n$ facets of $P'\cs P''$ do not intersect, it
follows that $w_1\cdots w_n=0$ in $H^{2n}(M'\scs M'')$, hence
\[
  (p'^{\,*}v'_1+p''^{\,*}v''_1)\cdots
  (p'^{\,*}v'_n+p''^{\,*}v''_n)
  =p'^{\,*}(v'_1\cdots v'_n)+p''^{\,*}(v''_1\cdots v''_n)=0.
\]
For any choice of an orientation with the corresponding
fundamental class $\langle M'\scs M''\rangle\in H_{2n}(M'\scs
M'')$, we obtain
\[
  v'_1\cdots v'_n\bigl(p'_*\langle M'\scs M''\rangle\bigr) +
  v''_1\cdots v''_n\bigl(p''_*\langle M'\scs M''\rangle\bigr)= 0.
\]
An orientation of $M'\scs M''$ is compatible with the orientations
of $M'$ and $M''$ if and only if $p'_*\langle M'\scs
M''\rangle=\langle M'\rangle$ and $p''_*\langle M'\scs
M''\rangle=\langle M''\rangle$. Substituting this to the identity
above, we obtain $\sigma(x')+\sigma(x'')=0$.
\end{proof}

\begin{proposition}\label{scscs}
Let $M'=M(P',\varLambda'')$ and $M''=M(P'',\varLambda'')$ be two
omnioriented quasitoric manifolds, and assume that
$\sigma(x')=-\sigma(x'')$. Then the stably complex structure
defined on the equivariant connected sum $M'\scs_{x',x''} M''$ by
the characteristic matrix~\eqref{cmcs} and the compatible
orientation is equivalent to the sum of the canonical stably
complex structures on $M'$ and~$M''$. In particular, the
corresponding complex bordism classes satisfy
\[
  [M'\scs M'']=[M']+[M''].
\]
\end{proposition}
\begin{proof}
The connected sum of the two canonical stably complex structures
on $M'$ and $M''$ is defined by the isomorphism
\begin{equation}\label{csscs}
  \mathcal T(M'\scs M'')\oplus\underline{\R}^{2(m'+m''-n)}
  \stackrel{\cong}\longrightarrow
  p'^{\,*}(\rho'_1\oplus\cdots\oplus\rho'_{m'})\oplus
  p''^{\,*}(\rho''_1\oplus\cdots\oplus\rho''_{m''})
\end{equation}
(see Construction~\ref{cobcs}). We have $p'^{\,*}(\rho'_j)=\xi'_j$
for $n+1\le j\le m'$ and $p''^{\,*}(\rho''_k)=\xi''_k$ for $n+1\le
k\le m''$. Furthermore, we claim that $p'^{\,*}(\rho'_i)\oplus
p''^{\,*}(\rho''_i)=\xi_i\oplus\underline{\C}$ for $1\le i\le n$.
Indeed, relations~\eqref{viref} for $M'$ imply
\[
  p'^{\,*}(\rho'_i)=(\xi'_{n+1})^{-\lambda'_{i,n+1}}\otimes\cdots
  \otimes(\xi'_{m'})^{-\lambda'_{i,m'}},\quad1\le i\le n,
\]
and the same relations for $M''$ imply
\[
  p''^{\,*}(\rho''_i)=(\xi''_{n+1})^{-\lambda''_{i,n+1}}\otimes\cdots
  \otimes(\xi''_{m''})^{-\lambda''_{i,m''}},\quad1\le i\le n.
\]
The line bundle $\xi'_j$ over $M'\scs M''$ has a section whose
zero set is precisely the characteristic submanifold corresponding
to the facet $F'_j\subset P'\cs P''$, for $n+1\le j\le m'$. There
is an analogous property of the bundles $\xi''_k$, for $n+1\le
k\le m''$. Now, since the facets $F'_j$ and $F''_k$ do not
intersect in $P'\cs P''$ for any $j$ and~$k$, the bundle
$p'^{\,*}(\rho'_i)\oplus p''^{\,*}(\rho''_i)$ has a nowhere
vanishing section, for $1\le i\le n$. Therefore,
$p'^{\,*}(\rho'_i)\oplus
p''^{\,*}(\rho''_i)=\eta\oplus\underline{\C}$ for some line
bundle~$\eta$. By comparing the first Chern classes and
using~\eqref{wi}, we obtain $\eta=\xi_i$, as needed.

Then the stably complex structure~\eqref{csscs} takes the form
\begin{multline*}
  \mathcal T(M'\scs M'')\oplus\underline{\R}^{2(m'+m''-n)}\stackrel\cong\longrightarrow\\
  \stackrel\cong\longrightarrow\xi_1\oplus\cdots\oplus\xi_n\oplus\xi'_{n+1}\oplus
    \cdots\oplus\xi'_{m'}
    \oplus\xi''_{n+1}\oplus\cdots\oplus\xi''_{m''}\oplus\underline{\C}^n.
\end{multline*}
This differs by a trivial summand~$\underline{\C}^n$ from the
stably complex structure defined by matrix~\eqref{cmcs}.
\end{proof}

\begin{corollary}
The complex bordism class $[M'\scs M'']$ does not depend on the
choice of initial vertices and the ordering of facets of $P'$ and
$P''$.
%(although the combinatorial type of $P'\cs P''$ and the
%topological type of $M'\scs M''$ may depend on these choices).
\end{corollary}

The relationship between the equivariant connected sum
$M'\scs_{x',x''} M''$ of omnioriented quasitoric manifolds and the
standard connected sum $M'\cs M''$ of oriented (or stably complex)
manifolds is now clear: the two operations produce the same
manifold if and only if $\sigma(x')=-\sigma(x'')$. Otherwise the
equivariant connected sum gives $M'\cs \overline{M''}$ or
$\overline{M'}\cs M''$ depending on the choice of orientation.
This implies that the equivariant connected sum cannot always be
used to obtain the sum of bordism classes. If the sign of every
vertex of~$P$ is positive, for example, then it is impossible to
obtain the bordism class $2[M]$ directly from $M\scs M$. This is
the case when $M$ is a toric manifold.

\subsection*{Proof of the main result}
We start with an example.
\begin{example}\label{prodtwos}
Consider the standard cube $I^n$ with the orientation induced
from~$\R^n$. The quasitoric manifold over $I^n$ corresponding to
the characteristic $n\times 2n$-matrix $(I|\,-\!I)$ (where $I$ is
the unit $n\times n$-matrix) is the product $(\C P^1)^n$ with the
standard complex structure. It represents a nontrivial complex
bordism class, and the signs of all vertices of the cube are
positive.

On the other hand, we can consider the omnioriented quasitoric
manifold over~$I^n$ corresponding to the matrix $(I|I)$. It is
easy to see that the corresponding stably complex structure on
$(\C P^1)^n\cong(S^2)^n$ is the product of $n$ copies of the
trivial structure on $\C P^1$ from Example~\ref{2cp1}. We denote
this omnioriented quasitoric manifold by~$S$. The bordism class
$[S]$ is zero, and the sign of a vertex
$(\varepsilon_1,\ldots,\varepsilon_n)\in I^n$ where
$\varepsilon_i=0$ or $1$ is given by
\[
  \sigma(\varepsilon_1,\ldots,\varepsilon_n)=
  (-1)^{\varepsilon_1}\cdots(-1)^{\varepsilon_n}.
\]
So adjacent vertices of $I^n$ have opposite signs.
\end{example}

We are now in a position to prove the next key lemma which
emphasises an important principle; however unsuitable a quasitoric
manifold $M$ may be for the formation of connected sums, a good
alternative representative always exists within the complex
bordism class~$[M]$.

\begin{lemma}\label{difsi}
Let $M$ be an omnioriented quasitoric manifold of dimension $>2$
over a polytope~$P$. Then there exists an omnioriented $M'$ over a
polytope~$P'$ such that $[M']=[M]$ and $P'$ has at least two
vertices of opposite sign.
\end{lemma}
\begin{proof}
Suppose that $x$ is the initial vertex of~$P$. Let $S$ be the
omnioriented product of 2-spheres of Example~\ref{prodtwos}, with
initial vertex~$w=(0,\ldots,0)$.

If $\sigma(x)=-1$, define $M'$ to be $S\scs_{w,x}M$ over
$P'=I^n\cs_{w,x}P$. Then $[M']=[M]$, because $S$ bounds. Moreover,
there is a pair of adjacent vertices of $I^n$ which survive under
formation of connected sum~$I^n\cs_{w,x}P$ (because~$n>1$). These
two vertices have opposite signs, as sought.

If $\sigma(x)=-1$, we make the same construction using the
opposite orientation of $I^n$ (and therefore of~$S$). Since $-S$
also bounds, the same conclusions hold.
\end{proof}

We may now complete the proof of the main result.
\begin{theorem}[\cite{b-p-r07}]\label{6.11}
In dimensions $>2$, every complex bordism class contains a
quasitoric manifold, necessarily connected, whose stably complex
structure is induced by an omniorientation, and is therefore
compatible with the torus action.
\end{theorem}
\begin{proof}
Consider bordism classes $[M_1]$ and $[M_2]$ in~$\varOmega_n^U$,
represented by omnioriented quasitoric manifolds over polytopes
$P_1$ and $P_2$ respectively. It then suffices to construct a
quasitoric manifold $M$ such that $[M]=[M_1]+[M_2]$, because
Theorem~\ref{ducob} gives an additive basis of~$\varOmega_n^U$
represented by quasitoric manifolds.

Firstly, we follow Lemma~\ref{difsi} and replace $M_2$ by $M_2'$
over $P_2'=\I^n\cs P_2$. Then we choose the initial vertex of
$P_2'$ so as to ensure that it has the opposite sign to the
initial vertex of~$P_1$, thereby guaranteeing the construction of
$M_1\scs M_2'$ over $P_1\cs P_2'$. The resulting omniorientation
defines the required bordism class, by Proposition~\ref{scscs} and
Lemma~\ref{difsi}.
\end{proof}

Combining Theorem~\ref{6.11} with Proposition~\ref{kact} and the
quadratic description~\eqref{zpqua} of the moment-angle
manifold~$\zp$ leads to another interesting conclusion relating
toric topology to complex cobordism:

\begin{theorem}\label{repquad}
Every complex bordism class may be represented by a stably complex
manifold obtained as the quotient of a free torus action on a
complete intersection of real quadrics.
\end{theorem}

One further deduction from Theorem~\ref{6.11} is the following
result of~Ray:

\begin{theorem}[{\cite{ray86}}]
Every complex bordism class contains a representative whose stable
tangent bundle is a sum of line bundles.
\end{theorem}

\subsection*{Examples}
Here we consider some examples of 4-dimensional quasitoric
manifolds (i.e.~$n=2$) illustrating the constructions of this
section.

\begin{example}
When $n=2$, the complex bordism class $[\C P^2]$ of the standard
complex structure of Example~\ref{cpnoo} is an additive generator
of the group $\varOmega_4^U\cong\Z^2$, with $c_2(\C P^2)=3$ and
all signs of the vertices of the quotient 2-simplex $\varDelta^2$
being positive.
\end{example}

The question then arises of representing the bordism class $2[\C
P^2]$ by an omnioriented quasitoric manifold~$M$. We cannot expect
to use $\C P^2\cs\C P^2$ for~$M$, because no vertices of sign $-1$
are available in~$\varDelta^2$, as required by Lemma~\ref{cssign}.
Moreover, $M$ must satisfy $c_2(M)=6$, by additivity, so the
quotient polytope $P$ has 6 or more vertices (see
Exercise~\ref{cnsign}). It follows that $P$ cannot be
$\varDelta^2\cs\varDelta^2$, which is a square! So we proceed to
appealing to Lemma~\ref{difsi}, and replace the second copy of $\C
P^2$ by the omnioriented quasitoric manifold $(-S)\scs\,\C P^2$
over $P'=I^2\cs\varDelta^2$. Of course $(-S)\scs\,\C P^2$ is
bordant to $\C P^2$, and $P'$ is a pentagon. These observations
lead naturally to our second example:

\begin{example}\label{cptboxcpt}
The quasitoric manifold $M=\C P^2\scs\,(-S)\scs\,\C P^2$
represents the bordism class $2[\C P^2]$, and its quotient is
polytope is $\varDelta^2\cs I^2\cs\varDelta^2$, which is a
hexagon. Fig.~\ref{figecs1} illustrates the procedure
diagrammatically, in terms of characteristic functions and
orientations. Every vertex of the hexagon has sign 1, so $M$
admits an equivariant almost complex structure by
Theorem~\ref{kusth}; in fact it coincides with the manifold from
Exercise~\ref{qtnt}.
\begin{figure}[h]
\begin{picture}(120,35)
  %lines:
  \put(0,5){\line(0,1){25}}
  \put(0,5){\line(5,3){21}}
  \put(0,30){\line(5,-3){21}}
  \put(40,5){\line(-1,1){12.5}}
  \put(40,5){\line(1,1){12.5}}
  \put(40,30){\line(-1,-1){12.5}}
  \put(40,30){\line(1,-1){12.5}}
  \put(80,5){\line(0,1){25}}
  \put(80,5){\line(-5,3){21}}
  \put(80,30){\line(-5,-3){21}}
  \put(95,10){\line(3,-2){12.5}}
  \put(95,10){\line(0,1){15}}
  \put(95,25){\line(3,2){12.5}}
  \put(120,10){\line(-3,-2){12.5}}
  \put(120,10){\line(0,1){15}}
  \put(120,25){\line(-3,2){12.5}}
  %orientations:
  \put(7.5,17.5){\oval(8,8)[b]}
  \put(7.5,17.5){\oval(8,8)[tr]}
  \put(8.5,21.5){\vector(-1,0){2}}
  \put(40,17.5){\oval(8,8)[b]}
  \put(40,17.5){\oval(8,8)[tr]}
  \put(41,21.5){\vector(-1,0){2}}
  \put(72.5,17.5){\oval(8,8)[b]}
  \put(72.5,17.5){\oval(8,8)[tr]}
  \put(73.5,21.5){\vector(-1,0){2}}
  \put(107.5,17.5){\oval(8,8)[b]}
  \put(107.5,17.5){\oval(8,8)[tr]}
  \put(108.5,21.5){\vector(-1,0){2}}
  %tuples:
%  \put(0.5,19){$\scriptscriptstyle(-1,-1)$}
  \put(0.5,18.5){$\scriptscriptstyle(-1,-1)$}
  \put(9.5,25){$\scriptscriptstyle(0,1)$}
  \put(9.5,9){$\scriptscriptstyle(1,0)$}
  \put(29,25){$\scriptscriptstyle(0,1)$}
  \put(29,9){$\scriptscriptstyle(1,0)$}
  \put(46,25){$\scriptscriptstyle(1,0)$}
  \put(46,9){$\scriptscriptstyle(0,1)$}
  \put(65,25){$\scriptscriptstyle(1,0)$}
  \put(65,9){$\scriptscriptstyle(0,1)$}
  \put(80.5,23){$\scriptscriptstyle(-1,-1)$}
  \put(85.5,12){$\scriptscriptstyle(-1,-1)$}
  \put(96.5,30.2){$\scriptscriptstyle(0,1)$}
  \put(113,30.2){$\scriptscriptstyle(1,0)$}
  \put(96.5,3.7){$\scriptscriptstyle(1,0)$}
  \put(113,3.7){$\scriptscriptstyle(0,1)$}
%  \put(111,21){$\scriptscriptstyle(-1,-1)$}
  \put(110.5,21.5){$\scriptscriptstyle(-1,-1)$}
  %operations:
  \put(23,17){$\scriptstyle\scs$}
  \put(55,17){$\scriptstyle\scs$}
  \put(86,17){$\scriptstyle=$}
\end{picture}
\caption{Equivariant connected sum $\C
P^2{\scriptstyle\mathbin{\widetilde\#}}\,(-S)\,
{\scriptstyle\mathbin{\widetilde\#}}\,\C P^2$.}
\label{figecs1}
\end{figure}
\end{example}

Our third example shows a related 4-dimensional situation in which
the connected sum of the quotient polytopes does support a
compatible orientation.

\begin{example}\label{cptmincpt}
Let $\overline{\C P^2\!}$ denote the quasitoric manifold obtained
by reverting the standard orientation of $\C P^2$ (equivalently,
reverting the orientation of the standard simplex~$\varDelta^2)$.
Every vertex has sign~$-1$, and we may construct $\C
P^2\scs\,\overline{\C P^2\!}$ as an omnioriented quasitoric
manifold over $\varDelta^2\cs\varDelta^2$. The corresponding
characteristic functions and orientations are shown in
Fig.~\ref{cp2cs}.
\begin{figure}[h]
\begin{picture}(120,35)
  \put(5,5){\line(0,1){30}}
  \put(5,5){\line(5,3){25}}
  \put(5,35){\line(5,-3){25}}
  \put(40,20){\line(5,3){25}}
  \put(40,20){\line(5,-3){25}}
  \put(65,35){\line(0,-1){30}}
  \put(85,5){\line(0,1){30}}
  \put(85,5){\line(1,0){30}}
  \put(85,35){\line(1,0){30}}
  \put(115,5){\line(0,1){30}}
  \put(15,20){\oval(10,10)[b]}
  \put(15,20){\oval(10,10)[tr]}
  \put(15.7,25){\vector(-1,0){2}}
  \put(55,20){\oval(10,10)[b]}
  \put(55,20){\oval(10,10)[tr]}
  \put(55.7,25){\vector(-1,0){2}}
  \put(100,20){\oval(10,10)[b]}
  \put(100,20){\oval(10,10)[tr]}
  \put(100.7,25){\vector(-1,0){2}}
  %
%  \put(-6,18){$\scriptstyle(-1,-1)$}
\put(5.3,21.5){$\scriptstyle(-1,-1)$}
  \put(16,29){$\scriptstyle(0,1)$}
  \put(16,9.6){$\scriptstyle(1,0)$}
  \put(47,29){$\scriptstyle(0,1)$}
  \put(47,10){$\scriptstyle(1,0)$}
  \put(65.4,28){$\scriptstyle(-1,-1)$}
  \put(74,10){$\scriptstyle(-1,-1)$}
  \put(98,2){$\scriptstyle(1,0)$}
  \put(98,36){$\scriptstyle(0,1)$}
%  \put(116,18){$\scriptstyle(-1,-1)$}
\put(104,28){$\scriptstyle(-1,-1)$}
  \put(33,19){$\scs$}
  \put(73,19){$=$}
\end{picture}
\vspace{-3mm}
\caption{Equivariant connected sum $\C
P^2{\scriptstyle\mathbin{\widetilde\#}}\,\overline{\C P^2\!}$.}
\label{cp2cs}
\end{figure}
Of course $[\overline{\C P^2\!}]=-[\C P^2]$. So $[\C
P^2]+[\overline{\C P^2\!}]=0$ in $\varOmega_4^U$, and the manifold
$\C P^2\scs\,\overline{\C P^2\!}$ bounds by
Proposition~\ref{scscs}.
\end{example}

A situation similar to that of Example~\ref{cptboxcpt} arises in
higher dimensions, when we consider the problem of representing
complex bordism classes by toric manifolds. For any such~$V$, the
top Chern number $c_n[V]$ coincides with the Euler characteristic,
and is therefore equal to the number of vertices of the quotient
polytope~$P$.

Suppose that toric manifolds $V_1$ and $V_2$ are of dimension
$\ge4$, and have quotient polytopes $P_1$ and $P_2$ respectively.
Then $c_n[V_1]=f_0(P_1)$ and $c_n[V_2]=f_0(P_2)$, yet $f_0(P_1\cs
P_2)=f_0(P_1)+f_0(P_2)-2$, where $f_0(\cdot)$ denotes the number
of vertices. Since $c_n([V_1]+[V_2])=c_n[V_1]+c_n[V_2]$, no
omnioriented quasitoric manifold over $P_1\cs P_2$ can represent
$[V_1]+[V_2]$ (see Exercise~\ref{cnsign}). This objection vanishes
for $P_1\cs I^n\cs P_2$, because it enjoys additional $2^n-2$
vertices with opposite signs.

\subsection*{Exercises}
\begin{exercise}\label{hijef}
Construct an effective action of a $j$-dimensional torus $T^j$ on
a Milnor hypersurface~$H_{ij}$ and a representation of $T^j$ in
$\C^{(i+1)(j+1)}$ such that the composition $H_{ij}\to\C
P^i\times\C P^j\to\C P^{(i+1)(j+1)-1}$ with the Segre embedding
becomes equivariant. Describe the fixed points of this action.
\end{exercise}

\begin{exercise}\label{hijnoef}
A torus $T^{i+j-1}$ cannot act effectively with isolated fixed
points on a Milnor hypersurface $H_{ij}$ with $i>1$. (Hint: use
Theorem~\ref{theo:stcoh} and other results of
Section~\ref{torusman}.)
\end{exercise}

\begin{exercise}
Show that the Milnor hypersurface $H_{11}$ is isomorphic (as a
complex manifold) to the Hirzebruch surface $F_1$ from
Example~\ref{hirzebruch}. In particular, $H_{11}$ is not
homeomorphic to~$F_0=\C P^1\times\C P^1$ (see
Exercise~\ref{hstop}).
\end{exercise}

\begin{exercise}
Show by comparing the cohomology rings that $H_{1j}$ is not
homeomorphic to $\C P^1\times \C P^{j-1}$. On the other hand, the
two manifolds are complex bordant by Exercise~\ref{h1jcob}.
\end{exercise}

\begin{exercise}
The procedure described in Example~\ref{cptboxcpt} allows one to
construct an almost complex quasitoric manifold $M$ (with all
signs positive) representing the sum of cobordism classes
$[V_1]+[V_2]$ of any two projective toric manifolds of real
dimension~4. Describe a similar procedure giving almost complex
quasitoric representatives for $[V_1]+[V_2]$ in dimension~6.
(Hint: modify the intermediate zero-cobordant manifold~$S$.) What
about higher dimensions?
\end{exercise}

\section{The universal toric genus}\label{utg}
A theory of equivariant genera for stably complex manifolds
equipped with compatible actions of a torus~$T^k$ was developed
in~\cite{b-p-r10}. This theory focuses on the notion of
\emph{universal toric genus} $\varPhi$, defined on stably complex
$T^k$-manifolds and taking values in the complex cobordism ring
$U^*(BT^k)$ of the classifying space. The construction of
$\varPhi$ goes back to the works of tom Dieck, Krichever and
L\"offler from the 1970s. The universal toric genus $\varPhi$ is
an equivariant analogue of the universal Hirzebruch genus
(Example~\ref{unigenus}) corresponding to the identity
homomorphism from the complex cobordism ring $\varOmega_U$ to
itself.

Here is an idea behind the construction of~$\varPhi$; details are
provided below. We start by defining a composite transformation of
$T^k$-equivariant cohomology functors
\begin{equation}\label{PhiX}
  \varPhi_X\colon U^*_{T^k}(X)\stackrel{\nu}\longrightarrow
  \MU^*_{T^k}(X)\stackrel{\alpha}\longrightarrow
  U^*(ET^k\times_{T^k}X).
\end{equation}
Here $U^*_{T^k}(X)$ (respectively $\MU^*_{T^k}(X)$) denotes the
\emph{geometric} (respectively the \emph{homotopic})
\emph{$T^k$-equivariant complex cobordism ring} of a
$T^k$-manifold~$X$, and $U^*(ET^k\times_{T^k}X)$ denotes the
ordinary complex cobordism of the Borel construction. (Note that
the geometric and homotopical versions of equivariant cobordism
are different, because of the lack of equivariant transversality.)

By restricting~\eqref{PhiX} to the case $X=\pt$ we get a
homomorphism of $\varOmega_U$-modules
\[
  \varPhi\colon \varOmega_{U:T^k}\to\varOmega_U[[u_1,\ldots,u_k]]
\]
from the geometric $T^k$-(co)bordism ring
$U^{*}_{T^k}(\pt)=\varOmega_{U:T^k}^{*}=\varOmega^{U:T^k}_{-*}$ to
the ring $U^{*}(BT^k)=\varOmega_U[[u_1,\ldots,u_k]]$.  Here $u_j$
is the cobordism Chern class $c_1^{U}(\bar\eta_j)$ of the
canonical line bundle (the conjugate of the Hopf bundle) over the
$j$th factor of~$BT^k=(\C P^\infty)^k$, for $1\le j\le k$. We
refer to $\varPhi$ as the \emph{universal toric
genus}\label{defutge}. It assigns to a bordism class
$[M,c_{\mathcal T}]\in\varOmega_{2n}^{U:T^k}$ of a
$2n$-dimensional stably complex $T^k$-manifold~$M$ the `cobordism
class' of the map $ET^k\times_{T^k}M\to BT^k$. The value
$\varPhi(M)$ is a power series in $u_1,\ldots,u_k$ with
coefficients in~$\varOmega_U$ and constant term~$[M]$.

\smallskip

We now proceed to providing the details of the construction. All
our $T^k$-spaces $X$ have homotopy type of cell complexes.
%are compactly generated and weakly Hausdorff~\cite{vogt71}, and
%underlie the model category of $T^k$-spaces and $T^k$-maps.
It is often important to take account of basepoints, in which case
we insist that they be fixed by~$T^k$. If $X$ itself does not have
a $T^k$-fixed point, then a disjoint fixed basepoint can be added;
the result is denoted by~$X_+$.

\subsection*{Homotopic equivariant cobordism}
The homotopic version of equivariant cobordism is defined via the
Thom $T^k$-spectrum $MU_{T^k}$, whose spaces are indexed by the
inclusion poset of complex representations $V$ of $T^k$ (of
complex dimension~$|V|$). Each $MU_{T^k}(V)$ is the Thom
$T^k$-space of the universal $|V|$-dimensional complex
$T^k$-equivariant vector bundle $\gamma_{V}$ over~$BU_{T^k}(V)$,
and each spectrum map $\varSigma^{2(|W|-|V|)}MU_{T^k}(V)\to
MU_{T^k}(W)$ is induced by the inclusion $V\subset W$ of a
$T^k$-submodule. The \emph{homotopic $T^k$-equivariant complex
cobordism group}\label{homocobor} $MU^n_{T^k}(X)$ of a pointed
$T^k$-space $X$ is defined by stabilising the pointed
$T^k$-homotopy sets:
\[
  MU^n_{T^k}(X)=\lim\limits_{\longrightarrow}\bigl[
  \varSigma^{2|V|-n}(X_+),MU_{T^k}(V)\bigr]_{T^k},
\]
The details of this construction can be found
in~\cite[Chapters~XXV--XXVIII]{may96}.

Applying the Borel construction to $\gamma_{V}$ yields a complex
$|V|$-dimensional bundle $ET^k\times_{T^k}\gamma_{V}$ over
$ET^k\times_{T^k}BU_{T^k}(V)$, whose Thom space is
$ET^k_+\wedge_{T^k}MU_{T^k}(V)$. The classifying map for the
bundle $ET^k\times_{T^k}\gamma_{V}$ induces a map of Thom spaces
$ET^k_+\wedge_{T^k}MU_{T^k}(V)\to MU(|V|)$.  Now consider a
$T^k$-map $\varSigma^{2|V|-n}(X_+)\to MU_{T^k}(V)$ representing a
homotopic cobordism class in $MU^n_{T^k}(X)$. By applying the
Borel construction and composing with the classifying map above,
we obtain a composite map of Thom spaces
\[
  \varSigma^{2|V|-n}(ET^k\times_{T^k} X)_+\longrightarrow
  ET^k_+\wedge_{T^k}MU_{T^k}(V)\to MU(|V|).
\]
This construction is homotopy invariant, so we get a map
\[
  \bigl[\varSigma^{2|V|-n}(X_+),MU_{T^k}(V)\bigr]_{T^k}
  \longrightarrow
  \bigl[\varSigma^{2|V|-n}(ET^k\times_{T^k} X)_+,MU(|V|)\bigr].
\]
Furthermore, it preserves stabilisation and therefore yields the
transformation
\[
  \alpha\colon MU^*_{T^k}(X)
  \longrightarrow U^*(ET^k\times_{T^k}X),
\]
which is multiplicative and preserves Thom classes
\cite[Proposition~1.2]{tomd70}.

The construction of $\alpha$ may be also interpreted using the
homomorphism $MU^*_{T^k}(X)\to MU^*_{T^k}(ET^k\times X)$ induced
by the $T^k$-projection $ET^k\times X\to X$; since $T^k$ acts
freely on $ET^k\times X$, the target may be replaced by
$U^*(ET^k\times_{T^k}X)$. Moreover, $\alpha$ is an isomorphism
whenever $X$ is compact and $T^k$ acts freely.

\begin{remark} According to the result of L\"offler~\cite[Chapter~XXVII]{may96}, $\alpha$ is the
homomorphism of completion with respect to the augmentation ideal
in $MU^*_{T^k}(X)$.
\end{remark}

\subsection*{Geometric equivariant cobordism}
The geometric version of equivariant cobordism can be defined
naturally by providing an equivariant version of Quillen's
geometric approach~\cite{quil71} to complex cobordism via
\emph{complex oriented maps} (see Construction~\ref{defgeomcob}).
However, this approach relies on normal complex structures,
whereas many of our examples present themselves most readily in
terms of tangential information. In the non-equivariant situation,
the two forms of data are, of course, interchangeable; but the
same does not hold equivariantly. This fact was often ignored in
early literature, and appears only to have been made explicit in
1995, by Comeza\~na~\cite[Chapter~XXVIII,~\S3]{may96}. As we shall
see below, tangential structures may be converted to normal, but
the procedure is not reversible.

We recall from Construction~\ref{defgeomcob} that elements in the
cobordism group $U^{-d}(X)$ of a manifold $X$ can be represented
by \emph{stably tangentially complex} bundles $\pi\colon E\to X$
with $d$-dimensional fibre~$F$, i.e. by those $\pi$ for which the
bundle $\mathcal T_F(E)$ of tangents along the fibre is equipped
with a stably complex structure $c_{\mathcal T}(\pi)$.
%Two such bundles $\pi$ and $\pi'$ are \emph{equivalent} when the
%stably complex structures $c_{\mathcal T}(\pi)$ and $c_{\mathcal
%T}(\pi')$ are equivalent in the standard sense (see
%Appendix~\ref{orcobbord}). Equivalence classes of bundles
%$\pi\colon E\to X$ and $\pi\colon E'\to X'$ are \emph{cobordant}
%when there exists a third bundle $\rho\colon L\to X$ whose
%boundary $\partial L\to X$ may be identified with $E\sqcup E'\to
%X$ in the standard fashion.

If $\pi$ is $T^k$-equivariant bundle, then it is stably
tangentially complex \emph{as a $T^k$-equivariant bundle} when
$c_{\mathcal T}(\pi)$ is also $T^k$-equivariant. The notions of
\emph{equivariant equivalence}  and \emph{equivariant cobordism}
apply to such bundles accordingly.

The \emph{geometric $T^k$-equivariant complex cobordism
group}\label{geomequcobor} $U^{-d}_{T^k}(X)$ consists of
equivariant cobordism classes of $d$-dimensional stably
tangentially complex $T^k$-equivariant bundles over~$X$. If
$X=\pt$, then we may identify both $F$ and $E$ with some
$d$-dimensional smooth $T^k$-manifold $M$, and $\mathcal T_F(E)$
with its tangent bundle $\mathcal T(M)$. So $c_{\mathcal T}(\pi)$
reduces to a $T^k$-equivariant stably tangentially complex
structure $c_\mathcal T$ on~$M$, and its cobordism class belongs
to the group $\varOmega^{-d}_{U:T^k}=U^{-d}_{T^k}(\pt)$ (the
cobordism group of point). The bordism group of point is given by
$\varOmega_d^{U:T^k}=\varOmega^{-d}_{U:T^k}$. The direct sums
$\varOmega^{U:T^k}=\varOmega^{U:T^k}_*=\bigoplus_d\varOmega^{U:T^k}_d$
and $\varOmega_{U:T^k}=\bigoplus_d\varOmega_{U:T^k}^d$ are the
\emph{geometric $T^k$-equivariant bordism} and \emph{cobordism
rings} respectively, and $U^*_{T^k}(X)$ is a graded
$\varOmega_{U:T^k}$--module under cartesian product. Furthermore,
$U^*_{T^k}(\,\cdot)$ is functorial with respect to pullback along
smooth $T^k$-maps $Y\to X$.

\begin{proposition}\label{cttcn}
Given any smooth compact $T^k$-manifold $X$, there are canonical
homomorphisms
\[
  \nu\colon U^{-d}_{T^k}(X)\to MU^{-d}_{T^k}(X), \quad d\ge0.
\]
\end{proposition}
\begin{proof}
The idea is to convert the tangential structure used in the
definition of geometric cobordism group~$U^{-d}_{T^k}(X)$ to the
normal structure required for the Pontryagin--Thom collapse map in
the homotopical approach.

Let $\pi$ in \smash{$\varOmega^{-d}_{U:T^k}(X)$} denote the
cobordism class of a stably tangentially complex $T^k$-bundle
$\pi\colon E\to X$. Choose a $T^k$-equivariant embedding $i\colon
E\to V$ into a complex $T^k$-representation space $V$ (see
Theorem~\ref{eqembth}) and consider the embedding $(\pi,i)\colon
E\to X\times V$. It is a map of vector bundles over $X$ which is
$T^k$-equivariant with respect to the diagonal action on $X\times
V$. There is an equivariant isomorphism
$c\colon\tau_F(E)\oplus\nu(\pi,i)\to\underline{V}= E\times V$ of
bundles over~$E$, where $\nu(\pi,i)$ is the normal bundle
of~$(\pi,i)$. Now combine $c$ with the stably complex structure
isomorphism $c_\mathcal T(\pi)\colon\mathcal
T_F(E)\oplus\underline{\R}^{2l-d}\to\xi$ to obtain an equivariant
isomorphism
\begin{equation}\label{stnocplxtk}
  \underline W\oplus\nu(\pi,i)\longrightarrow\xi^{\perp}\oplus
  \underline V\oplus\underline{\R}^{2l-d},
\end{equation}
where $\underline W=\xi^\perp\oplus\xi$ is a $T^k$-decomposition
for some complex representation space~$W$.

If $d$ is even, \eqref{stnocplxtk} determines a complex
$T^k$-structure on an equivariant stabilisation of $\nu(\pi,i)$;
if $d$ is odd, a further summand $\R$ must be added. For
notational convenience, assume the former, and write $d=2n$. We
compose~\eqref{stnocplxtk} with the classifying map
$\xi^{\perp}\oplus \underline
V\oplus\underline{\C}^{l-n}\to\gamma_R$, where $R$ is a
$T^k$-representation of complex dimension $|R|=|V|+|W|-n$ and then
pass to the Thom spaces to get a sequence of $T^k$-equivariant
maps
\[
  \varSigma^{2|W|}\Th(\nu(\pi,i))\longrightarrow \Th(\xi^\perp\oplus\underline V\oplus\R^{2l-d})
  \longrightarrow\Th(\gamma_R)=MU_{T^k}(R).
\]
We compose this with the Pontryagin--Thom collapse map on
$\nu(\pi,i)$ to obtain
\[
  f(\pi)\colon \varSigma^{2|V|+2|W|}X_+=\varSigma^{2|R|-d}X_+\longrightarrow MU_{T^k}(R).
\]
If $\pi$ and $\pi'$ are equivalent, then $f(\pi)$ and $f(\pi')$
differ only by suspension; if they are cobordant, then $f(\pi)$
and $f(\pi')$ are stably $T^k$-homotopic. So we define $\nu(\pi)$
to be the $T^k$-homotopy class of $f(\pi)$, as an element of
\smash{$MU^{-d}_{T^k}(X)$}.

The linearity of $\nu$ follows immediately from the fact that
addition in $U^{-d}_{T^k}(X)$ is induced by disjoint union.
\end{proof}

The proof of Proposition \ref{cttcn} also shows that $\nu$ factors
through the geometric cobordism group of stably normally complex
$T^k$-manifolds over~$X$.

\subsection*{The universal toric genus}
For any smooth compact $T^k$-manifold $X$, we define the
homomorphism
\[
  \varPhi_X\colon U^*_{T^k}(X)\stackrel{\nu}\longrightarrow
  \MU^*_{T^k}(X)\stackrel{\alpha}\longrightarrow
  U^*(ET^k\times_{T^k}X).
\]

\begin{definition}\label{definitutg}
The homomorphism
\[
  \varPhi\colon\varOmega_{U:T^k}\longrightarrow\varOmega_U[[u_1,\dots,u_k]]
\]
corresponding to the case $X=\pt$ above is called the
\emph{universal toric genus}.
\end{definition}

The genus $\varPhi$ is a multiplicative cobordism invariant of
stably complex $T^k$-manifolds, and it takes values in the ring
$\varOmega_U[[u_1,\dots,u_k]]$. As such it is an equivariant
extension of Hirzebruch's original notion of genus (see
Appendix~\ref{aphige}), and is closely related to the theory of
formal group laws. We explore this relation and study other
equivariant genera in the next sections.

\begin{remark}
By the result of Hanke~\cite{hank05} and L\"offler~\cite[(3.1)
Satz]{loff74}, when $X=\pt$, both homomorphisms $\nu$ and $\alpha$
are monic; therefore so is $\varPhi$.

On the other hand, there are two important reasons why $\nu$
cannot be epic. Firstly, it is defined on stably tangential
structures by converting them into stably normal information; this
procedure cannot be reversed equivariantly, because the former are
stabilised only by trivial representations of $T^k$, whereas the
latter are stabilised by arbitrary representations~$V$. Secondly,
homotopical equivariant cobordism groups are periodic, and each
$T^k$-representation $W$ gives rise to an invertible \emph{Euler
class}\label{eulerclassec} $e(W)$ in $MU^{2|W|}_{T^k}(\pt)$, while
$U^d_{T^k}(\pt)=\varOmega^d_{U:T^k}$ is zero for positive~$d$;
this phenomenon exemplifies the failure of equivariant
transversality.
\end{remark}

In geometric terms, the universal toric genus~$\varPhi$ assigns to
a geometric bordism class $[M,c_{\mathcal
T}]\in\varOmega_d^{U:T^k}$ of a $d$-dimensional stably complex
$T^k$-manifold~$M$ the `cobordism class' of the map
$ET^k\times_{T^k}M\to BT^k$. Since both $ET^k\times_{T^k}M$ and
$BT^k$ are infinite-dimensional, one needs to use their finite
approximations to define the cobordism class $\varPhi(M)$ purely
in terms of stably complex structures. Here is a conceptual way to
make this precise.

\begin{proposition}\label{utggysin}
Let $[M]\in \varOmega_{U:T^k}$ be a geometric equivariant
cobordism class represented by a $d$-dimensional $T^k$-manifold
$M$. Then
\[
  \varPhi(M)=(\id\times_{T^k}\pi)_{\,!}1,
\]
where
\[
  (\id\times_{T^k}\pi)_{\,!}\colon U^*(ET^k\times_{T^k}M)\longrightarrow
  U^{*-d}(ET^k\times_{T^k}\pt)=U^{*-d}(BT^k)
\]
is the Gysin homomorphism in cobordism induced by the projection
$\pi\colon M\to\pt$.
\end{proposition}
\begin{proof}
Choose an equivariant embedding $M\to V$ into a complex
representation space. Then the projection
$\id\times_{T^k}\pi\colon ET^k\times_{T^k}M\to BT^k$ factorises as
\[
  ET^k\times_{T^k}M\stackrel i\longrightarrow ET^k\times_{T^k}V
  \longrightarrow BT^k.
\]
We can approximate $ET^k$ by the products $(S^{2q+1})^k$ with the
diagonal action of~$T^k$. Then the above sequence is approximated
by the appropriate factorisations of smooth bundles over
finite-dimensional manifolds,
\begin{equation}\label{finiappr}
  (S^{2q+1})^k\times_{T^k}M\stackrel {i_q}\longrightarrow (S^{2q+1})^k\times_{T^k}V
  \longrightarrow (\C P^q)^k.
\end{equation}
This defines a complex orientation for the map
$(S^{2q+1})^k\times_{T^k}M\to (\C P^q)^k$ (see
Construction~\ref{defgeomcob}) and therefore a complex cobordism
class in $U^{-d}((\C P^q)^k)$, which we denote by~$\varPhi_q(M)$.
By definition of the Gysin homomorphism
(Construction~\ref{gysincob}),
$\varPhi_q(M)=(\id_q\times_{T^k}\pi)_{\,!}1$, where $\id_q$ is the
identity map of~$(S^{2q+1})^k$. As $q$ increases, the classes
$\varPhi_q(M)$ form an inverse system, whose limit is
$\varPhi(M)\in U^{-d}(BT^k)$. In particular,
$\varPhi(M)=(\id\times_{T^k}\pi)_{\,!}1$.
\end{proof}

\subsection*{Coefficients of the expansion of~$\varPhi(M)$}
Recall that the cobordism ring $U^*(\C P^\infty)$ is isomorphic to
$\varOmega_U[[u]]$, where $u=c_1^U(\bar\eta)\in U^2(\C P^\infty)$
is the generator represented geometrically by a codimension-one
complex projective subspace $\C P^{\infty-1}\subset\C P^\infty$
(viewed as the direct limit of inclusions $\C P^{N-1}\subset\C
P^N$).

A basis for the $\varOmega_U$-module
$\varOmega_U[[u_1,\dots,u_k]]=U^*(BT^k)$ in dimension $2|\omega|$
is given by the monomials $u^\omega=u_1^{\omega_1}\cdots
u_k^{\omega_k}$, where $\omega$ ranges over nonnegative integral
vectors $(\omega_1,\ldots,\omega_k)$, and
$|\omega|=\sum_j\omega_j$. A monomial $u^\omega$ is represented
geometrically by a $k$-fold product of complex projective
subspaces of codimension $(\omega_1,\ldots,\omega_k)$ in $(\C
P^\infty)^k$. If we write
\begin{equation}\label{defgomega}
  \varPhi(M)\;=\;\sum_\omega g_\omega(M)\,u^\omega
\end{equation}
in $\varOmega_U[[u_1,\dots,u_k]]$, then the coefficients
$g_\omega(M)$ lie in $\varOmega_U^{-2(|\omega|+n)}$. We shall
describe their representatives geometrically as universal
operations on~$M$. The cobordism class $g_\omega(M)$ will be
represented by the total space $G_\omega(M)$ of a bundle with
fibre $M$ over the product
$B_\omega=B_{\omega_1}\times\cdots\times B_{\omega_k}$, where each
$B_{\omega_i}$ is the \emph{bounded flag
manifold}\label{boundedfmcob} (Section~\ref{bfman}), albeit with
the stably complex structure representing zero
in~$\varOmega^U_{2\omega_i}$ and therefore not equivalent to the
standard complex structure on~$\BF_{\omega_i}$.

We start by describing stably complex structures on~$\BF_n$.
\begin{construction}
We denote by $\xi_n$ the `tautological' line bundle over the
bounded flag manifold $\BF_n$, whose fibre over $\mathcal
U\in\BF_n$ is the first space~$U_1$. We recall from
Proposition~\ref{charbfmbt} that $\BF_n$ has a structure of a Bott
tower in which each stage $B_k=\BF_k$ is the projectivisation $\C
P(\xi_{k-1}\oplus\underline\C)$. We shall denote the pullback of
$\xi_i$ to the top stage $\BF_n$ by the same symbol~$\xi_i$, so
that we have $n$ line bundles $\xi_1,\ldots,\xi_n$ over~$\BF_n$.

Using the matrix $A$ corresponding to the Bott manifold $\BF_n$
(described in Example~\ref{btbfn}), we identify $\BF_n$ with the
quotient of the product of 3-spheres
\begin{equation}\label{nsprodu}
  (S^3)^n=
  \{(z_1,\ldots,z_{2n})\in\C^{2n}\colon|z_k|^2+|z_{k+n}|^2=1,\ 1\le
  k\le n\}
\end{equation}
by the action of $T^n$ given by
\begin{equation}\label{tabfm1}
  (z_1,\ldots,z_{2n})\mapsto (t_1z_1,\:t_1^{-1}t_2z_2,\:\ldots,\:
  t_{n-1}^{-1}t_nz_n,\:t_1z_{n+1},\:t_2z_{n+2},\:\ldots,\:t_nz_{2n})
\end{equation}
(see the proof of Theorem~\ref{almat}). The manifold $\BF_n$ is a
complex algebraic variety (a toric manifold). The corresponding
stably complex structure can be described by viewing $\BF_n$ as a
quasitoric manifold and applying Theorem~\ref{taum}, which gives
\begin{equation}\label{stabbfm1}
  \mathcal T(\BF_n)\oplus\underline\C^n\cong\rho_1\oplus\cdots\oplus\rho_m,
\end{equation}
where the $\rho_i$ are the line bundles~\eqref{rhoi} corresponding
to characteristic submanifolds. We have $c_1(\rho_i)=v_i\in
H^2(\BF_n)$, \ $1\le i\le 2n$, the canonical ring generators of
the cohomology ring of $\BF_n$ viewed as a quasitoric manifold. On
the other hand, we have the cohomology ring generators
$u_k=c_1(\xi_k)$, \ $1\le k\le n$, for the Bott manifold~$\BF_n$.
The two sets are related by the identities $u_k=-v_{k+n}$ (see
Exercise~\ref{reluivi}). Then identities~\eqref{viref} in
$H^*(\BF_n)$ imply that $v_1=-u_1$ and $v_k=u_{k-1}-u_k$, or
equivalently, $\rho_1=\bar\xi_1$ and $\rho_k=\xi_{k-1}\bar\xi_k$
(where we dropped the sign of tensor product of line bundles), for
$2\le k\le n$. The stably complex structure~\eqref{stabbfm1}
therefore becomes
\begin{equation}\label{stabbfm2}
  \mathcal T(\BF_n)\oplus\underline\C^n\cong\bar\xi_1\oplus\xi_1\bar\xi_2\oplus\cdots\oplus
  \xi_{n-1}\bar\xi_n\oplus\bar\xi_1\oplus\bar\xi_2\oplus\cdots\oplus\bar\xi_n
\end{equation}
(note that when $n=1$ we obtain the standard isomorphism $\mathcal
T\C P^1\oplus\underline\C\cong\bar\eta\oplus\bar\eta$, as
$\xi_1=\eta$ is the tautological line bundle).

We shall change the stably complex structure on $\BF_n$ so that
the resulting bordism class in $\varOmega^U_{2n}$ will be zero. To
see that this is possible, we regard $\BF_n$ as a \emph{sphere
bundle} over $\BF_{n-1}$ rather that the complex projectivisation
$\C P(\xi_{n-1}\oplus\underline\C)$. If a stably complex structure
$c_{\mathcal T}$ on $\BF_n$ restricts to a trivial stably complex
structure on each fibre~$S^2$ (see Example~\ref{2cp1}), then
$c_{\mathcal T}$ extends over the associated 3-disk bundle, so it
is cobordant to zero.

So we need to change the stably complex structure~\eqref{stabbfm2}
so that the new structure restricts to a trivial one on each
fibre~$S^2$. We decompose $(S^3)^n$ as $(S^3)^{n-1}\times S^3$ and
$T^n$ as $T^{n-1}\times T^1$, then $T^1$ acts trivially on
$(S^3)^{n-1}$ by~\eqref{tabfm1}, and we obtain
\begin{equation}\label{bfndec}
  \BF_n=(S^3)^n/T^n=\bigl((S^3)^{n-1}\times(S^3/T^1)\bigr)/T^{n-1}=
  \BF_{n-1}\times_{T^{n-1}}S^2.
\end{equation}
Here $T^1$ acts on~$S^3$ diagonally, so the stably complex
structure on~$S^2$ is the standard structure of~$\C P^1$ (see
Example~\ref{2cp1qt}).

Now change the torus action~\eqref{tabfm1} to the following:
\begin{equation}\label{tabfm2}
  (z_1,\ldots,z_{2n})\mapsto (t_1z_1,\:t_1^{-1}t_2z_2,\:\ldots,\:
  t_{n-1}^{-1}t_nz_n,\:t^{-1}_1z_{n+1},\:t^{-1}_2z_{n+2},\:\ldots,\:t^{-1}_nz_{2n}).
\end{equation}
Then decomposition~\eqref{bfndec} is still valid, but now $T^1$
acts on $S^3$ antidiagonally, so the stably complex structure
on~$S^2$ is trivial, as needed. The resulting stably complex
structure on $\BF_n$ is given by the isomorphism
\begin{equation}\label{stabbfm3}
  \mathcal T(\BF_n)\oplus\underline\R^{2n}\cong\bar\xi_1\oplus\xi_1\bar\xi_2\oplus\cdots\oplus
  \xi_{n-1}\bar\xi_n\oplus\xi_1\oplus\xi_2\oplus\cdots\oplus\xi_n,
\end{equation}
and its cobordism class in~$\varOmega_{2n}^U$ is zero.
\end{construction}

\begin{definition}\label{raybasis}
We denote by $B_n$ the manifold~$\BF_n$ with zero-cobordant stably
complex structure~\eqref{stabbfm3}. The `tautological' line bundle
$\xi_n$ is classified by a map $B_n\to\C P^\infty$ and therefore
defines a bordism class $\beta_n\in U_{2n}(\C P^\infty)$. We also
set $\beta_0=1$. The set of bordism classes $\{\beta_n\colon
n\ge0\}$ is called \emph{Ray's basis} of the $\varOmega_U$-module
$U_*(\C P^\infty)$.
\end{definition}

\begin{proposition}[\cite{ray86}]
The bordism classes $\{\beta_n\colon n\ge0\}$ form a basis of the
free $\varOmega_U$-module $U_*(\C P^\infty)$ which is dual to the
basis $\{u^k\colon k\ge0\}$ of the $\varOmega_U$-module $U^*(\C
P^\infty)=\varOmega_U[[u]]$. Here $u=c^U_1(\bar\eta)$ is the
cobordism first Chern class of the canonical line bundle over~$\C
P^\infty$ (represented by $\C P^{\infty-1}\subset\C P^\infty$).
\end{proposition}
\begin{proof}
Since $[B_n]=0$ in $\varOmega_U$, the bordism class $\beta_n$ lies
in the reduced bordism module $\widetilde U_{2n}(\C P^\infty)$
for~$n>0$. Therefore, to show that $\{\beta_n\colon n\ge0\}$ and
$\{u^k\colon k\ge0\}$ are dual bases it is enough to verify the
property
\[
  u\frown\beta_n=\beta_{n-1}
\]
(the $\frown$-product is defined in
Construction~\ref{cobproducts}). The bordism class
$u\frown\beta_n$ is obtained by making the map $B_n\to\C P^\infty$
transverse to the zero section $\C P^{\infty-1}$ of $\C
P^\infty=MU(1)$ or, equivalently, by restricting the map $B_n\to\C
P^\infty$ to the zero set of a transverse section of the line
bundle $\xi_n$. As is clear from the proof of
Proposition~\ref{linebunbt}, the zero set of a transverse section
of $\xi_n=\rho^0_1$ is obtained by setting $z_1=0$
in~\eqref{nsprodu} and~\eqref{tabfm2}, which gives
precisely~$B_{n-1}$.
\end{proof}

For a nonnegative integer vector
$\omega=(\omega_1,\ldots,\omega_k)$, define the manifold
$B_\omega=B_{\omega_1}\times\cdots\times B_{\omega_k}$ and the
corresponding product bordism class $\beta_\omega\in
U_{2|\omega|}(BT^k)$.

\begin{corollary}\label{betabasis}
The set $\{\beta_{\omega}\}$ is a basis of the free
$\varOmega_U$-module $U_*(BT^k)$; this basis is dual to the basis
$\{u^\omega\}$ of~$U^*(BT^k)=\varOmega_U[[u_1,\ldots,u_k]]$.
\end{corollary}

\begin{definition}\label{Gomegadef}
Let $M$ be a tangentially stably complex $T^k$-manifold~$M$. Let
$T^\omega$ be the torus $T^{\omega_1}\times\cdots\times
T^{\omega_k}$ and $(S^3)^\omega$ the product
$(S^3)^{\omega_1}\times\cdots\times (S^3)^{\omega_k}$, on which
$T^\omega$ acts coordinatewise by~\eqref{tabfm2}. Define the
manifold
\[
  G_\omega(M)=(S^3)^\omega\times_{T^{\omega}}M,
\]
where $T^\omega$ acts on $M$ via the representation
\[
%\begin{equation}\label{tomegatk}
  (t_{1,1},\ldots,t_{1,\omega_1};\,\dots\,;t_{k,1},\ldots,t_{k,\omega_k})
  \;\longmapsto\;(t_{1,\omega_1},\dots,t_{k,\omega_k}).
%\end{equation}
\]
The stably complex structure on $G_\omega(M)$ is induced by the
tangential structures on the base and fibre of the bundle $M\to
G_\omega(M)\to B_\omega$.
\end{definition}
\begin{theorem}[\cite{b-p-r10}]\label{Gmrepsgm}
The manifold $G_\omega(M)$ represents the bordism class of the
coefficient $g_\omega(M)\in\varOmega_U^{-2(|\omega|+n)}$
of~\eqref{defgomega}. In particular, the constant term of
$\varPhi(M)\in\varOmega^U[[u_1,\ldots,u_k]]$ is
$[M]\in\varOmega^{-2n}_U$.
\end{theorem}
\begin{proof}
By Corollary~\ref{betabasis}, the coefficient $g_\omega(M)$ is
identified with the Kronecker product
$\langle\varPhi(M),\beta_\omega\rangle$ (see
Construction~\ref{cobproducts}). In terms of~\eqref{finiappr}, it
is represented on the pullback of the diagram
\[
  B_\omega\longrightarrow (\C P^q)^k\longleftarrow(S^{2q+1})^k\times_{T^k}M
\]
for suitably large $q$, and therefore on the pullback of the
diagram
\[
  B_\omega\longrightarrow BT^k\longleftarrow ET^k\times_{T^k}M
\]
of direct limits. The latter pullback is exactly $G_\omega(M)$.
\end{proof}

\begin{remark}
There is also a similar description of the coefficients in the
expansion of $\varPhi_X(\pi)$ for the transformation~\eqref{PhiX}
in the case when $U^*(ET^k\times_{T^k}X)$ is a finitely generated
free $U^*(BT^k)$-module, see~\cite[Theorem~3.15]{b-p-r10}.
\end{remark}

\subsection*{Exercises}
\begin{exercise}
Proposition~\ref{utggysin} can be generalised to the following
description of homomorphism~$\varPhi_X$ given by~\eqref{PhiX}. Let
$\pi\in U^{-d}_{T^k}(X)$ be a geometric cobordism class
represented by a $T^k$-equivariant bundle $\pi\colon E\to X$. Then
\[
  \varPhi_X(\pi)=(1\times_{T^k}\pi)_{\,!}1,
\]
where
\[
  (1\times_{T^k}\pi)_{\,!}\colon U^*(ET^k\times_{T^k}E)\longrightarrow
  U^{*-d}(ET^k\times_{T^k}X)
\]
is the Gysin homomorphism in cobordism.
\end{exercise}

\begin{exercise}
The original approach of~\cite{ray86} to define a zero-cobordant
stably tangent structure on~$\BF_n$ is as follows. One can
identify $\BF_n$ with the sphere bundle
$S(\xi_{n-1}\oplus\underline\R)$ (rather than with $\C
P(\xi_{n-1}\oplus\underline\C)$). Let $\pi_n\colon
S(\xi_{n-1}\oplus\underline\R)\to\BF_{n-1}$ be the projection;
show that the tangent bundle of $\BF_n$ satisfies
\[
  \mathcal T\BF_n\oplus\underline\R\cong\pi_n^*(\mathcal T\BF_{n-1}
  \oplus\xi_{n-1}\oplus\underline\R).
\]
By identifying $\pi_n^*\xi_{n-1}$ with $\xi_{n-1}$ (as bundles
over~$\BF_n$), we obtain inductively
\[
  \mathcal T\BF_n\oplus\underline\R\cong
  \xi_1\oplus\xi_2\oplus\cdots\oplus\xi_{n-1}\oplus\underline\R.
\]
Show that this stably complex structure is equivalent to that
of~\eqref{stabbfm3}. (Hint: calculate the total Chern classes.)
\end{exercise}

\begin{exercise}
There is a canonical $T^k$-action on $G_\omega(M)$ and $B_\omega$
making $G_\omega(M)\to B_\omega$ into a $T^k$-equivariant bundle.
\end{exercise}

\section{Equivariant genera, rigidity and fibre multiplicativity}
Recall that a \emph{genus}\label{degenus} is a multiplicative
cobordism invariant of stably complex manifolds, i.e. a ring
homomorphism $\varphi\colon\varOmega_U\to R$ to a commutative ring
with unit. We only consider genera taking values in torsion-free
rings~$R$; such $\varphi$ are uniquely determined by series
$f(x)=x+\cdots\in R\otimes\Q[[x]]$ via \emph{Hirzebruch's
correspondence} (Construction~\ref{hirzcor}).

Historically, equivariant extensions of genera were first
considered  by Atiyah and Hirzebruch \cite{at-hi70}, who
established the \emph{rigidity} property of the $\chi_y$-genus and
$\widehat{A}$-genus of~$S^1$-manifolds. The origins of these
concepts lie in the Atiyah--Bott fixed point formula
\cite{at-bo67}, which also acted as a catalyst for the development
of equivariant index theory. This development culminated in the
celebrated result of Taubes~\cite{taub89} establishing the
rigidity of the Ochanine--Witten \emph{elliptic
genus}\label{ellitpigenus} on spin $S^1$-manifolds.

Here we develop an approach to equivariant genera and rigidity
based solely on the complex cobordism theory. It allows us to
define an equivariant extension and the appropriate concept of
rigidity for an arbitrary Hirzebruch genus, and agrees with the
classical index-theoretical approach when the genus is the index
of an elliptic complex.

\subsection{Equivariant genera}
Our definition of an equivariant genus uses a universal
transformation of cohomology theories, studied
%in~\cite[\S6]{co-fl66} and~\cite{buch70}:
in~\cite{buch70}:

\begin{construction}[Chern--Dold character]\label{cherndold}
Consider multiplicative transformations of cohomology theories
\[
  h\colon U^*(X)\to H^*(X;\varOmega_U\otimes\Q).
\]
Such $h$ is determined uniquely by a series $h(u)\in \smash{H^2(\C
P^\infty;\varOmega_U\otimes\Q)}=\varOmega_U\otimes\Q[[x]]$, where
$u=c_1^U(\bar\eta)\in U^2(\C P^\infty)$ and $x=c^H_1(\bar\eta)\in
H^2(\C P^\infty)$.

The \emph{Chern--Dold character} is the unique multiplicative
transformation
\[
  \ch_U\colon U^*(X)\to H^*(X;\varOmega_U\otimes\Q)
\]
which reduces to the canonical inclusion
$\varOmega_U\to\varOmega_U\otimes\Q$ in the case $X=\pt$.

\begin{proposition}
The Chern--Dold character satisfies
\[
  \ch_U(u)=f_U(x),
\]
where $f_U(x)$ is the exponential of the formal group law $F_U$ of
geometric cobordisms.
\end{proposition}
\begin{proof}
Since $\ch_U$ acts identically on $\varOmega_U$, it follows that
\begin{equation}\label{chiufu}
  \ch_U F_U(v_1,v_2)=F_U(\ch_U(v_1),\ch_U(v_2))
\end{equation}
for any $v_1,v_2\in U^2(\C P^\infty)$. Let $f(x)$ denote the
series $\ch_U(u)\in\varOmega_U\otimes\Q[[x]]$. Let
$v_i=c_1^U(\xi_i)$ and $x_i=c_1^H(\xi_i)$ for $i=1,2$. Then
\begin{gather*}
  \ch_U F_U(v_1,v_2)=
  \ch_U\bigl(c_1^U(\xi_1\otimes\xi_2)\bigr)
  =f\bigl(c_1^H(\xi_1\otimes\xi_2)\bigr)=f(x_1+x_2),\\
  \ch_U(v_1)=f(x_1),\quad\ch_U(v_2)=f(x_2).
\end{gather*}
Substituting these expressions in~\eqref{chiufu} we get
$f(x_1+x_2)=F_U(f(x_1),f(x_2))$, which means that $f(x)$ is the
exponential of~$F_U$.
\end{proof}

Given a Hirzebruch genus $\varphi\colon\varOmega\to R\otimes\Q$
corresponding to $f(x)\in R\otimes\Q[[x]]$, we define a
multiplicative transformation
\[
  h_\varphi\colon U^*(X)\stackrel{\ch_U}\longrightarrow H^*(X;\varOmega_U\otimes\Q)
  \stackrel{\varphi}\longrightarrow H^*(X;R\otimes\Q)
\]
where the second homomorphism acts by~$\varphi$ on the
coefficients only. In the case $X=BT^k$ we obtain a homomorphism
\[
  h_{\varphi}\colon\varOmega_U[[u_1,\ldots,u_k]]\to
  R\otimes\Q[[x_1,\ldots,x_k]]
\]
which acts on the coefficients as $\varphi$ and sends $u_i$
to~$f(x_i)$ for $1\le i\le k$.
\end{construction}

\begin{definition}[equivariant genus]\label{defeqgenus}
%Let $\varphi\colon\varOmega_U\to R$ be a genus taking value in a
%torsion-free ring~$R$ with the corresponding $R$-series
%$f(x)=x+\cdots$. By analogy with the above construction we define
%a homomorphism
%\[
%  h_{\bar\varphi}\colon\varOmega_U[[u_1,\ldots,u_k]]\to R[[x_1,\ldots,x_k]]
%\]
%which acts on the coefficients as $\varphi\colon\varOmega_U\to R$
%and sends $u_i$ to~$f(x_i)$ for $1\le i\le k$.
%
The \emph{$T^k$-equivariant extension} of $\varphi$ is the ring
homomorphism
\[
  \varphi^T\colon\varOmega_{U:T^k}\to R\otimes\Q[[x_1,\ldots,x_k]]
\]
defined as the composition $h_{\varphi}\cdot\varPhi$ with the
universal toric genus.
\end{definition}

\begin{remark}
In the definition of equivariant genus from~\cite{b-p-r10}, the
generator $u_i$ was sent to $x_i$ instead of~$f(x_i)$. Although
this does not affect the notion of rigidity defined below,
Definition~\ref{defeqgenus} is a more natural extension of
Krichever's and Atiyah--Hirzebruch's approaches.
\end{remark}

%In particular, the $T^k$-equivariant extension of the universal
%genus $\varphi_U\colon\varOmega^U\to\varOmega^U$ of
%Example~\ref{unigenus} is~$\varPhi$; hence the name
%\emph{universal toric genus}.

\subsection*{Rigidity}
Krichever \cite{kric74}, \cite{kric90} considers rational valued
genera $\varphi\colon\varOmega_U\to\Q$, and equivariant extensions
$\varphi^K\colon\varOmega_{U:T^k}\to K^*(BT^k)\otimes\Q$. The
equivariant genus~$\varphi^K$ is related to ours $\varphi^T$ via a
natural transformation $U^*(X)\to K^*(X)\otimes\Q$ defined
by~$\varphi$, as explained below. If $\varphi(M)$ can be realised
as the index $\ind(\mathcal E)$ of an elliptic complex $\mathcal
E$ of complex vector bundles over~$M$ (see e.g.~\cite{h-b-j94}),
then for any $T^k$-manifold $M$ this index has a natural
$T^k$-equivariant extension $\ind^{T^k}(\mathcal E)$ which is an
element of the complex representation ring $R_U(T^k)$, and hence
in its completion~$K^0(BT^k)$. Krichever's interpretation of
rigidity is to require that $\varphi^K$ should lie the subring of
constants $\Q$ for every~$M$. In the case of an index, this
amounts to insisting that the corresponding $T^k$-representation
is always trivial, and therefore conforms to Atiyah and
Hirzebruch's original notion~\cite{at-hi70}.

The composition of $h_{\varphi}\colon U^{even}(X)\to
H^{even}(X;\Q)$ with the inverse of the Chern character
$\mathrm{ch}\colon K^0(X)\otimes\Q\stackrel\cong\to
H^{even}(X;\Q)$ gives the transformation
\[
  h^K_{\varphi}=\mathrm{ch}^{-1}\cdot h_{\varphi}\colon U^{even}(X)\to K^0(X)\otimes\Q,
\]
considered by Krichever in~\cite{kric74}. The transformation
$U^{odd}(X)\to K^1(X)\otimes\Q$ is defined similarly; together
they give a transformation $h^K_\varphi\colon U^*(X)\to
K^*(X)\otimes\Q$.

\begin{remark} In the case of the Todd genus
$\td\colon\varOmega_U\to\Z$ this construction gives the
Conner--Floyd transformation $U^*\to K^*$ described after
Example~\ref{xyexa}.
\end{remark}

Krichever~\cite{kric90} referred to a genus
$\varphi\colon\varOmega_U\to\Q$ as \emph{rigid} if the composition
\[
  \varphi^K\colon\varOmega_{U:T^k}\stackrel{\varPhi}\longrightarrow U^*(BT^k)
  \stackrel{h^K_\varphi}\longrightarrow K^0(BT^k)\otimes\Q
\]
belongs to the subring $\Q\subset K^0(BT^k)\otimes\Q$, i.e.
satisfies $\varphi^K(M)=\varphi(M)$ for any $[M]\in
\varOmega_{U:T^k}$. Definition~\ref{defeqgenus} of equivariant
genera based on the notion of the universal toric genus leads
naturally to the following version of rigidity, which subsumes the
approaches of Atiyah--Hirzebruch and Krichever:

\begin{definition}\label{defrigid}
A genus $\varphi\colon\varOmega_U\to R$ is \emph{$T^k$-rigid} on a
stably complex $T^k$-manifold $M$ whenever
$\varphi^T\colon\varOmega_{U:T^k}\to R\otimes\Q[[u_1,\ldots,u_k]]$
satisfies $\varphi^T(M)=\varphi(M)$; if this holds for every $M$,
then $\varphi$ is \emph{$T^k$-rigid}.
\end{definition}

Since the Chern character $\mathrm{ch}\colon K^0(X)\otimes\Q\to
H^{even}(X;\Q)$ is an isomorphism, a rational genus $\varphi$ is
$T^k$-rigid in the sense of Definition~{\rm\ref{defrigid}} if and
only if it is rigid in the sense of Krichever (and therefore in
the original index-theoretical sense of Atiyah--Hirzebruch if
$\varphi$ is an index).

%\begin{proposition}
%A genus $\varphi\colon \varOmega_U\to\Q$ is $T^k$-rigid  if and
%only if it satisfies $\varphi^K(M)=\varphi(M)$ for any $[M]\in
%\varOmega_{U:T^k}$.
%\end{proposition}
%\begin{proof}
%By definition, $\varphi$ is rigid (i.e. satisfies
%$\varphi^T(M)=\varphi(M)$) if and only if the image of $[M]$ under
%the composite map
%\[
%  \varPhi_{U:T^k}\stackrel{\varPhi}\longrightarrow
%  \varOmega_U[[u_1,\ldots,u_k]]\stackrel{\varphi}\longrightarrow
%  \Q[[u_1,\ldots,u_k]]
%\]
%is $\varphi(M)$ (a constant). On the other hand, because the Chern
%character is an isomorphism, $\varphi$ satisfies
%$\varphi^K(M)=\varphi(M)$ if and only if the image of $[M]$ under
%\[
%  \varPhi_{U:T^k}\stackrel{\varPhi}\longrightarrow
%  \varOmega_U[[u_1,\ldots,u_k]]\stackrel{h_{\bar\varphi}}\longrightarrow
%  \Q[[x_1,\ldots,x_k]]
%\]
%is a constant. The map $h_{\bar\varphi}$ is a homomorphism of
%$\varOmega_U$-modules sending each $[N]\in\varOmega_U$ to
%$\varphi(N)$ and sending $u_i$ to $f(x_i)$ by
%Proposition~\ref{prophphi}. Since $f(x)=x+\cdots$, it follows that
%$\varphi\cdot\varPhi[M]$ is a constant if and only if
%$h_{\bar\varphi}\cdot\varPhi[M]$ is a constant.
%\end{proof}

In Section~\ref{qtgenera} we shall describe how toric methods can
be applied to establish the rigidity property for several
fundamental Hirzebruch genera. Now we consider another important
property of genera.

\subsection*{Fibre multiplicativity}
The following definition extends that of
Hirzebruch~\cite[Chapter~4]{h-b-j94} for the oriented case. It
applies to fibre bundles of the form $M\to E\times_G
M\stackrel{\pi}{\to} B$, where $M$ and $B$ are closed, connected
and stably tangentially complex, $G$ is a compact Lie group of
positive rank whose action preserves the stably complex structure
on~$M$, and $E\to B$ is a principal $G$-bundle. In these
circumstances, the bundle $\pi$ is stably tangentially complex,
and $N=E\times_G M$ inherits a canonical stably complex structure.

\begin{definition}\label{deffibmul}
A genus $\varphi\colon\varOmega_U\to R$ is \emph{fibre
multiplicative with respect to} the stably complex manifold~$M$
whenever $\varphi(N)=\varphi(M)\varphi(B)$ for any such bundle
$\pi$ with fibre~$M$; if this holds for every $M$, then $\varphi$
is \emph{fibre multiplicative}.
\end{definition}

For rational genera in the oriented category,
Ochanine~\cite[Proposition~1]{ocha87} proved that rigidity is
equivalent to fibre multiplicativity (see
also~\cite[Chapter~4]{h-b-j94}). In the toric case, we have the
following stably complex analogue, whose conclusions are integral.
It refers to bundles $E\times_GM\stackrel{\pi}{\longrightarrow}B$
of the form required by Definition~\ref{deffibmul}, where $G$ has
maximal torus $T^k$ with $k\ge1$.

\begin{theorem}[\cite{b-p-r10}]\label{multandrig}
If the genus $\varphi$ is $T^k$-rigid on~$M$, then it is fibre
multiplicative with respect to~$M$ for bundles whose structure
group $G$ has the property that $U^*(BG)$ is torsion-free.

On the other hand, if $\varphi$ is fibre multiplicative with
respect to a stably tangentially complex $T^k$-manifold~$M$, then
it is $T^k$-rigid on~$M$.
\end{theorem}
\begin{proof}
Let $\varphi$ be $T^k$-rigid on~$M$, and consider the pullback
squares
\[
%begin{equation}\label{pulsqs}
\begin{CD}
  E\times_GM@>f'>>EG\times_GM@<i'<<ET^k\times_{T^k}M\\
  @V\pi VV@V{\pi^G}VV@V{\pi^{T^k}}VV\\
  B@>f>>BG@<i<<BT^k,
\end{CD}
\]
%end{equation}
where $\pi^G$ is universal, $i$ is induced by inclusion, and $f$
classifies $\pi$. By Proposition~\ref{utggysin} and commutativity
of the right square, $\pi^G_{\,!}1=[M]\cdot 1+\beta$, where $1\in
U^0(EG\times_GM)$ and $\beta\in\widetilde U^{-2n}(BG)$. By
commutativity of the left square, $\pi_{\,!}1=[M]\cdot1+f^*\beta$
in $U^{-2n}(B)$. Applying the Gysin homomorphism associated with
the augmentation map $\varepsilon^B\colon B\to\pt$ yields
\begin{equation}\label{mult}
[E\times_GM]\;=\;\varepsilon^B_{\,!}\pi_{\,!}1\;=\;
[M][B]+\varepsilon^B_{\,!}f^*\beta
\end{equation}
in $\varOmega^U_{2(n+b)}$, where $\dim B=2b$; so
$\varphi(E\times_GM)=\varphi(M)\varphi(B)+
\varphi(\varepsilon^B_{\,!}f^*\beta)$. Moreover,
$i^*\beta=\sum_{|\omega|>0}g_\omega(M)u^\omega$ in
$U^{-2n}(BT^k)$, so $\varphi(i^*\beta)=0$ because $\varphi$ is
$T^k$-rigid. The assumptions on $G$ ensure that $i^*$ is
injective, which implies that $\varphi(\beta)=0$ in
$U^*(BG)\otimes_\varphi R$. Fibre multiplicativity then follows
from \eqref{mult}.

Conversely, suppose that $\varphi$ is fibre multiplicative with
respect to~$M$, and consider the manifold $G_\omega(M)$ of
Theorem~\ref{Gmrepsgm}. By Definition \ref{Gomegadef}, it is the
total space of the bundle $(S^3)^\omega\times_{T^\omega}\!M\to
B_\omega$,
%$((S^3)^\omega\times_{T^\omega}\!T^k)\times_{T^k}\!M\to B_\omega$,
which has structure group $T^k$; therefore
$\varphi(G_\omega(M))=0$, because $B_\omega$ bounds for every
$|\omega|>0$. So $\varphi$ is $T^k$-rigid on~$M$.
\end{proof}

\begin{remark}
We may define $\varphi$ to be \emph{$G$-rigid} when
$\varphi(\beta)=0$, as in the proof of Theorem~\ref{multandrig}.
It follows that $T$-rigidity implies $G$-rigidity for any $G$ such
that $\varOmega^*_U(BG)$ is torsion-free.
\end{remark}

\begin{example}\label{rigidex}
The signature (or the $L$-genus, see Example~\ref{xyexa}.2) is
fibre multiplicative over any simply connected
base~\cite[Chapter~4]{h-b-j94}, and so is rigid.
\end{example}

\section{Isolated fixed points: localisation formulae}
In this section we focus on stably tangentially complex
$T^k$-manifolds $(M^{2n},c_{\mathcal T})$ for which the fixed
points $p$ are isolated; in other words, the fixed point set $M^T$
is finite. We proceed by deducing a localisation formula for
$\varPhi(M)$ in terms of fixed point data. We give several
illustrative examples, and describe the consequences for certain
non-equivariant genera and their $T^k$-equivariant extensions.

Localisation theorems in equivariant generalised cohomology
theories appear in the works of tom Dieck \cite{tomd70}, Quillen
\cite{quil71}, Krichever~\cite{kric74}, Kawakubo~\cite{kawa80},
and elsewhere. We prove our Theorem~\ref{evutg} by interpreting
their results in the case of isolated fixed points, and
identifying the signs explicitly.

Each integer vector $\mb w=(w_1,\ldots,w_k)\in\Z^k$ determines a
line bundle
\[
  \bar\eta^{\mb w}=
  \bar\eta_1^{w_1}\otimes\cdots\otimes\bar\eta_k^{w_k}
\]
over $BT^k=(\C P^\infty)^k$, where $\bar\eta_j$ is the canonical
line bundle over the $j$th factor. Let
\[
  [\mb w,\mb u]=c_1^U(\bar\eta^{\mb w})
\]
denote the cobordism first Chern class of $\bar\eta^{\mb w}$. It
is given by the power series
\[
  F_U(\underbrace{u_1,\ldots,u_1}_{w_1},\ldots,\underbrace{u_k,\ldots,u_k}_{w_k})
  \in U^2(BT^k),
\]
where $F_U(u_1,\ldots,u_k)$ is the iterated substitution
$F_U(\cdots F_U(F_U(u_1,u_2),u_3),\ldots,u_k)$ in the formal group
law of geometric cobordisms, see Section~\ref{secfglgc}.
%(the associativity implies that the result does not depend on bracketing.
Modulo decomposables we have that
\begin{equation}\label{wjxmoddec}
  [\mb w,\mb u]\;\equiv\;w_1u_1+\cdots +w_ku_k,
\end{equation}
and it is convenient to rewrite the right hand side as a scalar
product $\langle \mb w,\mb u\rangle$.

Let $p$ be an isolated fixed point for the $T^k$-action on~$M$. We
recall from Section~\ref{secinvscs} that the weights $\mb
w_j(p)\in\Z^k$ and the sign $\sigma(p)=\pm1$ are defined, and
refer to $\{\mb w_j(p),\;\sigma(p)\colon1\le j\le n,\;p\in M^T\}$
as the \emph{fixed point data}\label{fpdata} of~$(M,c_{\mathcal
T})$.

\begin{theorem}[localisation formula]\label{evutg}
For any stably tangentially complex $2n$-dimensional
$T^k$-manifold $M$ with isolated fixed points~$M^T$, the equation
\begin{equation}\label{lfx}
  \varPhi(M)\;=\;\sum_{p\in M^T}\!
  \sigma(p) \prod_{j=1}^n\frac1{[\mb w_j(p),\mb u\:]}
\end{equation}
is satisfied in $U^{-2n}(BT^k)$.
\end{theorem}

\begin{remark}
The summands on the right hand side of~\eqref{lfx} formally belong
to the localised ring $S^{-1}U^*(BT^k)$ where $S$ is the set of
equivariant Euler classes of nontrivial representation of~$T^k$.
\end{remark}

\begin{proof}
Choose an equivariant embedding $i\colon M\to V$ into a complex
$N$-dimensional representation space~$V$ and consider the
commutative diagram
\begin{equation}\label{fpCD}
\begin{CD}
  M^T @>r_M>> M\\
  @Vi^TVV @VViV\\
  V^T @>r_V>> V
\end{CD}
\end{equation}
where $V^T\subset V$ is the $T^k$-fixed subspace, $r_M$ and $r_V$
denote the inclusions of fixed points, and $i^T$ is the
restriction of $i$ to~$M^T$.

We restrict~\eqref{fpCD} to a tubular neighbourhood of $M^T$ in
$V$, which can be identified with the total space $E$ of the
normal bundle $\nu=\nu(M^T\to V)$:
\begin{equation}\label{fpCDvb1}
\begin{CD}
  M^T @>r_M>> E_1\\
  @Vi^TVV @VViV\\
  E_2 @>r_V>> E
\end{CD}
\end{equation}
where $E_1=\nu(r_M)$ and $E_2=\nu(i^T)$. The normal bundle $\nu$
decomposes as
\begin{equation}\label{nuzeta}
  \nu=\nu(r_M)\oplus
  r_M^*(\nu(i))=\nu(r_M)\oplus\nu(i^T)\oplus\zeta,
\end{equation}
where $\zeta$ is the `excess' bundle over $M^T$, whose fibres are
the non-trivial parts of the $T^k$-representations in the fibres
of~$r_M^*\nu(i)$. We therefore rewrite~\eqref{fpCDvb1} as
\begin{equation}\label{fpCDvb2}
\begin{CD}
  M^T @>r_M>> E_1\\
  @Vi^TVV @VViV\\
  E_2 @>r_V>> E_1\oplus E_2\oplus F,
\end{CD}
\end{equation}
where $F$ is the total space of~$\zeta$. Since all relevant
bundles are complex, we have Gysin--Thom isomorphisms (see
Construction~\ref{gysincob} and Exercise~\ref{gysin-thom})
\begin{align*}
  i_{\,!}\colon& U^*(M)\stackrel\cong\longrightarrow U^{*+p}\bigl(\Th(\nu(i))\bigr)
  =U^{*+p}(V,V\setminus M)=U^{*+p}(E,E\setminus E_1),\\
  i^T_{\,!}\colon& U^*(M^T)\stackrel\cong\longrightarrow U^{*+q}\bigl(\Th(\nu(i^T))\bigr)
  =U^{*+q}(V^T,V^T\setminus M^T)=U^{*+q}(E_2,E_2\setminus M^T),
\end{align*}
where $p=\dim V-\dim M$ and $q=\dim V^T-\dim M^T$. Let $i_1\colon
E_1\to E_1\oplus E_2$, $i_2\colon E_2\to E_1\oplus E_2$ and
$k\colon E_1\oplus E_2\to E$ be the inclusion maps. Then for $x\in
U^*(M)$,
\[
  r_V^*i_{\,!}x=i_2^*k^*k_{\,!}i_{1!}x=i_2^*\bigl(e(\nu(k))\cdot
  i_{1!}x\bigr).
\]
Since $i_2^*(\nu(k))=\pi^*(\zeta)$, where $\pi\colon E_2\to M$ is
the projection, the last term above can be written as
\begin{align*}
  \pi^*(e(\zeta))\cdot i_2^*i_{1!}x&
  =\pi^*(e(\zeta))\cdot i^T_{\,!}r_M^*x&&\text{by
  Proposition~\ref{gysinprop}~(e)}\\
  &=i^T_{\,!}\bigl(i^{T*}\pi^*(e(\zeta))\cdot r_M^*x\bigr)
  &&\text{by Proposition~\ref{gysinprop}~(d)}\\
  &=i^T_{\,!}\bigl(e(\zeta)\cdot r_M^*x\bigr).
\end{align*}
We therefore obtain
\begin{equation}\label{excess}
  r_V^*i_{\,!}x=i^T_{\,!}\bigl(e(\zeta)\cdot r_M^*x\bigr)
\end{equation}
in $U^{*+p}(V^T,V^T\setminus M^T)$.

Given a $T^k$-equivariant map $f\colon M\to N$ we denote by
$\widehat f$ its `Borelification', i.e. the map
$ET^k\times_{T^k}M\to ET^k\times_{T^k}N$. Applying this procedure
to~\eqref{fpCD} we obtain a commutative diagram
\begin{equation}\label{BfpCD}
\begin{CD}
  BT^k\times M^T @>\widehat r_M>> ET^k\times_{T^k}M\\
  @V\widehat i^TVV @VV\widehat iV\\
  BT^k\times V^T @>\widehat r_V>> ET^k\times_{T^k}V.
\end{CD}
\end{equation}
Using the finite-dimensional approximation of the above diagram
(as in the proof of Proposition~\ref{utggysin}) we can view it as
a diagram of proper maps of smooth manifolds and therefore apply
Gysin homomorphisms in cobordism. By analogy with~\eqref{excess}
we obtain for $x\in U^*(ET^k\times_{T^k}M)$,
\begin{equation}\label{excess1}
  \widehat r_V^*\widehat i_{\,!}x=
  \widehat i^T_{\,!}\bigl(e(\widehat\zeta)\cdot \widehat
  r_M^*x\bigr),
\end{equation}
where $\widehat\zeta$ is the `excess' bundle over $BT^k\times M^T$
defined similarly to~\eqref{nuzeta}.

Similarly, by considering the diagram
\[
\begin{CD}
  \mathbf 0 @= \mathbf0\\
  @Vj^TVV @VVjV\\
  V^T @>r_V>> V
\end{CD}
\]
we obtain for $y\in U^*(BT^k)$,
\begin{equation}\label{excess2}
  \widehat r_V^*\widehat j_{\,!}y=
  \widehat j^T_{\,!}\bigl(e(\widehat\zeta_V)\cdot y\bigr),
\end{equation}
where $\zeta_V$ is the nontrivial part of the
$T^k$-representation~$V$, i.e. $V=V^T\oplus\zeta_V$. Let
$\pi\colon M\to\pt$ and $\pi^T\colon M^T\to\pt$ be the
projections. Then $\widehat i_{\,!}=\widehat
j_{\,!}\cdot\widehat\pi_{\,!}$, $\widehat i^T_{\,!}=\widehat
j^T_{\,!}\cdot\widehat\pi^T_{\,!}$. Substituting $y=\widehat
\pi_{\,!}x$ in~\eqref{excess2} we obtain
\[
  \widehat r_V^*\widehat i_{\,!}x=
  \widehat j^T_{\,!}\bigl(e(\widehat\zeta_V)\cdot\widehat\pi_{\,!}x\bigr)
\]
Comparing this to~\eqref{excess1} and using the fact that
\[
  \widehat j^T_{\,!}\colon U^*(BT^k)\to U^{*+r}(\Th(\nu(\widehat
  j^T)))%=U^{*+r}(BT^k\times V^T,BT^k\times(V^T\setminus\mathbf0))
  =U^{*+r}(\varSigma^r BT^k)
\]
is an isomorphism (here $r=\dim V^T$), we finally obtain
\[
  e(\widehat\zeta_V)\cdot\widehat\pi_{\,!}x= \widehat
  \pi^T_{\,!}\bigl(e(\widehat\zeta)\cdot \widehat
  r_M^*x\bigr).
\]
Now set $x=1$. Then $\widehat\pi_{\,!}1=\varPhi(M)$ by
Proposition~\ref{utggysin} and $\widehat r_M^*1=1$. We get
\begin{equation}\label{utgnonis}
  e(\widehat\zeta_V)\cdot\varPhi(M)= \widehat
  \pi^T_{\,!}\bigl(e(\widehat\zeta)\bigr).
\end{equation}
This formula is valid without restrictions on the fixed point set.
Now, if $M^T$ is finite, then
$\pi^T_{\,!}(e(\widehat\zeta))=\sum_{p\in
M^T}e(\widehat\zeta|_p)$. Recall that $\zeta$ is defined from the
decomposition $\nu=\nu(r_M)\oplus\nu(i^T)\oplus\zeta$, in which
$\nu|_p=\nu(M^T\to V)|_p$ can be identified with $V$ and
$\nu(i^T)|_p$ can be identified with~$V^T$ (because~$p$ is
isolated). Since $V=V^T\oplus\zeta_V$, it follows that
$e(\widehat\zeta_V)=e(\nu(\widehat r_M)|_p)e(\widehat\zeta|_p)$
for any $p\in M^T$. We therefore can rewrite~\eqref{utgnonis} as
\[
  \varPhi(M)= \sum_{p\in M^T}\frac 1{e\bigl(\nu(\widehat r_M)|_p\bigr)}.
\]
It remains to note that $\nu(r_M)|_p=\nu(p\to M)$ is the
tangential $T^k$-representation $\mathcal T_p M$, so $\nu(\widehat
r_M)|_p$ is the bundle $ET^k\times_{T^k}\mathcal T_p M$ over
$BT^k$, whose Euler class is
\[
  e\bigl(\nu(\widehat r_M)|_p\bigr)=\sigma(p)\prod_{j=1}^n[\mb w_j(p),\mb
  u\,]
\]
by the definition of sign $\sigma(p)$ and weights $\mb w_j(p)$.
\end{proof}

\begin{remark}\label{genkrich}
If the structure $c_{\mathcal T}$ is almost complex, then
$\sigma(x)=1$ for all fixed points $x$, and \eqref{lfx} reduces to
Krichever's formula~\cite[(2.7)]{kric90}).
\end{remark}

The left-hand side of \eqref{lfx} lies in
$\varOmega_U[[u_1,\ldots,u_k]]$, whereas the right-hand side
appears to belong to an appropriate localisation. It follows that
all terms of negative degree must cancel, thereby imposing
substantial restrictions on the fixed point data. These may be
made explicit by rewriting~\eqref{wjxmoddec} as
\begin{equation}\label{psit}
  [\mb w,t\mb u]\;\equiv\;(w_1u_1+\dots +w_ku_k)\,t \mod (t^2)
\end{equation}
in $\varOmega_U[[u_1,\dots,u_k,t]]$, and then defining the power
series
\begin{equation}\label{tseries}
\sum_l\cf_l\,t^l\;=\;t^n\varPhi(M)(t\mb u)\;=\; \sum_{p\in
M^T}\!\sigma(p)\prod_{j=1}^n\frac{t}{[\mb w_j(p),t\mb u]}
\end{equation}
over the localisation of $\varOmega_U[[u_1,\ldots,u_k]]$.

\begin{proposition}\label{cfrels}
The coefficients $\cf_l$ are zero for $\,0\le l<n$, and satisfy
\[
  \cf_{n+m}\;=\;\sum_{|\omega|=m}g_\omega(M)u^\omega
\]
for $m\ge 0$; in particular, $\cf_n=[M]$.
\end{proposition}
\begin{proof}
Combine the definitions of $\cf_l$ in \eqref{tseries} and
$g_\omega$ in \eqref{defgomega}.
\end{proof}

\begin{remark}
The equations $\cf_l=0$ for $0\le l<n$ are the $T^k$-analogues of
the \emph{Conner--Floyd relations} for
$\Z_p$-actions~\cite[Appendix~4]{novi67}; the extra equation
$\cf_n=[M]$ provides an expression for the cobordism class of $M$
in terms of fixed point data. This is important because, according
to Theorem~\ref{6.11}, every element of $\varOmega_U$ may be
represented by a stably tangentially complex $T^k$-manifold with
isolated fixed points. We explore on this in the next section.
\end{remark}

Now let $\varphi\colon\varOmega_U\to R$ be a genus taking values
in a torsion-free ring~$R$, with the corresponding series
$f(x)=x+\cdots\in R\otimes\Q[[x]]$. We may adapt
Theorem~\ref{evutg} to express $\varphi(M)$ in terms of fixed
point data. The resulting formula is much simpler, because the
formal group law $\varphi F_U$ may be linearised over
$R\otimes\Q$:

\begin{proposition}\label{evgenus}
Let $\varphi\colon\varOmega_U\to R$ be a genus with
torsion-free~$R$, and let $M$ be a stably tangentially complex
$2n$-dimensional $T^k$-manifold with isolated fixed points~$M^T$.
Then the equivariant genus $\varphi^T(M)=\varphi(M)+\cdots$ is
given by
\begin{equation}\label{lfgenus}
  \varphi^T(M)\;=\;\sum_{p\in M^T}\!
  \sigma(p) \prod_{j=1}^n\frac1{f\bigl(\langle \mb w_j(p),\mb
  x\rangle\bigr)},
\end{equation}
where $\langle \mb w,\mb x\rangle=w_1x_1+\cdots+w_kx_k$ for $\mb
w=(w_1,\ldots,w_k)$.
\end{proposition}
\begin{proof}
By Theorem~\ref{fexpser}, $f(x)$ is the exponential series of the
formal group law $\varphi F_U$, i.e. $\varphi
F_U(u_1,u_2)=f(f^{-1}(u_1)+f^{-1}(u_2))$, and therefore
$h_{\varphi} F_U(u_1,u_2)=f(x_1+x_2)$. An iterated application of
this formula gives $h_{\varphi}([\mb w_j(p),\mb u])=f(\langle \mb
w_j(p),\mb x\rangle)$. Since $\varphi^T=h_{\varphi}\cdot\varPhi$,
the result follows from Theorem~\ref{evutg}.
\end{proof}

\begin{example}\label{agcfrels}
The augmentation genus $\varepsilon\colon\varOmega_U\to\Z$
corresponds to the series $f(x)=x$; it vanishes on any $M^{2n}$
with $n>0$. Formula~\eqref{lfgenus} then gives
\begin{equation}\label{agzero}
\sum_{p\in M^T}\!
  \sigma(p) \prod_{j=1}^n\frac1{\langle \mb w_j(p),\mb
  x\rangle}\;=\;0.
\end{equation}
Let $M=\C P^n$ on which $T^{n+1}$ acts homogeneous coordinatewise.
There are $n+1$ fixed points $p_0$, \ldots, $p_n$, each having a
single nonzero coordinate. So the weight vector $\mb w_j(p_k)$ is
$\mb e_j-\mb e_k$ for $0\le j\le n$, $j\neq k$, and every
$\sigma(p_k)$ is positive; thus~\eqref{agzero} reduces to the
classical identity
\[
  \sum_{k=0}^n\;\;\mathop{\prod_{0\le j\le
  n}}\limits_{j\neq k}\!\frac1{x_j-x_k}\;=\;0\,.
\]
\end{example}

\begin{example}
Consider the $S^1$-action preserving the standard complex
structure on~$\C P^1$. It has two fixed points, both with signs~1
and weights $1$ and $-1$, respectively (see Example~\ref{2cp1eq}).
Theorem~\ref{evutg} gives the following expression for the
universal toric genus:
\[
  \varPhi(\C P^1)=\frac1u+\frac1{\bar u}
\]
in $U^{-2}(\C P^\infty)$ where $\bar u=[-1,u]$ is the inverse
series in the formal group law of geometric cobordisms.

By Theorem~\ref{evgenus}, a genus $\varphi\colon\varOmega_U\to R$
is rigid on $\C P^1$ only if its defining series $f(x)$ satisfies
the equation
\begin{equation}\label{cp1rig}
  \frac 1{f(x)}+\frac 1{f(-x)}=c
\end{equation}
in $R\otimes\Q[[x]]$. The general analytic solution is of the form
\begin{equation}\label{ansolcp1}
  f(x)=\frac{x}{q(x^2)+cx/2},\qquad\text{where}\quad q(0)=1
\end{equation}
(an exercise). The Todd genus of Example~\ref{xyexa}.3 is defined
by the series $f(x)=1-e^{-x}$, and \eqref{cp1rig} is satisfied
with $c=1$. So $\td$ is $T$-rigid on $\C P^1$, and
\[
  q(x^2)=\frac x2\cdot\frac{e^{x/2}+e^{-x/2}}
  {e^{x/2}-e^{-x/2}}
\]
in $\Q[[x]]$. In fact $\td$ is fibre multiplicative with respect
to $\C P^1$ by~\cite{hirz66}, so rigidity also follows from
Theorem~\ref{multandrig}.

We can also consider the $S^1$-action on~$M=\C P^1$ with trivial
stably complex structure. It has two fixed points of signs $1$ and
$-1$, both with weights~1. Theorem~\ref{evutg} gives the universal
toric genus $\varPhi(M)=\frac1u-\frac1{u}=0$, which also follows
from the fact that $M$ bounds equivariantly.
\end{example}

Another classical application of the localisation formula is the
Atiyah--Hirzebruch formula~\cite{at-hi70} expressing the
$\chi_y$-genus of a complex $S^1$-manifold in terms of the fixed
point data. We discuss a generalisation of this formula due to
Krichever~\cite{kric74}. It refers to the 2-parameter
\emph{$\chi_{a,b}$-genus} corresponding to the series
\begin{equation}\label{ftab}
  f(x)=\frac{e^{a x}-e^{b x}}{a e^{b x}-b e^{ax}}\in\Q[a,b].
\end{equation}
The $\chi_y$-genus (Example~\ref{chigenu}) corresponds to the
parameter values $a=y$, $b=-1$.

We choose a circle subgroup in $T^k$ defined by a primitive vector
$\nu\in\Z^k$:
%\begin{equation}\label{cirnu1}
\[
  S(\nu)=\{(e^{2\pi i\nu_1\varphi},\ldots,e^{2\pi i\nu_k\varphi})\in\mathbb T^k\colon
  \varphi\in\R\}.
\]
%\end{equation}
We have $M^{S(\nu)}=M^T$ for a generic circle $S(\nu)$ (see
Lemma~\ref{sa}). The weights of the tangential representation of
$S(\nu)$ at~$p$ are $\langle \mb w_j(p),\nu\rangle$, $1\le j\le
n$. If fixed points $M^T$ are isolated and $M^{S(\nu)}=M^T$, then
\begin{equation}\label{gennu}
  \langle \mb w_j(p),\nu\rangle\ne0\qquad\text{ for $1\le j\le n$ and any $p\in M^T$.}
\end{equation}
We define the \emph{index} $\mathop{\mathrm{ind}}_\nu p$ as the
number of negative weights at~$p$, i.e.
\[
  \ind_\nu p=\#\{j\colon \langle \mb
  w_j(p),\nu\rangle<0\}.
\]

\begin{theorem}[generalised Atiyah--Hirzebruch formula~\cite{kric74}]
The $\chi_{a,b}$-genus is $T^k$-rigid. Furthermore, the $\chi_{ a,
b}$-genus of a stably tangentially complex $2n$-dimensional
$T^k$-manifold~$M$ with finite $M^T$ is given by
\begin{equation}\label{tablocal}
  \chi_{ a, b}(M)=\sum_{p\in M^T}\sigma(p)(- a)^{\mathop{\mathrm{ind}}_\nu p}
  (- b)^{n-\mathop{\mathrm{ind}}_\nu p}
\end{equation}
for any $\nu\subset\mathbb Z^k$ satisfying $M^{S(\nu)}=M^T$.
%, with the convention~$0^0=1$.
\end{theorem}

\begin{remark}
The original formula of Atiyah--Hirzebruch~\cite{at-hi70} was
given for complex manifolds. Krichever implicitly assumed
manifolds to be almost complex when deducing his formula, as no
signs were mentioned in~\cite{kric74}. However his proof,
presented below, automatically extends to the stably complex
situation by incorporating the signs of fixed points.

Also, both Atiyah--Hirzebruch's and Krichever's formulae are valid
without assuming the fixed points to be isolated, see
Exercise~\ref{chiabgen}. We give the formula in the case of
isolated fixed points to emphasise the role of signs.
\end{remark}

\begin{proof}[Proof of Theorem~\ref{tablocal}]
We establish the rigidity and prove the formula for $M$ with
isolated fixed points only; the proof in the general case is
similar. By Proposition~\ref{evgenus}, to prove the $T^k$-rigidity
on $M$ it suffices to prove that $\chi_{a,b}^{S(\nu)}$ is constant
for any $S(\nu)$ satisfying $M^{S(\nu)}=M^T$.
Formula~\eqref{lfgenus} gives
\begin{equation}\label{eetab}
  \chi_{ a, b}^{S(\nu)}(M)=\sum_{p\in M^T}\sigma(p)\prod_{j=1}^n
  \frac{ a e^{ b\langle\mb w_j(p),\nu\rangle x}- b e^{ a\langle\mb w_j(p),\nu\rangle x}}
  {e^{ a\langle\mb w_j(p),\nu\rangle x}-e^{ b\langle\mb w_j(p),\nu\rangle
  x}}.
\end{equation}
This expression belongs to $\mathbb Z[ a, b][[x]]$ (that is, it is
non-singular at zero) and its constant term is $\chi_{ a, b}(M)$.
We denote $\omega_j=\langle\mb w_j(p),\nu\rangle$ and $e^{(\!a-
b)x}=q$; then we may rewrite each factor in the product above as
\begin{equation}\label{factor}
  \frac{ a e^{ b\omega_jx}- b e^{ a\omega_j x}}
  {e^{ a\omega_j x}-e^{ b\omega_j x}}=
  \frac{ a- b\, q^{\omega_j}}{q^{\omega_j}-1}.
\end{equation}
Then~\eqref{eetab} takes the following form:
\begin{equation}\label{eetab1}
  \chi_{ a, b}^{S(\nu)}(M)=\sum_{p\in M^T}\sigma(p)\prod_{j=1}^n
  \frac{ a- b\,q^{\omega_j}}{q^{\omega_j}-1}.
\end{equation}
%The idea is to view $\chi_{ a, b}^{S(\nu)}(M)$ as a function in
%complex variable $z$ and study its analytic properties.
Now we let $q\to\infty$. Then~\eqref{factor} has limit $- b$ if
$\omega_j>0$ and limit $- a$ if $\omega_j<0$. Therefore, the limit
of~\eqref{eetab1} is
\[
  \sum_{p\in M^T}\sigma(p)(- a)^{\mathop{\mathrm{ind}}_\nu p}
  (- b)^{n-\mathop{\mathrm{ind}}_\nu p}.
\]
Similarly, the limit of~\eqref{eetab1} as $q\to 0$ is
\[
  \sum_{p\in M^T}\sigma(p)(- a)^{n-\mathop{\mathrm{ind}}_\nu p}
  (- b)^{\mathop{\mathrm{ind}}_\nu p}.
\]
Therefore, $\chi_{ a, b}^{S(\nu)}(M)$ has value for $q=0$ as well
as for $q=\infty$. Since it is a finite Laurent series in~$q$, it
must be constant in~$q$, and its value coincides with either of
the limits above.
\end{proof}

\begin{remark}
It follows also that the right hand side of~\eqref{tablocal} is
independent of~$\nu$; the two limits above are taken to each other
by substitution $\nu\to-\nu$.
\end{remark}

\subsection*{Exercises}

\begin{exercise}
The general analytic solution of~\eqref{cp1rig} is given
by~\eqref{ansolcp1}.
\end{exercise}

\begin{exercise}\label{fgchiab}
The formal group law corresponding to the $\chi_{a,b}$-genus is
given by
\[
  F(u,v)=\frac{u+v+(a+b)uv}{1-ab\cdot uv}.
\]
\end{exercise}

\begin{exercise}\label{chiabgen}
Let $M$ be a stably tangentially complex $T^k$-manifold~$M$ with
fixed point set~$M^T$. Then
\[
  \chi_{ a, b}(M)=\sum_{F\subset M^T}\chi_{ a, b}(F)(- a)^{\mathop{\mathrm{ind}}_\nu F}
  (- b)^{\ell-\mathop{\mathrm{ind}}_\nu F}
\]
for any $\nu\subset\mathbb Z^k$ satisfying $M^{S(\nu)}=M^T$, where
the sum is taken over connected fixed submanifolds $F\subset M^T$,
$2\ell=\dim M-\dim F$, and $\ind_\nu F$ is the number of negative
weights of the $S(\nu)$-action in the normal bundle of~$F$.
\end{exercise}

\section{Quasitoric manifolds and genera}\label{qtgenera}
In the case of quasitoric manifolds, the combinatorial description
of signs of fixed points and weights obtained in
Section~\ref{qtman} opens a way to effective calculation of
characteristic numbers and Hirzebruch genera using localisation
techniques. We illustrate this approach by presenting formulae
expressing the $\chi_{a,b}$-genus (in particular, the signature
and the Todd genus) of a quasitoric manifold as a sum of
contributions depending only on the `local combinatorics' near the
vertices of the quotient polytope. These formulae were obtained
in~\cite{pano01}; they can also be deduced from the results
of~\cite{masu99} in the more general context of torus manifolds.

Localisation formulae for genera on quasitoric manifolds can be
interpreted as functional equations on the series $f(x)$. By
resolving these equations for particular examples of quasitoric
manifolds one may derive different `universality theorems' for
rigid genera. For example, according to a result of
Musin~\cite{musi11}, the $\chi_{a,b}$-genus is universal for
$T^k$-rigid genera (this implies that any $T^k$-rigid rational
genus is $\chi_{a,b}$ for some rational parameters~$a,b$). We
prove this result as Theorem~\ref{musin} by resolving the
functional equation coming from localisation
formula~\eqref{lfgenus} on $\C P^2$ with a nonstandard
omniorientation.

\medskip

We assume given a combinatorial quasitoric pair $(P,\varLambda)$
(see Definition~\ref{cqtp}) and the corresponding omnioriented
quasitoric manifold $M=M(P,\varLambda)$. This fixes a
$T^n$-invariant stably complex structure on~$M$ and the
corresponding bordism class $[M]\in\varOmega^U_{2n}$, as described
in Corollary~\ref{qtscs}.
%The orientation of $M$ defines the fundamental homology class
%$\langle M\rangle\in H_{2n}(M)$.

Any fixed point $v\in M$ is given by the intersection of $n$
characteristic submanifolds $v=M_{j_1}\cap\cdots\cap M_{j_n}$ and
corresponds to a vertex of the polytope $P$, which we also denote
by~$v$. The expressions for the weights $\mb w_j(v)$ and the sign
$\sigma(v)$\label{weigsign} in terms of the quasitoric pair
$(P,\varLambda)$ are given by Proposition~\ref{qtw} and
Lemma~\ref{qts}.

In the quasitoric case the condition~\eqref{gennu} guarantees that
the circle $S(\nu)$ satisfies $M^T=M^{S(\nu)}$ (an exercise).

\begin{example}[{Chern number $c_n[M]$}]
The series~\eqref{ftab} defining the $\chi_{a,b}$-genus has limit
as $a-b\to0$, which can be calculated as follows:
\[
  \frac{e^{a x}-e^{b x}}{a e^{b x}-b e^{ax}}=\frac{1-e^{(b-a)x}}{ae^{(b-a)x}-b}=
  \frac{(a-b)x+\cdots}{(a-b)-a(a-b)x+\cdots}=\frac x{1-ax}.
\]
For $a=b=-1$ we obtain the defining series for the top Chern
number~$c_n[M]$ (see Example~\ref{xyexa}.1). Plugging these values
into~\eqref{tablocal} we obtain
\[
  c_n[M]=\sum_{v\in M^T}\sigma(M),
\]
which we already know from Exercise~\ref{cnsign}. When all signs
are positive, we obtain $c_n[M]=\chi(M)=f_0(P)$, i.e. the Euler
characteristic\label{eulercha} of $M$ is equal to the number of
vertices of~$P$.
\end{example}

\begin{example}[signature]
Substituting $a=1$, $b=-1$ in~\eqref{ftab} we obtain the series
$\tanh(x)$ defining the $L$-genus or the signature (see
Example~\ref{xyexa}.2). Being an oriented cobordism invariant, the
signature $\sign(M)$ does not depend on the stably complex
structure, i.e. only the global orientation part of the
omniorientation data affects the signature. The following
statement gives a formula for $\sign(M)$ which depends only on the
orientation:

\begin{proposition}
\label{signor} For an oriented quasitoric manifold $M$,
\[
  \sign(M)=\sum_{v\in M^T}
  \det\bigl(\widetilde{\mb w}_1(v),\ldots,\widetilde{\mb w}_n(v)\bigr),
\]
where $\widetilde{\mb w}_j(v)$ are the vectors defined by the
conditions
\[
  \widetilde{\mb w}_j(v)=\pm \mb w_j(v)\quad
  \text{ and }\quad\bigl(\widetilde{\mb w}_j(v),\nu\bigr)>0,\qquad
  1\le j\le n.
\]
\end{proposition}
\begin{proof}
Plugging $a=1$, $b=-1$ into~\eqref{tablocal} we obtain
\begin{equation}\label{signat}
  \sign(M)=\sum_{v\in M^T}(-1)^{\ind_\nu(v)}\sigma(v).
\end{equation}
Using expression~\eqref{wsimpl} for the sign~$\sigma(v)$ we
calculate
\[
  (-1)^{\ind_\nu(v)}\sigma(v)=
  (-1)^{\ind_\nu(v)}\det\bigl(\mb w_1(v),\ldots,\mb w_n(v)\bigr)
  =\det\bigl(\widetilde{\mb w}_1(v),\ldots,\widetilde{\mb
  w}_n(v)\bigr),
\]
which implies the required formula.
\end{proof}

If $M$ is a projective toric manifold $V_P$, then $\sigma(v)=1$
for any~$v$ and formula~\eqref{signat} gives
\[
  \sign(V_P)=\sum_{v}(-1)^{\ind_\nu(v)}.
\]
Furthermore, in this case the weights $\mb w_1(v),\ldots,\mb
w_n(v)$ are the primitive vectors along the edges of $P$ pointing
out of~$v$ (see Example~\ref{torvarsig}). It follows that the
index $\ind_\nu(v)$ coincides with the index defined in the proof
of Dehn--Sommerville equations (Theorem~\ref{ds}), and we obtain
the formula known in toric geometry
(see~\cite[Theorem~3.12]{oda88}):
\[
  \sign(V_P)=\sum_{k=0}^n(-1)^kh_k(P).
\]
Note that if $n$ is odd then the sum vanishes by the
Dehn--Sommerville equations.
\end{example}

\begin{example}[Todd genus]
Substituting $a=0$, $b=-1$ in~\eqref{ftab} we obtain the series
$1-e^{-x}$ defining the Todd genus (see Example~\ref{xyexa}.3). We
cannot plug $a=0$ directly into~\eqref{tablocal}, but it is clear
from the proof that when $a=0$, only vertices of index~0
contribute $(-b)^n$ to the sum. This gives the following formula
for the Todd genus of a quasitoric manifold:
\begin{equation}\label{tdqtm}
  \td(M)=\sum_{v\colon\ind_\nu(v)=0}\sigma(v).
\end{equation}

When $M$ is projective toric manifold~$V_P$, there is only one
vertex of index~0. It is the `bottom' vertex, which has all
incident edges pointing out (in the notation used in the proof of
Theorem~\ref{ds}). Since $\sigma(v)=1$ for every $v\in P$,
formula~\eqref{tdqtm} gives $\td(V_P)=1$, which is well-known
(see, e.g.~\cite[\S5.3]{fult93}).

In the almost complex case we have the following result:

\begin{proposition}
If a quasitoric manifold $M$ admits an equivariant almost complex
structure, then $\td(M)>0$.
\end{proposition}
\begin{proof}
We choose a compatible omniorientation, so that $\sigma(v)>0$ for
any~$v$. Then~\eqref{tdqtm} implies $\td(M)\ge0$, and we need to
show that there is at least one vertex of index~0. Let $v$ be any
vertex. Since $\mb w_1(v),\ldots,\mb w_n(v)$ are linearly
independent, we may choose $\nu$ so that $\langle \mb
w_j(v),\nu\rangle>0$ for any~$j$. Then $\ind_\nu v=0$. The result
follows by observation that $\td(M)$ is independent of~$\nu$.
\end{proof}

A description of $\td(M)$ which is independent of $\nu$ is
outlined in Exercise~\ref{toddindep}.
\end{example}

The following result of Musin~\cite{musi80} can be proved by
application of localisation formula for quasitoric manifolds:

\begin{theorem}\label{musin}
The 2-parameter genus $\chi_{a,b}$ is universal for $T^k$-rigid
genera. In particular, any $T^k$-rigid rational genus is
$\chi_{a,b}$ for some rational parameters~$a,b$.
\end{theorem}
\begin{proof}
The rigidity of $\chi_{a,b}$ is established by
Theorem~\ref{tablocal}. To see that any $T^k$-rigid genus is
$\chi_{a,b}$ we solve the functional equation arising from the
localisation formula for one particular example of $T^k$-manifold.

We consider the quasitoric manifold $M=\C P^2$ with a nonstandard
omniorientation defined by the characteristic matrix
\[
  \varLambda=\begin{pmatrix}1&0&1\\0&1&-1\end{pmatrix}
\]
It has three fixed points $v_1,v_2,v_3$, whose corresponding
minors $\varLambda_{v_i}$ are obtained by deleting the $i$th
column of~$\varLambda$. The weights are given by
$\{(1,0),(1,1)\}$, $\{(0,-1),(1,1)\}$ and $\{(0,1),(1,0)\}$,
respectively, see Fig.~\ref{cp2nstd}.
\begin{figure}[h]
\begin{center}
\begin{picture}(100,60)
  \put(20,10){\line(0,1){45}}
  \put(20,10){\line(1,0){45}}
  \put(20,55){\line(1,-1){45}}
  \put(33,23){\oval(13,13)[b]}
  \put(33,23){\oval(13,13)[tr]}
  \put(34,29.5){\vector(-1,0){2}}
  \put(16,6){$v_3$}
  \put(22,6){\small $(1,0)$}
  \put(54,6){\small $(1,0)$}
  \put(67,6){$v_1$}
  \put(32,0){$\lambda_2=(0,1)$}
  \put(11,13){\small $(0,1)$}
  \put(0,30){$\lambda_1=(1,0)$}
  \put(9,50){\small $(0,-1)$}
  \put(16,57){$v_2$}
  \put(25,52){\small $(1,1)$}
  \put(62,15){\small $(1,1)$}
  \put(47,32){$\lambda_3=(1,-1)$}
\end{picture}%
\caption{$\C P^2$ with nonstandard omniorientation}\label{cp2nstd}
\end{center}
\end{figure}
The signs are calculated using formula~\eqref{wsimpl}:
\[
  \sigma(v_1)=\begin{vmatrix} 1&1\\1&0 \end{vmatrix}=-1,\quad
  \sigma(v_2)=\begin{vmatrix} 0&1\\-1&1 \end{vmatrix}=1,\quad
  \sigma(v_3)=\begin{vmatrix} 1&0\\0&1 \end{vmatrix}=1.
\]
Plugging these data into formula~\eqref{lfgenus} we obtain that a
genus $\varphi$ is rigid on $M$ only if its defining series $f(x)$
satisfies the equation
\[
  -\frac{1}{f(x_1)f(x_1+x_2)}+\frac{1}{f(-x_2)f(x_1+x_2)}+
  \frac{1}{f(x_1)f(x_2)}\;=\;c.
\]
Interchanging $x_1$ and $x_2$ gives
\begin{equation}\label{help1}
  -\frac{1}{f(x_2)f(x_1+x_2)}+\frac{1}{f(-x_1)f(x_1+x_2)}+
  \frac{1}{f(x_2)f(x_1)}\;=\;c\,,
\end{equation}
and subtraction yields
\[
  \left(\frac{1}{f(x_1)}+\frac{1}{f(-x_1)}\right)\frac{1}{f(x_1+x_2)}
  \;=\;
  \left(\frac{1}{f(x_2)}+\frac{1}{f(-x_2)}\right)\frac{1}{f(x_1+x_2)}\,.
  \]
It follows that
\[
  \frac{1}{f(x)}+\frac{1}{f(-x)}\;=\;c'\quad\text{ and }\quad
  \frac{1}{f(-x)}\;=\; c'-\frac{1}{f(x)}
\]
for some constant $c'$. Substituting in \eqref{help1} gives
\[
  \left(\frac{1}{f(x_1)}+\frac{1}{f(x_2)}-c'\right)\frac{1}{f(x_1+x_2)}\;=\;
  \frac{1}{f(x_1)f(x_2)}-c\,,
\]
which rearranges to
\[
  f(x_1+x_2)\;=\;
  \frac{f(x_1)+f(x_2)-c'f(x_1)f(x_2)}{1-cf(x_1)f(x_2)}\,.
\]
So $f$ is the exponential series of the formal group law
$F(x_1,x_2)$ corresponding to $\chi_{a,b}$, with $c'=-a-b$ and
$c=ab$ (see Exercise~\ref{fgchiab}).
\end{proof}

\subsection*{Exercises}
\begin{exercise}
If $\nu\in\Z^n$ satisfies~\eqref{gennu}, then the circle subgroup
$S(\nu)\subset T^n$ acts on the quasitoric manifold $M$ with
isolated fixed points corresponding to the vertices of the
quotient polytope~$P$.
\end{exercise}

\begin{exercise}\label{toddindep}
Let $M=M(P,\varLambda)$ be a quasitoric manifold. We realise the
dual complex $\mathcal K_P$ as a triangulated sphere with vertices
at the unit vectors $\frac{\mb a_i}{|\mb a_i|}$, $i=1,\ldots,m$.
Then one can define a continuous piecewise smooth map $f\colon
S^{n-1}\to S^{n-1}$ by sending $\frac{\mb a_i}{|\mb a_i|}$ to
$\frac{\lambda_i}{|\lambda_i|}$ (here $\lambda_i$ is the $i$th
column of $\varLambda$) and extending the map smoothly on the
spherical simplices corresponding to $I\in K_P$. Such an extension
is well defined because the vectors $\{\lambda_i\colon i\in I\}$
are linearly independent, so one always chooses the smallest
spherical simplex spanned by them.

Show that $\td(M)=\deg f$. (Hint: the number of preimages of
$\frac{\nu}{|\nu|}\in S^{n-1}$ under $f$ is equal to the number of
maximal simplices $I_v\in\sK$ such that all coefficients in the
decomposition of $\nu$ via $\{\lambda_i\colon i\in I_v\}$ are
positive; these coefficients are $\langle\mb w_i(v),\nu\rangle$.)

More generally, the Todd genus of a torus manifold can be
calculated in this way as the \emph{degree} of the corresponding
multi-fan, see~\cite[\S3]{ha-ma03}.
\end{exercise}

\begin{exercise}
Calculate $c_n[M]$, the signature and the Todd genus for
quasitoric manifolds of Example~\ref{cpnoo}, Example~\ref{cpquas}
and Exercise~\ref{qtnt}.
\end{exercise}

% Более подробно написать про две стаб. комплексные стр-ры на CP^1
% (как фактор C^2 и H^1)

%T главу о формальнvх группах привести конструкциі универс. форм. группv и доказательство теоремv Tазара.

\begin{appendix}
%\renewcommand{\thechapter}{\Roman{chapter}}
%\renewcommand{\thetheorem}{\Roman{chapter}.\arabic{theorem}}
%\numberwithin{theorem}{chapter}
%\renewcommand{\theexercise}{\Roman{chapter}.\arabic{exercise}}
%\numberwithin{exercise}{chapter}

%\chapter{Classification of simple 3-polytopes with few facets}

\chapter{Commutative and homological algebra}\label{hab}
Here we review some basic algebraic notions and results in a way
suited for topological applications. In order to make algebraic
constructions compatible with topological ones we sometimes use a
notation which may seem unusual to a reader with an algebraic
background. This in particular concerns the way we treat gradings
and resolutions.

We fix a ground ring~$\k$, which is always assumed to be a field
or the ring $\Z$ of integers. In the latter case by a `$\k$-vector
space of dimension~$d$' we mean an abelian group of rank~$d$.

\section{Algebras and modules}\label{algmod}
A $\k$-\emph{algebra} (or simply \emph{algebra}) $A$ is a ring
which is also a $\k$-vector space, and whose multiplication
$A\times A\to A$ is $\k$-bilinear. (The latter condition is void
if $\k=\Z$, so $\Z$-algebras are ordinary rings.) All our algebras
will be commutative and with unit~$1$, unless explicitly stated
otherwise. The basic example is $A=\k[v_1,\ldots,v_m]$, the
\emph{polynomial algebra} in $m$ generators, for which we shall
often use a shortened notation $\k[m]$.

An algebra $A$ is \emph{finitely generated} if there are finitely
many elements $a_1,\ldots,a_n$ of $A$ such that every element of
$A$ can be written as a polynomial in $a_1,\ldots,a_n$ with
coefficients in~$\k$. Therefore, a finitely generated algebra is
the quotient of a polynomial algebra by an ideal.

An \emph{$A$-module}\label{dmodule} is a $\k$-vector space $M$ on
which $A$ acts linearly, that is, there is a map $A\times M\to M$
which is $\k$-linear in each argument and satisfies $1m=m$ and
$(ab)m=a(bm)$ for all $a,b\in A$, $m\in M$. Any ideal $I$ of $A$
is an $A$-module. If $A=\k$, then an $A$-module is a $\k$-vector
space.

An $A$-module $M$ is \emph{finitely generated}\label{fgmodule} if
there exist $x_1,\ldots,x_n$ in $M$ such that every element $x$ of
$M$ can be written (not necessarily uniquely) as
$x=a_1x_1+\cdots+a_nx_n$, \ $a_i\in A$.

\smallskip

An algebra $A$ is \emph{$\Z$-graded}\label{gradedalgebra} (or
simply \emph{graded}) if it is represented as a direct sum
$A=\bigoplus_{i\in\Z}A^i$ such that $A^i\cdot A^j\subset A^{i+j}$.
Elements $a\in A^i$ are said to be \emph{homogeneous} of
degree~$i$, denoted $\deg a=i$. The set of homogeneous elements of
$A$ is denoted by $\mathcal H(A)=\bigcup_{i}A^i$. An ideal $I$ of
$A$ is \emph{homogeneous} if it is generated by homogeneous
elements. In most cases our graded algebras will be either
\emph{nonpositively graded} (i.e. $A^i=0$ for $i>0$) or
\emph{nonnegatively graded} (i.e. $A^i=0$ for $i<0$); the latter
is also called an $\N$-graded algebra. A nonnegatively graded
algebra $A$ is \emph{connected} if $A^0=\k$. For a nonnegatively
graded algebra~$A$, define the \emph{positive ideal} by
$A^+=\bigoplus_{i>0}A^i$; if $A$ is connected and $\k$ is a field
then $A^+$ is a maximal ideal.

If $A$ is a graded algebra, then an $A$-module $M$ is
\emph{graded}\label{gradedmodule} if $M=\bigoplus_{i\in\Z}M^i$
such that $A^i\cdot M^j\subset M^{i+j}$. An $A$-module map
$f\colon M\to N$ between two graded modules is
\emph{degree-preserving} (or of \emph{degree~$0$}) if
$f(M^i)\subset N^i$, and is of \emph{degree~$k$} if $f(M^i)\subset
N^{i+k}$ for all~$i$.

Graded algebras arising in topology are often
\emph{graded-commutative}\label{gradedcommu} (or
\emph{skew-commutative}) rather than commutative in the usual
sense. This means that
\[
  ab=(-1)^{ij}\,ba\quad\text{for any }a\in A^i, b\in A^j.
\]
If the characteristic of $\k$ is not 2, then the square of an
odd-degree element in a graded-commutative algebra is zero. To
avoid confusion we double the grading in commutative algebras $A$;
the resulting graded algebras $A=\bigoplus_{i\in\Z}A^{2i}$ are
commutative in either sense.

For example, we make the polynomial algebra $\k[v_1,\ldots,v_m]$
graded by setting $\deg v_i=2$. It then becomes a
\emph{free}\label{freegco} graded commutative algebra on $m$
generators of degree two (free means no relations apart from the
graded commutativity). The \emph{exterior algebra}
$\Lambda[u_1,\ldots,u_m]$ has relations $u_i^2=0$ and
$u_iu_j=-u_ju_i$. We shall assume $\deg u_i=1$ unless otherwise
specified. An exterior algebra is a free graded commutative
algebra if the characteristic of $\k$ is not~$2$.
%Given a subset $I=\{i_1,\ldots,i_k\}\subset[m]$ we denote by $v_I$
%the square-free monomial $v_{i_1}\cdots v_{i_k}$ in~$\k[m]$. We
%also denote by $u_I$ the exterior monomial $u_{i_1}\cdots u_{i_k}$
%where $i_1<\cdots<i_k$.

\emph{Bigraded}\label{bigradedalge} (i.e. $\Z\oplus\Z$-graded) and
\emph{multigraded} ($\Z^m$-graded) algebras $A$ are defined
similarly; their homogeneous elements $a\in A$ have
\emph{bidegree} $\bideg a=(i,j)\in\Z\oplus\Z$ or
\emph{multidegree} $\mathop{\mathrm{mdeg}}a=\mb i\in\Z^m$
respectively.

\smallskip

We continue assuming $A$ to be (graded) commutative. The
\emph{tensor product}\label{tensoramod} $M\otimes_A N$ of
$A$-modules $M$ and $N$ is the quotient of a free $A$-module with
generator set $M\times N$ by the submodule generated by all
elements of the following types:
\begin{gather*}
  (x+x',y)-(x,y)-(x',y),\quad (x,y+y')-(x,y)-(x,y'),\\
  (ax,y)-a(x,y),\quad (x,ay)-a(x,y),
\end{gather*}
where $x,x'\in M$, $y,y'\in N$, $a\in A$. For each basis element
$(x,y)$, its image in $M\otimes_A N$ is denoted by $x\otimes y$.

We shall denote the tensor product $M\otimes_\k N$ of $\k$-vector
spaces by simply $M\otimes N$. For example, if $M=N=\k[v]$, then
$M\otimes N=\k[v_1,v_2]$.

The tensor product $A\otimes B$ of graded-commutative algebras $A$
and $B$ is a graded commutative algebra, with the multiplication
defined on homogeneous elements by
\[
  (a\otimes b)\cdot(a'\otimes b')=(-1)^{\deg b\deg a'}aa'\otimes
  bb'.
\]

An $A$-module $F$ is \emph{free}\label{freemodule} if it is
isomorphic to a direct sum $\bigoplus_{i\in I}F_i$, where each
$F_i$ is isomorphic to $A$ as an $A$-module. If both $A$ and $F$
are graded then every $F_i$ is isomorphic to a $j$-fold
\emph{suspension} $s^jA$ for some~$j$, where $s^jA$ is the graded
$A$-module with $(s^jA)^k=A^{k-j}$. A \emph{basis} of a free
$A$-module $F$ is a set $\mathcal S$ of elements of $F$ such that
each $x\in F$ can be uniquely written as a finite linear
combination of elements of $\mathcal S$ with coefficients in~$A$.
If $A$ is finitely generated then all bases have the same
cardinality (an exercise), called the \emph{rank} of~$F$. If
$\mathcal S$ is a basis of a free $A$-module $F$, then for any
$A$-module $M$ a set map $\mathcal S\to M$ extends uniquely to an
$A$-module homomorphism $F\to M$.

A module $P$ is \emph{projective}\label{projectivemod} if for any
epimorphism of modules $p\colon M\to N$ and homomorphism $f\colon
P\to N$, there is a homomorphism $f'\colon P\to M$ such that
$pf'=f$. This is described by the following commutative diagram:
\[
\xymatrix{
  M\ar[r]^p & N\ar[r] & 0\\
            & P\ar[u]_f\ar@{-->}[ul]^{f'}
}
\]
Equivalently $P$ is projective if it is a direct summand of a free
module (an exercise). In particular, free modules are projective.

\smallskip

A sequence of homomorphisms of $A$-modules
\[
  \cdots\longrightarrow M_1\stackrel{f_1}\longrightarrow M_2\stackrel{f_2}\longrightarrow
  M_3\stackrel{f_3}\longrightarrow M_4\longrightarrow\cdots
\]
is called an \emph{exact sequence}\label{exaseque} if $\Im
f_i=\Ker f_{i+1}$ for all~$i$.

A \emph{chain complex}\label{chaincomplex} is a sequence
$C_*=\{C_i,\partial_i\}$ of $A$-modules $C_i$ and homomorphisms
$\partial_i\colon C_i\to C_{i-1}$ such that
$\partial_i\partial_{i+1}=0$. This is usually written as
\[
 \cdots\longrightarrow
 C_{i+1}\stackrel{\partial_{i+1}}\longrightarrow C_i
 \stackrel{\partial_i}\longrightarrow
 C_{i-1}\longrightarrow\cdots
 %\longrightarrow C_1
 %\stackrel{\partial_1}\longrightarrow C_0\longrightarrow 0.
\]
The condition $\partial_i\partial_{i+1}=0$ implies that
$\Im\partial_{i+1}\subset\Ker\partial_i$. The elements of
$\Ker\partial$ are called \emph{cycles}, and the elements of
$\Im\partial$ are \emph{boundaries}. The \emph{$i$th homology
group} (or \emph{homology module}) of $C_*$ is defined by
\[
  H_i(C_*)=\Ker\partial_i/\Im\partial_{i+1}.
\]

A \emph{cochain complex}\label{cochaincomplex} is a sequence
$C^*=\{C^i,d^i\}$ of $A$-modules $C^i$ and homomorphisms
$d^i\colon C^i\to C^{i+1}$ such that $d^id^{i-1}=0$. This is
usually written as
\[
  %0\longrightarrow C^0\stackrel{d^0}\longrightarrow C^1\longrightarrow
  \cdots\longrightarrow
  C^{i-1}\stackrel{d^{i-1}}\longrightarrow
  C^i\stackrel{d^i}\longrightarrow C^{i+1}\longrightarrow\cdots.
\]
The elements of $\Ker d$ are called \emph{cocycles}, and the
elements of $\Im d$ are \emph{coboundaries}. The \emph{$i$th
cohomology group} (or \emph{cohomology module}) of $C^*$ is
defined by
\[
  H^i(C^*)=\Ker d^i/\Im d^{i-1}.
\]
A cochain complex may be also viewed as a graded $\k$-vector space
$C^*=\bigoplus_i C^i$ in which every graded component $C^i$ is an
$A$-module, together with an $A$-linear map $d\colon C^*\to C^*$
raising the degree by~1 and satisfying the condition $d^2=0$.

Note that a chain complex may be turned to a cochain complex by
inverting the grading (i.e. turning the $i$th graded component
into the $(-i)$th).

A \emph{map of cochain complexes} is a graded $A$-module map
$f\colon C^*\to D^*$ which commutes with the differentials. Such a
map induces a map in cohomology $\widetilde f\colon H(C^*)\to
H(D^*)$, which is also an $A$-module map.

Let $f,g\colon C^*\to D^*$ be two maps of cochain complexes. A
\emph{cochain homotopy}\label{cochainhomotopy} between $f$ and $g$
is a set of maps $s=\{s^i\colon C^i\to D^{i-1}\}$ satisfying the
identities
\[
  ds+sd=f-g
\]
(more precisely, $d^{i-1}s^i+s^{i+1}d^i=f^i-g^i$). This is
described by the following commutative diagram
\[
\xymatrix{
  \cdots \ar[r] &
  C^{i-1} \ar[r]^{d^{i-1}} \ar[ld]^{s^{i-1}} \ar[d]
  & C^i \ar[r]^{d^i} \ar[ld]^{s^i} \ar[d]^(0.3){f^i-g^i} &
  C^{i+1} \ar[r] \ar[ld]^{s^{i+1}} \ar[d] & \cdots\\
  \cdots \ar[r] &
  D^{i-1} \ar[r]_{d^{i-1}} & D^i \ar[r]_{d^i} &
  D^{i+1} \ar[r] & \cdots
}
\]
If there is a cochain homotopy between $f$ and~$g$, then $f$ and
$g$ induce the same map in cohomology (an exercise). A \emph{chain
homotopy} between maps of chain complexes is defined similarly.

\smallskip

A \emph{differential graded algebra} (a \emph{dg-algebra} for
short) is a graded algebra $A$ together with a $\k$-linear map
$d\colon A\to A$, called the \emph{differential}, which raises the
degree by one, and satisfies the identity $d^2=0$ (so that
$\{A^i,d^i\}$ is a cochain complex) and the \emph{Leibniz
identity}
\begin{equation}\label{leibnitz}
  d(a\cdot b)=da\cdot b+(-1)^ia\cdot db\quad\text{for }a\in A^i, b\in
  A.
\end{equation}
In order to emphasise the differential, we may display a
dg-algebra $A$ as $(A,d)$. Its cohomology $H(A,d)=\Ker d/\Im d$ is
a graded algebra (an exercise). Differential graded algebras whose
differential lowers the degree by one are also considered, in
which case homology is a graded algebra.

A \emph{quasi-isomorphism}\label{dquism} between dg-algebras is a
homomorphism $f\colon A\to B$ which induces an isomorphism in
cohomology, $\widetilde f\colon H(A)\stackrel\cong\longrightarrow
H(B)$.

\subsection*{Exercises.}
\begin{exercise}
If $A$ is a finitely generated algebra, then all bases of a free
$A$-module have the same cardinality.
\end{exercise}

\begin{exercise}
A module is projective if and only if it is a direct summand of a
free module.
\end{exercise}

\begin{exercise}\label{5lemma}
Prove the following extended version of the \emph{5-lemma}. Let
\[
\xymatrix{
  C^1 \ar[r] \ar[d]^(0.4){f^1} & C^2 \ar[r] \ar[d]^(0.4){f^2}
  & C^3 \ar[r] \ar[d]^(0.4){f^3} & C^4 \ar[r] \ar[d]^(0.4){f^4}
  & C^5 \ar[d]^(0.4){f^5}\\
  D^1 \ar[r] & D^2 \ar[r] & D^3 \ar[r]
  & D^4 \ar[r] & D^5
}
\]
be a commutative diagram with exact rows. Then
\begin{itemize}
\item[(a)] if $f^2$ and $f^4$ are monomorphisms and $f^1$ is an
epimorphism, then $f^3$ is a monomorphism;
\item[(b)] if $f^2$ and $f^4$ are epimorphisms and $f^5$ is a
monomorphism, then $f^3$ is an epimorphism.
\end{itemize}
Therefore, is $f^1,f^2,f^4,f^5$ are isomorphisms, then $f^3$ is
also an isomorphism.
\end{exercise}

\begin{exercise}
Cochain homotopic maps between cochain complexes induce the same
maps in cohomology.
\end{exercise}

\begin{exercise}
The cohomology of a differential graded algebra is a graded
algebra.
\end{exercise}

\section{Homological theory of graded rings and modules}\label{torapdx}
From now on we assume that $A$ is a commutative finitely generated
$\k$-algebra with unit, graded by nonnegative even numbers (i.e.
$A=\bigoplus_{i\ge0}A^i$) and connected (i.e. $A^0=\k$). The basic
example to keep in mind is $A=\k[m]=\k[v_1,\ldots,v_m]$ with $\deg
v_i=2$, however we shall need a greater generality occasionally.
We also assume that all $A$-modules $M$ are nonnegatively graded
and finitely generated, and all module maps are degree-preserving,
unless the contrary is explicitly stated.

A \emph{free} (respectively, \emph{projective}) \emph{resolution}
of $M$ is an exact sequence of $A$-modules
\begin{equation}
\label{resol}
  \cdots \stackrel{d}{\longrightarrow} R^{-i} \stackrel{d}{\longrightarrow}
  \cdots \stackrel{d}{\longrightarrow} R^{-1}
  \stackrel{d}{\longrightarrow} R^0 \to M \to 0
\end{equation}
in which all $R^{-i}$ are free (respectively, projective)
$A$-modules. A free resolution exists for every~$M$ (an exercise,
or see constructions below). The minimal number $p$ for which
there exists a projective resolution~(\ref{resol}) with $R^{-i}=0$
for $i>p$ is called the \emph{projective} (or \emph{homological})
\emph{dimension} of the module $M$; we shall denote it by
$\pdim_AM$ or simply $\pdim M$. If such $p$ does not exist, we set
$\pdim M=\infty$. The module $M_i=\Ker[d\colon R^{-i+1}\to
R^{-i+2}]$ is called the \emph{$i$th syzygy module} for~$M$.
%Tаким образом, $\pdim M$ есть наименьшее $h$, для которого модуль
%$M_h$ проективен.

If $\k$ is a field and $A$ is as above, then an $A$-module is
projective if and only if it is free (see
Exercise~\ref{exprojfree}), and we therefore need not to
distinguish between free and projective resolutions in this case.

We can convert a resolution~\eqref{resol} into a bigraded
$\k$-vector space $R=\bigoplus_{i,j} R^{-i,j}$ where
$R^{-i,j}=(R^{-i})^j$ is the $j$th graded component of the
module~$R^{-i}$, and the $(-i,j)$th component of $d$ acts as
$d^{-i,j}\colon R^{-i,j}\to R^{-i+1,j}$. We refer to the first
grading of~$R$ as \emph{external}\label{extdegree}; it comes from
the indexing of the terms in the resolution and is therefore
nonpositive by our convention. The second, \emph{internal},
grading of $R$ comes from the grading in the modules $R^{-i}$ and
is therefore even and nonnegative. The \emph{total} degree of an
element of $R$ is defined as the sum of its external and internal
degrees. We can view $A$ as a bigraded algebra with trivial first
grading (i.e. $A^{i,j}=0$ for $i\ne0$ and $A^{0,j}=A^j$); then $R$
becomes a bigraded $A$-module.

If we drop the term $M$ in resolution~\eqref{resol}, then the
resulting cochain complex is exactly $R$ (with respect to its
external grading), and we have
\begin{align*}
  H^{-i,j}(R,d)&=\Ker d^{-i,j}/\Im d^{-i-1,j}=0\quad\text{for }
  i>0,\\
  H^{0,j}(R,d)&=M^j.
\end{align*}
We may view $M$ as a trivial cochain complex $0\to M\to 0$, or as
a bigraded module with trivial external grading, i.e. $M^{i,j}=0$
for $i\ne0$ and $M^{0,j}=M^j$. Then resolution~\eqref{resol} can
be interpreted as a map of cochain complexes of $A$-modules:
\begin{equation}\label{quism}
\begin{CD}
  \cdots @>>> R^{-i} @>d>> \cdots @>d>> R^{-1} @>d>> R^0 @>>> 0\\
  @.          @VVV         @.           @VVV         @VVV\\
  \cdots @>>> 0      @>>>  \cdots @>>>  0      @>>>  M   @>>> 0
\end{CD}
\end{equation}
or simply as a map $(R,d)\to(M,0)$ inducing an isomorphism in
cohomology.

The \emph{Poincar\'e series}\label{poiser} of a graded $\k$-vector
space $V=\bigoplus_i V^i$ whose graded components are
finite-dimensional is given by
$$
  F(V;\lambda)=\sum_i(\dim_\k V^i)\lambda^i.
$$

%The Poincar\'e series of $M$ (viewed as a graded $\k$-vector
%space) can be calculated from any free resolution~(\ref{resol}):

\begin{proposition}
\label{psresol} Let~\eqref{resol} be a finite free resolution of
an $A$-module, in which $R^{-i}$ is a free module of rank $q_i$ on
generators of degrees $d_{1i},\ldots,d_{q_ii}$. Then
$$
  F(M;\lambda)=F(A;\lambda)\sum_{i\ge0}(-1)^i(\lambda^{d_{1i}}+\cdots+\lambda^{d_{q_ii}}).
$$
\end{proposition}
\begin{proof}
Since $H^{-i,j}(R,d)=0$ for $i>0$ and $H^{0,j}(R,d)=M^j$, we
obtain
$$
  \sum_{i\ge0}(-1)^i\dim_\k R^{-i,j}=\dim_\k M^j
$$
by a basic property of the Euler characteristic. Multiplying by
$\lambda^j$ and summing up over~$j$ we obtain
$$
  \sum_{i\ge0}(-1)^iF(R^{-i};\lambda)=F(M;\lambda).
$$
Since each $R^{-i}$ is a free $A$-module, its Poincar\'e series is
given by
$F(R^{-i};\lambda)=F(A;\lambda)(\lambda^{d_{1i}}+\cdots+\lambda^{d_{q_ii}})$,
which implies the required formula.
\end{proof}

\begin{construction}[minimal resolution]
\label{minimal} Let $\k$ be a field, and let
$M=\bigoplus_{i\ge0}M^i$ be a graded $A$-module, which is not
necessarily finitely generated, but for which every graded
component $M^i$ is finite-dimensional as a $\k$-vector space.
There is the following canonical way to construct a free
resolution for~$M$.

Take the lowest degree $i$ in which $M^i\ne0$ and choose a
$\k$-vector space basis in~$M^i$. Span an $A$-submodule $M_1$ by
this basis and then take the lowest degree in which $M\ne M_1$. In
this degree choose a $\k$-vector space basis in the complement
of~$M_1$, and span a module $M_2$ by this basis and~$M_1$.
Continuing this process we obtain a system of generators for~$M$
which has a finite number of elements in each degree, and has the
property that images of the generators form a basis of the
$\k$-vector space $M\otimes_A\k=M/(A^+\cdot M)$. A system of
generators of $M$ obtained in this way is referred to as
\emph{minimal} (or as a \emph{minimal basis}).

Now choose a minimal generating set in $M$ and span by its
elements a free $A$-module $R_{\min}^0$. Then we have an
epimorphism $R_{\min}^0\to M$. Next we choose a minimal basis in
the kernel of this epimorphism, and span by it a free module
$R_{\min}^{-1}$. Then choose a minimal basis in the kernel of the
map $R_{\min}^{-1}\to R_{\min}^0$, and so on. At the $i$th step we
choose a minimal basis in the kernel of the map $d\colon
R_{\min}^{-i+1}\to R_{\min}^{-i+2}$ constructed in the previous
step, and span a free module $R_{\min}^{-i}$ by this basis. As a
result we obtain a free resolution of~$M$, which is referred to as
\emph{minimal}. A minimal resolution is unique up to an
isomorphism.
\end{construction}

\begin{proposition}\label{mintensor}
For a minimal resolution of $M$, the induced maps
\[
  R_{\min}^{-i}\otimes_A\k\to R_{\min}^{-i+1}\otimes_A\k
\]
are zero for $i\ge1$.
\end{proposition}
\begin{proof} By construction, the map $R^0_{\min}\otimes_\k A\to
M\otimes_\k A$ is an isomorphism, which implies that the kernel of
the map $R_{\min}^0\to M$ is contained in $A^+\cdot R_{\min}^0$.
Similarly, for each $i\ge1$ the kernel of the map $d\colon
R_{\min}^{-i}\to R_{\min}^{-i+1}$ is contained in $A^+\cdot
R_{\min}^{-i}$, and therefore the image of the same map is
contained in $A^+\cdot R_{\min}^{-i+1}$. This implies that the
induced maps $R_{\min}^{-i}\otimes_A\k\to
R_{\min}^{-i+1}\otimes_A\k$ are zero.
\end{proof}

\begin{remark}
If $\k=\Z$ then the above described inductive procedure still
gives a minimal basis for an $A$-module~$M$, but the kernel of the
map $d\colon R^0_{\min}\to M$ may not be contained in $A^+\cdot
R^0_{\min}$, and the induced map $R_{\min}^{-1}\otimes_A\Z\to
R_{\min}^0\otimes_A\Z$ may be nonzero.
\end{remark}

\begin{construction}[Koszul resolution]
\label{koszul} Let $A=\k[v_1,\ldots,v_m]$ and $M=\k$ with the
$A$-module structure given by the augmentation map sending each
$v_i$ to zero.
%Consider the exterior $\k$-algebra
%$\Lambda[u_1,\ldots,u_m]$ on degree-one generators, which is
%determined by the relations $u_i^2=0$ and $u_iu_j+u_ju_i=0$.
We turn the tensor product
\[
  E=E_m=\Lambda[u_1,\ldots,u_m]\otimes\k[v_1,\ldots,v_m]
\]
into a \emph{bigraded differential algebra} by setting
\begin{equation}\label{kosdiff}
\begin{gathered}
  \bideg u_i=(-1,2),\quad\bideg v_i=(0,2),\\
  du_i=v_i,\quad dv_i=0
\end{gathered}
\end{equation}
and requiring $d$ to satisfy the Leibniz
identity~\eqref{leibnitz}. Then $(E,d)$ together with the
augmentation map $\varepsilon\colon E\to\k$ defines a cochain
complex of $\k[m]$-modules
\begin{multline}\label{kosz}
  0\to\Lambda^m[u_1,\ldots,u_m]\otimes\k[v_1,\ldots,v_m]
  \stackrel{d}{\longrightarrow}\cdots\\
  \stackrel{d}{\longrightarrow}\Lambda^1[u_1,\ldots,u_m]\otimes\k[v_1,\ldots,v_m]
  \stackrel{d}{\longrightarrow}
  \k[v_1,\ldots,v_m]\stackrel{\varepsilon}{\longrightarrow}\k\to0,
\end{multline}
where $\Lambda^i[u_1,\ldots,u_m]$ is the subspace of
$\Lambda[u_1,\ldots,u_m]$ generated by monomials of length~$i$. We
shall show that the complex above is an exact sequence, or
equivalently, that $\varepsilon\colon(E,d)\to[\k,0]$ is a
quasi-isomorphism. There is an obvious inclusion $\eta\colon\k\to
E$ such that $\varepsilon\eta=\id$. To finish the proof we shall
construct a cochain homotopy between $\id$ and $\eta\varepsilon$,
that is, a set of $\k$-linear maps $s=\{s^{-i,2j}\colon
E^{-i,2j}\to E^{-i-1,2j}\}$ satisfying the identity
\begin{equation}\label{chaineq}
  ds+sd=\id-\eta\varepsilon.
\end{equation}
For $m=1$ we define the map $s_1\colon E_1^{0,*}=\k[v]\to
E_1^{-1,*}$ by the formula
$$
  s_1(a_0+a_1v+\cdots+a_jv^j)=u(a_1+a_2v+\cdots+a_jv^{j-1}).
$$
Then for $f=a_0+a_1v+\cdots+a_jv^j\in E_1^{0,*}$ we have
$ds_1f=f-a_0=f-\eta\varepsilon f$ and $s_1df=0$. On the other
hand, for $uf\in E_1^{-1,*}$ we have $s_1d(uf)=uf$ and
$ds_1(uf)=0$. In any case~\eqref{chaineq} holds. Now we may assume
by induction that for $m=k-1$ the required cochain homotopy
$s_{k-1}\colon E_{k-1}\to E_{k-1}$ is already constructed. Since
$E_k=E_{k-1}\otimes E_1$, \
$\varepsilon_k=\varepsilon_{k-1}\otimes\varepsilon_1$ and
$\eta_k=\eta_{k-1}\otimes\eta_1$, a direct calculation shows that
the map
$$
  s_k=s_{k-1}\otimes\id+\eta_{k-1}\varepsilon_{k-1}\otimes s_1
$$
is a cochain homotopy between $\id$ and $\eta_k\varepsilon_k$.

Since $\Lambda^i[u_1,\ldots,u_m]\otimes\k[m]$ is a free
$\k[m]$-module, \eqref{kosz} is a free resolution for the
$\k[m]$-module~$\k$. It is known as the \emph{Koszul resolution}.
It can be shown to be minimal (an exercise).
\end{construction}

Let~\eqref{resol} be a projective resolution of an $A$-module $M$,
and $N$ is another $A$-module. Applying the functor $\ \otimes_A
N$ to~\eqref{quism} we obtain a homomorphism of cochain complexes
$$
  (R\otimes_AN,d)\to(M\otimes_{A}N,0),
$$
which does not induce a cohomology isomorphism in general. The
$(-i)$th graded cohomology module of the cochain complex
\begin{equation}\label{rotimesn}
  \cdots\to R^{-i}\otimes_{A}N\to \cdots \to
  R^{-1}\otimes_{A}N \to R^0\otimes_{A}N \to 0
\end{equation}
is denoted by $\Tor^{-i}_{A}(M,N)$. We shall also consider the
bigraded $A$-module
$$
  \Tor_{A}(M,N)=\bigoplus_{i,j\ge0}\Tor^{-i,j}_{A}(M,N)
$$
where $\Tor^{-i,j}_{A}(M,N)$ is the $j$th graded component of
$\Tor_{A}^{-i}(M,N)$

The following properties of $\Tor^{-i}_{A}(M,N)$ are well-known
(see e.g.~\cite{macl63}).

\begin{theorem}\label{torprop}\
\begin{itemize}
\item[(a)] The module $\Tor^{-i}_{A}(M,N)$ does not
depend, up to isomorphism, on a choice of
resolution~{\rm(\ref{resol})};

\item[(b)] $\Tor^{-i}_{A}(\:\cdot\:,N)$ and
$\Tor^{-i}_{A}(M,\:\cdot\:)$ are covariant functors;

\item[(c)] $\Tor^0_{A}(M,N)=M\otimes_{A}N$;

\item[(d)] $\Tor^{-i}_{A}(M,N)\cong \Tor^{-i}_{A}(N,M)$;

\item[(e)] An exact sequence of $A$-modules
$$
  0\longrightarrow M_1 \longrightarrow M_2
  \longrightarrow M_3 \longrightarrow 0
$$
induces the following long exact sequence:
\begin{align*}
  \cdots&\longrightarrow\Tor^{-i}_A(M_1,N)\longrightarrow\Tor^{-i}_A(M_2,N)
    \longrightarrow\Tor^{-i}_A(M_3,N)\longrightarrow\cdots\\
  \cdots&\longrightarrow\Tor^{-1}_A(M_1,N)\longrightarrow\Tor^{-1}_A(M_2,N)
    \longrightarrow\Tor^{-1}_A(M_3,N)\\
        &\longrightarrow\Tor^0_A(M_1,N)\longrightarrow\Tor^0_A(M_2,N)
    \longrightarrow\Tor^0_A(M_3,N)\longrightarrow0.
\end{align*}
\end{itemize}
\end{theorem}

In the case $N=\k$ the $\Tor$-modules can be read from a minimal
resolution of $M$ as follows:

\begin{proposition}\label{mintor}
Let $\k$ be a field, and let \eqref{resol} be a minimal resolution
of an $A$-module~$M$. Then
\[
\begin{aligned}
  \Tor^{-i}_A(M,\k)&\cong R_{\min}^{-i}\otimes_A\k,\\
  \dim_\k\Tor^{-i}_A(M,\k)&=\mathop{\mathrm{rank}}R_{\min}^{-i}.
\end{aligned}
\]
\end{proposition}
\begin{proof}
Indeed, the differentials in the cochain complex
$$
  \cdots\longrightarrow R^{-i}_{\min}\otimes_A\k \longrightarrow \cdots
  \longrightarrow R_{\min}^{-1}\otimes_A\k \longrightarrow
  R_{\min}^0\otimes_A\k \longrightarrow 0
$$
are all trivial by Proposition~\ref{mintensor}.
\end{proof}

\begin{corollary}\label{hdtor}
Let $\k$ be a field, and let $M$ be a $A$-module. Then
$$
  \pdim M=\max\bigl\{ i\colon \Tor_A^{-i}(M,\k)\ne0 \bigr\}.
$$
\end{corollary}

\begin{corollary}%[Hilbert Syzygy Theorem]
If $\k$ is a field, then $\pdim M\le m$ for any
$\k[v_1,\ldots,v_m]$-module~$M$.
\end{corollary}
\begin{proof}
By the previous corollary and Theorem~\ref{torprop}~(d),
\[
  \pdim M=\max\bigl\{ i\colon \Tor_{\k[m]}^{-i}(M,\k)\ne0 \bigr\}=\max\bigl\{ i\colon
  \Tor_{\k[m]}^{-i}(\k,M)\ne0\bigr\}.
\]
Using the Koszul resolution for the $\k[m]$-module $\k$ we obtain
\begin{multline*}
  \Tor_{\k[m]}^{-i}(\k,M)=H^{-i}\bigl(\Lambda[u_1,\ldots,u_m]\otimes\k[m]
  \otimes_{\k[m]}M,d\bigr)\\=H^{-i}\bigl(\Lambda[u_1,\ldots,u_m]\otimes
  M,d\bigr).
\end{multline*}
Therefore,
\[
  \pdim M=\max\bigl\{ i\colon \Tor_{\k[m]}^{-i}(\k,M)\ne0\bigr\}\le
  \max\bigl\{ i\colon \Lambda^i[u_1,\ldots,u_m]\otimes
  M\ne0\bigr\}=m.
\]
\end{proof}

\begin{example} Let $A=\k[v_1,\ldots,v_m]$ and $M=N=\k$. By the minimality of the Koszul resolution,
$$
  \Tor_{\k[v_1,\ldots,v_m]}(\k,\k)=\Lambda[u_1,\ldots,u_m]
$$
and $\pdim_A\k=m$.
\end{example}

We note that in general there is no canonical way to define a
multiplication in $\Tor_A(M,N)$, even if both $M$ and $N$ are
$A$-algebras rather than just $A$-modules. However, in the
particular case when $A=\k[m]$, $M$ is an algebra with a unit and
$N=\k$ there is the following canonical way to define a product in
$\Tor_A(M,N)$, extending the previous example. We consider the
differential bigraded algebra $(\Lambda[u_1,\ldots,u_m]\otimes
M,d)$ whose bigrading and differential are defined similarly
to~\eqref{kosdiff}:
\begin{equation}\label{koszuldiff}
\begin{gathered}
  \bideg u_i=(-1,2),\quad\bideg x=(0,\deg x)\quad\text{for }x\in M,\\
  du_i=v_i\cdot 1,\quad dx=0
\end{gathered}
\end{equation}
(here $v_i\cdot 1$ is the element of $M$ obtained by applying
$v_i\in\k[m]$ to $1\in M$, and we identify $u_i$ with
$u_i\otimes1$ and $x$ with $1\otimes x$ for simplicity). Using the
fact that the cohomology of a differential graded algebra is a
graded algebra we obtain:

\begin{lemma}\label{koscom}
Let $M$ be a graded $\k[v_1,\ldots,v_m]$-algebra. Then
$\Tor_{\k[m]}(M,\k)$ is a bigraded $\k$-algebra whose product is
defined via the isomorphism
$$
  \Tor_{\k[v_1,\ldots,v_m]}\bigl(M,\k\bigr)\cong
  H\bigl(\Lambda[u_1,\ldots,u_m]\otimes M,d\bigr).
$$
\end{lemma}
\begin{proof}
Using the Koszul resolution in the definition of
$\Tor_{\k[m]}(\k,M)$ and Theorem~\ref{torprop}~(d) we calculate
\begin{multline*}
  \Tor_{\k[m]}(M,\k)
  \cong\Tor_{\k[m]}(\k,M)\\
  =H\bigl( \Lambda[u_1,\ldots,u_m]\otimes\k[m]
  \otimes_{\k[m]}M,d\bigr)\cong
  H\bigl( \Lambda[u_1,\ldots,u_m]\otimes M,d \bigr).\qedhere
\end{multline*}
\end{proof}

The algebra $\bigl(\Lambda[u_1,\ldots,u_m]\otimes M,d\bigr)$ is
known as the \emph{Koszul algebra}\label{koszulal} (or the
\emph{Koszul complex}) of~$M$.

\begin{remark}
Lemma~\ref{koscom} holds also in the case when $M$ does not have
unit (e.g., when it is a graded ideal in $\k[m]$). In this case
formula~\eqref{koszuldiff} for the differential needs to be
modified as follows:
\[
  d(u_ix)=v_i\cdot x,\quad dx=0\quad\text{for }x\in M.
\]
\end{remark}

We finish this discussion of $\Tor_{\k[m]}(M,\k)$  by mentioning
an important conjecture of commutative homological algebra.

\begin{conjecture}[{Horrocks, see~\cite[p.~453]{bu-ei77}}]
Let $M$ be a graded $\k[v_1,\ldots,v_m]$-module  such that
$\dim_{\k}M<\infty$, where $\k$ is a field. Then
\[
  \dim_{\k}\Tor^{-i}_{\k[v_1,\ldots,v_m]}(M,\k)\ge\bin mi.
\]
\end{conjecture}

It is sometimes formulated in a weaker form:
\begin{conjecture}[weak Horrocks' Conjecture]\label{weakhorrocks}
Let $M$ be a graded $\k[v_1,\ldots,v_m]$-module such that
$\dim_{\k}M<\infty$, where $\k$ is a field. Then
\[
  \dim_{\k}\Tor_{\k[v_1,\ldots,v_m]}(M,\k)\ge 2^m.
\]
\end{conjecture}

If the algebra $A$ is not necessarily commutative, then
$\Tor_A(M,N)$ is defined for a right $A$-module $M$ and a left
$A$-module~$N$ in the same way as above. However, in this case
$\Tor_A(M,N)$ is no longer an $A$-module, and is just a
$\k$-vector space. If both $M$ and $N$ are $A$-bimodules, then
$\Tor_A(M,N)$ is an $A$-bimodule itself.

The construction of $\Tor$ can also be extended to the case of
differential graded modules and algebras, see
Section~\ref{apemss}.

In the standard notation adopted in the algebraic literature, the
modules in a resolution~\eqref{resol} are numbered by nonnegative
rather than nonpositive integers:
\[
  \cdots \stackrel{d}{\longrightarrow} R^{i} \stackrel{d}{\longrightarrow} \cdots
  \stackrel{d}{\longrightarrow} R^{1} \stackrel{d}{\longrightarrow} R^0 \to M \to
  0.
\]
In this notation, the $i$th $\Tor$-module is denoted by
$\Tor^A_i(M,N)$, \ $i\ge0$ (note that~\eqref{rotimesn} becomes a
chain complex, and $\Tor^A_*(M,N)$ is its homology). Therefore,
the two notations are related by
\[
  \Tor_A^{-i}(M,N)=\Tor^A_i(M,N).
\]
Applying the functor $\Hom_A(\ \ ,N)$ to~\eqref{quism} (with
$R^{-i}$ replaced by~$R^i$) we obtain the cochain complex
$$
  0\to\Hom_A(R^0,N)\to\Hom_A(R^1,N)\to\cdots\to\Hom_A(R^i,N)\to\cdots.
$$
Its $i$th cohomology module is denoted by $\Ext^{i}_A(M,N)$.

The properties of the functor $\Ext$ are similar to those given by
Theorem~\ref{torprop} for $\Tor$, with the exception of~(d):

{\samepage
\begin{theorem}
\label{extprop}\
\begin{itemize}
\item[(a)] The module $\Ext^i_{A}(M,N)$ does not depend, up to
isomorphism, on a choice of resolution~{\rm(\ref{resol})};

\item[(b)] $\Ext^i_{A}(\:\cdot\:,N)$ is a contravariant functor,
and $\Ext^i_{A}(M,\:\cdot\:)$ is a covariant functor;

\item[(c)] $\Ext^0_{A}(M,N)=\Hom_A(M,N)$;

\item[(d)] An exact sequence of $A$-modules
$$
  0\longrightarrow M_1 \longrightarrow M_2
  \longrightarrow M_3 \longrightarrow 0
$$
induces the following long exact sequence:
\begin{align*}
  0&\longrightarrow\Ext^0_A(M_3,N)\longrightarrow\Ext^0_A(M_2,N)
    \longrightarrow\Ext^0_A(M_1,N)\\
   &\longrightarrow\Ext^1_A(M_3,N)\longrightarrow\Ext^1_A(M_2,N)
    \longrightarrow\Ext^1_A(M_1,N)\longrightarrow\cdots\\
   \cdots &\longrightarrow\Ext^i_A(M_3,N)\longrightarrow\Ext^i_A(M_2,N)
    \longrightarrow\Ext^i_A(M_1,N)\longrightarrow\cdots;
\end{align*}

\item[(e)] An exact sequence of $A$-modules
$$
  0\longrightarrow N_1 \longrightarrow N_2
  \longrightarrow N_3 \longrightarrow 0
$$
induces the following long exact sequence:
\begin{align*}
  0&\longrightarrow\Ext^0_A(M,N_1)\longrightarrow\Ext^0_A(M,N_2)
    \longrightarrow\Ext^0_A(M,N_3)\\
   &\longrightarrow\Ext^1_A(M,N_1)\longrightarrow\Ext^1_A(M,N_2)
    \longrightarrow\Ext^1_A(M,N_3)\longrightarrow\cdots\\
   \cdots &\longrightarrow\Ext^i_A(M,N_1)\longrightarrow\Ext^i_A(M,N_2)
    \longrightarrow\Ext^i_A(M,N_3)\longrightarrow\cdots.
\end{align*}
\end{itemize}
\end{theorem}
}

\subsection*{Exercises.}
\begin{exercise}
Show that a free resolution exists for every $A$-module~$M$.
(Hint: use the fact that every module is the quotient of a free
module.)
\end{exercise}

\begin{exercise}\label{exprojfree}
If $\k$ is a field and $A=\k[m]$, then every projective graded
$A$-module is free (hint: see \cite[Lemma~VII.6.2]{macl63}). This
is also true in the ungraded case, but is much harder to prove (a
theorem of Quillen and Suslin, settling the famous problem of
Serre). More generally, if $A$ is a finitely generated
nonnegatively graded commutative connected algebra over a
field~$\k$, then every projective $A$-module is free
(see~\cite[Theorem~A3.2]{eise95}). Give an example of a projective
module over a ring which is not free.
\end{exercise}

\begin{exercise}
The Koszul resolution is minimal.
\end{exercise}

\section[Cohen--Macaulay algebras]{Regular sequences and Cohen--Macaulay algebras}\label{cma}
Cohen--Macaulay algebras and modules play an important role in
commutative algebra, algebraic geometry and combinatorics.
%A detailed account of their theory and applications is given
%in~\cite{br-he98}.
Their definition uses the notion of a regular sequence (see
Definition~\ref{regseq} below), which also plays an important role
in algebraic topology, namely in the construction of new
cohomology theories (see~\cite{land76} and Appendix,
Section~\ref{aphige}). In the case of finitely generated algebras
over a field~$\k$, an algebra is Cohen--Macaulay if and only if it
is a free module of finite rank over a polynomial subalgebra.

Here we consider nonnegatively evenly graded finitely generated
commutative connected algebras~$A$ over a field $\k$ and finitely
generated nonnegatively graded $A$-modules $M$ (the case $\k=\Z$
requires extra care, and is treated separately in some particular
cases in the main chapters of the book). The positive part $A^+$
is the unique homogeneous maximal ideal of~$A$, and the results we
discuss here are parallel to those from the homological theory of
Noetherian local rings (we refer to~\cite[Chapters~1--2]{br-he98}
or~\cite[Chapter~19]{eise95} for the details).

Given a sequence of elements $\mb t=(t_1,\ldots,t_k)$ of $A$, we
denote by $A/\mb t$ the quotient of $A$ by the ideal generated by
$\mb t$, and denote by $M/\mb t M$ the quotient of $M$ by the
submodule $t_1M+\cdots+t_k M$. An element $t\in A$ is called a
\emph{zero divisor} on $M$ if $tx=0$ for some nonzero $x\in M$. An
element $t\in A$ is not a zero divisor on $M$ if and only if the
map $M\stackrel{t}\longrightarrow M$ given by multiplication by
$t$ is injective.

\begin{definition}\label{regseq}
Let $M$ be an $A$-module. A homogeneous sequence $\mb
t=(t_1,\ldots,t_k)\in\mathcal H(A^+)$ is called an
$M$-\emph{regular sequence} if $t_{i+1}$ is not a zero divisor on
$M/(t_1M+\cdots+t_iM)$ for $0\le i< k$. We often refer to
$A$-regular sequences simply as \emph{regular}.
\end{definition}

The importance of regular sequences in homological algebra builds
on the fundamental fact that an exact sequence of modules remains
exact after taking quotients by a regular sequence:

\begin{proposition}\label{regex}
Assume given an exact sequence of $A$-modules:
$$
\begin{array}{ccccccccc}
  \cdots \longrightarrow S^i & \stackrel{f_i}\longrightarrow & S^{i-1}
  & \stackrel{f_{i-1}}\longrightarrow & \cdots &
  \stackrel{f_1}\longrightarrow & S^0 & \stackrel{f_0}\longrightarrow & M\to0
\end{array}
$$
If $\mb t$ is an $M$-regular and $S^i$-regular sequence for all
$i\ge0$, then
$$
\begin{array}{ccccccccc}
  \cdots \to S^i/\mb t S^i &
  \stackrel{\overline{f}_i}\longrightarrow & S^{i-1}/\mb t S^{i-1}
  & \stackrel{\overline{f}_{i-1}}\longrightarrow & \cdots &
  \stackrel{\overline{f}_1}\longrightarrow & S^0/\mb t S^0 &
  \stackrel{\overline{f}_0}\longrightarrow & M/\mb t M \to0
\end{array}
$$
is an exact sequence of $A/\mb t$-modules.
\end{proposition}
\begin{proof}
Using induction we reduce the statement to the case when $\mb t$
consists of a single element~$t$. Since
$$
  S^i/t S^i=S^i\otimes_A(A/t),
$$
and $\ \ \otimes_A(A/t)$ is a right exact functor, it is enough to
verify exactness of the quotient sequence starting from the term
$S^1/t S^1$.

Consider the following fragment of the quotient sequence
($i\ge1$):
$$
  S^{i+1}/tS^{i+1}
  \stackrel{\overline{f}_{i+1}}\longrightarrow S^i/tS^i
  \stackrel{\overline{f}_i}\longrightarrow S^{i-1}/tS^{i-1}
  \stackrel{\overline{f}_{i-1}}\longrightarrow S^{i-2}/tS^{i-2}
$$
(where we denote $S^{-1}=M$). For any element $x\in S^i$ we denote
by $\overline{x}$ its residue class in~$S^i/tS^i$. Let
$\overline{f}_i(\overline{x})=0$, then $f_i(x)=ty$ for some $y\in
S^{i-1}$ and $tf_{i-1}(y)=0$. Since $t$ is $S^{i-2}$-regular, we
have $f_{i-1}(y)=0$. Hence, there is $x'\in S^{i}$ such that
$y=f_i(x')$. This implies that $f_i(x-tx')=0$. Therefore,
$x-tx'\in f_{i+1}(S^{i+1})$ and
$\overline{x}\in\overline{f}_{i+1}(S^{i+1}/tS^{i+1})$. Thus, the
quotient sequence is exact.
\end{proof}

The following proposition is often used as the definition of
regular sequences:

\begin{proposition}\label{reg-free}
A sequence $t_1,\ldots,t_k\in\mathcal H(A^+)$ is $M$-regular if
and only if $M$ is a free (not necessarily finitely generated)
$\k[t_1,\ldots,t_k]$-module.
\end{proposition}
\begin{proof}
If $M$ is a free $\k[t_1,\ldots,t_k]$-module, then
$M/(t_1M+\cdots+t_{i-1}M)$ is a free $\k[t_i,\ldots,t_k]$-module,
which implies that $t_i$ is $M/(t_1M+\cdots+t_{i-1}M)$-regular for
$1\le i\le k$. Therefore, $t_1,\ldots,t_k$ is an $M$-regular
sequence.

Conversely, let $\mb t=(t_1,\ldots,t_k)$ be an $M$-regular
sequence. Consider a minimal resolution $(R_{\min},d)$ for the
$\k[\mb t]$-module~$M$. Then, by Proposition~\ref{regex}, the
sequence of $\k$-modules
$$
  \cdots \longrightarrow
  R^{-1}_{\min}/\mb t R^{-1}_{\min} \longrightarrow
  R^0_{\min}/\mb t R^0_{\min} \longrightarrow
  M/\mb t M \longrightarrow 0
$$
is exact. Note that $R^{-i}_{\min}/\mb t R^{-i}_{\min}=
R^{-i}_{\min}\otimes_{\k[\mb t]}\k$. Since the resolution is
minimal, the map $R^0_{\min}\otimes_{\k[\mb t]}\k\to
M\otimes_{\k[\mb t]}\k$ is an isomorphism. Hence,
$R^{-i}_{\min}\otimes_{\k[\mb t]}\k=0$ for $i>0$, which implies
that $R_{\min}^{-i}=0$. Thus, $R_{\min}^0\to M$ is an isomorphism,
i.e. $M$ is a free $\k[\mb t]$-module.
\end{proof}

The following is a direct corollary of Proposition~\ref{reg-free}:

\begin{proposition}\label{rsorder}
The property of being a regular sequence does not depend on the
order of elements in~$\mb t=(t_1,\ldots,t_k)$.
\end{proposition}

%T качестве следствия мv получаем утверждение, которое полезно при
%вvчислении $\Tor$-алгебр симплициальнvх комплексов.

\begin{lemma}
\label{tortor} Let $\mb t$ be a sequence of elements of~$A$ which
is $A$-regular and $M$-regular. Then
$$
  \Tor^{-i}_A(M,\k)=\Tor_{A/\!\mb t\,}^{-i}(M/\mb t M,\k).
$$
\end{lemma}
\begin{proof}
Applying Proposition~\ref{regex} to a minimal resolution of~$M$,
we obtain a minimal resolution of the $A/\mb t$-module $M/\mb t
M$. The rest follows from Proposition~\ref{mintor}.
\end{proof}

An $M$-regular sequence is \emph{maximal} if it is not contained
in an $M$-regular sequence of greater length.

\begin{theorem}[D.~Rees]\label{reesthm}
All maximal $M$-regular sequences in $A$ have the same length
given by
\begin{equation}\label{depth}
  \depth_A M=\min\bigl\{ i\colon\Ext^i_A(\k,M)\ne0 \bigr\}.
\end{equation}
\end{theorem}
This number given by \eqref{depth} is referred to as the
\emph{depth} of~$M$; the simplified notation $\depth M$ will be
used whenever it creates no confusion. The proof of
Theorem~\ref{reesthm} uses the following fact:

\begin{lemma}\label{reghomext}
Let $\mb t=(t_1,\ldots,t_n)\in \mathcal H(A^+)$ be a $M$-regular
sequence. Then
\[
  \Ext^n_A(\k,M)\cong\Hom_A(\k,M/\mb t M).
\]
\end{lemma}
\begin{proof}
We use induction on~$n$. The case $n=0$ is tautological. It
follows from Lemma~\ref{rsorder} that $t_n$ is an $M$-regular
element, so we have the exact sequence
\[
  0\longrightarrow M\stackrel{\cdot t_n}\longrightarrow M
  \longrightarrow M/t_n M\longrightarrow 0.
\]
The map $\Ext^i_A(\k,M)\to\Ext^i_A(\k,M)$ induced by
multiplication by~$t_n$ is zero (an exercise). Therefore, the
second long exact sequence for $\Ext$
%(Theorem~\ref{extprop}~(e))
induced by the short exact sequence above splits into short exact
sequences of the form
\[
  0\longrightarrow \Ext^{n-1}_A(\k,M)\longrightarrow
  \Ext^{n-1}_A(\k,M/t_nM) \longrightarrow
  \Ext^n_A(\k,M)\longrightarrow 0.
\]
Let $\mb t'=(t_1,\ldots,t_{n-1})$. By induction,
\[
  \Ext^{n-1}_A(\k,M)\cong\Hom_A(\k,M/\mb t'M)=0,
\]
where the latter identity follows from Exercise~\ref{0depth},
since $t_n$ is $M/\mb t'M$-regular. Now the exact sequence above
implies that
\[
  \Ext^n_A(\k,M)\cong\Ext^{n-1}_A(\k,M/t_nM)
  \cong\Hom_A(\k,M/\mb t M),
\]
where the latter identity follows by induction.
\end{proof}

\begin{proof}[Proof of Theorem~\ref{reesthm}] Let $\mb
t=(t_1,\ldots,t_n)$ be a maximal $M$-regular sequence. Then, by
Lemma~\ref{reghomext} and Exercise~\ref{0depth},
\[
  \Ext^n_A(\k,M)\cong\Hom_A(\k,M/\mb t M)\ne0,
\]
as $A$ does not contain an $M/\mb t M$-regular element. On the
other hand,
\[
  \Ext^i_A(\k,M)\cong\Hom_A\bigl(\k,M/(t_1 M+\cdots+t_i M)\bigr)=0
\]
for $i<n$, since $t_{i+1}$ is $M/(t_1 M+\cdots+t_i M)$-regular.
\end{proof}

The following fundamental result relates the depth to the
projective dimension:
% of a module.

\begin{theorem}[Auslander--Buchsbaum]\label{abthm}
Let $M\ne0$ be an $A$-module such that $\pdim M<\infty$. Then
$$
  \pdim M+\depth M=\depth A.
$$
\end{theorem}
\begin{proof}
First let $\depth A=0$. Assume that $\pdim M=p>0$. Consider the
minimal resolution for $M$ (which is finite by hypothesis):
$$
\begin{CD}
  0 \to R_{\min}^{-p} @>d_p>> R_{\min}^{-p+1} @>>>
  \cdots @>>> R_{\min}^0 @>>> M \to 0.
\end{CD}
$$
Since $\depth A=0$, we have $\Hom_A(\k,A)=\Ext^0_A(\k,A)\ne0$ by
Theorem~\ref{reesthm}, i.e.
%. Therefore,
there is a monomorphism of
$A$-modules $i\colon\k\to A$. In the commutative diagram
$$
\begin{CD}
  R^{-p}\otimes_A\k @>d_p\otimes_A\k>> R^{-p+1}\otimes_A\k\\
  @V\id\otimes_AiVV                     @VV\id\otimes_AiV\\
  R^{-p}            @>d_p>>            R^{-p+1}
\end{CD}
$$
the maps $d_p$ and $\id\otimes_Ai$ are injective (the latter
because the module $R^{-p}$ is free). Hence $d_p\otimes_A\k$ is
also injective, which contradicts minimality of the resolution. We
obtaine $\pdim M=0$, i.e. $M$ is a free $A$-module and $\depth
M=\depth A=0$.

Now let $\depth A>0$. Assume that $\depth M=0$. Consider the first
syzygy module $M_1=\Ker[R^0\to M]$ for~$M$. It follows
from~\eqref{depth} and the exact sequence for $\Ext$ that $\depth
M_1=1$. Since $\pdim M_1=\pdim M-1$, it is enough to prove the
Auslender--Buchsbaum formula for the module~$M_1$. Hence, we may
assume that $\depth M>0$. This implies that there is an element
$t\in A$ which is $A$-regular and $M$-regular (an exercise). Then
$$
  \depth_{A/t}A/t=\depth_AA-1,\quad
  \depth_{A/t}M/tM=\depth_AM-1
$$
by the definition of depth, and
$$
  \pdim_{A/t}M/tM=\pdim_AM
$$
by Corollary~\ref{hdtor} and Lemma~\ref{tortor}. Now we finish by
induction on~$\depth A$.
\end{proof}

The
%(Krull)
\emph{dimension}\label{dimmodule} of $A$, denoted $\dim A$, is the
maximal number of elements of $A$ algebraically independent
over~$\k$. The \emph{dimension} of an $A$-module $M$ is $\dim
M=\dim(A/\mathop{\mathrm{Ann}}M)$, where
$\mathop{\mathrm{Ann}}M=\{a\in A\colon aM=0\}$ is the
\emph{annihilator} of $M$.

\begin{definition}\label{krdim}
A sequence $t_1,\ldots,t_n$ of algebraically independent
homogeneous elements of $A$ is called a \emph{homogeneous system
of parameters} (briefly \emph{hsop}) for $M$ if $\dim
M/(t_1M+\cdots+t_nM)=0$. Equivalently, $t_1,\ldots,t_n$ is an hsop
if $n=\dim M$ and $M$ is a finitely-generated
$\k[t_1,\ldots,t_n]$-module.

An hsop consisting of linear elements (i.e. elements of lowest
positive degree~2) is referred to as a \emph{linear system of
parameters} (briefly \emph{lsop}).
\end{definition}

The following result (due to Hilbert) is a graded version of the
well-known \emph{Noether Normalisation Lemma}:

\begin{theorem}[{\cite[Theorem~1.5.17]{br-he98}}]
\label{noether} An hsop exists for any $A$-module~$M$. If $\k$ is
an infinite field and $A$ is generated by degree-two elements,
then a lsop can be chosen for~$M$.
\end{theorem}

It is easy to see that a regular sequence consists of
algebraically independent elements, which implies that
$\mathop{\rm depth}M\le\dim M$.

\begin{definition}\label{CM}
$M$ is a \emph{Cohen--Macaulay $A$-module} if $\depth M=\dim M$,
that is, if $A$ contains an $M$-regular sequence $t_1,\ldots,t_n$
of length $n=\dim M$. If $A$ is a Cohen--Macaulay $A$-module, then
it is called a \emph{Cohen--Macaulay algebra}.

By Proposition~\ref{cm-fm}, $A$ is a Cohen--Macaulay algebra if
and only if it is a free finitely generated module over a
polynomial subalgebra.
\end{definition}

\begin{proposition}
\label{rlsop} Let $M$ be a Cohen--Macaulay $A$-module. Then a
sequence $\mb t=(t_1,\ldots,t_k)\in \mathcal H(A^+)$ is
$M$-regular if and only if it is a part of an hsop for~$M$.
\end{proposition}
\begin{proof}
Let $\dim M=n$. Assume that $\mb t$ is an $M$-regular sequence.
The fact that $t_i$ is an $M/(t_1M+\cdots+t_{i-1}M)$-regular
element implies that
$$
  \dim M/(t_1M+\cdots+t_iM)=\dim M/(t_1M+\cdots+t_{i-1}M)-1,\quad
  i=1,\ldots,k
$$
(an exercise). Therefore, $\dim M/\mb tM=n-k$, i.e. $\mb t$ is a
part of an hsop for~$M$.

For the other direction, see~{\cite[Theorem~2.1.2~(c)]{br-he98}}.
\end{proof}

In particular, any hsop in a Cohen--Macaulay algebra~$A$ is
regular.

\begin{proposition}\label{cm-fm}
An algebra is Cohen--Macaulay if and only if it is a free finitely
generated module over a polynomial subalgebra.
\end{proposition}
\begin{proof}
Assume that $A$ is Cohen--Macaulay and $\dim A=n$. Then there is a
regular sequence $\mb t=(t_1,\ldots,t_n)$ in~$A$. By the previous
proposition, $\mb t$ is an hsop, so that $\dim A/\mb t=0$ and
therefore $A$ is a finitely generated $\k[t_1,\ldots,t_n]$-module.
This module is also free by Proposition~\ref{reg-free}.

On the other hand, if $A$ is free finitely generated over
$\k[t_1,\ldots,t_n]$ where $t_1,\ldots,t_n\in A$, then $\dim A=n$
and $t_1,\ldots,t_n$ is a regular sequence by
Proposition~\ref{reg-free}, so that $A$ is Cohen--Macaulay.
\end{proof}

\begin{proposition}\label{pscma}
If $A$ is Cohen--Macaulay with an lsop~$\mb t=(t_1,\ldots,t_n)$,
then there is the following formula for the Poincar\'e series
of~$A$:
\[
  F(A;\lambda)=\frac{F\bigl( A/(t_1,\ldots,t_n);\lambda \bigr)}{(1-\lambda^2)^n},
\]
where $F(A/(t_1,\ldots,t_n);\lambda)$ is a polynomial with
nonnegative integer coefficients.
\end{proposition}
\begin{proof}
Since $A$ is a free finitely generated module over
$\k[t_1,\ldots,t_n]$, we have an isomorphism of $\k$-vector spaces
$A\cong(A/\mb t)\otimes\k[t_1,\ldots,t_n]$. Calculating the
Poincar\'e series of both sides yields the required formula.
\end{proof}

\begin{remark}
If $A$ is generated by its elements $a_1,\ldots,a_n$ of positive
degrees $d_1,\ldots,d_n$ respectively, then it may be shown that
the Poincar\'e series of $A$ is a rational function of the form
\[
  F(A;\lambda)=\frac{P(\lambda)}{(1-\lambda^{d_1})(1-\lambda^{d_2})\cdots(1-\lambda^{d_n})},
\]
where $P(\lambda)$ is a polynomial with integer coefficients.
However, in general the polynomial $P(\lambda)$ cannot be given
explicitly, and some of its coefficients may be negative.
\end{remark}

\subsection*{Exercises}
\begin{exercise}
An $M$-regular sequence consists of algebraically independent
elements.
% over~$\k$.
\end{exercise}

\begin{exercise}
The map $\Ext^i_A(\k,M)\to\Ext^i_A(\k,M)$ induced by
multiplication by an element $x\in\mathcal H(A^+)$ is zero.
\end{exercise}

\begin{exercise}\label{0depth}
The following conditions are equivalent for an $A$-module~$M$:
\begin{itemize}
\item[(a)] Every element of $\mathcal H(A^+)$ is a zero divisor on~$M$,
i.e. $\depth M=0$;
\item[(b)] $\Hom_A(\k,M)\ne0$.
\end{itemize}
(Hint: show that if $\mathcal H(A^+)$ consists of zero divisors
on~$M$ then the ideal $A^+$ annihilates a homogeneous element
of~$M$, see~\cite[Corollary~3.2]{eise95}.)
\end{exercise}

\begin{exercise}
Let $\depth A>0$, let $M$ be an $A$-module with $\depth M=0$, and
let $M_1=\Ker[R^0\to M]$ be the first syzygy module for~$M$. Then
$\depth M_1=1$.
\end{exercise}

\begin{exercise}
If $\depth A>0$ and $\depth M>0$, then there exists an element
$t\in A$ which is $A$-regular and $M$-regular.
\end{exercise}

\begin{exercise}
%The Auslander--Buchsbaum formula (
Theorem~\ref{abthm} does not hold if~$\pdim M=\infty$.
\end{exercise}

\begin{exercise}
Show that $\dim A=0$ if and only if $A$ is finite-dimensional as a
$\k$-vector space. Is it true that $\depth A=0$ implies that
$\dim_\k A$ is finite?
\end{exercise}

\begin{exercise}
Give an example of an algebra $A$ over a field $\k$ of finite
characteristic which is generated by linear elements, but does not
have an lsop.
\end{exercise}

\begin{exercise}
A regular sequence consists of algebraically independent elements.
\end{exercise}

\begin{exercise}
If $t\in\mathcal H(A^+)$ is an $M$-regular element, then $\dim
M/tM=\dim M-1$.
\end{exercise}

\begin{exercise}
Let $\k=\Z$. Show that if $A$ is a free finitely generated module
over a polynomial subalgebra~$\Z[t_1,\ldots,t_k]$ then $t_{i+1}$
is not a zero divisor on $A/(t_1,\ldots,t_i)$ for $0\le i< k$, but
the converse is not true. Therefore, the two possible definitions
of a regular sequence over $\Z$ do not agree. (The reason why
Proposition~\ref{cm-fm} fails over $\Z$ is that minimal
resolutions do not have the required good properties, see the
remark after Construction~\ref{minimal}.)
%Indeed, let $A=\Z[v_1,v_2]/(2v_2)$. Then $v_1$ is a regular
%element, but $A$ is not a free $\Z[v_1]$-module.
\end{exercise}

\section{Formality and Massey products}\label{dgaap}
Here we develop the algebraic formalism used in rational homotopy
theory. We work with dg-algebras $A=\bigoplus_{i\ge0}A^i$ over a
field $\k$ of zero characteristic (usually $\R$ or~$\Q$).
%, and refer to such $A$ as \emph{dg-algebras} for short.
We do not assume $A$ to be finitely generated. Commutativity of
dg-algebras is always understood in the graded sense.

A dg-algebra $A$ is called \emph{homologically
connected}\label{homoconne} if $H^0(A,d)=\k$.
%(Note that a homologically connected dg-algebra is
%not necessarily connected as a graded algebra.)
%A connected dg-algebra $A$ is \emph{simply
%connected} if $H^1[A,d]=0$.

Recall that a homomorphism between dg-algebras $(A,d_A)$ and
$(B,d_B)$ is a $\k$-linear map $f\colon A\to B$ which preserves
degrees, i.e. $f(A^i)\subset B^i$, and satisfies $f(ab)= f(a)f(b)$
and $d_Bf(a)=f(d_Aa)$ for all $a,b\in A$. Such a homomorphism
induces a homomorphism $\widetilde f\colon H(A,d_A)\to H(B,d_B)$
of cohomology algebras. We refer to $f$ as a
\emph{quasi-isomorphism}\label{defiquism} if $\widetilde f$ is an
isomorphism. The equivalence relation generated by
quasi-isomorphisms of dg-algebras is referred to as \emph{weak
equivalence}. Since quasi-isomorphisms are often not invertible, a
weak equivalence between $A$ and $B$ implies only the existence of
a zigzag of quasi-isomorphisms of the form
$$
  A\leftarrow A_1\to A_2\leftarrow A_3\to\cdots\leftarrow A_k\to B.
$$
A dg-algebra $B$ weakly equivalent to $A$ is called a \emph{model}
of $A$. The above `long' zigzag of quasi-isomorphisms can be
reduced to a `short' zigzag $A\leftarrow M\to B$ using the notion
of a minimal model.

\begin{definition}\label{midga}
A commutative dg-algebra $M=\bigoplus_{i\ge0}M^i$ is called
\emph{minimal} (in the sense of Sullivan) if the following three
conditions are satisfied:
\begin{itemize}
\item[(a)] $M^0=\k$ and $d(M^0)=0$;

\item[(b)] $M$ is a free commutative dg-algebra, i.e.
$$
  M=\Lambda[x_k\colon\deg x_k\text{ is odd}]\otimes\k[x_k\colon\deg
  x_k\text{ is even}],
$$
there are only finitely many generators in each degree, and
$$
  \deg x_k\le\deg x_l \quad\text{for }k\le l;
$$

\item[(c)] $dx_k$ is a polynomial
on generators $x_1,\ldots,x_{k-1}$, for each $k\ge1$ (this is
called the \emph{nilpotence condition} on~$d$).
\end{itemize}
\end{definition}

Clearly, a minimal dg-algebra $M$ is \emph{simply connected} (i.e.
$H^1(M,d)=0$)\label{simplycondga} if and only if $M^1=0$. In this
case $\deg x_k\ge2$ for each $k$, and the nilpotence condition is
equivalent to the \emph{decomposability} of~$d$, i.e.
$$
  d(M)\subset M^+\cdot M^+,
$$
where $M^+$ is the subspace generated by elements of positive
degree. For non-simply connected dg-algebras decomposability does
not imply nilpotence: the algebra $\Lambda[x,y]$, $\deg x=\deg
y=1$, with $dx=0$, $dy=xy$ is not minimal.

\begin{definition}\label{minmo}
A minimal dg-algebra $M$ is called a \emph{minimal model} for a
commutative dg-algebra $A$ if there is a quasi-isomorphism
$h\colon M\to A$.
\end{definition}

\begin{theorem}\label{thminmod}\

\begin{itemize}

\item[(a)] For each homologically connected commutative dg-algebra~$A$ satisfying the condition $\dim
H^i[A]<\infty$ for all~$i$, there exists a minimal model~$M_A$,
which is unique up to isomorphism.

\item[(b)] A homomorphism of commutative dg-algebras $f\colon A\to B$ lifts to a
homomorphism $\widehat{f}\colon M_A\to M_B$
% of their minimal models
closing the commutative diagram
$$
\begin{CD}
  M_A @>\widehat{f}>> M_B\\
  @Vh_A VV @VVh_B V\\
  A @>f>> B.
\end{CD}
$$

\item[(c)] If $f$ is a quasi-isomorphism, then $\widehat{f}$ is an isomorphism.
\end{itemize}
\end{theorem}

\begin{remark}
Minimal models can be also defined for dg-algebras~$A$ which do
not satisfy the finiteness condition $\dim H^i[A]<\infty$, but we
shall not need this.
\end{remark}

This theorem is due to Sullivan (simply connected case) and
Halperin (general). A proof can be found
in~\cite[Theorem~II.6]{lehm77} or~\cite[\S12]{f-h-t01}.

\begin{corollary}
A weak equivalence between two commutative dg-algebras $A,B$
satisfying the condition of Theorem~\ref{thminmod} can be
represented by a `short' zigzag $A\gets M\to B$ of
quasi-isomorphisms, where $M$ is the minimal model for $A$
(or~$B$).
\end{corollary}

\begin{definition}\label{formalg}
A dg-algebra is $A$ called \emph{formal} if it is weakly
equivalent to its cohomology $H[A]$ (viewed as a dg-algebra with
zero differential).
\end{definition}

\begin{corollary}
A commutative dg-algebra $A$ with the minimal model $M_A$ is
formal if and only if $M_A$ in formal. In this case there is a
zigzag of quasi-isomorphisms $A\gets M_A\to H[A]$.
\end{corollary}

If $A$ is formal with minimal model $M_A$, then the minimal model
can be recovered from the cohomology algebra $H[A]$ using an
inductive procedure.

\begin{remark}
Even if $A$ is formal, the zigzag $A\gets M_A\to H[A]$ usually
cannot be reduced to a single quasi-isomorphism $A\to H[A]$ or
$H[A]\to A$.
\end{remark}

\begin{example}\label{nonformalg}
Let $M$ be a dg-algebra with three generators $a_1,a_2,a_3$ of
degree~1 and differential given by
$$
  da_1=da_2=0,\quad da_3=a_1a_2.
$$
This $M$ is minimal, but nor formal. Indeed the first degree
cohomology $H^1[M]$ is generated by the classes $\alpha_1$,
$\alpha_2$ corresponding to the cocycles $a_1$,~$a_2$, and we have
$\alpha_1\alpha_2=0$. Assume there is a quasi-isomorphism $f\colon
M\to H[M]$; then we have $f(a_3)=k_1\alpha_1+k_2\alpha_2$ for some
$k_1,k_2\in\mathbf k$. This implies that $f(a_1a_3)=0$, which is
impossible since $a_1a_3$ represents a nontrivial cohomology
class.
\end{example}

We next review Massey products, which provide a simple and
effective tool for establishing nonformality of a dg-algebra.
Massey products constitute a series of higher-order operations (or
\emph{brackets}) in the cohomology of a dg-algebra, with the
second-order operation coinciding with the cohomology
multiplication, while the higher-order brackets are only defined
for certain tuples of cohomology classes. We shall only consider
triple (third-order) Massey products here.

\begin{construction}[triple Massey product]\label{massey}
Let $A$ be a dg-algebra, and let $\alpha_1,\alpha_2,\alpha_3$ be
three cohomology classes such that
$\alpha_1\alpha_2=\alpha_2\alpha_3=0$ in~$H[A]$. Choose their
representing cocycles $a_i\in A^{k_i}$, \ $i=1,2,3$. Since the
pairwise cohomology products vanish, there are elements $a_{12}\in
A^{k_1+k_2-1}$ and $a_{23}\in A^{k_2+k_3-1}$ such that
\[
  da_{12}=a_1a_2 \quad\text{and}\quad da_{23}=a_2a_3.
\]
Then one easily checks that
\[
  (-1)^{k_1+1}a_1a_{23}+a_{12}a_3
\]
is a cocycle in $A^{k_1+k_2+k_3-1}$. Its cohomology class is
called a (triple) \emph{Massey product} of $\alpha_1$, $\alpha_2$
and $\alpha_3$, and denoted by
$\langle\alpha_1,\alpha_2,\alpha_3\rangle$.

More precisely, the Massey product
$\langle\alpha_1,\alpha_2,\alpha_3\rangle$ is the set of all
elements in $H^{k_1+k_2+k_3-1}[A]$ obtained by the above
procedure. Since there are choices of $a_{12}$ and $a_{23}$
involved, the set $\langle\alpha_1,\alpha_2,\alpha_3\rangle$ may
consist of more than one element. In fact, $a_{12}$ is defined up
to addition of a cocycle in $A^{k_1+k_2-1}$, and $a_{23}$ is
defined up to a cocycle in $A^{k_2+k_3-1}$. Therefore, any two
elements in $\langle\alpha_1,\alpha_2,\alpha_3\rangle$ differ by
an element of the subset
\[
  \alpha_1\cdot H^{k_2+k_3-1}[A]+\alpha_3\cdot H^{k_1+k_2-1}[A]
  \subset H^{k_1+k_2+k_3-1}[A],
\]
which is called the \emph{indeterminacy} of the Massey product
$\langle\alpha_1,\alpha_2,\alpha_3\rangle$.

A Massey product $\langle\alpha_1,\alpha_2,\alpha_3\rangle$ is
called \emph{trivial}\label{trivialmp} (or \emph{vanishing}) if it
contains zero. Clearly, a Massey product
$\langle\alpha_1,\alpha_2,\alpha_3\rangle$ is trivial if and only
if its image in the quotient algebra $H[A]/(\alpha_1,\alpha_3)$ is
zero.
\end{construction}

\begin{proposition}\label{mquism}
Let $f\colon A\to B$ be a quasi-isomorphism of dg-algebras. Then
all Massey products in $H[A]$ are trivial if and only they are all
trivial in~$H[B]$.
\end{proposition}
\begin{proof}
Assume that all Massey products in $H[B]$ are trivial. Let
$\langle\alpha_1,\alpha_2,\alpha_3\rangle$ be a Massey product
in~$H[A]$. We define elements $a_1,a_2,a_3,a_{12},a_{23}\in A$ as
in Construction~\ref{massey}, and set $b_i=f(a_i)$,
$b_{ij}=f(a_{ij})$. Let $\beta_i$ denote the cohomology class
corresponding to the cocycle~$b_i$. Since
$\beta_1\beta_2=\widetilde f(\alpha_1\alpha_2)=0$ and
$\beta_2\beta_3=0$, the Massey product
$\langle\beta_1,\beta_2,\beta_3\rangle$ is defined. By the
assumption, $0\in\langle\beta_1,\beta_2,\beta_3\rangle$. This
means that we may choose $b'_{12},b'_{23}\in B$ in such a way that
\[
  b_1b_2=db'_{12},\quad b_2b_3=db'_{23},\quad\text{and}\quad
  (-1)^{k_1+1}b_1b'_{23}+b'_{12}b_3=db_{123}
\]
for some $b_{123}\in B^{k_1+k_2+k_3-2}$. Since $db_{12}=b_1b_2$,
we have $d(b'_{12}-b_{12})=0$. Since $f\colon A\to B$ is a
quasi-isomorphism, there is a cocycle $c_{12}\in A$ such that
$f(c_{12})=b'_{12}-b_{12}$, and similarly there is a cocycle
$c_{23}\in A$ such that $f(c_{23})=b'_{23}-b_{23}$. Set
\[
  a'_{12}=a_{12}+c_{12},\quad a'_{23}=a_{23}+c_{23}.
\]
Then $da'_{12}=da_{12}=a_1a_2$ and similarly $da'_{23}=a_2a_3$.
Therefore, the cohomology class of
$c=(-1)^{k_1+1}a_1a'_{23}+a'_{12}a_3$ is a Massey product of
$\alpha_1,\alpha_2,\alpha_3$. Then
\[
  f(c)=(-1)^{k_1+1}b_1b'_{23}+b'_{12}b_3=db_{123}.
\]
%Using again the fact that $f\colon A\to B$ is a quasi-isomorphism,
Since $f$ is a quasi-isomorphism, $c$ is a coboundary, i.e.
$c=da_{123}$ for some $a_{123}\in A$. Therefore, the Massey
product $\langle\alpha_1,\alpha_2,\alpha_3\rangle$ is trivial.

The fact that the triviality of Massey products in $H[A]$ implies
their triviality in $H[B]$ is proved similarly (an exercise).
\end{proof}

\begin{corollary}
If a dg-algebra $A$ is formal, then all Massey products in $H[A]$
are trivial.
\end{corollary}
\begin{proof}
Apply Proposition~\ref{mquism} to a zigzag $A\gets\cdots\to H[A]$
of quasi-isomorphisms and use that fact that all Massey products
for the dg-algebra $H[A]$ with zero differential are trivial.
\end{proof}

\subsection*{Exercises}
\begin{exercise}
Let $A\to B$ be a quasi-isomorphism of dg-algebras. Show that the
triviality of Massey products in $H[A]$ implies their triviality
in $H[B]$.
\end{exercise}

\begin{exercise}
Find a nontrivial Massey product in the cohomology of the
dg-algebra of Example~\ref{nonformalg}.
\end{exercise}

\chapter{Algebraic topology}\label{algtop}

\section{Homotopy and homology}\label{homothomol}
Here we collect the basic constructions and facts, with almost no
proofs. For the details the reader is referred to standard sources
on algebraic topology, such as~\cite{fo-fu89},~\cite{hatc02}
or~\cite{novi96}. Unlike the other appendices, we do not include
exercises here, as there would be too many of them.

All spaces here are topological spaces, and all maps are
continuous. We denote by $\k$ a commutative ring with unit
(usually $\Z$ or a field; in fact an abelian group is enough until
we start considering the multiplication in cohomology).

\subsection*{Basic homotopy theory}
Two maps $f,g\colon X\to Y$ of spaces are
\emph{homotopic}\label{dehomotopy} (denoted by $f\simeq g$) if
there is a map $F\colon X\times\I\to Y$ (where $\I=[0,1]$ is the
unit interval) such that $F|_{X\!\times0}=f$ and
$F|_{X\!\times1}=g$. We denote the map $F|_{X\!\times\, t}\colon
X\to Y$ by~$F_t$, for~$t\in\I$. Homotopy is an equivalence
relation, and we denote by $[X,Y]$ the set of homotopy classes of
maps from $X$ to~$Y$.

Two spaces $X$ and $Y$ are \emph{homotopy
equivalent}\label{dehoequ} if there are maps $f\colon X\to Y$ and
$g\colon Y\to X$ such that $g\circ f$ and $f\circ g$ are homotopic
to the identity maps of $X$ and $Y$ respectively. The
\emph{homotopy type} of a space $X$ is the class of spaces
homotopy equivalent to~$X$.

A space $X$ is \emph{contractible} if it is homotopy equivalent to
a point.

A \emph{pair}\label{depairsp} $(X,A)$ of spaces consists of a
space $X$ and its subspace~$A$. A \emph{map of pairs}
$f\colon(X,A)\to(Y,B)$ is a continuous map $f\colon X\to Y$ such
that $f(A)\subset B$.

A \emph{pointed space}\label{pointedspace} (or \emph{based space})
is a pair $(X,\pt)$ where $\pt$ is a point of~$X$, called the
\emph{basepoint}. We denote by $X_+$ the pointed space
$(X\sqcup\pt,\pt)$, where $X\sqcup\pt$ is $X$ with a disjoint
point added. A \emph{map of pointed spaces} (or \emph{pointed
map}) is a map of pairs $(X,\pt)\to(Y,\pt)$. We denote a pointed
space $(X,\pt)$ simply by $X$ whenever the choice of the basepoint
is clear or irrelevant. Given two pointed spaces $(X,\pt)$ and
$(Y,\pt)$, their \emph{wedge} (or \emph{bouquet}) is defined as
the pointed space $X\vee Y$ obtained by attaching $X$ and $Y$ at
the basepoints. Then $X\vee Y$ is contained as a pointed subspace
in the product $X\times Y$, and the quotient $X\wedge Y=(X\times
Y)/(X\vee Y)$ is called the \emph{smash product} of $X$ and~$Y$.

For any pointed space $(X,\pt)$ the set of homotopy classes of
pointed maps $(S^k,\pt)\to(X,\pt)$ (where $S^k$ is a
$k$-dimensional sphere, $k\ge0$) is a group for $k>0$, which is
called the \emph{$k$th homotopy group}\label{homotopygroup}
of~$(X,\pt)$ and denoted by $\pi_k(X,\pt)$ or simply~$\pi_k(X)$.
We have that $\pi_0(X)$ is the set of path connected components
of~$X$. The group $\pi_1(X)$ is called the \emph{fundamental
group} of~$X$. The groups $\pi_k(X)$ are abelian for $k>1$. A
pointed map $f\colon X\to Y$ induces a homomorphism
$f_*\colon\pi_k(X)\to\pi_k(Y)$ for each~$k$, and the homomorphisms
induced by homotopic maps are the same.

\medskip

A \emph{locally trivial fibration}\label{ltfibra} (or a
\emph{fibre bundle}) is a quadruple $(E,B,F,p)$ where $E$, $B$,
$F$ are spaces and $p$ is a map $E\to B$ such that for any point
$x\in B$ there is a neighbourhood $U\subset B$ and a homeomorphism
$\varphi\colon p^{-1}(U)\stackrel{\cong}{\longrightarrow}U\times
F$ closing the commutative diagram
\[
  \xymatrix{
  p^{-1}(U) \ar[rr]^{\varphi}\ar[dr]_p && U\times F \ar[dl]\\
  &B
  }
\]
The space $E$ is called the \emph{total space}, $B$ the
\emph{base}, and $F$ the \emph{fibre} of the fibre bundle. The
terms `locally trivial fibration' and `fibre bundle' are also
often used for the map $p\colon E\to B$.

A \emph{cell complex}\label{cellcomplex} (or a \emph{CW-complex})
is a Hausdorff topological space $X$ represented as a union
$\bigcup e_i^q$ of pairwise nonintersecting subsets $e_i^q$,
called \emph{cells}, in such a way that for each cell $e_i^q$
there is a map of a closed $q$-disk $D^q$ to $X$ (the
\emph{characteristic map} of~$e_i^q$) whose restriction to the
interior of $D^q$ is a homeomorphism onto~$e_i^q$. Furthermore,
the following two conditions are assumed:
\begin{itemize}
\item[(C)] the boundary $\bar e_i^q\setminus e_i^q$ of a cell is
contained in a union of finitely many cells of dimensions $<q$;
\item[(W)] a subset $Y\subset X$ is closed if and only the
intersection $Y\cap \bar e_i^q$ is closed for every cell~$e_i^q$
(i.e. the topology of $X$ is the weakest topology in which all
characteristic maps are continuous).
\end{itemize}
The union of cells of $X$ of dimension $\le n$ is called the
\emph{$n$th skeleton} of~$X$ and denoted by~$\sk^n X$ or by~$X^n$.
A cell complex $X$ can be obtained from its 0th skeleton $\sk^0X$
(which is a discrete set) by iterating the operation of
\emph{attaching a cell}: a space $Z$ is obtained from $Y$ by
attaching an $n$-cell along a map $f\colon S^{n-1}\to Y$ if $Z$ is
the pushout of the form
\begin{equation}\label{celliter}
\begin{CD}
  S^{n-1} @>f>> Y\\
  @VVV @VVV\\
  D^n @>>> Z
\end{CD}
\end{equation}
We use the notation $Z=Y\cup_f D^n$.

A \emph{cell subcomplex} of a cell complex $X$ is a closed
subspace which is a union of cells of~$X$. Each skeleton of $X$ is
a cell subcomplex.

A map $f\colon X\to Y$ between cell complexes is called a
\emph{cellular map}\label{cellumap} if $f(\sk^n X)\subset\sk^n Y$
for all~$n$.

\begin{theorem}[Cellular approximation]\label{celluappr}
A map between cell complexes is homotopic to a cellular map.
\end{theorem}

For any integer $n>0$ and any group $\pi$ (which is assumed to be
abelian if $n>1$) there exists a connected cell complex $K(\pi,n)$
such that $\pi_n\bigl(K(\pi,n)\bigr)\cong\pi$ and
$\pi_k\bigl(K(\pi,n)\bigr)=0$ for $k\ne n$. (A cell complex
$K(\pi,n)$ can be constructed by taking the wedge of $n$-spheres
corresponding to a set of generators of~$\pi$ and then killing the
higher homotopy groups by attaching cells.) The space $K(\pi,n)$
is called the \emph{Eilenberg--Mac\,Lane space}\label{Eilespac}
(corresponding to $n$ and~$\pi$), and it is unique up to homotopy
equivalence. Examples of Eilenberg--Mac\,Lane spaces include the
circle $S^1=K(\Z,1)$, the infinite-dimensional real projective
space $\R P^\infty=K(\Z_2,1)$ and the infinite-dimensional complex
projective space $\C P^\infty=K(\Z,2)$.

A space $X$ which is homotopy equivalent to $K(\pi,1)$ for some
$\pi$ is called \emph{aspherical}\label{aspheri}. Any surface (a
closed 2-dimensional manifold) which not a 2-sphere or a real
projective plane is aspherical.

\medskip

A locally trivial fibration $p\colon E\to B$ satisfies the
following \emph{covering homotopy property} (\emph{CHP} for short)
with respect to maps of cell complexes~$Y$: for any homotopy
$F\colon Y\times\I\to B$ and any map $f\colon Y\to E$ such that
$p\circ f=F_0$ there is a covering homotopy $\widetilde F\colon
Y\times\I\to E$, satisfying $\widetilde F_0=f$ and
$p\circ\widetilde F=F$. This is described by the commutative
diagram
\begin{equation}\label{chpdiag}
\xymatrix{
  Y \ar[r]^f \ar[d]_{i_0} & E \ar[d]^p\\
  Y\times\I \ar[r]^{F} \ar@{-->}[ur]^{\widetilde F} & B.
}
\end{equation}

A map $p\colon E\to B$ satisfying the CHP with respect to maps of
cell complexes $Y\to B$ is called a \emph{Serre
fibration}\label{defibrat}. A map $p\colon E\to B$ satisfying CHP
with respect to all maps $Y\to B$ is called a \emph{Hurewitz
fibration}. The difference between Serre and Hurewitz fibrations
is not important for our constructions and will be ignored; we
shall use the term \emph{fibration} for both of them. A fibre
bundle is a fibration, but not every fibration is a fibre bundle.
All fibres $p^{-1}(b)$, $b\in B$, of a fibration are homotopy
equivalent if $B$ is connected.

\begin{theorem}[Exact sequence of fibration]
For a fibration $p\colon E\to B$ with fibre $F$ there exists a
long exact sequence
\[
  \cdots\to \pi_k(F)\stackrel{i_*}\longrightarrow
  \pi_k(E)\stackrel{p_*}\longrightarrow
  \pi_k(B)\stackrel{\partial}\longrightarrow \pi_{k-1}(F)\to\cdots
\]
where the map $i_*$ is induced by the inclusion of the fibre
$i\colon F\to E$, the map $p_*$ is induced by the projection
$p\colon E\to B$, and $\partial$ is  the connecting homomorphism.
\end{theorem}
The connecting homomorphism $\partial$ is defined as follows. Take
an element $\gamma\in\pi_k(B)$ and choose a representative
$g\colon S^k\to B$ in the homotopy class~$\gamma$. By considering
the composition $S^{k-1}\times\I\to S^k\stackrel{g}\longrightarrow
B$ (where the first map contracts the top and bottom bases of the
cylinder to the north and south poles of the sphere), we may view
the map $g$ as a homotopy $F\colon S^{k-1}\times\I\to B$ of the
trivial map $F_0\colon S^{k-1}\to\pt$ with $F_1$ also being
trivial. Using the CHP we lift $F$ to a homotopy $\widetilde
F\colon S^{k-1}\times\I\to E$ with $\widetilde F_0$ still trivial,
but $\widetilde F_1\colon S^{k-1}\to E$ being trivial only after
projecting onto~$B$. The latter condition means that $\widetilde
F_1$ is in fact a map $S^{k-1}\to F$, homotopy class of which we
take for~$\partial\gamma$.

The \emph{path space}\label{pathloopsp} of a space $X$ is the
space $PX$ of pointed maps $(\I,0)\to(X,\pt)$. The \emph{loop
space} $\varOmega X$ is the space of pointed maps
$f\colon(\I,0)\to(X,\pt)$ such that $f(1)=\pt$. The map $p\colon
PX\to X$ given by $p(f)=f(1)$ is a fibration, with fibre (homotopy
equivalent to)~$\varOmega X$.

\begin{proposition}\label{prophf}
For any map $f\colon X\to Y$ there exist a homotopy equivalence
$h\colon X\to \widetilde X$ and a fibration $p\colon\widetilde
X\to Y$ such that $f=p\circ h$. Furthermore, this decomposition is
functorial in the sense that a commutative diagram of maps
\[
  \xymatrix{
  X \ar[r] \ar[d]_f & X' \ar[d]_{f'}\\
  Y \ar[r] & Y'
  }
\]
induces a commutative diagram
\[
  \xymatrix{
  X \ar[rr] \ar[rd]^h \ar[dd]_(0.25)f && X' \ar[dd]_(0.25){f'} \ar[rd]^{h'}\\
  & \widetilde X \ar[ld]^{p} \ar[rr] && \widetilde{X'} \ar[ld]^{p'}\\
  Y \ar[rr] && Y'
  }
\]

\end{proposition}
\begin{proof}
Let $\widetilde X$ be the set of pairs $(x,g)$ consisting of a
point $x\in X$ and a path $g\colon\I\to Y$ with $g(0)=f(x)$.This
is described by the pullback diagram
\[
\xymatrix{
  \widetilde X \ar[r] \ar[d] & Y^{\I} \ar[d]^{p_0}\\
  X \ar[r]^{f} & Y
}
\]
where $Y^{\I}$ is the space of all paths $g\colon\I\to Y$ and the
map $p_0$ takes $g$ to $g(0)$.

Then the homotopy equivalence $h\colon X\to\widetilde X$ is given
by $h(x)=(x,c_{f(x)})$, where $c_{f(x)}\colon\I\to Y$ is the
constant path $t\mapsto f(x)$, and the fibration
$p\colon\widetilde X\to Y$ is given by $p(x,h)=h(1)$. The
functoriality property follows by inspection.
\end{proof}

A space homotopy equivalent to the fibre of the fibration $p\colon
\widetilde X\to Y$ from Proposition~\ref{prophf} is referred to as
the \emph{homotopy fibre}\label{homotopyfibre} of the map~$f\colon
X\to Y$, and denoted by $\mathop{\mathrm{hofib}}f$. The
functoriality of the construction of $\widetilde X$ implies that
the homotopy fibre is well-defined: for any other decomposition
$f=p'\circ h'$ into a composition of a fibration $p'$ and a
homotopy equivalence~$h'$ the homotopy fibre of $f$ is homotopy
equivalent to the fibre of~$p'$. The homotopy fibre of the
inclusion $\pt\to X$ is the loop space~$\varOmega X$.

\medskip

An inclusion $i\colon A\to X$ of a cell subcomplex in a cell
complex~$X$ satisfies the following \emph{homotopy extension
property} (\emph{HEP} for short): for any map $f\colon X\to Y$, a
homotopy $F\colon A\times\I\to Y$ such that $F_0=f|_A$ can be
extended to a homotopy $\widehat F\colon X\times\I\to Y$ such that
$\widehat F_0=\widehat f$ and $\widehat F|_{A\times\I}=F$. This is
described by the commutative diagram
\begin{equation}\label{hepdiag}
\xymatrix{
  A \ar[r]^{F'} \ar[d]_{i} & Y^{\I} \ar[d]^{p_0}\\
  X \ar[r]^{f} \ar@{-->}[ur]^{\widehat F'} & Y
}
\end{equation}
where $F'$ is the \emph{adjoint} of $F$ (i.e. $F'(a)=h$ where
$h(t)=F(a,t)\in Y$), the map $p_0$ takes $h$ to $h(0)$, and
$\widehat F'$ is the adjoint of $\widehat F$.

A map $i\colon A\to X$ satisfying the HEP is called
a\label{decofibration} \emph{cofibration}, and $X/i(A)$ is its
\emph{cofibre}. A pair $(X,A)$ for which $i\colon A\to X$ is a
cofibration is called a \emph{Borsuk pair}.

Diagram~\eqref{chpdiag} expresses the fact that fibrations obey
the \emph{right lifting property}\label{RLPLLP} with respect to
particular cofibrations $Y\to Y\times\I$, which are also homotopy
equivalences. Similarly, diagram~\eqref{hepdiag} expresses the
fact that cofibrations obey the \emph{left lifting property} with
respect to particular fibrations $Y^{\I}\to Y$, which are also
homotopy equivalences. This will be important for axiomatising
homotopy theory via the concept of \emph{model category}, see
Section~\ref{secmc}.

The \emph{cone}\label{conesuspspace} over a space $X$ is the
quotient space $\cone X=(X\times\I)/(X\times1)$. The
\emph{suspension} $\varSigma X$ is the quotient $(\cone
X)/(X\times 0)$. There are inclusions $X\hookrightarrow\cone
X\hookrightarrow\varSigma X$ of closed subspaces, where the first
map is given by $x\mapsto(x,0)$ and the second by
$(x,t)\mapsto\bigl(x,(t+1)/2\bigr)$ for $x\in X$, $t\in\I$. The
inclusion $i\colon X\to\cone X$ is a cofibration, with cofibre
$\varSigma X$.

\begin{proposition}\label{prophcf}
For any map $f\colon X\!\to Y$ there exist a cofibration $i\colon
X\!\to\widehat Y$ and a homotopy equivalence $h\colon\widehat Y\to
Y$ such that $f=h\circ i$. This decomposition is functorial.
\end{proposition}
\begin{proof}
Let $\widehat Y$ be the quotient of $(X\times\I)\sqcup Y$ obtained
by identifying $x\times 0\in X\times\I$ with $f(x)\in Y$. This is
described by the pushout diagram
\[
\xymatrix{
  X \ar[r]^f \ar[d]_{i_0} & Y \ar[d]\\
  X\times\I \ar[r] & \widehat Y
}
\]
where the map $i_0$ takes $x$ to $x\times0$.

The cofibration $i\colon X\to\widehat Y$ is given by $i(x)=x\times
1$, and the homotopy equivalence $h\colon\widehat Y\to Y$ is given
by $h(x\times t)=f(x)$ and $h(y)=y$.
\end{proof}

The space $\widehat Y$\label{dehocofibre} from
Proposition~\ref{prophcf} is known as the \emph{mapping cylinder}
of the map~$f\colon X\to Y$. The space $\widehat Y/i(X)=\widehat
Y/(X\!\times\!1)$ is known as the \emph{mapping cone} of $f\colon
X\to Y$. A space homotopy equivalent to the mapping cone of $f$ is
called the \emph{homotopy cofibre} of the map~$f\colon X\to Y$.
The homotopy cofibre of the projection $X\to\pt$ is the
suspension~$\varSigma X$.

\subsection*{Simplicial homology and cohomology}
Let $\sK$ be a simplicial complex on the set $[m]=\{1,\ldots,m\}$.
An \emph{oriented $q$-simplex} $\sigma$ of~$\sK$ is a $k$-simplex
$I=\{i_1,\ldots,i_{q+1}\}\in\sK$ together with an equivalence
class of total orderings of the set~$I$, two orderings being
equivalent if they differ by an even permutation. Denote by $[I]$
the oriented $q$-simplex with the equivalence class of orderings
containing $i_1<\cdots<i_{q+1}$.

%Here we assume that $\k$ is the group $\Z$ or a field; by a
%$\k$-module we mean an abelian group or a $\k$-vector space, and
%the dimension of a $\k$-module is the rank of an abelian group or
%the dimension of a vector space.

Define the \emph{$q$th simplicial chain group (or module)
$C_q(\sK;\k)$ with coefficients in~$\k$}\label{desimplchain} as
the free $\k$-module with basis consisting of oriented
$q$-simplices of $\sK$, modulo the relations $\sigma+\bar\sigma=0$
whenever $\sigma$ and $\bar\sigma$ are differently oriented
$q$-simplices corresponding to the same simplex of~$\sK$.
Therefore, $C_q(\sK,\k)$ is a free $\k$-module of rank~$f_q(\sK)$
(the number of $q$-simplices of~$\sK$) for $q\ge0$, and
$C_q(\sK;\k)=0$ for $q<0$.

Define the \emph{simplicial chain boundary homomorphisms}
\begin{align*}
  &\partial_q\colon C_q(\sK;\k)\to C_{q-1}(\sK;\k),\quad q\ge1,\\
  &\partial_q[i_1,\ldots,i_{q+1}]=
  \sum_{j=1}^{q+1}(-1)^{j-1}[i_1,\ldots,\widehat
  i_j,\ldots,i_{q+1}],
\end{align*}
where $\widehat i_j$ denotes that $i_j$ is missing. It is easily
checked that $\partial_q\partial_{q+1}=0$, so that
$C_*(\sK;\k)=\{C_q(\sK;\k),\partial_q\}$ is a chain complex,
called the \emph{simplicial chain complex} of~$\sK$ (with
coefficients in~$\k$).

The \emph{$q$th simplicial homology group of $\sK$ with
coefficients in~$\k$}\label{simphomology}, denoted by
$H_q(\sK;\k)$, is defined as the $q$th homology group of the
simplicial chain complex $C_*(\sK;\k)$.

The \emph{Euler characteristic}\label{Eulerchar} $\chi(\sK)$ of
$\sK$ is defined by
\[
  \chi(\sK)=\sum_{q\ge0}(-1)^q\rank H_q(\sK;\k).
\]
It is also given by
\[
  \chi(\sK)=f_0-f_1+f_2-\cdots,
\]
where $f_i$ is the number of $i$-simplices of~$\sK$, and therefore
$\chi(\sK)$ is independent of $\k$.

The \emph{augmented simplicial chain complex}\label{augmcs}
$\widetilde C_*(\sK;\k)$ is obtained by taking into account the
empty simplex $\varnothing\in\sK$. That is, $\widetilde
C_q(\sK;\k)=C_q(\sK;\k)$ for $q\ne-1$ and $\widetilde
C_{-1}(\sK;\k)\cong\k$ is a free $\k$-module with
basis~$[\varnothing]$, so that $\widetilde C_*(\sK;\k)$ is written
as
\[
  \cdots\to C_q(\sK;\k)\stackrel{\partial_q}\longrightarrow\cdots
  \longrightarrow C_1(\sK;\k)\stackrel{\partial_1}\longrightarrow
  C_0(\sK;\k)\stackrel{\varepsilon}\longrightarrow
  \widetilde C_{-1}(\sK;\k)\to0,
\]
where $\varepsilon=\partial_0$ is the \emph{augmentation} taking
each vertex $[i]$ to~$[\varnothing]$.

The $q$th homology group of $\widetilde C_*(\sK;\k)$ is called the
\emph{$q$th reduced simplicial homology group of~$\sK$ with
coefficients in~$\k$}\label{reducedhomology} and is denoted by
$\widetilde H_q(\sK;\k)$. We have that $H_q(\sK;\k)=\widetilde
H_q(\sK;\k)$ for $q\ge1$, and $H_0(\sK;\k)\cong\widetilde
H_0(\sK;\k)\oplus\k$ unless $\sK$ consists of $\varnothing$ only,
in which case $\widetilde H_{-1}(\varnothing;\k)\cong\k$.

The \emph{simplicial cochain complex}\label{simplcochain} of $\sK$
with coefficients in~$\k$ is defined to be
\[
  C^*(\sK;\k)=\Hom_{\k}\bigl(C_*(\sK;\k),\k\bigr).
\]
In explicit terms, $C^*(\sK;\k)=\{C^q(\sK;\k),d_q\}$. Here
$C^q(\sK;\k)$ is the \emph{$q$th simplicial cochain group (or
module) of~$\sK$ with coefficients in~$\k$}; it is a free
$\k$-module with basis consisting of simplicial cochains
$\alpha_I$ corresponding to $q$-simplices $I\in\sK$; the cochain
$\alpha_I$ takes value 1 on the oriented simplex $[I]$ and
vanishes on all other oriented simplices. The value of the
\emph{cochain differential}
\[
  d_q=\partial_{q+1}^*\colon C^q(\sK;\k)\to C^{q+1}(\sK;\k)
\]
on the basis elements is given by
\[
  d\alpha_I=\sum_{j\in[m]\setminus I,\,j\cup I\in\sK}
  \varepsilon(j,j\cup I)\alpha_{j\cup I},
\]
where the sign is given by $\varepsilon(j,j\cup I)=(-1)^{r-1}$ if
$j$ is the $r$th element of the set $j\cup I\subset[m]$, written
in increasing order.

The \emph{$q$th simplicial cohomology group of $\sK$ with
coefficients in~$\k$}\label{simpcohomo}, denoted by $H^q(\sK;\k)$,
is defined as the $q$th cohomology group of the cochain complex
$C^*(\sK;\k)$.

The \emph{$q$th reduced simplicial cohomology group of $\sK$ with
coefficients in~$\k$}\label{reduscoho}, denoted by $\widetilde
H^q(\sK;\k)$, is defined as the $q$th cohomology group of the
cochain complex
\[
  0\to\k\stackrel{d_{-1}}\longrightarrow
  C^0(\sK;\k)\stackrel{d_0}\longrightarrow
  C^1(\sK;\k)\stackrel{d_1}\longrightarrow\cdots
  \longrightarrow
  C^q(\sK;\k)\stackrel{d_q}\longrightarrow\cdots
\]
obtained by applying the functor $\Hom$ to the augmented chain
complex $\widetilde C_*(\sK;\k)$. The map $d_{-1}$ takes
$1=\alpha_{\varnothing}$ to the sum of $\alpha_{\{i\}}$
corresponding to all vertices $\{i\}\in\sK$. We have that
$H^q(\sK;\k)=\widetilde H^q(\sK;\k)$ for $q\ge1$, and
$H^0(\sK;\k)\cong\widetilde H^0(\sK;\k)\oplus\k$ unless $\sK$
consists of $\varnothing$ only, in which case $\widetilde
H^{-1}(\varnothing;\k)\cong\k$.

%We often denote $H_q(X;\Z)$ and $H^q(X;\Z)$ by simply $H_q(X)$ and $H^q(X)$.
%The first nontrivial cohomology group $H^q(X;\Z)$ with $q>0$ is
%torsion-free.

\subsection*{Singular homology and cohomology}
Let $X$ be a topological space. A \emph{singular $q$-simplex}
of~$X$ is a continuous map $f\colon\Delta^q\to X$. The \emph{$q$th
singular chain group $C_q(X;\k)$ of $X$ with coefficients
in~$\k$}\label{singularchain} is the free $\k$-module generated by
all singular $q$-simplices.

For each $i=1,\ldots,q+1$ there is a linear map
$\varphi_q^i\colon\Delta^{q-1}\to\Delta^q$ which sends
$\Delta^{q-1}$ to the face of $\Delta^q$ opposite its $i$th
vertex, preserving the order of vertices. The \emph{$i$th face} of
a singular simplex $f\colon\Delta^q\to X$, denoted by $f^{(i)}$ is
defined to be the singular $(q-1)$-simplex given by the
composition
\[
  f^{(i)}=f\circ\varphi_q^i\colon\Delta^{q-1}\to\Delta^q\to X.
\]
The \emph{singular chain boundary homomorphisms} are given by
\begin{align*}
  &\partial_q\colon C_q(X;\k)\to C_{q-1}(X;\k),\quad q\ge1,\\
  &\partial_q f=
  \sum_{i=1}^{q+1}(-1)^{i-1}f^{(i)}.
\end{align*}
Then $\partial_q\partial_{q+1}=0$, so that
$C_*(X;\k)=\{C_q(X;\k),\partial_q\}$ is a chain complex, called
the \emph{singular chain complex} of~$X$ (with coefficients
in~$\k$).

The \emph{$q$th singular homology group of $X$ with coefficients
in~$\k$}\label{singuhomo}, denoted by $H_q(X;\k)$, is the $q$th
homology group of the singular chain complex $C_*(X;\k)$.

Assume that a space $X$ has only finite number of nontrivial
homology groups (with $\Z$ coefficients), and each of these groups
has finite rank. Then the \emph{Euler
characteristic}\label{eulechara} $\chi(X)$ of $X$ is defined by
\[
  \chi(X)=\sum_{q\ge0}(-1)^q\rank H_q(X;\Z).
\]

\begin{proposition}
The group $H_0(X;\k)$ is a free $\k$-module of rank equal to the
number of path connected components of~$X$.
\end{proposition}

The \emph{augmented singular chain complex}\label{augme3}
$\widetilde C_*(X;\k)$ is defined by $\widetilde
C_q(X;\k)=C_q(X;\k)$ for $q\ne-1$ and $\widetilde
C_{-1}(X;\k)=\k$; the \emph{augmentation}
$\varepsilon=\partial_0\colon C_0(X;\k)\to\k$ is given by
$\varepsilon(f)=1$ for all singular $0$-simplices~$f$.

The $q$th homology group of $\widetilde C_*(X;\k)$ is called the
\emph{$q$th reduced singular homology group of~$\sK$ with
coefficients in~$\k$}\label{redusingu} and is denoted by
$\widetilde H_q(X;\k)$. We have $H_q(X;\k)=\widetilde H_q(X;\k)$
for $q\ge1$, and $H_0(X;\k)\cong\widetilde H_0(X;\k)\oplus\k$ if
$X$ is nonempty.

A nonempty space $X$ is \emph{acyclic}\label{acyclspac}
if~$\widetilde H_q(X;\Z)=0$ for all~$q$.

We shall drop the coefficient group $\k$ in the notation of chains
and homology occasionally.

For any simplicial complex $\sK$, there is an obvious inclusion
$i\colon C_q(\sK)\to C_q(|\sK|)$ of the simplicial chain groups
into the singular chain groups of the geometric realisation~$\sK$,
defining an inclusion of chain complexes.

\begin{theorem}
The map $i\colon C_q(\sK)\to C_q(|\sK|)$ induces an isomorphism
$H_q(\sK)\stackrel{\cong}\longrightarrow H_q(|\sK|)$ between the
simplicial homology groups of $\sK$ and the singular homology
groups of~$|\sK|$.
\end{theorem}

If $(X,A)$ is pair of spaces, then $C_*(A)$ is a chain subcomplex
in $C_*(X)$ (that is, $C_q(A)$ is a submodule of $C_q(X)$ and
$\partial_q C_q(A)\subset C_{q-1}(A)$). The quotient complex
\[
  C_*(X,A)=\bigl\{C_q(X,A)=C_q(X)/C_q(A),\;\;
  \bar\partial_q\colon C_q(X,A)\to C_{q-1}(X,A)\bigr\}
\]
is called the \emph{singular chain complex of the pair $(X,A)$}.
Its $q$th homology group $H_q(X,A)$ is called the \emph{$q$th
singular homology group of~$(X,A)$}\label{homopair}, or
\emph{$q$th relative singular homology group of $X$ modulo~$A$}.
Note that $H_q(X)=H_q(X,\varnothing)$ and $\widetilde
H_q(X)=H_q(X,\pt)$.

The \emph{singular cochain complex}\label{singucochain}
$C^*(X;\k)=\{C^q(X;\k),d_q\}$ is defined to be
\[
  C^*(X;\k)=\Hom_{\k}\bigl(C_*(X;\k),\k\bigr).
\]
Its $q$th cohomology group is called the \emph{$q$th singular
cohomology group of~$X$} and denoted by $H^q(X;\k)$. The
\emph{reduced singular cohomology groups} $\widetilde H^q(X;\k)$
and the \emph{relative cohomology groups} $H^q(X,A)$ are defined
similarly.

The ranks of the groups $H_q(X;\Z)$ and $H^q(X;\Z)$ coincide, and
the number $b^q(X)=\rank H^q(X;\Z)$ is called the \emph{$q$th
Betti number}\label{bettinumb} of~$X$. The Betti numbers with
coefficients in a field $\k$ are defined similarly.

The (co)homology groups have the following fundamental properties.

\begin{theorem}[Functoriality and homotopy invariance] A map $f\colon X\to
Y$ induces homomorphisms $f_*\colon H_q(X;\k)\to H_q(Y;\k)$ and
$f^*\colon H^q(Y;\k)\to H^q(X;\k)$. If two maps $f,g\colon X\to Y$
are homotopic, then $f_*=g_*$ and $f^*=g^*$.
\end{theorem}

We omit the coefficient group $\k$ in the notation for the rest of
this subsection.

\begin{theorem}[Exact sequences of pairs] For any pair $(X,A)$ there are long
exact sequences
\begin{align*}
  &\cdots\to H_q(A)\stackrel{i_*}\longrightarrow
  H_q(X)\stackrel{j_*}\longrightarrow
  H_q(X,A)\stackrel{\partial}\longrightarrow H_{q-1}(A)\to\cdots,\\
  &\cdots\to H^q(X,A)\stackrel{j^*}\longrightarrow
  H^q(X)\stackrel{i^*}\longrightarrow
  H^q(A)\stackrel{d}\longrightarrow H^{q+1}(X,A)\to\cdots
\end{align*}
Furthermore, if the inclusion $A\to X$ is a cofibration (e.g., if
it is an inclusion of a cell subcomplex), then the quotient
projection $(X,A)\to(X/A,\pt)$ induces isomorphisms
\[
  H_q(X,A)\cong H_q(X/A,\pt)=\widetilde H_q(X/A),\quad
  H^q(X,A)\cong H^q(X/A,\pt)=\widetilde H^q(X/A).
\]
\end{theorem}

In the homology exact sequence of a pair, the map $i_*$ is induced
by the inclusion $A\to X$, the map $j_*$ is induced by the map of
pairs $(X,\varnothing)\to(X,A)$, and the \emph{connecting
homomorphism} $\partial$ is the homology homomorphism
corresponding the boundary homomorphism sending a relative cocycle
$c\in C_q(X,A)$ to $\partial c\in C_q(A)$. The maps in the
cohomology exact sequence are defined dually.

Let $(X,A)$ be a pair and $B\subset A$ be a subspace. The
inclusion of pairs $(X\backslash B,A\backslash B)\to(X,A)$ induces
maps
\[
  H_q(X\backslash B,A\backslash B)\to H_q(X,A),
\]
which are referred to as the \emph{excision homomorphisms}.

\begin{theorem}[Excision]\label{excisionthm}
If the closure of $B$ is contained in the interior of~$A$, then
the excision homomorphisms $H_q(X\backslash B,A\backslash B)\to
H_q(X,A)$ are isomorphisms.
\end{theorem}

\begin{theorem}[Mayer--Vietoris exact sequences]\label{maviseq}
Let $X$ be a space, $A\subset X$, $B\subset X$ and $A\cup B=X$.
Assume that the excision homomorphisms
\[
  H_q(B,A\cap B)\to H_q(X,A)\quad\text{and}\quad
  H_q(A,A\cap B)\to H_q(X,B)
\]
are isomorphisms for all~$q$. Then there are exact sequences
\begin{align*}
  &\cdots\to H_q(A\cap B)\stackrel{\beta}\longrightarrow
  H_q(A)\oplus H_q(B)\stackrel{\alpha}\longrightarrow
  H_q(X)\stackrel{\delta}\longrightarrow H_{q-1}(A\cap B)\to\cdots,\\
  &\cdots\to H^q(X)\stackrel{\alpha}\longrightarrow
  H^q(A)\oplus H^q(B)\stackrel{\beta}\longrightarrow
  H^q(A\cap B)\stackrel{\delta}\longrightarrow H^{q+1}(X)\to\cdots.
\end{align*}
\end{theorem}

In the homology Mayer--Vietoris sequence, the map $\beta$ is the
difference of the homomorphism induced by the inclusions $A\cap
B\to A$ and $A\cap B\to B$, the map $\alpha$ is the sum of the
homomorphism induced by the inclusions $A\to X$ and $B\to X$, and
the connecting homomorphism $\delta$ is the composite
\[
  H_q(X)\stackrel{j_*}\longrightarrow H_q(X,A)=H_q(B,A\cap B)
  \stackrel{\partial}\longrightarrow H_{q-1}(A\cap B).
\]

The first key calculation of homology groups follows directly from
the general properties above:

\begin{proposition}
The $n$-disk $D^n$ is acyclic, and the reduced homology of an
$n$-sphere~$S^n$ is given by $\widetilde H_n(S^n;\k)\cong\k$ and
$\widetilde H_i(S^n;\k)=0$ for $i\ne n$.
\end{proposition}

Here is another direct corollary:

\begin{theorem}[Suspension isomorphism]\label{suspiso}
For any space $X$ and any $q>0$ there are isomorphisms
\[
  \widetilde H_q(\varSigma X)\cong\widetilde
  H_{q-1}(X)\quad\text{and}\quad
  \widetilde H^q(\varSigma X)\cong\widetilde
  H^{q-1}(X).
\]
\end{theorem}

A closed connected topological $n$-dimensional manifold $X$ is
\emph{orientable} over~$\k$ if $H_n(X;\k)\cong\k$. Every compact
connected manifold is oriented over $\Z_2$ or any field of
characteristic two.

\begin{theorem}[Poincar\'e duality]\label{pdtheorem} If a closed connected
$n$-manifold is orientable over~$\k$, then $H_q(X;\k)\cong
H^{n-q}(X;\k)$.
\end{theorem}

\subsection*{Cellular homology, and cohomology multiplication}
Let $X$ be a cell complex with skeletons
$\mathop{\mathrm{sk}}^nX=X^n$ for $n=0,1,2,\ldots$

We start with the case $\k=\Z$, and omit this coefficient group in
the notation. The group $\mathcal C_q(X)=H_q(X^q,X^{q-1})$ is
called the \emph{$q$th cellular chain group}\label{celluchain}.
This is a free abelian group with basis the $q$-dimensional cells
of~$X$, and therefore a cellular chain may be viewed as an
integral combination of cells of~$X$.

The \emph{cellular boundary homomorphism} $\partial_q\colon
\mathcal C_q(X)\to\mathcal C_{q-1}(X)$ is defined as the
composition
\[
  \partial_q\colon\mathcal C_q(X)=H_q(X^q,X^{q-1})\to
  H_{q-1}(X^{q-1})\to H_{q-1}(X^{q-1},X^{q-2})=
  \mathcal C_{q-1}(X)
\]
of homomorphisms from the homology exact sequences of pairs
$(X^q,X^{q-1})$ and $(X^{q-1},X^{q-2})$. The resulting chain
complex
\[
  \cdots\to\mathcal C_q(X)\stackrel{\partial_q}\longrightarrow
  \cdots\longrightarrow\mathcal C_1(X)
  \stackrel{\partial_1}\longrightarrow
  \mathcal C_0(X)\to0
\]
is called the \emph{cellular chain complex} of~$X$. Its $q$th
homology group, which we denote by $\mathcal H_q(X)$ for the
moment, is called the \emph{$q$th cellular homology
group}\label{celluhomol} of~$X$.

\begin{theorem}
For any cell complex $X$, there is a canonical isomorphism
$\mathcal H_q(X)\cong H_q(X)$ between the cellular and singular
homology groups.
\end{theorem}

The cellular homology groups $\mathcal H_q(X;\k)$ with
coefficients in~$\k$ and the cohomology groups $\mathcal
H^q(X;\k)$ are defined similarly; they are canonically isomorphic
to the appropriate singular homology and cohomology groups.

We shall therefore not distinguish between the singular and
cellular homology and cohomology groups of cell complexes.

\medskip

The product $X_1\times X_2$ of cell complexes is a cell complex
with cells of the form $e_1\times e_2$, where $e_1$ is a cell
of~$X_1$ and $e_2$ is a cell of~$X_2$.

Given two cellular cochains $c_1\in\mathcal C^{q_1}(X)$ and
$c_2\in\mathcal C^{q_2}(X)$, define their \emph{product} as the
cochain $c_1\times c_2\in\mathcal C^{q_1+q_2}(X\times X)$ whose
value on a cell $e_1\times e_2$ of $X\times X$ is given by
$(-1)^{q_1q_2}c_1(e_1)c_2(e_2)$. This product satisfies the
identity
\[
  \delta(c_1\times c_2)=\delta c_1\times
  c_2+(-1)^{q_1}c_1\times\delta c_2,
\]
where $\delta$ is the cellular cochain differential, and therefore
defines a map of cochain complexes
\[
  \mathcal C^*(X)\otimes\mathcal C^*(X)\to\mathcal C^*(X\times X)
\]
and induces a cohomology map
\[
  \times\colon H^{q_1}(X)\otimes H^{q_2}(X)\to
  H^{q_1+q_2}(X\times X),
\]
which is called the (cohomology) \emph{$\times$-product}.

The composition of the $\times$-product with the cohomology map
induced by the \emph{diagonal map}\label{diagomap}
$\varDelta\colon X\to X\times X$, \ $\varDelta(x)=(x,x)$, defines
a product
\[
  \smallsmile\colon H^{q_1}(X)\otimes H^{q_2}(X)
  \stackrel{\times}\longrightarrow H^{q_1+q_2}(X\times X)
  \stackrel{\varDelta^*}\longrightarrow H^{q_1+q_2}(X),
\]
which is called the \emph{cup product}, or \emph{cohomology
product}. It turns $H^*(X)=\bigoplus_{q\ge0}H^q(X)$ into an
associative and graded commutative ring. This ring structure on
$H^*(X)$ is a homotopy invariant of~$X$; it does not depend on the
cell complex structure. It is also functorial, in the sense that a
map $f\colon X\to Y$ induces a ring homomorphism $f^*\colon
H^*(Y)\to H^*(X)$.

There is also a relative version of the cohomology product, given
by the map
\[
  H^{q_1}(X;A)\otimes H^{q_2}(X,B)
  \longrightarrow H^{q_1+q_2}(X,A\cup B).
\]

\subsection*{Whitehead product, Samelson product and Pontryagin product}
Let $w\colon S^{k+l-1}\to S^k\vee S^l$ be the attaching map of the
$(k+l)$-cell of $S^k\times S^l$ with the standard cell structure.
Explicitly, the map $w$ can be defined as the composition
\[
  S^{k+l-1}=\partial(D^k\times D^l)=
  D^k\times S^{l-1}\cup_{S^{k-1}\times S^{l-1}}S^{k-1}\times
  D^l\to S^k\vee S^l,
\]
where the last map consists of two projections
\begin{align*}
  &D^k\times S^{l-1}\to D^k\to D^k/S^{k-1}=S^k\hookrightarrow
  S^k\vee S^l\quad\text{and}\\
  &S^{k-1}\times D^l\to D^l\to D^l/S^{l-1}=S^l\hookrightarrow
  S^k\vee S^l
\end{align*}
and maps $S^{k-1}\times S^{l-1}$ to the basepoint.

Given two pointed maps $f\colon S^k\to X$ and $g\colon S^l\to X$,
their \emph{Whitehead product}\label{dWhitehead} is defined as the
composition
\[
  [f,g]_w\colon S^{k+l-1}\stackrel w\longrightarrow S^k\vee S^l
  \stackrel{f\vee g}\longrightarrow X.
\]
It gives rise to a well-defined product
\[
  [\;\cdot\,,\,\cdot\,]_w\colon
  \pi_k(X)\times\pi_l(X)\to\pi_{k+l-1}(X),
\]
which is also called the Whitehead product. When $k=l=1$, the
Whitehead product is the commutator product in~$\pi_1(X)$. We have
$[f,g]_w=0$ in $\pi_{k+l-1}(X)$ whenever the map $f\vee g\colon
S^k\vee S^l\to X$ extends to a map $S^k\times S^l\to X$.

\begin{theorem}\
\begin{itemize}
\item[(a)] If $\alpha\in\pi_k(X)$ and $\beta,\gamma\in\pi_l(X)$ with $l>1$, then
\[
  [\alpha,\beta+\gamma]_w=[\alpha,\beta]_w+[\alpha,\gamma]_w.
\]
\item[(b)] If $\alpha\in\pi_k(X)$ and $\beta\in\pi_l(X)$ with $k,l>1$, then
\[
  [\alpha,\beta]_w=(-1)^{kl}[\beta,\alpha]_w.
\]
\item[(c)] If $\alpha\in\pi_k(X)$, $\beta\in\pi_l(X)$ and $\gamma\in\pi_m(X)$ with $k,l,m>1$, then
\[
  (-1)^{km}[[\alpha,\beta]_w,\gamma]_w+(-1)^{lk}[[\beta,\gamma]_w,\alpha]_w+
  (-1)^{ml}[[\gamma,\alpha]_w,\beta]_w=0.
\]
\end{itemize}
\end{theorem}

Now consider the loop space $\varOmega X$. The commutator of
loops, $(x,y)\mapsto xyx^{-1}y^{-1}$, induces a map $c\colon
\varOmega X\wedge\varOmega X\to\varOmega X$. Given two pointed
maps $f\colon S^p\to\varOmega X$ and $g\colon S^q\to\varOmega X$,
their \emph{Samelson product}\label{samproduc} is defined as
\[
  [f,g]_s\colon S^{p+q}=S^p\wedge S^q\stackrel{f\wedge g}\longrightarrow
  \varOmega X\wedge\varOmega X
  \stackrel c\longrightarrow \varOmega X.
\]
It gives rise to a well-defined product
\[
  [\;\cdot\,,\,\cdot\,]_s\colon
  \pi_p(\varOmega X)\times\pi_q(\varOmega X)\to\pi_{p+q}(\varOmega X),
\]
which is also called the Samelson product.

\begin{theorem}\
\begin{itemize}
\item[(a)] If $\varphi\in\pi_p(\varOmega X)$ and $\psi,\eta\in\pi_q(\varOmega X)$, then
\[
  [\varphi,\psi+\eta]_s=[\varphi,\psi]_s+[\varphi,\eta]_s.
\]
\item[(b)] If $\varphi\in\pi_p(\varOmega X)$ and $\psi\in\pi_q(\varOmega X)$, then
\[
  [\varphi,\psi]_s=-(-1)^{pq}[\psi,\varphi]_s.
\]
\item[(c)] If $\varphi\in\pi_p(\varOmega X)$, $\psi\in\pi_q(\varOmega X)$ and $\eta\in\pi_r(\varOmega X)$, then
\[
  [\varphi,[\psi,\eta]_s]_s=[[\varphi,\psi]_s,\eta]_s+(-1)^{pq}[\psi,[\varphi,\eta]_s]_s.
\]
\end{itemize}
\end{theorem}
The Samelson bracket makes the rational vector space
$\pi_*(\varOmega X)\otimes\Q$ into a graded Lie algebra,
see~\eqref{glaide}. It is called the \emph{rational homotopy Lie
algebra}\label{homotLie} of~$X$.

The \emph{Pontryagin product}\label{pontproduct} is defined as the
composition
\[
  H_*(\varOmega X;\k)\otimes H_*(\varOmega X;\k)
  \stackrel\times\longrightarrow H_*(\varOmega X\times\varOmega X;\k)
  \stackrel{m_*}\longrightarrow H_*(\varOmega X;\k),
\]
where $\k$ is a commutative ring with unit, $\times$ is the
homology cross-product, and $m\colon \varOmega X\times\varOmega
X\to\varOmega X$ is the loop multiplication. Pontryagin product
makes $H_*(\varOmega X;\k)$ into an associative (but not generally
commutative) algebra with unit.

Whitehead, Samelson and Pontryagin products are related by the
following classical result of Samelson:

\begin{theorem}[\cite{same53}]
There is a choice of adjunction isomorphism
$t\colon\pi_n(X)\to\pi_{n-1}(\varOmega X)$ such that
\[
  t[\alpha,\beta]_w=(-1)^{k-1}[t\alpha,t\beta]_s
\]
for $\alpha\in\pi_k(X)$ and $\beta\in\pi_l(X)$. Furthermore, if
$h\colon\pi_n(\varOmega X)\to H_n(\varOmega X)$ is the Hurewicz
homomorphism and $\star$ denotes the Pontryagin product, then
\[
  h[\varphi,\psi]_s=h(\varphi)\mathbin\star h(\psi)-(-1)^{pq}
  h(\psi)\mathbin\star h(\varphi)
\]
for $\varphi\in\pi_p(\varOmega X)$ and $\psi\in\pi_q(\varOmega
X)$.
\end{theorem}

It follows that the Pontryagin algebra $H_*(\varOmega X,\Q)$ is
the \emph{universal enveloping algebra} of the graded Lie algebra
$\pi_*(\varOmega X)\otimes\Q$, and $\pi_*(\varOmega X)\otimes\Q$
is the Lie algebra of \emph{primitive elements} in the Hopf
algebra $H_*(\varOmega X,\Q)$.

More details on the relationship between the three products,
including their generalisations to $H$-spaces, can be found in the
monograph by Neisendorfer~\cite{neis10}.

\section{Elements of rational homotopy theory}\label{aprht}
We only review some basis notions and results here. For a detailed
account of this much elaborated theory we refer to~\cite{bo-gu76},
\cite{lehm77} and~\cite{f-h-t01}.

\begin{definition}\label{ratioequ}
A map $f\colon X\to Y$ between spaces is called a \emph{rational
equivalence} if it induces isomorphisms in all rational homotopy
groups, that is, $f_*\colon\pi_i(X)\otimes\Q\to\pi_i(Y)\otimes\Q$
is an isomorphism for all~$i$.

The \emph{rational homotopy type} of $X$ is its equivalence class
in the equivalence relation generated by rational equivalences.
\end{definition}

\begin{proposition}
If both $X$ and $Y$ are simply connected, then $f\colon X\to Y$ is
a rational equivalence if and only if $f_*\colon H_i(X;\Q)\to
H_i(Y;\Q)$ (or equivalently, $f^*\colon H^i(Y;\Q)\to H^i(X;\Q)$)
is an isomorphism for each~$i$.
\end{proposition}

A space $X$ is \emph{nilpotent}\label{nilpospace} if its
fundamental group $\pi_1(X)$ is nilpotent and $\pi_1(X)$ acts
nilpotently on higher homotopy groups $\pi_n(X)$. Simply connected
spaces are obviously nilpotent.

Rational homotopy theory, whose foundation was laid in the work of
Quillen \cite{quil69R} and Sullivan~\cite{sull77}, translates the
study of the rational homotopy type of nilpotent spaces $X$ into
the algebraic language of dg-algebras and minimal models (see
Section~\ref{dgaap}). This translation is made via Sullivan's
\emph{algebra of piecewise polynomial differential
forms}~$\APL(X)$\label{sullivanalg}, whose properties we briefly
discuss below. Further remarks relating rational homotopy theory
to the theory of model categories are given in
Section~\ref{secmc}.

The most basic dg-algebra model of a space $X$ is its singular
cochains $C^*(X;\k)$. However, this dg-algebra is non-commutative,
and therefore is difficult to handle. If $X$ is a smooth manifold,
and $\k=\R$, then we may consider the dg-algebra $\Omega^*(X)$ of
de~Rham differential forms instead of $C^*(X;\R)$. It provides a
functorial (with respect to smooth maps) and commutative
dg-algebra model for~$X$, with the same cohomology $H^*(X;\R)$. It
is therefore natural to ask whether a functorial commutative
dg-algebra model exists for arbitrary cell complexes $X$ over a
field $\k$ of characteristic zero (over finite fields there are
secondary cohomological operations obstructing such a
construction). A first construction of a commutative dg-algebra
model which worked for simplicial complexes was suggested by Thom
in the end of the 1950s. Later, in the mid-1970s, Sullivan
provided a functorial construction of a dg-algebra model $\APL(X)$
whose cohomology is $H^*(X;\Q)$, using similar ideas as those of
Thom (a combinatorial version of differential forms).

The algebra $\APL(X)$ has the following two important properties.

\begin{theorem}[PL de Rham Theorem]\label{plderham}\

\begin{itemize}
\item[(a)] $\APL(X)$ is weakly equivalent to $C^*(X;\Q)$ via a short zigzag of
the form
\[
  \APL(X)\to D(X)\gets C^*(X;\Q),
\]
where $D(X)$ is another naturally defined dg-algebra;

\item[(b)] there is a natural map of cochain complexes $\APL(X)\to C^*(X;\Q)$,
the `Stokes map', which induces an isomorphism in cohomology.

\item[(c)] if $X$ is a smooth manifold, then the dg-algebra
$\Omega^*(X)$ of de~Rham forms is weakly equivalent to
$\APL(X)\otimes_{\Q}\R$.
\end{itemize}
\end{theorem}
The proof of (a) can be found in~\cite[\S{}III.3]{lehm77} or
in~\cite[Corollary~10.10]{f-h-t01}. For~(b),
see~\cite[\S2]{bo-gu76}, and for~(c),
see~\cite[Theorem~V.2]{lehm77} or~\cite[Theorem~11.4]{f-h-t01}.

\begin{definition}\label{minmodelspace}
Given a cell complex $X$, we refer to any commutative dg-algebra
$A$ weakly equivalent to $\APL(X)$ as a (rational) \emph{model}
of~$X$. The minimal model of $\APL(X)$ is called the \emph{minimal
model} of~$X$, and is denoted by~$M_X$.
\end{definition}

In the case of nilpotent spaces, the rational homotopy type of $X$
is fully determined by the commutative dg-algebra $\APL(X)$ or its
minimal model~$M_X$. More precisely, there is the following
fundamental result.

\begin{theorem}\label{rathomgp}\
There is a bijective correspondence between the set of rational
homotopy types of nilpotent spaces and the set of classes of
isomorphic minimal dg-algebras over~$\Q$. Under this
correspondence, there is a natural isomorphism
$$
  \Hom\bigl(\pi_i(X),\Q\bigr)\cong M^i_X\bigl/(M_X^+\cdot M_X^+),
$$
i.e. the rank of the $i$th rational homotopy group of $X$ equals
the number of generators of degree $i$ in the minimal model
of~$X$.
\end{theorem}

For the proof, see~\cite[Theorem~IV.8]{lehm77}.

\begin{definition}\label{formsp}
A space $X$ is \emph{formal} if $C^*(X;\Q)$ is a formal dg-algebra
(equivalently, if $\APL(X)$ is a formal commutative dg-algebra,
see Definition~\ref{formalg}).
\end{definition}

If $X$ in nilpotent, then it is formal if and only if there is a
quasi-isomorphism $M_X\to H^*(X;\Q)$. If $X$ is a smooth manifold,
then it is formal if and only if the de Rham algebra $\Omega^*(X)$
is formal; in this case we can define $M_X$ as the minimal model
of $\Omega^*(X)$ instead of~$\APL(X)$.

\begin{example}\label{formalex}\

1. Let $X=S^{2n+1}$ be an odd sphere, $n\ge1$. Then
\[
  M_X=\Lambda[x],\quad \deg x=2n+1,\quad dx=0.
\]
There is a quasi-isomorphism $M_X\to\Omega^*(X)$ sending $x$ to
the volume form of~$S^{2n+1}$.

2. Let $X=S^{2n}$, $n\ge1$. Then
\[
  M_X=\Lambda[y]\otimes\R[x],\quad \deg x=2n,\,\deg y=4n-1,\quad
  dx=0,\,dy=x^2.
\]
The map $M_X\to\Omega^*(X)$ sends $x$ to the volume form of
$S^{2n}$ and $y$ to zero.

3. Let $X=\C P^n$, $n\ge1$. Then
\[
  M_X=\Lambda[y]\otimes\R[x],\quad  \deg x=2,\,\deg y=2n+1,\quad
  dx=0,\,dy=x^{n+1}.
\]
There is a quasi-isomorphism $M_X\to\Omega^*(X)$ sending $x$ to
the Fubini--Study 2-form
$\omega=\frac{i}{2\pi}\partial\overline\partial\log|\mb
z|^2\in\Omega^2(\C P^n)$ and $y$ to zero.

4. Let $X=\C P^\infty$. Then
\[
  M_X=\Q[v]=H^*(X;\Q),\quad \deg v=2,\quad dv=0.
\]

5. Let $X=\C P^\infty\vee\C P^\infty$. Then
\begin{gather*}
  H^*(X;\Q)=\Q[v_1,v_2]/(v_1v_2),\quad\deg
  v_1=\deg v_2=2,\\
  M_X=\Q[v_1,v_2,w],\quad \deg w=3,\quad dv_1=dv_2=0,\,dw=v_1v_2,
\end{gather*}
and the map $M_X\to H^*(X;\Q)$ sends $w$ to zero.
Theorem~\ref{rathomgp} gives the following nontrivial rational
homotopy groups:
\[
  \pi_2(X)\otimes\Q=\Q^2,\quad\pi_3(X)\otimes\Q=\Q.
\]
This conforms with the homotopy fibration
\[
  S^3\to\C P^\infty\vee\C P^\infty\to\C P^\infty\times\C P^\infty
\]
(see Example~\ref{zk2wedge}).
\end{example}

In all examples above the space $X$ is formal (an exercise). Here
is an example of a nonformal manifold.

\begin{example}
Let $G$ be the \emph{Heisenberg group}\label{heisenb} consisting
of matrices
$$
  \begin{pmatrix} 1&x&z\\0&1&y\\0&0&1\end{pmatrix},\quad x,y,z\in\R,
$$
and let $\Gamma=G_\Z$ be the subgroup consisting of matrices with
$x,y,z\in\Z$. The quotient manifold $X=G/\Gamma$ is the
classifying space for the nilpotent group $\Gamma=\pi_1(X)$, i.e.
$X=K(\Gamma,1)$. The minimal model $M_X$ is generated by the
left-invariant forms
$$
  \omega_1=dx,\quad \omega_2=dy, \quad\omega_3=xdy-dz.
$$
This dg-algebra is isomorphic to the dg-algebra of
Example~\ref{nonformalg}. Therefore, the manifold $X$ is not
formal.
\end{example}

\subsection*{Exercises}

\begin{exercise}
Show that all spaces of Example~\ref{formalex} are formal.
\end{exercise}

\section{Group actions and equivariant cohomology}\label{gractions}
There is a vast literature available on this classical subject; we
mention the monographs of Bredon~\cite{bred72},
Hsiang~\cite{hsia75}, Allday--Puppe~\cite{al-pu93} and
Guillemin--Ginzburg--Karshon~\cite{g-g-k02}, among others. We
briefly review some basic concepts and results used in the main
part of the book.

Let $X$ be a Hausdorff space and $G$ a Hausdorff topological
group. One says that $G$ \emph{acts} on $X$ if for any element
$g\in G$ there is a homeomorphism $\phi_g\colon X\to X$, and the
assignment $g\mapsto\phi_g$ respects the algebraic and topological
structure. In more precise terms, a (left)
\emph{action}\label{defiaction} of $G$ on $X$ is given by a
continuous map
\[
  G\times X \to X,\quad
  (g,x) \mapsto gx
\]
such that $g(hx)=(gh)x$ for any $g,h\in G$, \ $x\in X$, and
$ex=x$, where $e$ is the unit of~$G$. The space $X$ is called a
(left) \emph{$G$-space}. Right actions and right $G$-spaces are
defined similarly. In the case when $G$ is an abelian group, the
notions of left and right action coincide.

A continuous map $f\colon X\to Y$ of $G$-spaces is
\emph{equivariant}\label{equivmap} if it commutes with the group
actions, i.e. $f(gx)=g\cdot f(x)$ for all $g\in G$ and $x\in X$. A
map $f$ is \emph{weakly equivariant} if there is an automorphism
$\psi\colon G\to G$ such that $f(gx)=\psi(g)(f(x))$ for all $g\in
G$ and $x\in X$. A weakly equivariant map corresponding to an
automorphism $\psi$ is also referred to as
\emph{$\psi$-equivariant}.

Let $x\in X$. The set
\[
  G_x=\{g\in G\colon gx=x\}
\]
of elements of $G$ fixing the point $x$ is a closed subgroup
in~$G$, called the \emph{stationary subgroup}, or the
\emph{stabiliser}\label{stasg} of~$x$. The subspace
$$
  Gx=\{gx\in X\colon g\in G\}\subset X
$$
is called the \emph{orbit} of $x$ with respect to the action
of~$G$ (or the \emph{$G$-orbit} for short). If points $x$ and $y$
are in the same orbit, then their stabilisers $G_x$ and $G_y$ are
conjugate subgroups in~$G$. The \emph{type} of an orbit $Gx$ is
the conjugation class of stabiliser subgroups of points in~$Gx$.

The set of all orbits is denoted by~$X/G$, and we have the
canonical projection $\pi\colon X\to X/G$. The space $X/G$ with
the standard quotient topology (a subset $U\subset X/G$ is open if
and only if $\pi^{-1}(U)$ is open) is referred to as the
\emph{orbit space}, or the \emph{quotient
space}\label{orbitspace}. If $G$ is a compact group, then the
quotient $X/G$ is Hausdorff, and the projection $\pi\colon X\to
X/G$ is a closed and proper map (i.e. the image of a closed subset
is closed, and the preimage of a compact subset is compact).

A point $x\in X$ is \emph{fixed}\label{defifixedp} if $Gx=x$, i.e.
$G_x=G$. The set of all fixed points of a $G$-space $X$ will be
denoted by~$X^G$.
%If $X$ is a smooth manifold, $G$ is a compact Lie group, and the
%$G$-action on $X$ is smooth, then each connected component of
%$X^G$ is a submanifold. Furthermore, if $X$ is closed (i.e.
%compact), then $X^G$ consists of finitely many submanifolds.
A $G$-action on $X$ is
\begin{itemize}
\item[--]\emph{effective} if the trivially acting subgroup
$\{g\in G\colon gx=x\text{ for all }x\in X\}$ is trivial (consists
of the single element $e\in G$);

\item[--]\emph{free} if all stabilisers $G_x$ are trivial;

\item[--]\emph{almost free} if all stabilisers $G_x$ are
finite subgroups of~$G$;

\item[--]\emph{semifree} if any stabiliser $G_x$ is either trivial
or is the whole~$G$;

\item[--]\emph{transitive} if for any two points $x,y\in X$ there
is an element $g\in G$ such that $gx=y$ (i.e. $X$ is a single
orbit of the $G$-action).
\end{itemize}

A \emph{principal $G$-bundle}\label{pribu} is a locally trivial
bundle $p\colon X\to B$ such that $G$ acts on $X$ %from the left
preserving fibres, and the induced $G$-action on each fibre is
free and transitive. It follows that each fibre is homeomorphic
to~$G$, the $G$-action on $X$ is free, the $G$-orbits are
precisely the fibres, and the projection $p\colon X\to B$ induces
a homeomorphism between the quotient $X/G$ and the base~$B$.
Therefore $p$ can be regarded as the projection onto the orbit
space of a free $G$-action. If the group $G$ is compact, the
converse is also true  under some mild topological assumptions
on~$X$: a free $G$-space $X$ is the total space of a principal
$G$-bundle (see~\cite[Chapter~II]{bred72}). Therefore, in this
case the notions of a principal $G$-bundle and a free $G$-action
are equivalent.

Now let $G$ be a compact Lie group. Then there exists a principal
$G$-bundle $EG\to BG$ whose total space $EG$ is contractible. This
bundle has the following universality property. Let $E\to B$ be
another principal $G$-bundle over a cell complex~$B$. Then there
is a unique up to homotopy map $f\colon B\to BG$ such that the
pullback of the bundle $EG\to BG$ along $f$ is the bundle $E\to
B$. The space $EG$ is referred to as the \emph{universal
$G$-space}, and the space $BG$ is the \emph{classifying
space}\label{classispa} for free $G$-actions (or simply the
classifying space for~$G$).

Let $X$ be a $G$-space. The \emph{diagonal} $G$-action on
$EG\times X$, given by
$$
  g(e,x)=(ge,gx), \qquad g\in G,\; e\in EG,\; x\in X,
$$
is free. Its orbit space is denoted by $EG\times_G X$ (we shall
also use the notation $B_GX$) and is called the \emph{Borel
construction}, or the \emph{homotopy quotient} of $X$ by~$G$. (The
latter term is used since the free $G$-space $EG\times X$ is
homotopy equivalent to the $G$-space~$X$.) There are two canonical
projections
\begin{equation}\label{2projections}
\begin{aligned}
  EG\times X &\to EG\\
  (e,x) &\mapsto e
\end{aligned}
\qquad\text{and}\qquad
\begin{aligned}
  EG\times X &\to X\\
  (e,x) &\mapsto x.
\end{aligned}
\end{equation}
After taking quotient by the $G$-actions, the second projection
above induces a map $EG\times_G X\to X/G$ between the homotopy and
ordinary quotients, which is a homotopy equivalence when the
$G$-action is free. The first projection gives rise to a bundle
$EG\times_G X\to BG$ with fibre $X$ and structure group~$G$,
called the bundle \emph{associated}\label{associbundle} with the
$G$-space~$X$.

More generally, if $E\to B$ is a principal $G$-bundle (i.e. $E$ is
a free $G$-space) and $X$ is a $G$-space, then we have a bundle
$E\times_G X\to B$ over $B$ with fibre~$X$. When $X$ is an
$n$-dimensional $G$-representation space, we obtain a \emph{vector
bundle}\label{devectorbundle} over~$B$ with structure group~$G$.
\emph{Real}, \emph{oriented} and \emph{complex} $n$-dimensional
(\emph{$n$-plane}) vector bundles correspond to the cases
$G=\mathop{\it GL\/}(n,\R)$, $\mathop{\it SL\/}(n,\R)$ and
$\mathop{\it GL\/}(n,\C)$, respectively. By introducing
appropriate metrics, their structure groups can be reduced to
$\mathop{\it O\/}(n)$, $\mathop{\it SO\/}(n)$ and $\mathop{\it
U\/}(n)$, respectively.

%The second projection induces a map $EG\times_G X\to X/G$. There
%is also a principal $G$-bundle $EG\times X\to EG\times_G X$.

Now let a compact Lie group $G$ act on a smooth manifold $M$ by
diffeomorphisms. Given a point $x\in M$ with stabiliser $G_x$ and
orbit $Gx$, the differential of the action of an element $g\in
G_x$ is a linear transformation of the tangent space $\mathcal
T_xM$ which is the identity on the tangent space to the orbit
of~$x$. Therefore, we obtain the induced representation of $G_x$
in the space $\mathcal T_xM/\mathcal T_x(Gx)$ of the orbit, which
is called the \emph{isotropy representation}\label{isotropyrepr}.
In particular, when $x$ is a fixed point (i.e. $G_x=G$ and
$Gx=x$), we obtain a representation of $G$ in $\mathcal T_xM$,
which is also called the \emph{tangential representation} of $G$
at a fixed point~$x$.

\begin{theorem}[Slice Theorem]\label{slicethm}
Let a compact Lie group $G$ act on a smooth manifold $M$. Then,
for each point $x\in M$, the orbit $Gx$ has a $G$-invariant
neighbourhood $G$-equivariantly diffeomorphic to
$G\times_{G_x}\bigl(\mathcal T_xM/\mathcal T_x(Gx)\bigr)$. The
latter is a vector bundle with fibre $\mathcal T_xM/\mathcal
T_x(Gx)$ over $G/G_x\cong Gx$ and the diffeomorphism takes the
orbit $Gx$ to the zero section of this bundle.
\end{theorem}

The slice theorem was proved by Koszul in the beginning of 1950s.
Its more general version for proper actions of noncompact Lie
groups was proved by Palais; this proof can be found e.g.
in~\cite[Theorem~B.23]{g-g-k02} The slice theorem has many
important consequences. One is that the union of orbits of a given
type is a (smooth, but possibly disconnected) submanifold of~$M$.
In particular, the fixed point set $M^G$ is a submanifold. Another
consequence is that the quotient $M/G$ by a free action of~$G$ is
a smooth manifold.

The second fundamental result is the equivariant embedding
theorem:

\begin{theorem}\label{eqembth}
Let a compact Lie group $G$ act on a compact smooth manifold~$M$.
Then there exist an equivariant embedding of $M$ into a
finite-dimensional linear representation space of~$G$.
\end{theorem}

This theorem was proved by Mostow and Palais in 1957 (instead of
compactness of $M$ they only assumed it to have a finite number of
$G$-orbit types). A simple proof in the case of compact~$M$ (due
to Mostow) can be found in~\cite[Theorem~B.50]{g-g-k02}.

%For each fixed point $x\in X$ there is a section $BG\to EG\times_G
%X$ of the bundle $EG\times_G X\to BG$, and therefore we have an
%inclusion
%\[
%  BG\times X^G\to EG\times_G X.
%\]
%We therefore have an inclusion $X^G\to B_GX$ whose composition
%with the projection $B_GX\to X/G$ coincides with the canonical
%inclusion $X^G\to X/G$.
%
\medskip

The \emph{equivariant cohomology}\label{equcohomolde} of $X$ with
coefficients in a ring~$\k$ is defined by
$$
  H^*_G(X;\k)=H^*(EG\times_G X;\k).
$$
Hence, $H^*_G(\pt;\k)=H^*(BG;\k)$, and the projection $EG\times_G
X\to BG$ turns $H^*_G(X;\k)$ into a $H^*_G(pt;\k)$-module.

For a pair of $G$-spaces $(X,A)$ (where the inclusion $A\subset X$
is an equivariant map), there is a long exact sequence
\[
  \cdots\to H_G^q(X,A)\stackrel{j^*}\longrightarrow
  H_G^q(X)\stackrel{i^*}\longrightarrow
  H_G^q(A)\stackrel{d}\longrightarrow H_G^{q+1}(X,A)\to\cdots
\]
Its existence follows from the exact sequence in ordinary
cohomology and the equivariant homeomorphism $(EG\times
X)/(EG\times A)\cong(EG\times(X/A))/(EG\times\pt)$.

\medskip

A bundle $\pi\colon E\to X$ with fibre $F$ is called a
\emph{$G$-equivariant bundle}\label{equbundle} if $\pi$ is an
equivariant map of $G$-spaces. By applying the Borel construction
to a $G$-equivariant bundle $\pi\colon E\to X$ we obtain a bundle
$B_GE\to B_GX$ with the same fibre~$F$. The \emph{equivariant
characteristic classes} of a $G$-equivariant vector bundle
$\pi\colon E\to X$ are defined as the ordinary characteristic
classes of the corresponding bundle $B_GE\to B_GX$. For example,
the \emph{equivariant Stiefel--Whitney classes} of a
$G$-equivariant vector bundle $\pi\colon E\to X$ belong to
$H^*_G(X;\Z_2)$ and are denoted by $w^G_i(E)$. If $E\to X$ is an
oriented $G$-equivariant vector bundle, then the \emph{equivariant
Euler class}\label{eulereqcoh} $e^G(E)\in H_G^*(X;\Z)$ is defined.
If $E\to X$ is a complex bundle and the $G$-action preserves the
fibrewise complex structure, then the \emph{equivariant Chern
classes} $c^G_i(E)\in H_G^{2i}(X;\Z)$ are defined.

Now let $M$ be a smooth oriented $G$-manifold of dimension~$n$,
where $G$ is a compact Lie group. Let $N\subset M$ be a
$G$-invariant (e.g., fixed) oriented submanifold of
codimension~$k$. We can identify the $G$-equivariant normal bundle
$\nu(N\subset M)$ with a $G$-invariant tubular neighbourhood $U$
of $N$ in~$M$ by means of a $G$-equivariant diffeomorphism. The
same diffeomorphism identifies the Thom space $\Th\nu$ (see
Section~\ref{thomspaces}) of the bundle $\nu$ with the quotient
space $\overline{U}/\partial\overline{U}$. We have the embedding
$i\colon N\subset M$, the projection $\pi\colon U\to N$, and the
\emph{Pontryagin--Thom map} $p\colon M\to\Th\nu$ contracting the
complement $M\setminus U$ to a point. In equivariant cohomology,
the \emph{Thom class}\label{thomcla} $\tau_N\in H^k_G(\Th\nu)$ is
uniquely determined by the identity
$$
  \bigl(\alpha\cdot p^*(\tau_N),\langle M\rangle\bigr)
  =\bigl(i^*\alpha,\langle N\rangle\bigr),
$$
for any $\alpha\in H^{n-k}_G(M)$. Here $\langle M\rangle\in
H^G_n(M)$ denotes the fundamental class of $M$ in equivariant
homology. The \emph{Gysin homomorphism}\label{defingysin} in
equivariant cohomology is defined by the composition
$$
\begin{CD}
  H^{*-k}_G(N) @>\pi^*>> H^{*-k}_G(U) @>\cdot\tau_N>>
  H^*_G(\Th\nu) @>p^*>> H^*_G(M),
\end{CD}
$$
and is denoted by~ $i_*$. Then $i^*(i_*(1))\in H^k_G(N)$ is the
equivariant Euler class of the normal bundle $\nu(N\subset M)$.

\subsection*{Exercises}

\begin{exercise}
Let $G$ be a compact group acting on a Hausdorff space~$X$. Show
that the $G$-orbits are closed, the orbit space $X/G$ is
Hausdorff, and the projection $\pi\colon X\to X/G$ is a closed and
proper map.
\end{exercise}

\begin{exercise}\label{hausclosed}
If the quotient $X/G$ is Hausdorff, then all $G$-orbits are
closed.
\end{exercise}

\begin{exercise}
Give an example of a Hausdorff space $X$ with a $G$-action whose
orbits are closed, but the quotient $X/G$ is not Hausdorff.
\end{exercise}

\begin{exercise}\label{hopfexercise}
Show than any complex line bundle over the complex projective
space $\C P^n$ has the form $S^{2n+1}\times_{S^1}\C$. Here the
$S^1$-action on~$S^{2n+1}$ is standard (the diagonal action on a
unit sphere in~$\C^{n+1}$, with quotient~$\C P^n$), and $\C$ is a
certain 1-dimensional $S^1$-representation space. The same line
bundle can be also given as $(\C^{n+1}\setminus\{\bf
0\})\times_{\C^\times}\C$, where $\C^\times=\mathop{\it
GL\,}(1,\C)$.

In particular, when $\C$ is the standard (weight~1)
$S^1$-representation space (given by $S^1\times\C\to\C$,
$(g,z)\mapsto gz$), the bundle $S^{2n+1}\times_{S^1}\C$ is the
\emph{canonical}\label{canolinebu} line bundle over~$\C P^n$ (it
is also known as the line bundle of hyperplane section and denoted
by~$\mathcal O(1)$ in algebraic geometry literature). The
\emph{Hopf} (or \emph{tautological}) line bundle, whose fibre over
a line $\ell\in\C P^n$ is~$\ell$ itself, corresponds to the
$S^1$-representation of weight~$-1$, given by $(g,z)\mapsto
g^{-1}z$.
\end{exercise}

%\begin{exercise}
%A smooth surjective map of manifolds $M\to N$ whose differential
%is surjective at each point is a locally trivial fibration.
%\end{exercise}
%
\begin{exercise}
Deduce from the slice theorem that if a compact Lie group $G$ acts
smoothly on a smooth manifold $M$, then the union of orbits of a
given type (in particular, the fixed point set $M^G$) is a
submanifold of~$M$. Also, deduce that if the action is free then
the quotient $M/G$ is a smooth manifold.
\end{exercise}

\begin{exercise}\label{faisre}
Show that if a compact Lie group $G$ acts on $M$ effectively and
$M$ is connected, then the isotropy representation
$G_x\to\mathit{GL}\bigl(\mathcal T_xM/\mathcal T_x(Gx)\bigr)$ is
faithful.
\end{exercise}

\section{Eilenberg--Moore spectral sequences}\label{apemss}
In their paper~\cite{ei-mo66} of 1966, Eilenberg and Moore
constructed a spectral sequence, which became one of the important
computational tools of algebraic topology. It particular, it
provides a method for calculation of the cohomology of the fibre
of a bundle $E\to B$ using the canonical $H^*(B)$-module structure
on~$H^*(E)$. This spectral sequence can be considered as an
extension of Adams' approach to calculating cohomology of loop
spaces~\cite{adam56}. In the 1960--70s applications of the
Eilenberg--Moore spectral sequence led to many important results
on cohomology of homogeneous spaces for Lie groups. More recently
it has been used for different calculations with toric spaces.
This section contains the necessary information about the spectral
sequence; we mainly follow L.~Smith's paper~\cite{smit67} in this
description. For a detailed account of differential homological
algebra and the Eilenberg--Moore spectral sequence, as well as its
applications which go beyond the scope of this book, we refer to
McCleary's book~\cite{mccl01}.

Here we assume that $\k$ is a field. The following theorem
provides an algebraic setup for the Eilenberg--Moore spectral
sequence.

\begin{theorem}[{Eilenberg--Moore \cite[Theorem~1.2]{smit67}}]
\label{algemss} Let $A$ be a differential graded $\k$-algebra, and
let $M$, $N$ be differential graded $A$-modules. Then there exists
a spectral sequence $\{E_r,d_r\}$ converging to $\Tor_A(M,N)$ and
whose $E_2$-term is
$$
  E_2^{-i,j}=\Tor^{-i,j}_{H[A]}\bigl(H[M],H[N]\bigr),\quad i,j\ge0,
$$
where $H[\,\cdot\,]$ denotes the algebra or module of cohomology.
\end{theorem}

\begin{remark}
The construction of $\Tor$ for differential graded objects
requires some additional considerations (see e.g.~\cite{smit67}
or~\cite[Chapter~XII]{macl63}).
\end{remark}

The spectral sequence of Theorem~\ref{algemss} lives in the second
quadrant and its differentials $d_r$ add $(r,1-r)$  to the
bidegree, for $r\ge1$. We shall refer to it as the \emph{algebraic
Eilenberg--Moore spectral sequence}\label{Eilemooress}. Its
$E_\infty$-term is expressed via a certain decreasing filtration
$\{F^{-p}\Tor_A(M,N)\}$ in $\Tor_A(M,N)$ by the formula
$$
  E_\infty^{-p,n+p}=F^{-p}\Bigl( \sum_{-i+j=n}\Tor^{-i,j}_A(M,N) \Bigr)
  \Bigl/ F^{-p+1}\Bigl( \sum_{-i+j=n}\Tor^{-i,j}_A(M,N) \Bigr).
$$

Topological applications of Theorem~\ref{algemss} arise in the
case when $A,M,N$ are cochain algebras of topological spaces. The
classical situation is described by the commutative diagram
\begin{equation}
\begin{CD}
  E @>>> E_0\\
  @VVV @VVV\\
  B @>>> B_0,
\end{CD}
\end{equation}
where $E_0\to B_0$ is a Serre fibre bundle with fibre $F$ over a
simply connected base~$B_0$, and $E\to B$ is the pullback along a
continuous map $B\to B_0$. For any space~$X$, let $C^*(X)$ denote
the singular $\k$-cochain algebra of~$X$. Then $C^*(E_0)$ and
$C^*(B)$ are $C^*(B_0)$-modules. Under these assumptions the
following statement holds.

\begin{lemma}[{\cite[Proposition~3.4]{smit67}}]\label{gentoralg}
$\Tor_{C^*(B_0)}(C^*(E_0),C^*(B))$ is a $\k$-algebra in a natural
way, and there is a canonical isomorphism of algebras
$$
  \Tor_{C^*(B_0)}\bigl(C^*(E_0),C^*(B)\bigr)\to H^*(E).
$$
\end{lemma}

Applying Theorem~\ref{algemss} in the case $A=C^*(B_0)$,
$M=C^*(E_0)$, $N=C^*(B)$ and taking into account
Lemma~\ref{gentoralg}, we come to the following statement.

\begin{theorem}[{Eilenberg--Moore}]
\label{topemss} There exists a spectral sequence $\{E_r,d_r\}$ of
commutative algebras converging to $H^*(E)$ with
\[
  E_2^{-i,j}=\Tor^{-i,j}_{H^*(B_0)}\bigl(H^*(E_0),H^*(B)\bigr).
\]
\end{theorem}

The spectral sequence of Theorem~\ref{topemss} is known as the
(topological) \emph{Eilenberg--Moore spectral sequence}. The case
when $B$ is a point is of particular importance, and we state the
corresponding result separately.

\begin{corollary}
\label{onefib} Let $E\to B$ be a fibration over a simply connected
space $B$ with fibre~$F$. Then there exists a spectral sequence
$\{E_r,d_r\}$ of commutative algebras converging to $H^*(F)$ with
\[
  E_2=\Tor_{H^*(B_0)}\bigl(H^*(E_0),\k\bigr).
\]
\end{corollary}

We refer to the spectral sequence of Corollary~\ref{onefib} as the
\emph{Eilenberg--Moore spectral sequence of the fibration $E\to
B$}. In the case when $E_0$ is a contractible we obtain a spectral
sequence converging to the cohomology of the loop space $\Omega
B_0$.

\begin{remark}
As we outlined in Section~\ref{aprht}, the Sullivan algebra
$\APL(X)$ provides a commutative rational model for $X$. It can be
proved~\cite[\S3]{bo-gu76} that the above results on the
Eilenberg--Moore spectral sequence hold over $\Q$ with $C^*$
replaced by~$\APL$. This result is not a direct corollary of
algebraic properties of $\Tor$, since the integration map
$\APL(X)\to C^*(X,\Q)$ is not multiplicative.
\end{remark}

\chapter{Categorical constructions}\label{catconstr}
In this appendix we introduce aspects of category theory that are
directly relevant to the study of toric spaces. Our exposition
follows closely the introductory sections of the work of Panov and
Ray~\cite{pa-ra08}.

\section{Diagrams and model categories}\label{secmc}
We use small capitals to denote categories. The set of morphisms
between objects $c$ and $d$ in a category $\cat{c}$ will be
denoted by $\mathop{\mathrm{Mor}}_{\scat{c}}(c,d)$, or simply by
$\cat{c}(c,d)$. The \emph{opposite category} $\cat{c}^{op}$ has
the same objects with morphisms reverted.
%i.e. $\cat c^{op}(c,d)=\cat{c}(d,c)$.

We shall work with the following categories of combinatorial
origin:
\begin{itemize}
\item[$\cdot$] $\cat{set}$: sets and set maps;
\item[$\cdot$] $\del$: finite ordered sets $[n]$ and nondecreasing maps;
\item[$\cdot$] $\ca(\sK)$: simplices of a finite simplicial complex $\sK$ and
their inclusions (the \emph{face category} of~$\sK$).
\end{itemize}

Here $\ca(\sK)$ is an example of a more general \emph{poset
category}\label{posetcateg} $\cat{p}$, whose objects are elements
$\sigma$ of a poset $(\mathcal P,\le)$ and there is a morphism
$\sigma\to\tau$ whenever $\sigma\le\tau$.

A category is \emph{small} if both its objects and morphisms are
sets. Among the three basic categories above, $\del$ and
$\ca(\sK)$ are small, while $\cat{set}$ is not. Furthermore,
$\ca(\sK)$ is finite.

Given a small category $\cat{s}$ and an arbitrary category
$\cat{c}$, a covariant functor $\mathcal D\colon\cat{s}\to\cat{c}$
is known as an \emph{$\cat{s}$-diagram}\label{dediagram}
in~$\cat{c}$. The source $\cat{s}$ is referred to as the
\emph{indexing category} of the diagram~$\mathcal D$. Such
diagrams are themselves the objects of a \emph{diagram category}
$[\cat{s},\cat{c}]$, whose morphisms are natural transformations.
When $\cat{s}$ is $\del^{op}$, the diagrams are known as
\emph{simplicial objects} in $\cat{c}$, and are written as
$\mathcal D_\bullet$; the object $\mathcal D[n]$ is abbreviated to
$\mathcal D_n$ for every $n\ge 0$, and forms the
\emph{$n$-simplices} of~$\mathcal D_\bullet$. A \emph{simplicial
set} is therefore a simplicial object in $\cat{set}$, i.e. an
object in the diagram category $[\del^{op},\cat{set}]$.
%Motivated by the example $\cat{sset}$ of simplicial sets, we may
%abbreviate the diagram category to $\cat{sr}$ in this case only.

We may interpret every object $c$ of $\cat{c}$ as a \emph{constant
$\cat{s}$-diagram}, and so define the constant functor
$\kappa\colon\cat{c}\rightarrow[\cat{s},\cat{c}]$. Whenever
$\kappa$ admits a right or left adjoint
$[\cat{s},\cat{c}]\rightarrow\cat{c}$, it is known as the
\emph{limit} or \emph{colimit} functor
respectively\label{dfinlimcolim}. In more detail, the limit of a
diagram $\mathcal D\colon\cat{s}\to\cat{c}$ is an object $\lim
\mathcal D$ in $\cat{c}$ for which there are natural
identifications of the morphism sets,
\[
  \mathop{\mathrm{Mor}}\nolimits_{[\scat{s},\scat{c}]}(\kappa(c),\mathcal D)\cong
  \mathop{\mathrm{Mor}}\nolimits_{\scat{c}}(c,\lim \mathcal D)
\]
for any object $c$ of~$\cat{c}$. Similarly, $\colim \mathcal D$
satisfies
\[
  \mathop{\mathrm{Mor}}\nolimits_{\scat{c}}(\colim \mathcal D,c)\cong
  \mathop{\mathrm{Mor}}\nolimits_{[\scat{s},\scat{c}]}(\mathcal D,\kappa(c)).
\]
Simple examples are \emph{products} and more general
\emph{pullbacks}, which are limits over the indexing category
$\bullet\to\bullet\gets\bullet$. Similarly, \emph{coproducts} and
more general \emph{pushouts}\label{depushout} are colimits over
the indexing category $\bullet\gets\bullet\to\bullet$.

For any object $c$ of $\cat{c}$, the objects of the
\emph{overcategory}\label{overundercategory} $\cat{c}\under c$ are
morphisms $f\colon b\rightarrow c$ with fixed target, and the
morphisms are the corresponding commutative triangles; the full
subcategory $\cat{c}\!\Downarrow\!c$ is given by restricting
attention to non-identities~$f$. Similarly, the objects of the
\emph{undercategory} $c\under\cat{c}$ are morphisms $f\colon
c\rightarrow d$ with fixed source, and the morphisms are the
corresponding triangles; $c\!\Downarrow\!\cat{c}$ is given by
restriction to the non-identities. In $\cat{cat}(\sK)$ for
example, we have
\begin{gather*}
\cat{cat}(\sK)\under I=\cat{cat}(\varDelta(I)),\quad
\cat{cat}(\sK)\!\Downarrow\!I=\cat{cat}(\partial\varDelta(I)),\\
I\under\cat{cat}(\sK)=\cat{cat}(\st_{\sK}I),\quad
I\!\Downarrow\!\cat{cat}(\sK)=\cat{cat}(\lk_{\sK}I),
\end{gather*}
for any $I\in\sK$, where $\varDelta(I)$ and $\partial\varDelta(I)$
denote the simplex with vertices~$I$ and its boundary, and star
and link are given in Definition~\ref{deflink}.

A \emph{model category} $\cat{mc}$\label{modelcategory} is a
category which is closed with respect to formation of small limits
and colimits, and contains three \emph{distinguished} classes of
morphisms: \emph{weak equivalencies} $w$, \emph{fibrations} $p$,
and \emph{cofibrations}~$i$. Unless otherwise stated, these
letters denote such morphisms henceforth. A fibration or
cofibration is \emph{acyclic} whenever it is also a weak
equivalence. The three distinguished morphisms are required to
satisfy the following axioms (see
Hirschhorn~\cite[Definition~7.1.3]{hirs03}):
\begin{itemize}
\item[(a)] a retract of a distinguished morphism is a distinguished morphism of the same class;
\item[(b)] if $f'\cdot f$ is a composition of morphisms $f$ and $f'$,
and two of the three morphisms $f$, $f'$ and $f'\cdot f$ is a weak
equivalence, then so is the third;
\item[(c)] acyclic cofibrations obey the left lifting property with respect to fibrations,
and cofibrations obey the left lifting property with respect to
acyclic fibrations;
\item[(d)] every morphism $f$ factorises functorially as
\begin{equation}\label{mcfacts}
h\;=\;p\cdot i\;=\;p'\cdot i',
\end{equation}
for some acyclic $p$ and $i'$.
\end{itemize}
These strengthen Quillen's original axioms for a closed model
category~\cite{quil67} in two minor but significant ways. Quillen
demanded only closure with respect to \emph{finite} limits and
colimits, and only existence of factorisation~\eqref{mcfacts}
rather than its functoriality. When using results of pioneering
authors such as Bousfield and Gugenheim \cite{bo-gu76} and Quillen
\cite{quil69R}, we must take account of these differences.

The axioms for a model category are actually self-dual, in the
sense that any general statement concerning fibrations,
cofibrations, limits, and colimits is equivalent to the statement
in which they are replaced by cofibrations, fibrations, colimits,
and limits respectively. In particular, $\cat{mc}^{op}$ always
admits a dual model structure.

The axioms imply that initial and terminal objects $\circ$ and $*$
exist in $\cat{mc}$, and that $\cat{mc}\!\downarrow\!M$ and
$M\!\downarrow\!\cat{mc}$ inherit model structures for any object
$M$.

An object of $\cat{mc}$ is \emph{cofibrant}\label{cofibobject}
when the natural morphism $\circ\to M$ is a cofibration, and is
\emph{fibrant} when the natural morphism $M\to *$ is a fibration.
A \emph{cofibrant approximation} to an object $N$ is a weak
equivalence $N'\to N$ with cofibrant source, and a \emph{fibrant
approximation} is a weak equivalence $N\to N''$ with fibrant
target. The full subcategories $\cat{mc}_c$, $\cat{mc}_{\!f}$ and
$\cat{mc}_{cf}$ are defined by restricting attention to those
objects of $\cat{mc}$ that are respectively cofibrant, fibrant,
and both. When applied to $\circ\to N$ and $N\to *$, the
factorisations \eqref{mcfacts} determine a \emph{cofibrant
replacement} functor $\omega\colon\cat{mc}\to\cat{mc}_c$, and a
\emph{fibrant replacement} functor
$\phi\colon\cat{mc}\to\cat{mc}_{\!f}$. It follows from the
definitions that $\omega$ and $\phi$ preserve weak equivalences,
and that the associated acyclic fibrations $\omega(N)\to N$ and
acyclic cofibrations $N\to \phi(N)$ form cofibrant and fibrant
approximations respectively.  These ideas are central to the
definition of homotopy limits and colimits given in
Section~\ref{holico} below.

Weak equivalences need not be invertible, so objects $M$ and $N$
are deemed to be \emph{weakly equivalent} if they are linked by a
zigzag
$M\stackrel{e_1}{\longleftarrow}\dots\stackrel{e_n}{\longrightarrow}N$
in $\cat{mc}$; this is the smallest equivalence relation generated
by the weak equivalences. An important consequence of the axioms
is the existence of a localisation functor
$\gamma\colon\cat{mc}\to\Ho(\cat{mc})$, such that $\gamma(w)$ is
an isomorphism in the \emph{homotopy category}\label{homotcatego}
$\Ho(\cat{mc})$ for every weak equivalence~$w$ (i.e.
$\Ho(\cat{mc})$ is obtained from $\cat{mc}$ by inverting all weak
equivalences). Here $\Ho(\cat{mc})$ has the same objects as
$\cat{mc}$, and is equivalent to a category whose objects are
those of $\cat{mc}_{cf}$, but whose morphisms are homotopy classes
of morphisms between them.

Any functor $F$ of model categories that preserves weak
equivalences induces a functor $\Ho(F)$ on their homotopy
categories. Examples include
\begin{equation}\label{hoomfi}
\Ho(\omega)\colon\Ho(\cat{mc})\to\Ho(\cat{mc}_c)\quad\text{and}\quad
\Ho(\phi)\colon\Ho(\cat{mc})\to\Ho(\cat{mc}_{\!f}).
\end{equation}
Such functors often occur as adjoint pairs
\begin{equation}\label{adjpair}
  F\colon\cat{mb}\rightleftarrows\cat{mc}:\! G\;,
\end{equation}
where $F$ is \emph{left Quillen} if it preserves cofibrations and
acyclic cofibrations, and $G$ is \emph{right Quillen} if it
preserves fibrations and acyclic fibrations. Either of these
implies the other, leading to the notion of a \emph{Quillen pair}
$(F,G)$; then Ken Brown's Lemma \cite[Lemma 7.7.1]{hirs03} applies
to show that $F$ and $G$ preserve all weak equivalences on
$\cat{mb}_c$ and $\cat{mc}_{\!f}$ respectively. So they may be
combined with \eqref{hoomfi} to produce an adjoint pair of
\emph{derived functors}
\[
  LF\colon\Ho(\cat{mb})\rightleftarrows\Ho(\cat{mc}):\! RG,
\]
which are equivalences of the homotopy categories (or certain of
their full subcategories) in favourable cases.

Our first examples of model categories are of topological origin,
as follows:
\begin{itemize}
\item[$\cdot$] $\cat{top}$: pointed topological spaces and  pointed continuous maps;
\item[$\cdot$]
$\cat{tmon}$: topological monoids and continuous homomorphisms;
\item[$\cdot$] $\cat{ssets}$: simplicial sets.
\end{itemize}

Homotopy theorists often impose restrictions on topological spaces
defining the category~$\cat{top}$, ensuring that it behaves nicely
with respect to formation of limits and colimits, etc. For
example,~$\cat{top}$ is often defined to consist of
\emph{compactly generated} Hausdorff spaces (a space $X$ is
compactly generated if a subset $A$ is closed whenever the
intersection of~$A$ with any compact subset of $X$ is closed), or
\emph{$k$-spaces}~\cite{vogt71}. We ignore this subtlety however,
as spaces we work with will be nice enough anyway.

There is a model structure on $\cat{top}$ in which fibrations are
Hurewicz fibrations, cofibrations are defined by the homotopy
extension property~\eqref{hepdiag}, and weak equivalences are
homotopy equivalences. However, in the more convenient and
standard model structure on $\cat{top}$, weak equivalences are
maps inducing isomorphisms of homotopy groups, fibrations are
Serre fibrations, and cofibrations obey the left lifting property
with respect to acyclic fibrations (this is a narrower class than
maps obeying the the HEP~\eqref{hepdiag}). The axioms for this
model structure on $\cat{top}$ are verified in~\cite[Theorem
2.4.23]{hove99}, for example.

In either of the model structures above, cell complexes are
cofibrant objects in~$\cat{top}$. Recall that a cell complex can
be defined as a result of iterating the operation of attaching a
cell, i.e. pushing out the standard cofibration $S^{n-1}\to D^n$,
see~\eqref{celliter}. In the second (standard) model structure
on~$\cat{top}$, every space $X$ has a \emph{cellular model}, i.e.
there is a weak equivalence $W\to X$ with $W$ a cell complex,
providing a cofibrant approximation. Two weakly equivalent
topological spaces have homotopy equivalent cellular models.
Cellular models are not functorial, however. A genuine cofibrant
replacement functor $\omega(X)\to X$ must be constructed with
care, and is defined in \cite[\S 98]{dw-sp95}, for example.

A \emph{topological monoid}\label{defntopmonoid} is a space with a
continuous associative product and identity element. We assume
that topological monoids
are pointed by their identities, so that %$\cat{tgp}$ and
$\cat{tmon}$ is a subcategory of $\cat{top}$.  The model structure
for $\cat{tmon}$ is originally due to Schw\"{a}nzl and Vogt
\cite{sc-vo91}, and may also be deduced from Schwede and Shipley's
theory \cite{sc-sh00} of monoids in monoidal model categories;
weak equivalences and fibrations are those homomorphisms which are
weak equivalences and fibrations in \cat{top}, and cofibrations
obey the appropriate lifting property.

In the standard model structure on~$\cat{sset}$, weak equivalences
are maps of simplicial sets whose realisations are weak
equivalences of spaces, fibrations are \emph{Kan fibrations}
(whose realisations are Serre fibrations), and cofibrations are
monomorphisms of simplicial sets. There is a Quillen equivalence
\[
  |\,\cdot\,|\colon\cat{sset}\rightleftarrows\cat{top}:\! S_\bullet\;,
\]
where $|\,\cdot\,|$ denotes the \emph{geometric
realisation}\label{geomrealss} of a simplicial set, and
$S_\bullet$ is the total singular complex of a space. It induces
an equivalence of homotopy categories of simplicial sets and
topological spaces.

\smallskip

Our algebraic categories are defined over arbitrary commutative
rings $\k$, but tend only to acquire model structures when $\k$ is
a field of characteristic zero.
%If $R=\Q$, and in this case only, we omit the subscript from the
%notation.
\begin{itemize}
\item[$\cdot$] $\cat{ch}_\k$ and $\cat{coch}_\k$: augmented chain and cochain
complexes;
\item[$\cdot$] $\cat{cdga}_\k$: commutative augmented differential graded
algebras, with cohomology differential (raising the degree by~1);
\item[$\cdot$]
$\cat{cdgc}_\k$: cocommutative coaugmented differential graded
coalgebras, with homology differential (lowering the degree by~1);
\item[$\cdot$]
$\cat{dga}_\k$: augmented differential graded algebras, with
homology differential;
\item[$\cdot$]
$\cat{dgc}_\k$: coaugmented differential graded coalgebras, with
homology differential;
\item[$\cdot$]
$\cat{dgl}$: differential graded Lie algebras over $\Q$, with
homology differential.
\end{itemize}
For any model structure on these categories, weak equivalences are
the \emph{quasi-isomorphisms}\label{quismmodel}, which induce
isomorphisms in homology or cohomology. The fibrations and
cofibrations are described in Section~\ref{almoca} below. The
augmentations and coaugmentations act as algebraic analogues of
basepoints.

We reserve the notation $\cat{amc}$ for any of the algebraic model
categories above, and assume that objects are graded over the
nonnegative integers. We denote the full subcategory of connected
objects by $\cat{amc}_0$.
%; for $i\geq 0$, we denote the full subcategory of $i$-connected
%objects by $\cat{amc}_i$.
In order to emphasise the differential, we may display an object
$M$ as $(M,d)$. The (co)homology group $H(M,d)$ is also an
$R$-module, and inherits all structure on $M$ except for the
differential. Nevertheless, we may interpret any graded algebra,
coalgebra or Lie algebra as an object of the corresponding
differential category, by imposing $d=0$.

Extending Definition~\ref{formalg}, we refer to an object $(M,d)$
is {\it formal in\/} $\cat{amc}$ whenever there exists a zigzag of
quasi-isomorphisms
\begin{equation}\label{zigzag}
  (M,d)=M_1\stackrel\simeq\longleftarrow\cdots\stackrel\simeq\longrightarrow M_k=(H(M),0).
\end{equation}
Formality only has meaning in an algebraic model category.

There is special class of cofibrant objects in $\cat{cdga}_\Q$,
which are analogous to cell complexes in~$\cat{top}$. These are
minimal dg-algebras, see Definition~\ref{midga}. Any minimal
dg-algebra can be constructed by successive pushouts of the form
\begin{equation}\label{mmiter}
\begin{CD}
  S_\Q(x) @>f>> (A,d_A)\\
  @VVjV @VVV\\
  S_\Q(w,dw) @>>> B=(A\otimes S(w),d_B)
\end{CD}
\end{equation}
in a way similar to constructing cell complexes by
pushouts~\eqref{celliter}. Here $S_\Q(x)$ denotes a free
commutative dg-algebra with one generator $x$ of positive degree
and zero differential, so that $S_\Q(x)$ is the exterior algebra
$\Lambda_{\Q}[x]$ when $\deg x$ is odd and the polynomial algebra
$\Q[x]$ when $\deg x$ is even. The dg-algebra $S_\Q(x,dx)$ has
zero cohomology. The map $j$ is defined by $j(x)=dw$. The
differential in the pushout dg-algebra $B\cong A\otimes S_\Q(w)$
is given by
\[
  d_B(a\otimes1)=d_A a\otimes1, \quad d_B(1\otimes w)=f(x)\otimes 1.
\]
Theorem~\ref{thminmod} asserts the existence of a cofibrant
approximation $f\colon M_A\to A$ for a homologically connected
dg-algebra~$A$, where $M_A$ is a \emph{minimal
model}\label{minmodelmodel} for~$A$; any two minimal models for
$A$ are necessarily isomorphic, and $M_A$ and $M_B$ are isomorphic
for quasi-isomorphic $A$ and $B$. The advantage of $M_A$ is that
it simplifies many calculations concerning $A$; disadvantages
include the fact that it may be difficult to describe for
relatively straightforward objects $A$, and that it cannot be
chosen functorially. A genuine cofibrant replacement functor
requires additional care, and seems first to have been made
explicit in~\cite[\S4.7]{bo-gu76}.

Sullivan's approach to rational homotopy theory is based on the
PL-cochain functor $\APL\colon\cat{top}\rightarrow\cat{cdga}_\Q$.
Basic results related to this approach are given in
Section~\ref{aprht}. Following \cite{f-h-t01}, $\APL(X)$ is
defined as $A^*(S_\bullet X)$, where $S_\bullet(X)$ denotes the
total singular complex of $X$ and
$A^*\colon\cat{sset}\rightarrow\cat{cdga}_\Q$ is the polynomial de
Rham functor of~\cite{bo-gu76}. The PL-de Rham Theorem
(Theorem~\ref{plderham}) yields a natural isomorphism
$H(\APL(X))\rightarrow H^*(X,\Q)$, so $\APL(X)$ provides a
commutative replacement for rational singular cochains, and $\APL$
descends to homotopy categories. Bousfield and Gugenheim prove
that it restricts to an equivalence of appropriate full
subcategories of $\Ho(\cat{top})$ and $\Ho(\cat{cdga}_\Q)$. In
other words, it provides a contravariant algebraic model for the
rational homotopy theory of well-behaved spaces.

Quillen's approach involves the homotopy groups $\pi_*(\varOmega
X)\otimes_\Z\Q$, which form the \emph{rational homotopy Lie
algebra of $X$} under the \emph{Samelson
product}\label{quilsamel}. He constructs a covariant functor
$Q\colon\cat{top}\rightarrow\cat{dgl}$, and a natural isomorphism
\[
H[Q(X)]\stackrel{\cong}{\longrightarrow}\pi_*(\varOmega
X)\otimes_{\Z}\Q.
\]
for any simply connected $X$. He concludes that $Q$ passes to an
equivalence of homotopy categories; in other words, its derived
functor provides a covariant algebraic model for the rational
homotopy theory of simply connected spaces.

%The two approaches are Eckmann--Hilton dual, but the details are
%subtle. Each has enabled important calculations, leading to the
%solution of significant geometric problems. For examples, and
%further details, we refer readers to \cite{f-h-t01}.

We recall from Definition~\ref{formsp} that a space $X$ is
\emph{formal} when $\APL(X)$ is formal in $\cat{cdga}_\Q$. A space
$X$ is referred to as \emph{coformal}\label{defncoforma} when
$Q(X)$ is formal in~$\cat{dgl}$.

\smallskip

The importance of categories of simplicial objects is due in part
to the structure of the indexing category $\del^{op}$. Every
object $[n]$ has \emph{degree} $n$, and every morphism may be
factored uniquely as a composition of morphisms that raise and
lower degree. These properties are formalised in the notion of a
\emph{Reedy category}\label{defnreedy} $\cat{a}$, which admits
generating subcategories $\cat{a}_+$ and $\cat{a}_-$ whose
non-identity morphisms raise and lower degree respectively.  The
diagram category $\fcat{a}{mc}$ then supports a canonical model
structure of its own \cite[Theorem 15.3.4]{hirs03}. By duality,
$\cat{a}^{op}$ is also Reedy, with
$(\cat{a}^{op})_+=(\cat{a}_-)^{op}$ and vice-versa. A simple
example is provided by $\cat{cat}(\sK)$, whose degree function
assigns the dimension $|I|-1$ to each simplex $I$ of~$\sK$. So
$\cat{cat}_+(\sK)$ is the same as $\cat{cat}(\sK)$, and
$\cat{cat}_-(\sK)$ consists entirely of identities.

In the Reedy model structure on $\fcat{cat$(\sK)$}{mc}$, weak
equivalences $w\colon \mathcal C\rightarrow \mathcal D$ are given
\emph{objectwise}, in the sense that $w(I)\colon \mathcal
C(I)\rightarrow \mathcal D(I)$ is a weak equivalence in $\cat{mc}$
for every $I\in\sK$. Fibrations are also objectwise. To describe
the cofibrations, we restrict $\mathcal C$ and $\mathcal D$ to the
overcategories
$\cat{cat}(\sK)\!\Downarrow\!I=\cat{cat}(\partial\varDelta(I))$,
and write $L_I \mathcal C$ and $L_I \mathcal D$ for their
respective colimits. So $L_I$ is the \emph{latching
functor}\label{delatch} of \cite{hove99}, and $g\colon \mathcal
C\rightarrow \mathcal D$ is a cofibration precisely when the
induced maps
\begin{equation}\label{diagcof}
\mathcal C(I)\amalg_{L_I \mathcal C}L_I \mathcal D\longrightarrow
\mathcal D(I)
\end{equation}
are cofibrations in $\cat{mc}$ for all $I\in\sK$. Thus $\mathcal
D\colon\cat{cat}(\sK)\to\cat{mc}$ is cofibrant when every
canonical map $\colim \mathcal
D|_{\scat{cat}(\partial\varDelta(I))}\rightarrow \mathcal D(I)$ is
a cofibration.

In the dual model structure on $\fcat{cat$^{op}(\sK)$}{mc}$, weak
equivalences and cofibrations are given objectwise. To describe
the fibrations, we restrict $\mathcal C$ and $\mathcal D$ to the
undercategories $\cat{cat}^{op}(\partial\varDelta(I))$, and write
$M_I \mathcal C$ and $M_I \mathcal D$ for their respective limits.
So $M_I$ is the \emph{matching functor}\label{dematch} of
\cite{hove99}, and $f\colon \mathcal C\rightarrow \mathcal D$ is a
fibration precisely when the induced maps
\begin{equation}\label{diagf}
\mathcal C(I)\longrightarrow \mathcal D(I)\times_{M_I \mathcal
D}M_I \mathcal C
\end{equation}
are fibrations in $\cat{mc}$ for all $I\in\sK$. Thus $\mathcal
C\colon\cat{cat}^{op}(\sK)\to\cat{mc}$ is fibrant when every
canonical map $\mathcal C(I)\rightarrow\lim \mathcal
C|_{\scat{cat}^{op}(\partial\varDelta(I))}$ is a fibration.

\section{Algebraic model categories}\label{almoca}

Here we give further details of the algebraic model categories
introduced in the previous section. We describe the fibrations and
cofibrations in each category, comment on the status of the
strengthened axioms, and give simple examples in less familiar
cases. We also discuss two important adjoint pairs.

So far as general algebraic notation is concerned, we work over an
arbitrary commutative ring~$\k$. We indicate the coefficient ring
$\k$ by means of a subscript when necessary, but often omit it. In
some situations $\k$ is restricted to the rational numbers~$\Q$,
in this case we always clearly indicate it in the notation.

We consider finite sets $\mathcal W$ of generators $w_1$, \dots,
$w_m$. We write the graded tensor $\k$-algebra on $\mathcal W$ as
$T_\k(w_1,\dots,w_m)$, and use the abbreviation $T_\k(\mathcal W)$
whenever possible. Its symmetrisation $S_\k(\mathcal W)$ is the
graded commutative $\k$-algebra generated by $\mathcal W$. If
$\mathcal U$, $\mathcal V\subset\mathcal W$ are the subsets of odd
and even grading respectively, then $S_\k(\mathcal W)$ is the
tensor product of the exterior algebra $\Lambda_\k[\mathcal U]$
and the polynomial algebra $\k[\mathcal V]$. When $\k$ is a field
of characteristic zero, it is also convenient to denote the free
graded Lie algebra on $\mathcal W$ and its commutative counterpart
by $\fl_\k(\mathcal W)$ and $\cl_\k(\mathcal W)$ respectively; the
latter is nothing more than a free $\k$-module.

Almost all of our graded algebras have finite type, leading to a
natural coalgebraic structure on their duals. We write the free
tensor coalgebra on $\mathcal W$ as $T_\k\langle\mathcal
W\rangle$; it is isomorphic to $T_\k(\mathcal W)$ as $\k$-modules,
and its diagonal is given by
\[
  \Delta(w_{j_1}\otimes\dots\otimes w_{j_r})=
  \sum_{k=0}^r(w_{j_1}\otimes\dots\otimes w_{j_k})
  \otimes(w_{j_{k+1}}\otimes\dots\otimes w_{j_r}).
\]
The submodule $S_\k\langle\mathcal W\rangle$ of symmetric elements
($w_i\otimes w_j+(-1)^{\deg w_i\deg w_j}w_j\otimes w_i$, for
example) is the graded cocommutative $\k$-coalgebra cogenerated
by~$\mathcal W$.

Given $\mathcal W$, we may sometimes define a differential by
denoting the set of elements $dw_1,\ldots,dw_m$ by $d\mathcal W$.
For example, we write the free dg-algebra on a single generator
$w$ of positive dimension as $T_\k(w,dw)$; the notation is
designed to reinforce the fact that its underlying algebra is the
tensor $\k$-algebra on elements $w$ and~$dw$. Similarly,
$T_\k\langle w,dw\rangle$ is the free differential graded
coalgebra on~$w$. For further information on differential graded
coalgebras, \cite{h-m-s74} remains a valuable source.

\subsection*{Chain and cochain complexes}
The existence of a model structure on categories of chain
complexes was first proposed by Quillen~\cite{quil67}, whose view
of homological algebra as homotopy theory in $\cat{ch}_\k$ was a
crucial insight. Variations involving bounded and unbounded
complexes are studied by Hovey~\cite{hove99}, for example. In
$\cat{ch}_\k$, we assume that the fibrations are epimorphic in
positive degrees and the cofibrations are monomorphic with
degree-wise projective cokernel~\cite{dw-sp95}. In particular,
every object is fibrant.

The existence of limits and colimits is assured by working
dimensionwise, and functoriality of the
factorisations~\eqref{mcfacts} follows automatically from the fact
that $\cat{ch}_\k$ is \emph{cofibrantly generated} \cite[Chapter
2]{hove99}.

Model structures on $\cat{coch}_\k$ are established by analogous
techniques. It is usual to assume that the fibrations are
epimorphic with degree-wise injective kernel, and the cofibrations
are monomorphic in positive degrees. Then every object is
cofibrant. There is an alternative structure based on projectives,
but we shall only refer to the rational case so we ignore the
distinction.

Tensor product of (co)chain complexes invests $\cat{ch}_\k$ and
$\cat{coch}_\k$ with the structure of a monoidal model category,
as defined by Schwede and Shipley \cite{sc-sh00}.

\subsection*{Commutative differential graded algebras}
We consider commutative differential graded algebras over $\Q$
with cohomology differentials, so they are commutative monoids in
$\cat{coch}_\Q$. A model structure on $\cat{cdga}_\Q$ was first
defined in this context by Bousfield and Gugenheim~\cite{bo-gu76},
and has played a significant role in the theoretical development
of rational homotopy theory ever since. The fibrations are
epimorphic, and the cofibrations are determined by the appropriate
lifting property; some care is required to identify sufficiently
many explicit cofibrations.

Limits in $\cat{cdga}_\Q$ are created in the underlying category
$\cat{coch}_\Q$ and endowed with the natural algebra structure,
whereas colimits exist because $\cat{cdga}_\Q$ has finite
coproducts and filtered colimits. The proof of the factorisation
axioms in~\cite{bo-gu76} is already functorial.

By way of example, we note that the product of algebras $A$ and
$B$ is their augmented sum $A\oplus B$, defined by pulling back
the diagram of augmentations,
\[
\begin{CD}
A\oplus B @>>> A\\
@VVV @VV\varepsilon_A V\\
B @>\varepsilon_B>> \Q
\end{CD}
\]
in $\cat{coch}$ and imposing the standard multiplication on the
result. The coproduct is their tensor product $A\otimes B$ over
$\Q$. Examples of cofibrations include extensions $A\to(A\otimes
S(w),d)$ given by~\eqref{mmiter}; such an extension is determined
by a cocycle $z=f(x)$ in~$A$. This illustrates the fact that the
pushout of a cofibration is a cofibration. A larger class of
cofibrations $A\rightarrow A\otimes S(\mathcal W)$ is given by
iteration, for any set $\mathcal W$ of positive dimensional
generators corresponding to cocycles in~$A$.

The factorisations \eqref{mcfacts} are only valid over fields of
characteristic~$0$, so the model structure does not extend to
$\cat{cdga}_\k$ for arbitrary rings $\k$.

\subsection*{Differential graded algebras}
Our differential graded algebras have homology differentials, and
are the monoids in $\cat{ch}_\k$. A model category structure in
$\cat{dga}_\k$ is therefore induced by applying Quillen's path
object argument, as in \cite{sc-sh00}; a similar structure was
first proposed by Jardine~\cite{jard97} (albeit with cohomology
differentials), who proceeds by modifying the methods of
\cite{bo-gu76}. Fibrations are epimorphisms, and cofibrations are
determined by the appropriate lifting property.

Limits are created in $\cat{ch}_\k$, whereas colimits exist
because $\cat{dga}_\k$ has finite coproducts and filtered
colimits. Functoriality of the factorisations follows by adapting
the proofs of \cite{bo-gu76}, and works over arbitrary~$\k$.

For example, the coproduct of algebras $A$ and $B$ is the free
product \mbox{$A\star B$}, formed by factoring out an appropriate
differential graded ideal \cite{jard97} from the free (tensor)
algebra $T_\k(A\otimes B)$ on the chain complex $A\otimes B$.
Examples of cofibrations include the extensions $A\to(A\star
T_\k(w),d)$, determined by cycles $z$ in $A$. By analogy with the
commutative case, such an extension is defined by the pushout
diagram
%\begin{equation}\label{frexdga}
\[
\begin{CD}
  T_\k(x) @>f>> A\\
  @VVjV @VVV\\
  T_\k(w,dw) @>>> A\star T_\k(w)
\end{CD}
\]
%\end{equation}
where $f(x)=z$ and $j(x)=dw$. The differential on $A\star T_\k(w)$
is given by
\[
  d(a\star 1)=d_A a\star 1, \quad d(1\star w)=f(x)\star 1.
\]
Further cofibrations $A\rightarrow A\mathbin{\star}T_\k(\mathcal
W)$ arise by iteration, for any set $\mathcal W$ of positive
dimensional generators corresponding to cycles in $A$.

\subsection*{Cocommutative differential graded coalgebras}
The cocommutative comonoids in $\cat{ch}_\k$ are the objects of
$\cat{cdgc}_\k$, and the morphisms preserve comultiplication. The
model structure is defined only over fields of characteristic~$0$;
in view of our applications, we shall restrict attention to the
case~$\Q$. In practice, we interpret $\cat{cdgc}_\Q$ as the full
subcategory $\cat{cdgc}_{0,\Q}$ of connected objects $C$, which
are necessarily coaugmented. Model structure was first defined on
the category of simply connected rational cocommutative coalgebras
by Quillen~\cite{quil69R}, and refined to $\cat{cdgc}_{0,\Q}$ by
Neisendorfer~\cite{neis78}. The cofibrations are monomorphisms,
and the fibrations are determined by the appropriate lifting
property.

Limits exist because $\cat{cdgc}_\Q$ has finite products and
filtered limits, whereas colimits are created in $\cat{ch}_\Q$,
and endowed with the natural coalgebra structure. Functoriality of
the factorisations again follows by adapting the proofs
of~\cite{bo-gu76}.

For example, the product of coalgebras $C$ and $D$ is their tensor
product $C\otimes D$ over~$\Q$. The coproduct is their coaugmented
sum, given by pushing out the diagram of coaugmentations
\[
\begin{CD}
\Q @>\delta_C>> C\\
@V\delta_DVV @VVV\\
D @>>> C\oplus D
\end{CD}
\]
in $\cat{ch}_\Q$ and imposing the standard comultiplication on the
result. Examples of fibrations include the projections $(C\otimes
S_\Q\langle dt\rangle,d)\rightarrow C$, which are determined by
cycles $z$ in $C$ and defined by the pullback diagram
\[
\begin{CD}
  C\otimes
  S\langle dt\rangle @>>> S\langle t,dt\rangle\\
  @VVV @VVqV\\
  C@>h>> S\langle x\rangle
\end{CD}
\]
where $q(t)=x$, $q(dt)=0$ and $h(z)=x$. The differential on
$C\otimes S_\Q\langle dt\rangle$ satisfies
\[
  d(z\otimes 1)=1\otimes dt,\quad d(1\otimes dt)=0.
\]
This illustrates the fact that the pullback of a fibration is a
fibration. Further fibrations $C\otimes S_\Q\langle d\mathcal
T\rangle\to C$ are given by iteration, for any set $\mathcal T$ of
generators corresponding to elements of degree $\ge 2$ in~$C$.

\subsection*{Differential graded coalgebras}
Model structures on more general categories of differential graded
coalgebras have been publicised by Getzler and
Goerss~\cite{ge-go99}, who also work over a field. Once again, we
restrict attention~to $\Q$. The objects of $\cat{dgc}_\Q$ are
comonoids in $\cat{ch}_\Q$, and the morphisms preserve
comultiplication. The cofibrations are monomorphisms, and the
fibrations are determined by the appropriate lifting property.

Limits exist because $\cat{dgc}_\Q$ has finite products and
filtered limits, and colimits are created in $\cat{ch}_\Q$.
Functoriality of factorisations follows from the fact that the
model structure is cofibrantly generated.

For example, the product of coalgebras $C$ and $D$ is the cofree
product ${C\star D}$ \cite{ge-go99}. Their coproduct is the
coaugmented sum, as in the case of $\cat{cdgc}_\Q$. Examples of
fibrations include the projections $[C\star T_\Q\langle
dt\rangle,d]\rightarrow C$, which are determined by cycles $z$ in
$C$ and defined by the pullback diagram
\[
\begin{CD}
  C\star
  T_\Q\langle dt\rangle @>>> T_\Q\langle t,dt\rangle\\
  @VVV @VVqV\\
  C@>h>> T_\Q\langle x\rangle
\end{CD}
\]
where $q(t)=x$, $q(dt)=0$ and $h(z)=x$.

\subsection*{Differential graded Lie algebras}

A (rational) differential graded Lie algebra $L$ is a chain
complex in $\cat{ch}_\Q$, equipped with a bracket morphism
$[\:,\,]\colon L\otimes L\rightarrow L$ satisfying signed versions
of the antisymmetry and Jacobi identity:
\begin{equation}\label{glaide}
%\begin{gather*}
  [x,y]=-(-1)^{\deg x\deg y}[y,x],\quad
  [x,[y,z]]=[[x,y],z]+(-1)^{\deg x\deg y}[y,[x,z]].
%\end{gather*}
\end{equation}
Differential graded Lie algebras over $\Q$ are the objects of the
cagegory $\cat{dgl}$. Quillen \cite{quil69R} originally defined a
model structure on the subcategory of \emph{reduced} objects,
which was extended to $\cat{dgl}$ by Neisendorfer \cite{neis78}.
Fibrations are epimorphisms, and cofibrations are determined by
the appropriate lifting property.

Limits are created in $\cat{ch}_\Q$, whereas colimits exist
because $\cat{dgl}$ has finite coproducts and filtered colimits.
Functoriality of the factorisations follows by adapting the proofs
of~\cite{neis78}.

For example, the product of Lie algebras $L$ and $M$ is their
product $L\oplus M$ as chain complexes, with the induced bracket
structure. Their coproduct is the free product $L\star M$,
obtained by factoring out an appropriate differential graded ideal
from the free Lie algebra $\fl(L\otimes M)$ on the chain complex
$L\otimes M$. Examples of cofibrations include the extensions
$L\to(L\star\fl(w),d)$, which are determined by cycles $z$ in $L$
and defined by the pushout diagram
\[
\begin{CD}
  \fl(x) @>f>> L\\
  @VVjV @VVV\\
  \fl(w,dw) @>>> L\star\fl(w)
\end{CD}
\]
%\end{equation}
where $f(x)=z$ and $j(x)=dw$. The differential on $L\star\fl(w)$
is given by
\[
  d(l\star 1)=d_Ll\star 1,\quad d(1\star w)=z\star 1.
\]
For historical reasons, a differential graded Lie algebra $L$ is
said to be \emph{coformal}\label{coformLie} whenever it is formal
in $\cat{dgl}$.

\subsection*{Adjoint pairs}
Following Moore \cite{moor71}, \cite{h-m-s74}, we consider the
algebraic \emph{classifying functor} $B_*$ and the \emph{loop
functor} $\varOmega_*$ as an adjoint pair
\begin{equation}\label{barcobar}
  \varOmega_*\colon\cat{dgc}_{0,\k}\rightleftarrows\cat{dga}_\k:\! B_*.
\end{equation}
For any object $A$ of $\cat{dga}_\k$, the classifying coalgebra
$B_*A$ agrees with Eilenberg and Mac Lane's normalised \emph{bar
construction} as objects of $\cat{ch}_\k$. For any object $C$ of
$\cat{dgc}_{0,\k}$, the loop algebra $\varOmega_*C$ is given by
the tensor algebra $T_\k(s^{-1}\overline{C})$ on the desuspended
$\k$-module $\overline{C}=\Ker(\varepsilon\colon C\to \k)$, and
agrees with Adams' \emph{cobar construction}~\cite{adam56} as
objects of $\cat{ch}_\k$.

The classical result of Adams links the Moore loop functor
$\varOmega\colon\cat{top}\to\cat{tmon}$ with its algebraic
analogue~$\varOmega_*$:

\begin{theorem}[\cite{adam56}]\label{adamscobar}
For a simply connected pointed space $X$ and a commutative
ring~$\k$, there is a natural isomorphism of graded algebras
\[
  H(\varOmega_*C_*(X;\k))\cong H_*(\varOmega X;\k),
\]
where $C_*(X;\k)$ denotes the suitably reduced singular chain
complex of~$X$.
\end{theorem}

The isomorphism of Theorem~\ref{adamscobar} is induced by a
natural homomorphism
\[
  \varOmega_*C_*(X;\k)\longrightarrow \mathit{CU}_*(\varOmega
  X;\k)
\]
of $\cat{dga}_\k$, where $\mathit{CU}_*(\varOmega X;\k)$ denotes
the suitably reduced cubical chains on $\varOmega X$ with the
dg-algebra structure induced from composition of Moore loops.

The graded homology algebra $H(\varOmega_*C)$ is denoted by
$\mathop{\mathrm{Cotor}}_C(\k,\k)$. When $\k$ is a field, there is
an isomorphism
\begin{equation}\label{cotorext}
\mathop{\mathrm{Cotor}}\nolimits_C(\k,\k)
\;\cong\;\Ext_{C^*}(\k,\k)
\end{equation}
of graded algebras \cite[page~41]{prid70}, where $C^*$ is the
graded algebra dual to $C$ and $\Ext_{C^*}(\k,\k)$ is the
\emph{Yoneda algebra}\label{yonedaalg} of $C^*$~\cite{macl63}.

\begin{proposition}\label{bcqpa}
The loop functor $\varOmega_*$ preserves cofibrations of connected
coalgebras and weak equivalences of simply connected coalgebras;
the classifying functor $B_*$ preserves fibrations of connected
algebras and all weak equivalences.
\end{proposition}
\begin{proof}
The fact that $B_*$ and $\varOmega_*$ preserve weak equivalences
of algebras and simply connected coalgebras respectively is proved
by standard arguments with the Eilenberg--Moore spectral sequence
\cite[page~538]{f-h-t92}. The additional assumption for coalgebras
is necessary to ensure that the cobar spectral sequence converges,
because the relevant filtration is decreasing.

Given any cofibration $i\colon C_1\to C_2$ of connected
coalgebras, we must check that $\varOmega_*i\colon\varOmega_*
C_1\to\varOmega_* C_2$ satisfies the left lifting property with
respect to any acyclic fibration $p\colon A_1\to A_2$ in
$\cat{dga}_\k$. This involves finding lifts $\varOmega_*C_2\to
A_1$ and $C_2\to B_*A_1$ in the respective diagrams
\[
\begin{CD}
  \varOmega_*C_1 @>>> A_1\\
  @V{\varOmega_*i}VV @VVpV\\
  \varOmega_*C_2 @>>> A_2
\end{CD}\qquad\text{and}\qquad
\begin{CD}
  C_1 @>>> B_*A_1\\
  @ViVV @VV{B_*p}V\\
  C_2 @>>> B_*A_2
\end{CD}\quad\quad;
\]
each lift implies the other, by adjointness. Since $p$ is an
acyclic fibration, its kernel $A$ satisfies $H(A)\cong \k$. In
this circumstances, the projection $B_*p$ splits
by~\cite[Theorem~IV.2.5]{h-m-s74}, so $B_*A_1$ is isomorphic to
the cofree product $B_*A_2\star B_*A$. Therefore, $B_*p$ is an
acyclic fibration in $\cat{dgc}_{0,\k}$, and our lift is assured.

A second application of adjointness shows that $B_*$ preserves all
fibrations of connected algebras.
\end{proof}

\begin{remark}
It follows from Proposition \ref{bcqpa} that the restriction
of~\eqref{barcobar} to simply connected coalgebras and connected
algebras respectively,
\[
  \varOmega_*\colon\cat{dgc}_{1,\k}\rightleftarrows\cat{dga}_{0,\k}:\!
  B_*,
\]
acts as a Quillen pair, and induces an adjoint pair of
equivalences on appropriate full subcategories of the homotopy
categories. An example is given in~\cite[p.~538]{f-h-t92} which
shows that $\varOmega_*$ fails to preserve quasi-isomorphisms (or
even acyclic cofibrations) if the coalgebras are not simply
connected.
\end{remark}

Over $\Q$, the adjunction maps $C\mapsto B_*\varOmega_*C$ and
$\varOmega_*B_*A\mapsto A$ are quasi-\-iso\-mor\-phisms for any
objects $A$ and~$C$.

Following Neisendorfer \cite[Proposition 7.2]{neis78}, we consider
a second pair of adjoint functors
\begin{equation}\label{lcadj}
L_*\colon\cat{cdgc}_{0,\Q}\rightleftarrows\cat{dgl}:\! M_*,
\end{equation}
whose derived functors induce an equivalence between
$\Ho(\cat{cdgc}_{0,\Q})$ and a certain full subcategory of
$\Ho(\cat{dgl})$. This extends Quillen's original results
\cite{quil69R} for $L_*$ and $M_*$, which apply only to simply
connected coalgebras and connected Lie algebras. Given a connected
cocommutative coalgebra $C$, the underlying graded Lie algebra of
$L_*C$ is the free Lie algebra $\fl(s^{-1}\overline{C})\subset
T(s^{-1}\overline{C})$. This is preserved by the differential in
$\varOmega_*C$ because $C$ is cocommutative, thereby identifying
$L_*C$ as the differential graded Lie algebra of primitives in
$\varOmega_*C$. The right adjoint functor $M_*$ may be regarded as
a generalisation to differential graded objects of the standard
complex for calculating the cohomology of Lie algebras. Given any
$L$ in $\cat{dgl}$ the underlying cocommutative coalgebra of
$M_*L$ is the symmetric coalgebra $C(sL)$ on the suspended vector
space $L$.

The ordinary (topological) classifying space functor
$B\colon\cat{top}\to\cat{tmon}$ and the Moore loop functor
$\varOmega\colon\cat{tmon}\to\cat{top}$ are not formally adjoint,
because $\varOmega$ does not preserve products. However, as it was
shown by Vogt~\cite{vogt}, after passing to appropriate
localisations, $\varOmega$ becomes \emph{right} adjoint to $B$ in
the homotopy categories.

There is also a similar result in simplicial category: the loop
functor from simplicial sets to simplicial groups is \emph{left}
adjoint to the classifying functor, as in the case of algebraic
functors~$\varOmega_*$ and~$B_*$.

\section{Homotopy limits and colimits}\label{holico}

The $\lim$ and $\colim$ functors $\fcat{a}{mc}\to\cat{mc}$ do not
generally preserve weak equivalences, and the theory of homotopy
limits and colimits has been developed to remedy this deficiency.
We outline their construction in this section, and discuss basic
properties.
%The literature is still in a state of considerable
%flux, and we refer to Recke's thesis~\cite{reck} for a comparison
%of several alternative treatments.
%Here we focus mainly on those of Recke's statements that are
%inspired by Hirschhorn, and make detailed appeal to \cite{hirs03}
%as necessary.

With $\cat{cat}(\sK)$ and $\cat{cat}^{op}(\sK)$ in mind as primary
examples, we assume throughout that $\cat{a}$ is a finite Reedy
category.

A Reedy category $\cat{a}$ has \emph{cofibrant constants} if the
constant $\cat{a}$-diagram $M$ is cofibrant in $\fcat{a}{mc}$, for
any cofibrant object $M$ of an arbitrary model category
$\cat{mc}$. Similarly, $\cat{a}$ has \emph{fibrant constants} if
the constant $\cat{a}$-diagram $N$ is fibrant for any fibrant
object $N$ of $\cat{mc}$.

%Note that the initial and terminal objects of $\fcat{a}{mc}$ are
%the constant diagrams $\circ$ and $*$ respectively.

As shown in \cite[Theorem 15.10.8]{hirs03}, a Reedy category
$\cat{a}$ has fibrant constants if and only if the first pair of
adjoint functors
\begin{equation}\label{quilpairs}
\colim\colon\fcat{a}{mc}\rightleftarrows\cat{mc}:\!\kappa,\qquad
\kappa\colon\cat{mc}\rightleftarrows\fcat{a}{mc}:\!\lim
\end{equation}
is a Quillen pair (i.e. $\colim$ is left Quillen) for every model
category $\cat{mc}$. Similarly, $\cat{a}$ has cofibrant constants
if and only if the second pair above is a Quillen pair, i.e.
$\lim$ is right Quillen.
%In these circumstances, it follows from
%Ken Brown's Lemma that $\lim$ and $\colim$ preserve weak
%equivalences on Reedy fibrant and Reedy cofibrant diagrams
%respectively.

We now apply the fibrant and cofibrant replacement functors
associated to the Reedy model structure on $\fcat{a}{mc}$, and
their homotopy functors \eqref{hoomfi}.
\begin{definition}\label{defhocolhol}
For any Reedy category $\cat{a}$ with fibrant and cofibrant
constants, and any model category $\cat{mc}$:
\begin{enumerate}
\item[(a)] the \emph{homotopy colimit} functor is the composition
\[
\hocolim\colon\Ho\fcat{a}{mc}
\stackrel{\Ho(\omega)}{\llllongrightarrow} \Ho\fcat{a}{mc}_c
\stackrel{\Ho(\colim)}{\lllllongrightarrow}\Ho(\cat{mc})\;;
\]
\item[(b)] the \emph{homotopy limit} functor is the composition
\[
\holim\colon\Ho\fcat{a}{mc}
\stackrel{\Ho(\phi)}{\llllongrightarrow} \Ho\fcat{a}{mc}_{\!f}
\stackrel{\Ho(\lim)}{\lllllongrightarrow}\Ho(\cat{mc}).
\]
\end{enumerate}
\end{definition}

\begin{remark}\label{htpyinvce}
Definition \ref{defhocolhol} incorporates the fact that $\holim$
and $\hocolim$ map objectwise weak equivalences of diagrams to
weak equivalences in $\cat{mc}$.
\end{remark}

The Reedy categories $\cat{cat}(\sK)$ and $\cat{cat}^{op}(\sK)$
satisfy the criteria of \cite[Proposition 15.10.2]{hirs03} and
therefore have fibrant and cofibrant constants, for every
simplicial complex $\sK$. This implies that $\holim$ and
$\hocolim\colon\Ho\fcat{cat$(\sK)$}{mc}\rightarrow\Ho(\cat{mc})$
are defined; and similarly for $\cat{cat}^{op}(\sK)$.

Describing explicit models for homotopy limits and colimits has
been a major objective for homotopy theorists since their study
was initiated by Bousfield and Kan~\cite{bo-ka72} and
Vogt~\cite{vogt73}. In terms of Definition \ref{defhocolhol}, the
issue is to choose fibrant and cofibrant replacement functors
$\phi$ and $\omega$. Many alternatives exist, including those
defined by the two-sided bar and cobar constructions of
\cite{p-r-v04} or the frames of \cite[\S16.6]{hirs03}, but no
single description yet appears to be convenient in all cases.
Instead, we accept a variety of possibilities, which are often
implicit; the next few results ensure that they are as compatible
and well-behaved as we need.

\begin{proposition}\label{hclhlpreswe}
Any cofibrant approximation $\mathcal
D'\stackrel{\simeq}{\longrightarrow} \mathcal D$ of diagrams
induces a weak equivalence $\colim \mathcal
D'\stackrel{\simeq}{\longrightarrow}\hocolim \mathcal D$ in
$\cat{mc}$; and any fibrant approximation $\mathcal
D\stackrel{\simeq}{\longrightarrow} \mathcal D''$ induces a weak
equivalence $\holim \mathcal
D\stackrel{\simeq}{\longrightarrow}\lim \mathcal D''$.
\end{proposition}
\begin{proof}
Using the left lifting property (axiom~(c)) of the cofibration
$\circ\to \mathcal D'$ with respect to the acyclic fibration
$\omega(\mathcal D)\to \mathcal D$ we obtain a factorisation
$\mathcal D'\to\omega(\mathcal D)\to \mathcal D$, in which the
left hand map is a weak equivalence by axiom~(b). But $\mathcal
D'$ and $\omega(\mathcal D)$ are cofibrant, and $\colim$ is left
Quillen, so the induced map $\colim \mathcal
D'\to\colim\omega(\mathcal D)$ is a weak equivalence, as required.
The proof for $\lim$ is dual.
\end{proof}

\begin{remark}\label{unirepl}
Such arguments may be strengthened to include uniqueness
statements, and show that the replacements $\phi(\mathcal D)$ and
$\omega(\mathcal D)$ are themselves unique up to homotopy
equivalence over $\mathcal D$, see \cite[Proposition
8.1.8]{hirs03}.
\end{remark}

\begin{proposition}\label{hotolim}
For any cofibrant diagram $\mathcal D$ and fibrant diagram
$\mathcal E$, there are natural weak equivalences $\hocolim
\mathcal D\stackrel{\simeq}{\longrightarrow}\colim \mathcal D$ and
$\lim \mathcal E\stackrel{\simeq}{\longrightarrow}\holim \mathcal
E$.
\end{proposition}
\begin{proof}
For $\mathcal D$, it suffices to apply the left Quillen functor
$\colim$ to the acyclic fibration $\omega(\mathcal D)\to \mathcal
D$. The proof for $\mathcal E$ is dual.
\end{proof}

\begin{proposition}\label{wecolim}
A weak equivalence $\mathcal
D'\stackrel{\simeq}{\longrightarrow}\mathcal D$ of cofibrant
diagrams induces a weak equivalence $\colim\mathcal
D'\stackrel{\simeq}{\longrightarrow}\colim\mathcal D$, and a weak
equivalence $\mathcal E\stackrel{\simeq}{\longrightarrow}\mathcal
E'$ of fibrant diagrams induces a weak equivalence $\lim\mathcal
E\stackrel{\simeq}{\longrightarrow}\lim\mathcal E'$.
\end{proposition}
\begin{proof}
This follows from Propositions~\ref{hclhlpreswe}
and~\ref{hotolim}.
\end{proof}

\begin{proposition}\label{hopu}
In any model category $\cat{mc}$:
\begin{enumerate}

\item[(a)]
if all three objects of a pushout diagram $\mathcal D\colon
L\leftarrow M\to N$ are cofibrant, and either of the maps is a
cofibration, then there exists a weak equivalence $\hocolim
\mathcal D\stackrel{\simeq}{\longrightarrow}\colim \mathcal D$;
\item[(b)] if all three objects of a pullback diagram $\mathcal E\colon P\to
Q\leftarrow R$ are fibrant, and either of the maps is a fibration,
then there exists a weak equivalence $\lim \mathcal
E\stackrel{\simeq}{\longrightarrow}\holim \mathcal E$.
\end{enumerate}
\end{proposition}

\begin{proof} For (a), assume that $M\to N$ is a cofibration, and that
the indexing category $\cat{b}$ for $\mathcal D$ has non-identity
morphisms $\lambda\leftarrow\mu\to\nu$. The degree function
$\deg(\lambda)=0$, $\deg(\mu)=1$, and $\deg(\nu)=2$ turns
$\cat{b}$ into a Reedy category with fibrant constants, and
ensures that $\mathcal D$ is cofibrant. So Proposition
\ref{hotolim} applies. If $M\to L$ is a cofibration, the
corresponding argument holds by symmetry.

For (b), the proofs are dual.
\end{proof}

%Further details  may be found in \cite[Proposition~19.9.4]{hirs03}.

There is an important situation when the homotopy colimit over a
poset category can be described explicitly:

\begin{lemma}[{Wedge
Lemma~\cite[Lemma~4.9]{w-z-z99}}]\label{wedgelemma}
Let $(\mathcal P,\le)$ be a poset with initial element~$\hatzero$,
and let $\cat{p}$ be the corresponding poset category. Suppose
there is diagram $\mathcal D\colon\cat{p}\to\cat{top}$ of spaces
so that $\mathcal D(\hatzero)=\pt$ and $\mathcal
D(\sigma)\to\mathcal D(\tau)$ is the constant map to the basepoint
for all $\sigma<\tau$ in~$\mathcal P$.
%and there exist points $p_\sigma\in\mathcal D(\sigma)$ for all
%$\sigma\in\mathcal P$ such that $\mathcal D(\sigma)\to\mathcal
%D(\tau)$ is the constant map to $p_\tau$ for all $\sigma<\tau$.
Then there is a homotopy equivalence
\[
  \hocolim\mathcal D\stackrel{\simeq}{\longrightarrow}
  \bigvee_{\sigma\in\mathcal P}
  \bigl(|\ord(\mathcal P_{>\sigma})|\mathop{*}\mathcal
  D(\sigma)\bigr),
\]
where $|\ord(\mathcal P_{>\sigma})|$ is the geometric realisation
of the order complex of the upper semi-interval
$P_{>\sigma}=\{\tau\in\mathcal P\colon\tau>\sigma\}$.
\end{lemma}

According to a result of Panov, Ray and Vogt, the classifying
space functor $B\colon\cat{tmon}\to\cat{top}$ commutes with
homotopy colimits (of topological monoids and topological spaces,
respectively) in the following sense:

\begin{theorem}[{\cite[Theorem 7.12, Proposition 7.15]{p-r-v04}}]\label{PRVthm}
For any diagram $\mathcal D\colon\cat{a}\rightarrow\cat{tmon}$ of
well-pointed topological monoids with the homotopy types of cell
complexes, there is a natural homotopy equivalence
\[
  g\colon\hocolim^{\scat{top}}B\mathcal D\stackrel\simeq\longrightarrow
  B\hocolim^{\scat{tmon}}\mathcal D.
\]
Furthermore, there is a commutative square
\begin{equation}\label{bhocolimdiag}
\begin{CD}
  \hocolim^{\scat{top}}B\mathcal D@>\simeq>>B\hocolim^{\scat{tmon}}\mathcal D\\
  @VVp^{\scat{top}}V@VV Bp^{\scat{tmon}}V\\
  \colim^{\scat{top}}B\mathcal D@>>>B\colim^{\scat{tmon}}\mathcal D
\end{CD}
\end{equation}\\
where $p^{\scat{top}}$ and $p^{\scat{tmon}}$ are the natural
projections.
\end{theorem}

A weaker version of the theorem above can be stated for the Moore
loop functor:

\begin{corollary}\label{loopshocolim}
For $\mathcal D\colon\cat{a}\rightarrow\cat{tmon}$ as above, there
is a commutative square
\[
\begin{CD}
  \varOmega\hocolim^{\scat{top}}B\mathcal D@>\simeq>>\hocolim^{\scat{tmon}}\mathcal D\\
  @VV\varOmega p^{\scat{top}}V@VV p^{\scat{tmon}}V\\
  \varOmega\colim^{\scat{top}}B\mathcal D@>>>\colim^{\scat{tmon}}\mathcal D
\end{CD}
\]\\
in $Ho(\cat{tmon})$, where the upper homomorphism is a homotopy
equivalence.
\end{corollary}
\begin{proof}
This follows by applying $\varOmega$ to~\eqref{bhocolimdiag} and
then composing the horizontal maps with the canonical weak
equivalence $\varOmega BG\rightarrow G$ in $\cat{tmon}$, where
$G=\hocolim^{\scat{tmon}}\mathcal D$ and
$\colim^{\scat{tmon}}\mathcal D$ respectively.
\end{proof}

The lower map in the diagram above is not a homotopy equivalence
in general, although $\varOmega$ preserves coproducts. Appropriate
examples are given in Section~\ref{loops}.

\chapter{Bordism and cobordism}\label{cobor}
Here we summarise the required facts from the theory of bordism
and cobordism, with the most attention given to complex
(co)bordism. Cobordism theory is one of the deepest and most
influential parts of algebraic topology, which experienced a
spectacular development in the 1960s. Proofs of results presented
here would require a separate monograph and a substantial
background in algebraic topology. There are exceptions where the
proofs are concise and included here. For the rest an interested
reader is referred to the works of Novikov \cite{novi67},
\cite{novi96}, and monographs of
Conner--Floyd~\cite{co-fl64},~\cite{co-fl66} and
Stong~\cite{ston68}.

We consider topological spaces which have the homotopy type of
cell complexes. All manifolds are assumed to be smooth, compact
and closed (without boundary), unless otherwise specified.

\section{Bordism of manifolds}\label{bordismman}
Given two $n$-dimensional manifolds $M_0$ and $M_1$, a
\emph{bordism} between them is an $(n+1)$-dimensional manifold $W$
with boundary, whose boundary is the disjoint union of $M_0$ and
$M_1$, that is, $\partial W=M_0\sqcup M_1$. If such a $W$ exists,
$M_0$ and $M_1$ are called \emph{bordant}. The bordism relation
splits
%MA
the set of manifolds into equivalence classes (see
Fig.~\ref{trcob}.1), which are called \emph{bordism classes}.
\begin{figure}
\label{trcob} \vspace*{-5mm}
\begin{center}
\end{center}
\caption{Transitivity of the bordism relation.}
\end{figure}

We denote the bordism class of $M$ by $[M]$, and denote by
$\varOmega_n^O$ the set of bordism classes of $n$-dimensional
manifolds. Then $\varOmega_n^O$ is an Abelian group with respect
to the disjoint union operation: $[M_1]+[M_2]=[M_1\sqcup M_2]$.
Zero is represented by the bordism class of the empty set (which
is counted as a manifold in any dimension), or by the bordism
class of any manifold which bounds. We also have $\partial(M\times
I)=M\sqcup M$. Hence, $2[M]=0$ and $\varOmega_n^O$ is a 2-torsion
group.

Set $\varOmega^O=\bigoplus _{n \ge 0}\varOmega _n^O$. The direct
product of manifolds induces a multiplication of bordism classes,
namely $[M_1]\times [M_2]=[M_1 \times M_2]$. It makes
$\varOmega^O$ a graded commutative ring, the \emph{unoriented
bordism ring}\label{unoriebordi}.

For any space $X$ the bordism relation can be extended to maps of
manifolds to~$X$: two maps $M_1\to X$ and $M_2\to X$ are
\emph{bordant} if there is a bordism $W$ between $M_1$ and $M_2$
and the map $M_1\sqcup M_2\to X$ extends to a map $W\to X$. The
set of bordism classes of maps $M\to X$ with $\dim M=n$ forms an
abelian group called the \emph{$n$-dimensional unoriented bordism
group of~$X$} and denoted $O_n(X)$ (other notations: $N_n(X)$,
$MO_n(X)$). We note that $O_n(pt)=\varOmega_n^O$, where $pt$ is a
point.

The bordism group $O_n(X,A)$ of a pair $A\subset X$ is defined as
the set of bordism classes of maps of manifolds with boundary,
$(M,\partial M)\to (X,A)$, where $\dim M=n$.
%MA
(Two such maps $f_0\colon(M_0,\partial M_0)\to(X,A)$ and
$f_1\colon(M_1,\partial M_1)\to(X,A)$ are \emph{bordant} if there is
$W$ such that $\partial W=M_0\cup M_1\cup M$, where $M$ is a bordism
between $\partial M_0$ and $\partial M_1$, and a map $f\colon W\to
X$ such that $f|_{M_0}=f_0$, $f|_{M_1}=f_1$ and $f(M)\subset A$.) We
have $O_n(X,\varnothing)=O_n(X)$.

There is an obviously defined map $\varOmega_m^O\times O_n(X)\to
O_{m+n}(X)$ turning $O_*(X)=\bigoplus_{n\ge0}O_n(X)$ into a graded
$\varOmega^O$-module. The assignment $X\mapsto O_*(X)$ defines a
\emph{generalised homology theory}\label{genhomotheory}, that is,
it is functorial in $X$, homotopy
%MA
invariant, has the excision property and exact sequences of pairs.

\section{Thom spaces and cobordism functors}\label{thomspaces}
A remarkable geometric construction due to Pontryagin and Thom
reduces the calculation of the bordism groups $O_n(X)$ to a
homotopical problem. Here we assume known basic facts from the
theory of vector bundles.

Given an $n$-dimensional real Euclidean vector bundle $\xi$ with
total space $E=E\xi$ and Hausdorff compact base~$X$, the
\emph{Thom space}\label{thomspaced} of~$\xi$ is defined as the
quotient
\[
  \Th \xi=E/E_{\ge1},
\]
where $E_{\ge1}$ is the subspace consisting of vectors of length
$\ge1$ in the fibres of~$\xi$. Equivalently, $\Th\xi=BE/SE$, where
$BE$ is the total space of the $n$-ball bundle associated
with~$\xi$ and $SE=\partial BE$ is the $(n-1)$-sphere bundle.
Also, $\Th\xi$ is the one-point compactification of~$E$.
%Given two spaces $X$ and $Y$ with basepoints $\pt_X$ and $\pt_Y$
%respectively, their \emph{smash product} is defined as
%\[
%  X\wedge Y=X\times Y\big/\bigl((X\times \pt_Y)\cup(\pt_X\times Y)\bigr).
%\]
The Thom space $\Th\xi$ has a canonical basepoint, the image
of~$E_{\ge1}$.

\begin{proposition}\label{thomprod}
If $\xi$ and $\eta$ are vector bundles over $X$ and $Y$
respectively, and $\xi\times\eta$ is the product vector bundle over
$X\times Y$, then
\[
  \Th(\xi\times\eta)=\Th\xi\wedge \Th\eta.
\]
\end{proposition}
The proof is left as an exercise.

\begin{example}\

1. Regarding $\R^k$ as the total space of a $k$-plane bundle over
a point, we obtain that the corresponding Thom space $\Th(\R^k)$
is a $k$-sphere~$S^k$.

2. If $\xi$ is a 0-dimensional bundle over $X$, then
$\Th\xi=X_+=X\sqcup pt$.

3. Let $\underline{\R}^k$ denote the trivial $k$-plane bundle
over~$X$. The Whitney sum $\xi\oplus\underline{\R}^k$ can be
identified with the product bundle $\xi\times\R^k$, where $\R^k$ is
the $k$-plane bundle over a point. Then Proposition~\ref{thomprod}
implies that
\[
  \Th(\xi\oplus\underline{\R}^k)=\varSigma^k\Th\xi,
\]
where $\varSigma^k$ denote the $k$-fold suspension.

4. Combining the previous two examples, we obtain
\[
  \Th(\underline{\R}^k)=\varSigma^k X\vee S^k.
\]
\end{example}

\begin{construction}[Pontryagin--Thom
construction]\label{ponthomcol} Let $M$ be a submanifold in $\R^m$
with normal bundle $\nu=\nu(M\subset\R^m)$.  The
\emph{Pontryagin--Thom map}
\[
  S^m\to\Th\nu
\]
identifies the tubular neighbourhood of $M$ in $\R^m\subset S^m$
with the set of vectors of length $<1$ in the fibres of~$\nu$, and
collapses the complement of the tubular neighbourhood to the
basepoint of the Thom space~$\Th\nu$.

This construction can be generalised to submanifolds $M\subset
E\xi$ in the total space of an arbitrary $m$-plane bundle~$\xi$
over a manifold, giving the collapse map
\begin{equation}\label{xinucollapse}
  \Th\xi\to\Th\nu,
\end{equation}
where $\nu=\nu(M\subset E\xi)$. Note that the Pontryagin--Thom
collapse map is a particular case of~\eqref{xinucollapse}, as
$S^m$ is the Thom space of an $m$-plane bundle over a point.
\end{construction}

Recall that a smooth map $f\colon W\to Z$ of manifolds is called
\emph{transverse} along a submanifold $Y\subset Z$ if, for every
$w\in f^{-1}(Y)$, the image of the tangent space to $W$ at $w$
together with the tangent space to $Y$ at $f(w)$ spans the tangent
space to $Z$ at~$f(w)$:
\[
  f_*\mathcal T_wW+\mathcal T_{f(w)}Y=\mathcal T_{f(w)}Z.
\]
If $f\colon W\to Z$ is transverse along $Y\subset Z$, then
$f^{-1}(Y)$ is a submanifold in $W$ of codimension equal to the
codimension of $Y$ in~$Z$.

\begin{construction}[cobordism classes of $\eta$-submanifolds in
$E\xi$]\label{etaxicob}
Let $\xi$ be an $m$-plane bundle with total
space $E\xi$ over a manifold~$X$, and let $\eta$ be an $n$-plane
bundle over a manifold~$Y$. An \emph{$\eta$-submanifold} of $E\xi$
is a pair $(M,f)$ consisting of a submanifold $M$ in $E\xi$ and a
map
\[
  f\colon\nu(M\subset E\xi)\to\eta
\]
such that $f$ is an isomorphism on each fibre (so that the
codimension of $M$ in $E\xi$ is~$n$). Two $\eta$-submanifolds
$(M_0,f_0)$ and $(M_1,f_1)$ are \emph{cobordant} if there is an
$\eta$-submanifold with boundary $(W,f)$ in the cylinder
$E\xi\times I\subset E(\xi\oplus\underline\R)$ such that
\[
  \partial(W,f)=\bigl((M_0,f_0)\times0\bigr)\cup\bigl((M_1,f_1)\times1\bigr).
\]
\end{construction}

\begin{theorem}\label{etacob}
The set of cobordism classes of $\eta$-submanifolds in $E\xi$ is in
one-to-one correspondence with the set $[\Th\xi,\Th\eta]$ of
homotopy classes of based maps of Thom spaces.
\end{theorem}
\begin{proof}
Assume given a based map $g\colon\Th\xi\to\Th\eta$. By changing $g$
within its homotopy class we may achieve that $g$ is transverse
along the zero section $Y\subset\Th\eta$ (transversality is a local
condition, and both $\Th\xi$ and $\Th\eta$ are manifolds outside the
basepoints). Since $\eta$ is an $n$-plane bundle, $M=g^{-1}(Y)$ is a
submanifold of codimension~$n$ in $E\xi=\Th\xi\setminus pt$ such
that
\[
  \nu(M\subset E\xi)=g^*(\nu(Y\subset E\eta))=g^*\eta.
\]
That is, $M$ is an $\eta$-submanifold in $E\xi$.

Conversely, assume given an $\eta$-submanifold $M\subset E\xi$. We
therefore have the map of Thom spaces $\Th\nu\to\Th\eta$ (induced by
the map of $\nu$ to $\eta$), whose composition with the collapse
map~\eqref{xinucollapse} gives the required map $\Th\xi\to\Th\eta$.

The fact that homotopic based maps $\Th\xi\to\Th\eta$ correspond to
cobordant $\eta$-submanifolds, and vice versa, is left as an
exercise.
\end{proof}

\begin{construction}[cobordism groups]
Let $\eta_k$ be the universal vector $k$-plane bundle $\EO(k)\to
\BO(k)$. Following the original notation of Thom, we denote
$\MO(k)=\Th\eta_k$.

Every submanifold $M\subset\R^{n+k}$ of dimension $n$ is an
$\eta_k$-submanifold via the classifying map of the normal bundle
$\nu(M\subset\R^{n+k})$. Denote by $\varOmega_O^{-n,k}$ the set of
cobordism classes of $M\subset\R^{n+k}$. The base of~$\eta_k$ is not
a manifold, but it is a direct limit of Grassmannians, and a simple
limit argument shows that Theorem~\ref{etacob} still holds for
$\eta_k$-submanifolds. Hence,
\[
  \varOmega_O^{-n,k}=[S^{n+k},\Th\eta_k]=\pi_{n+k}(\MO(k)).
\]
There is the stabilisation map
$\varOmega_O^{-n,k}\to\varOmega_O^{-n,k+1}$ obtained by composing
the suspended map $S^{n+k+1}\to\varSigma\MO(k)$ with the map
$\varSigma\MO(k)\to\MO(k+1)$ induced by the bundle map
$\eta_{k}\oplus\underline{\R}\to\eta_{k+1}$. The \emph{$(-n)$th
cobordism group} is defined by
\begin{equation}\label{-ncobgr}
  \varOmega_O^{-n}=\lim_{k\to\infty}\varOmega_O^{-n,k}=\lim_{k\to\infty}\pi_{n+k}(\MO(k)).
\end{equation}
\end{construction}

\begin{proposition}\label{borcobor}
We have a canonical isomorphism
\[
  \varOmega_O^{-n}\cong\varOmega^O_n
\]
between the cobordism and bordism groups, for $n\ge0$. In other
words, two $n$-dimensional manifolds $M_0$ and $M_1$ are bordant if
and only if there exist embeddings of $M_0$ and $M_1$ in the same
$\R^{n+k}$ which are cobordant.
\end{proposition}
\begin{proof}
Forgetting the embedding $M\subset\R^{n+k}$ we obtain a map
$\varOmega_O^{-n,k}\to\varOmega^O_n$, which may be shown to be a
group homomorphism. It is compatible with the stabilisation maps,
and therefore defines a homomorphism
$\varOmega_O^{-n}\to\varOmega^O_n$. Since every manifold $M$ may
be embedded in some $\R^{n+k}$, it is an isomorphism.
\end{proof}

Together with~\eqref{-ncobgr}, Proposition~\ref{borcobor} gives a
homotopical interpretation for the (unoriented) bordism groups. This
also implies that the notions of the `bordism class' and `cobordism
class' of a manifold $M$ are interchangeable. Theorem~\ref{etacob}
may be applied further to obtain a homotopical interpretation for
the bordism groups $O_n(X)$ of a space $X$:

\begin{construction}[bordism and cobordism groups of a
space]\label{oborspace} Let $X$ be a space. We set $\xi=\R^{n+k}$
(an $(n+k)$-plane bundle over a point) and $\eta=X\times\eta_k$
(the product of a 0-plane bundle over $X$ and the universal
$k$-plane bundle~$\eta_k$ over $BO(k)$), and consider cobordism
classes of $\eta$-submanifolds in~$E\xi$. Such a submanifold is
described by a pair $(f,\iota)$ consisting of a map $f\colon M\to
X$ from an $n$-dimensional manifold to~$X$ and an embedding
$\iota\colon M\to \R^{n+k}$ (the bundle map from $\nu(\iota)$ to
$\eta$ is the product of the map $E\nu(\iota)\to M\to X$ and the
classifying map of $\nu(\iota)$). By Theorem~\ref{etacob}, the set
of cobordism classes of $\eta$-submanifolds in $E\xi$ is given by
\[
  [S^{n+k},\Th(X\times\eta_k)]=\pi_{n+k}\bigl((X_+)\wedge\MO(k)\bigr).
\]
As in the proof of Proposition~\ref{borcobor}, there is the map from
the above set of cobordism classes to the bordism group $O_n(X)$,
which forgets the embedding $M\subset\R^{n+k}$. A stabilisation
argument shows that
\begin{equation}\label{oborgr}
  O_n(X)=\lim_{k\to\infty}\pi_{k+n}\bigl((X_+)\wedge\MO(k)\bigr),
\end{equation}
providing a homotopical interpretation for the bordism groups
of~$X$.

We define the \emph{cobordism groups} of $X$ as
\begin{equation}\label{ocobgr}
  O^n(X)=\lim_{k\to\infty}\bigl[\varSigma^{k-n}(X_+),\MO(k)\bigr].
\end{equation}

If $X$ is a (not necessarily compact) manifold, then the groups
$O^n(X)$ may also be obtained by stabilising the set of cobordism
classes of $\eta$-submanifolds in~$E\xi$. Namely, we need to set
$\xi=\underline{\R}^{k-n}$ (the trivial $(k-n)$-plane bundle
over~$X$), and $\eta=\eta_k$.
%Then $\eta$-submanifolds in $E\xi$ are described by pairs $(g,j)$,
%where $g\colon M\to X$, $\dim X-\dim M=n$ and $j\colon M\to\R^{k-n}$
%(the latter map may be assumed to be an embedding after stabilisation).
In other words, a cobordism class in $O^n(X)$ is described by the
composition
\begin{equation}\label{geomonx}
  M\hookrightarrow X\times\R^{k-n}\longrightarrow X,
\end{equation}
where the first map is an embedding of codimension~$k$.

The maps of Thom spaces $\MO(k)\wedge\MO(l)\to\MO(k+l)$ (induced
by the classifying maps $\eta_k\times\eta_l\to\eta_{k+l}$) turn
$O^*(X)=\prod_{n\in\Z}O^n(X)$ into a graded ring, called the
\emph{unoriented cobordism ring of~$X$}\label{unorcobo}. (One
needs to consider direct product instead of direct sum to take
care of infinite complexes like~$\R P^\infty$.)
\end{construction}

\subsection*{Exercises.}
\begin{exercise}
Prove Proposition~\ref{thomprod}.
\end{exercise}

\begin{exercise}
If $\eta$ is the tautological line bundle over $\R P^n$
(respectively, $\C P^n$), then $\Th\eta$ can be identified with $\R
P^{n+1}$ (respectively, $\C P^{n+1}$).
\end{exercise}

\begin{exercise}
Prove that cobordism of $\eta$-submanifolds is an equivalence
relation.
\end{exercise}

\begin{exercise}
Complete the proof of Theorem~\ref{etacob}.
\end{exercise}

\begin{exercise}
The forgetful map $\varOmega_O^{-n,k}\to\varOmega^O_n$ is a
homomorphism of groups.
\end{exercise}

\begin{exercise}\label{oorient}
Given any $(k-n)$-plane bundle $\xi$ over a manifold~$X$ and an
embedding $M\hookrightarrow E\xi$ of codimension~$k$, the
composition
\[
  M\hookrightarrow E\xi\longrightarrow X
\]
determines a cobordism class in~$O^n(X)$. (Hint: reduce
to~\eqref{geomonx} by embedding $\xi$ into a trivial bundle over the
same~$X$.)
\end{exercise}

\begin{exercise}[{Poincar\'e--Atiyah duality~\cite{atiy61} in unoriented bordism}]
If $X$ is an $n$-dimensional manifold, then
\[
  O^{n-k}(X)=O_k(X)\quad\text{for any }k.
\]
In particular, for $X=pt$ we obtain the isomorphisms of
Proposition~\ref{borcobor}.
\end{exercise}

\section{Oriented and complex bordism}\label{orcobbord}
%MA
The bordism relation may be extended to manifolds endowed with some
additional structure, which leads to important bordism theories.
%To take account of the additional structure in the
%definition of bordism one requires that $\partial
%W=M_1\sqcup\overline{M}_2$, where the structure on $\partial W$ is
%induced from that on $W$, and $\overline{M}$ denotes the manifold
%with the opposite structure.

The simplest additional structure is an orientation. By
definition, two oriented $n$-dimensional manifolds $M_1$ and $M_2$
are \emph{oriented bordant}\label{defnorbor} if there is an
oriented $(n+1)$-dimensional manifold $W$ with boundary such that
$\partial W=M_1\sqcup\overline{M}_2$, where $\overline M_2$
denotes $M_2$ with the orientation reversed. The \emph{oriented
bordism groups} $\varOmega_n^{SO}$ and the \emph{oriented bordism
ring} $\varOmega^{SO}=\bigoplus_{n\ge0}\varOmega_n^{SO}$ are
defined accordingly. Given an oriented manifold $M$, the manifold
$M\times I$ has the canonical orientation such that
$\partial(M\times I)=M\sqcup\overline M$. Hence,
$-[M]=[\overline{M}]$ in $\varOmega_n^{SO}$. Unlike
$\varOmega_n^O$, elements of $\varOmega^{SO}$ generally do not
have order~2.

Complex structure gives another important example of an additional
structure on manifolds. However, a direct attempt to define the
bordism relation on complex manifolds fails because the manifold
$W$ must be odd-dimensional and therefore cannot be complex. This
can be remedied by considering \emph{stably complex} (also known
as \emph{stably almost complex} or \emph{quasicomplex})
structures.

%MA
Let ${\mathcal T}\!M$ denote the tangent bundle of~$M$. We say that
$M$ admits a \emph{tangential stably complex structure} if there is
an isomorphism of real vector bundles
\begin{equation}\label{scs}
  c_{\mathcal T}\colon {\mathcal T}\!M\oplus \underline{\R}^k\to \xi
\end{equation}
between the `stable' tangent bundle and a complex vector
bundle~$\xi$ over~$M$. Some of the choices of such isomorphisms are
deemed to be equivalent, that is, determine the same stably complex
structures. This equivalence relation is generated by
\begin{itemize}
\item[(a)] additions of trivial complex summands; that is, $c_{\mathcal
T}$ is equivalent to
\[
  {\mathcal T}\!M\oplus \underline{\R}^k\oplus\underline\C
  \stackrel{c_{\mathcal T}\oplus\mathop{\mathrm{id}}}\lllra
  \xi\oplus\underline\C,
\]
where $\underline\C$ in the left hand side is canonically
identified with $\underline{\R}^2$;
\item[(b)] compositions with isomorphisms of complex bundles; that
is, $c_{\mathcal T}$ is equivalent to $\varphi\cdot c_{\mathcal
T}$ for every $\C$-linear isomorphism $\varphi\colon\xi\to\zeta$.
\end{itemize}
The equivalence class of $c_{\mathcal T}$ may be described
homotopically as the equivalence class of lifts of the map $M\to
BO(2N)$ classifying the stable tangent bundle to a map $M\to
BU(N)$ up to homotopy and stabilisation
(see~\cite[Chapters~II,~VII]{ston68}). A \emph{tangential stably
complex manifold} is a pair consisting of $M$ and an equivalence
class of isomorphisms $c_{\mathcal T}$; we shall use a simplified
notation $(M,c_{\mathcal T})$ for such pairs. This notion is a
generalisation of complex and \emph{almost complex}\label{tscmani}
manifolds (where the latter means a manifold with a choice of a
complex structure on ${\mathcal T}\!M$, that is, a stably complex
structure~\eqref{scs} with $k=0$).

We say that $M$ admits a \emph{normal complex
structure}\label{normalcs} if there is an embedding $i\colon
M\hookrightarrow\mathbb R^N$ with the property that the normal
bundle $\nu(i)$ admits a structure of a complex vector bundle.
There is an appropriate notion of stable equivalence for such
embeddings~$i$, and a normal complex structure $c_\nu$ on $M$ is
defined as the corresponding equivalence class. Tangential and
normal stably complex structures on $M$ determine each other by
means of the canonical isomorphism $\mathcal
T\!M\oplus\nu(i)\cong\underline{\R}^N$.
%We therefore may restrict our attention to tangential structures only.

\begin{example}\label{2cp1}
Let $M=\mathbb{C}P^1$. The standard complex structure on $M$ is
equivalent to the stably complex structure determined by the
isomorphism
\[
  {\mathcal T}(\mathbb{C}P^1)\oplus\underline{\R}^2\stackrel{\cong}{\longrightarrow}
  \overline{\eta}\oplus \overline{\eta}
\]
where $\eta$ is the tautological line bundle. On the other hand,
one can view $\C P^1$ as $S^2$ embedded into $\R^4\cong\C^2$ with
trivial normal bundle. We therefore have an isomorphism
\[
  {\mathcal T}(\mathbb{C}P^1)\oplus\underline{\R}^2\stackrel{\cong}{\longrightarrow}
  \underline{\C}^2\cong\eta\oplus \overline{\eta}
\]
which determines a trivial stably complex structure on~$\mathbb C
P^1$.
\end{example}

The bordism relation can be defined between stably complex
manifolds by taking account of the stably complex structure in the
bordism relation. As in the case of unoriented bordism, the set of
bordism classes $[M,c_{\mathcal T}]$ of $n$-dimensional stably
complex manifolds is an Abelian group with respect to disjoint
union. This group is called the \emph{$n$-dimensional complex
bordism group}\label{defncomplbor} and denoted by $\varOmega^U_n$.
The sphere $S^n$ has the canonical normally complex structure
determined by a complex structure on the trivial normal bundle of
the embedding $S^n\hookrightarrow\R^{n+2}$. The corresponding
bordism class represents the zero element in~$\varOmega^U_n$.
%The zero element is
%represented by the bordism class of any manifold $M$ which bounds
%and whose stable tangent bundle is trivial (and therefore
%isomorphic to a product complex vector bundle $M\times\C^k$).
%The sphere $S^n$ provides an example of such a manifold.
The opposite element to the bordism class $[M,c_{\mathcal T}]$ in
the group $\varOmega^U_n$ may be represented by the same manifold
$M$ with the stably complex structure determined by the isomorphism
\[
  {\mathcal T}\!M\oplus\underline{\R}^k\oplus\underline{\C}\stackrel{c_{\mathcal
  T}\oplus\tau}{\lllra}\xi\oplus\underline{\C}
\]
where $\tau\colon\C\to\C$ is the complex conjugation. We shall use
the abbreviated notations $[M]$ and $[\overline M]$ for the complex
bordism class and its opposite whenever the stably complex structure
$c_{\mathcal T}$ is clear from the context. There is a stably
complex structure on $M\times I$ such that $\partial(M\times
I)=M\sqcup\overline M$.

The direct product of stably complex manifolds turns
$\varOmega^U=\bigoplus_{n\ge0}\varOmega_n^U$ into a graded ring,
called the \emph{complex bordism ring}.

\begin{construction}[homotopic approach to cobordism]\label{defccobord}
The \emph{complex bordism groups}
$U_n(X)$ and \emph{cobordism groups} $U^n(X)$ may be defined
homotopically similarly to~\eqref{oborgr} and~\eqref{ocobgr}:
\begin{equation}
\label{comcobgr}
\begin{aligned}
  U_n(X)&=\lim_{k\to\infty}\pi_{2k+n}\bigl((X_+)\wedge\MU(k)\bigr),\\
  U^n(X)&=\lim_{k\to\infty}\bigl[\varSigma^{2k-n}(X_+),\MU(k)\bigr],
\end{aligned}
\end{equation}
where $\MU(k)$ is the Thom space of the universal complex
$k$-plane bundle $\EU(k)\to\BU(k)$. Here the direct limit uses the
maps $\varSigma^2MU(k)\to MU(k+1)$.
\end{construction}

\begin{construction}[geometric approach to
cobordism]\label{defgeomcob} Both groups $U_n(X)$ and $U^n(X)$ may
also be defined geometrically, in a way similar to the geometric
construction of unoriented bordism and cobordism groups
(Construction~\ref{oborspace}). The complex bordism group $U_n(X)$
consists of bordism classes of maps $M\to X$ of stably complex
$n$-dimensional manifolds $M$ to~$X$.

The complex cobordism group $U^n(X)$ of a manifold $X$ may be
defined via cobordisms of $\eta$-submanifolds in $E\xi$, like in
the unoriented case. Let $\xi=X\times\R^{2k-n}$ (a trivial
bundle), and let $\eta$ be the canonical (universal) complex
$k$-plane bundle over~$BU(k)$. Then an $\eta$-submanifold $M$ in
$E\xi$ defines a composite map of manifolds
\[
  M\hookrightarrow X\times\R^{2k-n}\longrightarrow X
\]
(compare~\eqref{geomonx}), where the first map is an embedding
whose normal bundle has a structure of a complex $k$-plane bundle,
and the second map is the projection onto the first factor. A map
$M\to X$ between manifolds which can be decomposed as above is
said to be \emph{complex orientable of
codimension~$n$}\label{compormap}. A choice of this decomposition
together with a complex bundle structure in the normal bundle is
called a \emph{complex orientation} of the map $M\to X$. As usual,
the equivalence relation on the set of complex orientations of
$M\to X$ is generated by bundle isomorphisms and stabilisations.
The group $U^n(X)$ consists of cobordism classes of complex
oriented maps $M\to X$ of codimension~$n$.

Let $y\in U^n(Y)$ be a cobordism class represented by a complex
oriented map $M\to Y$, and let $f\colon X\to Y$ be a map of
manifolds. If these two maps are transverse, the cobordism class
$f^*(y)\in U^n(X)$ is represented by the pullback $X\times_Y M\to
X$ with the induced complex orientation.

When $M\to X$ is a fibre bundle with fibre~$F$, the normal
structure used in the definition of a complex orientation can be
converted to a tangential structure. Namely, an equivalence class
of complex orientations of the bundle projection $M\to X$ is
determined by a choice of stably complex structure for the bundle
$\mathcal T_F(M)$ of tangents along the fibres of $M\to X$ (an
exercise). Such a bundle $M\to X$ is called \emph{stably
tangentially complex}.
\end{construction}

The equivalence of the homotopic and geometric approaches to
cobordism is established using transversality arguments and the
Pontryagin--Thom construction, as in the unoriented case.

If $X=pt$, then we obtain
\[
  U^{-n}(pt)=U_n(pt)=\varOmega^U_n
\]
for $n\ge0$, from either the homotopic or geometric description of
the (co)bordism groups. We also set $\varOmega_U^{-n}=U^{-n}(pt)$
and $\varOmega_U=\bigoplus_{n\ge0}\varOmega_U^{-n}$.

\begin{construction}[pairing and products]\label{cobproducts}
The product operations in cobordism are defined using the maps of
Thom spaces $\MU(k)\wedge\MU(l)\to\MU(k+l)$ induced by the
classifying maps of the products of canonical bundles.

There is a canonical pairing (the \emph{Kronecker product})
\[
  \langle\;\,,\,\rangle\colon U^m(X)\otimes U_n(X)\to
  \varOmega^U_{n-m},
\]
the \emph{$\frown$-product}
\[
  \frown\colon U^m(X)\otimes U_n(X)\to
  U_{n-m}(X),
\]
and the \emph{$\smile$-product} (or simply \emph{product})
\[
  \smile\colon U^m(X)\otimes U^n(X)\to
  U^{m+n}(X),
\]
defined as follows. Assume given a cobordism class $x\in U^m(X)$
represented by a map $\varSigma^{2l-m}X_+\to MU(l)$ and a bordism
class $\alpha\in U_n(X)$ represented by a map $S^{2k+n}\to
X_+\wedge MU(k)$. Then $\langle
x,\alpha\rangle\in\varOmega^U_{n-m}$ is represented by the
composite map
\[
\begin{CD}
  S^{2k+2l+n-m} @>\varSigma^{2l-m}\alpha>>
  \varSigma^{2l-m}X_+\wedge MU(k)
  @>x\wedge\,\id>> MU(l)\wedge MU(k)\to MU(l{+}k)
\end{CD}
\]
If $\varDelta\colon X_+\to (X\times X)_+=X_+\wedge X_+$ is the
diagonal map, then $x\frown\alpha\in U_{n-m}(X)$ is represented by
the composite map
\begin{gather*}
\begin{CD}
  S^{2k+2l+n-m} @>\varSigma^{2l-m}\alpha>>
  \varSigma^{2l-m}X_+\wedge MU(k) @>\varSigma^{2l-m}\!\varDelta\,\wedge\,\id>>
  X_+\wedge\varSigma^{2l-m}X_+\wedge MU(k)
\end{CD}\\
\begin{CD}
  @>\id\wedge x\wedge\id>> X_+\wedge MU(l)\wedge MU(k)\to X_+\wedge MU(l+k)
\end{CD}
\end{gather*}
The $\smile$-product is defined similarly; it turns
$U^*(X)=\prod_{n\in\Z}U^n(X)$ into a graded ring, called the
\emph{complex cobordism ring of~$X$}. It is a module
over~$\varOmega_U$.

The products operations can be also interpreted geometrically.

For example, assume that $x\in U^m(X)$ is represented by an
embedding of manifolds $M^{k-m}\to X=X^k$ with a complex structure
in the normal bundle, and $\alpha\in U_n(X)$ is represented by an
embedding $N^n\to X^k$ of a tangentially stably complex
manifold~$N$. Assume further that $M$ and $N$ intersect
transversely in~$X$, i.e. $\dim M\cap N=n-m$. Then $\langle
x,\alpha\rangle$ is the bordism class of the intersection $M\cap
N$, and $x\frown\alpha$ is the bordism class of the embedding
$M\cap N\to X$. The tangential complex structure of $M\cap N$ is
defined by the tangential structure of $N$ and the complex
structure in the normal bundle of $M\cap N\to N$ induced from the
normal bundle of $M\to X$.

Similarly, if $x\in U^{-d}(X)$ is represented by a smooth fibre
bundle $E^{k+d}\to X^k$ and $\alpha\in U_n(X)$ is represented by a
smooth map $N\to X$, then $\langle
x,\alpha\rangle\in\varOmega^U_{n+d}$ is the bordism class of the
pull-back $E'$, and $x\frown\alpha\in U_{n+d}(X)$ is the bordism
class of the composite map $E'\to X$ in the pull-back diagram
\[
\begin{CD}
  E'@>>> E\\
  @VVV @VVV\\
  N @>>> X
\end{CD}
\]
\end{construction}

\begin{construction}[Poincar\'e--Atiyah duality in cobordism]\label{padual}
Let $X$ be a manifold of dimension $d$. The inclusion of a point
$pt\subset X$ defines the bordism class $1\in U_0(X)$ and the
\emph{fundamental cobordism class of $X$} in $U^d(X)$. (The normal
bundle of a point has a complex structure if $d$ is even,
otherwise the normal bundle of a point in $X\times\R$ has a
complex structure.)

The identity map $X\to X$ defines the cobordism class $1\in
U^0(X)$. It defines the \emph{fundamental bordism class of~$X$} in
$U_d(X)$ only when $X$ is stably complex.

Now let $X$ be stably complex manifold with fundamental bordism
class $[X]\in U_d(X)$. The map
\[
  D=\;\cdot\frown\![X]\colon U^k(X)\to U_{d-k}(X),\quad x\mapsto x\frown[X]
\]
is an isomorphism (an exercise); it is called the
\emph{Poincar\'e--Atiyah duality} map.
\end{construction}

\begin{construction}[Gysin homomorphism]\label{gysincob} Let $f\colon X^k\to Y^{k+d}$ be a complex
oriented map of codimension~$d$ between manifolds (manifolds may
be not compact, in which case $f$ is assumed to be proper). It
induces a covariant map
\[
  f_{\,!}\colon U^n(X)\to U^{n+d}(Y)
\]
called the \emph{Gysin homomorphism}, whose geometric definition
is as follows. Let $x\in U^n(X)$ be represented by a complex
oriented map $g\colon M^{k-n}\to X^k$. Then $f_{\,!}(x)$ is
represented by the composition~$fg$.
\end{construction}

\begin{proposition}\label{gysinprop}
The Gysin homomorphism has the following properties:
\begin{itemize}
\item[(a)] $f_{\,!}\colon U^*(X)\to U^{*+d}(Y)$ depends only on
the homotopy class of~$f$;
\item[(b)] $f_{\,!}$ is a homomorphism of $\varOmega_U$-modules;
\item[(c)] $(fg)_{\,!}=f_{\,!}g_{\,!}$;
\item[(d)] $f_{\,!}(x\cdot f^*(y))=f_{\,!}(x)\cdot y$ for any $x\in
U^n(X)$, $y\in U^m(Y)$;
\item[(e)] assume that\\[-15pt]
\[
\begin{CD}
  X\times_Y Z @>g'>> X\\
  @Vf'VV @VVfV\\
  Z @>g>> Y
\end{CD}
\]
is a pullback square of manifolds, where $g$ is transverse to $f$
and $f'$ is endowed with the pullback of the complex orientation
of~$f$. Then
\[
  g^*f_{\,!}=f'_{\,!}\,g'^*\colon U^*(X)\to U^{*+d}(Z).
\]
\end{itemize}
\end{proposition}
\begin{proof}
(a) is clear from the homotopic definition of cobordism, while the
other properties follow easily from the geometric definition. For
example, to prove~(d) one needs to choose maps $Z\to X$ and $W\to
Y$ representing $x$ and $y$, respectively, and consider the
commutative diagram
\[
\begin{CD}
  Z\times_Y W @>>> X\times_Y W @>>> W\\
  @VVV @VVV @VVV\\
  Z @>>> X @>f>> Y
\end{CD}
\]
Both sides of~(d) are represented by the composite map $Z\times_Y
W\to Y$ above.
\end{proof}

Let $\xi$ be a complex $n$-plane bundle with total space $E$ over
a manifold~$X$ and let $i\colon X\to E$ be the zero section. The
element $i^*i_{\,!}1\in U^{2n}(X)$ where $1\in U^0(X)$ is called
the \emph{Euler class}\label{eulerclcob} of~$\xi$ and is denoted
by~$e(\xi)$.

%Here is another example of situation when certain cobordism classes
%can be represented very explicitly by maps of manifolds.

\begin{construction}[geometric cobordisms]\label{geomcob}
For any cell complex $X$ the cohomology group $H^2(X)$ can be
identified with the set $[X,\C P^\infty]$ of homotopy classes of
maps into~$\C P^\infty$. Since $\C P^\infty=MU(1)$, it follows
from~\eqref{comcobgr} that every element $x\in H^2(X)$ determines
a cobordism class $u_x\in U^2(X)$. The elements of $U^2(X)$
obtained in this way are called \emph{geometric cobordisms}
of~$X$. We therefore may view $H^2(X)$ as a subset in $U^2(X)$,
however the group operation in $H^2(X)$ is not obtained by
restricting the group operation in $U^2(X)$ (the relationship
between the two operations is discussed in Appendix~\ref{genera}).

When $X$ is a manifold, geometric cobordisms may be described by
submanifolds $M\subset X$ of codimension~2 with a fixed complex
structure on the normal bundle.

Indeed, every $x\in H^2(X)$ corresponds to a homotopy class of
maps $f_x\colon X\to\C P^\infty$. The image $f_x(X)$ is contained
in some $\C P^N\subset\C P^\infty$, and we may assume that
$f_x(X)$ is transverse to a certain hyperplane $H\subset\C P^N$.
Then $M_x=f_x^{-1}(H)$ is a codimension-2 submanifold in~$X$ whose
normal bundle acquires a complex structure by restriction of the
complex structure on the normal bundle of $H\subset\C P^N$.
Changing the map $f_x$ within its homotopy class does not affect
the bordism class of the embedding $M_x\to X$.

Conversely, assume given a submanifold $M\subset X$ of codimension
2 whose normal bundle is endowed with a complex structure. Then
the composition
\[
  X\to\Th(\nu)\to MU(1)=\C P^\infty
\]
of the Pontryagin--Thom collapse map $X\to\Th(\nu)$ and the map of
Thom spaces corresponding to the classifying map $M\to BU(1)$
of~$\nu$ defines an element $x_M\in H^2(X)$, and therefore a
geometric cobordism.

If $X$ is an oriented manifold, then a choice of complex structure
on the normal bundle of a codimension-2 embedding $M\subset X$ is
equivalent to orienting~$M$. The image of the fundamental class of
$M$ in $H_*(X)$ is Poincar\'e dual to $x_M\in H^2(X)$.
\end{construction}

\begin{construction}[connected sum]\label{cobcs}
For manifolds of positive dimension the disjoint union $M_1\sqcup
M_2$ representing the sum of bordism classes $[M_1]+[M_2]$ may be
replaced by their \emph{connected sum}, which represents the same
bordism class.

The connected sum $M_1\cs M_2$ of manifolds $M_1$ and $M_2$ of the
same dimension $n$ is constructed as follows. Choose points
$v_1\in M_1$ and $v_2\in M_2$, and take closed $\varepsilon$-balls
$B_\varepsilon(v_1)$ and~$B_\varepsilon(v_2)$ around them (both
manifolds may be assumed to be endowed with a Riemannian metric).
Fix an isometric embedding $f$ of a pair of standard
$\varepsilon$-balls $D^n\times S^0$ (here $S^0=\{0,1\}$) into
$M_1\sqcup M_2$ which maps $D^n\times0$ onto $B_\varepsilon(v_1)$
and $D^n\times1$ onto $B_\varepsilon(v_2)$. If both $M_1$ and
$M_2$ are oriented we additionally require the embedding $f$ to
preserve the orientation on the first ball and reverse in on the
second. Now, using this embedding, replace in $M_1\sqcup M_2$ the
pair of balls $D^n\times S^0$ by a `pipe' $S^{n-1}\times D^1$.
After smoothing the angles in the standard way we obtain a smooth
manifold $M_1\cs M_2$.

If both $M_1$ and $M_2$ are connected the smooth structure on
$M_1\cs M_2$ does not depend on a choice of points $v_1$, $v_2$
and embedding $D^n\times S^0\hookrightarrow M_1\sqcup M_2$. It
does however depend on the orientations; $M_1\cs M_2$ and
$M_1\cs\overline{M_2}$ are not diffeomorphic in general. For
example, the manifolds $\mathbb CP^2\cs\mathbb CP^2$ and $\mathbb
CP^2\cs\overline{\C P^2\!}$ are not diffeomorphic (and not even
homotopy equivalent because they have different signatures).
%if at least one of the
%manifolds $M_1$ and $M_2$ is non-orientable or \emph{reversible}
%(i.e. admits an orientation-reversing diffeomorphism). Otherwise
%the \emph{connected sum of oriented manifolds} needs to be
%defined; in its definition one requires the embedding $f$ to
%preserve the orientation on the first ball and reverse in on the
%second.

There are smooth contraction maps $p_1\colon M_1\cs M_2\to M_1$
and $p_2\colon M_1\cs M_2\to M_2$. In the oriented case the
manifold $M_1\cs M_2$ can be oriented in such a way that both
contraction maps preserve the orientations.

A bordism between $M_1\sqcup M_2$ and $M_1\cs M_2$ may be
constructed as follows. Consider a cylinder $M_1\times I$, from
which we remove an $\varepsilon$-neighbourhood
$U_\varepsilon(v_1\times1)$ of the point $v_1\times 1$. Similarly,
remove the neighbourhood $U_\varepsilon(v_2\times1)$ from
$M_2\times I$ (each of these two neighbourhoods can be identified
with the half of a standard open $(n+1)$-ball). Now connect the
two
%MA
remainders of cylinders by a `half pipe' $S^n_\le\times I$ in such
a way that the half-sphere $S^n_\le\times 0$ is identified with
the half-sphere on the boundary of $U_\varepsilon(v_1\times1)$,
and $S^n_\le\times1$ is identified with the half-sphere on the
boundary of $U_\varepsilon(v_2\times 1)$.
%MA
Smoothing the angles we obtain a manifold with boundary $M_1\sqcup
M_2\sqcup(M_1\cs M_2)$ (or
$\overline{M_1}\sqcup\overline{M_2}\sqcup(M_1\cs M_2)$ in the
oriented case), see Fig.~\thechapter.2.

\begin{figure}
\label{figc}
\begin{center}
\end{center}
\caption{Disjoint union and connected sum.}
\end{figure}

Finally, if $M_1$ and $M_2$ are stably complex manifolds, then
there is a canonical stably complex structure on $M_1\cs M_2$,
which is constructed as follows. Assume the stably complex
structures on $M_1$ and $M_2$ are determined by isomorphisms
\[
  c_{\mathcal T\!,1}\colon{\mathcal T}\!M_1\oplus\underline{\R}^{k_1}\to\xi_1
  \quad\text{and}\quad
  c_{\mathcal T\!,2}\colon{\mathcal T}\!M_2\oplus\underline{\R}^{k_2}\to\xi_2.
\]
Using the isomorphism ${\mathcal T}(M_1\cs
M_2)\oplus\underline{\R}^n\cong p_1^*{\mathcal T}\!M_1\oplus
p_2^*{\mathcal T}\!M_2$, we define a stably complex structure on
$M_1\cs M_2$ by the isomorphism
\begin{multline*}
  {\mathcal T}(M_1\cs M_2)\oplus\underline{\R}^{n+k_1+k_2}\\
  \cong
  p_1^*{\mathcal T}\!M_1\oplus\underline{\R}^{k_1}\oplus
  p_2^*{\mathcal T}\!M_2\oplus\underline{\R}^{k_2}
  \xrightarrow{c_{{\mathcal T},1}\oplus c_{{\mathcal T},2}}
  p_1^*\xi_1\oplus p_2^*\xi_2.
\end{multline*}
We shall refer to this stably complex structure as the
\emph{connected sum of stably complex
structures}\label{consumstcom} on $M_1$ and $M_2$. The
corresponding complex bordism class is $[M_1]+[M_2]$.
\end{construction}

\subsection*{Exercises.}
\begin{exercise}
Assume given a complex $(k-l)$-plane bundle $\xi$ over a
manifold~$X$ and an embedding $M\hookrightarrow E\xi$ whose normal
bundle has a structure of a complex $k$-plane bundle. Then the
composition
\[
  M\hookrightarrow E\xi\longrightarrow X
\]
determines a complex orientation for the map $M\to X$ of
codimension~$2l$, and therefore a complex cobordism class
in~$U^{2l}(X)$. (Compare Exercise~\ref{oorient}.)

Similarly, an embedding $M\hookrightarrow E(\xi\oplus\underline\R)$
whose normal bundle has a structure of a complex $k$-plane bundle
determines a complex cobordism class in~$U^{2l-1}(X)$, via the
composition
\[
  M\hookrightarrow E(\xi\oplus\underline\R)\longrightarrow X.
\]
This is how complex orientations were defined in~\cite{quil71}.
\end{exercise}

\begin{exercise}
Let $\pi\colon E\to B$ be a bundle with fibre~$F$. The map $\pi$
is complex oriented if and only a stably complex structure is
chosen for the bundle $\mathcal T_F(E)$ of tangents along the
fibres of~$\pi$.
\end{exercise}

\begin{exercise}
The Poincar\'e--Atiyah duality map
\[
  D=\;\cdot\frown\![X]\colon U^k(X)\to U_{d-k}(X),\quad x\mapsto x\frown[X]
\]
is an isomorphism for any stably complex manifold~$X$ of
dimension~$d$.
\end{exercise}

\begin{exercise}
Let $f\colon X^d\to Y^{p+d}$ be a complex oriented map of
manifolds, and let $D_X\colon U^k(X)\to U_{d-k}(X)$, $D_Y\colon
U^{p+k}(Y)\to U_{d-k}(Y)$ be the duality isomorphisms for $X$,
$Y$. Then the Gysin homomorphism satisfies $f_{\,!}=D^{-1}_Y f_*
D_X$.
\end{exercise}

\begin{exercise}\label{gysin-thom}
Let $\xi$ be a complex $n$-plane bundle over a manifold $M$ with
total space~$E$, and let $i\colon M\to E$ be the inclusion of zero
section. Define the Gysin homomorphism
\[
  i_{\,!}\colon U^*(M)\to U^{*+2n}(E,E\setminus M)=U^{*+2n}(\Th(\xi))
\]
by analogy with Construction~\ref{gysincob} and show that
$i_{\,!}$ is an isomorphism. It is called the \emph{Gysin--Thom
isomorphism} corresponding to~$\xi$.
\end{exercise}

\section{Structure results}
Let $M$ be an $n$-dimensional manifold and let $f\colon M\to
BO(n)$ be the classifying map of the tangent bundle. Given a
universal Stiefel--Whitney characteristic class $w\in
H^*(BO(n);\Z_2)=\Z_2[w_1,\ldots,w_n]$ of $n$-plane bundles, the
corresponding corresponding \emph{tangential Stiefel--Whitney
characteristic number}\label{charnumbersd} $w[M]$ is defined as
the result of pairing of $f^*(w)\in H^n(M;\Z_2)$ with the
fundamental class $\langle M\rangle\in H_n(M;\Z_2)$. The number
$w[M]$ is an unoriented bordism invariant.

\emph{Tangential Chern characteristic numbers} $c[M]$ of stably
complex $2n$-manifolds and \emph{Pontryagin characteristic
numbers} $p[M]$ of oriented $4n$-manifolds are defined similarly;
they are complex and oriented bordism invariants, respectively.

Normal characteristic numbers are of equal importance in cobordism
theory. If $M\hookrightarrow\R^N$ is an embedding with a fixed
complex structure in the normal bundle~$\nu$, classified by the
map $g\colon\nu\to BU(n)$, then the \emph{normal Chern
characteristic number} $\bar c[M]$ corresponding to $c\in
H^*(BU(n))=\Z[c_1,\ldots,c_n]$ is defined as $(g^*c)\langle
M\rangle$. Normal Stiefel--Whitney and Pontryagin numbers are
defined similarly. Since $\mathcal T M\oplus\nu=\underline{\R}^N$,
the tangential and normal characteristic numbers determine each
other.

In what follows, all characteristic numbers will be tangential.

\medskip

The theory of unoriented (co)bordism was the first to be
completed: the coefficient ring $\varOmega^O$ was calculated by
Thom, and the bordism groups $O_*(X)$ of cell complexes $X$ were
reduced to homology groups of $X$ with coefficients
in~$\varOmega^O$. The corresponding results are summarised as
follows.

\begin{theorem}\label{uocob}\
\begin{itemize}
\item[(a)] Two manifolds are unorientedly bordant if and only if they have
identical sets of Stiefel--Whitney characteristic numbers.

\item[(b)] $\varOmega^O$ is a polynomial ring over $\Z_2$ with
one generator $a_i$ in every positive dimension $i\ne 2^k-1$.

\item[(c)] For every cell complex $X$ the module $O_*(X)$ is a free
graded $\varOmega^O$-module isomorphic to
$H_*(X;\Z_2)\otimes_{\Z_2}\varOmega^O$.
\end{itemize}
\end{theorem}

Parts (a) and (b) were done by Thom~\cite{thom54}. Part~(c) was
first formulated by Conner and Floyd~\cite{co-fl64}; it also
follows from the results of Thom.

%Describing the complex bordism ring $\varOmega^U$ turned out to be
%a more difficult problem:

\begin{theorem}\label{comcob}\
\begin{itemize}
\item[(a)] $\varOmega^U\otimes\Q$ is a polynomial ring over $\Q$ generated by the bordism
classes of complex projective spaces $\C P^i$, \ $i\ge1$.

\item[(b)] Two stably complex manifolds are bordant if and only if they have
identical sets of Chern characteristic numbers.

\item[(c)] $\varOmega^U$ is a polynomial ring over $\Z$ with one generator $a_i$
in every even dimension $2i$, where $i\ge1$.
\end{itemize}
\end{theorem}

Part (a) can be proved by the methods of Thom. Part~(b) follows
from the results of Milnor~\cite{miln60} and
Novikov~\cite{novi60}. Part~(c) is the most difficult one; it was
done by Novikov~\cite{novi60} using the Adams spectral sequence
and structure theory of Hopf algebras (see also~\cite{novi62} for
a more detailed account) and Milnor (unpublished\footnote{Milnor's
proof was announced in~\cite{hirz60}; it was intended to be
included in the second part of~\cite{miln60}, but has never been
published.}). Another more geometric proof was given by
Stong~\cite{ston68}.

Note that part (c) of Theorem~\ref{uocob} does not extend to
complex bordism; $U_*(X)$ is not a free $\varOmega^U$-module in
general (although it is a free $\varOmega^U$-module if $H_*(X;\Z)$
is free abelian). Unlike the case of unoriented bordism, the
calculation of complex bordism of a space $X$ does not reduce to
calculating the coefficient ring $\varOmega^U$ and homology groups
$H_*(X)$. The theory of complex (co)bordism is much richer than
its unoriented analogue, and at the same time is not as
complicated as oriented bordism or other bordism theories with
additional structure, since the coefficient ring does not have
torsion. Thanks to this, % richness and lack of torsion,
complex cobordism theory found many striking and important
applications in algebraic topology and beyond. Many of these
applications were outlined in the pioneering work of
Novikov~\cite{novi67}.

The calculation of the oriented bordism ring was completed by
Novikov~\cite{novi60} (ring structure modulo torsion and odd
torsion) and Wall~\cite{wall60} (even torsion), with important
earlier contributions made by Rokhlin, Averbuch, and Milnor.
Unlike complex bordism, the ring $\varOmega^{SO}$ has additive
torsion. We give only a partial result here (which does not fully
describe the torsion elements).

\begin{theorem}\
\begin{itemize}
\item[(a)] $\varOmega^{SO}\otimes\Q$ is a polynomial ring over $\Q$ generated by the bordism
classes of complex projective spaces $\C P^{2i}$, \ $i\ge1$.

\item[(b)] The subring $\mathrm{Tors}\subset\varOmega^{SO}$ of
torsion elements contains only elements of order~$2$. The quotient
$\varOmega^{SO}/\mathrm{Tors}$ is a polynomial ring over $\Z$ with
one generator $a_i$ in every dimension $4i$, where $i\ge1$.

\item[(c)] Two oriented manifolds are bordant if and only if they have
identical sets of Pontryagin and Stiefel--Whitney characteristic
numbers.
\end{itemize}
\end{theorem}

\section{Ring generators}
To describe a set of ring generators for $\varOmega^U$ we shall
need a special characteristic class of complex vector bundles. Let
$\xi$ be a complex $k$-plane bundle over a manifold~$M$. Write its
total Chern class formally as follows:
\[
  c(\xi)=1+c_1(\xi)+\cdots+c_k(\xi)=(1+x_1)\cdots(1+x_k),
\]
so that $c_i(\xi )=\sigma_i(x_1,\ldots,x_k)$ is the $i$th
elementary symmetric function in formal indeterminates. These
indeterminates acquire a geometric meaning if $\xi$ is a sum
$\xi_1\oplus\cdots\oplus\xi_k$ of line bundles; then
$x_j=c_1(\xi_j)$, \ $1\le j\le k$. Consider the polynomial
\[
  P_n(x_1,\ldots x_k)=x_1^n+\cdots +x_k^n
\]
and express it via the elementary symmetric functions:
\[
  P_n(x_1,\ldots ,x_k)=s_n(\sigma_1,\ldots ,\sigma_k).
\]
Substituting the Chern classes for the elementary symmetric
functions we obtain a certain characteristic class of~$\xi$:
\[
  s_n(\xi)=s_n(c_1(\xi),\ldots,c_k(\xi))\in H^{2n}(M).
\]
This characteristic class plays an important role in detecting the
polynomial generators of the complex bordism ring, because of the
following properties (which follow immediately from the
definition).

\begin{proposition}\label{snprop}The characteristic class $s_n$
satisfies
\begin{itemize}
\item[(a)] $s_n(\xi)=0$ if $\xi$ is a bundle over $M$ and $\dim M<2n$;

\item[(b)]  $s_n(\xi\oplus\eta)=s_n(\xi)+s_n(\eta)$;

\item[(c)] $s_n(\xi)=c_1(\xi)^n$ if $\xi$ is a line bundle.
\end{itemize}
\end{proposition}

Given a stably complex $2n$-manifold $(M,c_{\mathcal T})$, define
its characteristic number
\begin{equation}\label{defsn}
  s_n[M]=s_n(\mathcal T M)\langle M\rangle\in\Z.
\end{equation}
Here $s_n(\mathcal T M)$ is understood to be $s_n(\xi)$, where
$\xi$ is the complex bundle from~\eqref{scs}.

\begin{corollary}
If a bordism class $[M]\in\varOmega_{2n}^U$ decomposes as
$[M_1]\times[M_2]$ where $\dim M_1>0$ and $\dim M_2>0$, then
$s_n[M]=0$.
\end{corollary}

It follows that the characteristic number $s_n$ vanishes on
decomposable elements of $\varOmega^U_{2n}$. It also detects
indecomposables that may be chosen as polynomial generators. The
following result featured in the proof of Theorem~\ref{comcob}.

\begin{theorem}\label{mulgen}
A bordism class $[M]\in\varOmega_{2n}^U$ may be chosen as a
polynomial generator $a_n$ of the ring $\varOmega^U$ if and only
if
\[
  s_n[M]=\begin{cases}
  \pm1,  &\text{if $n\ne p^k-1$ for any prime $p$;}\\
  \pm p, &\text{if $n=p^k-1$ for some prime $p$.}
  \end{cases}
\]
\end{theorem}

%MA
There is no universal description of connected manifolds
representing the polynomial generators $a_n\in\varOmega^U$. On the
other hand, there is a particularly nice family of manifolds whose
bordism classes generate the whole ring $\varOmega^U$. This family
is redundant though, so there are algebraic relations between
their bordism classes.

\begin{construction}[Milnor hypersurfaces]
Fix a pair of integers $j\ge i\ge0$ and consider the product $\C
P^i\times\C P^j$. Its algebraic subvariety
\begin{equation}
\label{hij}
  H_{ij}=\{
  (z_0:\cdots :z_i)\times (w_0:\cdots :w_j)\in
    \mathbb{C}P^i\times \mathbb{C}P^j\colon z_0w_0+\cdots +z_iw_i=0\}
\end{equation}
is called a \emph{Milnor hypersurface}. Note that $H_{0j}\cong\C
P^{j-1}$.
\end{construction}

Denote by $p_1$ and $p_2$ the projections of $\C P^i\times\C P^j$
onto its factors. Let $\eta$ be the tautological line bundle over
a complex projective space and $\bar\eta$ its conjugate (the
canonical line bundle). We have
\[
  H^*(\C P^i\times\C P^j)=\Z[x,y]/(x^{i+1}=0,\;y^{j+1}=0)
\]
where $x=p_1^*c_1(\bar\eta)$, $y=p_2^*c_1(\bar\eta)$.

\begin{proposition}
The geometric cobordism in $\C P^i\times\C P^j$ corresponding to
the element $x+y\in H^2(\C P^i\times\C P^j)$ is represented by the
submanifold~$H_{ij}$. In particular, the image of the fundamental
class $\langle H_{ij}\rangle$ in $H_{2(i+j-1)}(\C P^i\times\C
P^j)$ is Poincar\'e dual to $x+y$.
\end{proposition}
\begin{proof}
We have $x+y=c_1(p_1^*(\bar\eta)\otimes p_2^*(\bar\eta))$. The
classifying map $f_{x+y}\colon \C P^i\times\C P^j\to\C P^\infty$
is the composition of the \emph{Segre embedding}\label{segreemb}
\begin{align*}
 \sigma\colon  \C P^i\times\C P^j&\to\C P^{(i+1)(j+1)-1},\\
  (z_0:\cdots :z_i)\times (w_0:\cdots :w_j)&\mapsto
  (z_0w_0:z_0w_1:\cdots:z_kw_l:\cdots:z_iw_j),
\end{align*}
and the embedding $\C P^{ij+i+j}\to\C P^\infty$. The codimension 2
submanifold in $\C P^i\times\C P^j$ corresponding to the
cohomology class $x+y$ is obtained as the preimage
$\sigma^{-1}(H)$ of a generally positioned hyperplane in $\C
P^{ij+i+j}$ (that is, a hyperplane $H$ transverse to the image of
the Segre embedding, see Construction~\ref{geomcob}).
By~\eqref{hij}, the Milnor hypersurface is exactly
$\sigma^{-1}(H)$ for one such hyperplane~$H$.
\end{proof}

\begin{lemma}\label{shij}
\[
  s_{i+j-1}[H_{ij}]=\begin{cases}
%MA
  j,&\text{if \ $i=0$};\\
  2,&\text{if \ $i=j=1$};\\
  0,&\text{if \ $i=1$, $j>1$};\\
  -\bin{i+j}i,&\text{if \ $i>1$}.
  \end{cases}
\]
\end{lemma}
\begin{proof}
Let $i=0$. Since the stably complex structure on $H_{0j}=\C
P^{j-1}$ is determined by the isomorphism ${\mathcal T}(\C
P^{j-1})\oplus\C\cong\bar\eta\oplus\cdots\oplus\bar\eta$ ($j$
summands) and $x=c_1(\bar\eta)$, we have
\[
  s_{j-1}[\C P^{j-1}]=jx^{j-1}\langle\C P^{j-1}\rangle=j.
\]

Now let $i>0$. Then
\begin{multline*}
  s_{i+j-1}({\mathcal T}(\C P^i\times\C P^j))=(i+1)x^{i+j-1}+(j+1)y^{i+j-1}\\
  =\begin{cases}
  2x^j+(j+1)y^j,&\text{if $i=1$};\\
  0,&\text{if $i>1$.}
  \end{cases}
\end{multline*}
Denote by $\nu$ the normal bundle of the embedding $\iota\colon
H_{ij}\to\C P^i\times\C P^j$. Then
\begin{equation}\label{tauhij}
  {\mathcal T}(H_{ij})\oplus\nu=\iota^*({\mathcal T}(\C P^i\times\C P^j)).
\end{equation}
Since $c_1(\nu)=\iota^*(x+y)$, we obtain
$s_{i+j-1}(\nu)=\iota^*(x+y)^{i+j-1}$.

Assume $i=1$. Then~\eqref{tauhij} and Proposition~\ref{snprop}
imply that
\begin{multline*}
  s_j[H_{1j}]=s_j\bigl({\mathcal T}(H_{1j})\bigr)\langle H_{1j}\rangle
  =\iota^*\bigl(2x^j+(j+1)y^j-(x+y)^j\bigr)\langle H_{1j}\rangle\\
  =(2x^j+(j+1)y^j-(x+y)^j)(x+y)\langle\C P^1\times\C P^j\rangle=
  \begin{cases}
  2,\text{ if $j=1$};\\
  0,\text{ if $j>1$}.
  \end{cases}
\end{multline*}

Assume now that $i>1$. Then $s_{i+j-1}({\mathcal T}(\C P^i\times\C
P^j))=0$, and we obtain from~\eqref{tauhij} and
Proposition~\ref{snprop} that
\begin{multline*}
  s_{i+j-1}[H_{ij}]=-s_{i+j-1}(\nu)\langle H_{ij}\rangle=
  -\iota^*(x+y)^{i+j-1}\langle H_{ij}\rangle\\
  =-(x+y)^{i+j}\langle\C P^i\times\C P^j\rangle=-\bin{i+j}i.
\end{multline*}
\end{proof}

\begin{remark}
Since $s_1[H_{11}]=2=s_1[\C P^1]$ the manifold $H_{11}$ is bordant
to $\C P^1$. In fact $H_{11}\cong\C P^1$ (an exercise).
\end{remark}

\begin{theorem}\label{hijgen}
The bordism classes $\{[H_{ij}],0\le i\le j\}$ multiplicatively
generate the ring $\varOmega^U$.
\end{theorem}
\begin{proof}
A simple calculation shows that
\[
  \mathop{\text{g.c.d.}}\Bigl(\bin{n+1}i,\;1\le i\le n\Bigr)=
        \begin{cases}
            p, & \text{if $n=p^k-1$,}\\
            1, & \text{otherwise.}
        \end{cases}
\]
Now Lemma~\ref{shij} implies that a certain integer linear
combination of bordism classes $[H_{ij}]$ with $i+j=n+1$ can be
taken as the polynomial generator $a_n$ of $\varOmega^U$, see
Theorem~\ref{mulgen}.
\end{proof}

\begin{remark}
There is no universal description for a linear combination of
bordism classes $[H_{ij}]$ with $i+j=n+1$ giving the polynomial
generator of~$\varOmega^U$.

All algebraic relations between the classes $[H_{ij}]$ arise from
the associativity of the formal group law of geometric cobordism
(see Corollary~\ref{aijgen}).
\end{remark}

\begin{example}
Since $s_1[\C P^1]=2$, $s_2[\C P^2]=3$, the bordism classes $[\C
P^1]$ and $[\C P^2]$ may be taken as polynomial generators $a_1$
and $a_2$ of $\varOmega^U$. However $[\C P^3]$ cannot be taken as
$a_3$, since $s_3[\C P^3]=4$, while $s_3(a_3)=\pm2$. The bordism
class $[H_{22}]+[\C P^3]$ may be taken as~$a_3$.
\end{example}

Theorem~\ref{hijgen} admits the following important addendum,
which is due to Milnor (see~\cite[Chapter~7]{ston68} for the
proof).

\begin{theorem}[Milnor]
Every bordism class $x\in\varOmega_n^U$ with $n>0$ contains a
nonsingular algebraic variety (not necessarily connected).
\end{theorem}

The proof of this fact uses a construction of a (possibly
disconnected) algebraic variety representing the class $-[M]$ for
any bordism class $[M]\in\varOmega_n^U$ of $2n$-dimensional
manifold. The following question is still open.

\begin{problem}[Hirzebruch]
\label{hirz} Describe the set of bordism classes in $\varOmega^U$
containing connected nonsingular algebraic varieties.
\end{problem}

\begin{example}
\label{2dalg} The group $\varOmega^U_2$ is isomorphic to $\Z$ and
is generated by~$[\C P^1]$. Every class $k[\C
P^1]\in\varOmega^U_2$ contains a nonsingular algebraic variety,
namely, a disjoint union of $k$ copies of $\C P^1$ for $k>0$ and a
Riemann surface of genus $(1-k)$ for~$k\le0$. Connected algebraic
varieties are contained only in the classes $k[\C P^1]$
with~$k\le1$.
\end{example}

\subsection*{Exercises.}
\begin{exercise}
Properties listed in Proposition~\ref{snprop} determine the
characteristic class $s_n$ uniquely.
\end{exercise}

\begin{exercise}
Show that $H_{11}\cong\C P^1$.
\end{exercise}

\begin{exercise}\label{h1jcob}
Show that $H_{1j}$ is complex bordant to $\C P^1\times \C
P^{j-1}$. (Hint: calculate the characteristic numbers; no
geometric construction of this bordism is known!)
\end{exercise}

\begin{exercise}
An alternative set of ring generators of $\varOmega^U$ can be
constructed as follows. Let $M^{2n}_k$ be a submanifold in $\C
P^{n+1}$ dual to $kx\in H^2(\C P^{n+1})$, where $k$ is a positive
integer, and $x$ is the first Chern class of the hyperplane
section bundle. For example, one can take $M^{2n}_k$ to be a
nonsingular hypersurface of degree~$k$. Then the set
$\{[M^{2n}_k]\colon n\ge1,\,k\ge1\}$ multiplicatively generates
the ring~$\varOmega^U$.
\end{exercise}

\section{Invariant stably complex structures}\label{secinvscs}
Let $M$ be a $2n$-dimensional manifold with a stably complex
structure determined by the isomorphism
\begin{equation}\label{ctau2}
  c_{\mathcal T}\colon{\mathcal T}\!M\oplus\underline{\R}^{2(l-n)}\to\xi.
\end{equation}
Assume that the torus $T^k$ acts on~$M$.

\begin{definition}\label{iscs}
A stably complex structure $c_{\mathcal T}$ is
\emph{$T^k$-invariant} if for every $\mb t\in T^k$ the composition
\begin{equation} \label{rt}
  r(\mb t)\colon\xi\stackrel{c_{\mathcal T}^{-1}}{\llra}{\mathcal T}\!M\oplus
  \underline{\R}^{2(l-n)}\stackrel{d\mb t\,\oplus\,\mathrm{id}}{\lllra}
  {\mathcal T}\!M\oplus\underline{\R}^{2(l-n)}\stackrel{c_{\mathcal T}}{\longrightarrow}\xi
\end{equation}
is a complex bundle map, where $d\mb t$ is the differential of the
action by~$\mb t$. In other words,~\eqref{rt} determines a
representation $r\colon T^k\to\Hom_\mathbb{C}(\xi,\xi)$.
\end{definition}

Let $x\in M$ be an isolated fixed point of the $T^k$-action
on~$M$. Then we have a representation $r_x\colon T^k\to
GL(l,\mathbb{C})$ in the fibre of $\xi$ over~$x$. This fibre
$\xi_x\cong\mathbb{C}^l$ decomposes as
$\mathbb{C}^n\oplus\mathbb{C}^{l-n}$, where $r_x$ has no trivial
summands on $\C^n$ and is trivial on $\mathbb{C}^{l-n}$. The
nontrivial part of $r_x$ decomposes into a sum
$r_{i_1}\oplus\cdots\oplus r_{i_n}$ of one-dimensional complex
$T^k$-representations. In the corresponding coordinates
$(z_1,\ldots,z_n)$, an element $\mb t=(e^{2\pi
i\varphi_1},\ldots,e^{2\pi i\varphi_k})\in T^k$ acts by
\[
  \mb t\cdot\!(z_1,\ldots,z_n)=(e^{2\pi i\langle\mb w_1,\,\varphi\rangle}z_1,
  \ldots,e^{2\pi i\langle\mb w_n,\,\varphi\rangle}z_n),
\]
where $\varphi=(\varphi_1,\ldots,\varphi_k)\in\R^k$ and $\mb
w_j\in\Z^k$, \ $1\le j\le n$, are the \emph{weights} of the
representation $r_x$ at the fixed point~$x$. Also, the isomorphism
$c_{{\mathcal T},x}$ of~\eqref{ctau2} induces an orientation of
the tangent space ${\mathcal T}_x(M)$.

\begin{definition}\label{defsign}
For any fixed point $x\in M$, the \emph{sign} $\sigma(x)$ is $+1$
if the isomorphism
\[
  {\mathcal T}_x(M)\stackrel{\mathrm{id}\,\oplus\,0}{\lllra}{\mathcal T}_x(M)\oplus
  \mathbb{R}^{2(l-n)}\stackrel{c_{{\mathcal T},x}}{\llra}\xi_x
  \cong\mathbb{C}^n\oplus\mathbb{C}^{l-n}\stackrel{p}{\longrightarrow}
  \mathbb{C}^n,
\]
respects the canonical orientations, and $-1$ if it does not; here
$p$ is the projection onto the first summand.
\end{definition}

So $\sigma(x)$ compares the orientations induced by $r_x$ and
$c_{{\mathcal T},x}$ on ${\mathcal T}_x(M)$. If $M$ is an almost
complex $T^k$-manifold (i.e. $l=n$) then $\sigma(x)=1$ for every
fixed point~$x$.

\begin{example}\label{2cp1eq}
We let the circle $T^1=S^1$ act on $\C P^1$ in homogeneous
coordinates by $t\cdot(z_0:z_1)=(z_0:tz_1)$. This action has two
fixed points $(0:1)$ and $(1:0)$. In the standard stably complex
structure (see Example~\ref{2cp1}) the signs of both vertices are
positive (this structure is complex and the map in
Definition~\ref{defsign} is a complex linear map $\C\to\C$).

On the other hand, the trivial stably complex structure on $\C
P^1\cong S^2$ may be thought of as induced from the embedding of
$S^2$ into $\R^3\oplus\R^1$ with trivial normal bundle, which can
be made complex by identifying it with~$\underline\C$. We
therefore may think of the circle action as the rotation of the
unit sphere in $\R^3$ around the vertical axis. The fixed points
are the north and south poles. The rotations induced on the
tangent planes to the two fixed points are in different
directions. Therefore the signs of the two fixed points are
different (which sign is positive depends on the convention one
uses to induce an orientation on $S^2$ from the orientation
of~$\R^3$).
\end{example}

\chapter{Formal group laws and Hirzebruch genera}\label{genera}
The theory of \emph{formal groups} originally appeared in
algebraic geometry and plays an important role in number theory
and cryptography. Formal groups laws were brought into bordism
theory in the pioneering work of Novikov~\cite{novi67}, and
provided a very powerful tool for the theory of group actions on
manifolds and generalised homology theories. Early applications of
formal group laws in cobordism concerned finite group actions on
manifolds, or `differentiable periodic maps'. Subsequent
developments included constructions of complex oriented cohomology
theories and applications to \emph{Hirzebruch genera}, one of the
most important class of invariants of manifolds.

\section{Elements of the theory of formal group laws}
Let $R$ be a commutative ring with unit.

A formal power series $F(u,v)\in R[[u,v]]$ is called a
(commutative one-dimensional) \emph{formal group law} over $R$ if
it satisfies the following equations:
\begin{itemize}
\item[(a)] $F(u,0)=u$, $F(0,v)=v$;

\item[(b)] $F(F(u,v),w)=F(u,F(v,w))$;

\item[(c)] $F(u,v)=F(v,u)$.
\end{itemize}

The original example of a formal group law over a field $\k$ is
provided by the expansion near the unit of the multiplication map
$G\times G\to G$ in a one-dimensional algebraic group over $\k$.
This also explains the terminology.

A formal group law $F$ over $R$ is called \emph{linearisable} if
there exists a coordinate change $u\mapsto
g_F(u)=u+\sum_{i>1}g_iu^i\in R[[u]]$ such that
\begin{equation}\label{logprop}
  g_F(F(u,v))=g_F(u)+g_F(v).
\end{equation}
Note that every formal group law over $R$ determines a formal
group law over $R\otimes\Q$.

\begin{theorem}\label{logth}
Every formal group law $F$ is linearisable over $R\otimes\Q$.
\end{theorem}
\begin{proof} Consider the series $\omega(u)=\frac{\partial F(u,w)}{\partial
w}\Bigl|_{w=0}$. Then
\[
  \omega(F(u,v))=\frac{\partial
  F(F(u,v),w)}{\partial w}\Bigl|_{w=0}=\frac{\partial
  F(F(u,w),v)}{\partial F(u,w)}\cdot\frac{\partial F(u,w)}{\partial
  w}\Bigl|_{w=0}=\frac{\partial F(u,v)}{\partial u}\omega(u).
\]
We therefore have
$\frac{du}{\omega(u)}=\frac{dF(u,v)}{\omega(F(u,v))}$. Set
\begin{equation}\label{log}
  g(u)=\int_0^u\frac{dv}{\omega(v)};
\end{equation}
then $dg(u)=dg(F(u,v))$. This implies that $g(F(u,v))=g(u)+C$.
Since $F(0,v)=v$ and $g(0)=0$, we get $C=g(v)$. Thus,
$g(F(u,v))=g(u)+g(v)$.
\end{proof}

A series $g_F(u)=u+\sum_{i>1}g_iu^i$ satisfying
equation~\eqref{logprop} is called a
\emph{logarithm}\label{logfglaw} of the formal group law~$F$;
Theorem~\ref{logth} shows that a formal group law over
$R\otimes\Q$ always has a logarithm. Its functional inverse series
$f_F(t)\in R\otimes\Q[[t]]$ is called an \emph{exponential} of the
formal group law, so that we have $F(u,v)=f_F(g_F(u)+g_F(v))$ over
$R\otimes\Q$. If $R$ does not have torsion (i.e. $R\to R\otimes\Q$
is monomorphic), the latter formula shows that a formal group law
(as a series with coefficients in $R$) is fully determined by its
logarithm (which is a series with coefficients in $R\otimes\Q$).

%MA
\begin{example}\label{multfgl}
An example of a formal group law is given by the series
\begin{equation}\label{mfgleq}
  F(u,v)=(1+u)(1+v)-1=u+v+uv,
\end{equation}
over $\Z$, called the \emph{multiplicative formal group law}.
Introducing a formal indeterminate $\beta$ of degree~$-2$, we may
consider the 1-parameter extension of the multiplicative formal
group law, given by $F_\beta(u,v)=u+v-\beta uv$, with coefficients
in~$\Z[\beta]$. Its logarithm and exponential series are given by
\[
  g(u)=-\frac{\ln(1-\beta u)}\beta,\quad f(x)=\frac{1-e^{-\beta
  x}}\beta.
\]
\end{example}

Let $F=\sum_{k,l}a_{kl}u^kv^l$ be a formal group law over a ring
$R$ and $r\colon R\to R'$ a ring homomorphism. Denote by $r(F)$
the formal series $\sum_{k,l}r(a_{kl})u^kv^l\in R'[[u,v]]$; then
$r(F)$ is a formal group law over~$R'$.

A formal group law $\mathcal F$ over a ring $A$ is
\emph{universal}\label{univfglaw} if for any formal group law $F$
over any ring $R$ there exists a unique homomorphism $r\colon A\to
R$ such that $F=r(\mathcal F)$.

\begin{proposition}
If a universal formal group law $\mathcal F$ over $A$ exists, then
\begin{itemize}
\item[(a)] the ring $A$ is multiplicatively generated by the coefficients of the series $\mathcal F$;

\item[(b)] $\mathcal F$ is unique: if
$\mathcal F'$ is another universal formal group law over $A'$, then
there is an isomorphism $r\colon A\to A'$ such that $\mathcal
F'=r(\mathcal F)$.
\end{itemize}
\end{proposition}
\begin{proof}
To prove the first statement, denote by $A'$ the subring in $A$
generated by the coefficients of~$\mathcal F$. Then there is a
monomorphism $i\colon A'\to A$ satisfying $i(\mathcal F)=\mathcal
F$. On the other hand, by universality there exists a homomorphism
$r\colon A\to A'$ satisfying $r(\mathcal F)=\mathcal F$. It follows
that $ir(\mathcal F)=\mathcal F$. This implies that
$ir=\mathrm{id}\colon A\to A$ by the uniqueness requirement in the
definition of~$\mathcal F$. Thus $A'=A$. The second statement is
proved similarly.
\end{proof}

\begin{theorem}[Lazard~\cite{laza55}]\label{lazardth}
The universal formal group law $\mathcal F$ exists, and its
coefficient ring $A$ is isomorphic to the polynomial ring
$\Z[a_1,a_2,\ldots]$ on an infinite number of generators.
\end{theorem}

\section{Formal group law of geometric cobordisms}\label{secfglgc}
The applications of formal group laws in cobordism theory build
upon the following fundamental construction.

\begin{construction}[Formal group law of geometric
cobordisms~\cite{novi67}]\label{cfglgc} Let $X$ be a cell complex
and $u,v\in U^2(X)$ two geometric cobordisms (see
Construction~\ref{geomcob}) corresponding to elements $x,y\in
H^2(X)$ respectively. Denote by $u+_{\!{}_H}\!v$ the geometric
cobordism corresponding to the cohomology class $x+y$.

\begin{proposition}\label{u2x}
The following relation holds in $U^2(X)$:
\begin{equation}\label{fglgc}
  u+_{\!{}_H}\!v=F_U(u,v)=u+v+\sum_{k\ge1,\,l\ge1}\alpha_{kl}\,u^kv^l,
\end{equation}
where the coefficients $\alpha_{kl}\in\varOmega_U^{-2(k+l-1)}$ do
not depend on~$X$. The series $F_U(u,v)$ given by~\eqref{fglgc} is
a formal group law over the complex cobordism ring~$\varOmega_U$.
\end{proposition}
\begin{proof}
We first do calculations with the universal example $X=\C
P^\infty\times\C P^\infty$. Then
\[
  U^*(\C P^\infty\times\C P^\infty)=\varOmega^*_U[[\underline u,\underline
  v]],
\]
where $\underline u,\underline v$ are canonical geometric
cobordisms given by the projections of $\C P^\infty\times\C
P^\infty$ onto its factors. We therefore have the following
relation in $U^2(\C P^\infty\times\C P^\infty)$:
\begin{equation}\label{univcalc}
  \underline u+_{\!{}_H}\!\underline v=
  \sum_{k,l\ge0} \alpha_{kl}\,\underline u^k\underline v^l,
\end{equation}
where $\alpha_{kl}\in\varOmega_U^{-2(k+l-1)}$.

Now let the geometric cobordisms $u,v\in U^2(X)$ be given by maps
$f_u,f_v\colon X\to\C P^\infty$ respectively. Then $u=(f_u\times
f_v)^*(\underline u)$, $v=(f_u\times f_v)^*(\underline v)$ and
$u+_{\!{}_H}\!v=(f_u\times f_v)^*(\underline
u+_{\!{}_H}\!\underline v)$, where $f_u\times f_v\colon X\to\C
P^\infty\times\C P^\infty$. Applying the $\varOmega_U$-module map
$(f_u\times f_v)^*$ to~\eqref{univcalc} we obtain the required
formula~\eqref{fglgc}. The identities (a) and (c) for $F_U(u,v)$
are obvious, and the associativity~(b) follows from the identity
$(u+_{\!{}_H}\!v)+_{\!{}_H}\!w=F_U(F_U(u,v),w)$ and the
associativity of~$+_{\!{}_H}$.
\end{proof}

Series~\eqref{fglgc} is called the \emph{formal group law of
geometric cobordisms}\label{fglawgc}; nowadays it is also usually
referred to as the \emph{complex cobordism formal group law}.

By definition, the geometric cobordism $u\in U^2(X)$ is the first
\emph{Conner--Floyd Chern class} $c_1^U(\xi)$ of the complex line
bundle $\xi$ over $X$ obtained by pulling back the canonical
bundle along the map $f_u\colon X\to\C P^\infty$ (it also
coincides with the Euler class $e(\xi)$ as defined in
Section~\ref{orcobbord}). It follows that the formal group law of
geometric cobordisms gives an expression of
$c_1^U(\xi\otimes\eta)\in U^2(X)$ in terms of the classes
$u=c_1^U(\xi)$ and $v=c_1^U(\eta)$ of the factors:
\[
  c_1^U(\xi\otimes\eta)=F_U(u,v).
\]
\end{construction}

The coefficients of the formal group law of geometric cobordisms
and its logarithms may be described geometrically by the following
results.

\begin{theorem}[{Buchstaber~\cite[Theorem~4.8]{buch70}}]\label{buchth}
\[
  F_U(u,v)=\frac{\sum_{i,j\ge0}[H_{ij}]u^iv^j}
  {\bigl(\sum_{r\ge0}[\C P^r]u^r\bigr)\bigl(\sum_{s\ge0}[\C
  P^s]v^s\bigr)},
\]
where $H_{ij}$ ($0\le i\le j$) are Milnor
hypersurfaces~\eqref{hij} and $H_{ji}=H_{ij}$.
\end{theorem}
\begin{proof}
Set $X=\C P^i\times\C P^j$ in Proposition~\ref{u2x}. Consider the
Poincar\'e--Atiyah duality map $D\colon U^2(\C P^i\times\C P^j)\to
U_{2(i+j)-2}(\C P^i\times\C P^j)$ (Construction~\ref{padual}) and
the map $\varepsilon\colon U_*(\C P^i\times\C P^j)\to
U_*(pt)=\varOmega^U$ induced by the projection $\C P^i\times\C
P^j\to pt$. Then the composition
\[
  \varepsilon D\colon U^2(\C P^i\times\C
  P^j)\to\varOmega_{2(i+j)-2}^U
\]
takes geometric cobordisms to the bordism classes of the
corresponding submanifolds. In particular, $\varepsilon
D(u+_{\!{}_H}\!v)=[H_{ij}]$, $\varepsilon D(u^kv^l)=[\C
P^{i-k}][\C P^{j-l}]$. Applying $\varepsilon D$ to~\eqref{fglgc}
we obtain
\[
  [H_{ij}]=\sum_{k,\,l}\alpha_{kl}[\C P^{i-k}][\C P^{j-l}].
\]
Therefore,
\[
  \sum_{i,j}[H_{ij}]u^iv^j=\Bigl(\sum_{k,\,l}\alpha_{kl}u^kv^l\Bigr)
  \Bigl(\sum_{i\ge k}[\C P^{i-k}]u^{i-k}\Bigr)
  \Bigl(\sum_{j\ge l}[\C P^{j-l}]v^{j-l}\Bigr),
\]
which implies the required formula.
\end{proof}

%\begin{corollary}\label{hijrel}
%There is a ring isomorphism
%\[
%  \varOmega_U=\Z\bigl[[H_{ij}]\colon0\le i\le j\bigr]\big/
%  \bigl(\text{\rm associativity relations for~$F_U(u,v)$}\bigr).
%\]
%\end{corollary}

\begin{corollary}\label{aijgen}
The coefficients of the formal group law of geometric cobordisms
generate the complex cobordism ring~$\varOmega_U$.
\end{corollary}
\begin{proof}
By Theorem~\ref{hijgen}, $\varOmega_U$ is generated by the
cobordism classes $[H_{ij}]$, which can be integrally expressed
via the coefficients of $F_U$ using Theorem~\ref{buchth}.
\end{proof}

\begin{theorem}[{Mishchenko~\cite[Appendix~1]{novi67}}]\label{mishth}
The logarithm of the formal group law of geometric cobordisms is
given by
\[
  g_U(u)=u+\sum_{k\ge1}\frac{[\C P^k]}{k+1}u^{k+1}
  \in\varOmega_U\otimes\Q[[u]].
\]
\end{theorem}
\begin{proof}
By~\eqref{log},
\[
  dg_U(u)=\frac{du}{\frac{\partial F_U(u,v)}{\partial
  v}\Bigl|_{v=0}}.
\]
Using the formula of Theorem~\ref{buchth} and the identity
$H_{i0}=\C P^{i-1}$, we calculate
\[
  \frac{dg_U(u)}{du}=\frac{1+\sum_{k>0}[\C P^k]u^k}
  {1+\sum_{i>0}([H_{i1}]-[\C P^1][\C P^{i-1}])u^i}.
\]
Now $[H_{i1}]=[\C P^1][\C P^{i-1}]$ (see Exercise~\ref{h1jcob}).
It follows that $\frac{dg_U(u)}{du}=1+\sum_{k>0}[\C P^k]u^k$,
which implies the required formula.
\end{proof}

Using these calculations the following most important property of
the formal group law $F_U$ can be easily established:

\begin{theorem}[{Quillen~\cite[Theorem~2]{quil69O}}]\label{quillenth}
The formal group law $F_U$ of geometric cobordisms is universal.
\end{theorem}
\begin{proof}
Let $\mathcal F$ be the universal formal group law over a
ring~$A$. Then there is a homomorphism $r\colon A\to\varOmega_U$
which takes $\mathcal F$ to $F_U$. The series $\mathcal F$, viewed
as a formal group law over the ring $A\otimes\Q$, has the
universality property for all formal group laws over
$\Q$-algebras. By Theorem~\ref{logth}, such a formal group law is
determined by its logarithm, which is a series with leading
term~$u$. It follows that if we write the logarithm of $\mathcal
F$ as $\sum b_k\frac{u^{k+1}}{k+1}$ then the ring $A\otimes\Q$ is
the polynomial ring $\Q[b_1,b_2,\ldots]$. By Theorem~\ref{mishth},
$r(b_k)=[\C P^k]\in\varOmega_U$. Since
$\varOmega_U\otimes\Q\cong\Q[[\C P^1],[\C P^2],\ldots]$, this
implies that $r\otimes\Q$ is an isomorphism.

By Theorem~\ref{lazardth} the ring $A$ does not have torsion, so
$r$ is a monomorphism. On the other hand, Theorem~\ref{buchth}
implies that the image $r(A)$ contains the bordism classes
$[H_{ij}]\in\varOmega_U$, $0\le i\le j$. Since these classes
generate the whole ring $\varOmega_U$ (Theorem~\ref{hijgen}), the
map $r$ is onto and thus an isomorphism.
\end{proof}

\begin{remark} Complex cobordism theory is the \emph{universal complex-oriented
cohomology theory}, i.e. it is universal in the category of
cohomology theories in which the first Chern class is defined for
complex bundles. This implies that the formal group law $F_U$ of
geometric cobordisms is universal among formal group laws realised
by complex-oriented cohomology theories. Theorem~\ref{quillenth}
establishes the universality of $F_U$ among all formal group laws.
It therefore opens a way for application of results from the
algebraic theory of formal groups in cobordism.
\end{remark}

\section{Hirzebruch genera}\label{aphige}
Every homomorphism $\varphi\colon\varOmega^U\to R$ from the
complex bordism ring to a commutative ring $R$ with unit can be
regarded as a multiplicative characteristic of manifolds which is
an invariant of bordism classes. Such a homomorphism is called a
(complex) \emph{$R$-genus}. (The term \emph{`multiplicative
genus'} is also used, to emphasise that such a genus is a ring
homomorphism.)
%in classical algebraicgeometry, there are instances of genera
%which are not multiplicative.)

Assume that the ring $R$ does not have additive torsion. Then
every $R$-genus $\varphi$ is fully determined by the corresponding
homomorphism $\varOmega^U\otimes\Q\to R\otimes\Q$, which we shall
also denote by~$\varphi$. A construction due to
Hirzebruch~\cite{hirz66} describes homomorphisms
$\varphi\colon\varOmega^U\otimes\Q\to R\otimes\Q$ by means of
universal $R$-valued characteristic classes of special type.

\begin{construction}[Hirzebruch genera]\label{hirzcor}
Let $BU=\lim\limits_{n\to\infty}BU(n)$. Then $H^*(BU)$ is
isomorphic to the ring of formal power series
$\Z[[c_1,c_2,\ldots]]$ in universal Chern classes, $\deg c_k=2k$.
The set of tangential Chern characteristic numbers of a given
manifold $M$ defines an element in $\Hom(H^*(BU),\Z)$, which
belongs to the subgroup $H_*(BU)$ in the latter group. We
therefore obtain a group homomorphism
\begin{equation}\label{oubu}
  \varOmega^U\to H_*(BU).
\end{equation}
Since the product in $H_*(BU)$ arises from the maps $BU(k)\times
BU(l)\to BU(k+l)$ corresponding to the product of vector bundles,
and the Chern classes have the appropriate multiplicative
property, the map~\eqref{oubu} is a ring homomorphism.

Part~2 of Theorem~\ref{comcob} says that~\eqref{oubu} is a
monomorphism, and Part~1 of the same theorem says that the map
$\varOmega^U\otimes\Q\to H_*(BU;\Q)$ is an isomorphism. It follows
that every homomorphism $\varphi\colon\varOmega^U\otimes\Q\to
R\otimes\Q$ can be interpreted as an element of
\[
  \Hom_{\Q}(H_*(BU;\Q),R\otimes\Q)=H^*(BU;\Q)\otimes R,
\]
or as a sequence of homogeneous polynomials
$\{K_i(c_1,\ldots,c_i),\;i\ge0\}$, $\deg K_i=2i$. This sequence of
polynomials cannot be chosen arbitrarily; the fact that $\varphi$
is a ring homomorphism imposes certain conditions. These
conditions may be described as follows: an identity
\[
  1+c_1+c_2+\cdots=(1+c'_1+c'_2+\cdots)\cdot(1+c''_1+c''_2+\cdots)
\]
implies the identity
\begin{equation}\label{mulseqeq}
  \sum_{n\ge0}K_n(c_1,\ldots,c_n)=
  \sum_{i\ge0}K_i(c'_1,\ldots,c'_i)\cdot
  \sum_{j\ge0}K_j(c''_1,\ldots,c''_j).
\end{equation}
A sequence of homogeneous polynomials $\mathcal
K=\{K_i(c_1,\ldots,c_i),i\ge0\}$ with $K_0=1$ satisfying the
identities~\eqref{mulseqeq} is called a \emph{multiplicative
Hirzebruch sequence}.

\begin{proposition}
A multiplicative sequence $\mathcal K$ is completely determined by
the series
\[
  Q(x)=1+q_1x+q_2x^2+\cdots\in R\otimes\Q[[x]],
\]
where $x=c_1$, and $q_i=K_i(1,0,\ldots,0)$; moreover, every series
$Q(x)$ as above determines a multiplicative sequence.
\end{proposition}
\begin{proof}
Indeed, by considering the identity
\begin{equation}\label{formalchern}
  1+c_1+\cdots+c_n=(1+x_1)\cdots(1+x_n)
\end{equation}
we obtain from \eqref{mulseqeq} that
\begin{multline*}
  Q(x_1)\cdots Q(x_n)=1+K_1(c_1)+K_2(c_1,c_2)+\cdots\\+
  K_n(c_1,\ldots,c_n)+K_{n+1}(c_1,\ldots,c_n,0)+\cdots.\qedhere
\end{multline*}
\end{proof}

Along with the series $Q(x)$ it is convenient to consider the
series $f(x)\in R\otimes\Q[[x]]$ with leading term $x$ given by
the identity
\[
  Q(x)=\frac x{f(x)}.
\]
It follows that the $n$th term $K_n(c_1,\ldots,c_n)$ in the
multiplicative Hirzebruch sequence corresponding to a genus
$\varphi\colon\varOmega^U\otimes\Q\to R\otimes\Q$ is the
degree-$2n$ part of the series $\prod^n_{i=1}\frac{x_i}{f(x_i)}\in
R\otimes\Q[[c_1,\ldots,c_n]]$. By some abuse of notation, we
regard $\prod^n_{i=1}\frac{x_i}{f(x_i)}$ as a characteristic class
of complex $n$-plane bundles.
%(although it is a series in $c_1,\ldots,c_n$, rather than a polynomial).
Then the value of $\varphi$ on an $2n$-dimensional stably complex
manifold $M$ is given by
\begin{equation}\label{phim}
  \varphi[M]=\Bigl(\prod^n_{i=1}\frac{x_i}{f(x_i)}(\mathcal TM)\Bigr)
  \langle M\rangle,
\end{equation}
where the Chern classes $c_1,\ldots,c_n$ are expressed via the
indeterminates $x_1,\ldots,x_n$ by the
relation~\eqref{formalchern}, and the degree-$2n$ part of
$\prod^n_{i=1}\frac{x_i}{f(x_i)}$ is taken to obtain a
characteristic class of~$\mathcal TM$. We shall also denote the
`characteristic class' $\prod^n_{i=1}\frac{x_i}{f(x_i)}$ of a
complex $n$-plane bundle $\xi$ by $\varphi(\xi)$; so that
$\varphi[M]=\varphi({\mathcal T}\!M)\langle M\rangle$.

We refer to the homomorphism $\varphi\colon\varOmega^U\to
R\otimes\Q$ given by~\eqref{phim} as the \emph{Hirzebruch genus}
corresponding to the series $f(x)=x+\cdots\in R\otimes\Q[[x]]$.
\end{construction}

%\begin{theorem}\label{hirztheo}
%Let $R$ be a torsion-free ring. Any $R$-genus
%$\varphi\colon\varOmega^U\to R$ (or $R\otimes\Q$-genus
%$\varphi\colon\varOmega^U\to R\otimes\Q$) defines a series
%$f(x)=x+\cdots\in R\otimes\Q[[x]]$, and any such series defines an
%$R\otimes\Q$-genus.
%\end{theorem}
%
%The characteristic class
%$\varphi=\prod^n_{i=1}\frac{x_i}{f(x_i)}\in H^*(BU(n);R\otimes\Q)$
%of $n$-dimensional complex bundles is called the \emph{Hirzebruch
%genus} corresponding to a series $f(x)=x+\cdots\in
%R\otimes\Q[[x]]$. The above construction can be summarised as
%follows:

\begin{remark}
A parallel theory of genera exists for oriented manifolds. These
genera are homomorphisms $\varOmega^{SO}\to R$ from the oriented
bordism ring, and the Hirzebruch construction expresses genera
over torsion-free rings via Pontryagin characteristic classes
(which replace the Chern classes).
\end{remark}

Given a torsion-free ring $R$ and a series $f(x)\in
R\otimes\Q[[x]]$ one may ask whether the corresponding Hirzebruch
genus $\varphi\colon\varOmega^U\to R\otimes\Q$ actually takes
values in~$R$. This constitutes the \emph{integrality problem}
for~$\varphi$. A solution to this problem can be obtained using
the theory of formal group laws.

Every genus $\varphi\colon\varOmega^U\to R$ gives rise to a formal
group law $\varphi(F_U)$ over~$R$, where $F_U$ is the formal group
of geometric cobordisms (Construction~\ref{cfglgc}).

\begin{theorem}[Novikov~\cite{novi68}]\label{fexpser}
For every genus $\varphi\colon\varOmega^U\to R$, the exponential
of the formal group law $\varphi(F_U)$ is the series $f(x)\in
R\otimes\Q[[x]]$ corresponding to~$\varphi$.
\end{theorem}
This can be proved either directly, by appealing to the
construction of geometric cobordisms, or indirectly, by
calculating the values of~$\varphi$ on projective spaces and
comparing to the formula for the logarithm of the formal group
law.
\begin{proof}[1st proof.] Let $X$ be a stably complex $d$-manifold
with $x,y\in H^2(X)$, and let $u,v\in U^2(X)$ be the
corresponding geometric cobordisms (see
Construction~\ref{geomcob}) represented by codimension-2
submanifolds $M_x\subset X$ and $M_y\subset X$  respectively. Then
$u^kv^l\in U^{2(k+l)}(X)$ is represented by a submanifold of
codimension $2(k+l)$, which we denote by~$M_{kl}$.
By~\eqref{fglgc}, we have the following relation in $U^2(X)$:
\[
  [M_{x+y}]=\sum_{k,l\ge0}\alpha_{kl}[M_{kl}]
\]
where $M_{x+y}\subset X$ is the codimension-2 submanifold dual to
$x+y\in H^2(X)$. We apply the composition $\varepsilon D$ of the
Poincar\'e--Atiyah duality map $D\colon U^2(X)\to U_{d-2}(X)$ and
the augmentation $U_{d-2}(X)\to\varOmega^U_{d-2}$ to the identity
above, and then apply the genus $\varphi$ to the resulting idenity
in~$\varOmega^U_{d-2}$ to obtain
\begin{equation}\label{mxy}
  \varphi[M_{x+y}]=\sum\varphi(\alpha_{kl})\varphi[M_{kl}].
\end{equation}
Let $\iota\colon M_{x+y}\subset X$ be the embedding. Considering
the decomposition
\[
  \iota^*({\mathcal T} X)={\mathcal T}\!M_{x+y}\oplus\nu(\iota)
\]
and using the multiplicativity of the characteristic class
$\varphi$ we obtain
\[
  \iota^*\varphi({\mathcal T} X)=\varphi({\mathcal T}
  M_{x+y})\cdot\iota^*({\textstyle\frac{x+y}{f(x+y)}}).
\]
Therefore,
\begin{equation}\label{f1}
  \varphi[M_{x+y}]=\iota^*\bigl(\varphi({\mathcal T}
  X)\cdot{\textstyle\frac{f(x+y)}{x+y}}\bigr)\langle M_{x+y}\rangle=
  \bigl(\varphi({\mathcal T} X)\cdot f(x+y)\bigr)\langle X\rangle.
\end{equation}
Similarly, by considering the embedding $M_{kl}\to X$ we obtain
\begin{equation}\label{f2}
  \varphi[M_{kl}]=\bigl(\varphi({\mathcal T} X)\cdot f(x)^kf(y)^l\bigr)\langle X\rangle.
\end{equation}
Plugging \eqref{f1} and \eqref{f2} into \eqref{mxy} we finally
obtain
\[
  f(x+y)=\sum_{k,l\ge0}\varphi(\alpha_{kl})f(x)^k f(y)^l.
\]
This implies, by definition, that $f$ is the exponential of
$\varphi(F_U)$.
\end{proof}

\begin{proof}[2nd proof.] The complex bundle isomorphism $\mathcal T(\C
P^k)\oplus\underline{\C}=\bar\eta\oplus\cdots\oplus\bar\eta$
($k+1$ summands) allows us to calculate the value of a genus on
$\C P^k$ explicitly. Let $x=c_1(\bar\eta)\in H^2(\C P^k)$ and let
$g$ be the series functionally inverse to $f$; then
\begin{align*}
\varphi[\C P^k]&=\Bigl(\frac
x{f(x)}\Bigr)^{k+1}\langle\C P^k\rangle\\
&=\text{coefficient of $x^k$ in }\Bigl(\frac
x{f(x)}\Bigr)^{k+1}=\mathop{\mathrm{res}_0}\Bigl(\frac
1{f(x)}\Bigr)^{k+1}\\
&=\frac1{2\pi i}\oint\Bigl(\frac 1{f(x)}\Bigr)^{k+1}dx=
\frac1{2\pi i}\oint\frac 1{u^{k+1}}g'(u)du\\
&=\mathop{\mathrm{res}_0}\Bigl(\frac{g'(u)}{u^{k+1}}\Bigr)=
\text{coefficient of $u^k$ in $g'(u)$.}
\end{align*}
(Integrating over a closed path around zero makes sense only for
convergent power series with coefficients in $\C$, however the
result holds for all power series with coefficients in
$R\otimes\Q$.) Therefore,
\[
  g'(u)=\sum_{k\ge0}\varphi[\C P^k]u^k.
\]
Theorem~\ref{mishth} then implies that $g$ is the logarithm of the
formal group law $\varphi(F_U)$, and thus $f$ is its exponential.
\end{proof}

Now we can formulate the following criterion for the integrality
of a genus:

\begin{theorem}\label{intgenfgl}
Let $R$ be a torsion-free ring and $f(x)\in R\otimes\Q[[x]]$. The
corresponding genus $\varphi\colon\varOmega^U\to R\otimes\Q$ takes
values in~$R$ if and only if the coefficients of the formal group
law $\varphi(F_U)(u,v)=f(f^{-1}(u)+f^{-1}(v))$ belong to~$R$.
\end{theorem}
\begin{proof}
This follows from Theorem~\ref{fexpser} and
Corollary~\ref{aijgen}.
\end{proof}

\begin{example}\label{unigenus}
The \emph{universal genus} maps a stably complex manifold $M$ to
its bordism class $[M]\in\varOmega^U$ and therefore corresponds to
the identity homomorphism
$\varphi_U\colon\varOmega^U\to\varOmega^U$. The corresponding
characteristic class with coefficients in $\varOmega^U\otimes\Q$,
the \emph{universal Hirzebruch genus}, was studied
in~\cite{buch70}; its corresponding series $f_U(x)$ is the
exponential of the universal formal group law of geometric
cobordisms.
\end{example}

\begin{example}\label{xyexa} We take $R=\Z$ in these examples.

1. The top Chern number $c_n[M]$ is a Hirzebruch genus, and its
corresponding $f$-series is $f(x)=\frac x{1+x}$. The value of this
genus on a stably complex manifold $(M,c_{\mathcal T})$ equals the
Euler characteristic of $M$ if $c_{\mathcal T}$ is an
\emph{almost} complex structure.

2. The \emph{$L$-genus} $L[M]$ corresponds to the series
$f(x)=\mathop{\mathrm{tanh}}(x)$ (the hyperbolic tangent). It is
equal to the \emph{signature} $\sign(M)$ by the classical
Hirzebruch formula~\cite{hirz66}.

3. The \emph{Todd genus} $\td[M]$ corresponds to the series
$f(x)=1-e^{-x}$. The corresponding formal group law is given by
$F(u,v)=u+v-uv$, compare Example~\eqref{mfgleq}. By
Theorem~\ref{intgenfgl}, the Todd genus is integral on any complex
bordism class. Furthermore, it takes value 1 on every complex
projective space~$\C P^k$.
\end{example}

The `trivial' genus $\varepsilon\colon\varOmega_U\to\Z$
corresponding to the series $f(x)=x$ gives rise to the
\emph{augmentation transformation} $U^*\to H^*$ from complex
cobordism to ordinary cohomology (also known as the \emph{Thom
homomorphism}). More generally, for every genus
$\varphi\colon\varOmega_U\to R$ and a space $X$ we may set
$h^*_\varphi(X)=U^*(X)\otimes_{\varOmega_U}R$. Under certain
conditions guaranteeing the exactness of the sequences of pairs
(known as the \emph{Landweber Exact Functor
Theorem}~\cite{land76})\label{Landefth} the functor
$h^*_\varphi(\cdot)$ gives rise to a complex-oriented cohomology
theory with the coefficient ring~$R$.

As an example of this procedure, consider the genus corresponding
to the 1-parameter extension of the multiplicative formal group
law, see Example~\ref{multfgl}. It is also usually called the
\emph{Todd genus}, and takes values in the ring $\Z[\beta]$,
$\deg\beta=-2$. By interpreting $\beta=1-\bar\eta$ as the
\emph{Bott element} in the complex $K$-group $\widetilde
K^0(S^2)=K^{-2}(pt)$ we obtain a homomorphism $\td\colon
\varOmega_U\to K^*(pt)$. It gives rise to a multiplicative
transformation $U^*\to K^*$ from complex cobordism to complex
$K$-theory introduced by Conner and Floyd~\cite{co-fl66}. In this
paper Conner and Floyd proved that complex cobordism determines
complex $K$-theory by means of the isomorphism
$K^*(X)=U^*(X)\otimes_{\varOmega_U}\mathbb Z[\beta]$, where the
$\varOmega_U$-module structure on $\mathbb Z[\beta]$ is given by
the Todd genus. Their proof makes use of the Conner--Floyd Chern
classes; several proofs were given subsequently, including one
which follows directly from the Landweber exact functor theorem.

\begin{example}\label{chigenu}
Another important example from the original work of Hirzebruch is
given by the \emph{$\chi_y$-genus}. It corresponds to the series
\[
  f(x)=\frac{1-e^{-x(1+y)}}{1+ye^{-x(1+y)}},
\]
where $y\in\R$ is a parameter. Setting $y=-1$, $y=0$ and $y=1$ we
get $c_n[M]$, the Todd genus $\td[M]$ and the $L$-genus
$L[M]=\sign(M)$ respectively.

If $M$ is a complex manifold then the value $\chi_y[M]$ can be
calculated in terms of the Euler characteristics of
\emph{Dolbeault complexes} on~$M$, see~\cite{hirz66}.
%This was the
%original Hirzebruch's motivation to study the $\chi_y$-genus.
\end{example}

\subsection*{Exercises.}
%\begin{exercise}\label{h1i}
%Show that $[H_{1i}]=[\C P^1][\C P^{i-1}]$ (as complex bordism
%classes) for $i\ge 1$ by calculating the Chern numbers. (The
%authors do not know of an explicit construction realising this
%cobordism.)
%\end{exercise}

\begin{exercise}
The formal group law corresponding to the $\chi_y$-genus is given
by
\[
  F(u,v)=\frac{u+v+(y-1)uv}{1+y\cdot uv}.
\]
\end{exercise}
\end{appendix}

\begin{theindex}
\item Action (of group), \pageref{defiaction}
  \subitem almost free, \pageref{defftra}, \pageref{defifixedp}
  \subitem effective, \pageref{defifixedp}
  \subitem free, \pageref{defifixedp}
  \subitem proper, \pageref{Gprop}, \pageref{propergad}
  \subitem semifree, \pageref{semifreeHam}, \pageref{defifixedp}
  \subitem transitive, \pageref{defifixedp}
\item Acyclic (space), \pageref{acysp}, \pageref{acyclspac}
\item Affine equivalence, \pageref{combeq}
\item Alexander duality, \pageref{simaldual}, \pageref{alexadua}
\item Algebra, \pageref{algmod}
  \subitem bigraded, \pageref{bigradedalge}
  \subitem connected, \pageref{gradedalgebra}
  \subitem exterior, \pageref{freegco}
  \subitem finitely generated, \pageref{algmod}
  \subitem graded, \pageref{gradedalgebra}
  \subitem graded-commutative, \pageref{gradedalgebra}
    \subsubitem free, \pageref{freegco}
  \subitem polynomial, \pageref{algmod}
  \subitem multigraded, \pageref{bigradedalge}
  \subitem with straightening law (ASL), \pageref{straighrel}
\item Almost complex structure, \pageref{tscmani}
   \subitem associated with omniorientation, \pageref{omniacs}
   \subitem $T$-invariant, \pageref{omniacs},
   \pageref{kusth}, \pageref{iscs}
\item Ample divisor, \pageref{amplediviso}
\item Annihilator (of module), \pageref{dimmodule}
\item Arrangement (of subspaces), \pageref{arrangem}
  \subitem coordinate, \pageref{casim}
\item Aspherical (space), \pageref{aspheri}, \pageref{asphHm}
\item Associahedron, \pageref{defassociah}
  \subitem generalised, \pageref{stasheptpe}
\item Atiyah--Hirzebruch formula, \pageref{tablocal}
\item Axial function, \pageref{torusaf}, \pageref{defiaxif}
  \subitem $n$-independent, \pageref{defiTgraph}
  \subitem 2-independent, \pageref{defiaxif}

\indexspace

\item Bar construction, \pageref{barcobar}
\item Barnette sphere, \pageref{barsph}, \pageref{barnstar},
\pageref{barnesph}
\item Barycentre (of simplex), \pageref{defbaryce}
\item Barycentric subdivision
  \subitem of polytope, \pageref{bcptope}
  \subitem of simplicial complex, \pageref{barsub}
  \subitem of simplicial poset, \pageref{baryc}, \pageref{prop:barst}
\item Base (of Schlegel diagram), \pageref{basesch}
\item Based map, \pageref{pointedspace}
\item Based space, \pageref{pointedspace}
\item Basepoint, \pageref{pointedspace}
\item Basis (of free module), \pageref{freemodule}
  \subitem minimal, \pageref{minimal}
\item Betti numbers (algebraic), \pageref{toralg},
\pageref{aBettisp}, \pageref{multibnsp}
\item Betti numbers (topological), \pageref{bettinumb}
  \subitem bigraded, \pageref{defbn}
\item Bidegree, \pageref{bigradedalge}
\item Bistellar equivalence, \pageref{bistequ}
\item Bistellar move, \pageref{bist}
\item Blow-down, \pageref{deblowdown}
\item Blow-up
  \subitem of $T$-graph, \pageref{blowupTgr}
  \subitem of $T$-manifold, \pageref{dblowup}
\item Boolean lattice, \pageref{defsp}
\item Bordism, \pageref{bordismman}
  \subitem complex, \pageref{defncomplbor}
  \subitem oriented, \pageref{defnorbor}
  \subitem unoriented, \pageref{unoriebordi}
\item Borel construction, \pageref{2projections}
\item Borel spectral sequence, \pageref{borelss}
\item Bott manifold, \pageref{defbottmani}
  \subitem generalised, \pageref{genbott}
\item Bott--Taubes polytope, \pageref{defcyclohe}
\item Bott tower, \pageref{defbotttower}
  \subitem generalised, \pageref{genbott}
  \subitem real, \pageref{bottsmalcov}
  \subitem topologically trivial, \pageref{btowetri}
\item Boundary, \pageref{chaincomplex}
\item Bounded flag, \pageref{bfm}
\item Bounded flag manifold, \pageref{bfm}, \pageref{btbfn},
\pageref{boundedfmcob}
\item Bouquet (of spaces), \pageref{pointedspace}
\item Br\"uckner sphere, \pageref{brueckner}
\item Buchsbaum complex, \pageref{buchsbaumcplx}
\item Buchstaber invariant, \pageref{buchinvari}
\item Building set, \pageref{defbs}
  \subitem graphical, \pageref{graphass}

\indexspace

\item Calabi--Eckmann manifold, \pageref{calabieck}
\item Canonical line bundle, \pageref{canolinebu}
\item Catalan number, \pageref{defassociah}
\item Categorical quotient, \pageref{catequot}
\item CAT(0) inequality, \pageref{Hopfconj}
\item Cell complex, \pageref{cellcomplex}
\item Cellular approximation, \pageref{celluappr}
\item Cellular chain, \pageref{celluchain}
\item Cellular map, \pageref{cellumap}
\item Chain complex, \pageref{chaincomplex}
  \subitem augmented, \pageref{augmcs}, \pageref{augme3}
\item Chain homotopy, \pageref{cochainhomotopy}
\item Characteristic function, \pageref{charf}
  \subitem directed, \pageref{dircharf}
\item Characteristic matrix, \pageref{charmat}
  \subitem refined, \pageref{lamat}
\item Characteristic number, \pageref{charnumbersd}
\item Characteristic pair, \pageref{charpair}
  \subitem combinatorial, \pageref{cqtp}
  \subitem equivalence of, \pageref{psieq}, \pageref{cqtp}
\item Characteristic submanifold, \pageref{charf},
\pageref{charm}, \pageref{charmttm}
\item Charney--Davis Conjecture, \pageref{defflp}
\item Chern--Dold character, \pageref{cherndold}
\item Chow ring, \pageref{cohtm}
\item CHP (covering homotopy property), \pageref{chpdiag}
\item Classifying functor (algebraic), \pageref{barcobar}
\item Classifying space, \pageref{classispa}
\item Clique, \pageref{defcliq}, \pageref{dclique}
\item Cobar construction, \pageref{adamsSR}, \pageref{barcobar}
\item Cobordism, \pageref{etaxicob}, \pageref{ocobgr}
  \subitem complex, \pageref{defccobord}
  \subitem equivariant
    \subsubitem geometric, \pageref{PhiX}, \pageref{geomequcobor}
    \subsubitem homotopic, \pageref{PhiX}, \pageref{homocobor}
  \subitem unoriented, \pageref{unorcobo}
\item Coboundary, \pageref{cochaincomplex}
\item Cochain complex, \pageref{cochaincomplex}
\item Cochain homotopy, \pageref{cochainhomotopy}
\item Cocycle, \pageref{cochaincomplex}
\item Cofibrant (object, replacement, approximation), \pageref{cofibobject}
\item Cofibration, \pageref{decofibration}
  \subitem in model category, \pageref{overundercategory}
\item Cofibre, \pageref{decofibration}
\item Coformality, \pageref{djcoform}, \pageref{defncoforma},
\pageref{coformLie}
\item Cohen--Macaulay
  \subitem algebra, \pageref{CM}
  \subitem module, \pageref{CM}
  \subitem simplicial complex, \pageref{CMsim}
  \subitem simplicial poset, \pageref{CMsimpose}
\item Cohomological rigidity, \pageref{defcohorig}
\item Cohomology
  \subitem cellular, \pageref{celluhomol}
  \subitem of cochain complex, \pageref{cochaincomplex}
  \subitem simplicial, \pageref{simpcohomo}
    \subsubitem reduced, \pageref{reduscoho}
  \subitem singular, \pageref{singucochain}
\item Cohomology product, \pageref{diagomap}
\item Cohomology product length, \pageref{cohoprole}
\item Colimit, \pageref{colimitmac}, \pageref{colimittv},
\pageref{facecate}, \pageref{dfinlimcolim}
\item Combinatorial equivalence
  \subitem of polyhedral complexes, \pageref{schldiag}
  \subitem of polytopes, \pageref{combeq}
\item Combinatorial neighbourhood, \pageref{combnei}
\item Complex orientable map, \pageref{compormap}
\item Complete intersection algebra, \pageref{frcomplint}
\item Cone
  \subitem over simplicial complex, \pageref{consusp}
  \subitem over space, \pageref{conesuspspace}
\item Cone (convex polyhedral), \pageref{cpcone}
  \subitem dual, \pageref{dualcone}
  \subitem rational, \pageref{cpcone}
  \subitem regular, \pageref{cpcone}
  \subitem simplicial, \pageref{cpcone}
  \subitem strongly convex, \pageref{cpcone}
\item Connected sum
  \subitem of manifolds, \pageref{connectesu}, \pageref{cobcs}
    \subsubitem equivariant, \pageref{equcs}
  \subitem of simple polytopes, \pageref{consum}
  \subitem of simplicial complexes, \pageref{simcs}
  \subitem of stably complex manifolds, \pageref{consumstcom}
\item Connection (on graph), \pageref{defconn}
\item Contraction (of building set), \pageref{contracti}
\item Coordinate subspace, \pageref{li}
\item Coproduct, \pageref{depushout}
\item Core (of simplicial complex), \pageref{defcore}
\item Cross-polytope, \pageref{cross}
\item Cube, \pageref{cube}
  \subitem standard, \pageref{cube}
  \subitem topological, \pageref{topcubcom}
\item Cubical complex
  \subitem abstract, \pageref{cubcom}
  \subitem polyhedral, \pageref{polyhedralcc}
  \subitem topological, \pageref{topcubcom}
\item Cubical subdivision, \pageref{cubicsubd}
\item Cup product, \pageref{diagomap}
\item CW complex, \pageref{cellcomplex}
\item Cycle, \pageref{chaincomplex}
\item Cyclohedron, \pageref{defcyclohe}

\indexspace

\item Davis--Januszkiewicz space, \pageref{ddjs}, \pageref{djsformalsp}
\item Deformation retraction, \pageref{deret}
\item Degree (grading), \pageref{gradedalgebra}
  \subitem external, \pageref{extdegree}
  \subitem internal, \pageref{extdegree}
  \subitem total, \pageref{extdegree}
\item Dehn--Sommerville relations
  \subitem for polytopes, \pageref{ds}, \pageref{DSgenli}
  \subitem for simplicial complexes, \pageref{dsgor}
  \subitem for simplicial posets, \pageref{theo:DS}
  \subitem for triangulated manifolds, \pageref{dstm}
\item Depth (of module), \pageref{depth}
\item Diagonal approximation, \pageref{diagappro}
\item Diagonal map, \pageref{diagomap}
\item Diagram (functor), \pageref{dediagram}
\item Diagram category, \pageref{dediagram}
\item Differential graded algebra (dg-algebra, dga),
\pageref{leibnitz}, \pageref{dgaap}
  \subitem formal, \pageref{formalg}
  \subitem homologically connected, \pageref{homoconne}
  \subitem minimal, \pageref{midga}
  \subitem simply connected, \pageref{simplycondga}
\item Dimension
  \subitem of module, \pageref{dimmodule}
  \subitem of polytope, \pageref{suphp}
  \subitem of simplicial complex, \pageref{defsubcomplex}
  \subitem of simplicial poset, \pageref{defsp}
\item Dolbeault cohomology, \pageref{dolbcohom}
\item Dolbeault complex, \pageref{dolbcohom}, \pageref{chigenu}
\item Double (simplicial), \pageref{simplicialwedge}

\indexspace

\item Edge, \pageref{facet}
\item Eilenberg--Mac\,Lane space, \pageref{Eilespac}
\item Eilenberg--Moore spectral sequence, \pageref{Eilemooress}
\item Elliptic genus, \pageref{ellitpigenus}
\item Equivariant bundle, \pageref{equbundle}
\item Equivariant characteristic class, \pageref{equbundle}
\item Equivariant cohomology, \pageref{equcohomolde}
  \subitem of $T$-graph, \pageref{eqcohTgr}
\item Equivariant map, \pageref{equivmap}
\item Euler class
  \subitem in cobordism, \pageref{eulerclcob}
  \subitem in equivariant cobordism, \pageref{eulerclassec}
  \subitem in equivariant cohomology, \pageref{gysin1}, \pageref{eulereqcoh}
\item Euler formula, \pageref{euler}
\item Euler characteristic, \pageref{eulerqtori}, \pageref{eulercha},
\pageref{Eulerchar}, \pageref{eulechara}, \pageref{xyexa}
\item Eulerian poset, \pageref{theo:DS}
\item Exact sequence, \pageref{exaseque}
\item Excision, \pageref{excisionthm}
\item Exponential (of formal group law), \pageref{logfglaw}

\indexspace

\item Face
  \subitem of cubical complex, \pageref{cubcom}
  \subitem of cone, \pageref{faceofcone}
  \subitem of manifold with corners, \pageref{facemwco}, \pageref{facesmwc}
  \subitem of manifold with faces, \pageref{mawfa}
  \subitem of polytope, \pageref{facep}
  \subitem of simplicial complex, \pageref{faceofsc}
  \subitem of simplicial poset, \pageref{faceofsp}
  \subitem of $T$-graph, \pageref{faceTgr}
\item Face category (of simplicial complex),
\pageref{facecategor}, \pageref{facecate}, \pageref{posetcateg}
\item Face coalgebra, \pageref{facecoalgebra}
\item Face poset, \pageref{faceposet}
\item Face ring (Stanley--Reisner ring)
  \subitem of manifold with corners, \pageref{qfrgen}
  \subitem of simple polytope, \pageref{frptope}
  \subitem of simplicial complex, \pageref{frsim}
    \subsubitem exterior, \pageref{extefr}
  \subitem of simplicial poset, \pageref{frspo}
\item Facet, \pageref{facet}, \pageref{mawfa}, \pageref{faceTgr}
\item Face truncation, \pageref{hypcut}
\item Fan, \pageref{deffan}, \pageref{defps}
  \subitem complete, \pageref{deffan}
  \subitem normal, \pageref{nf}
  \subitem rational, \pageref{deffan}
  \subitem regular, \pageref{deffan}
  \subitem simplicial, \pageref{deffan}
\item Fat wedge, \pageref{dfatwedge}
\item Fibrant (object, replacement, approximation), \pageref{cofibobject}
\item Fibration, \pageref{defibrat}
  \subitem in model category, \pageref{modelcategory}
  \subitem locally trivial, \pageref{ltfibra}
\item Fibre bundle, \pageref{ltfibra}
  \subitem associated (with $G$-space), \pageref{associbundle}
\item Fixed point, \pageref{defifixedp}
\item Fixed point data, \pageref{fpdata}
\item Flag complex (simplicial), \pageref{defflagc}, \pageref{homotflag}
\item Flagtope (flag polytope), \pageref{defflp}
\item Folding map, \pageref{fold}
\item Formal group law, \pageref{genera}
  \subitem of geometric cobordisms, \pageref{fglawgc}
  \subitem linearisable, \pageref{logprop}
  \subitem universal, \pageref{univfglaw}
\item Formality
  \subitem in model category, \pageref{defncoforma}
  \subitem integral, \pageref{intformality}
  \subitem of dg-algebra, \pageref{formalg}
  \subitem of space
  \pageref{Masseymac}, \pageref{xkform}, \pageref{formtvman},
  \pageref{formsp}
\item $F$-polynomial (of polytope), \pageref{fhvect}
\item Fr\"olicher spectral sequence, \pageref{borelss}
\item Full subcomplex (of simplicial complex), \pageref{defcore}
\item Fundamental group, \pageref{homotopygroup}
\item Fundamental homology class, \pageref{fundpseudo}
\item $f$-vector (face vector)
  \subitem of cubical complex, \pageref{fvectorcc}
  \subitem of polytope, \pageref{fhvect}
  \subitem of simplicial complex, \pageref{simfve}
  \subitem of simplicial poset, \pageref{fvsimpose}

\indexspace

\item Gal Conjecture, \pageref{galconj}
\item Gale diagram, \pageref{galediagra}, \pageref{ex3quadr},
\pageref{3quad}
  \subitem combinatorial, \pageref{combgalediag}, \pageref{malink}
\item Gale duality, \pageref{galeduality}, \pageref{polquad},
\pageref{propmmap}, \pageref{galeduadef}, \pageref{mumatrix}
\item Gale transform, \pageref{galeduality}
\item Ganea's Theorem, \pageref{GaneaT}
\item $g$-conjecture, \pageref{gconj}, \pageref{gconje}
\item Generalised homology theory, \pageref{genhomotheory}
\item Generalised Lower Bound Conjecture, \pageref{glbt},
\pageref{glbctor}
\item Generating series
  \subitem of face polynomials, \pageref{genesefp}
  \subitem of polytopes, \pageref{genserie}
\item Genus, \pageref{degenus}, \pageref{aphige}
\item Geometric cobordism, \pageref{geomcob}
\item Geometric quotient, \pageref{geomequot}
\item Geometric realisation
  \subitem of cubical complex, \pageref{geomrelcc}
  \subitem of simplicial complex, \pageref{defsubcomplex}
  \subitem of simplicial set, \pageref{geomrealss}
\item Ghost vertex, \pageref{ghostvertex}
\item GKM-graph, \pageref{dGKMmani}
\item GKM-manifold, \pageref{dGKMmani}
\item Golod
  \subitem ring, \pageref{defgolod}, \pageref{Golodrin}
  \subitem simplicial complex, \pageref{defgolod}
\item Gorenstein, Gorenstein*
  \subitem simplicial complex, \pageref{goren}, \pageref{pdcom}
  \subitem simplicial poset, \pageref{goresp}
\item Goresky--MacPherson formula, \pageref{gomaformu}
\item Graph
  \subitem chordal, \pageref{dchordal}
  \subitem of polytope, \pageref{graphptope}
  \subitem simple, \pageref{graphptope}, \pageref{graphass}
\item Graph-associahedron, \pageref{graphass}
\item Graph product, \pageref{graphprodu}
\item $g$-theorem, \pageref{gth}, \pageref{proofgth}
\item $g$-vector
  \subitem of polytope, \pageref{fhvect}
  \subitem of simplicial complex, \pageref{simfve}
\item Gysin homomorphism
  \subitem in equivariant cohomology, \pageref{gysin1}, \pageref{defingysin}
  \subitem in cobordism, \pageref{utggysin}, \pageref{gysincob}
\item Gysin--Thom isomorphism, \pageref{fpCDvb2}, \pageref{gysin-thom}

\indexspace

\item Half-smash product (left, right), \pageref{halfsmash}
\item Hamiltonian action, \pageref{hamiltac}, \pageref{semifreeHam}
\item Hamiltonian-minimal submanifold, \pageref{hminimalsubd}
\item Hamiltonian vector field, \pageref{lagrsubm}
\item Hard Lefschetz Theorem, \pageref{hardlef}, \pageref{hlth}
\item Hauptvermutung, \pageref{haupt}
\item Heisenberg group, \pageref{heisenb}
\item HEP (homotopy extension property), \pageref{hepdiag}
\item Hermitian quadric, \pageref{zpqua}
\item Hilton--Milnor Theorem, \pageref{HMTheorem}
\item Hirzebruch genus, \pageref{degenus}, \pageref{hirzcor}
  \subitem equivariant, \pageref{defeqgenus}
  \subitem fibre multiplicative, \pageref{deffibmul}
  \subitem rigid, \pageref{ellitpigenus}, \pageref{defrigid}
  \subitem universal, \pageref{unigenus}
\item Hirzebruch surface, \pageref{hirzebruch}
\item $H$-minimal submanifold, \pageref{hminimalsubd}
\item Hodge algebra, \pageref{straighrel}
\item Hodge number, \pageref{dolbcohom}
\item Homological dimension (of module), \pageref{resol}
\item Homology
  \subitem cellular, \pageref{celluhomol}
  \subitem of chain complex, \pageref{chaincomplex}
  \subitem simplicial, \pageref{simphomology}
    \subsubitem reduced, \pageref{reducedhomology}
  \subitem singular, \pageref{singuhomo}
    \subsubitem of pair, \pageref{homopair}
    \subsubitem reduced, \pageref{redusingu}
\item Homology polytope, \pageref{deffaceacy}
\item Homology sphere, \pageref{cannon}
\item Homotopy, \pageref{dehomotopy}
\item Homotopy category (of model category), \pageref{homotcatego}
\item Homotopy cofibre, \pageref{dehocofibre}
\item Homotopy colimit, \pageref{Dcofib}, \pageref{defhocolhol}
\item Homotopy equivalence, \pageref{dehoequ}
\item Homotopy fibre, \pageref{homotopyfibre}
\item Homotopy group, \pageref{homotopygroup}
\item Homotopy Lie algebra, \pageref{djcoform},
\pageref{homotLie}, \pageref{quilsamel}
\item Homotopy limit, \pageref{Dcofib}, \pageref{defhocolhol}
\item Homotopy quotient, \pageref{2projections}
\item Homotopy type, \pageref{dehoequ}
\item Hopf Conjecture, \pageref{Hopfconj}
\item Hopf equation, \pageref{Hopfequ}
\item Hopf line bundle, \pageref{canolinebu}
\item Hopf manifold, \pageref{hopf}
\item $H$-polynomial (of polytope), \pageref{fhvect}
\item hsop (homogeneous system of parameters), \pageref{krdim}
  \subitem in face rings, \pageref{lsopr}, \pageref{lsopscc}
\item $h$-vector
  \subitem of polytope, \pageref{fhvect}
  \subitem of simplicial complex, \pageref{simfve}
  \subitem of simplicial poset, \pageref{fvsimpose}
\item Hyperplane cut, \pageref{hypcut}

\indexspace

\item Intersection homology, \pageref{intehomo}, \pageref{intecohomol}
\item Intersection poset (of arrangement), \pageref{gomaformu}
\item Isotropy representation, \pageref{isotropyrepr}

\indexspace

\item Join (least common upper bound), \pageref{joinandmeet},
\pageref{joinmee}
\item Join (operation)
  \subitem of simplicial complexes, \pageref{join}
  \subitem of simplicial posets, \pageref{joinsimpose}
  \subitem of spaces, \pageref{formtvman}

\indexspace

\item Klein bottle, \pageref{Kleinbot}
\item Koszul algebra, \pageref{koszulal}
\item Koszul complex, \pageref{koszulal}
\item Koszul resolution, \pageref{koszul}

\indexspace

\item Lagrangian immersion, \pageref{lagrsubm}
\item Lagrangian submanifold, \pageref{lagrsubm}
\item Landweber Exact Functor Theorem, \pageref{Landefth}
\item Latching functor, \pageref{delatch}
\item Lattice, \pageref{tvff}, \pageref{defnlattice}
\item Lefschetz pair, \pageref{lefpair}
\item Left lifting property, \pageref{RLPLLP}
\item Leibniz identity, \pageref{leibnitz}
\item $L$-genus (signature), \pageref{xyexa}
\item Limit (of diagram), \pageref{dfinlimcolim}
\item Link (in simplicial complex), \pageref{deflink}, \pageref{starlinksp}
\item Link (of intersection of quadrics), \pageref{linkintq}
\item Localisation formula, \pageref{evutg}
\item Locally standard (torus action), \pageref{defls}, \pageref{locastatm}
\item Logarithm (of formal group law), \pageref{logfglaw}
\item Loop functor (algebraic), \pageref{barcobar}
\item Loop space, \pageref{loops}, \pageref{pathloopsp}
\item Lower Bound Theorem, \pageref{lbt}
\item lsop (linear system of parameters), \pageref{krdim}
  \subitem in face rings, \pageref{lsopr}, \pageref{lsopscc}
  \subitem integral, \pageref{integrallsop}
\item LVM-manifold, \pageref{lvmmanid}, \pageref{lvmmanicomp}

\indexspace

\item Manifold with corners, \pageref{maniwithcor}, \pageref{mawco}
  \subitem face-acyclic, \pageref{deffaceacy}
  \subitem nice, \pageref{nicemwco}
\item Manifold with faces, \pageref{mawfa}
\item Mapping cone, \pageref{dehocofibre}
\item Mapping cylinder, \pageref{dehocofibre}
\item Massey product, \pageref{baskconstr}, \pageref{massey}
  \subitem indecomposable, \pageref{indec3massey}
  \subitem indeterminacy of, \pageref{trivialmp}
  \subitem trivial (vanishing), \pageref{trivialmp}
\item Matching functor, \pageref{dematch}
\item Mayer--Vietoris sequence, \pageref{maviseq}
\item Maximal action (of torus), \pageref{maxitorusact}
\item Meet (greatest common lower bound), \pageref{joinandmeet},
\pageref{joinmee}
\item Milnor hypersurface, \pageref{hij}
\item Minimal model
  \subitem of dg-algebra, \pageref{minmo}, \pageref{minmodelmodel}
  \subitem of space, \pageref{minmodelspace}
\item Minimal submanifold, \pageref{hminimalsubd}
\item Minkowski sum, \pageref{minkosu}
\item Missing face, \pageref{defcliq}, \pageref{higherWP},
\pageref{homotflag}
\item Model (of dg-algebra), \pageref{defiquism}
\item Model category, \pageref{modelcategory}
\item Module, \pageref{dmodule}
  \subitem finitely generated, \pageref{fgmodule}
  \subitem free, \pageref{freemodule}
  \subitem graded, \pageref{gradedmodule}
  \subitem projective, \pageref{projectivemod}
\item Moment-angle complex, \pageref{constrmac}, \pageref{macsimpos}
  \subitem real, \pageref{realmac}
\item Moment-angle manifold, \pageref{mamanifo}, \pageref{momeanma}
  \subitem polytopal, \pageref{mamanifo}, \pageref{polymamanif}
  \subitem non-formal, \pageref{nonformalmac}
\item Moment map, \pageref{dmomentmap}, \pageref{symred1}
  \subitem proper, \pageref{constsymred}, \pageref{propmmap}
\item Monoid, \pageref{monoidpair}
\item Monomial ideal, \pageref{sfmi}
\item Moore loops, \pageref{mooreloop}
\item Multi-fan, \pageref{multifan}, \pageref{toddindep}
\item Multidegree, \pageref{bigradedalge}
\item Multigrading, \pageref{mgrad}, \pageref{multigrator}, \pageref{multibnsp}
\item $M$-vector, \pageref{defmvect}

\indexspace

\item Nerve complex, \pageref{polsph}
\item Nested set, \pageref{defnestohe}
\item Nestohedron, \pageref{defnestohe}
\item Nilpotent space, \pageref{nilpospace}
\item Noether Normalisation Lemma, \pageref{noether}
\item Non-$PL$ sphere, \pageref{nonPLsph}, \pageref{non-PL}
\item Normal complex structure, \pageref{normalcs}

\indexspace

\item Omniorientation, \pageref{defomni}, \pageref{omnitorma},
\pageref{omniorittm}, \pageref{omniorieorie}
\item Opposite category, \pageref{secmc}
\item Orbifold, \pageref{dorbifold}
\item Orbit (of group action), \pageref{stasg}
\item Orbit space, \pageref{orbitspace}
\item Order complex (of poset), \pageref{ordercom}
\item Overcategory, \pageref{overundercategory}

\indexspace

\item Path space, \pageref{pathloopsp}
\item Pair (of spaces), \pageref{depairsp}
\item Perfect elimination order, \pageref{peorder}
\item Permutahedron, \pageref{permutahe}
\item Picard group, \pageref{picardgro}
\item $PL$ map, \pageref{defplmap}
  \subitem homeomorphism, \pageref{defplmap}
\item $PL$ manifold, \pageref{simman1}
\item $PL$ sphere, \pageref{defts}
\item Poincar\'e algebra, \pageref{pdalgebra}
\item Poincar\'e--Atiyah duality, \pageref{padual}
\item Poincar\'e duality, \pageref{pdtheorem}
\item Poincar\'e duality space, \pageref{defPds}
\item Poincar\'e series \pageref{poiser}, \pageref{loopsps}
  \subitem of face ring, \pageref{psfr}
\item Poincar\'e sphere, \pageref{non-PL}
\item Pointed map, \pageref{pointedspace}
\item Pointed space, \pageref{pointedspace}
\item Polar set, \pageref{polarset}
\item Polyhedral complex, \pageref{schldiag}
\item Polyhedral product, \pageref{nsc}, \pageref{polyprosp},
\pageref{homotopy}
\item Polyhedral smash product, \pageref{phsp}
\item Polyhedron (convex), \pageref{pol2}
  \subitem Delzant, \pageref{delzpolyh}
\item Polyhedron (simplicial complex), \pageref{polyhed}
\item Polytopal sphere, \pageref{defps}
\item Polytope
  \subitem combinatorial, \pageref{combinptpe}
  \subitem convex, \pageref{pol1}
  \subitem cyclic, \pageref{cyclic}
  \subitem Delzant, \pageref{delzantptope}, \pageref{delzaptope},
  \pageref{hamitoma}
  \subitem dual, \pageref{dualpol}
  \subitem generic, \pageref{genpolyt}
  \subitem lattice, \pageref{constrproj}
  \subitem neighbourly, \pageref{neighb}
  \subitem nonrational, \pageref{nrptope}
  \subitem polar, \pageref{polarset}
  \subitem regular, \pageref{24cell}
  \subitem self-dual, \pageref{selfdual}
  \subitem simple, \pageref{simpleptpe}, \pageref{ksimpleptope}
  \subitem simplicial, \pageref{scpolyt}, \pageref{ksimplicial}
  \subitem stacked, \pageref{stacked}
  \subitem triangle-free, \pageref{trianfre}
\item Polytope algebra, \pageref{hardlef}
\item Pontryagin algebra, \pageref{paseq}
\item Pontryagin product, \pageref{pontproduct}
\item Pontryagin--Thom map, \pageref{thomcla}, \pageref{ponthomcol}
\item Poset (partially ordered set), \pageref{faceposet}
\item Poset category, \pageref{posetcateg}
\item Positive orthant, \pageref{simpleptpe}, \pageref{posiorth}
\item Presentation (of a polytope by inequalities), \pageref{presptpe}
  \subitem generic, \pageref{scpolyt}
  \subitem irredundant, \pageref{presptpe}
  \subitem rational, \pageref{rationprese}
\item Primitive (lattice vector), \pageref{cpcone}
\item Principal bundle, \pageref{pribu}
\item Principal minor (of matrix), \pageref{priminor}
\item Product
  \subitem of building sets, \pageref{prodbse}
  \subitem of polytopes, \pageref{prodsp}
\item Product (categorical), \pageref{depushout}
\item Product (in cobordism), \pageref{cobproducts}
\item Projective dimension (of module), \pageref{resol}
\item Pseudomanifold, \pageref{defpm}, \pageref{pseudosimpos}
  \subitem orientable, \pageref{fc}
\item Pullback, \pageref{depushout}
\item Pushout, \pageref{depushout}

\indexspace

\item Quadratic algebra, \pageref{quadralg}, \pageref{quadralge}
\item Quasi-isomorphism, \pageref{dquism}, \pageref{defiquism},
\pageref{quismmodel}
\item Quasitoric manifold, \pageref{qtm}, \pageref{qtmanmode},
\pageref{qtgenera}
  \subitem almost complex structure on, \pageref{kusth}
  \subitem canonical smooth structure on, \pageref{cansmstqt}
  \subitem complex structure on, \pageref{kusth}
  \subitem equivalence of, \pageref{psieq}, \pageref{cqtp}
\item Quillen pair (of functors), \pageref{adjpair}
\item Quotient space, \pageref{orbitspace}

\indexspace

\item Rank function, \pageref{defsp}
\item Rational equivalence, \pageref{ratioequ}
\item Rational homotopy type, \pageref{ratioequ}
\item Ray's basis, \pageref{raybasis}
\item Redundant inequality, \pageref{presptpe}
\item Reedy category, \pageref{defnreedy}, \pageref{defhocolhol}
\item Regular sequence, \pageref{regseq}
\item Regular subdivision, \pageref{regsubsp}
\item Regular value (of moment map), \pageref{constsymred},
\pageref{propmmap}
\item Reisner Theorem, \pageref{reisner}
\item Resolution (of module)
  \subitem free, \pageref{resol}
  \subitem minimal, \pageref{minimal}
  \subitem projective, \pageref{resol}
\item Restriction (of a building set), \pageref{restricti}
\item Restriction map
  \subitem algebraic, \pageref{rmap}, \pageref{lsopr},
  \pageref{restr}, \pageref{rmapsp}
  \subitem in equivariant cohomology, \pageref{rmevc}
\item Right-angled Artin group, \pageref{raCoxeter}
\item Right-angled Coxeter group, \pageref{raCoxeter}
\item Right lifting property, \pageref{RLPLLP}
\item Ring of polytopes, \pageref{ringofptopes}

\indexspace

\item Samelson product, \pageref{samelprodu}, \pageref{samproduc}
  \subitem higher, \pageref{hsamelprodu}, \pageref{quilsamel}
\item Schlegel diagram, \pageref{schldiag}
\item Segre embedding, \pageref{segreemb}
\item Seifert fibration, \pageref{seifertfib}
\item Sign (of fixed point), \pageref{qts}, \pageref{weigsign},
\pageref{defsign}
\item Signature, \pageref{signor}, \pageref{xyexa}
\item Simplex, \pageref{simcub}, \pageref{faceofsp}
  \subitem abstract, \pageref{absimcom}
  \subitem regular, \pageref{stsim}
  \subitem standard, \pageref{stsim}
\item Simplicial cell complex, \pageref{simpcellcomp}
\item Simplicial chain, \pageref{desimplchain}
\item Simplicial cochain, \pageref{simplcochain}
\item Simplicial complex
  \subitem abstract, \pageref{absimcom}
  \subitem Alexander dual, \pageref{dual}
  \subitem geometric, \pageref{polyhed}
  \subitem pure, \pageref{defsubcomplex}
  \subitem shifted, \pageref{defshiftedco}, \pageref{zkshifted}
\item Simplicial manifold, \pageref{simman1}
\item Simplicial map, \pageref{simmap}
  \subitem isomorphism, \pageref{simmap}
  \subitem nondegenerate, \pageref{simmap}
\item Simplicial object (in category), \pageref{dediagram}
\item Simplicial poset, \pageref{defsp}
  \subitem dual (of manifold with corners), \pageref{dualsp}
\item Simplicial set, \pageref{dediagram}
\item Simplicial sphere, \pageref{defts}
\item Simplicial subdivistion, \pageref{simpsubdi}
\item Simplicial wedge, \pageref{simplicialwedge}
\item Singular chain, \pageref{singularchain}
  \subitem of pair, \pageref{homopair}
\item Singular cochain, \pageref{singucochain}
\item Skeleton (of cell complex), \pageref{cellcomplex}
\item Slice Theorem, \pageref{slicethm}
\item Small category, \pageref{posetcateg}
\item Small cover, \pageref{smalcove}, \pageref{realtoriv}, \pageref{bottsmalcov}
\item Smash product, \pageref{halfsmash}, \pageref{pointedspace}
\item Stabiliser (of group action), \pageref{stasg}
\item Stably complex structure, \pageref{qtscs}, \pageref{scs}
  \subitem $T$-invariant, \pageref{iscs}
\item Stanley--Reisner ideal, \pageref{frsim}
\item Star, \pageref{deflink}, \pageref{starlinksp}, \pageref{starinfan}
\item Starshaped sphere, \pageref{defps}, \pageref{starshsphd}
\item Stasheff polytope, \pageref{stasheptpe}
\item Stationary subgroup, \pageref{stasg}
\item Steinitz Problem, \pageref{steinproblem}
\item Steinitz Theorem, \pageref{stein}
\item Stellahedron, \pageref{defstellahe}
\item Stellar subdivision
  \subitem of simplicial complex, \pageref{bist}
  \subitem of simplicial poset, \pageref{prop:barst}
\item Straightening relation, \pageref{straighrel}
\item Subcomplex
  of simplicial complex, \pageref{defsubcomplex}
\item Sullivan algebra (of piecewise polynomial differential
forms), \pageref{sullivanalg}
\item Supporting hyperplane, \pageref{suphp}
\item Suspension
  \subitem of module, \pageref{freemodule}
  \subitem of simplicial complex, \pageref{consusp}
  \subitem of space, \pageref{conesuspspace}
\item Suspension isomorphism, \pageref{suspiso}
\item Symplectic manifold, \pageref{sympmani}
\item Symplectic reduction, \pageref{constsymred},
\pageref{symred1}, \pageref{Hminimaltoricd}
\item Symplectic quotient, \pageref{sympquoti}, \pageref{symred1}
\item Sysygy, \pageref{resol}

\indexspace

\item Tangential representation, \pageref{isotropyrepr}
\item Tautological line bundle, \pageref{canolinebu}
\item Tensor product (of modules), \pageref{tensoramod}
\item $T$-graph (torus graph), \pageref{defiTgraph}
\item Thom class, \pageref{thomclasstm}, \pageref{ThomTgr},
\pageref{thomcla}
\item Thom space, \pageref{thomspaced}
\item $T$-manifold, \pageref{torus}
\item Todd genus, \pageref{tdqtm}, \pageref{xyexa}
\item Topological fan, \pageref{deftopfan}
  \subitem complete, \pageref{deftopfan}
\item Topological monoid, \pageref{defntopmonoid}
\item Topological toric manifold, \pageref{defittm}
\item $\Tor$-algebra, \pageref{toralg}
\item Toral rank, \pageref{defftra}
\item Toric manifold, \pageref{cohtm}, \pageref{Hminimaltoricd},
\pageref{tos}
  \subitem Hamiltonian, \pageref{hamitoma}, \pageref{hamitomamdef}
  \subitem over cube, \pageref{tmancube}
  \subitem projective, \pageref{cohtm}
\item Toric space, \pageref{toricspace}, \pageref{toricspac}
\item Toric variety, \pageref{hardlef}, \pageref{toric2}
  \subitem affine, \pageref{affinetv}
  \subitem projective, \pageref{constrproj}
  \subitem real, \pageref{realtoriv}, \pageref{smalcove}
\item Torus
  \subitem algebraic, \pageref{toricfan}
  \subitem standard, \pageref{toricspac}
\item Torus manifold, \pageref{deftoma}
\item Triangulated manifold, \pageref{simman1}
\item Triangulated sphere, \pageref{trisphe}, \pageref{defts}
\item Triangulation, \pageref{triangu}, \pageref{simpsubdi}
  \subitem neighbourly, \pageref{brueckner}
\item Tropical geometry, \pageref{realtoriv}
\item Truncated cube, \pageref{truncatedcube}
\item Truncated simplex, \pageref{trancasim}

\indexspace

\item Undercategory, \pageref{overundercategory}
\item Underlying simplicial complex (of fan), \pageref{undecompl}
\item Unipotent (upper triangular matrix), \pageref{dunipot}
\item Universal toric genus, \pageref{defutge}, \pageref{definitutg}
\item Upper Bound Theorem, \pageref{ubt}

\indexspace

\item Vector bundle, \pageref{devectorbundle}
\item Vertex
  \subitem of polytope, \pageref{facet}
  \subitem of simplicial complex, \pageref{ghostvertex}
  \subitem of simplicial poset, \pageref{defsp}
\item Vertex truncation, \pageref{hypcut}
\item Volume polynomial, \pageref{hardlef}

\indexspace

\item Weak equivalence
  \subitem of dg-algebras, \pageref{defiquism}
  \subitem in model category, \pageref{modelcategory}
\item Wedge (of spaces), \pageref{wedgedef}, \pageref{pointedspace}
  \subitem of spheres, \pageref{defshiftedco}, \pageref{zkglue},
  \pageref{flws}
\item Weight (of torus representation), \pageref{weighttoru},
\pageref{weigsign}, \pageref{defsign}
\item Weight graph, \pageref{weightgraphs}
\item Weight lattice (of torus), \pageref{defnlattice}
\item Whitehead product, \pageref{dWhitehead}
  \subitem higher, \pageref{higherWP}

\indexspace

\item Yoneda algebra, \pageref{yonedaalg}

\indexspace

\item Zonotope, \pageref{zonotop}

\indexspace

\item $\gamma$-polynomial, \pageref{defgamma}
\item $\gamma$-vector
  \subitem of polytope, \pageref{defgamma}
  \subitem of triangulated sphere, \pageref{Hgammas}
\item $\chi_{a,b}$-genus, \pageref{ftab}
\item $\chi_y$-genus, \pageref{chigenu}

\indexspace

\item $2$-truncated cube, \pageref{truncatedcube}
\item $24$-cell, \pageref{24cell}
\item $5$-lemma, \pageref{5lemma}
\end{theindex}
\end{document}